\documentclass[envcountsect]{svmono}

\usepackage{type1cm}

\usepackage{makeidx}         
\usepackage{graphicx}        
\usepackage{multicol}        
\usepackage[bottom]{footmisc}

\usepackage[nofontspec]{newtxtext}
\usepackage{newtxmath}       

\usepackage[unicode]{hyperref}
\usepackage[nameinlink,capitalize]{cleveref}
\crefname{appendix}{Appendix}{Appendices}
\Crefname{appendix}{Appendix}{Appendices}
\hypersetup{
colorlinks=true,
linktoc=all
}


\usepackage{enumitem}
\usepackage{mathtools}
\providecommand{\xlongrightarrow}[2][]{\xrightarrow[#1]{#2}}
\usepackage{tikz-cd}
\usetikzlibrary{decorations.pathreplacing,arrows.meta,calc,patterns}
\usepackage{dsfont}

\usepackage{fancyhdr}

\usepackage[final]{microtype}
\usepackage{inputenc}
\usepackage{textcomp}
\usepackage[T1]{fontenc}

\makeindex             

\graphicspath{{./}{images/}}

\newcommand{\NS}{\mathrm{NS}}
\newcommand{\lex}{\mathrm{lex}}
\newcommand{\NA}{\mathrm{NA}}
\newcommand{\SH}{\mathrm{SH}}
\newcommand{\PSH}{\mathrm{PSH}}
\newcommand{\QPSH}{\mathrm{QPSH}}
\newcommand{\BM}{\mathrm{BM}}
\newcommand{\Conv}{\mathrm{Conv}}
\newcommand{\FS}{\mathrm{FS}}
\newcommand{\An}{\mathrm{An}}
\newcommand{\val}{\mathrm{val}}
\newcommand{\Macron}{\mathrm{Defect}}
\newcommand{\Bl}{\mathrm{Bl}}
\newcommand{\Ele}{\mathrm{Ele}}
\newcommand{\triv}{\mathrm{triv}}
\newcommand{\gf}{\mathrm{gf}}
\newcommand{\tor}{\mathrm{tor}}
\newcommand{\prop}{\mathrm{prop}}
\newcommand{\Star}{\mathrm{Star}}
\newcommand{\Haus}{\mathrm{Haus}}
\newcommand{\sg}{\mathrm{sg}}
\newcommand{\TC}{\mathrm{TC}}
\newcommand{\Kah}{\textup{Käh}}

\newcommand{\uscsup}[1][]{\operatorname*{sup}_{#1}{}\!{}^{*}}
\newcommand{\ddc}{\mathrm{dd}^{\mathrm{c}}}
\newcommand{\ddt}{\frac{\mathrm{d}}{\mathrm{d}t}}

\newcommand{\loc}{\mathrm{loc}}

\newcommand{\qsup}{\operatorname*{q-sup}}
\newcommand{\qinf}{\operatorname*{q-inf}}

\DeclareMathOperator{\Nm}{Norm}
\DeclareMathOperator{\nn}{nn}
\DeclareMathOperator{\nK}{nK}
\DeclareMathOperator{\Capa}{Cap}
\DeclareMathOperator{\Real}{Re}
\DeclareMathOperator{\Spec}{Spec}
\DeclareMathOperator{\Sing}{Sing}

\DeclareMathOperator{\Rest}{Res}
\DeclareMathOperator{\DHm}{DH}
\DeclareMathOperator{\Fin}{Fin}
\DeclareMathOperator{\Reg}{Reg}
\DeclareMathOperator{\Tr}{Tr}
\DeclareMathOperator{\Trop}{Trop}
\DeclareMathOperator{\Div}{div}
\DeclareMathOperator{\ord}{ord}
\DeclareMathOperator{\Dom}{Dom}
\DeclareMathOperator{\epi}{epi}
\DeclareMathOperator{\LCE}{CE}
\DeclareMathOperator{\Supp}{Supp}
\DeclareMathOperator{\MA}{MA}
\DeclareMathOperator{\vol}{vol}

\DeclareMathOperator{\Int}{Int}
\DeclareMathOperator{\RelInt}{RelInt}
\DeclareMathOperator{\cl}{cl}
\DeclareMathOperator{\cor}{cor}
\DeclareMathOperator{\GL}{GL}
\DeclareMathOperator{\Aut}{Aut}

\DeclareMathOperator{\mult}{mult}
\DeclareMathOperator{\Bs}{Bs}

\setlist[enumerate,1]{label=\textup{(\arabic*)},ref=\textup{(\arabic*)}}

\makeatletter
\newcommand{\itemref}[1]{\itemref@split\cref#1:\@nil}
\newcommand{\Itemref}[1]{\itemref@split\Cref#1:\@nil}
\def\itemref@split#1itm:#2:#3:#4:\@nil{#1{#2:#3}\,\ref{itm:#2:#3:#4}}
\makeatother

\begin{document}

\author{Mingchen Xia}
\title{Singularities in global pluripotential theory}
\subtitle{-- Lectures at Zhejiang University --}
\date{Updated on \today.}
\pagenumbering{Alph}
\maketitle

\frontmatter

\preface

This book is an expanded version of my lecture notes at the Institute for Advanced Study in Mathematics (IASM) at Zhejiang University. My initial goal was modest: to write a self-contained reference for the participants in the lectures. But I soon realized that many results circulating in the area had never been rigorously proved anywhere in the literature. Each attempt to resolve these loose ends produced new ones, and the notes grew steadily until they became the present book.

\section*{Three pluripotential theories.}
The aim of this book is to present, from a unified point of view, three different but interrelated theories that all deserve to be called a \emph{global pluripotential theory}:
\begin{enumerate}
    \item\label{itm:preface:L7:1} the pluripotential theory on compact Kähler manifolds;
    \item\label{itm:preface:L7:2} the pluripotential theory on the Berkovich analytification of a projective variety;
    \item\label{itm:preface:L7:3} the toric pluripotential theory on a projective toric variety.
\end{enumerate}
The three theories are bound together by a common circle of ideas. Many results were originally proved in one setting and have analogues, sometimes evident and sometimes far from evident, in the others.

\section*{The central question: classifying singularities.}
Let us begin with the analytic side. Fix a compact Kähler manifold $X$. The central objects are the \emph{quasi-plurisubharmonic functions} on $X$ --- functions which are locally the sum of a smooth function and a plurisubharmonic function. We are mostly interested in their \emph{singularities}: the locus where $\varphi=-\infty$ and the rate at which $\varphi$ tends to $-\infty$ near this locus.

Singularities are not pathological exceptions; they are forced upon us by geometry. In typical applications, $X$ is the compactification of the moduli space of some class of geometric objects, with a Zariski-open locus $U\subseteq X$ parameterizing smooth objects. The natural metric on the associated polarizing line bundle is then smooth only on $U$, and the Grauert--Remmert extension theorem (\cref{thm:GRexten2}) extends it to all of $X$ at the cost of introducing singularities along $X\setminus U$. Understanding such singularities is therefore not an aesthetic refinement but a prerequisite for any moduli-theoretic application.

The local theory of these singularities is intricate, and a complete classification is essentially out of reach: locally, (quasi-)plurisubharmonic functions exhibit extremely complicated behavior, and the natural local invariants (Lelong numbers and multiplier ideal sheaves) capture only their rough features. For example, a function with log-log singularities has the same local invariants as a bounded one.

The picture changes drastically when one passes to the global setting. On a compact manifold, there are three natural equivalence relations on $\QPSH(X)$, of decreasing fineness:
\begin{enumerate}
    \item\label{itm:preface:L22:1} the \emph{singularity-type} relation $\sim$, identifying $\varphi,\psi$ whose difference $\varphi-\psi$ is bounded;
    \item\label{itm:preface:L22:2} the \emph{$P$-equivalence} $\sim_P$, two potentials with $\varphi\leq \psi$ are identified if the mixed non-pluripolar Monge--Amp\`ere masses
    \[
        \int_X \omega_\varphi^j\wedge \omega_\phi^{n-j} = \int_X \omega_\psi^j\wedge \omega_\phi^{n-j},\qquad j=0,\ldots,n,
    \]
    agree against every reference potential $\phi$. Here $\omega$ is a sufficiently large Kähler form. In other words, $\varphi$ and $\psi$ behave similarly with respect to the global non-pluripolar mass (\cref{def:Pmoresing});
    \item\label{itm:preface:L22:3} the \emph{$\mathcal{I}$-equivalence} $\sim_\mathcal{I}$, identifying $\varphi,\psi$ for which the multiplier ideal sheaves agree, $\mathcal{I}(k\varphi)=\mathcal{I}(k\psi)$ for every $k\in\mathbb{Z}_{>0}$ --- equivalently, $\varphi$ and $\psi$ have the same generic Lelong numbers along every prime divisor over every modification of $X$ (\cref{def:Imoresing}).
\end{enumerate}
The first relation records the finest information one can hope to capture by global invariants; $\sim_P$ retains the global \emph{mass} side, while $\sim_\mathcal{I}$ retains the algebraic-geometric side. They are studied systematically in \cref{chap:comp}.

A natural plan to classify singularities is the three-step decomposition:
\begin{enumerate}
    \item\label{itm:preface:L34:1} classify the $\mathcal{I}$-equivalence classes;
    \item\label{itm:preface:L34:2} classify the $P$-equivalence classes inside a given $\mathcal{I}$-class;
    \item\label{itm:preface:L34:3} classify the singularity types inside a given $P$-class.
\end{enumerate}
Step~3 has been the subject of much activity in the past decade under the name of \emph{pluripotential theory with prescribed singularities}, and is by now reasonably well understood. The first goal of this book is to give a complete answer to Step~1, and a satisfactory answer to Step~2 in the cases that arise in geometry.

\section*{$\mathcal{I}$-good singularities.}
Step~1 is the subject of \cref{chap:Igood}. We characterize $\mathcal{I}$-equivalence classes very explicitly through approximations by an explicit class of singularities, known as \emph{analytic singularities}; the analogous problem of approximating a plurisubharmonic function by ones with analytic singularities was studied in detail by Demailly.

The class of $\mathcal{I}$-good potentials is the central protagonist of this book. In a single equivalence theorem (\cref{thm:charIgoodasclosure}), $\mathcal{I}$-goodness admits three intuitive descriptions, each capturing a different facet of what it means for a singularity to be ``well-behaved.'' In the non-degenerate case, a quasi-plurisubharmonic function $\varphi$ with $\theta_{\varphi}=\theta+\ddc\varphi\geq 0$ for some closed smooth real $(1,1)$-form $\theta$ is $\mathcal{I}$-good if and only if any of the following hold:
\begin{enumerate}
    \item\label{itm:preface:L45:1} $\varphi$ is the limit, with respect to the $d_S$-pseudometric introduced in \cref{chap:comp}, of a sequence of potentials with analytic singularities. (\emph{Topological description}: $\mathcal{I}$-good potentials are precisely those reachable as limits of the explicit, well-understood class of analytic singularities.)
    \item\label{itm:preface:L45:2} The non-pluripolar Monge--Amp\`ere mass of $\varphi$ equals the algebraic volume introduced in \cref{def:volqpsh}:
    \[
        \int_X \theta_\varphi^n = \vol\theta_\varphi.
    \]
    (\emph{Mass--volume identity}: no mass is lost to pluripolar sets, so the analytic mass and the algebraic volume coincide.)
    \item\label{itm:preface:L45:3} The $P$-envelope and the $\mathcal{I}$-envelope of $\varphi$ introduced in \cref{chap:envelopes} agree:
    \[
        P_\theta[\varphi] = P_\theta[\varphi]_{\mathcal{I}}.
    \]
    (\emph{Analytic--algebraic envelope agreement}: the two natural least-singular representatives, one constructed from the $P$-equivalence relation and one from multiplier ideals, coincide.)
\end{enumerate}
The reader is invited to keep these three pictures in mind throughout the book; we will use whichever is most convenient at any given moment.

Step~2 is the more delicate one, and its difficulty turns out to be a feature, not a bug.

A given $\mathcal{I}$-equivalence class can consist of an enormous number of $P$-equivalence classes; constructing such examples is easy. By contrast, in the toric and non-Archimedean settings, Step~2 is essentially trivial: every $\mathcal{I}$-equivalence class consists of a single $P$-equivalence class. In the toric case, an $\mathcal{I}$- (or $P$-) equivalence class is a sub-convex body of the Newton body. In the non-Archimedean case, it is a homogeneous plurisubharmonic metric.

The discrepancy between the analytic and the toric/non-Archimedean theories suggests that, in the pluripotential theory on compact Kähler manifolds, a great many $\mathcal{I}$-classes are pathological, and that what one really wants is to single out within each $\mathcal{I}$-class a canonical $P$-class. The functions in this canonical class will be called \emph{$\mathcal{I}$-good}. The development of the theory of $\mathcal{I}$-good singularities, in \cref{chap:Igood}, is the technical core of this book. As the reader will see, almost every singularity that arises naturally in geometric applications is $\mathcal{I}$-good.

A useful guiding analogy is with measure theory. There exist many non-measurable functions, but unless one constructs one by hand using the axiom of choice, one only encounters measurable functions in practice. Similarly, although the space of non-$\mathcal{I}$-good singularities is non-empty, it is essentially invisible to geometric applications. Once one accepts the framework of $\mathcal{I}$-good singularities, Step~2 becomes essentially solved.

\section*{Tools, applications, conventions.}
The classification of singularities, in the form just described, is somewhat abstract. To put it to use, we develop two general techniques that allow induction on $\dim X$. For a prime divisor $Y$ in general position, the \emph{analytic Bertini theorem} relates quasi-plurisubharmonic functions on $X$ and on $Y$. For non-generic $Y$, the \emph{trace operator} associates with $\varphi$ a quasi-plurisubharmonic function on the normalization of $Y$, well-defined modulo $\mathcal{I}$-equivalence. These tools occupy \cref{chap:trace}.

The toric pluripotential theory on ample line bundles is mostly classical and is recalled, with full proofs, in \cref{chap:toric_ample}. The corresponding theory for big toric line bundles, which has never been written down in the literature, requires the partial Okounkov body machinery developed in \cref{chap:Okounkov} and is the content of \cref{chap:toric_big}.

The applications to non-Archimedean pluripotential theory occupy \cref{chap:NApotential}, building on the test-curve theory of \cref{chap:testcurve}. The convergence of partial Bergman kernels with arbitrary singular weight, extending the Berman--Boucksom--Witt Nyström theorem, is proved in \cref{chap:Bergman}.

\section*{Reading guide.}

This book is not an elementary introduction to pluripotential theory. It should not be the reader's first textbook in this field. The reader is expected to be comfortable with the basic pluripotential theory at the level of \cite{GZ17}, which is essentially the only prerequisite for the main body of the book.

My companion book \emph{Foundations of pluripotential theory} \cite{FPT} develops the analytic background for Part~I and gives detailed proofs of the foundational results used there.

This book is not ``a compilation of all the author's works''. Rather, it presents a new point of view on this theory, with new results and new insights. A reader who has gone through all my related papers is still encouraged to read at least part of this book.

This book is not an encyclopedia of pluripotential theory. For example, the reader would not find the definition of a pseudoconvex domain, nor the full details from the long research papers or books I depend on in the application part. These points irritated the referee, but I hope that these omissions are understandable, taking into account that this book is already over 500 pages long. Although I tried my best to make the book as self-contained as possible, certain proofs based on techniques far beyond the scope of this book are still omitted.

Results in this book are presented in the greatest possible generality. An expert reader may use this book as a convenient reference for the most general statements, and skip the proofs of the more general results if they are already familiar with them in a special case.

The reader is expected to read Part~II linearly, and consult the necessary preliminaries in Part~I when necessary.

A reader interested only in one application chapter may skip directly to it once the relevant Part~II prerequisites are in place. \cref{chap:toric_ample} is needed only for \cref{chap:toric_big}; it can otherwise be skipped without loss.

Some chapters in this book draw on more specialised background, recorded in the appendices: \cref{chap:convex} for convex analysis, \cref{chap:unibranch} for pluripotential theory on unibranch K\"ahler spaces, and \cref{chap:almostsg} for the almost-semigroup machinery underlying the partial Okounkov body construction.

In order to develop the theory of b-divisors as in \cref{chap:bdiv}, we need to take limits with respect to the net of all modifications of $X$. In general, this net does not admit a cofinal sequence.
As a consequence, nets are widely used in this book, and the reader is expected to be familiar with their basic properties. The books \cite{Wil70, Kel75} are good references for the basic theory.
It is essential to remember that the dominated convergence theorem and Fatou's lemma both fail for nets. A version of the monotone convergence theorem does hold, as in \cite[Proposition~7.12]{Fol99}.

I made a great effort to make sure that every single expression involving $(-\infty)-(-\infty)$ or $0\cdot \infty$ is explained. This point is often overlooked in the literature. In most cases, the fix is simple, but in some cases, it requires a non-trivial amount of work. Therefore a number of our proofs seem more complicated than in the original papers. Similarly, in the literature, many constants are only correct up to a constant involving $2$ or $\pi$. I tried to make sure that every constant in this book is correct, and hence there is an extensive discussion of the normalizations of various objects in \cref{chap:pshfunctions}, which again irritated the same referee. I hope that the reader will appreciate the effort and find the book more readable as a result.

All statements in this book have been formalized in Lean~4. The formalization is available in the \href{https://github.com/MingchenXia/SGPT-lean}{SGPT-lean repository}.

Comments and questions are always welcome, and I will maintain an \href{https://mingchenxia.github.io/Lectures/Errata_SGPT.pdf}{errata page} on my homepage as long as I am still alive, no matter whether I remain in academia or not.

\vspace{\baselineskip}
\begin{flushright}\noindent
\hfill \textit{Mingchen Xia}\\
in Hangzhou, March 2024\\
\end{flushright}

\enlargethispage{14\baselineskip}
\vfill

\begin{large}
Contact information:

\textit{Mingchen Xia}, \textsc{\href{https://igp.ustc.edu.cn/main.htm}{Institute of Geometry and Physics}, USTC}\,\href{https://www.ustc.edu.cn/}{\raisebox{-0.25ex}{\includegraphics[height=1.2em]{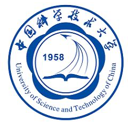}}}\par\nopagebreak
\textit{Email address}: \texttt{xiamingchen2008@gmail.com}\par\nopagebreak
\textit{Homepage}: \href{https://mingchenxia.github.io/}{\raisebox{-0.2ex}{\includegraphics[height=1.1em]{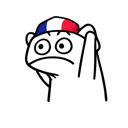}}}\,\url{https://mingchenxia.github.io/}.
\end{large}

\par\smallskip
\noindent\makebox[\textwidth][c]{%
    \includegraphics[width=0.75\textwidth]{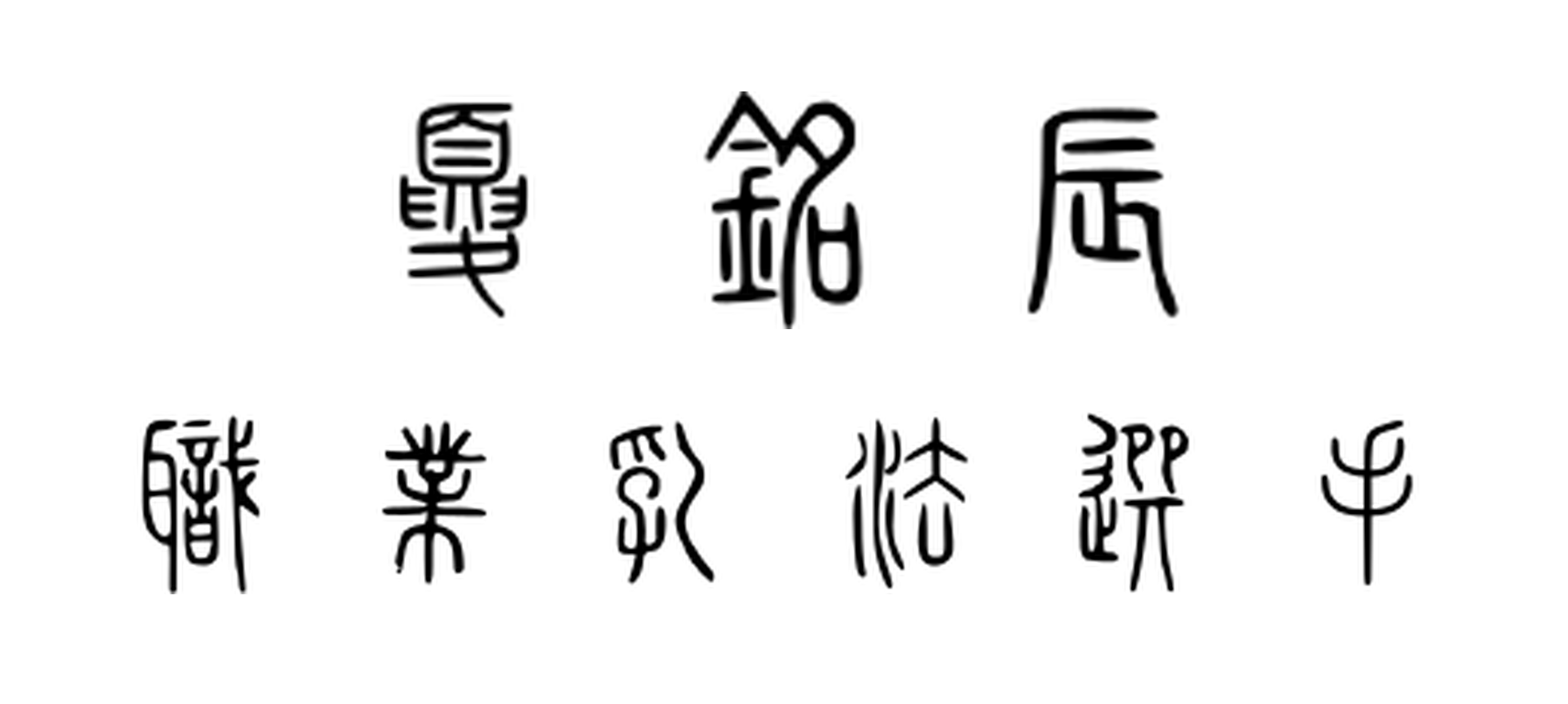}%
}
\par\vspace*{2.5\baselineskip}

\motto{The clerks are French, and, like most French people, are in a bad temper till they have eaten their lunch.

--- George Orwell (1903--1950)}

\extrachap{Acknowledgments}

I would like to express my gratitude to Bing Wang and Song Sun for their gracious invitations to China and for giving me the opportunity to give a series of lectures.

Furthermore, I am indebted to the dedicated researchers and administrative staff of the University of Science and Technology of China (USTC) and the Institute for Advanced Study in Mathematics (IASM) for creating an exceptional research environment. Their commitment to excellence has allowed me to immerse myself fully in mathematics during my time in China.

I am also immensely thankful to the participants in the course, including Song Sun, Mingyang Li, Xin Fu, Jiyuan Han, Junsheng Zhang, Yifan Chen, Yueqing Feng, Minghao Miao, and Federico Giust. Their active engagement and insightful discussions greatly enriched my lectures and deepened my understanding of the subject.

Special thanks go to Yi Yao, Kewei Zhang, and Robert Berman for their invaluable contributions to our discussions on toric geometry, which ultimately inspired the theory developed in \cref{chap:toric_big}.

Most results in this book were developed in collaboration with Tam\'as Darvas and Kewei Zhang, whose insights have been crucial to the development of the theory.
I would like to thank them for their collaboration over the years.

A substantial part of this book was already contained in my PhD thesis. I would like to thank my advisor, Robert Berman, for his guidance, and my colleagues in Göteborg and Paris for constant discussions, especially Bo Berndtsson, David Witt Nystr\"om, Sébastien Boucksom, and Elizabeth Wulcan.

I am especially grateful to Zetian Yan and Andreas Hohl. The former helped solve the problem of typesetting Polish, Romanian, and Swedish special characters within the framework of the SVMono class provided by Springer. Most of the book was typeset on the computer I inherited from my officemate Andreas at IMJ-PRG, which proved to be a worthy successor to the ``computeroid'' graciously provided by the university --- a machine on which coding was, shall we say, aspirational.

This work would not have been possible without the unwavering support and encouragement of all those mentioned above. Thank you for your generosity, guidance, and camaraderie throughout this endeavor.

I want to thank the following people and tools for pointing out mistakes or typos in the original version of the book: Vasanth Pidaparthy, Prakhar Gupta, Yi Yao, Tamás Darvas, Yan He, ChatGPT, Claude, and Codex.

I am very grateful to all but one of the referees for their detailed, professional, and helpful comments and criticism, which greatly improved the presentation of the whole book. Reading such a long and technical book is not an easy task, and their efforts are highly appreciated.

The book was written in the pre-AI era, but during the long publication process, AI tools developed rapidly, and I used them to check for typos and mistakes and to improve the writing. I am grateful to Claude, Codex, and ChatGPT for their help in improving the quality of the book. I am grateful to Xu Wang and Ragnar Sigurðsson for inviting me to Europe, where I had reliable access to the global internet and, in particular, to Claude.
At the same time, I am extremely disappointed by Anthropic's extreme hatred toward my homeland and by Dario Amodei's decision to block all users from mainland China, including myself.

I am very grateful to Jian Wang and Siqi He for their substantial help with the Lean~4 formalization of this book.

Enfin, je tiens à exprimer ma gratitude à Sébastien Boucksom et à Madame Natalia Hristic de Sorbonne Université, qui m'ont aidé à contacter le ministère de l'Intérieur en France. Sans leur intervention, je serais encore coincé en France, victime de l'extraordinaire efficacité du gouvernement français, en particulier de la préfecture de Créteil, et ce livre n'aurait jamais vu le jour. Merci également au comité du CNRS de m'avoir privé d'une audition deux années de suite, à savoir en 2023 et en 2024, me libérant ainsi à jamais de mon combat contre la bureaucratie française.

The author was supported initially by the Knut och Alice Wallenbergs Stiftelse KAW 2024.0273 during the early stages of this work, and later by the National Key R\&D Program of China 2025YFA1018200.

\tableofcontents

%
%

\extrachap{Conventions}

In the whole book, we adopt the following conventions:
\begin{itemize}
    \item A complex space is always assumed to be \emph{reduced}, \emph{paracompact} and \emph{Hausdorff}.
    \item A \emph{modification} of a complex space $X$ is a proper bimeromorphic morphism $\pi\colon Y\rightarrow X$ that is locally obtained from a finite composition of blow-ups with smooth centers. If $X$ is compact, we require $\pi$ to be a finite composition of blow-ups with smooth centers as in \cref{def:modif}.
    \item A \emph{subnet} of a net refers to a Willard subnet. More precisely, if $f\colon I\rightarrow X$ and $f'\colon I'\rightarrow X$ are nets, namely, maps from directed sets $I$ and $I'$ to the same set $X$, we say $f'$ is a subnet of $f$ if there is an \emph{increasing} final function $h\colon I'\rightarrow I$ such that $f\circ h=f'$.
    \item A \emph{domain} in $\mathbb{C}^n$ refers to a connected open subset.
    \item A \emph{complex manifold} is assumed to be paracompact.
    \item A \emph{submanifold} of a complex manifold means a closed complex submanifold.
    \item A \emph{neighborhood} is not necessarily open.
    \item The set $\mathbb{N}$ of natural numbers includes $0$.
    \item \emph{Increasing functions} and \emph{decreasing functions} are not necessarily strictly monotone.
    \item \emph{Positive forms} are not necessarily strictly positive.
    \item The function $\mathbb{R}\ni x \mapsto x^2$ is \emph{convex} instead of concave. This is different from the convention in \cite{CLS11}.
    \item We say a positive measure $\mu$ on a measurable space $X$ is \emph{carried} by a measurable subset $Y$ if $\mu(X\setminus Y)=0$.
    \item $\ddc$ means $(2\pi)^{-1}\mathrm{i}\partial \overline{\partial}$.
\end{itemize}

\mainmatter
\begin{partbacktext}
\part{Preliminaries}
\noindent
This part sets up the analytic and geometric language used throughout the book. It recalls the basic theory of quasi-plurisubharmonic functions, non-pluripolar products, envelope operators, geodesics of singular potentials and the toric-convex dictionary in the ample case. A reader already familiar with pluripotential theory on compact K\"ahler manifolds may skim parts of \cref{chap:pshfunctions,chap:npp}, but the later chapters use the precise conventions fixed here.

\cref{chap:pshfunctions} recalls plurisubharmonic and quasi-plurisubharmonic functions, pluripolar sets, the plurifine topology, Lelong numbers, multiplier ideals, analytic singularities, positivity of currents and singular Hermitian metrics. The presentation is concise, but we include proofs of several points whose conventions are important later, especially those involving Lelong numbers, multiplier ideals and quasi-equisingular approximation.

\cref{chap:npp} introduces the non-pluripolar product and the basic tools for manipulating Monge--Amp\`ere masses of singular currents.

\cref{chap:envelopes} develops the two envelope operators that organize the singularity theory of the book. The $P$-envelope records the singularity type relevant to non-pluripolar mass, while the $\mathcal I$-envelope records the singularity type detected by multiplier ideals. Their interaction is one of the main themes of the book.

\cref{chap:rays} studies subgeodesics, geodesics and geodesic rays in spaces of singular potentials. It also introduces the finite-energy space $\mathcal E^1(X,\theta;\phi)$ and the metric structures that will later be used to define convergence of singularity types.

Finally, \cref{chap:toric_ample} treats the ample toric case. There the pluripotential theory of toric metrics is translated into convex analysis on the associated polytope. This chapter is logically independent of much of the rest of the book, but it provides a model for the big toric theory developed later in \cref{chap:toric_big}.
\end{partbacktext}

\motto{Once Frigyes Riesz\footnote{Frigyes Riesz (1880--1956), known as Frédéric Riesz in French and Frederic Riesz in English, was the first mathematician to define the general notion of subharmonic functions, and he also gave these functions a Frenglish name from the very beginning --- \emph{fonctions subharmoniques}.} gave a brilliant explanation of why scientific work is easy. ``Everyone has ideas, both right ideas and wrong ideas,'' he said. ``Scientific work consists merely of separating them.''

--- Istvan Vincze}
\chapter{Plurisubharmonic functions}
\label{chap:pshfunctions}

The aim of this chapter is to recall the basic theory of plurisubharmonic functions and to fix the conventions used throughout the rest of the book. We assume only standard background in several complex variables and complex geometry, but we record a number of details whose normalization will matter later. General references for the material include \cite{GZ17}, \cite{Dem12} and \cite{DemBook}. My companion book \emph{Foundations of pluripotential theory} \cite[Chapters~7--8, 10--14]{FPT} develops this material in detail.

The chapter begins with plurisubharmonic functions on domains in $\mathbb C^n$ and on complex manifolds. We also include the $0$-dimensional case, which is occasionally convenient in inductive arguments. The basic permanence properties of plurisubharmonic functions are recalled: stability under upper envelopes and decreasing limits, Choquet's lemma, local integrability, gluing, extension across analytic subsets, and Kiselman's principle.

We then discuss the plurifine topology. This topology is the natural one for many local arguments in pluripotential theory, especially those involving contact sets and non-pluripolar products. We include the Bedford--Taylor base theorem in a form suited to later applications, following the approach of \cite{BT87} and \cite{EMW06}. See also \cite[Proposition~10.3.2]{FPT}.

The next part of the chapter fixes our conventions for Lelong numbers and multiplier ideal sheaves. These conventions differ by harmless but important factors from some common references; the conversions are explained in \cref{rmk:Lelongconv,rmk:misconv}. We also recall Siu's semicontinuity theorem, the strong openness theorem of Guan--Zhou, the Ohsawa--Takegoshi extension theorem, and the valuative characterization of multiplier ideals.

The remaining sections develop the global language used throughout the book. We introduce $\theta$-plurisubharmonic functions, quasi-plurisubharmonic functions, the preorder $\preceq$ of singularity types, and compactness properties in $L^1$. We then recall analytic singularities, log resolutions and quasi-equisingular approximation. Finally, we discuss closed positive $(1,1)$-currents, pseudo-effective and big classes, K\"ahler currents, Siu decomposition, non-K\"ahler and non-nef loci, modified nef classes, and singular Hermitian metrics on line bundles.

A reader familiar with this material may skim the chapter on a first pass, but the normalizations in \cref{sec:Lelongmis} and the approximation results in \cref{sec:anasing} will be used repeatedly later.

In this book, unlike in many references, we use the language of nets a lot. We remind the readers that subnets in the whole book mean \emph{Willard subnets}. In particular, a subnet of a monotone net is still monotone. The readers can find more details about nets or Moore--Smith convergence in Willard's textbook \cite{Wil70} or Kelley's \cite{Kel75}. The Willard convention and its use in compactness arguments are also discussed in \cite[Section~1.6]{FPT}. But one has to be careful that Willard's notion of subnets is stricter than Kelley's.

\section{The definition of plurisubharmonic functions}
\label{sec:pshdef}
In this section, we recall the notion of plurisubharmonic functions. We will also take care of the $0$-dimensional case, which makes a number of induction arguments easier to carry out.
None of our references treats the $0$-dimensional case, but the reader can easily verify that the results in this section hold in this exceptional case.

\subsection{The \texorpdfstring{$1$}{1}-dimensional case}
\label{subsec:psh_1d}
Let $\Omega$ be a domain (a connected open subset) in $\mathbb{C}$.

\begin{definition}\label{def:subhar1}
    A \emph{subharmonic function} \index{subharmonic function} on $\Omega$ is a function $\varphi\colon \Omega\rightarrow [-\infty,\infty)$ satisfying the following three conditions:
    \begin{enumerate}
        \item\label{itm:def:subhar1:1} $\varphi\not\equiv -\infty$;
        \item\label{itm:def:subhar1:2} $\varphi$ is upper semicontinuous;
        \item\label{itm:def:subhar1:3} $\varphi$ satisfies the \emph{sub-mean value inequality}: For any $a\in \Omega$ and $r>0$ such that $\overline{B_1(a,r)}\subseteq \Omega$, we have
        \[
        \varphi(a)\leq \frac{1}{2\pi}\int_0^{2\pi} \varphi(a+r\mathrm{e}^{\mathrm{i}\theta})\,\mathrm{d}\theta.\footnotemark
        \]
    \end{enumerate}
    \footnotetext{Condition~\ref{itm:def:subhar1:2} guarantees that $\varphi$ is measurable and locally bounded from above, and hence the integral in Condition~\ref{itm:def:subhar1:3} makes sense.}

    We will denote the set of subharmonic functions on $\Omega$ by $\SH(\Omega)$\index{$\SH(\Omega)$}.
\end{definition}
Here, $B_1(a,r)\coloneqq \{z\in \mathbb{C}:|z-a|<r\}$ denotes the open ball with center $a$ and radius $r$.

In fact, for each $a\in \Omega$, in \ref{itm:def:subhar1:3}, it suffices to require the sub-mean value inequality for all small enough $r>0$.

Intuitively, at a specific point $a\in \Omega$, Condition~\ref{itm:def:subhar1:2} says $\varphi(a)$ is no smaller than the upper limit of nearby values, while Condition~\ref{itm:def:subhar1:3} says $\varphi(a)$ is no larger than the mean. This intuition leads to the following rigidity theorem:
\begin{theorem}\label{thm:sh_rigid}
    Let $\varphi\colon \Omega\rightarrow [-\infty,\infty)$ be a measurable function. Then the following are equivalent:
    \begin{enumerate}
        \item\label{itm:thm:sh_rigid:1} $\varphi$ is locally integrable and $\Delta\varphi\geq 0$.
        \item\label{itm:thm:sh_rigid:2} $\varphi$ coincides almost everywhere with a subharmonic function $\psi$ on $\Omega$.
    \end{enumerate}
    Moreover, the subharmonic function $\psi$ in \ref{itm:thm:sh_rigid:2} is unique.
\end{theorem}
Here in Condition~\ref{itm:thm:sh_rigid:1}, $\Delta\varphi=2\mathrm{i}\partial\bar\partial\varphi$ is the Laplacian in the sense of currents. This is a special case of \cref{thm:psh_rigid} below.

This theorem gives a very useful way of constructing subharmonic functions.

\subsection{The higher dimensional case}
\label{subsec:psh_hd}
We will fix $n\in \mathbb{N}$ and a domain $\Omega$ (a connected open subset) in $\mathbb{C}^n$.

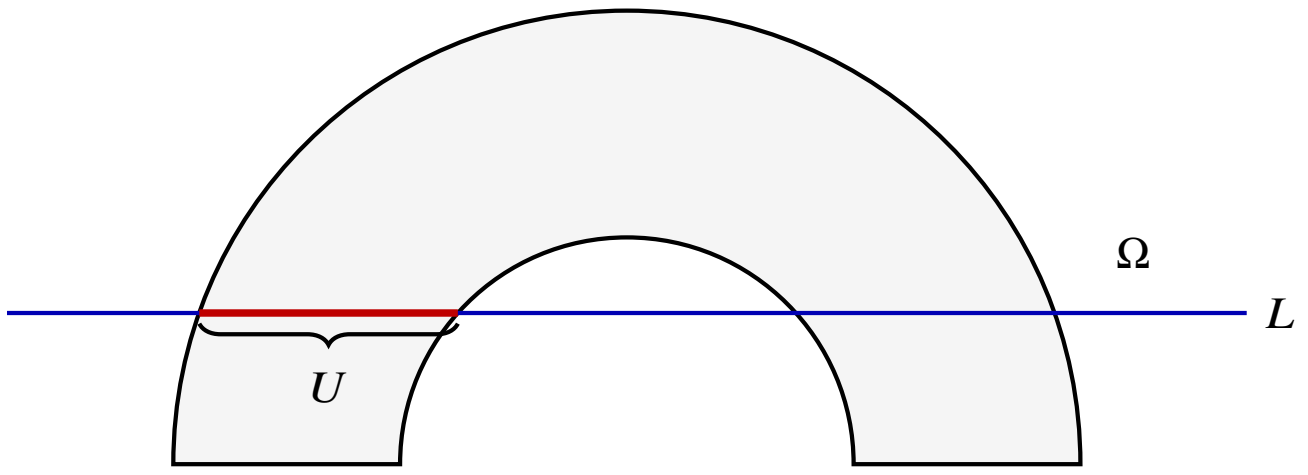
\begin{figure}[ht]
\centering
\begin{tikzpicture}[scale=1.0, line width=0.7pt]
  \draw[thick, fill=black!4]
    (-3,0) arc[start angle=180, end angle=0, radius=3]
    -- (1.5,0) arc[start angle=0, end angle=180, radius=1.5]
    -- cycle;
  \node at (3.35,1.4) {$\Omega$};
  \draw[blue!70!black, thick] (-4.1,1) -- (4.1,1) node[right, black] {$L$};
  \draw[red!75!black, line width=1.4pt] (-2.828,1) -- (-1.118,1);
  \draw[decorate, decoration={brace, mirror, amplitude=4pt, raise=2pt}, thick]
    (-2.828,1) -- (-1.118,1);
  \node[below=8pt] at (-1.973,1) {$U$};
\end{tikzpicture}
\caption{A domain $\Omega$ cut by a complex line $L$.}
\label{fig:psh}
\end{figure}

\begin{definition}\label{def:psh}
    When $n\geq 1$, a \emph{plurisubharmonic function}\index{plurisubharmonic function} on $\Omega$ is a function $\varphi\colon \Omega\rightarrow [-\infty,\infty)$ satisfying the following three conditions:
    \begin{enumerate}
        \item\label{itm:def:psh:1} $\varphi\not\equiv -\infty$;
        \item\label{itm:def:psh:2} $\varphi$ is upper semicontinuous;\footnote{This condition is redundant by \cite[Theorem~8.1.3]{FPT}.}
        \item\label{itm:def:psh:3} for any complex line $L\subseteq \mathbb{C}^n$ and any connected component $U$ of $L\cap \Omega$, the restriction $\varphi|_U$ is either subharmonic or constantly $-\infty$.\footnote{An extremely common mistake in the literature is to replace \ref{itm:def:psh:3} by the condition that $\varphi$ is locally integrable and $\ddc\varphi\geq 0$ in the sense of currents. For a concrete counterexample, consider a function $\varphi$ that takes a constant value $0$ at all but one single point, at which the value of $\varphi$ is $1$.}
    \end{enumerate}

    When $n=0$, the only domain $\Omega$ is the singleton. In this case, a \emph{plurisubharmonic function} on $\Omega$ is a real-valued function on $\Omega$.

    The set of plurisubharmonic functions on $\Omega$ is denoted by $\PSH(\Omega)$\index{$\PSH(\Omega)$}.
\end{definition}
A plurisubharmonic function is also called a psh function for short. The relevant notations are indicated in \cref{fig:psh}.\footnote{Recall that most figures in this book are somewhat misleading: We usually draw a complex dimension as a real dimension. The figures should not be read literally!}

\begin{example}\label{ex:pshfunctions:L104}
    When $n=0$, we have a canonical bijection $\PSH(\Omega)\cong \mathbb{R}$.
\end{example}

\begin{example}
    When $n=1$, we have $\PSH(\Omega)=\SH(\Omega)$.
\end{example}

Similar to \cref{thm:sh_rigid}, we have a rigidity theorem for plurisubharmonic functions as well.
\begin{theorem}\label{thm:psh_rigid}
    Let $\varphi\colon \Omega\rightarrow [-\infty,\infty)$ be a measurable function. Then the following are equivalent:
    \begin{enumerate}
        \item\label{itm:thm:psh_rigid:1} $\varphi$ is locally integrable and $\ddc\varphi\geq 0$;
        \item\label{itm:thm:psh_rigid:2} $\varphi$ coincides almost everywhere with a plurisubharmonic function $\psi$ on $\Omega$.
    \end{enumerate}
    Moreover, the plurisubharmonic function $\psi$ is unique.
\end{theorem}
Here, the operator $\ddc$ is normalized so that
\[
\ddc =\frac{\mathrm{i}}{2\pi}\partial \overline{\partial}.
\]
For the proof, we refer to \cite[Proposition~1.43]{GZ17} and \cite[Theorem~8.2.1]{FPT}.

Plurisubharmonic functions have nice functorialities:
\begin{proposition}\label{prop:func_domain}
    Let $n'\in \mathbb{N}$ and $\Omega'\subseteq \mathbb{C}^{n'}$ be a domain. Given any holomorphic map $f\colon \Omega\rightarrow \Omega'$ and any $\varphi\in \PSH(\Omega')$, exactly one of the following cases occurs:
    \begin{enumerate}
        \item\label{itm:prop:func_domain:1} $f^*\varphi\equiv -\infty$;
        \item\label{itm:prop:func_domain:2} $f^*\varphi\in \PSH(\Omega)$.
    \end{enumerate}
\end{proposition}
We refer to \cite[Proposition~1.44]{GZ17} for the proof\footnote{Recall that the statement of \cite[Proposition~1.44]{GZ17} is flawed. One has to reduce to the case where Case~\ref{itm:prop:func_domain:1} does not occur before following their proof.}. See also \cite[Proposition~8.2.3]{FPT}.

For each $n\in \mathbb{N}$, $a\in \mathbb{C}^n$ and $r>0$, we write
\begin{equation}\label{eq:Bnar}
B_n(a,r)=\left\{ z\in \mathbb{C}^n: |z-a|<r \right\}.
\end{equation}
\begin{proposition}\label{prop:ballpshconvex}
    Let $\varphi\in \PSH(B_n(a,r_0))$ for some $r_0>0$. Then the function
    \[
    (-\infty,\log r_0)\rightarrow \mathbb{R},\quad \log r\mapsto \sup_{B_n(a,r)}\varphi
    \]
    is convex and increasing.
\end{proposition}
See \cite[Theorem~4.1.13]{Hor07} for the case $n>1$ and \cite[Corollary~2.4]{Bou17} for the general case. See also \cite[Proposition~12.1.1]{FPT}.

\begin{proposition}\label{prop:subhimplyconv}
Let $a<b$ be two real numbers.
    Let $f\colon (a,b)\rightarrow [-\infty,\infty)$ be a function. Define
    \[
    g\colon \{z\in \mathbb{C}:a<\Real z<b\}\rightarrow [-\infty,\infty),\quad z\mapsto f(\Real z).
    \]
    Suppose that $g$ is subharmonic. Then $f$ is convex. In particular, $f$ takes real values only.
\end{proposition}

See \cite[Theorem~2.12]{HK76} for a more general result. For the convexity assertion, see also \cite[Proposition~18.1.5]{FPT}.

\subsection{The manifold case}
\label{subsec:psh_manifold}
Let $X$ be a complex manifold. In the whole book, complex manifolds are assumed to be paracompact, namely, all connected components have countable bases.

\begin{definition}\label{def:pshmfd}
    A \emph{plurisubharmonic function}\index{plurisubharmonic function} on $X$ is a function $\varphi\colon X\rightarrow [-\infty,\infty)$ such that for any $x\in X$, there exists an open neighborhood $U$ of $x$ in $X$, an integer $n\in \mathbb{N}$, a domain $\Omega\subseteq \mathbb{C}^n$ and a biholomorphic map $F\colon \Omega\rightarrow U$ such that $F^*(\varphi|_U)\in \PSH(\Omega)$.

    The set of plurisubharmonic functions on $X$ is denoted by $\PSH(X)$\index{$\PSH(X)$}.
\end{definition}

\begin{example}\label{ex:pshfunctions:L172}
    When $X$ is a domain in $\mathbb{C}^n$, the notions of plurisubharmonic functions in \cref{def:pshmfd} and in \cref{def:psh} coincide.
\end{example}

\begin{example}\label{ex:pshfunctions:L176}
    Write $\{X_i\}_{i\in I}$ for the set of connected components of $X$. Then we have a natural bijection
    \[
    \PSH(X) \cong \prod_{i\in I} \PSH(X_i).
    \]
    Here the product is in the category of sets.
    In particular, if $X=\varnothing$, then $\PSH(X)$ is a singleton.
\end{example}
This example allows us to reduce to the case of connected manifolds when studying general plurisubharmonic functions.

\begin{proposition}\label{prop:pullbackpsh}
    Let $Y$ be a complex manifold and $f\colon Y\rightarrow X$ be a holomorphic map. Then for any $\varphi\in \PSH(X)$, exactly one of the following cases occurs:
    \begin{enumerate}
        \item\label{itm:prop:pullbackpsh:1} $f^*\varphi$ is identically $-\infty$ on some connected component of $Y$;
        \item\label{itm:prop:pullbackpsh:2} $f^*\varphi\in \PSH(Y)$.
    \end{enumerate}
\end{proposition}
This proposition follows from \cref{prop:func_domain} applied separately on each connected component of $Y$.

\cref{thm:psh_rigid} implies immediately the general form of the rigidity theorem:
\begin{theorem}\label{thm:psh_rigid_gen}
    Let $\varphi\colon X\rightarrow [-\infty,\infty)$ be a measurable function. Then the following are equivalent:
    \begin{enumerate}
        \item\label{itm:thm:psh_rigid_gen:1} $\varphi$ is locally integrable and $\ddc \varphi\geq 0$;
        \item\label{itm:thm:psh_rigid_gen:2} $\varphi$ coincides almost everywhere with a plurisubharmonic function $\psi$ on $X$.
    \end{enumerate}
    Moreover, the plurisubharmonic function $\psi$ in \ref{itm:thm:psh_rigid_gen:2} is unique.
\end{theorem}

\begin{definition}\label{def:pluripolarsets}
    A subset $E\subseteq X$ is \emph{pluripolar}\index{set!pluripolar set} if for any $x\in X$, there is an open neighborhood $U$ of $x$ in $X$ and a function $\psi\in \PSH(U)$ such that
    \[
    \psi|_{E\cap U}\equiv -\infty.
    \]
    A subset $E\subseteq X$ is \emph{non-pluripolar}\index{set!non-pluripolar set} if $E$ is not pluripolar.

    A subset $F\subseteq X$ is \emph{co-pluripolar}\index{set!co-pluripolar set} if $X\setminus F$ is pluripolar.
\end{definition}
When $X$ has dimension $1$, a pluripolar set is called a \emph{polar set}. We say some property about objects on $X$ holds \emph{quasi-everywhere}\index{quasi-everywhere} if it holds outside a pluripolar set.

\begin{theorem}[Josefson's theorem]\label{thm:Josefson}\index{theorem!Josefson's theorem}
    Let $E\subseteq \mathbb{C}^n$ be a pluripolar set. Then there is $\varphi\in \PSH(\mathbb{C}^n)$ such that $\varphi|_E\equiv -\infty$.
\end{theorem}
See \cite[Corollary~4.41]{GZ17} for the proof of a more general result, and \cite[Theorem~10.2.1]{FPT} for the result stated here.

There is also a global version of Josefson's theorem:
\begin{theorem}\label{thm:gloJosefson}
    Assume that $X$ is a compact complex manifold and $E\subseteq X$ is a pluripolar set. Then there is a quasi-plurisubharmonic function $\varphi$ on $X$ with $\varphi|_E\equiv -\infty$.
\end{theorem}
See \cite{Vu19} and \cite[Theorem~14.8.7]{FPT} for proofs. We recall the definition of a quasi-plurisubharmonic function in \cref{def:qpsh}.

\section{Properties of plurisubharmonic functions}
\label{sec:psh_properties}
In this section, we explore the basic properties of plurisubharmonic functions.

Let $X$ be a complex manifold.

\begin{proposition}\label{prop:pshlocLp}
    Let $\varphi\in \PSH(X)$. Then for any $p\geq 1$, $\varphi\in L^p_{\loc}(X)$.
\end{proposition}
See \cite[Theorem~1.46, Theorem~1.48]{GZ17} and \cite[Theorem~8.3.7]{FPT}.\footnote{We remind the readers that the third displayed formula in the proof of \cite[Theorem~1.48]{GZ17} Step~2 is wrong. But this is irrelevant for our purpose. See also the discussion in \cite[Section~B.8]{FPT}.}

\begin{proposition}\label{prop:pshfunction_closedseq}
\leavevmode
\begin{enumerate}
    \item\label{itm:prop:pshfunction_closedseq:1} Assume that $(\varphi_i)_{i\in I}$ is a non-empty family in $\PSH(X)$ that is locally uniformly bounded from above. Then $\uscsup[i\in I] \varphi_i\in \PSH(X)$.
    \item\label{itm:prop:pshfunction_closedseq:2} Assume that $(\varphi_i)_{i\in I}$ is a decreasing net in $\PSH(X)$ such that $\varphi\coloneqq \inf_{i\in I}\varphi_i$ is not identically $-\infty$ on any connected component of $X$. Then $\varphi\in \PSH(X)$. Furthermore, $\varphi_i\xrightarrow{L^1_{\loc}} \varphi$.
\end{enumerate}
\end{proposition}
Here $\uscsup$ denotes the upper semicontinuous regularization of the supremum. When $I$ is a finite family, observe that
\[
\uscsup[i\in I]\varphi_i=  \sup_{i\in I}\varphi_i.
\]
When $I=\{1,\ldots,m\}$, we write
\[
    \varphi_1\lor \cdots\lor \varphi_m\coloneqq   \sup_{i\in I}\varphi_i.
\]
We refer to \cite[Proposition~1.28, Proposition~1.40]{GZ17} and \cite[Theorems~8.3.1--8.3.2]{FPT}\footnote{In \cite[Proposition~1.28]{GZ17}, the second part is only stated for sequences, the net version is obvious using the sub-mean value inequality and the monotone convergence theorem for nets as in \cite[Proposition~7.12]{Fol99}. For a proof for decreasing nets of psh functions, see \cite[Theorem~8.3.1]{FPT}.}\footnote{The proof of \cite[Proposition~1.40]{GZ17} contains a gap: The measurability of $u$ is not established during the proof, so the claim $u \ast \chi_{\epsilon}$ converges to $u$ in $L^1_{\loc}$ is not valid. One needs to choose the sequence $(v_j)_j$ carefully so that $v\leq u$, and use the fact that $v \ast \chi_{\epsilon}$ converges to $v$ in $L^1_{\loc}$ instead. See also \cite[Theorem~8.3.2]{FPT} for a proof by countable selection.}.

\begin{proposition}[Choquet's lemma]\label{prop:Choquet}\index{theorem!Choquet's lemma}
    Assume that $X$ has countably many connected components.
    Assume that $(\varphi_i)_{i\in I}$ is a non-empty family in $\PSH(X)$ that is locally uniformly bounded from above. There exists a countable subset $J\subseteq I$ such that
    \[
    \uscsup[i\in I] \varphi_i=\uscsup[j\in J] \varphi_j.
    \]
\end{proposition}
In the sequel, we say a subset $B\subseteq X$ is a \emph{holomorphic ball}\index{holomorphic ball} if $B$ has smooth boundary in $X$, and $B$ is biholomorphic to an open ball in $\mathbb{C}^n$.

\begin{proof}
    We may assume that $X$ is connected. Since by our convention, the complex manifold $X$ is paracompact, it can be covered by countably many holomorphic balls, so we can easily reduce to the case where $X$ is an open ball in $\mathbb{C}^n$. In this case, the result is proved in \cite[Lemma~4.31]{GZ17} and \cite[Lemma~3.4.1]{FPT}.
\end{proof}

\begin{proposition}\label{prop:ppsnull}
    A pluripolar set $E\subseteq X$ is a Lebesgue null set.
\end{proposition}
\begin{proof}
    This is a trivial consequence of \cref{prop:pshlocLp}.
\end{proof}

\begin{proposition}\label{prop:supsupstardiff}
    Let $(\varphi_i)_{i\in I}$ be a non-empty family in $\PSH(X)$ that is locally uniformly bounded from above. Then the set
    \[
        \left\{x\in X: \sup_{i\in I}\varphi_i<\uscsup[i\in I]\varphi_i \right\}
    \]
    is pluripolar and hence Lebesgue null.
\end{proposition}
See \cite[Corollary~4.28]{GZ17} and \cite[Corollary~10.2.2]{FPT}.

\begin{proposition}\label{prop:pshfuncdetdense}
    Suppose that $\varphi,\psi\in \PSH(X)$. If $\varphi\leq \psi$ almost everywhere, then $\varphi\leq \psi$ everywhere.

    In particular, if $\varphi=\psi$ almost everywhere, then $\varphi=\psi$ everywhere.
\end{proposition}
\begin{proof}
    Replace $\varphi$ by $\varphi\lor\psi$. Then $\varphi,\psi$ are two plurisubharmonic functions that coincide almost everywhere. The uniqueness statement in \cref{thm:psh_rigid_gen} gives $\varphi=\psi$.
\end{proof}

\begin{proposition}\label{prop:pluripolarunion}
    Let $(E_i)_{i\in \mathbb{Z}_{>0}}$ be a sequence of pluripolar sets in $X$. Then
    \[
    E\coloneqq \bigcup_{i=1}^{\infty}E_i
    \]
    is also pluripolar.
\end{proposition}
\begin{proof}
    The problem is local, so we may assume that $X\subseteq \mathbb{C}^n$ is a bounded domain. In this case, by \cref{thm:Josefson} for each $i\in \mathbb{Z}_{>0}$ we can choose $\psi_i\in \PSH(\mathbb{C}^n)$ such that
    \[
    \psi_i|_{E_i}\equiv -\infty,\quad \psi_i|_X\leq 0
    \]
    for all $i>0$. By \cref{prop:pshlocLp}, we have $\psi_i|_X\in L^1(X)$ for all $i>0$. After rescaling, we may also assume that $\|\psi_i\|_{L^1(X)}\leq 1$ for all $i>0$.

    We then define
    \[
    \psi=\sum_{i=1}^{\infty}2^{-i} \psi_i|_X.
    \]
    Then $\psi\in \PSH(X)$ according to \cref{prop:pshfunction_closedseq} and $\psi|_{E}=-\infty$.
\end{proof}

\begin{corollary}\label{cor:L1limipp}
    Let $(\varphi_j)_{j\in \mathbb{Z}_{>0}}$ be a sequence in $\PSH(X)$ such that $\varphi_j\xrightarrow{L^1_{\loc}}\varphi\in \PSH(X)$. Then
    \begin{equation}\label{eq:L1limsup1}
        \varphi=\mathop{\vphantom{\operatorname{sup}^{*}}\inf}\limits_{j>0} \uscsup[k\geq j]\varphi_k,
    \end{equation} and the set
    \[
        \left\{x\in X: \varphi(x)\neq \varlimsup_{j\to\infty}\varphi_j(x)\right\}
    \]
    is pluripolar.
\end{corollary}
\begin{proof}
    We may assume that $X$ is connected.

    We first observe that $(\varphi_j)_j$ is locally uniformly bounded from above. This follows from \cite[Exercise~1.20]{GZ17} and the proof of \cite[Theorem~8.3.7]{FPT}.

    For each $j\geq 1$, let
    \[
        \psi_j=\uscsup[k\geq j]\varphi_k.
    \]
    Then $\psi_j\in \PSH(X)$ by \cref{prop:pshfunction_closedseq}. Moreover, $(\psi_j)_j$ is a decreasing sequence and $\psi_j\geq \varphi_j$ for all $j$.
    After considering a subsequence of $(\varphi_j)_j$ which converges to $\varphi$ almost everywhere, we get $\varphi\leq \psi\coloneqq \inf_j \psi_j$ almost everywhere. In particular, $\psi\in \PSH(X)$ by \cref{prop:pshfunction_closedseq}.

    On the other hand, by \cref{prop:supsupstardiff}, there exist pluripolar sets $Z_j\subseteq X$ such that
    \[
    \psi_j=\sup_{k\geq j}\varphi_k
    \]
    on $X\setminus Z_j$. Let
    \[
    Z=\bigcup_{j=1}^{\infty} Z_j.
    \]
    Then $Z$ is a pluripolar set by \cref{prop:pluripolarunion}, and for any $x\in X\setminus Z$, we have $\psi(x)=\varlimsup_j \varphi_j(x)$.
    Moreover, $\psi=\varphi$ almost everywhere by \cite[Theorem~1.46]{GZ17} (see also Hartogs' lemma in \cite[Theorem~8.3.4]{FPT}), and hence everywhere by \cref{prop:pshfuncdetdense}. Hence we conclude \eqref{eq:L1limsup1}.
    Moreover, for all $x\in X\setminus Z$,
    \[
        \varphi(x)=\varlimsup_{j\to\infty}  \varphi_j(x).
    \]
    This completes the proof.
\end{proof}

\begin{theorem}[Brelot, Grauert--Remmert]\label{thm:GRexten}\index{theorem!Grauert--Remmert extension theorem}
    Let $Z$ be an analytic subset in $X$ and $\varphi\in \PSH(X\setminus Z)$. Then the function $\varphi$ admits an extension to $\PSH(X)$ in the following two cases:
    \begin{enumerate}
        \item\label{itm:thm:GRexten:1} The set $Z$ has codimension at least $2$ everywhere.
        \item\label{itm:thm:GRexten:2} The set $Z$ has codimension at least $1$ everywhere and $\varphi$ is locally bounded from above near each point of $Z$.
    \end{enumerate}
    In both cases, the extension is unique and is given by
    \begin{equation}\label{eq:GRextvarphiexp}
    \varphi(x)=\varlimsup_{X\setminus Z\ni y\to x}\varphi(y),\quad x\in Z.
    \end{equation}
\end{theorem}
\begin{figure}[ht]
\centering
\begin{tikzpicture}[scale=0.95, line width=0.6pt, every node/.style={font=\small}]
  \def\R{2.6}
  \def\Re{0.55}
  \draw[thick] (0,0) circle (\R);
  \draw[thick]         (-\R,0) arc[start angle=180, end angle=360, x radius=\R, y radius=\Re];
  \draw[thick, dashed] (-\R,0) arc[start angle=180, end angle=  0, x radius=\R, y radius=\Re];
  \draw[-{Stealth[length=2.5mm]}, thin] (3.6,1.85) -- (1.95,1.7);
  \node[anchor=west, inner sep=2pt] at (3.6,1.85) {$\partial B_n(z,r)$};
  \draw[thick] (-3.65, 2.05) .. controls (-3.0, 1.55) and (-2.75, 1.25) .. (-2.42, 0.95);
  \draw[thick, densely dashed]
    (-2.42, 0.95) .. controls (-1.5, 0.5) and (-0.6, 0.15) .. (0, 0)
    .. controls (0.6, -0.15) and (1.5, -0.5) .. (2.42, -0.95);
  \draw[thick] (2.42, -0.95) .. controls (2.75, -1.25) and (3.0, -1.55) .. (3.65, -2.05);
  \node at (3.85, -1.65) {$Z$};
  \fill (0,0) circle (1.7pt);
  \node[anchor=north east, inner sep=2pt] at (0,0) {$z$};
  \draw[gray, thin] ($ (1.5, 0.9) + 3.0*(0.2924, 0.9563) $) --
                    ($ (1.5, 0.9) - 3.0*(0.2924, 0.9563) $);
  \fill (1.5, 0.9) circle (1.3pt) node[right=1pt, black] {$x_1$};
  \node[inner sep=2pt] at ($ (1.5,0.9)+2.55*(0.2924,0.9563) $) {$L_1$};
  \draw[gray, thin] ($ (0.9, 0.4) + 3.0*(0.0872, 0.9962) $) --
                    ($ (0.9, 0.4) - 3.0*(0.0872, 0.9962) $);
  \fill (0.9, 0.4) circle (1.3pt) node[right=1pt, black] {$x_2$};
  \node[inner sep=2pt] at ($ (0.9,0.4)+2.65*(0.0872,0.9962) $) {$L_2$};
  \draw[gray, thin] ($ (0.4, 0.15) + 3.0*(-0.0349, 0.9994) $) --
                    ($ (0.4, 0.15) - 3.0*(-0.0349, 0.9994) $);
  \fill (0.4, 0.15) circle (1.3pt) node[right=1pt, black] {$x_3$};
  \node[inner sep=2pt] at ($ (0.4,0.15)+2.75*(-0.0349,0.9994) $) {$L_3$};
  \draw[blue!70!black, line width=1pt] (0,-2.95) -- (0,2.95);
  \node[above, black, inner sep=2pt] at (-0.18, 2.95) {$L$};
  \draw[red!75!black, thick, fill=red!10] (0, 2.60) circle (0.30);
  \draw[red!75!black, thick, fill=red!10] (0,-2.60) circle (0.30);
  \draw[-{Stealth[length=2.5mm]}, red!75!black, thin] (-1.4, 2.6) -- (-0.36, 2.6);
  \node[red!75!black, anchor=east, inner sep=2pt] at (-1.4, 2.6) {$U$};
  \draw[-{Stealth[length=2.5mm]}, red!75!black, thin] (-1.4,-2.6) -- (-0.36,-2.6);
  \node[red!75!black, anchor=east, inner sep=2pt] at (-1.4,-2.6) {$U$};
\end{tikzpicture}
\caption{Notation in the proof of \itemref{itm:thm:GRexten:1}.}
\label{fig:GR}
\end{figure}
\begin{proof}
    The extension is unique thanks to \cref{prop:pshfuncdetdense}. We shall prove the existence.

    \ref{itm:thm:GRexten:2} Thanks to the uniqueness of the extension, the problem is local, so we may assume that $X$ is a domain in $\mathbb{C}^n$ with $n>0$ and there is a non-zero holomorphic function $f$ on $X$ vanishing on $Z$. For each $\epsilon>0$, we claim that the function $\varphi_{\epsilon}$ defined by
    \[
    \varphi_{\epsilon}(x)\coloneqq \left\{
    \begin{aligned}
        \varphi(x)+\epsilon\log|f(x)|^2,&\quad x\in X\setminus Z;\\
        -\infty,&\quad x\in Z
    \end{aligned}
    \right.
    \]
    is plurisubharmonic on $X$. By \cref{def:psh}, it suffices to verify the case $n=1$. In this case, we may assume that $Z=\{0\}$.
    It is clear that $ \varphi_{\epsilon}\in \SH(X\setminus Z)$.
    It suffices to verify the sub-mean value inequality at $0$, which is immediate.

    Next observe that the net $\varphi_{\epsilon}$ is increasing as $\epsilon\searrow 0$ and $\varphi_{\epsilon}$ is locally uniformly bounded from above. It follows from \cref{prop:pshfunction_closedseq} that $\tilde{\varphi}\coloneqq \uscsup[\epsilon>0]\varphi_{\epsilon}\in \PSH(X)$. Moreover, $\tilde{\varphi}$ clearly extends $\varphi$. Note that \eqref{eq:GRextvarphiexp} follows from the construction.

    \ref{itm:thm:GRexten:1} We invite the reader to consult \cref{fig:GR} for the notations used in the proof.

    Thanks to \ref{itm:thm:GRexten:2}, it suffices to verify that $\varphi$ is locally bounded from above near each point of $Z$. The problem is local, so we may assume that $X$ is a domain in $\mathbb{C}^n$ with $n\geq 2$.

    Suppose the assertion fails. Then there exist $z\in Z$ and a sequence $(x_j)_{j\geq 1}$ in $X\setminus Z$ converging to $z$ such that
    \begin{equation}\label{eq:limphiinfty_temp1}
    \lim_{j\to\infty} \varphi(x_j)=\infty.
    \end{equation}
    Since $Z$ has codimension at least $2$\footnote{In fact, codimension at least~$1$ suffices for this step.}, we may choose a complex line $L$ through $z$ which meets $Z$ in a discrete set; after shrinking $X$ we may further assume that
    \[
        L\cap Z=\{z\}.
    \]
    Pick an open ball $B_n(z,r)\Subset X$. After adding a constant to $\varphi$, we may arrange $\varphi<0$ on $L\cap \partial B_n(z,r)$. Since $\varphi$ is upper semicontinuous, there is an open
    neighborhood $U\subseteq X\setminus Z$ of $L\cap \partial B_n(z,r)$ such that
    \[
        \varphi|_U<0.
    \]

    For each $j\geq 1$ we now choose a complex line $L_j$ through $x_j$, disjoint from $Z\cap B_n(z,r)$, with $L_j\to L$ as $j\to\infty$ (the convergence being in the Grassmannian of affine complex lines in $\mathbb{C}^n$). The existence of such $L_j$ is precisely where the assumption $\operatorname{codim} Z\geq 2$ enters: after replacing $B_n(z,r)$ by a slightly smaller ball if necessary, the closure of $Z\cap B_n(z,r)$ is analytic in a neighborhood of the closed ball. The projection to the space of directions is proper on this compact analytic set, so by Remmert's theorem its image is contained in a proper analytic subset of $\mathbb{P}^{n-1}$. Hence its complement contains directions arbitrarily close to that of $L$.

    For $j$ large enough we now have
    \[
        L_j\cap \partial B_n(z,r)\subseteq U,
    \]
    so applying the maximum principle to $\varphi|_{L_j\cap X}$ at $x_j$ yields $\varphi(x_j)<0$, contradicting \eqref{eq:limphiinfty_temp1}.
\end{proof}

\begin{lemma}\label{lma:invariantpshfunfinite}
    Let $\varphi\in \PSH((\Delta^*)^n)$ be an $(S^1)^n$-invariant plurisubharmonic function. Then $\varphi$ is finite everywhere.
\end{lemma}
Here
\[
\Delta^*\coloneqq \{z\in \mathbb{C}: 0<|z|<1\}
\]
denotes the punctured unit disk.
\begin{proof}
    It clearly suffices to handle the case $n=1$. In this case, by \cite[Theorem~2.12]{HK76} (see also \cite[Proposition~18.1.5]{FPT} applied in logarithmic coordinates), the map
    \[
    \log r\mapsto \int_{0}^{1} \varphi\left(r\mathrm{e}^{2\pi\mathrm{i}\theta}\right)\,\mathrm{d}\theta=\varphi(r)
    \]
    is a convex function of $\log r$.
    So, the set $\{r\in (0,1): \varphi(r)=-\infty\}$ is convex.
    But $\varphi$ is almost everywhere finite by \cref{prop:pshlocLp}. Since $\varphi$ is $S^1$-invariant, $\varphi|_{(0,1)}$ is almost everywhere finite.
    It follows from the convexity that it is everywhere finite. For another proof, see \cite[Corollary~18.1.6]{FPT}.
\end{proof}

\begin{proposition}[Kiselman's principle]\label{prop:Kis}\index{theorem!Kiselman's principle}
Let $\Omega\subseteq \mathbb{C}^m\times \mathbb{C}^n$ be a pseudoconvex domain, and let $\Omega'$ be its projection to $\mathbb{C}^m$. Assume that for each $z\in \Omega'$, the fiber
\[
    \Omega_z\coloneqq \bigl\{w\in\mathbb{C}^n: (z,w)\in \Omega \bigr\}
\]
has the form $E_z+\mathrm{i}\mathbb{R}^n$, where $E_z\subseteq \mathbb{R}^n$ is a nonempty open convex set. Let $\varphi\in \PSH(\Omega)$ be independent of the imaginary part of $w$, and define
\[
    \psi(z)=\inf_{t\in E_z}\varphi(z,t),\quad z\in \Omega'.
\]
Then either $\psi\equiv -\infty$ or $\psi\in \PSH(\Omega')$.
\end{proposition}
See \cite[Chapter~I, Section~7.B, Theorem~(7.5) and its proof, p.~57]{DemBook} for the proof as well as the notion of pseudoconvex domains. See also \cite[Corollary~11.2.6]{FPT} for the pseudoconvex tube case.

\begin{lemma}\label{lma:gl}
    Let $\Omega\subseteq \mathbb{C}^n$ be a domain and $\Omega'\subseteq \Omega$ be a subdomain. Consider $\varphi\in \PSH(\Omega)$ and $\psi\in \PSH(\Omega')$. Assume that
    \begin{equation}\label{eq:gluecond}
    \varliminf_{\substack{\Omega'\ni y\to x,\\ \psi(y)\neq -\infty}}\bigl(\varphi(y)-\psi(y)\bigr)\geq 0
    \end{equation}
    for any $x\in \Omega\cap \partial \Omega'$. Define
    \[
    \eta(z)=\left\{
    \begin{aligned}
        \varphi(z)\lor \psi(z),&\quad \text{if }z\in \Omega',\\
        \varphi(z),&\quad \text{if }z\in \Omega\setminus\Omega'.
    \end{aligned}
    \right.
    \]
    Then $\eta\in \PSH(\Omega)$.
\end{lemma}
This is essentially \cite[Proposition~1.30]{GZ17}. See also \cite[Proposition~8.2.5]{FPT}. Since the statement in the reference is slightly misleading for our purposes, we include the proof for clarification.
\begin{figure}[ht]
\centering
\begin{tikzpicture}[scale=0.85, line width=0.6pt]
  \draw[thick, fill=black!4] (0,0) ellipse[x radius=4.0, y radius=2.5];
  \node at (4.4,2.0) {$\Omega$};
  \draw[thick] (1.000, 2.421)
    .. controls (1.000, 0.800) and (0.400, -0.776) .. (0.400, -2.487);
  \node at (-2.6, 0.0) {$\Omega'$};
  \begin{scope}
    \clip (0.7, 0) circle (0.85);
    \fill[pattern=north east lines, pattern color=red!50] (-1, -1) rectangle (0.7, 1);
  \end{scope}
  \draw[red!75!black, thick] (0.7, 0) circle (0.85);
  \node[red!75!black, anchor=west, inner sep=2pt] at (1.7, 1.05) {$U_x$};
  \draw[red!75!black, thin] (1.7, 1.05) -- (1.30, 0.50);
  \node[red!75!black, anchor=east, inner sep=2pt] at (-0.6, -1.05) {$U_x\cap\Omega'$};
  \draw[red!75!black, thin] (-0.6, -1.05) -- (-0.05, -0.55);
  \fill (0.7, 0) circle (1.6pt);
  \node[anchor=west, inner sep=2pt] at (0.78, 0) {$x$};
  \fill (0.25, 0.40) circle (1.6pt);
  \draw[-{Stealth[length=2mm]}, thin] (-0.55, 1.45) -- (0.20, 0.50);
  \node[anchor=east, inner sep=2pt] at (-0.55, 1.45) {$y$};
\end{tikzpicture}
\caption{Gluing procedure}
\label{fig:glue}
\end{figure}
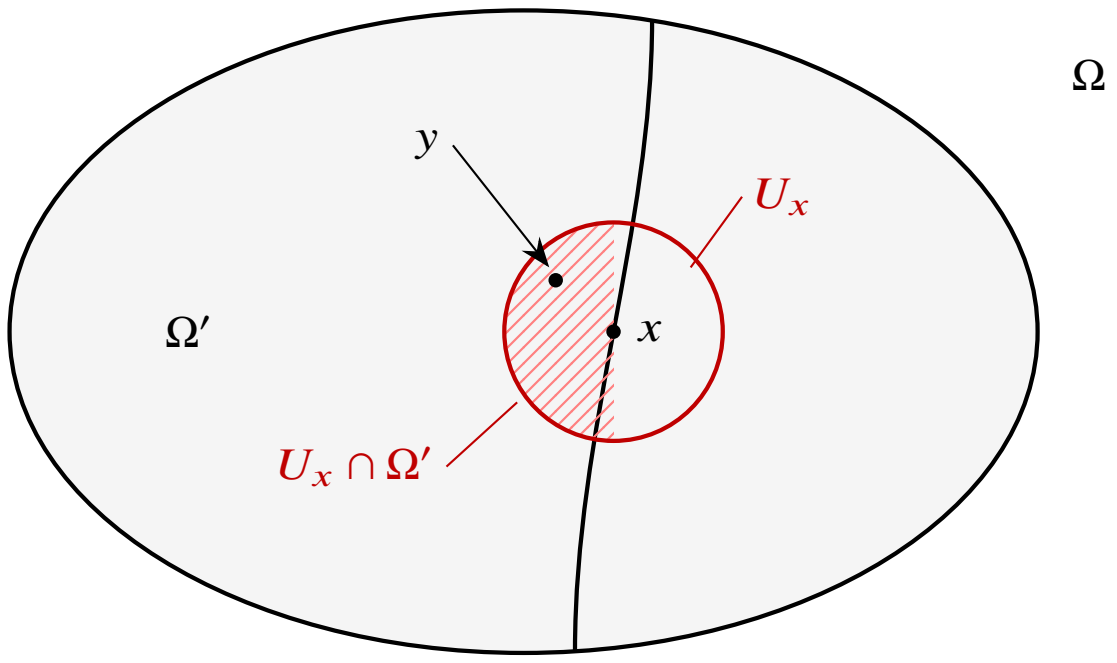
\begin{proof}
See \cref{fig:glue} for the notations used in the proof.

    Take $\epsilon>0$.
    We first define
    \[
    \eta_{\epsilon}(z)=\begin{cases}
        \varphi(z)\lor (\psi(z)-2\epsilon),& \textup{if }z\in \Omega';\\
        \varphi(z),& \textup{if }z\in \Omega\setminus\Omega'.
    \end{cases}
    \]
    We claim that
    \[
    \eta_{\epsilon}\in \PSH(\Omega).
    \]
    By our assumption, for each $x\in \Omega\cap \partial \Omega'$, we can find an open neighborhood $U_x\subseteq \Omega$ such that for any $y\in U_x\cap \Omega'$, we have $\varphi(y)\geq \psi(y)-\epsilon$.
    Therefore, by taking a union of the $U_x$'s for various $x$, there is an open neighborhood $U$ of $\Omega\cap \partial \Omega'$ such that
    \[
    \varphi(y)\geq \psi(y)-\epsilon,\quad \forall y\in U\cap \Omega'.
    \]
    Since $U$ contains $\Omega\cap\partial\Omega'$, the set $(\Omega\setminus \Omega')\cup U$ is open in $\Omega$. On this open set we have $\eta_{\epsilon}=\varphi$, hence $\eta_{\epsilon}$ is plurisubharmonic there. It is plurisubharmonic on $\Omega'$ by \cref{prop:pshfunction_closedseq}. These two descriptions agree on the overlap $U\cap\Omega'$, so our claim follows from the local nature of plurisubharmonicity.

    Next we observe that as $\epsilon$ decreases to $0$, the functions $\eta_{\epsilon}$ increase to $\eta$. The family $(\eta_\epsilon)_{0<\epsilon\leq 1}$ is locally uniformly bounded from above. This is clear away from $\Omega\cap\partial\Omega'$. Near $\Omega\cap\partial\Omega'$, applying the previous estimate with $\epsilon=1$ gives an open neighborhood $U_1$ of $\Omega\cap\partial\Omega'$ such that $\psi\leq \varphi+1$ on $U_1\cap\Omega'$, hence $\eta_\epsilon\leq \varphi+1$ on $U_1$ for all $0<\epsilon\leq 1$. Therefore, $\eta^*\in \PSH(\Omega)$ by \cref{prop:pshfunction_closedseq}.

    It remains to observe that $\eta$ is upper semicontinuous. This is only non-trivial on the boundary of $\Omega'$. Take $x\in \Omega\cap \partial \Omega'$ and let $(y_i)_{i>0}$ be a sequence in $\Omega'$ with limit $x$. We first show that
    \begin{equation}\label{eq:limsuppsiphitemp1}
    \varlimsup_{i\to\infty} \psi(y_i)\leq \varphi(x).
    \end{equation}
    We may assume that $\psi(y_i)\neq -\infty$ for all $i>0$ and the left-hand side of \eqref{eq:limsuppsiphitemp1} is not $-\infty$. Then we can compute
    \[
    \varlimsup_{i\to\infty} \psi(y_i)\leq \varlimsup_{i\to\infty} \psi(y_i)+\varliminf_{i\to\infty} \bigl(\varphi(y_i)-\psi(y_i)\bigr)\leq \varlimsup_{i\to\infty}\varphi(y_i)\leq \varphi(x).
    \]
    Combining this with the upper semicontinuity of $\varphi$, we get
    \[
        \varlimsup_{i\to\infty}\bigl(\varphi(y_i)\lor\psi(y_i)\bigr)\leq \varphi(x)=\eta(x).
    \]
    This proves upper semicontinuity across $\Omega\cap\partial\Omega'$.
    Therefore, $\eta=\eta^*\in \PSH(\Omega)$.
\end{proof}

Note that the boundary condition \eqref{eq:gluecond} enters the proof only through \eqref{eq:limsuppsiphitemp1}, therefore, we can formulate a variant of \cref{lma:gl}:
\begin{corollary}\label{cor:glbdry}
    Let $\Omega\subseteq \mathbb{C}^n$ be a domain and $\Omega'\subseteq \Omega$ be a subdomain\footnote{In \cite[Proposition~1.30]{GZ17}, Guedj--Zeriahi assumed in addition that $\Omega'$ is relatively compact in $\Omega$. This assumption is not necessary for the proof, and the applications of this proposition in \cite{GZ17} usually do not satisfy this condition. See also the gluing statement in \cite[Proposition~8.2.5]{FPT}.}. Consider $\varphi\in \PSH(\Omega)$ and $\psi\in \PSH(\Omega')$. Assume that
    \[
        \varlimsup_{\Omega'\ni y\to x}\psi(y)\leq \varphi(x)
    \]
    for any $x\in \Omega\cap\partial\Omega'$. Define
    \[
    \eta(z)=\begin{cases}
        \varphi(z)\lor \psi(z),&\quad \textup{if }z\in \Omega';\\
        \varphi(z),&\quad \textup{if }z\in \Omega\setminus\Omega'.
    \end{cases}
    \]
    Then $\eta\in \PSH(\Omega)$.
\end{corollary}

Using the ideas of \cref{lma:gl} and the solution of the homogeneous Monge--Ampère equations on a strong pseudoconvex domain, one can prove the following:
\begin{proposition}\label{prop:BT_balayage}
Let $\Omega \subseteq \mathbb{C}^n$ be a domain with $n\geq 1$, let $D \Subset \Omega$ be a strongly pseudoconvex domain, and let $\varphi \in \PSH(\Omega)$ be locally bounded.
Then there exists a unique locally bounded $\psi\in \PSH(\Omega)$
such that
\[
\begin{cases}
\left(\ddc\psi\right)^n = 0& \textup{on } D,\\
\psi = \varphi & \textup{on } \Omega\setminus D.
\end{cases}
\]
Moreover, $\psi\geq \varphi$.
\end{proposition}
See \cite[Proposition~9.1]{BT82} for the proof, and \cite[Corollary~9.6.4]{FPT} for balayage on balls. We sometimes say $\psi$ is the \emph{balayage}\index{balayage} of $\varphi$ on $D$.

Let us also recall the following domination principle:
\begin{proposition}\label{prop:dom_prin_loc}
    Let $\Omega\subseteq \mathbb{C}^n$ be a bounded domain.
    Let $\varphi,\psi\in \PSH(\Omega)$ be locally bounded. Assume that
    \begin{enumerate}
        \item\label{itm:prop:dom_prin_loc:1} for each $z\in \partial \Omega$, we have
            \[
                \varliminf_{\Omega\ni w\to z} \Bigl(\psi(w)-\varphi(w)\Bigr)\geq 0;
            \]
        \item\label{itm:prop:dom_prin_loc:2} $\varphi\leq \psi$ a.e.\ with respect to $(\ddc \psi)^n$,
    \end{enumerate}
    Then $\varphi\leq \psi$.
\end{proposition}
This is \cite[Corollary~3.31]{GZ17}.\footnote{The statement in the reference is flawed, but the proof is correct. See also \cite[Corollary~9.2.9]{FPT}.}

\section{Plurifine topology}\label{sec:plurifine}

In this section, we introduce the notion of plurifine topology. Unlike other sections in this chapter, this section contains full details. This is mainly due to the unfortunate omissions and numerous minor problems in the foundational paper \cite{BT87}. In particular, the key theorem \cref{thm:pf_basis} was claimed by Bedford--Taylor without proof. Our presentation largely follows \cite{EMW06}. See also \cite[Proposition~10.3.2 and Section~B.10]{FPT}. We also include enough details so that this section is readable for those who are not familiar with the classical potential theory.

Very recently, El Kadiri--Fuglede published a book \cite{EKF} about classical fine potential theory, containing essentially all results in this section. The interested reader may consult their book for further details. See also \cite[Chapters~5 and~10]{FPT}.

\subsection{Plurifine topology on domains}
\label{subsec:plurifine}
Let $\Omega\subseteq \mathbb{C}^n$ ($n\in \mathbb{N}$) be a domain.

\begin{definition}\label{def:pftopologydomain}
  The \emph{plurifine topology}\index{plurifine topology} on $\Omega$ is the weakest topology making all $\mathbb{R}$-valued plurisubharmonic functions on $\Omega$ continuous.
\end{definition}
We want to distinguish the Euclidean topology from the plurifine topology.
In the whole book, topological notions without adjectives refer to those with respect to the Euclidean topology. We include the symbol $\mathcal{F}$ in order to denote those with respect to the plurifine topology. For example, we will say $\mathcal{F}$-open subset, $\mathcal{F}$-neighborhood, $\mathcal{F}$-closure, etc.
The $\mathcal{F}$-closure of a set $E\subseteq \Omega$ will be denoted by $\overline{E}^{\mathcal{F}}$.

Recall that in the whole book, we follow Bourbaki's convention: a neighborhood is not necessarily open. Similarly, an $\mathcal{F}$-neighborhood is not necessarily $\mathcal{F}$-open.

A priori, we should include $\Omega$ into the notations as well, but as we will see shortly in \cref{cor:pf_compatible}, this is usually unnecessary.

When $n=0$, the domain $\Omega$ is a singleton and all assertions in this subsection are trivial. In the proofs below we therefore assume $n>0$.

\begin{proposition}\label{prop:pf_finer}
  The plurifine topology on $\Omega$ is finer than the Euclidean topology.
\end{proposition}
\begin{proof}
  It suffices to show that the unit ball $\{z\in \mathbb{C}^n:|z|<1\}$ is $\mathcal{F}$-open. This follows from the observation that this set can be written as
  \[
      \{\psi<0\} \quad\textup{with}\quad \psi(z)\coloneqq (\log |z|) \lor (-1).
  \]
  This completes the proof.
\end{proof}

\begin{example}\label{ex:pfopenpsh}
    Let $\varphi\in \PSH(\Omega)$ and $C\in \mathbb{R}$. Then the sets $\{\varphi>C\}$ and $\{\varphi<C\}$ are both $\mathcal{F}$-open.

    In fact, the latter case follows from \cref{prop:pf_finer}, while the former follows from the observation
    \[
    \{\varphi>C\}=\{\varphi\lor (C-1)>C\}.
    \]
\end{example}

\begin{definition}\label{def:pshfunctions:L691}
  A subset $E\subseteq \Omega$ is \emph{thin}\index{thin subset}\footnote{A more proper name would be \emph{plurithin}. But since we will never need the classical notion of thin sets à la Cartan in this book, we prefer omitting the prefix \emph{pluri-}.} at $x\in \Omega$ if one of the following conditions holds:
  \begin{enumerate}
    \item\label{itm:chapter-pshfunctions:L698:1} $x\not\in \overline{E}$;
    \item\label{itm:chapter-pshfunctions:L698:2} $x\in \overline{E}$ and there is an open neighborhood $U\subseteq \Omega$ of $x$ and $\varphi\in \PSH(U)$ such that
    \[
      \varlimsup_{(U\cap E)\setminus\{x\} \ni y\to x} \varphi(y)<\varphi(x).
    \]
  \end{enumerate}
  We say $E$ is \emph{thin} if it is thin at all $x\in \Omega$.
\end{definition}
\begin{remark}\label{rmk:cstthindef}
    In the second case, we can always arrange that
    \[
    \varphi|_{(E\setminus \{x\})\cap U}=0.
    \]
    In fact, we may assume that $\varphi\leq 0$ and $C<0$ is such that
    \[
    \varlimsup_{(U\cap E)\setminus\{x\} \ni y\to x} \varphi(y)<C<\varphi(x).
    \]
    We let
    \[
        \psi=(-C)^{-1}(\varphi\lor C)+1.
    \]
    Then $\psi$ satisfies our requirements for a smaller $U$.
\end{remark}

In the second case, the function $\varphi$ can be very much improved.
\begin{proposition}[Bedford--Taylor]\label{prop:BTthin}\index{theorem!Bedford--Taylor base theorem}
  Consider a set $E\subseteq \Omega$ and $x\in \overline{E}$. Assume that $E$ is thin at $x$. Then there is $\varphi\in \PSH(\mathbb{C}^n)$ with the following properties:
  \begin{enumerate}
    \item\label{itm:prop:BTthin:1} $\varphi$ is locally bounded outside $x$;
    \item\label{itm:prop:BTthin:2} $\varphi(x)>-\infty$;
    \item\label{itm:prop:BTthin:3}
    \[
    \lim_{E\setminus\{x\} \ni y\to x} \varphi(y)=-\infty.
    \]
  \end{enumerate}
\end{proposition}
\begin{proof}\footnote{The original argument in \cite[Proposition~10.2]{BT82} was quite intriguing: Neither the auxiliary functions $\varphi_j$'s nor the simple computations were correct. However, I believe that Bedford--Taylor had a correct proof in mind. Something more than a typo, but not yet a mistake, could be properly called a \emph{thinkpo}, a terminology invented by R.~Berman. For a proof of the local base theorem, see also \cite[Proposition~10.3.2]{FPT}.}
By \cref{rmk:cstthindef}, there is an open neighborhood $U\subseteq \Omega$ of $x$ and $\psi\in \PSH(U)$ such that
\[
    \psi|_{U\cap (E\setminus\{x\})}=-1 \quad\textup{and}\quad \psi(x)=0.
\]
Without loss of generality we may assume $x=0$ and that $U=B_n(0,1)$ is the open unit ball.

Since $\psi$ is upper semicontinuous at $0$ with $\psi(0)=0$, we may choose a strictly decreasing sequence $(\delta_j)_{j\geq 1}\subseteq (0,\mathrm{e}^{-1})$ such that
\begin{equation}\label{eq:BTthin_choosedelta}
    \psi(y)<2^{-j-2}\quad\textup{whenever } |y|<\delta_j.
\end{equation}
In particular $|\log\delta_j|>1$ for every $j\geq 1$, so the series of constants
\begin{equation}\label{eq:BTthin_summable}
    \sum_{j=1}^\infty \frac{2^{-j-1}}{|\log\delta_j|} \leq \sum_{j=1}^\infty 2^{-j-1}=\frac{1}{2}
\end{equation}
converges. We further set
\[
    \gamma_j\coloneqq \exp\bigl(2^{j+1}\log\delta_j\bigr)=\delta_j^{2^{j+1}}\in(0,1),
\]
which is also decreasing in $j$.

For each $j\geq 1$ define
\[
    \varphi_j(z)\coloneqq
    \begin{cases}
        \frac{2^{-j-1}}{|\log\delta_j|}\log|z|\lor\bigl(\psi(z)-2^{-j}\bigr), &\textup{if } |z|<\delta_j,\\
        \frac{2^{-j-1}}{|\log\delta_j|}\log|z|,&\textup{if } |z|\geq \delta_j.
    \end{cases}
\]
On the boundary $|z|=\delta_j$, the logarithmic term tends to $-2^{-j-1}$. By \eqref{eq:BTthin_choosedelta},
\[
    \varlimsup_{|z|<\delta_j,\ z\to z_0}(\psi(z)-2^{-j})\leq 2^{-j-2}-2^{-j}<-2^{-j-1}
\]
for every $z_0$ with $|z_0|=\delta_j$. Hence the boundary condition in \cref{lma:gl} is satisfied. We obtain that $\varphi_j\in \PSH(\mathbb{C}^n)$.

By construction, $\varphi_j\leq 0$ on $U$, $\varphi_j(0)=-2^{-j}$, and for $z\in E\setminus\{0\}$ with $|z|<\gamma_j$ we have $\varphi_j(z)\leq -1$, since
\[
    \frac{2^{-j-1}}{|\log\delta_j|}\log|z|<\frac{2^{-j-1}}{|\log\delta_j|}\,(2^{j+1}\log\delta_j)=-1
\]
and, on $E\setminus\{0\}$, $\psi(z)-2^{-j}=-1-2^{-j}<-1$.

We then define
\[
    \varphi\coloneqq \sum_{j=1}^{\infty}\varphi_j.
\]
On the unit ball $U$ each $\varphi_j\leq 0$, so the partial sums $S_N\coloneqq \sum_{j=1}^N\varphi_j$ form a decreasing sequence of plurisubharmonic functions with
\[
    S_N(0)=-\sum_{j=1}^N 2^{-j}\to -1>-\infty,\quad \textup{as } N\to\infty
\]
so $\varphi|_U\in \PSH(U)$ by \cref{prop:pshfunction_closedseq}. On $\mathbb{C}^n\setminus B_n(0,\delta_1)$ each $\varphi_j$ equals $\frac{2^{-j-1}}{|\log\delta_j|}\log|z|$, and \eqref{eq:BTthin_summable} shows $S_N$ converges uniformly on every compact subset to $\bigl(\sum_{j\geq 1}\tfrac{2^{-j-1}}{|\log\delta_j|}\bigr)\log|z|$, a plurisubharmonic limit. Both formulas agree on $U\setminus B_n(0,\delta_1)$, so altogether $\varphi\in \PSH(\mathbb{C}^n)$.

We now verify the three claims.

  \emph{Property \ref{itm:prop:BTthin:2}}. This is the equality $\varphi(0)=-\sum_{j\geq 1}2^{-j}=-1>-\infty$ obtained above.

  \emph{Property \ref{itm:prop:BTthin:3}}. Fix $j_0\geq 1$. For $z\in E\setminus\{0\}$ with $|z|<\gamma_{j_0}$, the bound $\varphi_j(z)\leq -1$ holds for every
  $j\leq j_0$ (since the $\gamma_j$ are decreasing), while $\varphi_j(z)\leq 0$ for $j>j_0$. Hence $\varphi(z)\leq -j_0$, and
  letting $j_0\to \infty$ yields
  \[
    \varlimsup_{E\setminus \{0\} \ni y\to 0}\varphi(y)=-\infty.
  \]

  \emph{Property \ref{itm:prop:BTthin:1}.} Fix any $z_0\in \mathbb{C}^n$ with $z_0\neq 0$ and choose $J\geq 1$ such that $\delta_J<|z_0|/2$. On the open set
  $W=\{w\in \mathbb{C}^n:|w|>|z_0|/2\}$, all $\varphi_j$ with $j\geq J$ lie in the outer regime, and \eqref{eq:BTthin_summable} shows $\sum_{j\geq J}\varphi_j(w)
   =\bigl(\sum_{j\geq J}\tfrac{2^{-j-1}}{|\log\delta_j|}\bigr)\log|w|$ is locally bounded on $W$. The remaining finitely many terms $\varphi_1,\ldots,\varphi_{J-1}$ are plurisubharmonic on $\mathbb{C}^n$, hence locally bounded above on $W$, while $\varphi_j(w)\geq \tfrac{2^{-j-1}}{|\log\delta_j|}\log|w|$ keeps
  them locally bounded below as well. Therefore $\varphi$ is locally bounded on $W$. Since $z_0\neq 0$ was arbitrary, $\varphi$ is locally bounded on $\mathbb{C}^n\setminus\{0\}$, which is property~\ref{itm:prop:BTthin:1}.
\end{proof}

\begin{lemma}\label{lma:unionthin}
    Let $E_1,E_2\subseteq \Omega$. Assume that $E_1$, $E_2$ are both thin at $x\in \Omega$. Then so is $E_1\cup E_2$.
\end{lemma}
\begin{proof}
    We may assume that $x\in \overline{E_1}\cap \overline{E_2}$. Take an open neighborhood $U\subseteq \Omega$ of $x$ and $\varphi_1,\varphi_2\in \PSH(U)$ such that
    \[
        \varlimsup_{y\to x,y\in E_i\setminus\{x\}} \varphi_i(y)<\varphi_i(x),\quad i=1,2.
    \]
    Then $\varphi_1+\varphi_2\in \PSH(U)$ and
    \[
    \varlimsup_{y\to x,y\in (E_1\cup E_2)\setminus\{x\}} (\varphi_1+\varphi_2)(y)<\varphi_1(x)+\varphi_2(x).
    \]
    In particular, $E_1\cup E_2$ is thin at $x$.
\end{proof}

\begin{theorem}[H.~Cartan]\label{thm:Cartan}\index{theorem!H.~Cartan's coherence theorem}
  Consider a set $E\subseteq \Omega$ and $x\in E$.
  Then the following are equivalent:
  \begin{enumerate}
    \item\label{itm:thm:Cartan:1} $E$ is an $\mathcal{F}$-neighborhood of $x$;
    \item\label{itm:thm:Cartan:2} $\Omega\setminus E$ is thin at $x$.
  \end{enumerate}
\end{theorem}
\begin{proof}
  \ref{itm:thm:Cartan:2} $\implies$ \ref{itm:thm:Cartan:1}. Assume \ref{itm:thm:Cartan:2}. We may assume that $x\in \overline{\Omega\setminus E}$. Otherwise, our assertion follows from \cref{prop:pf_finer}.

  By \cref{prop:BTthin}, there is an open neighborhood $U$ of $x$ in $\Omega$ and $\varphi\in \PSH(\mathbb{C}^n)$ such that
  \[
      \varphi(x)>\sup_{y\in U\cap (\Omega\setminus E)} \varphi(y)\eqqcolon \lambda.
  \]
  Let $F=\{y\in \Omega: \varphi(y)>\lambda\}$. Then $x\in F$ and $F$ is $\mathcal{F}$-open by \cref{ex:pfopenpsh}. Moreover, $U\cap F\subseteq E$. By \cref{prop:pf_finer}, we conclude \ref{itm:thm:Cartan:1}.

  \ref{itm:thm:Cartan:1} $\implies$ \ref{itm:thm:Cartan:2}. Assume \ref{itm:thm:Cartan:1}. We may always replace $E$ by smaller $\mathcal{F}$-neighborhoods of $x$. In particular, we may assume that $E$ has the following form
  \[
      \{y\in U:\varphi_1(y)>\lambda_1,\ldots,\varphi_m(y)>\lambda_m\},
  \]
  where $U\subseteq \Omega$ is an open neighborhood of $x$, and $\varphi_1,\ldots,\varphi_m$ are $\mathbb{R}$-valued psh functions on $\Omega$, and $\lambda_1,\ldots,\lambda_m\in \mathbb{R}$. Since $x\in U$, the set $\Omega\setminus U$ is thin at $x$. Moreover,
    \[
        \Omega\setminus E=(\Omega\setminus U)\cup \bigcup_{i=1}^m\{y\in U:\varphi_i(y)\leq\lambda_i\}.
    \]
    By \cref{lma:unionthin}, it is enough to prove that each set $\{y\in U:\varphi_i(y)\leq\lambda_i\}$ is thin at $x$. Thus we may assume $m=1$. In this case, $\Omega\setminus E$ is clearly thin at $x$.
\end{proof}

\begin{theorem}\label{thm:pf_basis}
  A base of the plurifine topology on $\Omega$ is given by sets of the following form:
  \begin{equation}\label{eq:basis_fine}
    \left\{x\in U: \varphi(x)>0\right\},
  \end{equation}
  where $U\subseteq \Omega$ is an open subset and $\varphi\in \PSH(U)$.
\end{theorem}

\begin{proof}
  Observe that sets of the form \eqref{eq:basis_fine} are $\mathcal{F}$-open.\footnote{This is not entirely obvious by definition, as $\varphi$ is not defined on the whole $\Omega$.} By \cref{thm:Cartan}, it suffices to show its complement in $\Omega$ is thin at each point of \eqref{eq:basis_fine}, which is obvious.

  Now consider $x\in \Omega$ and an $\mathcal{F}$-open neighborhood $V\subseteq \Omega$ of $x$. We want to find a set of the form \eqref{eq:basis_fine} contained in $V$ and containing $x$.

  Write $E=\Omega\setminus V$. In case $x\in \Int V$, there is nothing to prove. So we may assume that $x\in \overline{E}$. By \cref{thm:Cartan}, $E$ is thin at $x$. By definition, there is an open neighborhood $U\subseteq \Omega$ of $x$ and $\varphi\in \PSH(U)$ such that
  \[
    \varlimsup_{y\to x,y\in U\cap (E\setminus\{x\})}\varphi(y)<\varphi(x).
  \]
  We may assume that $\varphi|_{E\cap U}\leq 0<\varphi(x)$. Then the set $\{y\in U:\varphi(y)>0\}$ suffices for our purpose.
\end{proof}

\begin{remark}
    Recall that in general, an $\mathcal{F}$-open set is \emph{not} a \emph{countable} union of sets of the form \eqref{eq:basis_fine}. In fact, an $\mathcal{F}$-open set is not a Borel set in general.
    See \cite{EK23} and \cite[Example~10.3.5]{FPT} for concrete examples.
\end{remark}

\begin{corollary}\label{cor:pf_compatible}
  Let $\Omega_1\subseteq \Omega_2\subseteq \mathbb{C}^n$ be two non-empty open subsets. Then the plurifine topology on $\Omega_1$ is the same as the subspace topology induced from the plurifine topology on $\Omega_2$.
\end{corollary}
In particular, when we talk about an $\mathcal{F}$-open set $U$ in $\mathbb{C}^n$, we no longer have to specify the domain $\Omega\supseteq U$.

\begin{corollary}\label{cor:compactnhformbase}
  Let $\Omega\subseteq \mathbb{C}^n$ be an $\mathcal{F}$-open subset and $x\in \Omega$. Then $x$ has a compact $\mathcal{F}$-neighborhood contained in $\Omega$.
\end{corollary}
\begin{proof}
  By \cref{thm:pf_basis}, it suffices to prove the assertion after replacing $\Omega$ by a smaller $\mathcal F$-open neighborhood of $x$. Thus we may assume that
    \[
        \Omega=\{y\in U:\varphi(y)>0\}
    \]
    for some open set $U\subseteq\mathbb C^n$ and some $\varphi\in\PSH(U)$.

  Take a compact neighborhood $K$ of $x$ in $U$. Now $\{y\in K: \varphi(y)\geq \varphi(x)/2\}$ is a compact $\mathcal{F}$-neighborhood of $x$ contained in $\Omega$.
\end{proof}

\begin{corollary}\label{cor:holomappfcont}
    Let $\Omega\subseteq \mathbb{C}^n$, $\Omega'\subseteq \mathbb{C}^{n'}$ be two domains and $F\colon \Omega'\rightarrow \Omega$ be a holomorphic map. Then $F$ is $\mathcal{F}$-continuous.
\end{corollary}

\begin{proof}
    It suffices to show that the inverse image $F^{-1}(U)$ of each $\mathcal{F}$-open subset $U\subseteq \Omega$ is $\mathcal{F}$-open.
    By \cref{thm:pf_basis}, after possibly shrinking $\Omega$ and $\Omega'$, we may assume that $U$ has the form $\{x\in \Omega:\psi(x)>0\}$, where $\psi \in \PSH(\Omega)$. Clearly, we may assume that $F^{-1}(U)$ is non-empty. Then $F^*\psi\in \PSH(\Omega')$ by \cref{prop:pullbackpsh}, we find that
    \[
        F^{-1}(U)=\{y\in \Omega': F^*\psi(y)>0\}
    \]
    is $\mathcal{F}$-open.
\end{proof}

\subsection{Plurifine topology on manifolds}\label{subsec:pftmfd}
Let $X$ be a complex manifold.
\begin{definition}\label{def:pftopologygeneral}
    The \emph{plurifine topology}\index{plurifine topology} on $X$ is the topology with a base consisting of sets of the form $F^{-1}(V)$, where $U\subseteq X$ is an open subset, $F\colon U\rightarrow \Omega$ is a biholomorphic map, $\Omega\subseteq \mathbb{C}^n$ is a domain for some $n\in \mathbb{N}$, and $V\subseteq \Omega$ is $\mathcal{F}$-open.
\end{definition}
Note that these sets form a topological base thanks to \cref{cor:holomappfcont}. Moreover, it also follows from \cref{cor:holomappfcont} that the plurifine topologies on domains defined in \cref{def:pftopologygeneral} and in \cref{def:pftopologydomain} coincide.

We refer to \cref{def:qpsh} for the notion of quasi-plurisubharmonic functions.
\begin{proposition}\label{prop:pshfunFcont}
    Let $\varphi\in \QPSH(X)$. Then $\varphi$ is $\mathcal{F}$-continuous.
\end{proposition}
\begin{proof}
    It suffices to show that $\varphi|_{\{\varphi\neq -\infty\}}$ is $\mathcal{F}$-continuous, since $\{\varphi=-\infty\}$ is $\mathcal{F}$-closed.

    The problem is local, so we may assume that $X\subseteq \mathbb{C}^n$ is a domain and $\varphi=\psi+g$, where $\psi\in \PSH(X)$ and $g\in C^{\infty}(X)$ and $|g|\leq C$ for some $C>0$. For an open interval $(a,b)\subseteq \mathbb{R}$, it suffices to show that
    \[
        U\coloneqq \left\{ x\in X:a<\varphi(x)<b\right\}=\left\{ x\in X:a-g(x)<\psi(x)<b-g(x)\right\}
    \]
    is $\mathcal{F}$-open. Take $x\in U$. We can find an open neighborhood $V$ of $x$ in $U$ such that
    \[
    \sup_{y\in V} \left( a-g(y)\right)<\psi(x)< \inf_{y\in V} \left( b-g(y) \right).
    \]
    Therefore,
    \[
    \left\{z\in V:\sup_{y\in V} \left(a-g(y)\right)<\psi(z)<\inf_{y\in V} \left(b-g(y)\right)\right\}
    \]
    is an $\mathcal{F}$-open neighborhood of $x$ in $U$. We conclude that $U$ is $\mathcal{F}$-open.
\end{proof}

\begin{corollary}\label{cor:diffpshFopen}
    Let $\varphi,\psi\in \QPSH(X)$. Then the set
    \[
    \{x\in X: \varphi(x)>\psi(x)\}
    \]
    is $\mathcal{F}$-open.
\end{corollary}
\begin{proof}
    It suffices to show that for any $x\in X$ such that $\varphi(x)>\psi(x)$, the same holds on an $\mathcal{F}$-neighborhood $U$ of $x$. Observe that $\varphi(x)\neq -\infty$.
    If $\psi(x)\neq -\infty$, then it suffices to apply \cref{prop:pshfunFcont}. Otherwise, take
    \[
    U\coloneqq \{y\in X: \varphi(y)>\varphi(x)-1\}\cap \{y\in X: \psi(y)<\varphi(x)-1\}.
    \]
    This is an $\mathcal{F}$-open neighborhood of $x$ contained in $\{\varphi>\psi\}$.
\end{proof}

\begin{lemma}\label{lma:pshfunfinitelocuspfdense}
    Let $Z\subseteq X$ be a pluripolar subset. Then
    \[
    \overline{X\setminus Z}^{\mathcal{F}}=X.
    \]
\end{lemma}
\begin{proof}
    The problem is local, so we may assume that $X$ is a domain in $\mathbb{C}^n$ and that there is $\varphi\in \PSH(X)$ with $Z\subseteq\{\varphi=-\infty\}$. It suffices to show that $\{\varphi>-\infty\}$ is $\mathcal{F}$-dense.

    Let $x\in X$ be a point with $\varphi(x)=-\infty$ and $U\subseteq X$ be an $\mathcal{F}$-open neighborhood of $x$ in $X$. We need to show that $U\cap \{\varphi>-\infty\}\neq \varnothing$.

    Thanks to \cref{thm:pf_basis}, after shrinking $U$, we may assume that there is $\psi\in \PSH(X)$ such that $U=\{\psi>0\}$. Observe that $U$ is not a pluripolar set: Otherwise, $\psi\leq 0$ almost everywhere by \cref{prop:ppsnull}, and hence everywhere by \cref{prop:pshfuncdetdense}. So $\varphi|_U\not\equiv -\infty$. We conclude that $U\cap \{\varphi>-\infty\}\neq \varnothing$, as desired.
\end{proof}

\begin{corollary}\label{cor:diffsupinfindeppluripolar}
    Let $\varphi,\psi\in \QPSH(X)$. Set
    \[
    W=\left\{x\in X:\varphi(x)=-\infty\right\}\textup{ or }W=\left\{x\in X:\psi(x)=-\infty\right\}.
    \]
    Then for any pluripolar set $Z\subseteq X$, we have
    \[
    \sup_{X\setminus W}(\varphi-\psi) =\sup_{X\setminus (W\cup Z)}(\varphi-\psi),\quad \inf_{X\setminus W}(\varphi-\psi) =\inf_{X\setminus (W\cup Z)}(\varphi-\psi).
    \]
    In particular, taking $\psi=0$, we find that
    \[
    \sup_{X\setminus Z}\varphi=\sup_X \varphi.
    \]
\end{corollary}
\begin{proof}
    This is an immediate consequence of \cref{lma:pshfunfinitelocuspfdense} and \cref{prop:pshfunFcont}.
\end{proof}

In the literature about pluripotential theory, one often finds the careless expressions like $\sup_X (\varphi-\psi)$. The issue is that $\varphi-\psi$ is not defined everywhere, and hence this expression does not make sense if you read it literally. \cref{cor:diffsupinfindeppluripolar} tells you that you do not have to worry too much about the details on a pluripolar set. In other words, $\sup$ and $\inf$ could always be understood as a kind of essential supremum and essential infimum modulo pluripolar sets.

There is a convenient way to fix this issue in the literature --- just replace the suprema by quasi-suprema:
\begin{definition}\label{def:quasi_sup}
    Let $Z\subseteq X$ be a pluripolar set, and $f\colon X\setminus Z\rightarrow [-\infty,\infty]$ be a function. We define the \emph{quasi-supremum}\index{quasi-supremum} and the \emph{quasi-infimum}\index{quasi-infimum} of $f$ as follows:
    \[
    \begin{split}
    \qsup_X f\coloneqq &\inf\left\{\sup_{X\setminus W} f:W\supseteq Z, W \textup{ is pluripolar}\right\},\\
    \qinf_X f\coloneqq &\sup\left\{\inf_{X\setminus W} f:W\supseteq Z, W \textup{ is pluripolar}\right\}.
    \end{split}
    \]
    \index{$\qsup_X f$}\index{$\qinf_X f$}
\end{definition}
For two functions $f$ and $g$ equal quasi-everywhere, the quasi-suprema and quasi-infima of them are equal, as is clear from the definition.

For a quasi-plurisubharmonic function $\varphi$, we have
\[
\qsup_X \varphi=\sup_X \varphi,
\]
as a consequence of \cref{cor:diffsupinfindeppluripolar}.

\section{Lelong numbers and multiplier ideal sheaves}\label{sec:Lelongmis}
In this section, we briefly recall the notions of Lelong numbers and multiplier ideal sheaves. Our presentation is by no means intended to be complete. Readers are encouraged to read the textbooks \cite[Section~2.3]{GZ17}, \cite{Dem12} and \cite[Chapter~12]{FPT}.

Let $X$ be a complex manifold.

\begin{definition}\label{def:pshfunctions:L1002}
    Let $\varphi\in \PSH(X)$ and $x\in X$. The \emph{Lelong number}\index{Lelong number} $\nu(\varphi,x)$ of $\varphi$ at $x$ is defined as follows:
    If $x$ is a point on a $0$-dimensional connected component of $X$, we set
    \[
        \nu(\varphi,x)=0.
    \]
    Otherwise, take an open neighborhood $U$ of $x$ in $X$ and a biholomorphic map $F\colon U\rightarrow \Omega$, where $\Omega$ is a domain in $\mathbb{C}^n$. Then we define
    \begin{equation}\label{eq:nuvarphix}
    \nu(\varphi,x)\coloneqq \sup\left\{ \gamma\in \mathbb{R}_{\geq 0}: \varphi(F^{-1}(y))\leq \gamma \log|y-F(x)|^2+\mathcal{O}(1)\textup{ as }y\to F(x) \right\}.
    \end{equation}
    \index{$\nu(\varphi,x)$}
\end{definition}
Observe that $\nu(\varphi,x)$ does not depend on the choices of $U$ and $F$, see \cite[Page~167, Corollary~7.2]{DemBook} and \cite[Proposition~12.1.1 and Definition~12.1.2]{FPT}.
Furthermore, it follows from \cref{prop:Lelongreform} below that the supremum in \eqref{eq:nuvarphix} is a maximum.

\begin{remark}\label{rmk:Lelongconv}
Our definition of the Lelong number is not standard. It differs from the standard definition by a factor of $2$. As a mnemonic, just remember
\[
\nu\left(\log |z|^2,0\right)=1 \quad (\text{instead of } 2).
\]
Our convention of the Lelong numbers together with the convention of the multiplier ideal sheaves below make sure that \eqref{eq:LelongMIS} has no extra factors.
\end{remark}

\begin{proposition}\label{prop:Lelongreform}
    Let $\varphi\in \PSH(B_n(0,1))$. Then
    \begin{equation}\label{eq:Lelongnewdef}
    \nu(\varphi,0)=\lim_{r \to 0+} \frac{\sup_{B_n(0,r)}\varphi}{\log r^2}\in [0,\infty).
    \end{equation}
\end{proposition}
\begin{proof}
    Note that \eqref{eq:Lelongnewdef} is trivial when $n=0$, so we may assume that $n>0$ in the proof.

    It follows from \cref{prop:ballpshconvex} that the limit in \eqref{eq:Lelongnewdef} exists and is finite. We shall denote the limit by $\nu'(\varphi,0)$ for the time being.

    We first observe that by \cref{prop:ballpshconvex}, for any fixed $r_1\in (0,1)$ there is a constant $C_{r_1}$ such that
    \begin{equation}\label{eq:varphixlocalupperbd}
    \varphi(x)\leq  \nu'(\varphi,0)\log |x|^2+C_{r_1}
    \end{equation}
    when $x\in B_n(0,r_1)\setminus \{0\}$. In particular, $\nu(\varphi,0)\geq \nu'(\varphi,0)$.

    In order to argue the reverse inequality, we may assume that $\nu(\varphi,0)>0$.

    Next observe that by \eqref{eq:nuvarphix}, for each small enough $\epsilon>0$, we can find $r_0\in (0,1)$ and $C>0$ so that for all $x\in B_n(0,r_0)\setminus \{0\}$, we have
    \[
    \varphi(x)\leq (\nu(\varphi,0)-\epsilon)\log |x|^2+C.
    \]
    It follows that $\nu'(\varphi,0)\geq \nu(\varphi,0)-\epsilon$. Letting $\epsilon\to 0+$, we conclude that $\nu'(\varphi,0)\geq \nu(\varphi,0)$, which combined with the previous inequality gives \eqref{eq:Lelongnewdef}.
\end{proof}

We recall Siu's semicontinuity theorem.
\begin{theorem}[Siu]\label{thm:Siusemi}\index{theorem!Siu's semicontinuity theorem}
    Let $\varphi\in \PSH(X)$. Then the map $X\ni x\mapsto \nu(\varphi,x)$ is upper semicontinuous with respect to the Zariski topology.
\end{theorem}
For an elegant proof we refer to \cite[Theorem~2.10]{Dem12} and \cite[Theorem~12.3.4]{FPT}.

\begin{proposition}\label{prop:Lelongmax}
    Let $\varphi,\psi\in \PSH(X)$, $\lambda\in \mathbb{R}_{\geq 0}$ and $x\in X$. Then
    \[
        \begin{aligned}
        \nu(\varphi\lor \psi,x)=&\min\{\nu(\varphi,x),\nu(\psi,x)\},\\
        \nu(\varphi+\psi,x)=&\nu(\varphi,x)+\nu(\psi,x),\\
        \nu(\lambda\varphi,x)=&\lambda\nu(\varphi,x).
        \end{aligned}
    \]
\end{proposition}
\begin{proof}
    All properties are local, so we may assume that $X=B_n(0,1)$ for some $n\in \mathbb{N}$. All properties follow directly from \cref{prop:Lelongreform}.
\end{proof}

\begin{corollary}\label{cor:supsLelong}
    Let $(\varphi_i)_{i\in I}$ be a non-empty family in $\PSH(X)$ locally uniformly bounded from above and $x\in X$. Then
    \[
    \nu\left(\uscsup[i\in I]\varphi_i,x\right)=\inf_{i\in I}\nu(\varphi_i,x).
    \]
\end{corollary}
\begin{proof}
    We may assume that $X$ is connected. Write $\varphi=\uscsup[i\in I] \varphi_i$. Then $\varphi\in \PSH(X)$ by \cref{prop:pshfunction_closedseq}.

    We observe that the $\leq$ inequality is trivial. It remains to argue the reverse inequality.

    It follows from \cref{prop:Choquet} that we may assume that $I$ is countable. When $I$ is finite, this is already proved in \cref{prop:Lelongmax}. Otherwise, we may further assume that $I=\mathbb{Z}_{>0}$. Thanks to \cref{prop:Lelongmax}, we may further assume that $(\varphi_i)_{i\in \mathbb{Z}_{>0}}$ is an increasing sequence.
    Furthermore, since the problem is local, we may assume that $x=0$, that $X=B_n(0,1)$ for some $n\in \mathbb{N}$, and that $(\varphi_i)_i$ is uniformly bounded from above on some ball $B_n(0,r_1)$ with $0<r_1<1$. Since the sequence $(\varphi_i)_i$ is increasing, the sequence $\nu(\varphi_i,0)$ is decreasing. Thus the constants in \eqref{eq:varphixlocalupperbd} can be chosen independently of $i$, and we have
    \[
    \varphi_i(x)\leq \nu(\varphi_i,0)\log |x|^2+C
    \]
    for all $x\in B_n(0,r_1)\setminus \{0\}$ and all $i\geq 1$, where $C$ is a constant independent of $i$. In particular, thanks to \cref{prop:supsupstardiff}, for almost all $x\in B_n(0,r_1)$, we have
    \[
    \varphi(x)\leq \lim_{i\to\infty}\nu(\varphi_i,0)\log |x|^2+C.
    \]
    Thanks to \cref{prop:pshfuncdetdense}, the same holds for all $x\neq 0$ near $0$, and hence
    \[
    \nu(\varphi,0)\geq \lim_{i\to \infty} \nu(\varphi_i,0).
    \]
    This completes the proof.
\end{proof}

\begin{definition}\label{def:pshfunctions:L1101}
    Let $F\subseteq X$ be a non-empty irreducible analytic subset. Then we define the \emph{generic Lelong number}\index{generic Lelong number} of $\varphi$ along $F$ as
    \[
    \nu(\varphi,F)\coloneqq \min_{x\in F}\nu(\varphi,x).
    \]
    Note that the minimum is attained by Siu's semicontinuity theorem \cref{thm:Siusemi}\index{$\nu(\varphi,F)$}.
\end{definition}

\begin{definition}\label{def:pshfunctions:L1112}
    Let $\varphi\in \PSH(X)$.
    Let $E$ be a prime divisor over $X$ (see \cref{def:primdiv}). Take a proper bimeromorphic morphism $\pi\colon Y\rightarrow X$ from a complex manifold $Y$ such that $E$ is a prime divisor on $Y$, then we define the \emph{generic Lelong number} of $\varphi$ along $E$ as
    \[
        \nu(\varphi,E)\coloneqq \nu(\pi^*\varphi,E).
    \]
    \index{$\nu(\varphi,E)$}
\end{definition}
It follows from \cref{thm:Siusemi} that $\nu(\varphi,E)$ does not depend on the choice of $\pi$.

\begin{definition}\label{def:pshfunctions:L1122}
Let $\varphi\in \PSH(X)$, the \emph{multiplier ideal sheaf} \index{multiplier ideal sheaf} $\mathcal{I}(\varphi)$ \index{$\mathcal{I}(\varphi)$} of $\varphi$ is by definition the ideal sheaf given by
\[
\Gamma\left(U,\mathcal{I}(\varphi)\right)=\left\{ f\in \mathcal{O}_X(U) : |f|^2\exp(-\varphi)\in L^1_{\loc}(U) \right\}
\]
for any open subset $U\subseteq X$.
\end{definition}

\begin{remark}\label{rmk:misconv}
This definition is different from a few references, where instead of $\exp(-\varphi)$, they use $\exp(-2\varphi)$. The conventions adopted in the current book are the most convenient ones as far as I know. They simplify a number of formulae. As a mnemonic, for any real $\lambda>0$, we have
\[
\mathcal{I}(\lambda \log |z|^2)=\mathcal{O}_{\mathbb{C}}\left(-\lfloor \lambda \rfloor \{0\}\right)\quad \Bigl(\textup{instead of }\mathcal{O}_{\mathbb{C}}\left(-\lfloor 2\lambda \rfloor \{0\}\right)\Bigr),
\]
where $z$ is a variable in $\mathbb{C}$ and $\{0\}$ is the divisor defined by $0\in \mathbb{C}$.
\end{remark}

\begin{proposition}[Nadel]\index{theorem!Nadel coherence theorem}
    Let $\varphi\in \PSH(X)$. Then $\mathcal{I}(\varphi)$ is coherent.
\end{proposition}
See \cite[Proposition~5.7]{Dem12} and \cite[Theorem~12.2.7]{FPT}.

\begin{theorem}\label{thm:multipsubadd}
    Let $\varphi,\psi\in \PSH(X)$. Then
    \[
        \mathcal{I}(\varphi+\psi)\subseteq \mathcal{I}(\varphi)\cdot \mathcal{I}(\psi).
    \]
\end{theorem}
See \cite{DEL00}, \cite[Theorem~14.2]{Dem12} and \cite[Theorem~12.4.1]{FPT}.

The two invariants are related by the following simple result:
\begin{proposition}\label{prop:Lelongnumfrommis}
    Let $\varphi\in \PSH(X)$ and $E$ be a prime divisor over $X$. Then
    \begin{equation}\label{eq:LelongMIS}
        \nu(\varphi,E)=\lim_{k\to\infty}\frac{1}{k}\ord_E\mathcal{I}(k\varphi).
    \end{equation}
\end{proposition}
See \cite[Proposition~2.14]{DX21} and \cite[Corollary~12.6.8]{FPT}.

Recall that this particular form of the formula is compatible with our conventions of $\nu$ and $\mathcal{I}$. As a consistency check, consider $\varphi=\log |z|^2$ with $z\in \mathbb{C}$ and let $E$ be the divisor $\{0\}$. Then both sides of \eqref{eq:LelongMIS} are equal to $1$. See \cref{rmk:Lelongconv} and \cref{rmk:misconv}.

Also observe the following simple lemma:
\begin{lemma}\label{lma:blowupLelong}
    Let $x \in X$ lie on a positive-dimensional component and let $\varphi\in \PSH(X)$. Let $\pi\colon Y\rightarrow X$ be the blow-up of $X$ at $x$ and $E$ the divisor defined by the pullback of the ideal of $x$.
    Then
    \[
        \nu(\varphi,x)=\nu(\varphi,E),
    \]
\end{lemma}
See \cite[Corollaire~1.1.8]{Bou02} and \cite[Proposition~12.3.6]{FPT}.

Conversely, the information of the generic Lelong numbers determines the multiplier ideal sheaves:
\begin{theorem}\label{thm:valuativemulti}
    Let $\varphi\in \PSH(X)$.
    Let $x\in X$ and $f\in \mathcal{O}_{X,x}$, $f\neq 0$. Then the following are equivalent:
    \begin{enumerate}
        \item\label{itm:thm:valuativemulti:1} $f\in \mathcal{I}(\varphi)_x$;
        \item\label{itm:thm:valuativemulti:2} there exists $\epsilon>0$ such that for any prime divisor $E$ over $X$ such that $x$ is contained in the center of $E$ on $X$, we have
            \begin{equation}\label{eq:ordfLelonglogdis}
                \ord_E(f)\geq (1+\epsilon) \nu(\varphi,E)-A_X(E).
            \end{equation}
    \end{enumerate}
    If $\varphi$ has analytic singularities and $\pi\colon Y\rightarrow X$ is a log resolution of $\varphi$ (see \cref{def:logres}), let $\{E_i\}$ be the components of its resolved polar divisor that meet $\pi^{-1}(x)$. In this case \ref{itm:thm:valuativemulti:2} can be replaced by the finite set of inequalities
    \begin{equation}\label{eq:ordfLelonglogdis2}
    \ord_{E_i}(f)> \nu(\varphi,E_i)-A_X(E_i)\quad \forall i.
    \end{equation}
\end{theorem}
Here $\ord_E(f)$ denotes the minimum of the local orders of vanishing of $\pi^*f$ along $E$ at the points of $E\cap\pi^{-1}(x)$, on a model $\pi\colon Y\to X$ on which $E$ is realized. At each such point, the local order is the largest power of a local equation of $E$ dividing the pullback germ. It depends only on the germ of $f$ at $x$. On a common resolution, the strict transform of $E$ maps onto $E$ and preserves these local orders, since they are locally constant along the divisor and agree where the models are identified. Thus the minimum is independent of the model.

Here $A_X(\bullet)$ denotes the log discrepancy of prime divisors. In \ref{itm:thm:valuativemulti:2} only divisors with $\nu(\varphi,E)>0$ impose a condition. In the finite test, the components $E_i$ include strict transforms of polar divisors on $X$. Near a point of the fiber where $E_i=\{z_i=0\}$, the exponent of $|z_i|^2$ after change of variables is $m_i+A_X(E_i)-1-\nu(\varphi,E_i)$, where $m_i$ is the local order of $\pi^*f$ along $E_i$. Integrability near the whole fiber requires these exponents to be strictly greater than $-1$ on every local piece, which gives \eqref{eq:ordfLelonglogdis2} after taking the minimum for each $i$.
We refer to \cite[Corollary~10.18, Proposition~10.12]{Bou17} and \cite[Definition~12.6.3 and Theorem~12.6.9]{FPT} for the proof and the precise definition of $A_X$.

The notion of analytic singularities is recalled in \cref{sec:anasing}.

Recall the strong openness theorem:
\begin{theorem}[Guan--Zhou]\label{thm:strongopen}\index{theorem!Guan--Zhou's strong openness theorem}
    Let $\varphi,\psi_j\in \PSH(X)$ ($j\in \mathbb{Z}_{>0}$) such that $\psi_j$ is an increasing sequence converging to $\varphi$ almost everywhere. Then for any $x\in X$, the germs satisfy
    \[
        \mathcal{I}\left(\psi_j\right)_x=\mathcal{I}(\varphi)_x
    \]
    when $j$ is large enough.
\end{theorem}
See \cite{GZ15, Hiep14} and \cite[Corollary~12.5.4]{FPT} for the proof.

\begin{proposition}\label{prop:pull-backmis}
    Let $\pi\colon Y\rightarrow X$ be a smooth morphism between complex manifolds. Assume that $\varphi\in \PSH(X)$. Then
    \[
    \mathcal{I}(\pi^*\varphi)=\pi^*\mathcal{I}(\varphi).
    \]
\end{proposition}
\begin{proof}
    It follows from \cite[Théorème~3.10]{SHC6} that locally $\pi$ can be written as the composition of an étale morphism and a projection. It suffices to handle the two cases separately. See also the coordinate argument in \cite[Corollary~12.4.2]{FPT}.

    Recall that in the complex analytic setting, an étale morphism is locally biholomorphic, so there is nothing to prove in this case.

    Next, assume that $Y=X\times U$, where $U\subseteq \mathbb{C}^n$ is a domain and $\pi$ is the natural projection. It follows from Fubini's theorem that
    \[
    \mathcal{I}(\pi^*\varphi)\supseteq \pi^*\mathcal{I}(\varphi).
    \]
    The reverse inequality is proved in \cite[Proposition~14.3]{Dem12} and \cite[Corollary~12.4.2]{FPT}\footnote{In \cite[Proposition~14.3]{Dem12}, Demailly used the highly non-standard notation $f^*\mathcal{I}(\varphi)$ to denote the image of $f^*\mathcal{I}(\varphi)\rightarrow \mathcal{O}_X$, even when $f$ is not flat. See also the convention in \cite[Section~B.12]{FPT}.}.
\end{proof}

\begin{definition}\label{def:restidealsheaf}
    Given a coherent ideal sheaf $\mathcal{I}$ on $X$ and a closed complex subspace $Y\subseteq X$, the \emph{restriction} \index{restriction ideal} $\Rest_Y\mathcal{I}$ \index{$\Rest_Y\mathcal{I}$} is the coherent ideal sheaf given by
\begin{equation}\label{eq:RestI}
    \Rest_Y\mathcal{I}\coloneqq \mathcal{I}/(\mathcal{I}\cap \mathcal{I}_Y),
\end{equation}
where $\mathcal{I}_Y$ is the ideal sheaf defining $Y$.
\end{definition}
In the literature, it is common to denote this sheaf by the misleading notation $\mathcal{I}|_Y$.

There is a natural morphism
\begin{equation}\label{eq:pullbacktoinverimage}
    i_Y^*\mathcal{I}=  \mathcal{I}/(\mathcal{I}\cdot \mathcal{I}_Y)\rightarrow \Rest_Y \mathcal{I},
\end{equation}
where $i_Y\colon Y\rightarrow X$ is the inclusion.

\begin{theorem}[Ohsawa--Takegoshi]\label{thm:OT}\index{theorem!Ohsawa--Takegoshi extension theorem}
    Let $Y$ be a connected submanifold of $X$ and $\varphi\in \PSH(X)$. Assume that $\varphi|_Y\not\equiv -\infty$. Then
    \[
    \mathcal{I}\left(\varphi|_Y\right)\subseteq \Rest_Y\mathcal{I}(\varphi).
    \]
\end{theorem}
See \cite[Theorem~14.1]{Dem12} and \cite[Theorems~11.4.4 and~12.4.1]{FPT}.

\section{Quasi-plurisubharmonic functions}
\label{sec:qpsh}
In practice, it is important to consider a variant of plurisubharmonic functions. We will fix a complex manifold $X$.

\begin{definition}\label{def:qpsh} Let $\theta$ be a closed real smooth $(1,1)$-form on $X$.

    A \emph{$\theta$-plurisubharmonic function}\index{plurisubharmonic function!theta-plurisubharmonic function@$\theta$-plurisubharmonic function} or \emph{$\theta$-psh function}\index{plurisubharmonic function!theta-psh function@$\theta$-psh function} for short on $X$ is a function $\varphi\colon X\rightarrow [-\infty,\infty)$ such that for each $x\in X$ and each open neighborhood $U$ of $x$ in $X$ satisfying the condition that $\theta=\ddc g$ for some smooth function $g$ on $U$, we have $g+\varphi|_U\in \PSH(U)$. The set of $\theta$-psh functions on $X$ is denoted by $\PSH(X,\theta)$. \index{$\PSH(X,\theta)$}

    A \emph{quasi-plurisubharmonic function} \index{plurisubharmonic function!quasi-plurisubharmonic function} on $X$ is a function $\varphi\colon X\rightarrow [-\infty,\infty)$ such that there exists a smooth closed real $(1,1)$-form $\theta'$ on $X$ such that $\varphi\in \PSH(X,\theta')$. The set of quasi-plurisubharmonic functions on $X$ is denoted by $\QPSH(X)$.\index{$\QPSH(X)$}
\end{definition}

There is a natural preorder on $\QPSH(X)$ defined as follows:
\begin{definition}\label{def:parorder}Assume that $X$ is compact.
    Given $\varphi,\psi\in \QPSH(X)$, we say that $\varphi$ is \emph{more singular} than $\psi$ and write $\varphi\preceq \psi$\index{$\varphi\preceq \psi$}\footnote{Some people write $\psi\preceq \varphi$.} if there is $C\in \mathbb{R}$ such that $\varphi\leq \psi+C$. We also say $\psi$ is less singular than $\varphi$ and write $\psi\succeq \varphi$.

    In case $\varphi\preceq \psi$ and $\psi\preceq \varphi$, we say $\varphi$ and $\psi$ have the same \emph{singularity type}. We write $\varphi\sim \psi$ \index{$\varphi\sim \psi$} in this case.
\end{definition}
It is sometimes convenient to include the constant function $-\infty$ as well. Then we always write $-\infty\preceq \varphi$ for any $\varphi\in \QPSH(X)$. Formally, we write $-\infty\sim -\infty$ as well.

When $X$ is not compact, one can still define similar notions, but the generalization is not unique, and we shall not consider them in this book.
\begin{remark}
For a fixed background $\theta$, we will use the local results for plurisubharmonic functions after adding a smooth local potential of $\theta$.
\end{remark}

\begin{proposition}\label{prop:L1compa}Assume that $X$ is connected and compact.
    Let $\theta$ be a closed real smooth $(1,1)$-form on $X$. Then for any $a,b\in \mathbb{R}$, $a\leq b$, the set
    \[
        \left\{ \varphi\in \PSH(X,\theta): \sup_X \varphi \in [a,b]\right\}
    \]
    is compact with respect to the $L^1$-topology. Moreover, $\varphi\mapsto \sup_X \varphi$ is $L^1$-continuous for $\varphi\in \PSH(X,\theta)$.
\end{proposition}
This is an immediate consequence of \cite[Proposition~8.5, Exercise~1.20]{GZ17}\footnote{In the proof of \cite[Proposition~8.5]{GZ17}, there is a misleading sentence: \emph{This is a dcreasing sequence of $\omega$-psh functions}, which should be removed.}. See also \cite[Proposition~14.2.1]{FPT}.
\begin{remark}\label{rmk:L1cpt}
    Suppose that $X$ is connected and compact and $K\subseteq X$ is a closed non-pluripolar subset. Then
    \[
        \left\{ \varphi\in \PSH(X,\theta): \sup_K \varphi \in [a,b]\right\}
    \]
    is relatively compact with respect to the $L^1$-topology.

    In fact, $\sup_X\varphi-\sup_K\varphi$ is uniformly bounded on $\PSH(X,\theta)$. Otherwise we could find $u_j\in\PSH(X,\theta)$ with $\sup_X u_j=0$ and $\sup_K u_j\leq-2^j$. By \cref{prop:L1compa} the $L^1$-norms of the $u_j$ are uniformly bounded, so $\sum_{j\geq 1}2^{-j}u_j$ is locally integrable. On a relatively compact coordinate ball, add to each $u_j$ the same smooth local potential of $\theta$, normalized to be nonpositive. The weighted partial sums are then decreasing and psh. Their limit is psh by \cref{prop:pshfunction_closedseq} and is the original sum plus the local potential, since the weights sum to one. Thus the original sum is $\theta$-psh and equals $-\infty$ on $K$, contradicting non-pluripolarity.

    The bound shows that $\sup_K\varphi\leq C_0$ implies a uniform upper bound on $\sup_X\varphi$. If also $\sup_K\varphi\geq a$, then $\sup_X\varphi\geq a$, so relative compactness follows from \cref{prop:L1compa}.
\end{remark}

\begin{corollary}\label{cor:decreasingnetL1conv}
    Assume that $X$ is connected and compact. Let $\theta$ be a closed real smooth $(1,1)$-form on $X$, and $(\varphi_i)_{i\in I}$ be a decreasing net in $\PSH(X,\theta)$ such that $\sup_X \varphi_i$ is uniformly bounded from below. Then $\varphi=\inf_{i\in I} \varphi_i$ is a $\theta$-psh function and
    \begin{equation}\label{eq:decreasingnetL1conv}
    \adjustlimits\lim_{i\in I}\sup_X \varphi_i=\sup_X \varphi,\quad \varphi_i\xrightarrow{L^1}\varphi.
    \end{equation}
\end{corollary}
\begin{proof}
    Choose $a\in\mathbb R$ with $\sup_X\varphi_i\geq a$ for every $i$. The nonempty closed sets $\{\varphi_i\geq a\}$ have the finite intersection property by directedness, so compactness gives a point where $\varphi\geq a$. By \cref{prop:pshfunction_closedseq}, on each connected coordinate ball the limit plus a local potential of $\theta$ is either psh or identically $-\infty$. These alternatives agree on overlapping balls, so connectedness and the point where $\varphi\geq a$ give $\varphi\in\PSH(X,\theta)$. The same proposition gives local $L^1$-convergence, hence $L^1$-convergence on the compact manifold $X$. The continuity of the supremum follows from \cref{prop:L1compa}.
\end{proof}

\begin{proposition}\label{prop:Lelongnumberupperbound}
Assume that $X$ is compact.
    Let $\theta$ be a closed real smooth $(1,1)$-form on $X$ and $E$ be a prime divisor over $X$.
    Then
    \[
        \sup \Bigl\{ \nu(\varphi,E): \varphi\in \PSH(X,\theta) \Bigr\}<\infty.
    \]
\end{proposition}
\begin{proof}
    We may restrict to the connected component of $X$ containing the center of $E$, since the Lelong number along $E$ only depends on the restriction to this component.
    It follows from the proof of \cite[Exercise~2.7]{GZ17} that $\nu(\bullet,E)$ is upper semicontinuous with respect to the $L^1$-topology on $\PSH(X,\theta)$, as also shown in \cite[Proposition~14.2.8]{FPT}. Since adding a constant does not change Lelong numbers, we may restrict to the compact set $\{\varphi\in \PSH(X,\theta):\sup_X\varphi=0\}$ from \cref{prop:L1compa}. The desired upper bound follows.
\end{proof}

\begin{proposition}\label{prop:PSHpullbij}
    Let $\pi\colon Y\rightarrow X$ be a proper bimeromorphic morphism from a compact Kähler manifold $Y$. Let $\theta$ be a closed real smooth $(1,1)$-form on $X$. Then the pull-back gives a bijection
    \[
    \pi^*\colon \PSH(X,\theta)\xlongrightarrow{\sim}\PSH(Y,\pi^*\theta).
    \]
\end{proposition}
This follows from a more general result \cref{thm:Zariskimain}.

Since we shall constantly need the balayage argument of \cref{prop:BT_balayage}, we explain explicitly the global version of it.

\begin{proposition}\label{prop:bal_general}
    Fix a closed smooth real $(1,1)$-form $\theta$ on $X$.
    Let $\varphi\in \PSH(X,\theta)$, and $B\Subset B' \subseteq X$ be two holomorphic balls of complex dimension $n\geq 1$. Assume that $\varphi$ is locally bounded on $B'$. Then there is a unique $\psi \in \PSH(X,\theta)$ so that $\psi$ is bounded on an open neighborhood of $\overline{B}$, and
    \[
    \begin{cases}
        \left( \theta|_B+\ddc \psi|_B \right)^n = 0,\\
        \psi = \varphi \quad \textup{on } X\setminus B.
    \end{cases}
    \]
    Moreover, $\psi\geq \varphi$.
\end{proposition}
We shall call $\psi$ the \emph{balayage}\index{balayage} of $\varphi$ on $B$.
\begin{proof}
    Choose a smooth potential $g$ of $\theta$ on $B'$ by the Dolbeault--Grothendieck lemma. In the coordinates of $B'$, the domain $B$ is smoothly bounded and biholomorphic to a Euclidean ball, hence strongly pseudoconvex by \cite{Bell81}. Apply \cref{prop:BT_balayage} to $\varphi+g$ on $B'$ with $D=B$. Subtracting $g$ gives the desired function on $B'$, which agrees with $\varphi$ on $B'\setminus B$ and hence glues to $\varphi$ outside $\overline B$. Local boundedness gives boundedness near $\overline B$, and uniqueness follows from \cref{prop:BT_balayage}. See also \cite[Corollary~14.2.9]{FPT} for balayage on coordinate balls.
\end{proof}

\section{Analytic singularities}\label{sec:anasing}

The simplest type of plurisubharmonic singularities is given by the so-called \emph{analytic singularities}.
The notion is fairly subtle and there are several mutually \emph{incompatible} definitions in the literature, as we shall explain below.

Let $X$ be a complex manifold.
\begin{definition}\label{def:neatanasing}
    We say $\varphi\in \QPSH(X)$ has \emph{analytic singularities} \index{singularities!analytic singularities} if for each $x\in X$, we can find a connected open neighborhood $U$ of $x$ such that $\varphi|_U$ has the form:
    \begin{equation}\label{eq:anasinglocal}
    c\log (|f_1|^2+\cdots+|f_N|^2)+R,
    \end{equation}
    where $f_1,\ldots,f_N$ are holomorphic functions on $U$, $c\in \mathbb{Q}_{>0}$ and $R$ is a bounded function on $U$.\footnote{Note that since $\varphi|_U\not\equiv -\infty$ by assumption, it is implicit in the definition that $\{f_1=\cdots=f_N=0\}\subseteq U$ is a proper subset.}

    When $R$ can be taken to be smooth\footnote{The decomposition \eqref{eq:anasinglocal} is highly non-unique. Here we mean for any $x$, there is an open neighborhood $U$ and a decomposition of the form \eqref{eq:anasinglocal} with $R$ smooth. In the non-trivial cases, $R$ cannot be smooth for all decompositions \eqref{eq:anasinglocal}.}, we say $\varphi$ has \emph{neat analytic singularities}\index{singularities!analytic singularities!neat analytic singularities}.

    Suppose that there is a coherent ideal $\mathcal{I}\subseteq \mathcal{O}_X$ on $X$ such that we can choose $U$ so that the $f_1,\ldots,f_N$ can be chosen as the generators of $\Gamma(U,\mathcal{I})$ and $c$ is independent of the choice of $U$. Then we say $\varphi$ has analytic singularities of \emph{type} $(c,\mathcal{I})$.
\end{definition}
If $X$ is compact, each potential with analytic singularities has a type. The type is not uniquely determined.
We refer to \cite{Bou02} and \cite{Bou02b} for the details. See also \cite[Section~13.1]{FPT}.

Some authors allow $c\in \mathbb{R}_{>0}$ in \eqref{eq:anasinglocal}. We keep $c\in \mathbb{Q}_{>0}$ because the following closure property, which is essential in constructing Demailly approximations, would otherwise fail.
\begin{proposition}\label{prop:analysingclosed}
    Let $\varphi,\psi\in \QPSH(X)$ be potentials with analytic singularities. Then so are $\lambda\varphi$ for $\lambda\in \mathbb{Q}_{>0}$ and $\varphi+\psi$. If $\varphi$ and $\psi$ are $\theta$-psh for the same smooth closed real $(1,1)$-form $\theta$, then $\varphi\lor\psi$ is $\theta$-psh and has analytic singularities.
\end{proposition}
\begin{proof}
    Scaling and addition preserve quasi-plurisubharmonicity by scaling and adding the respective backgrounds. For the maximum, adding a common local potential of $\theta$ reduces $\theta$-plurisubharmonicity to the local maximum property.

    Multiplying the coefficient in \eqref{eq:anasinglocal} gives the assertion for $\lambda\varphi$. For the other two operations, use a common positive rational coefficient in the local presentations, by replacing the generating systems with generators of suitable powers of their ideals. For addition, take all pairwise products of the two systems. For the maximum, take their union, using that the logarithm of the sum of two nonnegative quantities differs from the maximum of their logarithms by at most $\log 2$ wherever the sum is nonzero.
\end{proof}

\begin{definition}\label{def:pshfunctions:L1393}
    Let $D$ be an effective $\mathbb{Q}$-divisor\footnote{Divisors and $\mathbb{Q}$-divisors are implicitly assumed to have locally finite coefficients as usual.} on $X$.
    We say $\varphi\in \QPSH(X)$ has \emph{log singularities}\index{singularities!log singularities} (along $D$) on $X$ if for each $x\in X$, there is an open neighborhood $U$ of $x$ such that
    \begin{enumerate}
        \item\label{itm:chapter-pshfunctions:L1415:1} $D|_U$ has finitely many irreducible components and can be written as
        \[
            D|_U=\sum_{i=1}^N a_i D_i
        \]
        with $D_i$ prime divisors on $U$, $a_i\in \mathbb{Q}_{>0}$ and $s_i$ holomorphic functions on $U$ defining $D_i$, and
        \item\label{itm:chapter-pshfunctions:L1415:2} we have
        \begin{equation}\label{eq:logsingreminder}
            \varphi|_U=\sum_{i=1}^N a_i\log |s_i|^2+R,
        \end{equation}
        where $R$ is a bounded function on $U$.
    \end{enumerate}
\end{definition}
By \cref{prop:analysingclosed}, $\varphi$ has analytic singularities.

\begin{lemma}\label{lma:logsingrem}
    Suppose that $X$ is a compact Kähler manifold, $\theta$ is a closed smooth real $(1,1)$-form on $X$, and $\varphi\in \PSH(X,\theta)$. Suppose that $\varphi$ has log singularities along an effective $\mathbb{Q}$-divisor $D$ on $X$. Then the cohomology class $\{\theta\}-\{D\}$ is nef.

    Moreover, if in addition $\theta_{\varphi}$ is a Kähler current\footnote{That is, there is a K\"ahler form $\omega$ on $X$ such that $\theta_{\varphi}\geq \omega$ in the sense of currents.}, then the cohomology class $\{\theta\}-\{D\}$ is Kähler.
\end{lemma}
Here and in the sequel, we write $\theta_{\varphi}$ for $\theta+\ddc\varphi$.

\begin{proof}
    The current $\theta_{\varphi}-[D]$ has bounded local potentials. Therefore, the class $\{\theta\}-\{D\}$ is nef.

    If $\theta_{\varphi}$ is a Kähler current, say $\theta_{\varphi}\geq \omega$ for a Kähler form $\omega$, then $\theta_{\varphi}-[D]-\omega$ is positive and has bounded local potentials, and therefore represents a nef class. We conclude that $\{\theta\}-\{D\}$ is Kähler.
\end{proof}

The following proposition follows immediately from the definitions:
\begin{proposition}\label{prop:pshfunctions:L1425}
    Let $\pi\colon Y\rightarrow X$ be a proper bimeromorphic morphism from a complex manifold $Y$.
    Suppose that $\varphi\in \QPSH(X)$ has analytic singularities (resp.\ has log singularities along an effective $\mathbb{Q}$-divisor $D$).
    Then $\pi^*\varphi$ has analytic singularities (resp.\ has log singularities along $\pi^*D$).
\end{proposition}

\begin{definition}\label{def:logres}
    Let $\varphi\in \QPSH(X)$ be a potential with analytic singularities.
    A \emph{log resolution}\index{log resolution} of $\varphi$ is a modification $\pi\colon Y\rightarrow X$ such that $\pi^*\varphi$ has log singularities along a simple normal crossing effective $\mathbb{Q}$-divisor.
\end{definition}

\begin{theorem}\label{thm:resolvelogsing}
    Assume that $X$ is compact.
    Suppose that $\varphi\in \QPSH(X)$ has analytic singularities. Then there is a log resolution of $\varphi$.
\end{theorem}
For a proof, we refer to the arguments on \cite[Page~104]{MM07} and to \cite[Proposition~13.1.10]{FPT}.

A general quasi-plurisubharmonic function can be nicely approximated by those with analytic singularities. We need a few preliminary definitions.
\begin{definition}\label{def:quasiequsing}
    Let $X$ be a compact Kähler manifold and $\theta$ be a closed real smooth $(1,1)$-form on $X$. Consider $\varphi\in \PSH(X,\theta)$. A sequence $(\varphi_j)_{j\in \mathbb{Z}_{>0}}$ in $\QPSH(X)$ is \emph{quasi-equisingular approximation} \index{quasi-equisingular approximation} of $\varphi$ if
    \begin{enumerate}
        \item\label{itm:def:quasiequsing:1} $\varphi_j$ has analytic singularities for each $j$;
        \item\label{itm:def:quasiequsing:2} $\varphi_j$ is decreasing with limit $\varphi$;
        \item\label{itm:def:quasiequsing:3} there is a decreasing sequence $\epsilon_j\geq 0$ with limit $0$ and a Kähler form $\omega$ on $X$ such that $\varphi_j\in \PSH(X,\theta+\epsilon_j\omega)$;
        \item\label{itm:def:quasiequsing:4} for each $\lambda'>\lambda>0$, there is $j>0$ such that
        \begin{equation}\label{eq:quasi_eq_prop}
        \mathcal{I}(\lambda' \varphi_j)\subseteq \mathcal{I}(\lambda \varphi).
        \end{equation}
    \end{enumerate}
    We also say $(\theta_{\varphi_j})_j$ is a quasi-equisingular approximation of $\theta_{\varphi}$.
\end{definition}

\begin{definition}\label{def:analy-sing}
    Let $\mathcal{I}\subseteq\mathcal O_X$ be a coherent ideal sheaf and $c\in\mathbb{Q}_{>0}$. A function $\varphi\in \QPSH(X)$ is said to have \emph{gentle analytic singularities} \index{singularities!analytic singularities!gentle analytic singularities} (of type $(c,\mathcal{I})$) if
\begin{enumerate}
    \item\label{itm:def:analy-sing:1} $\varphi$ has analytic singularities of type $(c,\mathcal{I})$;
    \item\label{itm:def:analy-sing:2} $\mathrm{e}^{{\varphi}/{c}}\colon X\to \mathbb{R}_{\geq 0}$ is a smooth function;
    \item\label{itm:def:analy-sing:3} there is a proper bimeromorphic morphism $\pi\colon \widetilde{X}\to X$ from a Kähler manifold $\widetilde{X}$ and an effective $\mathbb{Z}$-divisor $D$ on $\widetilde{X}$ such that one can write $\pi^*\varphi$ locally as
    \[
        \pi^*\varphi=c\log|g|^2+h,
    \]
    where $g$ is a local equation of the divisor $D$ and $h$ is smooth.
\end{enumerate}
\end{definition}

\begin{theorem}\label{thm:qequi}
    Let $X$ be a compact Kähler manifold and $\theta$ be a closed real smooth $(1,1)$-form on $X$. Then any $\varphi\in \PSH(X,\theta)$ admits a quasi-equisingular approximation $(\varphi_j)_{j\in \mathbb{Z}_{>0}}$.

    Moreover, we can guarantee that for all $j>0$, $\varphi_j$ has gentle analytic singularities of type $\left(2^{-j},\mathcal{I}(2^j\varphi)\right)$.
\end{theorem}
We refer to \cite{Dem92reg, DPS01} and \cite[Theorem~13.4.5]{FPT} for the proof.

Quasi-equisingular approximations are essentially unique in the following sense:
\begin{proposition}\label{prop:compqequi}
    Let $X$ be a compact Kähler manifold and $\theta$ be a closed real smooth $(1,1)$-form on $X$. Consider $\varphi\in \PSH(X,\theta)$. Let $(\varphi_j)_j$ and $(\psi_j)_j$ be two quasi-equisingular approximations of $\varphi$. Then for any $\epsilon\in (0,1)$ and any $j>0$, we can find $k_0>0$ such that for any $k\geq k_0$, we have
    \[
        \psi_k\preceq (1-\epsilon)\varphi_j.
    \]
\end{proposition}
See \cite[Corollary~4.1.7]{Dem15} and \cite[Proposition~13.5.2]{FPT}.

\begin{definition}\label{def:Iinfty}
    Let $\varphi\in \QPSH(X)$ be a potential with analytic singularities. Then we define $\mathcal{I}_{\infty}(\varphi)$ \index{$\mathcal{I}_{\infty}(\varphi)$} as the ideal sheaf consisting of germs $f$ of holomorphic functions such that $\log |f|^2\leq \varphi+C$ for some $C\in \mathbb{R}$.
\end{definition}
By definition, $\mathcal{I}_{\infty}(\varphi)\subseteq \mathcal{I}(\varphi)$.

\begin{lemma}\label{lma:pshfunctions:L1494}
    Assume that $X$ is compact. Let $\varphi\in \QPSH(X)$ be a potential with analytic singularities.
    The sheaf $\mathcal{I}_{\infty}(\varphi)$ is coherent.
\end{lemma}
\begin{proof}
    By \cref{thm:resolvelogsing}, we may find a modification $\pi\colon Y\rightarrow X$ such that $\pi^*\varphi$ has log singularities. Observe that
    \[
        \mathcal{I}_{\infty}(\varphi)=\pi_*\mathcal{I}_{\infty}(\pi^*\varphi),
    \]
    so we may replace $X$ and $\varphi$ by $Y$ and $\pi^*\varphi$ and assume that $\varphi$ has log singularities along an effective $\mathbb{Q}$-divisor $D$. We decompose $D$ into its irreducible components:
    \[
        D=\sum_{i=1}^N a_iD_i.
    \]
    In this case, observe that
    \[
        \mathcal{I}_{\infty}(\varphi)=\mathcal{O}_X\left(-\sum_{i=1}^N\left(\lceil a_i \rceil D_i\right)\right)
    \]
    is clearly coherent.
\end{proof}

The multiplier ideal sheaf $\mathcal{I}$ and the sheaf $\mathcal{I}_{\infty}$ are not very different in the asymptotic sense, as shown by the following lemma:
\begin{lemma}\label{lma:IandIinf}
    Assume that $X$ is compact. Let $\varphi\in \QPSH(X)$ be a potential with analytic singularities. Then for any $\epsilon>0$, we can find $k_0>0$ such that for each $k\geq k_0$, we have
    \begin{equation}\label{eq:IandIinf}
        \mathcal{I}\left(k(1+\epsilon)\varphi\right)\subseteq \mathcal{I}_{\infty}(k\varphi).
    \end{equation}
\end{lemma}
\begin{proof}
    We shall prove a more precise local result. Take $x\in X$, take an open neighborhood $U\subseteq X$ of $x$, on which
    \[
    \varphi=c\log \left( |g_1|^2 +\cdots+|g_N|^2\right)+\mathcal{O}(1),
    \]
    where $c\in \mathbb{Q}_{>0}$ and $g_1,\ldots,g_N$ are holomorphic functions on $U$. Then we claim that for any $\lambda>0$, we have
    \begin{equation}\label{eq:IIinf_temp1}
        \mathcal{I}(\lambda\varphi)_x\subseteq \mathcal{I}_{\infty}\left( \left( \lambda-c^{-1}n \right)_+\varphi \right)_x,
    \end{equation}
    where for any $a\in \mathbb{R}$, $a_+$ means $a\lor 0$. Note that \eqref{eq:IandIinf} is a straightforward consequence of \eqref{eq:IIinf_temp1}. The global result follows easily from the local result by the compactness of $X$.

    Fix a smooth volume form $\mathrm{d}V$ on $X$.
    Take $f\in \mathcal{I}(\lambda\varphi)$. By strong openness \cref{thm:strongopen}, we can find $\epsilon>0$ so that $f\in \mathcal{I}((\lambda+\epsilon/c)\varphi)$. Therefore, we can find an open neighborhood $W\subseteq U$ of $x$ so that
    \[
    \int_W |f|^2\left( |g_1|^2 +\cdots+|g_N|^2\right)^{-c\lambda-\epsilon}\,\mathrm{d}V<\infty.
    \]
    But it follows from the Briançon--Skoda division theorem (\cite{Sko72}, \cite{SB74}, \cite[Theorem~11.17]{Dem12}) that
    \[
    f_x\in \left(g_{1,x},\ldots,g_{N,x} \right)^\alpha,
    \]
    where $\alpha=\left(\lfloor c \lambda \rfloor-n+1\right)_+$. Here $\left(g_{1,x},\ldots,g_{N,x} \right)$ denotes the ideal in $\mathcal{O}_{X,x}$ generated by the germs of $g_1,\ldots,g_N$ at $x$. Note that $\alpha\geq (c\lambda-n)_+$.

    It follows that on a neighborhood of $x$, we have
    \[
    \log |f|^2\leq \alpha \log \left( |g_1|^2 +\cdots+|g_N|^2\right)+\mathcal{O}(1)\leq c^{-1}\alpha \varphi+\mathcal{O}(1).
    \]
    Hence \eqref{eq:IIinf_temp1} follows. For a proof of the asymptotic inclusion using log resolutions, see \cite[Corollary~13.6.3]{FPT}.
\end{proof}

\begin{theorem}\label{thm:CT-thm-refined'}Let $X$ be a connected compact Kähler manifold and $Y\subseteq X$ be a connected closed submanifold.
    Take a Kähler form $\omega$ on $X$ and $\varphi\in \PSH(Y,\omega|_Y)$ such that $\omega|_Y+\ddc\varphi$ is a Kähler current and $\mathrm{e}^{\varphi}$ is a H\"older continuous function on $Y$.
    Then there exists $\tilde{\varphi}\in \PSH(X,\omega)$ satisfying
    \begin{enumerate}
        \item\label{itm:thm:CT-thm-refined':1} $\tilde{\varphi}|_Y=\varphi$;
        \item\label{itm:thm:CT-thm-refined':2} $\omega_{\tilde{\varphi}}$ is a Kähler current.
    \end{enumerate}
        In addition, if $\varphi$ has analytic singularities, then so does $\tilde{\varphi}$.
\end{theorem}
See \cite[Theorem~6.1]{DRWNXZ} and \cite[Theorem~14.7.12]{FPT}.

\section{The space of currents}
\label{sec:currents}
Let $X$ be a connected compact Kähler manifold of dimension $n$ and $\alpha\in \mathrm{H}^{1,1}(X,\mathbb{R})$.

\begin{definition}\label{def:pos_curr}
    Let $Y$ be a complex manifold of dimension $N$ and $m\in \mathbb{N}$. We say an $(m,m)$-current $T$ on $Y$ is \emph{positive}\index{positive currents}\footnote{This notion is sometimes known as \emph{weak positivity}.} if either $m>N$ or for any smooth $(1,0)$-forms $\beta_1,\ldots,\beta_{N-m}$ on $Y$, the measure
    \[
    T\wedge \mathrm{i}\beta_1\wedge \overline{\beta_1}\wedge \cdots \wedge \mathrm{i}\beta_{N-m}\wedge \overline{\beta_{N-m}}
    \]
    is positive.
\end{definition}
The basic properties of positive currents can be found in \cite[Section~III.1]{DemBook} and \cite[Section~7.2]{FPT}. Recall that a positive current is necessarily real.

\begin{definition}\label{def:psefclass}
    We say $\alpha$ is \emph{pseudo-effective}\index{cohomology class!pseudo-effective class}\index{line bundle!pseudo-effective line bundle}\index{class!pseudo-effective class} if there is a closed positive $(1,1)$-current in $\alpha$.

    We say $\alpha$ is \emph{big}\index{cohomology class!big class}\index{line bundle!big line bundle}\index{class!big class} if there is a closed positive $(1,1)$-current $T$ in $\alpha$ dominating a Kähler form. Such currents are called \emph{Kähler currents}\index{Kähler currents}.
\end{definition}
Given classes $\alpha,\beta\in \mathrm{H}^{1,1}(X,\mathbb{R})$, we say $\alpha\leq \beta$ if $\beta-\alpha$ is pseudo-effective.

\begin{definition}\label{def:spaceofcurrents}
    We introduce the following notations:
    \begin{enumerate}
        \item\label{itm:def:spaceofcurrents:1} $\mathcal{Z}_+(X)$ \index{$\mathcal{Z}_+(X)$} denotes the space of closed positive $(1,1)$-currents on $X$;
        \item\label{itm:def:spaceofcurrents:2} given a pseudo-effective $(1,1)$-class $\alpha$ on $X$, we write $\mathcal{Z}_+(X,\alpha)$ \index{$\mathcal{Z}_+(X,\alpha)$} for the set of $T\in \mathcal{Z}_+(X)$ such that $\{T\}=\alpha$.
    \end{enumerate}
\end{definition}
Here $\{T\}$ denotes the cohomology class represented by $T$.

\cref{def:parorder} has a natural analogue for currents.

\begin{definition}\label{def:pshfunctions:L1594}
Given $T,T'\in \mathcal{Z}_+(X)$, we write $T\preceq T'$
and say $T$ is \emph{more singular} than $T'$ if when we write $T=\theta+\ddc\varphi$, $T'=\theta'+\ddc \varphi'$, we have $\varphi\preceq \varphi'$. We write $T\sim T'$
if $T\preceq T'$ and $T'\preceq T$. In this case, we say $T$ and $T'$ have the same \emph{singularity type}\index{singularity type}.
\end{definition}

\begin{remark}\label{rmk:qpshtocurrents}
Observe that
\[
\mathcal{Z}_+(X)/\sim\cong \QPSH(X)/\sim
\]
canonically. The correspondence sends the class of a closed positive current $\theta_{\varphi}=\theta+\ddc \varphi$ to the class of $\varphi$.

We will adopt the following convention: Whenever we have a notion for quasi-plurisubharmonic functions which depends only on the singularity type, we use the same notation and the same definition for closed positive $(1,1)$-currents.
\end{remark}

\begin{example}\label{ex:Lelongcurrent}
    As an important example of \cref{rmk:qpshtocurrents}, given $T=\theta+\ddc\varphi\in \mathcal{Z}_+(X)$ and $x\in X$, we define
\begin{equation}\label{eq:convnuTx}
\nu(T,x)=\nu(\varphi,x).
\end{equation}
Again, as in \cref{rmk:Lelongconv}, this differs from the definitions in some literature by a factor of $2$. But given our normalization
\[
\ddc=\frac{\mathrm{i}}{2\pi}\partial\overline{\partial},
\]
\eqref{eq:convnuTx} seems to be the most natural choice.

The key example to keep in mind is the following:
\[
\nu([0],0)=1,
\]
where $[0]$ is the current of integration at $0\in \mathbb{P}^1$. In fact, as a simple application of the Green's second identity, one can verify that
\[
\frac{\mathrm{i}}{2\pi}\partial\overline{\partial}\log |z|^2=\delta_0,
\]
where the right-hand side is the Dirac delta distribution at $0\in \mathbb{C}$.
\end{example}

\begin{definition}\label{def:polarlocus}
    Given $T\in \mathcal{Z}_+(X)$, represent $T$ as $\theta+\ddc\varphi$ for some closed smooth real $(1,1)$-form $\theta$ on $X$ and $\varphi\in \PSH(X,\theta)$. The \emph{polar locus}\index{polar locus} of $T$ is defined as the set $\{\varphi=-\infty\}$.
\end{definition}
It is clear that the polar locus of $T$ is independent of the choices of $\theta$ and $\varphi$. We also call the set $\{-\infty\}$ the polar locus of $\varphi$.

\begin{definition}\label{def:nKlocus}
    Assume that $\alpha$ is big. The \emph{non-Kähler locus}\index{non-Kähler locus} $\nK(\alpha)$\index{$\nK(\alpha)$} of $\alpha$ is the intersection of the polar loci of all Kähler currents with analytic singularities in $\alpha$.
\end{definition}
\begin{theorem}[Boucksom]\label{thm:pshfunctions:L1641}
Assume that $\alpha$ is big.
    There is a Kähler current $T\in \alpha$ with analytic singularities, such that the polar locus of $T$ is exactly $\nK(\alpha)$. In particular, $\nK(\alpha)$ is a proper Zariski closed subset of $X$.
\end{theorem}
See \cite[Théorème~2.1.20]{Bou02} and \cite[Theorem~14.4.4]{FPT}.

\begin{definition}\label{def:nnlocus}
     Assume that $\alpha$ is big. The \emph{non-nef locus}\index{non-nef locus} $\nn(\alpha)$\index{$\nn(\alpha)$} of $\alpha$ is the following set:
     \[
     \nn(\alpha)\coloneqq \{x\in X :\nu(T_{\min},x)>0\},
     \]
     where $T_{\min}$ is a current with minimal singularities in $\alpha$.
\end{definition}
Note that $\nn(\alpha)\subseteq \nK(\alpha)$. Thanks to \cref{thm:Siusemi}, $\nn(\alpha)$ is a countable union of proper Zariski closed subsets of $X$. The non-Kähler locus and non-nef locus are studied in detail in \cite{Bou02} and \cite[Section~14.4]{FPT}.

\begin{lemma}[Siu's decomposition]\label{lma:Siudec}\index{theorem!Siu's decomposition}
    Let $E$ be a prime divisor on $X$. Then for any closed positive $(1,1)$-current $T$ on $X$, the difference $T-\nu(T,E)[E]$ is a closed positive $(1,1)$-current.
\end{lemma}
Here $[E]$ is the current of integration associated with $E$, while $\{E\}$ denotes the cohomology class of $E$. See \cite[Page~386, Example~1]{GH14} and \cite[Section~7.4]{FPT} for the current of integration. See \cite[Lemma~2.17]{Dem12} and \cite[Theorem~14.3.7]{FPT} for the proof.

As a consequence, for each closed positive $(1,1)$-current $T$ on $X$, we can write
\begin{equation}\label{eq:Siudec}
    T=\Reg T+\sum_{i} c_i [E_i],
\end{equation}
where $\{E_i\}$ is a countable collection of prime divisors on $X$, $c_i>0$.

\begin{definition}\label{def:nd_curr}
    A closed positive $(1,1)$-current $T$ on $X$ is \emph{non-divisorial}\index{non-divisorial current} (resp. \emph{divisorial}\index{divisorial current}) if $T=\Reg T$ (resp. $\Reg T=0$).
\end{definition}

It is helpful to check that our conventions are always consistent: There is no extra factor of $2$ or $1/2$ anywhere. One can verify this using our favorite example as in \cref{ex:Lelongcurrent}.

Next we recall the notion of modified nef classes.

\begin{definition}\label{def:mnclass}
    A class $\alpha\in \mathrm{H}^{1,1}(X,\mathbb{R})$ is \emph{modified nef}\index{modified nef class} if the following condition holds: after fixing a reference Kähler metric $\omega$ on $X$, for any $\epsilon>0$, we can find a closed $(1,1)$-current $T\in \alpha$ such that
    \begin{enumerate}
        \item\label{itm:def:mnclass:1} $T+\epsilon\omega\geq 0$;
        \item\label{itm:def:mnclass:2} $\nu(T+\epsilon\omega,D)=0$ for any prime divisor $D$ on $X$.
    \end{enumerate}
    This definition is independent of the choice of $\omega$.
\end{definition}
These classes are called \emph{nef en codimension $1$} in Boucksom's thesis \cite{Bou02}, where they were introduced for the first time. They are also known as \emph{movable classes} in other references. See also \cite[Definition~14.3.8]{FPT}.
Modified nef classes form a closed convex cone in $\mathrm{H}^{1,1}(X,\mathbb{R})$. Note that a modified nef class is necessarily pseudo-effective. A nef class is obviously modified nef.

Recall the multiplicity of a cohomology class as defined in \cite[Section~2.1.3]{Bou02} and \cite[Definition~14.4.1 and Proposition~14.4.2]{FPT}:
\begin{definition}\label{def:mult_class}
    Let $\alpha\in \mathrm{H}^{1,1}(X,\mathbb{R})$ be a pseudo-effective class and $D$ be a prime divisor on $X$. We define the \emph{Lelong number}\index{Lelong number} $\nu(\alpha,D)$\index{$\nu(\alpha,D)$} as follows:
    \begin{enumerate}
        \item\label{itm:def:mult_class:1} When $\alpha$ is big, define $\nu(\alpha,D)=\nu(T,D)$ for any closed positive $(1,1)$-current $T\in \alpha$ with minimal singularities.
        \item\label{itm:def:mult_class:2} In general, define
        \[
        \nu(\alpha,D)\coloneqq \lim_{\epsilon\to 0+}\nu(\alpha+\epsilon \{\omega\},D).
        \]
    \end{enumerate}
\end{definition}
When $\alpha$ is big, \ref{itm:def:mult_class:2} is compatible with \ref{itm:def:mult_class:1} and the definition is independent of the choice of $\omega$.

By definition, a pseudo-effective class $\alpha$ is modified nef if and only if $\nu(\alpha,D)=0$ for all prime divisors $D$ on $X$.

Let us recall the behavior of several cones under modifications.
\begin{proposition}\label{prop:cones_modif}
    Let $\pi\colon Y\rightarrow X$ be a proper bimeromorphic morphism from a Kähler manifold $Y$.
    \begin{enumerate}
        \item\label{itm:prop:cones_modif:1} For any nef class $\alpha\in \mathrm{H}^{1,1}(X,\mathbb{R})$, $\pi^*\alpha$ is nef.
        \item\label{itm:prop:cones_modif:2} For any modified nef class $\beta\in \mathrm{H}^{1,1}(Y,\mathbb{R})$, $\pi_*\beta$ is modified nef.
    \end{enumerate}
\end{proposition}

\begin{proof}
    Only \ref{itm:prop:cones_modif:2} requires a proof. Fix a Kähler class $\gamma$ on $Y$. Since the modified nef cone is closed, replacing $\beta$ by $\beta+\epsilon\gamma$ and then letting $\epsilon\to 0+$ reduces the proof to the case where $\beta$ is big as well. Let $T$ (resp.\ $S$) be a current with minimal singularities in $\pi_*\beta$ (resp.\ in $\beta$), and let $D$ be a prime divisor on $X$. It suffices to show that
    \[
    \nu(T,D)=0.
    \]
    Since $\beta$ is modified nef, $S$ is non-divisorial. Hence by \cref{lma:ndvpush} below, $\nu(\pi_*S,D)=0$. Since $T$ has minimal singularities in $\pi_*\beta$, we have $T\succeq \pi_*S$, and the desired equality follows.
\end{proof}

Recall that non-divisorial currents are introduced in \cref{def:nd_curr}.
\begin{lemma}\label{lma:ndvpush}
    Let $\pi\colon Y\rightarrow X$ be a proper bimeromorphic morphism from a Kähler manifold $Y$. Let $T$ be a non-divisorial current on $Y$. Then $\pi_*T$ is non-divisorial.
\end{lemma}
Conversely, if $S$ is a non-divisorial current on $X$, $\pi^*S$ could have divisorial part. As a simple example, consider $S$ on $\mathbb{P}^2$, whose local potential near $0\in \mathbb{C}^2_{z,w}$ looks like $\log (|z|^2+|w|^2)$.

\begin{proof}
    Let $D$ be a prime divisor on $X$. It follows from Zariski's main theorem \cref{thm:Zariskimain} that $D$ is not contained in the exceptional locus of $\pi$.
    Let $D'$ be the strict transform of $D$. Thanks to Siu's semicontinuity theorem \cref{thm:Siusemi}, we have
    \[
    \nu(\pi_*T,D)=\nu(T,D')=0.
    \]
    Hence $\pi_*T$ is non-divisorial.
\end{proof}

\section{Plurisubharmonic metrics on line bundles}\label{sec:pshmetriconlb}
A natural source of quasi-plurisubharmonic functions is the metrics on line bundles.

Let $X$ be a connected Kähler manifold and $L$ be a holomorphic line bundle on $X$. Usually, we do not distinguish $L$ from the associated invertible sheaf $\mathcal{O}_X(L)$.

\begin{definition}\label{def:pshfunctions:L1743}
    Let $V$ be a $1$-dimensional complex linear space. A \emph{Hermitian form}\index{Hermitian form} $h$ on $V$ is a map $h\colon V\times V\rightarrow \mathbb{C}$ such that
    \begin{enumerate}
        \item\label{itm:chapter-pshfunctions:L1764:1} $h$ is $\mathbb{C}$-linear in the second variable and conjugate linear in the first, and
        \item\label{itm:chapter-pshfunctions:L1764:2}
        \[
        |v|_h^2\coloneqq h(v,v)\in \mathbb{R}_{>0}
        \]
        for each $v\in V\setminus\{0\}$.
    \end{enumerate}
    We usually identify $h$ with the quadratic form $V\rightarrow \mathbb{R}$ sending $v$ to $|v|_h^2$. We write $|v|_h=\sqrt{|v|_h^2}$ for any $v\in V$.

    The \emph{singular Hermitian form} \index{Hermitian form!singular Hermitian form} on $V$ is the map $V\rightarrow \{0,\infty\}$ sending $0$ to $0$ and other elements to $\infty$.
\end{definition}
\begin{definition}\label{def:1dimlinspacehermten}
    Let $V_1$ and $V_2$ be $1$-dimensional complex linear spaces. Given two maps $h_i\colon V_i\rightarrow [0,\infty]$ ($i=1,2$), each of which is either a Hermitian form or a singular Hermitian form, we define the \emph{tensor product} $h_1\otimes h_2\colon V_1\otimes V_2\rightarrow [0,\infty]$ as follows:
    \index{Hermitian form!tensor product}
    \begin{enumerate}
        \item\label{itm:def:1dimlinspacehermten:1} If either $h_1$ or $h_2$ is singular, we define $h_1\otimes h_2$ as the singular Hermitian form;
        \item\label{itm:def:1dimlinspacehermten:2} otherwise, define $h_1\otimes h_2$ as the usual tensor product: For any $v_1\in V_1$, $v_2\in V_2$, set
        \[
        h_1\otimes h_2(v_1\otimes v_2)=h_1(v_1)h_2(v_2).
        \]
    \end{enumerate}
\end{definition}

\begin{definition}\label{def:Hermmetric}
    A \emph{Hermitian metric}\index{Hermitian metric} $h$ on $L$ is a family of Hermitian forms $(h_x)_{x\in X}$, such that
    \begin{enumerate}
        \item\label{itm:def:Hermmetric:1} for each $x\in X$, $h_x$ is a Hermitian form on $L_x$, and
        \item\label{itm:def:Hermmetric:2} for each local holomorphic section $s$ of $\mathcal{O}_X(L)$, the map $x\mapsto |s(x)|^2_{h_x}$ is smooth.
    \end{enumerate}
    The pair $(L,h)$ is called a \emph{Hermitian line bundle}\index{Hermitian line bundle}.
    We shall write $\ddc h=c_1(L,h)$\footnote{The unusual notation $\ddc h$ is sometimes referred to as the \emph{Göteborg notation} because it is widely used by complex geometers in Göteborg (usually spelled as Gothenburg in English, the second largest (yet very poorly known) city in Sweden). As I identify myself as \emph{Göteborgare}, I do not feel guilty about this notation.} for the first Chern form of $h$\footnote{In the literature, people sometimes define the \emph{curvature form} of $(L,h)$ as $\Theta_{h}=-2\pi\mathrm{i}\ddc h$.}. With the normalization of $\ddc$ fixed in the conventions, $[c_1(L,h)]$ is the image of the topological first Chern class $c_1(L)\in H^2(X,\mathbb{Z})$ in real de Rham cohomology.
    The map $x\mapsto |s(x)|_{h_x}$ will be denoted by $|s|_h$.

    To be more precise, if $U\subseteq X$ is an open subset on which $L$ admits a nowhere vanishing holomorphic section $s$, then we define
    \[
    (\ddc h)|_U=\ddc \left(-\log |s|_h^2\right).
    \]
\end{definition}

\begin{proposition}[Lelong--Poincaré]\label{prop:LelongPoincare}\index{theorem!Lelong--Poincaré formula}
    Let $s\in \mathrm{H}^0(X,L)$ be non-zero and $h$ be a Hermitian metric on $L$. Then
    \begin{equation}\label{eq:LP}
        c_1(L,h)+\ddc \log|s|_h^2=[Z(s)],
    \end{equation}
    where $Z(s)$ is the zero divisor defined by $s$ and $[\bullet]$ denote the associated current of integration.
\end{proposition}
See \cite[(3.11)]{Dem12} and \cite[Proposition~14.5.4]{FPT}. Again, we want to check that our conventions are compatible by investigating the following simple example.
\begin{example}\label{ex:LPP1}
    Let $X=\mathbb{P}^1$ and $L=\mathcal{O}_{\mathbb{P}^1}(1)$. The homogeneous coordinates on $\mathbb{P}^1$ will be denoted by $[X_0:X_1]$. At a point $x=[X_0:X_1]\in \mathbb{P}^1$, the fiber $L_x$ is identified with the dual of $[x]$, where $[x]\subseteq \mathbb{C}^2$ is the line represented by $x$.

    In order to introduce the Hermitian metric $h$ on $L$, we fix the standard Hermitian norm $\|\bullet\|$ on $\mathbb{C}^2$. Then given $\lambda\in L_x=[x]^{\vee}$, we introduce
    \[
    |\lambda|_{h_x}=\frac{|\lambda (\tilde{x})|}{\|\tilde{x}\|},
    \]
    where $\tilde{x}$ is an arbitrary non-zero element in $[x]$. The reader can easily verify that $h$ is indeed a Hermitian metric on $L$.
    The Hermitian metric $h$ is known as the \emph{Fubini--Study metric}.

    A holomorphic section $s\in \mathrm{H}^0(X,L)$ can be formally identified with a linear form $a_0X_0+a_1X_1$: At $x\in X$, the corresponding linear form on $[x]$ is given by sending $(X_0,X_1)$ to $a_0X_0+a_1X_1$.

    Next we compute $\ddc h=c_1(L,h)$. For this purpose, we cover $\mathbb{P}^1$ by $\mathbb{C}=\mathbb{P}^1\setminus \{\infty\}$ and $\mathbb{P}^1\setminus \{0\}$.
    Both are holomorphic coordinate charts with coordinate function $z=X_0/X_1$ and $z^{-1}=X_1/X_0$ respectively.

    We claim that on $\mathbb{C}$,
    \begin{equation}\label{eq:curvFS}
        \ddc h=\ddc \log (1+|z|^2).
    \end{equation}
    In fact, let $t$ be the nowhere vanishing section of $L$ on $\mathbb{C}$ corresponding to $X_1$. Then for $z\in \mathbb{C}$, we have an obvious lift $(z,1)\in [z]$, so
    \[
    |t|_h^2(z)=\frac{1}{|z|^2+1}.
    \]
    So \eqref{eq:curvFS} follows.

    In order to obtain a non-trivial case of the Lelong--Poincaré formula, we need to consider a section which vanishes at some points in $\mathbb{C}$. Let $s$ be the holomorphic section of $L$ corresponding to $X_0$. Then
    \[
    \log |s|_h^2(z)=\log \frac{|z|^2}{|z|^2+1}
    \]
    for any $z\in \mathbb{C}$ using the same argument as above. Therefore, we find that restricted to $\mathbb{C}$, we have
    \[
    c_1(L,h)+\ddc \log |s|_h^2=\ddc f=[0],
    \]
    where $f(z)=\log|z|^2$. So the Lelong--Poincaré formula \eqref{eq:LP} is verified in this case.

    The K\"ahler form $\ddc h$ on $\mathbb{P}^1$ is also known as the \emph{Fubini--Study metric}.
\end{example}

\begin{definition}\label{def:pshfunctions:L1836}
    A (singular) \emph{plurisubharmonic metric}\index{plurisubharmonic metric} (or \emph{psh metric} for short)\footnote{In the literature, people usually refer to such metrics as \emph{positively curved singular Hermitian metrics}. I dislike this terminology, as having positive curvature only determines a plurisubharmonic metric almost everywhere, not everywhere.} $h$ on $L$ is a family $(h_x)_{x\in X}$ such that
    \begin{enumerate}
        \item\label{itm:chapter-pshfunctions:L1857:1} for each $x\in X$, $h_x$ is either a Hermitian form on $L_x$ or the singular Hermitian form on $L_x$, and
        \item\label{itm:chapter-pshfunctions:L1857:2} there is a Hermitian metric $h_0$ on $L$ and $\varphi\in \PSH(X,c_1(L,h_0))$ such that for each $x\in X$ and each $v\in L_x$, we have
        \begin{equation}\label{eq:htwist}
        |v|^2_{h_x}=
        \left\{
        \begin{aligned}
            0,&\quad \textup{if }v=0;\\
            |v|^2_{h_{0,x}}\mathrm{e}^{-\varphi(x)},&\quad \textup{if }v\neq 0.
        \end{aligned}
        \right.
        \end{equation}
    \end{enumerate}
    The \emph{(first) Chern current} of $h$ is by definition
    \[
    \ddc h=c_1(L,h)\coloneqq c_1(L,h_0)+\ddc\varphi.
    \]
    The set of psh metrics on $L$ will be denoted by $\PSH(L)$\index{$\PSH(L)$}.
\end{definition}
We shall write the plurisubharmonic metric defined by \eqref{eq:htwist} as $h_0\exp(-\varphi)$\footnote{Be careful, this is not $h_0^2\exp(-\varphi)$, as I prefer to think of $h_0$ as a quadratic form.}. With these conventions, $c_1(L,h)$ does not depend on the choice of $h_0$.

Given a local holomorphic section $s$ of $L$, we shall use the following notations interchangeably to denote the map $x\mapsto h_x(s)$:
\[
|s|^2_h=h(s)=h(s,s).
\]

\begin{remark}\label{rmk:pshfunctions:L1864}
    In the literature, some people prefer the convention that in \eqref{eq:htwist}, neither side has the square. Our choice seems to be the most natural one given our normalization of $\ddc$.
\end{remark}

Observe that once a Hermitian metric $h_0$ on $L$ is given, the construction in \ref{itm:chapter-pshfunctions:L1857:2} gives a bijection between $\PSH(X,c_1(L,h_0))$ and the set of plurisubharmonic metrics on $L$.

Based on this correspondence, plurisubharmonic metrics have the same functoriality property as quasi-plurisubharmonic functions in \cref{prop:pullbackpsh}. In particular, given a connected closed submanifold $Y$ of $X$, $h|_Y$ is a plurisubharmonic metric as long as it is not singular everywhere.

\begin{definition}\label{def:tensorprodpshmetric}
    Given two holomorphic line bundles $L_1$, $L_2$ on $X$ and plurisubharmonic metrics $h_1$ on $L_1$ and $h_2$ on $L_2$, we define the \emph{tensor product} plurisubharmonic metric $h_1\otimes h_2$ on $L_1\otimes L_2$ as follows: for each $x\in X$, define
    \index{plurisubharmonic metric!tensor product}
    \[
        (h_1\otimes h_2)_x=h_{1,x}\otimes h_{2,x}
    \]
    in the sense of \cref{def:1dimlinspacehermten}.
\end{definition}
We can easily verify that $h_1\otimes h_2$ is indeed a plurisubharmonic metric on $L_1\otimes L_2$.

\begin{example}\label{ex:logsingP1}
    We continue with our example \cref{ex:LPP1}. Let $X=\mathbb{P}^1$ and $L=\mathcal{O}_{\mathbb{P}^1}(1)$. Let $h^0$ denote the Fubini--Study metric on $L$ as defined in \cref{ex:LPP1}. Note that we have changed the notation from $h$ to $h^0$. Let $\omega=\ddc h^0$.

    We construct $\varphi\in \PSH(X,\omega)$ as follows: On $\mathbb{C}$, define
    \begin{equation}\label{eq:logsingP1ex}
    \varphi(z)=\log \frac{|z|^2}{1+|z|^2}.
    \end{equation}
    Then $\varphi\in \PSH(\mathbb{C},\omega|_{\mathbb{C}})$ by \eqref{eq:curvFS}. Setting $\varphi(\infty)=0$, we can easily verify that $\varphi\in \PSH(\mathbb{P}^1,\omega)$.\footnote{This can also be verified using the Grauert--Remmert extension theorem \cref{thm:GRexten}.}

    We then get a plurisubharmonic metric $h=h^0\exp(-\varphi)$. To be more explicit, the fiber metric $h$ at $0$ is singular, $h_{\infty}=h^0_{\infty}$, while for $z\in \mathbb{C}\setminus \{0\}$ and $\lambda\in [z]^{\vee}$, we have
    \[
    |\lambda|_{h_z}=\frac{|\lambda(z,1)|}{|z|}.
    \]
\end{example}

In the remainder of this section, we assume that $X$ is compact.

\begin{definition}\label{def:ArcFSmetric}Assume that $L$ is a pseudo-effective line bundle on $X$.
    A \emph{Fubini--Study metric}\index{Fubini--Study metric} on $L$ is a psh metric $h$ on $L$ of the following form: There exists $m\in \mathbb{Z}_{>0}$, finitely many sections $s_1,\ldots,s_N\in \mathrm{H}^0(X,L^m)$ and $\lambda_1,\ldots,\lambda_N\in \mathbb{Q}$ such that for any local nowhere vanishing holomorphic section $s$ of $L$, we have
    \[
    |s|_h^2=\min_{i=1,\ldots,N}\left| \frac{s^{\otimes m}}{\mathrm{e}^{\lambda_i/2}s_i}\right|^{2m^{-1}}.
    \]
    We write $\FS(L)$ for the set of Fubini--Study metrics on $L$.
\end{definition}
If we fix a reference smooth Hermitian metric $h_0$ on $L$ with $\theta=\ddc h_0$, we can write $h=h_0\exp(-\varphi)$ with
\[
\varphi=\frac{1}{m}\max_{i=1,\ldots,N}\left(\log |s_i|^2_{h_0^{m}}+\lambda_i\right).
\]
Similarly, we write $\FS(X,\theta)$ for the set of such functions.

\begin{definition}\label{def:gFSmetric}
Assume that $L$ is a pseudo-effective line bundle on $X$.
    The set $\widetilde{\FS}(L)$ of \emph{generalized Fubini--Study metrics}\index{generalized Fubini--Study metric} is the smallest subset of $\PSH(L)$ containing $\FS(L)$ which is closed under the following two operations:
    \begin{enumerate}
        \item\label{itm:def:gFSmetric:1} $\mathbb{Q}$-convex combinations: if $h_1,h_2\in \widetilde{\FS}(L)$ and $t\in (0,1)\cap\mathbb{Q}$, then
        \[
        h_1^t \otimes h_2^{1-t}\in \widetilde{\FS}(L);
        \]
        \item\label{itm:def:gFSmetric:2} minima: if $h_1,h_2\in \widetilde{\FS}(L)$, then
        \[
        \min\{h_1,h_2\}\in \widetilde{\FS}(L).
        \]
    \end{enumerate}
\end{definition}

We shall need the following Ohsawa--Takegoshi type extension theorem.
\begin{theorem}\label{thm:OT_ext} Assume that $L$ is big and $T$ is a holomorphic line bundle on $X$. Fix a Hermitian metric $h_T$ on $T$.
Fix a Hermitian metric $h$ on $L$ and let $\theta=c_1(L,h)$.
Take a Kähler form $\omega$ on $X$. Let $Y\subseteq X$ be a connected closed submanifold of dimension $m$.
Suppose that $\varphi \in \PSH(X,\theta - \delta \omega)$ for some $\delta>0$ and $\varphi|_Y\not\equiv -\infty$.
Then there exists $k_0(\delta,h_T)>0$ such that for all $k \geq k_0$ and $s \in \mathrm{H}^0(Y, T\otimes L|_Y^k \otimes \mathcal{I}(k\varphi|_Y))$\footnote{Here and in the sequel, we usually abbreviate $\otimes k$ in the superscript as $k$ to save space.}, there exists an extension $\tilde s \in \mathrm{H}^0(X, T\otimes L^k \otimes \mathcal{I}(k\varphi))$ such that
\[
\int_X |\tilde{s}|^2_{h^k\otimes h_T}\mathrm{e}^{-k\varphi}\,\omega^n \leq C \int_Y |s|^2_{(h^k\otimes h_T)|_Y} \mathrm{e}^{-k\varphi|_Y}\,\omega|_Y^m,
\]
where $C>0$ is an absolute constant, independent of the data $(\varphi,s,k)$.
\end{theorem}
This is a special case of \cite[Theorem~1.4]{His12}. See also \cite[Corollary~14.7.9]{FPT}.

Let us also recall a simple result about the Bergman kernel approximations:
\begin{proposition}\label{prop:Bergman_approx}
Let $(L,h)$ be a Hermitian line bundle on $X$ and set $\theta=c_1(L,h)$. Let $(T,h_T)$ be a Hermitian line bundle on $X$. Assume that $\varphi\in \PSH(X,\theta)$ is a potential with analytic singularities such that $\theta_{\varphi}$ is a K\"ahler current. Fix a K\"ahler form $\omega$ on $X$.
For each $k\geq 1$, we let
\begin{equation}\label{eq:Bergman_seq_def}
\varphi_k \coloneqq \frac{1}{k} \log \sup_{\substack{s \in \mathrm{H}^0(X,L^k \otimes T)\\ \int_X |s|^2_{h^k  \otimes h_T}\mathrm{e}^{-k\varphi} \omega^n \leq 1}} |s|^2_{h^k  \otimes h_T}.
\end{equation}
Then for all sufficiently large $k$,
\[
\varphi\preceq \varphi_k\preceq \alpha_k \varphi,
\]
where $\alpha_k\in (0,1)$ is an increasing sequence with limit $1$.
\end{proposition}
\begin{proof}

First we prove $\varphi\preceq \varphi_k$. Since $\theta_{\varphi}$ is a K\"ahler current, there is $\delta>0$ such that $\varphi\in \PSH(X,\theta-\delta\omega)$. Let $x\in X$ with $\varphi(x)>-\infty$. Applying \cref{thm:OT_ext} to the point $Y=\{x\}$, for all sufficiently large $k$ and every unit vector $e_x\in (L^k\otimes T)_x$, there is a section $\tilde{s}\in \mathrm{H}^0(X,L^k\otimes T\otimes \mathcal{I}(k\varphi))$ with $\tilde{s}(x)=e_x$ and
\[
    \int_X |\tilde{s}|^2_{h^k\otimes h_T}\mathrm{e}^{-k\varphi}\omega^n\leq C\mathrm{e}^{-k\varphi(x)}.
\]
After rescaling $\tilde{s}$, the definition of $\varphi_k$ gives
\[
    \varphi_k(x)\geq \varphi(x)-\frac{\log C}{k}.
\]
This proves $\varphi\preceq\varphi_k$.

The other part of the inequality follows from \cref{lma:IandIinf} and \eqref{eq:IIinf_temp1}.

More precisely, if $s_1,\ldots,s_N$ is an orthonormal basis of $\mathrm{H}^0(X,L^k\otimes T\otimes \mathcal{I}(k\varphi))$, with respect to the $L^2$-norm:
\[
s\mapsto \left( \int_X |s|^2_{h^k\otimes h_T}\mathrm{e}^{-k\varphi}\,\omega^n \right)^{1/2},
\]
then by \eqref{eq:IIinf_temp1},
\[
s_i\in \mathrm{H}^0\left(X,L^k\otimes T\otimes \mathcal{I}_{\infty}\left((k-c^{-1}n)_+\varphi\right)\right),
\]
where $c>0$ depends only on $\varphi$. So there is $C_k>0$ so that
\[
\sum_{i=1}^N |s_i|^2_{h^k\otimes h_T}\mathrm{e}^{-(k-c^{-1}n)_+\varphi}\leq C_k\quad \textup{quasi-everywhere}.
\]
In other words,
\[
\varphi_k\preceq \frac{1}{k}\left( k-c^{-1}n \right)_+ \varphi.
\]
It suffices to take $\alpha_k=(1-k^{-1}c^{-1}n)_+$.
\end{proof}

\motto{Pour exprimer d'une manière frappante que le monument que j'élève sera placé sous l'invocation de la Science, j'ai décidé d'inscrire en lettres d'or sur la grande frise du premier étage et à la place d'honneur, les noms des plus grands savants\footnote{Gaspard Monge, Comte de Péluse (1746--1818), known oddly by his family name instead of \emph{de Péluse}, is one of the 72 names scribed on the Eiffel tower. Apart from being a mathematician, he was also a politician, active mainly after the French Revolution.} qui ont honoré la France depuis 1789 jusqu'à nos jours.

---  Gustave Eiffel, 1889}
\chapter{Non-pluripolar products}
\label{chap:npp}

Let $X$ be a complex manifold and let $\varphi_1,\ldots,\varphi_p\in\PSH(X)$. When the functions $\varphi_i$ are smooth, the differential form
\begin{equation}\label{eq:mixedMAtype}
    \ddc\varphi_1\wedge\cdots\wedge\ddc\varphi_p
\end{equation}
is defined by the usual calculus. A central problem in pluripotential theory is to extend this product to singular plurisubharmonic functions in a reasonable way.

In this book we use the non-pluripolar product. This construction, originating in the work of Bedford--Taylor and developed in the global setting by Guedj--Zeriahi and Boucksom--Eyssidieux--Guedj--Zeriahi, has two features that are essential for us: it is defined for arbitrary quasi-plurisubharmonic singularities, and it satisfies the monotonicity theorem for non-pluripolar masses. These two properties make it possible to compare singularity types through Monge--Amp\`ere mass.

The chapter is organized as follows. We first recall the Bedford--Taylor product for locally bounded plurisubharmonic functions. We then introduce the non-pluripolar product and record its basic functorial and locality properties. After a short discussion of quasi-continuous functions and capacity, we prove the convergence and concentration results used later in the theory of envelopes. The key tools extracted from this chapter are the monotonicity theorem, semicontinuity theorem, the concentration theorem for envelopes and the domination principle.

The reader may consult \cite{BEGZ10} and \cite{DDNLsurv} for broader background. My companion book \emph{Foundations of pluripotential theory} \cite[Chapters~9--10, 15--16]{FPT} gives detailed proofs of the Bedford--Taylor, non-pluripolar product and capacity results used here. Our presentation is adapted to the later use of envelopes and singularity types, so some estimates are stated in a form slightly more specific than in the standard references.

\section{Bedford--Taylor theory}\label{sec:BTtheory}
Let $X$ be a complex manifold and $\varphi_1,\ldots,\varphi_p\in \PSH(X)$ ($p\in \mathbb{N}$) be locally bounded plurisubharmonic functions on $X$\footnote{In the literature, some people use $\PSH(X)\cap L^{\infty}_{\loc}(X)$ to denote the set of such functions, which is an abuse of notation. However, this is legitimate thanks to the rigidity \cref{thm:psh_rigid_gen}.}. In this case, there is a canonical definition of the Monge--Amp\`ere type product \eqref{eq:mixedMAtype}.
\begin{definition}\label{def:npp:L22}
    We define the closed positive $(p,p)$-current \eqref{eq:mixedMAtype} on $X$ as follows: We proceed by induction on $p\geq 0$. When $p=0$, we define \eqref{eq:mixedMAtype} as the $(0,0)$-current $[X]$. When $p>0$, we let
    \[
    \ddc \varphi_1\wedge \cdots \wedge \ddc\varphi_p\coloneqq \ddc \left( \varphi_1\,\ddc \varphi_2 \wedge \cdots \wedge \ddc\varphi_p \right).
    \]
    \index{$\ddc \varphi_1\wedge \cdots \wedge \ddc\varphi_p$}
    We call this product the \emph{Bedford--Taylor product} \index{product!Bedford--Taylor product}.
\end{definition}

\begin{remark}\label{rmk:extensionBT}
There is also a slightly more general version of this construction. Given a closed positive current $T$, one can also define the product
\[
\ddc \varphi_1\wedge \cdots \wedge \ddc\varphi_p \wedge T
\]
in a very similar way.
\end{remark}

\begin{proposition}\label{prop:npp:L40}
    The product $\ddc \varphi_1\wedge \cdots \wedge \ddc\varphi_p$ is a closed positive $(p,p)$-current on $X$. Moreover, the product is symmetric in the $\varphi_i$'s.
\end{proposition}
See \cite[Proposition~3.3, Corollary~3.12]{GZ17} and \cite[Theorem~9.2.5]{FPT}.\footnote{In \cite[Proposition~3.3]{GZ17}, the authors claim that $\mathrm{d}u\wedge \mathrm{d}^{\mathrm{c}}u\wedge T$ is closed, which is wrong already in the smooth case. See \cite[Lemma~9.13.1 and Section~B.15]{FPT} for the gradient construction.} The proof relies crucially on an important estimate, known as the \emph{Chern--Levine--Nirenberg inequality}. See \cite{CLN69}, \cite[Theorem~3.9]{GZ17} and \cite[Theorem~9.1.3]{FPT}.

The Bedford--Taylor theory has many satisfactory properties.

\begin{theorem}\label{thm:contMA}
    Let $(\varphi_i^j)_{j\in \mathbb{Z}_{>0}}$ be decreasing sequences (resp.\ increasing sequences) of locally bounded psh functions on $X$ converging (resp.\ converging a.e.) to locally bounded psh functions $\varphi_i$'s, where $i=0,\ldots,p$. Then
    \[
    \varphi_0^j\, \ddc \varphi_1^j\wedge \cdots \wedge \ddc\varphi_p^j\rightharpoonup \varphi_0\, \ddc \varphi_1\wedge \cdots \wedge \ddc\varphi_p
    \]
    as $j\to\infty$. In particular, if $(\varphi_0^j)_j$ is the constant sequence $1$, we have
    \[
    \ddc \varphi_1^j\wedge \cdots \wedge \ddc\varphi_p^j\rightharpoonup \ddc \varphi_1\wedge \cdots \wedge \ddc\varphi_p.
    \]
\end{theorem}
Here the notation $\rightharpoonup$ denotes the weak-* convergence of currents.

We refer to \cite[Theorem~3.18, Theorem~3.23]{GZ17} and \cite[Theorems~9.2.5 and~9.5.1]{FPT} for the proofs.

By contrast, we emphasize that the Bedford--Taylor product is not continuous with respect to the $L^1_{\loc}$-convergence in general. A simple example can be found in \cite[Example~3.25]{GZ17}. For the continuity theorem under convergence in capacity, see \cite[Theorem~9.4.8]{FPT}.

\begin{lemma}\label{lma:boundedfibrefubini}
    Let $Q\subseteq\mathbb{C}^n$ and $D\subseteq\mathbb{C}^m$ be domains, let $p_2\colon Q\times D\to D$ be the second projection, and let $\eta$ be a smooth real $(m,m)$-form on $D$. Let $u_1,\ldots,u_n$ be locally bounded psh functions on $Q\times D$, and put
    \[
        M=\ddc u_1\wedge\cdots\wedge\ddc u_n,\qquad
        M_z=\ddc_x(u_1)_z\wedge\cdots\wedge\ddc_x(u_n)_z.
    \]
    Here $(u_r)_z(x)=u_r(x,z)$.
    For every bounded Borel function $h$ with compact support in $Q\times D$,
    \begin{equation}\label{eq:boundedfibrefubini}
        \int_{Q\times D}h\,p_2^*\eta\wedge M
        =\int_D\left(\int_Q h_z\,M_z\right)\eta(z),
        \qquad h_z(x)=h(x,z).
    \end{equation}
\end{lemma}
\begin{proof}
    By linearity, it suffices to consider $\eta\geq0$.
    We first prove \eqref{eq:boundedfibrefubini} when $h$ is continuous.
    On a product neighborhood of $\Supp h$, choose smooth psh functions $u_r^{(l)}\searrow u_r$, $1\leq r\leq n$. For the smooth approximants, decompose each $\ddc u_r^{(l)}$ according to the $x$- and $z$-directions. After wedging with $p_2^*\eta$, every term containing a factor $\mathrm{d}z_\alpha$ or $\mathrm{d}\bar z_\beta$ vanishes; hence only the factors $\ddc_x(u_r^{(l)})_z$ remain, and \eqref{eq:boundedfibrefubini} follows from the usual Fubini theorem. As $l\to\infty$, the Bedford--Taylor monotone continuity theorem \cref{thm:contMA} gives convergence on $Q\times D$ and on every fibre. Fix $K\Subset Q$ and $L\Subset D$ containing the projections of $\Supp h$. The Chern--Levine--Nirenberg estimate \cite[Theorem~3.9]{GZ17} (see also \cite[Theorem~9.1.3]{FPT}), together with the local uniform bounds for $u_r^{(l)}$, bounds the masses of the fibre products on $K$ uniformly for $z\in L$ and $l$. Since $h_z=0$ for $z\notin L$, dominated convergence proves \eqref{eq:boundedfibrefubini}.

    It remains to pass from continuous to bounded Borel functions. Fix relatively compact open sets $U\Subset Q$ and $V\Subset D$. If $\chi\in C_c(U)$, then
    \[
        F_\chi(z)\coloneqq\int_U\chi\,M_z,
        \qquad z\in V,
    \]
    is Borel measurable. Indeed, choose the smooth decreasing approximations on a neighborhood of $\overline U\times\overline V$ and set
    \[
        M_z^{(l)}\coloneqq
        \ddc_x(u_1^{(l)})_z\wedge\cdots\wedge\ddc_x(u_n^{(l)})_z.
    \]
    The functions $z\mapsto\int_U\chi\,M_z^{(l)}$ are continuous on $V$ and converge pointwise to $F_\chi$ by \cref{thm:contMA}.

    Let $\mathcal{A}_{U,V}$ be the vector space of bounded Borel functions $\chi$ on $U$ for which $F_\chi$ is Borel measurable on $V$. It is closed under bounded increasing pointwise limits of non-negative functions by the monotone convergence theorem. Taking compactly supported continuous cut-off functions increasing to $1$ on $U$ shows that $1\in\mathcal{A}_{U,V}$. Thus $\mathcal{A}_{U,V}$ contains the algebra $\mathbb{R}+C_c(U)$, which generates the Borel $\sigma$-algebra of $U$. The functional monotone class theorem therefore shows that $\mathcal{A}_{U,V}$ contains every bounded Borel function on $U$.

    Next, let $\mathcal{B}_{U,V}$ be the vector space of bounded Borel functions $g$ on $U\times V$ such that
    \[
        z\longmapsto\int_U g_z\,M_z,
        \qquad g_z(x)=g(x,z),
    \]
    is Borel measurable on $V$. This space is again closed under bounded increasing pointwise limits of non-negative functions. It contains every product $g(x,z)=\chi(x)\psi(z)$ of bounded Borel functions, by the preceding paragraph. Finite sums of such products form an algebra containing the indicators of all Borel rectangles $A\times B\subseteq U\times V$; hence this algebra generates the Borel $\sigma$-algebra of $U\times V$. A second application of the functional monotone class theorem shows that $\mathcal{B}_{U,V}$ contains every bounded Borel function on $U\times V$.

    To pass from $U\times V$ to $Q\times D$, choose exhaustions $U_j\Subset Q$ and $V_j\Subset D$ by relatively compact open sets. For a Borel set $E\subseteq Q\times D$, the preceding paragraph, applied to the restriction of $\mathds{1}_E$ to $U_j\times V_j$, shows that
    \[
        G_j(z)\coloneqq
        \mathds{1}_{V_j}(z)\int_{U_j}\mathds{1}_E(x,z)\,M_z
    \]
    is Borel measurable on $D$. We may take both exhaustions increasing, so $G_j\nearrow M_z(E_z)$ pointwise, where
    \[
        E_z\coloneqq\{x\in Q:(x,z)\in E\}.
    \]
    Thus $z\mapsto M_z(E_z)$ is Borel measurable on $D$.
    We may therefore define
    \[
        \lambda(E)\coloneqq\int_D M_z(E_z)\,\eta(z).
    \]
    The monotone convergence theorem, first on each fibre and then on $D$, shows that $\lambda$ is a Borel measure on $Q\times D$. The local uniform mass bound above shows that $\lambda$ is locally finite and hence Radon. Approximation by non-negative simple functions gives
    \[
        \int_{Q\times D}h\,\mathrm{d}\lambda
        =\int_D\left(\int_Qh_z\,M_z\right)\eta(z)
    \]
    for every non-negative Borel function $h$. For continuous compactly supported $h$, the first part of the proof shows that $\lambda$ agrees with $p_2^*\eta\wedge M$. The uniqueness of Radon measures then gives the same equality for every bounded Borel function $h$ with compact support.
\end{proof}

\section{The non-pluripolar products}\label{sec:npp}
The proofs of the standard results recalled in this section can be found in \cite{BEGZ10} and \cite[Sections~15.1--15.5]{FPT}.

Let $X$ be a complex manifold.

\begin{definition}\label{def:npp:L69}
Let $\varphi_1,\ldots,\varphi_p\in \PSH(X)$. We set
\[
O_k \coloneqq \bigcap_{j=1}^p \{\varphi_j>-k\},\quad k\in \mathbb{N}.
\]
We say that $\ddc \varphi_1\wedge\cdots\wedge \ddc \varphi_p$ is \emph{well-defined} if for each connected open subset $U\subseteq X$, any smooth positive $(1,1)$-form $\omega$ on $U$, and each compact subset $K\subseteq U$, we have
\begin{equation}\label{eq:welldefinepluri}
\sup_{k\geq 0} \int_{K\cap O_k}\left.\left(\bigwedge_{j=1}^p \ddc (\varphi_j\lor (-k))\right)\right|_U\wedge \omega^{\dim U-p}<\infty.
\end{equation}

In this case, we define the \emph{non-pluripolar product} \index{product!non-pluripolar product} $\ddc \varphi_1\wedge\cdots\wedge \ddc \varphi_p$ by
\begin{equation}\label{eq:npp}
\mathds{1}_{O_k}\,\ddc \varphi_1\wedge\cdots\wedge \ddc \varphi_p=\mathds{1}_{O_k}\bigwedge_{j=1}^p \ddc  \left( \varphi_j\lor (-k)\right)
\end{equation}
on $\bigcup_{k\geq 0} O_k$ and extend it by zero to $X$.
\end{definition}

As recalled in \cref{sec:plurifine}, an $\mathcal{F}$-open subset means an open subset with respect to the plurifine topology.
\begin{proposition} \label{prop:npp1}
Let $\varphi_1,\ldots,\varphi_p\in \PSH(X)$.
\begin{enumerate}
\item\label{itm:prop:npp1:1} The product $\ddc \varphi_1\wedge\cdots\wedge \ddc \varphi_p$ is local with respect to the plurifine topology in the following sense: Let $O\subseteq X$ be an $\mathcal{F}$-open subset and $\psi_1,\ldots,\psi_p\in \PSH(X)$. Assume that
\[
\varphi_j|_{O}=\psi_j|_{O},\quad j=1,\ldots,p,
\]
and that
\[
\bigwedge_{j=1}^p \ddc \varphi_j\textup{ and } \bigwedge_{j=1}^p \ddc \psi_j
\]
are both well-defined, then
    \begin{equation}\label{eq:ppp1}
    \left.\left(\bigwedge_{j=1}^p \ddc \varphi_j\right)\right|_{O}=\left.\left(\bigwedge_{j=1}^p \ddc \psi_j\right)\right|_{O}.
    \end{equation}
If furthermore $O$ is open in the usual topology, then the product
\[
\bigwedge_{j=1}^p \ddc \varphi_j|_{O}
\]
on $O$ is well-defined and
\begin{equation}\label{eq:ppp2}
\left. \left(\bigwedge_{j=1}^p \ddc \varphi_j\right)\right|_{O}=\bigwedge_{j=1}^p \ddc \varphi_j|_{O}.
\end{equation}
Let $\mathcal{U}$ be an open covering of $X$. Then $\ddc \varphi_1\wedge\cdots\wedge \ddc \varphi_p$
is well-defined if and only if each product
\[
\bigwedge_{j=1}^p \ddc \varphi_j|_{U},\quad U\in \mathcal{U}.
\]
    \item\label{itm:prop:npp1:2} The current $\ddc \varphi_1\wedge\cdots\wedge \ddc \varphi_p$ and the fact that it is well-defined depend only on the currents $\ddc \varphi_j$, not on the choice of the $\varphi_j$'s nor on the ordering of the $\varphi_j$'s.
    \item\label{itm:prop:npp1:3} When $\varphi_1,\ldots,\varphi_p\in L^{\infty}_{\loc}(X)$, the non-pluripolar product $\ddc \varphi_1\wedge\cdots\wedge \ddc \varphi_p$ is well-defined and is equal to the Bedford--Taylor product. More generally, if $q\leq p$, and we only assume that $\varphi_1,\ldots,\varphi_q\in L^{\infty}_{\loc}(X)$, and that the non-pluripolar product $T\coloneqq \ddc \varphi_{q+1}\wedge \cdots \wedge \ddc \varphi_p$ is well-defined, then the non-pluripolar product $\ddc \varphi_1\wedge\cdots\wedge \ddc \varphi_p$ is also well-defined, and is equal to the Bedford--Taylor product $\ddc \varphi_1\wedge\cdots\wedge \ddc \varphi_q\wedge T$ as in \cref{rmk:extensionBT}.
    \item\label{itm:prop:npp1:4} Assume that $\ddc \varphi_1\wedge\cdots\wedge \ddc \varphi_p$ is well-defined. Then $\ddc \varphi_1\wedge\cdots\wedge \ddc \varphi_p$ puts no mass on pluripolar sets.
    \item\label{itm:prop:npp1:5} Assume that $\ddc \varphi_1\wedge\cdots\wedge \ddc \varphi_p$ is well-defined. Then $\bigwedge_{j=1}^p \ddc \varphi_j$ is a closed positive $(p,p)$-current on $X$.\footnote{This current is in fact strongly positive, as shown in \cite{Xiabdiv2} and \cite[Proposition~15.1.4]{FPT}.}
    \item\label{itm:prop:npp1:6} The product is multilinear: Let $\psi_1\in \PSH(X)$, $a,b>0$ then
\begin{equation}\label{eq:ppp6}
\ddc (a\varphi_1+b\psi_1)\wedge \bigwedge_{j=2}^p \ddc \varphi_j =a\,\ddc \varphi_1\wedge \bigwedge_{j=2}^p \ddc \varphi_j + b\,\ddc \psi_1\wedge \bigwedge_{j=2}^p \ddc \varphi_j
\end{equation}
in the sense that the left-hand side is well-defined if and only if both terms on the right-hand side are well-defined, and the equality holds in that case.
\end{enumerate}
\end{proposition}
In view of \ref{itm:prop:npp1:3}, we do not need to specify whether our product $\ddc \varphi_1\wedge\cdots\wedge \ddc \varphi_p$ is the Bedford--Taylor product or the non-pluripolar product when the $\varphi_i$'s are all locally bounded.

\begin{example}\label{ex:1Dnpp}
    Let $\varphi\in \PSH(X)$, then
    \begin{equation}\label{eq:1dimnpp}
        (\ddc \varphi)^1=\mathds{1}_{\{\varphi>-\infty\}}\,\ddc\varphi.
    \end{equation}
    Here $(\ddc \varphi)^1$ on the left-hand side denotes the non-pluripolar product, while $\ddc\varphi$ on the right-hand side denotes the usual derivative in the sense of current.

    The problem is local, so we may assume that $X=\Delta^n$, the unit polydisk in $\mathbb{C}^n$. Take $f\in C^{\infty}_c(\Delta^n)$ with $f\geq 0$. We need to show that
    \[
        \lim_{k\to\infty} \int_{\Delta^n} \mathds{1}_{\{\varphi>-k\}}f\, \ddc \bigl(\varphi\lor (-k)\bigr)\wedge \omega^{n-1}=\int_{\Delta^n} \mathds{1}_{\{\varphi>-\infty\}}f\,\ddc\varphi\wedge \omega^{n-1},
    \]
    where $\omega$ is the standard Kähler form on $\mathbb{C}^n$. By the dominated convergence theorem, this reduces immediately to the $1$-dimensional case:
    \[
        \lim_{k\to\infty} \int_{\Delta^1} \mathds{1}_{\{\varphi>-k\}}f\, \Delta \bigl(\varphi\lor (-k)\bigr)=\int_{\Delta^1} \mathds{1}_{\{\varphi>-\infty\}}f\,\Delta\varphi,
    \]
    for any subharmonic function $\varphi$ on the unit disk $\Delta^1$, which is a special case of \cite[Proposition~3.3]{EK24}. See also \cite[Proposition~15.3.6]{FPT}.
\end{example}

\begin{definition}\label{def:npp:L146}
Let $T_1,\ldots,T_p$ be closed positive $(1,1)$-currents on $X$. We say that $T_1\wedge\cdots\wedge T_p$ is \emph{well-defined} if there exists an open covering $\mathcal{U}$ of $X$ such that on each $U\in \mathcal{U}$, we can find $\varphi_j^U\in \PSH(U)$ ($j=1,\ldots,p$) such that
\[
\ddc \varphi_j^U=T_j,\quad j=1,\ldots,p
\]
and $\ddc \varphi_1^U\wedge \cdots\wedge \ddc \varphi_p^U$ is well-defined.
In this case, we define the \emph{non-pluripolar product} \index{product!non-pluripolar product} $T_1\wedge\cdots\wedge T_p$ as the closed positive $(p,p)$-current on $X$ defined by
\begin{equation}\label{eq:ppp5}
 \left(T_1\wedge\cdots\wedge T_p\right)|_U=\ddc \varphi_1^U\wedge \cdots\wedge \ddc \varphi_p^U,\quad U\in \mathcal{U}.
\end{equation}
\end{definition}
The product $T_1\wedge\cdots\wedge T_p$ is independent of the choices we made thanks to \itemref{itm:prop:npp1:1} and \ref{itm:prop:npp1:2}.

\cref{prop:npp1} can be formulated in terms of currents without any difficulty.

\begin{remark}\label{rmk:npp:L161}
    Similar to \cref{rmk:extensionBT}, there is also an extension of the non-pluripolar theory allowing us to define
    \[
    T_1\wedge \cdots \wedge T_p\cap T
    \]
    for any closed positive current $T$. This is the \emph{relative non-pluripolar product} introduced by Vu \cite{Vu20}. For the relative Bedford--Taylor construction, see \cite[Definition~9.1.1]{FPT}. Unlike the relative Bedford--Taylor products, the relative non-pluripolar products present some pathological behaviors. For example, they are not linear in general.
\end{remark}

\begin{remark}\label{rmk:npp:L169}
    Another possible generalization of the non-pluripolar products is motivated by \itemref{itm:prop:npp1:1}. One could begin by defining a generalized notion of plurisubharmonic functions on $\mathcal{F}$-open sets, called \emph{$\mathcal{F}$-plurisubharmonic functions} and define their non-pluripolar products. See \cite{EKFW11, EKW14}. For plurifine locality of the ordinary Bedford--Taylor product, see \cite[Theorem~10.4.3]{FPT}.
\end{remark}

\begin{proposition}\label{prop:nppwelldef}
    Let $X$ be a compact Kähler manifold and let $T_1,\ldots,T_p$ be closed positive $(1,1)$-currents on $X$. Then $T_1\wedge \cdots\wedge T_p$ is well-defined.
\end{proposition}
This proposition explains why we usually work in the setting of compact Kähler manifolds.

\begin{corollary}\label{cor:compactfibreslicing}
    Let $X$ and $Z$ be connected compact Kähler manifolds of dimensions $n$ and $m$, respectively, and let $p_1\colon X\times Z\to X$ and $p_2\colon X\times Z\to Z$ be the projections. Let $\theta$ be a closed real smooth $(1,1)$-form on $X$, let $\phi\in\PSH(X,\theta)$, and let $\gamma$ be a closed smooth positive $(1,1)$-form on $Z$. Put
    \[
        Y=X\times Z,\qquad
        \Theta=p_1^*\theta+p_2^*\gamma,\qquad
        \Phi=p_1^*\phi.
    \]
    Suppose that $W\in\PSH(Y,\Theta)$ and $\Phi-C\leq W\leq\Phi+C$ for some $C>0$. Set
    \[
        T=\Theta+\ddc W,\qquad S=\Theta+\ddc\Phi,
    \]
    and write $W_z(x)=W(x,z)$. Then $W_z\in\PSH(X,\theta)$ for every $z\in Z$. For $0\leq j\leq n$ and every smooth real $(m,m)$-form $\eta$ on $Z$,
    \begin{equation}\label{eq:compactrelativeenergyslicing}
        \begin{aligned}
        &\int_Y(W-\Phi)p_2^*\eta\wedge T^j\wedge S^{n-j}\\
        &\qquad=\int_Z\left(\int_X(W_z-\phi)\,
        \theta_{W_z}^j\wedge\theta_\phi^{n-j}\right)\eta(z),
        \end{aligned}
    \end{equation}
    and
    \begin{equation}\label{eq:compactfibreslicingconstant}
        \int_Y p_2^*\eta\wedge T^n
        =\int_Z\left(\int_X\theta_{W_z}^n\right)\eta(z).
    \end{equation}
\end{corollary}
\begin{proof}
    By linearity, it suffices to consider $\eta\geq0$.
    All the non-pluripolar products that occur below are well-defined by \cref{prop:nppwelldef}.
    We work on a relatively compact product chart $Q\times D\subseteq X\times Z$ and let $a$ and $b$ be smooth local potentials of $p_1^*\theta$ and $p_2^*\gamma$. After subtracting constants, put $w=W+a+b\leq0$ and $f=\Phi+a+b\leq0$. For $k\geq1$, set
    \[
        \begin{aligned}
        w_k&=w\lor(-k),& f_k&=f\lor(-k),\\
        O_k&=\{w>-k\}\cap\{f>-k\},&
        \mu_{j,k}&=\mathds{1}_{O_k}(\ddc w_k)^j\wedge(\ddc f_k)^{n-j}.
        \end{aligned}
    \]
    For $z\in D$, write $w_z(x)=w(x,z)$ and $f_z(x)=f(x,z)$, and set
    \[
        \mu_{j,k,z}
        =\mathds{1}_{\{w_z>-k\}\cap\{f_z>-k\}}
        \bigl(\ddc_x(w_z\lor(-k))\bigr)^j\wedge
        \bigl(\ddc_x(f_z\lor(-k))\bigr)^{n-j}.
    \]
    Apply \cref{lma:boundedfibrefubini} with $u_1=\cdots=u_j=w_k$, $u_{j+1}=\cdots=u_n=f_k$, and
    \[
        h=\rho G\mathds{1}_{O_k},\qquad
        \rho\in\mathcal{C}^{\infty}_c(Q\times D),\quad \rho\geq0,
    \]
    where $G$ is any bounded non-negative Borel function. This gives
    \begin{equation}\label{eq:finitelevelfubini}
        \int_{Q\times D}\rho G\,p_2^*\eta\wedge\mu_{j,k}
        =\int_D\left(\int_Q\rho_zG_z\,\mu_{j,k,z}\right)\eta(z).
    \end{equation}
    By plurifine locality, $\mu_{j,k}$ increases to the local expression of $T^j\wedge S^{n-j}$, while, for every $z\in D$, $\mu_{j,k,z}$ increases to the local expression of $\theta_{W_z}^j\wedge\theta_\phi^{n-j}$. Taking $G=\delta+C\geq0$ in \eqref{eq:finitelevelfubini}, applying monotone convergence, subtracting the corresponding constant-weight identity, and summing over a non-negative partition of unity proves \eqref{eq:compactrelativeenergyslicing}. The same argument with constant weight and degree $n$ proves \eqref{eq:compactfibreslicingconstant}.
\end{proof}

\section{Quasi-continuous functions}\label{sec:qc}

Let $X$ be a compact Kähler manifold of dimension $n$ and $\theta$ be a closed real smooth $(1,1)$-form on $X$ representing a big cohomology class.

\begin{definition}\label{def:cap}
    Let $A\subseteq X$ be a Borel subset. The \emph{$\theta$-capacity}\index{capacity!theta-capacity@$\theta$-capacity} $\Capa_{\theta}(A)$\index{$\Capa_{\theta}$} of $A$ is defined as
    \[
    \Capa_{\theta}(A)\coloneqq \sup\left\{ \int_A\theta_{\varphi}^n : \varphi\in \PSH(X,\theta),V_{\theta}-1\leq \varphi\leq V_{\theta} \right\}.
    \]
\end{definition}
Here $V_{\theta}$ is defined in \eqref{eq:Vtheta} below. As we will recall below, the non-pluripolar product $\theta_{\varphi}^n$ is well-defined.

The capacity is not very sensitive to the choice of $\theta$:
\begin{theorem}\label{thm:comp_capa}
    Let $\theta'$ be another closed real smooth $(1,1)$-form on $X$ representing a big cohomology class. Then there are continuous functions $f,g\colon [0,\infty)\rightarrow [0,\infty)$, with $f(0)=g(0)=0$, such that for any Borel subset $E\subseteq X$, we have
    \[
    \Capa_{\theta}(E)\leq f\left( \Capa_{\theta'}(E)\right),\quad \Capa_{\theta'}(E)\leq g\left( \Capa_{\theta}(E)\right).
    \]
\end{theorem}
A more general result is proved in \cite{Lu21}. See also \cite[Theorem~16.9.8]{FPT}. Similar comparison results hold between $\theta$-capacity and the classical Bedford--Taylor capacity, see \cite[Section~9.2]{GZ17} and \cite[Proposition~16.9.9]{FPT} for the proof. As a consequence, we can freely apply the results in \cite[Section~4.2]{GZ17} and \cite[Sections~9.3--9.4]{FPT}, even though capacity has a different meaning there.

\begin{definition}\label{def:qc_fun}
Let $U$ be an open subset of $X$.
    A function $f\colon U\rightarrow [-\infty,\infty]$ is \emph{quasi-continuous}\index{quasi-continuous function} if it is Borel measurable and for any $\epsilon>0$, there is an open subset $G\subseteq U$ such that
    \begin{enumerate}
        \item\label{itm:def:qc_fun:1} $\Capa_{\theta}(G)\leq \epsilon$;
        \item\label{itm:def:qc_fun:2} $f|_{U\setminus G}$ is real-valued and continuous.
    \end{enumerate}
\end{definition}
Thanks to \cref{thm:comp_capa}, the notion of quasi-continuous functions is independent of the choice of $\theta$. Note that if $f,g\colon U\rightarrow [-\infty,\infty]$ are two Borel measurable functions equal quasi-everywhere (see \cref{def:pluripolarsets}), then $f$ is quasi-continuous if and only if $g$ is.

\begin{theorem}\label{thm:pp_cap}
    Let $A\subseteq X$ be a Borel set. Then the following are equivalent:
    \begin{enumerate}
        \item\label{itm:thm:pp_cap:1} $A$ is pluripolar;
        \item\label{itm:thm:pp_cap:2} $\Capa_{\theta}(A)=0$.
    \end{enumerate}
\end{theorem}
See \cite[Theorem~4.40]{GZ17} for the local result, and \cite[Theorem~10.1.4 and Proposition~16.9.9]{FPT} for the local and global capacity comparisons.

\begin{example}\label{ex:psh_qc}
    A quasi-plurisubharmonic function on an open subset $U\subseteq X$ is always quasi-continuous. See \cite[Theorem~4.20]{GZ17} and \cite[Theorem~9.4.2 and Proposition~16.9.9]{FPT} for the proof.

    More generally, if $\varphi,\psi$ are two quasi-plurisubharmonic functions on $U$, then the following function
    \[
    f(x)\coloneqq \begin{cases}
        \varphi(x)-\psi(x),&\textup{if }(\varphi\lor \psi)(x)\neq -\infty;\\
        \infty,&\textup{otherwise}
    \end{cases}
    \]
    is quasi-continuous.\footnote{In the literature, people usually say carelessly that $\varphi-\psi$ is quasi-continuous.}
\end{example}

\begin{definition}\label{def:conv_cap}
Let $U$ be an open subset of $X$.
    Let $(f_i)_{i\in I}$ be a net of Borel measurable functions $f_i\colon U\rightarrow [-\infty,\infty]$, and $f\colon U\rightarrow [-\infty,\infty]$ be a Borel measurable function. We say $(f_i)_{i\in I}$ converges to $f$ \emph{in capacity}\index{convergence in capacity} if for any $\delta>0$, we have
    \[
        \lim_{i\in I} \Capa_{\theta}\Bigl( \{f_i>f+\delta\} \Bigr)=0,\quad \lim_{i\in I} \Capa_{\theta}\Bigl( \{f_i<f-\delta\} \Bigr)=0\footnotemark.
    \]
    \footnotetext{It is very tempting to write $\lim_{i\in I} \Capa_{\theta}\left( \{|f_i-f|>\delta\} \right)=0$, as in \cite[Definition~4.23]{GZ17} for example. But the set where $f_i-f$ is not defined is not a pluripolar set in general. Hence this abuse of notation is not acceptable. The two gap sets are also used in \cite[Definition~9.4.5]{FPT}.}
    We sometimes write $f_i\xrightarrow{\mathrm{C}} f$\index{$\xrightarrow{\mathrm{C}}$}. We say that the convergence is \emph{locally in capacity} if the two limits hold after intersecting the gap sets with each compact subset of $U$. When $U=X$, the two notions agree.
\end{definition}
Note that $f$ is not uniquely determined by the net $(f_i)_{i\in I}$. Thanks to \cref{thm:comp_capa}, the notion of convergence in capacity is independent of the choice of $\theta$.

\begin{proposition}\label{prop:mono_psh_cap}
Let $U$ be an open subset of $X$.
    Let $(\varphi_i)_{i\in I}$ be a net in $\PSH(U,\theta)$ and $\varphi\in \PSH(U,\theta)$. Assume one of the following conditions holds:
    \begin{enumerate}
        \item\label{itm:prop:mono_psh_cap:1} $(\varphi_i)_{i\in I}$ is decreasing and $\varphi$ is the limit of the net;
        \item\label{itm:prop:mono_psh_cap:2} $(\varphi_i)_{i\in I}$ is increasing and converges almost everywhere to $\varphi$.
    \end{enumerate}
    Then $(\varphi_i)_{i\in I}$ converges to $\varphi$ locally in capacity on $U$.
\end{proposition}
\begin{proof}
    First suppose that the net is a sequence. On each coordinate ball relatively compact in $U$, add a smooth local potential of $\theta$ to reduce to the psh case. The sequence converges on compact subsets in the local Bedford--Taylor capacity by \cite[Theorem~7.2]{BT82}. See also \cite[Proposition~15.6.7]{FPT}. To pass to the global capacity on each compact core of the ball, first use a Kähler form $\omega$: its capacity tests lie between $-1$ and $0$, so after adding a bounded local potential of $\omega$ they can be rescaled and translated to admissible Bedford--Taylor tests, with a common scale. Thus local capacity convergence on the core implies convergence there in $\Capa_\omega$, and then in $\Capa_\theta$ by \cref{thm:comp_capa}. Cover a given compact subset of $U$ by finitely many such cores.

    For a decreasing net, for each rational $c$ use second countability to select a countable subcover of $\{\varphi<c\}=\bigcup_i\{\varphi_i<c\}$ from the sets on the right. The selected indices for all $c$ are dominated successively by an increasing sequence of indices, whose potentials still decrease to $\varphi$. For an increasing net, $\varphi_i\leq\varphi$ everywhere by \cref{prop:pshfuncdetdense}. Choquet's lemma \cref{prop:Choquet}, applied after adding smooth local potentials on a countable cover, gives a countable subfamily with the same regularized supremum $\varphi$. Taking successive upper bounds of its indices gives a sequence of potentials increasing to $\varphi$ almost everywhere by \cref{prop:supsupstardiff}. In either case, each gap set of the net is contained, on the tail after a selected index, in the corresponding gap set of that sequence. The sequence case proves the result.
\end{proof}

\begin{theorem}\label{thm:loc_conv_cap}
    Let $U$ be an open subset of $X$. Let $(\chi_i\colon U\rightarrow [0,\infty))_i$ be a sequence of quasi-continuous functions. Assume that the sequence is uniformly bounded and converges in capacity to $\chi\colon U\rightarrow [0,\infty)$. For each $j=1,\ldots,p$, let $(\varphi_j^i)_i$ be uniformly bounded sequences in $\PSH(U,\theta_j)$ converging in capacity to bounded $\varphi_j\in \PSH(U,\theta_j)$. Then
    \[
        \chi_i\,\bigwedge_{j=1}^p \theta_{j,\varphi_j^i}\rightharpoonup  \chi\,\bigwedge_{j=1}^p \theta_{j,\varphi_j}
    \]
    weakly as currents on $U$.
\end{theorem}
See \cite[Theorem~4.26]{GZ17} and \cite[Theorem~9.4.8]{FPT}\footnote{In the displayed formula in the statement of \cite[Theorem~4.26]{GZ17}, $1\leq j\leq p$ should be $1\leq i\leq p$ on both sides.}.

\section{Properties of non-pluripolar products}\label{sec:nppprop}
Let $X$ be a connected compact Kähler manifold of dimension $n$ and $\theta,\theta_1,\ldots,\theta_n$ be closed real smooth $(1,1)$-forms on $X$.

We write
\begin{equation}\label{eq:PSHpos}
\PSH(X,\theta)_{>0}=\left\{ \varphi\in \PSH(X,\theta): \int_X\theta_{\varphi}^n>0 \right\}.
\end{equation}
The non-pluripolar product $\theta_{\varphi}^n$ is well-defined thanks to \cref{prop:nppwelldef}.

\begin{remark}\label{rmk:npp:L275}
    Suppose that $X$ is a connected complex manifold of dimension $0$, namely, $X$ is a single point. In this case, by definition, the non-pluripolar product $\theta_{\varphi}^n$ is given by the current of integration at the unique point. So $\PSH(X,\theta)_{>0}=\PSH(X,\theta)\cong \mathbb{R}$ in this case and $\int_X \theta_{\varphi}^n=1$ for all $\varphi\in \PSH(X,\theta)$.
\end{remark}

Recall the following basic result:
\begin{proposition}\label{prop:nppposimpbig}
    Assume that $\PSH(X,\theta)_{>0}$ is non-empty. Then the cohomology class $\{\theta\}$ is big.
\end{proposition}
See \cite[Proposition~1.22]{BEGZ10} and \cite[Theorem~17.4.20]{FPT}.

We recall a few basic facts about the non-pluripolar masses.
\begin{proposition}\label{prop:nppmassinv}
    Let $\pi\colon Y\rightarrow X$ be a proper bimeromorphic morphism from a K\"ahler manifold $Y$ and $T_1,\ldots,T_n$ be closed positive $(1,1)$-currents on $X$. Then
    \[
        \int_Y (\pi^*T_1)\wedge \cdots \wedge (\pi^*T_n)=\int_X T_1\wedge \cdots \wedge T_n.
    \]
    Similarly, if $S_1,\ldots,S_n$ are closed positive $(1,1)$-currents on $Y$, then
    \[
        \int_Y S_1\wedge \cdots \wedge S_n=\int_X \left(\pi_*S_1\right)\wedge \cdots \wedge \left(\pi_*S_n\right).
    \]
\end{proposition}
Here $\pi^*T_i$ is defined as follows: Write $T_i=\theta_i+\ddc \varphi_i$ for some closed smooth real $(1,1)$-form $\theta_i$ on $X$ and $\varphi_i\in \PSH(X,\theta_i)$, then $\pi^*T_i=\pi^*\theta_i+\ddc \pi^*\varphi_i$. The pull-back $\pi^*T_i$ is clearly independent of the choices of $\theta_i$ and $\varphi_i$.
\begin{proof}
    The assertions follow from \itemref{itm:prop:npp1:1} on the Zariski open set where $\pi$ is an isomorphism. The exceptional locus is pluripolar, and the relevant non-pluripolar products put no mass on pluripolar sets by \itemref{itm:prop:npp1:4}.
\end{proof}

\begin{theorem}\label{thm:logconc}
    Let $\varphi_0,\varphi_1\in \PSH(X,\theta)$. Then the map
    \[
    [0,1] \ni t\mapsto \log\int_X \theta_{t\varphi_1+(1-t)\varphi_0}^n
    \]
    is concave. 
    
    More generally, given closed positive $(1,1)$-currents $T_1,\ldots,T_n$ on $X$, we have
    \[
        \int_X T_1\wedge \cdots \wedge T_n\geq \prod_{i=1}^n \left( \int_X T_i^n \right)^{1/n}.
    \]
\end{theorem}
See \cite{DDNL19log} and \cite[Theorem~17.4.23]{FPT} for the proof.
\begin{remark}\label{rmk:linearcompsh}
    Here and in the sequel, when we write expressions like $t\varphi+(1-t)\psi$ for $\varphi,\psi\in \QPSH(X)$, we will follow the convention that when $t=0$, the value is $\psi$ and when $t=1$, the value is $\varphi$.
\end{remark}

We shall write
\begin{equation}\label{eq:Vtheta}
    V_{\theta}=\sup\left\{ \varphi\in \PSH(X,\theta):\varphi\leq 0 \right\}.
\end{equation}
Note that the upper semicontinuous regularization $V_{\theta}^*$ of $V_{\theta}$ is itself a candidate of \eqref{eq:Vtheta} if $\PSH(X,\theta)\neq \varnothing$ as follows from \cref{prop:pshfunction_closedseq}, therefore, $V_{\theta}=V_{\theta}^*\in \PSH(X,\theta)$ if $\PSH(X,\theta)\neq \varnothing$.
\index{$V_{\theta}$}

The function $V_{\theta}$ should be regarded as a canonical representative of the least singular potentials in $\PSH(X,\theta)$. We recall the following result:

\begin{lemma}\label{lma:Vthetalocbdd}
    Suppose that $\{\theta\}$ is big.
    The potential $V_{\theta}$ is locally bounded on $X\setminus \nK(\{\theta\})$.
\end{lemma}
Recall that $\nK(\{\theta\})$ is the non-Kähler locus introduced in \cref{def:nKlocus}.
\begin{proof}
    By definition of $\nK(\{\theta\})$, it suffices to show that if $\varphi\in \PSH(X,\theta)$ has analytic singularities and $\theta_{\varphi}$ is a Kähler current, then $V_{\theta}$ is locally bounded outside $\{\varphi=-\infty\}$. After replacing $\varphi$ by $\varphi-\sup_X\varphi$, the function $\varphi$ is a candidate in \eqref{eq:Vtheta}; hence $V_{\theta}\geq \varphi$ and $V_{\theta}\leq 0$. Since $\varphi$ is locally bounded outside its pole set, so is $V_{\theta}$.
\end{proof}

The non-pluripolar product has a lower semicontinuity property.
\begin{theorem}[Semicontinuity theorem]\label{thm:semicon}\index{theorem!semicontinuity theorem (non-pluripolar product)}
    Let $\varphi_j, \varphi_j^k\in \PSH(X,\theta_j)$ ($k\in \mathbb{N}$, $j=1,\ldots,n$).
    Let $\chi_k,\chi$ ($k\in \mathbb{N}$) be non-negative uniformly bounded quasi-continuous functions on $X$ such that $(\chi_k)_k$ converges to $\chi$ in capacity.
    Assume that for any $j=1,\ldots,n$, as $k\to\infty$, either $\varphi_j^k$ decreases to $\varphi_j\in \PSH(X,\theta_j)$ or increases to $\varphi_j\in \PSH(X,\theta_j)$ almost everywhere.
    Then
    \begin{equation}\label{eq:semicon1}
    \varliminf_{k\to\infty}\int_X\chi_k \,\theta_{1,\varphi_1^k}\wedge \cdots \wedge \theta_{n,\varphi_n^k} \geq \int_X \chi \,\theta_{1,\varphi_1}\wedge \cdots \wedge \theta_{n,\varphi_n}.
\end{equation}

If in addition,
\[
\lim_{k\to\infty}\int_X\theta_{1,\varphi_1^k}\wedge \cdots \wedge \theta_{n,\varphi_n^k} =\int_X \theta_{1,\varphi_1}\wedge \cdots \wedge \theta_{n,\varphi_n},
\]
then
\[
\chi_k \,\theta_{1,\varphi_1^k}\wedge \cdots \wedge \theta_{n,\varphi_n^k}\rightharpoonup \chi \,\theta_{1,\varphi_1}\wedge \cdots \wedge \theta_{n,\varphi_n}
\]
weakly as measures\footnote{Recall that the weak convergence of Radon measures is stronger than the weak convergence as currents in general. When the Radon measures have uniformly bounded total variation, they are equivalent.}. In particular, the limit in \eqref{eq:semicon1} exists and equality holds in \eqref{eq:semicon1}.
\end{theorem}
See \cite[Theorem~2.6]{DDNLsurv} and \cite[Corollary~15.6.8]{FPT} for the proof.

The non-pluripolar mass is a monotone quantity with respect to the singularity type.
\begin{theorem}[Monotonicity theorem]\label{thm:mono}\index{theorem!monotonicity theorem (non-pluripolar product)}
 Let $\varphi_j,\psi_j\in \PSH(X,\theta_j)$ for $j=1,\ldots,n$. Assume that $\varphi_j \succeq \psi_j$\footnote{See \cref{def:parorder} for the notation.} for every $j$, then
 \[
 \int_{X} \theta_{1,\varphi_1}\wedge \cdots \wedge \theta_{n,\varphi_n} \geq \int_{X}  \theta_{1,\psi_1}\wedge \cdots \wedge  \theta_{n,\psi_n}.
 \]
 In particular, if $\varphi,\psi\in \PSH(X,\theta)$ with $\varphi\succeq \psi$, then
 \[
 \int_X \theta_{\varphi}^n\geq  \int_X \theta_{\psi}^n.
 \]
 In other words, the non-pluripolar mass increases as the potentials become less singular.
\end{theorem}
See \cite[Theorem~1.1]{DDNL18mono} and \cite[Theorem~15.7.4 and Corollary~15.7.5]{FPT}. We will prove a vast extension of this theorem in \cref{prop:mono2}.

By this theorem, the non-pluripolar mass $\int_X \theta_{\varphi}^n$ can be used as a rough measure of the singularities of $\varphi\in \PSH(X,\theta)$. In \cref{sec:Penv}, we shall refine this measure by defining the notion of $P$-envelope.

As a corollary, we obtain that
\begin{corollary}\label{cor:incseqnppcont}
Fix a directed set $I$.
    For each $j=1,\ldots,n$, take an increasing net $(\varphi_j^i)_{i\in I}$ in $\PSH(X,\theta_j)$, uniformly bounded from above. Set
    \[
    \varphi_j\coloneqq \uscsup[i\in I] \varphi_j^i,\quad j=1,\ldots,n.
    \]
    Then
    \begin{equation}\label{eq:increseqnppcont}
    \lim_{i\in I}\int_X \theta_{1,\varphi_1^i}\wedge \cdots \wedge \theta_{n,\varphi_n^i} = \int_X \theta_{1,\varphi_1}\wedge \cdots \wedge \theta_{n,\varphi_n}.
\end{equation}
In particular, if the nets are indeed sequences, we have
\[
\theta_{1,\varphi_1^i}\wedge \cdots \wedge \theta_{n,\varphi_n^i} \rightharpoonup \theta_{1,\varphi_1}\wedge \cdots \wedge \theta_{n,\varphi_n}.
\]
\end{corollary}
\begin{proof}
The final assertion follows from \eqref{eq:increseqnppcont} together with \cref{thm:semicon}. It remains to establish \eqref{eq:increseqnppcont}.

We may assume that $I$ is infinite as there is nothing to prove otherwise. Thanks to \cref{thm:mono}, we already know the $\leq$ inequality in \eqref{eq:increseqnppcont}. We prove the reverse inequality.
When $I\cong\mathbb{Z}_{>0}$ as directed sets, the reverse inequality follows from \cref{thm:semicon}. In general, by Choquet's lemma \cref{prop:Choquet}, we can find a countable infinite subset $R\subseteq I$ such that
\[
\uscsup[r\in R]\varphi^r_j=\uscsup[i\in I]\varphi^i_j
\]
for all $j=1,\ldots,n$.
We fix a bijection $R\cong \mathbb{Z}_{>0}$. For any $j=1,\ldots,n$, we will then denote elements $\varphi_j^r$ ($r\in R$) by $\varphi_j^1,\varphi_j^2,\dots$. We shall write
\[
\psi_j^a=\varphi_j^1\lor \cdots \lor \varphi_j^a
\]
for each $a\in \mathbb{Z}_{>0}$.

For every $a$ we may choose $i_a\in I$ dominating the finitely many indices corresponding to $\varphi_j^1,\ldots,\varphi_j^a$ for all $j$. It then follows from the monotonicity theorem that
\[
\lim_{i\in I}\int_X \theta_{1,\varphi_1^i}\wedge \cdots \wedge \theta_{n,\varphi_n^i}\geq \lim_{a\to \infty}\int_X \theta_{1,\psi_1^a}\wedge \cdots \wedge \theta_{n,\psi_n^a}.
\]
From the special case mentioned above, we know that the right-hand side is exactly the right-hand side of \eqref{eq:increseqnppcont}, so we conclude.
\end{proof}

We prove an interesting inequality about the Monge--Ampère measure of the maximum of two potentials.
\begin{lemma}\label{lma:lorMAineqe}
    Let $\varphi,\psi\in \PSH(X,\theta)$. Then
    \begin{equation}\label{eq:thetavarphilorpsi1}
    \theta_{\varphi\lor \psi}^n\geq \mathds{1}_{\{\varphi\geq \psi\}}\,\theta_{\varphi}^n+\mathds{1}_{\{\varphi<\psi\}}\,\theta_{\psi}^n.
    \end{equation}
    In particular, if $\varphi\leq \psi$, then
    \[
    \mathds{1}_{\{\varphi=\psi\}}\,\theta_{\varphi}^n\leq \mathds{1}_{\{\varphi=\psi\}}\,\theta_{\psi}^n.
    \]
\end{lemma}
At a first sight, \eqref{eq:thetavarphilorpsi1} might seem trivial, and it is if we replace $\{\varphi\geq \psi\}$ by $\{\varphi>\psi\}$. The difficulty really lies on understanding the contact set $\{\varphi=\psi\}$.

\begin{proof}
    For each $k\in \mathbb{N}$, we set
    \[
    \psi_k=\psi\lor (V_{\theta}-k),\quad \varphi_k=\varphi\lor (V_{\theta}-k).
    \]
    For each $t>0$ and each $k\in \mathbb{N}$, we have
    \[
    \begin{split}
    \theta_{\psi_k\lor (\varphi_k+t)}^n\geq & \mathds{1}_{\{\psi_k>\varphi_k+t\}}\,\theta_{\psi_k\lor (\varphi_k+t)}^n+\mathds{1}_{\{\psi_k<\varphi_k+t\}}\,\theta_{\psi_k\lor (\varphi _k+t)}^n\\
    =&\mathds{1}_{\{\psi_k>\varphi_k+t\}}\,\theta_{\psi_k}^n+\mathds{1}_{\{\psi_k<\varphi_k+t\}}\,\theta_{\varphi_k}^n,
    \end{split}
    \]
    where the equality follows from the plurifine locality of non-pluripolar products \itemref{itm:prop:npp1:1}. We observe that for any $t>0$, $\psi_k\lor (\varphi_k+t)\sim \psi_k\lor \varphi_k$. Therefore, the monotonicity theorem \cref{thm:mono} guarantees that
    \[
        \int_X \theta_{\psi_k\lor (\varphi_k+t)}^n=\int_X\theta_{\psi_k\lor \varphi_k}^n.
    \]
    As $t\to 0+$,
    \[
    \theta_{\psi_k\lor (\varphi_k+t)}^n\rightharpoonup \theta_{\psi_k\lor \varphi_k}^n
    \]
    by \cref{thm:semicon}.

    Letting $t\to 0+$, we conclude that
    \[
    \theta_{\psi_k\lor \varphi_k}^n\geq \mathds{1}_{\{\psi_k>\varphi_k\}}\,\theta_{\psi_k}^n+\mathds{1}_{\{\psi_k\leq \varphi_k\}}\,\theta_{\varphi_k}^n.
    \]
    In particular, multiplying both sides by $\mathds{1}_{\{\min\{\varphi,\psi\}>V_{\theta}-k\}}$ and applying \itemref{itm:prop:npp1:1} again, we find that
    \[
    \begin{split}
    \mathds{1}_{\{\min\{\varphi,\psi\}>V_{\theta}-k\}}\,\theta_{\psi\lor \varphi}^n\geq &\mathds{1}_{\{\min\{\varphi,\psi\}>V_{\theta}-k\}\cap \{\psi\leq \varphi\}}\,\theta_{\varphi}^n \\
    &+\mathds{1}_{\{\min\{\varphi,\psi\}>V_{\theta}-k\}\cap \{\psi>\varphi\}}\,\theta_{\psi}^n.
    \end{split}
    \]
    Letting $k\to\infty$, we conclude \eqref{eq:thetavarphilorpsi1}.
\end{proof}
\begin{corollary}\label{cor:intcontlimsup}
    Let $(\varphi_i)_{i\in I}$ be a net in $\PSH(X,\theta)$, and $\varphi\in \PSH(X,\theta)$ so that $\varphi_i\to \varphi$ in $L^1$. Assume that $\varphi_i\leq 0$ for all $i\in I$.
    Then
    \begin{equation}\label{eq:contact_int_usc}
    \varlimsup_{i\in I} \int_{\{\varphi_i=0\}}\theta_{\varphi_i}^n\leq \int_{\{\varphi=0\}} \theta_{\varphi}^n.
    \end{equation}
\end{corollary}
\begin{proof}
    Observe that $\varphi\leq 0$ by \cref{prop:L1compa}.

    \textbf{Step~1}. We first assume that $I\cong \mathbb{Z}_{>0}$ as directed sets. Namely, $(\varphi_i)_{i>0}$ is a sequence.

    For each $k>0$, let
    \[
    \psi_k=\uscsup[j\geq k] \varphi_j.
    \]
    Then $(\psi_k)_{k}$ is a decreasing sequence in $\PSH(X,\theta)$ with limit $\varphi$, see \cref{cor:L1limipp} for example.

    Thanks to \cref{lma:lorMAineqe}, for each $k>0$, we have
    \[
    \mathds{1}_{\{\varphi_k=\psi_k\}}\,\theta_{\varphi_k}^n\leq \mathds{1}_{\{\varphi_k=\psi_k\}}\,\theta_{\psi_k}^n.
    \]
    Multiplying both sides with $\mathds{1}_{\{\varphi_k=0\}}$, we find that
    \[
    \mathds{1}_{\{\varphi_k=0\}}\,\theta_{\varphi_k}^n\leq \mathds{1}_{\{\varphi_k=0\}}\,\theta_{\psi_k}^n.
    \]
    Take $b,C>0$. For each $k>0$, we have
    \[
    \begin{split}
    \int_{\{\varphi_k=0\}} \theta_{\varphi_k}^n\leq \int_{\{\varphi_k=0\}} \theta_{\psi_k}^n\leq \int_{\{\psi_k=0\}} \theta_{\psi_k}^n\\=\int_{\{\psi_k=0\}} \theta_{\psi_k\lor (V_{\theta}-C)}^n\leq \int_X \mathrm{e}^{b\psi_k}\,\theta_{\psi_k\lor (V_{\theta}-C)}^n,
    \end{split}
    \]
    where the equality part follows from \itemref{itm:prop:npp1:1}.

    Letting $k\to\infty$ with the help of \cref{thm:semicon}, we conclude that
    \[
    \varlimsup_{k\to\infty}\int_{\{\varphi_k=0\}} \theta_{\varphi_k}^n\leq \int_X \mathrm{e}^{b\varphi}\,\theta_{\varphi\lor (V_{\theta}-C)}^n.
    \]
    Letting $b\to\infty$ and then $C\to\infty$, we conclude \eqref{eq:contact_int_usc}.

    \textbf{Step~2}. We now handle the general net case.

    Fix a volume form on $X$, so that we can talk about the $L^1$-norm on $X$ in the sequel.

    Suppose that \eqref{eq:contact_int_usc} fails. Then there is $\delta>0$ such that, for every $i_0\in I$, one can find $i\geq i_0$ satisfying
    \[
        \int_{\{\varphi_i=0\}}\theta_{\varphi_i}^n>\int_{\{\varphi=0\}}\theta_{\varphi}^n+\delta.
    \]
    Since $\varphi_i\to\varphi$ in $L^1$, for every $k>0$ there is $a_k\in I$ such that
    \[
        \|\varphi_i-\varphi\|_{L^1}<k^{-1},\quad i\geq a_k.
    \]
    We can then recursively choose an increasing sequence $(i_k)_{k>0}$ in $I$ with $i_k\geq a_k$ and
    \[
        \int_{\{\varphi_{i_k}=0\}}\theta_{\varphi_{i_k}}^n>\int_{\{\varphi=0\}}\theta_{\varphi}^n+\delta.
    \]
    This sequence converges to $\varphi$ in $L^1$, contradicting Step~1.
\end{proof}

Next we introduce an envelope construction, which will be repeatedly used in the sequel.
\begin{definition}\label{def:Pthetaf}
    Given a function $f\colon X\rightarrow [-\infty,\infty]$, we define $P_{\theta}(f)$ as follows:
    \begin{equation}\label{eq:Pthetafdef}
    P_{\theta}(f)\coloneqq \uscsup \left\{ \varphi\in \PSH(X,\theta): \varphi\leq f\textup{ quasi-everywhere} \right\}\footnotemark.
    \end{equation}
    \index{$P_{\theta}(f)$}
    \footnotetext{In the original article \cite{DDNLsurv}, the authors required that $\varphi\leq f$ everywhere. But in the proof of their Theorem~2.7 (\cref{thm:concentration} below), they actually relied on the current definition. In the most relevant case, if $f$ is $\mathcal{F}$-upper semicontinuous (for example, if $f$ is upper semicontinuous), then these two definitions coincide. This is due to \cref{lma:pshfunfinitelocuspfdense} and \cref{prop:pshfunFcont}. See also \cite[Definition~16.1.3 and Proposition~16.1.4]{FPT} for the obstacle inequality and its exceptional sets.}
\end{definition}
The function $P_{\theta}(f)$ is either constantly $\pm\infty$ or lies in $\PSH(X,\theta)$. Moreover, given another function $g\colon X\rightarrow [-\infty,\infty]$, equal to $f$ quasi-everywhere, we have $P_{\theta}(f)=P_{\theta}(g)$. In particular, it makes sense to talk about $P_{\theta}(f)$ even if $f$ is only defined outside a pluripolar set.

We also observe that
\begin{equation}\label{eq:Pthetafleqfqe}
P_{\theta}(f)\leq f\quad\textup{quasi-everywhere}.
\end{equation}
This is a consequence of \cref{prop:Choquet} and \cref{prop:supsupstardiff}.

\begin{lemma}\label{lma:BT_freeenvelope}
    Fix a closed smooth real $(1,1)$-form $\theta$ on $X$ representing a big cohomology class.
    Let $f\colon X\rightarrow \mathbb{R}$ be a lower semicontinuous function.
    Assume that $P_{\theta}(f)\not\equiv \infty$. Then $P_{\theta}(f)$ is locally bounded outside $\nK(\{\theta\})$. Set $U\coloneqq \{P_{\theta}(f)<f\}\setminus \nK(\{\theta\})$, then
    \begin{equation}\label{eq:bal_ddcn0}
        \Bigl(\theta|_U+\ddc P_{\theta}(f)|_U\Bigr)^n=0.
    \end{equation}
\end{lemma}
\begin{proof}
    First observe that $V_{\theta}+\inf_X f$ is a candidate in the definition \eqref{eq:Pthetafdef}. So
    \[
    P_{\theta}(f)\geq V_{\theta}+\inf_X f,
    \]
    and in particular $P_{\theta}(f)\not\equiv -\infty$.

    In particular, $P_{\theta}(f)\in \PSH(X,\theta)$. It is therefore bounded from above. By \cref{lma:Vthetalocbdd}, it is locally bounded away from $\nK(\{\theta\})$.

    At an isolated point $x$ of $X$, the envelope equals $f(x)$, since the value of a candidate at $x$ can be raised to $f(x)$. Thus every $x\in U$ lies on a positive-dimensional component.

    Fix $x\in U$, put $a=P_{\theta}(f)(x)$, and choose $a<t<f(x)$. Let $\delta=(t-a)/2$. Choose concentric coordinate balls $B\Subset B'\Subset U$ centered at $x$, with a smooth local potential $g$ of $\theta$ on $B'$. Shrink $B$ so that on $\overline{B}$,
    \begin{equation}\label{eq:PthetafBp_temp1}
        P_{\theta}(f)+g<a+g(x)+\delta,\qquad g>g(x)-\delta,\qquad f>t.
    \end{equation}
    Let $\varphi$ be the $\theta$-balayage of $P_{\theta}(f)$ on $B$ from \cref{prop:bal_general}. Then $\varphi+g$ is psh on $B'$ and bounded by $a+g(x)+\delta$ on $\partial B$. The maximum principle gives the same bound on $B$. Subtracting $g$ gives
    \[
        \varphi|_B<t<f|_B.
    \]
    Since $\varphi=P_{\theta}(f)$ outside $B$, we have $\varphi\leq f$ quasi-everywhere. Thus $\varphi$ is a candidate for $P_{\theta}(f)$, and $\varphi\geq P_{\theta}(f)$ by balayage, so equality holds. This proves \eqref{eq:bal_ddcn0} on $B$, hence on $U$.
\end{proof}

\begin{theorem}\label{thm:concentration}
    Let $f\colon X\rightarrow [-\infty,\infty]$ be a quasi-continuous function.

    Assume that $P_{\theta}(f)\not\equiv \pm\infty$. Then $P_{\theta}(f)\in \PSH(X,\theta)$ and
    \begin{equation}\label{eq:Pthetafconc}
    \int_{\{P_{\theta}(f)<f\}} \theta_{P_{\theta}(f)}^n=0.
    \end{equation}
\end{theorem}

Thanks to \eqref{eq:Pthetafleqfqe}, we can rewrite \eqref{eq:Pthetafconc} as
\[
\int_{\{P_{\theta}(f)\neq f\}} \theta_{P_{\theta}(f)}^n=0.
\]
We shall say that $\theta_{P_{\theta}(f)}^n$ is carried by the contact set $\{P_{\theta}(f)=f\}$.
\footnote{It is tempting to say that $\theta_{P_{\theta}(f)}^n$ is supported on the contact set $\{P_{\theta}(f)=f\}$, as people are already doing in the literature. It should be mentioned that this is an abuse of the language, since the support of $\theta_{P_{\theta}(f)}^n$ (as a closed subset of $X$) could be much larger in general. One could probably introduce a notion of plurifine support, similar to what Fuglede did in the classical potential theory \cite{Fug72}. See also \cite[Theorem~16.2.5]{FPT} for the contact-set concentration statement.}

\begin{proof}
    \textbf{Step~1}. We first reduce to the case where $f$ is bounded.

    \textbf{Step~1.1}. Reduce to the case where $f\leq 0$.

    Take $C\in \mathbb{R}$ so that $P_{\theta}(f)\leq C$. Then
    \[
    P_{\theta}(f)=P_{\theta}\Bigl( \min\{f,C\} \Bigr).
    \]
    Thus we can replace $f$ by $\min \{f,C\}$. Replacing $f$ by $f-\sup_X f$, we reduce easily to the case where $f\leq 0$.

    \textbf{Step~1.2}. Reduce to the case where $f$ is bounded. Assume that in this case, the theorem is known. We prove \eqref{eq:Pthetafconc} for $f\leq 0$.

    For each $j>0$, we have
    \begin{equation}\label{eq:Pthetafjint0_temp1}
    \int_{\{P_{\theta}(f_j)<f_j\}} \theta_{P_{\theta}(f_j)}^n=0,
    \end{equation}
    where $f_j=f\lor (-j)$.

    Now fix $C>0$ and two open sets $G'\Subset G\Subset X\setminus \nK(\{\theta\})$. Fix a smooth function $\chi\colon X\rightarrow [0,1]$ so that $\chi|_{G'}\equiv 1$ and $\chi$ is supported in $G$.

    Define
    \[
    U_C\coloneqq G \cap \left\{P_{\theta}(f)>V_{\theta}-C \right\},\quad U'_C\coloneqq G' \cap \left\{P_{\theta}(f)>V_{\theta}-C \right\}.
    \]
    Then $U_C$ is $\mathcal{F}$-open. It follows from \itemref{itm:prop:npp1:1} that for any $j>0$,
    \[
    \mathds{1}_{U_C}\, \theta_{P_{\theta}(f_j)\lor (V_{\theta}-C)}^n=\mathds{1}_{U_C}\, \theta_{P_{\theta}(f_j)}^n.
    \]
    In particular, thanks to \eqref{eq:Pthetafjint0_temp1},
    \begin{equation}\label{eq:intU1mexpPthetaf_temp1}
    \int_{U_C} \left( 1-\mathrm{e}^{P_{\theta}(f_j)-f_j} \right)\,\theta_{P_{\theta}(f_j)\lor (V_{\theta}-C)}^n=0.
    \end{equation}
    Note that the $P_{\theta}(f_j)\lor (V_{\theta}-C)$'s for various $j$ are uniformly bounded from below on $U_C$, thanks to \cref{lma:Vthetalocbdd}.

    Next we claim that $(P_{\theta}(f_j))_j$ is decreasing and converges to $P_{\theta}(f)$.

    The sequence is decreasing; let us compute its limit.
    In fact, $P_{\theta}(f_j)\leq f_j$ quasi-everywhere. It follows that $\inf_j P_{\theta}(f_j)\leq f$ quasi-everywhere and hence
    \[
    \inf_{j>0} P_{\theta}(f_j)\leq P_{\theta}(f).
    \]
    The reverse inequality is trivial, and our assertion follows.

    Since $P_{\theta}(f)\not\equiv -\infty$, we know that the set $\{f=-\infty\}$ is pluripolar. It follows that $P_{\theta}(f_j)-f_j$ converges to the following function $g\colon X\rightarrow [-\infty,\infty)$ in capacity:
    \[
    g(x)=\begin{cases}
        P_{\theta}(f)(x)-f(x),&\textup{if } f(x)\neq -\infty;\\
        -\infty,& \textup{otherwise}.
    \end{cases}
    \]
    Therefore,
    \[
    1-\mathrm{e}^{P_{\theta}(f_j)-f_j}\xlongrightarrow{\mathrm{C}} 1-\mathrm{e}^{g}.
    \]
    Next we claim that
     \begin{equation}\label{eq:intU1mexpPthetaf_temp2}
    \int_{U'_C} \left( 1-\mathrm{e}^{g} \right)\,\theta_{P_{\theta}(f)}^n=0.
    \end{equation}
    Since $g\leq 0$ quasi-everywhere, the left-hand side of \eqref{eq:intU1mexpPthetaf_temp2} is non-negative. It suffices to prove the $\leq$ direction.

    We wish to let $j\to\infty$ directly in \eqref{eq:intU1mexpPthetaf_temp1}, but since $U_C$ is not open in general, this cannot be done directly. In the sequel, we shall slightly enlarge $U_C$ to get an open set and then take the limit.

    Fix $\epsilon>0$. We can find an open subset $W\Subset X\setminus \nK(\{\theta\})$ so that
    \begin{equation}\label{eq:CapCinvthetaWmusmall_temp1}
    \Capa_{C^{-1}\theta}(W\setminus U_C)<\epsilon.
    \end{equation}
    For example, we can construct $W$ as follows. By \cref{ex:psh_qc}, we can find an open set $A\subseteq G$ so that $\Capa_{C^{-1}\theta}(A)<\epsilon$ and $P_{\theta}(f)-V_{\theta}$ is continuous on $G\setminus A$. Then $U_C\cap (G\setminus A)$ is relatively open in $G\setminus A$, so we may choose an open set $W_0\subseteq G$ such that
    \[
        W_0\cap (G\setminus A)=U_C\cap (G\setminus A).
    \]
    Taking $W=W_0\cup A$, we get $U_C\subseteq W$ and $W\setminus U_C\subseteq A$, hence \eqref{eq:CapCinvthetaWmusmall_temp1}.

    Thanks to \eqref{eq:CapCinvthetaWmusmall_temp1} and \eqref{eq:intU1mexpPthetaf_temp1}, we have
    \[
    \int_W \chi\left( 1-\mathrm{e}^{P_{\theta}(f_j)-f_j} \right)\,\theta_{P_{\theta}(f_j)\lor (V_{\theta}-C)}^n\leq C^n\epsilon.
    \]
    Letting $j\to\infty$ and applying \cref{thm:loc_conv_cap}, we find that
    \[
    \int_{U'_C} \left( 1-\mathrm{e}^{g} \right)\,\theta_{P_{\theta}(f)\lor (V_{\theta}-C)}^n\leq \int_W \chi\left( 1-\mathrm{e}^{g} \right)\,\theta_{P_{\theta}(f)\lor (V_{\theta}-C)}^n\leq C^n\epsilon.
    \]
    Since $\epsilon>0$ is arbitrary, \eqref{eq:intU1mexpPthetaf_temp2} follows.

    Now letting $C\to\infty$ in \eqref{eq:intU1mexpPthetaf_temp2}, we find
    \[
    \int_{G'}\left( 1-\mathrm{e}^{g} \right)\,\theta_{P_{\theta}(f)}^n=0.
    \]
    Since $G'\Subset X\setminus \nK(\{\theta\})$ is arbitrary, and $\nK(\{\theta\})\subseteq X$ is a proper analytic subset, we finally conclude
    \[
    \int_X \left( 1-\mathrm{e}^{g} \right)\,\theta_{P_{\theta}(f)}^n=0.
    \]
    In other words, $P_{\theta}(f)=f$ almost everywhere with respect to $\theta_{P_{\theta}(f)}^n$. This proves \eqref{eq:Pthetafconc}.

    \textbf{Step~2}. We now assume that $-C'\leq f\leq 0$ for some $C'>0$.

    For each $j>0$, take an open subset $V_j\subseteq X$ so that
    \begin{enumerate}
        \item\label{itm:chapter-npp:L662:1} $\Capa_{C'^{-1}\theta}(V_j)\leq 2^{-j-1}$, and
        \item\label{itm:chapter-npp:L662:2} $f|_{X\setminus V_j}$ is continuous.
    \end{enumerate}
    Without loss of generality, we may assume that the $V_j$'s are decreasing. Take a continuous function $f_j\colon X\rightarrow [-C',0]$ extending $f|_{X\setminus V_j}$. This is always possible thanks to Tietze's extension theorem. For each $j>0$, we let
    \[
    h_j\coloneqq \sup_{k\geq j} f_k.
    \]
    Then $h_j$ is lower semicontinuous and $h_j$ agrees with $f$ outside $V_j$.
    Set
    \[
    h=\inf_{j>0} h_j.
    \]
    Then the set $\{h\neq f\}$ is contained in the intersection of the $V_j$'s and hence is a pluripolar set, thanks to \cref{thm:pp_cap}. In particular, $P_{\theta}(f)=P_{\theta}(h)$.

    Now we can apply \cref{lma:BT_freeenvelope} to conclude that
    \[
    \int_{X\setminus \nK(\{\theta\})} \left(1-\mathrm{e}^{P_{\theta}(h_j)-h_j} \right) \theta_{P_{\theta}(h_j)}^n=0
    \]
    for each $j>0$.

    Fix two open subsets $G'\Subset G\Subset X\setminus \nK(\{\theta\})$. Note that $-C'\leq h_j\leq 0$. In particular,
    \[
    V_{\theta}-C'\leq P_{\theta}(h_j)\leq V_{\theta}.
    \]
    Hence,
    \[
    \begin{aligned}
        &\int_G \left| 1-\mathrm{e}^{P_{\theta}(h_j)-f}\right|\,\theta_{P_{\theta}(h_j)}^n\\
        =&\int_{G\cap V_j} \left|1-\mathrm{e}^{P_{\theta}(h_j)-f}\right|\,\theta_{P_{\theta}(h_j)}^n+\int_{G\setminus V_j} \left|1-\mathrm{e}^{P_{\theta}(h_j)-h_j}\right|\,\theta_{P_{\theta}(h_j)}^n\\
        =& \int_{G\cap V_j} \left|1-\mathrm{e}^{P_{\theta}(h_j)-f}\right|\,\theta_{P_{\theta}(h_j)}^n\\
        \leq & \int_{G\cap V_j} \left( 1+\mathrm{e}^{C'}\right)\,\theta_{P_{\theta}(h_j)}^n\\
        \leq & \left( 1+\mathrm{e}^{C'}\right) \cdot 2^{-j-1} C'^n.
    \end{aligned}
    \]
    Now we can apply the same arguments as in Step~1.2 to conclude that
    \[
    \int_{G'}\left(1-\mathrm{e}^{P_{\theta}(f)-f}\right)\,\theta_{P_{\theta}(f)}^n=0.
    \]
    Since $G'\Subset X\setminus \nK(\{\theta\})$, \eqref{eq:Pthetafconc} follows.
\end{proof}

The following lemma is striking in that we begin only with an upper bound of $\varphi$, but at the end of the day, we get a lower bound almost for free. This powerful method will be employed again and again in the whole book. A more general version will be proved in \cref{thm:pathoenvelope}, after introducing the terminology of $P$-preorder in \cref{chap:comp}.
\begin{lemma}\label{lma:pathoenvelope}
    Let $\varphi,\psi\in \PSH(X,\theta)$, $\varphi\preceq \psi$ and $\int_X \theta_{\varphi}^n>0$. Then for any
    \begin{equation}\label{eq:arangetemp}
    a\in \left( 1, \left(\frac{\int_X \theta_{\psi}^n}{\int_X \theta_{\psi}^n-\int_X \theta_{\varphi}^n}\right)^{1/n} \right)\footnotemark,
    \end{equation}
    there is $\eta\in \PSH(X,\theta)_{>0}$ such that
    \begin{equation}\label{eq:path_env}
        a^{-1}\eta+\left(1-a^{-1}\right)\psi\leq \varphi.
    \end{equation}
\end{lemma}
\footnotetext{The fraction in \eqref{eq:arangetemp} is understood as $\infty$ if either $\int_X \theta_{\psi}^n=\int_X \theta_{\varphi}^n$ or $n=0$. Thanks to the monotonicity theorem \cref{thm:mono}, the interval \eqref{eq:arangetemp} is non-empty.}

We write
\begin{equation}\label{eq:perversePtheta}
    \begin{aligned}
    P_{\theta}\Bigl(a\varphi+(1-a)\psi\Bigr)\coloneqq  &\uscsup\left\{\eta\in \PSH(X,\theta): a^{-1}\eta+\left(1-a^{-1}\right)\psi  \leq  \varphi \right\}\\
    \in & \PSH(X,\theta).
    \end{aligned}
\end{equation}
Note that if we regard $a\varphi+(1-a)\psi$ as a function defined outside the pluripolar set $\{\varphi=\psi=-\infty\}$ on $X$, then \eqref{eq:perversePtheta} coincides with the envelope in the sense of \eqref{eq:Pthetafdef}.

Observe that
\begin{equation}\label{eq:P_path_everywhere_ineq}
a^{-1}P_{\theta}\Bigl(a\varphi+(1-a)\psi \Bigr)+(1-a^{-1})\psi\leq \varphi.
\end{equation}
In fact, this equation holds outside a pluripolar set by \cref{prop:supsupstardiff}, hence it holds everywhere by \cref{prop:pshfuncdetdense}.

Finally, observe that
\[
P_{\theta}\Bigl(a\varphi+(1-a)\psi\Bigr)\preceq P_{\theta}\Bigl(a\varphi+(1-a)\varphi\Bigr)=\varphi.
\]
If $\varphi\leq \psi$, the same argument gives
\[
P_{\theta}\Bigl(a\varphi+(1-a)\psi\Bigr)\leq\varphi.
\]

\begin{proof}
    Without loss of generality, we may assume that $\varphi\leq \psi\leq 0$.

    \textbf{Step~1}. We show the existence of $\eta\in \PSH(X,\theta)$ satisfying \eqref{eq:path_env}.

    \textbf{Step~1.1}. We first consider the case where $\psi\sim \varphi$. In this case, we can take $\eta=\varphi-C$ enough $C>0$. But since by \eqref{eq:P_path_everywhere_ineq}, we also have $P_{\theta}(a\varphi+(1-a)\psi)\leq \psi$, we find that
    \[
        P_{\theta}\Bigl(a\varphi+(1-a)\psi\Bigr)\sim \varphi\sim \psi
    \]
    in this case.

    \textbf{Step~1.2}. We make a first reduction of the problem.

    Now for each $\gamma\in \PSH(X,\theta)$, we let
    \[
    P_{\theta}[\gamma]\coloneqq \uscsup\Bigl\{\eta\in \PSH(X,\theta): \eta\leq \min\{\gamma+C,0\}\textup{ for some }C>0 \Bigr\}\footnotemark.
    \]
    \footnotetext{This is precisely the $P$-envelope to be introduced later in \cref{def:Penv}.}

    Observe that due to \cref{cor:incseqnppcont} and \cref{thm:mono}, we have
    \begin{equation}\label{eq:intthetaP_pr_mass_temp1}
    \int_X \theta_{P_{\theta}[\gamma]}^n=\int_X \theta_{\gamma}^n.
    \end{equation}
    Further, $0\geq P_{\theta}[\gamma]\geq \gamma$ if $\gamma\leq 0$.

    We now define an increasing transfinite sequence $(\psi_{\alpha})_{\alpha}$ consisting of non-positive functions in $\PSH(X,\theta)$ as follows:
    \begin{enumerate}
        \item\label{itm:chapter-npp:L766:1} $\psi_0=\psi$;
        \item\label{itm:chapter-npp:L766:2} $\psi_{\alpha+1}=P_{\theta}[\psi_{\alpha}]$;
        \item\label{itm:chapter-npp:L766:3} if $\alpha$ is a limit ordinal, we define
        \[
        \psi_{\alpha}=\uscsup[\beta<\alpha] \psi_{\beta}.
        \]
    \end{enumerate}
    Observe that
    \[
    \int_X \theta_{\psi_{\alpha}}^n=\int_X \theta_{\psi}^n
    \]
    by \cref{cor:incseqnppcont} and \eqref{eq:intthetaP_pr_mass_temp1}. Since $\PSH(X,\theta)$ is indeed a set, the transfinite sequence must stabilize, say at some $\phi\in \PSH(X,\theta)$.
    Note that $0 \geq \phi\geq \psi$, and $P_{\theta}[\phi]=\phi$.

    Note that the interval \eqref{eq:arangetemp} remains the same if we replace $\psi$ by $\phi$. Take $a$ in this interval, it suffices to prove the existence of $\eta\in \PSH(X,\theta)$ so that
    \begin{equation}\label{eq:path_env_temp1}
    a^{-1}\eta+\left(1-a^{-1}\right)\phi\leq \varphi.
    \end{equation}

    Let us record the following observation for later use. Suppose that $\tau\in \PSH(X,\theta)$ and $\tau\preceq \phi$. Then
    \begin{equation}\label{eq:modelpot_sup_temp1}
        \sup_{\{\phi\neq -\infty\}} (\tau-\phi)=\sup_X \tau.
    \end{equation}
    Observe that $\phi\leq 0$, so on the set $\{\phi\neq -\infty\}$, we have
    \[
    \tau-\phi\geq \tau.
    \]
    So the $\geq$ direction in \eqref{eq:modelpot_sup_temp1} follows from \cref{cor:diffsupinfindeppluripolar}.
    Conversely, by assumption, we can find a constant $C>0$ so that
    \[
    \tau-\sup_X \tau\leq \min\{\phi+C,0 \}.
    \]
    It follows that $\tau-\sup_X \tau\leq P_{\theta}[\phi]=\phi$. Therefore, the reverse inequality follows, and \eqref{eq:modelpot_sup_temp1} follows.

    For each $k>0$, we introduce
    \[
    \varphi_k=\varphi\lor \left(\phi-k\right),\quad \eta_k\coloneqq P_{\theta}\Bigl( a \varphi_k+(1-a)\phi \Bigr).
    \]
    Since $\varphi_k\sim \phi$, using Step~1.1, we have $\eta_k\in \PSH(X,\theta)$ and $\eta_k\sim \phi$ as well.

    Note that $(\eta_k)_k$ is decreasing in $k$. Let
    \[
    \eta\coloneqq \inf_{k\in \mathbb{N}} \eta_k.
    \]
    Note that $\eta$ automatically satisfies \eqref{eq:path_env_temp1}.
    It remains to show that $\eta\not\equiv -\infty$.

    \textbf{Step~1.3}. We prove that $\eta\not\equiv -\infty$.

    Suppose, for contradiction, that $\sup_X \eta_k\to -\infty$. For each $k>0$, let
    \[
    \gamma_k\coloneqq a^{-1}\eta_k+\left(1-a^{-1}\right)\phi,\quad D_k\coloneqq \left\{\gamma_k=\varphi_k\right\}.
    \]
    We claim that
    \begin{equation}\label{eq:1Dkan_temp1}
        \theta_{\eta_k}^n\leq a^n\mathds{1}_{D_k}\, \theta_{\varphi_k}^n.
    \end{equation}

    Since $\gamma_k\leq \varphi_k$ with equality on $D_k$, \cref{lma:lorMAineqe} gives
    \[
    \mathds{1}_{D_k}\, \theta_{\gamma_k}^n\leq \mathds{1}_{D_k}\, \theta_{\varphi_k}^n.
    \]
    On the other hand,
    \[
    \theta_{\gamma_k}=a^{-1}\theta_{\eta_k}+\left(1-a^{-1}\right)\theta_{\phi},
    \]
    so positivity and multilinearity of the non-pluripolar product imply
    \[
        a^{-n}\theta_{\eta_k}^n\leq \theta_{\gamma_k}^n.
    \]
    Finally, it follows from \cref{thm:concentration} that
    \[
    \mathds{1}_{D_k}\,\theta_{\eta_k}^n=\theta_{\eta_k}^n.
    \]
    Putting these results together, we conclude \eqref{eq:1Dkan_temp1}.

    Fix $j>k>0$. Note that
    \begin{equation}\label{eq:intthetavarphijn_temp1}
    \int_{\{\varphi\leq \phi-k \}} \theta_{\varphi_j}^n=\int_X \theta_{\varphi_j}^n-\int_{\{\varphi>\phi-k \}} \theta_{\varphi_j}^n=\int_X \theta_{\phi}^n-\int_{\{\varphi>\phi-k \}} \theta_{\varphi}^n,
    \end{equation}
    where we applied \cref{thm:mono} and \itemref{itm:prop:npp1:1} in the second equality.

    Next, we compute
    \[
    \begin{aligned}
        \int_{\left\{ \eta_j\leq \phi-ak \right\}} \theta_{\eta_j}^n\leq & a^n  \int_{\left\{ \eta_j\leq \phi-ak \right\}} \mathds{1}_{D_j}\, \theta_{\varphi_j}^n & \textup{by \eqref{eq:1Dkan_temp1}}\\
        \leq & a^n \int_{\left\{ a\varphi_j+(1-a)\phi\leq \phi-ak \right\}\cap \left\{\phi\neq -\infty\right\} } \theta_{\varphi_j}^n\\
        = & a^n \int_{\left\{ \varphi_j\leq \phi-k \right\}} \theta_{\varphi_j}^n\\
        = & a^n \int_{\left\{ \varphi\leq \phi-k \right\}} \theta_{\varphi_j}^n\\
        \leq & a^n\left( \int_X \theta_{\phi}^n- \int_{\left\{\varphi>\phi-k \right\}} \theta_{\varphi}^n \right)&\textup{ by \eqref{eq:intthetavarphijn_temp1}}.
    \end{aligned}
    \]
    Next thanks to \eqref{eq:modelpot_sup_temp1},
    \[
    \sup_{ \left\{\phi\neq -\infty \right\} } \left(\eta_j-\phi \right)=\sup_X \eta_j\to -\infty.
    \]
    In particular, for a fixed $k$, if $j$ is large enough, we have
    \[
    \left\{\eta_j\leq \phi-ak \right\}=X.
    \]
    Therefore, for a fixed $k$, for any large enough $j$,
    \[
    \int_X \theta_{\phi}^n=\int_X \theta_{\eta_j}^n\leq a^n\left( \int_X \theta_{\phi}^n- \int_{\{\varphi>\phi-k\}} \theta_{\varphi}^n\right).
    \]
    Letting $k\to\infty$, we find
    \[
    \int_X \theta_{\phi}^n\leq a^n \left(\int_X \theta_{\phi}^n-\int_X \theta_{\varphi}^n \right).
    \]
    This is a contradiction with our choice of $a$.

    \textbf{Step~2}.
    Next we argue that $P_{\theta}(a\varphi+(1-a)\psi)\in \PSH(X,\theta)_{>0}$. Choose
    \[
        a'\in \left( a, \left(\frac{\int_X \theta_{\psi}^n}{\int_X \theta_{\psi}^n-\int_X \theta_{\varphi}^n}\right)^{1/n} \right).
    \]
    It follows from \eqref{eq:perversePtheta} that
    \begin{equation}\label{eq:Pthetalowerbdtemp1}
        P_{\theta}\Bigl(a\varphi+(1-a)\psi \Bigr)\geq \frac{a}{a'} P_{\theta}\Bigl(a'\varphi+(1-a')\psi \Bigr)+\frac{a'-a}{a'}\varphi.
    \end{equation}
    Therefore, by \cref{thm:mono}, we have
    \begin{equation}\label{eq:Pthetalowerbdtemp2}
        \int_X \theta_{P_{\theta}\left(a\varphi+(1-a)\psi\right)}^n\geq \frac{(a'-a)^n}{a'^n}\int_X \theta_{\varphi}^n>0.
    \end{equation}
    We conclude the proof.
\end{proof}

When $\varphi$ and $\psi$ have the same mass, we can say more:
\begin{corollary}\label{cor:pathoenvelopeeqmass}
    Let $\varphi,\psi\in \PSH(X,\theta)_{>0}$, $\varphi\preceq \psi$. Assume that $\int_X \theta_{\varphi}^n=\int_X \theta_{\psi}^n$. Then for any $\epsilon\in (0,1)$, there is $\eta\in \PSH(X,\theta)$ such that
    \begin{enumerate}
        \item\label{itm:cor:pathoenvelopeeqmass:1} $\int_X \theta_{\eta}^n=\int_X \theta_{\varphi}^n$;
        \item\label{itm:cor:pathoenvelopeeqmass:2} $\epsilon\eta+(1-\epsilon)\psi\leq \varphi$.
    \end{enumerate}
    In fact, we can take
    \[
        \eta=P_{\theta}\left(\epsilon^{-1}\varphi+(1-\epsilon^{-1})\psi\right).
    \]
\end{corollary}

\begin{proof}
    Fix $\epsilon \in (0,1)$. We define $\eta$ as in the statement of the corollary. By \cref{lma:pathoenvelope}, $\eta\in \PSH(X,\theta)_{>0}$. Observe that \ref{itm:cor:pathoenvelopeeqmass:2} is satisfied by \eqref{eq:P_path_everywhere_ineq}.

    Thanks to \eqref{eq:Pthetalowerbdtemp2}, for each $a'>\epsilon^{-1}$, we have
    \[
    \int_X \theta_{\eta}^n>\left(\frac{a'-\epsilon^{-1}}{a'}\right)^n\int_X \theta_{\varphi}^n.
    \]
    Letting $a'\to\infty$, we conclude that
    \[
    \int_X \theta_{\eta}^n\geq \int_X \theta_{\varphi}^n.
    \]
    On the other hand, since $\eta\preceq \psi$, using the monotonicity theorem \cref{thm:mono} we find that
    \[
    \int_X \theta_{\eta}^n\leq \int_X \theta_{\psi}^n= \int_X \theta_{\varphi}^n.
    \]
    Hence,
    \[
    \int_X \theta_{\eta}^n= \int_X \theta_{\varphi}^n.
    \]
    Therefore, \ref{itm:cor:pathoenvelopeeqmass:1} holds as well.
\end{proof}

Next we prove a domination principle.
\begin{theorem}[Domination principle]\label{thm:dom_prin_1}\index{theorem!domination principle}
    Let $\varphi,\psi\in \PSH(X,\theta)_{>0}$. Assume that
    \begin{equation}\label{eq:varphipsiPequiv}
    \exists\phi\in \PSH(X,\theta)\quad \textup{s.t.} \quad \varphi\preceq \phi, \quad \psi\preceq \phi,\quad \int_X \theta_{\varphi}^n=\int_X \theta_{\psi}^n=\int_X \theta_{\phi}^n,
    \end{equation}
    and
    \begin{equation}\label{eq:intvarphilpsitheta0}
    \int_{\{\varphi<\psi\}} \theta_{\varphi}^n=0.
    \end{equation}
    Then $\varphi\geq \psi$.
\end{theorem}
\begin{remark}\label{rmk:npp:L982}
    Using the terminologies to be introduced in \cref{chap:envelopes}, we can reformulate a special case of \cref{thm:dom_prin_1} as follows: Suppose that $\phi\in \PSH(X,\theta)_{>0}$ is a model potential and $\varphi,\psi\in \mathcal{E}(X,\theta;\phi)$. Assume that \eqref{eq:intvarphilpsitheta0} holds then $\varphi\geq \psi$.

    Equivalently, using the terminologies of \cref{chap:comp}, the assumption \eqref{eq:varphipsiPequiv} can be reformulated as $\varphi\sim_P \psi$.
\end{remark}

\begin{proof}
    Thanks to the monotonicity theorem \cref{thm:mono},
    \[
    \int_X \theta_{\varphi\lor\psi}^n=\int_X \theta_{\phi}^n.
    \]
    We may replace $\psi$ by $\varphi\lor\psi$ and assume that $\varphi\leq \psi$. Without loss of generality, we may also assume that $\psi\leq 0$.

    Now fix $a>1$. We define
    \[
    \eta_a\coloneqq P_{\theta}\Bigl( a\varphi+(1-a)\psi\Bigr).
    \]
    Note that $\eta_a\in \PSH(X,\theta)$ by \cref{cor:pathoenvelopeeqmass} and
    \[
    \int_X \theta_{\eta_a}^n=\int_X \theta_{\phi}^n.
    \]
    Furthermore, $\eta_a\leq P_{\theta}(a\varphi+(1-a)\varphi)=\varphi$.

    Define
    \[
    \gamma_a\coloneqq a^{-1}\eta_a+(1-a^{-1})\psi,\quad D_a\coloneqq \{\gamma_a=\varphi\}.
    \]
    Then $\gamma_a\leq \varphi$ with equality on $D_a$. Therefore,
    \[
    \begin{aligned}
    a^{-n}\,\theta_{\eta_a}^n\leq & a^{-n}\,\theta_{\eta_a}^n+\mathds{1}_{D_a} (1-a^{-1})^n \,\theta_{\psi}^n&\\
    =& \mathds{1}_{D_a} a^{-n}\,\theta_{\eta_a}^n+\mathds{1}_{D_a} (1-a^{-1})^n \,\theta_{\psi}^n & \textup{ by \cref{thm:concentration}}\\
    \leq & \mathds{1}_{D_a} \theta_{\gamma_a}^n &\\
    \leq & \mathds{1}_{D_a} \theta_{\varphi}^n & \textup{ by \cref{lma:lorMAineqe}}\\
    =& \mathds{1}_{D_a\cap \{\varphi=\psi\}}\, \theta_{\varphi}^n & \textup{ by \eqref{eq:intvarphilpsitheta0}}.
    \end{aligned}
    \]
    Note that on $D_a\cap \{\varphi=\psi\}$, we have $\eta_a=\varphi$.
    We deduce that
    \[
    \theta_{\eta_a}^n=\mathds{1}_{\{\eta_a=\varphi\}}\, \theta_{\eta_a}^n \leq \theta_{\varphi}^n,
    \]
    where the inequality follows from \cref{lma:lorMAineqe}.

    But the two ends have the same mass, and hence
    \[
    \theta_{\eta_a}^n=\theta_{\varphi}^n=\mathds{1}_{\{\eta_a=\varphi\}}\,\theta_{\varphi}^n.
    \]
    Therefore,
    \[
    \int_X \mathrm{e}^{\eta_a}\,\theta_{\varphi}^n=\int_X \mathrm{e}^{\varphi}\,\theta_{\varphi}^n>0.
    \]
    Note that $\eta_a$ is decreasing in $a$. The above equation shows that
    \[
    \eta\coloneqq\inf_{a>1} \eta_a\not\equiv -\infty.
    \]
    On the other hand, if $x\in X$ is such that $\varphi(x)<\psi(x)$, we then have
    \[
    \eta_a(x)\leq a\varphi(x)-(a-1)\psi(x)= \psi(x)+a\Bigl( \varphi(x)-\psi(x) \Bigr).
    \]
    Letting $a\to \infty$, we find that $\eta(x)=-\infty$. Therefore, $\{\varphi\neq \psi\}$ is pluripolar, and hence empty by \cref{prop:pshfuncdetdense}. Our assertion follows.
\end{proof}

\begin{lemma}\label{lma:kahcurrentposmass}
    For any $\varphi\in \PSH(X,\theta)_{>0}$, there is $\psi\in \PSH(X,\theta)$ such that
    \begin{enumerate}
        \item\label{itm:lma:kahcurrentposmass:1} $\theta_{\psi}$ is a Kähler current, and
        \item\label{itm:lma:kahcurrentposmass:2} $\psi\leq \varphi$.
    \end{enumerate}
    In particular, there is an increasing sequence $(\varphi_i)_{i>0}$ in $\PSH(X,\theta)$ converging almost everywhere to $\varphi$ such that $\theta_{\varphi_i}$ is a K\"ahler current for all $i\geq 1$.
\end{lemma}
\begin{proof}
    Using \cref{lma:pathoenvelope}, we can find $\epsilon\in (0,1)$ and $\gamma\in \PSH(X,\theta)$ such that
    \[
    \epsilon V_{\theta}+(1-\epsilon)\gamma\leq \varphi.
    \]
    We observe that the cohomology class $\{\theta\}$ is big as a consequence of \cref{prop:nppposimpbig}. Therefore, we can take $\eta\in \PSH(X,\theta)$ such that $\theta_{\eta}$ is a Kähler current and $\eta\leq 0$. Since $\eta\leq V_{\theta}$, we may take
    \[
    \psi\coloneqq \epsilon\eta+(1-\epsilon)\gamma.
    \]
    Then $\psi$ clearly satisfies \ref{itm:lma:kahcurrentposmass:1} and \ref{itm:lma:kahcurrentposmass:2}.

    For the latter claim, it suffices to take
    \[
    \varphi_i=\left(1-(i+1)^{-1}\right)\varphi +(i+1)^{-1}\psi
    \]
    for each $i>0$.
\end{proof}

\begin{lemma}\label{lma:existsecposmass}
    Let $L$ be a holomorphic line bundle on $X$ with $\theta\in c_1(L)$. Assume that $\varphi\in \PSH(X,\theta)_{>0}$. Then there exists $k_0>0$ such that for each $k\geq k_0$, we have
    \[
        \mathrm{H}^0\left(X,L^k\otimes \mathcal{I}(k\varphi)\right)\neq 0.
    \]
\end{lemma}

\begin{proof}
    By \cref{lma:kahcurrentposmass}, there is $\psi\in \PSH(X,\theta)$ with $\theta_{\psi}$ a Kähler current and $\psi\leq\varphi$. Since $\mathcal{I}(k\psi)\subseteq \mathcal{I}(k\varphi)$, it suffices to prove the assertion with $\psi$ in place of $\varphi$. In this case, the result follows from Hörmander's $L^2$-estimate (\cite{Hor66}), see \cite[Theorem~13.21]{Dem12} and \cite[Lemma~14.6.5 and Theorem~14.6.7]{FPT}.
\end{proof}

In the remainder of this section, we sharpen \cref{thm:concentration} for $C^{1,\overline{1}}$ barrier functions.

This part relies on various techniques not covered in this book.

\begin{theorem}[Tosatti, Chu--Zhou]\label{thm:Tos_C11}
Let $X$ be a compact Kähler manifold with a Kähler form $\omega$ and let $f\in C^{1,1}(X)$. Then
\[
P_\omega(f)\in C^{1,1}(X),
\]
and there is a constant $C=C(X,\omega,\|f\|_{C^{1,1}})$ such that
\[
    \|P_\omega(f)\|_{C^{1,1}}\leq C.
\]
\end{theorem}
\index{theorem!Tosatti--Chu--Zhou regularity}
This is the main result of \cite{Tos18,CZ19}. For a proof see \cite[Theorem~16.2.26]{FPT}. The proof of this theorem relies on some PDE techniques far beyond the scope of the current book. Weaker results were due to \cite{Ber09, Ber19}. See also \cite[Proposition~16.2.22]{FPT} for the approximation by penalized equations.

\begin{proposition}[Darvas--Rubinstein]\label{prop:DR16_C11bar}
Let $\omega$ be a K\"ahler form on $X$, and let $f_1,f_2\in C^{1,\overline{1}}(X)$. Then $P_\omega(\min\{f_1,f_2\})\in C^{1,\overline{1}}(X)$.
\end{proposition}
\index{theorem!Darvas--Rubinstein regularity}
This is \cite[Theorem~2.5]{DR16}. See also \cite[Theorem~16.2.24]{FPT}. The original proof relies unfortunately on \cite{BD12}, which can be replaced by \cref{thm:Tos_C11}. For a trace estimate for the penalized equation, see \cite[Lemma~16.2.15]{FPT}.

\begin{lemma}\label{lma:DNTroof}
Assume that $\omega$ is a K\"ahler form on $X$ and let $f_1,f_2\in C^{1,\overline{1}}(X)$. Then the functions $f_1$, $f_2$ and $P_\omega(\min\{f_1,f_2\})$ coincide up to second order at almost every point of the contact set $\{P_\omega(\min\{f_1,f_2\})=f_1=f_2\}$. In particular,
\[
    \mathds{1}_{\left\{P_\omega\left(\min\{f_1,f_2\}\right)=f_1=f_2\right\}}\,\omega_{f_1}^n=\mathds{1}_{\left\{P_\omega\left(\min\{f_1,f_2\}\right)=f_1=f_2\right\}}\,\omega_{f_2}^n.
\]
We also have for $j\in \{1,2\}$,
\begin{equation}\label{eq:DNTroof_Pomega}
    \begin{aligned}
    \omega^n_{P_\omega\left(\min\{f_1,f_2\}\right)} = & \mathds{1}_{\left\{P_\omega\left(\min\{f_1,f_2\}\right)=f_1\right\}}\,\omega^n_{f_1} + \mathds{1}_{\left\{P_\omega\left(\min\{f_1,f_2\}\right)=f_2\right\}}\,\omega^n_{f_2} \\
    &- \mathds{1}_{\left\{P_\omega(\min\{f_1,f_2\})=f_1=f_2\right\}}\,\omega^n_{f_j}\\
    \geq & \mathds{1}_{\left\{P_\omega\left(\min\{f_1,f_2\}\right)=f_j\right\}}\,\omega^n_{f_j}.
    \end{aligned}
\end{equation}
\end{lemma}

Here $\omega_{f_i}\coloneqq \omega+\ddc f_i\in L^n(X)$, since $f_i\in C^{1,\overline{1}}(X)$. Therefore, $\omega_{f_i}^n$ makes sense as an absolutely continuous Radon measure.

\begin{proof}
Set $\varphi\coloneqq P_\omega(\min\{f_1,f_2\})$. Since $\min\{f_1,f_2\}\in C^0(X)$, we have $\varphi\leq \min\{f_1,f_2\}$. By \cref{prop:DR16_C11bar}, $\varphi\in C^{1,\overline{1}}(X)$. As $f_1,f_2,\varphi\in C^{1,\overline{1}}(X)$, the measures $\omega^n_{f_1}$, $\omega^n_{f_2}$ and $\omega^n_\varphi$ are all absolutely continuous.

\textbf{Step~1}: We show that at almost every contact point, $\varphi$ and $f_i$ agree to order two.

For $i\in\{1,2,3\}$ set
\[
\Psi_1\coloneqq \varphi-f_1,\quad \Psi_2\coloneqq \varphi-f_2,\quad \Psi_3\coloneqq f_1-f_2.
\]
Each $\Psi_i$ lies in $C^{1,\overline{1}}(X)$.

We claim that for $i=1,2,3$, we have
\begin{equation}\label{eq:DNTroof_step1claim}
\Psi_i\textup{ vanishes to second order at almost every point of }\{\Psi_i=0\}.
\end{equation}
In local coordinates, the bounded Laplacian of $\Psi_i$ gives $\Psi_i\in W^{2,p}_{\loc}$ for every $1<p<\infty$ by interior elliptic regularity. For $p>2n$, Morrey's inequality applied to the gradient, after subtracting the linear part given by the weak Hessian at an $L^p$ Lebesgue point, shows that the gradient is classically differentiable there with derivative equal to that Hessian. By Stampacchia's lemma, $\mathrm{d}\Psi_i=0$ almost everywhere on $\{\Psi_i=0\}$, and its weak derivative is zero almost everywhere on $\{\mathrm{d}\Psi_i=0\}$. This proves \eqref{eq:DNTroof_step1claim}.

Specialising \eqref{eq:DNTroof_step1claim} to $\Psi_1,\Psi_2$, the functions $\varphi$ and $f_i$ agree to order two at almost every point of the contact set $\{\varphi=f_i\}$. Consequently $\omega+\ddc \varphi$ and $\omega+\ddc f_i$ have the same matrix of coefficients almost everywhere on $\{\varphi=f_i\}$, and therefore so do their $n$-fold wedge powers:
\begin{equation}\label{eq:DNTroof_step1eq}
\mathds{1}_{\{\varphi=f_i\}}\,\omega^n_\varphi = \mathds{1}_{\{\varphi=f_i\}}\,\omega^n_{f_i},\qquad i\in\{1,2\}.
\end{equation}
The analogous statement applied to $\Psi_3$ gives
\begin{equation}\label{eq:DNTroof_step1eqf12}
\mathds{1}_{\{f_1=f_2\}}\,\omega^n_{f_1} = \mathds{1}_{\{f_1=f_2\}}\,\omega^n_{f_2}.
\end{equation}

\textbf{Step~2}. The barrier $\min\{f_1,f_2\}$ is continuous, hence quasi-continuous, so \cref{thm:concentration} yields
\begin{equation}\label{eq:DNTroof_step2conc}
\int_{X\setminus \left(\{\varphi=f_1\}\cup\{\varphi=f_2\}\right)}\omega^n_\varphi
=
\int_{\left\{\varphi<\min\{f_1,f_2\}\right\}}\omega^n_\varphi
=0.
\end{equation}

Set $V_1\coloneqq \{f_1<f_2\}$ and $V_2\coloneqq \{f_1>f_2\}$, both open. We may decompose
\[
\{\varphi=f_1\}\cup\{\varphi=f_2\}
=
\bigl(V_1\cap\{\varphi=f_1\}\bigr)
\sqcup
\bigl(V_2\cap\{\varphi=f_2\}\bigr)
\sqcup
\{\varphi=f_1=f_2\}.
\]
This combined with \eqref{eq:DNTroof_step2conc} and \eqref{eq:DNTroof_step1eq} yields
\[
\begin{aligned}
\omega^n_\varphi
={}&
\mathds{1}_{V_1\cap\{\varphi=f_1\}}\,\omega^n_\varphi
+
\mathds{1}_{V_2\cap\{\varphi=f_2\}}\,\omega^n_\varphi
+
\mathds{1}_{\{\varphi=f_1=f_2\}}\,\omega^n_\varphi
\\
={}&
\mathds{1}_{V_1\cap\{\varphi=f_1\}}\,\omega^n_{f_1}
+
\mathds{1}_{V_2\cap\{\varphi=f_2\}}\,\omega^n_{f_2}
+
\mathds{1}_{\{\varphi=f_1=f_2\}}\,\omega^n_{f_j}.
\end{aligned}
\]
On the other hand,
\[
\mathds{1}_{V_j\cap\{\varphi=f_j\}}\,\omega^n_{f_j}
=
\mathds{1}_{\{\varphi=f_j\}}\,\omega^n_{f_j}
-
\mathds{1}_{\{\varphi=f_1=f_2\}}\,\omega^n_{f_j}.
\]
Therefore,
\[
\omega^n_\varphi
=
\mathds{1}_{\{\varphi=f_1\}}\,\omega^n_{f_1}
+
\mathds{1}_{\{\varphi=f_2\}}\,\omega^n_{f_2}
+
\mathds{1}_{\{\varphi=f_1=f_2\}}
\,\bigl(\omega^n_{f_j}-\omega^n_{f_1}-\omega^n_{f_2}\bigr).
\]
Together with \eqref{eq:DNTroof_step1eqf12} we obtain \eqref{eq:DNTroof_Pomega}.
\end{proof}

\begin{theorem}[Di Nezza--Trapani]\label{thm:DNTcontact}
Let $\theta_1,\dots,\theta_n$ be closed real smooth $(1,1)$-forms on $X$ representing pseudo-effective cohomology classes. For each $i=1,\dots,n$, let $\varphi_i\in\PSH(X,\theta_i)$ and $f_i\in C^{1,\overline{1}}(X)$ be functions with $\varphi_i\leq f_i$. Then
\begin{equation}\label{eq:DNTcontact}
\mathds{1}_{\bigcap_i\{\varphi_i=f_i\}}\, \theta_{1,\varphi_1}\wedge\cdots\wedge\theta_{n,\varphi_n} = \mathds{1}_{\bigcap_i\{\varphi_i=f_i\}}\, \theta_{1,f_1}\wedge\cdots\wedge\theta_{n,f_n}.
\end{equation}
\end{theorem}

\begin{proof}
    \textbf{Step~1}. We first reduce to the case where $\theta_1=\cdots=\theta_n$, $\varphi_1=\cdots=\varphi_n$ and $f_1=\cdots=f_n$.

    Assume for the moment that this special case is known. Let us handle the general case.

    For $\lambda=(\lambda_1,\dots,\lambda_n)\in(\mathbb{R}_{>0})^n$, we set
    \[
    \theta_\lambda\coloneqq \sum_{j=1}^n\lambda_j\theta_j,\quad \varphi_\lambda\coloneqq \sum_{j=1}^n\lambda_j\varphi_j,\quad f_\lambda\coloneqq \sum_{j=1}^n\lambda_j f_j.
    \]
    Then $\varphi_\lambda\in\PSH(X,\theta_\lambda)$, $f_\lambda\in C^{1,\overline{1}}(X)$, $\varphi_\lambda\leq f_\lambda$, and
    \[
        \bigcap_{j=1}^n\{\varphi_j=f_j\}\subseteq\{\varphi_\lambda=f_\lambda\}.
    \]
    So
    \[
        \mathds{1}_{\bigcap_j\{\varphi_j=f_j\}}(\theta_\lambda)^n_{\varphi_\lambda} = \mathds{1}_{\bigcap_j\{\varphi_j=f_j\}}(\theta_\lambda)^n_{f_\lambda}.
    \]
Both sides are polynomials in $\lambda_1,\dots,\lambda_n$; comparing the coefficient of $\lambda_1\cdots\lambda_n$ gives \eqref{eq:DNTcontact}.

From now on, we assume
\[
\theta=\theta_1=\dots=\theta_n,\quad
\varphi=\varphi_1=\cdots=\varphi_n\in \PSH(X,\theta),
\]
and
\[
f=f_1=\cdots=f_n\in C^{1,\overline{1}}(X)
\]
with
\[
\varphi\leq f.
\]
We wish to prove that
\begin{equation}\label{eq:DNT_temp1}
    \mathds{1}_{\left\{\varphi=f\right\}}\,\theta_{\varphi}^n=\mathds{1}_{\left\{\varphi=f\right\}}\,\theta_{f}^n.
\end{equation}

\textbf{Step~2}. We first handle the special case where $f\in \PSH(X,\theta)$.

\textbf{Step~2.1}. We assume that $\theta$ is a Kähler form, denoted by $\omega$ for clarity.

Pick a decreasing sequence $(\phi_j)_j$ in $C^{\infty}(X)$\footnote{The functions $\phi_j$ are not required to be $\omega$-psh.} with pointwise limit $\varphi$ and set
\[
\psi_j\coloneqq P_{\omega}\left(\min\{\phi_j,f\}\right).
\]
Then $\psi_j\in\PSH(X,\omega)$, and
\[
\varphi\leq\psi_j\leq\min\{\phi_j,f\} \leq f,
\]
and hence $\psi_j\searrow\varphi$.
Applying \cref{lma:DNTroof} with $f_1=\phi_j$ and $f_2=f$ gives
\begin{equation}\label{eq:omegapsij_temp1}
\omega^n_{\psi_j}\geq \mathds{1}_{\{\psi_j=f\}}\,\omega^n_f \geq \mathds{1}_{\{\varphi=f\}}\,\omega^n_f.
\end{equation}

Fix $C>0$ with $\min_X f>-C$ and $\epsilon>0$.
Let
\[
h_j^{C,\epsilon}\coloneqq \frac{(\psi_j+C)\lor 0}{\left((\psi_j+C)\lor 0 \right)+\epsilon},\qquad h^{C,\epsilon}\coloneqq \frac{(\varphi+C)\lor 0}{\left((\varphi+C)\lor 0\right)+\epsilon}.
\]
These functions are both quasi-continuous with values in $[0,1]$.
Observe that $h_j^{C,\epsilon}$ converges in capacity to $h^{C,\epsilon}$ as $j\to\infty$.

Now fix $g\in C^0(X)$, $g\geq 0$. Note that $h_j^{C,\epsilon}=0$ on the set $\{\psi_j\leq -C\}$ and $h^{C,\epsilon}=0$ on the set $\{\varphi\leq -C\}$.
Using locality of the non-pluripolar product \itemref{itm:prop:npp1:1}, we have
\[
h_j^{C,\epsilon}g \,\omega_{\psi_j}^n=h_j^{C,\epsilon}g \,\omega_{\psi_j\lor (-C)}^n,\quad h^{C,\epsilon}g \,\omega_{\varphi}^n=h^{C,\epsilon}g \,\omega_{\varphi\lor (-C)}^n.
\]
Therefore,
\[
\begin{aligned}
\int_X g\,\omega_{\varphi}^n\geq & \int_X h^{C,\epsilon}g \,\omega_{\varphi}^n & \\
=& \int_X h^{C,\epsilon}g \,\omega_{\varphi\lor (-C)}^n &\\
=& \lim_{j\to\infty} \int_X h_j^{C,\epsilon}g \,\omega_{\psi_j\lor (-C)}^n & \textup{by \cref{thm:loc_conv_cap}}\\
=&\lim_{j\to\infty} \int_X h_j^{C,\epsilon}g \,\omega_{\psi_j}^n &\\
\geq & \lim_{j\to\infty} \int_{\{\varphi=f\}} h_j^{C,\epsilon}g \,\omega_{f}^n & \textup{by }\eqref{eq:omegapsij_temp1}\\
=& \int_{\{\varphi=f\}}h^{C,\epsilon}g \,\omega_{f}^n,
\end{aligned}
\]
where the last equality follows from the dominated convergence theorem.
By our choice of $C$, we have $h^{C,\epsilon}>0$ on $\{\varphi=f\}$, and hence when $\epsilon\to 0+$, $h^{C,\epsilon}\to 1$ everywhere on this set. Letting $\epsilon\to 0+$, we find
\[
\int_X g\,\omega_{\varphi}^n\geq \int_{\{\varphi=f\}}g \,\omega_{f}^n
\]
by the dominated convergence theorem.

By the Riesz representation theorem, we conclude that
\[
\omega_{\varphi}^n\geq \mathds{1}_{\{\varphi=f\}}\,\omega_{f}^n.
\]
Then by \cref{lma:lorMAineqe}, the reverse inequality also holds, and \eqref{eq:DNT_temp1} follows in this case.

\textbf{Step~2.2}. We no longer assume that $\theta$ is a Kähler form.

Choose a Kähler form $\omega$ on $X$ so that $\theta+\omega$ is a Kähler form and $f\in \PSH(X,\theta+\omega)$.
Fix $g\in C^0(X)$ and consider, for $t\geq 0$,
\[
Q(t)\coloneqq \int_{\{\varphi=f\}}g\,(\theta+t\omega+\ddc\varphi)^n - \int_{\{\varphi=f\}}g\,(\theta+t\omega+\ddc f)^n.
\]
By multilinearity of the non-pluripolar product, $Q$ is a polynomial in $t$. Step~2.1 applied to $\theta+t\omega$ for $t>1$ gives $Q(t)\equiv 0$ for $t>1$, hence $Q\equiv 0$, and in particular $Q(0)=0$. Since $g$ is arbitrary, \eqref{eq:DNT_temp1} holds.

\textbf{Step~3}. We handle the general case. We no longer assume that $f\in \PSH(X,\theta)$.

Choose a Kähler form $\omega$ on $X$ so that $\theta+\omega$ is a Kähler form. Then
\[
\varphi\leq P_{\theta+\omega}(f)\leq f.
\]
By \cref{lma:DNTroof} applied to $f_1=f_2=f$, we find that $f$ and $P_{\theta+\omega}(f)$ agree to the second order on the contact set $\{ P_{\theta+\omega}(f)=f \}$, so
\[
\mathds{1}_{\left\{ P_{\theta+\omega}(f)=f \right\}}\,\theta_{P_{\theta+\omega}(f)}^n=\mathds{1}_{\left\{ P_{\theta+\omega}(f)=f \right\}}\,\theta_{f}^n.\footnotemark
\]
\footnotetext{This expression looks slightly unusual, as $P_{\theta+\omega}(f)$ is not $\theta$-psh in general, but it does make sense.}
However, Step~2 implies that
\[
\mathds{1}_{\left\{ P_{\theta+\omega}(f)=\varphi \right\}}\,\theta_{\varphi}^n=\mathds{1}_{\left\{ P_{\theta+\omega}(f)=\varphi \right\}}\,\theta_{P_{\theta+\omega}(f)}^n.
\]
Combining these results, we conclude \eqref{eq:DNT_temp1}.
\end{proof}

\begin{corollary}\label{cor:DNTenvelope1}
    Assume that $\PSH(X,\theta)\neq \varnothing$. For any $f\in C^{1,\overline{1}}(X)$, we have
    \[
    \theta^n_{P_\theta(f)} = \mathds{1}_{\{P_\theta(f)=f\}}\,\theta^n_f.
    \]
    In particular,
    \[
        \theta_{V_{\theta}}^n=\mathds{1}_{\{V_{\theta}=0\}}\,\theta^n.
    \]
\end{corollary}
\begin{proof}
    From \cref{thm:DNTcontact}, we find
    \[
    \mathds{1}_{\{P_\theta(f)=f\}}\,\theta^n_{P_\theta(f)}=\mathds{1}_{\{P_\theta(f)=f\}}\,\theta^n_f.
    \]
    By the concentration property \cref{thm:concentration}, our assertion follows.
\end{proof}

\begin{theorem}[Di Nezza--Trapani]\label{thm:DNT23_C11bar}
Let $\theta$ be a smooth closed real $(1,1)$-form on $X$ representing a big cohomology class. Let $f\in C^{1,\overline{1}}(X)$. Then
\[
P_\theta(f)\in C^{1,\overline{1}}\bigl(X\setminus\nK(\{\theta\})\bigr).
\]
In particular, $V_{\theta}\in C^{1,\overline{1}}\bigl(X\setminus\nK(\{\theta\})\bigr)$.
\end{theorem}
\index{theorem!Di Nezza--Trapani regularity}

We refer to \cite{DNT23} and \cite[Theorem~16.2.29]{FPT} for the proof.

\motto{Politiques et scientifiques ont le sens des réalités, mais ce ne sont pas les mêmes. Il en résulte --- et ce sera là un principe que le général de Gaulle fera sien que l'activité de recherche ne peut être évaluée, quant à sa qualité propre, que par des hommes qui la pratiquent eux-mêmes.

--- Pierre Lelong\footnote{Pierre Lelong (1912--2011) was the husband of another famous mathematician Jacqueline Ferrand. During their marriage (1947--1977), the latter published under the name of Jacqueline Lelong-Ferrand.}, 1999}
\chapter{The envelope operators}
\label{chap:envelopes}

This chapter develops the envelope formalism that will be used throughout the rest of the book. There are two parallel constructions. The first is the $P$-envelope, which records the singularity type relevant to non-pluripolar Monge--Amp\`ere mass. The second is the $\mathcal I$-envelope, which records the singularity type detected by multiplier ideal sheaves. The distinction between these two envelopes is one of the main organizing principles of the book.

The starting point is the rooftop operator. Given two singular potentials, the rooftop construction produces the largest potential lying below both, and it provides a flexible way to compare singularities. From this we define the $P$-envelope $P_\theta[\varphi]$. Its fixed points are the model potentials. Model potentials are the natural representatives of $P$-singularity types: they preserve non-pluripolar mass and satisfy strong contact-set properties. The first part of the chapter establishes the basic algebra of these envelopes, their behavior under monotone limits, and the domination and contact principles needed later.

We then introduce the Monge--Amp\`ere energy relative to a model potential. This leads to the finite-energy classes $\mathcal E^1(X,\theta;\phi)$ and to stability results for rooftop envelopes. These results will be used in the metric theory of singularities and in the study of geodesic rays.

The second part of the chapter develops the $\mathcal I$-envelope $P_\theta[\varphi]_{\mathcal I}$. This envelope is defined by imposing multiplier ideal conditions rather than ordinary singularity-type conditions. It is therefore sensitive to valuative data in a way that the ordinary $P$-envelope is not. We prove its basic functorial properties, its independence of the chosen representative of the cohomology class, its bimeromorphic behavior, and its compatibility with analytic singularities.

The comparison between $P_\theta[\varphi]$ and $P_\theta[\varphi]_{\mathcal I}$ is delicate. They agree in important cases, for instance for analytic singularities, but they differ in general. Later chapters use this comparison to define $\mathcal I$-good singularities, trace operators, partial Okounkov bodies and non-Archimedean potentials. Thus the technical results of this chapter should be regarded as the basic calculus of singularity types used in the sequel.

Several results about $P$-envelopes are developed by Darvas--Di Nezza--Lu in various papers and summarized in their survey \cite{DDNLsurv}. The reader is invited to read their survey for more details and further results, and may consult \cite[Chapters~16--17]{FPT} for proofs.

\section{The \texorpdfstring{$P$}{P}-envelope}\label{sec:Penv}
In this section, $X$ will denote a connected compact Kähler manifold of dimension $n$.

\subsection{Rooftop operator and the definition of the \texorpdfstring{$P$}{P}-envelope}
\label{subsec:rooftop_Penv}

Fix a smooth closed real $(1,1)$-form $\theta$ on $X$.

\begin{definition}\label{def:rooftop}
    Given $\varphi,\psi\in \PSH(X,\theta)$, we define their \emph{rooftop operator}\index{rooftop operator} as follows:
    \begin{equation}\label{eq:rooftopdef1}
    \varphi\land\psi=\sup\left\{\eta\in \PSH(X,\theta):\eta\leq \varphi,\eta\leq \psi\right\}.
    \end{equation}
    \index{$\varphi\land \psi$}
    For simplicity of notation, we extend the definition to the case where $\varphi$ or $\psi$ is constantly $-\infty$, in which case we simply set
    \[
    \varphi\land \psi=-\infty.
    \]
\end{definition}
When we want to be more specific, we may also write $\varphi\land_{\theta}\psi$. Note that it is possible that the set in \eqref{eq:rooftopdef1} is empty, in which case $\varphi\land \psi=-\infty$.

\begin{proposition}\label{prop:land_well_def}
    The operator $\land$ is a well-defined commutative, associative binary operator
    \[
    \bigl(\PSH(X,\theta)\cup\{-\infty\}\bigr)\times\bigl(\PSH(X,\theta)\cup\{-\infty\}\bigr)\rightarrow \PSH(X,\theta)\cup\{-\infty\}.
    \]
\end{proposition}
\begin{proof}
    We first show that the map is well-defined. For this purpose, take $\varphi,\psi\in \PSH(X,\theta)$. When the set in \eqref{eq:rooftopdef1} is empty, there is nothing to prove. So let us assume that the set is not empty.

    Define
    \[
    \gamma=\uscsup \left\{\eta\in \PSH(X,\theta):\eta\leq \varphi,\eta\leq \psi\right\}.
    \]
    Then by \cref{prop:pshfunction_closedseq}, we find that $\gamma\in \PSH(X,\theta)$ and hence $\gamma$ is a candidate for the supremum in \eqref{eq:rooftopdef1}. Therefore, $\gamma\leq \varphi\land \psi$. The reverse inequality follows from the definition of $\gamma$, so
    \[
    \varphi\land \psi=\gamma\in \PSH(X,\theta).
    \]
    The commutativity and the associativity of $\land$ follow directly from the defining maximality property.
\end{proof}
Note that from the proof, we have
\[
\varphi\land \psi=P_{\theta}\Bigl( \min\{\varphi,\psi\} \Bigr).
\]
See \cref{def:Pthetaf} for the notation $P_{\theta}(\bullet)$.

\begin{lemma}\label{lma:rooftopMA}
    Let $\varphi,\psi\in \PSH(X,\theta)$. Assume that $\varphi\land\psi\in \PSH(X,\theta)$. Then
    \begin{equation}\label{eq:rooftopMAineq1}
    \theta_{\varphi\land \psi}^n\leq \mathds{1}_{\{\varphi\land\psi=\varphi\}}\, \theta_{\varphi}^n+\mathds{1}_{\{\varphi\land\psi=\psi\}} \,\theta_{\psi}^n.
    \end{equation}
\end{lemma}
In particular, $\theta_{\varphi\land \psi}^n$ is carried by the union $\{\varphi\land \psi=\varphi\}\cup \{\varphi\land \psi=\psi\}$.

\begin{proof}
We first observe that as a consequence of \cref{lma:lorMAineqe}, we have
\[
\mathds{1}_{\{\varphi\land \psi=\varphi\}} \,\theta_{\varphi\land \psi}^n\leq \mathds{1}_{\{\varphi\land \psi=\varphi\}} \,\theta_{\varphi}^n,\quad \mathds{1}_{\{\varphi\land \psi=\psi\}} \,\theta_{\varphi\land \psi}^n\leq \mathds{1}_{\{\varphi\land \psi=\psi\}} \,\theta_{\psi}^n.
\]
Applying \cref{thm:concentration} to $\min\{\varphi,\psi\}$, we conclude that
\[
\theta_{\varphi\land \psi}^n\leq \mathds{1}_{\{\varphi\land \psi=\varphi\}} \,\theta_{\varphi\land \psi}^n+\mathds{1}_{\{\varphi\land \psi=\psi\}}\,\theta_{\varphi\land \psi}^n\leq \mathds{1}_{\{\varphi\land\psi=\varphi\}} \,\theta_{\varphi}^n+\mathds{1}_{\{\varphi\land\psi=\psi\}} \,\theta_{\psi}^n,
\]
and \eqref{eq:rooftopMAineq1} is established.
\end{proof}

We recall that the relations $\preceq$ and $\sim$ are introduced in \cref{def:parorder}.
\begin{definition}\label{def:Penv}
    Given $\varphi\in \PSH(X,\theta)$, we define its \emph{$P$-envelope}\index{envelope!P-envelope@$P$-envelope} as follows:
    \begin{equation}\label{eq:Pthetavarphi}
    P_{\theta}[\varphi]\coloneqq \uscsup \left\{ \psi\in \PSH(X,\theta): \psi\leq 0,\psi\preceq \varphi \right\}.
    \end{equation}
    \index{$P_{\theta}[\varphi]$}
\end{definition}
Note that $\varphi-\sup_X \varphi$ is always a candidate for the supremum. By \cref{prop:pshfunction_closedseq}, we have $P_{\theta}[\varphi]\in \PSH(X,\theta)$ and $P_{\theta}[\varphi]\leq 0$.
Moreover, the definition can be equivalently described as
\begin{equation}\label{eq:Penvsups}
    P_{\theta}[\varphi]=\uscsup[C\in \mathbb{Z}_{>0}](\varphi+C)\land V_{\theta}.
\end{equation}
Recall that $V_{\theta}$ is introduced in \eqref{eq:Vtheta}. For any $C\in \mathbb{R}$, we have $(\varphi+C)\land V_{\theta}\in \PSH(X,\theta)$ and
\[
(\varphi+C)\land V_{\theta}\sim \varphi.
\]
In other words, in \eqref{eq:Pthetavarphi}, we may replace the condition $\psi\preceq \varphi$ by $\psi\sim \varphi$. Also note that as a consequence,
\begin{equation}\label{eq:varphipreceqP}
    \varphi\preceq P_{\theta}[\varphi].
\end{equation}

The idea lying behind the definition of $P_{\theta}[\varphi]$ is that we choose the least singular element out of all potentials with the same singularity type as $\varphi$. As we shall see in \cref{ex:loglogsing} below, $P_{\theta}[\varphi]$ does not necessarily have the same singularity type as $\varphi$. This forces us to define a rougher equivalence relation in \cref{def:Pmoresing}.

The envelope depends on the choice of $\theta$, but the dependence is easy to understand:
\begin{proposition}\label{prop:Penvindeptheta}
    Let $\theta'=\theta+\ddc g$ for some $g\in C^{\infty}(X)$. Then for any $\varphi\in \PSH(X,\theta)$, we have $\varphi'\coloneqq \varphi-g\in \PSH(X,\theta')$ and
    \[
        P_{\theta}[\varphi]\sim P_{\theta'}[\varphi'].
    \]
\end{proposition}
\begin{proof}
    By symmetry, it suffices to show that
    \[
        P_{\theta}[\varphi]\preceq P_{\theta'}[\varphi'].
    \]
    We may assume that $g\geq 0$ (after replacing $g$ by $g+C$ with $C\geq -\inf_X g$, which leaves $\theta'$ unchanged and shifts $\varphi'$ by an additive constant).
    Then for any $\psi\in \PSH(X,\theta)$ with $\psi\preceq \varphi$ and $\psi\leq 0$, we set $\psi'\coloneqq \psi-g\in \PSH(X,\theta')$. Then $\psi'\preceq \varphi'$ and $\psi'\leq 0$, so $\psi'\leq P_{\theta'}[\varphi']$. Since $\psi$ is arbitrary, it follows that
    \[
        P_{\theta}[\varphi]-\sup_X g\leq P_{\theta}[\varphi]-g\leq P_{\theta'}[\varphi'],
    \]
    which gives $P_{\theta}[\varphi]\preceq P_{\theta'}[\varphi']$.
\end{proof}

The $P$-envelope preserves the non-pluripolar masses:
\begin{proposition}\label{prop:Ppresmass}
Suppose that $\theta_1,\ldots,\theta_n$ are smooth closed real $(1,1)$-forms on $X$. Let $\varphi_i\in \PSH(X,\theta_i)$ for each $i=1,\ldots,n$. Then
\begin{equation}\label{eq:Penvpremass}
\int_X \theta_{1,P_{\theta_1}[\varphi_1]}\wedge \cdots \wedge \theta_{n,P_{\theta_n}[\varphi_n]}=\int_X \theta_{1,\varphi_1}\wedge \cdots\wedge \theta_{n,\varphi_n}.
\end{equation}
\end{proposition}
This proposition together with \cref{thm:mono} will be combined to a common generalization in \cref{prop:mono2} after introducing the $P$-preorder.
\begin{proof}
For each $C\in \mathbb{Z}_{>0}$ and each $i=1,\ldots,n$, we have
\[
(\varphi_i+C)\land V_{\theta_i}\sim \varphi_i.
\]
It follows from the monotonicity theorem \cref{thm:mono} that
\[
\int_X \theta_{1,(\varphi_1+C)\land V_{\theta_1}}\wedge \cdots \wedge \theta_{n,(\varphi_n+C)\land V_{\theta_n}}=\int_X \theta_{1,\varphi_1}\wedge \cdots\wedge \theta_{n,\varphi_n}.
\]
So \eqref{eq:Penvpremass} follows from \eqref{eq:Penvsups} and \cref{cor:incseqnppcont}.
\end{proof}

Conversely, \cref{prop:Ppresmass} characterizes the $P$-envelope, as we will see in a moment in \cref{thm:Pvarphidiffdef}. We need some preparation for the proof. We first show that the $P$-envelope can be regarded as a concentration of the non-pluripolar mass:
\begin{theorem}\label{thm:Pvarphisupport}
    Let $\varphi\in \PSH(X,\theta)$. Then
    \begin{equation}\label{eq:thetaPthetaphin}
        \theta_{P_{\theta}[\varphi]}^n\leq \mathds{1}_{\{P_{\theta}[\varphi]=0\}}\,\theta^n.
    \end{equation}
\end{theorem}
In fact, equality holds as we will see in \cref{prop:DNTenvelope2} below.
\begin{proof}
    Thanks to \cref{lma:rooftopMA}, for each $C>0$, we have
    \[
    \begin{aligned}
    \theta_{(\varphi+C)\land V_{\theta}}^n\leq & \mathds{1}_{\{(\varphi+C)\land V_{\theta}=\varphi+C\} }\, \theta_{\varphi}^n+\mathds{1}_{\{(\varphi+C)\land V_{\theta}=V_{\theta}\} }\, \theta_{V_{\theta}}^n\\
    \leq & \mathds{1}_{\{\varphi+C\leq V_{\theta}\} }\, \theta_{\varphi}^n+\mathds{1}_{\{P_{\theta}[\varphi]=V_{\theta}\} } \,\theta_{V_{\theta}}^n.
    \end{aligned}
    \]
    We wish to let $C\to\infty$.
    The dominated convergence theorem and the semicontinuity theorem \cref{cor:incseqnppcont} give, as $C\to\infty$,
    \[
    \mathds{1}_{\{\varphi+C\leq V_{\theta}\} } \theta_{\varphi}^n\rightharpoonup 0,\quad
    \theta_{(\varphi+C)\land V_{\theta}}^n\rightharpoonup \theta_{P_{\theta}[\varphi]}^n.
    \]
    So we conclude that
    \[
    \theta_{P_{\theta}[\varphi]}^n\leq \mathds{1}_{\{P_{\theta}[\varphi]=V_{\theta}\}} \,\theta_{V_{\theta}}^n.
    \]
    Taking \cref{cor:DNTenvelope1} into consideration, we conclude \eqref{eq:thetaPthetaphin}.
\end{proof}

Using essentially the same proof, we arrive at the following conclusion:
\begin{corollary}\label{cor:Pvarphisupport}
    Let $\varphi,\psi\in \PSH(X,\theta)$. Let
    \[
    \eta\coloneqq \uscsup[C>0] (\varphi+C)\land \psi.\footnotemark
    \]
    \footnotetext{Using the terminology to be introduced in \cref{def:Pthetavarphif}, we can also write $\eta=P_{\theta}[\varphi](\psi)$.}
    Assume that $\eta\not\equiv -\infty$ (for example when $\varphi\preceq \psi$).
    Then
    \[
    \theta_{\eta}^n\leq \mathds{1}_{\{\eta=\psi\}}\,\theta_{\psi}^n.
    \]
\end{corollary}
\begin{proof}
    Thanks to \cref{lma:rooftopMA}, for each $C>0$, we have
    \[
    \begin{aligned}
    \theta_{(\varphi+C)\land \psi}^n\leq & \mathds{1}_{\{(\varphi+C)\land \psi=\varphi+C\} }\, \theta_{\varphi}^n+\mathds{1}_{\{(\varphi+C)\land \psi=\psi\} }\, \theta_{\psi}^n\\
    \leq & \mathds{1}_{\{\varphi+C\leq \psi\} }\, \theta_{\varphi}^n+\mathds{1}_{\{\eta=\psi\} } \,\theta_{\psi}^n.
    \end{aligned}
    \]
    We wish to let $C\to\infty$.
    The dominated convergence theorem and the semicontinuity theorem \cref{cor:incseqnppcont} give, as $C\to\infty$,
    \[
    \mathds{1}_{\{\varphi+C\leq \psi\} } \theta_{\varphi}^n\rightharpoonup 0,\quad
    \theta_{(\varphi+C)\land \psi}^n\rightharpoonup \theta_{\eta}^n.
    \]
    So we conclude that
    \[
    \theta_{\eta}^n\leq \mathds{1}_{\{\eta=\psi\}} \,\theta_{\psi}^n.
    \]
    This completes the proof.
\end{proof}

\begin{theorem}\label{thm:Pvarphidiffdef}
    Assume that $\varphi\in \PSH(X,\theta)_{>0}$. Then
    \begin{equation}\label{eq:Penvdef}
    P_{\theta}[\varphi]=\sup\left\{\psi\in \PSH(X,\theta): \psi\leq 0,\varphi\preceq \psi, \int_X \theta_{\varphi}^n=\int_X \theta_{\psi}^n\right\}.
    \end{equation}
    In particular, in this case,\footnote{Pang--Sun--Wang--Zhou \cite{PSWZ26} recently proved that, \eqref{eq:Penvprojop} continues to hold even when $\varphi$ has vanishing mass. See also \cite[Theorem~17.5.3]{FPT}.}
    \begin{equation}\label{eq:Penvprojop}
    P_{\theta}\Bigl[P_{\theta}[\varphi]\Bigr]=P_{\theta}[\varphi].
    \end{equation}
\end{theorem}
Note that in \eqref{eq:Penvdef} and \eqref{eq:Pthetavarphi}, the auxiliary function $\psi$ lies on different sides of $\varphi$.
\begin{proof}
     Let $\psi$ be a candidate of the right-hand side of \eqref{eq:Penvdef}. It follows from \cref{thm:Pvarphisupport} that
     \[
     \int_{\{P_{\theta}[\varphi]<\psi \}} \theta_{P_{\theta}[\varphi]}^n\leq \int_{\{P_{\theta}[\varphi]<\psi \}\cap \{P_{\theta}[\varphi]=0\} } \theta^n=0.
     \]
     On the other hand, we have
    \[
    P_{\theta}[\varphi]\preceq P_{\theta}[\psi],\quad \psi\preceq P_{\theta}[\psi]
    \]
    and all these potentials have the same mass due to \cref{prop:Ppresmass}. So the domination principle \cref{thm:dom_prin_1} is applicable and gives
    \[
    P_{\theta}[\varphi]\geq \psi.
    \]
    Hence we get the $\geq$ direction in \eqref{eq:Penvdef}.

    On the other hand, $P_{\theta}[\varphi]$ itself is a candidate of the right-hand side of \eqref{eq:Penvdef}, as a consequence of \cref{prop:Ppresmass}. Therefore, \eqref{eq:Penvdef} follows.

    As for \eqref{eq:Penvprojop}, the $\geq$ direction is trivial. For the other direction, let $\psi\in \PSH(X,\theta)$ be a potential satisfying
    \[
    \psi\leq 0,\quad P_{\theta}[\varphi]\preceq \psi,\quad \int_X \theta_{P_{\theta}[\varphi]}^n=\int_X \theta_{\psi}^n.
    \]
    Then it follows from \cref{prop:Ppresmass} that
    \[
    \psi\leq 0,\quad \varphi\preceq \psi,\quad \int_X \theta_{\varphi}^n=\int_X \theta_{\psi}^n.
    \]
    In view of \eqref{eq:Penvdef}, we conclude the $\leq$ direction of \eqref{eq:Penvprojop}.
\end{proof}

\begin{definition}\label{def:modelpot}
    If $\varphi=P_{\theta}[\varphi]$ and $\int_X \theta_{\varphi}^n>0$, we say $\varphi$ is a \emph{model potential}.\index{potential!model potential}
\end{definition}
Recall that the notion of model potentials depends heavily on the choice of $\theta$. When there is a risk of confusion, we also say $\varphi$ is a model potential in $\PSH(X,\theta)$.

\begin{remark}\label{rmk:envelopes:L267}
\cref{def:modelpot} is different from the common definition in the literature: We impose the extra condition $\int_X \theta_{\varphi}^n>0$. The author believes that this is the only case where this notion is natural. We sometimes emphasize this point by saying $\varphi\in \PSH(X,\theta)_{>0}$ is a model potential.
\end{remark}

There are plenty of model potentials:
\begin{corollary}\label{cor:Psendspotentialtomodel}
    Let $\varphi\in \PSH(X,\theta)_{>0}$. Then $P_{\theta}[\varphi]$ is a model potential in $\PSH(X,\theta)$. Moreover,
    \[
        \int_X \theta_{P_{\theta}[\varphi]}^n=\int_X \theta_{\varphi}^n.
    \]
\end{corollary}
\begin{proof}
    This follows immediately from \cref{thm:Pvarphidiffdef} and \cref{prop:Ppresmass}.
\end{proof}

As we have already seen in the proof of \cref{lma:pathoenvelope}, we have the following interesting property:
\begin{proposition}\label{prop:model_sup_prop}
    Let $\varphi\in \PSH(X,\theta)$. Consider a model potential $\phi\in \PSH(X,\theta)_{>0}$. Assume that $\varphi\preceq \phi$. Then
    \[
    \sup_{\{\phi\neq -\infty\}} (\varphi-\phi)=\sup_X \varphi.
    \]
\end{proposition}
In particular, for any $\varphi\in \PSH(X,\theta)$, $\varphi\preceq \phi$ and $\varphi\leq 0$ implies $\varphi\leq \phi$.

\begin{proof}
    Observe that $\phi\leq 0$, so on the set $\{\phi\neq -\infty\}$, we have
    \[
    \varphi-\phi\geq \varphi.
    \]
    So the $\geq$ direction follows from \cref{cor:diffsupinfindeppluripolar}.
    Conversely, by assumption, we can find a constant $C>0$ so that
    \[
    \varphi-\sup_X \varphi\leq \min\{\phi+C,0 \}.
    \]
    It follows that $\varphi-\sup_X \varphi\leq P_{\theta}[\phi]=\phi$. Therefore, the reverse inequality follows.
\end{proof}

\begin{example}\label{ex:loglogsing}
    We continue our favorite example \cref{ex:LPP1}. Let $X=\mathbb{P}^1$ and $\omega$ be the Fubini--Study metric. We define $\varphi\in \PSH(X,\omega)$ as follows: For $z\in \mathbb{C}$, we let
    \[
    \varphi(z)=
    \begin{cases}
        -\log \left(|z|^2+1\right)+\Bigl(-\log \left(-\log |z|^2\right)\Bigr) \lor \Bigl(2+\log |z|^2\Bigr),& \textup{if }|z|<1/\sqrt{2},\\
        2+\log \frac{|z|^2}{|z|^2+1},& \textup{otherwise},
    \end{cases}
    \]
    while $\varphi(\infty)=2$. Note that $\varphi$ is $\omega$-psh on $X=\mathbb{P}^1$ by the gluing lemma \cref{lma:gl}.
The singularity of $\varphi$ only occurs at $z=0$, close to which, $\varphi$ is approximately $-\log \left(-\log |z|^2\right)$. This type of singularity is therefore called the \emph{log-log type singularity}.

    We claim that
    \begin{equation}\label{eq:Pomegavarphi0temp1}
        P_{\omega}[\varphi]=0.
    \end{equation}
    In particular, we find that $\varphi$ and $P_{\omega}[\varphi]$ have different singularity types.

    By \cref{thm:Pvarphidiffdef}, in order to verify \eqref{eq:Pomegavarphi0temp1}, it suffices to verify that
    \begin{equation}\label{eq:Pomegavarphi0temp2}
    \int_X \omega_{\varphi}^1=1.
    \end{equation}
    Here $\omega_{\varphi}^1$ is taken in the non-pluripolar sense. Since $\{0,\infty\}\subseteq \mathbb{P}^1$ is pluripolar, this reduces to showing that
    \[
    \int_{\mathbb{C}^*} \ddc \psi=\frac{1}{4\pi}\int_{\mathbb{C}^*} (\Delta \psi)\,\mathrm{d}\mu=1,
    \]
    where $\psi(z)=\varphi(z)+\log (|z|^2+1)$ and $\mu$ is the standard Lebesgue measure on $\mathbb{C}$.

    Note that the Laplacian vanishes outside $\overline{B(0,0.7)}$ since $\psi(z)=2+\log |z|^2$ there, which is harmonic. Therefore,
    \[
    \int_{\mathbb{C}^*} \ddc \psi=\frac{1}{4\pi}\int_{|z|<1/\sqrt{2}} (\Delta \psi)(z)\,\mathrm{d}\mu.
    \]
    It is an elementary exercise to see that the right-hand side is exactly equal to $1$. If you are familiar with toric geometry, this is more or less trivial since
    \[
    \nabla_r \Bigl( \left(-\log (-r)\right)\lor (2+r) \Bigr)\bigl((-\infty,-\log 2)\bigr)= (0,1].
    \]
    Otherwise, just try to evaluate the integral using Green's identities. Therefore, \eqref{eq:Pomegavarphi0temp2} is proved and our assertion \eqref{eq:Pomegavarphi0temp1} follows.
\end{example}

Next we give a criterion on when the rooftop operator is not identically $-\infty$.
\begin{proposition}\label{prop:landfinitecond1}
    Assume that $\varphi,\psi\in \PSH(X,\theta)$ and
    \begin{equation}\label{eq:varphimassppsimass_greater}
        \int_X \theta_{\varphi}^n+\int_X \theta_{\psi}^n>\int_X \theta_{\varphi\lor\psi}^n.
    \end{equation}
    Then $\varphi\land \psi\in \PSH(X,\theta)$.
\end{proposition}
Thanks to the monotonicity theorem \cref{thm:mono}, we may also replace $\varphi\lor \psi$ on the right-hand side of \eqref{eq:varphimassppsimass_greater} by any $\gamma\in \PSH(X,\theta)$ such that $\varphi\lor \psi\preceq \gamma$.
\begin{proof}
    Without loss of generality, we may assume that $\varphi,\psi\leq 0$. For simplicity, we write
    \[
    \eta=P_{\theta}[\varphi\lor \psi].
    \]
    Note that due to our assumption \eqref{eq:varphimassppsimass_greater}, $\varphi,\psi\in \PSH(X,\theta)_{>0}$, and hence so is $\varphi\lor \psi$ thanks to the monotonicity theorem \cref{thm:mono}.
    By \cref{cor:Psendspotentialtomodel}, $\eta$ is a model potential, and $\int_X\theta_\eta^n=\int_X\theta_{\varphi\lor\psi}^n>0$.

    Take $C>0$ large enough, so that
    \begin{equation}\label{eq:intsumintbiggerinttemp1}
    \int_{\{\varphi>\eta-C\}} \theta_{\varphi}^n+\int_{\{\psi>\eta-C\}} \theta_{\psi}^n>\int_X \theta_{\eta}^n.
    \end{equation}
    This is possible thanks to \itemref{itm:prop:npp1:4}.
    Fix $C'>C$.
    Write
    \[
    \gamma_{C'}\coloneqq \Bigl(\varphi\lor (\eta-C')\Bigr)\land \Bigl(\psi\lor (\eta-C')\Bigr).
    \]
    Then observe that $\gamma_{C'}\sim \eta$, and
    \[
    \inf_{C'>C}\gamma_{C'}=\varphi\land \psi.
    \]
    Suppose, for contradiction, that $\varphi\land\psi\equiv -\infty$. Then
    \[
    \adjustlimits\lim_{C'\to \infty} \sup_X\gamma_{C'}=-\infty.
    \]
    Thanks to \cref{prop:model_sup_prop}, for each $C'>C$,
    \[
    \sup_{X} \gamma_{C'}=\sup_{\{\eta\neq -\infty\}} (\gamma_{C'}-\eta),
    \]
    since $\eta$ is a model potential.
    It follows that
    \begin{equation}\label{eq:limsupgammametatemp1}
        \adjustlimits \lim_{C'\to \infty} \sup_{\{\eta\neq -\infty\}}  (\gamma_{C'}-\eta)=-\infty.
    \end{equation}
    In particular, we can take $C'$ large enough so that
    \[
    \gamma_{C'}\leq \eta-C.
    \]
    For such $C'$, we have
    \[
    \begin{aligned}
        \int_X \theta_{\gamma_{C'}}^n=&\int_{\{\gamma_{C'}\leq \eta-C\}} \theta_{\gamma_{C'}}^n \\
        \leq & \int_{\left\{\varphi\lor (\eta-C')\leq \eta-C\right\}} \theta_{\varphi\lor (\eta-C')}^n +\int_{\left\{\psi\lor (\eta-C')\leq \eta-C \right\}} \theta_{\psi\lor (\eta-C')}^n\\
        =& 2\int_X \theta_{\eta}^n-\int_{\{\varphi>\eta-C\}}\theta_{\varphi}^n-\int_{\{\psi>\eta-C\}}\theta_{\psi}^n\\
        <& \int_X \theta_{\eta}^n,
    \end{aligned}
    \]
    where the second line follows from \cref{lma:rooftopMA}, the fourth line follows from \eqref{eq:intsumintbiggerinttemp1}. This contradicts the fact that $\gamma_{C'}\sim \eta$ in view of \cref{thm:mono}.
\end{proof}

\begin{proposition}\label{prop:landpresmodel}
    Let $\varphi,\psi\in \PSH(X,\theta)$. Assume that
    \begin{equation}\label{eq:varphipsimodel_landpsh}
        \varphi=P_{\theta}[\varphi],\quad \psi=P_{\theta}[\psi],\quad \varphi\land \psi\not\equiv -\infty.
    \end{equation}
    Then
    \begin{equation}\label{eq:Pthetaphilandpsi}
        P_{\theta}[\varphi\land \psi]=\varphi\land \psi.
    \end{equation}
\end{proposition}
\begin{proof}
    Observe that by \itemref{itm:prop:Pconc:3},
    \[
        P_{\theta}[\varphi\land \psi]\leq P_{\theta}[\varphi]=\varphi,\quad P_{\theta}[\varphi\land \psi]\leq P_{\theta}[\psi]=\psi.
    \]
    So the $\leq$ direction in \eqref{eq:Pthetaphilandpsi} holds. The reverse direction is trivial.
\end{proof}

\begin{lemma}\label{lma:landcorrectmass}
    Let $\varphi,\psi\in \PSH(X,\theta)$. Assume that \eqref{eq:varphipsimodel_landpsh} holds and
    \[
    \int_X \theta_{\varphi\land \psi}^n>0.
    \]
    Suppose that $\gamma,\eta\in \PSH(X,\theta)$ and satisfy
    \[
    \gamma\preceq \varphi,\quad \eta\preceq \psi,\quad \int_X \theta_{\gamma}^n=\int_X \theta_{\varphi}^n,\quad \int_X \theta_{\eta}^n=\int_X \theta_{\psi}^n.
    \]
    Then $\gamma\land \eta\not\equiv -\infty$ and
    \begin{equation}\label{eq:thetagammaeta_temp1}
    \int_X \theta_{\gamma\land \eta}^n=\int_X \theta_{\varphi\land \psi}^n.
    \end{equation}
\end{lemma}
Later on in \cref{chap:comp}, we shall reformulate and generalize this lemma in terms of the $P$-preorder, which is a more natural language to express the result. See \cref{prop:rooftopprePequiv}.
\begin{proof}
    Choose constants $A,B>0$ such that
    \[
        \gamma-A\leq \varphi,\quad \eta-B\leq \psi.
    \]
    Since
    \[
        (\gamma-A)\land(\eta-B)\sim \gamma\land\eta,
    \]
    the non-pluripolar masses of these two rooftop potentials are equal by the monotonicity theorem \cref{thm:mono}. Hence it suffices to prove the lemma after replacing $\gamma,\eta$ by $\gamma-A,\eta-B$. We can therefore assume that
    \[
        \gamma\leq \varphi,\quad \eta\leq \psi.
    \]

    \textbf{Step~1}. We first show that $\gamma\land \psi\not\equiv -\infty$ and
    \begin{equation}\label{eq:theta_gammalandpsi_mass_temp1}
    \int_X \theta_{\gamma\land \psi}^n=\int_X \theta_{\varphi\land \psi}^n.
    \end{equation}

    For any $a>1$, we define
    \[
    \gamma_a\coloneqq P_{\theta}\Bigl( a\gamma+(1-a)\varphi \Bigr).
    \]
    When $a=1$, we simply write $\gamma_1=\gamma$.

    Thanks to \cref{cor:pathoenvelopeeqmass}, for any $a\geq 1$, we have $\gamma_a\leq \varphi$ and
    \[
    \int_X \theta_{\gamma_a}^n=\int_X \theta_{\varphi}^n.
    \]
    By our assumption, for any $a\geq 1$, we have
    \[
    \int_X \theta_{\gamma_a}^n+\int_X \theta_{\varphi\land \psi}^n> \int_X \theta_{\varphi}^n,\quad \gamma_a\preceq \varphi,\quad \varphi\land \psi\preceq \varphi.
    \]
    So \cref{prop:landfinitecond1} gives
    \[
    \gamma_a\land \varphi\land \psi=\gamma_a\land \psi\in \PSH(X,\theta).
    \]
    In particular, for $a=1$, we find $\gamma\land \psi\not\equiv -\infty$. It remains to verify \eqref{eq:theta_gammalandpsi_mass_temp1}.

    For $a>1$, we have
    \[
    \gamma\geq a^{-1}\gamma_a+\left(1-a^{-1}\right)\varphi,
    \]
    thanks to the definition of $\gamma_a$ (cf. \eqref{eq:P_path_everywhere_ineq}). Hence
    \[
    \gamma\land \psi\geq a^{-1}\left( \gamma_a\land \psi\right)+\left(1-a^{-1}\right)\left(\varphi\land \psi\right).
    \]
    Therefore, by the monotonicity theorem \cref{thm:mono},
    \[
    \int_X \theta_{\gamma\land \psi}^n\geq \left(1-a^{-1}\right)^n \int_X \theta_{\varphi\land \psi}^n.
    \]
    Letting $a\to \infty$ and taking \cref{thm:mono} into account, \eqref{eq:theta_gammalandpsi_mass_temp1} follows.

    \textbf{Step~2}. We complete the proof.

    For any $a>1$, we let
    \[
    \eta_a\coloneqq P_{\theta}\Bigl( a\eta+(1-a)\psi \Bigr).
    \]
    When $a=1$, we set $\eta_1=\eta$.
    Thanks to \cref{cor:pathoenvelopeeqmass}, for each $a\geq 1$, we have $\eta_a\leq \psi$ and
    \[
    \int_X \theta_{\eta_a}^n=\int_X \theta_{\psi}^n.
    \]
    From Step~1, we have
    \[
    \int_X \theta_{\gamma\land\psi}^n+\int_X \theta_{\eta_a}^n=\int_X \theta_{\varphi\land \psi}^n+\int_X \theta_{\psi}^n>\int_X \theta_{\psi}^n.
    \]
    Applying \cref{prop:landfinitecond1} again, we find
    \[
    \gamma\land\psi\land \eta_a =\gamma\land \eta_a\in \PSH(X,\theta).
    \]
    When $a=1$, we find $\gamma\land \eta\in \PSH(X,\theta)$. It remains to prove \eqref{eq:thetagammaeta_temp1}.

    By definition of $\eta_a$, we have
    \[
    \eta\geq a^{-1}\eta_a+\left(1-a^{-1}\right)\psi.
    \]
    Therefore,
    \[
    \gamma\land \eta\geq a^{-1} (\gamma\land \eta_a)+\left(1-a^{-1}\right) \left( \gamma\land\psi \right).
    \]
    By the monotonicity theorem \cref{thm:mono},
    \[
    \int_X \theta_{\gamma\land \eta}^n\geq \left(1-a^{-1}\right)^n\int_X \theta_{\gamma\land \psi}^n.
    \]
    Letting $a\to \infty$ and applying \cref{thm:mono} again, we arrive at \eqref{eq:thetagammaeta_temp1}.
\end{proof}

There is an interesting diamond-like inequality regarding the non-pluripolar masses:
\begin{theorem}\label{thm:diamond}
    Assume that $\varphi,\psi\in \PSH(X,\theta)$ and $\varphi\land \psi\in \PSH(X,\theta)$. Then
    \begin{equation}\label{eq:diamond}
        \int_X \theta_{\varphi}^n+\int_X \theta_{\psi}^n\leq \int_X \theta_{\varphi\lor \psi}^n+\int_X \theta_{\varphi\land \psi}^n.
    \end{equation}
\end{theorem}
\begin{proof}
    We may assume that $\varphi,\psi\in \PSH(X,\theta)_{>0}$, as otherwise \eqref{eq:diamond} follows immediately from the monotonicity theorem \cref{thm:mono}.

    \textbf{Step~1}. We claim that it suffices to prove \eqref{eq:diamond} with $P_{\theta}[\varphi]$ and $P_{\theta}[\psi]$ in place of $\varphi$ and $\psi$.

    In fact, as $\varphi\land \psi\not \equiv -\infty$, we also have $P_{\theta}[\varphi]\land P_{\theta}[\psi]\not\equiv -\infty$. Moreover,
    \[
    \int_X \theta_{P_{\theta}[\varphi]}^n=\int_X \theta_{\varphi}^n,\quad \int_X \theta_{P_{\theta}[\psi]}^n=\int_X \theta_{\psi}^n
    \]
    by \cref{prop:Ppresmass}. On the other hand, as $C\to \infty$, the potentials
    \[
        \Bigl( (\varphi+C)\land V_{\theta} \Bigr) \lor \Bigl( (\psi+C)\land V_{\theta} \Bigr)
    \]
    converge to $P_{\theta}[\varphi]\lor P_{\theta}[\psi]$ almost everywhere. Therefore, thanks to \cref{cor:incseqnppcont} and \cref{thm:mono},
    \[
    \int_X \theta_{P_{\theta}[\varphi]\lor P_{\theta}[\psi]}^n=\int_X \theta_{\varphi\lor \psi}^n.
    \]
    Finally, we have
    \[
    \int_X \theta_{P_{\theta}[\varphi]\land P_{\theta}[\psi]}^n= \int_X \theta_{\varphi\land\psi}^n
    \]
    as well. When the left-hand side vanishes, this follows from \cref{thm:mono}. Otherwise, it follows from \cref{lma:landcorrectmass}.

    In particular, \eqref{eq:diamond} is equivalent to the corresponding result with $P_{\theta}[\varphi]$ and $P_{\theta}[\psi]$ in place of $\varphi$ and $\psi$.

    \textbf{Step~2}. We shall assume that $\varphi$ and $\psi$ are model potentials in the sequel.

    We write $\gamma=\varphi\lor \psi$ for simplicity. Note that $\gamma\leq 0$.

    For each $C>0$, we introduce
    \[
    \varphi_C=\varphi\lor (\gamma-C),\quad \psi_C=\psi\lor (\gamma-C).
    \]
    Note that $\varphi_C\sim \psi_C\sim \gamma$. Observe that \[
        \inf_{C>0} \varphi_C\land \psi_C=\varphi\land\psi.
    \]
    The $\geq$ direction is trivial, and the $\leq$ direction follows from the observation that the left-hand side is dominated by both $\varphi_C$ and $\psi_C$ for any $C>0$, and therefore by both $\varphi$ and $\psi$.

    For $C>0$, we compute
    \begin{equation}\label{eq:thetagammaminusvarphi_temp2}
    \begin{aligned}
    \theta_{\varphi_C}^n=& \mathds{1}_{\{\varphi>\gamma-C\}} \,\theta_{\varphi_C}^n+ \mathds{1}_{\{\varphi\leq \gamma-C\}} \,\theta_{\varphi_C}^n &\\
    =& \mathds{1}_{\{\varphi>\gamma-C\}} \,\theta_{\varphi}^n+ \mathds{1}_{\{\varphi\leq \gamma-C\}} \,\theta_{\varphi_C}^n & \textup{by \itemref{itm:prop:npp1:1}}\\
    =& \theta_{\varphi}^n+ \mathds{1}_{\{\varphi\leq \gamma-C\}} \,\theta_{\varphi_C}^n & \textup{by \cref{thm:Pvarphisupport}}.
    \end{aligned}
\end{equation}
    In the last line we use that $\varphi$ is model, so $\theta_\varphi^n$ is carried by $\{\varphi=0\}\subseteq\{\varphi>\gamma-C\}$.
    We also get a similar formula for $\theta_{\psi_C}^n$.

    Therefore, taking the monotonicity theorem \cref{thm:mono} into consideration, we find
    \begin{equation}\label{eq:intthetagammaminusvarphi_temp1}
    \int_X \theta_{\gamma}^n-\int_X \theta_{\varphi}^n=\int_{\{\varphi\leq \gamma-C\}} \theta_{\varphi_C}^n,\quad \int_X \theta_{\gamma}^n-\int_X \theta_{\psi}^n=\int_{\{\psi\leq \gamma-C\}} \theta_{\psi_C}^n.
    \end{equation}
    On the other hand, using \cref{lma:rooftopMA} and \eqref{eq:thetagammaminusvarphi_temp2}, we find
    \[
    \begin{aligned}
    \theta_{\varphi_C\land \psi_C}^n \leq & \mathds{1}_{\{\varphi_C\land \psi_C=\varphi_C\}} \,\theta_{\varphi_C}^n+\mathds{1}_{\{\varphi_C\land \psi_C=\psi_C\}} \,\theta_{\psi_C}^n \\
    \leq & \mathds{1}_{\{\varphi_C\land \psi_C=\varphi_C=0\}}\,\theta_{\varphi}^n+\mathds{1}_{\{\varphi_C\land \psi_C=\psi_C=0\}}\,\theta_{\psi}^n+\mathds{1}_{\{\varphi\leq \gamma-C\}}\,\theta_{\varphi_C}^n\\
    &+\mathds{1}_{\{\psi\leq \gamma-C\}}\,\theta_{\psi_C}^n.
    \end{aligned}
    \]
    In particular,
    \[
    \mathds{1}_{\{ \varphi_C\land \psi_C<0 \}}\, \theta_{\varphi_C\land \psi_C}^n \leq \mathds{1}_{\{\varphi\leq \gamma-C\}}\,\theta_{\varphi_C}^n+\mathds{1}_{\{\psi\leq \gamma-C\}}\,\theta_{\psi_C}^n.
    \]
    Integrating, we find
    \[
    \begin{aligned}
    \int_{\{ \varphi_C\land \psi_C=0 \}} \theta_{\varphi_C\land \psi_C}^n
    \geq & \int_X \theta_{\gamma}^n-\int_{\{\varphi\leq \gamma-C\}} \theta_{\varphi_C}^n-\int_{\{\psi\leq \gamma-C\}} \theta_{\psi_C}^n\\
    =& \int_X \theta_{\varphi}^n+\int_X \theta_{\psi}^n-\int_X \theta_{\gamma}^n,
    \end{aligned}
    \]
    where on the first line we used the fact that $\varphi_C\land \psi_C\sim \gamma$, while on the second line, we used \eqref{eq:intthetagammaminusvarphi_temp1}.

    Letting $C\to\infty$, and using \cref{cor:intcontlimsup}, we conclude \eqref{eq:diamond}.
\end{proof}

For later use, we also introduce a generalized envelope operator:
\begin{definition}\label{def:Pthetavarphif}
Given $\varphi\in\PSH(X,\theta)$ and a function $f\colon X\to[-\infty,\infty]$, we define
\[
P_\theta[\varphi](f)\coloneqq \uscsup \bigl\{\eta\in\PSH(X,\theta): \eta\leq f\textup{ quasi-everywhere},\ \eta\preceq\varphi\bigr\}.
\]
\index{$P_\theta[\varphi](f)$}
\end{definition}
In general, we cannot replace $\eta\leq f$ quasi-everywhere by everywhere.

\begin{example}\label{ex:envelopes:L625}
    When $f=0$, we have $P_\theta[\varphi](f)=P_{\theta}[\varphi]$.

    When $\varphi=V_{\theta}$, we have $P_\theta[V_{\theta}](f)=P_{\theta}(f)$.

    So \cref{def:Pthetavarphif} is a common generalization of \cref{def:Penv} and \cref{def:Pthetaf}.
\end{example}

\begin{proposition}\label{prop:DNTenvelope2}
    Let $\varphi\in\PSH(X,\theta)$ and $f\in C^{1,\overline{1}}(X)$ with $\varphi\leq f$. Then
    \[
    \theta^n_{P_\theta[\varphi](f)} = \mathds{1}_{\left\{P_\theta[\varphi](f)=f\right\}}\,\theta^n_f.
    \]
    In particular,
    \[
    \theta^n_{P_\theta[\varphi]} = \mathds{1}_{\left\{P_\theta[\varphi]=0\right\}}\,\theta^n.
    \]
\end{proposition}
This is a vast generalization of \cref{thm:Pvarphisupport}.
\begin{proof}
    By definition, we have
    \[
    P_\theta[\varphi](f)=P_\theta[\varphi]\Bigl(P_\theta(f)\Bigr)
    \]
    and $\varphi\leq P_\theta(f)$. By \cref{cor:Pvarphisupport} and \cref{cor:DNTenvelope1},
    \[
        \theta^n_{P_\theta[\varphi](f)} \leq \mathds{1}_{\left\{P_\theta[\varphi](f)=P_\theta(f)\right\}}\,\theta^n_{P_\theta(f)} = \mathds{1}_{\left\{P_\theta[\varphi](f)=f\right\}}\,\theta^n_f.
    \]
    The reverse inequality follows immediately from \cref{thm:DNTcontact}.
\end{proof}

\begin{proposition}\label{prop:DNTenvelope3}
    Let $\varphi\in\PSH(X,\theta)$ and $f\in C^{1,\overline{1}}(X)$ with $\varphi\leq f$. Then
    \[
    \begin{aligned}
        \mathds{1}_{\left\{\varphi=P_\theta(f)\right\}}\,\theta^n_{P_\theta(f)} =& \mathds{1}_{\left\{\varphi=P_\theta(f)\right\}}\,\theta^n_\varphi, \\
        \mathds{1}_{\left\{\varphi=P_\theta[\varphi](f)\right\}}\,\theta^n_{P_\theta[\varphi](f)} =& \mathds{1}_{\{\varphi=P_\theta[\varphi](f)\}}\,\theta^n_\varphi.
    \end{aligned}
    \]
\end{proposition}
\begin{proof}
    By \cref{cor:DNTenvelope1} and \cref{thm:DNTcontact}, we have
    \[
        \mathds{1}_{\{\varphi=P_\theta(f)\}}\,\theta^n_{P_\theta(f)}=\mathds{1}_{\{\varphi=P_\theta(f)=f\}}\,\theta^n_{f}=\mathds{1}_{\{\varphi=P_\theta(f)=f\}}\,\theta^n_{\varphi}\leq \mathds{1}_{\{\varphi=P_\theta(f)\}}\,\theta^n_{\varphi}.
    \]
    The reverse inequality follows from \cref{lma:lorMAineqe}. The first equality follows.

    As for the second, we compute using \cref{prop:DNTenvelope2} and \cref{thm:DNTcontact} that
    \[
    \begin{split}
    \mathds{1}_{\left\{\varphi=P_\theta[\varphi](f)\right\}}\,\theta^n_{P_\theta[\varphi](f)}=\mathds{1}_{\left\{\varphi=P_\theta[\varphi](f)=f\right\}}\,\theta^n_{f}\\
    =\mathds{1}_{\left\{\varphi=P_\theta[\varphi](f)=f\right\}}\,\theta_{\varphi}^n\leq \mathds{1}_{\{\varphi=P_\theta[\varphi](f)\}}\,\theta^n_\varphi.
    \end{split}
    \]
    Since $\varphi\leq P_{\theta}[\varphi](f)$, \cref{lma:lorMAineqe} gives the opposite inequality. The second equality follows.
\end{proof}

\subsection{Properties of the \texorpdfstring{$P$}{P}-envelope}
\label{subsec:Penv_prop}
Let $\theta,\theta_1,\theta_2$ be smooth closed real $(1,1)$-forms on $X$.

\begin{proposition}\label{prop:Penvbimero}
    Let $\pi\colon Y\rightarrow X$ be a proper bimeromorphic morphism from a Kähler manifold $Y$ to $X$. Then for any $\varphi\in \PSH(X,\theta)$, we have
    \[
        P_{\pi^*\theta}[\pi^*\varphi]=\pi^*P_{\theta}[\varphi].
    \]
    In particular, a potential $\varphi\in \PSH(X,\theta)_{>0}$ is model if and only if $\pi^*\varphi\in \PSH(Y,\pi^*\theta)_{>0}$ is model.
\end{proposition}
\begin{proof}
    By \cref{prop:PSHpullbij}, every $\chi\in\PSH(Y,\pi^*\theta)$ is of the form $\chi=\pi^*\psi$ for a unique $\psi\in\PSH(X,\theta)$. Since $\pi$ is surjective, $\chi\leq 0$ if and only if $\psi\leq 0$. Moreover,
    \[
        \chi\preceq \pi^*\varphi \Longleftrightarrow \psi\preceq \varphi.
    \]
    Thus the candidates defining $P_{\pi^*\theta}[\pi^*\varphi]$ are exactly the pull-backs of the candidates defining $P_\theta[\varphi]$, and hence
    \[
        P_{\pi^*\theta}[\pi^*\varphi]=\pi^*P_\theta[\varphi].
    \]
    Finally, the non-pluripolar mass is preserved under bimeromorphic pull-back by \cref{prop:nppmassinv}. Therefore $\varphi$ is model of positive mass if and only if $\pi^*\varphi$ is.
\end{proof}

We have the following concavity property of the $P$-envelope.
\begin{proposition}\label{prop:Pconc}
    \leavevmode
    \begin{enumerate}
        \item\label{itm:prop:Pconc:1} Suppose that $\varphi\in \PSH(X,\theta)$ and $\lambda\in \mathbb{R}_{>0}$. Then
            \[
                P_{\lambda\theta}[\lambda\varphi]=\lambda P_{\theta}[\varphi].
            \]
        \item\label{itm:prop:Pconc:2} Suppose that $\varphi_1\in \PSH(X,\theta_1)$ and $\varphi_2\in \PSH(X,\theta_2)$. Then
            \[
                P_{\theta_1+\theta_2}[\varphi_1+\varphi_2]\geq P_{\theta_1}[\varphi_1]+P_{\theta_2}[\varphi_2].
            \]
        \item\label{itm:prop:Pconc:3} Suppose that $\varphi,\psi\in \PSH(X,\theta)$ and $\varphi\preceq \psi$. Then
            \[
                P_{\theta}[\varphi]\leq P_{\theta}[\psi].
            \]
    \end{enumerate}
\end{proposition}
\begin{proof}
    \ref{itm:prop:Pconc:1} This follows by rescaling the candidates in the definition of the $P$-envelope.

    \ref{itm:prop:Pconc:2} Suppose that $\psi_1\in \PSH(X,\theta_1)$ and $\psi_2\in \PSH(X,\theta_2)$ satisfy
    \[
        \psi_i\leq 0,\quad \psi_i\preceq \varphi_i
    \]
    for $i=1,2$. Then
    \[
        \psi_1+\psi_2\leq 0,\quad \psi_1+\psi_2\preceq \varphi_1+\varphi_2.
    \]
    It follows from \eqref{eq:Pthetavarphi} that
    \[
        \psi_1+\psi_2\leq P_{\theta_1+\theta_2}[\varphi_1+\varphi_2].
    \]
    Since $\psi_1$ and $\psi_2$ are arbitrary, we conclude.

    \ref{itm:prop:Pconc:3} This follows by inclusion of the defining candidate sets.
\end{proof}

\begin{proposition}\label{prop:decseqmodel}
    Let $(\varphi_j)_{j\in I}$ be a decreasing net of potentials in $\PSH(X,\theta)$ satisfying $P_{\theta}[\varphi_j]=\varphi_j$ for each $j\in I$. Set $\varphi\coloneqq \inf_{j\in I} \varphi_j$. Then $P_{\theta}[\varphi]=\varphi$.
\end{proposition}
\begin{proof}
    Since $P_{\theta}[\varphi_j]=\varphi_j$, we have $\sup_X\varphi_j=0$. By \cref{cor:decreasingnetL1conv}, we have $\varphi\in \PSH(X,\theta)$, and $\varphi\leq 0$. Therefore, for each $j\in I$,
    \[
        \varphi\leq P_{\theta}[\varphi]\leq P_{\theta}[\varphi_j] =\varphi_j.
    \]
    Therefore, $\varphi=P_{\theta}[\varphi]$.
\end{proof}

\begin{proposition}\label{prop:vol_limit_model}
    Let $(\epsilon_i)_{i\in I}$ be a decreasing net in $\mathbb{R}_{\geq 0}$ with limit $0$. Take a Kähler form $\omega$ on $X$.
    Consider a decreasing net $(\varphi_i)_{i\in I}$ with $\varphi_i\in \PSH(X,\theta+\epsilon_i \omega)$ satisfying
    \begin{equation}\label{eq:Palmostmodeltemp}
        P_{\theta+\epsilon_i\omega}\left[\varphi_i\right]=\varphi_i
    \end{equation}
    for every $i\in I$. Set $\varphi\coloneqq\inf_{i\in I}\varphi_i$. Then $P_{\theta}[\varphi]=\varphi$, and
    \begin{equation}\label{eq:massmodeldec}
    \lim_{i\in I}\int_X \left(\theta + \epsilon_i \omega\right)^n_{\varphi_i}=\int_X \theta_{\varphi}^n.
    \end{equation}
    Moreover, if $\int_X \theta_{\varphi}^n>0$, then for any prime divisor $E$ over $X$, we have
    \begin{equation}\label{eq:Lelongcontdecseq}
    \lim_{i\in I} \nu\left(\varphi_i,E\right)=\nu(\varphi,E).
    \end{equation}
\end{proposition}
Note that both \eqref{eq:massmodeldec} and \eqref{eq:Lelongcontdecseq} fail without the assumption \eqref{eq:Palmostmodeltemp}.

\begin{proof}
        We first observe that $\varphi\in \PSH(X,\theta)$. In fact, for each $i_0\in I$, the tail $(\varphi_i)_{i\geq i_0}$ lies in $\PSH(X,\theta+\epsilon_{i_0}\omega)$ and decreases to $\varphi$. Since $\sup_X\varphi_i=0$ by \eqref{eq:Palmostmodeltemp}, the limit is not identically $-\infty$. Thus $\varphi\in \PSH(X,\theta+\epsilon_{i_0}\omega)$ for every $i_0\in I$. Since $\epsilon_{i_0}\to 0$, it follows that $\varphi\in \PSH(X,\theta)$.
        For every $i\in I$ we have
        \[
        P_{\theta}[\varphi]\leq P_{\theta+\epsilon_i\omega}[\varphi]\leq P_{\theta+\epsilon_i\omega}[\varphi_i]=\varphi_i.
        \]
        Therefore, $P_{\theta}[\varphi]\leq \varphi$. The reverse inequality is trivial. Hence $P_{\theta}[\varphi]=\varphi$.

        By the monotonicity theorem \cref{thm:mono} and the multi-linearity of the non-pluripolar product \itemref{itm:prop:npp1:6}, we have
        \[
            \varliminf_{i\in I} \int_X \left(\theta + \epsilon_i \omega\right)^n_{\varphi_i} \geq  \varliminf_{i\in I} \int_X \left(\theta + \epsilon_i \omega\right)^n_{\varphi}  = \int_X \theta_{\varphi}^n.
        \]
        We now argue the reverse inequality.

        Fix $i_0\in I$. Since passing to the tail $i\geq i_0$ does not change the limsup, we have
        \[
            \begin{aligned}
                \varlimsup_{i\in I} \int_X \left(\theta + \epsilon_i \omega\right)^n_{\varphi_i}=&\varlimsup_{i\geq i_0} \int_{\left\{\varphi_i =0\right\}} \left(\theta + \epsilon_i \omega\right)^n_{\varphi_i} \\
                \leq & \varlimsup_{i\geq i_0} \int_{\left\{\varphi_i =0\right\}} \left(\theta + \epsilon_{i_0} \omega\right)^n_{\varphi_i}\\
                \leq & \int_{\left\{\varphi =0\right\}} \left(\theta + \epsilon_{i_0} \omega\right)^n_{\varphi},
            \end{aligned}
        \]
    where in the first line we used \eqref{eq:Palmostmodeltemp} and \cref{thm:Pvarphisupport}, and in the last line we used the fact that $\varphi_i \searrow \varphi$ on the tail $i\geq i_0$, and on this tail, $\varphi_i \to \varphi$ in $L^1$ as well by \cref{cor:decreasingnetL1conv}, and therefore, \cref{cor:intcontlimsup} is applicable. Taking the limit with respect to $i_0$, we arrive at the desired conclusion:
    \[
    \varlimsup_{i\in I} \int_{X} \left(\theta + \epsilon_i \omega \right)^n_{\varphi_i} \leq \varliminf_{i_0\in I}\int_{\{\varphi =0\}} \left(\theta + \epsilon_{i_0} \omega \right)^n_{\varphi}  =  \int_{\{\varphi =0\}} \theta^n_{\varphi} \leq \int_{X} \theta^n_{\varphi},
    \]
    where the equality follows from the multi-linearity of the non-pluripolar product \itemref{itm:prop:npp1:6} and the last inequality follows from \cref{thm:Pvarphisupport}.
    This finishes the proof of \eqref{eq:massmodeldec}.

    It remains to argue \eqref{eq:Lelongcontdecseq}.
    Fix $i_0\in I$. By \cref{prop:Lelongnumberupperbound}, there is a constant $B>0$ such that
    \[
    \nu(\psi,E)\leq B,\quad \psi\in \PSH(X,\theta+\epsilon_{i_0}\omega).
    \]
    By \cref{lma:pathoenvelope} and \eqref{eq:massmodeldec}, for any $\delta \in (0,1)$ there exists $i_{\delta}\geq i_0$ such that for all $i\geq i_{\delta}$ there exists $\psi_i \in \PSH(X,\theta+\epsilon_i\omega)$ satisfying $(1-\delta) \varphi_i + \delta \psi_i \leq \varphi$. This implies that for $i\geq i_{\delta}$ we have
    \[
    (1-\delta)\nu(\varphi_i,E) + \delta\nu(\psi_i,E) \geq \nu(\varphi,E) \geq \nu(\varphi_i,E).
    \]
    Since $\psi_i\in \PSH(X,\theta+\epsilon_{i_0}\omega)$ for $i\geq i_{\delta}$, we get
    \[
        \varliminf_{i\in I}\nu(\varphi_i,E)\geq \frac{\nu(\varphi,E)-\delta B}{1-\delta}.
    \]
    Letting $\delta\searrow 0$ gives $\varliminf_{i\in I}\nu(\varphi_i,E)\geq \nu(\varphi,E)$. The reverse inequality follows from $\varphi\leq \varphi_i$, and hence we conclude \eqref{eq:Lelongcontdecseq}.
\end{proof}

\begin{corollary}\label{cor:varphiperturbtheta}
    Let $\varphi\in \PSH(X,\theta)_{>0}$ be a model potential. Let $\omega$ be a Kähler form on $X$.
    Then
    \[
    \varphi=\inf_{\epsilon>0} P_{\theta+\epsilon\omega}[\varphi].
    \]
\end{corollary}
\begin{proof}
    Set
    \[
        \eta\coloneqq\inf_{\epsilon>0}P_{\theta+\epsilon\omega}[\varphi].
    \]
    Since $\varphi$ is model, $\varphi=P_{\theta}[\varphi]\leq P_{\theta+\epsilon\omega}[\varphi]$ for every $\epsilon>0$, hence $\varphi\leq\eta\leq 0$.

    Next we prove the reverse direction.
    Applying \cref{prop:vol_limit_model} to the decreasing net $(P_{\theta+\epsilon\omega}[\varphi])_{\epsilon>0}$ and using \cref{prop:Ppresmass}, we get
    \[
        P_\theta[\eta]=\eta,\quad \int_X\theta_\eta^n=\int_X\theta_\varphi^n.
    \]
    Thus $\eta$ is a candidate in the characterization \cref{thm:Pvarphidiffdef} of $P_\theta[\varphi]$. Hence $\eta\leq P_\theta[\varphi]=\varphi$.
\end{proof}

\begin{proposition}\label{prop:incnetmodel}
    Let $(\varphi_i)_{i\in I}$ be an increasing net of potentials in $\PSH(X,\theta)_{>0}$ uniformly bounded from above. Let $\varphi\coloneqq \uscsup[i\in I]\varphi_i$. Then
    \[
    \uscsup[i\in I] P_{\theta}[\varphi_i]=P_{\theta}[\varphi].
    \]

    In particular, if $\varphi_i$ is model for all $i\in I$, then so is $\varphi$.
\end{proposition}
\begin{proof}
We may assume that $I$ is infinite since otherwise, there is nothing to prove.

    We write
    \[
    \eta\coloneqq \uscsup[i\in I] P_{\theta}[\varphi_i].
    \]
    Then it is clear that $\eta\leq P_{\theta}[\varphi]$.

    By \cref{cor:incseqnppcont}, we have
    \[
        \lim_{i\in I} \int_X\theta_{\varphi_i}^n=\int_X\theta_{\varphi}^n>0.
    \]
    We may assume that for each $i\in I$, $\int_X\theta_{\varphi_i}^n<\int_X\theta_{\varphi}^n$, since otherwise $P_{\theta}[\varphi_i]=P_{\theta}[\varphi]$ by \cref{thm:Pvarphidiffdef}, and the same holds for all $j\geq i$; in that case there is nothing to prove.

    So by \cref{lma:pathoenvelope}, we can find a decreasing net $\epsilon_i\searrow 0$ ($i\in I$) with $\epsilon_i\in (0,1)$ and $\psi_i\in \PSH(X,\theta)$ ($i\in I$) such that for all $i\in I$,
    \[
        (1-\epsilon_i) \varphi+\epsilon_i\psi_i\leq \varphi_i.
    \]
    We may, for example, take
    \begin{equation}\label{eq:epsilon_inetchoice}
        \epsilon_i^{-1}=\frac{1}{2}\left( 1+\left(\frac{\int_X \theta_{\varphi}^n}{\int_X \theta_{\varphi}^n-\int_X \theta_{\varphi_i}^n}\right)^{1/n} \right).
    \end{equation}

    By the concavity property of the envelope \cref{prop:Pconc}, we have
    \[
        P_{\theta}[\varphi]+\epsilon_i P_{\theta}[\psi_i]\leq (1-\epsilon_i) P_{\theta}[\varphi]+\epsilon_i P_{\theta}[\psi_i]\leq P_{\theta}[\varphi_i]\leq\eta.
    \]
    Each $P_{\theta}[\psi_i]$ has supremum $0$, hence the family is bounded in $L^1$ by \cref{prop:L1compa}. Thus $\epsilon_iP_{\theta}[\psi_i]\to 0$ in $L^1$. The previous inequality then gives $P_{\theta}[\varphi]\leq \eta$ almost everywhere, and hence everywhere by \cref{prop:pshfuncdetdense}.
\end{proof}

\subsection{Relative full mass classes}\label{subsec:fullmass}
Let $\theta$ be a smooth closed real $(1,1)$-form on $X$ representing a big cohomology class. Fix a model potential $\phi\in \PSH(X,\theta)_{>0}$.

\begin{definition}\label{def:envelopes:L906}
    We define
    \[
        \begin{aligned}
        \PSH(X,\theta;\phi)\coloneqq &\left\{\eta\in \PSH(X,\theta): \eta\preceq \phi\right\},\\
        \mathcal{E}^{\infty}(X,\theta;\phi)\coloneqq &\left\{\eta\in \PSH(X,\theta): \eta\sim \phi\right\},\\
        \mathcal{E}(X,\theta;\phi)\coloneqq &\left\{\eta\in \PSH(X,\theta;\phi): \int_X \theta_{\eta}^n=\int_X \theta_{\phi}^n\right\},\\
        \mathcal{E}^1(X,\theta;\phi)\coloneqq &\left\{\eta\in \mathcal{E}(X,\theta;\phi): \int_X |\phi-\eta|\,\theta_{\eta}^n<\infty\right\}.
        \end{aligned}
    \]
    \index{$\PSH(X,\theta;\phi)$}\index{$\mathcal{E}^{\infty}(X,\theta;\phi)$}\index{$\mathcal{E}(X,\theta;\phi)$} \index{$\mathcal{E}^1(X,\theta;\phi)$}
    Potentials in the last three classes are said to have \emph{relatively minimal singularities}\index{potential!potential with relatively minimal singularities}, \emph{full mass}\index{potential!potential with relative full mass} and \emph{finite energy}\index{potential!potential with relative finite energy} relative to $\phi$ respectively.
\end{definition}
We have the following inclusions:
\begin{equation}\label{eq:energyclassinc}
\mathcal{E}^{\infty}(X,\theta;\phi)\subseteq \mathcal{E}^1(X,\theta;\phi)\subseteq \mathcal{E}(X,\theta;\phi)\subseteq \PSH(X,\theta;\phi).
\end{equation}
The only non-trivial part is the first inclusion, which follows from the monotonicity theorem \cref{thm:mono}.

\begin{remark}\label{rmk:intwelldef}
    Note that this integral
    \[
    \int_X |\phi-\eta|\,\theta_{\eta}^n
    \]
    is defined: The locus where $\phi-\eta$ is undefined is a pluripolar set, while the product $\theta_{\eta}^n$ puts no mass on pluripolar sets (\cref{prop:npp1}).

    Similar remarks apply when we talk about similar integrals in the sequel.
\end{remark}

When $\phi=V_{\theta}$, we usually write
\[
    \begin{aligned}
        \mathcal{E}^{\infty}(X,\theta;V_{\theta})=& \mathcal{E}^{\infty}(X,\theta),\\
        \mathcal{E}(X,\theta;V_{\theta})=& \mathcal{E}(X,\theta),\\
        \mathcal{E}^{1}(X,\theta;V_{\theta})=& \mathcal{E}^{1}(X,\theta).
    \end{aligned}
\]
Potentials in the three classes are said to have \emph{minimal singularities}\index{potential!potential with minimal singularities}, \emph{full mass}\index{potential!potential with full mass} and \emph{finite energy}\index{potential!potential with finite energy} respectively. The relation \eqref{eq:energyclassinc} can be written as
\[
\mathcal{E}^{\infty}(X,\theta)\subseteq \mathcal{E}^{1}(X,\theta) \subseteq \mathcal{E}(X,\theta)
\]
in this case.

The $P$-envelope can be used to characterize the full mass classes:
\begin{proposition}\label{prop:fullmassP}
    Let $\varphi\in \PSH(X,\theta)$. Then the following are equivalent:
    \begin{enumerate}
        \item\label{itm:prop:fullmassP:1} $\varphi\in \mathcal{E}(X,\theta;\phi)$;
        \item\label{itm:prop:fullmassP:2} $P_{\theta}[\varphi]=\phi$.
    \end{enumerate}
\end{proposition}
\begin{proof}
    \ref{itm:prop:fullmassP:2} $\implies$ \ref{itm:prop:fullmassP:1}. Assume \ref{itm:prop:fullmassP:2}. By \cref{prop:Ppresmass},
    \[
    \int_X\theta_\varphi^n=\int_X\theta_{P_\theta[\varphi]}^n=\int_X\theta_\phi^n.
    \]
    Moreover, since the mass is positive, \cref{thm:Pvarphidiffdef} implies that
    \[
        \varphi\preceq P_\theta[\varphi]=\phi.
    \]
    Thus $\varphi\in\mathcal E(X,\theta;\phi)$.

    \ref{itm:prop:fullmassP:1} $\implies$ \ref{itm:prop:fullmassP:2}. Assume \ref{itm:prop:fullmassP:1}. Then $\varphi\in \PSH(X,\theta)_{>0}$.
    Note that $\phi$ is a candidate of $P_{\theta}[\varphi]$ as in \eqref{eq:Penvdef}. So $P_{\theta}[\varphi]\geq \phi$. The reverse inequality is trivial.
\end{proof}

We have the following comparison principle.
\begin{proposition}[Comparison principle]\label{prop:compa_princ}\index{theorem!comparison principle}
    Fix $j\in \{0,\ldots,n\}$.
    Let $\theta_1,\ldots,\theta_j$ be closed smooth real $(1,1)$-forms on $X$ and $\psi_i\in \PSH(X,\theta_i)$ for $i=1,\ldots,j$. Suppose that $\varphi,\psi\in \mathcal{E}(X,\theta;\phi)$. Then
    \begin{equation}\label{eq:comp_princi1}
    \int_{\{\varphi<\psi\}} \theta_{\varphi}^{n-j} \wedge \theta_{1,\psi_1}\wedge \cdots \wedge \theta_{j,\psi_j}\geq \int_{\{\varphi<\psi\}} \theta_{\psi}^{n-j} \wedge \theta_{1,\psi_1}\wedge \cdots \wedge \theta_{j,\psi_j}.
    \end{equation}
    As a consequence, we also have
    \begin{equation}\label{eq:comp_princi12}
    \int_{\{\varphi\leq \psi\}} \theta_{\varphi}^{n-j} \wedge \theta_{1,\psi_1}\wedge \cdots \wedge \theta_{j,\psi_j}\geq \int_{\{\varphi\leq\psi\}} \theta_{\psi}^{n-j} \wedge \theta_{1,\psi_1}\wedge \cdots \wedge \theta_{j,\psi_j}.
    \end{equation}
\end{proposition}
\begin{proof}
    Observe that $\varphi\lor \psi\in \mathcal{E}(X,\theta;\phi)$, as a consequence of the monotonicity theorem \cref{thm:mono}.
    Using the plurifine locality of the non-pluripolar product \itemref{itm:prop:npp1:1}, we compute
    \[
    \begin{aligned}
    &\int_X \theta_{\varphi\lor \psi}^{n-j}\wedge \theta_{1,\psi_1}\wedge \cdots \wedge \theta_{j,\psi_j}\\
    \geq & \int_{\{\varphi>\psi\} } \theta_{\varphi}^{n-j}\wedge \theta_{1,\psi_1}\wedge \cdots \wedge \theta_{j,\psi_j}+\int_{\{\varphi<\psi\} } \theta_{\psi}^{n-j}\wedge \theta_{1,\psi_1}\wedge \cdots \wedge \theta_{j,\psi_j}\\
    =& \int_X  \theta_{\varphi}^{n-j}\wedge \theta_{1,\psi_1}\wedge \cdots \wedge \theta_{j,\psi_j}-\int_{\{\varphi\leq \psi\} } \theta_{\varphi}^{n-j}\wedge \theta_{1,\psi_1}\wedge \cdots \wedge \theta_{j,\psi_j}\\
    &+\int_{\{\varphi<\psi\} } \theta_{\psi}^{n-j}\wedge \theta_{1,\psi_1}\wedge \cdots \wedge \theta_{j,\psi_j}\\
    =&\int_X  \theta_{\varphi\lor \psi}^{n-j}\wedge \theta_{1,\psi_1}\wedge \cdots \wedge \theta_{j,\psi_j}-\int_{\{\varphi\leq \psi\} } \theta_{\varphi}^{n-j}\wedge \theta_{1,\psi_1}\wedge \cdots \wedge \theta_{j,\psi_j}\\
    &+\int_{\{\varphi<\psi\} } \theta_{\psi}^{n-j}\wedge \theta_{1,\psi_1}\wedge \cdots \wedge \theta_{j,\psi_j},
    \end{aligned}
    \]
    where in the last step, we applied \cref{prop:Ppresmass} to $\varphi$ and $\varphi\lor \psi$. Therefore,
    \[
    \begin{aligned}
    \int_{\{\varphi\leq \psi\}}& \theta_{\varphi}^{n-j} \wedge \theta_{1,\psi_1}\wedge \cdots \wedge \theta_{j,\psi_j}\\
    &\geq \int_{\{\varphi<\psi\}} \theta_{\psi}^{n-j} \wedge \theta_{1,\psi_1}\wedge \cdots \wedge \theta_{j,\psi_j}.
    \end{aligned}
    \]
    Applying the preceding inequality to $\varphi+\epsilon$ in place of $\varphi$, and using $\theta_{\varphi+\epsilon}=\theta_\varphi$, we get
    \[
        \int_{\{\varphi+\epsilon\leq\psi\}}\theta_\varphi^{n-j}\wedge \theta_{1,\psi_1}\wedge\cdots\wedge\theta_{j,\psi_j} \geq \int_{\{\varphi+\epsilon<\psi\}}\theta_\psi^{n-j}\wedge \theta_{1,\psi_1}\wedge\cdots\wedge\theta_{j,\psi_j}.
    \]
    As $\epsilon\searrow 0$, by monotone convergence theorem, we obtain
    \[
    \int_{\{\varphi<\psi\}}\theta_\varphi^{n-j}\wedge \theta_{1,\psi_1}\wedge\cdots\wedge\theta_{j,\psi_j} \geq \int_{\{\varphi<\psi\}}\theta_\psi^{n-j}\wedge \theta_{1,\psi_1}\wedge\cdots\wedge\theta_{j,\psi_j},
    \]
    which proves \eqref{eq:comp_princi1}.

    Finally, the \eqref{eq:comp_princi12} follows by subtracting \eqref{eq:comp_princi1} (with the role of $\varphi$ and $\psi$ switched) from the total mixed mass equality
    \[
        \int_X\theta_\varphi^{n-j}\wedge \theta_{1,\psi_1}\wedge\cdots\wedge\theta_{j,\psi_j}= \int_X\theta_\psi^{n-j}\wedge \theta_{1,\psi_1}\wedge\cdots\wedge\theta_{j,\psi_j},
    \]
    which itself follows from \cref{prop:Ppresmass} and \cref{prop:fullmassP}.
\end{proof}

The full mass potentials are essential in resolving the Monge--Ampère equations. We recall the following two Aubin--Yau type theorems. Recall that a measure $\mu$ on $X$ is \emph{non-pluripolar}\index{non-pluripolar measure} if it is a Radon measure and $\mu(K)=0$ for each pluripolar set $K\subseteq X$.\footnote{Here we identify as usual a Radon measure with its completion.}
\begin{theorem}\label{thm:Yau_rel}
    Let $\mu$ be a non-pluripolar measure on $X$ with $\mu(X)=\int_X \theta_{\phi}^n$. Then there is a unique $\varphi\in \mathcal{E}(X,\theta;\phi)$ such that
    \[
    \theta_{\varphi}^n=\mu,\quad \sup_X \varphi=0.
    \]
\end{theorem}

\begin{theorem}\label{thm:AY_rel}
    Fix $\lambda>0$. Let $\mu$ be a non-pluripolar measure on $X$ with $\mu(X)>0$. Then there is a unique $\varphi\in \mathcal{E}(X,\theta;\phi)$ such that
    \[
    \theta_{\varphi}^n=\mathrm{e}^{\lambda\varphi}\mu.
    \]
    Furthermore, for any $\psi\in \mathcal{E}(X,\theta;\phi)$ satisfying
    \[
    \theta_{\psi}^n\geq \mathrm{e}^{\lambda\psi}\mu,
    \]
    we have $\varphi\geq \psi$.
\end{theorem}
The proofs are beyond the scope of this book; we refer to \cite{DDNL19log} and \cite[Proposition~17.4.6 and Theorems~17.4.11, 17.4.13--17.4.14]{FPT}.

In order to handle the finite energy classes, it is convenient to introduce the following quantity:
\begin{definition}\label{def:MAenergy}
    We define the \emph{Monge--Amp\`ere energy}\index{Monge--Amp\`ere energy} $E^{\phi}_{\theta}\colon \mathcal{E}^{\infty}(X,\theta;\phi)\rightarrow \mathbb{R}$\index{$E^{\phi}_{\theta}$} as follows
    \begin{equation}\label{eq:Edefbdd}
        E^{\phi}_{\theta}(\varphi)\coloneqq \frac{1}{n+1}\sum_{j=0}^n \int_X (\varphi-\phi)\, \theta_{\varphi}^j\wedge \theta_{\phi}^{n-j}.
    \end{equation}

    More generally, we extend $E_{\theta}^{\phi}$ to a functional $E^{\phi}_{\theta}\colon \PSH(X,\theta;\phi)\rightarrow [-\infty,\infty)$ as follows
    \begin{equation}\label{eq:Eextendgeneral}
        E^{\phi}_{\theta}(\varphi)\coloneqq \inf\left\{E^{\phi}_{\theta}(\psi):\psi\in \mathcal{E}^{\infty}(X,\theta;\phi),\varphi\leq \psi\right\}.
    \end{equation}

     We write $E_{\theta}$\index{$E_{\theta}$} instead of $E_{\theta}^{\phi}$ when $\phi=V_{\theta}$.
\end{definition}
The set in \eqref{eq:Eextendgeneral} is non-empty, since $\phi+\sup_X \varphi$ is a candidate. Therefore, $E^{\phi}_{\theta}(\varphi)<\infty$.

Note that
\begin{equation}\label{eq:Ephiaddcst}
E^{\phi}_{\theta}(\varphi+C)=E^{\phi}_{\theta}(\varphi)+C\int_X \theta_{\phi}^n
\end{equation}
for any $\varphi\in \PSH(X,\theta;\phi)$ and $C\in \mathbb{R}$.

\begin{lemma}\label{lma:EmonoEinfty}
    The functional $E_{\theta}^{\phi}\colon \mathcal{E}^{\infty}(X,\theta;\phi)\rightarrow \mathbb{R}$ is increasing. In particular, the extended definition \eqref{eq:Eextendgeneral} agrees with \eqref{eq:Edefbdd} on $\mathcal{E}^{\infty}(X,\theta;\phi)$, and is increasing as well.
\end{lemma}
\begin{proof}
    Let $\varphi,\psi\in \mathcal{E}^{\infty}(X,\theta;\phi)$, we have
    \[
    \begin{aligned}
    &(n+1)E^{\phi}_{\theta}(\varphi)-(n+1)E^{\phi}_{\theta}(\psi)-\sum_{j=0}^n \int_X (\varphi-\psi)\,\theta_{\varphi}^j\wedge \theta_{\psi}^{n-j} \\
    =& \sum_{j=0}^n \int_X (\varphi-\phi)\, \theta_{\varphi}^j\wedge \theta_{\phi}^{n-j}-\sum_{j=0}^n \int_X (\psi-\phi)\, \theta_{\psi}^j\wedge \theta_{\phi}^{n-j}-\sum_{j=0}^n \int_X (\varphi-\psi)\,\theta_{\varphi}^j\wedge \theta_{\psi}^{n-j}\\
    =& \sum_{j=0}^n \int_X (\varphi-\psi)\,\left( \theta_{\varphi}^j\wedge \theta_{\phi}^{n-j}-\theta_{\varphi}^j\wedge \theta_{\psi}^{n-j} \right)\\
    &+\sum_{j=0}^n \int_X (\psi-\phi)\,\left( \theta_{\varphi}^j\wedge \theta_{\phi}^{n-j}-\theta_{\psi}^j\wedge \theta_{\phi}^{n-j} \right)\\
    =& \sum_{\substack{j+a+b=n-1,\\ j,a,b\geq 0}} (n-j)\int_X (\psi-\phi)\,\left( \theta_{\phi}^a\wedge \theta_{\psi}^{b+1}\wedge \theta_{\varphi}^j- \theta_{\phi}^a\wedge \theta_{\psi}^{b}\wedge \theta_{\varphi}^{j+1}\right)\\
    &+\sum_{\substack{j+a+b=n-1,\\ j,a,b\geq 0}} (n-j)\int_X (\psi-\phi)\,\left( \theta_{\phi}^j\wedge \theta_{\varphi}^{a+1}\wedge \theta_{\psi}^b- \theta_{\phi}^j\wedge \theta_{\varphi}^{a}\wedge \theta_{\psi}^{b+1}\right)\\
    =&0,
    \end{aligned}
    \]
    where the third equality follows from the integration by parts formula. See \cite{Xia19b, Lu21} and \cite[Theorem~17.1.3]{FPT} for the proof in the context of non-pluripolar products.

    In other words,
    \begin{equation}\label{eq:int_by_part_E_temp1}
    E^{\phi}_{\theta}(\varphi)-E^{\phi}_{\theta}(\psi)=\frac{1}{n+1}\sum_{j=0}^n \int_X(\varphi-\psi)\,\theta_{\varphi}^j\wedge \theta_{\psi}^{n-j}.
    \end{equation}
    The monotonicity follows.
\end{proof}

\begin{lemma}\label{lma:comp_int_ene}
    Let $\varphi,\psi\in \mathcal{E}^{\infty}(X,\theta;\phi)$. Then
    \begin{equation}\label{eq:comp_int_ene}
    \int_X (\varphi-\psi)\,\theta_{\psi}^n \geq E_{\theta}^{\phi}(\varphi)-E_{\theta}^{\phi}(\psi)\geq \int_X (\varphi-\psi)\,\theta_{\varphi}^n.
    \end{equation}
\end{lemma}
\begin{proof}
    Thanks to \eqref{eq:Ephiaddcst}, we may assume that $\varphi\geq \psi$.

       As we have seen in \eqref{eq:int_by_part_E_temp1}, the middle term can be written as
       \[
       \frac{1}{n+1}\sum_{j=0}^n \int_X(\varphi-\psi)\,\theta_{\varphi}^j\wedge \theta_{\psi}^{n-j}.
       \]
       We claim that for each individual $j$, the inequality \eqref{eq:comp_int_ene} holds:
       \[
       \int_X (\varphi-\psi)\,\theta_{\psi}^n \geq \int_X(\varphi-\psi)\,\theta_{\varphi}^j\wedge \theta_{\psi}^{n-j} \geq \int_X (\varphi-\psi)\,\theta_{\varphi}^n.
       \]
       Thanks to \cref{prop:compa_princ}, we have
       \[
       \begin{aligned}
       \int_X(\varphi-\psi)\,\theta_{\varphi}^j\wedge \theta_{\psi}^{n-j}=&\int_0^{\infty} \int_{\{\varphi>\psi+t\}} \theta_{\varphi}^j\wedge \theta_{\psi}^{n-j}\,\mathrm{d}t\\
       \geq & \int_0^{\infty} \int_{\{\varphi>\psi+t\}} \theta_{\varphi}^n\,\mathrm{d}t \\
       =& \int_X (\varphi-\psi)\,\theta_{\varphi}^n.
       \end{aligned}
       \]
       The other inequality is similar.
\end{proof}

\begin{lemma}\label{lma:Econcave1}
    The function $E^{\phi}_{\theta}\colon \mathcal{E}^{\infty}(X,\theta;\phi)\rightarrow \mathbb{R}$ is concave.
\end{lemma}
\begin{proof}
    Let $\varphi_0,\varphi_1\in \mathcal{E}^{\infty}(X,\theta;\phi)$. We write $\varphi_t=t\varphi_1+(1-t)\varphi_0$ for each $t\in [0,1]$.

    We claim that for each $t\in [0,1]$, we have
    \begin{equation}\label{eq:ddtEphi_temp1}
        \ddt E^{\phi}_{\theta}(\varphi_t) =\int_X \left(\varphi_1-\varphi_0\right)\,\theta_{\varphi_t}^n.
    \end{equation}
    At the boundary $t=0,1$, the derivative is understood as the one-sided derivative.

    Since $\varphi_0,\varphi_1\in \mathcal{E}^{\infty}(X,\theta;\phi)$, the difference $\varphi_1-\varphi_0$ is quasi-everywhere bounded. Hence all mixed integrals below are finite. We compute the right derivative first. Fix $t\in [0,1)$ and take $s>0$ small enough so that $t+s\leq 1$. By the cocycle formula \eqref{eq:int_by_part_E_temp1},
    \[
        \begin{aligned}
            \frac{E^{\phi}_{\theta}(\varphi_{t+s})-E^{\phi}_{\theta}(\varphi_t)}{s}
            =&\frac{1}{(n+1)s}\sum_{j=0}^n
            \int_X(\varphi_{t+s}-\varphi_t)\,
            \theta_{\varphi_{t+s}}^j\wedge\theta_{\varphi_t}^{n-j}\\
            =&\frac{1}{n+1}\sum_{j=0}^n
            \int_X(\varphi_1-\varphi_0)\,
            \theta_{\varphi_{t+s}}^j\wedge\theta_{\varphi_t}^{n-j}.
        \end{aligned}
    \]
    For each fixed $j$, the last integral is a polynomial in $s$ by the multi-linearity of the non-pluripolar product \itemref{itm:prop:npp1:6}. Letting $s\to0+$, it converges to
    \[
        \int_X(\varphi_1-\varphi_0)\,\theta_{\varphi_t}^n.
    \]
    Thus the right derivative is given by the right-hand side of \eqref{eq:ddtEphi_temp1}. The proof for the left derivative is identical, using $s<0$ with $t+s\in[0,1]$. Hence the left and right derivatives agree for $t\in(0,1)$, and the endpoint derivatives are the corresponding one-sided derivatives.

    Now it suffices to show that the derivative \eqref{eq:ddtEphi_temp1} is decreasing in $t$.
    Take $0\leq s<t \leq 1$. We wish to show that
    \[
        \int_X (\varphi_1-\varphi_0)\,\theta_{\varphi_t}^n\leq \int_X (\varphi_1-\varphi_0)\,\theta_{\varphi_s}^n.
    \]
    This can be equivalently written as
    \[
    \int_X (\varphi_t-\varphi_s)\,\theta_{\varphi_t}^n\leq \int_X (\varphi_t-\varphi_s)\,\theta_{\varphi_s}^n.
    \]
    In particular, replacing $\varphi_0$ by $\varphi_s$ and $\varphi_1$ by $\varphi_t$, we are reduced to proving the following
    \[
        \int_X (\varphi_1-\varphi_0)\,\theta_{\varphi_1}^n\leq \int_X (\varphi_1-\varphi_0)\,\theta_{\varphi_0}^n.
    \]
    But this is just a consequence of \cref{lma:comp_int_ene}.
\end{proof}

\begin{proposition}\label{prop:E_cont_decnet}
    Let $(\varphi_i)_{i\in I}$ be a decreasing net in $\PSH(X,\theta;\phi)$ with
    \[
        \varphi=\inf_{i\in I}\varphi_i\not\equiv -\infty.
    \]
    Then
    \begin{equation}\label{eq:E_cont_decnet}
        \lim_{i\in I} E_{\theta}^{\phi}(\varphi_i)=E_{\theta}^{\phi}(\varphi).
    \end{equation}
\end{proposition}
\begin{proof}
    We first handle the case where $I=\mathbb{Z}_{>0}$. Namely, $(\varphi_i)_i$ is a decreasing sequence.
    Thanks to \cref{lma:EmonoEinfty}, we know that $E_{\theta}^{\phi}(\varphi_i)$ is decreasing in $i$. So the limit in \eqref{eq:E_cont_decnet} exists and the $\geq$ inequality holds. Conversely, let $\psi\in \mathcal{E}^{\infty}(X,\theta;\phi)$ and $\varphi\leq \psi$. We need to show that
    \[
    E_{\theta}^{\phi}(\psi)\geq \lim_{i\to\infty} E_{\theta}^{\phi}\left(\varphi_i\lor \psi\right).
    \]
    Replacing the sequence $(\varphi_i)_i$ by $(\varphi_i\lor\psi)_i$, we are reduced to proving the $\leq$ direction of \eqref{eq:E_cont_decnet} when $\varphi\in \mathcal{E}^{\infty}(X,\theta;\phi)$.

    In this case, thanks to \cref{lma:comp_int_ene}, we have
    \[
    E_{\theta}^{\phi}(\varphi_i)-E_{\theta}^{\phi}(\varphi)\leq \int_X (\varphi_i-\varphi)\,\theta_{\varphi}^n.
    \]
    Our assertion follows from the monotone convergence theorem.

    We now treat the general net case. By \cref{lma:EmonoEinfty}, the net $E_{\theta}^{\phi}(\varphi_i)$ is decreasing and bounded below by $E_{\theta}^{\phi}(\varphi)$. Thus
    \[
        L\coloneqq \lim_{i\in I}E_{\theta}^{\phi}(\varphi_i)
    \]
    exists and satisfies $L\geq E_{\theta}^{\phi}(\varphi)$. It remains to prove the reverse inequality.

    Since $\varphi\leq \varphi_i$ for every $i\in I$, the hypotheses of \cref{cor:decreasingnetL1conv} are satisfied, and hence $\varphi_i\xrightarrow{L^1}\varphi$. For each $k\geq 1$, choose $a_k\in I$ such that
    \[
        \|\varphi_i-\varphi\|_{L^1}<k^{-1},\quad i\geq a_k.
    \]
    We can recursively choose an increasing sequence $(i_k)_{k\geq 1}$ in $I$ such that $i_k\geq a_k$ for all $k$. Then $(\varphi_{i_k})_k$ is a decreasing sequence and $\varphi_{i_k}\to\varphi$ in $L^1$. Its pointwise decreasing limit is therefore equal to $\varphi$ almost everywhere, hence everywhere by \cref{prop:pshfuncdetdense}. Applying the sequence case to $(\varphi_{i_k})_k$, we get
    \[
        \lim_{k\to\infty}E_{\theta}^{\phi}(\varphi_{i_k})=E_{\theta}^{\phi}(\varphi).
    \]
    Since $L\leq E_{\theta}^{\phi}(\varphi_{i_k})$ for every $k$, we conclude $L\leq E_{\theta}^{\phi}(\varphi)$.
\end{proof}

\begin{proposition}\label{prop:E1charac}
    Let $\varphi\in \mathcal{E}(X,\theta;\phi)$. The following are equivalent:
    \begin{enumerate}
        \item\label{itm:prop:E1charac:1} $\varphi\in \mathcal{E}^1(X,\theta;\phi)$;
        \item\label{itm:prop:E1charac:2} $E^{\phi}_{\theta}(\varphi)>-\infty$.
    \end{enumerate}
\end{proposition}
\begin{proof}
    Fix $\varphi\in \mathcal{E}(X,\theta;\phi)$. Without loss of generality, we may assume that $\varphi\leq 0$. In particular, $\varphi\leq P_{\theta}[\varphi]\leq P_{\theta}[\phi]=\phi$.

    \ref{itm:prop:E1charac:2} $\implies$ \ref{itm:prop:E1charac:1}.
    Assume \ref{itm:prop:E1charac:2}. Using the plurifine locality of the non-pluripolar product \itemref{itm:prop:npp1:1}, we have that $\mathds{1}_{\{\varphi>\phi-C\}} \theta_{\varphi\lor (\phi-C)}^n=\mathds{1}_{\{\varphi>\phi-C\}} \theta_{\varphi}^n$ converges \emph{strongly}\footnote{Namely, the convergence holds after integrating against any $L^{\infty}$ function.} to $\theta_{\varphi}^n$ as $C\to \infty$ by the dominated convergence theorem. Therefore, for a fixed $D>0$,
    \[
    \begin{aligned}
    \int_X \Bigl( \varphi\lor (\phi-D)-\phi \Bigr)\,\theta_{\varphi}^n=&\lim_{C\to\infty} \int_{\{\varphi>\phi-C\} } \Bigl( \varphi\lor (\phi-D)-\phi \Bigr)\,\theta_{\varphi\lor (\phi-C)}^n\\
    \geq & \varlimsup_{C\to\infty} \int_{\{\varphi>\phi-C\} } \Bigl( \varphi\lor (\phi-C)-\phi \Bigr)\,\theta_{\varphi\lor (\phi-C)}^n\\
    \geq & \varlimsup_{C\to\infty} \int_X \Bigl( \varphi\lor (\phi-C)-\phi \Bigr)\,\theta_{\varphi\lor (\phi-C)}^n\\
    \geq & \varlimsup_{C\to\infty} (n+1)E_{\theta}^{\phi}\Bigl(\varphi \lor (\phi-C) \Bigr)\\
    \geq & (n+1)E_{\theta}^{\phi}(\varphi).
    \end{aligned}
    \]
    Here the last but one inequality follows from \eqref{eq:Edefbdd} and the fact that $\varphi\lor (\phi-C)\leq \phi$.
    Letting $D\to\infty$, we conclude \ref{itm:prop:E1charac:1}.

    \ref{itm:prop:E1charac:1} $\implies$ \ref{itm:prop:E1charac:2}. Assume \ref{itm:prop:E1charac:1}. Thanks to \cref{lma:comp_int_ene}, we know that for each $C>0$,
    \[
    \begin{aligned}
    E_{\theta}^{\phi}\Bigl(\varphi\lor (\phi-C)\Bigr)\geq & \int_X \Bigl(\varphi\lor (\phi-C)-\phi \Bigr)\,\theta_{\varphi\lor (\phi-C)}^n\\
    =& \int_{\left\{ \varphi>\phi-C \right\}} \Bigl(\varphi-\phi \Bigr)\,\theta_{\varphi \lor (\phi-C)}^n+ \int_{\left\{ \varphi\leq \phi-C \right\}} (-C)\,\theta_{\varphi\lor (\phi-C)}^n \\
    =& \int_{\{\varphi+C>\phi\}} \left(\varphi+C-\phi \right)\,\theta_{\varphi}^n-C\int_X \theta_{\phi}^n\\
    \geq & \int_X \left(\varphi+C-\phi \right)\,\theta_{\varphi}^n-C\int_X \theta_{\phi}^n\\
    =& \int_X (\varphi-\phi)\,\theta_{\varphi}^n.
    \end{aligned}
    \]
    Due to \cref{prop:E_cont_decnet}, we have
    \[
    \lim_{C\to\infty}E_{\theta}^{\phi}\Bigl(\varphi\lor (\phi-C)\Bigr)=E_{\theta}^{\phi}(\varphi),
    \]
    so \ref{itm:prop:E1charac:2} holds.
\end{proof}

\begin{corollary}\label{cor:E_concave2}
    The function $E^{\phi}_{\theta}\colon \PSH(X,\theta;\phi)\rightarrow [-\infty,\infty)$ is concave.
\end{corollary}
\begin{proof}
    Let $\varphi,\psi\in \PSH(X,\theta;\phi)$ and $t\in (0,1)$. We need to show that
    \begin{equation}\label{eq:Econcave_temp1}
    E^{\phi}_{\theta}\Bigl( t\varphi+(1-t)\psi \Bigr)\geq tE^{\phi}_{\theta}(\varphi)+(1-t)E^{\phi}_{\theta}(\psi).
    \end{equation}

    For this purpose, observe that when $C\to \infty$, we have
    \[
        t\Bigl(\varphi\lor (\phi-C)\Bigr)+(1-t)\Bigl(\psi\lor (\phi-C)\Bigr)\searrow t\varphi+(1-t)\psi.
    \]
    Observe that
    \[
        \varphi\lor (\phi-C),\psi\lor (\phi-C)\in \mathcal{E}^{\infty}(X,\theta;\phi).
    \]
    By \cref{lma:Econcave1}, for each $C>0$, we have
    \[
        E_{\theta}^{\phi}\Biggl( t\Bigl(\varphi\lor (\phi-C)\Bigr)+(1-t)\Bigl(\psi\lor (\phi-C)\Bigr) \Biggr)\geq tE_{\theta}^{\phi}\Bigl(\varphi\lor (\phi-C)\Bigr)+(1-t)E_{\theta}^{\phi}\Bigl(\psi\lor (\phi-C)\Bigr).
    \]
    Letting $C\to \infty$, we conclude \eqref{eq:Econcave_temp1} using \cref{prop:E_cont_decnet}.
\end{proof}

\begin{corollary}\label{cor:E_classes_convex}
    The sets $\PSH(X,\theta;\phi)$, $\mathcal{E}(X,\theta;\phi)$, $\mathcal{E}^1(X,\theta;\phi)$, $\mathcal{E}^{\infty}(X,\theta;\phi)$ are all convex.
\end{corollary}
\begin{proof}
    The sets $\PSH(X,\theta;\phi)$ and $\mathcal{E}^{\infty}(X,\theta;\phi)$ are trivially convex.

    We next handle $\mathcal{E}(X,\theta;\phi)$. In this case, suppose that $\varphi,\psi \in \mathcal{E}(X,\theta;\phi)$ and $t\in (0,1)$. Then we compute
    \[
    \begin{split}
        \int_X \theta_{t\varphi+(1-t)\psi}^n=\sum_{j=0}^n \binom{n}{j}t^j (1-t)^{n-j} \int_X \theta_{\varphi}^j\wedge \theta_{\psi}^{n-j}\\
        =\sum_{j=0}^n \binom{n}{j}t^j (1-t)^{n-j} \int_X \theta_{\phi}^n=\int_X \theta_{\phi}^n,
    \end{split}
    \]
    where we applied \cref{prop:Ppresmass} in the second equality. Therefore, $t\varphi+(1-t)\psi \in \mathcal{E}(X,\theta;\phi)$.

    Finally, let us consider $\mathcal{E}^1(X,\theta;\phi)$. In this case, suppose that $\varphi,\psi \in \mathcal{E}^1(X,\theta;\phi)$ and $t\in (0,1)$. We have already shown that $t\varphi+(1-t)\psi \in \mathcal{E}(X,\theta;\phi)$. Furthermore, the energy of this potential is finite due to \cref{cor:E_concave2}. Therefore, $t\varphi+(1-t)\psi \in \mathcal{E}^1(X,\theta;\phi)$ in view of \cref{prop:E1charac}.
\end{proof}

\begin{lemma}\label{lma:intdom1}
	Let $\varphi,\psi,\gamma\in \mathcal{E}(X,\theta;\phi)$. Assume that $\gamma\geq \varphi\lor\psi$.
	Then
	\[
	\int_X (\gamma-\varphi)\,\theta_{\psi}^n\leq 2\int_X (\gamma-\varphi)\,\theta_{\varphi}^n+2\int_X (\gamma-\psi)\,\theta_{\psi}^n.
	\]
\end{lemma}
\begin{proof}
	Observe that
	\begin{equation}\label{eq:gvp2}
		\{\gamma>\varphi+2t\}\subseteq \{\gamma>\psi+t \}\cup \{\psi>\varphi+t\}\,.
	\end{equation}
	So
	\[
	\begin{aligned}
		\int_X (\gamma-\varphi)\,\theta_{\psi}^n
		=& 2\int_0^{\infty} \int_{\{\gamma>\varphi+2t\}}\,\theta_{\psi}^n\,\mathrm{d}t \\
		\leq& 2\int_0^{\infty} \int_{\{\gamma>\psi+t\}}\,\theta_{\psi}^n\,\mathrm{d}t+
		2\int_0^{\infty} \int_{\{\psi>\varphi+t\}}\,\theta_{\psi}^n\,\mathrm{d}t\\
		\leq & 2\int_X (\gamma-\psi) \,\theta_{\psi}^n+2\int_0^{\infty} \int_{\{\psi>\varphi+t\}}\,\theta_{\varphi}^n\,\mathrm{d}t\\
		\leq & 2\int_X (\gamma-\psi) \,\theta_{\psi}^n+2\int_0^{\infty} \int_{\{\gamma>\varphi+t\}}\,\theta_{\varphi}^n\,\mathrm{d}t\\
		=&2\int_X (\gamma-\psi) \,\theta_{\psi}^n+2\int_X (\gamma-\varphi) \,\theta_{\varphi}^n,
	\end{aligned}
	\]
    where the second line follows from \eqref{eq:gvp2}, while the third line follows from \cref{prop:compa_princ}.
\end{proof}

\begin{corollary}\label{cor:mix_MA_energy_finite}
    Let $\varphi,\psi,\gamma_1,\ldots,\gamma_n\in \mathcal{E}^1(X,\theta;\phi)$. Then
    \begin{equation}\label{eq:mix_MA_energy}
        \int_X (\psi-\varphi)\,\theta_{\gamma_1}\wedge \cdots \wedge \theta_{\gamma_n}
    \end{equation}
    is well-defined and is in $\mathbb{R}$.
\end{corollary}
\begin{proof}
    \textbf{Step~1}. We first consider the case $\psi=\phi$.

    We can find $C>0$ so that $\varphi\leq \phi+C$. Therefore, $\phi-\varphi \geq -C$ almost everywhere with respect to $\theta_{\gamma_1}\wedge \cdots \wedge \theta_{\gamma_n}$. So the integral \eqref{eq:mix_MA_energy} is well-defined and takes value in $\mathbb{R}\cup \{\infty\}$. We show that the integral is finite. For this purpose, we may replace $\varphi$ by $\varphi-C$ and assume that $\varphi\leq \phi$. Put
    \[
        \eta=n^{-1}(\gamma_1+\cdots+\gamma_n).
    \]
    Then by the multi-linearity of the non-pluripolar product \itemref{itm:prop:npp1:6} and positivity, we have
    \[
        \int_X (\phi-\varphi)\,\theta_{\gamma_1}\wedge \cdots \wedge \theta_{\gamma_n}\leq n^n\int_X (\phi-\varphi)\,\theta_{\eta}^n.
    \]
    Thanks to \cref{cor:E_classes_convex}, $\eta\in \mathcal{E}^1(X,\theta;\phi)$. Choose $A>0$ so that $\eta\leq \phi+A$. Applying \cref{lma:intdom1} to $\varphi,\eta$ and $\phi+A$, we get
    \[
        \int_X(\phi+A-\varphi)\,\theta_{\eta}^n
        \leq 2\int_X(\phi+A-\varphi)\,\theta_{\varphi}^n
        +2\int_X(\phi+A-\eta)\,\theta_{\eta}^n<\infty.
    \]
    Hence $\int_X(\phi-\varphi)\,\theta_{\eta}^n<\infty$, and Step~1 follows.

    \textbf{Step~2}. We handle the general case. Take $C>0$ so that
    \[
        \psi\leq \phi+C,\quad \varphi\leq \phi+C.
    \]
    Then
    \[
        (\psi-C)-\phi \leq \psi-\varphi\leq \phi-(\varphi-C).
    \]
    Step~1 applied to $\varphi-C$ and to $\psi-C$ shows that $\phi-(\varphi-C)$ and $\phi-(\psi-C)$ are integrable with respect to $\theta_{\gamma_1}\wedge \cdots \wedge \theta_{\gamma_n}$. Hence the positive and negative parts of $\psi-\varphi$ are both integrable, so \eqref{eq:mix_MA_energy} is well-defined and finite.
\end{proof}

\begin{proposition}\label{prop:cocycE1}
    Given $\varphi,\psi\in \mathcal{E}^1(X,\theta;\phi)$, we have the following cocycle equality
    \begin{equation}\label{eq:Ecocyc}
        E^{\phi}_{\theta}(\psi)-E^{\phi}_{\theta}(\varphi)=\frac{1}{n+1}\sum_{j=0}^n \int_X (\psi-\varphi)\, \theta_{\psi}^j\wedge \theta_{\varphi}^{n-j}.
    \end{equation}
\end{proposition}
As a consequence of \eqref{eq:Ecocyc}, for $\varphi\in \mathcal{E}^1(X,\theta;\phi)$, \eqref{eq:Edefbdd} continues to hold. Note that each term on the right-hand side of \eqref{eq:Ecocyc} is finite due to \cref{cor:mix_MA_energy_finite}.
\begin{proof}
    Since the conclusion is known when $\varphi,\psi\in \mathcal{E}^{\infty}(X,\theta;\phi)$ (see \eqref{eq:int_by_part_E_temp1}), it suffices to prove that for each $k=0,\ldots,n$, we have
    \begin{equation}\label{eq:canonic_app_energy_conv}
    \lim_{C\to\infty} \int_X \Bigl( \varphi\lor (\phi-C)-\psi\lor (\phi-C) \Bigr) \,\theta_{\varphi\lor (\phi-C)}^k\wedge \theta_{\psi\lor (\phi-C)}^{n-k}=\int_X (\varphi-\psi)\,\theta_{\varphi}^k\wedge \theta_{\psi}^{n-k}.
    \end{equation}
    For this purpose, we may assume that $\varphi,\psi\leq \phi$ and it suffices to establish the following:
    \begin{equation}\label{eq:intvarphilorphimCmphi_temp1}
        \lim_{C\to\infty} \int_X \Bigl( \varphi\lor (\phi-C)-\phi \Bigr) \,\theta_{\varphi\lor (\phi-C)}^k\wedge \theta_{\psi\lor (\phi-C)}^{n-k}=\int_X (\varphi-\phi)\,\theta_{\varphi}^k\wedge \theta_{\psi}^{n-k}.
    \end{equation}
    Fix $C>0$. We compute the difference as follows:
    \[
    \begin{aligned}
    &\int_X \Bigl( \varphi\lor (\phi-C)-\phi \Bigr) \,\theta_{\varphi\lor (\phi-C)}^k\wedge \theta_{\psi\lor (\phi-C)}^{n-k}-\int_X (\varphi-\phi)\,\theta_{\varphi}^k\wedge \theta_{\psi}^{n-k}\\
    =& \int_{\{\min\{\varphi,\psi\}>\phi-C\}}\left( \varphi-\phi \right) \,\theta_{\varphi}^k\wedge \theta_{\psi}^{n-k}-\int_X (\varphi-\phi)\,\theta_{\varphi}^k\wedge \theta_{\psi}^{n-k}\\
    &+\int_{\{\min\{\varphi,\psi\}\leq \phi-C\}}\Bigl( \varphi\lor (\phi-C)-\phi \Bigr) \,\theta_{\varphi\lor (\phi-C)}^k\wedge \theta_{\psi\lor (\phi-C)}^{n-k}\\
    =& -\int_{\{\min\{\varphi,\psi\}\leq \phi-C\}}\left( \varphi-\phi \right) \,\theta_{\varphi}^k\wedge \theta_{\psi}^{n-k}+\\
     &+\int_{\{\min\{\varphi,\psi\}\leq \phi-C\}}\Bigl( \varphi\lor (\phi-C)-\phi \Bigr) \,\theta_{\varphi\lor (\phi-C)}^k\wedge \theta_{\psi\lor (\phi-C)}^{n-k}
    \end{aligned}
    \]
    We first handle the first term. In order to show that it converges to $0$ when $C\to \infty$, it suffices to show that
    \[
    \int_X\left( \phi-\varphi \right) \,\theta_{\varphi}^k\wedge \theta_{\psi}^{n-k}<\infty.
    \]
    For this purpose, it suffices to prove a stronger inequality:
    \begin{equation}\label{eq:intXphimvarphifinite_temp1}
    \int_X\left( \phi-\varphi \right) \,\theta_{(\varphi+\psi)/2}^n<\infty.
    \end{equation}
    From \cref{cor:E_classes_convex}, we know that $(\varphi+\psi)/2\in \mathcal{E}^1(X,\theta;\phi)$. Therefore, \eqref{eq:intXphimvarphifinite_temp1} follows from \cref{lma:intdom1}.

    It remains to handle the second term. It then suffices to prove the following:
    \begin{equation}\label{eq:limCCintvarphiphiC_temp1}
    \begin{split}
    \lim_{C\to\infty} C\int_{\{\varphi\leq \phi-C\}} \theta_{\varphi\lor (\phi-C)}^k\wedge \theta_{\psi\lor (\phi-C)}^{n-k}=&0,\\
    \lim_{C\to\infty} C\int_{\{\psi\leq \phi-C\}} \theta_{\varphi\lor (\phi-C)}^k\wedge \theta_{\psi\lor (\phi-C)}^{n-k}=&0.
    \end{split}
    \end{equation}
    By symmetry, it suffices to handle the first. Observe the following inclusions:
    \[
    \{\varphi\leq \phi-C\}\subseteq \left\{ \varphi\lor (\phi-C)\leq \frac{1}{2}\Bigl( \psi\lor (\phi-C)+\phi-C \Bigr) \right\}\subseteq \{\varphi\leq \phi-C/2\}.
    \]
    By \cref{prop:compa_princ}, we find
    \[
    \begin{aligned}
        &C\int_{\{\varphi\leq \phi-C\}} \theta_{\varphi\lor (\phi-C)}^k\wedge \theta_{\psi\lor (\phi-C)}^{n-k}\\
        \leq & C\int_{\left\{\varphi\lor (\phi-C)\leq \frac{1}{2}\bigl(\psi\lor (\phi-C)+\phi-C\bigr)\right\}} \theta_{\varphi\lor (\phi-C)}^k\wedge \theta_{\psi\lor (\phi-C)}^{n-k}\\
        \leq & 2^{n-k}C\int_{\left\{\varphi\lor (\phi-C)\leq \frac{1}{2}\bigl(\psi\lor (\phi-C)+\phi-C\bigr)\right\}} \theta_{\varphi\lor (\phi-C)}^k\wedge \theta_{\frac{1}{2}\bigl(\psi\lor (\phi-C)+\phi-C\bigr)}^{n-k}\\
        \leq & 2^{n-k}C\int_{\{\varphi\lor (\phi-C)\leq \frac{1}{2}\bigl(\psi\lor (\phi-C)+\phi-C\bigr)\}} \theta_{\varphi\lor (\phi-C)}^n\\
        \leq & 2^{n-k}C\int_{\{\varphi\leq \phi-C/2\}} \theta_{\varphi\lor (\phi-C)}^n\\
        = & 2^{n-k}C\int_{\{\varphi\leq \phi-C/2\}} \theta_{\varphi}^n.
    \end{aligned}
    \]
    Here on the final step, we have applied the following argument:
    \[
    \begin{aligned}
    &\int_{\left\{\varphi\leq \phi-C/2 \right\}} \theta_{\varphi\lor (\phi-C)}^n &\\
    = &\int_X \theta_{\varphi\lor (\phi-C)}^n-\int_{\left\{\varphi> \phi-C/2 \right\}} \theta_{\varphi\lor (\phi-C)}^n &\\
    = & \int_X \theta_{\varphi}^n-\int_{\left\{\varphi> \phi-C/2 \right\}} \theta_{\varphi}^n & \textup{by \itemref{itm:prop:npp1:1}}\\
    =& \int_{\left\{\varphi\leq \phi-C/2 \right\}} \theta_{\varphi}^n. &
    \end{aligned}
    \]

    On the other hand,
    \[
    \varlimsup_{C\to\infty} C\int_{\{\varphi\leq \phi-C/2 \}} \theta_{\varphi}^n\leq 2 \varlimsup_{C\to\infty} \int_{\{\varphi\leq \phi-C \}}(\phi-\varphi)\, \theta_{\varphi}^n=0,
    \]
    since $\varphi\in \mathcal{E}^1(X,\theta;\phi)$. We conclude \eqref{eq:limCCintvarphiphiC_temp1}.
\end{proof}
\begin{proposition}\label{prop:Ediffbdd}
Let $\varphi,\psi\in \mathcal{E}^1(X,\theta;\phi)$. Then
    \begin{equation}
        \int_X \left(\psi-\varphi\right)\theta^n_{\psi}\leq E_{\theta}^{\phi}(\psi)-E_{\theta}^{\phi}(\varphi)\leq \int_X \left(\psi-\varphi\right)\theta^n_{\varphi}.
    \end{equation}
\end{proposition}
\begin{proof}
    Thanks to \eqref{eq:canonic_app_energy_conv}, this can be reduced to the case where $\varphi,\psi\in \mathcal{E}^{\infty}(X,\theta;\phi)$, which is established in \cref{lma:comp_int_ene}.
\end{proof}

\begin{proposition}\label{prop:Econt_mono}
    Let $(\varphi_i)_{i\in I}$ be a net in $\mathcal{E}^1(X,\theta;\phi)$. Assume that $\varphi\in \mathcal{E}^1(X,\theta;\phi)$ and either one of the following conditions holds:
    \begin{enumerate}
        \item\label{itm:prop:Econt_mono:1} $(\varphi_i)_i$ is a decreasing net with $\varphi=\inf_{i\in I}\varphi_i$;
        \item\label{itm:prop:Econt_mono:2} $(\varphi_i)_i$ is an increasing net with almost everywhere limit $\varphi$.
    \end{enumerate}
    Then
    \begin{equation}\label{eq:Iconv_mono}
    \lim_{i\in I} E_{\theta}^{\phi}(\varphi_i)=E_{\theta}^{\phi}(\varphi).
    \end{equation}
\end{proposition}
\begin{proof}
    The decreasing case is already proved in \cref{prop:E_cont_decnet}, so let us focus on the case where $(\varphi_i)_i$ is increasing.

    We first handle the case where $I=\mathbb{Z}_{>0}$. Namely, $(\varphi_i)_i$ is an increasing sequence.

    \textbf{Step~1}. We first prove the result when $\varphi_j,\varphi\in \mathcal{E}^{\infty}(X,\theta;\phi)$ for each $j\geq 1$.

    Due to \cref{cor:incseqnppcont}, we have
    \[
        \lim_{j\to\infty} \int_X \theta_{\varphi}^k\wedge \theta_{\varphi_j}^{n-k}=\int_X \theta_{\varphi}^n
    \]
    for each $k=0,\ldots,n$.
    By \eqref{eq:int_by_part_E_temp1}, we have
    \[
    \begin{aligned}
    0\leq & E_{\theta}^{\phi}(\varphi)-E_{\theta}^{\phi}(\varphi_j)\\
    =& \frac{1}{n+1}\sum_{k=0}^n \int_X(\varphi-\varphi_j)\,\theta_{\varphi}^k\wedge \theta_{\varphi_j}^{n-k}.
    \end{aligned}
    \]
    We then apply \cref{thm:semicon} to conclude that the right-hand side converges to $0$.

    \textbf{Step~2}. We prove the general case. The $\leq$ direction in \eqref{eq:Iconv_mono} follows from the monotonicity of $E_{\theta}^{\phi}$, as proved in \cref{lma:EmonoEinfty}. It remains to prove the $\geq$ direction.

    After subtracting a constant from all $\varphi_j$ and from $\varphi$, which shifts all the energies by the same constant in view of \eqref{eq:Ephiaddcst}, we may assume that $\varphi\leq \phi$. For each $C>0$ and $j\geq 1$, we let
    \[
    \varphi_C\coloneqq \varphi\lor (\phi-C),\quad \varphi_{j,C}\coloneqq \varphi_j\lor (\phi-C).
    \]
    By Step~1, for each $C>0$, we have
    \[
    \lim_{j\to\infty} E_{\theta}^{\phi}(\varphi_{j,C})= E_{\theta}^{\phi}(\varphi_{C})\geq E_{\theta}^{\phi}(\varphi).
    \]
    It suffices therefore to show that
    \begin{equation}\label{eq:limCCEphi_temp1}
    \lim_{C\to \infty} E_{\theta}^{\phi}(\varphi_{j,C})-E_{\theta}^{\phi}(\varphi_{j})=0\quad \textup{uniformly for $j$ in a tail}.
    \end{equation}
    Fix $i\geq 1$.
    For $j\geq i$, we compute
    \[
    \begin{aligned}
        & E_{\theta}^{\phi}(\varphi_{j,C})-E_{\theta}^{\phi}(\varphi_{j})\\
        \leq & \int_X \left( \varphi_{j,C}-\varphi_{j}\right)\,\theta_{\varphi_{j}}^n &\textup{by \cref{prop:Ediffbdd}}\\
        \leq & \int_{\{\varphi_j\leq \phi-C\}} (\phi-C-\varphi_j)\,\theta_{\varphi_j}^n\\
        =& \int_C^{\infty}\int_{\{\varphi_j\leq \phi-t\}} \theta_{\varphi_j}^n\,\mathrm{d}t\\
        \leq & \int_C^{\infty}\int_{\{\varphi_i\leq (\varphi_j+\phi-t)/2\}} \theta_{\varphi_j}^n\,\mathrm{d}t\\
        \leq & 2^n\int_C^{\infty}\int_{\{\varphi_i\leq (\varphi_j+\phi-t)/2\}} \theta_{(\varphi_j+\phi-t)/2}^n\,\mathrm{d}t\\
        \leq & 2^n\int_C^{\infty}\int_{\{\varphi_i\leq (\varphi_j+\phi-t)/2\}} \theta_{\varphi_i}^n\,\mathrm{d}t & \textup{by \cref{prop:compa_princ}}\\
        \leq &2^n\int_C^{\infty}\int_{\{\varphi_i\leq \phi-t/2\}} \theta_{\varphi_i}^n\,\mathrm{d}t\\
        = & 2^{n+1}\int_{\{\varphi_i< \phi-C/2\}} (\phi-\varphi_i-C/2)\,\theta_{\varphi_i}^n\\
        = & 2^{n+1}\int_{\{\varphi_i< \phi-C/2\}} (\phi-\varphi_i)\,\theta_{\varphi_i}^n-C 2^n\int_{\{\varphi_i< \phi-C/2\}}\theta_{\varphi_i}^n.
    \end{aligned}
    \]
    As $C\to \infty$, the first term converges to $0$ since $\varphi_i\in \mathcal{E}^1(X,\theta;\phi)$. As for the second term,
    \[
    \lim_{C\to\infty}C\int_{\{\varphi_i< \phi-C/2 \}} \theta_{\varphi_i}^n\leq 2 \lim_{C\to\infty}\int_{\{\varphi_i< \phi-C/2 \}}(\phi-\varphi_i)\, \theta_{\varphi_i}^n = 0.
    \]
    We conclude \eqref{eq:limCCEphi_temp1}.

    We now pass from sequences to nets. By \cref{lma:EmonoEinfty}, the net $(E_{\theta}^{\phi}(\varphi_i))_{i\in I}$ is increasing and bounded above by $E_{\theta}^{\phi}(\varphi)$. Let
    \[
        L\coloneqq \lim_{i\in I}E_{\theta}^{\phi}(\varphi_i).
    \]
    Then $L\leq E_{\theta}^{\phi}(\varphi)$.

    By \cref{prop:pshfuncdetdense}, we have
    \[
        \varphi=\uscsup[i\in I]\varphi_i.
    \]
    Moreover, $\varphi_i\leq \varphi$ for all $i\in I$, so the family is locally uniformly bounded above.
    Applying Choquet's lemma \cref{prop:Choquet}, we can find a countable subset $J\subseteq I$ such that
    \[
        \varphi=\uscsup[j\in J]\varphi_j.
    \]
    Enumerate $J=\{j_1,j_2,\ldots\}$, repeating elements if necessary. Since $I$ is directed, we can choose an increasing sequence $(i_k)_{k\geq 1}$ in $I$ such that $i_k$ dominates $j_1,\ldots,j_k$. It follows that $(\varphi_{i_k})_k$ increases almost everywhere to $\varphi$. Applying the sequence case to this sequence gives
    \[
        \lim_{k\to\infty}E_{\theta}^{\phi}(\varphi_{i_k})=E_{\theta}^{\phi}(\varphi).
    \]
    Since $E_{\theta}^{\phi}(\varphi_{i_k})\leq L$ for every $k$, we conclude $E_{\theta}^{\phi}(\varphi)\leq L$.
\end{proof}

Next we want to prove that $\mathcal{E}^1(X,\theta;\phi)$ is closed under the rooftop operator. For this purpose, we shall need a few preliminary results.

We shall approximate the rooftop operator by solutions of certain Monge--Ampère equations.
\begin{lemma}\label{lma:rooftop_MAequation}
	Let $\varphi,\psi\in \mathcal{E}^{\infty}(X,\theta;\phi)$. Then there is $\gamma\in \mathcal{E}^{\infty}(X,\theta;\phi)$ such that
	\begin{equation}\label{eq:solution_thetavarphiephiu}
		\theta_{\gamma}^n=\mathrm{e}^{\gamma-\varphi}\,\theta_{\varphi}^n+\mathrm{e}^{\gamma-\psi}\,\theta_{\psi}^n.
	\end{equation}
\end{lemma}
It is not clear if $\int_X\mathrm{e}^{-\varphi}\theta_{\varphi}^n<\infty$, hence we cannot say that $\mathrm{e}^{-\varphi}\theta_{\varphi}^n+\mathrm{e}^{-\psi}\theta_{\psi}^n$ is a Radon measure, so \cref{thm:AY_rel} is not directly applicable.
\begin{proof}
	For each $j\geq 1$, let $\varphi_j\coloneqq \varphi\lor(-j)$, $\psi_j\coloneqq \psi\lor(-j)$. These functions are not necessarily $\theta$-psh.
	Let
	\[
	\mu_j=\mathrm{e}^{-\varphi_j}\theta_{\varphi}^n+\mathrm{e}^{-\psi_j}\theta_{\psi}^n.
	\]
	By \cref{thm:AY_rel}, we can find $\gamma_j\in \mathcal{E}(X,\theta;\phi)$ such that
	\begin{equation}\label{eq:thetagammaj_temp1}
		\theta_{\gamma_j}^n=\mathrm{e}^{\gamma_j}\mu_j=\mathrm{e}^{\gamma_j-\varphi_j}\theta_{\varphi}^n+\mathrm{e}^{\gamma_j-\psi_j}\theta_{\psi}^n.
	\end{equation}
	Take a constant $C>0$ so that $\psi-2C\leq \varphi\leq \psi+2C$.
	Let
	\[
	\eta\coloneqq \frac{\varphi+\psi}{2}-C-n\log 2.
	\]
	Then $\eta\in \mathcal{E}^{\infty}(X,\theta;\phi)$ and we have the following inequalities:
    \[
    \begin{aligned}
    \theta_{\eta}^n =& \theta_{(\varphi+\psi)/2}^n\\
    \geq & 2^{-n}\theta_{\varphi}^n+2^{-n}\theta_{\psi}^n\\
    \geq & 2^{-n}\mathrm{e}^{(\varphi+\psi)/2}\mathrm{e}^{-C}\left( \mathrm{e}^{-\varphi}\theta_{\varphi}^n+\mathrm{e}^{-\psi}\theta_{\psi}^n \right)\\
    \geq & 2^{-n}\mathrm{e}^{(\varphi+\psi)/2}\mathrm{e}^{-C}\left( \mathrm{e}^{-\varphi_j}\theta_{\varphi}^n+\mathrm{e}^{-\psi_j}\theta_{\psi}^n \right)\\
    = &
    \mathrm{e}^{\eta}\mu_j.
    \end{aligned}
    \]
	Hence, $\gamma_j\geq \eta$ by \cref{thm:AY_rel}. By \cref{thm:AY_rel} again, $\gamma_j$ is decreasing in $j$; let $\gamma=\inf_{j>0}\gamma_j$.
	Then $\gamma\geq \eta$, hence $\gamma\in \mathcal{E}^{\infty}(X,\theta;\phi)$.

    Now observe that
    \[
    \theta_{\gamma_j}^n\rightharpoonup \theta_{\gamma}^n
    \]
    as $j\to \infty$, by \cref{thm:semicon}.
    Finally observe that there is a constant $C'>0$ so that
    \[
    \gamma_1\leq \varphi+C',\quad \gamma_1\leq \psi+C',
    \]
    since $\varphi,\psi\in \mathcal{E}^{\infty}(X,\theta;\phi)$.
    Therefore, for each $j>0$,
    \[
    \gamma_j\leq \varphi_j+C',\quad \gamma_j\leq \psi_j+C'.
    \]
    Hence, \eqref{eq:solution_thetavarphiephiu} follows from \eqref{eq:thetagammaj_temp1} and the dominated convergence theorem by letting $j\to\infty$.
\end{proof}

We also need a few integral estimates.

\begin{lemma}\label{lma:fundam_estimate}
	Let $\varphi,\psi,\gamma\in \mathcal{E}(X,\theta;\phi)$. Assume that $\varphi\leq \psi\leq \gamma$. Then
	\[
	\int_X (\gamma-\psi)\,\theta_{\psi}^n\leq 2^{n+1}\int_X (\gamma-\varphi)\,\theta_{\varphi}^n.
	\]
\end{lemma}
\begin{proof}
	Observe that for any $t\geq 0$,
	\begin{equation}\label{eq:ppgcomp}
        \begin{split}
		\{\gamma>\psi+2t\}\setminus \{\psi=-\infty\}\subseteq \left\{(\gamma+\psi)/2>\psi+t\right\}\\
        \subseteq \left\{(\gamma+\psi)/2>\varphi+t\right\}\subseteq \{\gamma>\varphi+t\}.
        \end{split}
	\end{equation}
	So
	\[
	\begin{aligned}
		\int_X (\gamma-\psi)\,\theta_{\psi}^n
		=& 2\int_0^{\infty}  \int_{\{\gamma-\psi>2t\}}\,\theta_{\psi}^n\,\mathrm{d}t &\\
		\leq & 2\int_0^{\infty} \int_{\left\{(\gamma+\psi)/2>\psi+t\right\}}\,\theta_{\psi}^n\,\mathrm{d}t & \textup{by \eqref{eq:ppgcomp}}\\
		\leq  & 2^{n+1}\int_0^{\infty} \int_{\left\{(\gamma+\psi)/2>\varphi+t\right\}}\,\theta_{(\gamma+\psi)/2}^n\,\mathrm{d}t &\\
		\leq & 2^{n+1}\int_0^{\infty} \int_{\left\{(\gamma+\psi)/2>\varphi+t\right\}}\,\theta_{\varphi}^n\,\mathrm{d}t & \textup{by \cref{prop:compa_princ}}\\
		\leq &2^{n+1}\int_0^{\infty} \int_{\left\{\gamma>\varphi+t\right\}}\,\theta_{\varphi}^n\,\mathrm{d}t &\textup{by \eqref{eq:ppgcomp}}\\
		=& 2^{n+1} \int_X (\gamma-\varphi)\,\theta_{\varphi}^n &.
	\end{aligned}
	\]
	This completes the proof.
\end{proof}

\begin{lemma}\label{lma:stab}
	Let $\varphi_j\in \mathcal{E}(X,\theta;\phi)$ ($j\in \mathbb{Z}_{>0}$) and $\gamma\in \mathcal{E}(X,\theta;\phi)$. Assume that $\varphi_j\leq \gamma$ for each $j$ and that $\varphi_j$ converges to $\varphi\in \PSH(X,\theta)$ with respect to the $L^1$-topology.
	Assume that there is a constant $A>0$ such that for any $j>0$,
	\[
	\int_X (\varphi_j-\gamma)\,\theta_{\varphi_j}^n\geq -A.
	\]
	Then $\varphi\in \mathcal{E}(X,\theta;\phi)$ and
	\begin{equation}\label{eq:intvarphigammageqm2np3A}
	    \int_X (\varphi-\gamma)\,\theta_{\varphi}^n\geq - 2^{n+3}A.
	\end{equation}
\end{lemma}
\begin{proof}
    Replacing simultaneously $\gamma,\varphi_j,\varphi$ by $\gamma-B,\varphi_j-B,\varphi-B$ for $B\gg 1$, we may assume that $\gamma\leq 0$. The assumptions and the conclusion are invariant under this translation. Since $\gamma\in\mathcal E(X,\theta;\phi)$ and $P_\theta[\phi]=\phi$, we then have $\gamma\leq\phi$.

	\textbf{Step~1}. Assume that $(\varphi_j)_j$ is decreasing. In this case, we prove
	\begin{equation}\label{eq:intvarphigammageqm4A_temp1}
	\int_X (\varphi-\gamma)\,\theta_{\varphi}^n\geq -4A.
    \end{equation}
	By \cref{lma:intdom1}, for any $j, k>0$,
	\[
	\int_X (\varphi_j-\gamma)\,\theta_{\varphi_k}^n\geq -4A.
	\]
	For any $C>0$,
	\[
	\int_X \Bigl(\varphi_j \lor (\gamma-C)-\gamma\Bigr)\,\theta_{\varphi_k}^n\geq \int_X (\varphi_j-\gamma) \,\theta_{\varphi_k}^n\geq -4A.
	\]
	Letting $k\to \infty$, by \cref{thm:semicon}, we find
	\[
	\int_X \Bigl(\varphi_j \lor (\gamma-C)-\gamma\Bigr)\,\theta_{\varphi}^n\geq -4A.
	\]
	Letting $j\to \infty$, by the monotone convergence theorem, we get
	\[
	\int_X \Bigl(\varphi \lor (\gamma-C)-\gamma\Bigr)\,\theta_{\varphi}^n\geq -4A.
	\]
	Then we let $C\to \infty$, again by the monotone convergence theorem, we conclude \eqref{eq:intvarphigammageqm4A_temp1}.

	\textbf{Step~2}. In general, let
    \[
    \psi_j=\uscsup[k\geq j] \varphi_k.
    \]
    Then $\varphi=\inf_{j>0}\psi_j$ due to \cref{cor:L1limipp}.

	For each $C>0$, let
	\[
	\psi_{j,C}=\psi_j \lor (\gamma-C),\quad \varphi_C=\varphi \lor (\gamma-C).
	\]
    Then $\psi_{j,C},\varphi_C\in \mathcal{E}(X,\theta;\phi)$ since $\gamma-C\in \mathcal{E}(X,\theta;\phi)$ thanks to the monotonicity theorem \cref{thm:mono}.
	Observe that $\psi_{j,C}$ decreases to $\varphi_C$ as $j\to \infty$. Moreover,
    \[
    \gamma\geq \psi_{j,C} \geq \psi_j\geq \varphi_j.
    \]
	By \cref{lma:fundam_estimate},
	\[
	\int_X (\psi_{j,C}-\gamma)\,\theta_{\psi_{j,C}}^n\geq -2^{n+1}A.
	\]
	By Step~1,
	\begin{equation}\label{eq:inproof1}
		\int_X(\varphi_C-\gamma)\,\theta_{\varphi_C}^n\geq -2^{n+3}A.
	\end{equation}
	In particular, by the plurifine locality of the non-pluripolar product \itemref{itm:prop:npp1:1},
	\[
	    \int_{\{\varphi>\gamma-C\}}(\varphi-\gamma)\,\theta_{\varphi}^n\geq -2^{n+3}A.
	\]
	Letting $C\to \infty$, we conclude \eqref{eq:intvarphigammageqm2np3A}. As a consequence,
    \[
        \int_X (\gamma-\varphi)\,\theta_{\varphi}^n<\infty.
    \]
	We still have to prove that $\varphi\in \mathcal{E}(X,\theta;\phi)$. In fact, by \eqref{eq:inproof1}, for any $C>0$,
	\[
	\int_{\{\varphi\leq \gamma-C\}}\,\theta_{\varphi_C}^n\leq \frac{1}{C}\int_X (\gamma-\varphi_C)\,\theta_{\varphi_C}^n\leq 2^{n+3}C^{-1}A.
	\]
	Using the monotonicity theorem \cref{thm:mono} and the plurifine locality of the non-pluripolar product \itemref{itm:prop:npp1:1}, we find
    \[
    \int_X \theta_{\phi}^n=\int_X \theta_{\varphi_C}^n=\int_{\{\varphi\leq \gamma-C\}}\,\theta_{\varphi_C}^n+\int_{\{\varphi>\gamma-C\}}\,\theta_{\varphi}^n.
    \]
    Letting $C\to\infty$, we conclude that
    \[
    \int_X \theta_{\varphi}^n=\int_X \theta_{\phi}^n.
    \]
    This completes the proof.
\end{proof}

\begin{proposition}\label{prop:relrooftopclosed}
    Assume that $\varphi,\psi\in \mathcal{E}(X,\theta;\phi)$ (resp. $\mathcal{E}^1(X,\theta;\phi)$, $\mathcal{E}^{\infty}(X,\theta;\phi)$), then so is $\varphi\land \psi$.
\end{proposition}
\begin{proof}
    The case of $\mathcal{E}^{\infty}(X,\theta;\phi)$ is trivial.

    We consider the case $\mathcal{E}(X,\theta;\phi)$. It follows from \cref{prop:landfinitecond1} that $\varphi\land \psi\in \PSH(X,\theta)$. By \cref{thm:diamond}, we have
    \[
        \int_X \theta_{\varphi\land\psi}^n\geq \int_X \theta_{\phi}^n.
    \]
    By the monotonicity theorem \cref{thm:mono}, equality holds. By \cref{thm:Pvarphidiffdef}, we conclude that
    \[
        P_{\theta}[\varphi\land \psi]=\phi.
    \]

    Finally, assume that $\varphi,\psi\in \mathcal{E}^1(X,\theta;\phi)$.
    We may assume that $\varphi,\psi\leq \phi$. For each $j\geq 1$, consider the approximations:
	\[
	\varphi_j\coloneqq \varphi\lor(\phi-j)\,,\quad \psi_j\coloneqq \psi\lor(\phi-j).
	\]
	By \cref{lma:rooftop_MAequation}, we can take $\gamma_j\in \mathcal{E}^{\infty}(X,\theta;\phi)$ solving the following equation:
	\[
	\theta_{\gamma_j}^n=\mathrm{e}^{\gamma_j-\varphi_j}\,\theta_{\varphi_j}^n+\mathrm{e}^{\gamma_j-\psi_j}\,\theta_{\psi_j}^n.
	\]
    In particular,
    \[
        \theta_{\gamma_j}^n\geq \mathrm{e}^{\gamma_j-\varphi_j}\,\theta_{\varphi_j}^n.
    \]
   We claim that this implies $\gamma_j\leq \varphi_j$. Indeed, by \cref{prop:compa_princ},
    \[
        \int_{\{\varphi_j<\gamma_j\}}\theta_{\varphi_j}^n
        \geq
        \int_{\{\varphi_j<\gamma_j\}}\theta_{\gamma_j}^n
        \geq
        \int_{\{\varphi_j<\gamma_j\}}\mathrm{e}^{\gamma_j-\varphi_j}\,\theta_{\varphi_j}^n.
    \]
    Hence
    \[
        \int_{\{\varphi_j<\gamma_j\}}\left(\mathrm{e}^{\gamma_j-\varphi_j}-1\right)\theta_{\varphi_j}^n=0.
    \]
    It follows that $\int_{\{\varphi_j<\gamma_j\}}\theta_{\varphi_j}^n=0$, and the domination principle \cref{thm:dom_prin_1} gives $\gamma_j\leq \varphi_j$. Similarly, $\gamma_j\leq \psi_j$. In particular, $\gamma_j\leq \varphi_j\land \psi_j$.
	We claim that
	\begin{equation}\label{eq:intgammajphi_unif_lower}
		\int_X (\gamma_j-\phi)\,\theta_{\gamma_j}^n >-C
	\end{equation}
    for some $C$ independent of $j$.

	Assume the claim is true for now.
	We get immediately that
    \[
    \sup_{X}\gamma_j=\sup_{X\setminus \{\phi=-\infty\}}(\gamma_j-\phi)\geq -\frac{C}{\int_X \theta_{\phi}^n},
    \]
    where the first equality follows from \cref{prop:model_sup_prop}.
    By \cref{prop:L1compa}, after possibly extracting a subsequence, we may assume that $\gamma_j\to \gamma\in \PSH(X,\theta)$ in $L^1$-topology.
    Then $\gamma\in \mathcal{E}^1(X,\theta;\phi)$ by \cref{lma:stab} and \cref{prop:E1charac}. Since $\gamma\leq\varphi\land\psi\leq\phi$, the monotonicity theorem \cref{thm:mono} gives
    \[
        \int_X\theta_{\varphi\land\psi}^n=\int_X\theta_\phi^n.
    \]
    Moreover, by the monotonicity of $E^\phi_\theta$,
    \[
        E^\phi_\theta(\varphi\land\psi)\geq E^\phi_\theta(\gamma)>-\infty.
    \]
    Thus $\varphi\land\psi\in\mathcal E^1(X,\theta;\phi)$ by \cref{prop:E1charac}.

	Now we prove the claim \eqref{eq:intgammajphi_unif_lower}.
	By symmetry, it suffices to prove
	\[
	\int_X (\phi-\gamma_j) \mathrm{e}^{\gamma_j-\varphi_j}\,\theta_{\varphi_j}^n\leq C.
	\]
	But note that
	\[
	\int_X (\phi-\gamma_j) \mathrm{e}^{\gamma_j-\varphi_j}\,\theta_{\varphi_j}^n=\int_X (\phi-\varphi_j) \mathrm{e}^{\gamma_j-\varphi_j}\,\theta_{\varphi_j}^n+\int_X (\varphi_j-\gamma_j) \mathrm{e}^{\gamma_j-\varphi_j}\,\theta_{\varphi_j}^n.
	\]
	But $x \mathrm{e}^{-x}\leq \mathrm{e}^{-1}$ for $x\geq 0$ (with maximum at $x=1$), so the second term is bounded; it remains to prove
	\[
	\int_X (\phi-\varphi_j) \mathrm{e}^{\gamma_j-\varphi_j}\,\theta_{\varphi_j}^n\leq C.
	\]
	As $\gamma_j\leq \varphi_j$, it suffices to prove
	\begin{equation}\label{eq:inproof2}
		\int_X (\phi-\varphi_j)\,\theta_{\varphi_j}^n\leq C.
	\end{equation}
    We compute, using \eqref{eq:Edefbdd} and the fact that the $k<n$ terms
    $\int_X(\varphi_j-\phi)\theta_{\varphi_j}^k\wedge\theta_\phi^{n-k}$ are all non-positive (since $\varphi_j\leq \phi$),
    \[
        \int_X (\varphi_j-\phi)\,\theta_{\varphi_j}^n\geq (n+1)E_{\theta}^{\phi}(\varphi_j).
    \]
	Since $\varphi_j\searrow\varphi$ and $E_\theta^{\phi}$ is continuous along decreasing sequences (\cref{prop:E_cont_decnet}), $E_\theta^{\phi}(\varphi_j)$ is bounded below uniformly in $j$. Thus, \eqref{eq:inproof2} follows.
\end{proof}

\begin{proposition}\label{prop:relativeEupperclosed}
    Let $\varphi,\psi\in \PSH(X,\theta)$ be potentials such that $\varphi\preceq \psi \preceq \phi$. Assume that $\varphi\in \mathcal{E}(X,\theta;\phi)$ (resp. $\mathcal{E}^1(X,\theta;\phi)$, $\mathcal{E}^{\infty}(X,\theta;\phi)$), then so is $\psi$.
\end{proposition}
\begin{proof}
    We may assume that $\varphi\leq \psi$.

    The case $\mathcal{E}^{\infty}(X,\theta;\phi)$ is trivial. The case $\mathcal{E}(X,\theta;\phi)$ follows from \cref{thm:mono}. The case $\mathcal{E}^1(X,\theta;\phi)$ follows from the characterization of $\mathcal{E}^1(X,\theta;\phi)$ in \cref{prop:E1charac} and the monotonicity of $E_{\theta}^{\phi}$ proved in \cref{lma:EmonoEinfty}.
\end{proof}

\begin{proposition}\label{prop:supsEE1}
    Let $(\varphi_i)_{i\in I}$ be a uniformly bounded from above non-empty family in $\mathcal{E}(X,\theta;\phi)$ (resp. $\mathcal{E}^1(X,\theta;\phi)$, $\mathcal{E}^{\infty}(X,\theta;\phi)$), then so is $\uscsup[i\in I]\varphi_i$.
\end{proposition}
\begin{proof}
Thanks to \cref{prop:relativeEupperclosed}, it suffices to show that
\[
\uscsup[i\in I]\varphi_i\preceq \phi.
\]
Since $\phi$ is model and $\varphi_i\preceq \phi$, we know that
\[
\varphi_i-\sup_X \varphi_i \leq \phi
\]
for any $i\in I$.
Since $(\varphi_i)_{i\in I}$ is uniformly bounded from above by assumption, our assertion follows.
\end{proof}

\begin{proposition}\label{prop:envrelfullmass}
    Let $\varphi,\psi\in \mathcal{E}(X,\theta;\phi)$. Then
    \[
        P_{\theta}[\varphi](\psi)=\uscsup[C > 0] (\varphi+C)\land \psi=\psi.
    \]
\end{proposition}
\begin{proof}
    The first equality follows from the definition of $P_{\theta}[\varphi](\psi)$. We focus on the second equality.

    Since for each $C\geq 0$,
    \[
        (\varphi\land \psi+C)\land \psi\leq (\varphi+C)\land \psi\leq \psi,
    \]
    we may replace $\varphi$ by $\varphi\land \psi$ (cf. \cref{prop:relrooftopclosed}) and assume that $\varphi\leq \psi$.

    Let
    \[
        \gamma\coloneqq \uscsup[C > 0] (\varphi+C)\land \psi.
    \]
    Observe that $\gamma\leq \psi$, and $\gamma\in \mathcal{E}(X,\theta;\phi)$ due to \cref{prop:relrooftopclosed} and \cref{prop:supsEE1}.
    Then \cref{cor:Pvarphisupport} guarantees that
    \[
    \int_{\{\gamma<\psi\} } \theta_{\gamma}^n=0.
    \]
    Therefore, we can apply the domination principle \cref{thm:dom_prin_1} to conclude that $\gamma=\psi$.
\end{proof}

\begin{lemma}\label{lma:Ediff}
    Let $\varphi,\psi\in \mathcal{E}^{\infty}(X,\theta;\phi)$. Define
    \[
    \eta_t\coloneqq \Bigl( (1-t)\varphi+t\psi \Bigr)\land \psi,\quad t\in [0,1].
    \]
    Then $E_{\theta}^{\phi}(\eta_t)$ is differentiable for $t\in [0,1]$ and
    \begin{equation}\label{eq:Iphit_diff}
    \frac{\mathrm{d}}{\mathrm{d}t}E_{\theta}^{\phi}(\eta_t)=\int_{X} \Bigl( \psi-\min\{\varphi,\psi\} \Bigr)\,\theta_{\eta_t}^n.
    \end{equation}
\end{lemma}
    The derivatives at $t=0$ or $1$ are understood as one-sided derivatives.
\begin{proof}
    Let us prove \eqref{eq:Iphit_diff} with right-derivative instead of derivative, and $t\in [0,1)$. The left-derivative case is completely parallel.

    For each $t\in [0,1]$, we let
    \[
    f_t=\min \{ (1-t)\varphi+t\psi, \psi\}.
    \]

    Fix $t\in [0,1)$ and $s>0$ small enough so that $t+s<1$. Thanks to \cref{prop:Ediffbdd}, we have
    \[
    \begin{aligned}
    E_{\theta}^{\phi}(\eta_{t+s})-E_{\theta}^{\phi}(\eta_{t})\leq & \int_X \left(\eta_{t+s}-\eta_t\right)\,\theta_{\eta_t}^n \\
    =& \int_X \left(\eta_{t+s}-f_t\right)\,\theta_{\eta_t}^n &\textup{by \cref{lma:rooftopMA}}\\
    \leq & \int_X \left(f_{t+s}-f_t \right)\,\theta_{\eta_t}^n\\
    =& s\int_X  \Bigl(\psi-\min\{\varphi,\psi\}\Bigr)\,\theta_{\eta_t}^n.
    \end{aligned}
    \]
    Similarly, we get
    \[
    E_{\theta}^{\phi}(\eta_{t+s})-E_{\theta}^{\phi}(\eta_{t})\geq s\int_X  \Bigl(\psi-\min\{\varphi,\psi\}\Bigr)\,\theta_{\eta_{t+s}}^n.
    \]
    Observe that $\eta_{t+s}$ converges in capacity (in fact, even uniformly) to $\eta_t$ as $s\to 0+$. The desired result \eqref{eq:Iphit_diff} is then just a consequence of \cref{thm:loc_conv_cap}.
\end{proof}

\begin{lemma}\label{lma:Eaffinediff}
    Let $\varphi,\psi\in \mathcal{E}^{1}(X,\theta;\phi)$, and set
    \[
        \varphi_t\coloneqq (1-t)\varphi+t\psi,\qquad t\in[0,1].
    \]
    Then the function
    \[
        t\mapsto E_{\theta}^{\phi}(\varphi_t)
    \]
    is differentiable on $[0,1]$, where the derivative at the endpoints is understood as the one-sided derivative. Moreover,
    \begin{equation}\label{eq:Eaffinediff}
        \frac{\mathrm d}{\mathrm dt}E_{\theta}^{\phi}(\varphi_t)
        =
        \int_X(\psi-\varphi)\,\theta_{\varphi_t}^n.
    \end{equation}
\end{lemma}
\begin{proof}
    By \cref{cor:E_classes_convex}, $\varphi_t\in \mathcal E^1(X,\theta;\phi)$ for all $t\in[0,1]$.
    For $a=0,\ldots,n$, put
    \[
        M_a=\int_X(\psi-\varphi)\,\theta_{\psi}^a\wedge\theta_{\varphi}^{n-a}.
    \]
    Note that each $M_a$ is finite by \cref{cor:mix_MA_energy_finite}. By the multi-linearity of the non-pluripolar product \itemref{itm:prop:npp1:6}, the quantity
    \[
        \int_X(\psi-\varphi)\,\theta_{\varphi_t}^n
    \]
    is therefore a polynomial in $t$ with coefficients given by the $M_a$'s. In particular, the right-hand side of \eqref{eq:Eaffinediff} is finite and continuous in $t$.

    We prove the formula for the right derivative at a point $t\in[0,1)$. Let $s>0$ be small enough so that $t+s\leq 1$. By the cocycle formula \cref{prop:cocycE1}, we have
    \[
    \begin{aligned}
        \frac{E_{\theta}^{\phi}(\varphi_{t+s})-E_{\theta}^{\phi}(\varphi_t)}{s}
        =&\frac{1}{(n+1)s}\sum_{j=0}^n
        \int_X(\varphi_{t+s}-\varphi_t)\,
        \theta_{\varphi_{t+s}}^j\wedge\theta_{\varphi_t}^{n-j} \\
        =&\frac{1}{n+1}\sum_{j=0}^n
        \int_X(\psi-\varphi)\,
        \theta_{\varphi_{t+s}}^j\wedge\theta_{\varphi_t}^{n-j}.
    \end{aligned}
    \]
    Again by \itemref{itm:prop:npp1:6}, for each fixed $j$ the last integral is a polynomial in $s$, expressed as a linear combination of the $M_a$'s with coefficients depending on $t$. Letting $s\to 0+$, it converges to
    \[
        \int_X(\psi-\varphi)\,\theta_{\varphi_t}^n.
    \]
    Therefore,
    \[
        \lim_{s\to0+}
        \frac{E_{\theta}^{\phi}(\varphi_{t+s})-E_{\theta}^{\phi}(\varphi_t)}{s}
        =
        \int_X(\psi-\varphi)\,\theta_{\varphi_t}^n.
    \]

    The proof for the left derivative is identical, using $s<0$ with $t+s\in[0,1]$. Thus the left and right derivatives agree for $t\in(0,1)$, and the endpoint derivatives are given by the same formula.
\end{proof}

In particular, we have a diamond-like inequality.
\begin{corollary}\label{cor:diam_E}
    Let $\varphi,\psi\in \mathcal{E}^1(X,\theta;\phi)$. Then
    \begin{equation}\label{eq:diam_E}
    E_{\theta}^{\phi}(\varphi\lor \psi)+E_{\theta}^{\phi}(\varphi\land \psi)\geq E_{\theta}^{\phi}(\psi)+E_{\theta}^{\phi}(\varphi).
    \end{equation}
\end{corollary}
\begin{proof}
    \textbf{Step~1}. We first reduce to the case where $\varphi,\psi \in \mathcal{E}^{\infty}(X,\theta;\phi)$.

    Assume that \eqref{eq:diam_E} is known in that case. For each $C>0$, we let
    \[
    \varphi_C=\varphi\lor (\phi-C),\quad \psi_C=\psi\lor (\phi-C).
    \]
    Then
    \[
    E_{\theta}^{\phi}(\varphi_C\lor \psi_C)+E_{\theta}^{\phi}(\varphi_C\land \psi_C)\geq E_{\theta}^{\phi}(\psi_C)+E_{\theta}^{\phi}(\varphi_C).
    \]
    As $C\to\infty$, we have
    \[
        \varphi_C\searrow\varphi,\quad \psi_C\searrow\psi,\quad \varphi_C\lor\psi_C\searrow\varphi\lor\psi,
    \]
    and, by the monotonicity and maximality of the rooftop,
    \[
        \varphi_C\land\psi_C\searrow\varphi\land\psi.
    \]
    Applying \cref{prop:E_cont_decnet} to these four decreasing sequences gives \eqref{eq:diam_E}.

    \textbf{Step~2}. We assume that $\varphi,\psi \in \mathcal{E}^{\infty}(X,\theta;\phi)$.

    Thanks to \cref{lma:Eaffinediff} and the plurifine locality of the non-pluripolar product \itemref{itm:prop:npp1:1}, we have
    \[
    \begin{aligned}
    E_{\theta}^{\phi}(\varphi\lor \psi)-E_{\theta}^{\phi}(\varphi)=&\int_0^1\int_X (\varphi\lor \psi-\varphi)\,\theta_{(1-t)\varphi+t(\varphi\lor \psi)}^n\,\mathrm{d}t\\
    =& \int_0^1 \int_{\{\psi>\varphi\}} (\psi-\varphi)\,\theta_{(1-t)\varphi+t\psi}^n\,\mathrm{d}t.
    \end{aligned}
    \]
    On the other hand,
    \[
    \begin{aligned}
        &E_{\theta}^{\phi}(\psi)-E_{\theta}^{\phi}(\varphi\land \psi)\\
        =& \int_0^1 \int_X (\psi-\min\{\varphi,\psi\})\,\theta_{((1-t)\varphi+t\psi)\land \psi}^n\,\mathrm{d}t\\
        \leq & \int_0^1 \int_{\{\varphi<\psi\}} (\psi-\varphi)\,\theta_{(1-t)\varphi+t\psi}^n\,\mathrm{d}t& \textup{by \cref{lma:rooftopMA}}.
    \end{aligned}
    \]
    Adding these equations together, we conclude \eqref{eq:diam_E}.
\end{proof}

\section{The \texorpdfstring{$\mathcal{I}$}{I}-envelope}\label{sec:Ienv}
From the algebraic point of view, a more natural envelope operator is given by the $\mathcal{I}$-envelope.

In this section, $X$ will denote a connected compact Kähler manifold of dimension $n$.

\subsection{\texorpdfstring{$\mathcal{I}$}{I}-equivalence}
\label{subsec:Iequiv}

\begin{proposition}\label{prop:Iequivchar}
    Given $\varphi,\psi\in \QPSH(X)$, the following are equivalent:
    \begin{enumerate}
        \item\label{itm:prop:Iequivchar:1} For any $k\in \mathbb{Z}_{>0}$, we have
        \[
            \mathcal{I}(k\varphi)=\mathcal{I}(k\psi);
        \]
        \item\label{itm:prop:Iequivchar:2} for any $\lambda\in \mathbb{R}_{>0}$, we have
            \[
            \mathcal{I}(\lambda\varphi)=\mathcal{I}(\lambda\psi);
            \]
        \item\label{itm:prop:Iequivchar:3} for any modification $\pi\colon Y\rightarrow X$ and any $y\in Y$, we have
            \[
                \nu(\pi^*\varphi,y)=\nu(\pi^*\psi,y);
            \]
        \item\label{itm:prop:Iequivchar:4} for any proper bimeromorphic morphism $\pi\colon Y\rightarrow X$ from a Kähler manifold and any $y\in Y$, we have
            \[
                \nu(\pi^*\varphi,y)=\nu(\pi^*\psi,y);
            \]
        \item\label{itm:prop:Iequivchar:5} for any prime divisor $E$ over $X$, we have
            \[
                \nu(\varphi,E)=\nu(\psi,E).
            \]
    \end{enumerate}
\end{proposition}
See \cref{def:primdiv} for the definition of prime divisors over $X$. Recall that in the whole book, a \emph{modification} of a compact complex space means a finite composition of blow-ups with smooth centers. This terminology is highly non-standard.

For simplicity, we usually write
\[
\nu(\varphi,y)=\nu(\pi^*\varphi,y)
\]
when $\pi$ is obvious from the context.

\begin{proof}
    \ref{itm:prop:Iequivchar:4} $\iff$ \ref{itm:prop:Iequivchar:5}. For the forward direction, use a smooth Kähler model realizing $E$ by \cref{cor:primerealization}. For the converse, use \cref{lma:blowupLelong} at each positive-dimensional point. At isolated points the Lelong numbers are zero.

    \ref{itm:prop:Iequivchar:3} $\iff$ \ref{itm:prop:Iequivchar:5}. This follows from \cref{cor:primerealization}.

    \ref{itm:prop:Iequivchar:1} $\implies$ \ref{itm:prop:Iequivchar:5}. This follows from \cref{prop:Lelongnumfrommis}.

    \ref{itm:prop:Iequivchar:5} $\implies$ \ref{itm:prop:Iequivchar:2}. This follows from \cref{thm:valuativemulti}.

    \ref{itm:prop:Iequivchar:2} $\implies$ \ref{itm:prop:Iequivchar:1}. This is the restriction to positive integral values of $\lambda$.
\end{proof}

\begin{definition}\label{def:Iequiv}
    Given $\varphi,\psi\in \QPSH(X)$, we say they are \emph{$\mathcal{I}$-equivalent}\index{I-equivalent@$\mathcal{I}$-equivalence} and write $\varphi\sim_{\mathcal{I}}\psi$\index{$\sim_{\mathcal{I}}$} if the equivalent conditions in \cref{prop:Iequivchar} are satisfied.
\end{definition}
Clearly, $\sim_{\mathcal{I}}$ is an equivalence relation on $\QPSH(X)$.

\begin{proposition}\label{prop:Ienvbimero}
    Let $\pi\colon Y\rightarrow X$ be a proper bimeromorphic morphism from a Kähler manifold $Y$ to $X$. For $\varphi,\psi\in \QPSH(X)$, the following are equivalent:
    \begin{enumerate}
        \item\label{itm:prop:Ienvbimero:1} $\varphi\sim_{\mathcal{I}}\psi$;
        \item\label{itm:prop:Ienvbimero:2} $\pi^*\varphi\sim_{\mathcal{I}}\pi^*\psi$.
    \end{enumerate}
\end{proposition}
\begin{proof}
    \ref{itm:prop:Ienvbimero:1} $\implies$ \ref{itm:prop:Ienvbimero:2}. This follows from \itemref{itm:prop:Iequivchar:4}.

    \ref{itm:prop:Ienvbimero:2} $\implies$ \ref{itm:prop:Ienvbimero:1}. This follows from the simple fact that
    \[
        \mathcal{I}(k\varphi)=\pi_* \left( \omega_{Y/X}\otimes \mathcal{I}(k\pi^*\varphi)  \right),\quad \mathcal{I}(k\psi)=\pi_* \left( \omega_{Y/X}\otimes \mathcal{I}(k\pi^*\psi)  \right)
    \]
    for any $k\in \mathbb{Z}_{>0}$, where $\omega_{Y/X}$ is the relative dualizing sheaf. See \cite[Proposition~5.8]{Dem12} and \cite[Proposition~12.6.4]{FPT}.
\end{proof}

\begin{proposition}\label{prop:Iequivmax}
    Let $\varphi,\varphi',\psi,\psi'\in \QPSH(X)$ and $\lambda>0$. Assume that $\varphi\sim_{\mathcal{I}}\psi$ and $\varphi'\sim_{\mathcal{I}}\psi'$. Then
    \[
        \varphi\lor \varphi'\sim_{\mathcal{I}} \psi\lor\psi',\quad \varphi+\varphi'\sim_{\mathcal{I}} \psi+\psi',\quad \lambda\varphi\sim_{\mathcal{I}}\lambda\psi.
    \]
    Similarly, if $(\varphi_i)_{i\in I}$, $(\psi_i)_{i\in I}$ are two non-empty uniformly bounded from above families in $\PSH(X,\theta)$ for some closed smooth real $(1,1)$-form $\theta$ on $X$ such that $\varphi_i\sim_{\mathcal{I}}\psi_i$ for all $i\in I$, then
    \[
    \uscsup[i\in I] \varphi_i\sim_{\mathcal{I}}\uscsup[i\in I] \psi_i.
    \]
\end{proposition}
\begin{proof}
    This follows from \cref{prop:Lelongmax} and \cref{cor:supsLelong}.
\end{proof}

\subsection{The definition of the \texorpdfstring{$\mathcal{I}$}{I}-envelope}
\label{subsec:Ienv_def}
Fix a smooth closed real $(1,1)$-form $\theta$ on $X$.
\begin{definition}\label{def:Ienv}
    Given $\varphi\in \PSH(X,\theta)$, we define its \emph{$\mathcal{I}$-envelope}\index{envelope!I-envelope@$\mathcal{I}$-envelope} as follows:
    \begin{equation}\label{eq:Ienvelopedef}
        P_{\theta}[\varphi]_{\mathcal{I}}\coloneqq \uscsup\{\psi\in \PSH(X,\theta): \psi\leq 0,\psi\sim_{\mathcal{I}} \varphi\}.
    \end{equation}
    \index{$P_{\theta}[\varphi]_{\mathcal{I}}$}
    If $\varphi=P_{\theta}[\varphi]_{\mathcal{I}}$, we say $\varphi$ is an \emph{$\mathcal{I}$-model potential}\index{potential!I-model potential@$\mathcal{I}$-model potential} (in $\PSH(X,\theta)$).
\end{definition}
Note that by \cref{prop:pshfunction_closedseq}, $P_{\theta}[\varphi]_{\mathcal{I}}\in \PSH(X,\theta)$.

\begin{proposition}\label{prop:Ienvindeptheta}
    Let $\theta'=\theta+\ddc g$ for some $g\in C^{\infty}(X)$. Then for any $\varphi\in \PSH(X,\theta)$, we have $\varphi'\coloneqq \varphi-g\in \PSH(X,\theta')$ and
    \[
        P_{\theta}[\varphi]_{\mathcal{I}}\sim P_{\theta'}[\varphi']_{\mathcal{I}}.
    \]
\end{proposition}
\begin{proof}
    By symmetry, it suffices to show that
    \[
        P_{\theta}[\varphi]_{\mathcal{I}}\preceq P_{\theta'}[\varphi']_{\mathcal{I}}.
    \]
    After replacing $g$ by $g+C$ for a suitable constant $C$, which only changes $\varphi'$ by an additive constant, we may assume that $g\geq 0$.
    Let $\psi\in \PSH(X,\theta)$ be such that $\psi\leq 0$ and $\psi\sim_{\mathcal{I}}\varphi$. Then $\psi'\coloneqq \psi-g$ lies in $\PSH(X,\theta')$, satisfies $\psi'\leq 0$, and is $\mathcal{I}$-equivalent to $\varphi'$ since adding a smooth function does not change multiplier ideals. Hence
    \[
        \psi-g=\psi'\leq P_{\theta'}[\varphi']_{\mathcal{I}}.
    \]
    It follows that
    \[
        \psi-\sup_X g\leq P_{\theta'}[\varphi']_{\mathcal{I}}.
    \]
    Taking the upper envelope over all such $\psi$ gives
    \[
        P_{\theta}[\varphi]_{\mathcal{I}}-\sup_X g\leq P_{\theta'}[\varphi']_{\mathcal{I}},
    \]
    which proves the desired comparison.
\end{proof}

\begin{proposition}\label{prop:Ienvelopebimero}
    Let $\pi\colon Y\rightarrow X$ be a proper bimeromorphic morphism from a Kähler manifold $Y$ to $X$. Then for $\varphi \in \PSH(X,\theta)$, we have
    \[
        P_{\pi^*\theta}[\pi^*\varphi]_{\mathcal{I}}=\pi^* P_{\theta}[\varphi]_{\mathcal{I}}.
    \]
\end{proposition}
\begin{proof}
    Let
    \[
        \eta=P_{\pi^*\theta}[\pi^*\varphi]_{\mathcal{I}}.
    \]
    First take a candidate $\psi\in \PSH(X,\theta)$ in \eqref{eq:Ienvelopedef}, so that $\psi\leq 0$ and $\psi\sim_{\mathcal{I}}\varphi$. Then $\pi^*\psi\leq 0$ and $\pi^*\psi\sim_{\mathcal{I}}\pi^*\varphi$ by \cref{prop:Ienvbimero}; hence $\pi^*\psi\leq \eta$. Taking the upper envelope over all such $\psi$ gives
    \[
        \pi^* P_{\theta}[\varphi]_{\mathcal{I}}\leq \eta.
    \]

    Conversely, let $\chi\in \PSH(Y,\pi^*\theta)$ be a candidate in the definition of $\eta$, so that $\chi\leq 0$ and $\chi\sim_{\mathcal{I}}\pi^*\varphi$. By \cref{prop:PSHpullbij}, there is a unique $\psi\in \PSH(X,\theta)$ such that $\chi=\pi^*\psi$. Since $\pi$ is surjective, $\psi\leq 0$. By \cref{prop:Ienvbimero}, the relation $\chi\sim_{\mathcal{I}}\pi^*\varphi$ implies $\psi\sim_{\mathcal{I}}\varphi$. Thus $\psi\leq P_{\theta}[\varphi]_{\mathcal{I}}$, and hence
    \[
        \chi\leq \pi^* P_{\theta}[\varphi]_{\mathcal{I}}.
    \]
    Taking the upper envelope over all such $\chi$ gives the reverse inequality.
\end{proof}

\begin{proposition}\label{prop:Ienvprojection}
    Let $\varphi\in \PSH(X,\theta)$. Then
    \[
        \varphi\sim_{\mathcal{I}}  P_{\theta}[\varphi]_{\mathcal{I}}.
    \]
    In particular,
    \[
        P_{\theta}\Bigl[ P_{\theta}[\varphi]_{\mathcal{I}} \Bigr]_{\mathcal{I}}=P_{\theta}[\varphi]_{\mathcal{I}}
    \]
    and the upper semicontinuous regularization in \eqref{eq:Ienvelopedef} is not necessary.
\end{proposition}
\begin{proof}
    In view of \cref{prop:Iequivchar}, it suffices to show that for $k \in \mathbb{Z}_{>0}$, we have
    \begin{equation}\label{eq:IenvelopepreservLelong}
        \mathcal{I}(k\varphi)=\mathcal{I}(kP_{\theta}[\varphi]_{\mathcal{I}}).
    \end{equation}

    By \cref{prop:Choquet}, we can find $\psi_i\in \PSH(X,\theta)$ ($i\in \mathbb{Z}_{>0}$) such that $\psi_i\leq 0$, $\psi_i\sim_{\mathcal{I}}\varphi$ for all $i\geq 1$ and
    \[
        \uscsup[i>0]\psi_i=P_{\theta}[\varphi]_{\mathcal{I}}.
    \]
    By \cref{prop:Iequivmax}, we may replace $\psi_i$ by $\psi_1\lor\cdots\lor\psi_i$ and assume that the sequence $\psi_i$ is increasing. In this case, it follows from the strong openness theorem \cref{thm:strongopen} that for each $k\in \mathbb{Z}_{>0}$, we have
    \[
        \mathcal{I}(k\varphi)=\mathcal{I}(k\psi_j)=\mathcal{I} (kP_{\theta}[\varphi]_{\mathcal{I}})
    \]
    for $j$ large enough.
\end{proof}

\begin{definition}\label{def:volqpsh}
    Let $\varphi\in \PSH(X,\theta)$, we define the \emph{volume}\index{volume}\footnote{We choose to call this quantity the \emph{volume} instead of the \emph{$\mathcal{I}$-volume} so that the terminology is consistent with the line bundle case.} $\vol(\theta,\varphi)$ as
    \[
        \vol(\theta,\varphi)=\int_X \Bigl(\theta+\ddc P_{\theta}[\varphi]_{\mathcal{I}}  \Bigr)^n.
    \]
    \index{$\vol(\theta,\varphi)$}
\end{definition}

\begin{proposition}\label{prop:voldeponlyoncurr}
    Let $\theta'=\theta+\ddc g$ for some $g\in C^{\infty}(X)$. Then for any $\varphi\in \PSH(X,\theta)$, we have $\varphi'=\varphi-g\in \PSH(X,\theta')$ and
    \[
        \vol(\theta,\varphi)=  \vol\left(\theta',\varphi'\right).
    \]
\end{proposition}
\begin{proof}
    By \cref{prop:Ienvindeptheta}, the currents $\theta+\ddc P_\theta[\varphi]_{\mathcal I}$ and $\theta'+\ddc P_{\theta'}[\varphi']_{\mathcal I}$ have the same singularity type in the same cohomology class. Applying the monotonicity theorem \cref{thm:mono} gives equality of their non-pluripolar masses.
\end{proof}

In view of \cref{prop:voldeponlyoncurr}, the volume $\vol(\theta,\varphi)$ depends only on the current $\theta_{\varphi}$, and we may write
\begin{equation}\label{eq:volcurrdef}
    \vol \theta_{\varphi}=\vol(\theta,\varphi).
\end{equation}
\index{$\vol \theta_{\varphi}$}

\begin{definition}\label{def:volcohclass}
    Let $\alpha\in \mathrm{H}^{1,1}(X,\mathbb{R})$ be a pseudo-effective class. The \emph{volume}\index{volume!volume of a pseudo-effective class} $\vol \alpha$ of $\alpha$ is defined as
    \[
    \vol \alpha=\vol T_{\min},
    \]
    where $T_{\min}$ is a current with minimal singularities in $\alpha$. Note that $\vol \alpha$ is independent of the choice of $T_{\min}$, thanks to \cref{thm:mono}.
\end{definition}

Let us recall the following elementary result for later use.
\begin{proposition}\label{prop:cones_modif2}
    Let $\pi\colon Y\rightarrow X$ be a proper bimeromorphic morphism from a Kähler manifold $Y$.
    \begin{enumerate}
        \item\label{itm:prop:cones_modif2:1} For any big class $\alpha\in \mathrm{H}^{1,1}(X,\mathbb{R})$, $\pi^*\alpha$ is big. Moreover, $\vol \pi^*\alpha=\vol \alpha$.
        \item\label{itm:prop:cones_modif2:2} For any big class $\beta\in \mathrm{H}^{1,1}(Y,\mathbb{R})$, $\pi_*\beta$ is big. Moreover, $\vol \pi_*\beta\geq \vol \beta$.
    \end{enumerate}
\end{proposition}
\begin{proof}
    \ref{itm:prop:cones_modif2:1} Take a current $T_{\min}$ with minimal singularities in $\alpha$. Then $\pi^*T_{\min}$ is a current with minimal singularities in $\pi^*\alpha$. Our assertion follows from \cref{prop:nppmassinv}.

    \ref{itm:prop:cones_modif2:2} Let $S_{\min}$ be a current with minimal singularities in $\beta$. Then $\pi_*S_{\min}$ is a positive current in the cohomology class $\pi_*\beta$. A current with minimal singularities in $\pi_*\beta$ is less singular than $\pi_*S_{\min}$, hence by the monotonicity theorem \cref{thm:mono}
    \[
        \vol(\pi_*\beta)\geq \int_X\left(\pi_*S_{\min}\right)^n.
    \]
    By \cref{prop:nppmassinv},
    \[
        \int_X \left(\pi_*S_{\min}\right)^n\geq \int_Y S_{\min}^n=\vol\beta.
    \]
    We conclude $\vol \pi_*\beta\geq \vol \beta$.
\end{proof}

The $\mathcal{I}$-envelope and the $P$-envelope are related in a simple manner.
\begin{proposition}\label{prop:PandPI}
    Let $\varphi\in \PSH(X,\theta)$. Then
    \[
        P_{\theta}[\varphi]\leq P_{\theta}[\varphi]_{\mathcal{I}},\quad \varphi \sim_{\mathcal{I}}P_{\theta}[\varphi].
    \]
\end{proposition}
As a consequence,
\begin{equation}\label{eq:PthetaIP_comp}
    P_{\theta}\Bigl[P_{\theta}[\varphi]\Bigr]_{\mathcal{I}}=P_{\theta}[\varphi]_{\mathcal{I}}.
\end{equation}
\begin{proof}
    It suffices to show that $\varphi\sim_{\mathcal{I}} P_{\theta}[\varphi]$. Namely, for each $k\in \mathbb{Z}_{>0}$, we have
    \begin{equation}\label{eq:IkvarphiIkP}
        \mathcal{I}(k\varphi)=\mathcal{I}\Bigl(kP_{\theta}[\varphi]\Bigr).
    \end{equation}
    Fix $k$ for now.
    It follows from \eqref{eq:Penvsups} and the strong openness theorem \cref{thm:strongopen} that
    \[
        \mathcal{I}\Bigl(kP_{\theta}[\varphi]\Bigr)=\mathcal{I}\Bigl((k\varphi+C)\land kV_{\theta}\Bigr),
    \]
    when $C$ is large enough. Since $(k\varphi+C)\land kV_{\theta}\sim k\varphi$, we have
    \[
        \mathcal{I}\Bigl((k\varphi+C)\land kV_{\theta}\Bigr)=\mathcal{I}(k\varphi)
    \]
    and \eqref{eq:IkvarphiIkP} follows.
\end{proof}

In particular, we obtain an interesting relation between the non-pluripolar mass and the volume.
\begin{corollary}\label{cor:compnppmassandvol}
    Let $\varphi\in \PSH(X,\theta)$. Then
    \[
        \int_X \theta_{\varphi}^n\leq \vol \theta_{\varphi}.
    \]
\end{corollary}
The reverse inequality fails in general, see \cref{ex:BBJ}.
\begin{proof}
    This follows from \cref{prop:PandPI}, \cref{thm:mono} and \cref{prop:Ppresmass}.
\end{proof}

We note the following special case:
\begin{proposition}\label{prop:analysingcompPandPI}
    Let $\varphi\in \PSH(X,\theta)$. Assume that $\varphi$ has analytic singularities. Then
    \begin{equation}\label{eq:ana_sing_phi_PI_eq}
        \varphi\sim P_{\theta}[\varphi] \sim P_{\theta}[\varphi]_{\mathcal{I}}.
    \end{equation}
    Hence, $P_{\theta}[\varphi]$ and $P_{\theta}[\varphi]_{\mathcal{I}}$ both have analytic singularities.

    In particular,
    \begin{equation}\label{eq:ana_sing_mass_vol}
        \int_X \theta_{\varphi}^n= \vol \theta_{\varphi}.
    \end{equation}
\end{proposition}
\begin{proof}
    First observe that \eqref{eq:ana_sing_mass_vol} follows from \eqref{eq:ana_sing_phi_PI_eq} and \cref{thm:mono}. It remains to establish \eqref{eq:ana_sing_phi_PI_eq}.

    In view of \cref{prop:PandPI}, it suffices to show that
    \begin{equation}\label{eq:Pprecvarphitemp1}
        P_{\theta}[\varphi]_{\mathcal{I}}\preceq \varphi.
    \end{equation}
    By \cref{prop:Ienvelopebimero,prop:Penvbimero} and \cref{thm:resolvelogsing}, after passing to a log resolution we may assume that $\varphi$ has log singularities along an effective $\mathbb Q$-divisor $D$. This reduction does not change the assertion, since both the $P$- and the $\mathcal I$-envelope are compatible with bimeromorphic pull-back.

    After multiplying $\theta$ and $\varphi$ by a sufficiently divisible positive integer, and then dividing the final comparison by the same integer, we may assume that $D$ is integral. In this case the Bergman-type quasi-equisingular approximants of $\varphi$ associated to the ideals $\mathcal I(2^j\varphi)$ have the same singularity type as $\varphi$.

    Since
    \[
        P_\theta[\varphi]_{\mathcal I}\sim_{\mathcal I}\varphi
    \]
    by \cref{prop:Ienvprojection}, we have
    \[
        \mathcal I(2^jP_\theta[\varphi]_{\mathcal I})=\mathcal I(2^j\varphi).
    \]
    Thus the quasi-equisingular approximants $\eta_j$ of $P_\theta[\varphi]_{\mathcal I}$ and $\varphi_j$ of $\varphi$ may be chosen with the same singularity type $(2^{-j},\mathcal I(2^j\varphi))$.
    Consequently $\eta_j\sim\varphi_j\sim\varphi$ for $j\gg1$. Since $\eta_j\searrow P_\theta[\varphi]_{\mathcal I}$, choosing one such large $j$ gives
    \[
        P_\theta[\varphi]_{\mathcal I}\leq \eta_j\leq \varphi+C_j,
    \]
    and hence \eqref{eq:Pprecvarphitemp1} holds.
\end{proof}

\subsection{Properties of the \texorpdfstring{$\mathcal{I}$}{I}-envelope}
\label{subsec:Ienv_prop}
Let $\theta,\theta_1,\theta_2$ be smooth closed real $(1,1)$-forms on $X$.

We have the following properties of the $\mathcal{I}$-envelope.
\begin{proposition}\label{prop:PIconc}
    \leavevmode
    \begin{enumerate}
        \item\label{itm:prop:PIconc:1} Suppose that $\varphi\in \PSH(X,\theta)$ and $\lambda\in \mathbb{R}_{>0}$. Then
            \[
                P_{\lambda\theta}[\lambda\varphi]_{\mathcal{I}}=\lambda P_{\theta}[\varphi]_{\mathcal{I}}.
            \]
        \item\label{itm:prop:PIconc:2} Suppose that $\varphi_1\in \PSH(X,\theta_1)$ and $\varphi_2\in \PSH(X,\theta_2)$. Then
            \[
                P_{\theta_1+\theta_2}[\varphi_1+\varphi_2]_{\mathcal{I}}\geq P_{\theta_1}[\varphi_1]_{\mathcal{I}}+P_{\theta_2}[\varphi_2]_{\mathcal{I}}.
            \]
        \item\label{itm:prop:PIconc:3}
        Suppose that $\varphi_1\in \PSH(X,\theta_1)$ and $\varphi_2\in \PSH(X,\theta_2)$. Then
            \[
                P_{\theta_1+\theta_2}[\varphi_1+\varphi_2]_{\mathcal{I}}\sim_{\mathcal{I}} P_{\theta_1}[\varphi_1]_{\mathcal{I}}+P_{\theta_2}[\varphi_2]_{\mathcal{I}}.
            \]
        \item\label{itm:prop:PIconc:4}
            Suppose that $\varphi_1,\varphi_2\in \PSH(X,\theta)$. Then
                \[
                    P_{\theta}[\varphi_1\lor \varphi_2]_{\mathcal{I}}\sim_{\mathcal{I}} P_{\theta}[\varphi_1]_{\mathcal{I}}\lor P_{\theta}[\varphi_2]_{\mathcal{I}}.
                \]
    \end{enumerate}
\end{proposition}
\begin{proof}
    \ref{itm:prop:PIconc:1} This follows by rescaling the candidates in the definition of the $\mathcal{I}$-envelope.

    \ref{itm:prop:PIconc:2} Suppose that $\psi_1\in \PSH(X,\theta_1)$ and $\psi_2\in \PSH(X,\theta_2)$ satisfy
    \[
        \psi_i\leq 0,\quad \psi_i \sim_{\mathcal{I}} \varphi_i
    \]
    for $i=1,2$. Then thanks to \cref{prop:Iequivmax},
    \[
        \psi_1+\psi_2\leq 0,\quad \psi_1+\psi_2 \sim_{\mathcal{I}} \varphi_1+\varphi_2.
    \]
    It follows that
    \[
        \psi_1+\psi_2\leq P_{\theta_1+\theta_2}[\varphi_1+\varphi_2]_{\mathcal{I}}.
    \]
    Since $\psi_1$ and $\psi_2$ are arbitrary, we conclude.

    \ref{itm:prop:PIconc:3} and \ref{itm:prop:PIconc:4}: These follow easily from \cref{prop:Ienvprojection} and \cref{prop:Iequivmax}.
\end{proof}

\begin{lemma}\label{lma:PIenvmono1}
    Let $\varphi,\psi\in \PSH(X,\theta)$. Assume that $\varphi\preceq \psi$. Then
    \[
        P_{\theta}[\varphi]_{\mathcal{I}}\leq P_{\theta}[\psi]_{\mathcal{I}}.
    \]
\end{lemma}
\begin{proof}
    We may assume that $\psi\leq 0$.
    Observe that $P_{\theta}[\varphi]_{\mathcal{I}}\lor \psi\sim_{\mathcal{I}} \psi$ as a consequence of \cref{prop:Lelongmax} and \cref{prop:Ienvprojection}. So
    \[
    P_{\theta}[\varphi]_{\mathcal{I}}\lor \psi\leq P_{\theta}[\psi]_{\mathcal{I}}.
    \]
    The desired inequality follows.
\end{proof}

\begin{proposition}\label{prop:decnetmodelPI}
    Consider a decreasing net $(\varphi_i)_{i\in I}$ of model potentials, with $\varphi_i\in \PSH(X,\theta)_{>0}$ for all $i\in I$. Set
    \[
    \varphi\coloneqq \inf_{i\in I}\varphi_i.
    \]
    Suppose that $\varphi\not\equiv -\infty$ and $\int_X \theta_{\varphi}^n>0$. Then
    \[
        \inf_{i\in I} P_{\theta}[\varphi_i]_{\mathcal{I}}=P_{\theta}[\varphi]_{\mathcal{I}}.
    \]
    In particular, if the $\varphi_i$'s are all $\mathcal{I}$-model, then so is $\varphi$.
\end{proposition}
\begin{proof}
    Let $\eta=\inf_{i\in I} P_{\theta}[\varphi_i]_{\mathcal{I}}$. We have $\eta\geq P_{\theta}[\varphi]_{\mathcal{I}}$ as a consequence of \cref{lma:PIenvmono1}. In particular, $\eta\in \PSH(X,\theta)$.

    By \cref{prop:vol_limit_model}, we have
    \[
        \lim_{i\in I} \int_X\theta_{\varphi_i}^n=\int_X\theta_{\varphi}^n>0.
    \]
    We may assume that for each $i\in I$, $\int_X\theta_{\varphi_i}^n>\int_X\theta_{\varphi}^n$, since otherwise $\varphi_i=\varphi$ by \cref{thm:Pvarphidiffdef}, and hence $P_{\theta}[\varphi_i]_{\mathcal{I}}=P_{\theta}[\varphi]_{\mathcal{I}}$ due to \eqref{eq:PthetaIP_comp}; the same holds for all $j\geq i$, and there is nothing to prove.

    So by \cref{lma:pathoenvelope}, we can find a decreasing net $\epsilon_i\searrow 0$ ($i\in I$) with $\epsilon_i\in (0,1)$ and $\psi_i\in \PSH(X,\theta)$ such that for all $i\in I$,
    \[
        (1-\epsilon_i) \varphi_i+\epsilon_i\psi_i\leq \varphi.
    \]
    We can, for example, take
    \[
        \epsilon_i^{-1}=\frac{1}{2}\left( 1+\left(\frac{\int_X \theta_{\varphi_i}^n}{\int_X \theta_{\varphi_i}^n-\int_X \theta_{\varphi}^n}\right)^{1/n} \right).
    \]

    By \cref{prop:PIconc} and \cref{lma:PIenvmono1}, we have
    \[
        \eta+\epsilon_i P_{\theta}[\psi_i]_{\mathcal{I}}\leq (1-\epsilon_i) \eta+\epsilon_i P_{\theta}[\psi_i]_{\mathcal{I}}\leq (1-\epsilon_i) P_{\theta}[\varphi_i]_{\mathcal{I}}+\epsilon_i P_{\theta}[\psi_i]_{\mathcal{I}}\leq P_{\theta}[\varphi]_{\mathcal{I}}.
    \]
    Each $P_{\theta}[\psi_i]_{\mathcal{I}}$ has supremum $0$, so \cref{prop:L1compa} shows that the family is bounded in $L^1$. Hence $\epsilon_iP_{\theta}[\psi_i]_{\mathcal{I}}\to 0$ in $L^1$, and the preceding inequality gives $\eta\leq P_{\theta}[\varphi]_{\mathcal{I}}$ almost everywhere. By \cref{prop:pshfuncdetdense}, the inequality holds everywhere.
\end{proof}

\begin{proposition}\label{prop:incnetmodelPI}
    Let $(\varphi_i)_{i\in I}$ be an increasing net in $\PSH(X,\theta)_{>0}$ uniformly bounded from above. Let $\varphi\coloneqq \uscsup[i\in I]\varphi_i$. Then
    \[
        \uscsup[i\in I] P_{\theta}[\varphi_i]_{\mathcal{I}}=P_{\theta}[\varphi]_{\mathcal{I}}.
    \]
    In particular, if the $\varphi_i$'s are all $\mathcal{I}$-model, then so is $\varphi$.
\end{proposition}
\begin{proof}
    Let $\eta=\uscsup[i\in I] P_{\theta}[\varphi_i]_{\mathcal{I}}$. Then $\eta\leq P_{\theta}[\varphi]_{\mathcal{I}}$ as a consequence of \cref{lma:PIenvmono1}.

    By \cref{cor:incseqnppcont}, we have
    \[
        \lim_{i\in I} \int_X\theta_{\varphi_i}^n=\int_X\theta_{\varphi}^n>0.
    \]
    We may assume that for each $i\in I$, $\int_X\theta_{\varphi_i}^n<\int_X\theta_{\varphi}^n$, since otherwise $P_{\theta}[\varphi_i]=P_{\theta}[\varphi]$ by \cref{thm:Pvarphidiffdef}, and hence $P_{\theta}[\varphi_i]_{\mathcal{I}}=P_{\theta}[\varphi]_{\mathcal{I}}$ due to \eqref{eq:PthetaIP_comp}; the same holds for all $j\geq i$, and there is nothing to prove.

    So by \cref{lma:pathoenvelope}, we can find a decreasing net $\epsilon_i\searrow 0$ ($i\in I$) as in \eqref{eq:epsilon_inetchoice} with $\epsilon_i\in (0,1)$ and $\psi_i\in \PSH(X,\theta)$ such that for all $i\in I$,
    \[
        (1-\epsilon_i) \varphi+\epsilon_i\psi_i\leq \varphi_i.
    \]
    By \cref{prop:PIconc} and \cref{lma:PIenvmono1}, we have
    \[
        P_{\theta}[\varphi]_{\mathcal{I}}+\epsilon_i P_{\theta}[\psi_i]_{\mathcal{I}}\leq (1-\epsilon_i) P_{\theta}[\varphi]_{\mathcal{I}}+\epsilon_i P_{\theta}[\psi_i]_{\mathcal{I}}\leq P_{\theta}[\varphi_i]_{\mathcal{I}}\leq \eta.
    \]
    As above, $\epsilon_iP_{\theta}[\psi_i]_{\mathcal{I}}\to 0$ in $L^1$. It follows that $P_{\theta}[\varphi]_{\mathcal{I}}\leq \eta$ almost everywhere, hence everywhere by \cref{prop:pshfuncdetdense}.
\end{proof}

\begin{remark}\label{rmk:envelopes:L2443}
    One could also define the following interpolation between the $\mathcal{I}$-envelope and the $P$-envelope: Suppose $\varphi\in \PSH(X,\theta)_{>0}$, $j\in \{0,\ldots,n\}$. Then we let
    \[
    \begin{split}
    P_{\theta,j}[\varphi]\coloneqq \uscsup \left\{ \psi\in \PSH(X,\theta): \psi\leq 0,\varphi\preceq \psi, \int_X \theta_{\varphi}^j\wedge \theta_{P_{\theta}[\varphi]_{\mathcal{I}}}^{n-j}\right.\\
    \left.=\int_X \theta_{\psi}^j\wedge \theta_{P_{\theta}[\psi]_{\mathcal{I}}}^{n-j} \right\}.
    \end{split}
    \]
    The cases $j=0$ and $j=n$ recover the $\mathcal{I}$- and $P$-envelopes, respectively. If $\varphi$ is $\mathcal{I}$-good, then $P_{\theta,j}[\varphi]=P_\theta[\varphi]$ for every $j$. Indeed, the mixed mass in the definition for $\varphi$ then equals $\int_X\theta_\varphi^n$. For a candidate $\psi$, its mixed mass is at least $\int_X\theta_\psi^n$ by monotonicity, and $\int_X\theta_\psi^n\geq\int_X\theta_\varphi^n$. Thus equality of mixed masses implies equality of top masses, so $\psi\leq P_\theta[\varphi]$ by \cref{thm:Pvarphidiffdef}. The potential $P_\theta[\varphi]$ is itself a candidate by mass preservation and \eqref{eq:PthetaIP_comp}.

    Is $P_{\theta,j}[\bullet]$ a projection operator for $0<j<n$ in general?
\end{remark}

\makeatletter
\newcommand*{\germanmotto}[1]{{%
    \font\germanmottofont="[TeXGyreTermesX-Italic.otf]" at \f@size pt\relax
    \germanmottofont #1%
}}
\makeatother

\motto{\germanmotto{In den Dreißiger Jahren besuchte ich regelmäßig die Schweiz, teils um mich auch auf den Viertausendern zu tummeln, zum großen Teil aber auch, um Emigrantenblätter zu lesen und mich mit Kollegen über Naziverbrechen zu unterhalten. Aber auch die Schweizer schauten sich, wenn sie offen reden wollten, ebenso ängstlich um wie das bei uns üblich war.}\footnote{The \href{https://www.universityworldnews.com/post.php?story=20241107142108558}{recent policy of ETH} against Chinese students makes me feel that nothing has changed in Switzerland after the collapsing of Nazi for almost 80 years.}

--- Oskar Perron\footnote{Oskar Perron (1880--1975), after earning an \emph{Eisernes Kreuz} during WWI, obtained a position in München in 1922, initiating the glorious period there. Among his colleagues were Carathéodory, Tietze and Sommerfield.}}
\chapter{Geodesic rays in the space of potentials}
\label{chap:rays}

When $\omega$ is a K\"ahler form, the space of smooth strictly $\omega$-plurisubharmonic potentials carries the Mabuchi $L^2$-metric. Its geodesic equation is the homogeneous complex Monge--Amp\`ere equation on $X\times\{0<\Real z<1\}$. Subsolutions of this equation are subgeodesics, and the geodesic between two boundary potentials is obtained as a Perron envelope. This construction continues to make sense for singular potentials, even when no smooth geodesic exists.

The goal of this chapter is to further develop the corresponding theory for $\theta$-plurisubharmonic potentials, where $\theta$ is a smooth closed real $(1,1)$-form representing a big class. We define subgeodesics by requiring the complexification $\Phi(x,z)=\varphi_{\Real z}(x)$ to be plurisubharmonic. Geodesics are then defined as Perron envelopes of subgeodesics with prescribed boundary values. This point of view is flexible enough to handle unbounded potentials and singularity constraints.

The first section sets up the basic properties of subgeodesics. We then construct geodesics with prescribed singular endpoints and study geodesic rays. The results are formulated relative to a model potential $\phi$, so that the endpoints and rays live in a fixed relative singularity class.

The final part of the chapter studies the finite-energy space $\mathcal E^1(X,\theta;\phi)$ and the Darvas $d_1$-metric. We prove the basic metric properties needed later, including completeness and monotone convergence criteria. The radial version of the metric on geodesic rays will be used in the study of test curves and singularity types.

This chapter provides the analytic bridge between envelopes and the later non-Archimedean theory. The space $\mathcal E^1(X,\theta;\phi)$ and the $d_1$-geometry developed here feed into the pseudometric $d_S$ of \cref{chap:comp}, while geodesic rays will later be related to test curves and non-Archimedean potentials.

\section{Subgeodesics}\label{sec:subgeo}
Let $X$ be a connected compact Kähler manifold of dimension $n$ and $\theta$ be a smooth closed real $(1,1)$-form on $X$ representing a big cohomology class.

\begin{definition}\label{def:subgeod}
    Let us fix $\varphi_0,\varphi_1\in \PSH(X,\theta)$. A \emph{subgeodesic}\index{subgeodesic} from $\varphi_0$ to $\varphi_1$ is a family $(\varphi_t)_{t\in (0,1)}$ in $\PSH(X,\theta)$ such that
    \begin{enumerate}
        \item\label{itm:def:subgeod:1} if we define
            \[
                \Phi\colon X\times \Bigl\{z\in \mathbb{C}:\Real z\in (0,1) \Bigr\}\rightarrow [-\infty,\infty),\quad (x,z)\mapsto \varphi_{\Real z}(x),
            \]
            then $\Phi$ is $p_1^*\theta$-psh, where $p_1\colon X\times \{z\in \mathbb{C}:\Real z\in (0,1)\}\rightarrow X$ is the natural projection;
        \item\label{itm:def:subgeod:2} when $t\to 0+$ (resp. $t\to 1-$), $\varphi_t$ converges to $\varphi_0$ (resp. $\varphi_1$) with respect to the $L^1$-topology.
    \end{enumerate}
    We also say $(\varphi_t)_{t\in [0,1]}$ is a subgeodesic.

    We call $\Phi$ the \emph{complexification}\index{complexification} of the subgeodesic $(\varphi_t)_t$.
\end{definition}
When we do not want to specify $\varphi_0$ and $\varphi_1$, we shall say $(\varphi_t)_{t\in (0,1)}$ is a subgeodesic\footnote{Note that it is implicit in the definition that $\varphi_t$ as $t\to 0+$ or $t\to 1-$ admits an $L^1$-limit.}.

More generally, a family $(\psi_t)_{t\in [a,b]}$ (resp. $(\psi_t)_{t\in (a,b)}$) in $\PSH(X,\theta)$ for some $a\leq b$ is called a \emph{subgeodesic}\index{subgeodesic} if $(\psi_{tb+(1-t)a})_{t\in [0,1]}$ (resp. $(\psi_{tb+(1-t)a})_{t\in (0,1)}$) is a subgeodesic.

\begin{remark}
    If $a<b$, then a subgeodesic $(\psi_t)_{t\in (a,b)}$ uniquely determines its endpoints $\psi_a$ and $\psi_b$: they are the $L^1$-limits at $a$ and $b$, and \cref{prop:pshfuncdetdense} identifies their plurisubharmonic representatives.
\end{remark}

\begin{remark}\label{rmk:complexify1}
In the literature, people sometimes regard $\Phi$ as a function defined on $X\times\{z\in \mathbb{C}: \mathrm{e}^{-1}<|z|<1\}$, with $\Phi(x,z)=\varphi_{-\log |z|}(x)$. We sometimes also use this convention; the ambient strip or annulus will make the convention clear.
\end{remark}

In general, there are no subgeodesics from $\varphi_0$ to $\varphi_1$. In fact, by \cref{prop:convexsubgeod} below, the existence of a subgeodesic from $\varphi_0$ to $\varphi_1$ forces $\varphi_0\land\varphi_1\not\equiv-\infty$, a condition that may fail (see \cref{ex:landminfty}).

We first note that the subgeodesics are well-behaved under the change of $\theta$:
\begin{proposition}\label{prop:subgeodindeptheta}
    Let $g$ be a smooth real function on $X$.
    Let $\theta'=\theta+\ddc g$. Suppose that $(\varphi_t)_{t\in [0,1]}$ is a subgeodesic in $\PSH(X,\theta)$. Then $(\varphi_t-g)_{t\in [0,1]}$ is a subgeodesic in $\PSH(X,\theta')$.
\end{proposition}
\begin{proof}
    This follows directly from the definition of the complexification.
\end{proof}

\begin{example}\label{ex:linearsubgeod}
    Let $\varphi_0\in \PSH(X,\theta)$, $C\in \mathbb{R}$. Let
    \[
    \varphi_t=\varphi_0+tC,\quad t\in [0,1].
    \]
    Then $(\varphi_t)_{t\in [0,1]}$ is a subgeodesic.
    Indeed, since $\Real z$ is harmonic in $z$, the complexification $\Phi(x,z)=\varphi_0(x)+C\,\Real z$ is $p_1^*\theta$-psh.

    As a consequence, the constant $(\varphi_0)_{t\in [0,1]}$ is a subgeodesic, called the \emph{constant subgeodesic}\index{subgeodesic!constant subgeodesic} at $\varphi_0$.

    A more general version is as follows: Suppose that $(\varphi_t)_{t\in [0,1]}$ is a subgeodesic in $\PSH(X,\theta)$, $C_1,C_2\in \mathbb{R}$. Then $(\varphi_t+C_1t+C_2)_{t\in [0,1]}$ is also a subgeodesic.
\end{example}

\begin{proposition}\label{prop:convexsubgeod}
    Let $\varphi_0,\varphi_1\in \PSH(X,\theta)$ and $(\varphi_t)_{t\in (0,1)}$ be a subgeodesic from $\varphi_0$ to $\varphi_1$. Then for each $x\in X$, $[0,1] \ni t\mapsto \varphi_t(x)$ is a convex function. Furthermore,
    \[
    \inf_{t\in (0,1)}\varphi_t\in \PSH(X,\theta),\quad \inf_{t\in (0,1)}\varphi_t\leq \varphi_0\land \varphi_1.
    \]
\end{proposition}

\begin{proof}
    Let $\Phi$ be the complexification of $(\varphi_t)_{t\in (0,1)}$.

    For each $x\in X$, the map
    \[
    \Bigl\{z\in \mathbb{C}: \Real z\in (0,1) \Bigr\}\rightarrow [-\infty,\infty),\quad z\mapsto \Phi(x,z)
    \]
    is either subharmonic or constantly $-\infty$, as follows from \itemref{itm:def:subgeod:1} and \cref{prop:pullbackpsh}.
    In the latter case, the convexity of $[0,1] \ni t\mapsto \varphi_t(x)$ is trivial. In the former case, the convexity on the interval $(0,1)$ follows from
    \cref{prop:subhimplyconv}.

    In order to verify the convexity at the boundary, let us fix $s\in (0,1)$. We need to show that
    \begin{equation}\label{eq:varphisconvextemp1}
    \varphi_s(x)\leq s\varphi_1(x)+(1-s)\varphi_0(x)
    \end{equation}
    for all $x\in X$. Thanks to \cref{prop:pshfuncdetdense}, it suffices to prove this for almost all $x$.

    Choose sequences $a_j\searrow0$ and $b_j\nearrow1$ such that $\varphi_{a_j}\to\varphi_0$ and $\varphi_{b_j}\to\varphi_1$ almost everywhere. Outside a set $Z\subseteq X$ of zero Lebesgue measure, all these convergences hold and all the values involved are finite. For $x\in X\setminus Z$ and $j$ large enough that $a_j<s<b_j$, the convexity on $(0,1)$ gives
    \[
    \varphi_s(x)\leq \frac{b_j-s}{b_j-a_j}\varphi_{a_j}(x)+\frac{s-a_j}{b_j-a_j}\varphi_{b_j}(x).
    \]
    Letting $j\to\infty$ gives \eqref{eq:varphisconvextemp1} almost everywhere.

    Let us prove the last assertion. Let
    \[
    \varphi\coloneqq \inf_{t\in (0,1)}\varphi_t.
    \]
    By Kiselman's principle \cref{prop:Kis}, we know that $\varphi\in \PSH(X,\theta)\cup \{-\infty\}$. Since the complexification $\Phi$ is not identically $-\infty$, Fubini's theorem shows that for almost every $x\in X$, the function $t\mapsto \varphi_t(x)$ is finite for almost every $t\in (0,1)$. For such an $x$, the convexity just proved implies that $t\mapsto \varphi_t(x)$ is a finite convex function on $(0,1)$, hence it is bounded from below on $(0,1)$. Thus $\varphi\not\equiv -\infty$, and so $\varphi\in \PSH(X,\theta)$.

    For any $t\in (0,1)$, using the convexity established above, we have
    \[
    \varphi\leq (1-t)\varphi_0+t\varphi_1.
    \]
    It follows that $\varphi\leq \varphi_0$ and $\varphi\leq \varphi_1$ almost everywhere and hence everywhere by \cref{prop:pshfuncdetdense}. Our assertion follows.
\end{proof}

\begin{corollary}\label{cor:subgeodboundary}
    Let $\varphi_0,\varphi_1\in \PSH(X,\theta)$ and $(\varphi_t)_{t\in (0,1)}$ be a subgeodesic from $\varphi_0$ to $\varphi_1$. Then
    \begin{equation}\label{eq:varphitconvexbdry}
    \lim_{t\to 0+}\varphi_t\leq \varphi_0,\quad   \lim_{t\to 1-}\varphi_t\leq \varphi_1.
    \end{equation}
    Moreover,
    \[
    \lim_{t\to 0+}\varphi_t=\varphi_0,\quad   \lim_{t\to 1-}\varphi_t=\varphi_1\quad \textup{quasi-everywhere}.
    \]
\end{corollary}
The inequalities in \eqref{eq:varphitconvexbdry} are not equalities everywhere in general. For example, consider an arbitrary subgeodesic from an unbounded $\varphi_0\in \mathcal{E}(X,\omega)$ to a bounded $\varphi_1\in \PSH(X,\omega)$ for some Kähler form $\omega$ on $X$.
\begin{proof}
    The inequalities in \eqref{eq:varphitconvexbdry} and the existence of the limits follow from the convexity in \cref{prop:convexsubgeod}. Let us prove the quasi-everywhere equality at $t=0$. Take a decreasing sequence $t_j\searrow 0$ with $t_j>0$. Then $\varphi_{t_j}\to \varphi_0$ in $L^1(X)$ as $j\to\infty$. So by \cref{cor:L1limipp}, we have
    \[
    \varlimsup_{j\to\infty}\varphi_{t_j}= \varphi_0\quad \textup{quasi-everywhere}.
    \]
    Therefore, we have
    \[
    \varphi_0=\varlimsup_{j\to\infty}\varphi_{t_j}=\lim_{t\to 0+}\varphi_t\leq \varphi_0\quad \textup{quasi-everywhere}.
    \]
    The assertion at $t=1$ is similar.
\end{proof}

\begin{corollary}\label{cor:convexsubgeod}
    Let $(\varphi_t)_{t\in [0,1]}$ be a subgeodesic in $\PSH(X,\theta)$. Then
    \begin{equation}\label{eq:supconvexsubgeod}
        \adjustlimits\sup_X \inf_{s\in (0,1)}\varphi_s \leq \sup_X \varphi_t\leq \sup_X \varphi_0\lor \sup_X \varphi_1
    \end{equation}
    for all $t\in (0,1)$.

    Furthermore, $[0,1] \ni t\mapsto \varphi_t$ is continuous with respect to the $L^1$-topology.
\end{corollary}
In particular, in view of \cref{prop:L1compa}, the family $\{\varphi_t\}_{t\in [0,1]}$ is relatively compact with respect to the $L^1$-topology on $\PSH(X,\theta)$.
\begin{proof}
    \eqref{eq:supconvexsubgeod} follows immediately from \cref{prop:convexsubgeod}.
    Let us verify the continuity of $t\mapsto \varphi_t$. The continuity at $t=0$ or $t=1$ follows from the definition of subgeodesic. While the continuity on $(0,1)$ follows from the convexity of $t\mapsto \varphi_t$, and the fact that $\varphi_t\in L^1(X)$ for each $t\in (0,1)$.
\end{proof}

\begin{proposition}\label{prop:maxsubgeod}
    Let $(\varphi^i_0)_{i\in I}$, $(\varphi^i_1)_{i\in I}$ be two non-empty uniformly bounded from above families in $\PSH(X,\theta)$. Let $(\varphi^i_t)_{t\in (0,1)}$ be subgeodesics from $\varphi_0^i$ to $\varphi_1^i$ for each $i\in I$. Then
    \[
    \left(\uscsup[i\in I]\varphi^i_t\right)_{t\in (0,1)}
    \]
    is a subgeodesic from $\uscsup[i] \varphi^i_0$ to $\uscsup[i] \varphi^i_1$.
\end{proposition}
\begin{proof}
    We may assume that $\varphi_0^i,\varphi_1^i\leq 0$ for all $i\in I$. Then it follows that $\varphi_t^i\leq 0$ for all $t\in (0,1)$ and all $i\in I$ by \cref{prop:convexsubgeod}.

    We define
    \[
    \varphi_t\coloneqq \uscsup[i\in I] \varphi^i_t\in \PSH(X,\theta)
    \]
    for all $t\in [0,1]$.
    Observe that $[0,1] \ni t\mapsto \varphi_t$ is convex by the same argument leading to \eqref{eq:varphisconvextemp1}.

    Let $(\psi_t)_{t\in (0,1)}$ be the subgeodesic whose complexification $\Phi_{\psi}$ corresponds to $\uscsup[i] \Phi_{\varphi^i}$, where $\Phi_{\varphi^i}$ is the complexification of $(\varphi^i_t)_{t\in (0,1)}$. Then clearly, $\varphi_t\leq \psi_t$ for each $t\in (0,1)$. On the other hand, by \cref{prop:supsupstardiff},
    \[
    \psi_t=\sup_{i\in I}\varphi^i_t=\varphi_t\quad \textup{almost everywhere}
    \]
    for almost all $t\in (0,1)$. Therefore, using \cref{prop:pshfuncdetdense}, we find $\psi_t=\varphi_t$ for almost all $t\in (0,1)$. Since both functions are convex in $t$, we conclude that $\psi_t=\varphi_t$ for all $t\in (0,1)$.

    It remains to argue that $\varphi_t\xrightarrow{L^1}\varphi_0$ as $t\to 0+$ and $\varphi_t\xrightarrow{L^1}\varphi_1$ as $t\to 1-$.
    By symmetry, it suffices to argue the former.

    By Choquet's lemma \cref{prop:Choquet}, applied to the uniformly bounded above family of complexifications $(\Phi_{\varphi^i})_{i\in I}$ on the strip, there is a countable subset $I'\subset I$ such that
    \[
    \uscsup[i\in I]\Phi_{\varphi^i}=\uscsup[i\in I']\Phi_{\varphi^i}.
    \]
    Applying Choquet's lemma \cref{prop:Choquet} once more to the endpoint families if necessary, and enlarging $I'$, we may also assume
    \[
    \uscsup[i\in I]\varphi_0^i=\uscsup[i\in I']\varphi_0^i,\quad \uscsup[i\in I]\varphi_1^i=\uscsup[i\in I']\varphi_1^i.
    \]
    Thus, after replacing $I$ by $I'$, we may assume that $I$ is countable.

    We know that for any $t\in (0,1)$ and any $j\in I$,
    \[
    \varphi^j_t\leq \varphi_t\leq t\varphi_1+(1-t)\varphi_0.
    \]
    Letting $t\to 0+$, we find that
    \[
    \varphi_0^j\leq \varlimsup_{t\to 0+}\varphi_t\leq \varphi_0
    \]
    almost everywhere. Since $I$ is countable, we conclude that
    \begin{equation}\label{eq:varphi0limsuptemp1}
    \varphi_0=\varlimsup_{t\to 0+}\varphi_t
    \end{equation}
    almost everywhere.

    Fix $i_0\in I$.
    Recall that by \cref{prop:convexsubgeod}, for each $t\in (0,1)$, we have
    \[
    \adjustlimits \inf_{t\in (0,1)}\sup_X \varphi_t\geq \adjustlimits\inf_{t\in (0,1)}\sup_X \varphi^{i_0}_t\geq \adjustlimits\sup_X \inf_{t\in (0,1)}\varphi^{i_0}_t>-\infty,
    \]
    so the set $\{\varphi_t\}_{t\in (0,1)}$ is relatively compact with respect to the $L^1$-topology by \cref{prop:L1compa}.
    Let $\psi$ be a cluster point as $t\to 0+$. It suffices to show that $\psi=\varphi_0$. By \cref{cor:L1limipp} and \eqref{eq:varphi0limsuptemp1}, this holds almost everywhere. Therefore, it holds everywhere by \cref{prop:pshfuncdetdense}.
\end{proof}

\begin{proposition}\label{prop:subgeodrestsubgeod}
    Let $(\varphi_t)_{t\in [0,1]}$ be a subgeodesic in $\PSH(X,\theta)$. Then for any $0\leq a\leq b\leq 1$, the segment $(\varphi_{t})_{t\in [a,b]}$ is a subgeodesic.
\end{proposition}

\begin{proof}
    It suffices to show that
    \[
    \varphi_{tb+(1-t)a}\xrightarrow{L^1}\varphi_a,\quad \varphi_{tb+(1-t)a}\xrightarrow{L^1}\varphi_b
    \]
    as $t\to 0+$ and $t\to 1-$ respectively. In other words, we need to show that for any $c\in (0,1)$, we have
    \[
    \varphi_t\xrightarrow{L^1}\varphi_c
    \]
    as $t\to c$. For this purpose, observe that by \cref{prop:convexsubgeod},
    \[
    \adjustlimits\sup_X \inf_{s\in (0,1)}\varphi_s\leq \sup_{X}\varphi_t\leq \left( \sup_X \varphi_0 \right) \lor \left( \sup_X\varphi_1 \right)
    \]
    for any $t\in (0,1)$. Therefore, $\{\varphi_t\}_{t\in (0,1)}$ is a relatively compact family with respect to the $L^1$-topology on $\PSH(X,\theta)$ by \cref{prop:L1compa}. It suffices to show that any cluster point $\psi$ of $\varphi_t$ as $t\to c$ is equal to $\varphi_c$.
    By \cref{cor:L1limipp} and the convexity \cref{prop:convexsubgeod}, we have $\varphi_c=\psi$ almost everywhere and hence everywhere by \cref{prop:pshfuncdetdense}.
\end{proof}

For the purpose of this chapter, it is useful to reformulate the gluing lemma \cref{cor:glbdry} in the following way, which allows us to glue subgeodesics directly:
\begin{lemma}\label{lma:subgeodgluing}
    Let $(\varphi_t)_{t\in [0,1]}$ be a subgeodesic in $\PSH(X,\theta)$. Let $0\leq a<b\leq 1$ and let $(\psi_t)_{t\in [a,b]}$ be a subgeodesic such that $\psi_a=\varphi_a$, $\psi_b=\varphi_b$, and
    \[
        \psi_t\geq \varphi_t,\quad t\in (a,b).
    \]
    Define
    \[
    \gamma_t=
    \begin{cases}
        \psi_t,&\textup{if }t\in (a,b),\\
        \varphi_t,&\textup{if }t\in [0,1]\setminus (a,b).
    \end{cases}
    \]
    Then $(\gamma_t)_{t\in [0,1]}$ is a subgeodesic.
\end{lemma}
\begin{proof}
    The $L^1$-convergence of $\gamma_t$ to $\varphi_0$ at $0$ and $\varphi_1$ at $1$ follows immediately from the corresponding property for $(\varphi_t)_t$ and $(\psi_t)_t$.
    It remains to prove that the complexification of $(\gamma_t)_t$ is $p_1^*\theta$-psh.

    Let
    \[
        S=\{z\in\mathbb{C}:0<\Real z<1\},\quad S'=\{z\in\mathbb{C}:a<\Real z<b\}.
    \]
    Let $\Phi$ and $\Psi$ be the complexifications of $(\varphi_t)_t$ and $(\psi_t)_t$ on the corresponding strips. Since $\Psi\geq \Phi$ on $X\times S'$, the complexification of $(\gamma_t)_t$ is equal to $\Phi\lor\Psi$ on $X\times S'$ and to $\Phi$ on $X\times (S\setminus S')$.
    By the gluing lemma \cref{cor:glbdry}, it is enough to verify that
    \[
        \varlimsup_{\substack{S'\ni z'\to z\\ x'\to x}}\Psi(x',z')\leq \Phi(x,z)
    \]
    at every point $(x,z)\in X\times S$ with $\Real z=a$ or $\Real z=b$.

    Suppose for instance that $a>0$, and let $(x_j,z_j)$ tend to $(x,z)$ with $\Real z=a$ and $a<\Real z_j<b$. Set $t_j=\Real z_j$. Since $(\psi_t)_t$ is a subgeodesic from $\varphi_a$ to $\varphi_b$, the sequence $\psi_{t_j}$ converges to $\varphi_a$ in $L^1$. Applying \cref{cor:L1limipp} to this sequence gives
    \[
        \varphi_a=\inf_{m>0}\left(\uscsup[j\geq m]\psi_{t_j}\right).
    \]
    For fixed $m$, the upper semicontinuity of $\uscsup[j\geq m]\psi_{t_j}$ gives
    \[
        \varlimsup_{j\to\infty}\psi_{t_j}(x_j)\leq \left(\uscsup[j\geq m]\psi_{t_j}\right)(x).
    \]
    Letting $m\to\infty$, we get
    \[
        \varlimsup_{j\to\infty}\psi_{t_j}(x_j)\leq \varphi_a(x).
    \]
    Since $\Phi(x,z)=\varphi_a(x)$ when $\Real z=a$, this is the desired boundary inequality along $\{\Real z=a\}$. The same argument applies along $\{\Real z=b\}$ when $b<1$. Hence \cref{cor:glbdry} shows that the pasted complexification is $p_1^*\theta$-psh.
\end{proof}

\begin{definition}\label{def:subgeoray}
    A ray $\ell=(\ell_t)_{t\geq 0}$ is a \emph{subgeodesic ray}\index{subgeodesic ray} in $\PSH(X,\theta)$ if for any $0\leq a\leq b$, the segment $(\ell_{t})_{t\in [a,b]}$ is a subgeodesic in $\PSH(X,\theta)$. We say $\ell$ \emph{emanates} from $\ell_0$.

    The \emph{complexification}\index{complexification} of a subgeodesic ray $\ell$ is defined as the potential
    \[
    \Phi\colon X\times \{z\in \mathbb{C}:\Real z>0\}\rightarrow [-\infty,\infty),\quad (x,z)\mapsto \ell_{\Real z}(x).
    \]
    Note that $\Phi$ is $p_1^*\theta$-psh, where $p_1\colon X\times \{z\in \mathbb{C}:\Real z>0\}\rightarrow X$ is the natural projection.
\end{definition}
\begin{remark}
    As in \cref{rmk:complexify1}, one may also define the complexification as a function $X\times \{z\in \mathbb{C}: 0<|z|<1\}\rightarrow [-\infty,\infty)$.
\end{remark}

\section{Geodesics in the space of potentials}\label{sec:relativeray}
Let $X$ be a connected compact Kähler manifold of dimension $n$ and $\theta$ be a smooth closed real $(1,1)$-form on $X$ representing a big cohomology class.
Fix a model potential $\phi\in \PSH(X,\theta)_{>0}$. See \cref{def:modelpot} for the definition.

\begin{definition}\label{def:geo}
    Let $\varphi_0, \varphi_1 \in \PSH(X,\theta)$. The \emph{geodesic}\index{geodesic} $(\varphi_t)_{t\in (0,1)}$ from $\varphi_0$ to $\varphi_1$ is the family of potentials $\varphi_t\in \PSH(X,\theta)\cup \{-\infty\}$ such that
    \begin{equation}\label{eq:Perron2}
        \begin{split}
        \varphi_t=\uscsup & \left\{  \psi_t: (\psi_{s})_s \textup{ is a subgeodesic from }\psi_0\textup{ to }\psi_1,\right.\\
        &\psi_0,\psi_1\in \PSH(X,\theta),\left.\psi_0\leq \varphi_0,\psi_1\leq \varphi_1\right\}.
        \end{split}
    \end{equation}
    More generally, let $(\varphi_t)_{t\in [a,b]}$ ($a,b\in \mathbb{R}$, $a\leq b$) be a curve in $\PSH(X,\theta)$. We say $(\varphi_t)_{t\in [a,b]}$ is a \emph{geodesic}\index{geodesic} if the curve $(\varphi_{tb+(1-t)a})_{t\in (0,1)}$ is a geodesic from $\varphi_a$ to $\varphi_b$.

    We also say $(\varphi_t)_{t\in (a,b)}$ is a \emph{geodesic} in $\PSH(X,\theta)$ from $\varphi_a$ to $\varphi_b$.
\end{definition}

The envelopes of the form \eqref{eq:Perron2} are usually referred to as the \emph{Perron envelopes}. In general, the geodesic defined by \eqref{eq:Perron2} fails to have the correct limit when $t\to 0+$ or $t\to 1-$. Therefore, \emph{a priori} it is not a subgeodesic from $\varphi_0$ to $\varphi_1$. In the sequel, we use the word \emph{geodesic} only in situations where the endpoint behavior has been verified.

\begin{lemma}\label{lma:geodsubgeod}
    Let $\varphi_0,\varphi_1\in \PSH(X,\theta)$. Assume that there is a subgeodesic from $\varphi_0$ to $\varphi_1$. Then the geodesic $(\varphi_t)_{t\in (0,1)}$ from $\varphi_0$ to $\varphi_1$ is a subgeodesic from $\varphi_0$ to $\varphi_1$.
\end{lemma}
\begin{proof}
    This is an immediate consequence of \cref{prop:maxsubgeod}.
\end{proof}

\begin{example}\label{ex:lineargeod}
    Let $\varphi_0\in \PSH(X,\theta)$ and $C\in \mathbb{R}$. Then the subgeodesic $(\varphi_0+tC)_{t\in [0,1]}$ studied in \cref{ex:linearsubgeod} is a geodesic. This follows easily from \cref{prop:convexsubgeod}.

    In particular, when $C=0$, we find that the constant subgeodesic at $\varphi_0$ is indeed a geodesic, which we call the \emph{constant geodesic}\index{geodesic!constant geodesic} at $\varphi_0$.

    More generally, suppose that $(\varphi_t)_{t\in [0,1]}$ is a geodesic and $C_1,C_2\in \mathbb{R}$, then $(\varphi_t+C_1t+C_2)_{t\in [0,1]}$ is also a geodesic. This follows immediately from \cref{ex:linearsubgeod}.
\end{example}

Next we want to show that under mild assumptions, there exist subgeodesics between two potentials. The assumption below turns out to be necessary as well, as we shall prove in \cref{thm:geodexi} later.
\begin{proposition}\label{prop:perronenvissubgeod}
Given $\varphi_0, \varphi_1 \in \mathcal{E}(X,\theta;\phi)$, the geodesic $(\varphi_t)_{t\in (0,1)}$ from $\varphi_0$ to $\varphi_1$ defined by \eqref{eq:Perron2} is a subgeodesic from $\varphi_0$ to $\varphi_1$ and $\varphi_t\in \mathcal{E}(X,\theta;\phi)$ for each $t\in (0,1)$.

Moreover, for any $0\leq a\leq b\leq 1$, the restriction $(\varphi_t)_{t\in [a,b]}$ is a geodesic.

If furthermore $\varphi_0,\varphi_1\in \mathcal{E}^1(X,\theta;\phi)$ (resp. $\mathcal{E}^{\infty}(X,\theta;\phi)$), then $\varphi_t\in \mathcal{E}^1(X,\theta;\phi)$ (resp. $\mathcal{E}^{\infty}(X,\theta;\phi)$) for all $t\in (0,1)$.
\end{proposition}
We refer to \cref{subsec:fullmass} for the definition of $\mathcal{E}(X,\theta;\phi)$. Our assumption means that $P_{\theta}[\varphi_0]=P_{\theta}[\varphi_1]=\phi$.

As a consequence of this proposition, the restriction $(\varphi_t)_{t\in [a,b]}$ is also a subgeodesic since $\varphi_a,\varphi_b\in \mathcal{E}(X,\theta;\phi)$ as well.

\begin{proof}
    After subtracting suitable constants from $\varphi_0$ and $\varphi_1$ and adding back the corresponding affine function of $t$, we may assume that $\varphi_0,\varphi_1\leq \phi$. It follows from \cref{prop:convexsubgeod} that $\varphi_t\leq \phi$ for all $t\in (0,1)$. In fact, we have the stronger estimate
    \begin{equation}\label{eq:geodesicconvextemp1}
        \varphi_t\leq t\varphi_1+(1-t)\varphi_0,\quad t\in (0,1).
    \end{equation}

    We first observe that when $\varphi_0,\varphi_1\in \mathcal{E}(X,\theta;\phi)$, so is $\varphi_0\land \varphi_1$, see \cref{prop:relrooftopclosed}. In particular, the constant subgeodesic $t\mapsto \varphi_0\land \varphi_1$ is a candidate in the Perron envelope \eqref{eq:Perron2}. So
    \begin{equation}\label{eq:varphitgeqlandtemp1}
        \varphi_t\geq \varphi_0\land \varphi_1,\quad t\in (0,1).
    \end{equation}
    By \cref{prop:maxsubgeod}, $(\varphi_t)_{t\in (0,1)}$ is a subgeodesic.\footnote{Be careful, here $t\in (0,1)$ instead of $[0,1]$.}
    It follows from \cref{prop:relativeEupperclosed} that $\varphi_t\in \mathcal{E}(X,\theta;\phi)$ for all $t\in (0,1)$.

    Next, we show that as $t\to 0+$, we have $\varphi_t\xrightarrow{L^1}\varphi_0$. The corresponding result at $t=1$ is similar.

    We first argue the special case where $\varphi_0\preceq \varphi_1$. Take a constant $C>0$ such that
    \[
        \varphi_0-C\leq \varphi_1.
    \]
    Then $(\varphi_0-Ct)_{t\in (0,1)}$ is clearly a candidate in \eqref{eq:Perron2}, see \cref{ex:linearsubgeod}. Therefore, for all $t\in (0,1)$,
    \begin{equation}\label{eq:varphi0andvarphit}
        \varphi_0-Ct\leq \varphi_t\leq t\varphi_1+(1-t)\varphi_0.
    \end{equation}
    It follows that $\varphi_t\xrightarrow{L^1}\varphi_0$ as $t\to 0+$.

    Let us come back to the general case.
    By \cref{cor:convexsubgeod}, $\{\varphi_t:t\in (0,1)\}$ is a relatively compact subset of $\PSH(X,\theta)$ with respect to the $L^1$-topology.

    Let $\psi$ be an $L^1$-cluster point of $\varphi_t$ as $t\searrow 0$, and it suffices to show that $\psi=\varphi_0$.

    For each $M\in \mathbb{N}$, we write
    \[
        \varphi_0^M=\varphi_0\land (\varphi_1+M).
    \]
    Observe that $\varphi_0^M\in \mathcal{E}(X,\theta;\phi)$ by \cref{prop:relrooftopclosed}.
    Let $(\varphi_t^M)_{t\in (0,1)}$ be the geodesic from $\varphi_0^M$ to $\varphi_1$. Since $\varphi_0^M\preceq\varphi_1$, the special case proved above applies to the geodesic $(\varphi_t^M)_t$, and
    \[
        \varphi_t^M\xrightarrow{L^1}\varphi_0^M \quad\textup{as }t\to0+.
    \]
    Moreover, every candidate for the geodesic from $\varphi_0^M$ to $\varphi_1$ is also a candidate for the geodesic from $\varphi_0$ to $\varphi_1$, so
    \[
        \varphi_t^M\leq\varphi_t,\quad t\in(0,1).
    \]
    Letting $t\to 0+$ along a sequence for which $\varphi_t\to\psi$ in $L^1$, and using \cref{cor:L1limipp}, we obtain
    \[
        \psi\geq\varphi_0^M
    \]
    almost everywhere, hence everywhere by \cref{prop:pshfuncdetdense}.

    On the other hand, by \eqref{eq:geodesicconvextemp1}, $\psi\leq \varphi_0$. So it suffices to show that
    \[
        \varphi_0\land (\varphi_1+M)\xrightarrow{L^1}   \varphi_0
    \]
    as $M\to \infty$, which is shown in \cref{prop:envrelfullmass}.

    Now we have shown that $(\varphi_t)_{t\in [0,1]}$ is a subgeodesic.

    Next, take $0\leq a\leq b\leq 1$. We want to show that the restriction $(\varphi_t)_{t\in [a,b]}$ is the geodesic from $\varphi_a$ to $\varphi_b$. We may assume that $a<b$. The argument is the standard \emph{balayage} argument.

    Let $(\psi_t)_{t\in (a,b)}$ be the geodesic from $\varphi_a$ to $\varphi_b$. Since $(\varphi_t)_{t\in [a,b]}$ is a subgeodesic by \cref{prop:subgeodrestsubgeod}, we have $\psi_t\geq \varphi_t$ for all $t\in (a,b)$.

    We define
    \[
    \eta_t=
    \left\{
    \begin{aligned}
        \psi_t,&\quad\text{if }t\in (a,b),\\
        \varphi_t,&\quad\text{if }t\in (0,1)\setminus (a,b).
    \end{aligned}
    \right.
    \]
    Extending it by setting $\psi_a=\varphi_a$ and $\psi_b=\varphi_b$, the first part of the proof applied to the endpoints $\varphi_a,\varphi_b$ shows that $(\psi_t)_{t\in [a,b]}$ is a subgeodesic from $\varphi_a$ to $\varphi_b$. Since $\psi_t\geq \varphi_t$ for $t\in (a,b)$, the gluing lemma \cref{lma:subgeodgluing} shows that $(\eta_t)_{t\in [0,1]}$ is a subgeodesic from $\varphi_0$ to $\varphi_1$.

    Therefore, for all $t\in (0,1)$, we have
    \[
    \varphi_t\geq \eta_t.
    \]
    In particular, for $t\in (a,b)$, we have
    \[
    \varphi_t\geq \eta_t=\psi_t\geq \varphi_t.
    \]
    In other words, $(\varphi_t)_{t\in (a,b)}=(\psi_t)_{t\in (a,b)}$ is the geodesic from $\varphi_a$ to $\varphi_b$.

    Finally, assume furthermore that $\varphi_0,\varphi_1\in \mathcal{E}^1(X,\theta;\phi)$ (resp. $\mathcal{E}^{\infty}(X,\theta;\phi)$). Thanks to \eqref{eq:varphitgeqlandtemp1}, \cref{prop:relrooftopclosed} and \cref{prop:relativeEupperclosed}, we find $\varphi_t\in \mathcal{E}^1(X,\theta;\phi)$ (resp. $\mathcal{E}^{\infty}(X,\theta;\phi)$) for all $t\in (0,1)$.
\end{proof}

\begin{lemma}\label{lma:geoddecreasing}
    Let $\varphi_0,\varphi_1\in \PSH(X,\theta)$ with $\varphi_1\leq \varphi_0$. Let $(\varphi_t)_{t\in (0,1)}$ be the geodesic from $\varphi_0$ to $\varphi_1$. Assume that $(\varphi_t)_{t\in [0,1]}$ is also a subgeodesic. Then $\varphi_t$ is decreasing in $t\in [0,1]$.
\end{lemma}
\begin{proof}
    Fix $0< r<s\leq 1$ and put $\delta=(s-r)/(1-r)$. Then $(\varphi_{\delta+(1-\delta)t})_{t\in [0,1]}$ is a subgeodesic from $\varphi_\delta$ to $\varphi_1$. Since
    \[
    \varphi_\delta\leq (1-\delta)\varphi_0+\delta\varphi_1\leq \varphi_0,
    \]
    we have $\varphi_{\delta+(1-\delta)t}\leq \varphi_t$ for all $t\in (0,1)$ by the definition of geodesic. In particular, $\varphi_s=\varphi_{\delta+(1-\delta)r}\leq \varphi_r$.

    It remains to prove that for any $t\in (0,1)$, we have $\varphi_t\leq \varphi_0$, but this is an immediate consequence of \cref{prop:convexsubgeod}:
    \[
        \varphi_t\leq t\varphi_1+(1-t)\varphi_0\leq \varphi_0.
    \]
\end{proof}

\begin{proposition}\label{prop:geodsupsublinear}
    Let $\varphi_1,\varphi_0\in \mathcal{E}(X,\theta;\phi)$ with $\varphi_1\preceq \varphi_0$. Let $(\varphi_t)_{t\in (0,1)}$ be the geodesic from $\varphi_0$ to $\varphi_1$. Then
    \begin{equation}\label{eq:tsupsuptemp1}
        s\sup_{\{\varphi_0\neq -\infty\}}(\varphi_1-\varphi_0)=\sup_{\{\varphi_0\neq -\infty\}}(\varphi_s-\varphi_0)
    \end{equation}
    for all $s\in [0,1]$.
\end{proposition}
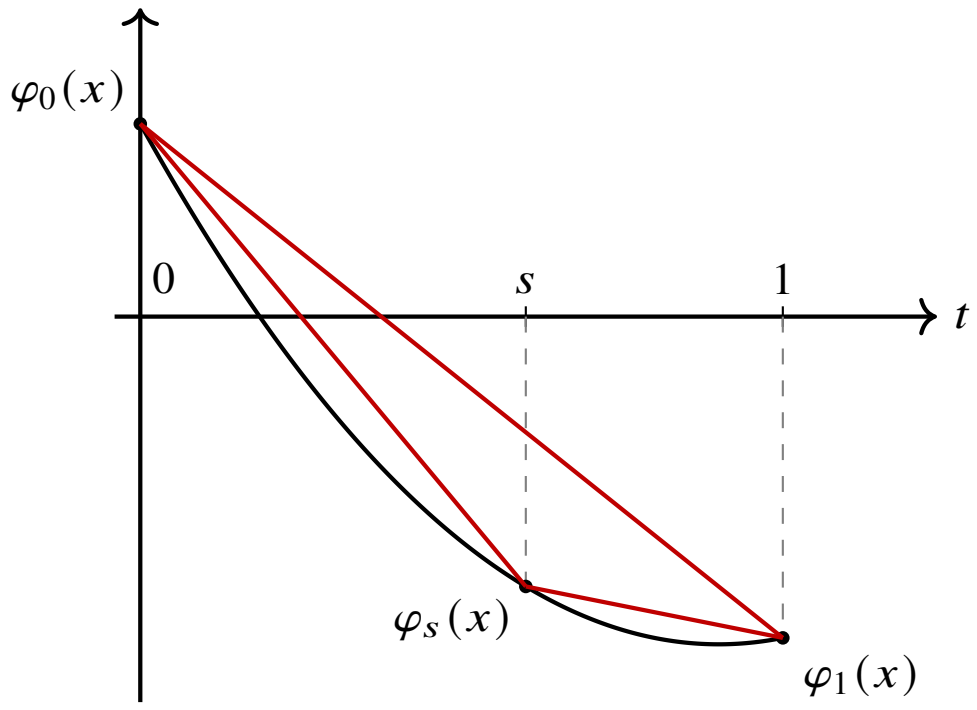
\begin{figure}[ht]
\centering
\begin{tikzpicture}[scale=0.85, line width=0.6pt, every node/.style={font=\small}]
  \draw[->, thick] (-0.2,0) -- (6.2,0) node[right] {$t$};
  \draw[->, thick] (0,-3.0) -- (0,2.4);
  \draw[thin] (0,0.08) -- (0,-0.08);
  \draw[thin] (3,0.08) -- (3,-0.08);
  \draw[thin] (5,0.08) -- (5,-0.08);
  \node[anchor=south west, inner sep=2pt] at (0,0.08) {$0$};
  \node[anchor=south, inner sep=2pt] at (3,0.08) {$s$};
  \node[anchor=south, inner sep=2pt] at (5,0.08) {$1$};
  \draw[gray, dashed, thin] (3, 0) -- (3,-2.10);
  \draw[gray, dashed, thin] (5, 0) -- (5,-2.50);
  \draw[thick] plot[smooth, samples=80, domain=0:5] (\x, {1.5 - 1.8*\x + 0.20*\x*\x});
  \fill (0,1.5) circle (1.5pt);
  \fill (3,-2.10) circle (1.5pt);
  \fill (5,-2.50) circle (1.5pt);
  \node[anchor=south east, inner sep=2pt] at (0,1.5) {$\varphi_0(x)$};
  \node[anchor=north east, inner sep=2pt] at (3,-2.10) {$\varphi_s(x)$};
  \node[anchor=north west, inner sep=3pt] at (5,-2.50) {$\varphi_1(x)$};
  \draw[red!75!black, thick] (0,1.5) -- (3,-2.10);
  \draw[red!75!black, thick] (0,1.5) -- (5,-2.50);
  \draw[red!75!black, thick] (3,-2.10) -- (5,-2.50);
\end{tikzpicture}
\caption{Typical behaviour of $t\mapsto\varphi_t(x)$}
\label{fig:convphi}
\end{figure}
\begin{proof}
The notations in the proof are indicated in \cref{fig:convphi}.\footnote{When dealing with convex functions, drawing a picture is the easiest way to keep track of the directions of inequalities.}

    We may assume that $s\in (0,1)$ since there is nothing to prove when $s=0$ or $s=1$.

    After replacing $\varphi_t$ by $\varphi_t-C't$ for some large enough $C'>0$, we may assume that $\varphi_1\leq \varphi_0$. This procedure preserves the geodesic property by \cref{ex:lineargeod}. Then $\varphi_t$ is decreasing in $t$ by \cref{lma:geoddecreasing}.

    Let
    \begin{equation}\label{eq:Cdefsupvarphi1m0temp1}
        C=\sup_{\{\varphi_1\neq -\infty\}} \left(\varphi_1-\varphi_0\right)=\sup_{\{\varphi_0\neq -\infty\}} \left(\varphi_1-\varphi_0\right)\leq 0,
    \end{equation}
    where the equality is due to \cref{cor:diffsupinfindeppluripolar}.
    Then by \cref{prop:pshfuncdetdense}, we have
    \[
        \varphi_1\leq \varphi_0+C.
    \]
    So $(\varphi_1-C(1-t))_{t\in (0,1)}$ is a candidate in \eqref{eq:Perron2} and hence
    \begin{equation}\label{eq:varphi1leqvarphittemp}
        \varphi_1-C(1-t)\leq \varphi_t,\quad t\in (0,1).
    \end{equation}

    By \cref{prop:perronenvissubgeod}, we have $\varphi_t\xrightarrow{L^1}\varphi_1$ as $t\to 1-$.
    Since $\varphi_t$ is decreasing in $t\in (0,1)$, it follows that $\varphi_1=\inf_{t\in (0,1)}\varphi_t$.
    Therefore, we can find a pluripolar set $Z\subseteq X$ such that $\varphi_t(x)\rightarrow \varphi_1(x)>-\infty$ as $t\to 1-$ for all $x\in X\setminus Z$.

    Similarly, since $\varphi_0=\uscsup[t\in (0,1)] \varphi_t$, after enlarging $Z$, we may also guarantee that $\varphi_t(x)\rightarrow \varphi_0(x)>-\infty$ as $t\to 0+$ for all $x\in X\setminus Z$ by \cref{prop:supsupstardiff}.

    For any such $x\in X\setminus Z$, the function $t\mapsto \varphi_t(x)$ is a real-valued continuous convex function on $[0,1]$. In particular, $t\mapsto \varphi_t(x)$ is absolutely continuous on $[0,1]$.
    Hence, for any $s\in [0,1)$, we have
    \begin{equation}\label{eq:varphi1minusvarphi0temp1}
        \varphi_1(x)-\varphi_s(x)=\int_s^1 \ddt \varphi_t(x)\,\mathrm{d}t\leq (1-s)\lim_{t\to 1-}\frac{\varphi_1(x)-\varphi_t(x)}{1-t}\leq (1-s)C,
    \end{equation}
    where the second inequality follows from \eqref{eq:varphi1leqvarphittemp}.

    Taking supremum in \eqref{eq:varphi1minusvarphi0temp1}, we find that
    \begin{equation}\label{eq:supvarphi1sleqleqtemp1}
    \sup_{X\setminus Z} (\varphi_1-\varphi_s)\leq (1-s)\adjustlimits\sup_{x\in X\setminus Z}\lim_{t\to 1-}\frac{\varphi_1(x)-\varphi_t(x)}{1-t}\leq (1-s)C.
    \end{equation}
    When $s=0$, we deduce from \cref{cor:diffsupinfindeppluripolar} and \eqref{eq:Cdefsupvarphi1m0temp1} that
    \[
    C=\sup_{\{\varphi_1\neq -\infty\}} \left(\varphi_1-\varphi_0\right)=\adjustlimits\sup_{x\in X\setminus Z}\lim_{t\to 1-}\frac{\varphi_1(x)-\varphi_t(x)}{1-t}.
    \]
    But this equality works equally well for the geodesic $(\varphi_{(1-s)t+s})_{t\in [0,1]}$ with the same argument. It follows that
    \[
    \sup_{\{\varphi_1\neq -\infty\}} \left(\varphi_1-\varphi_s\right)=(1-s)\adjustlimits\sup_{x\in X\setminus Z}\lim_{t\to 1-}\frac{\varphi_1(x)-\varphi_t(x)}{1-t}=(1-s)C.
    \]
    Therefore, invoking \cref{cor:diffsupinfindeppluripolar} again, we deduce that all inequalities in \eqref{eq:supvarphi1sleqleqtemp1} are in fact equalities. In other words,
    \begin{equation}\label{eq:supvarphi1mivarphi0temp1}
        \sup_{\{\varphi_1\neq -\infty\}}(\varphi_1-\varphi_0)=   \adjustlimits \sup_{x\in X\setminus Z}\lim_{t\to 1-}\frac{\varphi_1(x)-\varphi_t(x)}{1-t}
        =\sup_{\{\varphi_1\neq -\infty\}}\frac{\varphi_1-\varphi_s}{1-s}.
    \end{equation}
    On the other hand, we have the trivial inequality when $s\in (0,1)$:
    \[
    \sup_{\{\varphi_1\neq -\infty\}}(\varphi_1-\varphi_0)\leq s \sup_{\{\varphi_1\neq -\infty\}}\frac{\varphi_s-\varphi_0}{s}+   (1-s) \sup_{\{\varphi_1\neq -\infty\}}\frac{\varphi_1-\varphi_s}{1-s}.
    \]
    Together with \eqref{eq:supvarphi1mivarphi0temp1}, we find that for $s\in (0,1)$,
    \[
    \sup_{\{\varphi_1\neq -\infty\}}(\varphi_1-\varphi_0)\leq \sup_{\{\varphi_1\neq -\infty\}}\frac{\varphi_s-\varphi_0}{s}.
    \]
    The reverse inequality follows from the convexity, and hence for $s\in (0,1)$, we have
    \[
        \sup_{\{\varphi_1\neq -\infty\}}\frac{\varphi_s-\varphi_0}{s}=\sup_{\{\varphi_1\neq -\infty\}}(\varphi_1-\varphi_0).
    \]
    Using \cref{cor:diffsupinfindeppluripolar}, we conclude \eqref{eq:tsupsuptemp1}.
\end{proof}

With an almost identical proof, we find
\begin{proposition}\label{prop:geodinfsublinear}
    Let $\varphi_1,\varphi_0\in \mathcal{E}^{\infty}(X,\theta;\phi)$. Let $(\varphi_t)_{t\in (0,1)}$ be the geodesic from $\varphi_0$ to $\varphi_1$. Then
    \[
        t\inf_{\{\phi\neq -\infty\}}(\varphi_1-\varphi_0)=\inf_{\{\phi\neq -\infty\}}(\varphi_t-\varphi_0)
    \]
    for all $t\in (0,1]$.
\end{proposition}

\begin{definition}\label{def:geodraydef}
    Let $\ell=(\ell_t)_{t\geq 0}$ be a family in $\mathcal{E}(X,\theta;\phi)$. We say $\ell$ is a \emph{geodesic ray}\index{geodesic ray} in $\mathcal{E}(X,\theta;\phi)$ emanating from $\ell_0$ if for each $0\leq a\leq b$, the restriction $(\ell_t)_{t\in [a,b]}$ is a geodesic.

    The set of geodesic rays in $\mathcal{E}(X,\theta;\phi)$ emanating from $\phi$ is denoted by $\mathcal{R}(X,\theta;\phi)$\index{$\mathcal{R}(X,\theta;\phi)$}.

    We say a geodesic ray $\ell\in \mathcal{R}(X,\theta;\phi)$ has \emph{finite energy}\index{geodesic ray!geodesic ray with finite energy} if $\ell_t\in \mathcal{E}^1(X,\theta;\phi)$ for all $t>0$. The set of geodesic rays with finite energy is denoted by $\mathcal{R}^1(X,\theta;\phi)$\index{$\mathcal{R}^1(X,\theta;\phi)$}.

    We say a geodesic ray $\ell\in \mathcal{R}(X,\theta;\phi)$ is \emph{bounded}\index{geodesic ray!bounded geodesic ray} if $\ell_t\in \mathcal{E}^{\infty}(X,\theta;\phi)$ for all $t\geq 0$. The set of bounded geodesic rays is denoted by $\mathcal{R}^{\infty}(X,\theta;\phi)$\index{$\mathcal{R}^{\infty}(X,\theta;\phi)$}.

    Given $\ell,\ell'\in \mathcal{R}(X,\theta;\phi)$, we write $\ell\leq \ell'$ if $\ell_t\leq \ell'_t$ for each $t\geq 0$.

    When $\phi=V_{\theta}$, we usually omit it from the notations and write $\mathcal{R}(X,\theta)$\index{$\mathcal{R}(X,\theta)$}, $\mathcal{R}^1(X,\theta)$\index{$\mathcal{R}^1(X,\theta)$} and $\mathcal{R}^{\infty}(X,\theta)$\index{$\mathcal{R}^{\infty}(X,\theta)$} respectively.
\end{definition}

\begin{proposition}\label{prop:raysuplinear}
    Let $\ell\in \mathcal{R}(X,\theta;\phi)$. Then there is a constant $C\in \mathbb{R}$ such that
    \[
        \sup_{X}\ell_t=Ct,\quad t\geq 0.
    \]
\end{proposition}
\begin{proof}
    It follows from \cref{prop:geodsupsublinear} that
    \[
    \sup_{\{\phi\neq -\infty\}}(\ell_t-\phi)=t\sup_{\{\phi\neq -\infty\}}(\ell_1-\phi)
    \]
    for all $t\geq 0$.

    It suffices to show that for any $t\geq 0$,
    \[
    \sup_{\{\phi\neq -\infty\}}(\ell_t-\phi)=\sup_X \ell_t.
    \]
    This was already proved in \cref{prop:model_sup_prop}.
\end{proof}

\section{The metrics on the spaces of potentials and rays}\label{sec:metE1}

Let $X$ be a connected compact Kähler manifold of dimension $n$ and $\theta$ be a smooth closed real $(1,1)$-form on $X$ representing a big cohomology class.
Fix a model potential $\phi\in \PSH(X,\theta)_{>0}$.

We first study a natural metric on $\mathcal{E}^1(X,\theta;\phi)$.

\begin{definition}\label{def:d1onE12}
    Let $\varphi,\psi\in \mathcal{E}^1(X,\theta;\phi)$, we define the \emph{Darvas metric}\index{Darvas metric} $d_1$ as follows:
    \[
        d_1(\varphi,\psi)=E_{\theta}^{\phi}(\varphi)+E_{\theta}^{\phi}(\psi)-2E_{\theta}^{\phi}(\varphi\land \psi).
    \]
    Note that by \cref{prop:relrooftopclosed}, $\varphi\land \psi\in \mathcal{E}^1(X,\theta;\phi)$.
\end{definition}
Recall that $E_{\theta}^{\phi}$ is defined in \cref{def:MAenergy}.

In particular, if $\varphi\leq \psi$, we have
\begin{equation}\label{eq:d1asEdiff}
    d_1(\varphi,\psi)=E_{\theta}^{\phi}(\psi)-E_{\theta}^{\phi}(\varphi).
\end{equation}

We wish to show that $d_1$ is a complete metric. We first prove a contraction property:
\begin{proposition}\label{prop:d1contract}
    Let $\varphi,\psi,\gamma\in \mathcal{E}^1(X,\theta;\phi)$. Then
    \begin{equation}\label{eq:d1contract}
    d_1(\varphi,\psi)\geq d_1(\varphi\land \gamma,\psi\land \gamma).
    \end{equation}
\end{proposition}
\begin{proof}
    \textbf{Step~1}. We first assume that $\varphi\geq \psi$. Then
    \[
    \begin{aligned}
        d_1(\varphi,\psi)=& E_{\theta}^{\phi}(\varphi)-E_{\theta}^{\phi}(\psi)\\
        \geq & E_{\theta}^{\phi}\Bigl((\varphi\land \gamma)\lor \psi\Bigr)-E_{\theta}^{\phi}(\psi) & \textup{by \cref{lma:EmonoEinfty}}\\
        \geq & E_{\theta}^{\phi}(\varphi\land \gamma)-E_{\theta}^{\phi}(\psi\land \gamma)&\textup{by \cref{cor:diam_E}}.
    \end{aligned}
    \]

    \textbf{Step~2}. We prove the general case.

    By Step~1, we have
    \[
    d_1(\varphi,\varphi\land \psi)\geq d_1(\varphi\land \gamma,\varphi\land \psi\land \gamma),\quad d_1(\psi,\varphi\land \psi)\geq d_1(\psi\land \gamma,\varphi\land \psi\land \gamma).
    \]
    Adding the two inequalities together, we conclude \eqref{eq:d1contract}.
\end{proof}

\begin{lemma}\label{lma:d1metric}
    The function $d_1$ is a metric on $\mathcal{E}^1(X,\theta;\phi)$.
\end{lemma}
\begin{proof}
    There are two facts to prove, as in the two steps below.

    \textbf{Step~1}. Let $\varphi,\psi\in \mathcal{E}^1(X,\theta;\phi)$. Assume that $d_1(\varphi,\psi)=0$. Then we will show that $\varphi=\psi$.

    We first treat the case $\psi\leq \varphi$. Then $d_1(\varphi,\psi)=E_\theta^\phi(\varphi)-E_\theta^\phi(\psi)=0$, and \cref{prop:cocycE1} gives
    \[
    \int_X (\varphi-\psi)\,\theta_{\psi}^n=0.
    \]
    We conclude that $\psi=\varphi$ using the domination principle \cref{thm:dom_prin_1}.

    In general, put $\gamma=\varphi\land\psi$. Since
    \[
    d_1(\varphi,\psi)=d_1(\varphi,\gamma)+d_1(\psi,\gamma),
    \]
    both terms on the right vanish. Applying the comparable case to the pairs $(\varphi,\gamma)$ and $(\psi,\gamma)$ gives $\varphi=\gamma=\psi$.

    \textbf{Step~2}. Let $\varphi,\psi,\gamma\in \mathcal{E}^1(X,\theta;\phi)$. We prove the triangle inequality:
    \[
    d_1(\varphi,\psi)\leq d_1(\varphi,\gamma)+d_1(\psi,\gamma).
    \]
    This can be translated to
    \[
    E_{\theta}^{\phi}(\varphi\land \gamma)-E_{\theta}^{\phi}(\varphi\land \psi)\leq E_{\theta}^{\phi}(\gamma)-E_{\theta}^{\phi}(\psi\land \gamma).
    \]
    We just have to compute using \cref{prop:d1contract}:
    \[
    \begin{aligned}
        E_{\theta}^{\phi}(\gamma)-E_{\theta}^{\phi}(\psi\land \gamma)\geq & E_{\theta}^{\phi}(\varphi\land \gamma)-E_{\theta}^{\phi}(\varphi\land \psi\land \gamma)\\
        \geq &E_{\theta}^{\phi}(\varphi\land \gamma)-E_{\theta}^{\phi}(\varphi\land \psi),
    \end{aligned}
    \]
    where the second line follows from \cref{lma:EmonoEinfty}.
\end{proof}

We introduce an auxiliary functional.
\begin{definition}\label{def:Itheta}
    Define $I_{\theta}\colon \mathcal{E}^1(X,\theta;\phi)\times \mathcal{E}^1(X,\theta;\phi)\rightarrow \mathbb{R}\cup \{\infty\}$ as follows:
    \begin{equation}
        I_{\theta}(\varphi,\psi)\coloneqq \int_X |\varphi-\psi|\left( \theta_{\varphi}^n+\theta_{\psi}^n \right).
    \end{equation}
    \index{$I_{\theta}$}
\end{definition}
We will soon see that $I_{\theta}$ is always finite, see \cref{thm:I_d1_comp}.

We observe the following elementary equality.
\begin{lemma}\label{lma:I_equal_Imaxsum}
    Let $\varphi,\psi\in \mathcal{E}^1(X,\theta;\phi)$. Then
    \[
    I_{\theta}(\varphi,\psi)=I_{\theta}(\varphi\lor \psi,\varphi)+I_{\theta}(\varphi\lor \psi,\psi).
    \]
\end{lemma}
\begin{proof}
    By the plurifine locality of non-pluripolar products \itemref{itm:prop:npp1:1}, we can write
    \[
    \begin{aligned}
         I_{\theta}(\varphi,\psi)=&\int_{ \{\varphi<\psi\}} (\psi-\varphi) \left( \theta_{\varphi}^n+\theta_{\psi}^n \right)+\int_{ \{\varphi>\psi\}} (\varphi-\psi) \left( \theta_{\varphi}^n+\theta_{\psi}^n \right),\\
         I_{\theta}(\varphi\lor \psi,\varphi)=& \int_{\{\psi>\varphi\}} (\psi-\varphi)\left(\theta_{\psi}^n+\theta_{\varphi}^n \right),\\
         I_{\theta}(\varphi\lor \psi,\psi)=& \int_{\{\psi<\varphi\}} (\varphi-\psi)\left(\theta_{\psi}^n+\theta_{\varphi}^n \right).
    \end{aligned}
    \]
    The desired equality follows by comparing these terms.
\end{proof}

We have an interesting relation between $d_1$ and $I_{\theta}$ defined in \cref{def:Itheta}.
\begin{theorem}\label{thm:I_d1_comp}
    Let $\varphi,\psi\in \mathcal{E}^1(X,\theta;\phi)$. Then
    \begin{equation}\label{eq:I_d1_comp}
    \frac{1}{C_n} I_{\theta}(\varphi,\psi)\leq d_1(\varphi,\psi)\leq I_{\theta}(\varphi,\psi),
    \end{equation}
    where $C_n=3\cdot 2^{n+2}(n+1)$.
\end{theorem}
\begin{proof}
    \textbf{Step~1}. We first prove the right-hand part of \eqref{eq:I_d1_comp}.

    Thanks to \cref{prop:Ediffbdd} and \cref{lma:rooftopMA}, we have
    \[
    \begin{aligned}
    E_{\theta}^{\phi}(\varphi)-E_{\theta}^{\phi}(\varphi\land \psi)\leq & \int_X (\varphi-\varphi\land\psi)\,\theta_{\varphi\land \psi}^n \\
    \leq & \int_{\{\psi=\varphi\land \psi \}} (\varphi-\psi)\,\theta_{\psi}^n\\
    \leq & \int_X |\varphi-\psi|\,\theta_{\psi}^n.
    \end{aligned}
    \]
    Similarly,
    \[
    E_{\theta}^{\phi}(\psi)-E_{\theta}^{\phi}(\varphi\land \psi)\leq \int_X |\varphi-\psi|\,\theta_{\varphi}^n.
    \]
    Adding these inequalities up, we find
    \[
    d_1(\varphi,\psi)\leq I_{\theta}(\varphi,\psi).
    \]

    \textbf{Step~2}. We prove the left-hand part of \eqref{eq:I_d1_comp}.

    We claim that
    \begin{equation}\label{eq:d_1varphivarphipsimean_temp1}
        d_1\left(\varphi,\frac{\varphi+\psi}{2}\right)\leq \frac{3(n+1)}{2} d_1(\varphi,\psi).
    \end{equation}
    For this purpose, we compute directly
    \[
    \begin{aligned}
        &d_1\left(\varphi,\frac{\varphi+\psi}{2}\right)\\
        =&d_1\left(\varphi,\varphi\land \frac{\varphi+\psi}{2}\right)+d_1\left(\frac{\varphi+\psi}{2},\varphi\land \frac{\varphi+\psi}{2}\right)\\
        \leq & d_1\left(\varphi,\varphi\land \psi\right)+d_1\left(\frac{\varphi+\psi}{2},\varphi\land\psi\right)&\textup{by \cref{lma:EmonoEinfty}}\\
        \leq & \int_X \left( \varphi-\varphi\land\psi\right)\,\theta_{\varphi\land\psi}^n+ \int_X \left( \frac{\varphi+\psi}{2}-\varphi\land\psi\right)\,\theta_{\varphi\land \psi}^n &\textup{by \cref{prop:Ediffbdd}}\\
        =& \frac{3}{2}\int_X \left( \varphi-\varphi\land\psi\right)\,\theta_{\varphi\land\psi}^n+\frac{1}{2}\int_X \left( \psi-\varphi\land\psi\right)\,\theta_{\varphi\land\psi}^n \\
        \leq & \frac{3(n+1)}{2}d_1\left(\varphi,\varphi\land \psi\right) + \frac{n+1}{2}d_1(\psi,\varphi\land \psi)& \textup{by \eqref{eq:Ecocyc}}\\
        \leq & \frac{3(n+1)}{2}d_1(\varphi,\psi),
    \end{aligned}
    \]
    and \eqref{eq:d_1varphivarphipsimean_temp1} follows.

    We now estimate the left-hand side of \eqref{eq:d_1varphivarphipsimean_temp1}:
    \[
    \begin{aligned}
    d_1\left(\varphi,\frac{\varphi+\psi}{2}\right)\geq & d_1\left(\varphi,\varphi\land \frac{\varphi+\psi}{2}\right)\\
    \geq & \int_X \left(\varphi-\varphi\land \frac{\varphi+\psi}{2}\right)\,\theta_{\varphi}^n\\
    \geq & \int_{\{\varphi>\psi\}}\left(\varphi-\varphi\land \frac{\varphi+\psi}{2}\right)\,\theta_{\varphi}^n \\
    \geq & \int_{\{\varphi>\psi\}}\left(\varphi-\frac{\varphi+\psi}{2}\right)\,\theta_{\varphi}^n\\
    =& \frac{1}{2}\int_{\{\varphi>\psi\}}\left(\varphi-\psi\right)\,\theta_{\varphi}^n,
    \end{aligned}
    \]
    where the second inequality follows from \cref{prop:Ediffbdd}.

    Similarly,
    \[
    \begin{aligned}
        d_1\left(\varphi,\frac{\varphi+\psi}{2}\right)\geq & d_1\left(\frac{\varphi+\psi}{2},\varphi\land \frac{\varphi+\psi}{2}\right)\\
        \geq & \int_X \left(\frac{\varphi+\psi}{2}-\varphi\land \frac{\varphi+\psi}{2}\right)\,\theta_{(\varphi+\psi)/2}^n\\
        \geq & 2^{-n} \int_X \left(\frac{\varphi+\psi}{2}-\varphi\land \frac{\varphi+\psi}{2}\right)\,\theta_{\varphi}^n\\
        \geq & 2^{-n} \int_{\{\varphi<\psi\}} \left(\frac{\varphi+\psi}{2}-\varphi\land \frac{\varphi+\psi}{2}\right)\,\theta_{\varphi}^n\\
        \geq & 2^{-n} \int_{\{\varphi<\psi\}} \left(\frac{\varphi+\psi}{2}-\varphi\right)\,\theta_{\varphi}^n\\
        =& 2^{-n-1} \int_{\{\varphi<\psi\}} \left(\psi-\varphi\right)\,\theta_{\varphi}^n,
    \end{aligned}
    \]
    where the third inequality follows from the multi-linearity of the non-pluripolar product \itemref{itm:prop:npp1:6}.

    Adding these estimates up, we find
    \[
        2^{n+2} d_1\left(\varphi,\frac{\varphi+\psi}{2}\right)\geq \left(2^{n+1}+2\right) d_1\left(\varphi,\frac{\varphi+\psi}{2}\right)\geq \int_X |\varphi-\psi|\,\theta_{\varphi}^n.
    \]
    Therefore, combined with \eqref{eq:d_1varphivarphipsimean_temp1} we conclude
    \[
    3\cdot 2^{n+1}(n+1)d_1(\varphi,\psi)\geq \int_X |\varphi-\psi|\,\theta_{\varphi}^n.
    \]
    By symmetry, we also get a similar expression after exchanging $\varphi$ and $\psi$. Adding these inequalities together, we find
    \[
    3\cdot 2^{n+2}(n+1)d_1(\varphi,\psi)\geq I_{\theta}(\varphi,\psi).
    \]
    The left-hand part of \eqref{eq:I_d1_comp} then follows.
\end{proof}

\begin{corollary}\label{cor:d1maxonearg}
    Let $\varphi,\psi\in \mathcal{E}^1(X,\theta;\phi)$. Then
    \begin{equation}\label{eq:d1maxonearg}
        d_1(\varphi\lor \psi,\varphi)\leq C_n d_1(\varphi,\psi),
    \end{equation}
    where $C_n$ is the constant in \cref{thm:I_d1_comp}.
\end{corollary}
\begin{proof}
    By \cref{thm:I_d1_comp} and \cref{lma:I_equal_Imaxsum}, we have
    \[
        d_1(\varphi\lor\psi,\varphi)\leq I_{\theta}(\varphi\lor\psi,\varphi)\leq I_{\theta}(\varphi,\psi)\leq C_n d_1(\varphi,\psi).
    \]
    This completes the proof.
\end{proof}

\begin{lemma}\label{lma:d1bdsupX}
    There is $A,B>0$ so that for any $\varphi\in \mathcal{E}^1(X,\theta;\phi)$,  we have
    \begin{equation}\label{eq:d1bdsupX}
    d_1(\varphi,\phi)\geq -\int_X \theta_{\phi}^n\cdot \sup_X \varphi\geq -Ad_1(\varphi,\phi)-B.
    \end{equation}
\end{lemma}
\begin{proof}
    When $\sup_X \varphi\leq 0$, the right-hand part of \eqref{eq:d1bdsupX} is trivial. Since $\varphi\leq \phi$, we find
    \[
    d_1(\phi,\varphi)=-E_{\theta}^{\phi}(\varphi)\geq -\int_X \theta_{\phi}^n\cdot\sup_{\{\phi\neq -\infty\}} (\varphi-\phi)=-\int_X \theta_{\phi}^n\cdot\sup_X \varphi,
    \]
    where the last equality follows from \cref{prop:model_sup_prop}.

    We can therefore assume that $\sup_X \varphi> 0$. In this case, the left-hand part of \eqref{eq:d1bdsupX} is trivial. Take a Kähler form $\omega\geq \theta$ on $X$. Then thanks to \cref{thm:Pvarphisupport}, we can find a constant $C>0$ so that
    \[
    \theta_{\phi}^n\leq C\omega^n.
    \]
    Thanks to \cref{prop:L1compa}, there is a constant $C'>0$, independent of the choice of $\varphi$, so that
    \[
    \int_X \left( \phi-\varphi+\sup_X \varphi \right)\,\theta_{\phi}^n\leq C'.
    \]
    We estimate
    \[
    \begin{aligned}
    I_{\theta}(\phi,\varphi)\geq & \int_X |\varphi-\phi|\,\theta_{\phi}^n\\
    \geq & \left(\sup_X \varphi\right)\int_X \theta_{\phi}^n-\int_X \left( \phi-\varphi+\sup_X \varphi \right)\,\theta_{\phi}^n \\
    \geq & \left(\sup_X \varphi\right)\int_X \theta_{\phi}^n-C'.
    \end{aligned}
    \]
    The right-hand part of \eqref{eq:d1bdsupX} then follows from \cref{thm:I_d1_comp}.
\end{proof}

Next we handle the completeness of $d_1$. The completeness can be proved in a more general framework.
\begin{definition}\label{def:rays:L858}
	Let $E$ be a set. A \emph{pre-rooftop structure} on $E$ is a binary operator $\land:E\times E\rightarrow E$, satisfying the following axioms: For $x,y,z\in E$,
	\begin{enumerate}
		\item\label{itm:chapter-rays:L841:1} $x\land y=y\land x$.
		\item\label{itm:chapter-rays:L841:2} $(x\land y)\land z=x\land (y\land z)$.
		\item\label{itm:chapter-rays:L841:3} $x\land x=x$.
	\end{enumerate}
	We call $(E,\land)$ a \emph{pre-rooftop space}\index{pre-rooftop space}.
\end{definition}

A pre-rooftop structure $\land$ defines a partial order $\leq$ on $E$ as follows:
\[
x\leq y\quad \text{if and only if } \quad x\land y=x.
\]
Here by abuse of notation, we use $\leq$ to denote the partial order.

In particular, it makes sense to talk about increasing and decreasing sequences in $E$.

\begin{definition}\label{def:roof}
	Let $(E,d)$ be a metric space. A \emph{pre-rooftop structure}\index{pre-rooftop structure} on $(E,d)$ is a pre-rooftop structure $\land$ on $E$. We say $(E,d,\land)$ is a \emph{pre-rooftop metric space}\index{pre-rooftop metric space}.

	A \emph{rooftop structure}\index{rooftop structure} on $(E,d)$ is a pre-rooftop structure $\land$ on $E$ such that
	\begin{equation}\label{eq:rtdef}
		d(x\land z,y\land z)\leq d(x,y),\quad \forall x,y,z\in E.
	\end{equation}
	We call $(E,d,\land)$ a \emph{rooftop metric space}\index{rooftop metric space}.
\end{definition}

\begin{lemma}\label{lma:landdis}
	Let $(E,d,\land)$ be a rooftop metric space. Let $x,y,x',y'\in E$. Then
	\begin{equation}
		d(x\land y,x'\land y')\leq d(x,x')+d(y,y').
	\end{equation}
\end{lemma}
In particular, $\land\colon E\times E\rightarrow E$ is continuous.
\begin{proof}
	We compute
	\[
	d(x\land y,x'\land y')\leq d(x\land y,x\land y')+d(x\land y',x'\land y')\leq d(x,x')+d(y,y').
	\]
	This completes the proof.
\end{proof}

\begin{proposition}\label{prop:rooftopcomp}
	Let $(E,d,\land)$ be a rooftop metric space. Then $(E,d)$ is complete if and only if both of the following hold:
	\begin{enumerate}
		\item\label{itm:prop:rooftopcomp:1} Each increasing Cauchy sequence converges.
		\item\label{itm:prop:rooftopcomp:2} Each decreasing Cauchy sequence converges.
	\end{enumerate}
\end{proposition}
\begin{proof}
	The direct implication is trivial.

	Conversely, assume that both conditions are true.
	Let $(x_j)_{j>0}$ be a Cauchy sequence in $E$. We want to prove that $(x_j)_j$ converges. By passing to a subsequence, we may assume that
	\[
	d\left(x_j,x_{j+1}\right)\leq 2^{-j}.
	\]
	For $k,j\geq 1$, let
	\[
	y_j^k\coloneqq x_k\land \cdots \land x_{k+j}.
	\]
	Then $(y_j^k)_j$ is decreasing, and
	\[
	d\left(y_j^k,y_{j+1}^k\right)\leq d\left(x_{k+j},x_{k+j+1}\right)\leq 2^{-k-j}.
	\]
	So $(y_j^k)_j$ is a decreasing Cauchy sequence. Define
	\[
	y^k\coloneqq \lim_{j\to\infty}y_j^k.
	\]
	Then
	\[
	d\left(y^k,y^{k+1}\right)=\lim_{j\to\infty}d\left(y_{j+1}^k,y_j^{k+1}\right)\leq d\left(x_k,x_{k+1}\right)\leq 2^{-k}.
	\]
    We also need to check that the sequence $(y^k)_k$ is increasing. For each $j\geq 1$, we have
    \[
        y_j^k=x_k\land y_{j-1}^{k+1}\leq y_{j-1}^{k+1}.
    \]
    Equivalently,
    \[
        y_j^k\land y_{j-1}^{k+1}=y_j^k.
    \]
    Letting $j\to\infty$ and using the continuity of $\land$ from \cref{lma:landdis}, we obtain
    \[
        y^k\land y^{k+1}=y^k.
    \]
    Thus $y^k\leq y^{k+1}$.
	So $(y^k)_k$ is an increasing Cauchy sequence. Let
	\[
	y\coloneqq \lim_{k\to\infty}y^k.
	\]
	For every $j\geq 0$, we have
	\[
	    d\left(y_j^k,x_k\right)\leq \sum_{r=k}^{k+j-1}d\left(x_r,x_{r+1}\right),
	\]
	with the empty sum understood as $0$. Indeed, this is trivial for $j=0$, and the induction step follows from
	\[
	    y_j^k=x_k\land y_{j-1}^{k+1}
	\]
	and
	\[
	    d\left(y_j^k,x_k\right)\leq d\left(y_{j-1}^{k+1},x_k\right)\leq d\left(y_{j-1}^{k+1},x_{k+1}\right)+d\left(x_{k+1},x_k\right).
	\]
	Letting $j\to\infty$, we get
	\[
	d\left(y^k,x_k\right)\leq \sum_{r=k}^{\infty}d\left(x_r,x_{r+1}\right)\leq 2^{1-k}.
	\]
	So $(x_k)_k$ converges to $y$.
\end{proof}

\begin{theorem}\label{thm:d1_complete}
    The metric space $(\mathcal{E}^1(X,\theta;\phi),d_1)$ is complete.
\end{theorem}
It follows from the proof below that if $(\varphi_j)_{j>0}$ is a monotone sequence in $\mathcal{E}^1(X,\theta;\phi)$ with $d_1$-limit $\varphi$, then $\varphi$ is the almost everywhere limit of $(\varphi_j)_{j>0}$. We will use this result without further explanation in the sequel.
\begin{proof}
    As we have seen in \cref{lma:d1metric} and \cref{prop:d1contract}, the triple
    \[
    \left(\mathcal{E}^1(X,\theta;\phi),d_1,\land \right)
    \]
    is a rooftop metric space. Hence thanks to \cref{prop:rooftopcomp}, it remains to show that each increasing or decreasing Cauchy sequence in $\mathcal{E}^1(X,\theta;\phi)$ converges.

    We first consider an increasing sequence $(\varphi_i)_{i>0}$ in $\mathcal{E}^1(X,\theta;\phi)$. The Cauchy property simply means that
    \[
    \sup_{i>0} E_{\theta}^{\phi}(\varphi_i)<\infty.
    \]
    We claim that
    \[
    \varphi\coloneqq \uscsup[i>0]\varphi_i
    \]
    is the $d_1$-limit of the sequence. We first observe that $(\sup_X \varphi_i)_{i>0}$ is bounded, as a consequence of \cref{lma:d1bdsupX}. Therefore, \cref{prop:supsEE1} guarantees that $\varphi\in \mathcal{E}^1(X,\theta;\phi)$, and hence our assertion follows from \cref{prop:Econt_mono}.

    Next we consider a decreasing sequence $(\varphi_i)_{i>0}$ in $\mathcal{E}^1(X,\theta;\phi)$. The Cauchy property simply means that
    \[
    \inf_{i>0} E_{\theta}^{\phi}(\varphi_i)>-\infty.
    \]
    We claim that
    \[
     \varphi\coloneqq \inf_{i>0}\varphi_i
    \]
    is the $d_1$-limit of the sequence. Equivalently, we need to show that $\varphi\in \mathcal{E}^1(X,\theta;\phi)$ and
    \begin{equation}\label{eq:E_cont_decnet_temp1}
        E_{\theta}^{\phi}(\varphi)=\lim_{i\to\infty} E_{\theta}^{\phi}(\varphi_i).
    \end{equation}

    We first observe that $(\sup_X \varphi_i)_{i>0}$ is bounded, as a consequence of \cref{lma:d1bdsupX}. Therefore, $\varphi\in \PSH(X,\theta;\phi)$ and $\varphi_i\xlongrightarrow{L^1} \varphi$, see \cref{cor:decreasingnetL1conv}. We claim that $\varphi\in \mathcal{E}^1(X,\theta;\phi)$. For this purpose, we may assume that $\varphi_1\leq 0$.
    Since $(\varphi_i)_i$ is decreasing, $\varphi_i\leq \varphi_1\leq 0$ for all $i>0$. As $\varphi_i\preceq\phi$ and $\phi$ is model, each $\varphi_i$ is a candidate for $P_\theta[\phi]=\phi$, so $\varphi_i\leq \phi$. Thus, using the energy formula \eqref{eq:Edefbdd} and the fact that all mixed terms are non-positive, we get for each $i>0$,
    \[
        \frac{1}{n+1}\int_X \left(\varphi_i-\phi\right)\,\theta_{\varphi_i}^n\geq E_{\theta}^{\phi}(\varphi_i)\geq \inf_{j>0} E_{\theta}^{\phi}(\varphi_j)>-\infty.
    \]
    We can then apply \cref{lma:stab}, with $\gamma=\phi$, to conclude that $\varphi\in \mathcal{E}^1(X,\theta;\phi)$.
    Our assertion \eqref{eq:E_cont_decnet_temp1} then follows from \cref{prop:E_cont_decnet}. We conclude the proof.
\end{proof}

\begin{lemma}\label{lma:d1conv_inbetween_mono}
    Let $(\varphi_i)_{i>0}$ be a sequence in $\mathcal{E}^1(X,\theta;\phi)$ with $d_1$-limit $\varphi\in \mathcal{E}^1(X,\theta;\phi)$. Then after replacing $(\varphi_i)_{i>0}$ by a subsequence, we can find two sequences $(\psi_i)_{i>0}$ and $(\eta_i)_{i>0}$ in $\mathcal{E}^1(X,\theta;\phi)$ such that
    \begin{enumerate}
        \item\label{itm:lma:d1conv_inbetween_mono:1} $(\psi_i)_{i>0}$ is decreasing with $d_1$ and pointwise limit $\varphi$;
        \item\label{itm:lma:d1conv_inbetween_mono:2} $(\eta_i)_{i>0}$ is increasing with $d_1$ and almost everywhere limit $\varphi$.
        \item\label{itm:lma:d1conv_inbetween_mono:3} $\eta_i\leq \varphi_i\leq \psi_i$ for all $i>0$.
    \end{enumerate}
\end{lemma}
\begin{proof}
    We first note that as a consequence of \cref{lma:d1bdsupX}, $(\sup_X \varphi_i)_{i>0}$ is bounded.
    Let $C_n$ be the constant in \cref{thm:I_d1_comp}. After replacing $(\varphi_i)_{i>0}$ by a subsequence, we may assume that
    \[
    d_1(\varphi_i,\varphi_{i+1})\leq \min\{C_n^{-2i},2^{-i}\}.
    \]

    We first construct $(\psi_i)_{i>0}$. For this purpose, it suffices to define
    \[
    \psi_i\coloneqq \uscsup[j\geq i] \varphi_j
    \]
    for each $i>0$.
    It follows from \cref{prop:supsEE1} that $\psi_i\in \mathcal{E}^1(X,\theta;\phi)$ for each $i>0$.
    For $k\geq 0$, set
    \[
        \psi_i^k\coloneqq \varphi_i\lor\varphi_{i+1}\lor\cdots\lor\varphi_{i+k}.
    \]
    Iterating \cref{cor:d1maxonearg} gives
    \[
    d_1(\varphi_i,\psi_i^k)\leq \sum_{a=i}^{i+k-1}C_n^{a+1-i}d_1(\varphi_a,\varphi_{a+1})\leq \sum_{a=i}^{i+k-1}C_n^{a+1-i}C_n^{-2a}.
    \]
    Since $\psi_i^k$ increases almost everywhere to $\psi_i$, letting $k\to\infty$ and using \cref{prop:Econt_mono}, we get
    \[
        d_1(\varphi_i,\psi_i)\leq \frac{C_n^{1-2i}}{1-C_n^{-1}}.
    \]
    Thus $\psi_i\to\varphi$ in $d_1$. Since $(\psi_i)_i$ is decreasing, the monotone case in the proof of \cref{thm:d1_complete} shows that $\psi_i$ converges almost everywhere to $\varphi$, hence pointwise by \cref{prop:pshfuncdetdense}.

    Next we construct $(\eta_i)_{i>0}$.
    For each $i>0$ and $k\geq 0$, we let
    \[
    \eta_i^k\coloneqq \varphi_i\land \varphi_{i+1}\land \cdots\land \varphi_{i+k}.
    \]
    Then
    \[
    \begin{aligned}
    d_1(\varphi_i,\eta_i^k)\leq & \sum_{j=0}^{k-1} d_1\left(\eta_i^j,\eta_i^{j+1}\right)\\
    \leq &\sum_{j=0}^{k-1} d_1\left(\varphi_{i+j},\varphi_{i+j+1}\right) &\textup{by \cref{prop:d1contract}}\\
    \leq & \sum_{j=0}^{k-1} 2^{-i-j}\\
    \leq & 2^{1-i}.
    \end{aligned}
    \]
    Since $\eta_i^k\leq \varphi_i$, the preceding estimate gives
    \[
        E_{\theta}^{\phi}(\eta_i^k)\geq E_{\theta}^{\phi}(\varphi_i)-2^{1-i},\quad k\geq 0.
    \]
    The sequence $(\eta_i^k)_{k\geq 0}$ is decreasing, and the monotonicity of $E_{\theta}^{\phi}$ shows that the sequence $\left(E_{\theta}^{\phi}(\eta_i^k)\right)_{k\geq 0}$ is decreasing and bounded from below, hence convergent. Thus, for $k\leq l$,
    \[
        d_1(\eta_i^k,\eta_i^l)=E_{\theta}^{\phi}(\eta_i^k)-E_{\theta}^{\phi}(\eta_i^l),
    \]
    and $(\eta_i^k)_{k\geq 0}$ is a decreasing Cauchy sequence. Therefore, the monotone case in the proof of \cref{thm:d1_complete} shows that
    \[
    \eta_i\coloneqq \lim_{k\to\infty}\eta_i^k=\inf_{k\geq 0} \eta_i^k
    \]
    is the $d_1$-limit of $(\eta_i^k)_{k\geq 0}$ and \cref{prop:Econt_mono} shows that
    \[
    d_1(\varphi_i,\eta_i)\leq 2^{1-i}.
    \]
    By construction, $\eta_i\leq \varphi_i$ and $(\eta_i)_i$ is increasing. Therefore, $\varphi$ is the $d_1$-limit of the increasing sequence $(\eta_i)_{i>0}$. As we have seen in the proof of \cref{thm:d1_complete}, this implies that $\varphi$ is also the almost everywhere limit of $(\eta_i)_{i>0}$.
\end{proof}

\begin{theorem}\label{thm:Efuncont}
    The functional $E_{\theta}^{\phi}\colon \mathcal{E}^1(X,\theta;\phi) \rightarrow \mathbb{R}$ is continuous.
\end{theorem}
\begin{proof}
    Let $(\varphi_j)_{j>0}$ be a sequence in $\mathcal{E}^1(X,\theta;\phi)$ with $d_1$-limit $\varphi\in \mathcal{E}^1(X,\theta;\phi)$. We wish to show that
    \begin{equation}\label{eq:limjE_temp1}
        \lim_{j\to\infty} E_{\theta}^{\phi}(\varphi_j)=E_{\theta}^{\phi}(\varphi).
    \end{equation}
    It suffices to prove that every subsequence has a further subsequence along which \eqref{eq:limjE_temp1} holds. Thus, after replacing $(\varphi_j)_{j>0}$ by a subsequence, let $(\eta_j)_j$ and $(\psi_j)_j$ be the monotone sequences supplied by \cref{lma:d1conv_inbetween_mono}. From the construction, $\eta_j\leq \varphi_j\leq \psi_j$. By \cref{prop:Econt_mono},
    \[
        E_{\theta}^{\phi}(\eta_j)\to E_{\theta}^{\phi}(\varphi),\quad E_{\theta}^{\phi}(\psi_j)\to E_{\theta}^{\phi}(\varphi).
    \]
    The monotonicity of $E_{\theta}^{\phi}$ then squeezes $E_{\theta}^{\phi}(\varphi_j)$, proving \eqref{eq:limjE_temp1}.
\end{proof}

\begin{proposition}\label{prop:energylinear}
    Let $(\varphi_t)_{t\in [a,b]}$ be a geodesic in $\mathcal{E}^1(X,\theta;\phi)$. Then $t\mapsto E_{\theta}^{\phi}(\varphi_t)$ is a linear function on $[a,b]$.
\end{proposition}
We recall the terminology used in the proof. Let $\Omega\subset\mathbb{C}^m$ be Euclidean open. A function $u\in\PSH(\Omega)$ is \emph{maximal} if, for every open $G\Subset\Omega$ and every upper semicontinuous function $v$ on $\overline G$ such that $v\in\PSH(G)$ and $v\leq u$ on $\partial G$, one has $v\leq u$ on $G$. Let $\mathcal F$ denote the plurifine topology. If $U$ is $\mathcal F$-open and $u$ is $\mathcal F$-psh on $U$, then $u$ is \emph{$\mathcal F$-maximal} if, for every bounded $\mathcal F$-open set $G$ with $G\subset\overline G^{\mathcal F}\subset U$ and every bounded-above $\mathcal F$-upper semicontinuous function $v$ on $\overline G^{\mathcal F}$ that is $\mathcal F$-psh on $G$, one has
\[
    v\leq u\textup{ on }\partial_{\mathcal F}G
    \quad\Longrightarrow\quad
    v\leq u\textup{ on }G.
\]
See \cite[Definition~5.9.1]{EKF}. For Euclidean maximality, see also \cite[Definition~9.6.5]{FPT}.
\begin{proof}
    By an affine change of parameter, we may assume that $[a,b]=[0,1]$.

    \textbf{Step~1}. Assume first that $\varphi_0,\varphi_1\in \mathcal{E}^{\infty}(X,\theta;\phi)$. Let $Y=X\times\mathbb{P}^1$, and let $p_1,p_2$ be the two projections. On
    \[
        A\coloneqq\{z\in\mathbb{C}:\mathrm{e}^{-1}<|z|<1\},
        \qquad U(x,z)\coloneqq\varphi_{-\log|z|}(x),
    \]
    the function $U$ is $p_1^*\theta$-psh. Choose $C>0$ such that
    \[
        \phi-C\leq \varphi_0,\varphi_1\leq\phi+C.
    \]
    The constant subgeodesic at $\phi-C$ and \cref{prop:convexsubgeod} give
    \begin{equation}\label{eq:relgeoduniformbound}
        \phi-C\leq\varphi_t\leq(1-t)\varphi_0+t\varphi_1\leq\phi+C,
        \qquad 0<t<1.
    \end{equation}

    Fix concentric annuli $A_0\Subset A_1\Subset A$ and a unit-mass Fubini--Study form $\omega_{\mathrm{FS}}$ on $\mathbb{P}^1$. Put
    \[
        \beta=p_2^*\omega_{\mathrm{FS}},\qquad
        \Theta_B=p_1^*\theta+B\beta,\qquad
        \Phi=p_1^*\phi.
    \]
    The class represented by $\Theta_B$ is big. Since $\phi=P_\theta[\phi]$, we have $\phi\leq0$. For every $B>0$, $\Phi=P_{\Theta_B}[\Phi]$. Indeed, if $\Psi\in\PSH(Y,\Theta_B)$ satisfies $\Psi\leq0$ and $\Psi\preceq\Phi$, then on each fibre $X\times\{z\}$, either $\Psi_z\equiv-\infty$, or $\Psi_z\in\PSH(X,\theta)$ is a candidate in the definition of $P_\theta[\phi]=\phi$. Thus $\Psi_z\leq\phi$ in both cases. Conversely, $\Phi$ itself is a candidate. Moreover, \itemref{itm:prop:npp1:6} gives
    \begin{equation}\label{eq:productmodelmass}
        \int_Y(\Theta_B+\ddc\Phi)^{n+1}
        =(n+1)B\int_X\theta_\phi^n>0.
    \end{equation}

    Choose an $S^1$-invariant smooth function $q$ on $\mathbb{P}^1$ such that
    \[
        q<-C-1\ \textup{ on }A_0,
        \qquad q>C+1\ \textup{ near }\partial A_1.
    \]
    After increasing $B$, we may assume that $B\omega_{\mathrm{FS}}+\ddc q\geq0$. Identifying $q$ with $p_2^*q$, define
    \[
        W=
        \begin{cases}
            U\lor(\Phi+q),&\textup{on }X\times A_1,\\
            \Phi+q,&\textup{on }Y\setminus(X\times A_1).
        \end{cases}
    \]
    The gluing lemma \cref{cor:glbdry} shows that $W\in\PSH(Y,\Theta_B)$. By construction,
    \begin{equation}\label{eq:compactgluedpotential}
        W=U\ \textup{ on }X\times A_0,
        \qquad \Phi-C'\leq W\leq\Phi+C'
    \end{equation}
    for some $C'>0$. Hence $W\in\mathcal{E}^{\infty}(Y,\Theta_B;\Phi)$ by \eqref{eq:productmodelmass}.

    Set
    \[
        T=\Theta_B+\ddc W,\qquad S=\Theta_B+\ddc\Phi.
    \]
    Define the order-zero current on $\mathbb{P}^1$
    \begin{equation}\label{eq:compactrelativeenergycurrent}
        \mathcal{H}_W\coloneqq\frac{1}{n+1}(p_2)_*\left((W-\Phi)\sum_{j=0}^n T^j\wedge S^{n-j}\right).
    \end{equation}
    The integration-by-parts formula \cite{Xia19b} (see also \cite[Theorem~17.1.3 and the proof of Theorem~17.7.3]{FPT}), applied against a smooth test function on $\mathbb{P}^1$, gives
    \begin{equation}\label{eq:compactrelativeenergyddc}
        \ddc\mathcal{H}_W
        =\frac{1}{n+1}(p_2)_*\left(T^{n+1}-S^{n+1}\right).
    \end{equation}

    By \eqref{eq:Edefbdd}, \eqref{eq:compactgluedpotential}, and \cref{cor:compactfibreslicing}, on $A_0$ the current $\mathcal{H}_W$ is represented almost everywhere by
    \[
        z\longmapsto E_\theta^\phi(\varphi_{-\log|z|}).
    \]
    This representative is continuous. Indeed, on every smaller annulus, the convex functions $t\mapsto\varphi_t(x)-\phi(x)$, considered on the common locus where both values are finite, are uniformly bounded by \eqref{eq:relgeoduniformbound}, hence uniformly Lipschitz. The resulting inequalities between the potentials hold almost everywhere and therefore everywhere by \cref{prop:pshfuncdetdense}. The monotonicity in \cref{lma:EmonoEinfty} and \eqref{eq:Ephiaddcst} then give the local Lipschitz continuity of $t\mapsto E_\theta^\phi(\varphi_t)$.

    On $X\times A_0$, define the signed Radon measure
    \begin{equation}\label{eq:horizontaldefectrelative}
        \mu_U\coloneqq
        \left(T^{n+1}-(n+1)B\beta\wedge T^n\right)|_{X\times A_0}.
    \end{equation}
    This measure is positive. To see this, take a relatively compact product coordinate ball $Q\subset X\times A_0$. Let $a$ be a smooth local potential of $p_1^*\theta$, put $h=U+a\leq0$ and $h_k=h\lor(-k)$, and choose a smooth local potential $b$ of $B\beta$ with $|b|\leq C_Q$. For $p=n,n+1$, set
    \[
        \begin{aligned}
        \nu_k^{b,p}&=\mathds{1}_{\{h+b>-k\}}
        \bigl(\ddc((h+b)\lor(-k))\bigr)^p,\\
        \nu_k^{h,p}&=\mathds{1}_{\{h>-k\}}
        \bigl(\ddc h_k+B\beta\bigr)^p.
        \end{aligned}
    \]
    Plurifine locality and $|b|\leq C_Q$ give, for $k>C_Q$,
    \[
        \nu_{k-C_Q}^{h,p}\leq\nu_k^{b,p}\leq\nu_{k+C_Q}^{h,p}.
    \]
    Thus the two cofinal truncation families have the same limit. By \itemref{itm:prop:npp1:6} and $\beta^2=0$,
    \begin{equation}\label{eq:horizontalfinitelevelrelative}
        \begin{aligned}
        &\mathds{1}_{\{h>-k\}}\left[
        (\ddc h_k+B\beta)^{n+1}
        -(n+1)B\beta\wedge(\ddc h_k+B\beta)^n\right]\\
        &\qquad=\mathds{1}_{\{h>-k\}}(\ddc h_k)^{n+1}\geq0.
        \end{aligned}
    \end{equation}
    Passing to the increasing truncation limit yields
    \begin{equation}\label{eq:horizontaltruncationrelative}
        \mu_U|_Q=\lim_{k\to\infty}
        \mathds{1}_{\{h>-k\}}(\ddc h_k)^{n+1}\geq0.
    \end{equation}

    By \eqref{eq:compactfibreslicingconstant},
    \[
        (p_2)_*(\beta\wedge T^n)
        =\left(\int_X\theta_\phi^n\right)\omega_{\mathrm{FS}},
    \]
    since $W_z\sim\phi$ on every fibre and \cref{prop:Ppresmass} applies. On the other hand,
    \[
        (p_2)_*S^{n+1}
        =(n+1)B\left(\int_X\theta_\phi^n\right)\omega_{\mathrm{FS}}
    \]
    by \itemref{itm:prop:npp1:6}. Together with \eqref{eq:compactrelativeenergyddc} and \eqref{eq:horizontaldefectrelative}, these identities give on $A_0$
    \begin{equation}\label{eq:relativeenergycurvature}
        \ddc_z E_\theta^\phi(\varphi_{-\log|z|})
        =\frac{1}{n+1}(p_2)_*\mu_U.
    \end{equation}

    \textbf{Step~2}. We prove that $\mu_U=0$. Let $D\Subset X\times A$ be a coordinate domain, let $a$ be a smooth local potential of $p_1^*\theta$ near $\overline D$, and set $h=U+a$. If $v$ is a psh function on $D$ with the boundary inequality required in the definition of maximality, max-pasting gives
    \[
        \widetilde U=U\lor(v-a)\ \textup{ on }D,
        \qquad \widetilde U=U\ \textup{ outside }D.
    \]
    The function $\widetilde U$ need not be radial. Its $S^1$-symmetrization
    \[
        V(x,z)=\uscsup[|\lambda|=1]\widetilde U(x,\lambda z)
    \]
    is radial and $p_1^*\theta$-psh. The $S^1$-saturation of $D$ is relatively compact in $X\times A$, so $V=U$ near $X\times\partial A$. Hence the family $V_t(x)=V(x,\mathrm{e}^{-t})$ is a subgeodesic from $\varphi_0$ to $\varphi_1$ by \cref{def:subgeod}. The definition \eqref{eq:Perron2} therefore gives
    \[
        \widetilde U\leq V\leq U.
    \]
    It follows that $v\leq h$ on $D$, so $h$ is a maximal psh function.

    On a product coordinate ball $Q\Subset X\times A_0$, put $O_k=Q\cap\{h>-k\}$. By \cite[Proposition~5.9.3]{EKF}, $h$ is $\mathcal F$-maximal. The function $h_k=h\lor(-k)$ is locally bounded and agrees with $h$ on the $\mathcal F$-open set $O_k$, hence is $\mathcal F$-maximal there. Therefore \cite[Theorem~5.10.1]{EKF} gives
    \[
        \mathds{1}_{O_k}(\ddc h_k)^{n+1}=0.
    \]
    See also \cite[Lemma~17.7.1]{FPT} for a direct vanishing proof. By \eqref{eq:horizontaltruncationrelative}, $\mu_U=0$. Since $A_0\Subset A$ was arbitrary, \eqref{eq:relativeenergycurvature} shows that $z\mapsto E_\theta^\phi(\varphi_{-\log|z|})$ is radial and harmonic on $A$. Thus $t\mapsto E_\theta^\phi(\varphi_t)$ is affine on $(0,1)$.

    Choose $C_0>0$ such that $|\varphi_1-\varphi_0|\leq C_0$ on their common finite locus; equivalently, the two resulting inequalities between the potentials hold everywhere. Affine barriers and \cref{prop:convexsubgeod} give
    \[
        \begin{aligned}
        \varphi_0-C_0t&\leq\varphi_t\leq\varphi_0+C_0t,\\
        \varphi_1-C_0(1-t)&\leq\varphi_t\leq\varphi_1+C_0(1-t).
        \end{aligned}
    \]
    By \cref{lma:EmonoEinfty} and \eqref{eq:Ephiaddcst}, the energy converges to the corresponding endpoint energy as $t\to0+$ and $t\to1-$. Hence
    \begin{equation}\label{eq:relativeenergylinearbounded}
        E_\theta^\phi(\varphi_t)
        =(1-t)E_\theta^\phi(\varphi_0)+tE_\theta^\phi(\varphi_1)
    \end{equation}
    when the endpoints belong to $\mathcal E^\infty(X,\theta;\phi)$.

    \textbf{Step~3}. Let now $\varphi_0,\varphi_1\in\mathcal E^1(X,\theta;\phi)$. By \cref{ex:lineargeod} and \eqref{eq:Ephiaddcst}, we may add different constants to the endpoints and assume that $\varphi_0,\varphi_1\leq\phi$. Put
    \[
        \varphi_j^{(k)}=\varphi_j\lor(\phi-k),\qquad j=0,1,
    \]
    and let $(\varphi_t^{(k)})_{t\in[0,1]}$ be the corresponding geodesic. The endpoint truncations decrease to $\varphi_j$, and Perron monotonicity gives
    \[
        \varphi_t\leq\varphi_t^{(k+1)}\leq\varphi_t^{(k)}.
    \]
    Let $h_t=\lim_{k\to\infty}\varphi_t^{(k)}$. The decreasing limit of the complexifications is $p_1^*\theta$-psh and is not identically $-\infty$. Taking the limit in the convex upper bound gives, first almost everywhere and then everywhere by \cref{prop:pshfuncdetdense},
    \begin{equation}\label{eq:finiteenergygeodesicsandwich}
        \varphi_t\leq h_t\leq(1-t)\varphi_0+t\varphi_1.
    \end{equation}

    Note that
    \[
        |h_t-\varphi_0|
        \leq|\varphi_t-\varphi_0|+t|\varphi_1-\varphi_0|\quad \textup{a.e.}.
    \]
    Since $\varphi_t\to\varphi_0$ in $L^1$ as $t\to0+$ and $\varphi_0,\varphi_1\in L^1(X)$, it follows that $h_t\to\varphi_0$ in $L^1$. The same argument at $t=1$ gives $h_t\to\varphi_1$ in $L^1$. Together with the $p_1^*\theta$-psh complexification constructed above, these endpoint limits show by \cref{def:subgeod} that $(h_t)$ is a subgeodesic from $\varphi_0$ to $\varphi_1$. Hence \eqref{eq:Perron2} gives $h_t\leq\varphi_t$. Together with \eqref{eq:finiteenergygeodesicsandwich}, this proves
    \[
        \varphi_t^{(k)}\searrow\varphi_t,\qquad 0<t<1.
    \]
    Applying \cref{prop:Econt_mono} at the endpoints and at each interior slice, and then using \eqref{eq:relativeenergylinearbounded}, gives
    \[
        \begin{aligned}
        E_\theta^\phi(\varphi_t)
        &=\lim_{k\to\infty}E_\theta^\phi(\varphi_t^{(k)})\\
        &=\lim_{k\to\infty}\left((1-t)E_\theta^\phi(\varphi_0^{(k)})
        +tE_\theta^\phi(\varphi_1^{(k)})\right)\\
        &=(1-t)E_\theta^\phi(\varphi_0)+tE_\theta^\phi(\varphi_1).
        \end{aligned}
    \]
    This concludes the proof.
\end{proof}

\begin{proposition}\label{prop:supsgeod}
    Let $(\varphi^i_0)_{i\in I}$, $(\varphi^i_1)_{i\in I}$ be two uniformly bounded from above increasing nets in $\mathcal{E}^{1}(X,\theta;\phi)$. Let $(\varphi^i_t)_{t\in (0,1)}$ be the geodesic from $\varphi_0^i$ to $\varphi_1^i$ for each $i\in I$. Then
    \[
    \left(\uscsup[i\in I]\varphi^i_t\right)_{t\in (0,1)}
    \]
    is the geodesic from $\uscsup[i] \varphi^i_0$ to $\uscsup[i] \varphi^i_1$.
\end{proposition}
\begin{proof}
    Put
    \[
        \varphi_j=\uscsup[i\in I]\varphi_j^i\quad(j=0,1),
        \qquad h_t=\uscsup[i\in I]\varphi_t^i,
    \]
    and let $(g_t)_{t\in(0,1)}$ be the geodesic from $\varphi_0$ to $\varphi_1$. By \cref{prop:supsEE1}, $\varphi_0,\varphi_1\in\mathcal E^1(X,\theta;\phi)$. The slice families are uniformly bounded from above by \cref{prop:convexsubgeod}; hence $h_t\in\mathcal E^1(X,\theta;\phi)$ by \cref{prop:supsEE1}. Moreover, \cref{prop:maxsubgeod} shows that $(h_t)$ is a subgeodesic from $\varphi_0$ to $\varphi_1$. Thus
    \begin{equation}\label{eq:suprelativegeodesicorder}
        h_t\leq g_t.
    \end{equation}

    For each $t$, the net $(\varphi_t^i)_{i\in I}$ is increasing. By \cref{prop:supsupstardiff}, its pointwise supremum agrees with $h_t$ almost everywhere, and the analogous assertion holds at the endpoints. Therefore \cref{prop:Econt_mono,prop:energylinear} gives
    \[
        \begin{aligned}
        E_\theta^\phi(h_t)
        &=\lim_{i\in I}E_\theta^\phi(\varphi_t^i)\\
        &=\lim_{i\in I}\left((1-t)E_\theta^\phi(\varphi_0^i)
        +tE_\theta^\phi(\varphi_1^i)\right)\\
        &=(1-t)E_\theta^\phi(\varphi_0)+tE_\theta^\phi(\varphi_1)
        =E_\theta^\phi(g_t).
        \end{aligned}
    \]
    By \eqref{eq:d1asEdiff}, \eqref{eq:suprelativegeodesicorder}, and \cref{lma:d1metric},
    \[
        d_1(h_t,g_t)=E_\theta^\phi(g_t)-E_\theta^\phi(h_t)=0,
    \]
    so $h_t=g_t$ for all $t\in(0,1)$.
\end{proof}

\begin{proposition}\label{prop:d1convex}
    Let $(\varphi_t)_{t\in [0,1]}$, $(\psi_t)_{t\in [0,1]}$ be geodesics in $\mathcal{E}^1(X,\theta;\phi)$. Then the distance $d_1(\varphi_t,\psi_t)$ is a convex function of $t\in [0,1]$.
\end{proposition}
\begin{proof}
    By definition of $d_1$ and \cref{prop:energylinear}, it suffices to show the concavity of
    \[
    [0,1]\ni t\mapsto E_{\theta}^{\phi}(\varphi_t\land \psi_t).
    \]
    Let $(\eta_t)_{t\in [0,1]}$ be the geodesic from $\varphi_0\land \psi_0$ to $\varphi_1\land \psi_1$. Then for each $t\in [0,1]$, we have $\eta_t\leq \varphi_t\land \psi_t$. Hence, by the monotonicity of $E_\theta^\phi$ and \cref{prop:energylinear},
    \[
    E_{\theta}^{\phi}( \varphi_t\land \psi_t)\geq E_{\theta}^{\phi}(\eta_t)=(1-t)E_{\theta}^{\phi}(\varphi_0\land \psi_0)+tE_{\theta}^{\phi}(\varphi_1\land \psi_1)
    \]
    for each $t\in [0,1]$. Applying the same argument to the restrictions of $(\varphi_t)_t$ and $(\psi_t)_t$ to any subinterval of $[0,1]$, which are again geodesics by \cref{prop:perronenvissubgeod}, gives a similar inequality on every subinterval. Hence $t\mapsto E_\theta^\phi(\varphi_t\land\psi_t)$ is concave, and our assertion follows.
\end{proof}

In particular, we can introduce:
\begin{definition}\label{def:d1rays}
    Let $\ell,\ell'\in \mathcal{R}^1(X,\theta;\phi)$. We define
    \index{Darvas metric!radial Darvas metric}
    \index{$d_1$!on geodesic rays}
    \[
        d_1(\ell,\ell')\coloneqq \lim_{t\to\infty}\frac{1}{t}d_1(\ell_t,\ell'_t).
    \]
    The limit exists in $[0,\infty]$ by \cref{prop:d1convex}.
\end{definition}

\begin{theorem}\label{thm:d1raycomplete}
    The function $d_1$ defined in \cref{def:d1rays} is a metric on $\mathcal{R}^1(X,\theta;\phi)$.
\end{theorem}
\begin{proof}
    Symmetry and the triangle inequality follow from the corresponding pointwise properties of $d_1$ on $\mathcal{E}^1(X,\theta;\phi)$.

    We next show that $d_1(\ell,\ell')<\infty$ for any $\ell,\ell'\in \mathcal{R}^1(X,\theta;\phi)$. Let $C_\ell$ be the constant in \cref{prop:raysuplinear}, and set $\tilde{\ell}_t=\ell_t-C_\ell t$. Then $\tilde{\ell}$ is again a geodesic ray, $\sup_X\tilde{\ell}_t=0$, and hence $\tilde{\ell}_t\leq \phi$ by \cref{prop:model_sup_prop}. Thus, by \cref{prop:energylinear},
    \[
    d_1(\tilde{\ell}_t,\phi)=-E_\theta^\phi(\tilde{\ell}_t)=-tE_\theta^\phi(\tilde{\ell}_1).
    \]
    Moreover, by \eqref{eq:Ephiaddcst}, $d_1(\ell_t,\tilde{\ell}_t)=|C_\ell|t\int_X\theta_{\phi}^n$. Hence $d_1(\ell_t,\phi)$ grows at most linearly in $t$. The same argument applies to $\ell'$, and the triangle inequality gives $d_1(\ell,\ell')<\infty$.

It remains to prove separation: suppose that $\ell,\ell'\in \mathcal{R}^1(X,\theta;\phi)$ and $d_1(\ell,\ell')=0$.

    Fix $s>0$. Then it follows from \cref{prop:d1convex} that
    \[
    \frac{d_1(\ell_s,\ell'_s)}{s}\leq \lim_{t\to\infty}\frac{d_1(\ell_t,\ell'_t)}{t}=d_1(\ell,\ell')=0.
    \]
    Therefore, $\ell_s=\ell'_s$ and hence $\ell=\ell'$.
\end{proof}

\begin{definition}\label{def:radialMAenergy2}
    We define the \emph{radial Monge--Amp\`ere energy}\index{radial Monge--Amp\`ere energy} $\mathbf{E}^{\phi}\colon \mathcal{R}(X,\theta;\phi)\rightarrow \mathbb{R}\cup \{\pm\infty\}$\index{$\mathbf{E}^{\phi}$} as follows:
    \[
    \mathbf{E}^{\phi}(\ell)\coloneqq \varlimsup_{t\to\infty}\frac{E_{\theta}^{\phi}(\ell_t)}{t}.
    \]
    When $\phi=V_{\theta}$, we write $\mathbf{E}$\index{$\mathbf{E}$} instead of $\mathbf{E}^{V_{\theta}}$.
\end{definition}

\begin{proposition}\label{prop:d1geod_diff_E}
    Let $\ell,\ell'\in \mathcal{R}^1(X,\theta;\phi)$ and $\ell\leq \ell'$. Then
    \begin{equation}\label{eq:d1rayscompa}
        d_1(\ell,\ell')=\mathbf{E}^{\phi}(\ell')-\mathbf{E}^{\phi}(\ell).
    \end{equation}
\end{proposition}
\begin{proof}
    By \cref{prop:energylinear}, the energy is linear along finite energy rays. Hence
    \[
        E_{\theta}^{\phi}(\ell_t)=tE_{\theta}^{\phi}(\ell_1),\quad E_{\theta}^{\phi}(\ell'_t)=tE_{\theta}^{\phi}(\ell'_1).
    \]
    Since $\ell_t\leq \ell'_t$ for all $t\geq 0$, \eqref{eq:d1asEdiff} gives
    \[
        \frac{1}{t}d_1(\ell_t,\ell'_t)=\frac{1}{t}\left(E_{\theta}^{\phi}(\ell'_t)-E_{\theta}^{\phi}(\ell_t)\right).
    \]
    Letting $t\to\infty$ proves the assertion.
\end{proof}

Next we recall the $\lor$ operator at the level of geodesic rays.
\begin{definition}\label{def:lorray1}
    Let $\ell,\ell'\in \mathcal{R}^{\infty}(X,\theta;\phi)$. We define $\ell\lor \ell'$\index{$\ell\lor \ell'$} as the minimal ray in $\mathcal{R}^{\infty}(X,\theta;\phi)$ lying above both $\ell$ and $\ell'$.
\end{definition}
\begin{proposition}\label{prop:lorrays}
    Given $\ell,\ell'\in \mathcal{R}^{\infty}(X,\theta;\phi)$. Then $\ell\lor \ell'\in \mathcal{R}^{\infty}(X,\theta;\phi)$ exists, and
    \begin{equation}\label{eq:Elor}
    \mathbf{E}^{\phi}(\ell\lor \ell')=\lim_{t\to\infty}\frac{1}{t}E_{\theta}^{\phi}\left(\ell_t\lor \ell'_t\right).
    \end{equation}
\end{proposition}
\begin{proof}
    For each $t>0$, let $(\ell''^t_{s})_{s\in [0,t]}$ be the geodesic from $\phi$ to $\ell_t\lor \ell'_t$.

    \textbf{Step~1}.
    We first show that for each fixed $s\geq 0$, $\ell''^t_{s}$ is increasing in $t\in [s,\infty)$.

    To see this, fix $s\geq 0$ and choose $t'>t\geq s$. We need to show that
    \begin{equation}\label{eq:ellppdombyellpptemp1}
    \ell''^{t'}_s\geq \ell''^{t}_s.
    \end{equation}
    Since $(\ell''^{t'}_a)_{a\in [0,t]}$ is a geodesic, it suffices to show that $(\ell''^t_a)_{a\in [0,t]}$ is a candidate in the Perron envelope defining the former geodesic.
    In other words, in verifying \eqref{eq:ellppdombyellpptemp1}, we may assume that either $s=0$ or $s=t$. The case $s=0$ is of course trivial. So it remains to prove the following:
    \[
    \ell''^{t'}_t\geq \ell_t\lor \ell'_t.
    \]
    By symmetry, it suffices to prove
    \[
    \ell''^{t'}_t\geq \ell_t.
    \]
    But since $(\ell_a)_{a\in [0,t']}$ is a candidate in the Perron envelope defining $\ell''^{t'}$, this inequality follows.

    \textbf{Step~2}.
    Next, observe that for a fixed $s\geq 0$, we have
    \[
    \sup_X \ell''^t_s\leq \frac{s}{t}\sup_X\ell''^t_t+ \frac{t-s}{t}\sup_X \ell''^t_0=\frac{s}{t}\left( \sup_X\ell_t \right)\lor \left( \sup_X\ell'_t \right)
    \]
    for all $t>s$. The right-hand side is bounded from above by a constant independent of $t\geq s$ by \cref{prop:raysuplinear}.
    Let
    \begin{equation}\label{eq:llorlp}
    (\ell\lor \ell')_s\coloneqq \uscsup[t>s] \ell''^t_{s}.
    \end{equation}
    Then \cref{prop:supsgeod} guarantees that $\ell\lor \ell'\in \mathcal{R}^{\infty}(X,\theta;\phi)$.

    \textbf{Step~3}.
    We need to show that $\ell\lor \ell'$ defined in this way is indeed the minimal ray lying above $\ell$ and $\ell'$.

    First, by Step~1, we have
    \[
    \ell''^t_s\geq \ell''^s_s\geq \ell_s
    \]
    for any $t\geq s\geq 0$. Therefore,
    \[
    (\ell\lor \ell')_s\geq \ell_s
    \]
    for all $s\geq 0$. In other words, $\ell\lor \ell'\geq \ell$.
    Similarly, $\ell\lor \ell'\geq \ell'$.

    Next, let $L\in \mathcal{R}^{\infty}(X,\theta;\phi)$ be a ray lying above both $\ell$ and $\ell'$. Then we have
    \[
    L_t\geq \ell_t\lor \ell'_t
    \]
    for all $t\geq 0$. In particular,
    \[
    L_s\geq \ell''^t_s
    \]
    for all $t\geq s\geq 0$. It follows that
    \[
    L_s\geq (\ell\lor \ell')_s
    \]
    for all $s\geq 0$.

    \textbf{Step~4}. It remains to argue \eqref{eq:Elor}. Since the energy is linear along bounded rays,
    \[
    \mathbf{E}^{\phi}(\ell\lor \ell')=E_{\theta}^{\phi}\left((\ell\lor \ell')_1\right).
    \]
    Moreover, $\ell''^t_1$ increases to $(\ell\lor \ell')_1$ as $t\to\infty$. By \cref{prop:Econt_mono},
    \[
    E_{\theta}^{\phi}\left((\ell\lor \ell')_1\right)=\lim_{t\to \infty} E_{\theta}^{\phi}\left(\ell''^t_1\right).
    \]
    Finally, by \cref{prop:energylinear} applied to $(\ell''^t_s)_{s\in [0,t]}$,
    \[
    E_{\theta}^{\phi}\left(\ell''^t_1\right)=\frac{1}{t}E_{\theta}^{\phi}\left(\ell_t\lor \ell'_t\right).
    \]
    Combining these identities gives \eqref{eq:Elor}.
\end{proof}

\begin{lemma}\label{lma:d1rayineq}
For any $\ell,\ell'\in \mathcal{R}^{\infty}(X,\theta;\phi)$, we have
    \begin{equation}\label{eq:d1maxineq}
    d_1(\ell,\ell')\leq d_1(\ell,\ell\lor\ell')+d_1(\ell',\ell\lor\ell')\leq C_nd_1(\ell,\ell'),
    \end{equation}
    where $C_n=3(n+1)2^{n+2}$.
\end{lemma}
\begin{proof}
    The first inequality is trivial. As for the second, we estimate
    \[
    \begin{aligned}
    d_1(\ell,\ell\lor\ell')=& \mathbf{E}^{\phi}(\ell\lor\ell')-\mathbf{E}^{\phi}(\ell)&\textup{by \cref{prop:d1geod_diff_E}}\\
    =&\lim_{t\to\infty} \frac{1}{t}E_{\theta}^{\phi}\left(\ell_t\lor \ell'_t\right)-\mathbf{E}^{\phi}(\ell)&\textup{by \eqref{eq:Elor}}\\
    =& \lim_{t\to\infty} \frac{1}{t} d_1\left(\ell_t\lor \ell'_t,\ell_t\right)&\textup{by \eqref{eq:d1asEdiff}}.
    \end{aligned}
    \]
    By symmetry, we find
    \[
    d_1(\ell,\ell\lor\ell')+d_1(\ell',\ell\lor\ell')\leq \lim_{t\to\infty} \frac{1}{t} \Bigl(d_1\left(\ell_t\lor \ell'_t,\ell_t\right)+d_1\left(\ell_t\lor \ell'_t,\ell'_t\right)\Bigr).
    \]
    By \cref{thm:I_d1_comp} and \cref{lma:I_equal_Imaxsum}, for each $t>0$,
    \[
    d_1\left(\ell_t\lor \ell'_t,\ell_t\right)+d_1\left(\ell_t\lor \ell'_t,\ell'_t\right)\leq 3(n+1)2^{n+2} d_1(\ell_t,\ell'_t).
    \]
    Now \eqref{eq:d1maxineq} follows.
\end{proof}

\begin{example}\label{ex:rayasspsh}
    Let $\varphi\in \PSH(X,\theta)$. Assume that $\varphi\leq 0$.
    For each $C>0$, let $(\ell^{\varphi,C}_t)_{t\in [0,C]}$ be the geodesic from $V_{\theta}$ to $(V_{\theta}-C)\lor \varphi$. For each $t\geq 0$, there is $\ell^{\varphi}_t\in \mathcal{E}^{\infty}(X,\theta)$ such that
    \begin{equation}\label{eq:ellvarphiraydef}
        \ell^{\varphi,C}_t\xrightarrow{d_1}\ell^{\varphi}_t
    \end{equation}
    as $C\to \infty$.
    Then $\ell^{\varphi}\in \mathcal{R}^{\infty}(X,\theta)$ and
\begin{equation}\label{eq:Elphi}
\mathbf{E}(\ell^{\varphi})=\frac{1}{n+1} \sum_{j=0}^n \left(\int_X \theta_{\varphi}^j\wedge \theta_{V_{\theta}}^{n-j}-\int_X \theta_{V_{\theta}}^n\right).
\end{equation}
\end{example}
Comparing the above approximants with cutoff parameters shifted by a bounded amount shows that $\ell^{\varphi+C}=\ell^{\varphi}$ for any $C\in \mathbb{R}$. Therefore, for a general $\varphi\in \PSH(X,\theta)$, we may simply define
\[
\ell^{\varphi}\coloneqq \ell^{\varphi-\sup_X\varphi}.
\]
Then the conclusions of this example continue to hold.
\begin{proof}
    We first show that for each fixed $t\geq 0$, $\ell_t^{\varphi,C}$ is increasing in $C\geq t$.

    To see this, choose $t\leq C_1<C_2$. We need to show that
    \[
    \ell^{\varphi,C_1}_t\leq \ell^{\varphi,C_2}_t.
    \]
    Since both sides are geodesics for $t\in [0,C_1]$, it suffices to show that
    \begin{equation}\label{eq:VthetaminusC1temp1}
    (V_{\theta}-C_1)\lor \varphi \leq \ell_{C_1}^{\varphi,C_2}.
    \end{equation}
    Now $((V_{\theta}-t)\lor \varphi)_{t\in [0,C_2]}$ is a subgeodesic from $V_{\theta}$ to $(V_{\theta}-C_2)\lor \varphi$ by \cref{prop:maxsubgeod}.\footnote{Here we need $\varphi\leq 0$.} It is therefore a candidate in the Perron envelope defining $\ell^{\varphi,C_2}$, and evaluating at $t=C_1$ gives exactly \eqref{eq:VthetaminusC1temp1}.

    From \cref{prop:convexsubgeod}, we know that for any $C>t>0$, we have
    \[
    \ell^{\varphi,C}_t\leq \frac{t}{C}\Bigl((V_{\theta}-C)\lor \varphi\Bigr)+\frac{C-t}{C}\cdot V_{\theta}\leq 0,
    \]
    while the linear subgeodesic $V_\theta-t$ is a candidate in the same Perron envelope, so $V_\theta-t\leq \ell^{\varphi,C}_t$. Therefore, by \cref{prop:pshfunction_closedseq},
    \begin{equation}\label{eq:ellphitexp}
    \ell_t^{\varphi}\coloneqq \uscsup[C>t]\ell^{\varphi,C}_t\in \mathcal{E}^{\infty}(X,\theta)
    \end{equation}
    for all $t\geq 0$. Thanks to \cref{prop:Econt_mono}, we have
    \[
    \ell^{\varphi,C}_t\xrightarrow{d_1} \ell_t^{\varphi}
    \]
    as $C\to\infty$ for all $t\geq 0$. Applying \cref{prop:supsgeod} on each interval $[0,T]$ to the increasing family $(\ell^{\varphi,C}_t)_{t\in [0,T]}$, $C\geq T$, shows that $\ell^{\varphi}\in \mathcal{R}^{\infty}(X,\theta)$.

    It remains to compute the energy of $\ell^{\varphi}$. We first fix $C\geq t>0$ and compute using \cref{prop:energylinear}:
    \[
        E_{\theta}\left(\ell^{\varphi,C}_t\right)=\frac{t}{C}E_{\theta}\Bigl((V_{\theta}-C)\lor \varphi\Bigr).
    \]
    Letting $C\to\infty$ and applying \cref{prop:Econt_mono}, we find that
    \[
        E_{\theta}(\ell^{\varphi}_t)=\lim_{C\to\infty}\frac{t}{C}E_{\theta}\Bigl((V_{\theta}-C)\lor \varphi\Bigr)
    \]
    for any $t\geq 0$.
    It follows that
    \[
        \mathbf{E}(\ell^{\varphi})=\lim_{C\to\infty}\frac{1}{C}E_{\theta}\Bigl((V_{\theta}-C)\lor \varphi\Bigr).
    \]
    Using the definition of $E_{\theta}$, in order to obtain \eqref{eq:Elphi}, it suffices to show that for each $j=0,\ldots,n$, we have
    \begin{equation}\label{eq:limCintXtemp1}
    \lim_{C\to\infty}\int_X \frac{(V_{\theta}-C)\lor\varphi -V_{\theta}}{C} \,\theta_{(V_{\theta}-C)\lor \varphi}^j\wedge \theta_{V_{\theta}}^{n-j}=\int_X \theta_{\varphi}^{j}\wedge \theta_{V_{\theta}}^{n-j}-\int_X \theta_{V_{\theta}}^n.
    \end{equation}

    For this purpose, for each $C>0$, we decompose $X$ as $\{\varphi> V_{\theta}-C\}$ and $\{\varphi\leq V_{\theta}-C\}$.
    We have
    \[
    \begin{aligned}
    &\int_{\{\varphi> V_{\theta}-C\}} \frac{(V_{\theta}-C)\lor\varphi -V_{\theta}}{C} \,\theta_{(V_{\theta}-C)\lor \varphi}^j\wedge \theta_{V_{\theta}}^{n-j}\\
    =& \int_{\{\varphi> V_{\theta}-C\}} \frac{\varphi -V_{\theta}}{C} \,\theta_{\varphi}^j\wedge \theta_{V_{\theta}}^{n-j}.
    \end{aligned}
    \]
    On the other hand,
    \[
    \begin{aligned}
    &\int_{\{\varphi\leq V_{\theta}-C\}} \frac{(V_{\theta}-C)\lor\varphi -V_{\theta}}{C} \,\theta_{(V_{\theta}-C)\lor \varphi}^j\wedge \theta_{V_{\theta}}^{n-j}\\
    =& -\int_{\{\varphi\leq V_{\theta}-C\}} \theta_{(V_{\theta}-C)\lor\varphi}^j\wedge \theta_{V_{\theta}}^{n-j}\\
    =& -\int_X \theta_{V_{\theta}}^n+\int_{\{\varphi> V_{\theta}-C\}} \theta_{\varphi}^j\wedge \theta_{V_{\theta}}^{n-j}.
    \end{aligned}
    \]
    Observe that for $C>0$, the function $\mathds{1}_{\{\varphi>V_{\theta}-C\}}C^{-1}(\varphi-V_{\theta})$ is defined quasi-everywhere. Since $\varphi\leq V_\theta$, it takes values in $[-1,0]$ outside a pluripolar set.
    When $C\to \infty$, these functions converge to $0$ quasi-everywhere, so the first integral above converges to $0$ by dominated convergence. Similarly, the indicators $\mathds{1}_{\{\varphi>V_\theta-C\}}$ increase to $1$ outside a pluripolar set, hence
    \[
        \lim_{C\to\infty}\int_{\{\varphi>V_{\theta}-C\}} \theta_{\varphi}^j\wedge \theta_{V_{\theta}}^{n-j} = \int_X \theta_{\varphi}^j\wedge \theta_{V_{\theta}}^{n-j}.
    \]
    Therefore, \eqref{eq:limCintXtemp1} follows.
\end{proof}

\motto{There are two principal ways to formulate mathematical assertions (problems, conjectures, theorems, \dots ): Russian and French. The Russian way is to
choose the most simple and specific case (so that nobody could simplify the formulation preserving the main point). The French way is to generalize the statement
as far as nobody could generalize it further.

--- Vladimir Arnold\footnote{Vladimir Igorevich Arnold (1937--2010), who became a professor at l'Université Paris IX after the dissolution of USSR, was always sick of France (so am I!). In the public lecture entitled ``Sur l'éducation mathématique'' in 1997, he invented the famous joke ``Combien font $2+3$?'' to question the French education system.}}
\chapter{Toric pluripotential theory on ample line bundles}
\label{chap:toric_ample}

This chapter develops the toric form of pluripotential theory in the case of an ample line bundle. Throughout the chapter, $P\subseteq M_{\mathbb R}$ is a full-dimensional smooth lattice polytope, $\Sigma$ is the corresponding fan, $X=X_\Sigma$ is the associated smooth projective toric variety, and $L=\mathcal O_X(D_P)$ is the ample toric line bundle defined by $P$.

The basic principle is that toric plurisubharmonic metrics on $L$ can be described by convex functions on $N_{\mathbb R}$. Under this dictionary, singularity type, rooftop envelopes, Monge--Amp\`ere mass and energy all admit convex-analytic interpretations. The Newton body of a toric potential is the main invariant in this chapter: it records the effective convex domain seen by multiplier ideals and non-pluripolar mass.

The central result is Yao's theorem, stated as \cref{thm:Yao}. It characterizes the torus characters belonging to multiplier-ideal spaces in terms of the Newton body. As a consequence, the asymptotic dimension of
\[
    \mathrm{H}^0\left(X,L^k\otimes\mathcal I(k\varphi)\right)
\]
is governed by the Euclidean volume of the Newton body. This gives a toric model for the general volume theorem proved later in \cref{chap:Igood} and shows, in the ample toric setting, why multiplier ideals and non-pluripolar mass should be related.

The chapter also proves the real Monge--Amp\`ere formula for toric metrics, the toric energy formula, and the expression of the $d_1$-distance in convex terms. These results will be used again in \cref{chap:toric_big}, where the ample assumption is removed and the Newton body is compared with partial Okounkov bodies.

We assume familiarity with basic toric geometry at the level of \cite{CLS11}. For a construction of the ample toric setting, see also \cite[Sections~18.7--18.9]{FPT}. Some convex-analytic facts used below are recalled in \cref{chap:convex}. The chapter can be skipped on a first reading if the reader is only interested in the non-toric theory, but it provides useful examples and a guiding model for several later constructions.

Finally, throughout the book convex functions are understood in the usual sense: the function $\mathbb R\ni x\mapsto x^2$ is convex. This convention is opposite to the convention often used in toric geometry, including parts of \cite{CLS11}; care is therefore needed when comparing formulas. See also \cite[Definition~18.1.1]{FPT} for our convexity convention.

\section{Toric setup}\label{sec:toricsetup}
Let $T$ be a complex torus of dimension $n$\footnote{Namely, an algebraic group defined over $\mathbb{C}$, which is isomorphic to $\mathbb{G}_{\mathrm{m}}^n$.} and let $T_c\subseteq T(\mathbb{C})$ denote the corresponding compact torus.

Write $M$ for the character lattice of $T$, which is a free Abelian group of rank $n$. Similarly, let $N$ be the cocharacter lattice of $T$, which is the dual lattice of $M$. Given $m\in M$, the corresponding character of $T$ is denoted by $\chi^m$.
Write $M_{\mathbb{R}}=M\otimes_{\mathbb{Z}} \mathbb{R}$ and $N_{\mathbb{R}}=N\otimes_{\mathbb{Z}}\mathbb{R}$. The pairing between $M_{\mathbb{R}}$ and $N_{\mathbb{R}}$ is denoted by $\langle \bullet,\bullet \rangle$.

Let $P\subseteq M_{\mathbb{R}}$ be a full-dimensional \emph{smooth}\footnote{Recall that \emph{smooth} means that for every vertex $v\in P$, if $w_E$ denotes the first lattice point after $v$ as one traverses each edge $E$ of $P$ containing $v$, then $\{w_E-v\}_E$ forms a basis of $M$. See \cite[Definition~2.4.2]{CLS11} and \cite[Definition~18.7.1]{FPT}. We also say $P$ is a \emph{Delzant polytope}\index{polytope!Delzant polytope} in this case.} lattice polytope\footnote{A \emph{lattice polytope}\index{polytope!lattice polytope} in $M_{\mathbb{R}}$ is the convex hull of finitely many points in $M$.}.

Given any (closed) facet $F$ of $P$, let $u_F\in N$ denote the unique ray generator (the first non-zero integral element) of the inward normal ray of $F$. Then $P$ can be represented as
\begin{equation}\label{eq:Pfacetrep}
P=\left\{m\in M_{\mathbb{R}}: \langle m,u_F \rangle\geq -a_F\textup{ for all facets }F\textup{ of } P\right\}
\end{equation}
for some uniquely determined integers $a_F$. The presentation is called the \emph{facet presentation}\index{polytope!facet presentation} of $P$.

Given any (closed) face $Q$ of $P$, we let $\sigma_Q\subseteq N_{\mathbb{R}}$ be the closed convex cone generated by the $u_F$'s, where $F$ runs over all facets of $P$ containing $Q$. When $Q=P$, $\sigma_P$ is understood as $\{0\}$.

Let $\Sigma$ be the \emph{(inner) normal fan}\index{fan!normal fan of a polytope} of $P$. Namely,
\[
\Sigma=\left\{\sigma_Q:Q\textup{ is a face of }P\right\}.
\]
The notation $\Sigma(1)$ denotes the set of rays in $\Sigma$. Note that $\Sigma(1)$ is in bijective correspondence with the set of facets of $P$. In fact, given any facet $F$ of $P$, the cone $\sigma_F$ is just the ray generated by $u_F$, namely, the inward normal ray of $F$.

For any $\rho\in \Sigma(1)$, let $u_{\rho}\in N$ denote the ray generator of $\rho$, namely the first non-zero element in $N\cap \rho$. If $\rho=\sigma_F$ for some facet $F$ of $P$, then $u_{\rho}=u_F$.

Now the facet presentation \eqref{eq:Pfacetrep} can be equivalently rewritten as
\[
P=\left\{m\in M_{\mathbb{R}} : \langle m,u_{\rho} \rangle\geq -a_{\rho}\textup{ for all }\rho\in \Sigma(1) \right\}.
\]
Let $\Supp_P\colon N_{\mathbb{R}} \rightarrow \mathbb{R}$ denote the \emph{support function} of $P$. Recall that the support function (\cref{ex:suppfun}) of $P$ is defined as
\[
\Supp_P(n)=\max\left\{ \langle m,n \rangle: m\in P \right\}.
\]
Note that our support function differs from \cite[Proposition~4.2.14]{CLS11}, where instead of a maximum, they took the minimum. The maximum convention is also used in \cite[Definition~18.2.1]{FPT}.

Recall that the \emph{characteristic function} $\chi_P\colon M_{\mathbb{R}}\rightarrow \{0,\infty\}$ of $P$ is defined as in \cref{ex:charconvex}:
\[
\chi_P(m)\coloneqq
    \begin{cases}
        0,&\quad m\in P;\\
        \infty,&\quad m\not\in P.
    \end{cases}
\]

Let $X=X_{\Sigma}$ be the smooth projective toric variety corresponding to $\Sigma$. See \cite[Theorem~3.1.5]{CLS11} for the construction of $X$ and \cite[Theorem~3.1.19]{CLS11} for the smoothness of $X$. See also \cite[Theorem~18.7.3]{FPT} for both assertions. There is a canonical embedding $T\subseteq X$ as a dense Zariski open subset.

Let $D$ be the Cartier divisor on $X$ defined by $P$:
\[
D=\sum_{\rho\in \Sigma(1)} a_{\rho}D_{\rho},
\]
where $D_{\rho}$ is the toric prime divisor defined by $\rho$ under the orbit-cone correspondence \cite[Theorem~3.2.6]{CLS11}. See also \cite[Proposition~18.8.2]{FPT}.

Let $L$ be the toric line bundle induced by $P$, namely $L=\mathcal{O}_X(D)$. Since $P$ is a smooth (Delzant) lattice polytope, $L$ is very ample by \cite[Proposition~2.4.4, Proposition~6.1.10]{CLS11} and \cite[Theorems~18.7.3 and~18.8.5]{FPT}; in particular, $L$ is ample.

We will choose the base $\mathrm{e}$ for the logarithm map
\begin{equation}\label{eq:log}
\mathbb{C}^*\rightarrow \mathbb{R},\quad z\mapsto \log|z|^2.\footnotemark
\end{equation}
\footnotetext{Be careful when you compare with other references, some people prefer $\log |z|$, $-\log |z|$ or $-\log |z|^2$ instead.}
This choice is fixed throughout the book. We have a canonical identification
\[
T(\mathbb{C})=\mathrm{Hom}(M,\mathbb{C}^*)\xlongrightarrow{\sim} N\otimes_{\mathbb{Z}}\mathbb{C}^*.
\]
The logarithm map then induces a tropicalization map after tensoring with $N$:
\begin{equation}
\Trop\colon T(\mathbb{C})\rightarrow N_{\mathbb{R}},\quad p\mapsto \left( M \ni m\mapsto \log \left|\chi^m(p)\right|^2 \right).
\end{equation}
\index{$\Trop$}

Before proceeding, it is always helpful to understand everything in our favorite example.
\begin{example}\label{ex:P1toric}
    We take $n=1$ and $P=[0,1]\subseteq M_{\mathbb{R}}=\mathbb{R}$. In this case, the facet representation \eqref{eq:Pfacetrep} becomes
    \[
    P=\bigl\{m\in \mathbb{R}: \langle m,1 \rangle \geq 0, \langle m,-1 \rangle \geq -1\bigr\},
    \]
    with $u_{\{0\}}=1$, $u_{\{1\}}=-1$, $a_{\{0\}}=0$ and $a_{\{1\}}=1$. The normal fan $\Sigma$ is
    \[
    \Sigma=\bigl\{(-\infty,0], \{0\}, [0,\infty) \bigr\}.
    \]
    See \cref{fig:fanP1}.
        \begin{figure}[ht]
        \centering
        \begin{tikzpicture}[scale=1.0, line width=0.6pt, every node/.style={font=\small}]
          \draw[red!75!black, line width=1.4pt] (-4.0, 0) -- (-0.05, 0);
          \draw[blue!70!black, line width=1.4pt] (0.05, 0) -- (4.0, 0);
          \draw[thick, ->] (4.0, 0) -- (4.4, 0) node[right, black] {$N_{\mathbb{R}}$};
          \fill (0,0) circle (1.6pt);
          \node[above=2pt] at (-2.0, 0) {$(-\infty,0]$};
          \node[above=4pt] at (0, 0) {$\{0\}$};
          \node[above=2pt] at (2.0, 0) {$[0,\infty)$};
        \end{tikzpicture}
        \caption{The fan $\Sigma$ of $\mathbb{P}^1$ in $N_{\mathbb{R}}=\mathbb{R}$, consisting of the two rays $(-\infty,0]$, $[0,\infty)$ and the cone $\{0\}$.}
        \label{fig:fanP1}
    \end{figure}
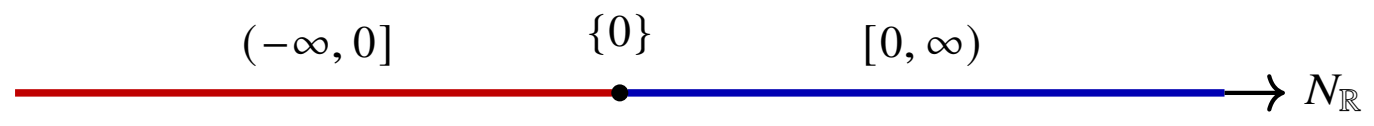

    The corresponding toric variety is just $X=\mathbb{P}^1$. Under the orbit-cone correspondence, we have
    \[
    D_{[0,\infty)}=[0],\quad D_{(-\infty,0]}=[\infty].
    \]
    The associated divisor $D=[\infty]$ and therefore,
    \[
    L=\mathcal{O}_X(D)=\mathcal{O}_{\mathbb{P}^1}(1).
    \]
\end{example}

\section{Toric plurisubharmonic functions}\label{sec:tpf}

We continue to use the notations of \cref{sec:toricsetup}.

\begin{lemma}\label{lma:convextopsh}
    Let $F\colon N_{\mathbb{R}}\rightarrow [-\infty,\infty]$ be a function. Then the following are equivalent:
    \begin{enumerate}
        \item\label{itm:lma:convextopsh:1} $F$ is convex and takes values in $\mathbb{R}$, and
        \item\label{itm:lma:convextopsh:2} $\Trop^*F$ is plurisubharmonic on $T(\mathbb{C})$.
    \end{enumerate}
\end{lemma}
\begin{proof}
We may choose an identification $N\cong \mathbb{Z}^n$ so that we have an identification $T(\mathbb{C})\cong \mathbb{C}^{*n}$.
Then $\Trop$ is identified with the map
\[
\Trop\colon \mathbb{C}^{*n}\rightarrow \mathbb{R}^n,\quad (z_1,\ldots,z_n)\mapsto \left(\log|z_1|^2,\ldots,\log|z_n|^2\right).
\]

\ref{itm:lma:convextopsh:1} $\implies$ \ref{itm:lma:convextopsh:2}. By the standard convolution regularization of finite convex functions, let $F_k\in C^{\infty}(\mathbb{R}^n)\cap \Conv(\mathbb{R}^n)$ be a decreasing sequence with limit $F$. It follows from a straightforward computation that
    \begin{equation}\label{eq:ddctrop}
    \begin{split}
    \ddc \Trop^*F_k(z_1,\ldots,z_n)
    =&\frac{\mathrm{i}}{2\pi} \sum_{i,j=1}^n \partial_{ij}F_k\left(\log|z_1|^2,\ldots,\log|z_n|^2\right)\\
    &\times z_i^{-1}\overline{z_j}^{-1}\mathrm{d}z_i\wedge \mathrm{d}\overline{z_j}.
    \end{split}
    \end{equation}
    So $\Trop^* F_k$ is plurisubharmonic. It follows from \cref{prop:pshfunction_closedseq} that $\Trop^* F$ is plurisubharmonic.

\ref{itm:lma:convextopsh:2} $\implies$ \ref{itm:lma:convextopsh:1}. It follows from \cref{lma:invariantpshfunfinite} that $F$ is finite. Moreover, by radial regularization, we may find a decreasing sequence $\varphi_k$ of $(S^1)^n$-invariant smooth psh functions on $\mathbb{C}^{*n}$ with limit $\Trop^*F$. Write $\varphi_k=\Trop^* F_k$ for some function $F_k\colon \mathbb{R}^n\rightarrow \mathbb{R}$; it follows from \eqref{eq:ddctrop} that $F_k$ is convex for all $k$. Therefore, $F$ is convex by \cref{lma:convdecnet}.
\end{proof}

Next we define a canonical K\"ahler form in $c_1(L)$.

Let $G_0\colon M_{\mathbb{R}}\rightarrow (-\infty,\infty]$ be defined as
\begin{equation}\label{eq:G0def}
G_0(m)\coloneqq
\begin{cases}
\sum\limits_{\rho \in \Sigma(1)} \left(\langle m, u_{\rho}\rangle +a_{\rho}\right)\log  \left(\langle m, u_{\rho}\rangle +a_{\rho}\right)\footnotemark,&\textup{if }m\in P,\\
\infty,& \textup{otherwise}.
\end{cases}
\end{equation}
\footnotetext{We understand that $0\log 0=0$ in this expression.}

This is a closed proper convex function and $G_0\sim \chi_P$, where $\sim$ is the relation defined in \cref{def:partialorderconv}.

Let
\begin{equation}\label{eq:F0def}
F_0=G_0^* \in \mathcal{E}^{\infty}(N_{\mathbb{R}},P).
\end{equation}
Here $G_0^*$ is the Legendre transform of $G_0$, as recalled in \cref{def:Legendregeneral}. The set $\mathcal{E}^{\infty}(N_{\mathbb{R}},P)$ is defined in \cref{def:convexPfunctions}.

By Guillemin's theorem \cite{Gui94, CDG03} (see also \cite[Theorem~18.9.1]{FPT}), $\ddc\Trop^* F_0$ can be extended to a unique K\"ahler form $\omega$ in $c_1(L)$. The K\"ahler form $\omega$ is clearly $T_c$-invariant.

For each $\rho\in \Sigma(1)$, we write
\[
    r_{\rho}(m)=\log \left(\langle m,u_{\rho}\rangle +a_{\rho}\right)+1,\quad m\in P.
\]
It follows from \eqref{eq:G0def} that
\begin{equation}\label{eq:nablaG0}
    \nabla G_0(m)=\sum_{\rho\in \Sigma(1)} r_{\rho}(m)u_{\rho},\quad m\in \Int P.
\end{equation}

\begin{example}\label{ex:P1toric1}
Let us move on with our favorite example \cref{ex:P1toric}. We continue to use the same notations. In this case,
\[
G_0(m)=
\begin{cases}
    m\log m+(1-m)\log (1-m),& \textup{if } m\in [0,1],\\
    \infty,&\textup{otherwise}.
\end{cases}
\]
The Legendre transform is given\footnote{While reading an advanced mathematical textbook/paper, I usually tend to trust the authors for their elementary computations. A few years ago, I was asked to present the result of a landmark paper written by two respected mathematicians on a conference. After spending a few days on the elementary integrals, I found out that all non-trivial constants in that paper were wrong. So I ask the readers to really verify this expression, if it is not obvious to you.} by
\[
F_0(n)=\log\left(1+\mathrm{e}^n\right).
\]
Composing with the tropicalization map, we find that
\[
\omega|_{\mathbb{C}^*}(z)=\ddc\log\left(1+|z|^2\right).
\]
This is exactly the Fubini--Study metric as we have seen in \cref{ex:LPP1}.

We can now explain one subtlety: in our expression \eqref{eq:G0def}, there is no factor $1/2$ before the sum; the absence of the factor $1/2$ matches the choice $\Trop(z)=\log|z|^2$ rather than $\log|z|$ in our tropicalization map \eqref{eq:log}.
\end{example}

Let $\PSH_{\tor}(X,\omega)$\index{$\PSH_{\tor}(X,\omega)$} denote the set of $T_c$-invariant $\omega$-psh functions.
\begin{theorem}\label{thm:toricpsh}
    There are canonical bijections between the following three sets:
    \begin{enumerate}
        \item\label{itm:thm:toricpsh:1} The set of $\varphi\in \PSH_{\tor}(X,\omega)$,
        \item\label{itm:thm:toricpsh:2} the set $\mathcal{P}(N_{\mathbb{R}},P)$ in \cref{def:convexPfunctions}, namely, the set of convex functions $F\colon N_{\mathbb{R}}\rightarrow \mathbb{R}$ satisfying $F\preceq \Supp_P$, and
        \item\label{itm:thm:toricpsh:3} the set of closed proper convex functions $G\in \Conv(M_{\mathbb{R}})$ satisfying
        \[
        G|_{M_{\mathbb{R}}\setminus P}\equiv \infty.
        \]
    \end{enumerate}
\end{theorem}
For the notions of closedness and properness, we refer to \cref{def:convexfun} and \cref{def:convexclosure}.

\begin{proof}
The bijection between \ref{itm:thm:toricpsh:2} and \ref{itm:thm:toricpsh:3} is the classical Legendre duality. Given $F$ as in \ref{itm:thm:toricpsh:2}, we construct $G=F^*$ and \emph{vice versa}, see \cref{prop:Fsuppchar}.

The map from \ref{itm:thm:toricpsh:1} to \ref{itm:thm:toricpsh:2} is given as follows: Given $\varphi\in \PSH_{\tor}(X,\omega)$, since $\varphi$ is $T_c$-invariant, we can find $f\colon N_{\mathbb{R}}\rightarrow [-\infty,\infty)$ such that
\begin{equation}\label{ex:varphitropf}
\varphi|_{T(\mathbb{C})}=\Trop^*f.
\end{equation}
We then define $F=f+F_0$. Then $\Trop^* F\in \PSH(T(\mathbb{C}))$.
By \cref{lma:convextopsh}, $F(n)$ is finite for any $n\in N_{\mathbb{R}}$ and $F$ is convex. Moreover, $F\preceq \Supp_P$ since $\varphi$ is bounded from above on the compact space $X$ and $F_0\preceq \Supp_P$.

Conversely, given a map $F\in \mathcal{P}(N_{\mathbb{R}},P)$, then
\[
\Trop^*(F-F_0)\in \PSH\left(T(\mathbb{C}),\omega|_{T(\mathbb{C})}\right).
\]
Since $F\preceq \Supp_P$ and $F_0\sim \Supp_P$, the difference $F-F_0$ is bounded from above. Hence the corresponding local psh potentials are locally bounded from above near $X\setminus T(\mathbb{C})$. It follows from \cref{thm:GRexten} that this function can be extended uniquely to an $\omega$-psh function on $X$. The uniqueness of the extension guarantees its $T_c$-invariance.

The two maps are inverse to each other by construction.
\end{proof}

Given $\varphi\in \PSH_{\tor}(X,\omega)$, we will write $F_{\varphi}$ and $G_{\varphi}$ for the convex functions given by \cref{thm:toricpsh}. From the proof, we have the following relations:
\begin{equation}\label{eq:FvarphidefGdef}
\varphi|_{T(\mathbb{C})}=\Trop^*\left(F_{\varphi}-F_0\right),\quad G_{\varphi}=F_{\varphi}^*.
\end{equation}
In particular, $\varphi|_{T(\mathbb{C})}$ is locally bounded. We will frequently use this result without further explanation.

\begin{example}\label{ex:FGsimplevarphi}
    Let us take our favorite example \cref{ex:P1toric1} again. We will continue to use the same notations.

    Recall that in \cref{ex:logsingP1} and \cref{ex:loglogsing}, we constructed two $S^1$-invariant functions in $\PSH(X,\omega)$.

    We begin with the function $\varphi$ in \cref{ex:logsingP1}. Recall that
    \[
    \varphi(z)=\log \frac{|z|^2}{|z|^2+1}
    \]
    for $z\in \mathbb{C}$. The function $f\colon \mathbb{R}\rightarrow \mathbb{R}$ in \eqref{ex:varphitropf} is therefore
    \[
    f(n)=\log \frac{\mathrm{e}^{n}}{1+\mathrm{e}^{n}}.
    \]
    Therefore, $F_{\varphi}\colon \mathbb{R}\rightarrow \mathbb{R}$ is
    \[
    F_{\varphi}(n)=n.
    \]
    Correspondingly, $G_{\varphi}\colon \mathbb{R}\rightarrow \mathbb{R}$ is
    \[
    G_{\varphi}(m)=
    \begin{cases}
        0,&\textup{if } m=1;\\
        \infty,&\textup{otherwise}.
    \end{cases}
    \]

    Similarly, if $\psi$ denotes the function in \cref{ex:loglogsing}, then the function $f$ in \eqref{ex:varphitropf} is
    \[
    f(n)=
    \begin{cases}
        -\log \left(\mathrm{e}^n+1 \right)+(-\log (-n))\lor (n+2),&\textup{if } n<-\log 2;\\
        2+\log \frac{\mathrm{e}^{n}}{1+\mathrm{e}^{n}},&\textup{otherwise}.
    \end{cases}
    \]
    Therefore,
    \[
    F_{\psi}(n)=
    \begin{cases}
        (-\log (-n))\lor (n+2),&\textup{if } n<-\log 2;\\
        2+n,&\textup{otherwise}.
    \end{cases}
    \]
    The Legendre transform is tricky to compute. Let $\lambda$ be the large solution of $\log x=x-2$. So $\lambda\approx 3.146$. The smaller solution is around $0.159<\log 2\approx 0.693$. It might be helpful to have a look at the poorly drawn picture \cref{fig:logx}.

    \begin{figure}[ht]
\centering
\begin{tikzpicture}[xscale=0.95, yscale=0.7]
  \draw[step=1, gray!30, dashed, very thin] (0,-3) grid (10,5);
  \draw[gray!60] (0,-3) rectangle (10,5);
  \draw[thin] (0,0) -- (10,0);
  \foreach \tx in {0,2,4,6,8,10}
    \node[below] at (\tx,-3) {\scriptsize $\tx$};
  \foreach \ty in {-3,-2,-1,0,1,2,3,4,5}
    \node[left] at (0,\ty) {\scriptsize $\ty$};
  \node[below] at (5,-3.6) {$x$};
  \node at (-0.9, 1) {$y$};
  \draw[orange, thick, smooth] plot[domain=0.1:10, samples=200] (\x, {ln(\x)});
  \draw[cyan!70!blue, thick] (0,-2) -- (7,5);
  \draw[red, thick] plot[only marks, mark=x, mark size=3pt]
    coordinates {(3.15,1.15) (0.16,-1.84)};
  \node[red, anchor=south, font=\footnotesize, inner sep=1pt]
    at (3.15,1.25) {$(3.15,\,1.15)$};
  \node[red, anchor=east, font=\footnotesize, inner sep=2pt]
    at (0.1,-1.84) {$(0.16,\,-1.84)$};
  \begin{scope}[shift={(7.1, 4.0)}]
    \draw[gray!50, fill=white, rounded corners] (-0.1,-0.7) rectangle (2.7,0.7);
    \draw[orange, thick] (0.05, 0.45) -- (0.55, 0.45);
    \node[anchor=west, font=\scriptsize] at (0.6, 0.45) {$y = \log x$};
    \draw[cyan!70!blue, thick] (0.05, 0.0) -- (0.55, 0.0);
    \node[anchor=west, font=\scriptsize] at (0.6, 0.0) {$y = x-2$};
    \draw[red, thick] plot[only marks, mark=x, mark size=2pt]
      coordinates {(0.3,-0.45)};
    \node[anchor=west, font=\scriptsize] at (0.6, -0.45) {Intersections};
  \end{scope}
\end{tikzpicture}
\caption{The graphs of $\log x$ and $x-2$.}
\label{fig:logx}
\end{figure}
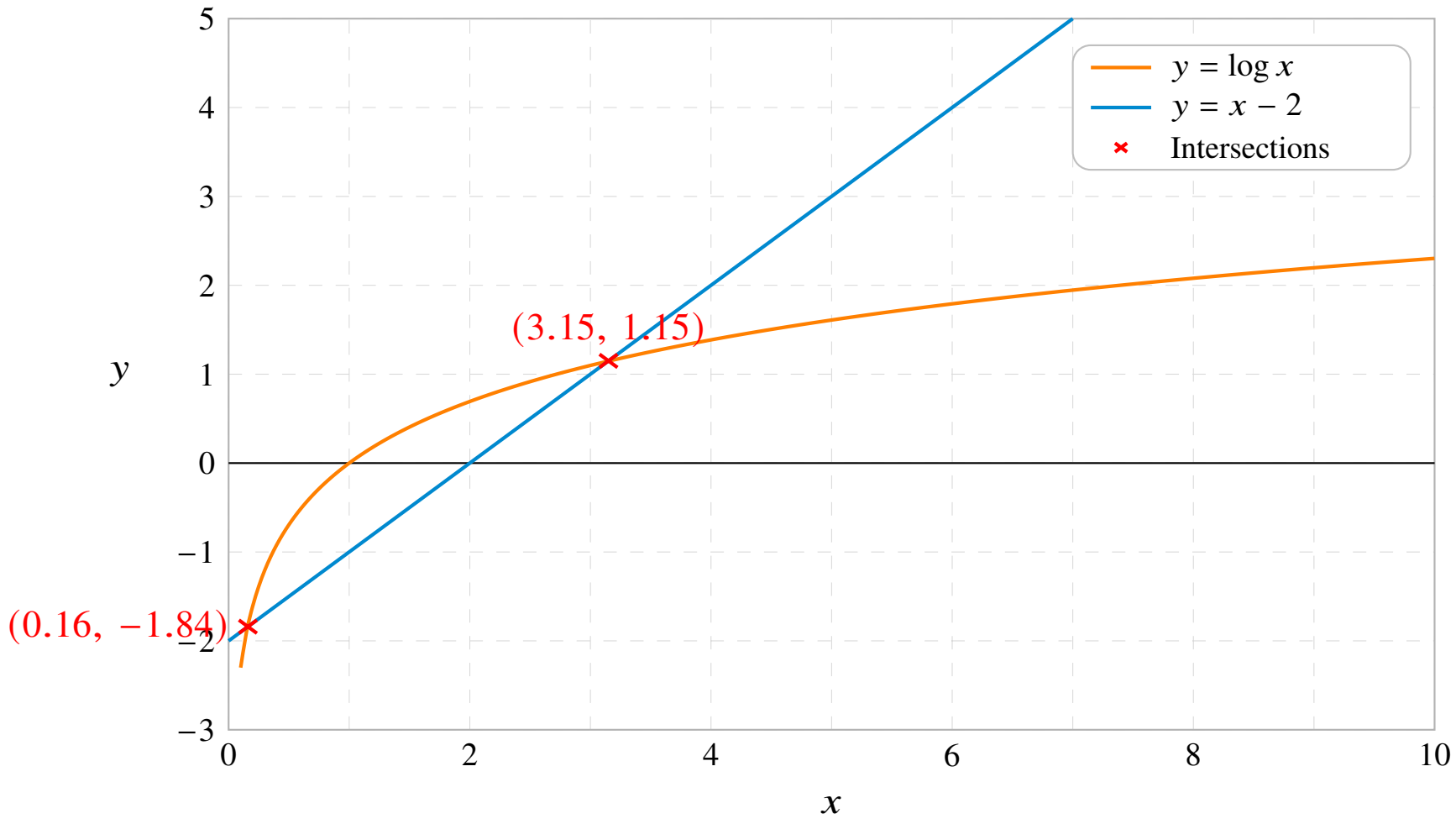

It is immediate that $G_{\psi}(m)=\infty$ unless $m\in [0,1]$. Let us assume that $m\in [0,1]$. Then
\[
\begin{aligned}
G_{\psi}(m)=&\sup_{n\in \mathbb{R}} \left( mn-F_{\psi}(n)\right)\\
=& \sup_{n<-\log 2} \left( mn-(-\log (-n))\lor (n+2)\right)  \lor \sup_{n\geq -\log 2} (mn-n-2)\\
=& \sup_{n>\log 2} \left( -mn+(\log n)\land (n-2)\right)  \lor \left((1-m)\log 2-2\right).
\end{aligned}
\]
Let us focus on the first part, which can be decomposed further into
\[
\begin{aligned}
& \sup_{n>\log 2} \left( -mn+(\log n)\land (n-2)\right)\\
=& \sup_{n\in \left(\log 2,\lambda\right]} (n-2-mn) \lor \sup_{n>\lambda} (\log n-mn)\\
=& \left( (1-m)\lambda-2 \right) \lor \sup_{n>\lambda} \left(\log n-mn\right).
\end{aligned}
\]
The latter part can be computed easily:
\[
\sup_{n>\lambda} \left(\log n-mn\right)=
\begin{cases}
    -\log m-1,&\textup{if } m\in \left(0,\lambda^{-1}\right];\\
    \log \lambda -m\lambda,& \textup{if }m\in \left(\lambda^{-1},1\right].
\end{cases}
\]
Putting everything together, we find
\[
G_{\psi}(m)=
\begin{cases}
    (-\log m-1)\lor \left( (1-m)\lambda-2 \right) ,&\textup{if } m\in \left(0,\lambda^{-1}\right];\\
    (\log \lambda -m\lambda)\lor \left( (1-m)\lambda-2 \right),& \textup{if }m\in \left(\lambda^{-1},1\right].
\end{cases}
\]
This can be further simplified, the final result is
\[
G_{\psi}(m)=
\begin{cases}
    -\log m-1 ,&\textup{if } m\in \left(0,\lambda^{-1}\right];\\
    (1-m)\lambda-2 ,& \textup{if }m\in \left(\lambda^{-1},1\right];\\
    \infty,&\textup{otherwise}.
\end{cases}
\]
The graph of $G_{\psi}$ on $(0,1]$ is sketched in \cref{fig:Gpsi}.
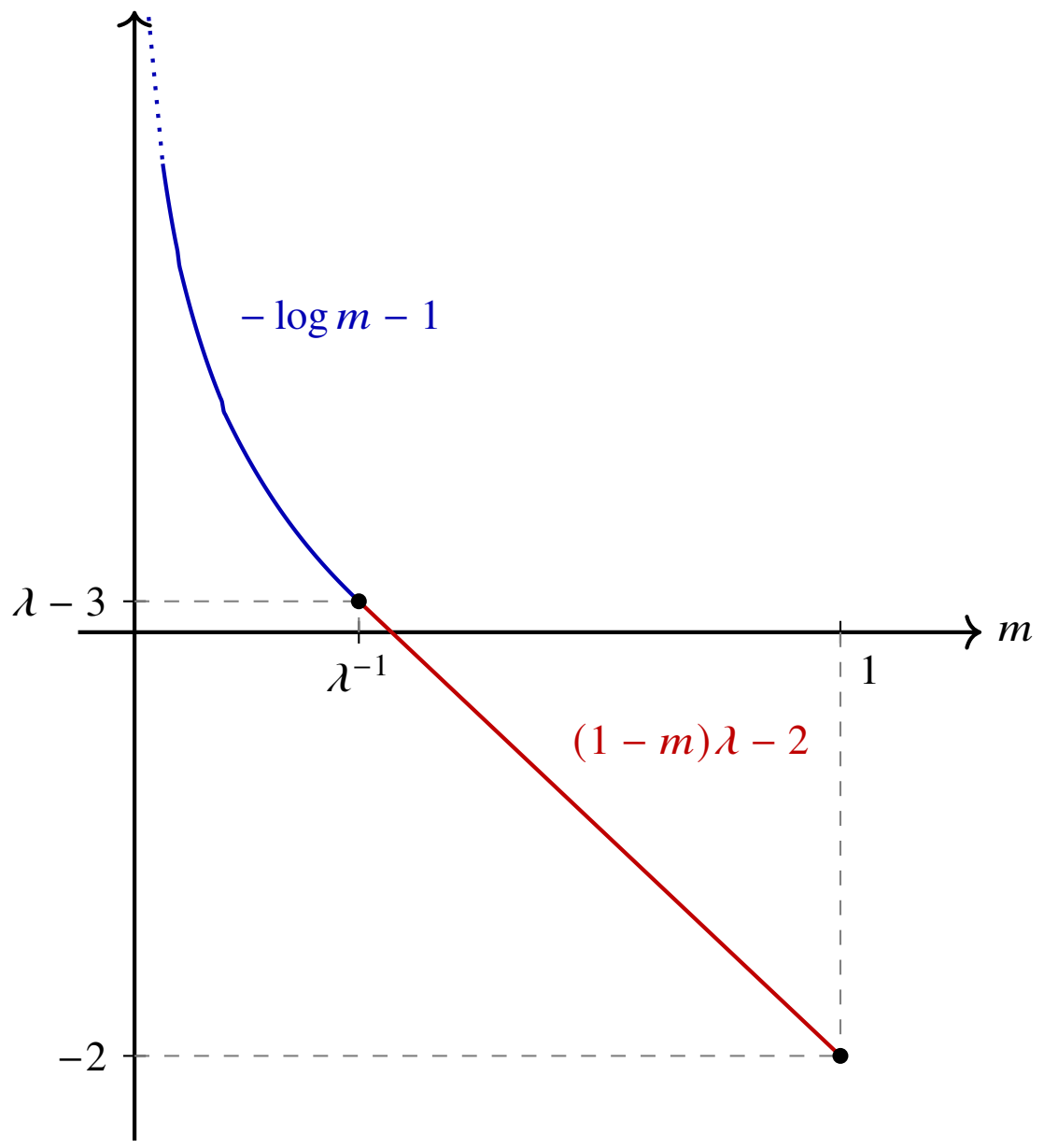
\begin{figure}[ht]
\centering
\begin{tikzpicture}[scale=1.0, line width=0.6pt, every node/.style={font=\small}]
  \draw[->, thick] (-0.4,0) -- (6.0,0) node[right] {$m$};
  \draw[->, thick] (0,-3.6) -- (0,4.4);
  \draw[thin] (1.589,0.08) -- (1.589,-0.08);
  \node[anchor=north, inner sep=2pt] at (1.589,-0.08) {$\lambda^{-1}$};
  \draw[thin] (5,0.08) -- (5,-0.08);
  \node[anchor=north west, inner sep=2pt] at (5.05,-0.08) {$1$};
  \draw[thin] (0.08,0.219) -- (-0.08,0.219);
  \node[anchor=east, inner sep=3pt] at (-0.08,0.219) {$\lambda-3$};
  \draw[thin] (0.08,-3) -- (-0.08,-3);
  \node[anchor=east, inner sep=3pt] at (-0.08,-3) {$-2$};
  \draw[gray, dashed, thin] (1.589,0)
    -- (1.589,0.219);
  \draw[gray, dashed, thin] (0,0.219)
    -- (1.589,0.219);
  \draw[gray, dashed, thin] (5,0) -- (5,-3);
  \draw[gray, dashed, thin] (0,-3) -- (5,-3);
  \draw[thick, blue!70!black]
    plot[smooth, samples=120, domain=0.20:1.589]
      (\x, {-1.5*ln(\x/5) - 1.5});
  \draw[thick, blue!70!black, dotted] (0.10, {-1.5*ln(0.02) - 1.5})
    -- (0.20, {-1.5*ln(0.04) - 1.5});
  \draw[thick, red!75!black] (1.589,0.219)
    -- (5,-3);
  \fill (1.589,0.219) circle (1.6pt);
  \fill (5,-3) circle (1.6pt);
  \node[blue!70!black, anchor=south west, inner sep=2pt] at (0.65,2.0)
    {$-\log m - 1$};
  \node[red!75!black, anchor=south west, inner sep=2pt] at (3.0,-1.0)
    {$(1-m)\lambda - 2$};
\end{tikzpicture}
\caption{The graph of $G_{\psi}$ on $(0,1]$.}
\label{fig:Gpsi}
\end{figure}
\end{example}

We observe a few elementary facts.
\begin{proposition}\label{prop:toricpshcomp}
    Given $\varphi,\psi\in \PSH_{\tor}(X,\omega)$. The following are equivalent:
    \begin{enumerate}
        \item\label{itm:prop:toricpshcomp:1} $\varphi\preceq \psi$,
        \item\label{itm:prop:toricpshcomp:2} $F_{\varphi}\preceq F_{\psi}$, and
        \item\label{itm:prop:toricpshcomp:3} $G_{\psi}\preceq G_{\varphi}$.
    \end{enumerate}
    The same holds if we replace all $\preceq$'s by $\leq$.
\end{proposition}
\begin{proof}
    The equivalence between \ref{itm:prop:toricpshcomp:1} and \ref{itm:prop:toricpshcomp:2} follows from the definition \eqref{eq:FvarphidefGdef}. The equivalence between \ref{itm:prop:toricpshcomp:2} and \ref{itm:prop:toricpshcomp:3} follows from the definition of the Legendre transform.
\end{proof}

Similarly, we have
\begin{proposition}\label{prop:toricpluscst}
    Given $\varphi\in \PSH_{\tor}(X,\omega)$ and $C\in \mathbb{R}$. We have
    \[
    F_{\varphi+C}=F_{\varphi}+C,\quad G_{\varphi+C}=G_{\varphi}-C.
    \]
\end{proposition}
\begin{proof}
    The first equality follows from \eqref{eq:FvarphidefGdef}. The second is obtained by taking the Legendre transform.
\end{proof}

\begin{proposition}\label{prop:toricrooftop}
    Given $\varphi,\psi\in \PSH_{\tor}(X,\omega)$ with $\varphi\land \psi\not\equiv -\infty$, then $\varphi\land \psi\in \PSH_{\tor}(X,\omega)$ and
    \[
    F_{\varphi\land \psi}=F_{\varphi}\land F_{\psi},\quad G_{\varphi\land \psi}=G_{\varphi}\lor G_{\psi}.
    \]
\end{proposition}
The operators $\land$ and $\lor$ are defined in \cref{def:convwedge} and \cref{def:convlor}.
\begin{proof}
By \cref{prop:land_well_def}, we have $\varphi\land \psi\in \PSH(X,\omega)$. Since $\varphi$ and $\psi$ are $T_c$-invariant, the set of candidates in the definition of $\varphi\land\psi$ is preserved by the action of $T_c$. Hence $\varphi\land\psi$ is $T_c$-invariant as well.
So $\varphi\land \psi$ is the biggest element in $\PSH_{\tor}(X,\omega)$ which is dominated by both $\varphi$ and $\psi$. In view of \cref{thm:toricpsh} and \cref{prop:toricpshcomp}, $G_{\varphi\land \psi}$ is the smallest closed proper convex function $G$ on $M_{\mathbb{R}}$ dominating both $G_{\varphi}$ and $G_{\psi}$, which is just $G_{\varphi}\lor G_{\psi}$.

For the statement about $F$, we use \cref{prop:Legendretranssup}:
\[
    F_{\varphi\land\psi}=(G_\varphi\lor G_\psi)^*=\cl(F_\varphi\land F_\psi).
\]
Since $\varphi\land\psi\not\equiv-\infty$, $F_{\varphi\land\psi}=(G_\varphi\lor G_\psi)^*$ is finite on all of $N_{\mathbb R}$. It follows that $F_\varphi\land F_\psi>-\infty$, and hence is finite on all of $N_{\mathbb R}$ as well. Being finite and convex on $N_{\mathbb R}$, it is continuous, hence closed. Therefore
\[
    F_{\varphi\land\psi}=F_\varphi\land F_\psi.
\]
\end{proof}

\begin{example}\label{ex:landminfty}
    Now we can give an example of $\varphi,\psi\in \PSH_{\tor}(X,\omega)$ with $\varphi\land \psi\equiv -\infty$.

    We take $P=[0,1]$ so that $X=\mathbb{P}^1$ and $\omega$ is the Fubini--Study metric. Let $\varphi\in \PSH(X,\omega)$ be such that
    \[
    \varphi(z)=\log \frac{|z|^2}{|z|^2+1}
    \]
    for $z\in \mathbb{C}$. We have computed that $G_{\varphi}$ in \cref{ex:FGsimplevarphi}:
    \[
    G_{\varphi}(m)=
    \begin{cases}
        0,&\textup{if } m=1,\\
        \infty,&\textup{otherwise}.
    \end{cases}
    \]
    Now we define $\psi\in \PSH_{\tor}(X,\omega)$ as the unique function such that
    \[
    \psi(z)=\log \frac{1}{|z|^2+1}
    \]
    for $z\in \mathbb{C}$. Then a similar computation shows that
    \[
    G_{\psi}(m)=
    \begin{cases}
        0,&\textup{if } m=0,\\
        \infty,&\textup{otherwise}.
    \end{cases}
    \]
    Now we claim that $\varphi\land \psi\equiv -\infty$. Otherwise, we would have
    \[
    G_{\varphi\land \psi}=G_{\varphi}\lor G_{\psi}\equiv\infty,
    \]
    which is not proper.
\end{example}

\begin{proposition}\label{prop:toricseq}
    Let $\{\varphi_i\}_{i\in I}$ be a non-empty family in $\PSH_{\tor}(X,\omega)$ uniformly bounded from above. Then $\uscsup[i\in I] \varphi_i\in \PSH_{\tor}(X,\omega)$ and
    \[
    F_{\uscsup[i\in I]\varphi_i}=\bigvee_{i\in I}F_{\varphi_i},\quad G_{\uscsup[i\in I]\varphi_i}=\cl \bigwedge_{i\in I} G_{\varphi_i}.
    \]
    Moreover, if $I$ is finite, then
    \[
        G_{\max_{i\in I}\varphi_i}=\bigwedge_{i\in I} G_{\varphi_i}.
    \]

    Similarly, if $\{\varphi_i\}_{i\in I}$ is a decreasing net in $\PSH_{\tor}(X,\omega)$ such that $\inf_{i\in I}\varphi_i\not\equiv -\infty$, then $\inf_{i\in I}\varphi_i\in \PSH_{\tor}(X,\omega)$ and
    \[
    F_{\inf_{i\in I}\varphi_i}=\inf_{i\in I}F_{\varphi_i},\quad G_{\inf_{i\in I}\varphi_i}=\bigvee_{i\in I}G_{\varphi_i}.
    \]
\end{proposition}
Recall that the closure $\cl$ is defined in \cref{def:convexclosure}.
\begin{proof}
Thanks to \cref{lma:convdecnet} and \cref{prop:supconv}, in both cases, the statement for $F$ is clear. The corresponding statement for $G$ is obtained via \cref{prop:Legendretranssup}. When $I$ is finite, the convex function $\bigvee_i F_{\varphi_i}$ is finite on all $N_{\mathbb R}$, hence closed. Therefore no additional closure is needed.
\end{proof}

The complex Monge--Ampère operator is closely related to the real one:
\begin{proposition}\label{prop:toricMAandrealMA}
    Let $\varphi\in \PSH_{\tor}(X,\omega)$. Then
    \begin{equation}\label{eq:tropMAmea}
    \Trop_* \left(\omega|_{T(\mathbb{C})}+\ddc\varphi|_{T(\mathbb{C})} \right)^n=\MA_{\mathbb{R}}(F_{\varphi}).
    \end{equation}
    In particular,
    \[
        \int_X \omega_{\varphi}^n =\int_{N_{\mathbb{R}}}\MA_{\mathbb{R}}(F_{\varphi})=n! \vol \overline{\{G_{\varphi}<\infty\}}
    \]
    and
    \[
        \int_X  \omega^n=n!\vol P.
    \]
\end{proposition}
Here the real Monge--Ampère operator is defined in \cref{def:realMA}. The normalization of the Lebesgue measure $\vol$ on $M_{\mathbb{R}}$ is such that the fundamental lattice cube has measure $1$.\footnote{In some references like \cite{GKZ08}, the normalization is so that the fundamental lattice cube has measure $n!$. Be careful when making comparisons with these references. See \cite[Theorem~18.4.1]{FPT} for the real Monge--Amp\`ere measure with our convention.}
\begin{proof}
    We only need to prove \eqref{eq:tropMAmea}. By \cref{prop:regularizationconvex}, we can find a decreasing sequence
    \[
        F_j\in \mathcal E^\infty(N_{\mathbb R},P)\cap C^\infty(N_{\mathbb R})
    \]
    with limit $F_\varphi$. By \cref{thm:toricpsh}, each $F_j$ corresponds to a bounded toric potential $\varphi_j\in\PSH_{\tor}(X,\omega)$.

    By \cref{thm:contMA} and \cref{thm:realMAcont}, it suffices to establish \eqref{eq:tropMAmea} for the $\varphi_j$'s. We may therefore reduce to the case where $F_{\varphi}$ is smooth. We write $F=F_{\varphi}$ to simplify the notations. The notations $a_i=\log |z_i|^2$ will be used, where $i=1,\ldots,n$.

    Next we fix an identification $N=\mathbb{Z}^n$. Fix a test function $f\in C^{0}_c(N_{\mathbb{R}})$. We need to show that
    \[
    \int_{\mathbb{C}^{*n}}f(a_1,\ldots,a_n) \left(\ddc \Trop^* F(z_1,\ldots,z_n)\right)^n=\int_{\mathbb{R}^n} f\,\MA_{\mathbb{R}}(F).
    \]
    Using \cref{prop:smoothMAreal} and \eqref{eq:ddctrop}, this reduces to
    \begin{equation}\label{eq:realcplxMAtemp1}
    \begin{split}
        \left( \frac{\mathrm{i}}{2\pi}\right)^n\int_{\mathbb{C}^{*n}}f(a_1,\ldots,a_n)
        \left(\sum_{i,j=1}^n \partial_{i,j}F(a_1,\ldots,a_n)z_i^{-1}\overline{z_j}^{-1}\mathrm{d}z_i\wedge \mathrm{d}\overline{z_j}\right)^n=\\
        n!\int_{\mathbb{R}^n} f \det \nabla^2 F\,\mathrm{d}\vol.
    \end{split}
    \end{equation}
    Expanding the bracket, we get
    \[
    \begin{aligned}
        \left(\sum_{i,j=1}^n \partial_{i,j}Fz_i^{-1}\overline{z_j}^{-1}\mathrm{d}z_i\wedge \mathrm{d}\overline{z_j}\right)^n
        =& \sum_{i_1,\ldots,i_n=1}^n \sum_{j_1,\ldots,j_n=1}^n \partial_{i_1j_1}F\cdots \partial_{i_nj_n}F\cdot \\
        &\mathrm{d}\log z_{i_1}\wedge \mathrm{d}\log \overline{z_{j_1}}\wedge \cdots\wedge \mathrm{d}\log z_{i_n}\wedge \mathrm{d}\log \overline{z_{j_n}},
    \end{aligned}
    \]
    where $\mathrm{d}\log z_i=z_i^{-1}\mathrm{d}z_i$ and $\mathrm{d}\log \overline{z_i}=\overline{z_i}^{-1}\mathrm{d}\overline{z_i}$ are understood.

    Using the apparent symmetry, the expression on the right-hand side becomes
    \[
    \begin{aligned}
    &\sum_{\sigma,\tau\in \mathfrak{S}_n}\prod_{k=1}^n \partial_{\sigma(k)\tau(k)}F\,\mathrm{d}\log z_{\sigma(1)}\wedge \mathrm{d}\log \overline{z_{\tau(1)}}\wedge \cdots\wedge \mathrm{d}\log z_{\sigma(n)}\wedge \mathrm{d}\log \overline{z_{\tau(n)}},\\
    =& n!\sum_{\tau\in \mathfrak{S}_n}\prod_{k=1}^n \partial_{k\tau(k)}F\,\mathrm{d}\log z_{1}\wedge \mathrm{d}\log \overline{z_{\tau(1)}}\wedge \cdots\wedge \mathrm{d}\log z_{n}\wedge \mathrm{d}\log \overline{z_{\tau(n)}}\\
    =& n!\sum_{\tau\in \mathfrak{S}_n}\operatorname{sgn}(\tau)\prod_{k=1}^n \partial_{k\tau(k)}F\,\mathrm{d}\log z_{1}\wedge \mathrm{d}\log \overline{z_1}\wedge \cdots\wedge \mathrm{d}\log z_{n}\wedge \mathrm{d}\log \overline{z_n}\\
    =& n! \det \nabla^2 F\mathrm{d}\log z_{1}\wedge \mathrm{d}\log \overline{z_1}\wedge \cdots\wedge \mathrm{d}\log z_{n}\wedge \mathrm{d}\log \overline{z_n},
    \end{aligned}
    \]
    where $\mathfrak{S}_n$ is the permutation group on $\{1,\ldots,n\}$ and $\operatorname{sgn}(\tau)$ is the sign of $\tau$.

    Next, switch to polar coordinates for each $z_i$: Let $z_i=r_i\exp(\mathrm{i}\theta_i)$ and recall that $r_i=\exp(a_i/2)$. Then the left-hand side of \eqref{eq:realcplxMAtemp1} becomes
    \[
    \begin{split}
    &\frac{n!}{(2\pi)^n}\int_{\mathbb{R}^n\times [0,2\pi)^n} f \det \nabla^2 F \,\mathrm{d}a_1\wedge \mathrm{d}\theta_1\wedge \cdots\wedge \mathrm{d}a_n\wedge \mathrm{d}\theta_n\\
    =& n!\int_{\mathbb{R}^n}f\det \nabla^2 F \,\mathrm{d}a_1\wedge \cdots \wedge\mathrm{d}a_n,
    \end{split}
    \]
    which is exactly what we have expected.
\end{proof}

Next we study the envelope operators developed in \cref{chap:envelopes} in the toric setting.
\begin{definition}\label{def:Newtonbody1}
    Let $\varphi\in \PSH_{\tor}(X,\omega)$. We define its \emph{Newton body}\index{Newton body} as
    \[
    \Delta(\omega,\varphi)\coloneqq \overline{\{G_{\varphi}<\infty\}}\subseteq P.
    \]
    \index{$\Delta(\omega,\varphi)$}
    Note that $\Delta(\omega,\varphi)$ is a convex body.
\end{definition}
By \cref{prop:gradDomFstar}, we have
\[
\Delta(\omega,\varphi)=\overline{\nabla F_{\varphi}(N_{\mathbb{R}})}.
\]

\begin{example}\label{ex:toric-ample:L621}
    By \eqref{eq:G0def}, we have
    \[
    \Delta(\omega,0)=P.
    \]
    In the case of \cref{ex:FGsimplevarphi}, we have
    \[
    \Delta(\omega,\varphi)=\{1\},\quad \Delta(\omega,\psi)=[0,1].
    \]
    Observe that in the latter case,
    \[
    \{G_{\psi}<\infty\}\subsetneq P.
    \]
\end{example}

\begin{proposition}\label{prop:GPenvelope}
    Let $\varphi\in \PSH_{\tor}(X,\omega)$. Then $P_{\omega}[\varphi]\in \PSH_{\tor}(X,\omega)$ and
    \begin{equation}\label{eq:toricPenv}
    G_{P_{\omega}[\varphi]}(x)=
    \begin{cases}
        G_0(x),&\textup{if } x\in \Delta(\omega,\varphi);\\
        \infty,&\textup{otherwise}.
    \end{cases}
    \end{equation}
\end{proposition}
\begin{proof}
    By \eqref{eq:Penvsups}, we have
    \[
    P_{\omega}[\varphi]=\uscsup[C\in \mathbb{R}]\bigl((\varphi+C)\land 0\bigr).
    \]
    It follows from \cref{prop:toricpluscst}, \cref{prop:toricrooftop} and \cref{prop:toricseq} that $P_{\omega}[\varphi]\in \PSH_{\tor}(X,\omega)$. Moreover, by the same propositions, we have
    \[
    G_{P_{\omega}[\varphi]}=\cl \inf_{C\in \mathbb{R}} \left( G_{0}\lor (G_{\varphi}-C) \right).
    \]
    We spell out the last closure. Set
    \[
        K=\inf_{C\in \mathbb{R}} \left( G_{0}\lor (G_{\varphi}-C) \right).
    \]
    If $x\in \Dom G_{\varphi}$, then $G_{\varphi}(x)-C\leq G_0(x)$ for $C$ large enough, hence $K(x)=G_0(x)$. If $x\notin \Dom G_{\varphi}$, then $K(x)=\infty$. Since $\Delta(\omega,\varphi)=\overline{\Dom G_{\varphi}}$ and $G_0$ is continuous on $P$, the lower semicontinuous regularization $\cl K$ is equal to $G_0$ on $\Delta(\omega,\varphi)$ and equal to $\infty$ outside $\Delta(\omega,\varphi)$. This is exactly \eqref{eq:toricPenv}.

\end{proof}

Recall that $\mathrm{H}^0(X,L)$ can be identified with the vector space generated by $\chi^m$ for all $m\in P\cap M$, see \cite[Proposition~4.3.3]{CLS11} and \cite[Theorem~18.8.5]{FPT}. In other words, a character $\chi^m$ of $T$ gives a global section of $L$ if and only if $m\in P$. This gives a beautiful characterization of the lattice points in $P$. The following theorem of Yi Yao gives an analogous characterization of the lattice points in the Newton body.

\begin{theorem}[Yao]\label{thm:Yao}\index{theorem!Yao's theorem (toric Riemann--Roch)}
    Let $\varphi\in \PSH_{\tor}(X,\omega)$. Let $m\in M$.
    \begin{enumerate}
        \item\label{itm:thm:Yao:1} Suppose that $m\in \Delta(\omega,\varphi)$. Then $\chi^m\in \mathrm{H}^0(X,L\otimes \mathcal{I}(\varphi))$.
        \item\label{itm:thm:Yao:2} Let $\sigma\in \Sigma$ be a cone of dimension $p$, with rays $\rho_1,\ldots,\rho_p$. Let $c_1,\ldots,c_p>0$ satisfy $\sum_{i=1}^p c_i=1$, and set $u=\sum_{i=1}^p c_i u_{\rho_i}$. If
        \begin{equation}\label{eq:mawayC1}
            \langle m,-u\rangle -\Supp_{\Delta(\omega,\varphi)}(-u)>1,
        \end{equation}
        then $\chi^m\not\in \mathrm{H}^0(X,L\otimes \mathcal{I}(\varphi))$.
    \end{enumerate}
\end{theorem}
In \ref{itm:thm:Yao:2}, there are exactly $p$ rays in $\sigma(1)$ since $\sigma$ is smooth.
\begin{proof}
    Put $\Delta=\Delta(\omega,\varphi)$. By \cref{prop:PandPI,prop:GPenvelope}, we may replace $\varphi$ by $P_{\omega}[\varphi]$ without changing its multiplier ideal or Newton body. We may therefore assume that $F_{\varphi}-\Supp_{\Delta}$ is bounded on $N_{\mathbb{R}}$. Set
    \[
        H(\nu)=\langle m,\nu\rangle-\Supp_{\Delta}(\nu),\quad \nu\in N_{\mathbb{R}}.
    \]
    Since $\Delta$ is compact, $H$ is Lipschitz and positively homogeneous.

    Observe also that $\mathrm{H}^0(X,L\otimes\mathcal{I}(\varphi))$ is a $T_c$-invariant subspace of $\mathrm{H}^0(X,L)$. Since the $T_c$-weight spaces in $\mathrm{H}^0(X,L)$ are the lines $\mathbb{C}\chi^m$, it is the direct sum of suitable $\mathbb{C}\chi^m$'s.

    If $m\notin P$, then $\chi^m$ is not a section of $L$. Suppose henceforth that $m\in M\cap P$. Membership in $\mathrm{H}^0(X,L\otimes\mathcal{I}(\varphi))$ is equivalent to the finiteness of
    \begin{equation}\label{eq:torictobefinite1}
        \int_{T(\mathbb{C})}|\chi^m|^2\exp(-\Trop^*F_0-\varphi)\,\omega^n.
    \end{equation}
    Since $\Trop^*F_0+\varphi=\Trop^*F_{\varphi}$, the boundedness of $F_{\varphi}-\Supp_{\Delta}$ and \cref{prop:toricMAandrealMA} show that this is equivalent to the finiteness of
    \[
        \int_{N_{\mathbb{R}}}\exp(H(\nu))\,\MA_{\mathbb{R}}(F_0)(\nu).
    \]
    Changing variables by $\nu=\nabla G_0(m')$, we obtain the equivalent integral
    \begin{equation}\label{eq:torictobefinite2}
        \int_P\exp\left(\langle m,\nabla G_0(m')\rangle-\Supp_{\Delta}(\nabla G_0(m'))\right)\,\mathrm{d}m'.
    \end{equation}
    If $m\in\Delta$, its integrand is at most $1$, proving \ref{itm:thm:Yao:1}.

    To prove \ref{itm:thm:Yao:2}, extend $\sigma$ to a maximal cone $\tau\in\Sigma$. Its primitive ray generators $u_1,\ldots,u_n$ form a basis of $N$, and we may arrange that $u_i=u_{\rho_i}$ for $1\leq i\leq p$. If $e_1,\ldots,e_n$ is the dual basis of $M$, then $z_i=\chi^{e_i}$ are coordinates on the affine chart $U_{\tau}\cong\mathbb{C}^n$.

    Let $v\in P\cap M$ be the vertex corresponding to $\tau$. The character $\chi^v$ is a local frame of $L$ on $U_{\tau}$. Writing
    \[
        \nu=\Trop(z)=\sum_{i=1}^n\log|z_i|^2\,u_i,
    \]
    the weight of the singular metric in this frame is $F_{\varphi}(\nu)-\langle v,\nu\rangle$. The local representative of $\chi^m$ is $\chi^{m-v}$, so its squared norm is
    \[
        \exp\left(\langle m,\nu\rangle-F_{\varphi}(\nu)\right).
    \]
    On the unit polydisc, $\omega^n$ is comparable to Euclidean volume. Set $t_i=-\log|z_i|^2$. Integrating out the angular variables and using the boundedness of $F_{\varphi}-\Supp_{\Delta}$, the weighted integral over this polydisc is finite if and only if the following integral is finite:
    \[
        \int_{[0,\infty)^n}\exp\left(H\left(-\sum_{i=1}^n t_i u_i\right)-\sum_{i=1}^n t_i\right)\,\mathrm{d}t.
    \]

    Restrict this integral to the tube parametrized by $s>1$ and $0<y_i<1$ for $2\leq i\leq n$, with
    \[
        t_1=c_1s,\qquad
        t_i=c_is+y_i\quad(2\leq i\leq p),\qquad
        t_i=y_i\quad(p<i\leq n).
    \]
    The Jacobian is $c_1>0$. Since $H$ is Lipschitz and positively homogeneous and $\sum_{i=1}^p c_i=1$, the exponent on this tube is
    \[
        H\left(-\sum_{i=1}^n t_i u_i\right)-\sum_{i=1}^n t_i
        =s\bigl(H(-u)-1\bigr)+\mathcal{O}(1),
    \]
    with a bounded error uniform in $s$ and the $y_i$'s. The local integral is therefore bounded from below by a positive multiple of
    \[
        \int_1^{\infty}\exp\left(s\bigl(H(-u)-1\bigr)\right)\,\mathrm{d}s,
    \]
    which diverges under \eqref{eq:mawayC1}. This proves \ref{itm:thm:Yao:2}.
\end{proof}

\begin{corollary}\label{cor:DXmaintoric}
    Let $\varphi\in \PSH_{\tor}(X,\omega)$. Then
    \[
    \lim_{k\to\infty}\frac{n!}{k^n}h^0\left(X,L^k\otimes \mathcal{I}(k\varphi)\right)= n!\vol \Delta(\omega,\varphi).
    \]
\end{corollary}
\begin{proof}
    Take $k>0$. Let $F_{0,k}$ denote the Guillemin potential associated with the polytope $kP$. On $T(\mathbb C)$, set
    \[
        \widetilde\varphi_k=\Trop^*(kF_\varphi-F_{0,k}).
    \]
    Then $\widetilde\varphi_k$ is the toric potential, with respect to the reference metric associated with $kP$, corresponding to the singular metric $h^k\exp(-k\varphi)$ on $L^k$. It differs from $k\varphi$ by the smooth toric function $\Trop^*(kF_0-F_{0,k})$. Hence
    \[
        \mathcal I(\widetilde\varphi_k)=\mathcal I(k\varphi).
    \]
    Moreover, the convex function associated with $\widetilde\varphi_k$ is $kF_\varphi$, so its Legendre transform is
    \[
        (kF_\varphi)^*(m)=kG_\varphi(m/k).
    \]
    Thus the Newton body of $\widetilde\varphi_k$ is $k\Delta(\omega,\varphi)$.

    Applying \itemref{itm:thm:Yao:1} to the line bundle $L^k$ and the potential $k\varphi$, whose Newton body is $k\Delta(\omega,\varphi)$, we have
    \[
    \# \left( k\Delta(\omega,\varphi)\cap M\right)\leq h^0\left(X,L^k\otimes \mathcal{I}(k\varphi)\right).
    \]
    Here $\#$ denotes the number of elements. Therefore,
    \[
    \lim_{k\to\infty}k^{-n}\#\left( k\Delta(\omega,\varphi)\cap M\right)\leq \varliminf_{k\to\infty} k^{-n} h^0\left(X,L^k\otimes \mathcal{I}(k\varphi)\right).
    \]
    The left-hand side is just $\vol \Delta(\omega,\varphi)$, as follows from \cref{thm:KK12growth}. So we conclude that
    \[
    \varliminf_{k\to\infty}\frac{n!}{k^n}h^0\left(X,L^k\otimes \mathcal{I}(k\varphi)\right)\geq n!\vol \Delta(\omega,\varphi).
    \]
    It remains to prove the reverse inequality.

    We choose an identification $M\cong \mathbb{Z}^n$. In particular, there is a natural distance $\mathrm{d}$ on $M_{\mathbb{R}}\cong \mathbb{R}^n$.
    Given a convex body $K$ and $\epsilon>0$, we let
    \[
    K^{\epsilon}\coloneqq \left\{x\in \mathbb{R}^n: \mathrm{d}(x,K)\leq \epsilon \right\}.
    \]
    This is again a convex body, since it can be realized as the Minkowski sum of $K$ with a closed ball of radius $\epsilon$.

    Applying \itemref{itm:thm:Yao:2} to the polytope $kP$ and the potential $\widetilde\varphi_k$, we can find $C>0$ uniform in $k$ such that for each $k\in \mathbb{Z}_{>0}$,
    \[
    \bigl\{ m\in M:\chi^m\in\mathrm{H}^0(X,L^k\otimes\mathcal{I}(k\varphi))\bigr\}
    \subseteq \left(k\Delta(\omega,\varphi)\right)^C.
    \]
    Indeed, if a point $m\in kP\cap M$ is at distance at least $C$ from $k\Delta(\omega,\varphi)$, the separating hyperplane theorem gives a unit vector $v\in N_{\mathbb R}$ such that
    \[
        \langle m,v\rangle-\Supp_{k\Delta(\omega,\varphi)}(v)\geq C.
    \]
    Since the fan $\Sigma$ is finite and covers $N_{\mathbb R}$, there is a constant $\Lambda>0$ with the following property: every unit vector $v\in N_{\mathbb R}$ can be written in the form
    \[
        v=-\lambda\sum_i c_i u_{\rho_i},\quad 0<\lambda\leq\Lambda,\quad c_i>0,\quad \sum_i c_i=1,
    \]
    where the rays $\rho_i$ are contained in a single cone of $\Sigma$.
    Put
    \[
        u\coloneqq \sum_i c_i u_{\rho_i}.
    \]
    Since the support function is positively homogeneous, we have
    \[
    \begin{aligned}
    \langle m,-u\rangle-\Supp_{k\Delta(\omega,\varphi)}(-u)&=\lambda^{-1}\left(\langle m,v\rangle-\Supp_{k\Delta(\omega,\varphi)}(v)\right)\\
    &\geq C/\Lambda.
    \end{aligned}
    \]
    Choosing $C>\Lambda$ makes \eqref{eq:mawayC1} applicable, uniformly in $k$.

    Fix $\epsilon>0$. Then for large enough $k$, we have
    \[
    \left(k\Delta(\omega,\varphi)\right)^{C}\subseteq k\left(\Delta(\omega,\varphi)^{\epsilon}\right).
    \]
    It follows that
    \[
    h^0\left(X,L^k\otimes \mathcal{I}(k\varphi)\right)\leq \# \left( k\left(\Delta(\omega,\varphi)^{\epsilon}\right) \cap M\right).
    \]
    Applying \cref{thm:KK12growth} again, we conclude that
    \[
    \varlimsup_{k\to\infty} \frac{n!}{k^n} h^0\left(X,L^k\otimes \mathcal{I}(k\varphi)\right)\leq n! \vol\left(\Delta(\omega,\varphi)^{\epsilon}\right).
    \]
    Letting $\epsilon\to 0+$ and applying \cref{thm:contvol}, we find
    \[
    \varlimsup_{k\to\infty} \frac{n!}{k^n} h^0\left(X,L^k\otimes \mathcal{I}(k\varphi)\right)\leq n! \vol \Delta(\omega,\varphi).
    \]
    Our assertion follows.
\end{proof}

\begin{example}\label{ex:toric-ample:L966}
    In general, in the setup of \cref{thm:Yao}, there exists $m\in M\cap (P\setminus \Delta(\omega,\varphi))$ such that $\chi^m\in \mathrm{H}^0(X,L\otimes \mathcal{I}(\varphi))$.

    As a concrete example, let us take $P=[0,1]$. Take $\varphi$ so that $\Delta(\omega,\varphi)=[0,1/2]$. We claim that $\chi^1$ is $L^2$-integrable.

    It suffices to verify the convergence of \eqref{eq:torictobefinite2}. Recall that
    \[
    \nabla G_0(m')=\log \frac{m'}{1-m'},\quad m'\in [0,1],
    \]
    while
    \[
    \Supp_{[0,1/2]}(a)=
    \begin{cases}
        a/2,&\textup{if }a>0;\\
        0,&\textup{otherwise}.
    \end{cases}
    \]
    Therefore, \eqref{eq:torictobefinite2} becomes
    \[
    \int_0^{1/2} \frac{m'}{1-m'}\,\mathrm{d}m'+\int_{1/2}^1  \left(\frac{m'}{1-m'}\right)^{1/2}\,\mathrm{d}m'<\infty.
    \]
\end{example}

We interpret various classes of potentials studied in \cref{subsec:fullmass} in the toric setting.
\begin{proposition}\label{prop:toric-ample:L990}
    Let $\varphi\in \PSH_{\tor}(X,\omega)$. Then the following are equivalent:
    \begin{enumerate}
        \item\label{itm:chapter-toric:L976:1} $\varphi\in \mathcal{E}^{\infty}(X,\omega)$;
        \item\label{itm:chapter-toric:L976:2} $F_{\varphi}\in \mathcal{E}^{\infty}(N_{\mathbb{R}},P)$;
        \item\label{itm:chapter-toric:L976:3} $G_{\varphi}\sim G_0$.
    \end{enumerate}
\end{proposition}
The notation $\mathcal{E}^{\infty}(N_{\mathbb{R}},P)$ is defined in \cref{def:convexPfunctions}.
\begin{proof}
    This follows immediately from \cref{prop:toricpshcomp}.
\end{proof}

\begin{proposition}\label{prop:toric-ample:L1003}
    Let $\varphi\in \PSH_{\tor}(X,\omega)$. Then the following are equivalent:
    \begin{enumerate}
        \item\label{itm:chapter-toric:L989:1} $\varphi\in \mathcal{E}(X,\omega)$;
        \item\label{itm:chapter-toric:L989:2} $F_{\varphi}\in \mathcal{E}(N_{\mathbb{R}},P)$;
        \item\label{itm:chapter-toric:L989:3} $\overline{\Dom G_{\varphi}}=P$.
    \end{enumerate}
\end{proposition}
The notation $\mathcal{E}(N_{\mathbb{R}},P)$ is defined in \cref{def:convexPfunctions}.
\begin{proof}
\ref{itm:chapter-toric:L989:1} $\iff$ \ref{itm:chapter-toric:L989:3}. By \cref{prop:toricMAandrealMA},
    \[
        \int_X \omega_{\varphi}^n=\int_{T(\mathbb{C})} \left(\omega|_{T(\mathbb{C})}+\ddc \varphi|_{T(\mathbb{C})} \right)^n= n! \vol\overline {\Dom G_{\varphi}},\quad \int_X \omega^n=n!\vol P,
    \]
    where in the first equality, we use the fact that $\omega_{\varphi}^n$ does not charge pluripolar sets \itemref{itm:prop:npp1:4}.

    Therefore, \ref{itm:chapter-toric:L989:1} and \ref{itm:chapter-toric:L989:3} are equivalent.

\ref{itm:chapter-toric:L989:2} $\iff$ \ref{itm:chapter-toric:L989:3}. This follows from \cref{prop:gradDomFstar}.
\end{proof}

\begin{proposition}\label{prop:Etoric}
    Let $\varphi\in \PSH_{\tor}(X,\omega)$. Then
    \begin{equation}\label{eq:Eomega_toric}
        E_{\omega}(\varphi)=n!\int_P \left(G_{0}-G_{\varphi}\right)\,\mathrm{d}\vol.
    \end{equation}
\end{proposition}
\begin{proof}
    For $\psi\in \PSH_{\tor}(X,\omega)$, write
    \[
        \mathcal C(\psi)\coloneqq n!\int_P\left(G_0-G_\psi\right)\,\mathrm{d}\vol.
    \]

    \textbf{Step~1}. We assume that $\varphi$ is bounded.

    Put $h=F_{\varphi}-F_0$, $\varphi_t=t\varphi$ and
    \[
        F_t\coloneqq F_0+th
    \]
    for $t\in[0,1]$. Then $F_t=F_{\varphi_t}$ and $G_t\coloneqq F_t^*=G_{\varphi_t}$. Since $\varphi$ is bounded, $h$ is bounded and continuous; in particular, $F_t\in \mathcal E^\infty(N_{\mathbb R},P)$ for every $t\in[0,1]$.

    We claim that
    \begin{equation}\label{eq:Etoric_convex_derivative}
        \frac{\mathrm d}{\mathrm dt}\mathcal C(\varphi_t)=\int_{N_{\mathbb R}}h\,\MA_{\mathbb R}(F_t).
    \end{equation}
    Let $s>0$ and $t,t+s\in[0,1]$. At a point $m\in \Int P$ where both $G_t$ and $G_{t+s}$ are differentiable, we have
    \[
        \begin{aligned}
            G_{t+s}(m)&\geq G_t(m)-s h\bigl(\nabla G_t(m)\bigr),\\
            G_t(m)&\geq G_{t+s}(m)+s h\bigl(\nabla G_{t+s}(m)\bigr).
        \end{aligned}
    \]
    Therefore,
    \[
        h\bigl(\nabla G_{t+s}(m)\bigr)\leq \frac{G_t(m)-G_{t+s}(m)}{s}\leq h\bigl(\nabla G_t(m)\bigr).
    \]
    Integrating the last inequality and using \cref{lma:realMA_Legendre_pushforward} gives
    \[
        \int_{N_{\mathbb R}}h\,\MA_{\mathbb R}(F_{t+s})\leq \frac{\mathcal C(\varphi_{t+s})-\mathcal C(\varphi_t)}{s}\leq \int_{N_{\mathbb R}}h\,\MA_{\mathbb R}(F_t).
    \]
    Letting $s\to 0+$ and using \cref{thm:realMAcont} for $F_{t+s}\to F_t$, we obtain the right derivative in \eqref{eq:Etoric_convex_derivative}. The left derivative follows in the same way. Hence \eqref{eq:Etoric_convex_derivative} holds.

    By \cref{prop:toricMAandrealMA}, the right-hand side of \eqref{eq:Etoric_convex_derivative} is
    \[
        \int_X\varphi\,\omega_{\varphi_t}^n.
    \]
    On the other hand, \cref{lma:Eaffinediff} gives
    \[
        \frac{\mathrm d}{\mathrm dt}E_{\omega}(\varphi_t)=\int_X\varphi\,\omega_{\varphi_t}^n.
    \]
    Thus $E_{\omega}(\varphi_t)-\mathcal C(\varphi_t)$ is constant in $t$. It vanishes at $t=0$, hence \eqref{eq:Eomega_toric} follows by taking $t=1$.

    \textbf{Step~2}. We handle the general case.

    Both sides of \eqref{eq:Eomega_toric} change by the same additive constant when we replace $\varphi$ by $\varphi+A$, namely by
    \[
        A\int_X\omega^n=A\,n!\vol P.
    \]
    We may therefore assume that $\varphi\leq 0$. For $C>0$, put
    \[
        \varphi_C\coloneqq \varphi\lor(-C).
    \]
    Then $\varphi_C\in \PSH_{\tor}(X,\omega)\cap L^{\infty}(X)$ and $\varphi_C\searrow \varphi$ as $C\to\infty$. By \cref{prop:toricseq}, $G_{\varphi_C}\nearrow G_{\varphi}$. Applying Step~1 to $\varphi_C$, we get
    \[
        E_{\omega}(\varphi_C)=n!\int_P \left(G_0-G_{\varphi_C}\right)\,\mathrm{d}\vol.
    \]
    Letting $C\to\infty$, the left-hand side converges to $E_{\omega}(\varphi)$ by \cref{prop:E_cont_decnet}. On the right-hand side, since $\varphi_C\leq 0$, we have $G_{\varphi_C}\geq G_0$, and the monotone convergence theorem applied to $G_{\varphi_C}-G_0\nearrow G_{\varphi}-G_0$ gives \eqref{eq:Eomega_toric}.
\end{proof}

\begin{corollary}\label{cor:toric-ample:L1092}
    Let $\varphi\in \PSH_{\tor}(X,\omega)$. Then the following are equivalent:
    \begin{enumerate}
        \item\label{itm:chapter-toric:L1032:1} $\varphi\in \mathcal{E}^1(X,\omega)$;
        \item\label{itm:chapter-toric:L1032:2} $F_{\varphi}\in \mathcal{E}^1(N_{\mathbb{R}},P)$;
        \item\label{itm:chapter-toric:L1032:3} $G_{\varphi}\in L^1(P)$.
    \end{enumerate}
\end{corollary}
The notation $\mathcal{E}^1(N_{\mathbb{R}},P)$ is defined in \cref{def:convexPfunctions}.
\begin{proof}
    The equivalence between \ref{itm:chapter-toric:L1032:1} and \ref{itm:chapter-toric:L1032:3} follows from \cref{prop:Etoric}. The equivalence between \ref{itm:chapter-toric:L1032:2} and \ref{itm:chapter-toric:L1032:3} is the definition of $\mathcal{E}^1(N_{\mathbb{R}},P)$, since $F_{\varphi}^*=G_{\varphi}$.
\end{proof}
\begin{definition}\label{def:toric-ample:L1104}
    We define
    \index{$\mathcal{E}^{\infty}_{\tor}(X,\omega)$}
    \index{$\mathcal{E}^{1}_{\tor}(X,\omega)$}
    \index{$\mathcal{E}_{\tor}(X,\omega)$}
    \[
        \begin{aligned}
        \mathcal{E}^{\infty}_{\tor}(X,\omega)= & \mathcal{E}^{\infty}(X,\omega) \cap \PSH_{\tor}(X,\omega),\\
        \mathcal{E}^{1}_{\tor}(X,\omega)= & \mathcal{E}^{1}(X,\omega) \cap \PSH_{\tor}(X,\omega),\\
        \mathcal{E}_{\tor}(X,\omega)= & \mathcal{E}(X,\omega) \cap \PSH_{\tor}(X,\omega).
        \end{aligned}
    \]
\end{definition}

\begin{corollary}\label{cor:toricd1}
    Let $\varphi,\psi\in \mathcal{E}^1_{\tor}(X,\omega)$. Then
    \[
        d_1(\varphi,\psi)= n!\int_P \left( 2\left(G_{\varphi}\lor G_{\psi}\right)-G_{\varphi}-G_{\psi}\right)\,\mathrm{d}\vol.
    \]
\end{corollary}
\begin{proof}
    This follows from \cref{prop:Etoric}, \cref{prop:toricrooftop} and \cref{def:d1onE12}.
\end{proof}

\begin{partbacktext}
\part{The theory of \texorpdfstring{$\mathcal{I}$}{I}-good singularities}
\noindent
This part develops the comparison theory of singularities that forms the technical core of the book. The basic objects are positive-mass $\theta$-psh potentials, their $P$- and $\mathcal I$-singularity types, and the envelopes introduced in \cref{chap:envelopes}. The central theme is that the ordinary pluripotential-theoretic equivalence relation and the equivalence relation detected by multiplier ideals are different in general, but coincide on a large and geometrically meaningful class of singularities: the $\mathcal I$-good ones.

The first chapter of the part, \cref{chap:comp}, introduces the general language for comparing singularities. It studies the $P$- and $\mathcal I$-preorders and defines the pseudometric $d_S$ on positive-mass potentials. This pseudometric is not restricted to $\mathcal I$-good potentials; rather, it is the ambient convergence tool used to compare singularity types. The chapter proves the main continuity and convergence properties of $d_S$, including its compatibility with envelopes and with changes of the background form.

In \cref{chap:Igood} we introduce $\mathcal I$-good singularities. The main structural result is the characterization theorem \cref{thm:charIgoodasclosure}, which relates $\mathcal I$-goodness to approximation by analytic singularities, equality between non-pluripolar mass and volume, and the coincidence of the $P$- and $\mathcal I$-envelopes. The chapter also develops the mixed volume formalism needed later and proves the asymptotic Riemann--Roch formula \cref{thm:DXmain1}, which identifies the non-pluripolar volume of an $\mathcal I$-good current with the asymptotic growth of multiplier-ideal sections.

\cref{chap:trace} develops the trace operator. Its purpose is to give a meaningful restriction of a quasi-psh singularity to a possibly singular subvariety, in a way compatible with multiplier ideals and birational modifications. The trace operator is accompanied by analytic Bertini-type results and by restricted-volume formulae. These results allow the singularity theory developed in the previous chapters to be used inductively on subvarieties.

\cref{chap:testcurve} studies test curves, namely concave decreasing families of model potentials. It proves the Ross--Witt Nystr\"om correspondence between test curves and geodesic rays and develops the operations on test curves needed later. The chapter also introduces the $\mathcal I$-model version of the theory, which will be used in the non-Archimedean applications.

\cref{chap:Okounkov} develops partial Okounkov bodies for singular metrics. These convex bodies extend the classical Newton--Okounkov construction to the setting of psh metrics and record the valuative data of multiplier ideals. The chapter treats both algebraic and transcendental versions and proves the continuity and volume formulae needed in the big toric and non-Archimedean chapters.

Finally, \cref{chap:bdiv} translates the preceding singularity theory into the language of b-divisors. It associates with a positive current a nef b-divisor recording the regular parts of its pullbacks on all modifications, and it compares $\mathcal I$-equivalence with b-divisorial data. The mixed intersection theory of nef b-divisors developed there provides a geometric counterpart to the mixed volume theory of $\mathcal I$-good currents.

The chapters are arranged linearly. The comparison theory of \cref{chap:comp} is used throughout the part; \cref{chap:Igood,chap:trace} provide the main analytic tools; \cref{chap:testcurve,chap:Okounkov,chap:bdiv} develop three complementary languages for the same singularity theory: geodesic/\allowbreak test-curve language, Okounkov-body language and b-divisorial language.

\newpage
\noindent These languages are then used separately in the application chapters \cref{chap:toric_big,chap:NApotential,chap:Bergman}.
\newpage
\end{partbacktext}

\motto{Algebra is the offer made by the devil to the mathematician. The devil says: ``I will give you this powerful machine, it will answer any question you like. All you need to do is give me your soul: give up geometry and you will have this marvelous machine.''

--- Michael Atiyah\footnote{Sir Michael Francis Atiyah (1929--2019) wrote the influential \emph{Introduction to commutative algebra} together with I.~G.~MacDonald, whose name is often omitted or misspelled.}}
\chapter{Comparison of singularities}
\label{chap:comp}

A central theme of this book is that singularities of quasi-plurisubharmonic functions can be compared at several levels of precision. The preorder $\preceq$ of \cref{def:parorder} compares potentials up to additive constants. It is the most concrete notion, but it is often too fine for global pluripotential theory: the $P$-envelope identifies potentials with the same model representative, while multiplier ideals identify singularities that have the same integrability behavior. This chapter makes these two coarser comparisons systematic.

We first introduce the $P$-preorder and the $\mathcal I$-preorder on $\QPSH(X)$. The $P$-preorder is defined by comparing $P$-envelopes, and its equivalence classes are the $P$-singularity types. The $\mathcal I$-preorder is defined by comparing multiplier ideals $\mathcal I(k\varphi)$ for all $k\in\mathbb Z_{>0}$. These two preorders are compatible with ordinary singularity comparison, but they capture different information: the former is governed by non-pluripolar mass, while the latter is governed by valuative and multiplier-ideal data.

The second main object of the chapter is the pseudometric $d_S$. For positive-mass potentials in a fixed class, $d_S$ is defined through the Darvas $d_1$-distance between the associated geodesic rays. Its zero locus is precisely $P$-equivalence. Thus $d_S$ should be viewed as a convergence structure on $P$-singularity types. The chapter proves the basic properties of $d_S$, including monotone convergence, sandwich criteria and compatibility with changing the background form.

The final results establish the continuity properties that will be used throughout the rest of the book. We prove continuity of the $\mathcal I$-envelope, continuity of non-pluripolar mass and volume, and a general equisingularity criterion. These results allow one to pass from analytic approximations to limiting singularities while controlling both $P$-types and multiplier ideals. They are used repeatedly in the theory of $\mathcal I$-good singularities, trace operators, test curves, partial Okounkov bodies and non-Archimedean potentials.

\section{The \texorpdfstring{$P$}{P} and \texorpdfstring{$\mathcal{I}$}{I}-preorders}\label{sec:PIpartialorder}
Let $X$ be a connected compact Kähler manifold of dimension $n$.

Recall that we defined a preorder on $\QPSH(X)$ in \cref{def:parorder} to compare the singularity types of quasi-plurisubharmonic functions. The problem with this preorder is that it is too fine. For our purposes, it is useful to consider rougher relations.

\subsection{The definitions of the preorders}
\label{subsec:preorder_def}

Recall that the $P$-envelope is defined in \cref{def:Penv}.
\begin{definition}\label{def:Pmoresing}
    For $\varphi,\psi\in \QPSH(X)$, we say $\varphi$ is \emph{$P$-more singular than $\psi$}\index{P-more singular@$P$-more singular} and write $\varphi\preceq_P \psi$\index{$\varphi\preceq_P \psi$} if for some (equivalently, every, due to \cref{lma:Pproj_insens_omega} below) closed smooth real $(1,1)$-form $\theta$ on $X$ such that $\varphi,\psi\in \PSH(X,\theta)_{>0}$, we have
        \begin{equation}\label{eq:Penvleq}
        P_{\theta}[\varphi]\leq P_{\theta}[\psi].
        \end{equation}
    If $\varphi\preceq_P \psi$ and $\psi \preceq_P \varphi$, we write $\varphi\sim_P \psi$ and say that $\varphi$ and $\psi$ have the same \emph{$P$-singularity type}\index{P-singularity type@$P$-singularity type}.
\end{definition}
Recall that $\PSH(X,\theta)_{>0}$ was introduced in \eqref{eq:PSHpos}.
Note that if $\varphi\preceq \psi$, then $\varphi\preceq_P \psi$. So the $P$-preorder is \emph{coarser} than $\preceq$.

The $P$-equivalence relation can also be formulated using the language of relative full mass classes developed in \cref{subsec:fullmass}.
In fact, $\varphi,\psi\in \PSH(X,\theta)_{>0}$ satisfies $\varphi\sim_P \psi$ if and only if $\varphi,\psi \in \mathcal{E}(X,\theta;\phi)$ for a common model potential $\phi\in \PSH(X,\theta)_{>0}$.

The condition \eqref{eq:Penvleq} is independent of the choice of $\theta$:
\begin{lemma}\label{lma:Pproj_insens_omega}
    Let $\varphi,\psi\in \PSH(X,\theta)_{>0}$. For any Kähler form $\omega$ on $X$, the following are equivalent:
    \begin{enumerate}
        \item\label{itm:lma:Pproj_insens_omega:1} $P_{\theta}[\varphi]\leq P_{\theta}[\psi]$;
        \item\label{itm:lma:Pproj_insens_omega:2} $P_{\theta+\omega}[\varphi]\leq P_{\theta+\omega}[\psi]$.
    \end{enumerate}
\end{lemma}
In particular, $\preceq_P$ defines a preorder on $\QPSH(X)$.
\begin{proof}
    \ref{itm:lma:Pproj_insens_omega:1} $\implies$ \ref{itm:lma:Pproj_insens_omega:2}. Assume \ref{itm:lma:Pproj_insens_omega:1}. Observe that
    \[
    \varphi\preceq P_{\theta}[\varphi]\leq P_{\theta+\omega}[\varphi].
    \]
    Taking the $P_{\theta+\omega}$-envelope, we find
    \[
    P_{\theta+\omega}[\varphi] \leq P_{\theta+\omega}\Bigl[P_{\theta}[\varphi]\Bigr] \leq P_{\theta+\omega}\Bigl[P_{\theta+\omega}[\varphi]\Bigr]=P_{\theta+\omega}[\varphi],
    \]
    where the equality follows from \cref{thm:Pvarphidiffdef}. In particular,
    \begin{equation}\label{eq:doubleP}
    P_{\theta+\omega}[\varphi]=P_{\theta+\omega}\Bigl[P_{\theta}[\varphi]\Bigr].
    \end{equation}
    A similar formula holds for $\psi$. So we see that \ref{itm:lma:Pproj_insens_omega:2} holds.

    \ref{itm:lma:Pproj_insens_omega:2} $\implies$ \ref{itm:lma:Pproj_insens_omega:1}. Assume \ref{itm:lma:Pproj_insens_omega:2}. By \eqref{eq:doubleP}, we may assume that $\varphi$ and $\psi$ are both model potentials in $\PSH(X,\theta)_{>0}$.

    Observe that $\varphi\lor \psi\preceq P_{\theta+\omega}[\psi]$.
    It follows that $P_{\theta+\omega}[\varphi\lor \psi]\leq P_{\theta+\omega}[\psi]$.
    The reverse inequality is trivial, so
    \[
    P_{\theta+\omega}[\varphi\lor \psi]= P_{\theta+\omega}[\psi].
    \]
    From the direction we have proved, for any $C\geq 1$,
    \[
    P_{\theta+C\omega}[\varphi\lor \psi]= P_{\theta+C\omega}[\psi].
    \]
    So by \cref{prop:Ppresmass},
    \[
    \int_X \Bigl(\theta+C\omega+\ddc(\varphi\lor \psi)\Bigr)^n=\int_X\left(\theta+C\omega+\ddc\psi\right)^n.
    \]
    Since both sides are polynomials in $C$, the equality extends to $C=0$, namely,
    \[
    \int_X \theta_{\varphi\lor \psi}^n=\int_X \theta_{\psi}^n.
    \]
    In particular, $\varphi\lor \psi\leq P_{\theta}[\psi]=\psi$ by \eqref{eq:Penvdef}. So \ref{itm:lma:Pproj_insens_omega:1} follows.
\end{proof}
As a consequence of \cref{lma:Pproj_insens_omega}, we can define the $P$-preorder at the level of currents.
Given closed positive $(1,1)$-currents $T=\theta_{\varphi}$, $S=\theta'_{\psi}$, we write $T\preceq_P S$ (resp. $T\sim_P S$) if $\varphi\preceq_P \psi$ (resp. $\varphi\sim_P \psi$). This definition is independent of the decompositions of $T$ and $S$.

As a first example of $P$-equivalence, we have:
\begin{example}\label{ex:Pequiv}
    Let $\theta$ be a closed smooth real $(1,1)$-form on $X$ and $\varphi\in \PSH(X,\theta)_{>0}$. Then
    \[
        \varphi\sim_P P_{\theta}[\varphi]\footnotemark.
    \]
    \footnotetext{It is not known to us whether the same holds when $\varphi$ has vanishing mass.}
    This follows immediately from \cref{thm:Pvarphidiffdef}.
\end{example}

We give a very useful criterion of the $P$-equivalence in terms of the non-pluripolar masses.
\begin{proposition}\label{prop:Pequivchar2}
    Let $\theta$ be a closed smooth real $(1,1)$-form on $X$ representing a big cohomology class.
    Let $\varphi,\psi\in \PSH(X,\theta)$ and $\varphi\preceq_P \psi$. Then the following are equivalent:
    \begin{enumerate}
        \item\label{itm:prop:Pequivchar2:1} $\varphi\sim_P \psi$;
        \item\label{itm:prop:Pequivchar2:2} for each $j=0,\ldots,n$, we have
            \begin{equation}\label{eq:mixedmassequal}
                \int_X \theta_{\varphi}^j\wedge \theta_{V_{\theta}}^{n-j}=\int_X \theta_{\psi}^j\wedge \theta_{V_{\theta}}^{n-j}.
            \end{equation}
    \end{enumerate}
    Assume furthermore that $\varphi,\psi\in \PSH(X,\theta)_{>0}$, then these conditions are equivalent to the following:
    \begin{enumerate}[resume]
        \item\label{itm:prop:Pequivchar2:3} We have
            \[
                \int_X \theta_{\varphi}^n=\int_X \theta_{\psi}^n.
            \]
    \end{enumerate}
\end{proposition}
Recall that $V_{\theta}$ is introduced in \eqref{eq:Vtheta}.

\begin{proof}
    We first prove the equivalence between \ref{itm:prop:Pequivchar2:1} and \ref{itm:prop:Pequivchar2:3} when $\varphi,\psi\in \PSH(X,\theta)_{>0}$.

    \ref{itm:prop:Pequivchar2:1} $\implies$ \ref{itm:prop:Pequivchar2:3}. Assume that $\varphi\sim_P \psi$. By \cref{lma:Pproj_insens_omega}, we have
    \[
        P_{\theta}[\varphi]=P_{\theta}[\psi].
    \]
    So \ref{itm:prop:Pequivchar2:3} follows from \cref{prop:Ppresmass}.

    \ref{itm:prop:Pequivchar2:3} $\implies$ \ref{itm:prop:Pequivchar2:1}. Since $\varphi\preceq_P\psi$, we have $P_{\theta}[\varphi]\leq P_{\theta}[\psi]$.
    Moreover, \cref{prop:Ppresmass} shows that $P_{\theta}[\varphi]$ and $P_{\theta}[\psi]$ have the same positive mass.
    By \cref{thm:Pvarphidiffdef}, $P_{\theta}[\psi]$ is a candidate in the characterization of $P_{\theta}\bigl[P_{\theta}[\varphi]\bigr]=P_{\theta}[\varphi]$.
    Hence $P_{\theta}[\psi]\leq P_{\theta}[\varphi]$, and therefore $P_{\theta}[\varphi]=P_{\theta}[\psi]$.
    So \ref{itm:prop:Pequivchar2:1} follows.

    Let us come back to the general case with $\varphi,\psi\in \PSH(X,\theta)$.

    \ref{itm:prop:Pequivchar2:1} $\implies$ \ref{itm:prop:Pequivchar2:2}.
    Fix $j\in \{0,\ldots,n\}$. We argue \eqref{eq:mixedmassequal}.

    Take a Kähler form $\omega$ on $X$. By \cref{lma:Pproj_insens_omega}, for each $\epsilon>0$, we have
    \[
        P_{\theta+\epsilon\omega}[\varphi]=  P_{\theta+\epsilon\omega}[\psi].
    \]
    It follows from \cref{prop:Ppresmass} that
    \[
        \begin{aligned}
        \int_X \left(\theta+\epsilon\omega+\ddc \psi  \right)^j\wedge \theta_{V_{\theta}}^{n-j}=&\int_X \Bigl(\theta+\epsilon\omega+\ddc P_{\theta+\epsilon\omega}[\psi]  \Bigr)^j\wedge \theta_{V_{\theta}}^{n-j}\\
        =&\int_X \Bigl(\theta+\epsilon\omega+\ddc P_{\theta+\epsilon\omega}[\varphi]  \Bigr)^j\wedge \theta_{V_{\theta}}^{n-j} \\
        =& \int_X \left(\theta+\epsilon\omega+\ddc \varphi  \right)^j\wedge \theta_{V_{\theta}}^{n-j}.
        \end{aligned}
    \]
    Since the two extremes are both polynomials in $\epsilon$ by the multi-linearity of the non-pluripolar products \itemref{itm:prop:npp1:6}, we conclude that the same holds when $\epsilon=0$, that is, \eqref{eq:mixedmassequal} holds.

    \ref{itm:prop:Pequivchar2:2} $\implies$ \ref{itm:prop:Pequivchar2:1}. Assume \eqref{eq:mixedmassequal} holds for all $j=0,\ldots,n$. Then the same holds with $P_{\theta}[\psi]$ in place of $\psi$, thanks to \cref{prop:Ppresmass}.
    Therefore, for each $t\in (0,1)$, we have
    \[
        \int_X \theta_{t\varphi+(1-t)V_{\theta}}^n=  \int_X \theta_{tP_{\theta}[\psi]+(1-t)V_{\theta}}^n
    \]
    by the binomial expansion. Furthermore, $t\varphi+(1-t)V_{\theta}\preceq tP_{\theta}[\psi]+(1-t)V_{\theta}$.
    By the implication \ref{itm:prop:Pequivchar2:3} $\implies$ \ref{itm:prop:Pequivchar2:1}, we have
    \[
        t\varphi+(1-t)V_{\theta}\sim_P t P_{\theta}[\psi]+(1-t)V_{\theta}
    \]
    for each $t\in (0,1)$.

    Fix a Kähler form $\omega$ on $X$.
    From the implication \ref{itm:prop:Pequivchar2:1} $\implies$ \ref{itm:prop:Pequivchar2:3}, we have
    \[
        \int_X (\theta+\omega)_{t\varphi+(1-t)V_{\theta}}^n=  \int_X (\theta+\omega)_{tP_{\theta}[\psi]+(1-t)V_{\theta}}^n.
    \]
    Since both sides are polynomials in $t$ by the multi-linearity of the non-pluripolar products \itemref{itm:prop:npp1:6}, the same holds when $t=1$. From the implication \ref{itm:prop:Pequivchar2:3} $\implies$ \ref{itm:prop:Pequivchar2:1} again, we have $\varphi\sim_P P_{\theta}[\psi]\sim_P \psi$.
\end{proof}

Next we introduce a different preorder.
\begin{proposition}\label{prop:Iequivchar2}
    Given $\varphi,\psi\in \QPSH(X)$, the following are equivalent:
    \begin{enumerate}
        \item\label{itm:prop:Iequivchar2:1} For any $k\in \mathbb{Z}_{>0}$, we have
        \[
            \mathcal{I}(k\varphi)\subseteq\mathcal{I}(k\psi);
        \]
        \item\label{itm:prop:Iequivchar2:2} for any $\lambda\in \mathbb{R}_{>0}$, we have
            \[
            \mathcal{I}(\lambda\varphi)\subseteq\mathcal{I}(\lambda\psi);
            \]
        \item\label{itm:prop:Iequivchar2:3} for any modification $\pi\colon Y\rightarrow X$ and any $y\in Y$, we have
            \[
                \nu(\pi^*\varphi,y)\geq \nu(\pi^*\psi,y);
            \]
        \item\label{itm:prop:Iequivchar2:4} for any proper bimeromorphic morphism $\pi\colon Y\rightarrow X$ from a Kähler manifold and any $y\in Y$, we have
            \[
                \nu(\pi^*\varphi,y)\geq \nu(\pi^*\psi,y);
            \]
        \item\label{itm:prop:Iequivchar2:5} for any prime divisor $E$ over $X$, we have
            \[
                \nu(\varphi,E)\geq \nu(\psi,E).
            \]
    \end{enumerate}
\end{proposition}
\begin{proof}
    The proof is almost identical to that of \cref{prop:Iequivchar}.
\end{proof}

\begin{definition}\label{def:Imoresing}
    Let $\varphi,\psi\in \QPSH(X)$, we say $\varphi$ is \emph{$\mathcal{I}$-more singular than $\psi$}\index{I-more singular@$\mathcal{I}$-more singular} and write $\varphi\preceq_{\mathcal{I}} \psi$\index{$\varphi\preceq_{\mathcal{I}} \psi$} if the equivalent conditions in \cref{prop:Iequivchar2} are satisfied.
\end{definition}
It is clear that $\preceq_{\mathcal{I}}$ is a preorder on $\QPSH(X)$.

Note that $\varphi\preceq_{\mathcal{I}}\psi$ and $\psi\preceq_{\mathcal{I}}\varphi$ both hold if and only if $\varphi\sim_{\mathcal{I}}\psi$ in the sense of \cref{def:Iequiv}.

Given closed positive $(1,1)$-currents $T=\theta_{\varphi}$, $S=\theta'_{\psi}$, we write $T\preceq_{\mathcal{I}} S$ (resp. $T\sim_{\mathcal{I}} S$) if $\varphi\preceq_{\mathcal{I}} \psi$ (resp. $\varphi\sim_{\mathcal{I}} \psi$). This definition is independent of the decompositions of $T$ and $S$.

\begin{lemma}\label{lma:Plor1}
    Let $\theta$ be a closed smooth real $(1,1)$-form on $X$ representing a big cohomology class.
    Let $\varphi,\psi\in \PSH(X,\theta)_{>0}$. Then
    \begin{equation}\label{eq:Plor1}
    P_{\theta}[\varphi\lor \psi]=P_{\theta}\Bigl[ P_{\theta}[\varphi]\lor P_{\theta}[\psi] \Bigr].
    \end{equation}
\end{lemma}
\begin{proof}
    Since $\varphi\lor \psi\preceq P_{\theta}[\varphi]\lor P_{\theta}[\psi]$, the $\leq$ direction of \eqref{eq:Plor1} follows. Conversely, it suffices to show that
    \[
    P_{\theta}[\varphi\lor \psi]\geq P_{\theta}[\varphi]\lor P_{\theta}[\psi],
    \]
    which follows from $\varphi\preceq\varphi\lor\psi$, $\psi\preceq\varphi\lor\psi$, and the monotonicity of the $P$-envelope with respect to $\preceq$.
\end{proof}

\begin{lemma}\label{lma:reform_preceqP}
    Let $\varphi,\psi\in \QPSH(X)$. Then the following are equivalent:
    \begin{enumerate}
        \item\label{itm:lma:reform_preceqP:1} $\varphi\preceq_P \psi$ (resp. $\varphi\preceq_{\mathcal{I}} \psi$);
        \item\label{itm:lma:reform_preceqP:2} $\varphi\lor \psi\sim_P\psi$ (resp. $\varphi\lor \psi\sim_{\mathcal{I}}\psi$).
    \end{enumerate}
\end{lemma}
\begin{proof}
Take a closed real smooth $(1,1)$-form $\theta$ on $X$ such that $\varphi,\psi\in \PSH(X,\theta)_{>0}$. We only prove the $P$ case, the $\mathcal{I}$ case is similar.

\ref{itm:lma:reform_preceqP:2} $\implies$ \ref{itm:lma:reform_preceqP:1}. Assume \ref{itm:lma:reform_preceqP:2}. By \ref{itm:lma:reform_preceqP:2} and \cref{ex:Pequiv}, $P_{\theta}[\varphi\lor \psi]=P_{\theta}[\psi]\sim_P \psi$. But $\varphi\preceq P_{\theta}[\varphi\lor \psi]$, so \ref{itm:lma:reform_preceqP:1} follows.

\ref{itm:lma:reform_preceqP:1} $\implies$ \ref{itm:lma:reform_preceqP:2}. Assume \ref{itm:lma:reform_preceqP:1}. We may assume that $\varphi,\psi$ are both model by \cref{lma:Plor1}.
    Then $\varphi\leq \psi$ and \ref{itm:lma:reform_preceqP:2} follows.
\end{proof}

\begin{corollary}\label{cor:PimpliesI}
    Let $\varphi,\psi \in \QPSH(X)$. Assume that $\varphi\preceq_P \psi$. Then $\varphi\preceq_{\mathcal{I}}\psi$.
\end{corollary}
Therefore, the preorder $\preceq_{\mathcal{I}}$ is coarser than $\preceq_P$. Given $\varphi,\psi\in \QPSH(X)$, we therefore have the implications:
\[
\varphi\preceq \psi \implies \varphi\preceq_P \psi \implies \varphi\preceq_{\mathcal{I}} \psi.
\]
\begin{proof}
    Thanks to \cref{lma:reform_preceqP}, we have $\varphi\lor \psi\sim_P \psi$. Take a closed smooth real $(1,1)$-form $\theta$ on $X$ such that $\varphi,\psi\in \PSH(X,\theta)_{>0}$.
    Then by \cref{prop:PandPI}, we have
    \[
        \psi\sim_{\mathcal{I}} P_{\theta}[\psi]= P_{\theta}\left[\varphi\lor \psi\right]\sim_{\mathcal{I}} \varphi\lor \psi.
    \]
    Our assertion follows from the $\mathcal{I}$-version of \cref{lma:reform_preceqP}.
\end{proof}

Next we give a few extra characterizations of the $P$-envelope.
\begin{corollary}\label{cor:Pvarphidef3}
    Let $\theta$ be a closed smooth real $(1,1)$-form on $X$ representing a big cohomology class.
    Assume that $\varphi\in \PSH(X,\theta)_{>0}$. Then
    \[
        \begin{aligned}
        P_{\theta}[\varphi]=&\sup\left\{ \psi\in \PSH(X,\theta): \psi\leq 0, \psi\sim_P \varphi \right\}\\
        =&  \sup\left\{ \psi\in \PSH(X,\theta): \psi\leq 0, \psi\preceq_P \varphi \right\}.
        \end{aligned}
    \]
\end{corollary}
Just for comparison, let us recall a few other characterizations of the $P$-envelope for $\varphi\in \PSH(X,\theta)_{>0}$:
\[
\begin{aligned}
    P_{\theta}[\varphi]=&\uscsup\left\{ \psi\in \PSH(X,\theta):\psi\leq 0,\psi\preceq \varphi \right\}\\
    =&\uscsup\left\{ \psi\in \PSH(X,\theta):\psi\leq 0,\psi\sim \varphi \right\}\\
    =&\uscsup[C\in \mathbb{Z}_{>0}] (\varphi+C)\land V_{\theta}\\
    =& \sup\left\{\psi\in \PSH(X,\theta): \psi\leq 0,\varphi\preceq \psi, \int_X \theta_{\varphi}^n=\int_X \theta_{\psi}^n\right\}.
\end{aligned}
\]
\begin{proof}
    Note that $\psi\sim_P \varphi$ implies that $\psi\in \PSH(X,\theta)_{>0}$ by \cref{prop:Ppresmass}. We observe that
    \[
    \begin{split}
    &\sup\left\{ \psi\in \PSH(X,\theta): \psi\leq 0, \psi\sim_P \varphi \right\}\\=&\sup\left\{ \psi\in \PSH(X,\theta): \psi\leq 0, \varphi\preceq \psi, \psi\sim_P \varphi \right\}
    \end{split}
    \]
    since for any candidate $\psi$ for the first supremum, $\psi \lor (\varphi-\sup_X \varphi)$ is larger than $\psi$, and is a candidate for the second supremum thanks to \cref{lma:reform_preceqP}.
    So the first equality is a direct consequence of \cref{prop:Pequivchar2} and \cref{thm:Pvarphidiffdef}.

    Next we prove the second equality. We only need to show that for any $\psi\in \PSH(X,\theta)$ with $\psi\leq 0$ and $\psi\preceq_P \varphi$, we have $\psi\leq P_{\theta}[\varphi]$.

    By \cref{lma:reform_preceqP} and \cref{ex:Pequiv}, we know that $P_{\theta}[\varphi]\lor \psi\sim_P \varphi$ and $P_{\theta}[\varphi]\lor \psi\leq 0$. It follows from the first equality that $\psi\leq P_{\theta}[\varphi]$.
\end{proof}
Similarly, we have a new characterization of the $\mathcal{I}$-envelope.
\begin{corollary}\label{cor:Ienvelopedef2}
    Let $\theta$ be a closed smooth real $(1,1)$-form on $X$ representing a pseudo-effective cohomology class.
    Assume that $\varphi\in \PSH(X,\theta)$. Then
    \[
        P_{\theta}[\varphi]_{\mathcal{I}}=\sup\left\{ \psi\in \PSH(X,\theta): \psi\leq 0, \psi\preceq_{\mathcal{I}} \varphi \right\}.
    \]
\end{corollary}
\begin{proof}
    It suffices to show that for any $\psi\in \PSH(X,\theta)$ with $\psi\leq 0$ and $\psi\preceq_{\mathcal{I}}\varphi$, we have $\psi\leq P_{\theta}[\varphi]_{\mathcal{I}}$.
    By \cref{lma:reform_preceqP} and \cref{prop:Ienvprojection}, we know that $P_{\theta}[\varphi]_{\mathcal{I}}\lor \psi\sim_{\mathcal{I}}\varphi$. Therefore, $P_{\theta}[\varphi]_{\mathcal{I}}\lor \psi$ is a candidate in the definition of $P_{\theta}[\varphi]_{\mathcal{I}}$ (cf. \eqref{eq:Ienvelopedef}), so
    \[
    \psi\leq P_{\theta}[\varphi]_{\mathcal{I}}\lor \psi\leq P_{\theta}[\varphi]_{\mathcal{I}}.
    \]
    We conclude the proof.
\end{proof}

\begin{proposition}\label{prop:Icomparandenvelope}
    Let $\theta$ be a closed real smooth $(1,1)$-form on $X$ and $\varphi,\psi\in \PSH(X,\theta)$. Then the following are equivalent:
    \begin{enumerate}
        \item\label{itm:prop:Icomparandenvelope:1} $\varphi\preceq_{\mathcal{I}}\psi$;
        \item\label{itm:prop:Icomparandenvelope:2} $P_{\theta}[\varphi]_{\mathcal{I}}\leq P_{\theta}[\psi]_{\mathcal{I}}$.
    \end{enumerate}
\end{proposition}
\begin{proof}
    \ref{itm:prop:Icomparandenvelope:1} $\implies$ \ref{itm:prop:Icomparandenvelope:2}. This follows immediately from \cref{cor:Ienvelopedef2}.

    \ref{itm:prop:Icomparandenvelope:2} $\implies$ \ref{itm:prop:Icomparandenvelope:1}. It follows from \cref{prop:Ienvprojection} that
    \[
        \varphi\sim_{\mathcal{I}} P_{\theta}[\varphi]_{\mathcal{I}}\leq P_{\theta}[\psi]_{\mathcal{I}} \sim_{\mathcal{I}} \psi.
    \]
    We conclude \ref{itm:prop:Icomparandenvelope:1}.
\end{proof}

\begin{example}\label{ex:Pequivnotequivex}
    Let us continue our example \cref{ex:loglogsing}, where $X=\mathbb{P}^1$, $\omega$ is the Fubini--Study metric and $\varphi\in \PSH(X,\omega)$ has log-log singularity at $0$. We have shown that $P_{\omega}[\varphi]=0$ in \eqref{eq:Pomegavarphi0temp1}, so $\varphi\sim_{P} 0$ and hence $\varphi\sim_{\mathcal{I}} 0$. In particular, $P$-equivalence is not equivalent to the equivalence of singularity types.

    On the other hand, consider a potential $\psi\in \PSH(X,\omega)$ with log singularity at $0$, as in \cref{ex:logsingP1}. We know that $\nu(\psi,0)=1$ from the explicit expression \eqref{eq:logsingP1ex}. So $\psi\not\sim_{\mathcal{I}} 0$ and hence $\psi\not\sim_{P} 0$.

    Moreover, $\psi \preceq_{P} \varphi$ and hence $\psi\preceq_{\mathcal{I}}\varphi$.
\end{example}

We give an example showing that $P$-equivalence is not equivalent to $\mathcal{I}$-equivalence.
\begin{example}\label{ex:BBJ}
    Let $X=\mathbb{P}^1$ and $\omega$ be the Fubini--Study metric. Let $K\subseteq \mathbb{P}^1$ be a polar Cantor set carrying an atom free probability measure $\mu$ supported on $K$ (see \cite[Page~31]{Car83} and \cite[Proposition~4.2.4]{FPT}). Write $\mu=\omega+\ddc\varphi$ for some $\omega$-subharmonic function $\varphi$. Since $\mu$ is atom free, we know that all Lelong numbers of $\varphi$ are $0$. On the other hand, $\varphi$ has $0$ non-pluripolar mass since $K$ is pluripolar.

    Then observe that $\varphi\sim_{\mathcal{I}} 0$ while $\varphi\not\sim_P 0$.
\end{example}

For later use, we give the following definition.
\begin{definition}\label{def:elemetric}
    Let $L$ be a pseudo-effective line bundle on $X$. An \emph{elementary metric}\index{elementary metric} on $L$ is a psh metric $h$ on $L$ such that there is a generalized Fubini--Study metric $h'$ on $L$ such that
    \[
    \ddc h\sim_{P} \ddc h'.
    \]
    The set of elementary metrics on $L$ is denoted by $\Ele(L)$\index{$\Ele(L)$}.

    We also say $\ddc h$ is elementary. If we have fixed a Hermitian metric $h_0$ on $L$, and if we represent $h$ as $h_0\exp(-\varphi)$, we also say the quasi-psh function $\varphi$ is elementary.
\end{definition}
Recall that the generalized Fubini--Study metrics are defined in \cref{def:gFSmetric}.

\subsection{Properties of the preorders}
\label{subsec:preorder_prop}

Now we state a more natural version of the monotonicity theorem \cref{thm:mono}.
\begin{proposition}[Monotonicity theorem]\label{prop:mono2}
    Let $\theta_1,\ldots,\theta_n$ be closed real smooth $(1,1)$-forms on $X$.
    Let $\varphi_i,\psi_i\in \PSH(X,\theta_i)$ for $i=1,\ldots,n$. Assume that $\varphi_i\preceq_P \psi_i$ for each $i$. Then
    \[
        \int_X \theta_{1,\varphi_1}\wedge \cdots \wedge  \theta_{n,\varphi_n}\leq \int_X \theta_{1,\psi_1}\wedge \cdots \wedge  \theta_{n,\psi_n}.
    \]
    In other words, the non-pluripolar mass is a decreasing quantity with respect to the $P$-singularity types of the potentials.
\end{proposition}
\begin{proof}
    Fix a Kähler form $\omega$ on $X$. For each $i=1,\ldots,n$, since $\varphi_i\preceq_P\psi_i$, we have
    \[
        P_{\theta_i+\epsilon\omega}[\varphi_i]\leq P_{\theta_i+\epsilon\omega}[\psi_i]
    \]
    for all $\epsilon>0$. Therefore, by \cref{prop:Ppresmass} and our previous monotonicity theorem \cref{thm:mono}, we have
    \[
        \int_X (\theta_1+\epsilon\omega)_{\varphi_1}\wedge \cdots \wedge (\theta_n+\epsilon\omega)_{\varphi_n}\leq \int_X (\theta_1+\epsilon\omega)_{\psi_1}\wedge \cdots \wedge (\theta_n+\epsilon\omega)_{\psi_n}.
    \]
    Letting $\epsilon\to 0+$ and applying the multi-linearity of the non-pluripolar products \itemref{itm:prop:npp1:6}, we find the desired inequality.
\end{proof}

The $P$-preorder allows us to derive a natural generalization of \cref{lma:pathoenvelope}.
\begin{theorem}\label{thm:pathoenvelope}
    Let $\varphi,\psi\in \PSH(X,\theta)_{>0}$, $\varphi\preceq_P \psi$. Then for any
    \begin{equation}\label{eq:arangetemp_reproduce}
        a\in \left( 1, \left(\frac{\int_X \theta_{\psi}^n}{\int_X \theta_{\psi}^n-\int_X \theta_{\varphi}^n}\right)^{1/n} \right),
    \end{equation}
    the function
    \[
        \eta\coloneqq P_{\theta}\Bigl(a\varphi+(1-a)\psi\Bigr)
    \]
    lies in $\PSH(X,\theta)_{>0}$, satisfies $\eta\preceq_P \varphi$ and
    \begin{equation}\label{eq:path_env_reproduce}
        a^{-1}\eta+\left(1-a^{-1}\right)\psi\leq \varphi.
    \end{equation}

    If in addition, $\varphi\sim_P \psi$, then $\eta\sim_P \varphi$ as well.
\end{theorem}
In \eqref{eq:arangetemp_reproduce}, the fraction is understood as $\infty$ if either $\int_X \theta_{\psi}^n=\int_X \theta_{\varphi}^n$ or $n=0$, exactly as in \cref{lma:pathoenvelope}.
\begin{proof}
    By our assumption and \eqref{eq:varphipreceqP}, we have $\varphi\preceq P_{\theta}[\varphi]\leq P_{\theta}[\psi]$.
    Due to \cref{prop:Ppresmass}, the interval \eqref{eq:arangetemp_reproduce} remains unchanged after replacing $\psi$ by $P_{\theta}[\psi]$. Fix $a$ in this interval.
    So by \cref{lma:pathoenvelope} and the discussions afterwards,
    \[
        \eta'\coloneqq  P_{\theta}\Bigl(a\varphi+(1-a)P_{\theta}[\psi]\Bigr)\in \PSH(X,\theta)_{>0}.
    \]
    We first check that $\eta$ is a genuine $\theta$-psh function. Choose a closed non-pluripolar set $K\subseteq X$ on which $\varphi$ and $\psi$ are bounded. Then the obstacle $a\varphi+(1-a)\psi$ is bounded above on $K$.
    Every candidate in the definition of $\eta$ is therefore uniformly bounded above on $X$ by \cref{rmk:L1cpt}, so $\eta\not\equiv \infty$ and hence $\eta\in \PSH(X,\theta)\cup\{-\infty\}$.
    On the other hand, since $\psi\preceq P_{\theta}[\psi]$, we can find $C>0$ such that $\psi\leq P_{\theta}[\psi]+C$. Then
    \[
        \eta'-(a-1)C\leq a\varphi+(1-a)\psi
    \]
    quasi-everywhere, hence $\eta\geq \eta'-(a-1)C$. Thus $\eta'\preceq \eta$, and the monotonicity theorem \cref{thm:mono} gives $\eta\in \PSH(X,\theta)_{>0}$.

    By \eqref{eq:Pthetafleqfqe},
    \[
        \eta\leq a\varphi+(1-a)\psi
    \]
    quasi-everywhere. Equivalently, \eqref{eq:path_env_reproduce} holds quasi-everywhere. Since both sides of \eqref{eq:path_env_reproduce} are $\theta$-psh, \cref{prop:pshfuncdetdense} upgrades the inequality to an everywhere inequality.

    It remains to prove the assertion about the $P$-order.
    Using \cref{prop:Pconc}, we get
    \[
        a^{-1} P_{\theta}[\eta]+\left(1-a^{-1}\right)P_{\theta}[\psi] \leq P_{\theta}[\varphi].
    \]
    By the assumption $\varphi\preceq_P\psi$, we also have $P_{\theta}[\varphi]\leq P_{\theta}[\psi]$. Therefore,
    \[
        a^{-1} P_{\theta}[\eta]\leq P_{\theta}[\varphi]-\left(1-a^{-1}\right)P_{\theta}[\psi]\leq a^{-1} P_{\theta}[\varphi]
    \]
    quasi-everywhere and hence everywhere by \cref{prop:pshfuncdetdense}. It follows that $P_{\theta}[\eta]\leq P_{\theta}[\varphi]$, namely $\eta\preceq_P\varphi$.

    Finally assume that $\varphi\sim_P\psi$. Then $\int_X\theta_{P_{\theta}[\psi]}^n=\int_X\theta_{\varphi}^n$ by \cref{prop:Ppresmass}.
    Applying \cref{cor:pathoenvelopeeqmass} gives
    \[
        \int_X\theta_{\eta'}^n=\int_X\theta_{\varphi}^n.
    \]
    Since $\eta'\preceq\eta$, the monotonicity theorem \cref{thm:mono} gives
    \[
        \int_X\theta_{\eta}^n\geq \int_X\theta_{\varphi}^n.
    \]
    The reverse inequality follows from $\eta\preceq_P\varphi$ and \cref{prop:mono2}. Hence $\eta\sim_P\varphi$ by \cref{prop:Pequivchar2}.
\end{proof}

Next we show that the $P$ and $\mathcal{I}$-preorders are preserved by some natural operations.

\begin{lemma}\label{lma:Ipartorderinv}
Let $\pi\colon Y\rightarrow X$ be a proper bimeromorphic morphism from a Kähler manifold $Y$.
    Given two quasi-plurisubharmonic functions $\varphi,\psi$ on $X$, then the following are equivalent:
    \begin{enumerate}
        \item\label{itm:lma:Ipartorderinv:1} $\varphi\preceq_{P} \psi$;
        \item\label{itm:lma:Ipartorderinv:2} $\pi^*\varphi\preceq_{P} \pi^*\psi$.
    \end{enumerate}
    The same holds with $\mathcal{I}$ in place of $P$.
\end{lemma}
\begin{proof}
   In the $P$-case, this follows from \cref{prop:Penvbimero}, while in the $\mathcal{I}$-case, this follows from \cref{prop:Ienvelopebimero}.
\end{proof}

\begin{proposition}\label{prop:Ppartialsum}
    Let $\varphi,\psi,\varphi',\psi'\in \QPSH(X)$. Assume that
    \[
    \varphi\preceq_P \psi,\quad \varphi'\preceq_P \psi'.
    \]
    Then
    \[
    \varphi+\varphi'\preceq_P \psi+\psi'.
    \]
    The same holds with $\preceq_{\mathcal{I}}$ in place of $\preceq_P$.
\end{proposition}
\begin{proof}
    Take a Kähler form $\omega$ on $X$ such that $\varphi,\psi,\varphi',\psi'\in \PSH(X,\omega)_{>0}$. The statement for $\preceq_{\mathcal{I}}$ is a simple consequence of \cref{prop:Lelongmax}. We only need to handle the case of $\preceq_{P}$.

    \textbf{Step~1}.
    We first show that
    \[
    P_{\omega}[\varphi]+P_{\omega}[\varphi']\sim_P \varphi+\varphi'.
    \]
    In fact, we clearly have
    \[
    P_{\omega}[\varphi]+P_{\omega}[\varphi']\succeq \varphi+\varphi'.
    \]
    So by \cref{prop:Pequivchar2}, it suffices to show that they have the same mass. We compute
    \[
    \begin{aligned}
    &\int_X \Bigl(2\omega+\ddc P_{\omega}[\varphi]+\ddc P_{\omega}[\varphi']\Bigr)^n\\
    =& \sum_{j=0}^n \binom{n}{j} \int_X \Bigl(\omega+\ddc P_{\omega}[\varphi]\Bigr)^j\wedge \Bigl(\omega+\ddc P_{\omega}[\varphi']\Bigr)^{n-j}\\
    =& \sum_{j=0}^n \binom{n}{j} \int_X \omega_{\varphi}^j\wedge \omega_{\varphi'}^{n-j}\\
    =& \int_X \left(2\omega+\ddc\varphi+\ddc\varphi'\right)^n,
    \end{aligned}
    \]
    where we applied \cref{prop:Ppresmass} to the mixed products on the third line.

    \textbf{Step~2}. By Step~1, we may assume that $\varphi,\psi,\varphi',\psi'$ are all model potentials, after replacing them by their $P_{\omega}$-envelopes. So $\varphi\leq \psi$ and $\varphi'\leq \psi'$. Our assertion follows.
\end{proof}
\begin{proposition}\label{prop:Ppartialsup}
    Let $(\varphi_i)_{i\in I}$, $(\psi_i)_{i\in I}$ be uniformly bounded from above non-empty families in $\QPSH(X)$. Assume that there exists a closed smooth real $(1,1)$-form $\theta$ such that $\varphi_i,\psi_i\in \PSH(X,\theta)$ and $\varphi_i\preceq_P \psi_i$ for all $i\in I$. Then
    \[
    \uscsup[i\in I] \varphi_i\preceq_P \uscsup[i\in I]\psi_i.
    \]

    The same holds with $\preceq_{\mathcal{I}}$ in place of $\preceq_P$.
\end{proposition}
\begin{proof}
    By adding a Kähler form to $\theta$, we may assume that $\varphi_i,\psi_i\in \PSH(X,\theta)_{>0}$ for all $i\in I$. The statement for $\preceq_{\mathcal{I}}$ is a simple consequence of \cref{cor:supsLelong}, we only have to consider the statement for $\preceq_P$.

    \textbf{Step~1}. We first handle the case where $I$ is a directed set and $(\varphi_i)_{i\in I}$ and $(\psi_i)_{i\in I}$ are increasing nets.

   In this case,  our assertion follows simply from \cref{prop:incnetmodel}.

    \textbf{Step~2}. We handle the case where $I$ is finite. We may assume that $I=\{0,1\}$.
    It suffices to show that
    \[
        P_{\theta}[\varphi_0]\lor P_{\theta}[\varphi_1]\sim_P \varphi_0\lor\varphi_1,
    \]
    which follows from \cref{lma:Plor1}.

    \textbf{Step~3}. The general case can be reduced to the two cases handled in Step~1 and Step~2. More precisely,
    by applying \cref{prop:Choquet} to the two families and taking the union of the resulting countable subsets, we can find a countable subset $J\subseteq I$ such that
    \[
    \uscsup[j\in J] \varphi_j=\uscsup[i\in I] \varphi_i,\quad \uscsup[j\in J] \psi_j=\uscsup[i\in I] \psi_i.
    \]
    We may replace $I$ by $J$ and assume that $I$ is countable. If $I$ is finite, Step~2 applies directly. Otherwise, write $I=\mathbb{Z}_{>0}$ and set
    \[
        \tilde{\varphi}_m=\varphi_1\lor\cdots\lor\varphi_m,\quad
        \tilde{\psi}_m=\psi_1\lor\cdots\lor\psi_m.
    \]
    By Step~2, $\tilde{\varphi}_m\preceq_P\tilde{\psi}_m$ for every $m$. The two new sequences are increasing and have the same upper envelopes as the original families, so Step~1 applies.

\end{proof}

We now refine \cref{lma:landcorrectmass} using the $P$-equivalence relation.
\begin{proposition}\label{prop:rooftopprePequiv}
    Let $\varphi,\psi,\varphi',\psi'\in \PSH(X,\theta)$ for some closed smooth real $(1,1)$-form $\theta$ on $X$. Assume that
    \[
    \varphi\sim_P \varphi',\quad \psi\sim_P\psi',\quad  \varphi'\land \psi'\in \PSH(X,\theta)_{>0}.
    \]
    Then
    \[
    \varphi\land \psi\in \PSH(X,\theta)_{>0},\quad \varphi\land \psi\sim_P \varphi'\land \psi'.
    \]
\end{proposition}
\begin{proof}
    We first observe that $\varphi,\varphi',\psi,\psi'\in \PSH(X,\theta)_{>0}$ by assumption. Let
    \[
    \phi\coloneqq P_{\theta}[\varphi]=P_{\theta}[\varphi'],\quad \gamma=P_{\theta}[\psi]=P_{\theta}[\psi'].
    \]
    Since $\varphi'\preceq\phi$ and $\psi'\preceq\gamma$, after adding a common constant to $\varphi'$ and $\psi'$, we may assume that $\varphi'\leq\phi$ and $\psi'\leq\gamma$. Then $\varphi'\land\psi'\leq\phi\land\gamma$, so the assumption $\varphi'\land\psi'\in \PSH(X,\theta)_{>0}$ implies $\phi\land \gamma \in \PSH(X,\theta)_{>0}$ by monotonicity. It follows from \cref{lma:landcorrectmass} that
    \[
    \int_X \theta_{\varphi'\land \psi'}^n=\int_X \theta_{\phi\land \gamma}^n.
    \]
    Next, we apply \cref{lma:landcorrectmass} again to conclude that $\varphi\land \psi\in \PSH(X,\theta)$ and
    \[
    \int_X \theta_{\varphi\land \psi}^n=\int_X \theta_{\phi\land \gamma}^n=\int_X \theta_{\varphi'\land \psi'}^n>0.
    \]
    Since $\varphi\land\psi\preceq \phi\land\gamma$ and the two potentials have the same positive mass, \cref{thm:Pvarphidiffdef} gives $P_{\theta}[\varphi\land\psi]=\phi\land\gamma$. The same argument gives $P_{\theta}[\varphi'\land\psi']=\phi\land\gamma$, hence $\varphi\land\psi\sim_P\varphi'\land\psi'$.
\end{proof}

We can formulate a generalization of \cref{prop:envrelfullmass} using the $P$-partial order.
\begin{corollary}\label{cor:P_par_envelope_rel}
    Let $\varphi,\psi\in \PSH(X,\theta)_{>0}$. Assume that $\psi\preceq_P \varphi$. Then for any $C>0$, $\psi\land (\varphi+C)\in \PSH(X,\theta)_{>0}$ and
    \begin{equation}\label{eq:supvarphipsipC}
        P_{\theta}[\varphi](\psi)=\uscsup[C>0] \psi\land (\varphi+C)=\psi.
    \end{equation}
\end{corollary}
We refer to \cref{def:Pthetavarphif} for the definition of $P_{\theta}[\varphi](\psi)$.
\begin{proof}
    Since $\psi\preceq_P\varphi$, by \cref{prop:rooftopprePequiv}, we get
    \[
        \psi\land(\varphi+C)\in \PSH(X,\theta)_{>0}
    \]
    and $\psi\land(\varphi+C)\sim_P\psi$ for every $C>0$.

    Let
    \[
        \eta\coloneqq \uscsup[C>0]\psi\land(\varphi+C).
    \]
    Then $\eta\sim_P \psi$.
    By \cref{cor:Pvarphisupport}, we have
    \[
        \theta_\eta^n\leq \mathds{1}_{\{\eta=\psi\}}\,\theta_\psi^n.
    \]
    The domination principle \cref{thm:dom_prin_1} gives $\eta\geq \psi$. The reverse inequality is trivial, so $\eta=\psi$.
\end{proof}

The $P$-equivalence relation characterizes when a subgeodesic exists. See \cref{def:subgeod} for the notion of subgeodesics.
\begin{theorem}\label{thm:geodexi}
    Let $\varphi_0,\varphi_1\in \PSH(X,\theta)_{>0}$. Then the following are equivalent:
    \begin{enumerate}
        \item\label{itm:thm:geodexi:1} There is a subgeodesic from $\varphi_0$ to $\varphi_1$;
        \item\label{itm:thm:geodexi:2} $\varphi_0\sim_P\varphi_1$.
    \end{enumerate}
\end{theorem}
\begin{proof}
    \ref{itm:thm:geodexi:2} $\implies$ \ref{itm:thm:geodexi:1}. This follows from \cref{prop:perronenvissubgeod}.

    \ref{itm:thm:geodexi:1} $\implies$ \ref{itm:thm:geodexi:2}.
    Let $(\varphi_t)_{t\in (0,1)}$ be a subgeodesic from $\varphi_0$ to $\varphi_1$.

    \textbf{Step~1}. We assume that $\varphi_0\geq \varphi_1$.

    We first replace the given subgeodesic by the geodesic from $\varphi_0$ to $\varphi_1$. Note that the geodesic is necessarily a subgeodesic as well by \cref{lma:geodsubgeod}.

    Observe that $t\mapsto\varphi_t$ is decreasing due to \cref{lma:geoddecreasing}.

    Let $\varphi_t=\varphi_1$ for all $t>1$. Then by the gluing lemma \cref{lma:subgeodgluing}, we find that $(\varphi_t)_{t\geq 0}$ is a subgeodesic ray.

    Next, we consider the Legendre transform
    \[
    \Gamma_{\tau}\coloneqq \inf_{t\geq 0}(\varphi_t-t\tau),\quad \tau\in \mathbb{R}.
    \]
    It follows from Kiselman's principle \cref{prop:Kis} that $\Gamma_{\tau}\in \PSH(X,\theta)\cup \{-\infty\}$. Note that for $\tau>0$, we clearly have $\Gamma_{\tau}\equiv -\infty$. On the other hand, for $\tau\leq 0$,
    \[
    \Gamma_{\tau}=\inf_{t\in [0,1]}(\varphi_t-t\tau)\in \PSH(X,\theta).
    \]

    By Legendre inversion (cf. \cref{cor:Leginvpartdef}), for $t>0$,
    \[
    \varphi_t=\sup_{\tau<0} \left(\Gamma_{\tau}+t\tau\right).
    \]
    Fix a K\"ahler form $\omega$ on $X$.
    It follows from \cref{prop:Ppartialsup} that for each $t>0$,
    \[
    \varphi_t\sim_P \uscsup[\tau<0] P_{\theta+\omega}[\Gamma_{\tau}].
    \]
    The right-hand side is independent of $t$.
    Here by adding $\omega$, we no longer have to worry about the possibility where $\Gamma_{\tau}$ has vanishing mass.

    Write
    \[
    \varphi_0=\uscsup[t\in (0,1)]\varphi_t.
    \]
    By \cref{prop:Ppartialsup} again, we find that
    \[
    \varphi_0\sim_P \uscsup[\tau<0] P_{\theta+\omega}[\Gamma_{\tau}]
    \]
    as well. So $\varphi_0\sim_P \varphi_1$.

    \textbf{Step~2}. We prove the general case.

    Observe that $(\varphi_t\lor\varphi_1)_{t\in (0,1)}$ is a subgeodesic from $\varphi_0\lor \varphi_1$ to $\varphi_1$. By Step~1,
    \[
    \varphi_0\lor \varphi_1\sim_P \varphi_1.
    \]
    Hence $\varphi_0\preceq_P \varphi_1$ by \cref{lma:reform_preceqP}. The converse is proved similarly. Hence \ref{itm:thm:geodexi:2} follows.

\end{proof}

Restating the \ref{itm:thm:geodexi:1} $\implies$ \ref{itm:thm:geodexi:2} direction in the theorem in terms of the complexification, we find the following interesting result.

Let $S=\{z\in \mathbb{C}:0< \Real z <1\}$. We write $p_1\colon X\times S\rightarrow X$ for the natural projection.
\begin{corollary}\label{cor:comp:L653}
    Let $\Phi\in \PSH(X\times S,p_1^*\theta)$. Assume that for any $c\in \mathbb{R}$, $x\in X$ and $z\in S$, we have
    \[
    \Phi(x,z)=\Phi(x,z+\mathrm{i}c).
    \]
    Then $\int_X (\theta+\ddc \Phi_z)^n$ is independent of $z\in S$, where $\Phi_z \in \PSH(X,\theta)$ is given by $\Phi_z(x)=\Phi(x,z)$.
\end{corollary}

\begin{proof}
    Fix a Kähler form $\omega$ on $X$. We may regard $\Phi$ as an element of $\PSH(X\times S,p_1^*(\theta+\omega))$. For any $0<a<b<1$, the restriction $(\Phi_t)_{t\in [a,b]}$ is then a subgeodesic from $\Phi_a$ to $\Phi_b$ in the class $\theta+\omega$. The slices have positive mass in this class, so \cref{thm:geodexi} gives $\Phi_a\sim_P\Phi_b$. By \cref{prop:mono2}, this directly implies
    \[
        \int_X \theta_{\Phi_a}^n=\int_X \theta_{\Phi_b}^n.
    \]
    Since $\Phi$ is invariant under imaginary translations, $\Phi_z$ depends only on $\Real z$, so the mass is independent of $z\in S$.
\end{proof}

\section[The \texorpdfstring{$\lowercase{d}_S$}{dS}-pseudometric]{The \texorpdfstring{$d_S$}{dS}-pseudometric}\label{sec:dsdef}
Let $X$ be a connected compact Kähler manifold of dimension $n$ and $\theta$ be a closed real smooth $(1,1)$-form on $X$ representing a big cohomology class. The goal of this section is to study a pseudometric on the space $\PSH(X,\theta)$.

\subsection{The definition of the \texorpdfstring{$d_S$}{dS}-pseudometric}
\label{subsec:dS_def}
Recall that for any $\varphi\in \PSH(X,\theta)$, the geodesic ray $\ell^{\varphi}\in \mathcal{R}^1(X,\theta)$ is defined in \cref{ex:rayasspsh}.
\begin{definition}\label{def:dS}
    For $\varphi,\psi\in \PSH(X,\theta)$, we define
    \[
        d_S(\varphi,\psi)\coloneqq d_1(\ell^{\varphi},\ell^{\psi}).
    \]
    \index{$d_S$}
    When we want to be more specific, we write $d_{S,\theta}$ instead of $d_S$.
\end{definition}
The $d_1$ distance of geodesic rays is defined in \cref{def:d1rays}.

\begin{proposition}\label{prop:comp:L686}
    The function $d_S$ defined in \cref{def:dS} is a pseudometric on $\PSH(X,\theta)$.
\end{proposition}
\begin{proof}
    This follows immediately from \cref{thm:d1raycomplete}.
\end{proof}
When studying a pseudometric, the first thing is to understand when the distance between two elements vanishes.

We first prove a preparation:
\begin{lemma}\label{lma:dSalmostriang}
    Let $\varphi,\psi\in \PSH(X,\theta)$. Then
    \[
        d_S(\varphi,\psi)\leq d_S(\varphi,\varphi\lor \psi)+d_S(\psi,\varphi\lor \psi) \leq C_n d_S(\varphi,\psi),
    \]
    where $C_n=3(n+1)2^{n+2}$.
\end{lemma}
We shall use the notations introduced in \cref{ex:rayasspsh}.
\begin{proof}
    We claim that
    \begin{equation}\label{eq:elllorsingtype}
        \ell^{\varphi}\lor \ell^{\psi}=\ell^{\varphi\lor \psi}.
    \end{equation}
    Recall that $\lor$ is defined in \cref{def:lorray1}.
    Note that this assertion implies our desired inequality by \cref{lma:d1rayineq}.

    In proving this assertion, we may assume that $\varphi,\psi\leq 0$ since
    \[
    \ell^{\varphi+C}=\ell^{\varphi},\quad \ell^{\psi+C}=\ell^{\psi},\quad \ell^{(\varphi+C)\lor (\psi+C)}=\ell^{\varphi\lor \psi}
    \]
    for any $C\in \mathbb{R}$.

    In fact, it is clear that
    \[
        \ell^{\varphi}\leq   \ell^{\varphi\lor \psi},\quad \ell^{\psi}\leq   \ell^{\varphi\lor \psi},
    \]
    so the $\leq$ direction in \eqref{eq:elllorsingtype} holds.

    Conversely, if $\ell'\in \mathcal{R}^1(X,\theta)$ and $\ell'\geq \ell^{\varphi}\lor \ell^{\psi}$, then for each $t\geq 0$,
    \[
        \ell'_t\geq \Bigl((V_{\theta}-t)\lor \varphi\Bigr) \lor \Bigl((V_{\theta}-t)\lor \psi\Bigr)=(V_{\theta}-t)\lor (\varphi\lor \psi).
    \]
    See \eqref{eq:ellphitexp} for the first inequality.
    Therefore, thanks to the definition of the Perron envelopes (cf. \eqref{eq:Perron2}), we have
    \[
    \ell'_s\geq \ell^{\varphi\lor \psi,t}_s
    \]
    for any $0\leq s\leq t$.
    It follows from \eqref{eq:ellphitexp} that $\ell'_s\geq \ell^{\varphi\lor \psi}_s$ for any $s\geq 0$.
\end{proof}

\begin{proposition}\label{prop:ds0char}
    Let $\varphi,\psi\in \PSH(X,\theta)$. Then the following are equivalent:
    \begin{enumerate}
        \item\label{itm:prop:ds0char:1} $\varphi\sim_P \psi$;
        \item\label{itm:prop:ds0char:2} $d_S(\varphi,\psi)=0$.
    \end{enumerate}

    In particular, $d_S(\varphi,P_{\theta}[\varphi])=0$ for all $\varphi\in \PSH(X,\theta)_{>0}$.
\end{proposition}
\begin{proof}
    By \cref{lma:reform_preceqP}, we have $\varphi\sim_P \psi$ if and only if $\varphi\sim_P \varphi\lor \psi$ and $\psi\sim_P \varphi\lor \psi$. By \cref{lma:dSalmostriang}, $d_S(\varphi,\psi)=0$ if and only if $d_S(\varphi,\varphi\lor \psi)=0$ and $d_S(\psi,\varphi\lor \psi)=0$. So it suffices to prove the assertion when $\varphi\leq \psi$. Assuming this, by \cref{prop:d1geod_diff_E}, \ref{itm:prop:ds0char:2} holds if and only if
    \[
        \mathbf{E}\left(\ell^{\varphi}\right)=  \mathbf{E}\left(\ell^{\psi}\right),
    \]
    where $\mathbf{E}$ is introduced in \cref{def:radialMAenergy2}.
    But by \eqref{eq:Elphi}, this holds if and only if
    \[
        \sum_{j=0}^n \int_X \theta_{\varphi}^j\wedge \theta_{V_{\theta}}^{n-j}=\sum_{j=0}^n \int_X \theta_{\psi}^j\wedge \theta_{V_{\theta}}^{n-j}.
    \]
    Thanks to the monotonicity theorem \cref{thm:mono}, this holds if and only if for all $j=0,\ldots,n$,
    \[
        \int_X \theta_{\varphi}^j\wedge \theta_{V_{\theta}}^{n-j}=\int_X \theta_{\psi}^j\wedge \theta_{V_{\theta}}^{n-j},
    \]
    which is equivalent to \ref{itm:prop:ds0char:1} by \cref{prop:Pequivchar2}.
\end{proof}

\begin{lemma}\label{lma:varphileqpsi_metric}
    Suppose that $\varphi,\psi\in \PSH(X,\theta)$ and $\varphi\preceq_P \psi$. Then
    \[
    d_S(\varphi,\psi)=\frac{1}{n+1}\sum_{j=0}^n \left(\int_X \theta_{\psi}^j\wedge \theta_{V_{\theta}}^{n-j}-\int_X \theta_{\varphi}^j\wedge \theta_{V_{\theta}}^{n-j}\right).
    \]
\end{lemma}
\begin{proof}
    Set $\chi\coloneqq\varphi\lor\psi$. By \cref{lma:reform_preceqP}, we have $\chi\sim_P\psi$. Hence \cref{prop:ds0char} gives $d_S(\chi,\psi)=0$, and the pseudometric property gives
    \[
        d_S(\varphi,\psi)=d_S(\varphi,\chi).
    \]
    Since $\varphi\leq\chi$, \eqref{eq:Elphi} and \cref{prop:d1geod_diff_E} give
    \[
    d_S(\varphi,\chi)=\frac{1}{n+1}\sum_{j=0}^n \left(\int_X \theta_{\chi}^j\wedge \theta_{V_{\theta}}^{n-j}-\int_X \theta_{\varphi}^j\wedge \theta_{V_{\theta}}^{n-j}\right).
    \]
    Finally, \cref{prop:Pequivchar2} applied to $\psi\leq\chi$ gives
    \[
        \int_X \theta_{\chi}^j\wedge \theta_{V_{\theta}}^{n-j}=\int_X \theta_{\psi}^j\wedge \theta_{V_{\theta}}^{n-j}
    \]
    for all $j=0,\ldots,n$.
\end{proof}
\begin{corollary}\label{cor:dsthreeterm}
    Suppose that $\varphi,\psi,\eta\in \PSH(X,\theta)$ and $\varphi\preceq_P\psi\preceq_P \eta$. Then
    \[
        d_S(\varphi,\eta)\geq d_S(\varphi,\psi),\quad d_S(\varphi,\eta)\geq d_S(\psi,\eta).
    \]
\end{corollary}
\begin{proof}
    This is an immediate consequence of \cref{lma:varphileqpsi_metric} and the monotonicity theorem \cref{prop:mono2}.
\end{proof}

\begin{corollary}\label{cor:dsmetricdoubleineq}
    For any $\varphi,\psi\in \PSH(X,\theta)$, we have
    \begin{equation}\label{eq:ds_biineq}
        \begin{split}
            d_S(\varphi,\psi)\leq & \frac{1}{n+1}\sum_{j=0}^n \left( 2\int_X \theta_{\varphi\lor \psi}^j\wedge \theta_{V_{\theta}}^{n-j}-\int_X \theta_{\varphi}^j\wedge \theta_{V_{\theta}}^{n-j}-\int_X \theta_{\psi}^j\wedge \theta_{V_{\theta}}^{n-j} \right) \\
            \leq & C_n d_S(\varphi,\psi),
        \end{split}
\end{equation}
where $C_n=3(n+1)2^{n+2}$.

In particular, if $(\varphi_i)_{i\in I}$ is a net in $\PSH(X,\theta)$ with $d_S$-limit $\varphi$, then for each $j=0,\ldots,n$,
\[
    \lim_{i\in I}\int_X \theta_{\varphi_i}^j\wedge \theta_{V_{\theta}}^{n-j}=\int_X \theta_{\varphi}^j\wedge \theta_{V_{\theta}}^{n-j}=\lim_{i\in I}\int_X \theta_{\varphi_i\lor \varphi}^j\wedge \theta_{V_{\theta}}^{n-j}.
\]
\end{corollary}
\begin{proof}
    The estimates \eqref{eq:ds_biineq} follow from the combination of \cref{lma:varphileqpsi_metric} and \cref{lma:dSalmostriang}.

    Suppose that $\varphi_i\xrightarrow{d_S}\varphi$. Then $\varphi_i\lor \varphi\xrightarrow{d_S}\varphi$ by \cref{lma:dSalmostriang}. Therefore, \cref{thm:mono} and \cref{lma:varphileqpsi_metric} imply that
    \[
    \lim_{i\in I}\int_X \theta_{\varphi_i\lor \varphi}^j\wedge \theta_{V_{\theta}}^{n-j}= \int_X \theta_{\varphi}^j\wedge \theta_{V_{\theta}}^{n-j}
    \]
    for any $j=0,\ldots,n$.
    The last assertion now follows from \eqref{eq:ds_biineq} and \cref{thm:mono}.
\end{proof}

\begin{corollary}\label{cor:dSconv_bimero}
    Let $\pi\colon Y\rightarrow X$ be a proper bimeromorphic morphism from a Kähler manifold. Given a net $(\varphi_i)_{i\in I}$ in $\PSH(X,\theta)$ and $\varphi\in \PSH(X,\theta)$, the following are equivalent:
    \begin{enumerate}
        \item $\varphi_i\xrightarrow{d_S}\varphi$;
        \item $\pi^*\varphi_i\xrightarrow{d_S}\pi^*\varphi$.
    \end{enumerate}
\end{corollary}
\begin{proof}
    It suffices to observe that the middle term in \eqref{eq:ds_biineq} is invariant under $\pi^*$.
\end{proof}

\begin{corollary}\label{cor:incnetdSconv}
    Suppose that $(\varphi_i)_{i\in I}$ is an increasing net in $\PSH(X,\theta)$, uniformly bounded from above. Then
    \[
        \varphi_i\xrightarrow{d_S} \uscsup[i\in I]\varphi_i.
    \]
\end{corollary}
If the $\varphi_i$'s are all model potentials in $\PSH(X,\theta)_{>0}$, then so is $\uscsup[i\in I]\varphi_i$, as we have seen in \cref{prop:incnetmodel}.
\begin{proof}
    Write $\varphi=\uscsup[j\in I]\varphi_j$. Recall that by \cref{prop:pshfunction_closedseq}, $\varphi\in \PSH(X,\theta)$.
    By \cref{lma:varphileqpsi_metric}, it suffices to show that for each $k=0,\ldots,n$, we have
    \[
        \lim_{j\in I}\int_X \theta_{\varphi_j}^k\wedge \theta_{V_{\theta}}^{n-k} =  \int_X \theta_{\varphi}^k\wedge \theta_{V_{\theta}}^{n-k}.
    \]
    The latter follows from \cref{cor:incseqnppcont}.
\end{proof}

\begin{corollary}\label{cor:massdiffdS}
    Let $\varphi,\psi\in \PSH(X,\theta)$. Then
    \[
    \left| \int_X \theta_{\varphi}^n-\int_X \theta_{\psi}^n \right|\leq D_n d_S(\varphi,\psi),
    \]
    where $D_n=3(n+1)C_n$ with $C_n$ being the same constant as in \cref{lma:dSalmostriang}.
\end{corollary}
In other words, the non-pluripolar mass is uniformly continuous with respect to the $d_S$-pseudometric.
\begin{proof}
    We compute
    \[
    \begin{aligned}
        \left| \int_X \theta_{\varphi}^n-\int_X \theta_{\psi}^n \right|\leq & \left|2\int_X \theta_{\varphi\lor \psi}^n-\int_X \theta_{\varphi}^n-\int_X \theta_{\psi}^n\right| + 2\left| \int_X \theta_{\varphi\lor \psi}^n-\int_X \theta_{\varphi}^n\right|\\
        \leq &  (n+1)C_n d_S(\varphi,\psi)+2(n+1) d_S(\varphi,\varphi\lor \psi)\\
        \leq & (n+1)C_n d_S(\varphi,\psi)+2(n+1)C_n d_S(\varphi,\psi),
    \end{aligned}
    \]
    where the first line is just the triangle inequality, the second line follows from \cref{cor:dsmetricdoubleineq} and the third line follows from \cref{lma:dSalmostriang}.
\end{proof}

By contrast to the increasing case, for nets decreasing with respect to the $P$-preorder, the situation is different:
\begin{corollary}\label{cor:decnetdS}
    Suppose that $(\varphi_i)_{i\in I}$ is a net in $\PSH(X,\theta)$ such that
    \[
        \varphi_j\preceq_P\varphi_i
    \]
    whenever $j\geq i$. Let $\varphi\in\PSH(X,\theta)$ be such that $\varphi\preceq_P\varphi_i$ for all $i\in I$. Then the following are equivalent:
    \begin{enumerate}
        \item\label{itm:cor:decnetdS:1} We have
            \[
                \varphi_i\xrightarrow{d_S} \varphi;
            \]
        \item\label{itm:cor:decnetdS:2} for each $k=0,\ldots,n$, we have
            \begin{equation}\label{eq:mixedmasslim}
                \lim_{j\in I}\int_X \theta_{\varphi_j}^k\wedge \theta_{V_{\theta}}^{n-k} =  \int_X \theta_{\varphi}^k\wedge \theta_{V_{\theta}}^{n-k}.
            \end{equation}
    \end{enumerate}
    If we assume furthermore that $\int_X \theta_{\varphi}^n>0$, then the above conditions are equivalent to the following:
    \begin{enumerate}[resume]
        \item\label{itm:cor:decnetdS:3} We have
            \[
                \lim_{j\in I}\int_X \theta_{\varphi_j}^n=\int_X \theta_{\varphi}^n.
            \]
    \end{enumerate}
    When $\int_X \theta_{\varphi}^n>0$ and these equivalent conditions are met, we also have
    \begin{equation}\label{eq:Pcontdecseq}
        P_{\theta}[\varphi]=\inf_{j\in I}P_{\theta}[\varphi_j].
    \end{equation}
\end{corollary}
\begin{proof}
    \ref{itm:cor:decnetdS:1} $\iff$ \ref{itm:cor:decnetdS:2}. This follows immediately from \cref{lma:varphileqpsi_metric}.

    Assume that $\int_X \theta_{\varphi}^n>0$.

    \ref{itm:cor:decnetdS:2} $\implies$ \ref{itm:cor:decnetdS:3}. This is the case $k=n$ in \eqref{eq:mixedmasslim}.

    \ref{itm:cor:decnetdS:3} $\implies$ \ref{itm:cor:decnetdS:2}.
    Assume (3). If for some $i\in I$, $\int_X \theta_{\varphi_i}^n= \int_X \theta_{\varphi}^n$, then for all $j\geq i$ we have
    \[
        \varphi\preceq_P\varphi_j\preceq_P\varphi_i.
    \]
    Hence $\int_X\theta_{\varphi_j}^n=\int_X\theta_{\varphi}^n$ by \cref{prop:mono2}. It follows from \cref{prop:Pequivchar2} that $\varphi_j\sim_P\varphi$, and (2) follows from \cref{prop:Ppresmass}. So in the sequel, we shall assume that for all $i\in I$, we have $\int_X \theta_{\varphi_i}^n>\int_X \theta_{\varphi}^n$.
    In particular, each $\varphi_i$ has positive mass.

    Let $(b_j)_{j\in I}$ be a net converging to $\infty$ such that
    \[
    b_j\in \left( 1, \left(\frac{\int_X \theta_{\varphi_j}^n}{\int_X \theta_{\varphi_j}^n-\int_X \theta_{\varphi}^n}\right)^{1/n} \right).
    \]
    This can be achieved for example by taking $b_j$ as the midpoint of this interval.
    By \cref{thm:pathoenvelope}, for each $j\in I$, we can find $\eta_j\in \PSH(X,\theta)$ such that
    \[
    b_j^{-1}\eta_j +\left(1-b_j^{-1}\right)\varphi_j \leq \varphi.
    \]
    It follows from the monotonicity theorem \cref{thm:mono} that for any $k=0,\ldots,n$,
    \[
    \int_X \theta_{\varphi}^k \wedge \theta_{V_{\theta}}^{n-k}\geq \left(1-b_j^{-1}\right)^k\int_X \theta_{\varphi_j}^k \wedge \theta_{V_{\theta}}^{n-k}.
    \]
    Taking the limit, we conclude the $\leq$ direction in \eqref{eq:mixedmasslim}. The $\geq$ direction follows from \cref{prop:mono2}.

    Finally, we argue \eqref{eq:Pcontdecseq}.
    Let $\psi_j=P_{\theta}[\varphi_j]$. It follows from \cref{cor:Psendspotentialtomodel} that $\psi_j$ is a model potential. Since the net $(\varphi_j)_{j\in I}$ is decreasing with respect to the $P$-preorder, $(\psi_j)_{j\in I}$ is decreasing pointwise. Let
    \[
        \psi=\inf_{j\in I}\psi_j\geq P_{\theta}[\varphi].
    \]
    By \cref{prop:decseqmodel}, $\psi$ is a model potential. Since $P_{\theta}[\varphi]\leq \psi\leq \psi_j$ for all $j\in I$, the monotonicity theorem \cref{thm:mono} and \cref{prop:Ppresmass} give
    \[
        \int_X \theta_{\varphi}^n=\int_X\theta_{P_{\theta}[\varphi]}^n\leq \int_X \theta_{\psi}^n\leq \lim_{j\in I}\int_X \theta_{\psi_j}^n=\lim_{j\in I}\int_X \theta_{\varphi_j}^n=\int_X \theta_{\varphi}^n.
    \]
    Hence $\psi$ has the same mass as $P_{\theta}[\varphi]$. As $\psi$ is model and $\psi\geq P_{\theta}[\varphi]$, \cref{thm:Pvarphidiffdef} gives $\psi=P_{\theta}[\varphi]$.
\end{proof}

Having understood the increasing and decreasing cases, we shall handle general convergent sequences. In fact, since $d_S$ is a pseudometric, the topology is completely determined by convergent sequences, so we do not need to consider nets in general.

\begin{proposition}\label{prop:incanddec}
    Let $\varphi_j,\varphi\in \PSH(X,\theta)$ ($j\geq 1$), $\varphi_j\xrightarrow{d_S}\varphi$. Assume that there is $\delta>0$ such that
    \[
    \int_X \theta_{\varphi_j}^n\geq \delta
    \]
    for all $j$ and the $\varphi_j$'s are all model potentials. Then up to replacing $(\varphi_j)_j$ by a subsequence, there is a decreasing sequence $(\psi_j)_j$ and an increasing sequence $(\eta_j)_j$ in $\PSH(X,\theta)$ such that
    \begin{enumerate}
        \item\label{itm:prop:incanddec:1} $\psi_j\xrightarrow{d_S}\varphi$, $\eta_j\xrightarrow{d_S}\varphi$;
        \item\label{itm:prop:incanddec:2} $\psi_j\geq \varphi_j\geq \eta_j$ for all $j$.
    \end{enumerate}
    In fact, for any $j\geq 1$, we will take
    \[
    \eta_j=\inf_{k\in \mathbb{N}}\varphi_j\land \varphi_{j+1}\land \cdots \land \varphi_{j+k},\quad
    \psi_j=\uscsup[k\geq j]\varphi_k.
    \]
\end{proposition}
\begin{proof}
    First note that $\int_X \theta_{\varphi}^n\geq \delta$ as well, due to the $d_S$-continuity of non-pluripolar masses as proved in \cref{cor:massdiffdS}. Replacing $\varphi$ by $P_{\theta}[\varphi]$, we may then assume that $\varphi$ is a model potential in $\PSH(X,\theta)_{>0}$.

    We are free to replace $(\varphi_j)_j$ by a subsequence. So we may assume that
    \begin{equation}\label{eq:conditiononvarphijtemp1}
    d_S(\varphi_j,\varphi_{j+1})\leq \frac{\delta 2^{-j-4}}{D_n C_n^{j+1}},\quad d_S(\varphi,\varphi_j)\leq \frac{\delta 2^{-j-4}}{D_n},
    \end{equation}
    where $C_n$ is the constant in \cref{cor:dsmetricdoubleineq}, $D_n$ is the constant in \cref{cor:massdiffdS}.

    In particular, by \cref{cor:massdiffdS},
    \begin{equation}\label{eq:varphijvarphimassdiffbdd}
        \left|\int_X \theta_{\varphi_j}^n-\int_X \theta_{\varphi}^n\right|\leq \delta 2^{-j-4}.
    \end{equation}

    \textbf{Step~1}. We handle the $\psi_j$'s. For each $j\geq 1$ and $k\geq 1$, by \cref{lma:dSalmostriang} we have
    \[
    \begin{aligned}
    d_S(\varphi_j,\varphi_j\lor\varphi_{j+1}\lor\cdots\lor \varphi_{j+k})
    &\leq C_n\,d_S(\varphi_j,\varphi_{j+1}\lor\cdots\lor \varphi_{j+k})\\
    &\leq C_n\,d_S(\varphi_j,\varphi_{j+1})\\
    &\quad + C_n\,d_S(\varphi_{j+1},\varphi_{j+1}\lor\cdots\lor \varphi_{j+k}).
    \end{aligned}
    \]
    By iteration, we find
    \[
    \begin{split}
        d_S(\varphi_j,\varphi_j\lor\varphi_{j+1}\lor\cdots\lor \varphi_{j+k})\leq \sum_{a=j}^{j+k-1}C_n^{a+1-j}d_S(\varphi_a,\varphi_{a+1})\\
        \leq \frac{C_n^{-j}}{D_n}\sum_{a=j}^{j+k-1}\delta 2^{-a-4}\leq \frac{\delta 2^{-j-3}}{D_n}.
    \end{split}
    \]
    Using \cref{cor:incnetdSconv}, we have
    \[
        \varphi_j\lor\varphi_{j+1}\lor\cdots\lor \varphi_{j+k}\xrightarrow{d_S}\psi_j
    \]
    as $k\to\infty$. Hence
    \begin{equation}\label{eq:dsvarphijpsijesttemp1}
        d_S(\varphi_j,\psi_j)\leq \frac{\delta 2^{-j-3}}{D_n}.
    \end{equation}
    We conclude that $\psi_j\xrightarrow{d_S}\varphi$.

    Moreover, we observe that
    \begin{equation}\label{eq:varphiexpressiontemp1}
        \varphi=\inf_{j\geq 1} P_{\theta}[\psi_j]
    \end{equation}
    Indeed, set $\tilde{\psi}_j=P_{\theta}[\psi_j]$. Then $(\tilde{\psi}_j)_j$ is a decreasing sequence of model potentials, and $d_S(\tilde{\psi}_j,\psi_j)=0$ by \cref{prop:ds0char}. Hence $\tilde{\psi}_j\xrightarrow{d_S}\varphi$. If $\chi=\inf_j\tilde{\psi}_j$, then \cref{prop:vol_limit_model} and \cref{cor:decnetdS} give $\tilde{\psi}_j\xrightarrow{d_S}\chi$. Thus $d_S(\chi,\varphi)=0$. Since both $\chi$ and $\varphi$ are model potentials, \cref{prop:ds0char} gives $P_{\theta}[\chi]=P_{\theta}[\varphi]$, hence $\chi=\varphi$, proving \eqref{eq:varphiexpressiontemp1}.

    \textbf{Step~2}. We consider the $\eta_j$'s.

    For each $j\geq 1$ and $k\geq 0$, we let
    \[
        \eta_j^k\coloneqq \varphi_j\land \cdots \land \varphi_{j+k}.
    \]
    Using \eqref{eq:conditiononvarphijtemp1}, \eqref{eq:dsvarphijpsijesttemp1} and \cref{cor:massdiffdS}, we have
    \[
        \left|\int_X\theta_{\psi_j}^n-\int_X \theta_{\varphi}^n \right|\leq \delta 2^{-j-2}
    \]
    for all $j\geq 1$. Taking \eqref{eq:varphijvarphimassdiffbdd} into account, we have
    \begin{equation}\label{eq:varphijpsijmassdiffbddtemp1}
     \left|\int_X\theta_{\varphi_j}^n-\int_X \theta_{\psi_{j-1}}^n \right|\leq \delta 2^{-j}
    \end{equation}
    for $j>1$.

    \textbf{Step~2.1}.
     We claim that for a fixed $j\geq 1$, for any $k\in \mathbb{N}$, we have
    $\eta_j^k\in \PSH(X,\theta)$ and
    \begin{equation}\label{eq:etajkmassbddbelowtemp1}
        \int_X \theta_{\eta_j^k}^n\geq \int_X \theta_{\varphi_j}^n-\sum_{a=1}^{k} \delta 2^{-j-a}.
    \end{equation}

    We argue by induction on $k\geq 0$.
    The case $k=0$ is trivial. When $k>0$, assume that the case $k-1$ is known. Then
    \[
        \begin{aligned}
            \int_X \theta_{\eta_j^{k-1}}^n+\int_X \theta_{\varphi_{j+k}}^n
            \geq \int_X \theta_{\varphi_j}^n-\sum_{a=1}^{k-1}\delta 2^{-j-a}+\int_X \theta_{\psi_{j+k-1}}^n-\delta 2^{-j-k}\\
            \geq \int_X \theta_{\varphi_j}^n-\delta 2^{-j}+\int_X \theta_{\psi_{j+k-1}}^n>\int_X \theta_{\psi_{j+k-1}}^n,
        \end{aligned}
    \]
    where the first inequality follows from the inductive hypothesis and \eqref{eq:varphijpsijmassdiffbddtemp1}.

    Observe that
    \[
    \eta_j^{k-1}\lor \varphi_{j+k}\leq \psi_{j+k-1},
    \]
    it follows from \cref{prop:landfinitecond1} that $\eta_j^k\in \PSH(X,\theta)$. By \cref{thm:diamond}, we deduce that
    \[
    \begin{aligned}
        \int_X \theta_{\eta_j^k}^n\geq & \int_X \theta_{\varphi_{j+k}}^n+  \int_X \theta_{\eta_j^{k-1}}^n-\int_X \theta_{\psi_{j+k-1}}^n\\
        \geq & \int_X \theta_{\varphi_j}^n-\sum_{a=1}^{k} \delta 2^{-j-a},
    \end{aligned}
    \]
    where the second inequality follows from the inductive hypothesis and \eqref{eq:varphijpsijmassdiffbddtemp1}. Therefore, \eqref{eq:etajkmassbddbelowtemp1} follows.

    \textbf{Step~2.2}.
    It follows by induction from \cref{prop:landpresmodel} that for any $j\geq 1$, $k\geq 0$,
    \[
        P_{\theta}\left[\eta_j^k\right]=\eta_j^k.
    \]
    By \cref{prop:vol_limit_model}, we have
    \[
        \lim_{k\to\infty} \int_X \theta_{\eta_j^k}^n=\int_X \theta_{\eta_j}^n
    \]
    for any $j\geq 1$.
    Letting $k\to\infty$ in \eqref{eq:etajkmassbddbelowtemp1}, we find that
    \[
        \int_X \theta_{\eta_j}^n\geq \int_X \theta_{\varphi_j}^n-\delta 2^{-j}>0
    \]
    for $j\geq 1$.
    Observe that we also have
    \[
        \int_X \theta_{\eta_j}^n\leq \int_X \theta_{\varphi_j}^n \leq \int_X \theta_{\psi_j}^n
    \]
    for $j\geq 1$ by \cref{thm:mono}.
    It follows from these inequalities and \eqref{eq:varphijvarphimassdiffbdd} that
    \[
        \lim_{j\to\infty}\int_X\theta_{\eta_j}^n=\int_X\theta_{\varphi}^n.
    \]
    For $k\geq 1$, we have
    \[
        \eta_j^k=\varphi_j\land\eta_{j+1}^{k-1}\leq \eta_{j+1}^{k-1}.
    \]
    Letting $k\to\infty$ gives $\eta_j\leq\eta_{j+1}$. Hence $(\eta_j)_j$ is increasing. Hence \cref{cor:incseqnppcont} gives
    \[
        \int_X \theta_{\eta}^n=\lim_{j\to\infty}\int_X \theta_{\eta_j}^n=\int_X \theta_{\varphi}^n,
    \]
    where $\eta=\uscsup[j\geq 1] \eta_j$.
    For fixed $m$, if $j\geq m$, then $\eta_j\leq\varphi_j\leq\psi_m$, whereas if $j<m$, then $\eta_j\leq\eta_m\leq\psi_m$. Thus $\eta_j$ is a candidate in the definition of $P_{\theta}[\psi_m]$ for all $j,m$, and \eqref{eq:varphiexpressiontemp1} gives $\eta\leq \varphi$.
    It follows from \cref{prop:Pequivchar2} that $\eta\sim_P \varphi$. By \cref{cor:incnetdSconv} and \cref{prop:ds0char}, we have $\eta_j\xrightarrow{d_S}\varphi$.
\end{proof}

\begin{corollary}\label{cor:completenessdS}
    Let $(\varphi_j)_{j\in I}$ be a Cauchy net (with respect to $d_S$) in $\PSH(X,\theta)$. Assume that there is $\delta>0$ such that $\int_X \theta_{\varphi_j}^n\geq \delta$ for all $j\in I$. Then $(\varphi_j)_{j\in I}$ converges with respect to $d_S$.

    In particular, if $(\varphi_j)_{j\in I}$ is a decreasing net such that $\int_X \theta_{\varphi_j}^n\geq \delta>0$ for all $j\in I$, then $(\varphi_j)_{j\in I}$ converges with respect to $d_S$.
\end{corollary}
We can obviously relax the decreasing condition to the following: the $P$-singularity types of $(\varphi_j)_{j\in I}$ are decreasing.

\begin{proof}
    In a pseudometric space, sequential completeness implies completeness for Cauchy nets. So for the first assertion, we may assume that $(\varphi_j)_j$ is a Cauchy sequence.

    We first treat the case of a decreasing sequence $(\gamma_j)_j$ with $\int_X\theta_{\gamma_j}^n\geq\delta$ for all $j$. Put $\tilde \gamma_j=P_\theta[\gamma_j]$. Then $(\tilde \gamma_j)_j$ is a decreasing sequence of model potentials, $d_S(\gamma_j,\tilde \gamma_j)=0$, and
    \[
        \int_X\theta_{\tilde \gamma_j}^n=\int_X\theta_{\gamma_j}^n\geq \delta
    \]
    by \cref{prop:Ppresmass}. By \cref{prop:vol_limit_model}, the decreasing limit $\tilde \gamma=\inf_j\tilde \gamma_j$ is model and
    \[
        \lim_{j\to\infty}\int_X\theta_{\tilde \gamma_j}^n=\int_X\theta_{\tilde \gamma}^n\geq\delta.
    \]
    Hence $\tilde \gamma_j\to\tilde \gamma$ in $d_S$ by \cref{cor:decnetdS}, and therefore $\gamma_j$ also converges in $d_S$.

    Next let $(\varphi_j)_{j>0}$ be a $d_S$-Cauchy sequence with masses bounded from below by $\delta$. Put
    \[
        \tilde{\varphi}_j=P_\theta[\varphi_j].
    \]
    By \cref{cor:Psendspotentialtomodel,ex:Pequiv,prop:ds0char}, the $\tilde{\varphi}_j$'s are model potentials, $d_S(\varphi_j,\tilde{\varphi}_j)=0$, and
    \[
        \int_X\theta_{\tilde{\varphi}_j}^n=\int_X\theta_{\varphi_j}^n\geq\delta
    \]
    by \cref{prop:Ppresmass}. Replacing $\varphi_j$ by $\tilde{\varphi}_j$, we may therefore assume that the $\varphi_j$'s are model potentials. Passing to a subsequence, we may also assume that
    \[
        d_S(\varphi_j,\varphi_{j+1})\leq C_n^{-2j},
    \]
    where $C_n$ is the constant in \cref{cor:dsmetricdoubleineq}.
    Define
    \[
    \psi_j=\uscsup[k\geq j]\varphi_k
    \]
    for each $j>0$. As in the proof of \cref{prop:incanddec} Step~1, especially \eqref{eq:dsvarphijpsijesttemp1}, we know that
    \[
    \lim_{j\to\infty} d_S(\varphi_j,\psi_j)=0.
    \]
    The sequence $(\psi_j)_j$ is decreasing and satisfies $\int_X\theta_{\psi_j}^n\geq\delta$ by \cref{thm:mono}. By the first paragraph, $(\psi_j)_j$ converges in $d_S$. Hence the chosen subsequence of $(\varphi_j)_j$ converges, and the original Cauchy sequence converges as well.

    Finally, suppose that $(\varphi_j)_{j\in I}$ is a decreasing net with masses bounded from below by $\delta$. For each $k=0,\ldots,n$, the net
    \[
        \int_X\theta_{\varphi_j}^k\wedge\theta_{V_\theta}^{n-k}
    \]
    is decreasing and bounded from below, hence Cauchy as a real-valued net. If $j'\geq j$, then $\varphi_{j'}\leq\varphi_j$, and \cref{lma:varphileqpsi_metric} expresses $d_S(\varphi_{j'},\varphi_j)$ as a finite sum of the corresponding differences of these mixed masses. For two sufficiently large indices, we choose a common upper bound in $I$ and use the triangle inequality. It follows that $(\varphi_j)_{j\in I}$ is $d_S$-Cauchy. The first assertion gives its convergence.
\end{proof}

\begin{lemma}\label{lma:dSsmallmult}
    There is a constant $C>0$ depending only on $X$ and $\theta$ such that for any $\varphi\in \PSH(X,\theta)$ satisfying that $\theta_{\varphi}$ is a Kähler current, we have
        \[
        d_{S,\theta}\Bigl((1-\epsilon)\varphi,\varphi\Bigr)\leq C \epsilon
        \]
        for small enough $\epsilon>0$ such that $(1-\epsilon)\varphi\in \PSH(X,\theta)$.
    \end{lemma}
    \begin{proof}
        After adding a constant, we may assume that $\varphi\leq 0$. Since $\theta_{\varphi}$ is a Kähler current, $(1-\epsilon)\varphi\in \PSH(X,\theta)$ for all small enough $\epsilon>0$. Fix such an $\epsilon$. Take any Kähler form $\omega$ on $X$ so that $\omega-\theta$ is still a Kähler form. Note that the choice of $\omega$ depends only on $X$ and $\theta$.

        Since $\varphi\preceq (1-\epsilon)\varphi$, by \cref{lma:varphileqpsi_metric} and the multi-linearity of the non-pluripolar product \itemref{itm:prop:npp1:3}, we can compute
        \[
        \begin{aligned}
            d_{S,\theta}\Bigl((1-\epsilon)\varphi,\varphi\Bigr)=& \frac{1}{n+1}\sum_{j=0}^n \left(\int_X \theta_{(1-\epsilon)\varphi}^j \wedge \theta_{V_{\theta}}^{n-j}-\int_X \theta_{\varphi}^j \wedge \theta_{V_{\theta}}^{n-j}\right)\\
            \leq & \frac{1}{n+1}\sum_{j=0}^n \left(\int_X \left((1-\epsilon)\theta_{\varphi}+\epsilon\omega\right)^j \wedge \theta_{V_{\theta}}^{n-j}-\int_X \theta_{\varphi}^j \wedge \theta_{V_{\theta}}^{n-j}\right)\\
            =&\sum_{m=1}^n c_m(\varphi)\epsilon^m,
        \end{aligned}
        \]
        where, for $1\leq m\leq n$,
        \[
            c_m(\varphi)=\frac{1}{n+1}\sum_{j=m}^n \binom{j}{m}
            \sum_{\ell=0}^m (-1)^{m-\ell}\binom{m}{\ell}
            \int_X \omega^\ell\wedge\theta_{\varphi}^{j-\ell}\wedge \theta_{V_{\theta}}^{n-j}.
        \]
        Note that
        \[
        \begin{aligned}
            \left|c_m(\varphi)\right| \leq & \frac{1}{n+1}\sum_{j=m}^n \binom{j}{m}
            \sum_{\ell=0}^m \binom{m}{\ell}
            \int_X \omega^\ell\wedge\theta_{\varphi}^{j-\ell}\wedge \theta_{V_{\theta}}^{n-j}\\
            \leq &  \frac{1}{n+1}\sum_{j=m}^n \binom{j}{m}
            \sum_{\ell=0}^m \binom{m}{\ell}
            \int_X \omega^\ell\wedge \theta_{V_{\theta}}^{n-\ell}
        \end{aligned}
        \]
        where we used the monotonicity theorem \cref{thm:mono}. Hence the whole expression is of the order of $\mathcal{O}(\epsilon)$, and our assertion follows.
    \end{proof}

\subsection{Convergence theorems}
\label{subsec:dS_conv}
Next we establish some important convergence theorems, allowing us to effectively manipulate the $d_S$-convergence.

\begin{lemma}\label{lma:dsconvpertV}
    Let $(\varphi_i)_{i\in I}$ be a net in $\PSH(X,\theta)$ and $\varphi\in \PSH(X,\theta)$. Assume that $\varphi_i\xrightarrow{d_S}\varphi$. Then for any $t\in (0,1]$,
    \[
    (1-t)\varphi_i+ tV_{\theta}\xrightarrow{d_S}(1-t)\varphi+ tV_{\theta}.
    \]
\end{lemma}
 When $t=1$, the sum is understood as in \cref{rmk:linearcompsh}.
\begin{proof}
    Fix $t\in (0,1]$. We write
    \[
    \varphi_{i,t}=(1-t)\varphi_i+ tV_{\theta},\quad \varphi_t=(1-t)\varphi+ tV_{\theta}
    \]
    for any $i\in I$.

    By \cref{cor:dsmetricdoubleineq}, it suffices to show that for each $j=0,\ldots,n$,
    \begin{equation}\label{eq:massconvafterpert}
    2\int_X \theta_{\varphi_{i,t}\lor \varphi_t}^j\wedge \theta_{V_{\theta}}^{n-j}-\int_X \theta_{\varphi_{i,t}}^j\wedge \theta_{V_{\theta}}^{n-j}-\int_X \theta_{\varphi_t}^j\wedge \theta_{V_{\theta}}^{n-j}\to 0.
    \end{equation}
    Observe that
    \[
    \varphi_{i,t}\lor \varphi_t=(1-t)(\varphi\lor \varphi_i)+t V_{\theta}.
    \]
    So after binomial expansion, \eqref{eq:massconvafterpert} follows from \cref{cor:dsmetricdoubleineq}.
\end{proof}

\begin{lemma}\label{lma:linearpertbyVtheta}
    Let $\varphi\in \PSH(X,\theta)$. For each $t\in (0,1)$, let $\varphi_t=(1-t)\varphi+t V_{\theta}$. Then
    \[
        \varphi_t\xrightarrow{d_S}\varphi
    \]
    as $t\to 0+$.
\end{lemma}
\begin{proof}
    Since $\varphi\preceq\varphi_t$, by \cref{lma:varphileqpsi_metric}, we need to show that for each $j=1,\dots,n$, we have
    \[
        \lim_{t\to 0+} \int_X \theta_{\varphi_t}^j\wedge \theta_{V_{\theta}}^{n-j}=\int_X \theta_{\varphi}^j\wedge \theta_{V_{\theta}}^{n-j}.
    \]
    For this purpose, we compute
    \[
        \begin{aligned}
            & \int_X \theta_{\varphi_t}^j\wedge \theta_{V_{\theta}}^{n-j}- \int_X \theta_{\varphi}^j\wedge \theta_{V_{\theta}}^{n-j}\\
            =&\left((1-t)^j-1\right)\int_X \theta_{\varphi}^j\wedge \theta_{V_{\theta}}^{n-j}\\
            &+\sum_{i=0}^{j-1} \binom{j}{i}(1-t)^i t^{j-i}\int_X\theta_{\varphi}^i\wedge \theta_{V_{\theta}}^{n-i}.
        \end{aligned}
    \]
    As $t\to 0+$, the right-hand side clearly tends to $0$.
\end{proof}

The following convergence theorem lies at the heart of the whole theory.
\begin{theorem}\label{thm:convdS}
    Let $\theta_1,\ldots,\theta_n$ be smooth closed real $(1,1)$-forms on $X$ representing big cohomology classes. Suppose that $(\varphi_j^k)_{k\in I}$ are nets in $\PSH(X,\theta_j)$ and $\varphi_j\in \PSH(X,\theta_j)$ for $j=1,\ldots,n$. We assume that $\varphi_j^k\xrightarrow{d_S}\varphi_j$ for each $j=1,\ldots,n$. Then
    \begin{equation}\label{eq:convmixedmassds}
        \lim_{k\in I}\int_X\theta_{1,\varphi_1^k}\wedge\cdots\wedge \theta_{n,\varphi_n^k}=\int_X\theta_{1,\varphi_1}\wedge\cdots\wedge \theta_{n,\varphi_n}.
    \end{equation}
\end{theorem}
\begin{proof}
Since $d_S$ is a pseudometric, in order to establish the continuity of mixed masses, it suffices to consider sequences instead of nets. So we may assume that $I=\mathbb{Z}_{>0}$ as ordered sets.

\textbf{Step~1}. We reduce to the case where $\varphi_j^k,\varphi_j$ all have positive masses and there is a constant $\delta>0$ such that for all $j=1,\ldots,n$ and $k\in I$,
	\[
		\int_X \theta_{j,\varphi_j^k}^n>\delta.
	\]
    Assuming in this case, the theorem has been proved, let us handle the general case.
    Take $t\in (0,1)$.
    By \cref{lma:dsconvpertV}, we have
    \[
    (1-t)\varphi_j^k+tV_{\theta_j}\xrightarrow{d_S}(1-t)\varphi_j+tV_{\theta_j}
    \]
    as $k\to\infty$ for each $j=1,\ldots,n$.
    Note that for all $k\in I$,
    \[
        \int_X \theta^n_{j,(1-t)\varphi_j^k+tV_{\theta_j}}\geq t^n\int_X \theta_{V_{\theta_j}}^n>0.
    \]
    Assume that we have proved the special case of the theorem, we have
    \[
    \begin{split}
    &\lim_{k\to\infty}\int_X\theta_{1,(1-t)\varphi_1^k+tV_{\theta_1}}\wedge\cdots\wedge \theta_{n,(1-t)\varphi_n^k+tV_{\theta_n}}\\
    =&\int_X\theta_{1,(1-t)\varphi_1+tV_{\theta_1}}\wedge\cdots\wedge \theta_{n,(1-t)\varphi_n+tV_{\theta_n}}.
    \end{split}
    \]
    More explicitly, for fixed $k$ the left-hand side is the value at $t$ of a polynomial $Q_k(t)$ of degree at most $n$, and the right-hand side is the value of a polynomial $Q(t)$ of degree at most $n$. Choose $n+1$ distinct numbers $t_0,\ldots,t_n\in (0,1)$. Lagrange interpolation expresses $Q_k(0)-Q(0)$ as a fixed linear combination of the numbers $Q_k(t_a)-Q(t_a)$, which tend to $0$. Hence the same formula holds at $t=0$.
    From this, \eqref{eq:convmixedmassds} follows.

    \textbf{Step~2}. We now assume that the assumptions mentioned at the beginning of Step~1 are met. Note that we may assume that $\varphi_j^k$, $\varphi_j$ are model potentials for all $j=1,\ldots,n$, $k>0$ by \cref{prop:ds0char} and \cref{cor:Psendspotentialtomodel}.

    We first handle the following monotone situation. Assume that, for each fixed $i$, the sequence $(\varphi_i^k)_k$ is either increasing or decreasing. Let $\tilde{\varphi}_i$ denote its pointwise limit, namely
\[
\tilde{\varphi}_i=\inf_{k>0}\varphi_i^k
\]
in the decreasing case and
\[
\tilde{\varphi}_i=\uscsup[k>0]\varphi_i^k
\]
in the increasing case. We assume moreover that $\tilde{\varphi}_i\sim_P\varphi_i$ for every $i$. Since mixed non-pluripolar masses are invariant under replacing each entry by a $P$-equivalent potential, it suffices in this monotone situation to prove convergence to the mixed mass of the $\tilde{\varphi}_i$'s.

After reordering the indices, we may assume that there is $i_0\in \{0,\ldots,n\}$ such that $(\varphi_i^k)_k$ is decreasing for $i\leq i_0$ and increasing for $i>i_0$. Therefore, for each $k>0$, using \cref{thm:mono}, we have
\[
\int_X \theta_{1,\varphi_1^k}\wedge\cdots\wedge\theta_{n,\varphi_n^k}\geq \int_X \theta_{1,\tilde{\varphi}_1}\wedge\cdots\wedge \theta_{i_0,\tilde{\varphi}_{i_0}}\wedge \theta_{i_0+1,\varphi^k_{i_0+1}}\wedge \cdots  \wedge\theta_{n,\varphi_n^k}.
\]
Using \cref{cor:incseqnppcont}, we therefore conclude that
\[
\varliminf_{k\to\infty}\int_X \theta_{1,\varphi_1^k}\wedge\cdots\wedge\theta_{n,\varphi_n^k}\geq \int_X \theta_{1,\tilde{\varphi}_1}\wedge\cdots\wedge\theta_{n,\tilde{\varphi}_n}.
\]
It remains to prove
	\begin{equation}\label{eq:limsup}
		\varlimsup_{k\to\infty}\int_X \theta_{1,\varphi_1^k}\wedge\cdots\wedge\theta_{n,\varphi_n^k}\leq\int_X \theta_{1,\tilde{\varphi}_1}\wedge\cdots\wedge\theta_{n,\tilde{\varphi}_n}.
	\end{equation}
By \cref{thm:mono}, for each $k>0$, we have
\[
\int_X \theta_{1,\varphi_1^k}\wedge\cdots\wedge\theta_{n,\varphi_n^k}\leq \int_X \theta_{1,\varphi_1^k}\wedge \cdots \wedge \theta_{i_0,\varphi_{i_0}^k}\wedge \theta_{i_0+1,\tilde{\varphi}_{i_0+1}}\wedge \cdots\wedge \theta_{n,\tilde{\varphi}_{n}}.
\]
When proving \eqref{eq:limsup}, we may replace $\varphi_j^k$ by $\tilde{\varphi}_j$ whenever $j>i_0$, $k>0$.
Thus, we are reduced to the case where for all $i$, $(\varphi_i^k)_k$ is decreasing.

By \cref{prop:vol_limit_model}, the masses $\int_X\theta_{i,\varphi_i^k}^n$ decrease to $\int_X\theta_{i,\tilde{\varphi}_i}^n$. Hence, thanks to \cref{lma:pathoenvelope}, for each $i=1,\ldots,n$, we may take an increasing sequence $(b^k_i)_k$ tending to $\infty$ satisfying
\[
b_i^k\in \left(1, \left( \frac{\int_X \theta_{i,\varphi_i^k}^n}{\int_X \theta_{i,\varphi_i^k}^n-\int_X \theta_{i,\tilde{\varphi}_i}^n} \right)^{1/n} \right)
\]
where the upper endpoint is understood as $\infty$ if the denominator is $0$. For each such $b_i^k$, \cref{lma:pathoenvelope} gives $\rho_i^k\in \PSH(X,\theta_i)$ with
\[
	(b^k_i)^{-1}\rho_i^k+\left(1-(b^k_i)^{-1}\right)\varphi_i^k\leq \tilde{\varphi}_i.
\]

Then by \cref{thm:mono} again,
\[
	\prod_{i=1}^n\left(1-(b^k_i)^{-1}\right)\int_X \theta_{1,\varphi_1^k}\wedge\cdots\wedge\theta_{n,\varphi_n^k}\leq \int_X \theta_{1,\tilde{\varphi}_1}\wedge\cdots\wedge\theta_{n,\tilde{\varphi}_n}.
\]
Letting $k\to\infty$, we conclude \eqref{eq:limsup}. This proves the monotone situation.

We now return to the general case. It suffices to prove that any subsequence of
\[
\int_X \theta_{1,\varphi_1^k}\wedge\cdots\wedge\theta_{n,\varphi_n^k}
\]
has a further subsequence converging to the desired limit. Starting with an arbitrary subsequence, we pass to a further subsequence, still indexed by $k$, so that the construction in \cref{prop:incanddec} applies simultaneously to the finitely many sequences $(\varphi_i^k)_k$, $i=1,\ldots,n$.

Thus, after discarding finitely many terms and relabeling, for each $i$ there are model potentials $\eta_i^k,\psi_i^k\in \PSH(X,\theta_i)$ such that $(\eta_i^k)_k$ is increasing, $(\psi_i^k)_k$ is decreasing, and
\[
\eta_i^k\leq \varphi_i^k\leq \psi_i^k,\quad
\eta_i^k\xrightarrow{d_S}\varphi_i,\quad
\psi_i^k\xrightarrow{d_S}\varphi_i.
\]
For the upper sequence, we take $\psi_i^k=P_{\theta_i}[\uscsup[\ell\geq k]\varphi_i^\ell]$; as in \eqref{eq:varphiexpressiontemp1}, this gives $\inf_k\psi_i^k=\varphi_i$. For the lower sequence, the construction in \cref{prop:incanddec} together with \cref{prop:decseqmodel} shows that the $\eta_i^k$'s are model after discarding finitely many terms, and the proof of \cref{prop:incanddec} shows that their increasing limit is $P$-equivalent to $\varphi_i$.

By \cref{thm:mono}, for every $k$,
\[
\int_X\theta_{1,\eta_1^k}\wedge\cdots\wedge\theta_{n,\eta_n^k}
\leq
\int_X\theta_{1,\varphi_1^k}\wedge\cdots\wedge\theta_{n,\varphi_n^k}
\leq
\int_X\theta_{1,\psi_1^k}\wedge\cdots\wedge\theta_{n,\psi_n^k}.
\]
The monotone situation just proved applies to both the lower and upper sequences, and both limiting mixed masses are equal to
\[
\int_X\theta_{1,\varphi_1}\wedge\cdots\wedge\theta_{n,\varphi_n}.
\]
The squeeze theorem therefore gives the desired convergence along this further subsequence. This concludes the proof.
\end{proof}

\begin{corollary}\label{cor:dsconvcrit}
    Suppose that $(\varphi_i)_{i\in I}$ is a net in $\PSH(X,\theta)$ and $\varphi\in \PSH(X,\theta)$. Then the following are equivalent:
    \begin{enumerate}
        \item\label{itm:cor:dsconvcrit:1} $\varphi_i\xrightarrow{d_S}\varphi$;
        \item\label{itm:cor:dsconvcrit:2} $\varphi_i\lor \varphi\xrightarrow{d_S}\varphi$ and
        \begin{equation}\label{eq:massconv_varphii}
        \lim_{i\in I}\int_X \theta_{\varphi_i}^j\wedge \theta_{V_{\theta}}^{n-j}=\int_X \theta_{\varphi}^j\wedge \theta_{V_{\theta}}^{n-j}
        \end{equation}
        for each $j=0,\ldots,n$;
        \item\label{itm:cor:dsconvcrit:3} for each $j=0,\ldots,n$, \eqref{eq:massconv_varphii} holds and
        \begin{equation}\label{eq:massconv_varphii2}
        \lim_{i\in I}\int_X \theta_{\varphi_i\lor \varphi}^j\wedge \theta_{V_{\theta}}^{n-j}=\int_X \theta_{\varphi}^j\wedge \theta_{V_{\theta}}^{n-j}.
        \end{equation}
    \end{enumerate}
\end{corollary}
The corollary allows us to reduce a number of convergence problems related to $d_S$ to the case $\varphi_i\geq \varphi$.
This is the most handy way of establishing $d_S$-convergence in practice.
\begin{proof}
    The equivalence between \ref{itm:cor:dsconvcrit:2} and \ref{itm:cor:dsconvcrit:3} follows directly from \cref{lma:varphileqpsi_metric}.

    \ref{itm:cor:dsconvcrit:1} $\implies$ \ref{itm:cor:dsconvcrit:2}. That $\varphi_i\lor \varphi\xrightarrow{d_S}\varphi$ follows from \cref{cor:dsmetricdoubleineq}. While \eqref{eq:massconv_varphii} follows from \cref{thm:convdS}.

    \ref{itm:cor:dsconvcrit:2} $\implies$ \ref{itm:cor:dsconvcrit:1}. By \eqref{eq:ds_biineq}, we need to show that for each $j=0,\ldots,n$, we have
    \[
    2\int_X \theta_{\varphi_i\lor \varphi}^j\wedge \theta_{V_{\theta}}^{n-j}-\int_X \theta_{\varphi}^j\wedge \theta_{V_{\theta}}^{n-j}-\int_X \theta_{\varphi_i}^j\wedge \theta_{V_{\theta}}^{n-j}\to 0.
    \]
    This follows from \cref{thm:convdS} and \eqref{eq:massconv_varphii}.
\end{proof}

\begin{corollary}\label{cor:dSconv_changetheta}
    Let $(\varphi_i)_{i\in I}$ be a net in $\PSH(X,\theta)$ and $\varphi\in \PSH(X,\theta)$. Let $\omega$ be a closed smooth positive $(1,1)$-form on $X$. Then the following are equivalent:
    \begin{enumerate}
        \item\label{itm:cor:dSconv_changetheta:1} $\varphi_i\xrightarrow{d_{S,\theta}}\varphi$;
        \item\label{itm:cor:dSconv_changetheta:2} $\varphi_i\xrightarrow{d_{S,\theta+\omega}}\varphi$.
    \end{enumerate}
\end{corollary}
In particular, there is no risk when we simply write $\varphi_i\xrightarrow{d_S}\varphi$. Similarly, consider two nets $(\varphi_i)_{i\in I}$, $(\psi_i)_{i\in I}$ in $\PSH(X,\theta)$. If both admit $d_S$-limits in $\PSH(X,\theta)$, then we can write $d_S(\varphi_i,\psi_i)\to 0$ without ambiguity.

\begin{proof}
    \ref{itm:cor:dSconv_changetheta:1} $\implies$ \ref{itm:cor:dSconv_changetheta:2}. It suffices to show that for each $j=0,\ldots,n$, we have
    \[
    \begin{split}
        2\int_X (\theta+\omega)_{\varphi_i\lor \varphi}^j\wedge (\theta+\omega)_{V_{\theta+\omega}}^{n-j}-\int_X (\theta+\omega)_{\varphi_i}^j\wedge (\theta+\omega)_{V_{\theta+\omega}}^{n-j}\\
        -\int_X (\theta+\omega)_{\varphi}^j\wedge (\theta+\omega)_{V_{\theta+\omega}}^{n-j}\to 0.
    \end{split}
    \]
    Note that this quantity is a linear combination of terms of the following form:
    \[
    \begin{split}
        2\int_X \theta_{\varphi_i\lor \varphi}^r\wedge \omega^{j-r}\wedge (\theta+\omega)_{V_{\theta+\omega}}^{n-j}-\int_X \theta_{\varphi_i}^r\wedge \omega^{j-r}\wedge (\theta+\omega)_{V_{\theta+\omega}}^{n-j}\\
        -\int_X \theta_{\varphi}^r\wedge \omega^{j-r}\wedge (\theta+\omega)_{V_{\theta+\omega}}^{n-j},
    \end{split}
    \]
    where $r=0,\ldots,j$. By \cref{thm:convdS}, it suffices to show that $\varphi \lor \varphi_i \xrightarrow{d_S}\varphi$. But this follows from \cref{cor:dsconvcrit}.

    \ref{itm:cor:dSconv_changetheta:2} $\implies$ \ref{itm:cor:dSconv_changetheta:1}. From the direction we already proved, for each $C\geq 1$, we have
    \[
        \varphi_i\xrightarrow{d_{S,\theta+C\omega}}\varphi.
    \]
    By \cref{thm:convdS}, it follows that
    \[
        \lim_{i\in I}\int_X (\theta+C\omega)_{\varphi_i}^j\wedge \theta_{V_{\theta}}^{n-j}=\int_X (\theta+C\omega)_{\varphi}^j\wedge \theta_{V_{\theta}}^{n-j}
    \]
    for all $j=0,\ldots,n$.
    It follows that
    \begin{equation}\label{eq:varphijmass_limit}
        \lim_{i\in I}\int_X \theta_{\varphi_i}^j\wedge \theta_{V_{\theta}}^{n-j}=\int_X \theta_{\varphi}^j\wedge \theta_{V_{\theta}}^{n-j}.
    \end{equation}
    By \cref{cor:dsconvcrit}, it remains to show that $\varphi_i\lor \varphi\xrightarrow{d_{S,\theta}}\varphi$. By \cref{cor:dsconvcrit} again, we know that $\varphi_i\lor \varphi\xrightarrow{d_{S,\theta+\omega}}\varphi$. So it suffices to apply \eqref{eq:varphijmass_limit} to $\varphi_i\lor \varphi$ instead of $\varphi_i$, and we conclude by \cref{lma:varphileqpsi_metric}.
\end{proof}

\begin{corollary}\label{cor:lineardsconv}
    Let $\varphi_0,\varphi_1\in \PSH(X,\theta)$. Define $\varphi_t=t\varphi_1+(1-t)\varphi_0$ for $t\in (0,1)$.
    Then
    \[
    \varphi_t\xrightarrow{d_S}\varphi_0
    \]
    as $t\to 0+$.
\end{corollary}
\begin{proof}
First note that for each $j=0,\ldots,n$,
\[
\lim_{t\to 0+}\int_X \theta_{\varphi_t}^j\wedge \theta_{V_{\theta}}^{n-j}=\int_X \theta_{\varphi_0}^j\wedge \theta_{V_{\theta}}^{n-j}.
\]
So thanks to \cref{cor:dsconvcrit}, it remains to argue that for all $j=0,\ldots,n$,
\[
\lim_{t\to 0+}\int_X \theta_{\varphi_t\lor \varphi_0}^j\wedge \theta_{V_{\theta}}^{n-j}=\int_X \theta_{\varphi_0}^j\wedge \theta_{V_{\theta}}^{n-j}.
\]
Observe that for $t\in (0,1)$, we have
\[
\varphi_t\lor \varphi_0=t(\varphi_1\lor \varphi_0)+(1-t)\varphi_0,
\]
so the desired limit follows by the same binomial expansion.
\end{proof}

We have the following sandwich criterion:
\begin{corollary}\label{cor:dsconvupplower}
    Let $(\varphi_i)_{i\in I}, (\psi_i)_{i\in I}, (\eta_i)_{i\in I}$ be three nets in $\PSH(X,\theta)$ and $\varphi\in \PSH(X,\theta)$. Assume that
    \begin{enumerate}
        \item\label{itm:cor:dsconvupplower:1} $\psi_i\preceq_P \varphi_i\preceq_P \eta_i$ for each $i\in I$;
        \item\label{itm:cor:dsconvupplower:2} $\eta_i\xrightarrow{d_S} \varphi$, $\psi_i\xrightarrow{d_S} \varphi$.
    \end{enumerate}
    Then $\varphi_i\xrightarrow{d_S}\varphi$.
\end{corollary}
\begin{proof}
    By \cref{cor:dSconv_changetheta}, we may replace $\theta$ by $\theta+\omega$, where $\omega$ is a Kähler form on $X$. In particular, we may assume that $\varphi_i,\psi_i,\eta_i\in \PSH(X,\theta)_{>0}$ for all $i\in I$.
    By \cref{prop:ds0char}, we may assume that $\varphi_i,\psi_i,\eta_i$ are model potentials for all $i\in I$ and hence $\psi_i \leq \varphi_i\leq \eta_i$ for all $i\in I$.

    It follows from \cref{thm:mono} that for each $k=0,\ldots,n$, we have
    \[
    \int_X \theta_{\psi_i}^k\wedge \theta_{V_{\theta}}^{n-k}\leq \int_X \theta_{\varphi_i}^k\wedge \theta_{V_{\theta}}^{n-k}\leq \int_X \theta_{\eta_i}^k\wedge \theta_{V_{\theta}}^{n-k}
    \]
    for all $i\in I$. By \cref{thm:convdS}, the limits, with respect to $i\in I$, of the two outer terms are $\int_X \theta_{\varphi}^k\wedge \theta_{V_{\theta}}^{n-k}$. It follows that
    \begin{equation}\label{eq:thetak_conv}
    \lim_{i\in I}\int_X \theta_{\varphi_i}^k\wedge \theta_{V_{\theta}}^{n-k}=\int_X \theta_{\varphi}^k\wedge \theta_{V_{\theta}}^{n-k}.
    \end{equation}
    By \cref{cor:dsconvcrit}, it remains to prove that $\varphi_i \lor \varphi \xrightarrow{d_S}\varphi$. By \cref{cor:dsconvcrit} and \cref{prop:Ppartialsup}, up to replacing $\psi_i$ (resp. $\varphi_i$, $\eta_i$) by $\psi_i\lor \varphi$ (resp. $\varphi_i\lor\varphi$, $\eta_i\lor\varphi$), we may assume from the beginning that $\psi_i,\varphi_i,\eta_i\geq \varphi$. Now $\varphi_i\xrightarrow{d_S}\varphi$ by \eqref{eq:thetak_conv} and \cref{lma:varphileqpsi_metric}.
\end{proof}

\begin{lemma}\label{lma:dS_res}
    Let $(\varphi_i)_{i\in I}$ be a net in $\PSH(X,\theta)$, and $\varphi\in \PSH(X,\theta)$. Fix $\lambda>0$. Then the following are equivalent:
    \begin{enumerate}
        \item\label{itm:lma:dS_res:1} $\varphi_i\xrightarrow{d_{S,\theta}}\varphi$;
        \item\label{itm:lma:dS_res:2} $\lambda\varphi_i\xrightarrow{d_{S,\lambda\theta}}\lambda\varphi$.
    \end{enumerate}
\end{lemma}
\begin{proof}
    We have $V_{\lambda\theta}=\lambda V_{\theta}$. The middle expression in \eqref{eq:ds_biineq} for the pair $(\lambda\varphi_i,\lambda\varphi)$ in the class $\lambda\theta$ is therefore $\lambda^n$ times the corresponding expression for $(\varphi_i,\varphi)$ in the class $\theta$. The equivalence follows from \cref{cor:dsmetricdoubleineq}.
\end{proof}

\begin{proposition}\label{prop:lor_dS_conv}
    Let $(\varphi_i)_{i\in I}$ (resp. $(\psi_i)_{i\in I}$) be a net in $\PSH(X,\theta)$ such that $\varphi_i \xrightarrow{d_S}\varphi\in \PSH(X,\theta)$ (resp. $\psi_i \xrightarrow{d_S}\psi\in \PSH(X,\theta)$).
    Then
    \[
    \varphi_i\lor\psi_i\xrightarrow{d_S}\varphi\lor \psi.
    \]
\end{proposition}
\begin{proof}
    Since $d_S$ is a pseudometric, we may assume that both nets are actually sequences and $I=\mathbb{Z}_{>0}$. By \cref{cor:dSconv_changetheta}, we may assume that the masses $\int_X \theta_{\varphi}^n>0$, $\int_X \theta_{\psi}^n>0$.

    By \cref{prop:ds0char} and \cref{prop:Ppartialsup}, replacing all potentials by their $P$-envelopes changes neither the hypotheses nor the conclusion. Thus we may assume that the $\varphi_j$'s, the $\psi_j$'s, $\varphi$ and $\psi$ are all model. After discarding finitely many terms, the masses of the two sequences are bounded below by a positive constant. Applying \cref{prop:incanddec} to both sequences and using the sandwich criterion \cref{cor:dsconvupplower}, it is enough to treat the case where each of $(\varphi_j)_j$ and $(\psi_j)_j$ is monotone. In this reduced situation, $(\varphi_j)_j$ (resp. $(\psi_j)_j$) converges to $\varphi$ (resp. $\psi$) almost everywhere due to \cref{prop:decPorderPenvconve}, \cref{cor:incnetdSconv} and \cref{prop:incnetmodel}.

    We handle three cases separately.

\textbf{Step~1}. Assume that both sequences are increasing.

In this case, we have $\varphi_j\lor \psi_j\nearrow \varphi\lor \psi$ almost everywhere. Therefore, $\varphi_j\lor \psi_j\xrightarrow{d_S} \varphi\lor \psi$ by \cref{cor:incnetdSconv}.

\textbf{Step~2}. Assume that one sequence, say $(\varphi_j)_j$ is increasing while the other is decreasing. Then we have
\[
\varphi_j \lor \psi \leq \varphi_j \lor \psi_j \leq \varphi \lor \psi_j.
\]
Thanks to \cref{cor:dsconvupplower}, it suffices to show that both sides converge to $\varphi\lor \psi$ with respect to $d_S$. So we reduce to the case where both sequences are decreasing.

\textbf{Step~3}. Assume that both sequences are decreasing.

In this case, due to \cref{cor:decnetdS}, it suffices to show that
\begin{equation}\label{eq:masslortemp1}
\lim_{j\to\infty} \int_X \theta_{\varphi_j\lor \psi_j}^n=\int_X \theta_{\varphi\lor \psi}^n.
\end{equation}
The $\geq$ direction follows from \cref{thm:mono}; it remains to argue the $\leq$ direction.

Thanks to \cref{lma:pathoenvelope}, we may find a sequence $(\epsilon_j)_j$ in $(0,1)$ with limit $0$ and a sequence $(\eta_j)_j$ in $\PSH(X,\theta)_{>0}$ such that
\[
    (1-\epsilon_j)\varphi_j+\epsilon_j \eta_j \leq \varphi,\quad \eta_j \leq \varphi_j.
\]
It follows that for each $j\geq 1$, we have
\[
(1-\epsilon_j)(\varphi_j\lor \psi_j)+\epsilon_j \eta_j \leq \varphi \lor \psi_j.
\]
Therefore by \cref{thm:mono},
\[
(1-\epsilon_j)^n \int_X \theta_{\varphi_j\lor \psi_j}^n\leq \int_X \theta_{\varphi \lor \psi_j}^n.
\]
Letting $j\to\infty$, we find that
\[
\lim_{j\to\infty}\int_X \theta_{\varphi_j\lor \psi_j}^n\leq \lim_{j\to\infty}\int_X \theta_{\varphi\lor \psi_j}^n.
\]
Therefore, in order to prove \eqref{eq:masslortemp1}, we may assume that one of the sequences is constant; let us say $\psi_j=\psi$ for all $j$. Repeating the same argument as before and constructing $(\epsilon_j)_j$, $(\eta_j)_j$ as above, we get
\[
(1-\epsilon_j)^n \int_X \theta_{\varphi_j\lor \psi}^n\leq \int_X \theta_{\varphi \lor \psi}^n.
\]
Letting $j\to\infty$, we conclude \eqref{eq:masslortemp1}.
\end{proof}

\begin{proposition}\label{prop:dsconvpresorder}
    Let $(\varphi_i)_{i\in I}$, $(\psi_i)_{i\in I}$ be nets in $\PSH(X,\theta)$ such that $\varphi_i\xrightarrow{d_S} \varphi\in \PSH(X,\theta)$ and $\psi_i\xrightarrow{d_S}\psi\in \PSH(X,\theta)$. Assume that $\varphi_i\preceq_P \psi_i$ for all $i\in I$. Then $\varphi\preceq_P \psi$.
\end{proposition}
\begin{proof}
    It follows from \cref{prop:lor_dS_conv} that
    \[
        \varphi_i\lor\psi_i\xrightarrow{d_S}\varphi\lor \psi.
    \]
    By \cref{lma:reform_preceqP}, we have $\varphi_i\lor\psi_i\sim_P \psi_i$ for all $i\in I$. In particular, by \cref{prop:ds0char},
    \[
        \varphi_i\lor\psi_i\xrightarrow{d_S}\psi.
    \]
    By \cref{prop:ds0char} again, $\varphi\lor\psi\sim_P\psi$ and hence $\varphi\preceq_P \psi$ by \cref{lma:reform_preceqP}.
\end{proof}

\begin{theorem}\label{thm:dSadditivity}
    Let $\theta_1$, $\theta_2$ be smooth real closed $(1,1)$-forms on $X$ representing big cohomology classes. Suppose that $(\varphi_i)_{i\in I}$ (resp. $(\psi_i)_{i\in I}$) be a net in $\PSH(X,\theta_1)$ (resp. $\PSH(X,\theta_2)$) and $\varphi\in \PSH(X,\theta_1)$ (resp. $\psi\in \PSH(X,\theta_2)$).
    Consider the following three conditions:
    \begin{enumerate}
        \item\label{itm:thm:dSadditivity:1} $\varphi_i\xrightarrow{d_{S}}\varphi$;
        \item\label{itm:thm:dSadditivity:2} $\psi_i\xrightarrow{d_{S}}\psi$;
        \item\label{itm:thm:dSadditivity:3} $\varphi_i+\psi_i\xrightarrow{d_{S}}\varphi+\psi$.
    \end{enumerate}
    Then any two of these conditions imply the third.
\end{theorem}
\begin{proof}
By \cref{cor:dSconv_changetheta}, we may assume that $\theta_1,\theta_2$ are both Kähler forms. We denote them by $\omega_1,\omega_2$ instead. Let $\omega=\omega_1+\omega_2$.

\ref{itm:thm:dSadditivity:1}+\ref{itm:thm:dSadditivity:2} $\implies$ \ref{itm:thm:dSadditivity:3}.
	It suffices to show that for each $r=0,\ldots,n$,
	\[
		2\int_X \omega_{(\varphi_i+\psi_i)\lor(\varphi+\psi)}^r\wedge \omega^{n-r}-\int_X \omega_{\varphi_i+\psi_i}^r\wedge \omega^{n-r}-\int_X \omega_{\varphi+\psi}^r\wedge\omega^{n-r}\to 0.
	\]
	Observe that for each $i\in I$,
	\[
		(\varphi_i+\psi_i)\lor (\varphi+\psi)\leq (\varphi_i\lor \varphi)+(\psi_i \lor \psi).
	\]
	The expression with the actual rooftop is non-negative and is bounded above by the same expression with this larger potential in place of $(\varphi_i+\psi_i)\lor(\varphi+\psi)$. Thus, it suffices to show that
	\[
		2\int_X \omega_{(\varphi_i\lor \varphi)+(\psi_i\lor \psi)}^r\wedge \omega^{n-r}-\int_X \omega_{\varphi_i+\psi_i}^r\wedge \omega^{n-r}-\int_X \omega_{\varphi+\psi}^r\wedge\omega^{n-r}\to 0.
	\]
	The left-hand side is a linear combination of
	\[
		2\int_X \omega_{1,\varphi_i\lor\varphi}^a\wedge \omega_{2,\psi_i\lor\psi}^{r-a}\wedge \omega^{n-r}-\int_X \omega_{1,\varphi_i}^a\wedge \omega_{2,\psi_i}^{r-a} \wedge \omega^{n-r}-\int_X \omega_{1,\varphi}^a\wedge \omega_{2,\psi}^{r-a} \wedge \omega^{n-r}
	\]
	with $a=0,\ldots,r$. Observe that $\varphi_i\lor\varphi\xrightarrow{d_S}\varphi$ and $\psi_i\lor\psi\xrightarrow{d_S}\psi$ by \cref{cor:dsmetricdoubleineq}, each term tends to $0$ by \cref{thm:convdS}.

\ref{itm:thm:dSadditivity:1}+\ref{itm:thm:dSadditivity:3} $\implies$ \ref{itm:thm:dSadditivity:2}. For each $C\geq 1$, from the direction we already proved together with \cref{lma:dS_res},
    \[
    C\varphi_i+\psi_i\xrightarrow{d_{S}} C\varphi+\psi.
    \]
By \cref{thm:convdS}, for each $j=0,\ldots,n$,
\[
\begin{split}
&\lim_{i\in I} \int_X \Bigl(C\omega_1+\omega_2+\ddc(C\varphi_i+\psi_i)\Bigr)^j\wedge \omega_2^{n-j}\\
= &\int_X \Bigl(C\omega_1+\omega_2+\ddc(C\varphi+\psi)\Bigr)^j\wedge \omega_2^{n-j}.
\end{split}
\]
Both sides are polynomials in $C$. Comparing the constant terms gives
\begin{equation}\label{eq:psii_quant_conv}
\lim_{i\in I} \int_X \omega_{2,\psi_i}^j\wedge \omega_2^{n-j}=\int_X \omega_{2,\psi}^j\wedge \omega_2^{n-j}.
\end{equation}
Therefore, \ref{itm:thm:dSadditivity:2} follows if $\psi_i\geq \psi$ for each $i$ by \cref{lma:varphileqpsi_metric}.

Next we prove the general case. By the direction that we already proved, we know that $\varphi_i+\psi\xrightarrow{d_{S}} \varphi+\psi$. By \cref{prop:lor_dS_conv}, we have
\[
\varphi_i+(\psi_i\lor \psi)\xrightarrow{d_{S}} \varphi+\psi.
\]
It follows from the special case above that $\psi_i\lor \psi\xrightarrow{d_{S}}\psi$. It follows from \eqref{eq:psii_quant_conv} and \cref{cor:dsconvcrit} that \ref{itm:thm:dSadditivity:2} holds.

\ref{itm:thm:dSadditivity:2}+\ref{itm:thm:dSadditivity:3} $\implies$ \ref{itm:thm:dSadditivity:1}. This is similar.

\end{proof}

\begin{proposition}\label{prop:decPorderPenvconve}
    Let $(\varphi_i)_{i\in I}$ be a net in $\PSH(X,\theta)$, decreasing with respect to the $P$-partial order. Assume that there is $\varphi\in \PSH(X,\theta)_{>0}$ with $\varphi_i\xrightarrow{d_S}\varphi$, then
    \begin{equation}\label{eq:decPorderPenvconve}
        \inf_{i\in I} P_{\theta}[\varphi_i]=P_{\theta}[\varphi].
    \end{equation}
\end{proposition}
\begin{proof}
    First observe that by our assumption, $(P_{\theta}[\varphi_i])_{i\in I}$ is a decreasing net. By the continuity of non-pluripolar masses \cref{cor:massdiffdS}, we have
    \[
        \lim_{i\in I} \int_X \theta_{P_{\theta}[\varphi_i]}^n=\lim_{i\in I} \int_X \theta_{\varphi_i}^n= \int_X \theta_{\varphi}^n.
    \]
    Replacing $I$ by a tail, we may therefore assume that $\varphi_i\in \PSH(X,\theta)_{>0}$ for each $i\in I$.
    Since $\varphi_i\sim_P P_\theta[\varphi_i]$, replacing $\varphi_i$ by $P_\theta[\varphi_i]$ does not change its $d_S$-class by \cref{prop:ds0char}. We may assume that each $\varphi_i$ is model and $(\varphi_i)_{i\in I}$ is decreasing. In this case, by \cref{prop:vol_limit_model}, the potential
    \[
        \psi\coloneqq \inf_{i\in I} \varphi_i
    \]
    is model in $\PSH(X,\theta)_{>0}$ with
    \[
        \int_X \theta_{\psi}^n=\lim_{i\in I} \int_X \theta_{\varphi_i}^n= \int_X \theta_{\varphi}^n.
    \]
    Therefore, by \cref{cor:decnetdS}, we find $\varphi_i\xrightarrow{d_S}\psi$. In particular, $\varphi\sim_P \psi$ by \cref{prop:ds0char}. Therefore,
    \[
        \psi=P_{\theta}[\psi]=P_{\theta}[\varphi].
    \]
    Therefore, \eqref{eq:decPorderPenvconve} follows.
\end{proof}

\begin{theorem}\label{thm:contPI}
    The map
    \[
        P_{\theta}[\bullet]_{\mathcal{I}}\colon \PSH(X,\theta)_{>0}\rightarrow \PSH(X,\theta)_{>0}
    \]
    is continuous with respect to $d_S$.
   \end{theorem}
   \begin{proof}
       Let $(\varphi_i)_{i\in \mathbb{Z}_{>0}}$ be a sequence in $\PSH(X,\theta)_{>0}$ such that $\varphi_i\xrightarrow{d_S}\varphi\in \PSH(X,\theta)_{>0}$. We want to show that
       \begin{equation}
           P_{\theta}[\varphi_i]_{\mathcal{I}}\xrightarrow{d_S}P_{\theta}[\varphi]_{\mathcal{I}}.
       \end{equation}
       We may assume that the $\varphi_i$'s and $\varphi$ are all model potentials by \cref{prop:ds0char}.

       By \cref{prop:incanddec} and \cref{cor:dsconvupplower}, we may assume that $(\varphi_i)_i$ is either increasing or decreasing.
       In the increasing case, we apply \cref{prop:incnetmodelPI} and \cref{cor:incnetdSconv}, while in the decreasing case, we apply \cref{prop:decnetmodelPI}, \cref{prop:vol_limit_model} and \cref{cor:decnetdS}.
    \end{proof}

We record the following result for later use.
\begin{lemma}\label{lma:Vthetaepsconv}
Fix a Kähler form $\omega$ on $X$.
As $\epsilon\to 0$ (from both sides), we have
\begin{equation}\label{eq:Vthetaeps_dS}
    V_{\theta+\epsilon\omega}\xlongrightarrow{d_S}V_{\theta}.
\end{equation}
\end{lemma}
\begin{proof}
There are two assertions to prove, as detailed in the two steps.

\textbf{Step~1}. We first handle the case where $\epsilon\to 0+$.

In this case \eqref{eq:Vthetaeps_dS} means
\[
    V_{\theta+\epsilon\omega}\xlongrightarrow{d_{S,\theta+\omega}}V_{\theta}
\]
as $\epsilon\to 0+$. So thanks to \cref{cor:decnetdS} and \cref{prop:vol_limit_model}, it suffices to prove the following:
\begin{equation}\label{eq:epsil_P_temp1}
    \inf_{\epsilon>0} P_{\theta+\omega}\left[ V_{\theta+\epsilon\omega} \right]=P_{\theta+\omega}\left[ V_{\theta} \right].
\end{equation}
First observe that
\[
V_{\theta}=\inf_{\epsilon>0} V_{\theta+\epsilon\omega}.
\]
In fact, the $\leq$ direction is trivial. As for the reverse inequality, it suffices to observe that the right-hand side lies in $\PSH(X,\theta)$.
Therefore, due to \cref{prop:vol_limit_model},
\[
\lim_{\epsilon \to 0+} \int_X \left( \theta+\epsilon\omega+\ddc V_{\theta+\epsilon\omega} \right)^n=\int_X \theta_{V_{\theta}}^n.
\]
Therefore, for all $\epsilon>0$ small enough, we can find $\eta_{\epsilon}\in \PSH(X,\theta+\epsilon\omega)$ and $a_{\epsilon}>0$ decreasing to $0$ so that
\[
(1-a_\epsilon)V_{\theta+\epsilon \omega} + a_\epsilon \eta_\epsilon \leq V_{\theta}.
\]
Therefore, thanks to \cref{prop:Pconc},
\[
(1-a_\epsilon)P_{\theta+\omega}\left[ V_{\theta+\epsilon \omega} \right] + a_\epsilon P_{\theta+\omega}\left[\eta_\epsilon\right] \leq P_{\theta+\omega}\left[V_{\theta}\right].
\]
Since the envelopes $P_{\theta+\omega}[\eta_\epsilon]$ have supremum $0$, they are uniformly bounded in $L^1$ by \cref{prop:L1compa}; hence $a_\epsilon P_{\theta+\omega}[\eta_\epsilon]\to 0$ in $L^1$. Letting $\epsilon\to 0$ and using \cref{prop:pshfuncdetdense}, we conclude \eqref{eq:epsil_P_temp1}.

\textbf{Step~2}. We then handle the case where $\epsilon\to 0-$.

In this case, \eqref{eq:Vthetaeps_dS} simply means
\[
V_{\theta-\epsilon\omega}\xrightarrow{d_{S,\theta}}V_{\theta}
\]
as $\epsilon\to 0+$.
But this follows from \cref{cor:incnetdSconv} if we can prove
\begin{equation}\label{eq:Vthetamepsomega_temp1}
\uscsup[\epsilon>0] V_{\theta-\epsilon\omega}=V_{\theta},
\end{equation}
where we understand that $\epsilon$ is small enough so that $\theta-\epsilon\omega$ represents a big cohomology class. Since $\{\theta\}$ is big, we can find $\varphi\in \PSH(X,\theta)$ so that
\[
\varphi\leq 0,\quad \theta_{\varphi}\geq \delta\omega
\]
for some $\delta>0$. For a small $\epsilon>0$, we then have
\[
\theta+\ddc \left( \left(1-\frac{\epsilon}{\delta}\right)V_{\theta}+ \frac{\epsilon}{\delta}\varphi\right)\geq \epsilon\omega.
\]
Therefore,
\[
\left(1-\frac{\epsilon}{\delta}\right)V_{\theta}+ \frac{\epsilon}{\delta}\varphi\leq V_{\theta-\epsilon\omega}.
\]
Letting $\epsilon\to 0+$, we then find
\[
V_{\theta}\leq \uscsup[\epsilon>0] V_{\theta-\epsilon\omega}.
\]
The reverse inequality is trivial, hence \eqref{eq:Vthetamepsomega_temp1} is established.
\end{proof}

\subsection{Continuity of invariants}\label{subsec:continv}
In this section, we prove the continuity of a few invariants of the singularities with respect to $d_S$.

    \begin{theorem}\label{thm:Lelongcont}
        Let $(\varphi_j)_{j\in I}$ be a net in $\PSH(X,\theta)$ and $\varphi_j\xrightarrow{d_S}\varphi\in \PSH(X,\theta)$. Then for any prime divisor $E$ over $X$, we have
        \begin{equation}\label{eq:convnu}
        \lim_{j\in I}\nu(\varphi_j,E)=\nu(\varphi,E).
        \end{equation}
    \end{theorem}
\begin{proof}
    First observe that since $d_S$ is a pseudometric, it suffices to prove \eqref{eq:convnu} when $I=\mathbb{Z}_{>0}$ as partially ordered sets.

    By \cref{cor:dSconv_changetheta}, after adding a Kähler form to $\theta$, we may assume that the masses of $\varphi_j$ and of $\varphi$ are bounded from below by a positive constant.

    By \cref{thm:contPI}, we may assume that $\varphi_j$ and $\varphi$ are both $\mathcal{I}$-model and hence model. When proving \eqref{eq:convnu}, we are free to pass to subsequences.

    By \cref{prop:incanddec} and a sandwich argument, we may assume that the sequence $(\varphi_j)_j$ is either increasing or decreasing. In the increasing case, the Lelong numbers $\nu(\varphi_j,E)$ decrease to a limit at least $\nu(\varphi,E)$, while $L^1$-convergence gives the opposite inequality by upper semicontinuity of Lelong numbers. In the decreasing case, \eqref{eq:convnu} follows from \cref{prop:vol_limit_model}.
\end{proof}

\begin{theorem}\label{thm:contvolu}
    Let $(\varphi_j)_{j\in I}$ be a net in $\PSH(X,\theta)$ and $\varphi\in \PSH(X,\theta)_{>0}$. Assume that $\varphi_j\xrightarrow{d_S}\varphi$. Then
    \begin{equation}\label{eq:Ivolcont}
        \vol \theta_{\varphi_j}\to \vol \theta_{\varphi},\quad \int_X \theta_{\varphi_j}^n\to \int_X \theta_{\varphi}^n.
    \end{equation}
\end{theorem}
Recall the volume is defined in \cref{def:volqpsh}. In fact, we do not have to assume the positivity of the mass of $\varphi$. The proof of the general statement is slightly more involved. See \cref{cor:contvolu2} below.
\begin{proof}
    The latter part of \eqref{eq:Ivolcont} is just a special case of \cref{thm:convdS}. It remains to prove the former part.

    We may therefore assume that $\int_X \theta_{\varphi_j}^n>0$ for all $j\in I$ after passing to a tail. Then by \cref{thm:contPI}, we have
    \[
        P_{\theta}[\varphi_j]_{\mathcal{I}}\xrightarrow{d_S}P_{\theta}[\varphi]_{\mathcal{I}}.
    \]
    Therefore, the first part of \eqref{eq:Ivolcont} follows again from \cref{thm:convdS}.
\end{proof}

Next we show that $d_S$-convergent sequences have a sort of quasi-equisingular property (cf. \eqref{eq:quasi_eq_prop}).
\begin{theorem}\label{thm:equising_cond_general}
    Let $\varphi_j,\varphi\in \PSH(X,\theta)$ ($j\in \mathbb{Z}_{>0}$). Assume that $\varphi_j\xrightarrow{d_S}\varphi$. Then for each $\lambda'>\lambda>0$, there is $j_0>0$ so that for $j\geq j_0$,
    \begin{equation}\label{eq:quasi_equi_cond}
        \mathcal{I}(\lambda' \varphi_j)\subseteq \mathcal{I}(\lambda \varphi).
    \end{equation}
\end{theorem}
\begin{proof}
Fix $\lambda'>\lambda>0$. We want to find $j_0>0$ so that for $j\geq j_0$, \eqref{eq:quasi_equi_cond} holds.

\textbf{Step~1}. We first assume that $\varphi$ has analytic singularities.

Let $\pi\colon Y\rightarrow X$ be a log resolution of $\varphi$ and let $E_1,\ldots,E_N$ be all prime divisors in the polar locus of $\varphi$ on $Y$.
Fix $x\in X$. By \cref{thm:valuativemulti}, a non-zero germ $f\in\mathcal O_{X,x}$ belongs to the stalk of the right-hand side of \eqref{eq:quasi_equi_cond} if and only if for every $i$ such that $E_i$ meets $\pi^{-1}(x)$,
\begin{equation}\label{eq:ordEif}
\ord_{E_i}(f)>\lambda \nu(\varphi,E_i)-A_X(E_i),
\end{equation}
where $A_X$ denotes the log discrepancy.
Similarly, if $f$ belongs to the stalk of the left-hand side of \eqref{eq:quasi_equi_cond}, then for each of these $i$,
\[
\ord_{E_i}(f)> \lambda'\nu(\varphi_j,E_i)-A_X(E_i).
\]
As Lelong numbers are continuous with respect to $d_S$ by \cref{thm:Lelongcont}, we can find $j_0>0$ so that when $j\geq j_0$, $\lambda'\nu(\varphi_j,E_i)\geq \lambda \nu(\varphi,E_i)$ for all $i$. In particular, \eqref{eq:ordEif} follows.

\textbf{Step~2}. We handle the general case.

By \cref{cor:dSconv_changetheta}, we are free to increase $\theta$ and assume that $\theta_{\varphi}$ is a Kähler current.

Take a quasi-equisingular approximation $(\psi_k)_k$ of $\varphi$ in $\PSH(X,\theta)$. The existence is guaranteed by \cref{thm:qequi}.
Take $\lambda''\in (\lambda,\lambda')$. Then by definition of quasi-equisingular approximations (cf. \eqref{eq:quasi_eq_prop}), we can find $k>0$ so that
\[
    \mathcal{I}(\lambda''\psi_k)\subseteq \mathcal{I}(\lambda \varphi).
\]
Observe that $\psi_k\geq \varphi$ for all $k>0$. By \cref{prop:lor_dS_conv}, we then find that for any $k>0$, we have $\varphi_j\lor \psi_k\xrightarrow{d_S}\varphi\lor \psi_k=\psi_k$ as $j\to\infty$.
By Step~1, we can find $j_0>0$ so that for $j\geq j_0$,
\[
    \mathcal{I}\Bigl(\lambda' (\varphi_j\lor \psi_k)\Bigr)\subseteq \mathcal{I}(\lambda''\psi_k).
\]
It follows that for $j\geq j_0$,
\[
    \mathcal{I}(\lambda' \varphi_j)\subseteq  \mathcal{I}(\lambda \varphi).
\]
This completes the proof.
\end{proof}

\motto{Le but de cette thèse est de munir son auteur du titre de Docteur.\footnote{Similarly, the purpose of the current book is to make my complaints about France in the acknowledgments published.}

--- Adrien Douady\footnote{Adrien Douady (1935--2006) was a French mathematician known for his pioneering work in complex dynamics and fractal geometry. Along with John H. Hubbard, he proved important results about the Mandelbrot set and developed renormalization theory for polynomial mappings. He discovered the \emph{Douady Rabbit}, a famous fractal Julia set.

Douady studied at École Normale Supérieure (the place where I began to hate France, thanks to Claude Viterbo) and taught at several French universities. He was also a member of the Bourbaki group.

Tragically, he died in a swimming accident in 2006.
}, at the beginning of his thesis}
\chapter{\texorpdfstring{$\mathcal{I}$}{I}-good singularities}
\label{chap:Igood}

The class of $\mathcal I$-good singularities is the central object of this part of the book. These are the singularities for which the pluripotential-theoretic and multiplier-ideal notions of singularity type agree. Equivalently, they are the positive-mass potentials for which the ordinary model envelope and the $\mathcal I$-model envelope coincide, or for which the non-pluripolar mass equals the volume defined through multiplier ideals.

The starting point is the characterization theorem \cref{thm:charIgoodasclosure}. For $\varphi\in\PSH(X,\theta)_{>0}$, it shows that the following three conditions are equivalent: $\varphi$ lies in the $d_S$-closure of potentials with analytic singularities; the mass-volume identity
\[
    \int_X\theta_\varphi^n=\vol\theta_\varphi
\]
holds; and
\[
    P_\theta[\varphi]=P_\theta[\varphi]_{\mathcal I}.
\]

After establishing this characterization, the chapter develops the basic calculus of $\mathcal I$-good singularities. We prove stability under natural operations, including linear combinations and suprema, and record the main examples: analytic singularities, $\mathcal I$-model potentials, elementary metrics, and toric potentials. These closure properties are used repeatedly in later chapters, especially for traces, test curves and non-Archimedean potentials.

The next part of the chapter introduces mixed volumes for singular currents. The definition is first made in positive-volume classes by replacing each current with its $\mathcal I$-model envelope, and is then extended by K\"ahler perturbation. We prove symmetry, multilinearity, monotonicity, birational invariance, log-concavity and continuity with respect to $d_S$-convergence.

The final section proves the asymptotic Riemann--Roch formula \cref{thm:DXmain1}. For a Hermitian pseudo-effective line bundle $(L,h)$ on a connected compact K\"ahler manifold $X$ of dimension $n$, the theorem identifies
\[
    \lim_{k\to\infty}\frac{n!}{k^n}
    h^0\left(X,L^k\otimes\mathcal I(h^k)\right)
\]
with the volume of the Chern current of $h$. This result is the main algebraic manifestation of $\mathcal I$-goodness. It will be used later in the construction of partial Okounkov bodies and in the study of partial Bergman kernels.

\section{The notion of \texorpdfstring{$\mathcal{I}$}{I}-good singularities}\label{sec:Igoodnotion}

Let $X$ be a connected compact Kähler manifold of dimension $n$.

\begin{theorem}\label{thm:charIgoodasclosure}
    Let $\theta$ be a closed real smooth $(1,1)$-form on $X$ representing a big cohomology class, and $\varphi\in \PSH(X,\theta)_{>0}$. Then the following are equivalent:
    \begin{enumerate}
        \item\label{itm:thm:charIgoodasclosure:1} There exists a sequence $(\varphi_j)_{j>0}$ in $\PSH(X,\theta)$ with analytic singularities such that $\varphi_j\xrightarrow{d_S}\varphi$;
        \item\label{itm:thm:charIgoodasclosure:2} we have
        \begin{equation}\label{eq:nppmassequalvolume}
            \int_X \theta_{\varphi}^n=\vol \theta_{\varphi};
        \end{equation}
        \item\label{itm:thm:charIgoodasclosure:3} we have
            \begin{equation}\label{eq:PequpI}
                P_{\theta}[\varphi]=P_{\theta}[\varphi]_{\mathcal{I}}.
            \end{equation}
    \end{enumerate}
    In \ref{itm:thm:charIgoodasclosure:1}, we can in addition require that each $\theta_{\varphi_j}$ is a Kähler current.

    Moreover, if $\theta_{\varphi}$ is a Kähler current, one may take the sequence in \ref{itm:thm:charIgoodasclosure:1} to be any quasi-equisingular approximation of $\varphi$ in $\PSH(X,\theta)$.
\end{theorem}
Since $(\PSH(X,\theta),d_S)$ is a pseudometric space, in \ref{itm:thm:charIgoodasclosure:1} we can also replace the word \emph{sequence} by \emph{net}.

Recall that according to \cref{cor:compnppmassandvol} and \cref{prop:PandPI}, one direction of \eqref{eq:nppmassequalvolume} and \eqref{eq:PequpI} always holds:
\[
 \int_X \theta_{\varphi}^n\leq \vol \theta_{\varphi},\quad P_{\theta}[\varphi]\leq P_{\theta}[\varphi]_{\mathcal{I}}.
\]
\begin{proof}
    \ref{itm:thm:charIgoodasclosure:1} $\implies$ \ref{itm:thm:charIgoodasclosure:2}. By \cref{thm:convdS}, we have
    \[
    \lim_{j\to\infty} \int_X \theta_{\varphi_j}^n=\int_X \theta_{\varphi}^n>0.
    \]
    We may therefore assume that $\int_X \theta_{\varphi_j}^n>0$ for all $j \geq 1$.
    It follows from \cref{prop:analysingcompPandPI} that
    \[
        \int_X \theta_{\varphi_j}^n=\vol \theta_{\varphi_j}
    \]
    for any $j\geq 1$. Using \cref{thm:contvolu}, we conclude \eqref{eq:nppmassequalvolume}.

    \ref{itm:thm:charIgoodasclosure:2} $\iff$ \ref{itm:thm:charIgoodasclosure:3}. By definition,
    \[
        \vol\theta_{\varphi}= \int_X\theta_{P_{\theta}[\varphi]_{\mathcal I}}^n.
    \]
    Moreover, $\varphi-\sup_X\varphi\leq P_{\theta}[\varphi]_{\mathcal I}$.
    Thus, if \ref{itm:thm:charIgoodasclosure:2} holds, \cref{thm:Pvarphidiffdef} gives $P_{\theta}[\varphi]_{\mathcal I}\leq P_{\theta}[\varphi]$; the reverse inequality follows from \cref{prop:PandPI}. Conversely, if \ref{itm:thm:charIgoodasclosure:3} holds, then \ref{itm:thm:charIgoodasclosure:2} follows from \cref{prop:Ppresmass}.

    \ref{itm:thm:charIgoodasclosure:3} $\implies$ \ref{itm:thm:charIgoodasclosure:1}. Note that the condition in \ref{itm:thm:charIgoodasclosure:1} characterizes the closure of analytic singularities in $\PSH(X,\theta)$.

    \textbf{Step~1}. We first assume that $\theta_{\varphi}$ is a Kähler current. We will prove the following more general result in this case: Without assuming \ref{itm:thm:charIgoodasclosure:3}, $P_{\theta}[\varphi]_{\mathcal{I}}$ always lies in the closure of analytic singularities.

    Let $(\varphi_j)_j$ be a quasi-equisingular approximation of $\varphi$ in $\PSH(X,\theta)$. We will show that $\varphi_j\xrightarrow{d_S}P_{\theta}[\varphi]_{\mathcal{I}}$.
    Let
\[
    \psi=\inf_{j\in \mathbb{Z}_{>0}} P_{\theta}[\varphi_j].
\]
Since $\varphi\leq \varphi_j$ for all $j>0$, we have
\[
    \varphi-\sup_X\varphi\leq P_{\theta}[\varphi_j]
\]
for all $j$, and hence
\[
    \varphi-\sup_X\varphi\leq \psi.
\]
In particular, $\int_X\theta_\psi^n>0$.

Now \cref{prop:ds0char}, \cref{prop:vol_limit_model} and
\cref{cor:decnetdS} give
\[
    \varphi_j\xrightarrow{d_S}\psi.
\]
Moreover, observe that $\psi$ is $\mathcal{I}$-model by \cref{prop:decnetmodelPI} and \cref{prop:analysingcompPandPI}. Hence, if $\varphi\sim_{\mathcal I}\psi$, then
\[
    P_\theta[\varphi]_{\mathcal I}=P_\theta[\psi]_{\mathcal I}=\psi.
\]
So it suffices to show that $\varphi\sim_{\mathcal I}\psi$.
It only remains to argue that $\psi \preceq_{\mathcal{I}} \varphi$. For this purpose, take $\lambda>0$. We need to show that
    \[
        \mathcal{I}(\lambda \psi)\subseteq \mathcal{I}(\lambda \varphi).
    \]
    By the strong openness \cref{thm:strongopen}, we may take $\lambda'>\lambda$ such that $\mathcal{I}(\lambda \psi)=\mathcal{I}(\lambda' \psi)$. It follows from the definition of the quasi-equisingular approximation that
    \[
        \mathcal{I}(\lambda'\psi)\subseteq \mathcal{I}(\lambda'\varphi_j)\subseteq \mathcal{I}(\lambda\varphi)
    \]
    for large enough $j$, where the first inclusion uses $\psi\leq P_\theta[\varphi_j]$ and $P_\theta[\varphi_j]\sim_{\mathcal{I}}\varphi_j$. Our assertion follows.

    It follows from the proof that we may take $\varphi_j$ so that $\theta_{\varphi_j}$ is a K\"ahler current for all $j\geq 1$.

    \textbf{Step~2}. We handle the general case.

    Assume \ref{itm:thm:charIgoodasclosure:3} holds.
    By \cref{lma:kahcurrentposmass}, we can find $\psi\in \PSH(X,\theta)$ so that $\theta_{\psi}$ is a Kähler current and $\psi\leq \varphi$. We let
    \[
        \psi_j=(1-j^{-1})\varphi+j^{-1}\psi
    \]
    for each $j\in \mathbb{Z}_{>1}$. Then $(\psi_j)_j$ is an increasing sequence converging almost everywhere to $\varphi$.
    Then
    \[
    P_{\theta}[\psi_j]_{\mathcal{I}}\xrightarrow{d_S}P_{\theta}[\varphi]_{\mathcal{I}}=P_{\theta}[\varphi]
    \]
    by \cref{prop:incnetmodelPI}, \cref{cor:incnetdSconv}.
    From Step~1, we know that each $P_{\theta}[\psi_j]_{\mathcal{I}}$ lies in the closure of analytic singularities, hence so is $P_{\theta}[\varphi]\sim_P \varphi$. Therefore, \ref{itm:thm:charIgoodasclosure:1} follows.
\end{proof}

\begin{definition}\label{def:Igoodpot}
    We say a potential $\varphi\in \QPSH(X)$ is \emph{$\mathcal{I}$-good}\index{potential!I-good potential@$\mathcal{I}$-good potential} if for some smooth closed real $(1,1)$-form $\theta$ on $X$ such that $\varphi\in \PSH(X,\theta)_{>0}$, we have
    \begin{equation}\label{eq:envelopeeq}
    P_{\theta}[\varphi]=P_{\theta}[\varphi]_{\mathcal{I}}.
    \end{equation}
\end{definition}

\begin{remark}\label{rmk:Macron}
    In view of \cref{thm:charIgoodasclosure} and \cref{cor:compnppmassandvol}, the failure of $\mathcal{I}$-goodness of a given $\varphi\in \PSH(X,\theta)_{>0}$ can be characterized using the difference between the volume and the mass. We therefore introduce
    \[
    \Macron(\theta_{\varphi})\coloneqq \vol \theta_{\varphi}-\int_X \theta_{\varphi}^n.
    \]
    \index{$\Macron$}
\end{remark}

An immediate question is to verify that \cref{def:Igoodpot} is independent of the choice of $\theta$.
\begin{lemma}\label{lma:Igoodinsenspert}
Let $\varphi\in \PSH(X,\theta)_{>0}$ for some smooth closed real $(1,1)$-form $\theta$ on $X$. Take a Kähler form $\omega$ on $X$. Then the following are equivalent:
\begin{enumerate}
    \item\label{itm:lma:Igoodinsenspert:1} $P_{\theta}[\varphi]=P_{\theta}[\varphi]_{\mathcal{I}}$;
    \item\label{itm:lma:Igoodinsenspert:2} $P_{\theta+\omega}[\varphi]=P_{\theta+\omega}[\varphi]_{\mathcal{I}}$.
\end{enumerate}
\end{lemma}
\begin{proof}
    \ref{itm:lma:Igoodinsenspert:1} $\implies$ \ref{itm:lma:Igoodinsenspert:2}. By \cref{thm:charIgoodasclosure}, we can find a sequence $(\varphi_j)_j$ in $\PSH(X,\theta)$ with analytic singularities such that $\varphi_j\xrightarrow{d_{S,\theta}}\varphi$. By \cref{cor:dSconv_changetheta}, we have $\varphi_j\xrightarrow{d_{S,\theta+\omega}}\varphi$. Therefore, by \cref{thm:charIgoodasclosure} again, \ref{itm:lma:Igoodinsenspert:2} holds.

    \ref{itm:lma:Igoodinsenspert:2} $\implies$ \ref{itm:lma:Igoodinsenspert:1}. Suppose that \ref{itm:lma:Igoodinsenspert:1} fails, so that
\[
\int_X (\theta+\ddc \varphi)^n< \int_X (\theta+\ddc P_{\theta}[\varphi]_{\mathcal{I}})^n.
\]
Since $\varphi \preceq P_{\theta}[\varphi]_{\mathcal I}$, the monotonicity theorem \cref{thm:mono} gives
\[
    \int_X \theta_{\varphi}^i\wedge\omega^{n-i}\leq \int_X \theta_{P_{\theta}[\varphi]_{\mathcal I}}^i\wedge\omega^{n-i}
\]
for all $i=0,\ldots,n$, and the inequality is strict for $i=n$.
It follows that
\[
\begin{aligned}
\int_X (\theta+\omega+\ddc \varphi)^n
=& \sum_{i=0}^n\binom{n}{i}\int_X \theta_{\varphi}^i\wedge \omega^{n-i}\\
<& \sum_{i=0}^n\binom{n}{i}\int_X \theta_{P_{\theta}[\varphi]_{\mathcal{I}}}^i\wedge \omega^{n-i}\\
=& \int_X (\theta+\omega+\ddc P_{\theta}[\varphi]_{\mathcal{I}})^n
\\
\leq &\int_X (\theta+\omega+\ddc P_{\theta+\omega}[\varphi]_{\mathcal{I}})^n.
\end{aligned}
\]
So \ref{itm:lma:Igoodinsenspert:2} fails as well.
\end{proof}

\begin{corollary}\label{cor:Igoodclosed}
    Let $\theta$ be a closed real smooth $(1,1)$-form on $X$ representing a big cohomology class.
    Let $(\varphi_j)_{j\in I}$ be a net of $\mathcal{I}$-good potentials in $\PSH(X,\theta)$ such that $\varphi_j\xrightarrow{d_S}\varphi$. Then $\varphi$ is $\mathcal{I}$-good.
\end{corollary}
Note that we do not need to assume that $\varphi\in \PSH(X,\theta)_{>0}$.
\begin{proof}
    By adding a Kähler form to $\theta$, and using \cref{lma:Igoodinsenspert,cor:dSconv_changetheta}, we may assume that $\varphi_j,\varphi\in \PSH(X,\theta)_{>0}$ for all $j\in I$. It follows from \cref{thm:charIgoodasclosure} that
    \[
        \int_X \theta_{\varphi_j}^n=\vol \theta_{\varphi_j}
    \]
    for all $j\in I$. Taking limit with respect to $j$ with the help of \cref{thm:contvolu}, we conclude that
    \[
        \int_X \theta_{\varphi}^n=\vol \theta_{\varphi}.
    \]
    Therefore, by \cref{thm:charIgoodasclosure} again, we find that $\varphi$ is $\mathcal{I}$-good.
\end{proof}

\begin{example}\label{ex:analyIgood}
    Assume that $\varphi\in \QPSH(X)$ has analytic singularities. Then $\varphi$ is $\mathcal{I}$-good. This is proved in \cref{prop:analysingcompPandPI}.

    In particular, the potential in \cref{ex:logsingP1} is $\mathcal{I}$-good.
\end{example}

\begin{example}\label{ex:ImodelIgood}
Let $\theta$ be a closed real smooth $(1,1)$-form on $X$ representing a big cohomology class, and let $\varphi\in \PSH(X,\theta)_{>0}$\footnote{It is not known to us whether the same holds when $\varphi$ has vanishing mass.} be an $\mathcal{I}$-model potential. Then $\varphi$ is $\mathcal{I}$-good.

In particular, $V_{\theta}$ is $\mathcal{I}$-good.
\end{example}

\begin{example}\label{ex:fullmassIgood}
Let $\theta$ be a closed real smooth $(1,1)$-form on $X$ representing a big cohomology class, and let $\varphi\in \mathcal{E}(X,\theta)$. Then $\varphi$ is $\mathcal{I}$-good. In fact, since $P_{\theta}[\varphi]=V_{\theta}$, we deduce that $P_{\theta}[\varphi]_{\mathcal{I}}=V_{\theta}$ as well.

    In particular, the potential in \cref{ex:loglogsing} is $\mathcal{I}$-good.
\end{example}

A further class of examples of $\mathcal{I}$-good singularities will be given in \cref{ex:toricIgood} below.

There do exist, however, non-$\mathcal{I}$-good potentials.
\begin{example}\label{ex:Igood:L228}
    The potential in \cref{ex:BBJ} is not $\mathcal{I}$-good. In fact, since all Lelong numbers of $\varphi$ vanish, we know that $\varphi\sim_{\mathcal{I}}0$, hence
    \[
    P_{2\omega}[\varphi]_{\mathcal{I}}=0.
    \]
    On the other hand,
    \[
    \int_X (2\omega+\ddc \varphi)^1=\int_X \omega<\int_X (2\omega),
    \]
    where $(2\omega+\ddc \varphi)^1$ is understood in the non-pluripolar sense.
\end{example}

Quasi-equisingular approximations and $d_S$-convergent sequences are related in the following manner:
\begin{corollary}\label{cor:quasi-equichar}
    Let $\varphi\in \PSH(X,\theta)_{>0}$ and $(\epsilon_j)_j$ be a decreasing sequence in $\mathbb{R}_{\geq 0}$ with limit $0$. Fix a Kähler form $\omega$ on $X$.
    Consider a decreasing sequence $(\varphi_j)_{j>0}$, where each $\varphi_j\in \PSH(X,\theta+\epsilon_j\omega)$ is a potential with analytic singularities. Assume that $\varphi=\inf_j \varphi_j$. Then the following are equivalent:
    \begin{enumerate}
        \item\label{itm:cor:quasi-equichar:1} $\varphi_j\xrightarrow{d_S}P_{\theta}[\varphi]_{\mathcal{I}}$\footnote{Just to be sure, this means $\varphi_j\xrightarrow{d_{S,\theta+\epsilon \omega}}P_{\theta}[\varphi]_{\mathcal{I}}$ for any $\epsilon>0$. The choice of $\epsilon$ is irrelevant due to \cref{cor:dSconv_changetheta}.}, and
        \item\label{itm:cor:quasi-equichar:2} $(\varphi_j)_j$ is a quasi-equisingular approximation of $\varphi$.
    \end{enumerate}
\end{corollary}
\begin{proof}
    By \cref{cor:dSconv_changetheta} and \cref{ex:ImodelIgood}, we may replace $\theta$ by $\theta+C\omega$ for some large constant $C>0$ and assume that $\varphi,\varphi_j\in \PSH(X,\theta-\omega)$ for all $j\geq 1$.

    Here the limit in \ref{itm:cor:quasi-equichar:1} is unchanged by this replacement. Indeed, $P_{\theta}[\varphi]_{\mathcal I}$ is $\mathcal I$-model, hence $\mathcal I$-good by \cref{ex:ImodelIgood}. This implies
\[
    P_{\theta+C\omega}\bigl[P_{\theta}[\varphi]_{\mathcal I}\bigr]=P_{\theta+C\omega}\bigl[P_{\theta}[\varphi]_{\mathcal I}\bigr]_{\mathcal I}=P_{\theta+C\omega}[\varphi]_{\mathcal I}.
\]
Thus the old and new targets are $P$-equivalent in the enlarged class, and \cref{prop:ds0char,cor:dSconv_changetheta} shows that condition \ref{itm:cor:quasi-equichar:1} is
unchanged.

    \ref{itm:cor:quasi-equichar:2} $\implies$ \ref{itm:cor:quasi-equichar:1}. This is already proved in the proof of \cref{thm:charIgoodasclosure}.

    \ref{itm:cor:quasi-equichar:1} $\implies$ \ref{itm:cor:quasi-equichar:2}. This follows from \cref{thm:equising_cond_general}.
\end{proof}

\begin{corollary}\label{cor:quasi_equising_bimero}
    Let $\pi\colon Y\rightarrow X$ be a proper bimeromorphic morphism from a compact K\"ahler manifold $Y$. Let $\varphi\in \PSH(X,\theta)$ and let $(\varphi_j)_j$ be a quasi-equisingular approximation of $\varphi$. Then $(\pi^*\varphi_j)_j$ is a quasi-equisingular approximation of the $\pi^*\theta$-psh function $\pi^*\varphi$.
\end{corollary}
This result is by no means obvious from the definition of quasi-equisingular approximations.
\begin{proof}
    Let $(\epsilon_j)_j$ and $\omega$ be as in \cref{def:quasiequsing}, so that $\epsilon_j\searrow 0$ and $\varphi_j\in\PSH(X,\theta+\epsilon_j\omega)$ for all $j$. Choose a K\"ahler form $\tilde\omega$ on $Y$ and a constant $C>0$ such that $\pi^*\omega\leq C\tilde\omega$. Then
    \[
        \pi^*\varphi_j\in \PSH(Y,\pi^*\theta+C\epsilon_j\tilde\omega).
    \]
    Moreover, the sequence $(\pi^*\varphi_j)_j$ is decreasing, has analytic singularities, and decreases to $\pi^*\varphi$.

    It remains to verify the quasi-equisingular condition \itemref{itm:def:quasiequsing:4}. Fix a constant $A>0$. Then $(\varphi_j)_j$ is also a quasi-equisingular approximation of $\varphi$ with respect to $\theta+A\omega$. Since $\theta_{\varphi}+A\omega$ is a K\"ahler current, \cref{cor:quasi-equichar} gives
    \[
        \varphi_j\xrightarrow{d_S}P_{\theta+A\omega}[\varphi]_{\mathcal I}.
    \]
    Viewing the convergence in a fixed enlarged background class, pulling back and using \cref{cor:dSconv_bimero,prop:Ienvelopebimero}, we get
    \[
        \pi^*\varphi_j\xrightarrow{d_S}\pi^*P_{\theta+A\omega}[\varphi]_{\mathcal I}
        =P_{\pi^*(\theta+A\omega)}[\pi^*\varphi]_{\mathcal I}.
    \]
    The domination $\pi^*\omega\leq C\tilde\omega$ also gives
    \[
        \pi^*\varphi_j\in \PSH\bigl(Y,\pi^*(\theta+A\omega)+C\epsilon_j\tilde\omega\bigr),
    \]
    so the preceding convergence is in the setting of \cref{cor:quasi-equichar} on $Y$. The potential $\pi^*\varphi$ has positive mass with respect to $\pi^*(\theta+A\omega)$ by the bimeromorphic invariance of non-pluripolar products \cref{prop:nppmassinv}. Applying \cref{cor:quasi-equichar} again on $Y$ shows that $(\pi^*\varphi_j)_j$ satisfies the quasi-equisingular condition \itemref{itm:def:quasiequsing:4} for $\pi^*\varphi$. Combining this with the first paragraph proves the assertion with respect to $\pi^*\theta$.
\end{proof}

\section{Properties of \texorpdfstring{$\mathcal{I}$}{I}-good singularities}\label{sec:Igoodprop}
Let $X$ be a connected compact Kähler manifold of dimension $n$.

We show that $\mathcal{I}$-goodness is preserved by a number of natural operations.
\begin{proposition}\label{prop:Igoodlinear}
    Let $\varphi,\psi\in \QPSH(X)$ be $\mathcal{I}$-good and $\lambda>0$. Then the following potentials are all $\mathcal{I}$-good:
    \begin{enumerate}
        \item\label{itm:prop:Igoodlinear:1} $\varphi+\psi$;
        \item\label{itm:prop:Igoodlinear:2} $\varphi\lor \psi$;
        \item\label{itm:prop:Igoodlinear:3} $\lambda\varphi$.
    \end{enumerate}
\end{proposition}
\begin{proof}
    Take a closed real smooth $(1,1)$-form $\theta$ on $X$ such that $\varphi,\psi\in \PSH(X,\theta)_{>0}$. It follows from \cref{thm:charIgoodasclosure} that there are sequences $(\varphi_j)_j,(\psi_j)_j$ in $\PSH(X,\theta)$ with analytic singularities such that $\varphi_j\xrightarrow{d_S}\varphi$ and $\psi_j\xrightarrow{d_S}\psi$.

    By \cref{thm:dSadditivity}, \cref{prop:lor_dS_conv}, we have
    \[
        \varphi_j+\psi_j\xrightarrow{d_S}\varphi+\psi,\quad \varphi_j\lor\psi_j\xrightarrow{d_S}\varphi\lor\psi.
    \]
    By \cref{prop:analysingclosed}, the potentials $\varphi_j+\psi_j$ and $\varphi_j\lor\psi_j$ have analytic singularities. So $\varphi+\psi$ and $\varphi\lor\psi$ are $\mathcal{I}$-good by \cref{thm:charIgoodasclosure}.

    On the other hand, \cref{lma:dS_res} gives
    \[
        \lambda \varphi_j\xrightarrow{d_S}\lambda \varphi.
    \]
    Here the convergence is understood in the scaled class. So by \cref{thm:charIgoodasclosure}, in order to show that $\lambda\varphi$ is $\mathcal{I}$-good, it suffices to show that $\lambda \varphi_j$ is $\mathcal{I}$-good. In other words, we may assume that $\varphi$ has analytic singularities. We may further assume that $\varphi\leq 0$.

    In this case, take a decreasing sequence of rational numbers $\lambda_j\searrow \lambda$. Take a sufficiently large Kähler form $\omega$ on $X$ so that $\lambda_j\varphi,\lambda\varphi\in \PSH(X,\omega)_{>0}$ for any $j$. Then observe that
    \[
        \lambda\varphi=\uscsup[j>0] \lambda_j \varphi.
    \]
    Therefore, $\lambda_j\varphi\xrightarrow{d_S} \lambda \varphi$ by \cref{cor:incnetdSconv}. Each $\lambda_j \varphi$ has analytic singularities, and hence is $\mathcal{I}$-good by \cref{ex:analyIgood}. Therefore, $\lambda\varphi$ is $\mathcal{I}$-good by \cref{thm:charIgoodasclosure}.
\end{proof}

\begin{example}\label{ex:eleIgood}
    Let $L$ be a pseudo-effective line bundle on $X$. Elementary metrics on $L$ are defined in \cref{def:elemetric}. Let $h$ be an elementary metric on $L$. Then $\ddc h$ is $\mathcal{I}$-good.

    This is a direct consequence of \cref{prop:Igoodlinear} and \cref{ex:analyIgood}, and the fact that $\mathcal{I}$-goodness depends only on the $P$-singularity type.
\end{example}

\begin{proposition}\label{prop:Igoodsup}
    Let $(\varphi_j)_{j\in I}$ be a non-empty family of $\mathcal{I}$-good potentials in $\PSH(X,\theta)$ for some closed real smooth $(1,1)$-form $\theta$ on $X$. Assume that this family is locally uniformly bounded from above. Then $\uscsup[j\in I] \varphi_j$ is $\mathcal{I}$-good.
\end{proposition}
\begin{proof}
    After adding a K\"ahler form to $\theta$, we may assume that $\varphi_j\in \PSH(X,\theta)_{>0}$ for all $j\in I$.

    When $I$ is finite, this result follows from \cref{prop:Igoodlinear}.
    When $I$ is infinite, by Choquet's lemma \cref{prop:Choquet}, we may assume that $I=\mathbb{Z}_{>0}$. Replacing $\varphi_j$ by $\varphi_1\lor\cdots\lor\varphi_j$ and using \cref{prop:Igoodlinear}, we may assume that the sequence $(\varphi_j)_j$ is increasing. In this case, as shown in \cref{cor:incnetdSconv},
    \[
        \varphi_j\xrightarrow{d_S}\uscsup[ i\in \mathbb{Z}_{>0}] \varphi_i.
    \]
    Therefore, $\uscsup[i\in \mathbb{Z}_{>0}] \varphi_i$ is $\mathcal{I}$-good by \cref{cor:Igoodclosed}.
\end{proof}

\section{Mixed volumes}\label{sec:mix_vol}

We first extend the notion of volume in \cref{def:volqpsh} to the mixed case.
Let $\theta_1,\ldots,\theta_n$ be smooth closed real $(1,1)$-forms on $X$ representing pseudo-effective classes.

\begin{definition}\label{def:mix_vol2}
    Let $\varphi_i\in \PSH(X,\theta_i)$ for $i=1,\ldots,n$. Write $T_i=\theta_i+\ddc \varphi_i$ for each $i=1,\ldots,n$. Assume that $\vol T_i>0$ for all $i=1,\ldots,n$.
    We define the \emph{mixed volume}\index{mixed volume} $\vol(T_1,\ldots,T_n)$\index{$\vol(T_1,\ldots,T_n)$} as follows:
            \begin{equation}\label{eq:volmixed}
                \vol(T_1,\ldots,T_n)=\int_X \left(\theta_1+\ddc P_{\theta_1}[\varphi_1]_{\mathcal{I}}\right)\wedge \dots\wedge \left(\theta_n+\ddc P_{\theta_n}[\varphi_n]_{\mathcal{I}}\right).
            \end{equation}
\end{definition}
When $\vol T_i>0$ for each $i=1,\ldots,n$, the definition \eqref{eq:volmixed} does not depend on how we represent $T_i$ as $T_i=\theta_i+\ddc \varphi_i$, this is a consequence of \cref{thm:mono} and \cref{prop:Ienvindeptheta}.

\begin{lemma}\label{lma:vollinearlimit2}\leavevmode
    \begin{enumerate}
        \item\label{itm:lma:vollinearlimit2:1} The mixed volume is symmetric in its $n$-variables.
        \item\label{itm:lma:vollinearlimit2:2} For any Kähler form $\omega$ on $X$, we have
        \begin{equation}\label{eq:mix_vol_comp1}
            \vol\left(T_1,\ldots,T_n\right)=\lim_{\epsilon\to 0+} \vol\left(T_1+\epsilon\omega,\ldots,T_n+\epsilon\omega\right).
        \end{equation}
        \item\label{itm:lma:vollinearlimit2:3} Consider closed positive $(1,1)$-currents $T_1,\ldots,T_n,T_1'$ on $X$ with positive volumes, then
        \begin{equation}
            \vol (T_1+T_1',T_2,\ldots,T_n)=\vol (T_1,T_2,\ldots,T_n)+\vol (T_1',T_2,\ldots,T_n).
        \end{equation}
        \item\label{itm:lma:vollinearlimit2:4} Consider closed positive $(1,1)$-currents $T_1,\ldots,T_n$ on $X$ with positive volumes and $\lambda > 0$, then
        \[
        \vol (\lambda T_1,T_2,\ldots,T_n)=\lambda \vol (T_1,T_2,\ldots,T_n).
        \]
        \item\label{itm:lma:vollinearlimit2:5} Suppose that $T_1,\ldots,T_n,S_1,\ldots,S_n$ are closed positive $(1,1)$-currents on $X$ with positive volumes such that $T_i\preceq_{\mathcal{I}} S_i$ and $\{T_i\}=\{S_i\}$ for each $i=1,\ldots,n$. Then
        \begin{equation}
            \vol(T_1,\ldots,T_n)\leq \vol(S_1,\ldots,S_n).
        \end{equation}
    \end{enumerate}
\end{lemma}
The notation $\Reg$ is defined in \eqref{eq:Siudec}.
\begin{proof}
    \ref{itm:lma:vollinearlimit2:1} This follows from the symmetry of the mixed non-pluripolar product in \eqref{eq:volmixed}, cf. \itemref{itm:prop:npp1:2}.

    \ref{itm:lma:vollinearlimit2:2} For each $i$ and each $\epsilon>0$, we have
    \[
        P_{\theta_i}[\varphi_i]_{\mathcal{I}}\sim_P P_{\theta_i+\epsilon\omega}\left[P_{\theta_i}[\varphi_i]_{\mathcal{I}}\right]_{\mathcal{I}}=P_{\theta_i+\epsilon\omega}[\varphi_i]_{\mathcal{I}}
    \]
    as a consequence of \cref{ex:ImodelIgood}. Hence using \cref{prop:mono2} and the multi-linearity of the non-pluripolar product \itemref{itm:prop:npp1:6},
    \[
        \begin{aligned}
            &\lim_{\epsilon\to 0+}\vol(T_1+\epsilon\omega,\ldots,T_n+\epsilon\omega)\\
            =&\lim_{\epsilon\to 0+}\int_X \left(\theta_1+\epsilon\omega+\ddc P_{\theta_1}[\varphi_1]_{\mathcal{I}}\right)\wedge \cdots \wedge \left(\theta_n+\epsilon\omega+\ddc P_{\theta_n}[\varphi_n]_{\mathcal{I}}\right)\\
            =& \vol(T_1,\ldots,T_n).
        \end{aligned}
    \]

    \ref{itm:lma:vollinearlimit2:3}
    We write $T_i=\theta_i+\ddc \varphi_i$ as before for each $i$ and $T_1'=\theta_1'+\ddc \varphi_1'$.
    The potentials $P_{\theta_1}[\varphi_1]_{\mathcal I}$ and $P_{\theta_1'}[\varphi_1']_{\mathcal I}$ are $\mathcal I$-model with positive masses, hence $\mathcal I$-good. By \cref{prop:Igoodlinear}, their sum is $\mathcal I$-good. Moreover, it is $\mathcal I$-equivalent to $\varphi_1+\varphi_1'$. Then thanks to \cref{prop:Igoodlinear},
    \[
    P_{\theta_1}[\varphi_1]_{\mathcal{I}}+P_{\theta_1'}[\varphi_1']_{\mathcal{I}}\sim_P P_{\theta_1+\theta_1'}\left[P_{\theta_1}[\varphi_1]_{\mathcal{I}}+P_{\theta_1'}[\varphi_1']_{\mathcal{I}} \right]_{\mathcal{I}} = P_{\theta_1+\theta_1'}\left[\varphi_1+\varphi_1'\right]_{\mathcal{I}}
    \]
    Thus by \cref{prop:mono2}, we have
    \[
    \begin{aligned}
        &\vol (T_1+T_1',T_2,\ldots,T_n)\\
        =&\int_X \left( \theta_1+\theta_1'+\ddc P_{\theta_1+\theta_1'}[\varphi_1+\varphi_1']_{\mathcal{I}} \right)\wedge \left( \theta_2+\ddc P_{\theta_2}[\varphi_2]_{\mathcal{I}}\right)\wedge \cdots \\
        &\wedge  \left( \theta_n+\ddc P_{\theta_n}[\varphi_n]_{\mathcal{I}}\right)\\
        =&\int_X \left( \theta_1+\theta_1'+\ddc P_{\theta_1}[\varphi_1]_{\mathcal{I}}+\ddc P_{\theta_1'}[\varphi_1']_{\mathcal{I}} \right)\wedge \left( \theta_2+\ddc P_{\theta_2}[\varphi_2]_{\mathcal{I}}\right)\wedge \cdots \\
        &\wedge  \left( \theta_n+\ddc P_{\theta_n}[\varphi_n]_{\mathcal{I}}\right)\\
        =& \vol (T_1,T_2,\ldots,T_n)+\vol (T_1',T_2,\ldots,T_n).
    \end{aligned}
    \]

    \ref{itm:lma:vollinearlimit2:4} This follows from \eqref{eq:volmixed}, \itemref{itm:prop:PIconc:1} and \itemref{itm:prop:npp1:6}.

    \ref{itm:lma:vollinearlimit2:5} Our assertion follows from the monotonicity theorem \cref{thm:mono}.
\end{proof}

In view of \itemref{itm:lma:vollinearlimit2:2}, we may extend the definition of mixed volumes to the more general setting as follows:
\begin{definition}\label{def:mix_vol}
    Let $\varphi_i\in \PSH(X,\theta_i)$ for $i=1,\ldots,n$. Write $T_i=\theta_i+\ddc \varphi_i$ for each $i=1,\ldots,n$. Take a Kähler form $\omega$ on $X$. We define the \emph{mixed volume}\index{mixed volume} $\vol(T_1,\ldots,T_n)$\index{$\vol(T_1,\ldots,T_n)$} as
        \begin{equation}\label{eq:volmixedgeneral}
            \vol\left(T_1,\ldots,T_n\right)=\lim_{\epsilon\to 0+}\vol\left(T_1+\epsilon\omega,\ldots,T_n+\epsilon\omega\right).
        \end{equation}
    Note that $\vol(T_1+\epsilon\omega,\ldots,T_n+\epsilon\omega)$ is decreasing as $\epsilon \to 0+$ due to \itemref{itm:lma:vollinearlimit2:3}, so the limit exists.

    Furthermore, $\vol(T_1,\ldots,T_n)$ does not depend on the choice of $\omega$, as a consequence of \itemref{itm:lma:vollinearlimit2:3}.
\end{definition}

\begin{proposition}\label{prop:vollinearlimit}\leavevmode
    \begin{enumerate}
        \item\label{itm:prop:vollinearlimit:1} The mixed volume is symmetric in its $n$-variables.
        \item\label{itm:prop:vollinearlimit:2} Consider closed positive $(1,1)$-currents $T_1,\ldots,T_n,T_1'$ on $X$, then
        \begin{equation}
            \vol (T_1+T_1',T_2,\ldots,T_n)=\vol (T_1,T_2,\ldots,T_n)+\vol (T_1',T_2,\ldots,T_n).
        \end{equation}
        \item\label{itm:prop:vollinearlimit:3} Consider closed positive $(1,1)$-currents $T_1,\ldots,T_n$ on $X$ and $\lambda\geq 0$, then
        \[
        \vol (\lambda T_1,T_2,\ldots,T_n)=\lambda \vol (T_1,T_2,\ldots,T_n).
        \]
        \item\label{itm:prop:vollinearlimit:4} Suppose that $T_1,\ldots,T_n,S_1,\ldots,S_n$ are closed positive $(1,1)$-currents on $X$ such that $T_i\preceq_{\mathcal{I}} S_i$ and $\{T_i\}=\{S_i\}$ for each $i=1,\ldots,n$. Then
        \begin{equation}
            \vol(T_1,\ldots,T_n)\leq \vol(S_1,\ldots,S_n).
        \end{equation}
        \item\label{itm:prop:vollinearlimit:5} Suppose that $T_1,\ldots,T_n$ are closed positive $(1,1)$-currents on $X$. Then
        \begin{equation}
            \vol\left(T_1,\ldots,T_n \right)=\vol\left(\Reg T_1,\ldots,\Reg T_n \right).
        \end{equation}
        \item\label{itm:prop:vollinearlimit:6} For any closed positive $(1,1)$-current $T$ on $X$, we have
        \begin{equation}\label{eq:mix_vol_ext_pure}
            \vol(T,\ldots,T)=\vol T.
        \end{equation}
    \end{enumerate}
\end{proposition}

\begin{proof}
    \ref{itm:prop:vollinearlimit:1}, \ref{itm:prop:vollinearlimit:2}, \ref{itm:prop:vollinearlimit:3}, \ref{itm:prop:vollinearlimit:4} are immediate consequences of the corresponding properties proved in \cref{lma:vollinearlimit2}.

    \ref{itm:prop:vollinearlimit:5} Using \ref{itm:prop:vollinearlimit:3} and \ref{itm:prop:vollinearlimit:4}, it suffices to establish the following: Suppose that
    \begin{equation}\label{eq:T_count_sum_temp1}
        T=\sum_i c_i[E_i]
    \end{equation}
    is a closed positive $(1,1)$-current on $X$, where $\{E_i\}$ is a countable collection of prime divisors on $X$ and $c_i>0$. Then
    \begin{equation}\label{eq:volTT2Tn_temp1}
    \vol\left(T,T_2,\ldots,T_n \right)=0
    \end{equation}
    for any closed positive $(1,1)$-currents $T_2,\ldots,T_n$ on $X$.

    \textbf{Step~1}. We first assume that \eqref{eq:T_count_sum_temp1} is a finite sum. Fix a Kähler form $\omega$ on $X$.

    In this case, thanks to \ref{itm:lma:vollinearlimit2:3} and \ref{itm:lma:vollinearlimit2:4} again, we may assume that $T=[E]$ for some prime divisor $E$ on $X$. Write $[E]=\theta_{\varphi}$. Then $\varphi$ has analytic singularities thanks to \cref{prop:LelongPoincare}. Therefore, $P_{\theta+\epsilon\omega}[\varphi]_{\mathcal{I}}\sim \varphi$ for any $\epsilon>0$ due to \cref{prop:analysingcompPandPI}. Write $T_i=\theta_i+\ddc \varphi_i$ for $i=2,\ldots,n$, and choose $C>0$ so that $C\omega+\theta_i$ are Kähler forms for each $i=2,\ldots,n$. Then
    \[
    \begin{aligned}
    &\vol(T,T_2,\ldots,T_n)\\
    =& \lim_{\epsilon\to 0+} \int_X \left([E]+\epsilon\omega\right)\wedge \left( \theta_2+\epsilon\omega+\ddc P_{\theta_2+\epsilon\omega}[\varphi_2]_{\mathcal{I}}\right)\wedge \cdots \\
    &\wedge \left( \theta_n+\epsilon\omega+\ddc P_{\theta_n+\epsilon\omega}[\varphi_n]_{\mathcal{I}}\right)\\
    =&\lim_{\epsilon\to 0+} \epsilon\int_X \omega \wedge \left( \theta_2+\epsilon\omega+\ddc P_{\theta_2+\epsilon\omega}[\varphi_2]_{\mathcal{I}}\right)\wedge \cdots \\
    &\wedge \left( \theta_n+\epsilon\omega+\ddc P_{\theta_n+\epsilon\omega}[\varphi_n]_{\mathcal{I}}\right)\\
    \leq &\lim_{\epsilon\to 0+} \epsilon\int_X \omega \wedge \left( \theta_2+C\omega\right)\wedge \cdots \wedge \left( \theta_n+C\omega\right)\\
    =&0,
    \end{aligned}
    \]
    where the inequality step follows from the monotonicity theorem \cref{thm:mono} and the fact that $P_{\theta_i+\epsilon\omega}[\varphi_i]_{\mathcal{I}}\leq 0$ for all $i=2,\ldots,n$.

    \textbf{Step~2}. We prove the case where $\{E_i\}$ is infinite. We may assume that $i$ runs over $\mathbb{Z}_{>0}$.

    Write $T_i=\theta_{i}+\ddc\varphi_i$ for $i=2,\ldots,n$ as before.
    Fix a Kähler form $\omega$ on $X$ so that $\theta_i+\omega$ is a Kähler form for each $i=2,\ldots,n$. Fix $\epsilon>0$.

    We can find $N_0>0$ so that for any $N\geq N_0$, the class of
    \[
    \epsilon\omega+ \sum_{i=N+1}^{\infty}c_i[E_i]
    \]
    is Kähler. Take a Kähler form $\omega_N$ in this class.
    Then the currents
    \[
    \sum_{i=1}^N c_i [E_i]+\omega_N,\quad \sum_{i=1}^{\infty} c_i [E_i]+\epsilon\omega
    \]
    all lie in the same cohomology class.
    Moreover,
    \[
        \sum_{i=1}^{\infty} c_i [E_i]+\epsilon\omega\preceq_{\mathcal{I}}\sum_{i=1}^N c_i [E_i]+\omega_N.
    \]
    Therefore, by the monotonicity in \ref{itm:lma:vollinearlimit2:5}, followed by the finite divisorial case of Step~1,
    \[
    \begin{aligned}
    \vol\left(T,T_2,\ldots,T_n\right)\leq & \vol\left(T+\epsilon\omega,T_2,\ldots,T_n\right)\\
    \leq & \vol\left(\sum_{i=1}^N c_i [E_i]+\omega_N,T_2,\ldots,T_n\right)\\
    =&  \vol\left(\omega_N,T_2,\ldots,T_n\right)\\
    \leq & \int_X \omega_N\wedge (\theta_2+\omega)\wedge \cdots\wedge (\theta_n+\omega).
    \end{aligned}
    \]
    Letting $N\to\infty$, and using $\{\omega_N\}\to \{\epsilon\omega\}$, we get
    \[
    \vol\left(T,T_2,\ldots,T_n\right)\leq \epsilon \int_X \omega\wedge (\theta_2+\omega)\wedge \cdots\wedge (\theta_n+\omega).
    \]
    Since $\epsilon>0$ is arbitrary, we conclude \eqref{eq:volTT2Tn_temp1}.

    \ref{itm:prop:vollinearlimit:6} Write $T=\theta_{\varphi}$. In more concrete terms, we need to show that
    \[
    \lim_{\epsilon\to 0+}\int_X (\theta+\epsilon\omega+\ddc P_{\theta+\epsilon\omega}[\varphi]_{\mathcal{I}})^n=\int_X (\theta+\ddc P_{\theta}[\varphi]_{\mathcal{I}})^n.
    \]

    We may replace $\varphi$ by $P_{\theta}[\varphi]_{\mathcal{I}}$ and assume that $\varphi$ is $\mathcal{I}$-model in $\PSH(X,\theta)$. Then we claim that
    \[
    \varphi=\inf_{\epsilon>0}P_{\theta+\epsilon\omega}[\varphi]_{\mathcal{I}}.
    \]
    From this, our assertion follows from \cref{prop:vol_limit_model}.

    The $\leq$ direction is clear. For the converse, since $\varphi$ is $\mathcal{I}$-model, it suffices to show that for each prime divisor $E$ over $X$, we have
    \[
    \nu(\varphi,E)\leq \nu\left(\inf_{\epsilon>0}P_{\theta+\epsilon\omega}[\varphi]_{\mathcal{I}},E\right).
    \]
    We simply compute
    \[
    \nu\left(\inf_{\epsilon>0}P_{\theta+\epsilon\omega}[\varphi]_{\mathcal{I}},E\right)\geq \sup_{\epsilon>0}\nu\left(P_{\theta+\epsilon\omega}[\varphi]_{\mathcal{I}},E\right)=\nu(\varphi,E).
    \]
\end{proof}

\begin{lemma}\label{lma:perturbvol}
    Let $\omega$ be a Kähler form on $X$. Then there is a constant $C>0$ depending only on $X,\omega,\{T_1\},\ldots,\{T_n\}$ such that
    \[
    0\leq \vol(T_1+\epsilon\omega,\ldots, T_n+\epsilon\omega)-\vol(T_1,\ldots,T_n)\leq C\epsilon
    \]
    for any $\epsilon\in [0,1]$.
\end{lemma}
\begin{proof}
    By the multi-linearity established in \cref{prop:vollinearlimit}, we can write
    \[
        \vol\left(T_1+\epsilon\omega,\ldots, T_n+\epsilon\omega\right)-\vol(T_1,\ldots,T_n)
    \]
    as a linear combination of the $2^n-1$-terms:
    \[
        \epsilon^m\vol\left(\omega^m, T_{i_1},\ldots,T_{i_{n-m}} \right)
    \]
    where $J\subseteq \{1,\ldots,n\}$ is a non-empty subset of cardinality $m$, and $\{i_1,\ldots,i_{n-m}\}$ is the complement of $J$.
    It suffices to show that the coefficients are bounded by a constant as above. Take $D\geq 0$ as the minimal non-negative number so that $D\{\omega\}+\{T_i\}$ is nef for each $i=1,\ldots,n$. Then by the monotonicity theorem \cref{thm:mono},
    \[
        \vol\left(\omega^m, T_{i_1},\ldots,T_{i_{n-m}} \right)\leq \{\omega\}^m\cap \left( \{T_{i_1}\}+D\{\omega\} \right) \cap \cdots \cap \left( \{T_{i_{n-m}}\}+D\{\omega\} \right).
    \]
    Our assertion follows.
\end{proof}

Next we show that \cref{thm:contvolu} continues to hold even when $\varphi$ has vanishing mass. We prove a slightly more general result:
\begin{theorem}\label{thm:mvolcont}
    Let $\theta_1,\ldots,\theta_n$ be smooth closed real $(1,1)$-forms representing big cohomology classes. Suppose that $(\varphi_j^k)_{k\in I}$ are nets in $\PSH(X,\theta_j)$ and $\varphi_j\in \PSH(X,\theta_j)$ for $j=1,\ldots,n$. We assume that $\varphi_j^k\xrightarrow{d_S}\varphi_j$ for each $j=1,\ldots,n$. Then
    \begin{equation}\label{eq:convmixedvolds}
        \lim_{k\in I}\vol \left( \theta_{1,\varphi_1^k},\ldots, \theta_{n,\varphi_n^k}\right)=\vol \left( \theta_{1,\varphi_1},\ldots, \theta_{n,\varphi_n}\right).
    \end{equation}
\end{theorem}
\begin{proof}
    Fix a Kähler form $\omega$. Then for any $\epsilon>0$, we have
    \[
    \lim_{k\in I}\vol \left( (\theta_1+\epsilon\omega)_{\varphi_1^k},\ldots, (\theta_n+\epsilon\omega)_{\varphi_n^k}\right)=\vol \left( (\theta_1+\epsilon\omega)_{\varphi_1},\ldots, (\theta_n+\epsilon\omega)_{\varphi_n}\right)
    \]
    as a consequence of \cref{cor:dSconv_changetheta}, \cref{thm:convdS} and \cref{thm:contPI}.

    Now thanks to \cref{lma:perturbvol}, there is $C>0$ so that for each $k\in I$,
    \[
    \begin{aligned}
    0\leq & \vol \left( (\theta_1+\epsilon\omega)_{\varphi_1^k},\ldots, (\theta_n+\epsilon\omega)_{\varphi_n^k}\right)-\vol \left( \theta_{1,\varphi_1^k},\ldots, \theta_{n,\varphi_n^k}\right)\leq C\epsilon,\\
    0\leq & \vol \left( (\theta_1+\epsilon\omega)_{\varphi_1},\ldots, (\theta_n+\epsilon\omega)_{\varphi_n}\right)-\vol \left( \theta_{1,\varphi_1},\ldots, \theta_{n,\varphi_n}\right)\leq C\epsilon.
    \end{aligned}
    \]
    Therefore, \eqref{eq:convmixedvolds} follows.
\end{proof}

\begin{corollary}\label{cor:contvolu2}
    Let $\theta$ be a smooth closed real $(1,1)$-form representing a big cohomology class. Let $(\varphi_j)_{j\in I}$ be a net in $\PSH(X,\theta)$ and $\varphi\in \PSH(X,\theta)$. Assume that $\varphi_j\xrightarrow{d_S}\varphi$. Then
    \begin{equation}\label{eq:Ivolcont2}
        \vol \theta_{\varphi_j}\to \vol \theta_{\varphi},\quad \int_X \theta_{\varphi_j}^n\to \int_X \theta_{\varphi}^n.
    \end{equation}
\end{corollary}

\begin{proof}
    The first part of \eqref{eq:Ivolcont2} is a special case of \cref{thm:mvolcont}, while the second part of \eqref{eq:Ivolcont2} is a special case of \cref{thm:convdS}.
\end{proof}

The mixed volume has a log-concavity property:
\begin{proposition}\label{prop:logvol}
    Let $T_1,\ldots,T_n$ be closed positive $(1,1)$-currents on $X$. Then
    \[
    \vol(T_1,\ldots,T_n)\geq \prod_{i=1}^n(\vol T_i)^{1/n}.
    \]
\end{proposition}
\begin{proof}
    We may assume that $\vol T_i>0$ for each $i=1,\ldots,n$ since there is nothing to prove otherwise. In this case, we need to show that
    \[
    \begin{split}
    &\int_X \left(\theta_1+\ddc P_{\theta_1}[\varphi_1]_{\mathcal{I}}\right)\wedge \cdots \wedge \left(\theta_n+\ddc P_{\theta_n}[\varphi_n]_{\mathcal{I}}\right)\\
    \geq & \prod_{i=1}^n\left(\int_X \left(\theta_i+\ddc P_{\theta_i}[\varphi_i]_{\mathcal{I}} \right)^n \right)^{1/n}.
    \end{split}
    \]
    This is a special case of \cref{thm:logconc}.
\end{proof}

Next, we also study how the mixed volume behaves under bimeromorphic transforms.
\begin{proposition}\label{prop:biminv}
    Let $\pi\colon Y\rightarrow X$ be a proper bimeromorphic morphism from a Kähler manifold $Y$ to $X$. Then
    \[
    \vol\left(\pi^*T_1,\ldots,\pi^*T_n\right)=\vol\left(T_1,\ldots,T_n\right).
    \]
\end{proposition}
\begin{proof}
    Using the definition of the mixed volume, we may easily reduce to the case where $\vol T_i>0$ for each $i=1,\ldots,n$. By \cref{prop:Ienvelopebimero}, we know that if we write $T_i=\theta_i+\ddc\varphi_i$, then
    \[
    \pi^*P_{\theta_i}[\varphi_i]_{\mathcal{I}}=P_{\pi^*\theta_i}[\pi^*\varphi_i]_{\mathcal{I}}
    \]
    for each $i=1,\ldots,n$.
    In particular,
    \[
    \vol \pi^*T_i=\vol T_i>0.
    \]
    Our assertion follows from the bimeromorphic invariance of the non-pluripolar product \cref{prop:nppmassinv}.
\end{proof}

As for pushforward, we also have a similar result. We need a preliminary result:

\begin{lemma}\label{lma:Siudec_pullpush}
    Let $\pi\colon Y\rightarrow X$ be a proper bimeromorphic morphism from a Kähler manifold $Y$.
    Then for any non-divisorial closed positive $(1,1)$-current $T$ on $Y$, we have
    \[
    \pi^* \pi_* T=T+\sum_{i=1}^N c_i [E_i]
    \]
    for finitely many $\pi$-exceptional divisors $E_i$ and $c_i>0$.
\end{lemma}

\begin{proof}
    Let $E$ be the exceptional locus of $\pi$. Then
    \[
    T=\mathds{1}_{Y\setminus E}\pi^* \pi_* T.
    \]
    Therefore,
    \[
    \pi^* \pi_* T-T=\mathds{1}_E\pi^* \pi_* T,
    \]
    which has the stated form, due to the support theorems, see \cite[Section~8]{DemBook}.
\end{proof}

\begin{corollary}\label{cor:vol_inv_push}
    Let $\pi\colon X\rightarrow Z$ be a proper bimeromorphic morphism from $X$ to a Kähler manifold $Z$. Then
    \begin{equation}\label{eq:volpushinc}
    \vol(T_1,\ldots,T_n)=\vol(\pi_*T_1,\ldots,\pi_*T_n).
    \end{equation}
\end{corollary}
\begin{proof}
    Observe that we may assume that $T_i=\Reg T_i$ for all $i=1,\ldots,n$. In fact, since the pushforward of the divisorial part of $T_i$ is divisorial as well, hence by \itemref{itm:prop:vollinearlimit:5}, they do not contribute to the volumes.

    Now by \cref{prop:biminv}, it remains to show that
    \[
    \vol(T_1,\ldots,T_n)=\vol(\pi^*\pi_*T_1,\ldots,\pi^*\pi_*T_n).
    \]
    By \cref{lma:Siudec_pullpush}, the difference $\pi^*\pi_*T_i-T_i$ is divisorial, hence our desired equality follows from \itemref{itm:prop:vollinearlimit:5}.
\end{proof}

A particular corollary of \cref{cor:vol_inv_push} will be useful later.
\begin{corollary}\label{cor:Igoodpush}
     Let $\pi\colon X\rightarrow Z$ be a proper bimeromorphic morphism from $X$ to a Kähler manifold $Z$. Assume that $T$ is an $\mathcal{I}$-good closed positive $(1,1)$-current on $X$. Then so is $\pi_*T$.
\end{corollary}
\begin{proof}
    We may assume that $\int_X T^n>0$. Then by \cref{cor:vol_inv_push},
    \[
    \vol \pi_*T=\vol T>0
    \]
    as well.
    Since $T$ is $\mathcal{I}$-good, we have
    \[
    \vol T=\int_X T^n.
    \]
    But $\int_X T^n=\int_Z (\pi_*T)^n$, so
    \[
    \vol \pi_*T=\int_Z (\pi_*T)^n>0.
    \]
    It follows that $\pi_*T$ is $\mathcal{I}$-good.
\end{proof}

Finally we establish a semicontinuity property of the mixed volumes.
\begin{theorem}\label{thm:usc_mix_vol}
   Let $(\varphi_i^j)_{j\in J}$ ($i=1,\ldots,n$) be nets in $\PSH(X,\theta_i)$, and $\varphi_i\in \PSH(X,\theta_i)$ ($i=1,\ldots,n$). Assume that for each prime divisor $E$ over $X$, we have
   \[
   \lim_{j\in J} \nu\left(\varphi_i^j,E\right)=\nu\left(\varphi_i,E\right),\quad i=1,\ldots,n.
   \]
   Then
   \[
   \varlimsup_{j\in J} \vol \left(\theta_1+\ddc \varphi_1^j, \dots,\theta_n+\ddc \varphi_n^j \right) \leq  \vol\left(\theta_1+\ddc \varphi_1, \dots,\theta_n+\ddc \varphi_n\right).
   \]
\end{theorem}
\begin{proof}
    \textbf{Step~1}. We first assume that $\vol (\theta_i+\ddc \varphi_i^j)>0$ and $\vol (\theta_i+\ddc \varphi_i)>0$ for all $i=1,\ldots,n$ and $j\in J$.

    Without loss of generality, we may assume that the $\varphi_i^j$'s and the $\varphi_i$'s are $\mathcal{I}$-model for all $i=1,\ldots,n$ and $j\in J$. Our assertion becomes
    \begin{equation}\label{eq:volineq_1}
    \varlimsup_{j\in J} \int_X \left(\theta_1+\ddc \varphi_1^j\right)\wedge \dots \wedge \left(\theta_n+\ddc \varphi_n^j \right) \leq  \int_X \left(\theta_1+\ddc \varphi_1\right) \wedge \dots \wedge \left(\theta_n+\ddc \varphi_n\right).
    \end{equation}
    For each $j\in J$, define
    \[
    \psi_i^j\coloneqq \uscsup[k\geq j] \varphi_i^k,\quad i=1,\ldots,n.
    \]
    Observe that $\psi_i^j$ is $\mathcal{I}$-good thanks to \cref{prop:Igoodsup}, since after the preceding reduction the tail family is bounded above by $0$.
    It follows from \cref{cor:supsLelong} and our assumption that
    \[
    \lim_{j\in J}\nu\left(\psi_i^j,E\right)=\nu\left(\varphi_i,E\right),\quad i=1,\ldots,n.
    \]
    Moreover, for each $j\in J$, the same formula gives $\varphi_i\preceq_{\mathcal{I}}\psi_i^j$. Since $\varphi_i$ is $\mathcal{I}$-model and $\psi_i^j$ is $\mathcal{I}$-good, \cref{prop:Icomparandenvelope} gives
    \[
        \varphi_i\leq P_{\theta_i}\left[\psi_i^j\right].
    \]
    For each $i=1,\ldots,n$, we define
    \[
    \psi_i=\inf_{j\in J}P_{\theta_i}\left[\psi_i^j\right].
    \]
    In particular, $\psi_i\geq\varphi_i$, so $\int_X\theta_{i,\psi_i}^n>0$.
    Since $\psi_i^j$ is $\mathcal{I}$-good, $P_{\theta_i}[\psi_i^j]$ is both model and $\mathcal{I}$-model. Due to \cref{prop:decnetmodelPI}, $\psi_i$ is $\mathcal{I}$-model.
    Thanks to \cref{prop:vol_limit_model}, we know
    \[
    \nu(\psi_i,E)=\nu(\varphi_i,E)
    \]
    for any $i=1,\ldots,n$ and any prime divisor $E$ over $X$. In other words, $\psi_i\sim_{\mathcal{I}} \varphi_i$ for $i=1,\ldots,n$. But both $\varphi_i$ and $\psi_i$ are $\mathcal{I}$-good, therefore,
    \[
    \psi_i\sim_{P} \varphi_i,\quad i=1,\ldots,n.
    \]
    By \cref{prop:mono2}, we have
    \[
    \int_X (\theta_1+\ddc \psi_1)\wedge \cdots\wedge (\theta_n+\ddc \psi_n)=\int_X (\theta_1+\ddc \varphi_1)\wedge \cdots\wedge (\theta_n+\ddc \varphi_n).
    \]
    Next by \cref{prop:mono2} again,
    \[
    \begin{split}
    &\varlimsup_{j\in J} \int_X \left(\theta_1+\ddc \varphi_1^j\right)\wedge\dots\wedge  \left(\theta_n+\ddc \varphi_n^j \right)\\
    \leq & \varlimsup_{j\in J} \int_X \left(\theta_1+\ddc \psi_1^j\right)\wedge \dots \wedge \left(\theta_n+\ddc \psi_n^j \right).
    \end{split}
    \]
    On the other hand, due to \cref{prop:vol_limit_model} and \cref{cor:decnetdS}, for each $i=1,\ldots,n$, we have
    \[
    \psi_i^j\xrightarrow{d_S}\psi_i.
    \]
    We conclude from \cref{thm:convdS} that
    \[
    \varlimsup_{j\in J} \int_X \left(\theta_1+\ddc \psi_1^j\right)\wedge \dots \wedge \left(\theta_n+\ddc \psi_n^j \right)=\int_X (\theta_1+\ddc \psi_1)\wedge \cdots\wedge (\theta_n+\ddc \psi_n).
    \]
    Putting these equations together, \eqref{eq:volineq_1} follows.

    \textbf{Step~2}. Next we handle the general case.

    Fix a Kähler form $\omega$ on $X$. For any $\epsilon\in (0,1]$, from Step~1, we know that
    \[
    \begin{split}
   &\varlimsup_{j\in J} \vol \left(\theta_1+\epsilon\omega+\ddc \varphi_1^j, \dots,\theta_n+\epsilon\omega+\ddc \varphi_n^j \right) \\
   \leq & \vol\left(\theta_1+\epsilon\omega+\ddc \varphi_1, \dots,\theta_n+\epsilon\omega+\ddc \varphi_n\right).
   \end{split}
   \]
   Using \cref{lma:perturbvol}, we have
   \[
   \begin{aligned}
     &\varlimsup_{j\in J} \vol \left(\theta_1+\ddc \varphi_1^j, \dots,\theta_n+\ddc \varphi_n^j \right)\\
     \leq & \varlimsup_{j\in J} \vol \left(\theta_1+\epsilon\omega+\ddc \varphi_1^j, \dots,\theta_n+\epsilon\omega+\ddc \varphi_n^j \right)\\
     \leq &\vol\left(\theta_1+\epsilon\omega+\ddc \varphi_1, \dots,\theta_n+\epsilon\omega+\ddc \varphi_n\right)\\
     \leq &\vol\left(\theta_1+\ddc \varphi_1, \dots,\theta_n+\ddc \varphi_n\right)+C\epsilon.
    \end{aligned}
   \]
   But since $\epsilon$ is arbitrary, our assertion follows.
\end{proof}

\section{The volumes of Hermitian pseudo-effective line bundles}\label{sec:volHermitianbig}
Let $X$ be a connected compact Kähler manifold of dimension $n$.

\begin{definition}\label{def:Hermpsef}
    A \emph{Hermitian pseudo-effective line bundle}\index{line bundle!Hermitian pseudo-effective line bundle} $(L,h)$ on a complex manifold $Y$ consists of a holomorphic line bundle $L$ on $Y$ together with a plurisubharmonic metric $h$ on $L$.
\end{definition}

\begin{theorem}[Asymptotic Riemann--Roch formula]\label{thm:DXmain1}\index{theorem!asymptotic Riemann--Roch formula}
    Let $(L,h)$ be a Hermitian pseudo-effective line bundle on $X$ and $T$ be a holomorphic line bundle on $X$.
    We have
    \begin{equation}\label{eq:DXmain1}
        \lim_{k\to\infty}\frac{n!}{k^n}h^0\left(X,T\otimes L^k\otimes \mathcal{I}(h^k)\right)=\vol \left(\ddc h\right).
    \end{equation}
    In particular, the limit exists.
\end{theorem}

For the proof, let us fix a smooth Hermitian metric $h_0$ on $L$ with $\theta=c_1(L,h_0)$. We identify $h$ with $h_0\exp(-\varphi)$ for some $\varphi\in \PSH(X,\theta)$. See \cref{sec:pshmetriconlb} for the relevant notations.

Recall that when $X$ admits a big line bundle, it is necessarily projective. See \cite[Theorem~2.2.26]{MM07}.

We first handle the case where $\varphi$ has analytic singularities.
\begin{proposition}\label{prop:DXmainanalytic}
    Under the assumptions of \cref{thm:DXmain1}, assume furthermore that $\varphi$ has analytic singularities such that $\theta_{\varphi}$ is a Kähler current, then \eqref{eq:DXmain1} holds.
\end{proposition}
\begin{proof}
    \textbf{Step~1}. Reduce to the case of log singularities.

    Let $\pi\colon Y\rightarrow X$ be a log resolution of $\varphi$. In this case, for each $k\in \mathbb{Z}_{>0}$, we have
    \[
        h^0\Bigl(X,T\otimes L^k\otimes \mathcal{I}(k\varphi)\Bigr)=h^0\Bigl(Y,K_{Y/X}\otimes \pi^*T\otimes \pi^*L^k\otimes \mathcal{I}(k\pi^*\varphi)\Bigr).
    \]
    By \cref{prop:Ienvelopebimero}, we have
    \[
        \vol(\ddc h)=\vol (\ddc \pi^*h).
    \]
    Therefore, it suffices to argue \eqref{eq:DXmain1} with $K_{Y/X}\otimes \pi^*T$, $\pi^*L$ and $\pi^*h$ in place of $T$, $L$ and $h$.

    Note that $\theta_{\varphi}$ is no longer a Kähler current, but it still has positive mass.

    \textbf{Step~2}. Assume that $\varphi$ has log singularities along an effective $\mathbb{Q}$-divisor $D$. We may assume that $D$ is simple normal crossing by replacing $X$ by a modification and $D$ by the pull-back. We decompose $D$ into irreducible components, say
    \[
        D=\sum_{i=1}^N a_iD_i.
    \]
    In this case, we can easily compute
    \[
        \mathcal{I}(k\varphi)=\mathcal{O}_X\left(-\sum_{i=1}^N\lfloor k a_i \rfloor D_i\right)
    \]
    for each $k\in \mathbb{Z}_{>0}$.
    Observe that $L-D$ is nef (see \cref{lma:logsingrem}). Choose an integer $q>0$ so that $qD$ is integral. For each residue class $s=0,\ldots,q-1$, we have
    \[
        \left\lfloor (mq+s)D\right\rfloor=mqD+\lfloor sD\rfloor.
    \]
    Hence, on the residue class $k=mq+s$, the line bundle
    \[
        T\otimes L^k\otimes \mathcal{O}_X\left(-\sum_{i=1}^N\lfloor k a_i \rfloor D_i\right)
    \]
    is a fixed twist of $\left(q(L-D)\right)^m$. Applying the asymptotic Riemann--Roch theorem \cite[Corollary~1.4.41]{PAG1} to the nef line bundle $q(L-D)$ and then letting $m\to\infty$ in each residue class, we conclude that
    \[
        \lim_{k\to\infty}\frac{n!}{k^n}h^0\left(X,T\otimes L^k\otimes \mathcal{O}_X\left(-\sum_{i=1}^N\lfloor k a_i \rfloor D_i\right)\right)= (L-D)^n.
    \]
    Observe that by \cref{prop:LelongPoincare},
    \[
        \theta_{\varphi}=[D]+S,
    \]
    where $S$ is a closed positive $(1,1)$-current with bounded potential. Therefore,
    \[
        (L-D)^n=\int_X S^n=\int_X  \theta_{\varphi}^n.
    \]
    By \cref{ex:analyIgood}, we know that the right-hand side is exactly $\vol \theta_{\varphi}$.
\end{proof}

\begin{proof}[Proof of \cref{thm:DXmain1}]
    \textbf{Step~1}. We first handle the case where $\theta_{\varphi}$ is a Kähler current. Fix a Kähler form $\omega\geq \theta$ on $X$ such that $\theta_{\varphi}\geq 2\delta\omega$ for some $\delta\in (0,1)$.

    Let $(\varphi_j)_j$ be a quasi-equisingular approximation of $\varphi$ in $\PSH(X,\theta)$.
    We may assume that $\theta_{\varphi_j}\geq \delta\omega$ for all $j$.
    From \cref{prop:DXmainanalytic}, we know that for each $j\geq 1$,
    \[
        \varlimsup_{k\to\infty}\frac{n!}{k^n}h^0\left(X,T\otimes L^k\otimes  \mathcal{I}(k\varphi)\right)\leq   \lim_{k\to\infty}\frac{n!}{k^n}h^0\left(X,T\otimes L^k\otimes  \mathcal{I}(k\varphi_j)\right)=\vol \theta_{\varphi_j}.
    \]
    It follows from \cref{thm:charIgoodasclosure} and \cref{thm:contvolu} that the right-hand side converges to $\vol \theta_{\varphi}$ as $j\to\infty$. Therefore,
    \[
        \varlimsup_{k\to\infty}\frac{n!}{k^n}h^0\left(X,T\otimes L^k\otimes  \mathcal{I}(k\varphi)\right)\leq \vol \theta_{\varphi}.
    \]
    Conversely, fix an integer $N>\delta^{-1}$.
    From \cref{thm:charIgoodasclosure} and \cref{thm:convdS}, we know that
    \begin{equation}\label{eq:quasiequmassconvtemp1}
        \lim_{j\to\infty}\int_X \theta_{\varphi_j}^n=\int_X \theta_{P_{\theta}[\varphi]_{\mathcal{I}}}^n>0.
    \end{equation}
    Therefore, by \cref{lma:pathoenvelope}, we can find $j_0>0$ such that for $j\geq j_0$, there is $\psi\in \PSH(X,\theta)_{>0}$ (depending on $j$) with
    \begin{equation}\label{eq:linearlowerbdPItemp1}
        (1-N^{-1})\varphi_j+N^{-1}\psi \leq P_{\theta}[\varphi]_{\mathcal{I}}.
    \end{equation}
    For each $k>0$, we write $k=k'N-r$, where $k'\in \mathbb{N}$ and $r\in \{0,1,\ldots,N-1\}$.
    Then we compute for $j\geq j_0$ and large enough $k$ (to be specified shortly) that
    \[
        \begin{aligned}
            & h^0\left(X,T\otimes L^{k}\otimes \mathcal{I}(k\varphi)\right)\\
            \geq & h^0\left(X,T\otimes L^{-r}\otimes L^{k'N}\otimes \mathcal{I}(k'N\varphi)\right)\\
            \geq & h^0\left(X,T\otimes L^{-r}\otimes L^{k'N}\otimes \mathcal{I}\left(k'(\psi+(N-1)\varphi_j)\right)\right)\\
            \geq & h^0\left(X,T\otimes L^{-r}\otimes L^{k'(N-1)}\otimes \mathcal{I}\left(k'N\varphi_j\right)\right),
        \end{aligned}
    \]
    where the third line follows from \eqref{eq:linearlowerbdPItemp1}, the fourth line can be argued as follows: For large enough $k$, there is a non-zero section $s\in \mathrm{H}^0(X,L^{k'}\otimes \mathcal{I}(k'\psi))$ by \cref{lma:existsecposmass}. It follows from \cref{lma:IandIinf} that for large enough $k$,
    \[
        \mathcal{I}\left(k'N\varphi_j\right)\subseteq \mathcal{I}_{\infty}\left(k'(N-1)\varphi_j\right).
    \]
    It follows that multiplication by $s$ gives an injective map
    \[
        \begin{split}
            \mathrm{H}^0\left(X,T\otimes L^{-r}\otimes L^{k'(N-1)}\otimes \mathcal{I}\left(k'N\varphi_j\right)\right)\hookrightarrow \\
            \mathrm{H}^0\left(X,T\otimes L^{-r}\otimes L^{k'N}\otimes \mathcal{I}\left(k'\psi+k'(N-1)\varphi_j\right)\right).
        \end{split}
    \]
    Next observe that
\[
    (N-1)\theta+N\ddc\varphi_j= N\theta_{\varphi_j}-\theta \geq N\delta\omega-\theta \geq (N\delta-1)\omega.
\]
Since $N>\delta^{-1}$, this is a Kähler current. Hence \cref{prop:DXmainanalytic} is applicable. We let $k\to\infty$ to conclude that
    \[
        \begin{aligned}
        \varliminf_{k\to\infty}\frac{n!}{k^n}h^0\left(X,T\otimes L^{k}\otimes \mathcal{I}(k\varphi)\right)
        \geq & N^{-n}\int_{X} \left( (N-1)\theta+N\ddc \varphi_j\right)^n\\
        =& \int_X \left((1-N^{-1})\theta+\ddc\varphi_j \right)^n\\
        \geq & \int_X \left(\theta+\ddc\varphi_j \right)^n -C N^{-1},
        \end{aligned}
    \]
    where $C$ is a constant independent of $N$ and $j$.
    Letting $j\to\infty$ and then $N\to\infty$ and using \eqref{eq:quasiequmassconvtemp1}, we find that
    \[
        \varliminf_{k\to\infty}\frac{n!}{k^n}h^0\left(X,T\otimes L^{k}\otimes \mathcal{I}(k\varphi)\right)\geq \int_X \theta_{P_{\theta}[\varphi]_{\mathcal{I}}}^n.
    \]
    Therefore, \eqref{eq:DXmain1} follows.

    \textbf{Step~2}. We handle the case where $\vol \theta_{\varphi}>0$. We may assume that $\varphi$ is $\mathcal{I}$-model.

    First observe that $L$ is big by \cref{prop:nppposimpbig}. Hence $X$ is projective.
    Take a very ample line bundle $A$ on $X$ and a Kähler form $\omega$ in $c_1(A)$. Take a Hermitian metric $h_A$ on $A$ with $\ddc h_A=\omega$.

    Fix $N\in \mathbb{Z}_{>0}$. We decompose any $k>0$ as $k=k'N+r$ with $k'\in \mathbb{N}$ and $r\in \{0,1,\ldots,N-1\}$. Then
    \[
    h^0\left(X,T\otimes L^k\otimes  \mathcal{I}(k\varphi)\right)\leq h^0\left(X,T\otimes L^r \otimes L^{k'N}\otimes  \mathcal{I}(k'N\varphi)\right).
    \]
    Therefore,
    \[
    \begin{aligned}
     &\varlimsup_{k\to\infty}\frac{n!}{k^n}h^0\left(X,T\otimes L^k\otimes  \mathcal{I}(k\varphi)\right)\\
     \leq & \max_{r=0,\ldots,N-1}\varlimsup_{k'\to\infty}\frac{n!}{k'^nN^n} h^0\left(X,T\otimes L^r \otimes L^{k'N}\otimes  \mathcal{I}(k'N\varphi)\right)\\
     \leq & \max_{r=0,\ldots,N-1}\varlimsup_{k'\to\infty}\frac{n!}{k'^nN^n} h^0\left(X,T\otimes L^r \otimes L^{k'N}\otimes A^{k'}\otimes  \mathcal{I}(k'N\varphi)\right)\\
     =& \int_X \left(N^{-1}\omega+\theta+\ddc P_{\theta+N^{-1}\omega}[\varphi]_{\mathcal{I}}\right)^{n},
    \end{aligned}
    \]
    where the third line uses multiplication by a non-zero section of $A^{k'}$, and we have applied Step~1 to the Hermitian pseudo-effective line bundle $(L^N\otimes A,h^{\otimes N}\otimes h_A)$ on the fourth line.
    On the other hand, since $\varphi$ is $\mathcal{I}$-good by \cref{ex:ImodelIgood}, \cref{lma:Igoodinsenspert} gives
    \[
        P_{\theta+N^{-1}\omega}[\varphi]_{\mathcal{I}}=P_{\theta+N^{-1}\omega}[\varphi].
    \]
    It follows from \cref{prop:Ppresmass} that
    \[
        \varlimsup_{k\to\infty}\frac{n!}{k^n}h^0\left(X,T\otimes L^k\otimes  \mathcal{I}(k\varphi)\right)\leq \int_X   \left(N^{-1}\omega+\theta+\ddc \varphi\right)^n.
    \]
    Letting $N\to\infty$, we conclude
    \[
        \varlimsup_{k\to\infty}\frac{n!}{k^n}h^0\left(X,T\otimes L^k\otimes  \mathcal{I}(k\varphi)\right)\leq \int_X \theta_{\varphi}^n.
    \]
    It remains to argue the reverse inequality.

    Choose $\psi\in \PSH(X,\theta)$ such that $\theta_{\psi}$ is a Kähler current and $\psi\leq \varphi$. The existence of $\psi$ is guaranteed by \cref{lma:kahcurrentposmass}. Then for any $t\in (0,1)$, we set
    \[
        \varphi_t=(1-t)\varphi+t\psi.
    \]
    It follows again from Step~1 that
    \[
        \varliminf_{k\to\infty}\frac{n!}{k^n}h^0\left(X,T\otimes L^k\otimes  \mathcal{I}(k\varphi)\right)\geq \lim_{k\to\infty}\frac{n!}{k^n}h^0\left(X,T\otimes L^k\otimes  \mathcal{I}(k\varphi_t)\right)=\vol \theta_{\varphi_t}.
    \]
    On the other hand, by \cref{thm:mvolcont} and \cref{cor:lineardsconv},
    \[
        \lim_{t\to 0+}\vol \theta_{\varphi_t}=\vol \theta_{\varphi}.
    \]
    So we find
    \[
        \varliminf_{k\to\infty}\frac{n!}{k^n}h^0\left(X,T\otimes L^k\otimes  \mathcal{I}(k\varphi)\right)\geq \vol \theta_{\varphi}.
    \]
    We conclude \eqref{eq:DXmain1} in this case.

    \textbf{Step~3}. We finally handle the case where $\vol \theta_{\varphi}=0$. Replacing $\varphi$ by $P_{\theta}[\varphi]_{\mathcal{I}}$, we may assume that $\varphi$ is $\mathcal{I}$-model.

    Assume that \eqref{eq:DXmain1} fails. That is,
    \[
    \varlimsup_{k\to\infty}\frac{n!}{k^n}h^0\left(X,T\otimes L^k\right)\geq \varlimsup_{k\to\infty} \frac{n!}{k^n}h^0\left(X,T\otimes L^k\otimes \mathcal{I}(k\varphi)\right)>0,
    \]
    then $L$ is a big line bundle and hence $X$ is projective.

    Fix a very ample line bundle $A$ on $X$ and a K\"ahler form $\omega\in c_1(A)$. Take a decreasing sequence $(\epsilon_j)_j$ of rational numbers with limit $0$ and a quasi-equisingular approximation $(\varphi_j)_j$ of $\varphi$ with $\varphi_j\in \PSH(X,\theta+\epsilon_j\omega)_{>0}$.

    We claim that as $j\to\infty$, the sequence $P_{\theta+\epsilon_j\omega}[\varphi_j]$ is decreasing with limit $\varphi$.

    It is clear that this sequence is decreasing. Let $\psi$ denote its limit for the moment. It is also clear that $\psi\geq\varphi$. Since $\varphi$ is $\mathcal{I}$-model, it remains to show that $\psi\preceq_{\mathcal{I}}\varphi$. Fix $\lambda>0$. By strong openness, choose $\lambda'>\lambda$ such that $\mathcal{I}(\lambda\psi)=\mathcal{I}(\lambda'\psi)$. For $j$ large enough, the quasi-equisingularity of $(\varphi_j)_j$ gives
    \[
        \mathcal{I}(\lambda'\psi)\subseteq \mathcal{I}(\lambda'\varphi_j)\subseteq \mathcal{I}(\lambda\varphi),
    \]
    where the first inclusion uses $\psi\leq P_{\theta+\epsilon_j\omega}[\varphi_j]$ and $P_{\theta+\epsilon_j\omega}[\varphi_j]\sim_{\mathcal{I}}\varphi_j$. Thus $\psi\preceq_{\mathcal{I}}\varphi$, and the claim follows.

    By our claim, \cref{prop:Ppresmass} and \cref{prop:vol_limit_model}, we find that
    \begin{equation}\label{eq:volumeconv0temp1}
    \lim_{j\to\infty} \int_X \left(\theta+\epsilon_j\omega+\ddc \varphi_j\right)^n=\int_X \theta_{\varphi}^n=0.
    \end{equation}

    Fix $j>0$, take an integer $N>0$ so that $N\epsilon_j$ is an integer. Then we compute
    \[
    \begin{aligned}
    &\varlimsup_{k\to\infty} \frac{n!}{k^n}h^0\left(X,T\otimes L^k\otimes \mathcal{I}(k\varphi)\right)\\
    \leq & \varlimsup_{k\to\infty} \frac{n!}{k^n}h^0\left(X,T\otimes L^k\otimes \mathcal{I}(k\varphi_j)\right)\\
    \leq & \max_{a=0,\ldots,N-1}\varlimsup_{k'\to\infty} \frac{n!}{(k'N)^n}h^0\left(X,T\otimes L^a\otimes L^{Nk'}\otimes \mathcal{I}(Nk'\varphi_j)\right)\\
    \leq & \max_{a=0,\ldots,N-1}\varlimsup_{k'\to\infty} \frac{n!}{(k'N)^n}h^0\left(X,T\otimes L^a\otimes L^{Nk'}\otimes A^{k'N\epsilon_j}\otimes \mathcal{I}(Nk'\varphi_j)\right)\\
    =& \frac{1}{N^n}\int_X \left( N\theta+\epsilon_jN\omega+N\ddc \varphi_j \right)^n,
    \end{aligned}
    \]
    where the third line follows by writing $k=Nk'+a$ as before,
    we applied Step~2 on the last line. Dividing by $N^n$, we find that
    \[
    \varlimsup_{k\to\infty} \frac{n!}{k^n}h^0\left(X,T\otimes L^k\otimes \mathcal{I}(k\varphi)\right)\leq \int_X \left(\theta+\epsilon_j\omega+\ddc \varphi_j\right)^n.
    \]
    Since we know \eqref{eq:volumeconv0temp1}, letting $j\to\infty$, we conclude that
    \[
    \varlimsup_{k\to\infty} \frac{n!}{k^n}h^0\left(X,T\otimes L^k\otimes \mathcal{I}(k\varphi)\right)=0,
    \]
    which is a contradiction. Hence \eqref{eq:DXmain1} is established in full generality.
 \end{proof}

\begin{corollary}\label{cor:volbigL}
Let $L$ be a pseudo-effective line bundle on $X$, and let $h$ be a smooth Hermitian metric on $L$ with $\theta=c_1(L,h)$. Then we have
    \begin{equation}\label{eq:volbig}
        \lim_{k\to\infty}\frac{n!}{k^n}h^0\left(X,L^k\right)=\int_X \theta_{V_{\theta}}^n.
    \end{equation}
\end{corollary}
\begin{proof}
    It suffices to observe that
    \begin{equation}\label{eq:H0eqH0IVtheta}
        \mathrm{H}^0\left(X,L^k \otimes \mathcal{I}\left(kV_{\theta}\right) \right)=\mathrm{H}^0\left(X,L^k\right)
    \end{equation}
    for all $k>0$, since for any holomorphic section $s$ of $L^k$, the function $\log |s|_{h^k}^2$ is $\theta$-psh thanks to \cref{prop:LelongPoincare} and therefore, we have
    \[
        \frac{1}{k}\log |s|_{h^k}^2\preceq V_{\theta}.
    \]
    In particular, $s\in \mathrm{H}^0\left(X,L^k \otimes \mathcal{I}\left(kV_{\theta}\right) \right)$.
\end{proof}
This common quantity is the \emph{volume}\index{volume!volume of a pseudo-effective line bundle} of $L$, usually denoted by $\vol L$. In view of \cref{def:volcohclass}, we have
\begin{equation}\label{eq:volL_volc1L}
\vol L=\vol c_1(L).
\end{equation}

\begin{example}\label{ex:toricIgood}
    If $X$ is a smooth projective toric variety and $\theta$ is invariant under the action of the compact torus, then any $\varphi\in \PSH_{\tor}(X,\theta)$ is $\mathcal{I}$-good.
\end{example}
\begin{proof}
    Thanks to \cref{lma:Igoodinsenspert}, we may assume that $\theta\in c_1(L)$ for some toric invariant ample line bundle $L$. In this case, the result follows from \cref{thm:charIgoodasclosure}, \cref{thm:DXmain1} and \cref{thm:Yao}.
\end{proof}

\makeatletter
\providecommand{\russianmotto}[1]{{%
        \font\russianmottofont="[Tempora-Italic.otf]" at \f@size pt\relax
        \russianmottofont #1%
}}
\makeatother

\motto{Some of our philosophers made a serious mistake: without having grasped the substance
of the matter, they began to deny the significance of a new direction in science.

--- S.~L.~Sobolev\footnote{Sergei Lvovich Sobolev (1908--1989) was a Soviet mathematician whose work on weak
derivatives, generalized functions, and embedding theorems became foundational for modern
analysis and partial differential equations. He also contributed to mathematical physics,
numerical analysis, cubature formulas, and cybernetics, and was later central to the
mathematical life of Akademgorodok in Novosibirsk.}, A.~I.~Kitov and A.~A.~Lyapunov in \russianmotto{Вопросы философии}
}

\chapter{The trace operator}
\label{chap:trace}

The goal of this chapter is to give a meaning to the restriction of a quasi-psh singularity to a subvariety. If $\varphi\in\QPSH(X)$ and $Y\subseteq X$ is a subvariety, the naive restriction $\varphi|_Y$ is often unsatisfactory: it may be identically $-\infty$, or it may fail to preserve quasi-equisingular approximations.

There is a useful analogy with the trace operator in real analysis. A Sobolev function on a domain cannot, in general, be restricted pointwise to the boundary; instead, it admits a trace obtained by approximation using smooth functions, and the trace is well-defined modulo almost everywhere equality.
\cref{tab:trace_analogy} summarizes the corresponding dictionary.
\begin{center}
\refstepcounter{table}\label{tab:trace_analogy}
\begin{tabular}{|c|c|}
    \hline
    Real analysis & Pluripotential theory \\
    \hline
    smooth functions & potentials with analytic singularities \\
    \hline
    Sobolev functions & potentials with $\mathcal I$-good singularities \\
    \hline
    Sobolev convergence & $d_S$-convergence \\
    \hline
    equality almost everywhere & $P$-equivalence \\
    \hline
\end{tabular}

\smallskip
\textbf{Table~\thetable} A dictionary for the trace operator.
\end{center}
The situation here is similar. A singular quasi-psh function need not have a meaningful pointwise restriction to $Y$, but after replacing it by quasi-equisingular approximants and passing to the appropriate $d_S$-limit of the pointwise restriction of these approximants to $Y$, one obtains a well-defined $P$-singularity type. This is the trace $\Tr_Y(\varphi)$.

The trace operator is designed to be compatible with the two comparison theories developed earlier. It depends only on the $\mathcal{I}$-singularity type of $\varphi$, and itself is well-defined modulo $P$-equivalence.

The first part of the chapter defines the trace and proves its basic properties: independence of the quasi-equisingular approximation, monotonicity, linearity, compatibility with bimeromorphic pullback, and continuity under decreasing $d_S$-limits. These results provide the analytic foundation for later applications.

The second part relates the trace operator to restricted volumes. For a big class $\alpha$ and a subvariety $Y$ not contained in the non-K\"ahler locus, the classical restricted volume $\vol_{X|Y}(\alpha)$ can be expressed as the non-pluripolar volume of the trace of a potential with minimal singularities. We also prove an asymptotic Riemann--Roch formula for the restricted multiplier ideal sheaves.

The final part develops analytic Bertini-type statements. These results describe when multiplier ideals restrict correctly to a general hypersurface and when $\mathcal I$-goodness is preserved under restriction.

\section{The definition of the trace operator}\label{sec:trace_def}

Let $X$ be a connected compact Kähler manifold and $Y\subseteq X$ be an irreducible analytic subset.

The trace operator gives a way to restrict a quasi-plurisubharmonic function on $X$ to $\widetilde{Y}$, the normalization of $Y$.
It follows from \cite[Proposition~3.6]{GK20} that $\widetilde{Y}$ is a normal Kähler space.
We refer to \cref{chap:unibranch} for the pluripotential theory on unibranch Kähler spaces.

For later applications, we need this generality even if initially we are only interested in the smooth case.

Before diving into the trace operator, let us try to understand what goes wrong with the naive restriction procedure: simply take $\varphi|_{\widetilde{Y}}$. Both examples below are local, but can easily be globalized since the singularities are isolated.
\begin{example}\label{ex:naive_res1}
    Consider the case where $\varphi$ has log-log singularities as in \cref{ex:loglogsing}. Locally we can take $\varphi(z)=-\log (-\log |z|^2)$. We have $\nu(\varphi,0)=0$ but $\varphi|_0=-\infty$. So the naive restriction is not defined in this case.
\end{example}

Even if the naive restriction is defined, it does not behave well.
\begin{example}\label{ex:naive_res2}
    Consider a psh function $\varphi$ in two variables, say
    \[
    \varphi(z,w)=\left( -\log (-\log |z|^2)\right) \lor \left(\log |w|^2\right).
    \]
    Take an arbitrary quasi-equisingular approximation $(\varphi_j)_j$. Then each $\varphi_j$ is locally bounded since $\nu(\varphi,(0,0))=0$. Let us consider the restrictions to $H=\{z=0\}$. Then $\varphi_j|_H$ is still bounded, while
    \[
    \varphi|_H(w)=\log |w|^2.
    \]
    In other words, the restrictions $\varphi_j|_H$ fail to be a quasi-equisingular approximation of $\varphi|_H$.
\end{example}
Our trace operator gives an elegant way to solve both problems.

We first observe that if $\varphi\in \QPSH(X)$ has analytic singularities and $\nu(\varphi,Y)=0$, then $\varphi|_Y\not\equiv -\infty$.
Indeed, if $\varphi|_Y\equiv-\infty$, then in a local presentation
\[
    \varphi=c\log(|f_1|^2+\cdots+|f_N|^2)+\mathcal{O}(1)
\]
all generators $f_1,\ldots,f_N$ vanish along $Y$, which gives $\nu(\varphi,Y)>0$, a contradiction.

This observation will be crucial in the sequel.

\begin{proposition}\label{prop:traceindquasiequisingapp}
    Let $\varphi\in \QPSH(X)$ be a function such that $\nu(\varphi,Y)=0$. Let $(\varphi_i)_i$, $(\psi_i)_i$ be quasi-equisingular approximations of $\varphi$. Then
    \begin{equation}\label{eq:dsequivtemp1}
        \lim_{i\to\infty} d_S\left(\varphi_i|_{\widetilde{Y}},\psi_i|_{\widetilde{Y}}\right)=0.
    \end{equation}
\end{proposition}
Here in \eqref{eq:dsequivtemp1}, $d_S$ refers to any $d_{S,\theta}$, where $\theta$ is any closed smooth real $(1,1)$-form on $\widetilde{Y}$ representing a big cohomology class and satisfying $\varphi_i|_{\widetilde{Y}},\psi_i|_{\widetilde{Y}}\in \PSH(\widetilde{Y},\theta)$. For different choices of $\theta$, \eqref{eq:dsequivtemp1} remains equivalent, as explained in \cref{cor:dSconv_changetheta}.
\begin{proof}
After subtracting a common constant from $\varphi$, $\varphi_i$ and $\psi_i$, we may assume that $\varphi_i\leq 0$ for all $i\geq 1$.

Take a K\"ahler form $\omega$ on $X$ such that $\varphi_i,\psi_i\in \PSH(X,\omega/2)$ for all $i\geq 1$. By \cref{cor:dSconv_changetheta}, we need to show that
\[
\lim_{i\to\infty} d_{S,\omega|_{\widetilde{Y}}}\left(\varphi_i|_{\widetilde{Y}},\psi_i|_{\widetilde{Y}}\right)=0.
\]
Assume that this fails. Then, up to replacing the sequences by subsequences, we may assume that the following limit exists and
\[
\lim_{i\to\infty} d_{S,\omega|_{\widetilde{Y}}}\left(\varphi_i|_{\widetilde{Y}},\psi_i|_{\widetilde{Y}}\right)>0.
\]
Take a K\"ahler form $\tilde{\omega}$ on $\widetilde{Y}$. Then
\[
\lim_{i\to\infty} d_{S,\omega|_{\widetilde{Y}}+\tilde{\omega}}\left(\varphi_i|_{\widetilde{Y}},\psi_i|_{\widetilde{Y}}\right)>0
\]
by \cref{cor:dSconv_changetheta}.

By \cref{cor:quasi-equichar}, after changing the background class from $\omega/2$ to $\omega$ if necessary, we have
\[
    \varphi_i\xrightarrow{d_S}P_{\omega}[\varphi]_{\mathcal I},\quad \psi_i\xrightarrow{d_S}P_{\omega}[\varphi]_{\mathcal I}.
\]
Let
\[
    \chi_i\coloneqq \varphi_i\lor\psi_i.
\]
By \cref{prop:lor_dS_conv}, we have $\chi_i\xrightarrow{d_S} P_{\omega}[\varphi]_{\mathcal{I}}$. Moreover, $(\chi_i)_i$ is again a quasi-equisingular approximation of the original potential $\varphi$ as follows from \cref{cor:quasi-equichar}.

The triangle inequality gives
\[
d_S(\varphi_i|_{\widetilde{Y}},\psi_i|_{\widetilde{Y}}) \leq  d_S(\varphi_i|_{\widetilde{Y}},\chi_i|_{\widetilde{Y}})+ d_S(\psi_i|_{\widetilde{Y}},\chi_i|_{\widetilde{Y}}).
\]
Thus, after passing to a subsequence and interchanging the two quasi-equisingular approximations if necessary, we may assume that
\[
    \varliminf_{i\to\infty}d_S(\varphi_i|_{\widetilde{Y}},\chi_i|_{\widetilde{Y}})>0.
\]
Replacing $(\psi_i)_i$ by $(\chi_i)_i$, we may therefore assume that
\[
    \varphi_i\leq \psi_i
\]
for all $i$, while preserving the contradiction assumption.

Take a decreasing sequence $(\epsilon_j)_j$ in $\mathbb{R}_{>0}$ with limit $0$ such that $(1-\epsilon_j)\varphi_j\in \PSH(X,\omega)$. We first observe that
    \[
        \lim_{i\to\infty} d_{S,\omega|_{\widetilde{Y}}}\left(\varphi_i|_{\widetilde{Y}},(1-\epsilon_i)\varphi_i|_{\widetilde{Y}}\right)=0.
    \]
This is a consequence of \cref{lma:dSsmallmult}. Hence, by \cref{cor:dSconv_changetheta}, we find
\[
    \lim_{i\to\infty} d_{S,\omega|_{\widetilde{Y}}+\tilde{\omega}}\left(\varphi_i|_{\widetilde{Y}},(1-\epsilon_i)\varphi_i|_{\widetilde{Y}}\right)=0.
\]
Since $(\varphi_i)_i$ and $(\epsilon_i)_i$ are decreasing and $\varphi_i\leq 0$, the sequence $\bigl((1-\epsilon_i)\varphi_i|_{\widetilde{Y}}\bigr)_i$ is decreasing. It has uniformly positive mass in the class $\omega|_{\widetilde{Y}}+\tilde{\omega}$. Hence \cref{cor:completenessdS}, together with the preceding estimate, gives a potential $\psi\in \PSH(\widetilde{Y},\omega|_{\widetilde{Y}}+\tilde{\omega})$ such that
\[
\varphi_i|_{\widetilde{Y}}\xrightarrow{d_S}\psi.
\]
Hence,
\[
\lim_{i\to\infty} d_{S,\omega|_{\widetilde{Y}}+\tilde{\omega}}\left(\psi,(1-\epsilon_i)\varphi_i|_{\widetilde{Y}}\right)=0.
\]

By \cref{prop:compqequi}, after choosing $j_i$ inductively with $j_i\geq i$ and $j_i>j_{i-1}$, we may assume that
\[
    \psi_{j_i}\preceq (1-\epsilon_i)\varphi_i.
\]
Since we have arranged $\varphi_j\leq\psi_j$ for all $j$, it follows that
\[
    \varphi_{j_i}\leq\psi_{j_i}\preceq (1-\epsilon_i)\varphi_i.
\]
Hence,
\[
    \varphi_{j_i}|_{\widetilde{Y}}\leq \psi_{j_i}|_{\widetilde{Y}}\preceq (1-\epsilon_i)\varphi_i|_{\widetilde{Y}}.
\]
Therefore, by \cref{cor:dsthreeterm},
\[
    \begin{aligned}
       \varlimsup_{i\to\infty} d_{S,\omega|_{\widetilde{Y}}+\tilde{\omega}}\left(\varphi_{j_i}|_{\widetilde{Y}},\psi_{j_i}|_{\widetilde{Y}}\right)
     \leq & \varlimsup_{i\to\infty} d_{S,\omega|_{\widetilde{Y}}+\tilde{\omega}}\left(\varphi_{j_i}|_{\widetilde{Y}},(1-\epsilon_i)\varphi_i|_{\widetilde{Y}}\right)\\
     =& \varlimsup_{i\to\infty} d_{S,\omega|_{\widetilde{Y}}+\tilde{\omega}}\left(\psi,(1-\epsilon_i)\varphi_i|_{\widetilde{Y}}\right)\\
     = & 0,
     \end{aligned}
\]
which is a contradiction.
\end{proof}

\begin{definition}\label{def:traceop}
    Let $\varphi\in \QPSH(X)$ be a function such that $\nu(\varphi,Y)=0$. We say a potential $\psi\in \QPSH(\widetilde{Y})$ is a\footnote{Let us resume our analogy in the introduction of this chapter.
In real analysis, the trace of a Sobolev function $f$ is only defined up to almost-everywhere equality; we say a particular function on $\partial\Omega$ is \emph{a} trace of $f$. Likewise here, we say $\psi$ is \emph{a} trace operator of $\varphi$ rather than that the $P$-equivalence class of $\psi$ is \emph{the} trace operator of $\varphi$.} \emph{trace operator}\index{trace operator} of $\varphi$ along $Y$ if there is a quasi-equisingular approximation $(\varphi_j)_{j>0}$ of $\varphi$ such that
    \begin{equation}\label{eq:deftrace}
        \varphi_j|_{\widetilde{Y}}\xrightarrow{d_S}\psi\footnotemark.
    \end{equation}
    \footnotetext{To be more precise, what we mean is the following: We can find a closed smooth real $(1,1)$-form on $X$ such that $\varphi\in \PSH(X,\theta)$. Then there is a K\"ahler form such that $\omega+\theta+\ddc\varphi_j\geq 0$ for all $j\geq 1$. Take a K\"ahler form $\tilde{\omega}$ on $\widetilde{Y}$ so that $\tilde{\omega}\geq (\theta+\omega)|_{\widetilde{Y}}$ and $\psi\in \PSH(\widetilde{Y},\tilde{\omega})$. Then our condition means that $\varphi_j|_{\widetilde{Y}}\xrightarrow{d_{S,\tilde{\omega}}}\psi$. This condition is independent of the choices of $\theta$, $\omega$ and $\tilde{\omega}$ by \cref{cor:dSconv_changetheta}.}
\end{definition}
Indeed, choosing the background Kähler form on $\widetilde{Y}$ sufficiently positive, the restricted quasi-equisingular approximation is a decreasing sequence with masses bounded below by a positive constant; hence it converges by \cref{cor:completenessdS}.
Observe that by \cref{prop:traceindquasiequisingapp}, the condition \eqref{eq:deftrace} is independent of the choice of $(\varphi_j)_j$.

Later on in \cref{thm:tracetor}, we shall prove that the trace operator corresponds to the natural way of restricting convex bodies in the toric setting.

\begin{proposition}\label{prop:traceunique}
    Let $\varphi\in \QPSH(X)$ such that $\nu(\varphi,Y)=0$. Suppose that $\psi$ and $\psi'$ are trace operators of $\varphi$ along $Y$. Then $\psi$ and $\psi'$ are $\mathcal{I}$-good and $\psi\sim_P\psi'$.
\end{proposition}
\begin{proof}
    The restrictions of the analytic approximants have analytic singularities on $\widetilde{Y}$, so the claim follows from \cref{thm:charIgoodasclosure}. The fact that $\psi\sim_P\psi'$ follows from \cref{prop:traceindquasiequisingapp} and \cref{prop:ds0char}.
\end{proof}

\begin{example}
    As a trivial example, when $Y$ is just a single point, then $\QPSH(Y)$ is canonically identified with $\mathbb{R}$. Any constant $c\in \mathbb{R}$ is a trace operator of a function $\varphi\in \QPSH(X)$ satisfying $\nu(\varphi,Y)=0$.
\end{example}

\begin{definition}\label{def:traceop2}
    Let $\varphi\in \QPSH(X)$ such that $\nu(\varphi,Y)=0$. We write $\Tr_{Y}(\varphi)$\index{$\Tr_{Y}(\varphi)$} for any trace operator of $\varphi$ along $Y$.

    Given a closed smooth real $(1,1)$-form $\theta$ on $X$. When $\Tr_Y(\varphi)$ can be chosen to lie in $\PSH(\widetilde{Y},\theta|_{\widetilde{Y}})_{>0}$, we write
    \[
        \Tr_Y^{\theta}(\varphi)\coloneqq P_{\theta|_{\widetilde{Y}}}\bigl[ \Tr_Y(\varphi) \bigr]=P_{\theta|_{\widetilde{Y}}}\bigl[ \Tr_Y(\varphi) \bigr]_{\mathcal{I}}.
    \]
\end{definition}
The trace operator $\Tr_Y(\varphi)$ is therefore well-defined only up to $P$-equivalence by \cref{prop:traceunique}. Also observe that if $\varphi\in \PSH(X,\theta)$ for some smooth closed real $(1,1)$-form $\theta$ on $X$, then for any K\"ahler form $\omega$ on $X$, the trace operator $\Tr_Y^{\theta+\omega}(\varphi)$ is always defined. In particular, if $\theta_{\varphi}$ is a Kähler current, $\Tr^{\theta}_Y(\varphi)$ is always defined.

\begin{remark}\label{rmk:tracecurrent}
    As in \cref{rmk:qpshtocurrents}, the trace operator could also be applied to closed positive $(1,1)$-currents on $X$. If $T\in \mathcal{Z}_+(X,\alpha)$ for some pseudo-effective class $\alpha$ on $X$ (see \cref{def:spaceofcurrents}) and $\beta\in \mathrm{H}^{1,1}(\widetilde{Y},\mathbb{R})$, then we write
    \[
    \Tr^{\beta}_Y(T)
    \]
    for any (if exists) closed positive $(1,1)$-current in $\beta$ representing $\Tr_Y(T)$\index{$\Tr_{Y}(T)$} when $\nu(T,Y)=0$.
\end{remark}

\begin{proposition}\label{prop:Trdominarest}
    Let $\varphi\in \QPSH(X)$ such that $\nu(\varphi,Y)=0$.
    Assume that $\varphi|_Y\not\equiv -\infty$.
    Then
    \[
    \varphi|_{\widetilde{Y}}\preceq_P \Tr_Y(\varphi).
    \]
\end{proposition}
\begin{proof}
    Take a Kähler form $\omega$ such that $\omega_{\varphi}$ is a Kähler current.
    Let $(\varphi_j)_{j>0}$ be a quasi-equisingular approximation of $\varphi$ in $\PSH(X,\omega)_{>0}$.
    We may assume that $\varphi_j\leq 0$ for all $j\geq 1$.

    Then
    \begin{equation}\label{eq:varphijrestrleqPtemp}
    \varphi_j|_{\widetilde{Y}}\leq P_{\omega|_{\widetilde{Y}}}\left[\varphi_j|_{\widetilde{Y}}\right]
    \end{equation}
    for all $j\geq 1$. In particular,
    \[
    \varphi|_{\widetilde{Y}}\leq \inf_{j\geq 1}P_{\omega|_{\widetilde{Y}}}\left[\varphi_j|_{\widetilde{Y}}\right].
    \]
    Thanks to \cref{cor:decnetdS},
\begin{equation}\label{eq:TrYnewexpression}
\Tr_Y(\varphi)\sim_P \inf_{j\geq 1} P_{\omega|_{\widetilde{Y}}}[\varphi_j|_{\widetilde{Y}}].
\end{equation}
  We conclude our assertion.
\end{proof}

\begin{example}\label{ex:resanalyt}
    Let $\varphi\in \QPSH(X)$ such that $\nu(\varphi,Y)=0$. Assume that $\varphi$ has analytic singularities. Then
    \[
        \Tr_Y(\varphi)\sim_P \varphi|_{\widetilde{Y}}.
    \]
\end{example}

\begin{example}\label{ex:trace:L236}
    Let $\varphi\in \QPSH(X)$. Take a closed real smooth $(1,1)$-form $\theta$ on $X$ such that $\varphi\in \PSH(X,\theta)_{>0}$, then
    \[
        \Tr_X(\varphi)\sim_P P_{\theta}[\varphi]_{\mathcal{I}},\quad   \Tr^{\theta}_X(\varphi) = P_{\theta}[\varphi]_{\mathcal{I}}.
    \]
    In particular, the trace operator can be regarded as a generalization of the $\mathcal{I}$-envelope.
\end{example}

\begin{example}\label{ex:tracedefinedposmass}
    Assume that $\varphi\in \PSH(X,\theta)$ for some closed smooth real $(1,1)$-form $\theta$ on $X$, $\nu(\varphi,Y)=0$ and
    \begin{equation}\label{eq:traceposmasscond}
    \lim_{\epsilon\searrow 0} \int_{\widetilde{Y}} \left( \theta|_{\widetilde{Y}} +\epsilon\omega|_{\widetilde{Y}}+\ddc\Tr_Y^{\theta+\epsilon\omega}(\varphi)\right)^{\dim Y} >0
    \end{equation}
    for any arbitrary choice of a Kähler form $\omega$ on $X$.
    Then it follows from \cref{prop:vol_limit_model} that $\Tr_Y^{\theta}(\varphi)$ is defined, and its mass is exactly the above limit.

    As a consequence, we have the following formula:
    \begin{equation}
        \int_{\widetilde{Y}} \left(\theta|_{\widetilde{Y}}+\ddc\Tr_Y^{\theta}(\varphi)\right)^{\dim Y}=\lim_{\epsilon\searrow 0} \int_{\widetilde{Y}} \left( \theta|_{\widetilde{Y}} +\epsilon\omega|_{\widetilde{Y}}+\ddc\Tr_Y^{\theta+\epsilon\omega}(\varphi)\right)^{\dim Y},
    \end{equation}
    where the left-hand side is understood as $0$ if $\Tr_Y^{\theta}(\varphi)$ is not defined.
\end{example}

\begin{remark}\label{rmk:movingSes}
    The trace operator allows us to introduce the following extension of the moving Seshadri constant. For $T\in \mathcal{Z}_+(X,\alpha)$ and $x\in X$, define
    \[
    \epsilon(T,x)\coloneqq \inf_{V\ni x}\left(\frac{\vol \Tr_V^{\alpha|_{\widetilde{V}}} T}{\mult_x V} \right)^{\frac{1}{\dim V}},
    \]
    where $\vol \Tr_V^{\alpha|_{\widetilde{V}}} T=0$ if $\Tr_V^{\alpha|_{\widetilde{V}}} T$ is not defined. Here $V$ runs over all positive-dimensional closed irreducible analytic subsets of $X$ containing $x$.

    These moving Seshadri constants seem to be new. But since I do not have particularly good applications in mind, I will not study these objects in this book.
\end{remark}

\section{Properties of the trace operator}\label{sec:proptrace}
Let $X$ be a connected compact Kähler manifold and $Y\subseteq X$ be an irreducible analytic subset.

We prove a few elementary properties of the trace operator.
\begin{proposition}\label{prop:tracelinear}
    Let $\varphi,\psi\in \QPSH(X)$, $\lambda>0$. Assume that $\nu(\varphi,Y)=\nu(\psi,Y)=0$.
    Then we have the following:
    \begin{enumerate}
        \item\label{itm:prop:tracelinear:1} Suppose that $\varphi\preceq_{\mathcal{I}}\psi$. Then $\Tr_Y(\varphi)\preceq_P \Tr_Y(\psi)$.
        \item\label{itm:prop:tracelinear:2} We have
            \[
                \Tr_Y(\varphi+\psi)\sim_P \Tr_Y(\varphi)+\Tr_Y(\psi).
            \]
        \item\label{itm:prop:tracelinear:3} We have
            \[
                \Tr_Y(\lambda\varphi)\sim_P \lambda\Tr_Y(\varphi).
            \]
        \item\label{itm:prop:tracelinear:4} We have
            \[
                \Tr_Y(\varphi\lor \psi)\sim_P \Tr_Y(\varphi)\lor \Tr_Y(\psi).
            \]
    \end{enumerate}
\end{proposition}
As a consequence of \ref{itm:prop:tracelinear:1}, the $P$-equivalence class of $\Tr_Y(\varphi)$ depends only on the $\mathcal{I}$-equivalence class of $\varphi$.
\begin{proof}
    Take a closed smooth real $(1,1)$-form $\theta$ on $X$ such that $\theta_{\varphi}$, $\theta_{\psi}$ are both Kähler currents.
    Let $(\varphi_j)_j$ and $(\psi_j)_j$ be quasi-equisingular approximations of $\varphi$ and $\psi$ in $\PSH(X,\theta)$ respectively. We may assume that $\varphi_j\leq 0$ and $\psi_j\leq 0$ for all $j\geq 1$.

    \ref{itm:prop:tracelinear:1} We choose the quasi-equisingular approximations compatibly. By \cref{thm:qequi}, we may take $\varphi_j$ and $\psi_j$ with analytic singularities of types
    \[
    \left(2^{-j},\mathcal I(2^j\varphi)\right), \quad \left(2^{-j},\mathcal I(2^j\psi)\right),
    \]
    respectively. Since $\varphi\preceq_{\mathcal I}\psi$, we have
    \[
        \mathcal I(2^j\varphi)\subseteq \mathcal I(2^j\psi)
    \]
    for every $j$. Hence $\varphi_j\preceq\psi_j$ for every $j$. Restricting to $\widetilde{Y}$, we get $\varphi_j|_{\widetilde{Y}}\preceq\psi_j|_{\widetilde{Y}}$. Passing to the $d_S$-limits by \cref{prop:dsconvpresorder}, we obtain
    \[
        \Tr_Y(\varphi)\preceq_P\Tr_Y(\psi).
    \]
    Restricting to $\widetilde{Y}$ and passing to the $d_S$-limits by \cref{prop:dsconvpresorder}, we obtain
    \[
        \Tr_Y(\varphi)\preceq_P \Tr_Y(\psi).
    \]

    \ref{itm:prop:tracelinear:2} It follows from \cref{thm:dSadditivity} that $\varphi_j+\psi_j\xrightarrow{d_S}P_{\theta}[\varphi]_{\mathcal{I}}+P_{\theta}[\psi]_{\mathcal{I}}$. However, by \cref{prop:PIconc} and \cref{prop:Igoodlinear}, we have
    \[
        P_{\theta}[\varphi]_{\mathcal{I}}+P_{\theta}[\psi]_{\mathcal{I}}\sim_P P_{2\theta}[\varphi+\psi]_{\mathcal{I}}.
    \]
    Therefore, by \cref{prop:ds0char}, \cref{cor:quasi-equichar} and \cref{prop:analysingclosed}, $(\varphi_j+\psi_j)_j$ is a quasi-equisingular approximation of $\varphi+\psi$. We conclude using \cref{thm:dSadditivity}.

    \ref{itm:prop:tracelinear:3} Let $(\lambda_j)_j$ be an increasing sequence of positive rational numbers with limit $\lambda$. Then $(\lambda_j\varphi_j)_j$ is a quasi-equisingular approximation of $\lambda\varphi$.
    Since $\lambda_j/\lambda\to 1$, \cref{lma:dSsmallmult} gives
    \[
        d_S\left((\lambda_j/\lambda)\varphi_j|_{\widetilde{Y}},\varphi_j|_{\widetilde{Y}}\right)\to 0.
    \]
    After scaling the background class, \cref{lma:dS_res} shows that $\lambda_j\varphi_j|_{\widetilde{Y}}$ and $\lambda\varphi_j|_{\widetilde{Y}}$ have the same $d_S$-limit. The assertion follows.

    \ref{itm:prop:tracelinear:4} By \cref{prop:lor_dS_conv}, we have
    \[
        \varphi_j\lor \psi_j\xrightarrow{d_S} P_{\theta}[\varphi]_{\mathcal{I}}\lor P_{\theta}[\psi]_{\mathcal{I}}.
    \]
    By \cref{prop:PIconc} and \cref{prop:Igoodlinear}, we have
    \[
        P_{\theta}[\varphi]_{\mathcal{I}}\lor P_{\theta}[\psi]_{\mathcal{I}}\sim_P P_{\theta}[\varphi\lor \psi]_{\mathcal{I}}.
    \]
    Therefore, our assertion follows exactly as in the proof of \ref{itm:prop:tracelinear:2}.
\end{proof}

The trace operator is continuous along $d_S$-convergent decreasing sequences.
\begin{proposition}\label{prop:tracedeclimit}
    Let $(\varphi_j)_{j\in I}$ be a decreasing net in $\QPSH(X)$. Assume that there exists a closed real smooth $(1,1)$-form $\theta$ such that $\varphi_j\in \PSH(X,\theta)$ for each $j\in I$. Assume that $\varphi_j\xrightarrow{d_S}\varphi\in \QPSH(X)$ and $\nu(\varphi,Y)=0$. Then
    \[
        \Tr_Y(\varphi_j)\xrightarrow{d_S}\Tr_Y(\varphi).
    \]
\end{proposition}
In view of \cref{cor:quasi-equichar}, the trace operator preserves the property of being a quasi-equisingular approximation, thereby solving the problem in \cref{ex:naive_res1}.
\begin{proof}
    It is enough to prove the sequential assertion. Indeed, suppose that the net statement fails. Then there is $\delta>0$ such that for every $i_0\in I$ one can find $i\geq i_0$ with
    \[
        d_S(\Tr_Y(\varphi_i),\Tr_Y(\varphi))\geq\delta.
    \]
    Using also $\varphi_i\xrightarrow{d_S}\varphi$, we may choose an increasing sequence $i_j\in I$ such that
    \[
        d_S(\varphi_{i_j},\varphi)<j^{-1},\quad d_S(\Tr_Y(\varphi_{i_j}),\Tr_Y(\varphi))\geq\delta.
    \]
    The sequence $(\varphi_{i_j})_j$ is decreasing and converges to $\varphi$ in $d_S$, contradicting the sequential assertion.

    Thus we assume that $(\varphi_j)_{j\geq 1}$ is a decreasing sequence and that $\varphi_j\xrightarrow{d_S}\varphi$.

    Take a Kähler form $\omega$ on $X$ so that $\varphi_j,\varphi\in \PSH(X,\omega)_{>0}$.  Replacing $\varphi_j$ by $P_{\omega}[\varphi_j]$, we may then assume that $(\varphi_j)_j$ is a decreasing sequence of model potentials. By \cref{prop:decPorderPenvconve}, we have
    \[
        P_{\omega}[\varphi]=\inf_{j\geq 1}\varphi_j.
    \]
    By \cref{prop:ds0char}, we may assume from now on that
    \[
        \varphi=\inf_{j\geq 1}\varphi_j.
    \]
    After adding a Kähler form to $\omega$, we may assume that $\omega_{\varphi}$ is a Kähler current.

    For each $j\geq 1$, choose a quasi-equisingular approximation $(\varphi_j^k)_{k\geq 1}$ of $\varphi_j$ in $\PSH(X,\omega)$. Using the standard construction in \cref{thm:qequi}, we may choose these approximations compatibly in $j$, so that for each fixed $k$,
    \[
        \varphi_{j+1}^k\preceq \varphi_j^k
    \]
    We shall choose an increasing sequence of integers $k_j$ and put
    \[
        \eta_j\coloneqq \varphi_j^{k_j}.
    \]
    The sequence $(k_j)_j$ is chosen inductively: Choose $k_1>1$, for $j>1$, we choose $k_j$ large enough so that
    \begin{enumerate}
        \item\label{itm:chapter-trace:L389:1} $k_{j}\geq k_{j-1}$;
        \item\label{itm:chapter-trace:L389:2} $d_{S,\omega}(\eta_j,P_{\omega}[\varphi_j]_{\mathcal{I}})\leq 1/(j-1)$;
        \item\label{itm:chapter-trace:L389:3} $d_{S,\omega|_{\widetilde{Y}}}(\eta_j|_{\widetilde{Y}},\Tr_Y^{\omega}(\varphi_j))<1/(j-1)$.
    \end{enumerate}
    It follows from the compatibility in $j$ and from $k_j\geq k_{j-1}$ that $(\eta_j)_j$ is decreasing in singularity types. Indeed,
    \[
        \eta_{j+1}=\varphi_{j+1}^{k_{j+1}}\preceq \varphi_{j+1}^{k_j}\preceq \varphi_j^{k_j}=\eta_j.
    \]
    By \cref{thm:contPI}, we have
    \[
        P_{\omega}[\varphi_j]_{\mathcal{I}}\xlongrightarrow{d_S} P_{\omega}[\varphi]_{\mathcal{I}}.
    \]
    By \ref{itm:chapter-trace:L389:2} and \cref{cor:quasi-equichar}, we have $\eta_j\xrightarrow{d_S} P_{\omega}[\varphi]_{\mathcal{I}}$. Observe that by \cref{prop:decPorderPenvconve},
    \[
        P_{\omega}[\varphi]_{\mathcal{I}}=\inf_{j>0} P_{\omega}\left[\eta_j\right].
    \]
    By \cref{prop:analysingcompPandPI}, for each $j>0$, $P_{\omega}[\eta_j]$ has analytic singularities. By \cref{thm:equising_cond_general}, $(P_{\omega}[\eta_j])_j$ is a quasi-equisingular approximation of $P_{\omega}[\varphi]_{\mathcal{I}}$. Therefore, by definition of the trace operator and \itemref{itm:prop:tracelinear:1},
    \[
        \eta_j|_{\widetilde{Y}}\xrightarrow{d_S} \Tr_Y(\varphi).
    \]
    It follows from \ref{itm:chapter-trace:L389:3} that
    \[
        \Tr_Y(\varphi_j)\xrightarrow{d_S}\Tr_Y(\varphi).
    \]
    We conclude the proof.
\end{proof}

\begin{example}\label{ex:trace:L407}
    The trace operator is not continuous along increasing sequences. Let us consider the case $X=\mathbb{P}^2$ with coordinates $(z_1,z_2)$ on $\mathbb{C}^2\subseteq X$. Let $\omega_{\mathrm{{FS}}}$ denote the Fubini--Study metric.
    The subvariety $Y\cong \mathbb{P}^1$ is defined by $z_2=0$. Consider an increasing sequence $(\varphi_j)_j$ in $\PSH(X,\omega_{\mathrm{FS}})$, whose potentials near $(0,0)$ are given by
    \[
    \log |z_1|^2\lor \left( j^{-1}\log |z_2|^2\right)+\mathcal{O}(1).
    \]
    The pointwise restriction of these potentials to $Y$ are given locally by
    \[
    \log |z_1|^2+\mathcal{O}(1).
    \]
    On the other hand, locally
    \[
    \log |z_1|^2\lor \left( j^{-1}\log |z_2|^2\right)\to 0
    \]
    almost everywhere as $j\to\infty$. So the trace operator is not continuous along the sequence $(\varphi_j)_j$.
\end{example}

\begin{lemma}\label{lma:rescommpullback}
    Let $\pi\colon Z\rightarrow X$ be a proper bimeromorphic morphism, where $Z$ is a connected Kähler manifold. Assume that $W$ (resp. $Y$) is an irreducible analytic subset of $Z$ (resp. $X$) such that the restriction $\Pi \colon W\rightarrow Y$ of $\pi$ is defined and is bimeromorphic, so that we have the following commutative diagram
    \[
    \begin{tikzcd}
        \widetilde{W}\arrow[r]\arrow[d,"\widetilde{\Pi}"] & W \arrow[d, "\Pi"] \arrow[r, hook] & Z \arrow[d, "\pi"] \\
        \widetilde{Y}\arrow[r] & Y \arrow[r, hook]                  & X.
\end{tikzcd}
    \]
    Then for any $\varphi\in \QPSH(X)$ with $\nu(\varphi,Y)=0$, we have
\begin{equation}\label{eq:rescommpullback}
    \widetilde{\Pi}^*\Tr_{Y}(\varphi)\sim_P\Tr_W(\pi^*\varphi).
    \end{equation}
\end{lemma}
\begin{proof}
We first observe that by Zariski's main theorem, $\nu(\pi^*\varphi,W)=0$. So the right-hand side of \eqref{eq:rescommpullback} makes sense.

\textbf{Step~1}.
    Assume that $\varphi$ has analytic singularities.
    It suffices to apply \cref{ex:resanalyt} to reformulate \eqref{eq:rescommpullback} as
    \[
    \widetilde{\Pi}^*(\varphi|_{\widetilde{Y}})\sim_P(\pi^*\varphi)|_{\widetilde{W}}.
    \]
    In fact, the strict equality holds, which is nothing but the functoriality of pullbacks.

\textbf{Step~2}.
    Next we handle the general case. Choose a smooth closed real $(1,1)$-form $\theta$ such that $\theta_{\varphi}$ is a Kähler current. Take a quasi-equisingular approximation $(\varphi_j)_j$ of $\varphi$ in $\PSH(X,\theta)$. By \cref{cor:quasi_equising_bimero}, $(\pi^*\varphi_j)_j$ is a quasi-equisingular approximation of $\pi^*\varphi$.
    From Step~1, we know that for each $j$,
    \[
    \widetilde{\Pi}^*\Tr_{Y}(\varphi_j)\sim_P\Tr_W(\pi^*\varphi_j).
    \]
    Letting $j\to\infty$, the right-hand side converges to $\Tr_W(\pi^*\varphi)$ by the definition of the trace operator. On the left-hand side, the $d_S$-convergence of $\Tr_Y(\varphi_j)$ to $\Tr_Y(\varphi)$, together with \cref{lma:Ipartorderinv} and \cref{cor:dSconv_bimero}, implies \eqref{eq:rescommpullback}.
\end{proof}

\begin{proposition}\label{prop:OT2}
    Let $\varphi\in \QPSH(X)$ with $\nu(\varphi,Y)=0$. Assume that $Y$ is smooth.
    Then for any $\lambda>0$, we have
    \begin{equation}\label{eq:OT}
        \mathcal{I}(\lambda \Tr_Y(\varphi))\subseteq \Rest_Y\mathcal{I}(\lambda \varphi).
    \end{equation}
\end{proposition}
See \cref{def:restidealsheaf} for the definition of $\Rest_Y$.

\begin{proof}
Take a Kähler form $\omega$ on $X$ such that $\omega_{\varphi}$ is a Kähler current.

Let $(\varphi_j)_j$ be a quasi-equisingular approximation of $\varphi$ in $\PSH(X,\omega)$.

For each $j\geq 1$, we have
\[
\Tr_Y(\varphi) \preceq_P \varphi_j|_Y.
\]
Indeed, $\Tr_Y(\varphi)$ is $P$-equivalent to the decreasing $d_S$-limit of
the restrictions $\varphi_j|_Y$, hence it is more singular than each
$\varphi_j|_Y$.

For any $\lambda'>\lambda>0$, we can find $j>0$ so that
\[
    \mathcal{I}(\lambda' \varphi_j)\subseteq \mathcal{I}(\lambda \varphi).
\]
By \cref{cor:PimpliesI} and \cref{thm:OT}, we have
\[
    \mathcal{I}(\lambda' \Tr_Y(\varphi))\subseteq \mathcal{I}(\lambda' \varphi_j|_Y)\subseteq \Rest_Y\mathcal{I}(\lambda'\varphi_j)\subseteq \Rest_Y\mathcal{I}(\lambda \varphi).
\]
Thanks to the strong openness theorem \cref{thm:strongopen}, we conclude \eqref{eq:OT}.
\end{proof}

Lastly, we turn our attention to global sections. For this we will need the following global Ohsawa--Takegoshi extension theorem for the trace operator:

\begin{theorem}\label{thm:OT_ext_global}
    Let $L$ be a big line bundle on $X$ and let $\theta$ be a closed real smooth $(1,1)$-form on $X$ representing $c_1(L)$.
Suppose that $\varphi\in \PSH(X,\theta)$ and $\theta_{\varphi}$ is a Kähler current. Assume that $Y$ is smooth and $\nu(\varphi,Y)=0$.
Let $T$ be a holomorphic line bundle on $X$.
Then there exists $k_0$ such that for all $k \geq k_0$ and $s \in \mathrm{H}^0(Y, T|_Y\otimes L|_Y^k\otimes \mathcal{I}(k \Tr_Y^{\theta}(\varphi)))$, there exists an extension $\tilde s \in \mathrm{H}^0(X, T\otimes L^k \otimes \mathcal{I}(k\varphi))$.
\end{theorem}
It is of interest to know whether one can control the $L^2$-norm of $\tilde{s}$ in the above result.

\begin{proof}
Fix a Kähler form $\omega$ on $X$ and a Hermitian metric $h_T$ on $T$.
We may assume that $Y\neq X$ and that $\theta_{\varphi} \geq 3\delta \omega$ for some $\delta>0$. Let $(\varphi_j)_j$ be the decreasing quasi-equisingular approximation of $\varphi$ in $\PSH(X,\theta)$. We can assume that $\theta_{\varphi_j} \geq 2\delta \omega$ for all $j\geq 1$. Also, there exists $\epsilon_0 >0$ such that $\theta_{(1+\epsilon)\varphi_j} \geq \delta \omega$ for any $\epsilon \in (0,\epsilon_0)$. Take $k_0=k_0(\delta,h_T)$ as in \cref{thm:OT_ext}.

We fix $k \geq k_0$ and $s \in \mathrm{H}^0(Y, T|_Y\otimes L|_Y^k\otimes \mathcal{I}(k\Tr_Y^{\theta}(\varphi)))$. By \cref{thm:strongopen}, there exists $\epsilon \in (0,\epsilon_0)$ such that $s \in \mathrm{H}^0(Y, T|_Y\otimes L|_Y^k\otimes \mathcal{I}(k(1+\epsilon)\Tr_Y^{\theta}(\varphi)))$.

Since $\Tr_Y^{\theta}(\varphi) \preceq_P \varphi_j|_Y$, \cref{cor:PimpliesI} implies that $s \in \mathrm{H}^0(Y, T|_Y\otimes L|_Y^k\otimes \mathcal{I}(k(1+\epsilon)\varphi_j|_Y))$. By \cref{thm:OT_ext}, there exists $\tilde s_j \in \mathrm{H}^0(X, T\otimes L^k \otimes \mathcal{I}(k(1+\epsilon)\varphi_j))$ such that $\tilde s_j|_Y =s$, for all $j$.

But by definition of quasi-equisingular approximation, we obtain that for high enough $j$ the inclusion $\mathcal{I}(k(1+\epsilon)\varphi_j) \subseteq \mathcal{I}(k\varphi)$ holds. As a result, $\tilde s_j \in \mathrm{H}^0(X, T\otimes L^k \otimes \mathcal{I}(k\varphi))$ for high enough $j$, finishing the argument.
\end{proof}

\section{Relation to the classical restricted volumes}\label{sec:rv}
Let $X$ be a connected compact Kähler manifold of dimension $n$ and $Y$ be a connected submanifold of dimension $m$. Fix a big class $\alpha\in \mathrm{H}^{1,1}(X,\mathbb{R})$. Take a closed smooth real $(1,1)$-form $\theta\in \alpha$. Fix a Kähler form $\omega$ on $X$.

Recall that the notions of non-Kähler locus and non-nef locus were defined in \cref{def:nKlocus} and \cref{def:nnlocus}.

When $Y\not\subseteq \nK(\alpha)$, Matsumura (\cite[Definition~1.4]{Mat13}) defines the restricted volume of $\{\theta\}$ to $Y$ in the following manner:
\begin{equation}\label{eq:rest_vol_Mat_def}
        \vol_{X|Y}(\alpha)\coloneqq \sup_{\varphi} \int_Y \left(\theta|_Y+\ddc \varphi|_Y\right)^{m},
\end{equation}
where $\varphi$ runs over elements in $\PSH(X,\theta)$ with analytic singularities such that $\varphi|_Y\not\equiv -\infty$. This definition is independent of the choice of $\theta$.

In case $Y\not\subseteq \nn(\alpha)$, Collins--Tosatti \cite{CT22} extend the above definition of restricted volume:
\begin{equation}\label{eq:rest_vol_CT_def}
        \vol_{X|Y}(\alpha)\coloneqq \lim_{\epsilon \searrow 0}\sup_{\varphi} \int_Y \left(\theta|_Y + \epsilon \omega|_Y +\ddc \varphi|_Y\right)^m.
\end{equation}
where $\varphi$ runs over elements in $\PSH(X,\theta+\epsilon\omega)$ with analytic singularities such that $\varphi|_Y\not\equiv -\infty$. This definition is independent of the choice of $\theta$.

These definitions extend the more classical definition in the algebraic setting due to \cite{ELMNP09}.

\begin{proposition}\label{prop:restvol_voltrace}
Assume that $Y\not\subseteq \nK(\alpha)$. Then
    \begin{equation}\label{eq:volrest}
        \vol_{X|Y}(\alpha)= \int_Y \left( \theta|_Y+\ddc \Tr_Y^{\theta}(V_{\theta}) \right)^m =\int_Y (\theta|_Y + \ddc V_\theta|_Y)^m.
    \end{equation}
\end{proposition}

\begin{proof}
We start with the first equality of \eqref{eq:volrest}. Since $Y\not\subseteq \nK(\alpha)$, $V_{\theta}|_Y \not \equiv -\infty$ as a consequence of \cref{lma:Vthetalocbdd}, hence also $\nu(V_{\theta},Y)=0$.

Take a quasi-equisingular approximation $(\varphi_j)_{j>0}$ of $V_{\theta}$ with $\varphi_j\in \PSH(X,\theta+\epsilon_j\omega)$ for a decreasing sequence $(\epsilon_j)_{j>0}$ converging to $0$. By the monotonicity theorem \cref{thm:mono}, we have for each $j>0$,
\[
    \int_Y \left(\theta|_Y+\epsilon_j\omega|_Y+\ddc \varphi_j|_Y\right)^m\geq \int_Y \left(\theta|_Y+\ddc\chi|_Y\right)^m
\]
for every $\chi\in \PSH(X,\theta)$ with $\chi|_Y\not\equiv -\infty$. Since $\chi$ is arbitrary, we conclude that
\[
    \lim_{j\to\infty}\int_Y \left(\theta|_Y+\epsilon_j\omega|_Y+\ddc \varphi_j|_Y\right)^m\geq \vol_{X|Y}(\alpha).
\]
For any fixed $\epsilon>0$, we have
\[
\begin{split}
\lim_{j\to\infty} \int_Y \left(\theta|_Y+\epsilon_j\omega|_Y+\ddc \varphi_j|_Y\right)^m\leq & \lim_{j\to\infty} \int_Y \left(\theta|_Y+\epsilon\omega|_Y+\ddc \varphi_j|_Y\right)^m\\
=& \int_Y \left(\theta|_Y+\epsilon\omega|_Y+\ddc \Tr_Y^{\theta+\epsilon\omega}(V_{\theta})\right)^m.
\end{split}
\]
Letting $\epsilon\to 0+$ and applying \cref{ex:tracedefinedposmass}, we conclude that the $\leq$ direction in the first equality of \eqref{eq:volrest} holds.
Moreover, $Y\not\subseteq \nK(\alpha)$ implies $\vol_{X|Y}(\alpha)>0$ by \eqref{eq:rest_vol_Mat_def}. Hence the limiting mass in \cref{ex:tracedefinedposmass} is positive, and $\Tr_Y^\theta(V_\theta)$ is indeed defined in the present situation.

   For the reverse direction, by definition, for any fixed $\epsilon>0$, we have
    \[
    \begin{aligned}
    \int_Y \left(\theta|_Y+\epsilon_j\omega|_Y+\ddc\varphi_j|_Y\right)^m
    &\leq \int_Y \left(\theta|_Y+\epsilon\omega|_Y+\ddc\varphi_j|_Y\right)^m\\
    &\leq \vol_{X|Y}(\{\theta\}+\epsilon\{\omega\})
    \end{aligned}
    \]
    for all large enough $j$. Therefore, for large $j$, by the monotonicity theorem \cref{thm:mono}, we have
    \[
        \int_Y \left(\theta|_Y+\epsilon_j\omega|_Y+\ddc \Tr_Y^{\theta}(V_{\theta})\right)^m\leq \vol_{X|Y}(\{\theta\}+\epsilon\{\omega\}).
    \]
    Letting $j\to\infty$ and $\epsilon \searrow 0$, using the continuity of $\vol_{X|Y}$ (\cite[Corollary~4.11]{Mat13}) together with \cref{ex:tracedefinedposmass}, we conclude the first equality of \eqref{eq:volrest}.

Now we address the second equality. Due to \cref{thm:mono}, the defining formula \eqref{eq:rest_vol_Mat_def}, and the definition of $V_\theta$, we obtain that
\[
\vol_{X|Y}(\{\theta\}) \leq \int_Y (\theta|_Y + \ddc V_\theta|_Y)^m.
\]
The reverse inequality now follows from the first equality of \eqref{eq:volrest}, \cref{thm:mono} and the fact that $V_\theta|_Y \preceq_P \Tr_Y^\theta(V_\theta)$ as proved in \cref{prop:Trdominarest}.
\end{proof}

\begin{theorem}\label{thm:rest_volume_tr_volume} If $Y\not\subseteq \nn(\alpha)$, then
\begin{equation}\label{eq:non-nef_rest_vol}
\begin{aligned}
\vol_{X|Y}(\alpha) &= \lim_{\epsilon \to 0+}\int_Y \Bigl(\theta|_Y + \epsilon \omega|_Y +  \ddc V_{\theta+ \epsilon \omega}|_Y\Bigr)^m\\
=&\lim_{\epsilon \to 0+}\int_Y \left(\theta|_Y + \epsilon \omega|_Y +  \ddc \Tr_Y^{\theta+ \epsilon \omega}(V_{\theta+ \epsilon \omega}) \right)^m \\
&=\int_Y \left( \theta|_Y+\ddc \Tr_Y^{\theta}(V_{\theta}) \right)^m,
\end{aligned}
\end{equation}
where the last term is understood as $0$ if $\Tr_Y^{\theta}(V_{\theta})$ is not defined.
\end{theorem}

\begin{proof} Since $\nu(V_\theta, Y)=0$, we have $Y\not\subseteq \nK(\alpha+\epsilon\{\omega\})$ for all $\epsilon>0$. Therefore, by \eqref{eq:rest_vol_CT_def} and \eqref{eq:volrest}, only the last equality of \eqref{eq:non-nef_rest_vol} needs to be argued.

Thanks to \cref{lma:Vthetaepsconv}, we have
\[
    V_{\theta+\epsilon\omega}\xrightarrow{d_S}V_{\theta}
\]
as $\epsilon\to 0+$.

Therefore, using \cref{prop:tracedeclimit}, we find
\[
\Tr_Y\left( V_{\theta+\epsilon\omega} \right) \xrightarrow{d_S} \Tr_Y\left(V_{\theta}\right).
\]
Set
\[
M\coloneqq \lim_{\epsilon_0\to 0+}\int_Y \left(\theta|_Y+\epsilon_0\omega|_Y\right)^m_{\Tr_Y^{\theta+\epsilon_0\omega}(V_{\theta})}.
\]
By \cref{ex:tracedefinedposmass}, if $\Tr_Y^{\theta}(V_{\theta})$ is defined, then
\[
    M=\int_Y (\theta|_Y)^m_{\Tr_Y^{\theta}(V_{\theta})},
\]
whereas if $\Tr_Y^{\theta}(V_{\theta})$ is not defined, then $M=0$. Thus $M$ is the last term of \eqref{eq:non-nef_rest_vol}, with the convention made in the statement.

Thanks to \cref{thm:convdS}, for any $\epsilon_0>0$, we have
\[
\begin{aligned}
\lim_{\epsilon\to 0+}\int_Y \left(\theta|_Y + \epsilon_0 \omega|_Y\right)^m_{\Tr_Y^{\theta + \epsilon_0 \omega}(V_{\theta + \epsilon \omega})}
&= \int_Y (\theta|_Y + \epsilon_0 \omega|_Y)^m_{\Tr_Y^{\theta + \epsilon_0 \omega}(V_{\theta})}.
\end{aligned}
\]
By the definition of $M$, this implies the existence of $(\epsilon_j)_{j>0}$ such that $0<\epsilon_j<1/j$ and
\begin{equation}\label{eq:limit_1}
\lim_{j\to\infty}\int_Y \left(\theta|_Y + j^{-1}\omega|_Y\right)^m_{\Tr_Y^{\theta + j^{-1} \omega}(V_{\theta + \epsilon_j \omega}) } = M.
\end{equation}
Moreover, for each $j>0$,
\begin{equation}\label{eq:limit_2}
\begin{aligned}
\int_Y \left(\theta|_Y + j^{-1} \omega|_Y\right)^m_{\Tr_Y^{\theta + j^{-1} \omega}(V_{\theta + \epsilon_j \omega}) } \geq & \int_Y \left(\theta|_Y + \epsilon_j \omega|_Y\right)^m_{\Tr_Y^{\theta + \epsilon_j\omega}(V_{\theta + \epsilon_j \omega}) } \\
\geq & M.
\end{aligned}
\end{equation}
Indeed, the first inequality follows from the multi-linearity of the non-pluripolar products \itemref{itm:prop:npp1:6}. For the second one, if $M=0$ there is nothing to prove; if $M>0$, then \cref{ex:tracedefinedposmass} shows that $\Tr_Y^\theta(V_\theta)$ is defined, and the inequality follows from the multi-linearity of the non-pluripolar products \itemref{itm:prop:npp1:6} and the monotonicity theorem \cref{prop:mono2}.
Putting \eqref{eq:limit_1} and \eqref{eq:limit_2} together, it follows that
\[
\lim_{j\to\infty}\int_Y \left(\theta|_Y + \epsilon_j \omega|_Y\right)^m_{\Tr_Y^{\theta + \epsilon_j \omega}(V_{\theta + \epsilon_j \omega}) } = M.
\]
Finally, since $\int_Y \left(\theta|_Y + \epsilon \omega|_Y +  \ddc \Tr_Y^{\theta+ \epsilon \omega}(V_{\theta+ \epsilon \omega}) \right)^m$ depends monotonically on $\epsilon>0$, we conclude that
\[
\lim_{\epsilon\to 0+}\int_Y \left(\theta|_Y + \epsilon \omega|_Y +  \ddc \Tr_Y^{\theta+ \epsilon \omega}(V_{\theta+ \epsilon \omega}) \right)^m=M.
\]
This completes the proof.
\end{proof}

\section{Restricted volumes of line bundles}\label{sec:resvol}
Let $X$ be a connected projective manifold of dimension $n$ and $Y\subseteq X$ be a connected submanifold of dimension $m$. Consider a big line bundle $L$ on $X$, a Hermitian metric $h_0$ on $L$ with $\theta=c_1(L,h_0)$.
Let $A$ be a very ample line bundle on $X$. Take a Hermitian metric $h_A$ on $A$ such that $\omega=\ddc h_A$ is a Kähler form.

Using the trace operator, one could prove the following generalization of \cref{thm:DXmain1}:
\begin{theorem}\label{thm:rest_volume}
Let $h$ be a (singular) plurisubharmonic metric on $L$ with $\nu(\ddc h,Y)=0$. Assume that
\begin{equation}\label{eq:traceposmasscond2}
    \lim_{\epsilon\to 0+} \int_Y\left(\Tr_Y^{c_1(L|_Y)+\epsilon\omega|_Y}\bigl(c_1(L,h)\bigr) \right)^m>0.
\end{equation}
Then for any holomorphic line bundle $T$ on $X$ we have
\begin{equation}\label{eq:DXmainrelative}
    \int_Y \left(\Tr_Y^{c_1(L|_Y)}\bigl(c_1(L,h)\bigr)\right)^m = \lim_{k \to \infty} \frac{m!}{k^m} h^0\left(Y,T|_Y\otimes L|_Y^k \otimes \Rest_Y\left(\mathcal{I}(h^{k})\right) \right).
\end{equation}
\end{theorem}
Recall that $\Rest_Y$ is defined in \cref{def:restidealsheaf}. Observe that by \cref{ex:tracedefinedposmass}, \eqref{eq:traceposmasscond2} implies that $\Tr_Y^{c_1(L|_Y)}(c_1(L,h))$ is defined. So \eqref{eq:DXmainrelative} is defined.

We will identify $h$ with $\varphi\in \PSH(X,\theta)$ as in \eqref{eq:htwist}.

For the proof, recall that the ideal sheaf $\mathcal{I}_{\infty}$ is introduced in \cref{def:Iinfty}.
\begin{proof}
\textbf{Step~1}. We first assume that $\ddc h$ is a Kähler current.
Let $(\varphi_j)_j$ be a quasi-equisingular approximation of $\varphi$ in $\PSH(X,\theta)$.
After possibly replacing $(\varphi_j)_j$ by a subsequence, there exists $\epsilon_0 \in (0,1)\cap \mathbb{Q}$ such that the $\theta_{(1-\epsilon) \varphi_j}$'s are also Kähler currents for any $\epsilon \in (0,\epsilon_0)$.

For any rational $\epsilon\in (0,\epsilon_0)$, and sufficiently large $k$, we have
\[
\begin{aligned}
&\frac{1}{m!}\int_Y \left(\theta|_Y+\ddc\Tr_Y^{\theta}(\varphi)\right)^m \\
=& \lim_{k \to \infty} \frac{1}{k^m} h^0\left(Y, T|_Y\otimes L|_Y^k\otimes \mathcal{I}\left(k\Tr_Y^{\theta}(\varphi)\right)\right) \quad \textup{by \cref{thm:DXmain1}}\\
\leq & \varliminf_{k \to \infty} \frac{1}{k^m} h^0\left(Y, T|_Y\otimes L|_Y^k\otimes \Rest_Y\left(\mathcal{I}(k\varphi)\right)\right) \quad \textup{ by \cref{prop:OT2}}\\
\leq & \varlimsup_{k \to \infty} \frac{1}{k^m} h^0\left(Y, T|_Y\otimes L|_Y^k\otimes \Rest_Y\left(\mathcal{I}(k\varphi)\right)\right)\\
\leq & \varlimsup_{k \to \infty} \frac{1}{k^m} h^0\left(Y, T|_Y\otimes L|_Y^k\otimes \Rest_Y\mathcal{I}\left(k\varphi_j\right)\right) \\
\leq & \varlimsup_{k \to \infty} \frac{1}{k^m} h^0\left(Y, T|_Y\otimes L|_Y^k\otimes \Rest_Y\mathcal{I}_{\infty}\left((1-\epsilon )k\varphi_j\right)\right)  \quad \textup{by \cref{lma:IandIinf}}\\
\leq & \varlimsup_{k \to \infty} \frac{1}{k^m} h^0\left(Y, T|_Y\otimes L|_Y^k\otimes \mathcal{I}_{\infty}\left((1-\epsilon )k\varphi_j|_Y\right)\right)\\
\leq & \varlimsup_{k \to \infty} \frac{1}{k^m} h^0\left(Y, T|_Y\otimes L|_Y^k\otimes \mathcal{I}\left((1-\epsilon) k\varphi_j|_Y\right) \right)\\
= & \frac{1}{m!} \int_Y \left(\theta|_Y+(1-\epsilon)\ddc\varphi_j|_Y\right)^m\quad \textup{by \cref{thm:DXmain1}}.
\end{aligned}
\]
Here we used the observation that for any $\psi\in \QPSH(X)$ with analytic singularities with $\psi|_Y\not\equiv -\infty$, we have
\[
\Rest_Y \mathcal{I}_{\infty}(\psi)\subseteq \mathcal{I}_{\infty}\left(\psi|_Y\right).
\]

Letting $\epsilon\to0+$, \cref{lma:dSsmallmult} gives
\[
    (1-\epsilon)\varphi_j|_Y\xrightarrow{d_S}\varphi_j|_Y.
\]
Hence, by \cref{thm:convdS}, the last term converges to
\[
    \frac{1}{m!}\int_Y \left(\theta|_Y+\ddc\varphi_j|_Y\right)^m.
\]
We conclude that
\[
\begin{aligned}
\frac{1}{m!}\int_Y \left(\theta|_Y+\ddc\Tr_Y^{\theta}(\varphi)\right)^m\leq & \varliminf_{k \to \infty} \frac{1}{k^m} h^0\left(Y, T|_Y\otimes L|_Y^k\otimes \Rest_Y\left(\mathcal{I}(k\varphi)\right)\right) \\
\leq & \varlimsup_{k \to \infty} \frac{1}{k^m} h^0\left(Y, T|_Y\otimes L|_Y^k\otimes \Rest_Y\left(\mathcal{I}(k\varphi)\right)\right)\\
\leq & \frac{1}{m!} \int_Y \left(\theta|_Y+\ddc\varphi_j|_Y\right)^m.
\end{aligned}
\]
Letting $j\to \infty$, we conclude \eqref{eq:DXmainrelative}.

\textbf{Step~2}. We handle the general case. After adding a constant to $\varphi$, we may assume that $\varphi\leq 0$.

Using \cref{prop:OT2} and \cref{thm:DXmain1} we obtain that
\[
\begin{aligned}
    \int_Y \left(\theta|_Y+\ddc \Tr_Y^{\theta}(\varphi)\right)^m &= \lim_{k \to \infty} \frac{m!}{k^m} h^0\left(Y,T|_Y \otimes L|_Y^k \otimes \mathcal{I}\left(k\Tr_Y^{\theta}(\varphi)\right)\right)\\
    & \leq \varliminf_{k \to \infty} \frac{m!}{k^m} h^0\left(Y,T|_Y \otimes L|_Y^k \otimes \Rest_Y\left(\mathcal{I}(k\varphi)\right)\right).
\end{aligned}
\]
Now we address the other direction in \eqref{eq:DXmainrelative}. Let $\phi \in \mathrm{H}^0(X,A)$ be a section that does not vanish identically on $Y$. Such $\phi$ exists since $A$ is very ample.

We fix $k_0 \in \mathbb{N}$. For any $k \geq 0$, we have $k = q k_0 + r$ with $q,r \in \mathbb{N}$ and $r \in \{0,\ldots, k_0-1\}$. Also, we have an injective linear map
\[
\mathrm{H}^0\left(Y,T|_Y \otimes L|_Y^k \otimes \Rest_Y\mathcal{I}\left(k \varphi\right)\right)\xrightarrow{\cdot \phi^{\otimes q}} \mathrm{H}^0\left(Y,T|_Y\otimes L|_Y^k \otimes A|_Y^q \otimes \Rest_Y\mathcal{I}\left(k \varphi\right)\right).
\]
Therefore,
\[
\begin{aligned}
 &\varlimsup_{k\to\infty} \frac{m!}{k^m} h^0\left(Y,T|_Y \otimes L|_Y^k \otimes \Rest_Y\mathcal{I}\left(k\varphi\right)\right)  \\
 \leq & \varlimsup_{k\to\infty} \frac{m!}{k^m} h^0\left(Y,T|_Y \otimes L|_Y^k \otimes A|_Y^q\otimes \Rest_Y\mathcal{I}\left(k\varphi\right)\right) \\
 \leq & \max_{r=0,\ldots,k_0-1}\frac{1}{k_0^m} \varlimsup_{q\to\infty} \frac{m!}{q^m} h^0\left(Y,T|_Y\otimes L|_Y^{q k_0} \otimes A|_Y^q\otimes L|_Y^{r}\otimes \Rest_Y\mathcal{I}\left((qk_0+r)\varphi\right)\right)\\
 \leq & \max_{r=0,\ldots,k_0-1}\frac{1}{k_0^m} \varlimsup_{q\to\infty} \frac{m!}{q^m} h^0\left(Y,T|_Y\otimes L|_Y^{q k_0} \otimes A|_Y^q\otimes L|_Y^{r}\otimes \Rest_Y\mathcal{I}\left(k_0q \varphi\right)\right)\\
 = & \int_Y \left(\theta|_Y+k_0^{-1}\omega|_Y+\ddc \Tr_Y^{\theta+k_0^{-1}\omega}(\varphi)\right)^m\\
 =& \int_Y \left(\theta|_Y+k_0^{-1}\omega|_Y+\ddc \Tr_Y^{\theta}(\varphi)\right)^m,
 \end{aligned}
\]
where the maximum is taken over the finitely many residue classes $k=qk_0+r$, the following inequality uses $k_0q\leq k$, the first equality follows by applying Step~1 for each $r$ to the big line bundle $L^{k_0} \otimes A$, the Kähler current $k_0 \theta_{\varphi} + \omega$, and twisting bundle $T \otimes L^r$, and the second equality follows from \cref{prop:Ppresmass}. Letting $k_0\to\infty$ and using the multi-linearity of the non-pluripolar product \itemref{itm:prop:npp1:6}, we conclude that
\[
\varlimsup_{k\to\infty} \frac{m!}{k^m}  h^0\left(Y,T|_Y \otimes L|_Y^k \otimes \Rest_Y\mathcal{I}\left(k \varphi\right)\right)\leq \int_Y \left(\theta|_Y+\ddc \Tr_Y^{\theta}(\varphi)\right)^m.
\]
This completes the proof.
\end{proof}

\begin{theorem}\label{thm:rest_volume_2} Let $\varphi\in \PSH(X,\theta)$ such that $\nu(\varphi,Y) =0$. Assume that $\theta_{\varphi}$ is a Kähler current.
Then for any holomorphic line bundle $T$ on $X$ we have
\[
    \int_Y \left(\theta|_Y+\ddc \Tr_Y^{\theta}(\varphi)\right)^m  = \lim_{k \to \infty} \frac{m!}{k^m} \dim_{\mathbb C}\left\{ s|_Y:  s\in \mathrm{H}^0\left(X,T\otimes L^k \otimes \mathcal{I}(k\varphi)\right)\right\}.
\]
\end{theorem}

\begin{proof}
This is a consequence of \cref{thm:DXmain1}, \cref{thm:OT_ext_global} and \cref{thm:rest_volume}:
\[
\begin{aligned}
    \int_Y \left(\theta|_Y+\ddc \Tr_Y^{\theta}(\varphi)\right)^m  & = \lim_{k \to \infty} \frac{m!}{k^m} h^0\left(Y,T|_Y\otimes L|_Y^k\otimes \mathcal{I}\left(k \Tr_Y^{\theta}(\varphi)\right)\right)\\
    & \leq  \varliminf_{k \to \infty} \frac{m!}{k^m} \dim_{\mathbb C}\left\{ s|_Y:  s\in \mathrm{H}^0\left(X,T\otimes L^k \otimes \mathcal{I}(k\varphi)\right)\right\}\\
    & \leq  \varlimsup_{k \to \infty} \frac{m!}{k^m} \dim_{\mathbb C}\left\{ s|_Y:  s\in \mathrm{H}^0\left(X,T\otimes L^k \otimes \mathcal{I}(k\varphi)\right)\right\}\\
  & \leq  \lim_{k \to \infty} \frac{m!}{k^m} h^0\left(Y ,T|_Y\otimes L|_Y^k\otimes \Rest_Y\mathcal{I}(k\varphi)\right)\\
  & =  \int_Y \left(\theta|_Y+\ddc \Tr_Y^{\theta}(\varphi)\right)^m.
  \end{aligned}
\]
This completes the proof.
\end{proof}

\begin{remark}
    One could also show that when \eqref{eq:traceposmasscond2} fails, the right-hand side of \eqref{eq:DXmainrelative} is $0$. See \cite{DX24}.
\end{remark}

\section{Analytic Bertini theorems}\label{sec:ABT}
Let $X$ be a connected projective manifold of dimension $n\geq 1$.

The analytic Bertini theorem handles the restriction along a generic subvariety.

\begin{theorem}\label{thm:Bert}
Let $\varphi\in \QPSH(X)$. Let $p\colon X\rightarrow \mathbb{P}^N$ be a morphism ($N\geq 1$). Define
\[
\mathcal{G}\coloneqq \Bigl\{H\in |\mathcal{O}_{\mathbb{P}^N}(1)|: H'\coloneqq H\cap X \textup{ is smooth and } \mathcal{I}(\varphi|_{H'})=\Rest_{H'}\bigl(\mathcal{I}(\varphi)\bigr)\Bigr\}.
\]
Then $\mathcal{G}\subseteq |\mathcal{O}_{\mathbb{P}^N}(1)|$ is co-pluripolar.
\end{theorem}
Recall that co-pluripolar sets are defined in \cref{def:pluripolarsets}. We adopt the convention that $\mathcal{I}(-\infty)=0$.
\begin{remark}\label{rmk:trace:L787}
Here and in the sequel, we slightly abuse the notation by writing $H\cap X$ for $p^{-1}H$, the scheme-theoretic inverse image of $H$. In other words, $H\cap X\coloneqq H\times_{\mathbb{P}^N} X$.

By definition, any $H\in |\mathcal{O}_{\mathbb{P}^N}(1)|$ such that $p^{-1}H=\emptyset$ lies in $\mathcal{G}$.
\end{remark}
\begin{proof}

Take an ample line bundle $L$ with a smooth Hermitian metric $h$ such that $c_1(L,h)+\ddc\varphi\geq 0$, where $c_1(L,h)$ is the first Chern form of $(L,h)$, namely the curvature form of $h$.
We introduce $\Lambda\coloneqq |\mathcal{O}_{\mathbb{P}^N}(1)|$ to simplify our notations.

\textbf{Step~1}.
We prove that the following set is co-pluripolar:
\[
\begin{split}
\mathcal{G}_{L}\coloneqq \left\{H\in \Lambda: H\cap X \text{ is smooth and } \mathrm{H}^0\Bigl(H\cap X,\omega_{H\cap X}\otimes L|_{H\cap X}\otimes \mathcal{I}(\varphi|_{H\cap X})\Bigr)= \right.\\
\left. \mathrm{H}^0\Bigl(H\cap X,\omega_{H\cap X}\otimes L|_{H\cap X}\otimes \Rest_{H\cap X}\left(\mathcal{I}(\varphi)\right)\Bigr)\right\}.
\end{split}
\]
Here $\omega_{H\cap X}$ denotes the dualizing sheaf of $H\cap X$.

Let $U\subseteq \Lambda\times X$ be the closed subvariety whose $\mathbb{C}$-points correspond to pairs $(H,x)\in \Lambda\times X$ with $p(x)\in H$.
Let $\pi_1\colon U\rightarrow \Lambda$ be the natural projection. We may assume that $\pi_1$ is surjective, as otherwise there is nothing to prove.

Observe that $U$ is a local complete intersection scheme by \emph{Krull's Hauptidealsatz} and \emph{a fortiori} a Cohen--Macaulay scheme. It follows from miracle flatness \cite[Theorem~23.1]{Mat89} that the natural projection $\pi_2\colon U\rightarrow X$ is flat. As the fibers of $\pi_2$ over closed points of $X$ are isomorphic to $\mathbb{P}^{N-1}$, it follows that $\pi_2$ is smooth. Thus, $U$ is smooth as well. Moreover, observe that
\begin{equation}\label{eq:pi2pullvarphiItemp1}
    \mathcal{I}(\pi_2^*\varphi)=\pi_2^*\mathcal{I}(\varphi)
\end{equation}
by \cref{prop:pull-backmis}.

In the following, we will construct pluripolar sets $\Sigma_1\subseteq \Sigma_2 \subseteq \Sigma_3\subseteq \Sigma_4\subseteq \Lambda$ such that the behavior of $\pi_1$ is improved successively on the complement of $\Sigma_i$.

\textbf{Step~1.1}.
The usual Bertini theorem shows that there is a proper Zariski closed set $\Sigma_1\subseteq \Lambda$ such that $\pi_1$ has smooth fibers outside $\Sigma_1$. Enlarging $\Sigma_1$, we may guarantee that $\pi_1$ and $\mathcal{I}(\pi_2^*\varphi)$ are both flat outside $\Sigma_1$. See \cite[Théorème~6.9.1]{EGAIV2}.
Then after further enlarging $\Sigma_1$ we may assume that $p^{-1}(H)$ avoids all associated points of $\mathcal{O}_X/\mathcal{I}(\varphi)$, for all $H\in \Lambda\setminus \Sigma_1$. Let $\pi_{1,H}$ denote the fiber of $\pi_1$ at $H$ and write $i_H\colon \pi_{1,H}\rightarrow U$ for the inclusion morphism.
We arrive at
\[
\Rest_{\pi_{1,H}}(\mathcal{I}(\pi_2^*\varphi))=i_H^*\mathcal{I}(\pi_2^*\varphi)
\]
for all $H\in \Lambda\setminus \Sigma_1$.\footnote{This subtle point was overlooked in the proof of \cite{XiaBer}.}

\textbf{Step~1.2}.
By Grauert's coherence theorem,
\[
\mathcal{F}^i\coloneqq R^i\pi_{1*}\left(\omega_{U/\Lambda}\otimes \pi_2^*L\otimes \mathcal{I}(\pi_2^*\varphi) \right)
\]
is coherent for all $i$. Here $\omega_{U/\Lambda}$ denotes the relative dualizing sheaf of the morphism $U\rightarrow \Lambda$.
Thus, there is a proper Zariski closed set $\Sigma_2\subseteq \Lambda$ such that
\begin{enumerate}
    \item\label{itm:chapter-trace:L815:1} $\Sigma_2\supseteq \Sigma_1$.
    \item\label{itm:chapter-trace:L815:2} The $\mathcal{F}^i$'s are locally free outside $\Sigma_2$.
\end{enumerate}
We write $\mathcal{F}=\mathcal{F}^0$.
 By cohomology and base change \cite[Theorem~III.12.11]{Har}, for any $H\in \Lambda\setminus \Sigma_2$,
 the fiber $\mathcal{F}|_H$ of $\mathcal{F}$ is given by
 \[
\mathcal{F}|_H= \mathrm{H}^0\left(\pi_{1,H},\omega_{U/\Lambda}|_{\pi_{1,H}}\otimes \pi_2^*L|_{\pi_{1,H}}\otimes \Rest_{\pi_{1,H}}(\mathcal{I}(\pi_2^*\varphi)) \right).
\]

\textbf{Step~1.3}.
In order to proceed, we need to make use of the Hodge metric $h_{\mathcal{H}}$ on $\mathcal{F}$ defined in \cite{HPS18}. We briefly recall its definition in our setting.
By \cite[Section~22]{HPS18}, we can find a proper Zariski closed set $\Sigma_3\subseteq \Lambda$ such that
\begin{enumerate}
    \item\label{itm:chapter-trace:L830:1} $\Sigma_3\supseteq \Sigma_2$,
    \item\label{itm:chapter-trace:L830:2} $\pi_1$ is smooth outside $\Sigma_3$,
    \item\label{itm:chapter-trace:L830:3} both $\mathcal{F}$ and $\pi_{1*}\left(\omega_{U/\Lambda}\otimes \pi_2^*L \right)/\mathcal{F}$ are locally free outside $\Sigma_3$, and
    \item\label{itm:chapter-trace:L830:4} for each $i$,
    \[
    R^i\pi_{1*}\left(\omega_{U/\Lambda}\otimes \pi_2^*L \right)
    \]
    is locally free outside $\Sigma_3$.
\end{enumerate}
Then for any $H\in \Lambda\setminus \Sigma_3$,
\[
\begin{gathered}
\mathrm{H}^0\left(H\cap X,\omega_{H\cap X}\otimes L|_{H\cap X}\otimes \mathcal{I}(\varphi|_{H\cap X})\right)\\
\subseteq \mathcal{F}|_H\subseteq \mathrm{H}^0\left(H\cap X,\omega_{H\cap X}\otimes L|_{H\cap X}\right).
\end{gathered}
\]
See \cite[Lemma~22.1]{HPS18}.

Now we can give the definition of the Hodge metric on $\Lambda\setminus \Sigma_3$.
Given any $H\in \Lambda\setminus \Sigma_3$, any $\alpha\in \mathcal{F}|_H$, the Hodge metric is defined as
\[
h_{\mathcal{H}}(\alpha,\alpha):=\int_{X\cap H} |\alpha|^2_{h}\mathrm{e}^{-\varphi}\in [0,\infty].
\]
Observe that $h_{\mathcal{H}}(\alpha,\alpha)<\infty$ if and only if $\alpha\in \mathrm{H}^0(H\cap X,\omega_{H\cap X}\otimes L|_{H\cap X}\otimes \mathcal{I}(\varphi|_{H\cap X}))$. Moreover, $h_{\mathcal{H}}(\alpha,\alpha)>0$ if $\alpha\neq 0$.
It is shown in \cite{HPS18} (cf. \cite[Theorem~3.3.5]{PT18}) that $h_{\mathcal{H}}$ is indeed a singular Hermitian metric, and it extends to a positive metric on $\mathcal{F}$.

\textbf{Step~1.4}.
The determinant $\det h_{\mathcal{H}}$ is singular at all $H\in \Lambda\setminus \Sigma_3$ such that
\[
\mathrm{H}^0\Bigl(H\cap X,\omega_{H\cap X}\otimes L|_{H\cap X}\otimes \mathcal{I}(\varphi|_{H\cap X})\Bigr)\neq \mathcal{F}|_H.
\]
As the map $\pi_2$ is smooth, we have $\pi_2^*\mathcal{I}(\varphi)= \mathcal{I}(\pi_2^*\varphi)$ by \cref{prop:pull-backmis}.
Under the identification $\pi_{1,H}\cong H\cap X$, we have
\[
\Rest_{\pi_{1,H}}\Bigl(\pi_2^*\mathcal{I}(\varphi)\Bigr)\cong \Rest_{H\cap X}\Bigl(\mathcal{I}(\varphi)\Bigr).
\]
Thus, we have the following inclusions:
\[
\begin{split}
&\mathrm{H}^0\Bigl(H\cap X,\omega_{H\cap X}\otimes L|_{H\cap X}\otimes \mathcal{I}\left(\varphi|_{H\cap X}\right)\Bigr)\\
\subseteq & \mathrm{H}^0\Bigl(H\cap X,\omega_{H\cap X}\otimes L|_{H\cap X}\otimes \Rest_{H\cap X}\bigl(\mathcal{I}(\varphi)\bigr)\Bigr),
\end{split}
\]
the right-hand side being $\mathcal{F}|_H$.

Recall that the first inclusion follows from \cref{thm:OT}.
Hence, $\det h_{\mathcal{H}}$ is singular at all $H\in |\mathcal{O}_{\mathbb{P}^N}(1)|\setminus \Sigma_3$ such that
\[
\begin{split}
&\mathrm{H}^0\Bigl(H\cap X,\omega_{H\cap X}\otimes L|_{H\cap X}\otimes \mathcal{I}\left(\varphi|_{H\cap X}\right)\Bigr)\\
\neq & \mathrm{H}^0\Bigl(H\cap X,\omega_{H\cap X}\otimes L|_{H\cap X}\otimes \Rest_{H\cap X}\bigl(\mathcal{I}(\varphi)\bigr)\Bigr).
\end{split}
\]
Let $\Sigma_4$ be the union of $\Sigma_3$ and the set of all such $H$.
Since the Hodge metric $h_{\mathcal{H}}$ is positive (\cite[Theorem~3.3.5]{PT18} and \cite[Theorem~21.1]{HPS18}),
its determinant $\det h_{\mathcal{H}}$ is also positive (\cite[Proposition~1.3]{Rau15} and \cite[Proposition~25.1]{HPS18}), it follows that $\Sigma_4$ is pluripolar. As a consequence, $\mathcal{G}_{L}$ is co-pluripolar.

\textbf{Step~2}.

Fix an ample invertible sheaf $S$ on $X$. The same result holds with $L\otimes S^{\otimes a}$ in place of $L$. Thus, the set
\[
A\coloneqq \bigcap_{a=0}^{\infty}\mathcal{G}_{L\otimes S^{\otimes a}}
\]
is co-pluripolar.
For each $H\in \Lambda$ such that $X\cap H$ is smooth and $\mathcal{I}(\varphi|_{X\cap H})\neq \Rest_{H\cap X}(\mathcal{I}(\varphi))$, let $\mathcal{K}$ be the following cokernel:
\[
0\rightarrow \mathcal{I}\left(\varphi|_{X\cap H}\right)\rightarrow \Rest_{H\cap X}\bigl(\mathcal{I}(\varphi)\bigr)\rightarrow \mathcal{K}\rightarrow 0.
\]
By Serre vanishing theorem, taking $a$ large enough, we may guarantee that
\[
\mathrm{H}^1\Bigl(X\cap H,\omega_{X\cap H}\otimes (L\otimes S^{\otimes a})|_{X\cap H}\otimes \mathcal{I}(\varphi|_{X\cap H})\Bigr)=0
\]
and
\[
\mathrm{H}^0\Bigl(X\cap H,\omega_{X\cap H}\otimes (L\otimes S^{\otimes a})|_{X\cap H}\otimes \mathcal{K}\Bigr)\neq 0.
\]
Then
\[
\begin{split}
\mathrm{H}^0\Bigl(X\cap H,\omega_{X\cap H}\otimes(L\otimes S^{\otimes a})|_{X\cap H}\otimes\mathcal{I}\left(\varphi|_{X\cap H}\right)\Bigr)\neq \\
\mathrm{H}^0\Bigl(X\cap H,\omega_{X\cap H}\otimes(L\otimes S^{\otimes a})|_{X\cap H}\otimes \Rest_{H\cap X}\left(\mathcal{I}(\varphi)\right)\Bigr).
\end{split}
\]
Thus, $H\not\in A$. We conclude that $\mathcal{G}$ is co-pluripolar.
\end{proof}

In the sequel of this section, we fix a base-point free linear system $\Lambda$ on $X$. For the results below that assert connectedness, we assume in addition that the set of connected smooth members of $\Lambda$ is co-pluripolar in $\Lambda$.
\begin{corollary}\label{cor:qpshgeneralres}
    Let $\varphi\in \QPSH(X)$. Then for quasi-every $H\in \Lambda$, we have $\varphi|_H\not\equiv -\infty$.
\end{corollary}
\begin{proof}
    For quasi-every $H$, \cref{thm:Bert} gives
    \[
        \mathcal I(\varphi|_H)=\Rest_H\mathcal I(\varphi).
    \]
    After removing the hyperplanes meeting the associated points of $\mathcal O_X/\mathcal I(\varphi)$ improperly, the right-hand side is non-zero. Since $\mathcal I(-\infty)=0$, this implies $\varphi|_H\not\equiv-\infty$.
\end{proof}

\begin{corollary}\label{cor:ABTfortrace}
    Assume that $n\geq 2$.
    Let $\varphi\in \QPSH(X)$. Then quasi-every $H\in \Lambda$ is connected and smooth, satisfies $\nu(\varphi,H)=0$ and we have
    \[
    \Tr_H(\varphi)\sim_{\mathcal{I}} \varphi|_H.
    \]
\end{corollary}
The connectedness assumption on $H$ is used here because we developed most of the trace theory above for connected submanifolds.
\begin{proof}
    First observe that the set $\{x\in X:\nu(\varphi,x)>0\}$ is a countable union of proper analytic subsets by Siu's semicontinuity theorem \cref{thm:Siusemi}. It follows that a very general element in $\Lambda$ is not contained in this set.

    Fix an ample line bundle $L$ so that there is a smooth Hermitian metric $h_L$ such that $c_1(L,h_L)+\ddc \varphi$ is a Kähler current. Thanks to \cref{thm:Bert} and the connectedness assumption on $\Lambda$, we can find a co-pluripolar set $\Lambda'\subseteq \Lambda$ such that each $H\in \Lambda'$ satisfies the following:
    \begin{enumerate}
        \item\label{itm:chapter-trace:L946:1} $H$ is smooth and connected;
        \item\label{itm:chapter-trace:L946:2} $\nu(\varphi,H)=0$;
        \item\label{itm:chapter-trace:L946:3} $\mathcal{I}(k\varphi|_H)=\Rest_{H}(\mathcal{I}(k\varphi))$ for all $k\in \mathbb{Z}_{>0}$.
    \end{enumerate}
    Put $\alpha_H\coloneqq c_1(L,h_L)|_H$. It follows from \cref{thm:rest_volume} and \cref{thm:DXmain1} that
    \[
    \int_H \left( \alpha_H+\ddc \Tr_H^{c_1(L,h_L)}(\varphi)\right)^{n-1}
    =
    \int_H \left( \alpha_H+\ddc P_{c_1(L,h_L)|_H}\left[\varphi|_H\right]_{\mathcal{I}}  \right)^{n-1}.
    \]
    Since $\varphi|_H\preceq_P \Tr_H(\varphi)$ by \cref{prop:Trdominarest}, we also have $\varphi|_H\preceq_{\mathcal I}\Tr_H(\varphi)$. Hence
    \[
        P_{c_1(L,h_L)|_H}[\varphi|_H]_{\mathcal I}\leq P_{c_1(L,h_L)|_H}[\Tr_H(\varphi)]_{\mathcal I}.
    \]
    The trace is $\mathcal I$-good, so the right-hand side agrees with its $P$-envelope. The equality of the masses above and \cref{thm:Pvarphidiffdef} then imply that the two sides have the same $P$-singularity type, hence
    \[
        \Tr_H(\varphi)\sim_{\mathcal I}\varphi|_H.
    \]
\end{proof}

\begin{lemma}\label{lma:posmasscurrres}
    Assume that $n\geq 2$.
    Let $T$ be a closed positive $(1,1)$-current on $X$ with $\int_X T^n>0$. Then quasi-every $H\in \Lambda$ is connected and smooth, $T|_H$ is well-defined and satisfies
    \[
    \int_H T|_H^{n-1}>0.
    \]
\end{lemma}
\begin{proof}
    Write $T=\theta_{\varphi}$ for some smooth closed real $(1,1)$-form $\theta$ on $X$ and $\varphi\in \PSH(X,\theta)_{>0}$. Thanks to \cref{lma:kahcurrentposmass}, we can find $\psi\in \PSH(X,\theta)$ such that $\theta_{\psi}$ is a K\"ahler current and $\psi\leq \varphi$. By \cref{cor:qpshgeneralres} and the connectedness assumption on $\Lambda$, we can find a co-pluripolar set $\Lambda'\subseteq \Lambda$ such that each $H\in \Lambda'$ satisfies:
    \begin{enumerate}
        \item\label{itm:chapter-trace:L977:1} $H$ is smooth and connected;
        \item\label{itm:chapter-trace:L977:2} the restriction $\psi|_H$ is not identically $-\infty$.
    \end{enumerate}
    Since $\theta_\psi$ is a Kähler current, its restriction to such an $H$ has positive non-pluripolar mass. As $\psi|_H\leq\varphi|_H$, \cref{thm:mono} gives
    \[
        \int_H(\theta|_H+\ddc\varphi|_H)^{n-1} \geq \int_H(\theta|_H+\ddc\psi|_H)^{n-1}>0.
    \]
    We conclude the proof.
\end{proof}

\begin{corollary}\label{cor:tracegeneralwelldef}
    Assume that $n\geq 2$. Let $T$ be a closed positive $(1,1)$-current on $X$ with $\vol T>0$. Then quasi-every $H\in \Lambda$ is connected and smooth, and $\Tr_H^{\{T\}|_H}(T)$ is well-defined and has positive volume.
\end{corollary}
\begin{proof}
    Write $T=\theta+\ddc\varphi$ and set $\psi=P_{\theta}[\varphi]_{\mathcal{I}}$. Then $\psi\sim_{\mathcal{I}}\varphi$ and
    \[
        \int_X \theta_{\psi}^{n}=\vol T>0.
    \]
    By \cref{lma:posmasscurrres}, for quasi-every $H\in \Lambda$, the hypersurface $H$ is connected and smooth, $\psi|_H$ is well-defined, and
    \[
        \int_H \left(\theta|_H+\ddc\psi|_H\right)^{n-1}>0.
    \]
    For such $H$, \cref{cor:ABTfortrace} gives $\Tr_H(\psi)\sim_{\mathcal{I}}\psi|_H$. Since $\varphi\sim_{\mathcal{I}}\psi$, \cref{prop:tracelinear} also gives $\Tr_H(\varphi)\sim_P\Tr_H(\psi)$. The assertion now follows from \cref{ex:tracedefinedposmass}.
\end{proof}

\begin{proposition}\label{prop:varphiP_rest_gen}
    Assume that $n\geq 2$. Let $\varphi,\psi\in \QPSH(X)$. Assume that $\varphi\preceq_P \psi$. Then quasi-every $H\in \Lambda$ is connected and smooth, and $\varphi|_H\preceq_P \psi|_H$.
\end{proposition}
\begin{proof}
    Thanks to \cref{lma:reform_preceqP}, we may replace $\varphi$ by $\varphi\lor \psi$ and assume that $\varphi\sim_P \psi$.
    It suffices to show that $\varphi|_H\sim_P \psi|_H$ for quasi-every $H\in \Lambda$.

    Take a smooth closed real $(1,1)$-form $\theta$ on $X$ so that $\varphi,\psi\in \PSH(X,\theta)_{>0}$. It suffices to compare $\varphi$ and $\psi$ with $P_{\theta}[\varphi]$, so without loss of generality, we may assume that $\psi$ is a model potential in $\PSH(X,\theta)_{>0}$. Up to adding a constant to $\varphi$, we may then assume that $\varphi\leq \psi$. It follows from \cref{lma:pathoenvelope} that we can find a sequence $(\eta_j)_j$ in $\PSH(X,\theta)_{>0}$ such that
    \[
    j^{-1}\eta_j+\left(1-j^{-1}\right)\psi\leq \varphi
    \]
    for all $j\geq 2$. By \cref{cor:qpshgeneralres}, \cref{lma:posmasscurrres}, we can find a co-pluripolar set $\Lambda'\subseteq \Lambda$ such that any $H\in \Lambda'$ satisfies:
    \begin{enumerate}
        \item\label{itm:chapter-trace:L1015:1} $H$ is smooth and connected;
        \item\label{itm:chapter-trace:L1015:2} $\eta_j|_H\in \PSH(H,\theta|_H)_{>0}$ for all $j\geq 2$ and $\psi|_H\in \PSH(H,\theta|_H)_{>0}$.
    \end{enumerate}
    Therefore, taking \cref{prop:Pconc} into account, we arrive at
    \[
    j^{-1} P_{\theta|_H}\bigl[\eta_j|_H\bigr]+\left(1-j^{-1}\right)P_{\theta|_H}\bigl[\psi|_H\bigr]\leq P_{\theta|_H}\bigl[\varphi|_H\bigr]
    \]
    for all $j\geq 2$. The normalized envelopes $P_{\theta|_H}[\eta_j|_H]$ are uniformly bounded in $L^1$, hence
    \[
        j^{-1}P_{\theta|_H}[\eta_j|_H]\xrightarrow{L^1} 0.
    \]
    We conclude that
    \[
    P_{\theta|_H}[\psi|_H]\leq P_{\theta|_H}[\varphi|_H]
    \]
    and hence $\psi|_H\preceq_P \varphi|_H$.
\end{proof}

\begin{lemma}\label{lma:Igoodrest}
Assume that $n\geq 2$.
    Let $\theta$ be a closed smooth $(1,1)$-form on $X$ representing a big cohomology class and $(\varphi_j)_j$ be a decreasing sequence in $\PSH(X,\theta)$. Assume that $\varphi\in \PSH(X,\theta)$ and $\varphi_j\xrightarrow{d_S}\varphi$. Furthermore, assume that $\varphi\leq \varphi_j$ for each $j>0$.
    Then quasi-every $H\in \Lambda$ is connected and smooth, $\varphi_j|_H\not\equiv-\infty$ for all $j\geq 1$, $\varphi|_H\not\equiv -\infty$, and
    \[
    \varphi_j|_H\xrightarrow{d_S}\varphi|_H.
    \]
\end{lemma}
\begin{proof}
    By \cref{cor:dSconv_changetheta}, we may assume that $\varphi\in \PSH(X,\theta)_{>0}$.

    Using \cref{lma:pathoenvelope}, we can find a decreasing sequence $(\epsilon_j)_j$ in $(0,1)$ with limit $0$ and $\eta_j\in \PSH(X,\theta)_{>0}$ such that $\eta_j\leq \varphi_j$ and
    \[
    \epsilon_j\eta_j+(1-\epsilon_j)\varphi_j\leq \varphi.
    \]
    By \cref{cor:qpshgeneralres}, \cref{lma:posmasscurrres}, we can find a co-pluripolar set $\Lambda'\subseteq \Lambda$ such that any $H\in \Lambda'$ satisfies:
    \begin{enumerate}
        \item\label{itm:chapter-trace:L1048:1} $H$ is smooth and connected;
        \item\label{itm:chapter-trace:L1048:2} $\eta_j|_H\in \PSH(H,\theta|_H)_{>0}$ for all $j\geq 1$ and $\varphi|_H\in \PSH(H,\theta|_H)_{>0}$.
    \end{enumerate}
    Therefore, taking \cref{prop:Pconc} into account, we arrive at
    \[
    \epsilon_j P_{\theta|_H}\bigl[\eta_j|_H\bigr]+\left(1-\epsilon_j\right)P_{\theta|_H}\bigl[\varphi_j|_H\bigr]\leq P_{\theta|_H}\bigl[\varphi|_H\bigr].
    \]
    Since the normalized envelopes $P_{\theta|_H}[\eta_j|_H]$ are relatively compact in $L^1$, we have $\epsilon_jP_{\theta|_H}[\eta_j|_H]\to 0$ in $L^1$. Letting $j\to\infty$, we get
    \[
    \lim_{j\to\infty} P_{\theta|_H}\bigl[\varphi_j|_H\bigr]\leq P_{\theta|_H}\bigl[\varphi|_H\bigr].
    \]
    The reverse inequality is trivial, so we get an equality here.
    By \cref{thm:mono} and \cref{prop:vol_limit_model}, we conclude that
    \[
    \lim_{j\to\infty}\int_H \Bigl(\theta|_H+\ddc \varphi_j|_H\Bigr)^{n-1}=\int_H \Bigl(\theta|_H+\ddc \varphi|_H\Bigr)^{n-1}.
    \]
    Therefore, using \cref{cor:decnetdS}, we conclude that $\varphi_j|_H\xrightarrow{d_S}\varphi|_H$.
\end{proof}

\begin{corollary}\label{cor:rest_gen_very_good}
Assume that $n\geq 2$.
    Let $\varphi\in \QPSH(X)$ be an $\mathcal{I}$-good potential. Then quasi-every $H\in \Lambda$ satisfies:
    \begin{enumerate}
        \item\label{itm:chapter-trace:L1071:1} $H$ is connected and smooth;
        \item\label{itm:chapter-trace:L1071:2} $\varphi|_H$ is $\mathcal{I}$-good;
        \item\label{itm:chapter-trace:L1071:3} $\nu(\varphi,H)=0$;
        \item\label{itm:chapter-trace:L1071:4} $\Tr_H \varphi\sim_P \varphi|_H$.
    \end{enumerate}
    Furthermore, if $\theta$ is a closed smooth real $(1,1)$-form on $X$ such that $\varphi\in \PSH(X,\theta)_{>0}$, then we can further guarantee, on the same co-pluripolar set of $H\in \Lambda$, that $\Tr_H(\varphi)$ has a representative $\Tr_H(\varphi)\in \PSH(H,\theta|_H)_{>0}$.
\end{corollary}
\begin{proof}
    Choose a closed smooth real $(1,1)$-form $\theta_0$ representing a big class such that $\varphi\in \PSH(X,\theta_0)$ and $\theta_0+\ddc \varphi$ is a Kähler current. Let $(\varphi_j)_j$ be a quasi-equisingular approximation of $\varphi$ in $\PSH(X,\theta_0)_{>0}$. Since $\varphi$ is $\mathcal I$-good, \cref{cor:quasi-equichar,thm:charIgoodasclosure,prop:ds0char} gives
    \[
        \varphi_j\xrightarrow{d_S}\varphi.
    \]
    Applying \cref{lma:Igoodrest}, we find a co-pluripolar set of $H\in\Lambda$ such that $H$ is connected and smooth, the restrictions $\varphi_j|_H$ and $\varphi|_H$ are not identically $-\infty$, and
    \[
        \varphi_j|_H\xrightarrow{d_S}\varphi|_H.
    \]
    For such an $H$, the potentials $\varphi_j|_H$ have analytic singularities, hence are $\mathcal I$-good by \cref{ex:analyIgood}. It follows from \cref{cor:Igoodclosed} that $\varphi|_H$ is $\mathcal I$-good.

    On the other hand, \cref{cor:ABTfortrace} gives, after intersecting with another co-pluripolar set of parameters $H$, that $\nu(\varphi,H)=0$ and
    \[
        \Tr_H(\varphi)\sim_{\mathcal I}\varphi|_H.
    \]
    The trace $\Tr_H(\varphi)$ is $\mathcal I$-good by \cref{prop:traceunique}. Since both $\Tr_H(\varphi)$ and $\varphi|_H$ are $\mathcal I$-good, their $\mathcal I$-equivalence identifies their $\mathcal I$-envelopes, which agree with their $P$-envelopes. Hence
    \[
        \Tr_H(\varphi)\sim_P\varphi|_H.
    \]

    Finally \cref{cor:tracegeneralwelldef} shows, after intersecting with one more co-pluripolar set of $H\in\Lambda$, that $\Tr_H(\varphi)$ has a representative in $\PSH(H,\theta|_H)_{>0}$.
\end{proof}

For later use, let us also prove a reverse Bertini theorem here.
\begin{lemma}[Reverse Bertini theorem]\label{lma:fibIeq}\index{theorem!reverse Bertini theorem}
Let $X$ be a complex manifold, $f\colon X\rightarrow \Delta^*$ be a projective surjective morphism to the punctured unit disk $\Delta^*$. Let $(L,h)$, $(L,h')$ be Hermitian pseudo-effective line bundles on $X$ with the same underlying line bundle. Assume that $S^1$ acts smoothly on $(X,L)$ by biholomorphic bundle automorphisms, that $f$ is equivariant for the rotation action on $\Delta^*$, and that $h$ and $h'$ are invariant.
Assume that for quasi-every $z\in \Delta^*$, $X_z$ is smooth and $h|_{X_z}\sim_{\mathcal{I}}h'|_{X_z}$. Then $h\sim_{\mathcal{I}}h'$.
\end{lemma}
\begin{proof}
We need to show that $\mathcal{I}(kh)=\mathcal{I}(kh')$ for every positive integer $k$. It suffices to prove the case $k=1$, since our assumptions are invariant after replacing $h$ and $h'$ by $h^k$ and $h'^k$. We will therefore prove $\mathcal{I}(h)=\mathcal{I}(h')$. First observe that it suffices to prove that
\begin{equation}\label{eq:fstarcoin}
f_*\Bigl(K_{X}\otimes L\otimes \mathcal{I}(h)\Bigr)=f_*\Bigl(K_{X}\otimes L\otimes \mathcal{I}(h')\Bigr)
\end{equation}
as subsheaves of $f_*(K_{X}\otimes L)$. In fact, suppose that \eqref{eq:fstarcoin} holds. Let $H$ be a $f$-ample invertible sheaf on $X$. Then \eqref{eq:fstarcoin} also holds with $L\otimes H^m$ in place of $L$. It follows from Grauert--Remmert's version of Serre vanishing theorem \cite[Theorem~2.1(A)]{BS76} that $\mathcal{I}(h)=\mathcal{I}(h')$.

It remains to prove \eqref{eq:fstarcoin}. Observe that both sides of \eqref{eq:fstarcoin} are locally free by \cite[Corollary~1.5]{Mat16}. We claim that it suffices to show that
\begin{equation}\label{eq:fstarcoin2}
f_*\Bigl(K_{X}\otimes L\otimes \mathcal{I}(h)\Bigr) \otimes \mathbb{C}(z)=f_*\Bigl(K_{X}\otimes L\otimes \mathcal{I}(h')\Bigr) \otimes \mathbb{C}(z)
\end{equation}
for almost every $z\in \Delta^*$ with respect to Lebesgue measure.
Indeed, the two direct images are coherent subsheaves of $f_*(K_X\otimes L)$. If their fibers agree almost everywhere, they agree outside a discrete subset of $\Delta^*$. Their locus of disagreement is $S^1$-invariant, and no nonempty discrete subset of $\Delta^*$ is invariant under rotations. Hence they agree everywhere.

For almost every $z\in\Delta^*$, the fiber $X_z$ is smooth and
\[
\mathcal I(h|_{X_z})=\Rest_{X_z}\mathcal I(h),\qquad
\mathcal I(h'|_{X_z})=\Rest_{X_z}\mathcal I(h').
\]
To see this, take a countable cover of the smooth locus of $f$ by relatively compact product charts. Applying Fubini's theorem to the weighted squared norms of local generators of $\mathcal I(h)$ gives
$\Rest_{X_z}\mathcal I(h)\subseteq\mathcal I(h|_{X_z})$ for almost every $z$. The reverse inclusion follows locally from \cref{thm:OT}. The same argument applies to $h'$, and the critical values of $f$ form a discrete subset of $\Delta^*$.

By generic flatness of $\mathcal O_X/\mathcal I(h)$ and $\mathcal O_X/\mathcal I(h')$, the natural map \eqref{eq:pullbacktoinverimage} for each of these two ideal sheaves is an isomorphism on general fibers. Cohomology and base change therefore gives, for almost every $z\in \Delta^*$,
\[
\begin{split}
f_*\Bigl(K_{X}\otimes L\otimes \mathcal{I}(h)\Bigr) \otimes \mathbb{C}(z)=\mathrm{H}^0\Bigl(X_z,K_{X}|_{X_z}\otimes L|_{X_z}\otimes \mathcal{I}(h|_{X_z})\Bigr),\\
f_*\Bigl(K_{X}\otimes L\otimes \mathcal{I}(h')\Bigr) \otimes \mathbb{C}(z)=\mathrm{H}^0\Bigl(X_z,K_{X}|_{X_z}\otimes L|_{X_z}\otimes \mathcal{I}(h'|_{X_z})\Bigr).
\end{split}
\]
By hypothesis, $h|_{X_z}\sim_{\mathcal{I}}h'|_{X_z}$ for quasi-every $z$. Since polar subsets of $\Delta^*$ have Lebesgue measure zero, \eqref{eq:fstarcoin2} holds for almost every $z\in\Delta^*$. The proof is complete.
\end{proof}

\motto{Comment se fait-il que M.~Gauss ait osé vous faire dire que la plupart de vos théorèmes lui étaient connus et qu'il en avait fait la découverte dès 1808. Cet excès d'impudence n'est pas croyable de la part d'un homme qui a assez de mérite personnel pour n'avoir pas besoin de s'approprier les découvertes des autres.

--- Adrien-Marie Legendre\footnote{Adrien-Marie Legendre (1752--1833) was a French mathematician known for his foundational contributions to number theory, statistics, and mathematical analysis. Apart from his mathematical contributions, he also helped formalize the metric system during the French Revolution.}, in a letter to Jacobi in 1827}
\chapter{Test curves}
\label{chap:testcurve}

A geodesic ray $(\varphi_t)_{t\in[0,\infty)}$ in $\PSH(X,\theta)$ records the asymptotic behavior of singularities along a path of potentials. In this chapter we encode the same information in the dual variable $\tau$. Instead of recording the potential at time $t$, the Legendre transform records, for each slope $\tau$, the model potential determined by the part of the ray growing at least linearly with slope $\tau$. This produces a concave decreasing family $\tau\mapsto\Gamma_\tau$ of model potentials; more precisely, a test curve $\Gamma$ is a concave decreasing family $(\Gamma_\tau)_{\tau<\Gamma_{\max}}$ of model potentials in $\PSH(X,\theta)$.

The central result is a correspondence initiated by Ross and Witt Nystr\"om\footnote{Witt and Nystr\"om are both family names of a single person. Some Swedes have double family names. It should not be spelled as Witt-Nystr\"om as some people do in the literature.}. It is the pluripotential-theoretic analogue of Legendre duality: a test curve gives a geodesic ray by taking a Legendre transform, and a suitable geodesic ray gives back a test curve by the inverse transform. Thus test curves provide an intrinsic language for the singularity data carried by geodesic rays.

We first develop test curves as self-contained objects. In \cref{sec:tc} we introduce their basic properties, including independence of the auxiliary form, behavior under bimeromorphic morphisms, closedness in the parameter $\tau$, log-concavity of the mass function and the $P$-equivalence relation on test curves. The Ross--Witt Nystr\"om correspondence is proved in \cref{sec:RWN}, where the Legendre transform is shown to give a bijection between test curves and geodesic rays.

The later sections adapt this picture to the singularity types needed elsewhere in the book. In \cref{sec:Imodeltc} we develop the $\mathcal I$-model version of the theory, including the $\mathcal I$-envelope of a test curve, the transition from $P$-model to $\mathcal I$-model test curves, and the corresponding notion of maximal geodesic ray. Finally, \cref{sec:operationtc} studies sums, translations, maxima, suprema, infima and rescalings of test curves. These operations are needed in the non-Archimedean theory of \cref{chap:NApotential}, where test curves provide the analytic input for the construction of non-Archimedean potentials.

Throughout the chapter we freely use the convex-analytic results of \cref{chap:convex}. Those results are stated for convex functions; when they are applied to concave functions here, they are applied to the negatives.

\section{The notion of test curves}\label{sec:tc}

Let $X$ be a connected compact K\"ahler manifold of dimension $n$ and $\theta$ be a smooth closed real $(1,1)$-form on $X$ representing a big cohomology class.

Recall that the notion of model potentials is defined in \cref{def:modelpot}.
\begin{definition}\label{def:testcur}
    A \emph{test curve}\index{test curve} $\Gamma$ in $\PSH(X,\theta)$ consists of a real number $\Gamma_{\max}$ together with a map $(-\infty,\Gamma_{\max})\rightarrow \PSH(X,\theta)$ denoted by $\tau\mapsto \Gamma_{\tau}$ satisfying the following conditions:
    \begin{enumerate}
        \item\label{itm:def:testcur:1} The map $\tau\mapsto \Gamma_{\tau}$ is concave and decreasing;
        \item\label{itm:def:testcur:2} $P_{\theta}[\Gamma_{\tau}]=\Gamma_{\tau}$;
        \item\label{itm:def:testcur:3} the potential
        \begin{equation}\label{eq:Gammaminf}
        \Gamma_{-\infty}\coloneqq \uscsup[\tau<\Gamma_{\max}]\Gamma_{\tau}
        \end{equation}
        satisfies
        \[
        \int_X \left(\theta+\ddc \Gamma_{-\infty}\right)^n>0.
        \]
    \end{enumerate}
    Let $\phi\in \PSH(X,\theta)_{>0}$ be a model potential.
    The set of test curves $\Gamma$ with $\Gamma_{-\infty}=\phi$ is denoted by $\TC(X,\theta;\phi)$\index{$\TC(X,\theta;\phi)$}.

    The union of all $\TC(X,\theta;\phi)$'s for various model potentials $\phi\in \PSH(X,\theta)_{>0}$ is denoted by $\TC(X,\theta)_{>0}$\index{$\TC(X,\theta)_{>0}$}.
\end{definition}
By \ref{itm:def:testcur:2}, $\sup_{X}\Gamma_{\tau}=0$ for each $\tau<\Gamma_{\max}$. So $\Gamma_{-\infty}\in \PSH(X,\theta)$ by \cref{prop:pshfunction_closedseq}.

Thanks to \cref{cor:incseqnppcont} and \ref{itm:def:testcur:3}, there is $C<\Gamma_{\max}$ so that whenever $\tau<C$, we have $\Gamma_{\tau}\in \PSH(X,\theta)_{>0}$. Then it follows that $\Gamma_{-\infty}$ is a model potential by \cref{prop:incnetmodel}.

\begin{remark}\label{rmk:extendtestcur}
    Sometimes it is convenient to extend $\Gamma_{\tau}$ to $\tau\geq \Gamma_{\max}$ as well. This can be done as follows: for $\tau>\Gamma_{\max}$, we set $\Gamma_{\tau}\equiv -\infty$. For $\tau=\Gamma_{\max}$, we set
    \[
        \Gamma_{\tau}\coloneqq \inf_{\tau'<\Gamma_{\max}}\Gamma_{\tau'}\in \PSH(X,\theta).
    \]
    Note that $\Gamma_{\tau}\in \PSH(X,\theta)$ by \cref{cor:decreasingnetL1conv}.
\end{remark}

\begin{figure}[ht]
\centering
\resizebox{0.95\textwidth}{!}{%
\begin{tikzpicture}[scale=1.0, line width=0.6pt, line cap=round, line join=round, every node/.style={font=\small}]
  \begin{scope}
    \draw[->, thick] (-4.45,0) -- (0.95,0) node[right] {$\tau$};
    \draw[->, thick] (0,-3.25) -- (0,0.62) node[above] {$\Gamma_\tau(x)$};
    \draw[thick] plot[smooth, samples=60, domain=-4.05:-0.90] (\x,{-0.22-0.19*(\x+4.05)*(\x+4.05)});
    \draw[dashed] (-0.90,0) -- (-0.90,-2.13);
    \draw[dashed] (-0.90,-2.13) -- (-0.90,-2.85);
    \draw[thick] (-0.90,-2.85) -- (0.34,-2.85);
    \fill (-0.90,-2.13) circle (1.5pt);
    \draw[fill=white] (-0.90,-2.85) circle (1.7pt);
    \node[above left] at (-0.90,0) {$\Gamma_{\max}$};
    \node[fill=white, inner sep=0.8pt] at (0.18,-0.20) {$0$};
    \node[right] at (0.50,-2.85) {$-\infty$};
    \node[above] at (-2.55,0.40) {(2)};
  \end{scope}

  \begin{scope}[xshift=6.3cm]
    \draw[->, thick] (-4.45,0) -- (0.95,0) node[right] {$\tau$};
    \draw[->, thick] (0,-3.25) -- (0,0.62) node[above] {$\Gamma_\tau(x)$};
    \draw[thick]
      (-4.05,-0.22)
      .. controls (-3.30,-0.28) and (-2.55,-0.55) .. (-2.10,-1.00)
      .. controls (-1.75,-1.40) and (-1.55,-2.25) .. (-1.50,-2.85);
    \draw[dashed] (-1.50,0) -- (-1.50,-3.0);
    \draw[dashed] (-0.75,0) -- (-0.75,-2.85);
    \draw[thick] (-1.50,-2.85) -- (0.34,-2.85);
    \node[above left] at (-1.50,0) {$\tau(x)$};
    \node[above left] at (-0.75,0) {$\Gamma_{\max}$};
    \node[fill=white, inner sep=0.8pt] at (0.18,-0.20) {$0$};
    \node[right] at (0.50,-2.85) {$-\infty$};
    \node[above] at (-2.55,0.40) {(3)};
  \end{scope}
\end{tikzpicture}
}
\caption{The last two pointwise possibilities for $\tau\mapsto\Gamma_\tau(x)$.}
\label{fig:testcurve_pointwise}
\end{figure}
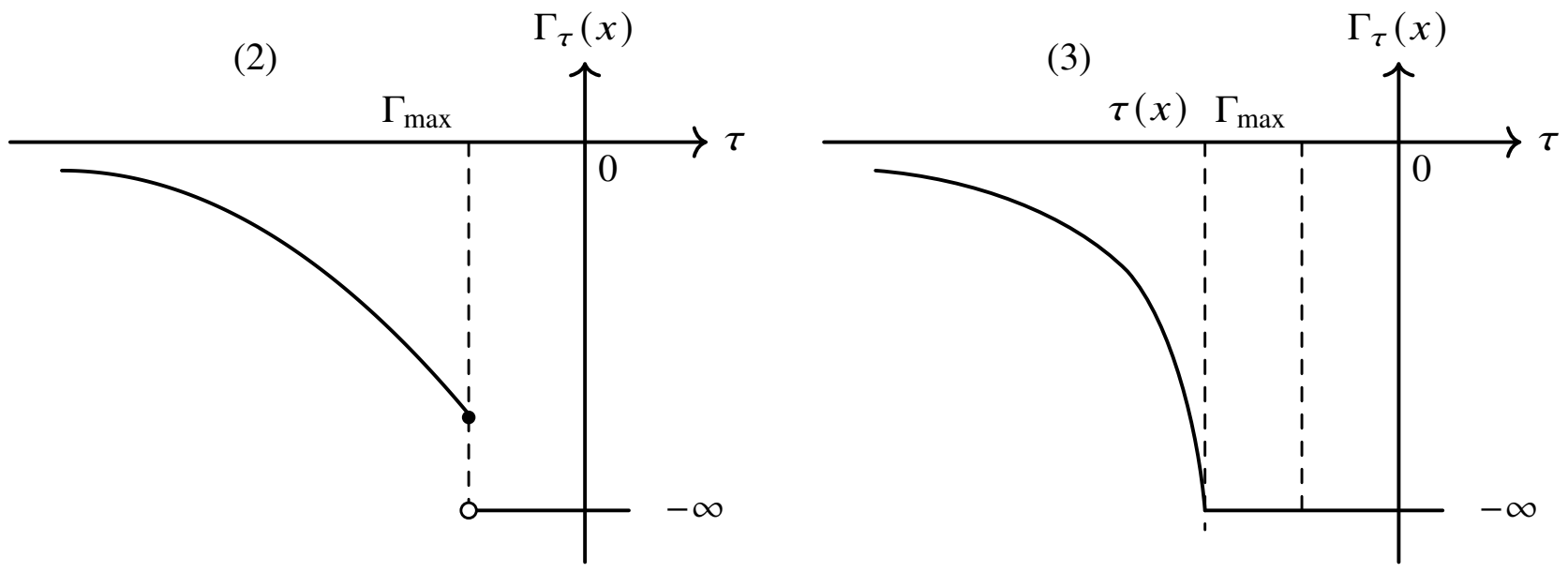

In general, for any $x\in X$, the map $\tau\mapsto \Gamma_{\tau}(x)$ falls into one of the following three categories:
\begin{enumerate}
    \item\label{itm:Gammataux_1} $\Gamma_{\tau}(x)\equiv -\infty$;
    \item\label{itm:Gammataux_2} $\Gamma_{\tau}(x)>-\infty$ for all $\tau\leq \Gamma_{\max}$;
    \item\label{itm:Gammataux_3} there is $\tau(x)\in (-\infty,\Gamma_{\max}]$ so that $\Gamma_{\tau}(x)>-\infty$ for all $\tau \in (-\infty,\tau(x))$, and $\Gamma_{\tau}(x)=-\infty$ for all $\tau\geq \tau(x)$.
\end{enumerate}
The last two possibilities are illustrated in \cref{fig:testcurve_pointwise}.

Note that \ref{itm:Gammataux_2} is the generic situation. It occurs whenever $x\notin \{\Gamma_{\Gamma_{\max}}=-\infty\}$. In particular, when this situation occurs, $(-\infty,\Gamma_{\max})\ni \tau\mapsto \Gamma_{\tau}(x)$ is finite concave, and hence continuous. As a consequence, for any $\tau \in (-\infty,\Gamma_{\max})$, we have
\begin{equation}\label{eq:Gammataucont}
    \Gamma_{\tau}=\inf_{\tau'<\tau}\Gamma_{\tau'},\quad \Gamma_{\tau}=\uscsup[\tau'\in (\tau,\Gamma_{\max})]\Gamma_{\tau'}
\end{equation}
almost everywhere and hence everywhere due to \cref{prop:pshfuncdetdense}.

\begin{figure}[ht]
\centering
\begin{tikzpicture}[scale=1.0, line width=0.6pt, every node/.style={font=\small}]
  \draw[->, thick] (-0.5,1.4) -- (8.0,1.4) node[right] {$\tau$};
  \draw[thick] plot[smooth, samples=80, domain=-0.5:6.5] (\x, {1 - 0.05*(\x+1)*(\x+1)});
  \draw[thick, dashed] (6.5,-1.8) -- (6.5,-3.2);
  \pgfmathsetmacro{\Yz}{1-0.05*4}      
  \pgfmathsetmacro{\Ymid}{1-0.05*16}   
  \pgfmathsetmacro{\Yp}{1-0.05*36}     
  \draw[red!75!black, dashed, thin] (1,1.4) -- (1,\Yz);
  \draw[red!75!black, dashed, thin] (3,1.4) -- (3,\Ymid);
  \draw[red!75!black, dashed, thin] (5,1.4) -- (5,\Yp);
  \draw[red!75!black, dashed, thin] (6.5,1.4) -- (6.5,-1.8);
  \node[red!75!black, above=2pt] at (1,1.4) {$\tau_0$};
  \node[red!75!black, above=2pt] at (3,1.4) {$\tau$};
  \node[red!75!black, above=2pt] at (5,1.4) {$\tau'$};
  \node[red!75!black, above=2pt] at (6.5,1.4) {$\tau_{\max}$};
  \draw[blue!70!black, thick] (1,\Yz) -- (5,\Yp);
  \fill[blue!70!black] (1,\Yz) circle (1.6pt);
  \fill[blue!70!black] (3,\Ymid) circle (1.6pt);
  \fill[blue!70!black] (5,\Yp) circle (1.6pt);
  \node[anchor=east] at (-0.5,0.95) {$\Gamma$};
  \node at (7.4,-2.6) {$-\infty$};
\end{tikzpicture}
\caption{The test curve $\Gamma$.}
\label{fig:tc}
\end{figure}
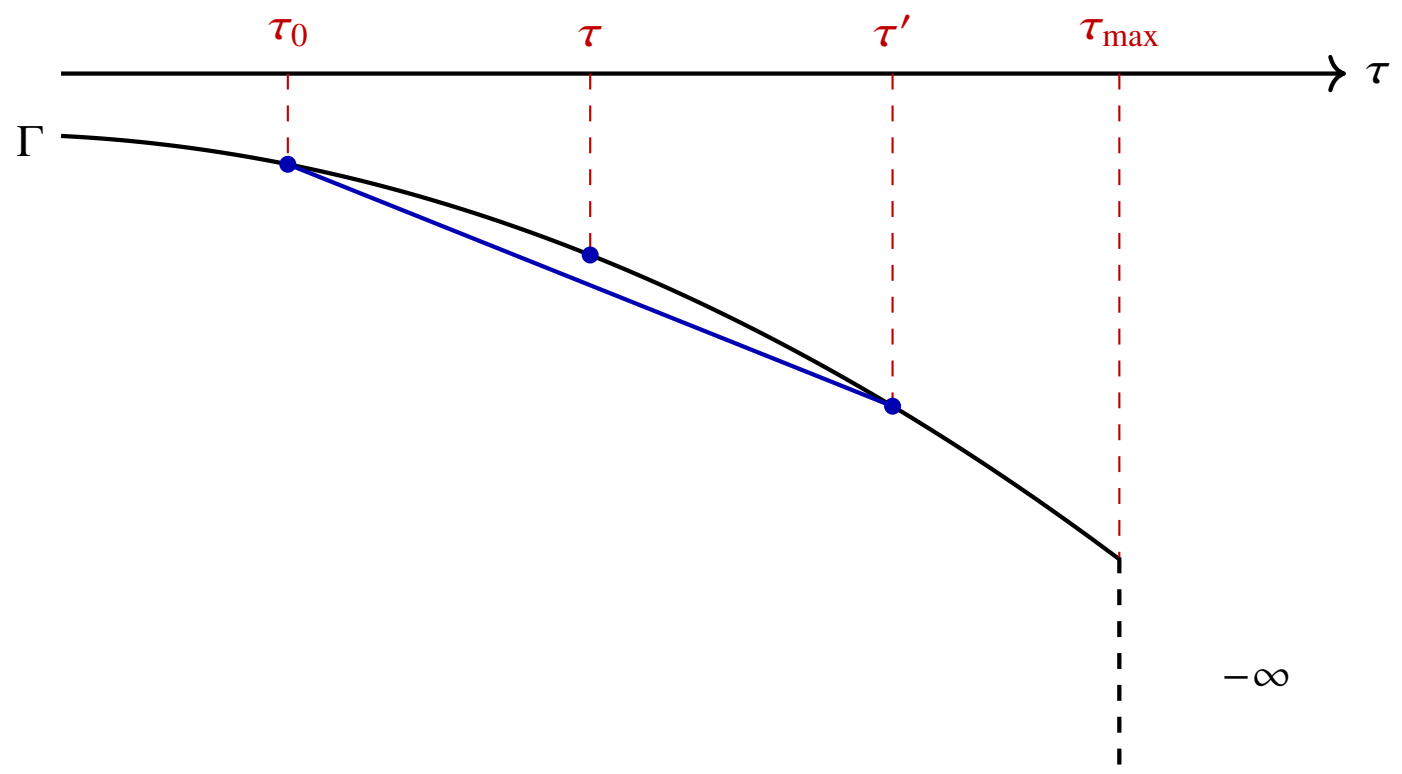

Recall that according to our general principle, we speak of model potentials only when a potential has positive mass. 
\begin{lemma}\label{lma:testcurvposmass}
Assume that $\Gamma\in \TC(X,\theta)_{>0}$. Then for each $\tau<\Gamma_{\max}$, we have
\begin{equation}\label{eq:dalethtauposmass}
\int_X \left(\theta+\ddc \Gamma_{\tau}\right)^n>0.
\end{equation}
In particular, $\Gamma_{\tau}$ is a model potential in $\PSH(X,\theta)_{>0}$.
\end{lemma}
The notations in the proof below are summarized in \cref{fig:tc}.

\begin{proof}

    Fix $\tau\in (-\infty,\Gamma_{\max})$.

    By assumption, $\Gamma_{-\infty}$ has positive mass. By \cref{cor:incseqnppcont}, we have
    \[
    \int_X \theta_{\Gamma_{-\infty}}^n=\lim_{\tau\to -\infty}\int_X \theta_{\Gamma_{\tau}}^n.
    \]
    In particular, for a sufficiently small $\tau_0<\tau$, we have
    \[
    \int_X \theta_{\Gamma_{\tau_0}}^n>0.
    \]
    Now take $\tau'\in (\tau, \Gamma_{\max})$ and $t\in (0,1)$ so that
    \[
    \tau=(1-t)\tau'+t\tau_0.
    \]
    From the concavity of $\Gamma$, we find that
    \[
    \Gamma_{\tau}\geq (1-t)\Gamma_{\tau'}+t\Gamma_{\tau_0}.
    \]
    By the monotonicity theorem \cref{thm:mono} and the multi-linearity of the non-pluripolar products \itemref{itm:prop:npp1:6},
    \[
    \int_X \theta_{\Gamma_{\tau}}^n\geq \int_X \theta_{(1-t)\Gamma_{\tau'}+t\Gamma_{\tau_0}}^n\geq t^n \int_X \theta_{\Gamma_{\tau_0}}^n>0
    \]
    and \eqref{eq:dalethtauposmass} follows.
\end{proof}

\begin{proposition}\label{prop:testcurvmasslogconc}
    Let $\Gamma\in \TC(X,\theta)_{>0}$. Then the map
    \[
    [-\infty,\Gamma_{\max})\rightarrow \mathbb{R},\quad \tau\mapsto \log \int_X \theta_{\Gamma_{\tau}}^n
    \]
    is continuous and concave on $(-\infty,\Gamma_{\max})$.
\end{proposition}
\begin{proof}
    The concavity of this function follows from \cref{thm:logconc} and \cref{thm:mono}. The continuity at $-\infty$ is a consequence of \cref{cor:incseqnppcont}.
    The continuity on $(-\infty,\Gamma_{\max})$ follows from the fact that finite concave functions on open intervals are continuous.
\end{proof}

\begin{definition}\label{def:relattestcurv}
    Let $\phi\in \PSH(X,\theta)_{>0}$ be a model potential.

    A test curve $\Gamma\in \TC(X,\theta;\phi)$ is said to be \emph{bounded}\index{test curve!bounded test curve} if for $\tau$ small enough, $\Gamma_{\tau}=\phi$. The subset of bounded test curves in $\TC(X,\theta;\phi)$ is denoted by $\TC^{\infty}(X,\theta;\phi)$\index{$\TC^{\infty}(X,\theta;\phi)$}.
    In this case, we write
    \begin{equation}\label{eq:Gammamin}
    \Gamma_{\min}\coloneqq \sup\{\tau\in \mathbb{R}:\Gamma_{\tau}=\phi\}.
    \end{equation}

    A test curve $\Gamma\in \TC(X,\theta;\phi)$ is said to have \emph{finite energy}\index{test curve!test curve with finite energy} if
    \begin{equation}\label{eq:tcfiniteenergy}
    \mathbf{E}^{\phi}(\Gamma)\coloneqq \Gamma_{\max} \int_X \theta_{\phi}^n+\int_{-\infty}^{\Gamma_{\max}} \left(\int_X\theta_{\Gamma_{\tau}}^n-\int_X \theta_{\phi}^n\right)\,\mathrm{d}\tau>-\infty.
    \end{equation}
    \index{$\mathbf{E}^{\phi}$}
    When $\phi=V_{\theta}$, we write $\mathbf{E}$ instead of $\mathbf{E}^{\phi}$.

    The subset of test curves with finite energy in $\TC(X,\theta;\phi)$ is denoted by $\TC^1(X,\theta;\phi)$\index{$\TC^1(X,\theta;\phi)$}.
\end{definition}

\begin{example}\label{ex:cantest}
    Given $\varphi\in \PSH(X,\theta)$, there is a canonically associated test curve $\Gamma^{\varphi}\in \TC^{\infty}(X,\theta;V_{\theta})$: Set $\Gamma^{\varphi}_{\max}=0$ and
    \[
    \Gamma^{\varphi}_{\tau}=
    \begin{cases}
        V_{\theta},&\textup{if }\tau\leq -1;\\
        P_{\theta}\left[ (1+\tau)\varphi-\tau V_{\theta} \right],& \textup{if }-1<\tau<0.
    \end{cases}
    \]
    Note that $\Gamma^{\varphi}$ is indeed a test curve, as follows from \cref{prop:Pconc}. Indeed, the concavity and monotonicity in $\tau$ follow from \cref{prop:Pconc}, and $\Gamma^\varphi_{-\infty}=V_\theta$ has positive mass.
\end{example}

We first observe that the notion of test curves does not really depend on the choice of $\theta$ within its cohomology class.
\begin{proposition}\label{prop:testcurveindeptheta}
    Let $\theta'$ be another smooth closed real $(1,1)$-form on $X$ representing the same cohomology class as $\theta$. Let $\phi\in \PSH(X,\theta)_{>0}$ be a model potential. Let $\phi'\in \PSH(X,\theta')_{>0}$ be the unique model potential satisfying $\phi\sim\phi'$.

    Then there is a canonical bijection
    \[
    \TC(X,\theta;\phi)\xrightarrow{\sim} \TC(X,\theta';\phi').
    \]
    This bijection induces the following bijections:
    \[
    \TC^1(X,\theta;\phi)\xrightarrow{\sim}\TC^1(X,\theta';\phi'),\quad \TC^{\infty}(X,\theta;\phi)\xrightarrow{\sim}\TC^{\infty}(X,\theta';\phi').
    \]
    These bijections satisfy the obvious cocycle conditions.
\end{proposition}
\begin{proof}
    Choose $g\in C^{\infty}(X)$ such that $\theta'=\theta+\ddc g$. Given any $\Gamma\in \TC(X,\theta;\phi)$, we observe that $\Gamma'\colon (-\infty,\Gamma_{\max})\rightarrow \PSH(X,\theta')$ defined as
    \[
    \tau\mapsto P_{\theta'}[\Gamma_{\tau}-g]
    \]
    lies in $\TC(X,\theta';\phi')$. In fact, the concavity and monotonicity of $\Gamma'$ follow from \cref{prop:Pconc}. Moreover, applying \cref{prop:incnetmodel} to the increasing net obtained as $\tau\to -\infty$ gives
    \[
        \Gamma'_{-\infty}=P_{\theta'}[\Gamma_{-\infty}-g]=P_{\theta'}[\phi-g]=\phi'.
    \]
    Thus $\Gamma'\in \TC(X,\theta';\phi')$.
    
    Moreover, the choice of $g$ is irrelevant since for any other choice of $g$, say $g'$, we have
    \[
    \Gamma_{\tau}-g\sim \Gamma_{\tau}-g'
    \]
    for all $\tau<\Gamma_{\max}$. The inverse construction is obtained by replacing $g$ with $-g$. The boundedness condition is clearly preserved. Finally, by \cref{prop:Ppresmass}, for every $\tau<\Gamma_{\max}$,
    \[
        \int_X(\theta'+\ddc \Gamma'_\tau)^n=\int_X(\theta+\ddc \Gamma_\tau)^n,
    \]
    and the same holds for $\phi$ and $\phi'$. Hence the finite energy condition is preserved. The cocycle condition is immediate from the uniqueness of the model representative in a fixed singularity type.
\end{proof}

\begin{proposition}\label{prop:ETCbimero}
    Let $\pi\colon Y\rightarrow X$ be a proper bimeromorphic morphism from a K\"ahler manifold $Y$. Then the pointwise pull-back induces a bijection
    \[
    \pi^*\colon \TC(X,\theta;\phi)\xlongrightarrow{\sim} \TC(Y,\pi^*\theta;\pi^*\phi).
    \]
\end{proposition}
\begin{proof}
     The bijectivity on the underlying potentials follows from \cref{prop:PSHpullbij}. The compatibility with model potentials follows from \cref{prop:Penvbimero}. These two facts also preserve concavity, monotonicity and the value of $\Gamma_{-\infty}$, hence give the asserted bijection of test curves.
\end{proof}

Next we verify the closedness of a test curve as a family of concave functions, so that no pathology arises in the Legendre transforms considered below. The notion of closedness is recalled in \cref{def:convexclosure}.
\begin{proposition}\label{prop:Gammaclosed}
    Let $\Gamma$ be a test curve in $\PSH(X,\theta)$. For each $x\in X$, the map $\mathbb{R} \ni\tau\mapsto \Gamma_{\tau}(x)$ is a closed concave function. Moreover, the map is proper as long as $\Gamma_{\Gamma_{\max}}(x)\neq -\infty$.
\end{proposition}

\begin{proof}
    We prove closedness.
    Fix $x\in X$. Assume that $\Gamma_{\tau}(x)\neq -\infty$ for some $\tau\in \mathbb{R}$. We only need to argue the upper semicontinuity of $\tau\mapsto \Gamma_{\tau}(x)$. The upper semicontinuity is clear at $\tau\geq \Gamma_{\max}$, so we are reduced to prove the following:
    \begin{equation}\label{eq:Gammatautemp1}
        \Gamma_{\tau}=\inf_{\tau'<\tau}\Gamma_{\tau'}
    \end{equation}
    for any $\tau<\Gamma_{\max}$. Take $\tau''\in (\tau,\Gamma_{\max})$. Outside the polar locus of $\Gamma_{\tau''}$, we know that \eqref{eq:Gammatautemp1} holds by continuity of real-valued concave functions. So \eqref{eq:Gammatautemp1} holds everywhere by \cref{prop:pshfuncdetdense}.

    The final assertion follows from properness of a finite closed concave function on its effective domain.
\end{proof}

\begin{definition}\label{def:Ptestcurve}
    Let $\Gamma\in \TC(X,\theta)_{>0}$ and $\omega$ be a smooth closed real positive $(1,1)$-form. Then we define
    $P_{\theta+\omega}[\Gamma]\in \TC(X,\theta+\omega)_{>0}$ as follows:
    \index{test curve!P-envelope@$P$-envelope of a test curve}
    \index{$P_{\theta+\omega}[\Gamma]$}
    \begin{enumerate}
        \item\label{itm:def:Ptestcurve:1} Define
        \[
        P_{\theta+\omega}[\Gamma]_{\max}=\Gamma_{\max};
        \]
        \item\label{itm:def:Ptestcurve:2} for each $\tau<\Gamma_{\max}$, define
        \[
        P_{\theta+\omega}[\Gamma]_{\tau}=P_{\theta+\omega}[\Gamma_{\tau}].
        \]
    \end{enumerate}
\end{definition}
It follows from \cref{prop:Pconc} that the family $P_{\theta+\omega}[\Gamma]$ is concave and decreasing. Moreover, since $\omega$ is positive,
\[
\int_X(\theta+\omega+\ddc\Gamma_\tau)^n \geq \int_X(\theta+\ddc\Gamma_\tau)^n>0
\]
for every $\tau<\Gamma_{\max}$, and hence $P_{\theta+\omega}[\Gamma]_\tau$ is a model potential by \cref{cor:Psendspotentialtomodel}. Thus $P_{\theta+\omega}[\Gamma]\in\TC(X,\theta+\omega)_{>0}$.

\begin{proposition}\label{prop:PpreGammainf}
    Let $\Gamma\in \TC(X,\theta)_{>0}$ and $\omega$ be a closed real smooth positive $(1,1)$-form on $X$. Then
    \[
    P_{\theta+\omega}[\Gamma]_{-\infty}=P_{\theta+\omega}[\Gamma_{-\infty}].
    \]
\end{proposition}
\begin{proof}
     This follows from \cref{prop:incnetmodel}.
\end{proof}

\section{Ross--Witt Nystr\"om correspondence}\label{sec:RWN}
Let $X$ be a connected compact K\"ahler manifold of dimension $n$ and $\theta$ be a smooth closed real $(1,1)$-form on $X$ representing a big cohomology class.
Fix a model potential $\phi\in \PSH(X,\theta)_{>0}$.

\cref{prop:Gammaclosed} allows us to talk about the Legendre transforms of test curves in the expected way.

The general definition of the Legendre transform \cref{def:Legendregeneral} can be translated as follows:
\begin{definition}\label{def:Legtrans}
    Let $\Gamma\in \TC(X,\theta;\phi)$. We define its \emph{Legendre transform}\index{Legendre transform} as the map $\Gamma^*\colon (0,\infty)\rightarrow \PSH(X,\theta)$ given by
    \begin{equation}\label{eq:testcurveLegtran}
        \Gamma^*_t=\sup_{\tau\in \mathbb{R}} \left(t\tau+ \Gamma_{\tau}\right)\footnotemark.
    \end{equation}
    \footnotetext{There is no usc regularization in the preceding formula. This is not a typo.}
\end{definition}
Thanks to \cref{rmk:extendtestcur}, \eqref{eq:testcurveLegtran} can be equivalently written as
\[
\Gamma^*_t=\sup_{\tau<\Gamma_{\max}} \left(t\tau+ \Gamma_{\tau}\right)=\sup_{\tau\leq \Gamma_{\max}} \left(t\tau+ \Gamma_{\tau}\right).
\]

It is sometimes handy to \emph{define}
\begin{equation}\label{eq:Gamma0star}
\Gamma^*_0\coloneqq \phi
\end{equation}
at $t=0$. But it is important to remember that, with this convention, \eqref{eq:testcurveLegtran} is not true at $t=0$ in general.

\begin{remark}\label{rmk:negativeray}
Here we do not discuss the case $t<0$ because its behavior is straightforward: take $x\in X$; if $\Gamma_{\tau}(x)=-\infty$ for all $\tau<\Gamma_{\max}$, then $\Gamma_t^*(x)=-\infty$; otherwise, $\Gamma_t^*(x)=\infty$.
\end{remark}

The information about $t>0$ suffices to characterize $\Gamma$.
\begin{proposition}\label{prop:Leginvtc}
    Let $\Gamma\in \TC(X,\theta;\phi)$. Then
    \begin{equation}\label{eq:Leginvtc}
    \Gamma_{\tau}=\inf_{t>0} \left(\Gamma^*_t-t\tau\right)
    \end{equation}
    for all $\tau\in \mathbb{R}$.
\end{proposition}
Due to our convention \eqref{eq:Gamma0star}, in \eqref{eq:Leginvtc} we may equivalently take $t\geq 0$.
\begin{proof}
    Fix $x\in X$. We want to establish \eqref{eq:Leginvtc} at $x$.
    We distinguish two cases. First suppose that $\Gamma_{\tau}(x)=-\infty$ for all $\tau<\Gamma_{\max}$ and hence all $\tau\in \mathbb{R}$. In this case, we have $\Gamma_t^*(x)=-\infty$ for all $t>0$. Therefore, \eqref{eq:Leginvtc} follows trivially.

    Otherwise, by \cref{rmk:negativeray}, we know that $\Gamma^*_t(x)=\infty$ for all $t<0$. The relative interior of the domain of $t\mapsto \Gamma^*_t(x)$ is contained in $(0,\infty)$.
    Therefore, \eqref{eq:Leginvtc} follows from \cref{thm:Legendredual}, \cref{prop:Gammaclosed}.
\end{proof}

In \cref{def:Legtrans}, we made the nontrivial claim that $\Gamma^*_t\in \PSH(X,\theta)$ for all $t>0$. Let us prove this.
\begin{lemma}\label{lma:testcurvelegusc}
    Let $\Gamma\in \TC(X,\theta;\phi)$. Then $\Gamma^*_t\in \PSH(X,\theta)$ for all $t>0$. In fact, $\Gamma^*$ is upper semicontinuous as a function of $X\times (0,\infty)$.
\end{lemma}
\begin{proof}
    We first observe that for each $x\in X$, we have
    \[
        \Gamma_t^*(x)\leq t\Gamma_{\max}<\infty.
    \]

    Let $R=\{a+\mathrm{i}b\in \mathbb{C}:a>0,b\in \mathbb{R}\}$. We consider
    \[
        F\colon X\times R\rightarrow [-\infty,\infty),\quad (x,a+\mathrm{i}b)\mapsto \Gamma_a^*(x).
    \]
    Let $\pi\colon X\times R\rightarrow X$ be the natural projection.
    Observe that the upper semicontinuous regularization $G$ of $F$ is $\pi^*\theta$-psh by \cref{prop:pshfunction_closedseq}.
    It suffices to show that $F=G$.
    We let
    \[
        E\coloneqq \left\{ (x,z)\in X\times R: F(x,z)<G(x,z)  \right\}.
    \]
    We want to argue that $E=\varnothing$. Since $F$ is invariant under translations in the imaginary direction of $R$, the same is true for its upper semicontinuous regularization $G$. So $E$ can be written as $B\times \mathrm{i}\mathbb{R}$ for some set $B\subseteq X\times (0,\infty)$. Since $E$ is a pluripolar set by \cref{prop:supsupstardiff}, it has zero Lebesgue measure. Hence, $B$ has zero Lebesgue measure. For each $x\in X$, write
    \[
        B_x=\left\{t\in (0,\infty): (x,t)\in B \right\}.
    \]
    By Fubini's theorem, $B_x$ has vanishing $1$-dimensional Lebesgue measure for all $x\in X\setminus Z$, where $Z\subseteq X$ is a subset of measure $0$. We may assume that $Z\supseteq \{\Gamma_{\Gamma_{\max}}=-\infty\}$ so that for $x\in X\setminus Z$, $\Gamma^*_t(x)\neq -\infty$ for all $t>0$.

    For any $x\in X\setminus Z$, both $t\mapsto F(x,t)$ and $G(x,t)$ are convex functions with values in $\mathbb{R}$ on $(0,\infty)$. They agree almost everywhere, hence everywhere by their continuity. It follows that for $x\in X\setminus Z$, we have $B_x=\varnothing$.

    By \cref{prop:Leginvtc}, for any $x\in X$, we have
    \[
        \Gamma_{\tau}(x)=\inf_{t>0} (F(x,t)-t\tau),\quad \tau<\Gamma_{\max}.
    \]
    On the other hand, let
    \begin{equation}\label{eq:chiGlegtemp1}
        \chi_{\tau}(x)=\inf_{t>0} (G(x,t)-t\tau),\quad \tau<  \Gamma_{\max},\quad x\in X.
    \end{equation}
    By Kiselman's principle \cref{prop:Kis}, $\chi_{\tau}\in \PSH(X,\theta)$. But on $X\setminus Z$, we already know that $\Gamma_{\tau}=\chi_{\tau}$ for all $\tau<\Gamma_{\max}$. By \cref{prop:pshfuncdetdense},
    \[
    \Gamma_{\tau}=\chi_{\tau},\quad \tau<\Gamma_{\max}.
    \]
    Now we conclude that $F(x,t)=G(x,t)$ by \cref{cor:Leginvpartdef}.
\end{proof}

\begin{corollary}\label{cor:GammastarP}
    Let $\Gamma \in \TC(X,\theta;\phi)$. Then $\Gamma^*_t\in \mathcal{E}(X,\theta;\phi)$ for all $t>0$.
\end{corollary}
\begin{proof}
Fix $t>0$.
    We already know that $\Gamma^*_t\in \PSH(X,\theta)$ by \cref{lma:testcurvelegusc}. It suffices to show that
    \[
    \Gamma^*_t\sim_P \phi.
    \]
    From \eqref{eq:testcurveLegtran} and \cref{prop:Ppartialsup}, we know that
    \[
    \Gamma^*_t\sim_P \uscsup[\tau<\Gamma_{\max}] \Gamma_{\tau}=\phi.
    \]
    This completes the proof.
\end{proof}

\begin{lemma}\label{lma:suplegenlinear}
    Let $\Gamma\in \TC(X,\theta;\phi)$. Then
    \[
        \sup_X \Gamma_t^*=t\Gamma_{\max}
    \]
    for all $t>0$.

    In particular, $t\mapsto \Gamma^*_t-t\Gamma_{\max}$ is a decreasing function in $t>0$.
\end{lemma}
\begin{proof}
    Since $\Gamma_{\tau}\leq 0$ for all $\tau<\Gamma_{\max}$, we have
    \[
        \sup_X \Gamma_t^*\leq t\Gamma_{\max}
    \]
    for all $t>0$.
    Conversely, for each $\tau<\Gamma_{\max}$,
    \[
        \sup_X\Gamma_t^*\geq \sup_X(t\tau+\Gamma_{\tau})=t\tau.
    \]
    Letting $\tau\to\Gamma_{\max}$ gives the reverse inequality.
    The final assertion follows directly from the formula
    \[
        \Gamma_t^*-t\Gamma_{\max}=\sup_{\tau<\Gamma_{\max}}\Bigl(\Gamma_\tau-t(\Gamma_{\max}-\tau)\Bigr).
    \]
    This completes the proof.
\end{proof}

\begin{lemma}\label{lma:LegsendsTCtoR}
    Given $\Gamma\in \TC(X,\theta;\phi)$, we have $\Gamma^*\in \mathcal{R}(X,\theta;\phi)$.
\end{lemma}
See \cref{def:geodraydef} for the notation $\mathcal{R}(X,\theta;\phi)$.
\begin{proof}
    It follows from \cref{lma:testcurvelegusc}, \eqref{eq:testcurveLegtran} and \cref{prop:pshfunction_closedseq} that the restriction of $\Gamma^*$ to each compact subinterval of $(0,\infty)$ is a subgeodesic. By \cref{cor:GammastarP}, for any $t>0$, $\Gamma^*_t\in \mathcal{E}(X,\theta;\phi)$.

    First observe that as $t\to 0+$, we have
    \begin{equation}\label{eq:GammatophiL1temp}
        \Gamma_t^*\xrightarrow{L^1}\phi.
    \end{equation}
    By \cref{lma:suplegenlinear} and \cref{prop:L1compa}, it suffices to show each $L^1$-cluster point $\psi\in \PSH(X,\theta)$ of $\Gamma_t^*$ as $t\to 0$ is equal to $\phi$.

    To see this, first observe that by \eqref{eq:testcurveLegtran}, for any fixed $t>0$,
    \[
        \Gamma_t^*\leq t\Gamma_{\max}+\phi.
    \]
    Therefore, $\psi\leq \phi$. On the other hand, for any fixed $\tau<\Gamma_{\max}$, by \eqref{eq:testcurveLegtran}, we have
    \[
        \Gamma^*_t\geq \Gamma_{\tau}+t\tau
    \]
    for any $t>0$. So $\psi\geq \Gamma_{\tau}$ almost everywhere and hence everywhere by \cref{prop:pshfuncdetdense}. It follows that $\psi\geq \phi$. Therefore, $\psi=\phi$.

     Assume that $\Gamma^*$ is not a geodesic ray. Then we can find $0\leq a<b$ such that $(\Gamma^*_t)_{t\in (a,b)}$ differs from the geodesic $(\eta_t)_{t\in (a,b)}$ from $\Gamma^*_a$ to $\Gamma^*_b$. The existence of $(\eta_t)_t$ is guaranteed by \cref{prop:perronenvissubgeod}.
     We consider the subgeodesic $(\ell_t)_{t>0}$ given by $\ell_t=\eta_t$ for $t\in (a,b)$ and $\ell_t=\Gamma^*_t$ otherwise. Note that $\ell$ is a subgeodesic due to the gluing lemma \cref{lma:subgeodgluing}.

     Consider the Legendre transform
    \[
        \Gamma'_{\tau}=\inf_{t>0} (\ell_t-t\tau),\quad \tau\in \mathbb{R}.
    \]
    Then $\Gamma'_{\tau}\geq \Gamma_{\tau}$ and $\Gamma'_{\tau}\in \PSH(X,\theta)\cup\{-\infty\}$ by \cref{prop:Kis} for all $\tau\in \mathbb{R}$.

    We claim that
    \begin{equation}\label{eq:GammapGammcomptemp1}
        \Gamma'_{\tau}\leq \Gamma_{\tau}+(b-a)(\Gamma_{\max}-\tau),\quad \tau\in \mathbb{R}.
    \end{equation}
    Observe that $\Gamma'_{\tau}\equiv -\infty$ when $\tau>\Gamma_{\max}$ by \cref{lma:suplegenlinear}. Thus it suffices to consider $\tau\leq \Gamma_{\max}$. Put
    \[
        C_{\tau}\coloneqq (b-a)(\Gamma_{\max}-\tau)\geq 0.
    \]
    For any $s\in [a,b]$, the monotonicity of $t\mapsto \Gamma_t^*-t\Gamma_{\max}$ from \cref{lma:suplegenlinear} gives
    \[
        \Gamma_b^*-b\tau \leq \Gamma_s^*-s\tau+C_{\tau}.
    \]
    Since $\ell_b=\Gamma_b^*$, we get
    \[
        \Gamma'_{\tau}\leq \Gamma_b^*-b\tau\leq \Gamma_s^*-s\tau+C_{\tau},\quad s\in [a,b].
    \]
    For $s\notin [a,b]$, we have $\ell_s=\Gamma_s^*$, hence
    \[
        \Gamma'_{\tau}\leq \Gamma_s^*-s\tau\leq\Gamma_s^*-s\tau+C_{\tau}.
    \]
    Taking the infimum over all $s>0$ gives
    \[
        \Gamma'_{\tau}\leq \Gamma_{\tau}+C_{\tau}.
    \]
    Therefore, \eqref{eq:GammapGammcomptemp1} follows.
    In particular, for any $\tau<\Gamma_{\max}$, we have $\Gamma'_{\tau}\sim \Gamma_{\tau}$.
    Moreover, $\sup_X\ell_t\leq t\Gamma_{\max}$ for all $t>0$: this is clear outside $(a,b)$ from \cref{lma:suplegenlinear}, and on $(a,b)$ it follows from the convexity of subgeodesics applied to the geodesic $(\eta_t)_{t\in(a,b)}$. Hence, for $\tau<\Gamma_{\max}$,
    \[
        \sup_X(\ell_t-t\tau)\leq t(\Gamma_{\max}-\tau).
    \]
    Letting $t\to 0+$ gives $\Gamma'_{\tau}\leq 0$. It follows from the fact that $\Gamma_{\tau}$ is a model potential that $\Gamma_{\tau}=\Gamma'_{\tau}$ for all $\tau<\Gamma_{\max}$. Therefore, by \cref{thm:Legendredual}, we have $\Gamma^*_t=\ell_t$ for all $t>0$, which is a contradiction.
\end{proof}
Given $\ell\in \mathcal{R}(X,\theta;\phi)$, define its Legendre transform
    \begin{equation}\label{eq:invLeg}
        \ell^*_{\tau}\coloneqq \inf_{t>0} (\ell_t-t\tau),\quad \tau \in \mathbb{R}.
    \end{equation}

\begin{lemma}\label{lma:LegendsRtoTC}
    Given $\ell\in \mathcal{R}(X,\theta;\phi)$, then $\ell^*=(\ell^*_\tau)_{\tau<\sup_X \ell_1}\in \TC(X,\theta;\phi)$.
\end{lemma}
\begin{proof}
    Note that it follows from \cref{prop:Kis} that $\ell^*_{\tau}\in \PSH(X,\theta)\cup \{-\infty\}$ for all $\tau\in \mathbb{R}$. It is clear that $\mathbb{R}\ni \tau\mapsto \ell^*_{\tau}$ is a decreasing and concave function.

    By \cref{prop:raysuplinear},
    \[
    \sup_X\ell_t= t\sup_X \ell_1 \quad \forall t\geq 0.
    \]
    Set $C=\sup_X\ell_1$. By \cref{prop:raysuplinear}, $\sup_X(\ell_t-Ct)=0$ for all $t>0$. For each $x$ outside a pluripolar set, the function $t\mapsto \ell_t(x)-Ct$ is convex on $(0,\infty)$ and bounded above by $0$, hence it is decreasing. By \cref{prop:pshfuncdetdense}, the same monotonicity holds as an inequality of $\theta$-psh functions. Thus $(\ell_t-Ct)_{t>0}$ is a decreasing net in $\PSH(X,\theta)$.

    Its infimum is not identically $-\infty$, since the family is normalized by $\sup_X(\ell_t-Ct)=0$ and is compact in $L^1$ by \cref{prop:L1compa}. Therefore
    \[
        \ell^*_{C}=\inf_{t>0}(\ell_t-Ct)\in \PSH(X,\theta)
    \]
    by \cref{prop:pshfunction_closedseq}.

    On the other hand, for $\tau>\sup_X \ell_1$, the same argument shows that
    \[
    \ell^*_{\tau}\equiv -\infty.
    \]
    Therefore, $\ell^*_{\tau}\in \PSH(X,\theta)$ if and only if $\tau\leq \ell^*_{\max}\coloneqq \sup_X \ell_1$.

    We claim that $(\ell^*_{\tau})_{\tau<\ell^*_{\max}}$ is a test curve. We first observe that for $\tau<\ell^*_{\max}$, we have
    \[
    \ell^*_{\tau}\leq \ell_1-\tau\sim_P \phi.
    \]
    Therefore,
    \begin{equation}\label{eq:elltaupreceqpphi}
        \ell^*_{\tau}\preceq_P \phi,\quad \forall \tau<\ell^*_{\max}.
    \end{equation}
    Also observe that for any $\tau\leq \ell^*_{\max}$ and any $t>0$, we have
    \[
    \sup_X \ell^*_{\tau}\leq \sup_X \ell_t-t\tau=\ell^*_{\max}t-t\tau.
    \]
    Letting $t\to 0+$, we find that for any $\tau\leq \ell^*_{\max}$, we have
    \begin{equation}\label{eq:supXelltaustarnegtemp1}
    \sup_X \ell^*_{\tau}\leq 0.
    \end{equation}

    Fix $\tau<\ell^*_{\max}$. We want to argue that
    \begin{equation}\label{eq:Pellstartemp1}
        P_{\theta}\left[\ell^*_{\tau}\right]=\ell^*_{\tau}.
    \end{equation}
    First we claim that for any $C>0$, we have
    \begin{equation}\label{eq:ellstartaupCtemp1}
    (\ell^*_{\tau}+C)\land \phi=(\ell^*_{\tau}+C)\land V_{\theta}.
    \end{equation}
    The $\leq$ direction is trivial. We argue the reverse inequality, which reduces to
    \[
     \phi\geq (\ell^*_{\tau}+C)\land V_{\theta}.
    \]
    Since $\phi$ is model and $(\ell^*_{\tau}+C)\land V_{\theta}\leq 0$, it suffices to show that
    \[
    \phi\succeq_P (\ell^*_{\tau}+C)\land V_{\theta},
    \]
    which follows from \eqref{eq:elltaupreceqpphi}. Therefore, \eqref{eq:ellstartaupCtemp1} is established. Thanks to \eqref{eq:supXelltaustarnegtemp1}, we have the obvious inequality
    \[
    (\ell^*_{\tau}+C)\land V_{\theta}\geq \ell^*_{\tau}
    \]
    for any $C>0$. Therefore, in order to prove \eqref{eq:Pellstartemp1}, it remains to argue that for any $C>0$,
    \begin{equation}\label{eq:ellstarleqetemp1}
        (\ell^*_{\tau}+C)\land \phi\leq  \ell^*_{\tau}.
    \end{equation}
    For this purpose, let us consider the following geodesics: For any $M>0$ and $t\in [0,1]$, let
    \[
        \ell^{1,M}_t= \ell_{tM}-tM\tau,\quad  \ell^{2,M}_t= (\ell^*_{\tau}+C)\land \phi-Ct.
    \]
    It is clear that at $t=0,1$, we have $\ell^{2,M}_t\leq \ell^{1,M}_t$. Hence, the same holds for all $t\in [0,1]$. In particular, for any fixed $s\in (0,1]$, we have
    \[
        (\ell^*_{\tau}+C)\land \phi-Cs\leq \ell_{sM}-sM\tau
    \]
    for all $M>0$.
    Taking infimum with respect to $M>0$, we find
    \[
        (\ell^*_{\tau}+C)\land \phi-Cs\leq \ell^*_{\tau}.
    \]
    Since $s\in (0,1]$ is arbitrary, we conclude \eqref{eq:ellstarleqetemp1}.

    Finally, we check the condition at $-\infty$. Set
    \[
        \psi\coloneqq \uscsup[\tau<\ell^*_{\max}]\ell^*_{\tau}.
    \]
    The family $(\ell^*_{\tau})_{\tau<\ell^*_{\max}}$ is locally uniformly bounded from above, so $\psi\in\PSH(X,\theta)$. Thanks to \cref{cor:subgeodboundary}, outside a pluripolar set, the function $t\mapsto \ell_t(x)$ is a finite closed convex function on $[0,\infty)$ with value $\phi(x)$ at $t=0$.
    Applying the one-variable Legendre duality gives
    \[
        \psi(x)=\adjustlimits \sup_{\tau<\ell^*_{\max}}\inf_{t>0}\bigl(\ell_t(x)-t\tau\bigr)=\ell_0(x)=\phi(x)\quad \textup{quasi-everywhere}.
    \]
    Hence $\psi=\phi$ by \cref{prop:pshfuncdetdense}. In particular, $\ell^*_{-\infty}=\phi$ has positive mass.
\end{proof}

\begin{theorem}\label{thm:Legenbij}
    The Legendre transform in \cref{def:Legtrans} is a bijection
    \begin{equation}\label{eq:RWNco}
        \TC(X,\theta;\phi)\xrightarrow{\sim} \mathcal{R}(X,\theta;\phi).
    \end{equation}
    Moreover, this bijection restricts to the following bijections:
    \begin{equation}\label{eq:RWNcotwocases}
    \TC^1(X,\theta;\phi)\xrightarrow{\sim} \mathcal{R}^1(X,\theta;\phi),\quad \TC^{\infty}(X,\theta;\phi)\xrightarrow{\sim} \mathcal{R}^{\infty}(X,\theta;\phi).
    \end{equation}
    For any $\Gamma\in \TC^1(X,\theta;\phi)$, we have
    \begin{equation}\label{eq:RWNenergy}
    \mathbf{E}^{\phi}(\Gamma)=\mathbf{E}^{\phi}(\Gamma^*).
    \end{equation}
\end{theorem}
Recall that the two energy functionals in \eqref{eq:RWNenergy} are defined in \eqref{eq:tcfiniteenergy} and \cref{def:radialMAenergy2} respectively.

The correspondence \eqref{eq:RWNco} will be referred to as the \emph{Ross--Witt Nyström correspondence}\index{Ross--Witt Nyström correspondence}.

To appreciate this result, just consider the simple case where $\theta$ is a Kähler form and $\phi=0$. In this case, elements in $\mathcal{R}^{\infty}(X,\theta)$ are rays of \emph{bounded} potentials, while elements in $\TC^{\infty}(X,\theta)$ are rays of \emph{singular} potentials. This result establishes a bridge between the pluripotential theory of regular potentials and that of singular potentials!

\begin{proof}
    \textbf{Step~1}. We first establish \eqref{eq:RWNco}.

    It follows from \cref{lma:LegsendsTCtoR} that the forward map is well-defined. The inverse map is given by \eqref{eq:invLeg}. We show that the inverse map is also well-defined. Given $\ell\in \mathcal{R}(X,\theta;\phi)$, we know from \cref{lma:LegendsRtoTC} that $\ell^*\in \TC(X,\theta;\phi)$.
    We need to show that $\ell^*\in \TC(X,\theta;\phi)$.

    The equality $\ell^*_{-\infty}=\phi$ was proved in \cref{lma:LegendsRtoTC}. Therefore, $\ell^*\in \TC(X,\theta;\phi)$ as expected.

    The two operations are inverse to each other thanks to \cref{cor:Leginvpartdef}. Hence, \eqref{eq:RWNco} is established.

    \textbf{Step~2}.
    Next we consider the bounded situation. Namely, we want to establish the second half of \eqref{eq:RWNcotwocases}.

    Suppose that $\Gamma\in \TC^{\infty}(X,\theta;\phi)$. Take $\tau_0\in \mathbb{R}$ so that $\Gamma_{\tau}=\phi$ for all $\tau\leq \tau_0$. It follows from \eqref{eq:testcurveLegtran} that
    \[
    \Gamma_t^*\geq \phi+t\tau_0
    \]
    for all $t>0$. Therefore, $\Gamma^*_t\sim \phi$ for all $t>0$ and hence $\Gamma^*\in \mathcal{R}^{\infty}(X,\theta;\phi)$.

    Conversely, suppose that $\ell\in \mathcal{R}^{\infty}(X,\theta;\phi)$. Thanks to \cref{prop:geodinfsublinear}, there is a constant $C>0$ such that
    \[
    \ell_t\geq \phi-Ct.
    \]
    Therefore, according to \eqref{eq:invLeg}, we have
    \[
    \ell^*_{\tau}\geq \inf_{t>0} \Bigl( \phi-(C+\tau)t \Bigr)=\phi
    \]
    if $\tau\leq -C$. Therefore, $\ell^*_{\tau}=\phi$ for all $\tau\leq -C$.

    \textbf{Step~3}. We establish \eqref{eq:RWNenergy} and the first half of \eqref{eq:RWNcotwocases}.

    \textbf{Step~3.1}. We reduce to the case where $\Gamma_{\max}=0$.

    Suppose that we define
    \[
    \Gamma'_{\tau}=\Gamma_{\tau+\Gamma_{\max}},\quad \forall \tau<0.
    \]
    Then $\Gamma'\in \TC(X,\theta;\phi)$ as well and for all $t>0$,
    \[
    \Gamma'^*_t=\sup_{\tau<0} \left( t\tau+\Gamma'_{\tau} \right)=\sup_{\tau<\Gamma_{\max}} \left( t\tau+\Gamma_{\tau} \right)-t\Gamma_{\max}=\Gamma^*_t-t\Gamma_{\max}.
    \]
    Therefore,
    \[
    \mathbf{E}^{\phi}\left(\Gamma'^*\right)=\mathbf{E}^{\phi}\left(\Gamma^*\right)-\Gamma_{\max}\int_X \theta_{\phi}^n.
    \]
    by \eqref{eq:Ephiaddcst}. Using \eqref{eq:tcfiniteenergy}, we also have
    \[
    \begin{aligned}
    \mathbf{E}^{\phi}(\Gamma')=&\int_{-\infty}^0 \left( \int_X \theta_{\Gamma'_{\tau}}^n-\int_X \theta_{\phi}^n \right)\,\mathrm{d}\tau\\
    =& \int_{-\infty}^{\Gamma_{\max}}\left( \int_X \theta_{\Gamma_{\tau}}^n-\int_X \theta_{\phi}^n \right)\,\mathrm{d}\tau\\
    =& \mathbf{E}^{\phi}(\Gamma)-\Gamma_{\max}\int_X \theta_{\phi}^n.
    \end{aligned}
    \]
    Therefore, it suffices to establish \eqref{eq:RWNenergy} for $\Gamma'$ in place of $\Gamma$.

    \textbf{Step~3.2}. We assume that $\Gamma_{\max}=0$ and $\Gamma\in \TC^{\infty}(X,\theta;\phi)$. We prove \eqref{eq:RWNenergy}.

    For $N\in \mathbb{Z}_{>0}$, $M \in \mathbb{Z}$, we introduce the following:
    \[
        \Gamma^{*,N,M}_t \coloneqq \max_{\substack{k \in \mathbb{Z} \\ k \leq M}} \left(\Gamma_{k/2^N} +  tk/2^N\right)\in \mathcal{E}^{\infty}(X,\theta;\phi),\quad t>0.
    \]
    We first claim that for all $t>0$, $N\in \mathbb{Z}_{>0}$ and $M\in \mathbb{Z}_{<0}$,
    \begin{equation}\label{eq:diff_eq_I}
        \frac{t}{2^N} \int_X \theta^n_{\Gamma_{(M+1)/2^N}}\leq E_{\theta}^{\phi}\left(\Gamma_t^{*,N,M+1}\right)-E_{\theta}^{\phi}\left(\Gamma_t^{*,N,M}\right)\leq \frac{t}{2^N} \int_X \theta^n_{\Gamma_{M/2^N}}.
    \end{equation}

    Assuming this, let us prove \eqref{eq:RWNenergy}.

    Fixing $N$, let $M=\lfloor 2^{N}\Gamma_{\min}\rfloor$. Recall that $\Gamma_{\min}$ is defined in \eqref{eq:Gammamin}.
    Then repeated applications of \eqref{eq:diff_eq_I} yield
    \[
        \sum_{j=M+1}^0 \frac{t}{2^N} \int_X \theta^n_{\Gamma_{j/2^N}}\leq E^{\phi}_{\theta}\left(\Gamma_t^{*,N,0}\right)-E^{\phi}_{\theta}\left(\Gamma_t^{*,N,M}\right)\leq \sum_{j=M}^{-1} \frac{t}{2^N} \int_X \theta^n_{\Gamma_{j/2^N}}.
    \]
    Since $M \leq 2^{N}\Gamma_{\min}$, we have
    \[
    \Gamma^{*,N,M}_t =\phi+t M/2^N.
    \]
    Using \eqref{eq:Ephiaddcst}, we can continue to write
    \[
        \sum_{j=M+1}^0 \frac{t}{2^N} \left(\int_X \theta^n_{\Gamma_{j/2^N}} - \int_X \theta^n_{\phi}\right)\leq E_{\theta}^{\phi}\left(\Gamma_t^{*,N,0}\right)\leq \sum_{j=M}^{-1} \frac{t}{2^N} \left(\int_X \theta^n_{\Gamma_{j/2^N}} - \int_X \theta^n_{\phi}\right).
    \]
We now notice that we have Riemann sums on both the left and right of the above inequality. Moreover, as $N\to\infty$, the potentials $\Gamma_t^{*,N,0}$ increase to $\Gamma_t^*$, hence their energies converge to $E^\phi_\theta(\Gamma_t^*)$ by \cref{prop:Econt_mono}. Using \cref{prop:testcurvmasslogconc}, it is possible to let $N \to \infty$ and obtain
\begin{equation}\label{eq:EGammastartt}
E_{\theta}^{\phi}(\Gamma_t^{*})=t\mathbf{E}^{\phi}(\Gamma)
\end{equation}
So \eqref{eq:RWNenergy} follows as desired.

It remains to argue \eqref{eq:diff_eq_I}. Fix $t>0$, $N\in \mathbb{Z}_{>0}$ and $M\in \mathbb{Z}$. By \cref{prop:Ediffbdd},
    \begin{equation}\label{eq:first_I_ineq}
    \begin{aligned}
        \int_X \left(\Gamma_t^{*,N,M+1}-\Gamma_t^{*,N,M}\right)\,\theta^n_{\Gamma_t^{*,N,M+1}}\leq & E_{\theta}^{\phi}\left(\Gamma_t^{*,N,M+1}\right)-E_{\theta}^{\phi}\left(\Gamma_t^{*,N,M}\right)\\
        \leq & \int_X \left(\Gamma_t^{*,N,M+1}-\Gamma_t^{*,N,M}\right)\,\theta^n_{\Gamma_t^{*,N,M}}.
    \end{aligned}
    \end{equation}
    Clearly $\Gamma_t^{*,N,M+1} \geq \Gamma_t^{*,N,M}$. Moreover, since $\mathbb{R}\ni \tau\mapsto \Gamma_{\tau}+t\tau$ is concave, we notice that
    \[
        U_t \coloneqq\left\{\Gamma^{*,N,M+1}_t>\Gamma^{*,N,M}_t\right\} = \left\{\Gamma_{(M+1)/2^N} +  2^{-N}t >\Gamma_{M/2^N}\right\},
    \]
    and on $U_t$ we have
    \begin{equation}\label{eq:GammstaronUntemp1}
        \Gamma_t^{*,N,M+1} = \Gamma_{(M+1)/2^N} + t (M+1)/2^{N},\quad \Gamma_t^{*,N,M} = \Gamma_{M/2^N} + t M /2^{N}.
    \end{equation}
     We also note that $U_t$ is $\mathcal{F}$-open by \cref{cor:diffpshFopen}.
     So from the lower bound in \eqref{eq:first_I_ineq}, we have
     \[
     \begin{aligned}
     E_{\theta}^{\phi}\left(\Gamma_t^{*,N,M+1}\right)-E_{\theta}^{\phi}\left(\Gamma_t^{*,N,M}\right)\geq & \int_{U_t} \left(\Gamma_t^{*,N,M+1}-\Gamma_t^{*,N,M}\right)\,\theta^n_{\Gamma_t^{*,N,M+1}}\\
     = & \int_{U_t} \left(\Gamma_{(M+1)/2^N}-\Gamma_{M/2^N}+t2^{-N}\right)\,\theta^n_{\Gamma_{(M+1)/2^N}}\\
     \geq & \int_{\left\{\Gamma_{(M+1)/2^N}=0\right\}}t2^{-N}\,\theta^n_{\Gamma_{(M+1)/2^N}},
     \end{aligned}
     \]
     where on the second line, we applied \eqref{eq:GammstaronUntemp1} and \cref{prop:npp1}, on the third line, we applied the fact that $\theta_{\Gamma_{(M+1)/2^N}}^n$ is supported on the set (cf. \cref{thm:Pvarphisupport})
     \[
         \left\{\Gamma_{(M+1)/2^N}=0\right\}\subseteq U_t\cap \left\{\Gamma_{M/2^N}=0\right\}.
     \]
     When $(M+1)/2^N=0$, we use \cref{prop:decseqmodel} to see that $P_\theta[\Gamma_0]=\Gamma_0$, since $\Gamma_0$ is the decreasing limit of the model potentials $\Gamma_\tau$ as $\tau\nearrow 0$. Thus \cref{thm:Pvarphisupport} applies to this endpoint as well.

     We have deduced the first inequality in \eqref{eq:diff_eq_I}. Next, we apply the upper bound part in \eqref{eq:first_I_ineq} and compute similarly
     \[
     \begin{aligned}
     E_{\theta}^{\phi}(\Gamma_t^{*,N,M+1})-E_{\theta}^{\phi}(\Gamma_t^{*,N,M})\leq & \int_X \left(\Gamma_t^{*,N,M+1}-\Gamma_t^{*,N,M}\right)\,\theta^n_{\Gamma_t^{*,N,M}}\\
     =& \int_{U_t} \left(\Gamma_{(M+1)/2^N}-\Gamma_{M/2^N}+t2^{-N}\right)\,\theta^n_{\Gamma_{M/2^N}}\\
     \leq & \int_{\{\Gamma_{M/2^N}=0\}\cap U_t}\left(\Gamma_{(M+1)/2^N}+t2^{-N}\right)\,\theta^n_{\Gamma_{M/2^N}} \\
     \leq &\int_{\{\Gamma_{M/2^N}=0\}\cap U_t} t2^{-N}\,\theta^n_{\Gamma_{M/2^N}}.
     \end{aligned}
     \]
     We conclude the latter half of \eqref{eq:diff_eq_I}.

    \textbf{Step~3.3}. We assume that $\Gamma_{\max}=0$. Now $\Gamma\in \TC(X,\theta;\phi)$ only.

    For each $\epsilon>0$, we introduce $\Gamma^{\epsilon}\in \TC^{\infty}(X,\theta;\phi)$ as follows:
    \begin{enumerate}
        \item\label{itm:chapter-testcurve:L718:1} Let $\Gamma^{\epsilon}_{\max}=0$, and
        \item\label{itm:chapter-testcurve:L718:2} we set
        \[
        \Gamma^{\epsilon}_{\tau}=
        \begin{cases}
            \phi,&\textup{if }\tau\leq -\epsilon^{-1};\\
            P_{\theta}\left[ (1+\epsilon\tau)\Gamma_{\tau}-\epsilon\tau \phi \right],&\textup{if }\tau\in \left(-\epsilon^{-1},0\right).
        \end{cases}
        \]
    \end{enumerate}
    The concavity and monotonicity in $\tau$ follow from \cref{prop:Pconc}. The value at $-\infty$ is $\phi$. Moreover, $\Gamma^\epsilon_\tau=\phi$ for $\tau\leq-\epsilon^{-1}$, so $\Gamma^\epsilon\in\TC^\infty(X,\theta;\phi)$.

    It follows from \cref{cor:lineardsconv} and \cref{cor:decnetdS} that for each $\tau<0$, the potentials $\Gamma^{\epsilon}_{\tau}$ decrease to $\Gamma_{\tau}$ as $\epsilon\searrow 0$. Applying \cref{prop:vol_limit_model}, we have
    \[
    \lim_{\epsilon \to 0+} \int_X \left(\theta+\ddc \Gamma^{\epsilon}_{\tau}\right)^n=\int_X \left(\theta+\ddc \Gamma_{\tau}\right)^n
    \]
    for all $\tau<0$. Hence, by the monotone convergence theorem and Step~3.2, we find
    \begin{equation}\label{eq:EphiGammatemp1}
        \mathbf{E}^{\phi}(\Gamma)=\lim_{\epsilon\to 0+}\mathbf{E}^{\phi}(\Gamma^{\epsilon})=\lim_{\epsilon\to 0+}\mathbf{E}^{\phi}(\Gamma^{\epsilon *})=\lim_{\epsilon\to 0+} E_{\theta}^{\phi}(\Gamma^{\epsilon *}_1),
    \end{equation}
    where the last equality follows from \eqref{eq:EGammastartt}.
    Furthermore, according to \cref{prop:Legendretranssup}, we have
    \[
    \Gamma^*_t=\inf_{\epsilon>0} \Gamma^{\epsilon *}_t
    \]
    for all $t>0$. Note that we do not have to take the closure of the right-hand side since it is automatically upper semicontinuous in $t$.

    Now suppose that $\Gamma\in \TC^1(X,\theta;\phi)$.
    Then by \eqref{eq:EphiGammatemp1}, as $\epsilon\to 0+$, $(\Gamma_t^{\epsilon *})_{\epsilon}$ is a decreasing Cauchy net in $\mathcal{E}^1(X,\theta;\phi)$ and hence by \cref{thm:Efuncont} for each $t>0$,
    \[
    E_{\theta}^{\phi}\left(\Gamma^*_t\right)=\lim_{\epsilon\to 0+}E_{\theta}^{\phi}\left(\Gamma^{\epsilon,*}_t\right)=t\mathbf{E}^{\phi}(\Gamma)>-\infty,
    \]
    where we have applied \eqref{eq:EGammastartt} and \eqref{eq:EphiGammatemp1}.
    Hence, $\Gamma^*\in \mathcal{R}^1(X,\theta;\phi)$. Moreover, \eqref{eq:RWNenergy} follows.

    Conversely, suppose that $\Gamma^*\in \mathcal{R}^1(X,\theta;\phi)$. Then \eqref{eq:EphiGammatemp1} implies that
    \[
    \mathbf{E}^{\phi}(\Gamma)=\lim_{\epsilon\to 0+} E_{\theta}^{\phi}\left(\Gamma^{\epsilon *}_1\right) \geq E_{\theta}^{\phi}\left(\Gamma^*_1\right)>-\infty.
    \]
    Hence, $\Gamma\in \TC^1(X,\theta;\phi)$.
\end{proof}

\begin{remark}\label{rmk:testcurve:L815}
    One could also consider geodesic rays emanating from another potential $\Phi \in \mathcal{E}(X,\theta;\phi)$. In this case, one can show that these geodesic rays are in bijection with $\Phi$-twisted test curves: In \cref{def:testcur}, we replace \ref{itm:def:testcur:2} by the following condition:
    \[
        P_{\theta}\left[\Gamma_{\tau}\right](\Phi)=\Gamma_{\tau}.
    \]
    Furthermore, we require that $\Gamma_{-\infty}=\Phi$.

    The Legendre transform and its inverse are defined by the same formulas. The proofs of \cref{lma:LegsendsTCtoR,lma:LegendsRtoTC,prop:Leginvtc} apply with base potential $\Phi$. For $t>0$, the normalization in \cref{lma:suplegenlinear} becomes
    \[
        \sup_{\{\Phi>-\infty\}}\bigl(\Gamma^*_t-\Phi\bigr)=t\Gamma_{\max}.
    \]
    Indeed, $\Gamma_\tau\leq\Phi$ for every $\tau<\Gamma_{\max}$, and $\sup_{\{\Phi>-\infty\}}(\Gamma_\tau-\Phi)=0$ for each such $\tau$ by the envelope condition above. Taking the supremum in the Legendre formula gives the identity.
\end{remark}

\begin{corollary}\label{cor:reltestcursuplinear}
    Let $\ell\in \mathcal{R}(X,\theta;\phi)$. Then $\sup_X \ell_t=\ell^*_{\max}t$ for any $t\geq 0$.

    In particular, $\ell_t-\ell^*_{\max}t$ is a decreasing function of $t\geq 0$.
\end{corollary}
\begin{proof}
    This follows from \cref{lma:suplegenlinear} and \cref{thm:Legenbij}.
\end{proof}

\begin{example}\label{ex:testcurve:L843}
    Let us see what the test curve in \cref{ex:cantest} corresponds to under the Ross--Witt Nyström correspondence. Fix $\varphi\in \PSH(X,\theta)$.
    We claim that
    \begin{equation}\label{eq:RWNsingtype}
    \ell^{\varphi}=\Gamma^{\varphi*},
    \end{equation}
    where $\ell^{\varphi}$ is as in \cref{ex:rayasspsh}.
    We may assume that $\varphi\leq 0$ since both sides are invariant after adding a constant to $\varphi$.

    We first prove the easy direction $\ell^{\varphi}\geq \Gamma^{\varphi*}$, which is equivalent to $\ell^{\varphi *}\geq \Gamma^{\varphi}$. Since $\ell^{\varphi *}$ is a test curve, the latter is equivalent to
    \[
    \ell^{\varphi*}_{\tau}\geq (1+\tau)\varphi-\tau V_{\theta}
    \]
    for all $\tau\in (-1,0)$, together with $\ell^{\varphi*}_{-1}\geq V_\theta$ for the constant branch of $\Gamma^\varphi$. By Legendre duality and monotonicity of $\ell^{\varphi*}$ in $\tau$, these inequalities are equivalent to
    \begin{equation}\label{eq:lvarphigeq_temp1}
    \ell^{\varphi}_t\geq (V_\theta-t)\lor \sup_{\tau\in (-1,0)} \Bigl((1+\tau)\varphi-\tau V_{\theta}+t\tau\Bigr)=\varphi\lor (V_{\theta}-t)
    \end{equation}
    for all $t>0$. Where $\varphi$ is finite, the expression in the supremum is affine in $\tau$ with endpoint values $V_\theta-t$ and $\varphi$. Where $\varphi=-\infty$, the term $V_\theta-t$ gives the displayed maximum.

    Using the notations of \cref{ex:rayasspsh}, we find easily that
    \[
    \ell^{\varphi,C}_t\geq \varphi\lor (V_{\theta}-t)
    \]
    for any $C>0$ and $t\in [0,C]$, since $\varphi\lor (V_\theta-t)$ is a candidate for the Perron envelope as in \cref{ex:rayasspsh}. Letting $C\to \infty$, we find \eqref{eq:lvarphigeq_temp1}.
    Therefore, $\ell^{\varphi}\geq \Gamma^{\varphi*}$ follows.

    In order to prove the equality in \eqref{eq:RWNsingtype}, it suffices to show that the two sides have the same energy, as a consequence of \eqref{eq:d1rayscompa}. So we compute
    \[
    \begin{aligned}
        \mathbf{E}\left(\Gamma^{\varphi*}\right)=& \mathbf{E}\left(\Gamma^{\varphi}\right)\\
        =& \int_{-1}^0 \left(\int_X \theta_{(1+\tau)\varphi-\tau V_{\theta}}^n-\int_X \theta_{V_{\theta}}^n\right)\,\mathrm{d}\tau\\
        =& \sum_{j=0}^{n} \binom{n}{j} \int_X \theta_{\varphi}^{j}\wedge \theta_{V_{\theta}}^{n-j} \int_0^1 \tau^j (1-\tau)^{n-j}\,\mathrm{d}\tau-\int_X \theta_{V_{\theta}}^n\\
        =& \sum_{j=0}^{n} \binom{n}{j}\frac{j!(n-j)!}{(n+1)!} \int_X \theta_{\varphi}^{j}\wedge \theta_{V_{\theta}}^{n-j}-\int_X \theta_{V_{\theta}}^n\\
        =& \mathbf{E}(\ell^{\varphi}),
    \end{aligned}
    \]
    where in the second line we used \cref{prop:Ppresmass}, in the fourth line we used the value of the $\beta$-function\footnote{Also known as Euler integral of the first kind.}, and the last line is just \eqref{eq:Elphi}.
\end{example}

The multiplier ideal sheaves of a test curve can be characterized using the corresponding geodesic ray in a very simple manner.
\begin{proposition}[He--Testorf--Wang]\label{prop:mis_RWN}\index{theorem!He--Testorf--Wang theorem}
    Let $\ell\in \mathcal{R}(X,\theta;\phi)$. Given any $\tau<\ell^*_{\max}$ and $x\in X$, we have
    \begin{equation}\label{eq:mis_RWN}
    \mathcal{I}\left(\ell^*_{\tau}\right)_{x}=\left\{ f\in \mathcal{O}_{X,x}: |f|^2\int_0^{\infty} \exp(-\ell_t+t\tau)\,\mathrm{d}t \textup{ is integrable near }x \right\}.
    \end{equation}
\end{proposition}
\begin{proof}
    Fix $x\in X$, $\tau<\ell^*_{\max}$ and $f\in \mathcal{O}_{X,x}$. Fix a K\"ahler form $\omega$ on $X$.

    \textbf{Step~1}. We first assume that $f$ lies in the right-hand side of \eqref{eq:mis_RWN}.

    Given any $y\in X$, it follows from \eqref{eq:invLeg} that there is $t_0>0$ with
    \[
    \ell^*_{\tau}(y)+1\geq \ell_{t_0}(y)-t_0\tau.
    \]
    Observe that $t\mapsto \ell_t-t\ell^*_{\max}$ is decreasing in $t$ by \cref{cor:reltestcursuplinear}, it follows that for $t\in [t_0,t_0+1]$, we have
    \[
    \ell^*_{\tau}(y)+1-t_0(\ell^*_{\max}-\tau)\geq \ell_{t_0}(y)-t_0\ell^*_{\max}\geq \ell_{t}(y)-t\ell^*_{\max}.
    \]
    Since $\tau< \ell^*_{\max}$ and $t_0\geq t-1$, we deduce that
    \begin{equation}\label{eq:mis_RWN_temp1}
    \ell^*_{\tau}(y)+1+\ell^*_{\max}-\tau \geq \ell_{t}(y)-t\tau,\quad t\in [t_0,t_0+1].
    \end{equation}
    Take a sufficiently small open neighborhood $U$ of $x$ such that
    \[
    \int_U |f|^2\int_0^{\infty} \exp(-\ell_t+t\tau)\,\mathrm{d}t\,\omega^n<\infty.
    \]
    Applying \eqref{eq:mis_RWN_temp1}, we deduce that
    \[
    \int_U |f|^2\exp\left(-\ell^*_{\tau}\right)\,\omega^n<\infty.
    \]
    Therefore, $f\in \mathcal{I}(\ell^*_{\tau})_x$.

    \textbf{Step~2}. Assume that $f\in \mathcal{I}(\ell^*_{\tau})_x$.

    It follows from the strong openness theorem \cref{thm:strongopen} (cf. \eqref{eq:Gammataucont}) that $f\in \mathcal{I}(\ell^*_{\tau+\epsilon})_x$ for some small enough $\epsilon>0$ with $\tau+\epsilon<\ell^*_{\max}$.
    Take a sufficiently small open neighborhood $U$ of $x$ such that
    \[
    \int_U |f|^2 \exp\left(-\ell^*_{\tau+\epsilon}\right)\,\omega^n<\infty.
    \]
    We compute
    \[
    \begin{aligned}
    \int_U |f|^2\int_0^{\infty} \exp\left(-\ell_t+t\tau\right)\,\mathrm{d}t\,\omega^n\leq & \int_U |f|^2\int_0^{\infty} \exp\left(-\ell^*_{\tau+\epsilon}-t\epsilon\right)\,\mathrm{d}t\,\omega^n\\
    =&\frac{1}{\epsilon}\int_U |f|^2\exp\left(-\ell^*_{\tau+\epsilon}\right)\,\omega^n\\
    <&\infty.
    \end{aligned}
    \]
    Therefore, $f$ lies in the right-hand side of \eqref{eq:mis_RWN}.
\end{proof}

The masses of potentials on a test curve can also be expressed in terms of the corresponding geodesic ray.
\begin{proposition}[Hisamoto]\label{prop:His_mass_tc}\index{theorem!Hisamoto's theorem}
        Let $\ell\in \mathcal{R}(X,\theta;\phi)$. Given any $\tau<\ell^*_{\max}$, we have
    \begin{equation}\label{eq:mass_tc}
    \int_X \left(\theta+\ddc \ell^*_{\tau}\right)^n=\int_{\{\dot{\ell}_0\geq \tau\}} \theta_{\phi}^n.
    \end{equation}
\end{proposition}
Here $\dot{\ell}_0$ denotes the right-derivative of $\ell_t$ with respect to $t$ at $t=0$. It is well-defined quasi-everywhere by \cref{cor:subgeodboundary}, and hence the right-hand side of \eqref{eq:mass_tc} makes sense, even though $\{\dot{\ell}_0\geq \tau\}$ is not well-defined.

\begin{proof}
    Fix $\tau<\ell^*_{\max}$.
    We first observe that it suffices to show that
    \begin{equation}\label{eq:mass_tc_temp1}
    \int_{\left\{\dot{\ell}_0\geq \tau\right\}} \theta_{\phi}^n=\int_{\left\{\ell^*_{\tau}=\phi\right\}} \theta_{\phi}^n.
    \end{equation}

    Assuming \eqref{eq:mass_tc_temp1} for the moment.
    By \cref{thm:Legenbij}, $\ell^*\in\TC(X,\theta;\phi)$. In particular, $\ell^*_\tau$ is a model potential and
    \[
        \ell^*_\tau\leq \ell^*_{-\infty}=\phi.
    \]
    Since $\phi$ is model as well, \cref{prop:DNTenvelope2} gives
    \[
        \theta^n_{\ell^*_{\tau}}=\mathds{1}_{\left\{\ell^*_{\tau}=0\right\}}\,\theta^n,\quad \theta^n_{\phi}=\mathds{1}_{\{\phi=0\}}\,\theta^n.
    \]
    Since $\ell^*_{\tau}\leq \phi \leq 0$, we have
    \[
        \mathds{1}_{\left\{\ell^*_{\tau}=\phi\right\}}\,\theta_{\phi}^n
        =\mathds{1}_{\left\{\ell^*_{\tau}=\phi=0\right\}}\,\theta^n
        =\mathds{1}_{\left\{\ell^*_{\tau}=0\right\}}\,\theta^n
        =\theta^n_{\ell^*_{\tau}}.
    \]
    Therefore, \eqref{eq:mass_tc} follows from \eqref{eq:mass_tc_temp1}.

    In order to prove \eqref{eq:mass_tc_temp1}, it suffices to show that
    \begin{equation}\label{eq:dotell_temp1}
    \left\{\dot{\ell}_0\geq \tau\right\}=\left\{\ell^*_{\tau}=\phi\right\}\quad \textup{outside a pluripolar set}.
    \end{equation}

    Take $x\in X$ so that $(\ell_t(x))_{t\geq 0}$ is finite and right-differentiable at $t=0$. Note that this condition holds quasi-everywhere. Suppose that $\ell^*_{\tau}(x)=\phi(x)$. Then
    \[
    \dot{\ell}_0(x)=\inf_{t>0} \frac{\ell_t(x)-\phi(x)}{t}\geq \inf_{t>0} \frac{\ell^*_{\tau}(x)+t\tau-\phi(x)}{t}=\tau.
    \]
    Therefore, the $\supseteq$ direction in \eqref{eq:dotell_temp1} follows.
    Conversely, suppose that $\dot{\ell}_0(x)\geq \tau$, then
    \[
    \ell^*_{\tau}(x)=\inf_{t>0} \left(\ell_t(x)-t\tau\right)\geq \inf_{t>0} \left(\phi(x)+t\dot{\ell}_0(x)-t\tau\right)\geq \phi(x).
    \]
    Together with $\ell^*_\tau\leq\phi$, this gives $\ell^*_\tau(x)=\phi(x)$.
    Therefore, the $\subseteq$ direction in \eqref{eq:dotell_temp1}  follows as well.
\end{proof}

\begin{corollary}\label{cor:ray_energy_slope}
    Let $\ell\in \mathcal{R}^1(X,\theta;\phi)$. Then $\dot{\ell}_0\in L^1(X,\theta_{\phi}^n)$ and
    \begin{equation}\label{eq:ray_energy_slope}
        \mathbf{E}^{\phi}(\ell)=\int_X \dot{\ell}_0\,\theta_{\phi}^n.
    \end{equation}
    In particular, for every $t\geq 0$,
    \[
        E_{\theta}^{\phi}(\ell_t)=t\int_X \dot{\ell}_0\,\theta_{\phi}^n.
    \]
\end{corollary}
\begin{proof}
    Set $\Gamma=\ell^*$. By \cref{thm:Legenbij}, we have $\Gamma\in \TC^1(X,\theta;\phi)$ and
    \[
        \mathbf{E}^{\phi}(\ell)=\mathbf{E}^{\phi}(\Gamma).
    \]
    By \cref{prop:His_mass_tc}, for every $\tau<\ell^*_{\max}$,
    \[
        \int_X \theta_{\Gamma_{\tau}}^n=\int_{\{\dot{\ell}_0\geq \tau\}} \theta_{\phi}^n.
    \]
    Hence \eqref{eq:tcfiniteenergy} gives
    \[
        \mathbf{E}^{\phi}(\ell)
        =\ell^*_{\max}\int_X \theta_{\phi}^n-\int_{-\infty}^{\ell^*_{\max}} \int_{\{\dot{\ell}_0<\tau\}}\theta_{\phi}^n\,\mathrm{d}\tau.
    \]
    For $\tau>\ell^*_{\max}$, we have $\ell^*_{\tau}\equiv -\infty$ by \cref{lma:LegendsRtoTC}. Then \eqref{eq:mass_tc_temp1} shows that 
    \[
        \int_{\{\dot{\ell}_0\geq \tau\}}\theta_{\phi}^n=0.
    \]
    The layer-cake formula therefore gives \eqref{eq:ray_energy_slope}. The final assertion follows from \cref{prop:energylinear}, since $E_{\theta}^{\phi}(\ell_0)=E_{\theta}^{\phi}(\phi)=0$.
\end{proof}

\section{\texorpdfstring{$\mathcal{I}$}{I}-model test curves}\label{sec:Imodeltc}
Let $X$ be a connected compact K\"ahler manifold of dimension $n$ and $\theta$ be a smooth closed real $(1,1)$-form on $X$ representing a big cohomology class.
Fix a model potential $\phi\in \PSH(X,\theta)_{>0}$.

\begin{definition}\label{def:Imodeltestcur}
    A test curve $\Gamma\in \TC(X,\theta;\phi)$ is \emph{$\mathcal{I}$-model}\index{test curve!I-model test curve@$\mathcal{I}$-model test curve} if for any $\tau<\Gamma_{\max}$, the potential $\Gamma_{\tau}$ is $\mathcal{I}$-model.

    The subset of $\mathcal{I}$-model test curves in $\TC(X,\theta;\phi)$ is denoted by $\mathcal{E}^{\NA}(X,\theta;\phi)$\index{$\mathcal{E}^{\NA}(X,\theta;\phi)$}. When $\phi=V_{\theta}$, we omit $\phi$ and write $\mathcal{E}^{\NA}(X,\theta)$\index{$\mathcal{E}^{\NA}(X,\theta)$} instead.

    The union of the sets of $\mathcal{I}$-model test curves in $\PSH(X,\theta)$ for all model potentials $\phi\in \PSH(X,\theta)_{>0}$ is denoted by $\PSH^{\NA}(X,\theta)_{>0}$\index{$\PSH^{\NA}(X,\theta)_{>0}$}.
\end{definition}
Note that $\Gamma_{\Gamma_{\max}}$ is automatically $\mathcal{I}$-model, since
\[
    \Gamma_{\Gamma_{\max}}\leq P_{\theta}\left[\Gamma_{\Gamma_{\max}}\right]_{\mathcal{I}}\leq \inf_{\tau<\Gamma_{\max}}P_{\theta}\left[\Gamma_{\tau}\right]_{\mathcal{I}}=\inf_{\tau<\Gamma_{\max}}\Gamma_{\tau}=\Gamma_{\Gamma_{\max}}.
\]

Here we write $\NA$ with non-Archimedean in mind. The precise relation with non-Archimedean pluripotential theory will be clear in \cref{chap:NApotential}.
This section and the next can be skipped on a first reading; the necessary results can be consulted when reading \cref{chap:NApotential}.

\begin{proposition}\label{prop:GammaminfImodel}
    Let $\Gamma\in \PSH^{\NA}(X,\theta)_{>0}$. Then $\Gamma_{-\infty}$ is an $\mathcal{I}$-model potential.
\end{proposition}
\begin{proof}
    This follows from \cref{prop:incnetmodelPI}.
\end{proof}

\begin{proposition}\label{prop:Imodeltestcurveindeptheta}
    Let $\theta'$ be another smooth closed real $(1,1)$-form on $X$ representing the same cohomology class as $\theta$. Then there is a canonical bijection
    \[
    \PSH^{\NA}(X,\theta)_{>0}\xrightarrow{\sim} \PSH^{\NA}(X,\theta')_{>0}.
    \]
    This bijection satisfies the obvious cocycle condition.
\end{proposition}
\begin{proof}
    Let $\theta'=\theta+\ddc g$. By \cref{prop:testcurveindeptheta}, the canonical map sends $\Gamma$ to the test curve with levels
    \[
        \Gamma'_\tau=P_{\theta'}[\Gamma_\tau-g].
    \]
    If $\Gamma_\tau$ is $\mathcal I$-model, then it is $\mathcal I$-good by \cref{ex:ImodelIgood}. Note that $\Gamma_\tau-g$ is $\mathcal I$-good since $\Gamma_{\tau}$ is. Hence
    \[
        P_{\theta'}[\Gamma_\tau-g]=P_{\theta'}[\Gamma_\tau-g]_{\mathcal I},
    \]
    so $\Gamma'_\tau$ is $\mathcal I$-model. The same argument applied to $-g$ gives the inverse map, and the cocycle condition follows from the one in \cref{prop:testcurveindeptheta}.
\end{proof}

\begin{proposition}\label{prop:ETCIbimero}
    Let $\pi\colon Y\rightarrow X$ be a proper bimeromorphic morphism from a K\"ahler manifold $Y$. Then the pointwise pull-back induces a bijection
    \[
    \pi^*\colon \mathcal{E}^{\NA}(X,\theta;\phi)\xrightarrow{\sim} \mathcal{E}^{\NA}(Y,\pi^*\theta;\pi^*\phi).
    \]
\end{proposition}
\begin{proof}
    This is an immediate consequence of \cref{prop:ETCbimero} and \cref{prop:Ienvelopebimero}.
\end{proof}

\begin{definition}\label{def:TCIenvelope}
    Given $\Gamma\in \TC(X,\theta;\phi)$, we define its \emph{$\mathcal{I}$-envelope}\index{envelope!I-envelope@$\mathcal{I}$-envelope} $P_{\theta}[\Gamma]_{\mathcal{I}}$\index{$P_{\theta}[\Gamma]_{\mathcal{I}}$} as the map
    \[
        (-\infty,\Gamma_{\max})\rightarrow \PSH(X,\theta),\quad \tau\mapsto P_{\theta}\left[\Gamma_{\tau}\right]_{\mathcal{I}}.
    \]
        More generally, for any closed real smooth positive $(1,1)$-form $\omega$ on $X$, we define $P_{\theta+\omega}[\Gamma]_{\mathcal{I}}$ as the map
        \[
        (-\infty,\Gamma_{\max})\rightarrow \PSH(X,\theta+\omega),\quad \tau\mapsto P_{\theta+\omega}\left[\Gamma_{\tau}\right]_{\mathcal{I}}.
    \]
\end{definition}

\begin{proposition}\label{prop:transitionPI}
    Let $\Gamma\in \TC(X,\theta;\phi)$. Then
    \[
        P_{\theta}[\Gamma]_{\mathcal{I}}\in \mathcal{E}^{\NA}(X,\theta;P_{\theta}[\phi]_{\mathcal{I}}).
    \]
    More generally, for any closed real smooth positive $(1,1)$-form $\omega$ on $X$, we have
    \[
        P_{\theta+\omega}[\Gamma]_{\mathcal{I}}\in \mathcal{E}^{\NA}(X,\theta+\omega;P_{\theta+\omega}[\phi]_{\mathcal{I}}).
    \]
\end{proposition}
\begin{proof}
     We prove the first assertion; the second is identical, replacing $\theta$ by $\theta+\omega$.

    Set
    \[
        \widetilde\Gamma_\tau
        \coloneqq P_\theta[\Gamma_\tau]_{\mathcal I}.
    \]
    The family is decreasing by \cref{lma:PIenvmono1}. If $\tau=(1-s)\tau_0+s\tau_1$, then the concavity of $\Gamma$ gives
    \[
        \Gamma_\tau\geq (1-s)\Gamma_{\tau_0}+s\Gamma_{\tau_1}.
    \]
    Hence, by \cref{lma:PIenvmono1} and \cref{prop:PIconc},
    \[
        \widetilde\Gamma_\tau \geq (1-s)\widetilde\Gamma_{\tau_0}+s\widetilde\Gamma_{\tau_1}.
    \]
    Thus the family is concave. Finally,
    \[
        \uscsup[\tau<\Gamma_{\max}]\widetilde\Gamma_\tau = P_\theta[\phi]_{\mathcal I}
    \]
    follows from \cref{prop:incnetmodelPI}. Since
    $P_\theta[\phi]_{\mathcal I}\geq\phi$, its mass is positive by the monotonicity theorem \cref{thm:mono}. This proves the assertion.
\end{proof}

\begin{definition}\label{def:maximalgeoray}
    Let $\phi\in \PSH(X,\theta)_{>0}$ be a model potential.
    A geodesic ray $\ell\in \mathcal{R}(X,\theta;\phi)$ is \emph{maximal}\index{geodesic ray!maximal geodesic ray} if $\ell^*$ is $\mathcal{I}$-model, namely, if $\ell^*\in\mathcal E^{\NA}(X,\theta;\phi)$.
\end{definition}

An important class of $\mathcal{I}$-model test curves is given by filtrations. We briefly recall the corresponding terminology.

\begin{definition}\label{def:testcurve:L1134}
    Let $L$ be a big line bundle. We write
    \[
    R(X,L)=\bigoplus_{k=0}^{\infty} \mathrm{H}^0(X,L^k)
    \]
    for the section ring\footnote{It is more natural to replace the direct sum by a disjoint union, leading to the notion of ringoids (\emph{annénoïdes} in French), introduced by Ducros in \cite{Duc21} in the context of Temkin's graded reduction of Berkovich germs.} of $L$.

    A \emph{filtration}\index{filtration} on $R(X,L)$ is a decreasing family of graded linear subspaces $(\mathcal{F}^{\lambda})_{\lambda\in \mathbb{R}}$ of $R(X,L)$ with graded pieces
    \[
    \mathcal{F}^{\lambda}=\bigoplus_{k=0}^{\infty} \mathcal{F}^{\lambda}_k,
    \]
    such that the following conditions are satisfied:
    \begin{itemize}
        \item The filtration is left-continuous: For any $\lambda\in \mathbb{R}$, we have
        \[
        \mathcal{F}^{\lambda}=\bigcap_{\lambda'<\lambda}\mathcal{F}^{\lambda'};
        \]
        \item the filtration is multiplicative: For any $\lambda,\lambda'\in \mathbb{R}$ and any $k,k'\in \mathbb{N}$, we have
        \[
        \mathcal{F}^{\lambda}_k\cdot \mathcal{F}^{\lambda'}_{k'}\subseteq \mathcal{F}^{\lambda+\lambda'}_{k+k'};
        \]
        \item there is an integer $C>0$ such that
        \begin{equation}\label{eq:Flinearlybdd}
        \mathcal{F}^{Cm}_{m}=0,\quad \mathcal{F}^{-Cm}_m=\mathrm{H}^0(X,L^m)
        \end{equation}
        for all $m\in \mathbb{Z}_{>0}$.
    \end{itemize}

    Given a filtration $\mathcal{F}$ on $R(X,L)$, we define
    \[
    \tau_k(\mathcal{F})=\sup\left\{\lambda\in \mathbb{R}: \mathcal{F}^{\lambda}_k\neq 0\right\}.
    \]
    By Fekete's lemma, we can introduce
    \[
    \tau(\mathcal{F})=\lim_{k\to\infty}\frac{1}{k}\tau_k(\mathcal{F})=\sup_{k\in \mathbb{Z}_{>0}}\frac{1}{k}\tau_k(\mathcal{F}).
    \]
\end{definition}
Note that $\tau(\mathcal{F})$ is bounded from above by the constant $C$ in \eqref{eq:Flinearlybdd}, hence finite.

\begin{example}\label{ex:RWN}
    Let $L$ be a big line bundle on $X$ and $\mathcal{F}$ be a filtration on $R(X,L)$. Fix a smooth Hermitian metric $h$ on $L$ and write $\theta=c_1(L,h)$.

    We introduce a few auxiliary functions. For each $k\in \mathbb{Z}_{>0}$, we introduce
    \[
    \Gamma^{\mathcal{F},k}_{\tau}\coloneqq \uscsup \left\{ \log |s|_{h^k}^2: s\in \mathcal{F}^{k\tau}_k, |s|_{h^k}^2\leq 1 \right\}.
    \]
    When $k\tau\leq \tau_k(\mathcal{F})$, we know that $\mathcal{F}_k^{k\tau}\neq 0$. Moreover, \cref{prop:LelongPoincare} and \cref{prop:pshfunction_closedseq} imply that
    \[
    \Gamma^{\mathcal{F},k}_{\tau}\in \PSH(X,k\theta),\quad \tau\leq k^{-1}\tau_k(\mathcal{F}).
    \]
    Observe that for $k,k'\in \mathbb{Z}_{>0}$, we have
    \[
    \Gamma^{\mathcal{F},k+k'}_{\tau}\geq \Gamma^{\mathcal{F},k}_{\tau}+\Gamma^{\mathcal{F},k'}_{\tau}.
    \]
    Consequently, for every $\tau<\tau(\mathcal F)$ and $x\in X$,
    \begin{equation}\label{eq:limkGammaFk}
    \varlimsup_{k\to\infty}\frac{1}{k}\Gamma^{\mathcal{F},k}_{\tau}(x)=\sup_{k\in \mathbb{Z}_{>0}}\frac{1}{k}\Gamma^{\mathcal{F},k}_{\tau}(x).
    \end{equation}
    Indeed, superadditivity gives the lower bound for the $\varlimsup$ along multiples of each degree at which the value is finite. If $\Gamma^{\mathcal{F},k}_{\tau}(x)>-\infty$ for all sufficiently large $k$, then the limit exists by Fekete's lemma: for each fixed degree with finite value, use a finite-valued representative of every residue class in the finite tail. At a point that belongs to the base locus in infinitely many degrees, the $\varlimsup$ is needed in general.

    We define $(\Gamma^{\mathcal{F}}_{\tau})_{\tau<\tau(\mathcal{F})}$ as follows:
    \[
    \Gamma^{\mathcal{F}}_{\tau}\coloneqq P_{\theta}\left[ \uscsup[k\in \mathbb{Z}_{>0}] \frac{1}{k} \Gamma^{\mathcal{F},k}_{\tau}\right]\footnotemark.
    \]
    \footnotetext{It is not clear if $P_{\theta}[\bullet]$ is necessary here. When $L$ is ample, it is shown in \cite[Proposition~7.11]{RWN14} that it is not necessary. The proof in the reference relies on a Skoda division theorem \cite[Theorem~7.10]{RWN14}, which is not known in the case of big line bundles.}
    We claim that $\Gamma^{\mathcal{F}}\in \mathcal{E}^{\NA}(X,\theta)$ and is bounded.

    It is clear that $(-\infty,\tau(\mathcal{F}))\ni \tau\mapsto \Gamma^{\mathcal{F}}_{\tau}$ is decreasing. We prove its concavity. By \cref{prop:Pconc}, it suffices to show that
    \[
    \Bigl(-\infty,\tau(\mathcal{F})\Bigr)\ni\tau\mapsto \uscsup[k\in \mathbb{Z}_{>0}] \frac{1}{k} \Gamma^{\mathcal{F},k}_{\tau}
    \]
    is concave. In other words, we need to prove the following: Given $\tau_0<\tau_1<\tau(\mathcal{F})$ and $t\in (0,1)$, we have
    \[
    \uscsup[k\in \mathbb{Z}_{>0}] \frac{1}{k} \Gamma^{\mathcal{F},k}_{t\tau_1+(1-t)\tau_0}\geq t\uscsup[k\in \mathbb{Z}_{>0}] \frac{1}{k} \Gamma^{\mathcal{F},k}_{\tau_1}+(1-t)\uscsup[k\in \mathbb{Z}_{>0}] \frac{1}{k} \Gamma^{\mathcal{F},k}_{\tau_0}.
    \]
    But thanks to \cref{prop:pshfuncdetdense} and \cref{prop:supsupstardiff}, it suffices to show that
    \[
    \sup_{k\in \mathbb{Z}_{>0}} \frac{1}{k}\Gamma^{\mathcal{F},k}_{t\tau_1+(1-t)\tau_0}\geq t\sup_{k\in \mathbb{Z}_{>0}} \frac{1}{k}\Gamma^{\mathcal{F},k}_{\tau_1}+(1-t)\sup_{k\in \mathbb{Z}_{>0}} \frac{1}{k}\Gamma^{\mathcal{F},k}_{\tau_0}
    \]
    for all $t\in (0,1)$.
    Take $s_i\in \mathcal{F}^{k_i\tau_i}_{k_i}$ for $i=0,1$ with $|s_i|_{h^{k_i}}^2\leq 1$, where $k_0,k_1\in \mathbb{Z}_{>0}$.
    We need to prove that
    \begin{equation}\label{eq:GammaFconctemp1}
    \sup_{k\in \mathbb{Z}_{>0}} \frac{1}{k}\Gamma^{\mathcal{F},k}_{t\tau_1+(1-t)\tau_0}\geq \frac{1-t}{k_0}\log |s_0|_{h^{k_0}}^2+\frac{t}{k_1}\log |s_1|_{h^{k_1}}^2.
    \end{equation}
    Approximating $t$ by rational numbers $t_m\searrow t$, and using the decreasingness of the filtration in the parameter $\lambda$, we reduce to the case $t\in\mathbb Q$. Write $t=p/q$ with $p,q\in \mathbb{Z}_{>0}$. Then
    \[
    s\coloneqq s_0^{k_1(q-p)}\otimes s_1^{k_0p}\in \mathcal{F}_{k_0k_1 q}^{k_0k_1\tau_0(q-p)+k_0k_1\tau_1 p},
    \]
    and
    \[
    \begin{aligned}
    &\frac{1}{k_0k_1 q} \log |s|_{h^{k_0k_1 q}}^2\\
    =&\frac{1}{k_0k_1 q}\left(k_1(q-p)\log |s_0|^2+ k_0p\log |s_1|^2\right)\\
    = &\frac{1-t}{k_0}\log |s_0|_{h^{k_0}}^2+\frac{t}{k_1}\log |s_1|_{h^{k_1}}^2.
    \end{aligned}
    \]
    So \eqref{eq:GammaFconctemp1} follows.

    Note that for each $k\in \mathbb{Z}_{>0}$, $\tau\leq k^{-1}\tau_k(\mathcal{F})$, we know that $k^{-1}\Gamma_{\tau}^{\mathcal{F},k}$ is $\mathcal{I}$-good by \cref{prop:Igoodsup}. The family is bounded above by $0$, so it follows from the same proposition that the obstacle in the definition of $\Gamma^{\mathcal{F}}_{\tau}$ is $\mathcal{I}$-good. Hence the $P$- and $\mathcal{I}$-envelopes of this obstacle agree, so $\Gamma^{\mathcal{F}}_{\tau}$ is $\mathcal{I}$-model for each $\tau<\tau(\mathcal{F})$.

    It remains to show that the test curve $\Gamma^{\mathcal{F}}$ is bounded and lies in $\mathcal{E}^{\NA}(X,\theta)$. Fix $\tau\leq -C$, where $C$ is as in \eqref{eq:Flinearlybdd}. We will show that
    \begin{equation}\label{eq:Gammatautheta}
    \Gamma^{\mathcal{F}}_{\tau}=V_{\theta}.
    \end{equation}
    This follows from the Bergman kernel technique. Based on the theory developed so far, we can also give an elementary argument.

    Fix $k\in\mathbb Z_{>0}$. Since $\tau\leq -C$, we have $\mathcal F^{k\tau}_k=\mathrm H^0(X,L^k)$. Let $s\in\mathrm H^0(X,L^k)$. After multiplying $s$ by a constant, we may assume that $|s|_{h^k}^2\leq1$. By the definition of $\Gamma^{\mathcal F,k}_\tau$,
    \[
        \Gamma^{\mathcal F,k}_\tau\geq \log|s|_{h^k}^2.
    \]
    Hence
    \[
        s\in \mathrm H^0\left(X,L^k\otimes\mathcal I(\Gamma^{\mathcal F,k}_\tau)\right).
    \]
    Since
    \[
    k\Gamma^{\mathcal F}_\tau \geq \Gamma^{\mathcal F,k}_\tau,
    \]
    we get
    \[
        s\in \mathrm H^0\left(X,L^k\otimes\mathcal I(k\Gamma^{\mathcal F}_\tau)\right).
    \]
    Thus
    \[
        \mathrm H^0\left(X,L^k\otimes\mathcal I(k\Gamma^{\mathcal F}_\tau)\right)=\mathrm H^0(X,L^k)
    \]
    for every $k\in\mathbb Z_{>0}$. By \cref{thm:DXmain1},
    \[
        \vol\left(\theta+\ddc\Gamma^{\mathcal F}_\tau\right)=\vol L.
    \]
    As $\Gamma^{\mathcal F}_\tau$ is $\mathcal I$-model, this volume is equal to $\int_X\left(\theta+\ddc \Gamma^{\mathcal F}_\tau\right)^n$. Hence $\Gamma^{\mathcal F}_\tau$ has full Monge--Amp\`ere mass. Since it is also model, it must be equal to $V_\theta$. This proves \eqref{eq:Gammatautheta}.
\end{example}

\begin{remark}\label{rmk:PhongSturm}
    There is an important special case of \cref{ex:RWN}: Suppose that $L$ is ample and $\mathcal{F}$ is the filtration induced by a smooth test configuration $(\mathcal{X},\mathcal{L})$ of $(X,L)$. Then the geodesic ray $\Gamma^{\mathcal{F}*}$ is exactly the Phong--Sturm geodesic ray associated with $(\mathcal{X},\mathcal{L})$. See \cite[Section~9]{RWN14}.
\end{remark}

\begin{remark}\label{rmk:max_filt}
    We deduce from \cref{ex:RWN} that the ray $\Gamma^{\mathcal{F}*}$ induced by a filtration $\mathcal{F}$ is maximal.
\end{remark}

\section{Operations on test curves}\label{sec:operationtc}
Let $X$ be a connected compact K\"ahler manifold of dimension $n$ and $\theta,\theta',\theta''$ be smooth closed real $(1,1)$-forms on $X$ representing big cohomology classes.

In this section, we develop several general operations on test curves, anticipating the applications in non-Archimedean geometry in \cref{chap:NApotential}. One may read \cref{chap:NApotential} first and consult this section when necessary.

\begin{definition}\label{def:potestcurve}
    Given $\Gamma\in \TC(X,\theta)_{>0}$, $\Gamma'\in \TC(X,\theta')_{>0}$, we say $\Gamma \leq \Gamma'$ if $\Gamma_{\max}\leq \Gamma'_{\max}$ and, for all $\tau<\Gamma_{\max}$, we have
    \index{test curve!partial order}
    \begin{equation}\label{eq:GammatauGammap}
    \Gamma_{\tau}\preceq_P \Gamma'_{\tau}.
    \end{equation}
\end{definition}
Observe that \eqref{eq:GammatauGammap} actually holds for all $\tau\in \mathbb{R}$ if $\theta=\theta'$. Here we use the convention that $-\infty \preceq_P \varphi$ for any $\varphi\in \QPSH(X)\cup \{-\infty\}$.
The relation $\leq$ defines a partial order on $\TC(X,\theta)_{>0}$.
\begin{lemma}\label{lma:testcurord1}
    Let $\Gamma,\Gamma'\in \TC(X,\theta)_{>0}$ and $\omega$ be a closed real smooth positive $(1,1)$-form on $X$. Then the following are equivalent:
    \begin{enumerate}
        \item\label{itm:lma:testcurord1:1} $\Gamma\leq \Gamma'$;
        \item\label{itm:lma:testcurord1:2} $P_{\theta+\omega}[\Gamma]\leq P_{\theta+\omega}[\Gamma']$.
    \end{enumerate}
\end{lemma}
\begin{proof}
    This follows from \cref{ex:Pequiv}.
\end{proof}

\begin{definition}\label{def:sumtestcur}
    Let $\Gamma \in \TC(X,\theta)_{>0}$ and $\Gamma' \in \TC(X,\theta')_{>0}$. Then we define $\Gamma+\Gamma'\in \TC(X,\theta+\theta')_{>0}$ as follows:
    \index{test curve!sum}
    \begin{enumerate}
        \item\label{itm:def:sumtestcur:1} We set
        \[
        (\Gamma+\Gamma')_{\max}\coloneqq \Gamma_{\max}+\Gamma'_{\max};
        \]
        \item\label{itm:def:sumtestcur:2} for any $\tau<(\Gamma+\Gamma')_{\max}$, we define\footnote{There is no usc regularization in the formula. It is not a typo.}
        \begin{equation}\label{eq:GammaGammapsum}
            (\Gamma+\Gamma')_{\tau }\coloneqq P_{\theta+\theta'}\left[\sup_{\delta\in \mathbb{R}} \left(\Gamma_\delta+\Gamma'_{\tau-\delta}\right)\right].
        \end{equation}
    \end{enumerate}
\end{definition}

\begin{lemma}\label{lma:testcurvplus}
    Let $\Gamma \in \TC(X,\theta)_{>0}$ and $\Gamma' \in \TC(X,\theta')_{>0}$. Then we have $\Gamma+\Gamma'\in \TC(X,\theta+\theta')_{>0}$. In particular, for any $\tau<(\Gamma+\Gamma')_{\max}$,
    \[
    \sup_{\delta\in \mathbb{R}} \left(\Gamma_\delta+\Gamma'_{\tau-\delta}\right)\in \PSH(X,\theta+\theta').
    \]
    This potential is $\mathcal{I}$-good if $\Gamma \in \PSH^{\NA}(X,\theta)_{>0}$ and $\Gamma' \in \PSH^{\NA}(X,\theta')_{>0}$.

    In particular, \eqref{eq:GammaGammapsum} in \cref{def:sumtestcur} makes sense.
\end{lemma}

\begin{proof}
    Let
    \[
    \eta_{\tau}=\sup_{\delta\in \mathbb{R}} \left(\Gamma_\delta+\Gamma'_{\tau-\delta}\right)=\sup_{\tau-\Gamma'_{\max}<\delta<\Gamma_{\max}}\left(\Gamma_\delta+\Gamma'_{\tau-\delta}\right)
    \]
    for all $\tau\in \mathbb{R}$.
    Set
    \[
    Z=\left\{x\in X: \Gamma_{\upsilon}(x)=-\infty \,\forall \upsilon\in \mathbb{R}\right\}\cup \left\{x\in X: \Gamma'_{\upsilon}(x)=-\infty \,\forall \upsilon\in \mathbb{R}\right\}.
    \]
    It follows from \cref{prop:sumLegendre} that for any $x\in X\setminus Z$, we have
    \[
    \eta^*_t(x)=\Gamma^*_t(x)+\Gamma'^*_t(x)
    \]
    for all $t>0$.
    The same trivially holds when $x\in Z$, so the equation holds everywhere. In particular, by \cref{cor:Leginvpartdef} and \cref{prop:Kis}, we have
    \[
    \eta_{\tau}=\left(\Gamma^*+\Gamma'^*\right)^*_{\tau}\in \PSH(X,\theta+\theta')
    \]
    when $\tau<\Gamma_{\max}+\Gamma'_{\max}$.

    We now check that $\Gamma+\Gamma'$ is a test curve. The family $\tau\mapsto\eta_\tau$ is decreasing, and the monotonicity of the $P$-envelope \itemref{itm:prop:Pconc:3} gives the decreasingness of $\Gamma+\Gamma'$. If $\tau=(1-s)\tau_0+s\tau_1$, then for arbitrary $\delta_0,\delta_1\in \mathbb{R}$ and $\delta=(1-s)\delta_0+s\delta_1$,
    \[
    \Gamma_\delta+\Gamma'_{\tau-\delta}\geq (1-s)(\Gamma_{\delta_0}+\Gamma'_{\tau_0-\delta_0})+s(\Gamma_{\delta_1}+\Gamma'_{\tau_1-\delta_1}).
    \]
    Taking suprema and using \cref{prop:Pconc} gives the concavity of $\Gamma+\Gamma'$.

    Fix $\tau<\Gamma_{\max}+\Gamma'_{\max}$. Choose $\delta\in \mathbb{R}$ with $\delta<\Gamma_{\max}$ and $\tau-\delta<\Gamma'_{\max}$. Then
    \[
    \eta_\tau\geq \Gamma_\delta+\Gamma'_{\tau-\delta}.
    \]
    By \cref{lma:testcurvposmass}, $\Gamma_\delta\in \PSH(X,\theta)_{>0}$ and $\Gamma'_{\tau-\delta}\in \PSH(X,\theta')_{>0}$. Hence $\Gamma_\delta+\Gamma'_{\tau-\delta}$ has positive $(\theta+\theta')$-mass by the multi-linearity of non-pluripolar products \itemref{itm:prop:npp1:6}. It follows from the monotonicity theorem \cref{thm:mono} that $\eta_\tau\in \PSH(X,\theta+\theta')_{>0}$. By \cref{cor:Psendspotentialtomodel},
    \[
        (\Gamma+\Gamma')_\tau=P_{\theta+\theta'}[\eta_\tau]
    \]
    is a model potential. Finally, $(\Gamma+\Gamma')_{-\infty}\in \PSH(X,\theta+\theta')$ by \cref{prop:pshfunction_closedseq}; since it dominates $(\Gamma+\Gamma')_\tau$ for any $\tau<\Gamma_{\max}+\Gamma'_{\max}$, it has positive mass by another application of the monotonicity theorem \cref{thm:mono}. Thus $\Gamma+\Gamma'\in \TC(X,\theta+\theta')_{>0}$.

    Next, assume that $\Gamma$ and $\Gamma'$ are $\mathcal{I}$-model. Fix $\tau<\Gamma_{\max}+\Gamma'_{\max}$. Then for each $\delta\in \mathbb{R}$ such that $\delta<\Gamma_{\max}$ and $\tau-\delta<\Gamma'_{\max}$, we know that $\Gamma_\delta\in \PSH(X,\theta)_{>0}$ and $\Gamma'_{\tau-\delta}\in \PSH(X,\theta')_{>0}$ by \cref{lma:testcurvposmass}.
    It follows from \cref{ex:ImodelIgood} that $\Gamma_\delta$ and $\Gamma'_{\tau-\delta}$ are both $\mathcal{I}$-good, hence so is $\Gamma_\delta+\Gamma'_{\tau-\delta}\in \PSH(X,\theta+\theta')_{>0}$ by \cref{prop:Igoodlinear}. These sums are bounded above by $0$, so $\eta_{\tau}$ is $\mathcal{I}$-good by \cref{prop:Igoodsup}. Since $(\Gamma+\Gamma')_\tau=P_{\theta+\theta'}[\eta_\tau]$, we conclude that $\Gamma+\Gamma'$ is $\mathcal{I}$-model.
\end{proof}

\begin{proposition}\label{prop:testcurvesumproperty}
    Let $\Gamma \in \TC(X,\theta)_{>0}$ and $\Gamma' \in \TC(X,\theta')_{>0}$. Then $\Gamma+\Gamma'\in \TC(X,\theta+\theta')_{>0}$. Moreover,
    \begin{equation}\label{eq:sumGammaGammap}
    (\Gamma+\Gamma')_{-\infty}=P_{\theta+\theta'}\left[\Gamma_{-\infty}+\Gamma'_{-\infty}\right].
    \end{equation}

    When $\Gamma\in \PSH^{\NA}(X,\theta)_{>0}$ and $\Gamma'\in \PSH^{\NA}(X,\theta')_{>0}$, we have $\Gamma+\Gamma'\in \PSH^{\NA}(X,\theta+\theta')_{>0}$.

    The operation $+$ is commutative and associative.
\end{proposition}
\begin{proof}
    It follows immediately from \cref{lma:testcurvplus} that $\Gamma+\Gamma'\in \TC(X,\theta+\theta')_{>0}$. It lies in $\PSH^{\NA}(X,\theta+\theta')_{>0}$ if $\Gamma\in \PSH^{\NA}(X,\theta)_{>0}$ and $\Gamma'\in \PSH^{\NA}(X,\theta')_{>0}$.

    We argue \eqref{eq:sumGammaGammap}. By definition, for any small enough $\tau$, we have
    \[
    (\Gamma+\Gamma')_{-\infty}\geq (\Gamma+\Gamma')_{2\tau}\succeq_P \Gamma_{\tau}+\Gamma'_{\tau}.
    \]
    Letting $\tau\to -\infty$ and applying \cref{prop:dsconvpresorder} and \cref{thm:dSadditivity}, we find that
    \[
    (\Gamma+\Gamma')_{-\infty}\succeq_P \Gamma_{-\infty}+\Gamma'_{-\infty}.
    \]
    On the other hand, for each small enough $\tau$, we have
    \[
    (\Gamma+\Gamma')_{\tau}\sim_P \sup_{\delta\in \mathbb{R}}\left( \Gamma_\delta+\Gamma'_{\tau-\delta}\right) \preceq_P \Gamma_{-\infty}+\Gamma'_{-\infty}
    \]
    by \cref{prop:Ppartialsum} and \cref{prop:dsconvpresorder}. We apply \cref{prop:dsconvpresorder} again, we conclude that
    \[
     (\Gamma+\Gamma')_{-\infty} \preceq_P \Gamma_{-\infty}+\Gamma'_{-\infty}.
    \]
    So \eqref{eq:sumGammaGammap} follows.

    Finally, let us show that $+$ is commutative and associative. Commutativity follows from the symmetry of the defining supremum. Let $\Gamma''\in \TC(X,\theta'')_{>0}$. Then we want to show that
    \[
    (\Gamma+\Gamma')+\Gamma''=\Gamma+(\Gamma'+\Gamma'').
    \]
    First observe that
    \[
    \left((\Gamma+\Gamma')+\Gamma'' \right)_{\max}=\left( \Gamma+(\Gamma'+\Gamma'')\right)_{\max}.
    \]
    Fix $\tau$ less than this common value. We compute that
    \[
    \begin{aligned}
    &\left( (\Gamma+\Gamma')+\Gamma''\right)_{\tau}\\
    =& P_{\theta+\theta'+\theta''}\left[\sup_{\delta_1\in \mathbb{R}} \left((\Gamma+\Gamma')_{\delta_1}+\Gamma''_{\tau-\delta_1}\right)\right]\\
    \sim_P & \sup_{\delta_1\in \mathbb{R}} \left((\Gamma+\Gamma')_{\delta_1}+\Gamma''_{\tau-\delta_1}\right)\\
    \sim_P & \sup_{\delta_1,\delta_2\in \mathbb{R}} \left(\Gamma_{\delta_2}+\Gamma'_{\delta_1-\delta_2}+\Gamma''_{\tau-\delta_1}\right),
    \end{aligned}
    \]
    where in the last line, we applied \cref{prop:dsconvpresorder} and \cref{prop:Ppartialsum}. Similarly, for $\left( \Gamma+(\Gamma'+\Gamma'')\right)_{\tau}$, we get the same expression. The associativity follows.
\end{proof}

\begin{lemma}\label{lma:testcursumcomp}
    Let $\Gamma \in \TC(X,\theta)_{>0}$ and $\Gamma' \in \TC(X,\theta')_{>0}$. Then for any closed smooth positive $(1,1)$-forms $\omega$ and $\omega'$ on $X$, we have
    \[
    P_{\theta+\theta'+\omega+\omega'}[\Gamma+\Gamma']=P_{\theta+\omega}[\Gamma]+P_{\theta'+\omega'}[\Gamma'].
    \]
\end{lemma}
\begin{proof}
    Observe that
    \[
    \begin{aligned}
    P_{\theta+\theta'+\omega+\omega'}[\Gamma+\Gamma']_{\max}=&\left(P_{\theta+\omega}[\Gamma]+P_{\theta'+\omega'}[\Gamma'] \right)_{\max}\\
    =&\Gamma_{\max}+\Gamma'_{\max}.
    \end{aligned}
    \]
    Take $\tau\in \mathbb{R}$ less than this common value. We need to verify that
    \[
    (\Gamma+\Gamma')_{\tau}\sim_P \Bigl( P_{\theta+\omega}[\Gamma]+P_{\theta'+\omega'}[\Gamma']\Bigr)_{\tau}.
    \]
    By definition, this means that
    \[
    \sup_{\tau-\Gamma'_{\max}<\delta<\Gamma_{\max}}\left(\Gamma_{\delta}+\Gamma'_{\tau-\delta}\right)\sim_P \sup_{\tau-\Gamma'_{\max}<\delta<\Gamma_{\max}}\left(P_{\theta+\omega}[\Gamma_{\delta}]+P_{\theta'+\omega'}[\Gamma'_{\tau-\delta}]\right).
    \]
    This is a consequence of \cref{prop:Ppartialsum} and \cref{prop:Ppartialsup}.
\end{proof}

\begin{definition}\label{def:testcurvplusC}
    Let $\Gamma\in \TC(X,\theta)_{>0}$ and $C\in \mathbb{R}$, we define $\Gamma+C\in \TC(X,\theta)_{>0}$ as follows:
    \index{test curve!translation}
    \begin{enumerate}
        \item\label{itm:def:testcurvplusC:1} We set
        \[
            (\Gamma+C)_{\max}\coloneqq \Gamma_{\max}+C;
        \]
        \item\label{itm:def:testcurvplusC:2} for any $\tau<(\Gamma+C)_{\max}$, we set
        \[
        (\Gamma+C)_{\tau}\coloneqq \Gamma_{\tau-C}.
        \]
    \end{enumerate}
\end{definition}
It is obvious that if $\Gamma\in \PSH^{\NA}(X,\theta)_{>0}$, then so is $\Gamma+C$.

\begin{proposition}\label{prop:testcurveplusC}
    Let $\Gamma\in \TC(X,\theta)_{>0}$, $\Gamma'\in \TC(X,\theta')_{>0}$ and $C,C'\in \mathbb{R}$. Then
    \begin{enumerate}
        \item\label{itm:prop:testcurveplusC:1} $(\Gamma+\Gamma')+C=\Gamma+(\Gamma'+C)=(\Gamma+C)+\Gamma'$;
        \item\label{itm:prop:testcurveplusC:2} $\Gamma+(C+C')=(\Gamma+C)+C'$.
    \end{enumerate}
\end{proposition}
\begin{proof}
    \ref{itm:prop:testcurveplusC:1} We first observe that
    \[
    \left((\Gamma+\Gamma')+C\right)_{\max}=\left( \Gamma+(\Gamma'+C) \right)_{\max}=\left( (\Gamma+C)+\Gamma' \right)_{\max}=\Gamma_{\max}+\Gamma'_{\max}+C.
    \]
    Take any $\tau\in \mathbb{R}$ less than this common value. We compute
    \[
        \begin{aligned}
        \left((\Gamma+\Gamma')+C\right)_{\tau}=&(\Gamma+\Gamma')_{\tau-C}
        = P_{\theta+\theta'}\left[ \sup_{\delta\in \mathbb{R}}\left( \Gamma_\delta+\Gamma'_{\tau-C-\delta} \right) \right],\\
        \left( \Gamma+(\Gamma'+C) \right)_{\tau}=& P_{\theta+\theta'}\left[ \sup_{\delta\in \mathbb{R}}\left( \Gamma_\delta+(\Gamma'+C)_{\tau-\delta} \right) \right]\\
        =& P_{\theta+\theta'}\left[ \sup_{\delta\in \mathbb{R}}\left( \Gamma_\delta+\Gamma'_{\tau-C-\delta} \right) \right],\\
        \left((\Gamma+C)+\Gamma'\right)_{\tau}=& P_{\theta+\theta'}\left[ \sup_{\delta\in \mathbb{R}}\left( (\Gamma+C)_{C+\delta}+\Gamma'_{\tau-C-\delta} \right) \right]\\
        =& P_{\theta+\theta'}\left[ \sup_{\delta\in \mathbb{R}}\left( \Gamma_\delta+\Gamma'_{\tau-C-\delta} \right) \right].
        \end{aligned}
    \]

    \ref{itm:prop:testcurveplusC:2} Observe that
    \[
        \left( \Gamma+(C+C') \right)_{\max}=\left( (\Gamma+C)+C' \right)_{\max}=\Gamma_{\max}+C+C'.
    \]
    For any $\tau\in \mathbb{R}$ less than this value, we have
    \[
        \left( \Gamma+(C+C') \right)_{\tau}=\Gamma_{\tau-C-C'}=\left( (\Gamma+C)+C' \right)_{\tau}.
    \]
    This completes the proof.
\end{proof}

\begin{definition}\label{def:testcurlor}
    Let $\Gamma,\Gamma'\in \TC(X,\theta)_{>0}$. We define $\Gamma\lor \Gamma'\in \TC(X,\theta)_{>0}$ as follows:
    \index{test curve!maximum}
    \begin{enumerate}
        \item\label{itm:def:testcurlor:1} We set
        \[
            (\Gamma\lor \Gamma')_{\max}\coloneqq \Gamma_{\max}\lor \Gamma'_{\max};
        \]
        \item\label{itm:def:testcurlor:2} for any $\tau<(\Gamma\lor \Gamma')_{\max}$, we define
        \begin{equation}\label{eq:testcurlordef}
            (\Gamma\lor \Gamma')_{\tau}\coloneqq P_{\theta}\left[ \LCE\left( \rho \mapsto \Gamma_{\rho}\lor \Gamma'_{\rho} \right)_{\tau} \right].
        \end{equation}
    \end{enumerate}
\end{definition}
Recall that the upper concave envelope $\LCE$ is defined in \cref{def:LCE}\footnote{In \cref{def:LCE}, we define the convex analogue, the \emph{lower convex envelope}. The concave version is obtained by applying this definition to negatives.}. By construction, we have $\Gamma\lor \Gamma'\geq \Gamma$ and $\Gamma\lor \Gamma'\geq \Gamma'$.

On the right-hand side of \eqref{eq:testcurlordef}, we first take the upper concave envelope $\LCE\left( \rho \mapsto \Gamma_{\rho}\lor \Gamma'_{\rho} \right)$, this gives a family of functions indexed by $\rho\in \mathbb{R}$. Then the notation $\LCE\left( \rho \mapsto \Gamma_{\rho}\lor \Gamma'_{\rho} \right)_{\tau}$ denotes the function obtained by evaluating this family at $\tau$.

\begin{lemma}\label{lma:testcurlor}
    Let $\Gamma,\Gamma'\in \TC(X,\theta)_{>0}$. Then for any $\tau<\Gamma_{\max}\lor \Gamma'_{\max}$, we have
    \[
        \LCE\left( \rho \mapsto \Gamma_{\rho}\lor \Gamma'_{\rho} \right)_{\tau}\in \PSH(X,\theta),
    \]
    and $\Gamma\lor \Gamma'\in \TC(X,\theta)_{>0}$.

    This potential is $\mathcal{I}$-good if $\Gamma,\Gamma'\in \PSH^{\NA}(X,\theta)_{>0}$.

    In particular, \eqref{eq:testcurlordef} in \cref{def:testcurlor} makes sense.
\end{lemma}

\begin{proof}
    To simplify the notation, we write
    \[
    \psi_{\tau}=\LCE\left( \rho \mapsto \Gamma_{\rho}\lor \Gamma'_{\rho} \right)_{\tau}
    \]
    for all $\tau\in \mathbb{R}$. Thanks to \cref{prop:Legendretranssup}, for any fixed $x\in X$, we have
    \begin{equation}\label{eq:psistartemp1}
    \psi^*_t(x)=\Gamma^*_t(x)\lor \Gamma'^*_t(x)
    \end{equation}
    for all $t>0$ as long as $\Gamma_{\tau}(x)\neq -\infty$ and $\Gamma'_{\tau}(x)\neq -\infty$ for some $\tau\in \mathbb{R}$. Otherwise, assume for instance that $\Gamma_{\tau}(x)=-\infty$ for all $\tau\in \mathbb{R}$; then by definition, $\psi_{\tau}(x)=\Gamma'_{\tau}(x)$ for all $\tau\in \mathbb{R}$. Therefore, $\Gamma^*_t(x)=-\infty$ for all $t>0$ and \eqref{eq:psistartemp1} continues to hold. The case where $\Gamma'_{\tau}(x)=-\infty$ for all $\tau$ is symmetric. Therefore, we have shown that
    \[
    \psi^*_t=\Gamma^*_t\lor \Gamma'^*_t\in \PSH(X,\theta).
    \]
    It follows from \cref{prop:maxsubgeod} that $(\psi^*_t)_{t\in [a,b]}$ is a subgeodesic for any $0< a<b$.

    Next we observe that $\psi_{\bullet}$ is closed by definition. So it follows from \cref{prop:Legendretranssup} and \cref{prop:Kis} that
    \[
    \psi_{\tau}=(\psi^*_{\bullet})^*_{\tau}\in \PSH(X,\theta)\cup \{-\infty\}.
    \]

    Due to \cref{prop:Gammaclosed} and \cref{prop:concavhull}, there is a pluripolar set $Z\subseteq X$ such that for $x\in X\setminus Z$, we have
    \[
    \begin{split}
    \psi_{\tau}(x)=\sup\left\{\lambda \Gamma_{\rho}(x)+(1-\lambda)\Gamma'_{\rho'}(x): \lambda \in (0,1),\rho,\rho'\in \mathbb{R},\lambda\rho+(1-\lambda)\rho'=\tau\right\}
    \end{split}
    \]
    for all $\tau<\Gamma_{\max}\lor \Gamma'_{\max}$. It follows from \cref{prop:pshfuncdetdense} that
    \begin{equation}\label{eq:psitausupslineartemp}
    \begin{split}
    \psi_{\tau}=\uscsup \left\{\lambda \Gamma_{\rho}+(1-\lambda)\Gamma'_{\rho'}: \lambda \in (0,1),\rho,\rho'\in \mathbb{R},\lambda\rho+(1-\lambda)\rho'=\tau \right\}
    \end{split}
    \end{equation}
    for all $\tau<\Gamma_{\max}\lor \Gamma'_{\max}$.

    We now check that $\Gamma\lor\Gamma'$ is a test curve. The function $\rho\mapsto\Gamma_\rho\lor\Gamma'_\rho$ is decreasing, hence its upper concave envelope $\psi_\bullet$ is decreasing. The monotonicity of the $P$-envelope \itemref{itm:prop:Pconc:3} gives the decreasingness of $\Gamma\lor\Gamma'$. Since $\psi_\bullet$ is concave by construction, \cref{prop:Pconc} gives the concavity of $\Gamma\lor\Gamma'$.

    Fix $\tau<\Gamma_{\max}\lor \Gamma'_{\max}$. Then either $\tau<\Gamma_{\max}$ or $\tau<\Gamma'_{\max}$. Since $\psi_\bullet$ is a concave majorant of $\rho\mapsto\Gamma_\rho\lor\Gamma'_\rho$, we have
    \[
    \psi_\tau\geq \Gamma_\tau
    \quad\text{if } \tau<\Gamma_{\max},
    \qquad
    \psi_\tau\geq \Gamma'_\tau
    \quad\text{if } \tau<\Gamma'_{\max}.
    \]
    Together with the fact that $\psi_\tau\in \PSH(X,\theta)\cup\{-\infty\}$, these inequalities show that $\psi_\tau\in \PSH(X,\theta)$. It follows from \cref{lma:testcurvposmass} and the monotonicity theorem \cref{thm:mono} that $\psi_\tau\in \PSH(X,\theta)_{>0}$. By \cref{cor:Psendspotentialtomodel},
    \[
        (\Gamma\lor\Gamma')_\tau=P_\theta[\psi_\tau]
    \]
    is a model potential. Finally, $(\Gamma\lor\Gamma')_{-\infty}\in \PSH(X,\theta)$ by \cref{prop:pshfunction_closedseq}; since it dominates $(\Gamma\lor\Gamma')_\tau$ for any $\tau<\Gamma_{\max}\lor\Gamma'_{\max}$, it has positive mass by another application of the monotonicity theorem \cref{thm:mono}. Thus $\Gamma\lor\Gamma'\in\TC(X,\theta)_{>0}$.

    It follows from \eqref{eq:psitausupslineartemp} that $\psi_{\tau}$ is $\mathcal{I}$-good if $\Gamma,\Gamma'\in \PSH^{\NA}(X,\theta)_{>0}$, thanks to \cref{prop:Igoodlinear} and \cref{prop:Igoodsup}; the family in the upper envelope is bounded above by $0$.
    Since $(\Gamma\lor\Gamma')_\tau=P_\theta[\psi_\tau]$, we conclude that $\Gamma\lor\Gamma'$ is $\mathcal{I}$-model in this case.
\end{proof}

\begin{corollary}\label{cor:testcurvlorprop}
    Let $\Gamma,\Gamma'\in \TC(X,\theta)_{>0}$. Then $\Gamma\lor \Gamma'\in\TC(X,\theta)_{>0}$ and
    \begin{equation}\label{eq:GammalorGammapminfty}
    (\Gamma\lor \Gamma')_{-\infty}=P_{\theta}\left[\Gamma_{-\infty}\lor \Gamma'_{-\infty}\right].
    \end{equation}
    If $\Gamma,\Gamma'\in \PSH^{\NA}(X,\theta)_{>0}$, then $\Gamma\lor \Gamma'\in\PSH^{\NA}(X,\theta)_{>0}$.

    For each $\Gamma''\in \TC(X,\theta)_{>0}$ and each $\Gamma''\geq \Gamma$ and $\Gamma''\geq \Gamma'$, we have $\Gamma''\geq \Gamma\lor \Gamma'$.

    Moreover, the operation $\lor$ is associative and commutative.
\end{corollary}
In particular, given a finite family $\{\Gamma_i\}_{i\in I}$ in $\TC(X,\theta)_{>0}$, we can define
\[
\bigvee_{i\in I}\Gamma_i
\]
without ambiguity.
\begin{proof}
    It follows immediately from \cref{lma:testcurlor} that $\Gamma\lor \Gamma'\in \TC(X,\theta)_{>0}$, and it lies in $\PSH^{\NA}(X,\theta)_{>0}$ if $\Gamma,\Gamma'\in \PSH^{\NA}(X,\theta)_{>0}$.

    We argue \eqref{eq:GammalorGammapminfty}. Set
    \[
        \eta\coloneqq \Gamma_{-\infty}\lor \Gamma'_{-\infty}.
    \]
    Since $\Gamma\lor\Gamma'\geq \Gamma$ and $\Gamma\lor\Gamma'\geq \Gamma'$, we have $(\Gamma\lor\Gamma')_{-\infty}\geq \eta$. As $(\Gamma\lor\Gamma')_{-\infty}$ is a model potential, it follows that
    \[
        (\Gamma\lor\Gamma')_{-\infty}\geq P_\theta[\eta].
    \]
    Conversely, for each $\rho\in \mathbb{R}$ we have $\Gamma_\rho\lor\Gamma'_\rho\leq \eta$. Hence the constant curve with value $\eta$ is a concave majorant of $\rho\mapsto \Gamma_\rho\lor\Gamma'_\rho$, and therefore
    \[
        \LCE\left(\rho\mapsto\Gamma_\rho\lor\Gamma'_\rho\right)_\tau\leq \eta,\quad \tau<\Gamma_{\max}\lor\Gamma'_{\max}.
    \]
    Applying the monotonicity of the $P$-envelope and then letting $\tau\to -\infty$ gives
    \[
        (\Gamma\lor\Gamma')_{-\infty}\leq P_\theta[\eta].
    \]
    This proves \eqref{eq:GammalorGammapminfty}.

    Take $\Gamma''$ as in the statement of the corollary.
    First observe that
    \[
        \Gamma''_{\max}\geq \Gamma_{\max} \lor \Gamma'_{\max}=(\Gamma\lor \Gamma')_{\max}.
    \]
    Take $\tau<(\Gamma\lor \Gamma')_{\max}$. We argue that
    \[
    \Gamma''_{\tau}\geq (\Gamma\lor \Gamma')_{\tau}.
    \]
    By the concavity of $\Gamma''$, the potential $\Gamma''_\tau$ dominates the concave envelope
    \[
    \LCE\left(\rho\mapsto \Gamma_{\rho}\lor \Gamma'_{\rho}\right)_\tau.
    \]
    Since $\Gamma''_\tau$ is model, the monotonicity of the $P$-envelope gives
    \[
    \Gamma''_{\tau}\geq (\Gamma\lor \Gamma')_{\tau}.
    \]
    Therefore,
    \[
    \Gamma''\geq \Gamma\lor \Gamma'.
    \]

    Commutativity is immediate from the definition. Associativity follows from the least-upper-bound property just proved: both $(\Gamma\lor\Gamma')\lor\Gamma''$ and $\Gamma\lor(\Gamma'\lor\Gamma'')$ are the least upper bound of $\Gamma,\Gamma',\Gamma''$.
\end{proof}

\begin{lemma}\label{lma:testcurlorPthetapomega}
    Let $\Gamma,\Gamma'\in \TC(X,\theta)_{>0}$ and $\omega$ be a closed smooth positive $(1,1)$-form on $X$. Then
    \[
    P_{\theta+\omega}[\Gamma\lor \Gamma']=P_{\theta+\omega}[\Gamma]\lor P_{\theta+\omega}[\Gamma'].
    \]
\end{lemma}
\begin{proof}
    We first observe that
    \[
    \left(P_{\theta+\omega}[\Gamma\lor \Gamma']\right)_{\max}=\left( P_{\theta+\omega}[\Gamma]\lor P_{\theta+\omega}[\Gamma']\right)_{\max}=\Gamma_{\max}\lor \Gamma'_{\max}.
    \]
    Let $\tau\in \mathbb{R}$ be less than this common value. We need to show that
    \[
    (\Gamma\lor \Gamma')_{\tau}\sim_P \left(P_{\theta+\omega}[\Gamma]\lor P_{\theta+\omega}[\Gamma'] \right)_{\tau}.
    \]

    We need the formula \eqref{eq:psitausupslineartemp} proved in the proof of \cref{lma:testcurlor}:
    \[
    \begin{split}
    (\Gamma\lor \Gamma')_{\tau}\sim_P\uscsup\left\{\lambda \Gamma_{\rho}+(1-\lambda)\Gamma'_{\rho'}: \lambda \in (0,1),\rho,\rho'\in \mathbb{R},\lambda\rho+(1-\lambda)\rho'=\tau \right\}.
    \end{split}
    \]
    A similar result holds with $P_{\theta+\omega}[\Gamma]$ and $P_{\theta+\omega}[\Gamma']$ in place of $\Gamma$ and $\Gamma'$. So our assertion is a direct consequence of \cref{prop:Ppartialsum} and \cref{prop:Ppartialsup}.
\end{proof}

\begin{definition}\label{def:testcursup}
    Let $(\Gamma^i)_{i\in I}$ be an increasing net in $\TC(X,\theta)_{>0}$. Assume that
    \index{test curve!usc supremum}
    \begin{equation}\label{eq:Gammaisupfinite1}
        \sup_{i\in I}\Gamma^{i}_{\max}<\infty.
    \end{equation}
    Then we define $\uscsup[i\in I] \Gamma^{i}\in \TC(X,\theta)_{>0}$ as follows:
    \begin{enumerate}
        \item\label{itm:def:testcursup:1} We set
        \[
            \left( \uscsup[i\in I] \Gamma^{i} \right)_{\max}=\sup_{i\in I}\Gamma^{i}_{\max};
        \]
        \item\label{itm:def:testcursup:2} for any $\tau< \sup_{i\in I}\Gamma^{i}_{\max}$, we let
        \[
        \left( \uscsup[i\in I] \Gamma^{i} \right)_{\tau}\coloneqq \uscsup[i\in I] \Gamma^{i}_{\tau}.
        \]
    \end{enumerate}
\end{definition}

\begin{proposition}\label{prop:supsincnetteestcur}
    Let $(\Gamma^i)_{i\in I}$ be an increasing net in $\TC(X,\theta)_{>0}$ satisfying \eqref{eq:Gammaisupfinite1}. Then $\uscsup[i\in I] \Gamma^{i}$ as defined in \cref{def:testcursup} lies in $\TC(X,\theta)_{>0}$. Moreover, if $\Gamma^i\in \PSH^{\NA}(X,\theta)_{>0}$ for all $i\in I$, then $\uscsup[i\in I] \Gamma^{i}$ lies in $\PSH^{\NA}(X,\theta)_{>0}$ as well.

    Moreover, we have
    \begin{equation}\label{eq:Gammiminf}
        \left( \uscsup[i\in I] \Gamma^{i} \right)_{-\infty}=\uscsup[i\in I]\Gamma^i_{-\infty}.
    \end{equation}
\end{proposition}
\begin{proof}
    Fix $\tau<\sup_i\Gamma^i_{\max}$. Choose $i_0$ such that $\tau<\Gamma^{i_0}_{\max}$. On the tail $i\geq i_0$, the potentials $\Gamma^i_\tau$ form an increasing net of model potentials of positive mass, uniformly bounded above by $0$. Hence \cref{prop:incnetmodel} shows that
    \[
        \uscsup[i\in I]\Gamma^i_\tau = \uscsup[i\geq i_0]\Gamma^i_\tau
    \]
    is model. Its mass is positive by monotonicity, since it dominates $\Gamma^{i_0}_\tau$. The decreasingness and concavity in $\tau$ follow from the directedness of the net and the corresponding properties of the curves $\Gamma^i$. Thus $\uscsup_i\Gamma^i$ is a test curve. If all $\Gamma^i$ are $\mathcal I$-model, the same tail argument and \cref{prop:incnetmodelPI} show that the supremum is $\mathcal I$-model.

    It remains to argue \eqref{eq:Gammiminf}. Fix any $i_0\in I$. Since the curves are decreasing and are extended as in \cref{rmk:extendtestcur}, the value at $-\infty$ can be computed by taking the supremum over rational $\tau<\Gamma^{i_0}_{\max}$.

    By \cref{prop:supsupstardiff}, there is a pluripolar set $Z\subseteq X$ such that for any $x\in X\setminus Z$,
    \[
    \left( \uscsup[i\in I] \Gamma^{i} \right)_{-\infty}(x)=\sup_{\mathbb{Q}\ni \tau< \Gamma^{i_0}_{\max}}\left(\uscsup[i\in I]\Gamma^i_{\tau}\right)(x)=\sup_{\mathbb{Q}\ni \tau< \Gamma^{i_0}_{\max},i\in I}\Gamma^i_{\tau}(x)=\sup_{i\in I}\Gamma^i_{-\infty}(x).
    \]
    So they are equal everywhere by \cref{prop:pshfuncdetdense}.
\end{proof}

\begin{lemma}\label{lma:suptestcurvcompatible}
    Let $(\Gamma^i)_{i\in I}$ be an increasing net in $\TC(X,\theta)_{>0}$ satisfying \eqref{eq:Gammaisupfinite1}. Assume that $\omega$ is a closed smooth positive $(1,1)$-form on $X$. Then
    \[
    P_{\theta+\omega}\left[ \uscsup[i\in I] \Gamma^{i} \right]= \uscsup[i\in I] P_{\theta+\omega}\left[ \Gamma^{i} \right].
    \]
\end{lemma}
\begin{proof}
    Observe that
    \[
    \left( P_{\theta+\omega}\left[ \uscsup[i\in I] \Gamma^{i} \right] \right)_{\max}= \left(\uscsup[i\in I] P_{\theta+\omega}\left[ \Gamma^{i} \right]\right)_{\max}=\sup_{i\in I}\Gamma^i_{\max}.
    \]
    Fix $\tau\in \mathbb{R}$ less than this common value.

    It suffices to show that
    \[
        \left( \uscsup[i\in I] \Gamma^{i} \right)_{\tau}\sim_P \left( \uscsup[i\in I] P_{\theta+\omega}\left[\Gamma^{i}\right] \right)_{\tau}.
    \]
    This is an immediate consequence of \cref{prop:Ppartialsup}.
\end{proof}

\begin{definition}\label{def:testcurvsupsgeneral}
    Let $(\Gamma^i)_{i\in I}$ be a non-empty family in $\TC(X,\theta)_{>0}$ satisfying \eqref{eq:Gammaisupfinite1}. Let $\Fin(I)$ denote the net of non-empty finite subsets of $I$, ordered by inclusion. Then we define
    \index{test curve!usc supremum}
    \begin{equation}\label{eq:generalsupstestcurv}
    \uscsup[i\in I]\Gamma^i\coloneqq \uscsup[J\in \Fin(I)] \left( \bigvee_{j\in J} \Gamma^j\right).
    \end{equation}
\end{definition}

Observe that by \cref{def:testcurlor}, we have
\[
\sup_{J\in \Fin(I)} \left( \bigvee_{j\in J} \Gamma^j\right)_{\max}=\sup_{i\in I}\Gamma^i_{\max}<\infty.
\]
So \eqref{eq:generalsupstestcurv} makes sense. In particular,
\begin{equation}\label{eq:testcursupmax}
    \left(\uscsup[i\in I]\Gamma^i\right)_{\max}=\sup_{i\in I}\Gamma^i_{\max}.
\end{equation}

It is clear that \cref{def:testcurvsupsgeneral} extends both \cref{def:testcursup} and \cref{def:testcurlor}.

\begin{proposition}\label{prop:suptcistc}
    Let $(\Gamma^i)_{i\in I}$ be a non-empty family in $\TC(X,\theta)_{>0}$ satisfying \eqref{eq:Gammaisupfinite1}. Then $\uscsup[i\in I]\Gamma^i\in \TC(X,\theta)_{>0}$. Moreover, if $\Gamma^i\in \PSH^{\NA}(X,\theta)_{>0}$ for all $i\in I$, then so is $\uscsup[i\in I]\Gamma^i$.

    Finally, we have
    \begin{equation}\label{eq:supsminfty}
        \left( \uscsup[i\in I]\Gamma^i \right)_{-\infty}=P_{\theta}\left[ \uscsup[i\in I]\Gamma^i_{-\infty} \right].
    \end{equation}
\end{proposition}
\begin{proof}
    The first assertion and the second follow from \cref{prop:supsincnetteestcur} and \cref{cor:testcurvlorprop}.

    It remains to argue \eqref{eq:supsminfty}. For this purpose, it suffices to show that
    \[
    \left( \uscsup[i\in I]\Gamma^i \right)_{-\infty}\sim_P \uscsup[i\in I]\Gamma^i_{-\infty}.
    \]
    For any $J\in \Fin(I)$, it follows from \cref{cor:testcurvlorprop} and \cref{prop:Ppartialsup} that
    \[
    \left(\bigvee_{j\in J}\Gamma^j\right)_{-\infty}\sim_P \bigvee_{j\in J}\Gamma^j_{-\infty}.
    \]
    From this, applying \cref{prop:incnetmodel}, \cref{prop:Ppartialsup} and \cref{prop:supsincnetteestcur}, we conclude our assertion.
\end{proof}

\begin{lemma}\label{lma:testcursupcompatible}
    Let $(\Gamma^i)_{i\in I}$ be a non-empty family in $\TC(X,\theta)_{>0}$ satisfying \eqref{eq:Gammaisupfinite1}. Assume that $\omega$ is a closed smooth positive $(1,1)$-form on $X$. Then
    \[
    P_{\theta+\omega}\left[ \uscsup[i\in I] \Gamma^{i} \right]= \uscsup[i\in I] P_{\theta+\omega}\left[ \Gamma^{i} \right].
    \]
\end{lemma}
\begin{proof}
    This is a direct consequence of \cref{lma:suptestcurvcompatible} and \cref{lma:testcurlorPthetapomega}.
\end{proof}

We prove a version of Choquet's lemma.
\begin{proposition}\label{prop:testcurvChoquet}
    Let $(\Gamma^i)_{i\in I}$ be a non-empty family in $\TC(X,\theta)_{>0}$ satisfying \eqref{eq:Gammaisupfinite1}. Then there is a countable subset $I'\subseteq I$ such that
    \[
    \uscsup[i\in I]\Gamma^i =\uscsup[i\in I']\Gamma^i.
    \]
\end{proposition}
\begin{proof}
    We may assume that $I$ is infinite.

    Set
    \[
        M\coloneqq \sup_{i\in I}\Gamma^i_{\max}.
    \]
    For each rational number $q<M$, apply Choquet's lemma \cref{prop:Choquet} to the locally uniformly bounded above family
    \[
        \left\{\left(\bigvee_{j\in J}\Gamma^j\right)_q:J\in\Fin(I)\right\}\subseteq \PSH(X,\theta).
    \]  
    We obtain a countable subfamily $\mathcal J_q\subseteq \Fin(I)$ such that
    \[
        \uscsup[J\in\Fin(I)]\left(\bigvee_{j\in J}\Gamma^j\right)_q = \uscsup[J\in\mathcal J_q]\left(\bigvee_{j\in J}\Gamma^j\right)_q.
    \]
    Let
    \[
        I'\coloneqq \bigcup_{q\in\mathbb Q, q<M} \bigcup_{J\in\mathcal J_q}J.
    \]
    Then $I'$ is countable and, for every rational $q<M$,
    \[
        \left(\uscsup[i\in I]\Gamma^i\right)_q= \left(\uscsup[i\in I']\Gamma^i\right)_q.
    \]
    Moreover, $\sup_{i\in I'}\Gamma^i_{\max}=M$, since the last equality is non-trivial for every rational $q<M$.
    Let $\Gamma'=\uscsup[i\in I']\Gamma^i$. Then clearly, $\Gamma'\leq \Gamma$. We claim that they are actually equal.

    Thanks to \cref{prop:Pequivchar2} and \cref{lma:testcurvposmass}, it suffices to show that for any $\tau<\uscsup[i\in I]\Gamma^i_{\max}$, we have
    \[
    \int_X \left(\theta+\ddc \Gamma'_{\tau}\right)^n=\int_X \left(\theta+\ddc \Gamma_{\tau}\right)^n.
    \]
    The equality is already known for all rational $\tau<M$. Since both
    \[
    \tau\mapsto \log\int_X\theta_{\Gamma'_\tau}^n,\quad \tau\mapsto \log\int_X\theta_{\Gamma_\tau}^n
    \]
    are continuous on $(-\infty,M)$ by \cref{prop:testcurvmasslogconc}, the equality of masses holds for every $\tau<M$.
\end{proof}

\begin{proposition}\label{prop:supGammiotherprop}
    Let $(\Gamma^i)_{i\in I}$ be a non-empty family in $\TC(X,\theta)_{>0}$ satisfying \eqref{eq:Gammaisupfinite1}. Let $C\in \mathbb{R}$. Then
    \[
    \uscsup[i\in I] (\Gamma^i+C)=\uscsup[i\in I] \Gamma^i+C.
    \]
    Suppose that $(\Gamma'^i)_{i\in I}$ is another family in $\TC(X,\theta')_{>0}$ satisfying \eqref{eq:Gammaisupfinite1}. Suppose that $\Gamma^i\leq \Gamma'^i$ for all $i\in I$. Then
    \[
    \uscsup[i\in I]\Gamma^i\leq \uscsup[i\in I]\Gamma'^i.
    \]
\end{proposition}
\begin{proof}
    This is immediate by definition.
\end{proof}

\begin{definition}\label{def:res}
    Let $\Gamma\in \TC(X,\theta)_{>0}$ and $\lambda>0$, we define $\lambda\Gamma\in \TC(X,\lambda\theta)_{>0}$
    \index{test curve!rescaling}
    as follows:
    \begin{enumerate}
        \item\label{itm:def:res:1} We set
        \[
        \left(\lambda\Gamma\right)_{\max}=\lambda \Gamma_{\max};
        \]
        \item\label{itm:def:res:2} for any $\tau< \left(\lambda\Gamma\right)_{\max}$, we set
        \[
        (\lambda\Gamma)_{\tau}=\lambda\Gamma_{\lambda^{-1}\tau}.
        \]
    \end{enumerate}
\end{definition}

\begin{proposition}\label{prop:testcurrescaling}
    Let $\Gamma\in \TC(X,\theta)_{>0}$ and $\lambda>0$. Then $\lambda\Gamma$ as defined in \cref{def:res} lies in $\TC(X,\lambda\theta)_{>0}$. Moreover, if $\Gamma\in \PSH^{\NA}(X,\theta)_{>0}$, then $\lambda\Gamma\in \PSH^{\NA}(X,\lambda\theta)_{>0}$.

    We have
    \begin{equation}
        (\lambda\Gamma)_{-\infty}=\lambda \Gamma_{-\infty}.
    \end{equation}
\end{proposition}
\begin{proof}
    This is immediate by definition.
\end{proof}

\begin{proposition}\label{prop:resclacompat}
    Let $\Gamma\in \TC(X,\theta)_{>0}$, $\Gamma'\in \TC(X,\theta')_{>0}$, $C\in \mathbb{R}$ and $\lambda,\lambda'>0$, we have
    \[
        \begin{aligned}
        \lambda(\Gamma+\Gamma')=&\lambda\Gamma+\lambda\Gamma',\\
        (\lambda\lambda')\Gamma=&\lambda(\lambda'\Gamma),\\
        \lambda(\Gamma+C)=&\lambda\Gamma+\lambda C.
        \end{aligned}
    \]
    Suppose that $(\Gamma^i)_{i\in I}$ is a non-empty family in $\TC(X,\theta)_{>0}$ satisfying \eqref{eq:Gammaisupfinite1}. Then
    \[
    \lambda \left( \uscsup[i\in I]\Gamma^i \right)=\uscsup[i\in I](\lambda \Gamma^i).
    \]

\end{proposition}
\begin{proof}
    This is immediate by definition.
\end{proof}

\begin{lemma}\label{lma:testcurvrescompatible}
    Let $\Gamma\in \TC(X,\theta)_{>0}$ and $\lambda>0$. Then for any closed smooth positive $(1,1)$-form $\omega$ on $X$, we have
    \[
    P_{\lambda\theta+\lambda\omega}[\lambda\Gamma]=\lambda P_{\theta+\omega}[\Gamma].
    \]
\end{lemma}
\begin{proof}
    This is clear by definition.
\end{proof}

\begin{definition}\label{def:testcurinf}
    Let $(\Gamma^i)_{i\in I}$ be a decreasing net in $\TC(X,\theta)_{>0}$. Assume that
    \index{test curve!infimum of a decreasing net}
    \begin{equation}\label{eq:Gammaiinffinite1}
        \inf_{i\in I}\Gamma^{i}_{\max}>-\infty,
    \end{equation}
    and
    \begin{equation}\label{eq:massnondeg}
    \inf_{i\in I} \int_X \left( \theta+\ddc \Gamma^i_{\tau} \right)^n>0,\quad \textup{for some }\tau< \inf_{i\in I}\Gamma^{i}_{\max}.
    \end{equation}
    Then we define $\inf_{i\in I} \Gamma^{i}\in \TC(X,\theta)_{>0}$ as follows:
    \begin{enumerate}
        \item\label{itm:def:testcurinf:1} We set
        \[
            \left( \inf_{i\in I} \Gamma^{i} \right)_{\max}=\inf_{i\in I}\Gamma^{i}_{\max};
        \]
        \item\label{itm:def:testcurinf:2} for any $\tau<  \left( \inf_{i\in I} \Gamma^{i} \right)_{\max}$, we let
        \[
        \left( \inf_{i\in I} \Gamma^{i} \right)_{\tau}\coloneqq \inf_{i\in I} \Gamma^{i}_{\tau}.
        \]
    \end{enumerate}
\end{definition}

\begin{proposition}\label{prop:inftcistc}
    Let $(\Gamma^i)_{i\in I}$ be a decreasing net in $\TC(X,\theta)_{>0}$ satisfying \eqref{eq:Gammaiinffinite1} and \eqref{eq:massnondeg}. Then $\inf_{i\in I} \Gamma^{i}\in \TC(X,\theta)_{>0}$.

    Moreover, if $\Gamma^i\in \PSH^{\NA}(X,\theta)_{>0}$ for all $i\in I$, then so is $\inf_{i\in I} \Gamma^{i}$.
\end{proposition}
\begin{proof}
    The first assertion is an immediate consequence of \cref{prop:decseqmodel}, except for the positivity of the mass at $-\infty$. Choose $\tau_0$ as in \eqref{eq:massnondeg}. Applying \cref{prop:vol_limit_model} gives
    \[
        \int_X \left(\theta+\ddc\inf_{i\in I}\Gamma^i_{\tau_0}\right)^n
        =\lim_{i\in I}\int_X\left(\theta+\ddc\Gamma^{i}_{\tau_0}\right)^n>0.
    \]
    Hence the value at $-\infty$ has positive mass by monotonicity. The last assertion follows from \cref{prop:decnetmodelPI}.
\end{proof}
In general, it is not true that
\[
\left( \inf_{i\in I}\Gamma^{i} \right)_{-\infty}=\inf_{i\in I}\Gamma^{i}_{-\infty}.
\]

\begin{lemma}\label{lma:testcurinfcompatible}
    Let $(\Gamma^i)_{i\in I}$ be a decreasing net in $\TC(X,\theta)_{>0}$ satisfying \eqref{eq:Gammaiinffinite1} and \eqref{eq:massnondeg}. Assume that $\omega$ is a closed smooth positive $(1,1)$-form on $X$. Then
    \[
    P_{\theta+\omega}\left[ \inf_{i\in I} \Gamma^{i} \right]= \inf_{i\in I} P_{\theta+\omega}\left[ \Gamma^{i} \right].
    \]
\end{lemma}
\begin{proof}
    First observe that
    \[
     \left( P_{\theta+\omega}\left[ \inf_{i\in I} \Gamma^{i} \right] \right)_{\max}= \left(\inf_{i\in I} P_{\theta+\omega}\left[ \Gamma^{i} \right]\right)_{\max}.
    \]
    Let $\tau\in \mathbb{R}$ be less than this common value. Then we need to show the following:
    \begin{equation}\label{eq:Pthetaomega_inf}
    P_{\theta+\omega}\left[ \inf_{i\in I} \Gamma^i_{\tau} \right]\sim_P \inf_{i\in I} P_{\theta+\omega}\left[  \Gamma^i_{\tau} \right].
    \end{equation}
    Applying \cref{prop:vol_limit_model} gives the mass convergence required in \cref{cor:decnetdS}. Hence $\Gamma^i_{\tau}\xrightarrow{d_S}\inf_{j\in I} \Gamma^j_{\tau} $. Thanks to \cref{cor:dSconv_changetheta} and \cref{cor:decnetdS}, we have
    \[
    P_{\theta+\omega}\left[ \Gamma^i_{\tau} \right]\xrightarrow{d_S} P_{\theta+\omega}\left[ \inf_{j\in I} \Gamma^j_{\tau} \right],\quad P_{\theta+\omega}\left[ \Gamma^i_{\tau} \right]\xrightarrow{d_S}  \inf_{j\in I} P_{\theta+\omega}\left[  \Gamma^j_{\tau} \right].
    \]
    Hence, \eqref{eq:Pthetaomega_inf} follows from \cref{prop:ds0char}.
\end{proof}

\motto{It is very fortunate that, unlike people who dig for gold, mathematicians can freely share their precious treasures with everybody. Once you understand something really well, it feels great to explain it to all.

--- Andrei Okounkov\footnote{Andrei Yuryevich Okounkov (1969--) is a Russian-American mathematician renowned for his contributions to representation theory. He was one of the key organizers of the ICM 2022 in St. Petersburg, which was unfortunately canceled under the indiscriminate discrimination against Russian citizens by the virtue signalers all over the western world after the war waged by the ruling class.}}
\chapter{The theory of Okounkov bodies}
\label{chap:Okounkov}

In this chapter, we apply our theory of singularities to the study of Okounkov bodies. We construct convex bodies $\Delta(L,\phi)$, the \emph{partial Okounkov bodies}, attached to a big line bundle $L$ on a smooth projective variety $X$ together with a plurisubharmonic metric $\phi$ on $L$. These convex bodies provide a universal numerical invariant of the $\mathcal{I}$-equivalence class of $\phi$, and reduce many problems in pluripotential theory to problems in convex geometry, which are often more tractable.

We develop two parallel constructions: an algebraic one in \cref{sec:PoB}, where $\phi$ is encoded by the multiplier ideals $\mathcal{I}(k\phi)$, and a transcendental one in \cref{sec:tpob}, where $\phi$ is treated directly as an analytic object via valuations on closed positive currents. The two constructions agree when both apply, but each has its own natural domain of definition and its own range of applications.

We assume some familiarity with the classical Okounkov bodies. Excellent references are the original papers of Lazarsfeld--Musta\c{t}\u{a} \cite{LM09} and Kaveh--Khovanskii \cite{KK12}; we give a rapid review here.

Let $X$ be an irreducible smooth projective variety of dimension $n$ and $L$ a big holomorphic line bundle on $X$. Given an admissible flag $X=Y_0\supseteq Y_1\supseteq \cdots\supseteq Y_n$ on $X$ (\cref{def:admfl}), one attaches to $L$ a natural $n$-dimensional convex body $\Delta(L)\subseteq \mathbb{R}^n$, generalizing the Newton polytope construction in toric geometry (\cref{def:Newtonbody1}). This construction was first considered by Okounkov \cite{Oko96, Oko03} and extended by Lazarsfeld--Musta\c{t}\u{a} \cite{LM09} and Kaveh--Khovanskii \cite{KK12}. The convex body $\Delta(L)$ is the \emph{Okounkov body} of $L$ with respect to the chosen flag. Concretely, the successive vanishing order along the flag defines a valuation $\nu\colon \mathbb{C}(X)^{\times}\to \mathbb{Z}^n$, and
\[
\Gamma(L)\coloneqq \left\{(\nu(s),k)\in \mathbb{Z}^{n+1}:k\in \mathbb{N},\, s\in \mathrm{H}^0(X,L^k)^{\times} \right\}
\]
is a semigroup whose closed convex cone, intersected with the slice $\{(x,1)\}$, is $\Delta(L)$. Lazarsfeld--Musta\c{t}\u{a} \cite{LM09} showed that $\Delta(L)$ depends only on the numerical class of $L$, and Jow \cite{Jow10} showed that the Okounkov bodies for all flags determine this class, so $\Delta(L)$ may be regarded as a universal numerical invariant of $L$.

The aim of this chapter is to perform the analogous construction at the level of singularity types of psh metrics. Setting
\[
\Gamma(L,\phi)\coloneqq \left\{(\nu(s),k)\in \mathbb{Z}^{n+1}:k\in \mathbb{N},\, s\in \mathrm{H}^0(X,L^k\otimes \mathcal{I}(k\phi))^{\times} \right\},
\]
one immediately encounters the obstruction that $\Gamma(L,\phi)$ is in general \emph{not} a semigroup, so the Lazarsfeld--Musta\c{t}\u{a} and Kaveh--Khovanskii constructions break down. The remedy is the notion of an \emph{almost semigroup}: $\Gamma(L,\phi)$ does satisfy a near-semigroup property (\cref{thm:Gammaasg}), and the convex-body construction $\Gamma\mapsto \Delta(\Gamma)$ extends continuously from semigroups to almost semigroups (\cref{thm:Okocont}). The partial Okounkov body of $(L,\phi)$ is then
\[
\Delta(L,\phi)\coloneqq \Delta\bigl(\Gamma(L,\phi)\bigr),
\]
and is a convex sub-body of $\Delta(L)$ --- explaining the name \emph{partial}.
The universality of partial Okounkov bodies is captured by \cref{cor:Deltacontimplyvarphi}: the $\mathcal{I}$-equivalence class of $\phi$ is determined by the Okounkov bodies of $(L,\phi)$ for sufficiently many flags.

The results in this chapter will be applied in \cref{chap:toric_big} to develop toric pluripotential theory on big line bundles.

\section{Almost semigroups}\label{sec:clo1}
We give a brief presentation of the theory of almost semigroups. The arguments are provided in \cref{chap:almostsg}.

Fix an integer $n\geq 0$. Fix a closed convex cone $C\subseteq \mathbb{R}^{n}\times \mathbb{R}_{\geq 0}$ such that $C\cap \{x_{n+1}=0\}=\{0\}$. Here $x_{n+1}$ is the last coordinate of $\mathbb{R}^{n+1}$.

Write $\hat{\mathcal{S}}(C)$ for the set of subsets of $C\cap \mathbb{Z}^{n+1}$ and $\mathcal{S}(C)$ for the set of subsemigroups $S\subseteq C\cap \mathbb{Z}^{n+1}$. For $S\in \hat{\mathcal{S}}(C)$ and each $k\in \mathbb{N}$ such that $S_k\neq\varnothing$, we write
\[
S_k\coloneqq \left\{x\in \mathbb{Z}^n: (x,k)\in S \right\}.
\]
Note that $S_k$ is a finite set by our assumption on $C$.

We introduce a pseudometric on $\hat{\mathcal{S}}(C)$ as follows:
\begin{equation}\label{eq:dsgdef}
d_{\sg}(S,S')\coloneqq \varlimsup_{k\to\infty} k^{-n} \left(|S_k|+|S'_k|-2|(S\cap S')_k|\right).
\end{equation}
Here $|\bullet|$ denotes the cardinality of a finite set. The quantity $d_{\sg}$ defined above is a pseudometric on $\hat{\mathcal{S}}(C)$.
Given $S,S'\in \hat{\mathcal{S}}(C)$, we say $S$ is equivalent to $S'$ and write $S\sim S'$ if $d_{\sg}(S,S')=0$. This is an equivalence relation.

The volume of $S\in \mathcal{S}(C)$ is defined as
\[
\vol S\coloneqq \lim_{k\to\infty} (ka)^{-n}|S_{ka}|=\varlimsup_{k\to\infty} k^{-n}|S_k|,
\]
where $a$ is a sufficiently divisible positive integer.
The existence of the limit and its independence from $a$ both follow from \cref{thm:KK12growth}.

We define $\overline{\mathcal{S}}(C)$ as the closure of $\mathcal{S}(C)$ in $\hat{\mathcal{S}}(C)$ with respect to the topology defined by the pseudometric $d_{\sg}$.
The function $\vol\colon \mathcal{S}(C)\rightarrow \mathbb{R}$ admits a unique $1$-Lipschitz extension to
\begin{equation}\label{eq:volex}
\vol\colon \overline{\mathcal{S}}(C)\rightarrow \mathbb{R}.
\end{equation}

Given $S\in \hat{\mathcal{S}}(C)$, we will write $C(S)\subseteq C$ for the closed convex cone generated by $S\cup \{0\}$. Moreover, for each $k\in \mathbb{Z}_{>0}$ such that $S_k\neq\varnothing$, we define
\[
\Delta_k(S)\coloneqq \Conv \left\{k^{-1}x\in \mathbb{R}^n:x\in S_k \right\}\subseteq \mathbb{R}^n.
\]
Here $\Conv$ denotes the convex hull.
\begin{definition}\label{def:Okounkov:L69}
    Let $\mathcal{S}'(C)$ be the subset of $\mathcal{S}(C)$ consisting of semigroups $S$ such that $S$ generates $\mathbb{Z}^{n+1}$ (as an Abelian group).
    \index{semigroup!$\mathcal{S}'(C)$}
    \index{$\mathcal{S}'(C)$}
\end{definition}
Note that for any $S\in \mathcal{S}'(C)$, the cone $C(S)$ has full dimension (i.e., the topological interior is non-empty). Given a full-dimensional subcone $C'\subseteq C$, we have $C'\cap \mathbb{Z}^{n+1}\in \mathcal{S}'(C)$.

We recall the following definition from \cite{KK12}.
\begin{definition}\label{def:Okokk}
 Given $S\in \mathcal{S}'(C)$, its \emph{Okounkov body} is defined by
\index{Okounkov body!Okounkov body of a semigroup}
\[
\Delta(S)\coloneqq \left\{x\in \mathbb{R}^n:(x,1)\in C(S) \right\}.
\]
\end{definition}

\begin{theorem}\label{thm:HausOkoun}
	For each $S\in \mathcal{S}'(C)$, we have
	\begin{equation}\label{eq:volWvolDelta}
		\vol S=\lim_{k\to\infty} k^{-n}|S_k|=\vol \Delta(S)>0.
	\end{equation}
    Moreover, as $k\to\infty$,
	\begin{equation}\label{eq:HausconvDeltaGLS}
		\Delta_{k}(S)\xrightarrow{d_{\Haus}}\Delta(S).
	\end{equation}
\end{theorem}
\begin{corollary}\label{cor:dist}
    Let $S,S'\in \mathcal{S}'(C)$. Assume that $\vol (S\cap S')>0$. Then we have
    \[
    d_{\sg}(S,S')=\vol(S)+\vol(S')-2\vol(S\cap S').
    \]
\end{corollary}

\begin{definition}\label{def:SpCbar}
    We define $\overline{\mathcal{S}'(C)}_{>0}$ as the set of elements in the closure of $\mathcal{S}'(C)$ in $\hat{\mathcal{S}}(C)$ with positive volume. An element in $\overline{\mathcal{S}'(C)}_{>0}$ is called an \emph{almost semigroup} in $C$.
    \index{$\overline{\mathcal{S}'(C)}_{>0}$}
\end{definition}
Recall that the volume here is defined in \eqref{eq:volex}.

\begin{theorem}\label{thm:Okocont}
	The Okounkov body map $\Delta\colon \mathcal{S}'(C) \rightarrow \mathcal{K}_n$
    as defined in \cref{def:Okokk} admits a unique continuous extension
	\begin{equation}\label{eq:Deltagensg}
		\Delta \colon \overline{\mathcal{S}'(C)}_{>0}\rightarrow \mathcal{K}_n.
	\end{equation}
     Moreover, for any $S\in \overline{\mathcal{S}'(C)}_{>0}$, we have
	\begin{equation}\label{eq:volWfinal}
		\vol S=\vol \Delta(S).
	\end{equation}
\end{theorem}

\begin{corollary}\label{cor:Okocomp}
    Suppose that $S,S'\in \overline{\mathcal{S}'(C)}_{>0}$ with $S\subseteq S'$. Then
    \begin{equation}\label{eq:Deltacontain}
    \Delta(S)\subseteq \Delta(S').
    \end{equation}
\end{corollary}

\section{Flags and valuations}\label{sec:flag}

\subsection{The algebraic setting}\label{subsec:flagvalalgebraic}
Let $X$ be an irreducible normal projective variety of dimension $n$.
\begin{definition}\label{def:admfl}
	An \emph{admissible flag}\index{admissible flag} $Y_{\bullet}$\index{$Y_{\bullet}$} on $X$ is a flag of subvarieties
	\[
		X=Y_0\supseteq Y_1\supseteq \cdots\supseteq Y_n
	\]
	such that $Y_i$ is irreducible of codimension $i$ and is smooth at the point $Y_n$.
\end{definition}

Given any admissible flag $Y_{\bullet}$, we can define a rank $n$ valuation $\nu_{Y_{\bullet}}\colon \mathbb{C}(X)^{\times}\rightarrow \mathbb{Z}^n$\index{$\nu_{Y_{\bullet}}$}. Here we consider $\mathbb{Z}^n$ as a totally ordered Abelian group with the lexicographic order.
We sometimes write $\mathbb{Z}^n_{\lex}$ to emphasize this point.

If we identify the elements in $\mathbb{Z}^n$ with a row vector, the automorphism group $\Aut(\mathbb{Z}^n_{\lex})$ of $\mathbb{Z}^n_{\lex}$ is then identified with the subgroup of $\GL(n,\mathbb{Z})$ consisting of matrices of the form $\mathrm{I}+U$, where $\mathrm{I}$ is the identity matrix and $U$ is a strictly upper triangular matrix with elements in $\mathbb{Z}$.

We recall the definition of $\nu_{Y_\bullet}$. Let $s^{(0)}=s\in\mathbb C(X)^\times$. For $i=1,\ldots,n$, after localizing around $Y_n$, choose a local defining equation $t_i$ of $Y_i$ in $Y_{i-1}$, and set
\[
    \nu_i=\ord_{Y_i}\left(s^{(i-1)}\right),\qquad s^{(i)}=\left(s^{(i-1)}t_i^{-\nu_i}\right)|_{Y_i} \in \mathbb C(Y_i)^\times.
\]
Then
\[
    \nu_{Y_\bullet}(s)=(\nu_1,\ldots,\nu_n)\in\mathbb Z^n.
\]
It is easy to verify that $\nu_{Y_{\bullet}}$ is indeed a rank $n$ valuation.

The same construction applies to a non-zero section $s\in \mathrm H^0(X,L)$ after choosing a local trivialization of $L$ near $Y_n$; the result is independent of this choice, since units have valuation $0$. We also write $\nu_{Y_\bullet}(D)$ for a non-zero effective Cartier divisor $D$, and extend the notation $\mathbb Q$-linearly to non-zero effective $\mathbb Q$-Cartier divisors. We also set $\nu_{Y_{\bullet}}(0)=0$.

\begin{remark}\label{rmk:Abhyankar}
Conversely, maximal-rank valuations of $\mathbb C(X)$ are closely related to flag valuations. More precisely, if a valuation has value group identified with $\mathbb Z^n_{\lex}$ and maximal rational rank $n$, then \cite[Theorem~2.9]{CFKLRS17} shows that it is equivalent\footnote{Two valuations $\nu,\nu'$ with values in $\mathbb Z^n$ are equivalent if one can find a matrix $G$ of the form $\mathrm I+N$, where $N$ is strictly upper triangular with integral entries, such that $\nu'=\nu G$.} to, but not necessarily equal to, a valuation induced by an admissible flag on a modification of $X$.
\end{remark}

\subsection{The transcendental setting}
\label{subsec:flagval_trans}
Let $X$ be a connected compact Kähler manifold of dimension $n$.

\begin{definition}\label{def:Okounkov:L170}
    A \emph{smooth flag}\index{smooth flag} $Y_{\bullet}$ on $X$ consists of a flag of connected closed submanifolds of $X$:
    \[
      X=Y_0\supseteq Y_1 \supseteq \cdots \supseteq Y_n,
    \]
    where $Y_i$ has dimension $n-i$.
\end{definition}
In this section, we will fix a smooth flag $Y_{\bullet}$ on $X$.

\begin{definition}\label{def:valcurr}
    Let $T$ be a closed positive $(1,1)$-current on $X$. We define the \emph{valuation}\index{valuation} of $T$ along $Y_{\bullet}$ as
    \[
        \nu_{Y_{\bullet}}(T)=\left(\nu_{Y_{\bullet}}(T)_1,\ldots,\nu_{Y_{\bullet}}(T)_n\right)\in \mathbb{R}_{\geq 0}^n
    \]
    \index{$\nu_{Y_{\bullet}}(T)$}
    by induction on $n$. When $n=0$, we define $\nu_{Y_{\bullet}}(T)$ as the unique point in $\mathbb{R}^0$. When $n\geq 1$, we define
    \[
        \nu_{Y_{\bullet}}(T)_1=\nu(T,Y_1);
    \]
    Then for $i=2,\ldots,n$, we define
    \[
        \nu_{Y_{\bullet}}(T)_i=\nu_{Y_1\supseteq \cdots \supseteq Y_n}\left(\Tr_{Y_1}(T-\nu(T,Y_1)[Y_1])\right)_{i-1}.
    \]
\end{definition}

\begin{proposition}\label{prop:valcurr_Iind}
    Let $T$ be a closed positive $(1,1)$-current on $X$. Then $\nu_{Y_{\bullet}}(T)\in \mathbb{R}_{\geq 0}^n$ defined in \cref{def:valcurr} is independent of the choices of the trace operators in the definition. Moreover, $\nu_{Y_{\bullet}}(T)$ depends only on the $\mathcal{I}$-equivalence class of $T$.
\end{proposition}
\begin{proof}
    We will prove both statements at the same time by induction on $n\geq 0$. The cases $n=0,1$ are trivial.

    Let us consider the case $n>1$ and assume that the result is known in dimension $n-1$.
    We first observe that $\nu_{Y_{\bullet}}(T)$ is independent of the choice of the trace operator: Different choices of $\Tr_{Y_1}(T-\nu(T,Y_1)[Y_1])$ are $\mathcal{I}$-equivalent by \cref{prop:traceunique}. Therefore, by induction, its valuation is well-defined.

    Next, let $T'$ be another closed positive $(1,1)$-current such that $T\sim_{\mathcal{I}}T'$.
    Using \cref{prop:Iequivchar}, we know that $\nu(T,Y_1)=\nu(T',Y_1)$. Therefore,
    \[
        T- \nu(T,Y_1)[Y_1]\sim_{\mathcal{I}} T'-\nu(T',Y_1)[Y_1].
    \]
    It follows by induction that
    \[
        \nu_{Y_1\supseteq \cdots \supseteq Y_n}\left(\Tr_{Y_1}(T-\nu(T,Y_1)[Y_1])\right)=\nu_{Y_1\supseteq \cdots \supseteq Y_n}\left(\Tr_{Y_1}(T'-\nu(T',Y_1)[Y_1])\right).
    \]
    This completes the proof.
\end{proof}

\begin{example}\label{ex:valuationdivcompatible}
    When $X$ is projective, and $D$ is a non-zero effective Cartier divisor, we have
    \[
    \nu_{Y_{\bullet}}([D])=\nu_{Y_{\bullet}}(D),
    \]
    where the right-hand side is defined in \cref{subsec:flagvalalgebraic}.
\end{example}

\begin{proposition}\label{prop:nuvaluationlinear}
    Let $T$, $S$ be closed positive $(1,1)$-currents on $X$, $\lambda\in \mathbb{R}_{\geq 0}$. Then
    \begin{enumerate}
     \item\label{itm:prop:nuvaluationlinear:1} if $T\preceq_{\mathcal{I}} S$, we have
        \begin{equation}\label{eq:nuTS}
            \nu_{Y_{\bullet}}(T)\geq_{\lex}\nu_{Y_{\bullet}}(S).
        \end{equation}
        \item\label{itm:prop:nuvaluationlinear:2} We have the following additivity property:
        \begin{equation}\label{eq:nuvaluationlinear}
            \nu_{Y_{\bullet}}(T+S)=\nu_{Y_{\bullet}}(T)+\nu_{Y_{\bullet}}(S),\quad \nu_{Y_{\bullet}}(\lambda T)=\lambda \nu_{Y_{\bullet}}(T).
        \end{equation}
    \end{enumerate}
\end{proposition}
\begin{proof}
\ref{itm:prop:nuvaluationlinear:1} We make an induction on $n\geq 0$. The case $n=0,1$ is trivial. Assume that $n\geq 2$ and the case $n-1$ is known. Observe that $\nu(T,Y_1)\geq \nu(S,Y_1)$, if the inequality is strict, we are done. So let us assume that $\nu(T,Y_1)=\nu(S,Y_1)$. By \cref{prop:tracelinear}, we find that
    \[
    \Tr_{Y_1}\Bigl(T-\nu(T,Y_1)[Y_1]\Bigr)\preceq_P\Tr_{Y_1}\Bigl(S-\nu(T,Y_1)[Y_1]\Bigr).
    \]
    Hence the same comparison holds for $\preceq_{\mathcal{I}}$ by \cref{cor:PimpliesI}.
    By the inductive hypothesis, we conclude \eqref{eq:nuTS}.

\ref{itm:prop:nuvaluationlinear:2} We make an induction on $n\geq 0$. The cases $n=0,1$ are trivial. Assume that $n\geq 2$ and the case $n-1$ is known. By \cref{prop:Lelongmax}, we have
    \[
    \nu(T+S,Y_1)=\nu(T,Y_1)+\nu(S,Y_1),\quad \nu(\lambda T,Y_1)=\lambda\nu(T,Y_1).
    \]
    By \cref{prop:tracelinear}, we have
    \[
    \begin{split}
    \Tr_{Y_1}\Bigl(T+S-\nu(T+S,Y_1)[Y_1]\Bigr)\sim_P&\Tr_{Y_1}\Bigl(T-\nu(T,Y_1)[Y_1]\Bigr)\\
    &+\Tr_{Y_1}\Bigl(S-\nu(S,Y_1)[Y_1]\Bigr),\\
    \Tr_{Y_1}\Bigl(\lambda T-\nu(\lambda T,Y_1)[Y_1]\Bigr)\sim_P& \lambda \Tr_{Y_1}\Bigl(T-\nu(T,Y_1)[Y_1]\Bigr).
    \end{split}
    \]
    Since $\sim_P$ implies $\sim_{\mathcal{I}}$ by \cref{cor:PimpliesI}, the inductive hypothesis and the independence of the valuation from the trace representative conclude \eqref{eq:nuvaluationlinear}.

\end{proof}

Assume that $n>0$ for the remainder of this section.
\begin{definition}\label{def:lifting}
    Let $\pi\colon Z\rightarrow X$ be a proper bimeromorphic morphism where $Z$ is a Kähler manifold. We say that a smooth flag $W_{\bullet}$ on $Z$ is a \emph{lifting} of $Y_{\bullet}$ to $Z$ if $\pi(W_i)\subseteq Y_i$ and the induced map $W_i\to Y_i$ is bimeromorphic for each $i=0,\ldots,n$.

    In this case, we define $\cor(Y_{\bullet},\pi)\in \Aut(\mathbb{Z}^n_{\lex})$ inductively as follows: When $n=1$, we define $\cor(Y_{\bullet},\pi)=[1]$; when $n>1$, we set
    \begin{equation}\label{eq:correcur}
    \cor(Y_{\bullet},\pi)\coloneqq
    \begin{bmatrix}
    1 & \nu_{W_1\supseteq \cdots \supseteq W_n}\Bigl((\pi^*[Y_1]-[W_1])|_{W_1}\Bigr)\\
    0 & \cor\Bigl(Y_1\supseteq \cdots \supseteq Y_n,\pi|_{W_1}\colon W_1\rightarrow Y_1\Bigr)
    \end{bmatrix}.
    \end{equation}
    \index{$\cor(Y_{\bullet},\pi)$}
\end{definition}
We observe that a lifting $W_{\bullet}$ of $Y_{\bullet}$ on $Z$ is unique if it exists: For each $i=0,\ldots,n-1$, the component $W_{i+1}$ is necessarily the strict transform of $Y_{i+1}$ with respect to the bimeromorphic morphism $W_i\rightarrow Y_i$. We shall also say that
$(W_{\bullet},\cor(Y_{\bullet},\pi))$ is \emph{the lifting} of $Y_{\bullet}$ to $Z$.

Note that in \eqref{eq:correcur}, $(\pi^*[Y_1]-[W_1])|_{W_1}$ can also be replaced by $\Tr_{W_1}((\pi^*[Y_1]-[W_1])|_{W_1})$, since $\pi^*[Y_1]-[W_1]$ has analytic singularities.

\begin{proposition}\label{prop:cormatrix}
    Let $\pi\colon Z\rightarrow X$ be a proper bimeromorphic morphism, where $Z$ is a Kähler manifold. Let $W_{\bullet}$ be a lifting of $Y_{\bullet}$. Then for any closed positive $(1,1)$-current $T$ on $X$, we have
    \begin{equation}
    \nu_{W_{\bullet}}(\pi^*T)= \nu_{Y_{\bullet}}(T)\cor(Y_{\bullet},\pi).
    \end{equation}
\end{proposition}
\begin{proof}
    We argue by induction on $n\geq 0$. The case $n=0$ is trivial. In general, assume that $n\geq 1$ and the result is proved in dimension $n-1$.

    For simplicity, we write $\nu=\nu_{Y_{\bullet}}$ and $\nu'=\nu_{W_{\bullet}}$.
    Let $\mu$ (resp. $\mu'$) be the valuation of currents defined by the truncated flag $Y_1\supseteq\cdots \supseteq Y_n$ (resp. $W_1\supseteq\cdots \supseteq W_n$). Then we need to show that
    \begin{equation}\label{eq:mubiration}
    \begin{split}
    &\begin{bmatrix}
     \nu'(\pi^*T)_1 &
     \mu'\Bigl( \Tr_{W_1}\bigl(\pi^*T-\nu'(\pi^*T)_1[W_1]\bigr)\Bigr)
    \end{bmatrix}\\
    =&
    \begin{bmatrix}
        \nu(T)_1 &
        \mu\Bigl(\Tr_{Y_1}\bigl(T-\nu(T)_1[Y_1]\bigr)\Bigr)
    \end{bmatrix}\cor(Y_{\bullet},\pi).
    \end{split}
    \end{equation}
    By Zariski's main theorem \cref{thm:Zariskimain}, $W_1\to Y_1$ is bimeromorphic and hence
    \[
        \nu'(\pi^*T)_1=\nu(T,Y_1)\eqqcolon c.
    \]
    By the inductive hypothesis, we have
    \begin{equation}\label{eq:ind_hypos}
        \mu'\Bigl(\Pi^{*} \Tr_{Y_1}(T-c[Y_1])\Bigr)=\mu\Bigl(\Tr_{Y_1}(T-c[Y_1])\Bigr)\cor\left(Y_1\supseteq \cdots \supseteq Y_n,\Pi\right),
    \end{equation}
    where $\Pi\colon W_1\rightarrow Y_1$ is the restriction of $\pi$. Since $\pi^*[Y_1]=[W_1]+(\pi^*[Y_1]-[W_1])$, \cref{lma:rescommpullback} and \cref{prop:tracelinear} give
    \[
    \begin{split}
        \Tr_{W_1}\bigl(\pi^*T-c[W_1]\bigr)\sim_P&\Tr_{W_1}\bigl(\pi^*(T-c[Y_1])\bigr)\\
        &+c \Tr_{W_1}\bigl(\pi^*[Y_1]-[W_1]\bigr)\\
        \sim_P&\Pi^{*}\Tr_{Y_1}\bigl(T-c[Y_1]\bigr)+c \Tr_{W_1}\bigl(\pi^*[Y_1]-[W_1]\bigr).
    \end{split}
    \]
    So
    \[
    \mu'\Bigl(\Tr_{W_1}\bigl(\pi^*T-c[W_1]\bigr)\Bigr)=\mu'\Bigl(\Pi^{*}\Tr_{Y_1}\bigl(T-c[Y_1]\bigr)\Bigr)+c\mu'\Bigl(\Tr_{W_1}\bigl(\pi^*[Y_1]-[W_1]\bigr)\Bigr).
    \]
Combining the above with \eqref{eq:ind_hypos}, we see that \eqref{eq:mubiration} follows.
\end{proof}

\begin{proposition}\label{prop:cormult}
    Let $\pi\colon Z\rightarrow X$, $p\colon Z'\rightarrow Z$ be proper bimeromorphic morphisms with $Z$ and $Z'$ being Kähler manifolds.
    Assume that $Y_{\bullet}$ admits a lifting $W_{\bullet}$ (resp. $W_{\bullet}'$) to $Z$ (resp. $Z'$). Then
    \begin{equation}\label{eq:cormul}
    \cor(Y_{\bullet},\pi\circ p)=\cor(Y_{\bullet},\pi)\cor(W_{\bullet},p).
    \end{equation}
\end{proposition}
\begin{proof}
    We let $\pi'=\pi\circ p$.
    Then $W'_{\bullet}$ is also the lifting of $W_{\bullet}$ with respect to $p$.
    \[
    \begin{tikzcd}
Z' \arrow[rr, "p"] \arrow[rd, "\pi'"] &   & Z \arrow[ld, "\pi"'] \\
                                      & X. &
\end{tikzcd}
    \]
    We argue by induction on $n\geq 1$. The case $n=1$ is trivial. Assume that $n\geq 2$ and the case $n-1$ has been solved. Then by \eqref{eq:correcur}, the desired formula \eqref{eq:cormul} can be reformulated as
    \[
    \begin{aligned}
    &\begin{bmatrix}
    1 & \nu_{W_1'\supseteq \cdots \supseteq W_n'}((\pi'^*[Y_1]-[W_1'])|_{W_1'})\\
    0 & \cor(Y_1\supseteq \cdots \supseteq Y_n,\pi'|_{W_1'}:W_1'\rightarrow Y_1)
    \end{bmatrix}=\\
    &\begin{bmatrix}
    1 & \nu_{W_1\supseteq \cdots \supseteq W_n}((\pi^*[Y_1]-[W_1])|_{W_1})\\
    0 & \cor(Y_1\supseteq \cdots \supseteq Y_n,\pi|_{W_1}:W_1\rightarrow Y_1)
    \end{bmatrix}\cdot\\
    &\begin{bmatrix}
    1 & \nu_{W_1'\supseteq \cdots \supseteq W_n'}((p^*[W_1]-[W_1'])|_{W_1'})\\
    0 & \cor(W_1\supseteq \cdots \supseteq W_n,p|_{W_1'}:W_1'\rightarrow W_1)
    \end{bmatrix}
    \end{aligned}.
    \]
    By the inductive hypothesis, this is equivalent to
    \[
    \begin{split}
    \nu_{W_1'\supseteq \cdots \supseteq W_n'}\left((\pi'^*[Y_1]-[W_1'])|_{W_1'}\right)=
        \nu_{W_1'\supseteq \cdots \supseteq W_n'}\left((p^*[W_1]-[W_1'])|_{W_1'}\right)+\\
        \nu_{W_1\supseteq \cdots \supseteq W_n}\left((\pi^*[Y_1]-[W_1])|_{W_1}\right)\cor\left(W_1\supseteq \cdots \supseteq W_n,p|_{W_1'}:W_1'\rightarrow W_1\right),
    \end{split}
    \]
    which, thanks to \cref{prop:cormatrix}, can be further rewritten as
    \[
    \begin{split}
    \nu_{W_1'\supseteq \cdots \supseteq W_n'}\left((\pi'^*[Y_1]-[W_1'])|_{W_1'}\right)=
        \nu_{W_1'\supseteq \cdots \supseteq W_n'}\left((p^*[W_1]-[W_1'])|_{W_1'}\right)+\\
        \nu_{W_1'\supseteq \cdots \supseteq W_n'}\left(p|_{W_1'}^*\left((\pi^*[Y_1]-[W_1])|_{W_1} \right)\right).
    \end{split}
    \]
    Indeed, as divisors on $W_1'$,
    \[
        (\pi'^*[Y_1]-[W_1'])|_{W_1'}=(p^*[W_1]-[W_1'])|_{W_1'}+p|_{W_1'}^*((\pi^*[Y_1]-[W_1])|_{W_1}).
    \]
    The desired equality is therefore exactly the additivity of the valuation, \cref{prop:nuvaluationlinear}.
\end{proof}

\begin{theorem}\label{thm:liftableflag}
    Let $\pi\colon Z\rightarrow X$ be a proper bimeromorphic morphism from a reduced complex space $Z$. Then there is a modification $W\rightarrow X$ dominating $Z\rightarrow X$ such that $Y_{\bullet}$ admits a lifting to $W$.
\end{theorem}
Recall that in this book, a modification means a finite composition of blowing-ups with smooth centers.
\begin{proof}
    By Hironaka's Chow lemma \cref{thm:HironakaChow}, we may assume that $\pi$ is a modification.

    We begin by setting $W_0=Z$. We will construct $W_i$ inductively for each $i$. Assume that for $0\leq i<n$ a smooth partial flag $W_0\supseteq\cdots \supseteq W_i$ has been constructed on a modification $\pi_i:Z_i\rightarrow Z$ so that $\pi\circ \pi_i$ restricts to bimeromorphic morphisms $W_j\rightarrow Y_j$ for each $j=0,\ldots,i$.

By Zariski's main theorem, $W_i\rightarrow Y_i$ is an isomorphism outside a codimension $2$ subset of $Y_i$.
We let $W_{i+1}$ be the strict transform of $Y_{i+1}$ in $W_i$.
The problem is that $W_{i+1}$ is not necessarily smooth.

We will further modify $Z_i$ and lift $W_1,\ldots, W_{i+1}$ in order to make the flag smooth. Take the embedded resolution of $(W_{j},W_{i+1})$, say $W_j'\rightarrow W_j$ for each $j=0,\ldots,i$.

We have canonical embeddings $W_{i}'\hookrightarrow W_{i-1}'\hookrightarrow\cdots \hookrightarrow W_0'$ making the following diagram commutative:
\[
\begin{tikzcd}
W_i' \arrow[d] \arrow[r, hook] & W_{i-1}' \arrow[d] \arrow[r, hook] & \cdots \arrow[r, hook] \arrow[d, "\vdots", phantom] & W_0' \arrow[d] \\
W_i \arrow[r, hook]            & W_{i-1} \arrow[r, hook]            & \cdots \arrow[r, hook]                              & W_0
\end{tikzcd}
\]
Let $W_{i+1}'$ be the strict transform of $W_{i+1}$ in $W_i'$. It suffices to define $\pi_{i+1}$ as the morphism $W_0'\rightarrow Z_i\rightarrow Z$ and
replace $W_0\supseteq \cdots\supseteq W_{i+1}$ by $W_0'\supseteq \cdots\supseteq W_{i+1}'$.
\end{proof}

We also include a more general version of \cref{def:valcurr} for currents on bimeromorphic models.
\begin{definition}\label{def:modelvalcurr}
    Let $\pi\colon Z\rightarrow X$ be a proper bimeromorphic morphism from a normal compact Kähler space, and let $S$ be a closed positive $(1,1)$-current on $Z$. Choose a proper bimeromorphic morphism
    \index{valuation!valuation of a current on a model}
    $h\colon \widehat Z\rightarrow Z$ from a compact Kähler manifold such that
    $Y_{\bullet}$ admits a lifting $(W_{\bullet},g)$ to $\widehat Z$ with
    respect to $\pi\circ h$, where
    \[
        g=\cor(Y_{\bullet},\pi\circ h).
    \]
    We define
    \begin{equation}\label{eq:nuYbulletS}
        \nu_{Y_{\bullet}}(S)\coloneqq \nu_{W_{\bullet}}(h^*S)g^{-1}\in \mathbb{R}^n.
    \end{equation}
\end{definition}
Note that it could happen that $\nu_{Y_{\bullet}}(S)\notin \mathbb{R}_{\geq 0}^n$.

\begin{lemma}\label{lma:modelvalcurrbir}
    The definition in \cref{def:modelvalcurr} is independent of the choice of $h$. Moreover, if $p\colon Z'\rightarrow Z$ is a proper bimeromorphic morphism from a normal compact Kähler space, then, with respect to the model $\pi\circ p\colon Z'\rightarrow X$,
    \[
        \nu_{Y_{\bullet}}(p^*S)=\nu_{Y_{\bullet}}(S).
    \]
\end{lemma}
\begin{proof}
    We first prove independence of the choice of $h$. Let $h_i\colon \widehat Z_i\rightarrow Z$, $i=1,2$, be two choices, and let $(W^i_{\bullet},g_i)$ be the corresponding liftings, with $g_i=\cor(Y_{\bullet},\pi\circ h_i)$. Choose a common modification over $Z$, $q_i\colon U\rightarrow \widehat Z_i$ such that $Y_{\bullet}$ admits a lifting $W_{\bullet}$ to $U$. Then $W_{\bullet}$ is also the lifting of $W^i_{\bullet}$ with respect to $q_i$. Put
    \[
        G_i=\cor(W^i_{\bullet},q_i).
    \]
    By \cref{prop:cormatrix,prop:cormult},
    \[
    \begin{aligned}
        \nu_{W_{\bullet}}\bigl(q_i^*h_i^*S\bigr)
        \cor(Y_{\bullet},\pi\circ h_i\circ q_i)^{-1}
        &=
        \nu_{W^i_{\bullet}}(h_i^*S)G_i(g_iG_i)^{-1}  \\
        &=
        \nu_{W^i_{\bullet}}(h_i^*S)g_i^{-1}.
    \end{aligned}
    \]
    Since $h_1\circ q_1=h_2\circ q_2$, the left-hand side is independent of $i$. This proves that \eqref{eq:nuYbulletS} is well-defined.

    For the birational invariance, choose a proper bimeromorphic morphism $q\colon \widehat Z\rightarrow Z'$ from a compact Kähler manifold such that $Y_{\bullet}$ admits a lifting to $\widehat Z$ with respect to $\pi\circ p\circ q$. Computing $\nu_{Y_{\bullet}}(p^*S)$ using $q$, and computing $\nu_{Y_{\bullet}}(S)$ using $p\circ q$, gives the same expression.
\end{proof}

\begin{remark}\label{rmk:valcurr}
    Suppose that $X$ is a normal projective variety. Consider a rank $n$ (surjective) valuation $\nu\colon \mathbb{C}(X)^{\times}\rightarrow \mathbb{Z}^n$ and a closed positive $(1,1)$-current $T$ on $X$. Then we can always define $\nu(T)\in \mathbb{R}^n$ as follows: Take a resolution $\pi\colon Y\rightarrow X$ such that there is a smooth flag $Y_{\bullet}$ on $Y$ and $g\in \Aut(\mathbb{Z}^n_{\lex})$ such that
    \[
    \nu=\nu_{Y_{\bullet}}g.
    \]
    Then we define
    \[
    \nu(T)\coloneqq \nu_{Y_{\bullet}}(\pi^*T)g.
    \]
    This definition does not depend on the choice of $\pi$, by passing to a common resolution and applying \cref{prop:cormatrix,prop:cormult}.
\end{remark}

\section{Algebraic partial Okounkov bodies}\label{sec:PoB}
Let $X$ be a connected smooth complex projective variety of dimension $n$ and $(L,h)$ be a Hermitian big line bundle on $X$.

Let $h_0$ be a smooth Hermitian metric on $L$. Let $\theta=c_1(L,h_0)$. Then we can identify $h$ with a function $\varphi\in \PSH(X,\theta)$.
We will use interchangeably the notations $(\theta,\varphi)$ and $(L,h)$. We assume that $\vol \theta_{\varphi}>0$ in this section.

Fix a rank $n$ valuation $\nu\colon \mathbb{C}(X)^{\times}\rightarrow \mathbb{Z}^n$, which without loss of generality can be assumed to be surjective. We shall frequently use the standard fact that such a full-rank valuation has one-dimensional leaves: for every finite-dimensional subspace $V\subset \mathbb C(X)$,
\[
    \#\nu(V^\times)=\dim V.
\]
See \cite[Lemma~1.3]{LM09} for instance.

We will adopt the notations of \cref{sec:clo1}.

\subsection{The spaces of sections}
\label{subsec:section_spaces}

\begin{definition}\label{def:Okounkov:L484}
    We will write
\[
\begin{aligned}
\Gamma(\theta,\varphi)\coloneqq & \left\{ (\nu(s),k) : k\in \mathbb{N}, s\in \mathrm{H}^0(X,L^k\otimes \mathcal{I}(k\varphi))^{\times} \right\},\\
\Delta_k(\theta,\varphi)\coloneqq & \Conv\left\{ k^{-1}\nu(s): s\in \mathrm{H}^0\left(X,L^k\otimes \mathcal{I}(k\varphi)\right)^{\times} \right\}\subseteq \mathbb{R}^n,
\end{aligned}
\]
\index{$\Gamma(\theta,\varphi)$}\index{$\Delta_k(\theta,\varphi)$}
where the second set is defined only for $k>0$ such that $\mathrm{H}^0(X,L^k\otimes \mathcal{I}(k\varphi))\neq 0$.

When $\varphi=V_{\theta}$, we simply write $\Gamma(L)$ and $\Delta_k(L)$ instead.
\end{definition}
Here and in the sequel, the cross notation means excluding $0$.
Here $\Conv$ denotes the convex hull. For large enough $k$, $\Delta_k(\theta,\varphi)$ is non-empty thanks to \cref{thm:DXmain1}.

\begin{definition}\label{def:Okounkov:L500}
    Assume that $\varphi$ has analytic singularities.
We define
    \begin{equation}\label{eq:Weps1}
        \Gamma^{\infty}(\theta,\varphi)\coloneqq \left\{ (\nu(s),k):k\in \mathbb{N}, s\in \mathrm{H}^0(X,L^k\otimes \mathcal{I}_{\infty}(k\varphi))^{\times} \right\}.
    \end{equation}\index{$\Gamma^{\infty}(\theta,\varphi)$}
\end{definition}
Recall that $\mathcal{I}_{\infty}$ is introduced in \cref{def:Iinfty}.

For later use, we introduce a twisted version as well.

\begin{definition}\label{def:DeltakT}
If $T$ is a holomorphic line bundle on $X$, we introduce
\[
\begin{aligned}
\Delta_{k,T}(\theta,\varphi)\coloneqq & \Conv\left\{ k^{-1}\nu(s): s\in \mathrm{H}^0(X,T\otimes L^k\otimes \mathcal{I}(k\varphi))^{\times} \right\}\subseteq \mathbb{R}^n,\\
\Delta_{k,T}(L)\coloneqq & \Conv\left\{ k^{-1}\nu(s): s\in \mathrm{H}^0(X,T\otimes L^k)^{\times} \right\}\subseteq \mathbb{R}^n
\end{aligned}
\]
\index{$\Delta_{k,T}(\theta,\varphi)$}\index{$\Delta_{k,T}(L)$}
for all $k\in \mathbb{Z}_{>0}$ such that the relevant spaces of sections are non-zero.
\end{definition}

\subsection{Algebraic Okounkov bodies}
\begin{proposition}\label{prop:Okounbiglbdl}
    There is a convex body $\Delta\in \mathcal{K}_n$ such that $\Gamma(L)\in \mathcal{S}'(\Delta)$.
\end{proposition}
Here and in the sequel, the notation $\mathcal{S}'(\Delta)$ means $\mathcal{S}'(C(\Delta))$, where $C(\Delta)$ is the cone over $\Delta$.
Recall that the notations $\mathcal{K}_n$ and $\mathcal{S}'(\Delta)$ are introduced in \cref{chap:almostsg}.

\begin{proof}
\textbf{Step~1}. We first show that there is $\Delta\in \mathcal{K}_n$ such that $\Delta_k(L)\subseteq \Delta$. For this purpose, using \cref{rmk:Abhyankar}, we may assume that $\nu$ is induced by an admissible flag $Y_{\bullet}$ on $X$. We may further assume that each $Y_i$ is smooth.

Fix $s\in \mathrm{H}^0(X,L^k)^{\times}$ for some $k\in \mathbb{Z}_{>0}$. Assume that $s\neq 0$. We need to show that for each $i=1,\ldots,n$, $\nu(s)_i\leq Ck$ for some constant $C>0$, independent of the choices of $k$ and $s$.

Fix an ample divisor $H$ on $X$. Take a large enough integer $b_1>0$ such that
\[
    (L-b_1 Y_1)\cdot H^{n-1}<0.
\]
Then $\nu(s)_1\leq b_1 k$. Next take a large enough integer $b_2$ such that, for every $a\in[0,b_1]$,
\[
\left((L-aY_1)|_{Y_1}-b_2Y_2\right) \cdot H^{n-2}<0.
\]
After removing the first vanishing order of $s$ along $Y_1$, the resulting section on $Y_1$ is a section of $k((L-aY_1)|_{Y_1})$ for some $a\in[0,b_1]$. Hence $\nu(s)_2\leq b_2k$. Continuing in this manner and choosing the constants uniformly over the compact ranges of the preceding normalized vanishing orders, we conclude that $\nu(s)_i/k$ is bounded for each $i$.

Here the expressions involving $L-aY_1$ are understood numerically. At the $i$-th step the previous normalized vanishing orders lie in a compact polytope, and the intersection number is continuous in these parameters; hence the constants $b_i$ can be chosen uniformly.

\textbf{Step~2}. Observe that $\Gamma(L)$ is clearly a semigroup. It remains to show that $\Gamma(L)$ generates $\mathbb{Z}^{n+1}$ as an Abelian group.

For this purpose, take two very ample divisors $A$ and $B$ so that $L=\mathcal{O}_X(A-B)$. After choosing $A$ and $B$ ample enough, we may guarantee that there exist sections $s_0\in \mathrm{H}^0(X,A)$, $t_i\in \mathrm{H}^0(X,B)$ for $i=0,\dots,n$ such that
\[
\nu(s_0)=\nu(t_0)=0
\]
and $\nu(t_i)$ is the $i$-th unit vector $e_i\in \mathbb{R}^n$ for $i=1,\ldots,n$.

Since $L$ is big, we can find $m_0>0$ such that for any $m\geq m_0$, we can find an effective divisor $F_m$ on $X$ linearly equivalent to $mL-B$. Let $f_m=\nu([F_m])$. Then
\[
(f_m,m),(f_m+e_1,m),\dots,(f_m+e_n,m)\in \Gamma(L).
\]
Since $(m+1)L$ is linearly equivalent to $A+F_m$, we have
\[
(f_m,m+1)\in \Gamma(L).
\]
It follows that $\Gamma(L)$ generates $\mathbb{Z}^{n+1}$.
\end{proof}

Thanks to \cref{prop:Okounbiglbdl}, we can introduce the next definition.
\begin{definition}\label{def:Okobd}
    We define the \emph{Okounkov body}\index{Okounkov body} of $L$ with respect to the valuation $\nu$ as
    \[
    \Delta_{\nu}(L)\coloneqq \Delta(\Gamma(L)).
    \]
    \index{$\Delta_{\nu}(L)$}
\end{definition}
When $\nu$ is induced by a smooth flag $Y_{\bullet}$ on $X$, we also write $\Delta_{Y_{\bullet}}(L)$ instead. The same convention applies to the partial Okounkov bodies studied below as well.

\begin{proposition}\label{prop:Okounonlydepnum}
    The Okounkov body $\Delta_{\nu}(L)$ depends only on the numerical class of $L$.
\end{proposition}
See \cite[Proposition~4.1]{LM09} for the elegant proof.

\begin{corollary}\label{cor:Okounvol}
    We have
    \begin{equation}
        \vol \Delta_{\nu}(L)=\frac{1}{n!}\vol L.
    \end{equation}
\end{corollary}
\begin{proof}
    This follows immediately from \cref{prop:Okounbiglbdl} and \cref{thm:HausOkoun}.
\end{proof}

\begin{proposition}\label{prop:GammaepsSp}
Assume that $\varphi$ has analytic singularities and $\theta_{\varphi}$ is a Kähler current. Then we have
    \[
    \Gamma^{\infty}(\theta,\varphi)\in \mathcal{S}'(\Delta_\nu(L))
    \]
    and
    \[
    \vol \Gamma^{\infty}(\theta,\varphi)=\frac{1}{n!}\int_X \theta_{\varphi}^n.
    \]
\end{proposition}
\begin{proof}
    By \cref{thm:resolvelogsing}, we can choose a modification $\pi\colon Y\to X$ such that $\pi^*\varphi$ has log singularities along an effective simple normal crossing $\mathbb Q$-divisor $D$. We regard $\nu$ as a valuation of $\mathbb C(Y)=\mathbb C(X)$. By the definition of $\mathcal I_\infty$, the section spaces defining $\Gamma^\infty(\theta,\varphi)$ are naturally identified with
    \[
        \mathrm{H}^0\left(Y,\pi^*L^k\otimes\mathcal O_Y(-\lceil kD\rceil)\right).
    \]
    Thus we may compute $\Gamma^\infty(\theta,\varphi)$ on $Y$.

    In this case,
    \[
    \Gamma^{\infty}(\theta,\varphi)=
    \Biggl\{
    (\nu(s),k):k\in \mathbb{N},  s\in \mathrm{H}^0\left(Y,\pi^*L^k\otimes \mathcal{O}_Y(-\lceil k D \rceil)\right)^{\times}
    \Biggr\}
    \]
    Since $\pi^*L-D$ is nef and big by \cref{lma:logsingrem}, our assertion follows from the same argument as \cref{prop:Okounbiglbdl}.
\end{proof}

We first extend \cref{thm:HausOkoun} to the twisted case.
\begin{proposition}\label{prop:Deltaconvtwisted}
    For any holomorphic line bundle $T$ on $X$, as $k\to\infty$
    \[
        \Delta_{k,T}(L)\xrightarrow{d_{\Haus}} \Delta_{\nu}(L).
    \]
\end{proposition}

\begin{proof}
    As $L$ is big, we can take $k_0\in \mathbb{Z}_{>0}$ so that
    \begin{enumerate}
        \item\label{itm:chapter-Okounkov:L657:1} $T^{-1}\otimes L^{k_0}$ admits a non-zero global holomorphic section $s_0$, and
        \item\label{itm:chapter-Okounkov:L657:2} $T\otimes L^{k_0}$ admits a non-zero global holomorphic section $s_1$.
    \end{enumerate}
    Then for $k\in \mathbb{Z}_{>k_0}$, we have injective linear maps
    \[
        \mathrm{H}^0(X,L^{k-k_0})\xrightarrow{\times s_1}\mathrm{H}^0(X,T\otimes L^{k})\xrightarrow{\times s_0}\mathrm{H}^0(X, L^{k+k_0}).
    \]
    It follows that
    \[
        (k-k_0)\Delta_{k-k_0}(L)+\nu(s_1)\subseteq k\Delta_{k,T}(L)\subseteq (k+k_0)\Delta_{k+k_0}(L)-\nu(s_0).
    \]
    Using \cref{thm:HausOkoun}, we conclude.
\end{proof}

\begin{proposition}\label{prop:subaddOkoun}
    Let $L'$ be another big line bundle on $X$. Then
    \[
    \Delta_{\nu}(L)+\Delta_{\nu}(L')\subseteq \Delta_{\nu}(L\otimes L').
    \]
\end{proposition}
\begin{proof}
    Observe that for each $k\in \mathbb{N}$, we have
    \[
    \Delta_k(L)+\Delta_k(L')\subseteq \Delta_k(L\otimes L').
    \]
    So our assertion follows immediately from \cref{thm:HausOkoun}.
\end{proof}

\begin{proposition}\label{prop:Okourescaling}
    For any $a\in \mathbb{Z}_{>0}$, we have
    \[
    \Delta_{\nu}(L^a)=a\Delta_{\nu}(L).
    \]
\end{proposition}
\begin{proof}
    This is an immediate consequence of \cref{thm:HausOkoun}.
\end{proof}

\subsection{Construction of partial Okounkov bodies}
\label{subsec:pob_construct}

\begin{theorem}\label{thm:Gammaasg}
    We have
    \[
    \Gamma(\theta,\varphi)\in \overline{\mathcal{S}'(\Delta_{\nu}(L))}_{>0}.
    \]
\end{theorem}
We refer to \cref{def:SpCbar} for the definition of $\overline{\mathcal{S}'(\Delta_{\nu}(L))}_{>0}$.

This theorem allows us to give the following definition:
\begin{definition}\label{def:Okounkov:L687}
The \emph{partial Okounkov body}\index{Okounkov body!partial Okounkov body} of $(L,h)$ is defined as
\begin{equation}\label{eq:Deltalbdef}
\Delta_{\nu}(L,h)=\Delta_{\nu}(\theta,\varphi)\coloneqq \Delta \left(\Gamma(\theta,\varphi)\right).
\end{equation}
\index{$\Delta_{\nu}(L,h)$}
\end{definition}
When $\nu$ is induced by an admissible flag $Y_{\bullet}$ on $X$ (see \cref{def:admfl}), we also say that $\Delta_{\nu}(\theta,\varphi)$ is the \emph{partial Okounkov body} of $(L,h)$ or of $(\theta,\varphi)$ with respect to $Y_{\bullet}$. In this case, we also write $\Delta_{Y_{\bullet}}$ instead of $\Delta_{\nu}$.

Note that when $h$ has minimal singularities, we have
\[
\Delta_{\nu}(L,h)=\Delta_{\nu}(L).
\]
So partial Okounkov bodies generalize Okounkov bodies.

\begin{corollary}\label{cor:POBvolume}
    We have
\begin{equation}\label{eq:Okov}
    \vol \Delta_{\nu}(\theta,\varphi)=\frac{1}{n!}\vol \theta_{\varphi}.
\end{equation}
\end{corollary}

We will prove \cref{thm:Gammaasg} and \cref{cor:POBvolume} at the same time. The proof relies on the pseudometric $d_{\sg}$ introduced in \eqref{eq:dsgdef}.

\begin{proof}
    \textbf{Step~1}. We first assume that $\varphi$ has analytic singularities and $\theta_{\varphi}$ is a Kähler current.

    We claim that
    \begin{equation}\label{eq:Gamma0Gammaanalytic}
        d_{\sg}(\Gamma^{\infty}(\theta,\varphi),\Gamma(\theta,\varphi))=0.
    \end{equation}
    Observe that for each $\epsilon\in \mathbb{Q}_{>0}$, we have
    \[
    \mathrm{H}^0\Bigl(X,L^k\otimes \mathcal{I}_{\infty}(k\varphi)\Bigr)\subseteq \mathrm{H}^0\Bigl(X,L^k\otimes \mathcal{I}(k\varphi)\Bigr)\subseteq \mathrm{H}^0\left(X,L^k\otimes \mathcal{I}_{\infty}(k(1-\epsilon)\varphi)\right)
    \]
    for all large enough $k$. This is a consequence of \cref{lma:IandIinf}. Therefore, it suffices to show that
    \[
    \lim_{\mathbb{Q}\ni \epsilon\to 0+}\vol \Gamma^{\infty}(\theta,(1-\epsilon)\varphi)=\vol \Gamma^{\infty}(\theta,\varphi).
    \]
    This follows from the explicit formula in \cref{prop:GammaepsSp}.

    \textbf{Step~2}. We next handle the case where $\theta_{\varphi}$ is a Kähler current.

    Let $(\varphi_j)_j$ be a quasi-equisingular approximation of $\varphi$ in $\PSH(X,\theta)$. Then $\varphi_j\xrightarrow{d_S}P_{\theta}[\varphi]_{\mathcal{I}}$ by \cref{cor:quasi-equichar}.

In this case, it suffices to prove that
\begin{equation}\label{eq:WtoWclaim}
	\Gamma(\theta,\varphi_j)\xrightarrow{d_{\sg}} \Gamma(\theta,\varphi).
\end{equation}
In fact, by \cref{thm:DXmain1}, we have
\[
	\begin{aligned}
		&d_{\sg}\left(\Gamma(\theta,\varphi_j),\Gamma(\theta,\varphi)\right) \\
  =& \varlimsup_{k\to\infty} k^{-n} \left( h^0\left(X,L^k\otimes \mathcal{I}(k\varphi_j)\right)-h^0\left(X,L^k\otimes \mathcal{I}(k\varphi)\right)\right)\\
  = & \lim_{k\to\infty} k^{-n}h^0\left(X,L^k\otimes \mathcal{I}(k\varphi_j)\right)-\lim_{k\to\infty} k^{-n}h^0\left(X,L^k\otimes \mathcal{I}(k\varphi)\right) \\
		 =&\frac{1}{n!}\vol \theta_{\varphi_j}-\frac{1}{n!}\vol \theta_{\varphi}.
	\end{aligned}
\]
Letting $j\to\infty$, we conclude \eqref{eq:WtoWclaim} by \cref{thm:contvolu}.

\textbf{Step~3}. Now we only assume that $\vol \theta_{\varphi}>0$.
We may replace $\varphi$ with $P_{\theta}[\varphi]_{\mathcal{I}}$ and then assume that $\varphi\in \PSH(X,\theta)_{>0}$.

Take a potential $\psi\in \PSH(X,\theta)$ such that $\psi \leq \varphi$ and $\theta_{\psi}$ is a Kähler current. The existence of $\psi$ is proved in \cref{lma:kahcurrentposmass}.
For each $\epsilon\in (0,1)$, let $\varphi_{\epsilon}=(1-\epsilon)\varphi+\epsilon\psi$.
It suffices to show that
\[
\Gamma(\theta,\varphi_{\epsilon})\xrightarrow{d_{\sg}}\Gamma(\theta,\varphi)
\]
as $\epsilon \to 0+$.
We compute using \cref{thm:DXmain1}:
\[
\begin{aligned}
    &d_{\sg}\left(\Gamma(\theta,\varphi_{\epsilon}),\Gamma(\theta,\varphi)\right)\\
    =& \varlimsup_{k\to\infty} k^{-n} \left( h^0\left(X,L^k\otimes \mathcal{I}(k\varphi)\right)- h^0\left(X,L^k\otimes \mathcal{I}(k\varphi_{\epsilon})\right) \right)\\
    =& \lim_{k\to\infty}k^{-n}h^0\left(X,L^k\otimes \mathcal{I}(k\varphi)\right)-\lim_{k\to\infty}k^{-n}h^0\left(X,L^k\otimes \mathcal{I}(k\varphi_{\epsilon})\right)\\
    =& \frac{1}{n!}\vol \theta_{\varphi}- \frac{1}{n!}\vol \theta_{\varphi_{\epsilon}}\\
    \to & 0
\end{aligned}
\]
by \cref{thm:contvolu}, as $\epsilon \to 0+$. Combining the three steps, $\Gamma(\theta,\varphi)$ lies in $\overline{\mathcal S'(\Delta_\nu(L))}_{>0}$. This proves \cref{thm:Gammaasg}. Moreover, by the continuity of the volume on $\overline{\mathcal S}(\Delta_\nu(L))$ and by \cref{thm:DXmain1},
\[
\vol \Gamma(\theta,\varphi)=\lim_{k\to\infty}k^{-n}h^0(X,L^k\otimes\mathcal I(k\varphi))=\frac{1}{n!}\vol\theta_\varphi.
\]
Since $\vol\Delta(\Gamma(\theta,\varphi))=\vol\Gamma(\theta,\varphi)$ by \cref{thm:Okocont}, \eqref{eq:Okov} follows.
\end{proof}

\begin{remark}\label{rmk:DeltaanaW0}
It follows from the proof that if $\varphi$ has analytic singularities and $\theta_{\varphi}$ is a Kähler current, then \eqref{eq:Gamma0Gammaanalytic} holds.

If we take a modification $\pi\colon Y\rightarrow X$ such that $\pi^*\varphi$ has log singularities along an effective $\mathbb{Q}$-divisor $D$ on $Y$, then
\begin{equation}\label{eq:DeltaanaW0}
\Delta_{\nu}(\theta,\varphi)=\Delta_{\nu}(\pi^*L-D)+\nu(D).
\end{equation}

This is a very special case of \cref{thm:pobmain}.
\end{remark}

\subsection{Basic properties of partial Okounkov bodies}
\label{subsec:pob_prop}

\begin{proposition}\label{prop:Okocurrent}
	The partial Okounkov body $\Delta_{\nu}(L,h)$ depends only on $\ddc h$, not on the explicit choices of $L,h_0,h$.
\end{proposition}
Thanks to this result, given a closed positive $(1,1)$-current $T\in c_1(L)$ on $X$ with $\int_X T^n>0$, we can write
\[
\Delta_{\nu}(T)\coloneqq \Delta_{\nu}(\theta,\varphi)
\]
if $T=\theta+\ddc\varphi$ for some $\varphi\in \PSH(X,\theta)$.
\begin{proof}
	There are two different claims to prove, as detailed in the two steps below.

\textbf{Step~1}. Let $h_0'$ be another Hermitian metric on $L$. Set $\theta'=c_1(L,h_0')$. Write $\ddc f=\theta-\theta'$.
	Let $\varphi'=\varphi+f\in \PSH(X,\theta')$. Then
	\begin{equation}\label{eq:DeltaDelta1}
		\Delta_{\nu}(\theta,\varphi)=\Delta_{\nu}(\theta',\varphi').
	\end{equation}

 This is obvious since $\Gamma(\theta,\varphi)=\Gamma(\theta',\varphi')$.

 \textbf{Step~2}. Let $L'$ be another big line bundle on $X$. By Step~1, we may assume that the reference Hermitian metric $h'_0$ on $L'$ is such that $c_1(L',h'_0)=\theta$.

 Let $h'$ be a plurisubharmonic metric on $L'$ with $c_1(L,h)=c_1(L',h')$. Then
 \[
 \Delta_{\nu}(L,h)=\Delta_{\nu}(L',h').
 \]
 By the approximation step in the construction of $\Delta_\nu(L,h)$ and the continuity of the extension in \cref{thm:Okocont}, it is enough to prove the assertion when $c_1(L,h)$ has analytic singularities. After taking a birational resolution, it suffices to deal with the case where $c_1(L,h)$ has analytic singularities along an effective $\mathbb{Q}$-divisor $D$. By rescaling, we may also assume that $D$ is a divisor. By \cref{rmk:DeltaanaW0}, we further reduce to the case where $c_1(L,h)$ is not singular.

In this case, the assertion is proved in \cref{prop:Okounonlydepnum}.
\end{proof}

\begin{proposition}\label{prop:IcompimplyDeltacomp}
	Let $\varphi,\psi\in \PSH(X,\theta)_{>0}$. Assume that $\varphi\preceq_{\mathcal{I}} \psi$. Then
	\begin{equation}\label{eq:Deltacomp}
		\Delta_{\nu}(\theta,\varphi)  \subseteq \Delta_{\nu}(\theta,\psi).
	\end{equation}
\end{proposition}
In particular, we always have
\[
\Delta_{\nu}(\theta,\varphi)\subseteq \Delta_{\nu}(L).
\]
\begin{proof}
	This follows from \cref{cor:Okocomp}.
\end{proof}

\begin{theorem}\label{thm:Okoucont}
	The Okounkov body map
	\[
		\Delta_{\nu}(\theta,\bullet):\left( \PSH(X,\theta)_{>0},d_S\right)\rightarrow \left(\mathcal{K}_n,d_{\Haus}\right)
	\]
	is continuous.
\end{theorem}

\begin{proof}
Let $\varphi_j\to \varphi$ be a $d_S$-convergent sequence in $\PSH(X,\theta)_{>0}$.
	We want to show that
	\begin{equation}\label{eq:Deltavjv}
		\Delta_{\nu}(\theta,\varphi_j)\xrightarrow{d_{\Haus}}\Delta_{\nu}(\theta,\varphi).
	\end{equation}
	Replacing the potentials by their $P$-envelopes does not change their partial Okounkov bodies, since $P$-equivalent potentials are $\mathcal{I}$-equivalent, nor does it change their $d_S$-classes. Thus we may assume that all $\varphi_j$'s and $\varphi$ are model potentials.

	By \cref{prop:incanddec}, after passing to a subsequence we can sandwich $\varphi_j$ between a decreasing sequence and an increasing sequence, both converging to $\varphi$ in $d_S$. Using \cref{prop:IcompimplyDeltacomp} and \cref{thm:Blaschke}, it is enough to prove \eqref{eq:Deltavjv} when $(\varphi_j)_j$ is monotone. By \cref{thm:contPI}, replacing $\varphi_j$ and $\varphi$ by their $\mathcal{I}$-envelopes preserves this reduction, so we may further assume that the $\varphi_j$'s and $\varphi$ are $\mathcal{I}$-model.
 In both cases, we claim that
 \[
 \Gamma(\theta,\varphi_j)\xrightarrow{d_{\sg}} \Gamma(\theta,\varphi)
 \]
 as $j\to\infty$.
 In fact, using \cref{thm:DXmain1}, we can compute
 \[
 \begin{aligned}
 d_{\sg}\left(\Gamma(\theta,\varphi_j), \Gamma(\theta,\varphi)\right) = &\varlimsup_{k\to\infty} k^{-n}\left| h^0\left(X,L^k\otimes \mathcal{I}(k\varphi_j)\right)-h^0\left(X,L^k\otimes \mathcal{I}(k\varphi)\right) \right|\\
 =& \frac{1}{n!} \left| \vol \theta_{\varphi_j}-\vol \theta_{\varphi}\right|,
 \end{aligned}
 \]
 which converges to $0$ by \cref{thm:contvolu}.
\end{proof}

\begin{proposition}\label{prop:birinvO}
	Let $\pi \colon Y\rightarrow X$ be a modification. Then
	\[
		\Delta_{\nu}(\pi^*L,\pi^*h)=\Delta_{\nu}(L,h).
	\]
\end{proposition}
\begin{proof}
    Thanks to \cref{prop:Ienvelopebimero}, we may assume that $\varphi$ is $\mathcal{I}$-model.
   By \cref{thm:charIgoodasclosure}, we can find a sequence $(\varphi_j)_j$ with analytic singularities in $\PSH(X,\theta)$ such that $\varphi_j\xrightarrow{d_S}\varphi$. It is clear that $\pi^*\varphi_j\xrightarrow{d_S}\pi^*\varphi$.
    By \cref{thm:Okoucont}, we may then reduce to the case where $\varphi$ has analytic singularities. In this case, it suffices to apply \cref{rmk:DeltaanaW0}.
\end{proof}

\begin{proposition}\label{prop:suba}
	Let $(L', h')$ be another Hermitian big line bundle on $X$. Then
	\[
		\Delta_{\nu}(L,h)+\Delta_{\nu}(L', h')\subseteq\Delta_{\nu}(L\otimes L',h\otimes h').
	\]
\end{proposition}
\begin{proof}
	Take a Hermitian metric $h_0'$ on $L'$ and let $\theta'=c_1(L',h_0')$. We identify $h'$ with $\varphi'\in \PSH(X,\theta')$. Then we need to show
	\begin{equation}\label{eq:suba}
		\Delta_{\nu}(\theta,\varphi)+\Delta_{\nu}(\theta',\varphi')\subseteq \Delta_{\nu}(\theta+\theta',\varphi+\varphi').
	\end{equation}
    We observe that
    \[
    P_{\theta}[\varphi]_{\mathcal{I}}+P_{\theta'}[\varphi']_{\mathcal{I}}\sim_{\mathcal{I}} \varphi+\varphi'.
    \]
    Thus, after replacing $\varphi$ and $\varphi'$ by their $\mathcal{I}$-envelopes, in view of \cref{prop:IcompimplyDeltacomp}, we may assume that $\varphi$ and $\varphi'$ are $\mathcal{I}$-good.

	By \cref{thm:charIgoodasclosure}, we can find sequences $(\varphi_j)_j$ and $(\varphi'_j)_j$ in $\PSH(X,\theta)_{>0}$ and $\PSH(X,\theta')_{>0}$ respectively such that
	\begin{enumerate}
		\item\label{itm:chapter-Okounkov:L930:1} $\varphi_j$ and $\varphi'_j$ both have analytic singularities for all $j\geq 1$, and
		\item\label{itm:chapter-Okounkov:L930:2} $\varphi_j\xrightarrow{d_S}\varphi$, $\varphi'_j\xrightarrow{d_S}\varphi'$.
	\end{enumerate}
	Then $\varphi_j+\varphi'_j\in \PSH(X,\theta+\theta')_{>0}$ and $\varphi_j+\varphi'_j\xrightarrow{d_S} \varphi+\varphi'$ by \cref{thm:dSadditivity}. Thus, by \cref{thm:Okoucont}, we may assume that $\varphi$ and $\varphi'$ both have analytic singularities. Taking a birational resolution, we may further assume that they have log singularities. By \cref{rmk:DeltaanaW0}, we reduce to the case without singularities, in which case the result is just \cref{prop:subaddOkoun}.
\end{proof}

\begin{theorem}\label{thm:concOko}
	Let $\varphi,\psi\in \PSH(X,\theta)_{>0}$. Then for any $t\in (0,1)$,
	\begin{equation}\label{eq:Deltaconcave}
		\Delta_{\nu}\left(\theta,t\varphi+(1-t)\psi\right)\supseteq t\Delta_{\nu}(\theta,\varphi)+(1-t)\Delta_{\nu}(\theta,\psi).
	\end{equation}
\end{theorem}
\begin{proof}
		We may assume that $t$ is rational as a consequence of \cref{thm:Okoucont}. Similarly, as in the proof of \cref{prop:suba},
		we may reduce to the case where both $\varphi$ and $\psi$ have analytic singularities.
	In this case, let $N>0$ be an integer such that $Nt$ is an integer.
    Then for any $s\in \mathrm{H}^0(X,L^k\otimes \mathcal{I}_{\infty}(k\varphi))$ and $r\in \mathrm{H}^0(X,L^k\otimes \mathcal{I}_{\infty}(k\psi))$, we have
    \[
    s^{tN}\otimes r^{N-tN}\in \mathrm{H}^0\left(X,L^{kN}\otimes \mathcal{I}_{\infty}\left(kN(t\varphi+(1-t)\psi)\right)\right).
    \]
	By \cref{thm:HausOkoun} and \cref{rmk:DeltaanaW0}, \eqref{eq:Deltaconcave} follows.
\end{proof}

\begin{proposition}\label{prop:res}
	For any $a\in \mathbb{Z}_{>0}$,
	\[
		\Delta_{\nu}(a\theta,a\varphi)=a\Delta_{\nu}(\theta,\varphi).
	\]
\end{proposition}
\begin{proof}
As in the proof of \cref{prop:suba}, we may assume that $\varphi$ has log singularities. Using \cref{rmk:DeltaanaW0}, we reduce to the case without the singularities, which is proved in
\cref{prop:Okourescaling}.
\end{proof}

In particular, if $S$ is a closed positive $(1,1)$-current on $X$ with $\int_X S^n>0$ and such that
\[
\{S\}\in \NS^1(X)_{\mathbb{Q}},
\]
we can define
\begin{equation}\label{eq:DeltanuTalgebraic1}
\Delta_{\nu}(S)\coloneqq a^{-1}\Delta_{\nu}(aS)
\end{equation}
for a sufficiently divisible positive integer $a$. This definition is independent of the choice of $a$ thanks to \cref{prop:res}.

We also need the following perturbation. Let $A$ be an ample line bundle on $X$. Fix a Hermitian metric $h_A$ on $A$ such that $\omega\coloneqq c_1(A,h_A)$ is a Kähler form on $X$.

\begin{proposition}\label{prop:Deltapert}
    As $\mathbb{Q}_{>0}\ni\delta\searrow 0$, the convex bodies $\Delta_{\nu}(\theta+\delta\omega+\ddc\varphi)$ are decreasing and
    \[
    \Delta_{\nu}(\theta+\delta\omega+\ddc\varphi)\xrightarrow{d_{\Haus}} \Delta_{\nu}(\theta_{\varphi}).
    \]
\end{proposition}

\begin{proof}
	Let $0\leq \delta<\delta'$ be two rational numbers. Take $C\in \mathbb{N}_{>0}$ divisible enough, so that $C\delta$ and $C\delta'$ are both integers. Then by \cref{prop:suba},
	\[
		\Delta_{\nu}(C\theta+C\delta\omega+C\ddc\varphi)\subseteq \Delta_{\nu}(C\theta+C\delta'\omega+C\ddc \varphi).
	\]
	It follows that
	\[
		\Delta_{\nu}(\theta+\delta\omega+\ddc\varphi)\subseteq \Delta_{\nu}(\theta+\delta'\omega+\ddc\varphi).
	\]
	On the other hand,
	\[
		\vol   \Delta_{\nu}(\theta+\delta\omega+\ddc\varphi)=\frac{1}{n!}\vol (\theta+\delta\omega)_{\varphi}.
	\]
	As $\delta\to 0+$, the right-hand side converges to
	\[
		\vol \Delta_{\nu}(\theta,\varphi)=\frac{1}{n!}\vol \theta_{\varphi},
	\]
    thanks to \cref{prop:vollinearlimit}.
	Our assertion therefore follows.
\end{proof}

\subsection{The Hausdorff convergence property}
\label{subsec:HCP}
Let $T$ be a holomorphic line bundle on $X$. The goal of this section is to prove the following:
\begin{theorem}\label{thm:HCP}
    As $k\to\infty$, we have
    \begin{equation}\label{eq:HCP1}
    \Delta_{k,T}(\theta,\varphi)\xrightarrow{d_{\Haus}} \Delta_{\nu}(\theta,\varphi).
    \end{equation}
\end{theorem}
Recall that $\Delta_{k,T}(\theta,\varphi)$ is defined in \cref{def:DeltakT}. Note that it is defined when $k$ is large enough.

Although we are only interested in the untwisted case, the proof given below requires the twisted case.

We first observe that the sequence $\Delta_{k,T}(\theta,\varphi)$ is uniformly bounded: This follows easily from \cref{prop:Deltaconvtwisted}. So Blaschke's selection theorem \cref{thm:Blaschke} is applicable. We will apply this observation without further comments.

\begin{lemma}\label{lma:twistedHcp}
    Assume that $\varphi$ has analytic singularities and $\theta_{\varphi}$ is a Kähler current. Then \eqref{eq:HCP1} holds.
\end{lemma}
\begin{proof}
    Up to replacing $X$ by a modification and twisting $T$ accordingly, we may assume that $\varphi$ has log singularities along an effective $\mathbb{Q}$-divisor $D$, see \cref{prop:birinvO} and \cref{thm:resolvelogsing}.

    Take $\epsilon\in \mathbb{Q}\cap (0,1)$.
    In this case, for large enough $k\in \mathbb{Z}_{>0}$ we have
    \[
    \begin{split}
        \mathrm{H}^0\left(X,T\otimes L^k\otimes \mathcal{I}_{\infty}(k\varphi)\right) \subseteq & \mathrm{H}^0\left(X,T\otimes L^k\otimes \mathcal{I}(k\varphi)\right) \\
        \subseteq & \mathrm{H}^0\left(X,T\otimes L^k\otimes \mathcal{I}_{\infty}(k(1-\epsilon)\varphi)\right).
    \end{split}
    \]
    Take an integer $N\in \mathbb{Z}_{>0}$ so that $ND$ is a divisor and $N\epsilon$ is an integer.

    Let $\Delta'$ be the limit of a subsequence of $(\Delta_{k,T}(\theta,\varphi))_k$, say the sequence defined by the indices $k_1,k_2,\ldots$. Thanks to \cref{thm:Blaschke}, it suffices to show that $\Delta'=\Delta_{\nu}(\theta,\varphi)$.

    There exists $t\in \{0,1,\ldots,N-1\}$ such that $k_i\equiv t$ modulo $N$ for infinitely many $i$, up to replacing $(k_i)_i$ by a subsequence, we may assume that $k_i\equiv t$ modulo $N$ for all $i$. Write $k_i=N g_i+t$.
    Then for large enough $i$, we have
    \[
    \begin{split}
        \mathrm{H}^0\left(X,T\otimes L^{-N+t} \otimes L^{N(g_i+1)}\otimes \mathcal{I}_{\infty}\left(N(g_i+1)\varphi\right)\right) \subseteq \mathrm{H}^0\left(X,T\otimes L^{k_i}\otimes \mathcal{I}(k_i\varphi)\right) \\
        \subseteq \mathrm{H}^0\left(X,T\otimes L^t \otimes L^{Ng_i}\otimes \mathcal{I}_{\infty}\left(g_i N(1-\epsilon)\varphi\right)\right).
    \end{split}
    \]
    So
    \[
    \begin{split}
    (g_i+1) \Delta_{g_i+1,T\otimes L^{-N+t}}(NL-ND)+N(g_i+1)\nu(D)\subseteq (Ng_i+t)\Delta_{k_i,T}(\theta,\varphi)\\
    \subseteq
    g_i\Delta_{g_i,T\otimes L^t}\left(NL-N(1-\epsilon)D\right)+Ng_i(1-\epsilon)\nu(D).
    \end{split}
    \]
    Letting $i\to\infty$, by \cref{prop:Deltaconvtwisted},
    \[
        \Delta_{\nu}(L-D)+\nu(D)\subseteq \Delta'\subseteq \Delta_{\nu}(L-(1-\epsilon)D)+(1-\epsilon)\nu(D).
    \]
    Letting $\epsilon\to 0+$, using the continuity of classical Okounkov bodies in the big cone recalled in \cref{thm:Okounkovtranmain}, we find that
    \[
        \Delta'=\Delta_{\nu}(L-D)+\nu(D)=\Delta_{\nu}(\theta,\varphi),
    \]
    where we applied \cref{rmk:DeltaanaW0} as well. Our assertion follows.
\end{proof}

\begin{lemma}\label{lma:Hausconvbetato0}
    Assume that $\theta_{\varphi}$ is a Kähler current. Then as $\mathbb{Q}\ni\beta\to 0+$, we have
    \[
        \Delta_{\nu}\Bigl((1-\beta)\theta,\varphi\Bigr)\xrightarrow{d_{\Haus}} \Delta_{\nu}(\theta,\varphi).
    \]
\end{lemma}
Here and in the sequel, $\Delta_{\nu}((1-\beta)\theta,\varphi)=\Delta_{\nu}((1-\beta)\theta+\ddc\varphi)$ is defined in \eqref{eq:DeltanuTalgebraic1}.
\begin{proof}
    By \cref{prop:suba} and \cref{prop:IcompimplyDeltacomp}, we have
    \[
        \Delta_{\nu}\Bigl((1-\beta)\theta,\varphi\Bigr)+\beta\Delta_{\nu}(L)\subseteq \Delta_\nu\left(\theta,\varphi+\beta V_\theta\right)\subseteq \Delta_{\nu}(\theta,\varphi)
    \]
    for any small enough $\beta\in \mathbb{Q}_{>0}$.
    In particular, if $\Delta'$ is the Hausdorff limit of a subnet of $(\Delta((1-\beta)\theta,\varphi))_{\beta}$, then $\Delta'\subseteq \Delta_{\nu}(\theta,\varphi)$.
    But
    \[
    \begin{split}
        \vol \Delta'=&\lim_{\beta\to 0+}\vol \Delta_{\nu}((1-\beta)\theta,\varphi)\\
        =&\frac{1}{n!}\lim_{\beta\to 0+}\int_X \left((1-\beta)\theta+\ddc P_{(1-\beta)\theta}[\varphi]_{\mathcal{I}}\right)^n\\
        =&\frac{1}{n!}\int_X (\theta+\ddc P_{\theta}[\varphi]_{\mathcal{I}})^n.
    \end{split}
    \]
    Since we have not developed the theory of nef b-divisors yet, we give a direct proof as well. Take a Kähler form $\omega$ so that $\theta_{\varphi}\geq \omega$. Let $\psi=P_{\theta-\omega}[\varphi]_{\mathcal{I}}$. Then $\varphi\sim_{\mathcal{I}}\psi$. In order to establish the last equality, we may replace $\varphi$ by $\psi$ and hence assume that $\varphi$ is $\mathcal{I}$-good. For small enough $\beta>0$, the form $(1-\beta)\theta-(\theta-\omega)=\omega-\beta\theta$ is Kähler, so this remains true in the class $(1-\beta)\{\theta\}$ by \cref{lma:Igoodinsenspert}. Thus the equality of masses in an $\mathcal{I}$-equivalence class reduces the desired equality to
    \[
    \lim_{\beta\to 0+}\int_X ((1-\beta)\theta+\ddc \varphi)^n=\int_X \theta_{\varphi}^n,
    \]
    which is obvious.

    It follows that $\Delta'=\Delta_{\nu}(\theta,\varphi)$. We conclude by \cref{thm:Blaschke}.
\end{proof}

\begin{proof}[Proof of \cref{thm:HCP}]
Fix a very ample line bundle $H$ on $X$ and a Kähler form $\omega\in c_1(H)$ such that $\omega\geq \theta$.

\textbf{Step~1}. We first handle the case where $\theta_{\varphi}$ is a Kähler current, say $\theta_{\varphi}\geq 2\delta\omega$ for some $\delta\in (0,1)$. Take a quasi-equisingular approximation $(\varphi_j)_j$ of $\varphi$ in $\PSH(X,\theta)$. We may assume that $\theta_{\varphi_j}\geq \delta\omega$ for all $j\geq 1$. After adding constants, we also assume that $\varphi\leq 0$ and $\varphi_j\leq 0$ for all $j$.

Let $\Delta'$ be a limit of a subsequence of $(\Delta_{k,T}(\theta,\varphi))_k$. Let us say the indices of the subsequence are $k_1<k_2<\cdots$. By \cref{thm:Blaschke}, it suffices to show that $\Delta'=\Delta_{\nu}(\theta,\varphi)$.

Observe that for each $j\geq 1$, we have $\Delta'\subseteq \Delta_{\nu}(\theta,\varphi_j)$ by \cref{lma:twistedHcp}. Letting $j\to \infty$, we find $\Delta'\subseteq \Delta_{\nu}(\theta,\varphi)$ as a consequence of \cref{thm:Okoucont}.
Therefore, it suffices to prove that
\[
    \vol \Delta'\geq \vol \Delta_{\nu}(\theta,\varphi).
\]
Fix an integer $N>\delta^{-1}$. Observe that for any $j\geq 1$, we have $\varphi_j\in \PSH(X,(1-N^{-1})\theta)$. Similarly, $\varphi\in \PSH(X,(1-N^{-1})\theta)$.
By \cref{lma:Hausconvbetato0}, it suffices to argue that
\begin{equation}\label{eq:volDeltatoprove}
\vol \Delta'\geq \vol \Delta_{\nu}\left((1-N^{-1})\theta,\varphi\right).
\end{equation}

\textbf{Step~1.1}. We first reduce to the case where $N|k_i$ for all $i$.

We are free to replace $(k_i)_i$ by a subsequence, so we may assume that $k_i\equiv a$ modulo $N$ for all $i\geq 1$, where $a\in \{0,1,\ldots,N-1\}$. We write $k_i=g_i N+a$.
Observe that for each $i\geq 1$,
\[
    \mathrm{H}^0\left(X,T\otimes L^{k_i}\otimes \mathcal{I}(k_i\varphi)\right)\supseteq \mathrm{H}^0\left(X,T\otimes L^{-N+a}\otimes L^{g_i N+ N}\otimes \mathcal{I}\left((g_iN+N)\varphi \right)\right).
\]
Up to replacing $T$ by $T\otimes L^{-N+a}$, we may therefore assume that $a=0$, so that $k_i=g_i N$.

\textbf{Step~1.2}. Write $k_i=g_i N$ for all $i$. We prove \eqref{eq:volDeltatoprove}.

Since
\[
    \int_X\theta_{\varphi_j}^n\to \int_X\theta_{P_\theta[\varphi]_{\mathcal I}}^n,
\]
we can apply \cref{lma:pathoenvelope} to find $j'\in \mathbb{Z}_{>0}$ such that for all $j\geq j'$, there is $\psi_j\in \PSH(X,\theta)_{>0}$ satisfying
\begin{equation}\label{eq:PvarphiIgeq1mNm1_temp1}
P_{\theta}[\varphi]_{\mathcal{I}}\geq \left(1-N^{-1}\right)\varphi_j+N^{-1} \psi_j.
\end{equation}
Fix $j\geq j'$. It suffices to show that
\begin{equation}\label{eq:DeltatransinDeltaprime}
    \Delta_{\nu}\left((1-N^{-1})\theta,\varphi_{j}\right)+v'\subseteq \Delta'
\end{equation}
for some $v'\in \mathbb{R}^n$. In fact, if this is true, we have
\[
\vol \Delta'\geq \vol \Delta_{\nu}\left((1-N^{-1})\theta,\varphi_j\right).
\]
Letting $j\to\infty$ and applying \cref{cor:dSconv_changetheta} and \cref{thm:Okoucont}, we conclude \eqref{eq:volDeltatoprove}.

It remains to prove \eqref{eq:DeltatransinDeltaprime}.

By \cref{lma:existsecposmass}, there is $g'>0$ such that for any $g\geq g'$, we can find a non-zero section
\[
    s_g\in \mathrm{H}^0\left(X,L^g\otimes\mathcal I(g\psi_j)\right).
\]
By \cref{lma:IandIinf}, after possibly increasing $g'$, for all $g\geq g'$,
\[
    \mathcal I\left(gN\varphi_j\right)\subseteq \mathcal I_\infty\left(g(N-1)\varphi_j\right).
\]
Moreover, the inequality \eqref{eq:PvarphiIgeq1mNm1_temp1} gives that for $g\geq g'$,
\[
\begin{aligned}
    \mathcal I(gN\varphi_j)\cdot\mathcal I(g\psi_j)\subseteq & \mathcal I_\infty\left(g(N-1)\varphi_j\right)\cdot\mathcal I(g\psi_j)\\
    \subseteq & \mathcal I\Bigl(g(N-1)\varphi_j+g\psi_j\Bigr)\\
    \subseteq & \mathcal I\Bigl(gN P_\theta[\varphi]_{\mathcal I}\Bigr)\\
    =&\mathcal I(gN\varphi),
\end{aligned}
\]
where the second inclusion follows directly from the definition of $\mathcal I_\infty$. Thus multiplication by $s_g$ gives an injective linear map
\[
    \mathrm{H}^0\left(X,T\otimes L^{g(N-1)}\otimes\mathcal I(gN\varphi_j)\right)
    \xrightarrow{\times s_g}
    \mathrm{H}^0\left(X,T\otimes L^{gN}\otimes\mathcal I(gN\varphi)\right).
\]
In particular, when $g=g_i$ for some $i$ large enough, we then find
\[
    \Delta_{g_i,T}\left((N-1)\theta,N\varphi_{j}\right)+(g_i)^{-1}\nu(s_{g_i})\subseteq N\Delta_{k_i,T}(\theta,\varphi).
\]
We observe that the $(g_i)^{-1}\nu(s_{g_i})$'s are bounded, since both convex bodies appearing in this equation remain bounded as $i$ varies.
After passing to a subsequence, we may assume that $(g_i)^{-1}\nu(s_{g_i})$ converges. Then \cref{lma:twistedHcp}, followed by division by $N$, gives a vector $v'\in \mathbb{R}^n$ such that \eqref{eq:DeltatransinDeltaprime} holds.

\textbf{Step~2}. Next we handle the general case.

Let $\Delta'$ be the Hausdorff limit of a subsequence of $(\Delta_{k,T}(\theta,\varphi))_k$, say the subsequence with indices $k_1<k_2<\cdots$.
By \cref{thm:Blaschke}, it suffices to prove that $\Delta'=\Delta_{\nu}(\theta,\varphi)$.

Take $\psi\in \PSH(X,\theta)$ such that $\theta_{\psi}$ is a Kähler current and $\psi\leq \varphi$. The existence of $\psi$ follows from \cref{lma:kahcurrentposmass}.

    Then for any $\epsilon\in \mathbb{Q}\cap (0,1)$,
    \[
        \Delta_{k,T}(\theta,\varphi)\supseteq \Delta_{k,T}\Bigl(\theta,(1-\epsilon)\varphi+\epsilon\psi \Bigr)
    \]
    for all $k\geq 1$. It follows from Step~1 that
    \[
        \Delta' \supseteq \Delta_{\nu}\Bigl(\theta,(1-\epsilon)\varphi+\epsilon\psi\Bigr).
    \]

    As $\epsilon\to 0+$, we have
    \[
        (1-\epsilon)\varphi+\epsilon\psi\xrightarrow{d_S}\varphi
    \]
    by \cref{lma:dSsmallmult}. Hence \cref{thm:Okoucont} gives
    \[
        \Delta_{\nu}\Bigl(\theta,(1-\epsilon)\varphi+\epsilon\psi\Bigr)\to\Delta_{\nu}(\theta,\varphi)
    \]
    in the Hausdorff metric. It follows that $\Delta' \supseteq \Delta_{\nu}(\theta,\varphi)$.

    It remains to establish that
    \begin{equation}\label{eq:Deltapvolumeupp}
        \vol \Delta'\leq \vol \Delta_{\nu}(\theta,\varphi).
    \end{equation}
    Fix an integer $N>0$, it suffices to argue that
    \begin{equation}\label{eq:Deltapvolumeupp2}
        \vol \Delta'\leq \frac{1}{n!} \int_X \left(N^{-1}\omega+\theta+\ddc P_{N^{-1}\omega+\theta}[\varphi]_{\mathcal{I}}\right)^n.
    \end{equation}
    Assuming this, letting $N\to\infty$ and using \cref{lma:perturbvol} together with \itemref{itm:prop:vollinearlimit:6}, we conclude \eqref{eq:Deltapvolumeupp}.

    After adding a constant to $\varphi$, which does not change the multiplier ideals or the current $\theta+\ddc\varphi$, we may assume that $\varphi\leq 0$.

    \textbf{Step~2.1}. We first reduce to the case $N|k_i$ for all $i$.

    For this purpose, we are free to replace $k_1<k_2<\cdots$ by a subsequence. We may then assume that $k_i\equiv a$ modulo $N$ for all $i\geq 1$ for some $a\in \{0,1,\ldots,N-1\}$. We write $k_i=g_iN +a$. Observe that
    \[
        \mathrm{H}^0\left(X,T\otimes L^{k_i}\otimes \mathcal{I}(k_i \varphi)\right)\subseteq \mathrm{H}^0\left(X, T\otimes L^a \otimes L^{g_i N}\otimes \mathcal{I}(g_i N \varphi) \right).
    \]
    Up to replacing $T$ by $T\otimes L^a$, we may assume that $a=0$.

    \enlargethispage*{15pt}
    \textbf{Step~2.2}. We write $k_i=g_i N$ for all $i$.

    Take a non-zero section $s\in \mathrm{H}^0(X,H)$.

    We have an injective linear map
    \[
        \mathrm{H}^0\left(X,T\otimes L^{gN}\otimes \mathcal{I}(gN\varphi)\right)\xrightarrow{\times s^g}\mathrm{H}^0\left(X,T\otimes H^g\otimes L^{gN}\otimes \mathcal{I}(gN\varphi)\right)
    \]
    for each $g\geq 1$.
    In particular, for each $i\geq 1$,
    \[
        k_i\Delta_{k_i,T}(\theta,\varphi)+g_i\nu(s)\subseteq g_i\Delta_{g_i,T}\left(\omega+N\theta,N\varphi\right).
    \]
    Letting $i\to\infty$, by Step~1, we have
    \[
        N\Delta' + \nu(s) \subseteq \Delta_{\nu}(\omega+N\theta,N\varphi).
    \]
    So
    \[
        \vol \Delta'\leq \vol \Delta_{\nu}\left(N^{-1}\omega+\theta,\varphi\right)=\frac{1}{n!}\int_X \left(N^{-1}\omega+\theta+\ddc P_{N^{-1}\omega+\theta}[\varphi]_{\mathcal{I}}\right)^n.
    \]
    This completes the proof.
\end{proof}

\subsection{Recover Lelong numbers from partial Okounkov bodies}
\label{subsec:pob_lelong}

\pagebreak[4]
\begin{theorem}\label{thm:nuOk}
	Let $E$ be a prime divisor on $X$. Let $Y_{\bullet}$ be an admissible flag with $E=Y_1$. Then
	\begin{equation}\label{eq:numinOk}
		\nu(\varphi,E)=\min_{x\in \Delta_{Y_{\bullet}}(\theta,\varphi)}x_1.
	\end{equation}
\end{theorem}
Here $x_1$ denotes the first component of $x$.

\begin{figure}[ht]
\centering
\begin{tikzpicture}[scale=0.85, line width=0.65pt, line cap=round, line join=round, every node/.style={font=\small}]
    \def\nuval{1.42}

    \draw[-{Latex[length=2.4mm]}] (0,0) -- (7.2,0) node[below right] {$x_1$};
    \draw[-{Latex[length=2.4mm]}] (0,0) -- (0,5.1) node[above] {$(x_2,\ldots,x_n)$};

    \path[fill=black!3, draw=black]
        (1.05,1.02) -- (2.55,4.20) -- (6.05,1.50) -- (4.15,0.18) -- cycle;

    \path[fill=red!7, draw=red!70!black]
        (\nuval,1.80) -- (2.55,4.20) -- (4.40,2.77) -- (4.22,1.85)
        .. controls (3.55,1.28) and (2.40,1.07) .. (\nuval,1.80) -- cycle;

    \begin{scope}
        \clip (\nuval,1.80) -- (2.55,4.20) -- (4.40,2.77) -- (4.22,1.85)
            .. controls (3.55,1.28) and (2.40,1.07) .. (\nuval,1.80) -- cycle;
        \foreach \x in {1.55,1.95,...,4.35} {
            \draw[red!70!black] (\x,0.85) -- ++(-0.45,3.05);
        }
    \end{scope}
    \draw[red!70!black]
        (\nuval,1.80) -- (2.55,4.20) -- (4.40,2.77) -- (4.22,1.85)
        .. controls (3.55,1.28) and (2.40,1.07) .. (\nuval,1.80) -- cycle;

    \draw[blue!70!black, densely dashed] (\nuval,0) -- (\nuval,1.80);
    \draw[blue!70!black, decorate, decoration={brace, mirror, amplitude=4pt}]
        (0,-0.28) -- (\nuval,-0.28) node[midway, below=5pt] {$\nu(\varphi,Y_1)$};

    \draw[red!70!black] (4.58,2.18) -- (4.18,2.14);
    \node[red!70!black, anchor=west] at (4.62,2.18) {$\Delta_{Y_{\bullet}}(\theta,\varphi)$};
    \node at (5.65,3.10) {$\Delta_{Y_{\bullet}}(L)$};
\end{tikzpicture}
\caption{Partial Okounkov body and Lelong number}
\label{fig:nuOk}
\end{figure}
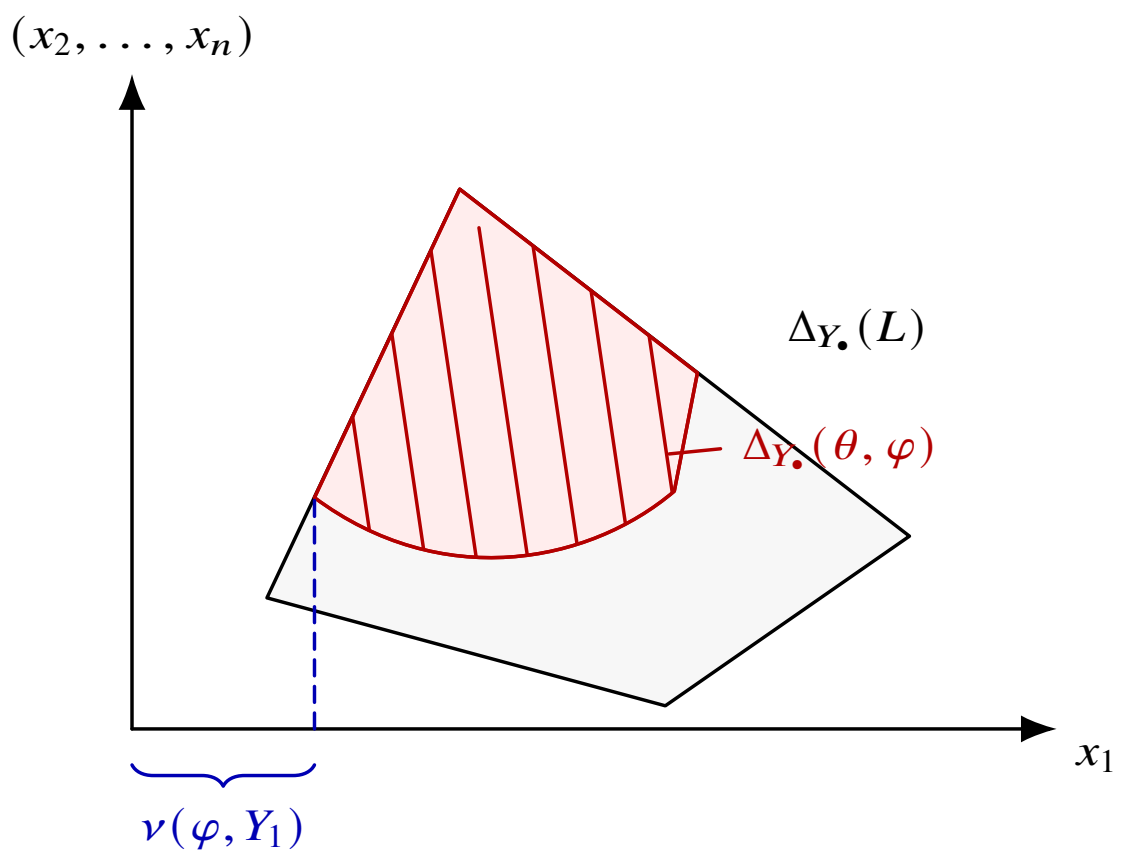

The picture in \cref{fig:nuOk} illustrates the formula: the left supporting hyperplane of $\Delta_{Y_{\bullet}}(\theta,\varphi)$ is given by $x_1=\nu(\varphi,Y_1)$.

\begin{proof}
Replacing $\varphi$ by $P_{\theta}[\varphi]_{\mathcal{I}}$, we may assume that $\varphi$ is $\mathcal{I}$-good.

\textbf{Step~1}. We first reduce to the case where $\varphi$ has analytic singularities.

By \cref{thm:charIgoodasclosure}, we can find a sequence $(\varphi_j)_j$ in $\PSH(X,\theta)_{>0}$ with analytic singularities such that $\varphi_j\xrightarrow{d_S}\varphi$. It follows from \cref{thm:Okoucont} that
\[
\Delta_{Y_{\bullet}}(\theta,\varphi_j)\xrightarrow{d_{\Haus}}\Delta_{Y_{\bullet}}(\theta,\varphi).
\]
Therefore, by \cref{lma:Hausdorffconvcoordbounds},
\[
\lim_{j\to\infty}\min_{x\in \Delta_{Y_{\bullet}}(\theta,\varphi_j)}x_1=\min_{x\in \Delta_{Y_{\bullet}}(\theta,\varphi)}x_1.
\]
In view of \cref{thm:Lelongcont}, it suffices to prove \eqref{eq:numinOk} with $\varphi_j$ in place of $\varphi$.

\textbf{Step~2}. Assume that $\varphi$ has analytic singularities. In view of \cref{prop:birinvO} and \cref{thm:resolvelogsing}, after replacing $X$ by a modification, we may assume that $\varphi$ has log singularities along an effective $\mathbb{Q}$-divisor $F$.

Perturbing $L$ by an ample $\mathbb{Q}$-line bundle by \cref{prop:Deltapert}, we may assume that $\theta_{\varphi}$ is a Kähler current. Therefore, $L-F$ is ample by \cref{lma:logsingrem}.
Finally, by rescaling, we may assume that $F$ is a divisor and $L$ is a line bundle.

	By \cref{thm:HCP}, we know that
	\[
		\min_{x\in \Delta_{Y_{\bullet}}(\theta,\varphi)}x_1=\adjustlimits\lim_{k\to\infty}\min_{x\in \Delta_k(\theta,\varphi)}x_1.
	\]
	By definition,
	\[
		\min_{x\in \Delta_k(\theta,\varphi)}x_1=k^{-1}\ord_E \mathrm{H}^0\left(X,L^k\otimes \mathcal{I}(k\varphi)\right).
	\]
	In view of \cref{prop:Lelongnumfrommis}, it remains to show that
	\begin{equation}\label{eq:temp1}
		\lim_{k\to\infty}k^{-1}\ord_E \mathrm{H}^0\left(X,L^k\otimes \mathcal{I}(k\varphi)\right)=\lim_{k\to\infty}k^{-1}\ord_E \mathcal{I}(k\varphi).
	\end{equation}
	The $\geq$ direction is trivial, we prove the converse. Observe that
	\[
		\mathrm{H}^0\left(X,L^k\otimes \mathcal{I}(k\varphi)\right)=\mathrm{H}^0\left(X,L^k\otimes \mathcal{O}_X(-kF)\right),\quad \mathcal{I}(k\varphi)=\mathcal{O}(-kF).
	\]
	As $L-F$ is ample, for large enough $k$, the line bundle $k(L-F)$ is globally generated. Hence we can choose a section of $k(L-F)$ which does not vanish at the generic point of $E$. Viewed as a section of $L^k\otimes\mathcal O_X(-kF)$, it has vanishing order exactly $\ord_E(kF)$ along $E$. Therefore,
	\[
		\ord_E  \mathrm{H}^0\left(X,L^k\otimes \mathcal{O}_X(-kF)\right)=\ord_E(kF).
	\]
	Thus, \eqref{eq:temp1} follows.
\end{proof}

\begin{corollary}\label{cor:Deltacontimplyvarphi}
	Let $\varphi,\psi\in \PSH(X,\theta)_{>0}$. If
	\[
		\Delta_{W_{\bullet}}(\pi^*\theta,\pi^*\varphi)\subseteq \Delta_{W_{\bullet}}(\pi^*\theta,\pi^*\psi)
	\]
	for all modifications $\pi\colon Y\rightarrow X$ and all admissible flags $W_{\bullet}$ on $Y$, then $\varphi\preceq_{\mathcal{I}} \psi$.
\end{corollary}
\begin{proof}
    By \cref{thm:nuOk}, the assumed inclusions imply that for every prime divisor $E$ on every modification $\pi$,
    \[
    \nu(\pi^*\varphi,E)\geq \nu(\pi^*\psi,E).
    \]
    Hence $\varphi\preceq_{\mathcal{I}}\psi$ by \cref{prop:Iequivchar2}.
\end{proof}
\begin{corollary}\label{cor:numin}
	Let $E$ be a prime divisor over $X$. Then
	\begin{equation}\label{eq:nuVthetaE}
		\nu(V_{\theta},E)=\lim_{k\to\infty}\frac{1}{k}\ord_E \mathrm{H}^0(X,L^k).
	\end{equation}
\end{corollary}

\begin{proof}
    After replacing $X$ by a modification, we may assume that $E$ is a prime divisor on $X$.
    Note that since $X$ is projective, we can always find an admissible flag with $E$ being the first component.

	Then \eqref{eq:nuVthetaE} follows from \cref{thm:nuOk} and the fact that $\Delta_{Y_{\bullet}}(\theta,V_{\theta})=\Delta_{Y_{\bullet}}(L)$ for any admissible flag $Y_{\bullet}$ on $X$.
\end{proof}

\section{Transcendental partial Okounkov bodies}\label{sec:tpob}
Let $X$ be a connected compact Kähler manifold of dimension $n>0$.
Fix a smooth flag $Y_{\bullet}$ on $X$. We will extend the theory of partial Okounkov bodies in the previous section to the transcendental setting.

\subsection{The traditional approach to the Okounkov body problem}
\label{subsec:trad_approach}

\begin{definition}\label{def:tradOk}
Let $\alpha$ be a big cohomology class on $X$. We define the \emph{Okounkov body}\index{Okounkov body} of $\alpha$ with respect to the flag $Y_{\bullet}$ as
    \begin{equation}\label{eq:twodefspob}
        \Delta_{Y_{\bullet}}(\alpha)\coloneqq \overline{\left\{\nu_{Y_{\bullet}}(S): S\in \mathcal{Z}_+(X,\alpha),S\textup{ has gentle analytic singularities} \right\}}.
    \end{equation}
    \index{$\Delta_{Y_{\bullet}}(\alpha)$}
\end{definition}
See \cref{def:analy-sing} for the definition of gentle analytic singularities.

The results of \cite{DRWNXZ} can be summarized as follows:
\begin{theorem}\label{thm:Okounkovtranmain}
    For any big cohomology class $\alpha$ on $X$, the set $\Delta_{Y_{\bullet}}(\alpha)\subseteq \mathbb{R}^n$ is a convex body satisfying the following properties:
    \begin{enumerate}
        \item\label{itm:thm:Okounkovtranmain:1} We have
        \[
        \vol \Delta_{Y_{\bullet}}(\alpha)=\frac{1}{n!}\vol \alpha;
        \]
        \item\label{itm:thm:Okounkovtranmain:2} given another big cohomology class $\alpha'$ on $X$, we have
        \[
        \Delta_{Y_{\bullet}}(\alpha)+\Delta_{Y_{\bullet}}(\alpha')\subseteq \Delta_{Y_{\bullet}}(\alpha+\alpha');
        \]
        \item\label{itm:thm:Okounkovtranmain:3} let $\pi\colon Y\rightarrow X$ be a proper bimeromorphic morphism where $Y$ is a Kähler manifold. Assume that $(W_{\bullet},g)$ is the lifting of $Y_{\bullet}$ to $Y$. Then
        \[
        \Delta_{W_{\bullet}}(\pi^*\alpha)=\Delta_{Y_{\bullet}}(\alpha)g;
        \]
        \item\label{itm:thm:Okounkovtranmain:4} the map $\alpha\mapsto \Delta_{Y_{\bullet}}(\alpha)$ is continuous in the big cone with respect to the Hausdorff metric;
        \item\label{itm:thm:Okounkovtranmain:5} if $\alpha$ is Kähler, then for all sufficiently small $t>0$ we have
        \[
            \left\{ y\in \mathbb{R}^{n-1}: (t,y)\in \Delta_{Y_{\bullet}}(\alpha) \right\}=\Delta_{Y_1\supseteq \cdots \supseteq Y_n}\Bigl((\alpha-t\{Y_1\})|_{Y_1}\Bigr).
        \]
    \end{enumerate}
\end{theorem}
See \cref{def:lifting} for the notion of lifting. The proofs of \ref{itm:thm:Okounkovtranmain:1}--\ref{itm:thm:Okounkovtranmain:4} require techniques not covered here.

For \ref{itm:thm:Okounkovtranmain:5}, choose a small rational $s>0$ such that $\alpha-t\{Y_1\}$ is Kähler for $0\leq t\leq s$. A Kähler form in $\alpha$ has first valuation coordinate $0$, and adding $s[Y_1]$ to a Kähler form in $\alpha-s\{Y_1\}$ gives first coordinate $s$. Thus $0<t<s$ lies in the interior of the projection of the body. In dimension one the two bodies in the slice formula are points.

In dimension at least two, take rational $t\in(0,s)$. By \cite[Lemma~7.4 in arXiv:2309.07584v3]{DRWNXZ}, the slice is the closure of the remaining coordinates of valuations of Kähler currents in $\alpha$ with analytic singularities and first coordinate $t$. Subtracting $t[Y_1]$ from such a current leaves a Kähler current with analytic singularities and generic Lelong number zero along $Y_1$. Its restriction to $Y_1$ gives the remaining coordinates. Conversely, a Kähler current with gentle analytic singularities on $Y_1$ in $(\alpha-t\{Y_1\})|_{Y_1}$ extends to a Kähler current with analytic singularities in $\alpha-t\{Y_1\}$ by \cref{thm:CT-thm-refined'}: the exponential of its potential is H\"older continuous by \cref{def:analy-sing}. Adding $t[Y_1]$ gives a current with the desired valuation. Here the closures of valuations for analytic and gentle analytic singularities agree, as in \cite[Section~3]{DRWNXZ}. The formula extends to real $t\in(0,s)$ by continuity of interior slices and \ref{itm:thm:Okounkovtranmain:4} on $Y_1$.

\subsection{Definitions of partial Okounkov bodies}
\label{subsec:tpob_def}
Let $\theta$ be a closed real smooth $(1,1)$-form on $X$ representing a big cohomology class $\alpha$.

Let $T=\theta_{\varphi}\in \mathcal{Z}_+(X,\alpha)$. We shall define a convex body $\Delta_{Y_{\bullet}}(T)\subseteq \mathbb{R}^n$\index{$\Delta_{Y_{\bullet}}(T)$}, which is also written as $\Delta_{Y_{\bullet}}(\theta,\varphi)$. This convex body is called the \emph{partial Okounkov body}\index{partial Okounkov body} of $T$ with respect to the flag $Y_{\bullet}$.

\subsubsection{The case of analytic singularities}
\label{subsubsec:tpob_analytic}

\begin{definition}\label{def:POBanalsing}
    When $T$ is a Kähler current with analytic singularities, we take a modification $\pi\colon Y\rightarrow X$ so that
    \index{partial Okounkov body!analytic singularities case}
    \begin{enumerate}
        \item\label{itm:def:POBanalsing:1}
        \begin{equation}\label{eq:resolveanalytic}
            \pi^*T=[D]+R,
        \end{equation}
        where $D$ is an effective $\mathbb{Q}$-divisor on $Y$ and $R$ is a closed positive $(1,1)$-current with bounded potential, and
        \item\label{itm:def:POBanalsing:2} the lifting $(Z_{\bullet},g)$ of $Y_{\bullet}$ to $Y$ exists.
    \end{enumerate}
Define
\begin{equation}\label{eq:POBanalsing}
\Delta_{Y_{\bullet}}(T)\coloneqq\Delta_{Z_{\bullet}}(\{R\})g^{-1}+\nu_{Z_{\bullet}}([D])g^{-1}.
\end{equation}
\end{definition}
The existence of $\pi$ is guaranteed by \cref{thm:resolvelogsing} and \cref{thm:liftableflag}.

\begin{lemma}\label{lma:Okounkov:L1426}
    The convex body $\Delta_{Y_{\bullet}}(T)$ defined in \cref{def:POBanalsing} is independent of the choice of $\pi$.
\end{lemma}
\begin{proof}
    Take another map $\pi'\colon Y'\rightarrow X$ with the same properties. We want to show that $\pi$ and $\pi'$ define the same $\Delta_{Y_{\bullet}}(T)$. We may assume that $\pi'$ dominates $\pi$ through $p\colon Y'\rightarrow Y$, so that we have a commutative diagram
    \[
    \begin{tikzcd}
Y' \arrow[rr, "p"] \arrow[rd, "\pi'"'] &   & Y \arrow[ld, "\pi"] \\
                                       & X. &
\end{tikzcd}
    \]
    We take $D$ and $R$ as in \eqref{eq:resolveanalytic}. Then
    \[
    \pi'^*T=[p^*D]+p^*R.
    \]
    Write $(Z_{\bullet},g)$ and $(Z_{\bullet}',g')$ for the liftings of $Y_{\bullet}$ to $Y$ and $Y'$ respectively.
    We need to prove that
    \[
    \Delta_{Z_{\bullet}}(\{R\})g^{-1}+\nu_{Z_{\bullet}}([D])g^{-1}=\Delta_{Z_{\bullet}'}(\{p^*R\})g'^{-1}+\nu_{Z_{\bullet}'}([p^*D])g'^{-1}.
    \]
    By \cref{thm:Okounkovtranmain}, applied to the class $\{R\}$, we have
    \[
        \Delta_{Z'_\bullet}(\{p^*R\})=\Delta_{Z_\bullet}(\{R\})\,\cor(Z_\bullet,p).
    \]
    By \cref{prop:cormatrix}, applied to the divisor current $[D]$, we also have
    \[
        \nu_{Z'_\bullet}([p^*D])=\nu_{Z_\bullet}([D])\,\cor(Z_\bullet,p).
    \]
    Finally, \cref{prop:cormult} gives
    \[
        g'=g\,\cor(Z_\bullet,p).
    \]
    Multiplying by $(g')^{-1}$ gives the desired equality.
\end{proof}

The above proof also describes the bimeromorphic behavior of $\Delta_{Y_{\bullet}}(T)$:
\begin{lemma}\label{lma:liftOkounana}
Let $T\in \mathcal{Z}_{+}(X,\alpha)$ be a Kähler current with analytic singularities.
    Let $\pi\colon Y\rightarrow X$ be a proper bimeromorphic morphism and $(W_{\bullet},g)$ be the lifting of $Y_{\bullet}$ to $Y$. Then
    \[
    \Delta_{W_{\bullet}}(\pi^*T)=\Delta_{Y_{\bullet}}(T)g.
    \]
\end{lemma}

\begin{lemma}\label{lma:Okounkovanalycomp}
    Assume that $T,S\in \mathcal{Z}_+(X,\alpha)$ are two Kähler currents with analytic singularities and $T\preceq S$. Then
    \[
    \Delta_{Y_{\bullet}}(T)\subseteq \Delta_{Y_{\bullet}}(S)\subseteq \Delta_{Y_{\bullet}}(\alpha).
    \]
    Moreover,
    \begin{equation}\label{eq:volpobanaly}
\vol \Delta_{Y_{\bullet}}(T)=\frac{1}{n!} \int_X T^n.
\end{equation}
\end{lemma}
\begin{proof}
    We first show that
    \[
    \Delta_{Y_{\bullet}}(T)\subseteq \Delta_{Y_{\bullet}}(S).
    \]
    Using \cref{lma:liftOkounana}, we may assume that $T$ and $S$ have log singularities along effective $\mathbb{Q}$-divisors $E$ and $F$ respectively. By assumption, $E\geq F$. Replacing $T$ and $S$ by $T-[F]$ and $S-[F]$ respectively, we may assume that $F=0$, cf. \eqref{eq:POBanalsing}.

    In this case, we need to show that
    \[
    \Delta_{Y_{\bullet}}(\alpha)\supseteq \Delta_{Y_{\bullet}}(\alpha-\{E\})+\nu_{Y_{\bullet}}([E]),
    \]
    which is obvious.

    Next we prove that
    \[
    \Delta_{Y_{\bullet}}(T)\subseteq \Delta_{Y_{\bullet}}(\alpha).
    \]
    By \cref{lma:liftOkounana} and \cref{thm:Okounkovtranmain} again, we may assume that $T$ has log singularities. We take $D$ and $R$ as in \eqref{eq:resolveanalytic}. We need to show that
    \[
    \Delta_{Y_{\bullet}}(\alpha-\{D\})+\nu_{Y_{\bullet}}([D])\subseteq \Delta_{Y_{\bullet}}(\alpha),
    \]
    which is again obvious.

    Finally, \eqref{eq:volpobanaly} follows immediately from \cref{thm:Okounkovtranmain}.
\end{proof}

\subsubsection{The case of Kähler currents}
\label{subsubsec:tpob_Kahler}

\begin{definition}\label{def:POBKahcurr}
     Let $T\in \mathcal{Z}_+(X,\alpha)$ be a Kähler current. Take a quasi-equisingular approximation $(T_j)_j$ of $T$ by Kähler currents with analytic singularities in $\mathcal{Z}_+(X,\alpha)$. Then we define
    \index{partial Okounkov body!Kähler current case}
    \[
        \Delta_{Y_{\bullet}}(T)\coloneqq \bigcap_{j=1}^{\infty}\Delta_{Y_{\bullet}}(T_j).
    \]
\end{definition}

To obtain such an approximation, write $T=\theta+\ddc\varphi$ and choose a Kähler form $\omega$ and $\epsilon>0$ with $T\geq 2\epsilon\omega$. Apply \cref{thm:qequi} to $\varphi$ with background $\theta-\epsilon\omega$. On the compact manifold $X$, the error forms are bounded by $\epsilon\omega/2$ for all large $j$. Hence the currents $\theta+\ddc\varphi_j$ are bounded below by $\epsilon\omega/2$. Passing to this tail preserves quasi-equisingularity.

\begin{lemma}\label{lma:Okounkov:L1517}
    The convex body $\Delta_{Y_{\bullet}}(T)$ in \cref{def:POBKahcurr} is independent of the chosen approximation by Kähler currents.
\end{lemma}
In particular, if $T$ also has analytic singularities, then the $\Delta_{Y_{\bullet}}(T)$'s defined in \cref{def:POBKahcurr} and in \cref{def:POBanalsing} coincide.
\begin{proof}
    Let $(S_k)_k$ and $(T_k)_k$ be two approximations as in \cref{def:POBKahcurr}, with potentials $\psi_k$ and $\varphi_k$ in $\PSH(X,\theta)$. By \cref{prop:compqequi}, for any small rational $\epsilon>0$ and $j>0$, we can find $k>0$ so that
    \[
        \psi_k\preceq (1-\epsilon)\varphi_j.
    \]
    By \cref{lma:Okounkovanalycomp},
    \[
    \bigcap_{k=1}^{\infty}\Delta_{Y_{\bullet}}(\theta,\psi_k)\subseteq \Delta_{Y_{\bullet}}(\theta,(1-\epsilon)\varphi_j).
    \]
    On the other hand, observe that
    \[
    \bigcap_{\epsilon\in \mathbb{Q}_{>0}\text{ small enough}}\Delta_{Y_{\bullet}}(\theta,(1-\epsilon)\varphi_j)=\Delta_{Y_{\bullet}}(\theta,\varphi_j).
    \]
    In fact, the $\supseteq$ direction follows from \cref{lma:Okounkovanalycomp}, so it suffices to show that the two sides have the same volume, which follows from \eqref{eq:volpobanaly}.

    It follows that
    \[
    \bigcap_{k=1}^{\infty}\Delta_{Y_{\bullet}}(\theta,\psi_k)\subseteq \bigcap_{j=1}^{\infty}\Delta_{Y_{\bullet}}(\theta,\varphi_j).
    \]
    The other inclusion follows by symmetry.
\end{proof}

The same argument shows that
\begin{corollary}\label{cor:Kahlercurrentcase}
    Suppose that $T,S\in \mathcal{Z}_+(X,\alpha)$ are two Kähler currents satisfying $T \preceq_{\mathcal{I}} S$. Then
    \[
    \Delta_{Y_{\bullet}}(T)\subseteq \Delta_{Y_{\bullet}}(S)\subseteq \Delta_{Y_{\bullet}}(\alpha).
    \]
\end{corollary}

\begin{proposition}\label{prop:Okounkov:L1551}
    Let $T\in \mathcal{Z}_+(X,\alpha)$ be a Kähler current. Then
    \begin{equation}\label{eq:volOkocur}
        \vol \Delta_{Y_{\bullet}}(T)=\frac{1}{n!}\vol T.
    \end{equation}
\end{proposition}
\begin{proof}
Take a quasi-equisingular approximation $(T_j)_j$ by Kähler currents as in \cref{def:POBKahcurr}. The bodies $\Delta_{Y_{\bullet}}(T_j)$ decrease by \cref{lma:Okounkovanalycomp}. By \cref{cor:quasi-equichar}, the potentials converge in $d_S$ to the $\mathcal I$-model representative of the potential of $T$. This representative has the same volume as $T$ by \cref{def:volqpsh}. The formula follows from \eqref{eq:volpobanaly} and \cref{thm:contvolu}.
\end{proof}

\begin{lemma}\label{lma:Okomonotone}
    Let $T\in \mathcal{Z}_+(X,\alpha)$ be a Kähler current and $\omega$ be a Kähler form on $X$. Then
    \begin{equation}\label{eq:DeltaTincreaseomegatemp1}
    \Delta_{Y_{\bullet}}(T)\subseteq \Delta_{Y_{\bullet}}(T+\omega).
    \end{equation}
    Moreover,
    \begin{equation}\label{eq:DeltaTincreaseomegatemp2}
    \Delta_{Y_{\bullet}}(T)=\bigcap_{\epsilon>0}\Delta_{Y_{\bullet}}(T+\epsilon\omega).
    \end{equation}
\end{lemma}
\begin{proof}
We first prove \eqref{eq:DeltaTincreaseomegatemp1}.
    Taking a quasi-equisingular approximation by Kähler currents as in \cref{def:POBKahcurr} (and adding $\omega$ to its terms for $T+\omega$), we reduce to the case where $T$ has analytic singularities. By \cref{lma:liftOkounana}, we may assume that $T$ has log singularities.
    Take $D$ and $R$ as in \eqref{eq:resolveanalytic}. By definition again, it suffices to show that
    \[
    \Delta_{Y_{\bullet}}(\{R\})\subseteq \Delta_{Y_{\bullet}}(\{R\}+\{\omega\}),
    \]
    which is clear by definition.

Next we prove \eqref{eq:DeltaTincreaseomegatemp2}. Thanks to \eqref{eq:DeltaTincreaseomegatemp1}, it remains to prove that both sides have the same volume:
\[
\lim_{\epsilon\to 0+}\vol (T+\epsilon\omega)=\vol T.
\]
This is proved in \cref{prop:vollinearlimit}.
\end{proof}

\subsubsection{The general case}
\label{subsubsec:tpob_general}

\begin{definition}\label{def:generalPOB}
    Let $T\in \mathcal{Z}_+(X,\alpha)$. Take a Kähler form $\omega$ on $X$. We define
    \index{partial Okounkov body!general case}
    \begin{equation}\label{eq:DeltaTgeneral}
        \Delta_{Y_{\bullet}}(T)=\bigcap_{j=1}^{\infty}\Delta_{Y_{\bullet}}(T+j^{-1}\omega).
    \end{equation}
    The same definition makes sense when $\alpha$ is only pseudo-effective.
\end{definition}
This definition is clearly independent of the choice of $\omega$ by \cref{lma:Okomonotone}. Moreover, it extends \cref{def:POBKahcurr} and \cref{def:POBanalsing} as a result of \cref{lma:Okomonotone}.

The main properties of $\Delta_{Y_{\bullet}}(T)$ are summarized as follows:
\begin{theorem}\label{thm:pobmain}
The convex bodies $\Delta_{Y_{\bullet}}(T)$ satisfy the following properties:
\begin{enumerate}
    \item\label{itm:thm:pobmain:1} Suppose that $T\in \mathcal{Z}_+(X,\alpha)_{>0}$. We have
    \begin{equation}\label{eq:volpobgeneral}
        \vol \Delta_{Y_{\bullet}}(T)=\frac{1}{n!}\vol T.
    \end{equation}
    \item\label{itm:thm:pobmain:2} For $T,S \in \mathcal{Z}_+(X,\alpha)$ satisfying $T\preceq_{\mathcal{I}}S$, we have
    \[
        \Delta_{Y_{\bullet}}(T)\subseteq \Delta_{Y_{\bullet}}(S)\subseteq \Delta_{Y_{\bullet}}(\alpha).
    \]
    \item\label{itm:thm:pobmain:3} For any current $T\in \mathcal{Z}_+(X,\alpha)$ with minimal singularities, we have
    \[
    \Delta_{Y_{\bullet}}(T)=\Delta_{Y_{\bullet}}(\alpha).
    \]
    \item\label{itm:thm:pobmain:4} The map $\mathcal{Z}_+(X,\alpha)_{>0}\rightarrow \mathcal{K}_n$ given by $T\mapsto \Delta_{Y_{\bullet}}(T)$ is continuous, where we endow the $d_S$-pseudometric on $\mathcal{Z}_+(X,\alpha)_{>0}$ and the Hausdorff topology on $\mathcal{K}_n$.
    \item\label{itm:thm:pobmain:5} Let $\pi\colon Y\rightarrow X$ be a proper bimeromorphic morphism, where $Y$ is a Kähler manifold. Assume that the lifting $(W_{\bullet},g)$ of $Y_{\bullet}$ to $Y$ exists. Then for any $T\in \mathcal{Z}_+(X,\alpha)_{>0}$, we have
    \[
        \Delta_{W_{\bullet}}(\pi^*T)=\Delta_{Y_{\bullet}}(T)g.
    \]
    \item\label{itm:thm:pobmain:6} For $T,S\in \mathcal{Z}_+(X,\alpha)$, we have
    \begin{equation}\label{eq:pobadditiv}
    \Delta_{Y_{\bullet}}(T)+\Delta_{Y_{\bullet}}(S)\subseteq \Delta_{Y_{\bullet}}(T+S).
    \end{equation}
\end{enumerate}
\end{theorem}

\begin{proof}
\ref{itm:thm:pobmain:1} By \eqref{eq:DeltaTgeneral} and \eqref{eq:volOkocur}, for any Kähler form $\omega$ on $X$,
\[
\vol \Delta_{Y_{\bullet}}(T)=\lim_{j\to\infty}\vol \Delta_{Y_{\bullet}}(T+j^{-1}\omega)=\frac{1}{n!}\lim_{j\to\infty}\vol(T+j^{-1}\omega).
\]
The right-hand side is computed in \cref{prop:vollinearlimit}. Hence, \eqref{eq:volpobgeneral} follows.

\ref{itm:thm:pobmain:2} Fix a Kähler form $\omega$ on $X$. By \cref{cor:Kahlercurrentcase}, for each $j\geq 1$,
\[
\Delta_{Y_{\bullet}}(T+j^{-1}\omega)\subseteq \Delta_{Y_{\bullet}}(S+j^{-1}\omega)\subseteq \Delta_{Y_{\bullet}}(\alpha+j^{-1}\{\omega\}).
\]
It remains to show that
\[
\Delta_{Y_{\bullet}}(\alpha)=\bigcap_{j=1}^{\infty}\Delta_{Y_{\bullet}}(\alpha+j^{-1}\{\omega\}).
\]
The $\subseteq$ direction is clear. Comparing the volumes using \cref{thm:Okounkovtranmain}, we conclude that equality holds.

\ref{itm:thm:pobmain:3} This follows from \ref{itm:thm:pobmain:1} and \ref{itm:thm:pobmain:2}.

\ref{itm:thm:pobmain:4} Let $(T_j)_j$ be a sequence in $\mathcal{Z}_+(X,\alpha)_{>0}$ converging to $T\in \mathcal{Z}_+(X,\alpha)_{>0}$ with respect to $d_S$. We want to show that $\Delta_{Y_{\bullet}}(T_j)\xrightarrow{d_{\Haus}} \Delta_{Y_{\bullet}}(T)$. By \ref{itm:thm:pobmain:2}, replacing the potential of each current by its $\mathcal I$-model representative does not change its partial Okounkov body. Thus we may represent the singularity types by model potentials of uniformly positive mass. Applying \cref{prop:incanddec}, after passing to a subsequence we may sandwich the sequence between an increasing and a decreasing sequence converging to $T$ in $d_S$. Hence it suffices to treat the monotone cases.

In either monotone case, \ref{itm:thm:pobmain:2} gives monotonicity of the convex bodies, while \ref{itm:thm:pobmain:1} and \cref{cor:massdiffdS} give convergence of their volumes to $\vol\Delta_{Y_{\bullet}}(T)$. Every Hausdorff cluster point is therefore a monotone limit with the same volume as $\Delta_{Y_{\bullet}}(T)$, hence is equal to $\Delta_{Y_{\bullet}}(T)$.

\ref{itm:thm:pobmain:5} We may assume that $T$ is $\mathcal{I}$-good. It follows from \ref{itm:thm:pobmain:4} and \cref{thm:charIgoodasclosure} that we may reduce to the case where $T$ has analytic singularities. Our assertion follows from \cref{lma:liftOkounana}.

\ref{itm:thm:pobmain:6} By \eqref{eq:DeltaTgeneral}, in order to prove \eqref{eq:pobadditiv}, we may assume that $T$ and $S$ are both Kähler currents. Choose quasi-equisingular approximations $(T_j)_j$ and $(S_j)_j$ by Kähler currents as in \cref{def:POBKahcurr}. By \cref{cor:quasi-equichar,thm:dSadditivity}, the sums of their potentials converge in $d_S$ to the sum of the two $\mathcal I$-model representatives. This sum is $\mathcal I$-equivalent to the potential of $T+S$ by \cref{prop:Iequivmax}. Therefore \ref{itm:thm:pobmain:2} and \ref{itm:thm:pobmain:4} reduce the inclusion to the case of analytic singularities. Replacing $X$ by a suitable modification, we may assume that $T$ and $S$ both have log singularities, say
\[
T=[D]+R,\quad S=[D']+R',
\]
where $D$ and $D'$ are $\mathbb{Q}$-divisors on $X$ and $R$ and $R'$ are closed positive $(1,1)$-currents with bounded potentials. We need to show that
\[
\Delta_{Y_{\bullet}}(\{R\})+\Delta_{Y_{\bullet}}(\{R'\})+\nu_{Y_{\bullet}}([D])+\nu_{Y_{\bullet}}([D'])\subseteq \Delta_{Y_{\bullet}}(\{R+R'\})+\nu_{Y_{\bullet}}([D+D']).
\]
By \cref{prop:nuvaluationlinear}, this is equivalent to
\[
\Delta_{Y_{\bullet}}(\{R\})+\Delta_{Y_{\bullet}}(\{R'\})\subseteq \Delta_{Y_{\bullet}}(\{R+R'\}),
\]
which is already proved in \cref{thm:Okounkovtranmain}.
\end{proof}

\begin{corollary}\label{cor:transalgPOB}
Assume that $L$ is a big line bundle on $X$ and $h$ is a plurisubharmonic metric on $L$ with positive volume. Then
\begin{equation}\label{eq:tranOkounandalgOkoun}
\Delta_{Y_{\bullet}}(\ddc h)=\Delta_{Y_{\bullet}}(L,h).
\end{equation}
\end{corollary}
Similarly, the definition \eqref{eq:DeltanuTalgebraic1} is compatible with the definition in \cref{def:generalPOB}.
\begin{proof}
    We may assume that $\ddc h$ has positive mass and is $\mathcal{I}$-good. By the $d_S$-continuity of both sides of \eqref{eq:tranOkounandalgOkoun} as proved in \cref{thm:pobmain} and \cref{thm:Okoucont}, together with \cref{thm:charIgoodasclosure}, we may assume that $\ddc h$ has analytic singularities.

    In this case, using the birational invariance of both sides of \eqref{eq:tranOkounandalgOkoun} as proved in \cref{prop:birinvO} and \cref{thm:pobmain}, we may assume that $\ddc h$ has log singularities. Finally, after all these reductions, the equality \eqref{eq:tranOkounandalgOkoun} holds by construction.
\end{proof}

\subsection{The valuative characterization}
\label{subsec:tpob_valuative}
In this section, we will characterize the partial Okounkov bodies using valuations of currents.

\begin{lemma}\label{lma:Kahlerclassokounrest}
    Let $\beta$ be a big and nef class on $X$. Then
    \begin{equation}\label{eq:Deltaresttox10}
    \left\{ y\in \mathbb{R}^{n-1}: (0,y)\in \Delta_{Y_{\bullet}}(\beta)\right\}=\Delta_{Y_1\supseteq \cdots \supseteq Y_n}(\beta|_{Y_1}).
    \end{equation}
\end{lemma}
Here the right-hand side is understood in the sense of \cref{def:generalPOB} when $\beta|_{Y_1}$ is not big.
\begin{proof}
    \textbf{Step~1}. We first reduce to the case where $\beta$ is a Kähler class.

    Take a Kähler class $\alpha$ on $X$. It follows from the volume formula in \cref{thm:Okounkovtranmain} that
    \[
        \Delta_{Y_{\bullet}}(\beta)=\bigcap_{\epsilon>0}\Delta_{Y_{\bullet}}(\beta+\epsilon\alpha).
    \]
    Since $\beta|_{Y_1}$ is nef, hence pseudo-effective, \cref{def:generalPOB} similarly gives
    \[
        \Delta_{Y_1\supseteq \cdots \supseteq Y_n}(\beta|_{Y_1})= \bigcap_{\epsilon>0}\Delta_{Y_1\supseteq \cdots \supseteq Y_n}(\beta|_{Y_1}+\epsilon \alpha|_{Y_1}).
    \]
    So it suffices to prove \eqref{eq:Deltaresttox10} with $\beta+\epsilon\alpha$ in place of $\beta$.

    \textbf{Step~2}. Assume that $\beta$ is a Kähler class.
    The $\supseteq$ direction in \eqref{eq:Deltaresttox10} follows from the extension theorem \cref{thm:CT-thm-refined'}. To prove the other direction, recall that by \cref{thm:Okounkovtranmain}, for $t>0$ small enough, we have
    \[
    \left\{ y\in \mathbb{R}^{n-1}: (t,y)\in \Delta_{Y_{\bullet}}(\beta) \right\}=\Delta_{Y_1\supseteq \cdots \supseteq Y_n}\left((\beta-t\{Y_1\})|_{Y_1}\right).
    \]
    As $t\to 0+$, the right-hand side converges to $\Delta_{Y_1\supseteq \cdots \supseteq Y_n}(\beta|_{Y_1})$ with respect to the Hausdorff metric as a consequence of \cref{thm:Okounkovtranmain}, while the left-hand side converges to
    \[
    \left\{ y\in \mathbb{R}^{n-1}: (0,y)\in \Delta_{Y_{\bullet}}(\beta) \right\}
    \]
    by \cref{lma:Hausdorffconvslice}. We conclude our assertion.
\end{proof}

\begin{lemma}\label{lma:slicepob}
    Let $T\in \mathcal{Z}_+(X,\alpha)$ be a Kähler current. Assume that $\nu(T,Y_1)=0$. Then
    \begin{equation}\label{eq:Deltaslice}
        \left\{ y\in \mathbb{R}^{n-1}: (0,y)\in \Delta_{Y_{\bullet}}(T)\right\}=\Delta_{Y_1\supseteq \cdots \supseteq Y_n}\left(\Tr^{\alpha|_{Y_1}}_{Y_1}(T)\right).
    \end{equation}

    More generally, if $T\in \mathcal{Z}_+(X,\alpha)$ and $\nu(T,Y_1)=0$, suppose in addition that $\Tr^{\alpha|_{Y_1}}_{Y_1}(T)$ is defined, then \eqref{eq:Deltaslice} still holds.
\end{lemma}
See \cref{rmk:tracecurrent} for the definition of $\Tr^{\alpha|_{Y_1}}_{Y_1}(T)$.
Note that $\Delta_{Y_1\supseteq \cdots \supseteq Y_n}\left(\Tr^{\alpha|_{Y_1}}_{Y_1}(T)\right)$ is independent of the choice of the representative $\Tr^{\alpha|_{Y_1}}_{Y_1}(T)$.
\begin{proof}
    Let $\omega$ be a K\"ahler form on $X$. The last assertion follows from the first by perturbing $\theta$ to $\theta+\epsilon\omega$.

    \textbf{Step~1}.
    We first handle the case where $T$ has analytic singularities. Let $\pi\colon Z\rightarrow X$ be a modification such that
    \begin{enumerate}
        \item\label{itm:chapter-Okounkov:L1728:1} $Y_{\bullet}$ admits a lifting $(W_{\bullet},g)$, and
        \item\label{itm:chapter-Okounkov:L1728:2} $\pi^*T=[D]+R$, where $D$ is an effective $\mathbb{Q}$-divisor on $Z$ and $R$ is a closed positive $(1,1)$-current with bounded potential.
    \end{enumerate}
    This is possible by \cref{thm:resolvelogsing} and \cref{thm:liftableflag}.

    By \cref{lma:rescommpullback},
    \[
    \Pi^*\Tr_{Y_1}(T)\sim_P \Tr_{W_1}(\pi^*T),
    \]
    where $\Pi\colon W_1\rightarrow Y_1$ is the restriction of $\pi$. It follows from \cref{thm:pobmain} that
    \[
    \begin{split}
    \Delta_{W_1\supseteq \cdots \supseteq W_n}(\Tr_{W_1}(\pi^*T))=& \Delta_{Y_1\supseteq \cdots \supseteq Y_n}(\Tr_{Y_1}(T))\cor(Y_1\supseteq \cdots \supseteq Y_n,\Pi),\\
    \Delta_{W_{\bullet}}(\pi^*T)=&\Delta_{Y_{\bullet}}(T)g.
    \end{split}
    \]
    Taking \eqref{eq:correcur} into account, we find that it suffices to show that
    \[
    \left\{y\in \mathbb{R}^{n-1}:(0,y)\in\Delta_{W_{\bullet}}(\pi^*T) \right\}=\Delta_{W_1\supseteq \cdots \supseteq W_n}(\Tr_{W_1}(\pi^*T)).
    \]
    We may assume that $\pi$ is the identity map. Since $\nu(T,Y_1)=0$, the divisor $D$ does not contain $Y_1$. Hence we have
    \[
    T=[D]+R,\quad T|_{Y_1}=[D]|_{Y_1}+R|_{Y_1}.
    \]
    Note that $[D]|_{Y_1}$ is the current of integration along an effective $\mathbb{Q}$-divisor on $Y_1$.

    In particular,
    \[
    \begin{split}
    \Delta_{Y_{\bullet}}(T)=& \Delta_{Y_{\bullet}}(\{R\})+\nu_{Y_{\bullet}}([D]),\\
    \Delta_{Y_1\supseteq \cdots \supseteq Y_n}(T|_{Y_1})=& \Delta_{Y_1\supseteq \cdots \supseteq Y_n}(\{R\}|_{Y_1})+\nu_{Y_1\supseteq \cdots \supseteq Y_n}([D]|_{Y_1}).
    \end{split}
    \]
    The first coordinate of $\nu_{Y_{\bullet}}([D])$ is $0$, and the remaining coordinates are exactly $\nu_{Y_1\supseteq \cdots \supseteq Y_n}([D]|_{Y_1})$. Thus it suffices to show that
    \[
    \left\{y\in \mathbb{R}^{n-1}: (0,y)\in \Delta_{Y_{\bullet}}(\{R\})\right\}=\Delta_{Y_1\supseteq \cdots \supseteq Y_n}(\{R\}|_{Y_1}).
    \]
    Since $R$ has bounded potential, its cohomology class is nef; it is big because $T$ is a Kähler current. Thus the last identity is exactly \cref{lma:Kahlerclassokounrest}.

    \textbf{Step~2}.
    Next we consider the case where $T$ is a Kähler current. Take a quasi-equisingular approximation $(T_j)_j$ by Kähler currents as in \cref{def:POBKahcurr}. From Step~1, we know that for large $j\geq 1$,
    \[
    \left\{y\in \mathbb{R}^{n-1}: (0,y)\in \Delta_{Y_{\bullet}}(T_j)\right\}=\Delta_{Y_1\supseteq \cdots \supseteq Y_n}(\Tr_{Y_1}(T_j)).
    \]
    Letting $j\to\infty$ and applying \cref{thm:pobmain} and \cref{prop:tracedeclimit}, we conclude \eqref{eq:Deltaslice}.
\end{proof}

\begin{remark}\label{rmk:slicepob}
    More generally, the same argument shows the following result: Let $k=0,\ldots,n$ and $T\in \mathcal{Z}_+(X,\alpha)$ such that $\nu(T,Y_k)=0$. Assume that $\Tr^{\alpha|_{Y_k}}_{Y_k}(T)$ is defined. Then
    \begin{equation}
        \left\{ y\in \mathbb{R}^{n-k}: (0,\ldots,0,y)\in \Delta_{Y_{\bullet}}(T)\right\}=\Delta_{Y_k\supseteq \cdots \supseteq Y_n}\left(\Tr^{\alpha|_{Y_k}}_{Y_k}(T)\right).
    \end{equation}
    Also note that this result extends \cite[Theorem~3.4]{Jow10} and hence gives simpler proofs of \cite[Theorem~A, Theorem~B]{Jow10}.
\end{remark}

\begin{theorem}\label{thm:KahcurrminOkoun}
Assume that $T\in \mathcal{Z}_+(X,\alpha)_{>0}$ is a Kähler current.
We have
\begin{equation}\label{eq:minOkounkov}
    \min_{\lex} \Delta_{Y_{\bullet}}(T)=\nu_{Y_{\bullet}}(T).
\end{equation}
\end{theorem}
Here the minimum is with respect to the lexicographic order.
\begin{proof}
We argue by induction on $n\geq 0$. The case $n=0$ is of course trivial. Let us assume that $n>0$ and the case $n-1$ has been proved.

We first observe that by \cref{thm:pobmain},
\[
\Delta_{Y_{\bullet}}(T-\nu(T,Y_1)[Y_1])+(\nu(T,Y_1),0,\ldots,0)\subseteq \Delta_{Y_{\bullet}}(T).
\]
By \itemref{itm:prop:vollinearlimit:5}, removing the divisorial part does not change the volume. Comparing the volumes of both sides using \cref{thm:pobmain} and \cref{prop:vollinearlimit}, we find that equality holds:
\[
\Delta_{Y_{\bullet}}(T-\nu(T,Y_1)[Y_1])+(\nu(T,Y_1),0,\ldots,0)=\Delta_{Y_{\bullet}}(T).
\]
Replacing $T$ by $T-\nu(T,Y_1)[Y_1]$, we may therefore assume that $\nu(T,Y_1)=0$. It suffices to apply \cref{lma:slicepob} and the inductive hypothesis.
\end{proof}

\begin{corollary}\label{cor:valuationcurrentinPOB}
For any $T\in \mathcal{Z}_+(X,\alpha)$,
\[
\nu_{Y_{\bullet}}(T)\in \Delta_{Y_{\bullet}}(T)\subseteq \Delta_{Y_{\bullet}}(\alpha).
\]
\end{corollary}
\begin{proof}
    When $T$ is a Kähler current, this follows from \cref{thm:KahcurrminOkoun}.

    In general, by definition, $\nu_{Y_{\bullet}}(T)=\nu_{Y_{\bullet}}(T+\omega)$ for any Kähler form $\omega$ on $X$. It follows that
    \[
    \nu_{Y_{\bullet}}(T)\in \Delta_{Y_{\bullet}}(T+\omega)
    \]
    for any Kähler form $\omega$. It follows that $\nu_{Y_{\bullet}}(T)\in \Delta_{Y_{\bullet}}(T)$.
\end{proof}

\begin{theorem}\label{thm:Deltapartialint}
    For any $T\in \mathcal{Z}_+(X,\alpha)_{>0}$,
    \begin{equation}\label{eq:DeltaTequalallval}
    \Delta_{Y_{\bullet}}(T)=\overline{\left\{\nu_{Y_{\bullet}}(S):S\in \mathcal{Z}_+(X,\alpha),S\preceq_{\mathcal{I}}T\right\}}.
    \end{equation}
    In particular,\footnote{According to Ya Deng, the definition \eqref{eq:TPOB_val} of $\Delta_{Y_{\bullet}}(\alpha)$ was what Demailly originally proposed for Deng's thesis. Due to the lack of the techniques of the trace operators, Deng had to work with analytic singularities instead. As a consequence, the transcendental analogue of \cref{prop:birinvO} is not obvious. This is one of the two main technical difficulties of \cref{thm:Okounkovtranmain}. This problem also led me to finally develop the theory of trace operators, a notion that had been on the author's mind for several years.}
    \begin{equation}\label{eq:TPOB_val}
    \Delta_{Y_{\bullet}}(\alpha)=\overline{\left\{\nu_{Y_{\bullet}}(S):S\in \mathcal{Z}_+(X,\alpha)\right\}}.
    \end{equation}
\end{theorem}

\begin{remark}\label{rmk:Okounkov:L1841}
We expect that the closure operation in \eqref{eq:DeltaTequalallval} is not necessary.
\end{remark}

\begin{proof}
    The $\supseteq$ direction in \eqref{eq:DeltaTequalallval} follows from \cref{cor:valuationcurrentinPOB} together with \itemref{itm:thm:pobmain:2}.

    Let us write
    \[
    D_{Y_{\bullet}}(T)=\left\{\nu_{Y_{\bullet}}(S):S\in \mathcal{Z}_+(X,\alpha),S\preceq_{\mathcal{I}}T\right\}
    \]
    for the time being.

    \textbf{Step~1}.
    Assume that $T$ has analytic singularities. We have
    \[
    \begin{split}
    \Delta_{Y_{\bullet}}(T)\supseteq & \overline{D_{Y_{\bullet}}(T)}\\
    \supseteq & \overline{\left\{\nu_{Y_{\bullet}}(S):\mathcal{Z}_+(X,\alpha)\ni S \textup{ has gentle analytic singularities}, S\preceq T\right\}}.
    \end{split}
    \]
    It follows easily from \cref{thm:Okounkovtranmain} that the volume of the right-hand side is equal to the volume of $\Delta_{Y_{\bullet}}(T)$, so \eqref{eq:DeltaTequalallval} holds.

    \textbf{Step~2}.
    Assume that $T$ is a Kähler current. Take a quasi-equisingular approximation $(T_j)_j$ by Kähler currents as in \cref{def:POBKahcurr}.
    Next we use the language of psh functions. Let $\varphi_j,\varphi\in \PSH(X,\theta)$ be the potentials corresponding to $T_j,T$ for each $j\geq 1$.

    Fix an integer $N>0$. For large enough $j\geq 1$, we can find $\psi_j\in \PSH(X,\theta)_{>0}$ such that
    \[
    P_{\theta}[\varphi]_{\mathcal{I}}\geq (1-N^{-1})\varphi_j+N^{-1} \psi_j.
    \]
    The existence of $\psi_j$ follows from \cref{lma:pathoenvelope}.
    It follows that
    \[
    \begin{aligned}
    D_{Y_{\bullet}}(T)\supseteq & D_{Y_{\bullet}}\left(\theta+\ddc \left((1-N^{-1})\varphi_j+N^{-1} \psi_j\right)\right)\\
    \supseteq & (1-N^{-1})D_{Y_{\bullet}}(T_j)+N^{-1} D_{Y_{\bullet}}(\theta+\ddc \psi_j).
    \end{aligned}
    \]
    We claim that
    \[
    \overline{D_{Y_{\bullet}}(T)}\supseteq \bigcap_{j=1}^{\infty}\overline{D_{Y_{\bullet}}(T_j)}.
    \]
    Let $x$ be a point in the right-hand side. Fix $N>0$. For all large $j$, choose
    \[
    y_j\in D_{Y_{\bullet}}(T_j),\quad z_j\in D_{Y_{\bullet}}(\theta+\ddc \psi_j)
    \]
    with $y_j\to x$. By the preceding inclusion,
    \[
    (1-N^{-1})y_j+N^{-1}z_j\in D_{Y_{\bullet}}(T).
    \]
    The points $x,y_j,z_j$ all lie in the compact set $\Delta_{Y_{\bullet}}(\alpha)$ after taking closures, by \cref{cor:valuationcurrentinPOB}. Hence, letting $j\to\infty$ and then $N\to\infty$, we get $x\in \overline{D_{Y_{\bullet}}(T)}$, proving the claim.
    Therefore, by Step~1, we conclude that
    \[
    \Delta_{Y_{\bullet}}(T)=\bigcap_{j=1}^{\infty}\overline{\Delta_{Y_{\bullet}}(T_j)}=\bigcap_{j=1}^{\infty}\overline{D_{Y_{\bullet}}(T_j)}\subseteq\overline{D_{Y_{\bullet}}(T)}.
    \]
The reverse direction is already known.

\textbf{Step~3}.
Finally, consider the general case. Take a Kähler current $T'\in \mathcal{Z}_+(X,\alpha)$ more singular than $T$ (the existence of which is given by \cref{lma:kahcurrentposmass}). For each $\epsilon\in (0,1)$, we know that
\[
\Delta_{Y_{\bullet}}((1-\epsilon)T+\epsilon T')=\overline{D_{Y_{\bullet}}((1-\epsilon)T+\epsilon T')}
\subseteq \overline{D_{Y_{\bullet}}(T)}.
\]
Letting $\epsilon\to 0+$ and using the $d_S$-continuity in \itemref{itm:thm:pobmain:4},
we find that
\[
\Delta_{Y_{\bullet}}(T)\subseteq \overline{D_{Y_{\bullet}}(T)}.
\]
As the other inclusion is already known, we conclude.
\end{proof}

\begin{corollary}\label{cor:KahcurrminOkoun}
Assume that $T\in \mathcal{Z}_+(X,\alpha)_{>0}$.
We have
\begin{equation}\label{eq:minOkounkov3}
    \min_{\lex} \Delta_{Y_{\bullet}}(T)=\nu_{Y_{\bullet}}(T).
\end{equation}
\end{corollary}
\begin{proof}
    By \cref{thm:Deltapartialint}, it is clear that
    \[
    \min_{\lex} \Delta_{Y_{\bullet}}(T)\leq_{\lex}\nu_{Y_{\bullet}}(T).
    \]
    On the other hand, we clearly have
    \[
    \Delta_{Y_{\bullet}}(T)\subseteq \Delta_{Y_{\bullet}}(T+\omega)
    \]
    for any Kähler form $\omega$ on $X$. It follows that
    \[
    \min_{\lex} \Delta_{Y_{\bullet}}(T)\geq_{\lex} \min_{\lex} \Delta_{Y_{\bullet}}(T+\omega).
    \]
    By \cref{thm:KahcurrminOkoun}, the right-hand side is just $\nu_{Y_{\bullet}}(T+\omega)=\nu_{Y_{\bullet}}(T)$. We conclude the proof.
\end{proof}

\section{Okounkov test curves}\label{sec:Okotc}
Fix $n\in \mathbb{N}$.
Let $\Delta,\Delta'\subseteq \mathbb{R}^n$ be convex bodies with positive volumes. The standard Lebesgue measure on $\mathbb{R}^n$ is denoted by $\vol$.

Recall that $\mathcal{K}_n$ denotes the set of convex bodies in $\mathbb{R}^n$ and $d_{\Haus}$ denotes the Hausdorff metric.
We refer to \cref{chap:almostsg} for the basic properties of these objects.

We study the notion of Okounkov test curves in this section; this leads to the definition of the Duistermaat--Heckman measure of a non-Archimedean metric in \cref{sec:DHmeasure} below. This section can be skipped on a first reading.

\begin{definition}\label{def:Otc}
	An \emph{Okounkov test curve}\index{test curve!Okounkov test curve} relative to $\Delta$ consists of
    \begin{enumerate}
        \item\label{itm:def:Otc:1} a number $\Delta_{\max}\in \mathbb{R}$ and
        \item\label{itm:def:Otc:2} an assignment $(-\infty,\Delta_{\max})\ni\tau \mapsto \Delta_{\tau}\in \mathcal{K}_n$ satisfying
        \begin{enumerate}
            \item\label{itm:def:Otc:2:a} the assignment $\tau\mapsto \Delta_{\tau}$ is decreasing and concave\footnote{Here concavity refers to the concavity with respect to the Minkowski sum.};
            \item\label{itm:def:Otc:2:b} we have $\Delta_{\tau}\xrightarrow{d_{\Haus}}\Delta$ as $\tau\to-\infty$.
        \end{enumerate}
    \end{enumerate}
    The set of Okounkov test curves relative to $\Delta$ is denoted by $\TC(\Delta)$\index{$\TC(\Delta)$}.

    An Okounkov test curve $\Delta_{\bullet}$ relative to $\Delta$ is \emph{bounded}\index{Okounkov test curve!bounded Okounkov test curve} if $\Delta_{\tau}=\Delta$ when $\tau$ is small enough. The subset of bounded Okounkov test curves is denoted by $\TC^{\infty}(\Delta)$\index{$\TC^{\infty}(\Delta)$}.

    An Okounkov test curve $\Delta_{\bullet}$ relative to $\Delta$ is said to have \emph{finite energy}\index{Okounkov test curve!Okounkov test curve with finite energy} if
    \begin{equation}\label{eq:Otestcurvenergy}
    \mathbf{E}(\Delta_{\bullet})\coloneqq n!\Delta_{\max} \vol \Delta+n!\int_{-\infty}^{\Delta_{\max}} \left(\vol \Delta_{\tau}-\vol \Delta\right)\,\mathrm{d}\tau>-\infty.
    \end{equation}
    \index{$\mathbf{E}(\Delta_{\bullet})$}
    The subset of Okounkov test curves with finite energy is denoted by $\TC^{1}(\Delta)$\index{$\TC^{1}(\Delta)$}.

    Given $\Delta_{\bullet}\in \TC(\Delta)$ and $\Delta'_{\bullet}\in \TC(\Delta')$, we say $\Delta_{\bullet}\leq \Delta'_{\bullet}$ if $\Delta_{\max}\leq \Delta'_{\max}$ and for any $\tau<\Delta_{\max}$, we have $\Delta_{\tau}\subseteq \Delta'_{\tau}$.
\end{definition}
Sometimes it is convenient to introduce
\begin{equation}\label{eq:DeltaDeltamax}
    \Delta_{\Delta_{\max}}=\bigcap_{\tau<\Delta_{\max}}\Delta_{\tau}\in \mathcal{K}_n.
\end{equation}
We shall always make this extension in the sequel when we talk about $\Delta_{\Delta_{\max}}$. Observe that $(-\infty,\Delta_{\max}]\ni\tau \mapsto \Delta_{\tau}$ is still concave.

\begin{proposition}\label{prop:Otccont}
	Any Okounkov test curve $(\Delta_{\tau})_{\tau< \Delta_{\max}}$ relative to $\Delta$ is continuous in $\tau$. Moreover, $\vol \Delta_{\tau}>0$ for all $\tau<\Delta_{\max}$.
\end{proposition}
\begin{proof}
    We first claim that $\vol\Delta_{\tau'}>0$ for all $\tau'<\Delta_{\max}$. By Condition~(2b) in \cref{def:Otc} and \cref{thm:contvol}, we know that $\vol\Delta_{\tau''}>0$ when $\tau''$ is small enough. Fix one such $\tau''$.
    We may assume that $\tau''\leq \tau'$ since otherwise there is nothing to prove. Next take $\tau'''\in (\tau',\Delta_{\max})$. Take $t\in (0,1)$ such that $\tau'=t\tau'''+(1-t)\tau''$. It follows that
    \[
    \vol \Delta_{\tau'}\geq \vol \left( t\Delta_{\tau'''}+(1-t)\Delta_{\tau''}\right)\geq (1-t)^n\vol \Delta_{\tau''}>0.
    \]

	Next we claim that $\vol\Delta_{\tau}$ is continuous for $\tau< \Delta_{\max}$. In fact, it follows from \cref{thm:BrunnMin} that $(-\infty,\Delta_{\max}) \ni \tau\mapsto \log \vol\Delta_{\tau}$ is concave, but we already know that it is finite, hence the continuity follows.

	Next we show that
	\[
		\Delta_{\tau}=\bigcap_{\tau'<\tau}\Delta_{\tau'}.
	\]
	The $\subseteq$ direction is obvious. By the continuity of the volume, both sides have the same volume and the volume is positive, we therefore obtain the equality.

	Similarly, we have
	\[
		\Delta_{\tau}=\overline{\bigcup_{\tau'>\tau}\Delta_{\tau'}}.
	\]
	The continuity of $\Delta_{\tau}$ at $\tau<\Delta_{\max}$ is proved.
\end{proof}

\begin{definition}\label{def:tf}
    A \emph{test function}\index{test function} on $\Delta$ is a function $G\colon \Delta\rightarrow [-\infty,\infty)$ such that
	\begin{enumerate}
		\item\label{itm:def:tf:1} $G$ is concave,
		\item\label{itm:def:tf:2} $G$ is finite on $\Int \Delta$, and
		\item\label{itm:def:tf:3} $G$ is upper semicontinuous.
	\end{enumerate}

    A test function $G$ is \emph{bounded}\index{test function!bounded test function} if $G$ is bounded from below.

    A test function $G$ has \emph{finite energy}\index{test function!test function with finite energy} if
 \begin{equation}\label{eq:EF}
	\mathbf{E}(G)\coloneqq n!\int_{\Delta} G\,\mathrm{d}\lambda>-\infty.
\end{equation}
\end{definition}

\begin{definition}\label{def:LegOkoun}
    Let $\Delta_{\bullet}\in \TC(\Delta)$. We define its \emph{Legendre transform}\index{Legendre transform} as
    \[
	   G[\Delta_{\bullet}] \colon \Delta\rightarrow [-\infty,\infty),\quad a\mapsto \sup\left\{\tau<\Delta_{\max}:a\in \Delta_{\tau}\right\}.
	\]
     Given a test function $G\colon \Delta\rightarrow [-\infty,\infty)$, we define its \emph{inverse Legendre transform}\index{Legendre transform!inverse Legendre transform} $\Delta[G]_{\bullet}$ as the Okounkov test curve relative to $\Delta$ defined as follows:
     \begin{enumerate}
         \item\label{itm:def:LegOkoun:1} $\Delta[G]_{\max}=\sup_{\Delta} G$, and
         \item\label{itm:def:LegOkoun:2} for each $\tau<\sup_{\Delta} G$, we set
         \[
         \Delta[G]_{\tau}=\{x\in \Delta:G\geq \tau\}.
         \]
     \end{enumerate}
\end{definition}
We observe that
\begin{equation}\label{eq:GDeltamax}
    G[\Delta_{\bullet}](a)=\max\left\{\tau\leq \Delta_{\max}: a\in \Delta_{\tau}\right\},\textup{ if } G[\Delta_{\bullet}](a)>-\infty.
\end{equation}

\begin{lemma}\label{lma:convbodyLegendre}
    Let $\Delta_{\bullet}\in \TC(\Delta)$. Then $G[\Delta_{\bullet}]$ defined in \cref{def:LegOkoun} is a test function.

    Similarly, if $G\colon \Delta\rightarrow [-\infty,\infty)$ is a test function, then $\Delta[G]_{\bullet}$ is an Okounkov test curve.
\end{lemma}
\begin{proof}
    First suppose that $\Delta_{\bullet}\in \TC(\Delta)$. We want to verify that $G[\Delta_{\bullet}]$ satisfies the conditions in \cref{def:tf}.

    We first verify the concavity.
    Take $a,b\in \Delta$. We want to prove that for any $t\in (0,1)$,
	\begin{equation}\label{eq:GDeltaconc}
		G[\Delta_{\bullet}](ta+(1-t)b)\geq tG[\Delta_{\bullet}](a)+(1-t)G[\Delta_{\bullet}](b).
	\end{equation}
	There is nothing to prove if $G[\Delta_{\bullet}](a)$ or $G[\Delta_{\bullet}](b)$ is $-\infty$. So we assume that both are finite. In this case, by \eqref{eq:GDeltamax},
    \[
    a\in \Delta_{G[\Delta_{\bullet}](a)},\quad b\in \Delta_{G[\Delta_{\bullet}](b)}.
    \]
    Thus,
	\[
		ta+(1-t)b\in t  \Delta_{G[\Delta_{\bullet}](a)}+(1-t)\Delta_{G[\Delta_{\bullet}](b)}\subseteq \Delta_{tG[\Delta_{\bullet}](a)+(1-t)G[\Delta_{\bullet}](b)}.
	\]
	We deduce that
	\[
		G[\Delta_{\bullet}](ta+(1-t)b)\geq tG[\Delta_{\bullet}](a)+(1-t)G[\Delta_{\bullet}](b).
	\]
	Therefore, \eqref{eq:GDeltaconc} follows.

 It is clear that $G[\Delta_{\bullet}]$ is finite on the interior of $\Delta$. It remains to argue that $G[\Delta_{\bullet}]$ is upper semicontinuous.

 Let $(a_i)_{i\geq 1}$ be a sequence in $\Delta$ with limit $a\in \Delta$. Define $\tau_i=G[\Delta_{\bullet}](a_i)$. Let $\tau=\varlimsup_i \tau_i$.
	We need to show that
	\begin{equation}\label{eq:ainDelta1}
		G[\Delta_{\bullet}](a)\geq \tau.
	\end{equation}
    There is nothing to prove if $\tau=-\infty$. We assume that this is not the case. Passing to a subsequence, we may assume that $\tau_i\to \tau$. In particular, we can assume that $\tau_i\neq -\infty$ for all $i\geq 1$.
    It follows from \eqref{eq:GDeltamax} that $a_i\in\Delta_{\tau_i}$ for all $i\geq 1$.
    If $\tau<\Delta_{\max}$, then $\Delta_{\tau_i}\xrightarrow{d_{\Haus}}\Delta_{\tau}$, so \cref{thm:Hausconvcond} implies that $a\in \Delta_{\tau}$. If $\tau=\Delta_{\max}$, then the same conclusion follows from the definition $\Delta_{\Delta_{\max}}=\bigcap_{\tau'<\Delta_{\max}}\Delta_{\tau'}$.
    Thus, \eqref{eq:ainDelta1} follows.

    Conversely, suppose that $G\colon \Delta\rightarrow [-\infty,\infty)$ is a test function. We argue that $\Delta[G]_{\bullet}$ is an Okounkov test curve. We verify the conditions in \cref{def:Otc}.

    Firstly, for each $\tau< \sup_{\Delta}G$, the set $\Delta[G]_{\tau}$ is a convex body as $G$ is concave and usc.
	Moreover, $\Delta[G]_{\tau}$ is clearly decreasing in $\tau$.

	Secondly, for each $a\in \Delta$, we can write $a=\lim_{i}a_i$ with $a_i\in \Int \Delta$. By assumption, $G$ is finite at $a_i$. Thus,
	\[
		a\in \overline{\{G>-\infty\}}=\overline{\bigcup_{\tau<\sup_{\Delta}G}\Delta[G]_{\tau}}.
	\]
	By \cref{thm:Hausconvcond}, $\Delta[G]_{\tau}\xrightarrow{d_{\Haus}}\Delta$ as $\tau\to -\infty$.

	Thirdly, $\Delta[G]$ is concave. To see this, take $\tau,\tau'< \Delta_{\max}$. We need to prove that for any $t\in (0,1)$,
	\begin{equation}\label{eq:Deconc}
		\Delta[G]_{t\tau+(1-t)\tau'}\supseteq t\Delta[G]_{\tau}+(1-t)\Delta[G]_{\tau'}.
	\end{equation}
	Let $a\in \Delta[G]_{\tau}$ and $b\in \Delta[G]_{\tau'}$. We have $G(a)\geq\tau$ and $G(b)\geq \tau'$. As $G$ is concave, we have $G(ta+(1-t)b)\geq t\tau+(1-t)\tau'$. Thus,
	\[
		ta+(1-t)b\in \Delta[G]_{t\tau+(1-t)\tau'}
	\]
	and \eqref{eq:Deconc} follows.
\end{proof}

\begin{theorem}\label{thm:Okotestcurve}
	The Legendre transform and inverse Legendre transform are inverse to each other,
	defining a bijection between $\TC(\Delta)$ and the set of test functions on $\Delta$.

    Under this bijection, $\TC^1(\Delta)$ corresponds to test functions on $\Delta$ with finite energy and $\TC^{\infty}(\Delta)$ corresponds to bounded test functions on $\Delta$.
\end{theorem}
\begin{proof}
    Thanks to \cref{lma:convbodyLegendre}, in order to prove the first assertion, it only remains to see that the Legendre transform and the inverse Legendre transform are inverse to each other, which is immediate by definition.

    It is obvious that $\TC^{\infty}(\Delta)$ corresponds to bounded test functions. Moreover, a direct computation shows that if $\Delta_{\bullet}\in \TC(\Delta)$, then
    \[
    \mathbf{E}(\Delta_{\bullet})=\mathbf{E}(G[\Delta_{\bullet}]),
    \]
    concluding the $\TC^{1}(\Delta)$ case.
\end{proof}

\begin{proposition}\label{prop:decnetLegend}
    Let $(\Delta^i)_{i\in I}$ be a decreasing net in $\mathcal{K}_n$ such that $\Delta\subseteq \bigcap_{i\in I}\Delta^i$. Consider a decreasing net $(\Delta^i_{\bullet})_{i\in I}$ with $\Delta^i_{\bullet}\in \TC(\Delta^i)$ for all $i\in I$ such that there is $\Delta_{\bullet}\in \TC(\Delta)$ satisfying the following properties:
    \begin{enumerate}
        \item\label{itm:prop:decnetLegend:1} $\Delta_{\max}= \lim_{i\in I}\Delta^i_{\max}$;
        \item\label{itm:prop:decnetLegend:2} for any $\tau<\Delta_{\max}$, we have $\Delta^i_{\tau}\xrightarrow{d_{\Haus}}\Delta_{\tau}$.
    \end{enumerate}
    Then for any $a\in \Delta$, we have
    \begin{equation}\label{eq:pwconvLegendre}
        \lim_{i\in I} G[\Delta^i_{\bullet}](a)=G[\Delta_{\bullet}](a).
    \end{equation}
\end{proposition}
Note that in general,
\[
\Delta\subsetneq \bigcap_{i\in I}\Delta^i.
\]
\begin{proof}
Fix $a\in \Delta$.
    It follows immediately from the definition of $G$ that the net $(G[\Delta^i_{\bullet}](a))_{i\in I}$ is decreasing and the $\geq$ direction in \eqref{eq:pwconvLegendre} holds. Let us prove the reverse inequality. Let $\tau$ denote the left-hand side of \eqref{eq:pwconvLegendre} for the moment. By definition, for any $\epsilon>0$ and any $i\in I$, we have $a\in \Delta^i_{\tau-\epsilon}$.
    It follows that
    \[
    a\in \Delta_{\tau-\epsilon}.
    \]
    Therefore,
    \[
    \tau\leq G[\Delta_{\bullet}](a).
    \]
    This completes the proof.
\end{proof}

Similarly, for increasing nets, we have:
\begin{proposition}\label{prop:incnetLegend}
    Let $(\Delta^i)_{i\in I}$ be an increasing net in $\mathcal{K}_n$ with Hausdorff limit $\Delta$ such that $\vol \Delta^i>0$ for all $i\in I$. Consider an increasing net $(\Delta^i_{\bullet})_{i\in I}$ with $\Delta^i_{\bullet}\in \TC(\Delta^i)$ for all $i\in I$.
    Let $\Delta_{\max}=\lim_{i\in I}\Delta^i_{\max}$. For any $\tau<\Delta_{\max}$, let $\Delta_{\tau}$ be the Hausdorff limit of $\Delta^i_{\tau}$. Then $\Delta_{\bullet}\in \TC(\Delta)$ and
    \begin{equation}\label{eq:apwconvLegendre}
    \lim_{i\in I}G[\Delta^i_{\bullet}](a)=G[\Delta_{\bullet}](a)
    \end{equation}
    for any $a\in \Int \Delta$.
\end{proposition}
\begin{proof}
    It is obvious that $\Delta_{\bullet}\in \TC(\Delta)$.

    Fix $a\in \Int \Delta$. Since $\Delta^i\xrightarrow{d_{\Haus}}\Delta$, we have $a\in \Delta^i$ for all large enough $i$.
    By definition, the net $(G[\Delta^i_{\bullet}](a))_{i\in I}$ is increasing and the $\leq$ direction in \eqref{eq:apwconvLegendre} holds.
    Let us write $\tau=G[\Delta_{\bullet}](a)$ for the time being.
    By definition of $G$, for any $\epsilon>0$, we have
    \[
    a\in \Delta_{\tau-\epsilon/2}.
    \]
    The concavity of $\Delta_{\bullet}$ guarantees that
    \[
    a\in \Int \Delta_{\tau-\epsilon}.
    \]
    Indeed, choose $\sigma\ll 0$ so that a small ball around $a$ is contained in $\Delta_\sigma$. Since $a\in\Delta_{\tau-\epsilon/2}$ and
    \[
        \tau-\epsilon=t\sigma+(1-t)(\tau-\epsilon/2)
    \]
    for a suitable $t\in(0,1)$, concavity gives a small ball around $a$ contained in $\Delta_{\tau-\epsilon}$.

    Since $\Delta^i_{\tau-\epsilon}\xrightarrow{d_{\Haus}}\Delta_{\tau-\epsilon}$ and $a\in \Int \Delta_{\tau-\epsilon}$, it follows that for all large enough $i\in I$,
    \[
    a\in \Delta^i_{\tau-\epsilon}.
    \]
    Therefore,
    \[
    \tau-\epsilon\leq G[\Delta^i_{\bullet}](a)
    \]
    for all large enough $i$. Taking the limit with respect to $i$ and then with respect to $\epsilon$, we conclude the desired inequality.
\end{proof}

\begin{definition}\label{def:DHmeasureOTC}
	Let $\Delta_{\bullet}$ be an Okounkov test curve relative to $\Delta$. We define the \emph{Duistermaat--Heckman measure}\index{Duistermaat--Heckman measure} $\DHm(\Delta_{\bullet})$\index{$\DHm(\Delta_{\bullet})$} as
	\[
		\DHm(\Delta_{\bullet})\coloneqq G[\Delta_{\bullet}]_*(\vol).
	\]
	It is a Radon measure on $\mathbb{R}$, since $G[\Delta_{\bullet}]$ is finite on $\Int \Delta$.
\end{definition}
In other words, $\DHm(\Delta_{\bullet})$ is the distribution of the random variable $G[\Delta_{\bullet}]$.

\begin{proposition}\label{prop:DHmoments}
    Let $\Delta_{\bullet}\in \TC(\Delta)$. Let $m\in \mathbb{Z}_{>0}$. Then the following extended moment identity holds:
    \begin{equation}\label{eq:momentcalc}
\int_{\mathbb{R}}x^m \DHm(\Delta_{\bullet})(x)=\Delta_{\max}^m\vol\Delta+m\int_{-\infty}^{\Delta_{\max}} \tau^{m-1}(\vol \Delta_{\tau}-\vol \Delta)\,\mathrm{d}\tau
    \end{equation}
and
\begin{equation}\label{eq:massDHm1}
\int_{\mathbb{R}}\DHm(\Delta_{\bullet})=\vol \Delta.
\end{equation}
\end{proposition}
\begin{proof}
In fact, \eqref{eq:massDHm1} follows immediately from the definition. Let $G=G[\Delta_{\bullet}]$. We first assume that $G$ is bounded from below. Then \eqref{eq:momentcalc} follows from a straightforward computation:
\[
\begin{aligned}
    &\int_{\mathbb{R}}x^m \DHm(\Delta_{\bullet})(x)\\
    =& \int_{\Delta} G[\Delta_{\bullet}](a)^m\,\mathrm{d}\vol(a)\\
    =& \int_{\Delta} \left( \Delta_{\max}^m-\int_{G[\Delta_{\bullet}](a)}^{\Delta_{\max}} m\tau^{m-1}\,\mathrm{d}\tau\right)\,\mathrm{d}\vol(a)\\
    =& \Delta_{\max}^m\vol \Delta-m\int_{\mathbb{R}}\int_{\Delta} \mathds{1}_{[G[\Delta_{\bullet}](a),\Delta_{\max}]}(\tau)\tau^{m-1}\,\mathrm{d}\vol(a)\,\mathrm{d}\tau\\
    =& \Delta_{\max}^m\vol \Delta-m\int_{-\infty}^{\Delta_{\max}}\int_{\Delta\setminus \Delta_{\tau}} \tau^{m-1}\,\mathrm{d}\vol(a)\,\mathrm{d}\tau\\
    =& \Delta_{\max}^m\vol \Delta-m\int_{-\infty}^{\Delta_{\max}} \tau^{m-1}\left(\vol \Delta-\vol \Delta_{\tau}\right)\,\mathrm{d}\tau.
\end{aligned}
\]
In general, set $G_A=G\lor (-A)$. The same pointwise identity applied to $G_A$ gives
\[
\int_{\Delta}G_A^m\,\mathrm{d}\vol=\Delta_{\max}^m\vol\Delta
+m\int_{-A}^{\Delta_{\max}}\tau^{m-1}(\vol \Delta_{\tau}-\vol\Delta)\,\mathrm{d}\tau.
\]
Letting $A\to\infty$ gives \eqref{eq:momentcalc}: if $m$ is even, this follows from monotone convergence applied to $G_A^m\nearrow G^m$; if $m$ is odd, apply monotone convergence to $\Delta_{\max}^m-G_A^m\nearrow \Delta_{\max}^m-G^m$. Thus both sides of \eqref{eq:momentcalc} are interpreted in the corresponding extended sense.
\end{proof}

\begin{lemma}\label{lma:DHmconv}
    Let $(\Delta^i)_{i\geq 1}$ be a decreasing sequence in $\mathcal{K}_n$ with limit $\Delta$.
    Suppose that $(\Delta_{\bullet}^i)_{i\geq 1}$ is a decreasing sequence with $\Delta_{\bullet}^i\in \TC(\Delta^i)$.
    Suppose that there is $\Delta_{\bullet}\in \TC(\Delta)$ such that
    \begin{enumerate}
        \item\label{itm:lma:DHmconv:1} $\Delta_{\max}=\lim_{i\to\infty}\Delta^i_{\max}$;
        \item\label{itm:lma:DHmconv:2} for any $\tau<\Delta_{\max}$, we have $\Delta^i_{\tau}\xrightarrow{d_{\Haus}}\Delta_{\tau}$.
    \end{enumerate}
    Then $\DHm(\Delta_{\bullet}^i)\rightharpoonup \DHm(\Delta_{\bullet})$.
\end{lemma}
\begin{proof}
    Let $f\in C_c(\mathbb{R})$. It follows from \cref{prop:decnetLegend} that
    \[
    G[\Delta_{\bullet}^i]\to G[\Delta_{\bullet}]
    \]
    pointwise on $\Delta$.
    Hence the dominated convergence theorem gives
    \[
        \int_{\Delta} f(G[\Delta_{\bullet}^i])\,\mathrm{d}\lambda\to \int_{\Delta} f(G[\Delta_{\bullet}])\,\mathrm{d}\lambda.
    \]
    Since $\Delta^i\xrightarrow{d_{\Haus}}\Delta$, \cref{thm:contvol} gives $\vol(\Delta^i\setminus \Delta)\to 0$. The contribution of $\Delta^i\setminus\Delta$ is bounded by $\|f\|_{\infty}\vol(\Delta^i\setminus\Delta)$, so
    \[
        \int_{\Delta^i} f(G[\Delta_{\bullet}^i])\,\mathrm{d}\lambda\to \int_{\Delta} f(G[\Delta_{\bullet}])\,\mathrm{d}\lambda.
    \]
    This is the weak convergence.
\end{proof}
Similarly, we have
\begin{lemma}\label{lma:DHmconv2}
    Let $(\Delta^i)_{i\geq 1}$ be an increasing sequence in $\mathcal{K}_n$ with Hausdorff limit $\Delta$ such that $\vol \Delta^i>0$ for all $i\geq 1$. Consider an increasing sequence $(\Delta^i_{\bullet})_{i\geq 1}$ with $\Delta^i_{\bullet}\in \TC(\Delta^i)$ for all $i\geq 1$.
    Let $\Delta_{\bullet}\in \TC(\Delta)$ be defined as
    \begin{enumerate}
        \item\label{itm:lma:DHmconv2:1} $\Delta_{\max}=\lim_{i\to\infty}\Delta^i_{\max}$;
        \item\label{itm:lma:DHmconv2:2} for any $\tau<\Delta_{\max}$, $\Delta_{\tau}$ is the Hausdorff limit of $\Delta^i_{\tau}$.
    \end{enumerate}
    Then we have
    \[
    \DHm(\Delta_{\bullet}^i)\rightharpoonup \DHm(\Delta_{\bullet}).
    \]
\end{lemma}
\begin{proof}
    Let $f\in C_c(\mathbb{R})$. It follows from \cref{prop:incnetLegend} that
    \[
    G[\Delta_{\bullet}^i]\to G[\Delta_{\bullet}]
    \]
    almost everywhere on $\Delta$.
    Extending $f(G[\Delta_{\bullet}^i])$ by $0$ outside $\Delta^i$, we also have $\mathds{1}_{\Delta^i}\to \mathds{1}_{\Delta}$ almost everywhere on $\Delta$, since the boundary of the convex body $\Delta$ has Lebesgue measure zero. The dominated convergence theorem gives
    \[
        \int_{\Delta^i} f(G[\Delta_{\bullet}^i])\,\mathrm{d}\lambda\to \int_{\Delta} f(G[\Delta_{\bullet}])\,\mathrm{d}\lambda.
    \]
    This is the weak convergence.
\end{proof}

The main source of Okounkov test curves is the following:
\begin{theorem}\label{thm:Okountescurvex}
Let $X$ be a connected compact Kähler manifold and $\theta$ be a closed smooth real $(1,1)$-form on $X$ representing a big cohomology class $\alpha$.
    Let $Y_{\bullet}$ be a smooth flag on $X$ and $\Gamma\in \TC(X,\theta)_{>0}$.
    Then the map
    \[
    (-\infty,\Gamma_{\max})\ni \tau\mapsto \Delta_{Y_{\bullet}}(\theta,\Gamma)_{\tau}\coloneqq \Delta_{Y_{\bullet}}(\theta,\Gamma_{\tau})
    \]
    defines an Okounkov test curve relative to $\Delta=\Delta_{Y_{\bullet}}(\theta,\phi)$, where $\phi=\Gamma_{-\infty}$.

    If $\Gamma\in\TC^{\infty}(X,\theta;\phi)$, then $\Delta_{Y_{\bullet}}(\theta,\Gamma)\in\TC^{\infty}(\Delta)$. If $\phi$ is $\mathcal I$-model and $\Gamma\in\TC^1(X,\theta;\phi)$, then $\Delta_{Y_{\bullet}}(\theta,\Gamma)\in\TC^1(\Delta)$. In the latter case, if all $\Gamma_\tau$ for $\tau<\Gamma_{\max}$ are also $\mathcal I$-model, then
    \[
        \mathbf{E}^{\phi}(\Gamma)=\mathbf{E}(\Delta_{Y_{\bullet}}(\theta,\Gamma)).
    \]
\end{theorem}
See \cref{def:testcur} and \cref{def:relattestcurv} for the relevant definitions.
\begin{proof}
    Consider $\Gamma\in \TC(X,\theta)_{>0}$. We need to verify that $\Delta_{Y_{\bullet}}(\theta,\Gamma)$ is an Okounkov test curve relative to $\Delta_{Y_{\bullet}}(\theta,\Gamma_{-\infty})$.

    First observe that $\tau\mapsto \Delta_{Y_{\bullet}}(\theta,\Gamma_{\tau})$ is concave and decreasing for $\tau<\Gamma_{\max}$. The decreasing property follows from \itemref{itm:thm:pobmain:2}. For concavity, use the valuative characterization \cref{thm:Deltapartialint} together with the additivity of valuations in \cref{prop:nuvaluationlinear}: if $x\in\Delta_{Y_\bullet}(\theta,\Gamma_\tau)$ and $y\in\Delta_{Y_\bullet}(\theta,\Gamma_{\tau'})$, then $tx+(1-t)y$ is obtained from the current $tS+(1-t)S'$, which is more singular than $\Gamma_{t\tau+(1-t)\tau'}$. Hence
    \[
        t\Delta_{Y_\bullet}\left(\theta,\Gamma_\tau\right)+(1-t)\Delta_{Y_\bullet}\left(\theta,\Gamma_{\tau'}\right)\subseteq\Delta_{Y_\bullet}\left(\theta,\Gamma_{t\tau+(1-t)\tau'}\right).
    \]

    Next we show that as $\tau\to -\infty$, we have
    \[
    \Delta_{Y_{\bullet}}(\theta,\Gamma_{\tau})\xrightarrow{d_{\Haus}}\Delta_{Y_{\bullet}}(\theta,\Gamma_{-\infty}).
    \]

    By \itemref{itm:thm:pobmain:2}, the bodies $\Delta_{Y_{\bullet}}(\theta,\Gamma_{\tau})$ are contained in $\Delta_{Y_{\bullet}}(\theta,\Gamma_{-\infty})$ and increase as $\tau\to-\infty$. Thus it suffices to compute their volumes:
    \[
    \begin{split}
    \lim_{\tau\to -\infty}\vol \Delta_{Y_{\bullet}}(\theta,\Gamma_{\tau})=\frac{1}{n!}\lim_{\tau\to -\infty}\vol (\theta+\ddc\Gamma_{\tau})=\frac{1}{n!}\vol (\theta+\ddc \Gamma_{-\infty}) \\
    =\vol \Delta_{Y_{\bullet}}(\theta,\Gamma_{-\infty}),
    \end{split}
    \]
    where we applied \cref{thm:pobmain} and \cref{thm:contvolu}.

    When $\Gamma\in \TC^{\infty}(X,\theta;\Gamma_{-\infty})$, it is clear that $\Delta_{Y_{\bullet}}(\theta,\Gamma)\in \TC^{\infty}(\Delta_{Y_{\bullet}}(\theta,\Gamma_{-\infty}))$.

    Put $\phi=\Gamma_{-\infty}$ and assume that $\phi$ is $\mathcal I$-model. By \itemref{itm:thm:pobmain:1} and \cref{def:volqpsh},
    \[
        n!\vol\Delta_{Y_{\bullet}}(\theta,\phi)=\int_X\theta_\phi^n,\qquad
        n!\vol\Delta_{Y_{\bullet}}(\theta,\Gamma_\tau)
        =\int_X\bigl(\theta+\ddc P_\theta[\Gamma_\tau]_{\mathcal I}\bigr)^n
        \geq \int_X\theta_{\Gamma_\tau}^n.
    \]
    The inequality follows from \cref{thm:mono}, since $\Gamma_\tau\preceq P_\theta[\Gamma_\tau]_{\mathcal I}$. Comparing \eqref{eq:tcfiniteenergy} and \eqref{eq:Otestcurvenergy}, the Okounkov energy is at least $\mathbf{E}^\phi(\Gamma)$, which proves the finite-energy assertion. If every $\Gamma_\tau$ is $\mathcal I$-model, the displayed inequality is an equality for each $\tau$, giving equality of the energies.
\end{proof}

\begin{remark}\label{rmk:Okounkov:L2335}
    As a special case of this construction, suppose that $\Gamma$ is the test curve induced by a test configuration as in \cref{ex:RWN} and \cref{rmk:PhongSturm}, then for any $\tau<\Gamma_{\max}$, $\Delta_{Y_{\bullet}}(\theta,\Gamma_{\tau})$ is the Okounkov body of a graded linear series
    \[
    \bigoplus_{k=0}^{\infty}\mathcal{F}^{k\tau}_k,
    \]
    where $\mathcal{F}$ is the filtration induced by the test configuration. See \cite[Theorem~5.28]{Xia21} for the details. In particular, in this case, our theory of partial Okounkov bodies recovers the Okounkov bodies of the filtered linear series in the sense of \cite{BC11}.
\end{remark}

\motto{The mathematician's patterns, like the painter's or the poet's must be beautiful; the ideas, like the colors or the words must fit together in a harmonious way. Beauty is the first test: There is no permanent place in this world for ugly mathematics.

--- Godfrey Harold Hardy\footnote{Godfrey Harold Hardy (1877--1947) was a British mathematician famous for his work in number theory and mathematical analysis. Apart from his research, Hardy was a strong advocate for pure mathematics and believed that mathematics should be pursued for its own beauty, not just for practical use.

He remained lifelong unmarried and dedicated much of his life entirely to mathematics, fitting into the common stereotype of a mathematician.}}
\chapter{The theory of b-divisors}
\label{chap:bdiv}

The purpose of this chapter is to translate the singularity theory developed in the preceding chapters into the language of b-divisors. A b-divisor records, on all proper bimeromorphic models of $X$, the cohomological data attached to a singular object.

Given a closed positive $(1,1)$-current $T$ on $X$, we associate to it a b-divisor $\mathbb D(T)$ by taking, on each modification $\pi\colon Y\to X$, the cohomology class of the regular part of $\pi^*T$ with respect to Siu's decomposition. This construction removes divisorial components and keeps track of the moving positivity of $T$ on every model. The resulting b-divisor is nef in a suitable sense, and its volume agrees with the volume of $T$.

The first part of the chapter develops this construction and studies its basic properties. We prove that $d_S$-convergence of currents implies convergence of the associated b-divisors, and we compare the $P$- and $\mathcal I$-orders through their b-divisorial data. This gives a birational way to read some of the singularity information encoded by envelopes and multiplier ideals.

The second part relates b-divisors to $\mathcal I$-good currents. In the big case, $\mathcal I$-good non-divisorial currents give rise to big and nef b-divisors, and conversely big and nef b-divisors can be represented by such currents.

The final part defines mixed intersection numbers of nef b-divisors. For big and nef b-divisors this is done through the corresponding $\mathcal I$-good currents and the mixed volume formalism of \cref{chap:Igood}; the general nef case is obtained by K\"ahler perturbation. We prove symmetry, multilinearity, monotonicity, continuity under decreasing nets, and upper semicontinuity. These results provide the b-divisorial counterpart of the mixed volume theory.

\section{The notion of b-divisors}\label{sec:bdivnotion}

The b-divisors defined in this section are sometimes called b-divisor classes. We always omit the word \emph{classes} to save space.

Let $X$ be a connected compact Kähler manifold of dimension $n$.

Let us recall the following elementary result regarding how cohomology behaves under blow-up.
\begin{proposition}\label{prop:blowup_cone}
    Let $\pi\colon Y\rightarrow X$ be a blow-up with \emph{connected} smooth center of codimension at least $2$ with exceptional divisor $E$. Then there is a natural identification
    \begin{equation}\label{eq:H11blowup}
    \mathrm{H}^{1,1}(Y,\mathbb{R})=\mathrm{H}^{1,1}(X,\mathbb{R})\oplus \mathbb{R}\{E\}.
    \end{equation}
\end{proposition}
See \cite{RYY19} for a much more general result.

\subsection{Nef b-divisors}
\label{subsec:nefbdiv}
\begin{definition}\label{def:bdiv:L38}
    A \emph{(Weil) b-divisor}\index{b-divisor} $\mathbb{D}$ over $X$ is an assignment $(\mathbb{D}_{\pi})_{\pi\colon Y\rightarrow X}$\index{$\mathbb{D}$}, where $\pi\colon Y\rightarrow X$ runs over all modifications of $X$ such that
    \begin{enumerate}
        \item\label{itm:chapter-bdiv:L43:1} $\mathbb{D}_{\pi}\in \mathrm{H}^{1,1}(Y,\mathbb{R})$;
        \item\label{itm:chapter-bdiv:L43:2} The classes are compatible under push-forwards: If $\pi'\colon Z\rightarrow X$ and $\pi\colon Y\rightarrow X$ are both modifications of $X$ and $\pi'$ dominates $\pi$ through $g\colon Z\rightarrow Y$, then $g_*\mathbb{D}_{\pi'}=\mathbb{D}_{\pi}$.
    \end{enumerate}
    We also write $\mathbb{D}_Y=\mathbb{D}_{\pi}$ if there is no risk of confusion.

    Given two Weil b-divisors $\mathbb{D}$ and $\mathbb{D}'$ over $X$, we say $\mathbb{D}\leq \mathbb{D}'$ if for each modification $\pi$ of $X$, we have $\mathbb{D}_{\pi}\leq \mathbb{D}'_{\pi}$. Recall that by definition, this means the class $\mathbb{D}'_{\pi}-\mathbb{D}_{\pi}$ is pseudo-effective.

    A Weil b-divisor with $\mathbb{D}\geq 0$ is called \emph{pseudo-effective}\index{b-divisor!pseudo-effective b-divisor}.
\end{definition}
The class $\mathbb{D}_X$ is called the \emph{root}\index{root!root of a b-divisor} of $\mathbb{D}$.
The set of Weil b-divisors over $X$ has the natural structure of a real vector space.

\begin{definition}\label{def:volb}
    The \emph{volume}\index{volume!volume of a b-divisor} of a pseudo-effective Weil b-divisor $\mathbb{D}$ over $X$ is
    \[
        \vol \mathbb{D}\coloneqq \lim_{\pi\colon Y\rightarrow X} \vol \mathbb{D}_Y.
    \]
    The right-hand side is a decreasing net due to \cref{prop:cones_modif2}, hence the limit always exists. The volume of a non-pseudo-effective Weil b-divisor is defined to be $0$.

    We say $\mathbb{D}$ is \emph{big}\index{b-divisor!big b-divisor} if $\vol \mathbb{D}>0$.
\end{definition}

\begin{lemma}\label{lma:volbdivnet}
    Let $(\mathbb{D}_i)_{i\in I}$ be a net of pseudo-effective b-divisors converging to a Weil b-divisor $\mathbb{D}$. Then
    \begin{equation}\label{eq:vollimsup}
    \varlimsup_{i\in I}\vol \mathbb{D}_i\leq \vol \mathbb{D}.
    \end{equation}
    If the net is decreasing, then
    \[
    \lim_{i\in I}\vol \mathbb{D}_i=\vol \mathbb{D}.
    \]
\end{lemma}
Here we say $(\mathbb{D}_i)_{i\in I}$ converges to $\mathbb{D}$ if for any modification $\pi\colon Y\rightarrow X$, we have $\mathbb{D}_{i,Y}\to \mathbb{D}_Y$ with respect to the Euclidean topology.
\begin{proof}
    First observe that $\mathbb{D}$ is pseudo-effective, since $\mathbb{D}_i\geq 0$ for all $i$.

    Let $\pi\colon Y\rightarrow X$ be a modification. Then
    \[
    \vol \mathbb{D}_Y=\lim_{i\in I} \vol \mathbb{D}_{i,Y}\geq \varlimsup_{i\in I} \vol \mathbb{D}_i.
    \]
    The inequality \eqref{eq:vollimsup} follows.
    As for the decreasing case, it suffices to observe that both sides of \eqref{eq:vollimsup} can be written as
    \[
    \inf_i \inf_{\pi\colon Y\rightarrow X} \vol \mathbb{D}_{i,Y}.
    \]
    This completes the proof.
\end{proof}

\begin{example}\label{ex:volbdivnetstrict}
    In general, equality need not hold in \eqref{eq:vollimsup}. Let $X=\mathbb{P}^2$, and let $H=c_1(\mathcal{O}_{\mathbb{P}^2}(1))$. Choose a sequence of points
    $p_i\in X$ converging to the generic point in the following sense: for every algebraic curve $C\subseteq X$, the set of indices $i$ with $p_i\in C$ is finite.

    Let
    \[
        \pi_i\colon X_i\rightarrow X
    \]
    be the blow-up of $X$ at $p_i$, with exceptional divisor $E_i$. Let $\mathbb{D}_i$ be the Cartier b-divisor realized on $X_i$ by the nef class
    \[
        \pi_i^*H-\{E_i\},
    \]
    and let $\mathbb{D}$ be the Cartier b-divisor realized on $X$ by $H$. Then $\mathbb{D}_i\to \mathbb{D}$.
    Indeed, for any fixed modification $\rho\colon Y\rightarrow X$, the morphism $\rho$ is an isomorphism over $p_i$ for all large enough $i$, hence the trace of $\mathbb{D}_i$ on $Y$ is $\rho^*H$ for all large enough $i$.

    On the other hand,
    \[
        \vol \mathbb{D}_i=\left(\pi_i^*H-\{E_i\}\right)^2=0
    \]
    for all $i$, whereas
    \[
        \vol \mathbb{D}=H^2=1.
    \]
    Thus the inequality in \eqref{eq:vollimsup} is strict for this sequence.
\end{example}

\begin{definition}\label{def:bdiv:L116}
    A \emph{Cartier b-divisor}\index{b-divisor!Cartier b-divisor} $\mathbb{D}$ over $X$ is a Weil b-divisor $\mathbb{D}$ over $X$ such that there exists a modification $\pi\colon Y\rightarrow X$ and a class $\alpha_Y\in \mathrm{H}^{1,1}(Y,\mathbb{R})$ so that for each $\pi'\colon Z\rightarrow X$ dominating $\pi$, the class $\mathbb{D}_Z$ is the pull-back of $\alpha_Y$. Any such $(\pi,\alpha_Y)$ is called a \emph{realization}\index{realization!realization of a Cartier b-divisor} of $\mathbb{D}$.
\end{definition}
By abuse of language, we also say $(Y,\alpha_Y)$ is a realization of $\mathbb{D}$.
The realization is not unique in general.

\begin{definition}\label{def:nefCar}
    A Cartier b-divisor $\mathbb{D}$ over $X$ is \emph{nef}\index{b-divisor!Cartier b-divisor!nef Cartier b-divisor} if there exists a realization $(\pi\colon Y\rightarrow X, \alpha_Y)$ of $\mathbb{D}$ such that $\alpha_Y$ is nef.
\end{definition}

\begin{definition}\label{def:nefb}
    A Weil b-divisor $\mathbb{D}$ over $X$ is \emph{nef}\index{b-divisor!nef b-divisor} if it is the limit of a net of nef Cartier b-divisors $(\mathbb{D}_i)_i$ over $X$.
\end{definition}
In other words, for each modification $\pi\colon Y\rightarrow X$, we have $\mathbb{D}_{i,Y}\to \mathbb{D}_Y$.

Note that thanks to \cref{prop:cones_modif}, each $\mathbb{D}_Y$ is necessarily modified nef, but it is not nef in general. The notion of modified nef classes is defined in \cref{def:mnclass}.

\emph{A priori}, for a Cartier b-divisor, nefness could mean either the notion defined by \cref{def:nefCar} or the one defined by \cref{def:nefb}. We will show in \cref{cor:nefbCequiv} that they are actually equivalent. Before that, by a nef Cartier b-divisor, we always mean nefness in the sense of \cref{def:nefCar}.

Our definition \cref{def:nefb} amounts to defining the set of Weil b-divisors as the closure of the set of Cartier b-divisors in $\varprojlim_{\pi}\mathrm{H}^{1,1}(Y,\mathbb{R})$ with respect to the projective limit topology. In particular, the limit of a converging net of nef b-divisors is still nef.

\subsection{The b-divisors of currents}\label{subsec:b-div_curr}

\begin{definition}\label{def:DT}
Let $T$ be a closed positive $(1,1)$-current on $X$. Given any modification $\pi\colon Y\rightarrow X$, we define
\begin{equation}
\mathbb{D}(T)_Y\coloneqq \left\{ \Reg \pi^*T \right\}\in \mathrm{H}^{1,1}(Y,\mathbb{R}).
\end{equation}
We call $\mathbb{D}(T)$\index{$\mathbb{D}(T)$} the \emph{b-divisor of $T$}\index{b-divisor!b-divisor of a current}.
\end{definition}
This assignment is compatible with push-forwards. Indeed, if $g\colon Z\to Y$ is a modification and $\pi\colon Y\to X$, then
\[
    g_*\Reg\bigl(g^*\pi^*T\bigr)=\Reg(\pi^*T)
\]
by Siu's decomposition and \cref{lma:Siudec_pullpush}. Hence
\[
    g_*\mathbb D(T)_Z=\mathbb D(T)_Y.
\]
Thus $\mathbb D(T)$ is a Weil b-divisor.

We observe that if $T'$ is another closed positive $(1,1)$-current on $X$ and $\lambda\geq 0$, then
\[
\mathbb{D}(T+T')=\mathbb{D}(T)+\mathbb{D}(T'),\quad \mathbb{D}(\lambda T)=\lambda \mathbb{D}(T).
\]
We shall use these identities implicitly in the sequel.

Note that when $T$ has analytic singularities, $\mathbb{D}(T)$ is Cartier.

\begin{theorem}\label{thm:DT}
    Let $T$ be a closed positive $(1,1)$-current on $X$. Then $\mathbb{D}(T)$ is nef. Moreover,
    \begin{equation}\label{eq:volDT}
    \vol T=\vol \mathbb{D}(T).
    \end{equation}
\end{theorem}

\begin{proof}
    Let $\omega$ be a Kähler form on $X$.

    \textbf{Step~1}. Reduce to the case where $T$ is a Kähler current.

    Note that $\mathbb{D}(\omega)$ is the Cartier b-divisor realized by $(X,\{\omega\})$.
    The b-divisors $\mathbb{D}(T+\epsilon\omega)=\mathbb{D}(T)+\epsilon \mathbb{D}(\omega)$ decrease to $\mathbb{D}(T)$ as $\epsilon\searrow 0$. Hence, by \cref{lma:volbdivnet} and \cref{prop:vollinearlimit}, it suffices to prove our assertion with $T+\epsilon\omega$ in place of $T$.

    \textbf{Step~2}. We prove the assertion under the additional assumption that $T$ has analytic singularities.

    Let $\pi\colon Y\rightarrow X$ be a modification so that
    \[
    \pi^*T=[D]+R,
    \]
    where $D$ is an effective $\mathbb{Q}$-divisor on $Y$ and $R$ is a closed positive $(1,1)$-current with locally bounded potentials. See \cref{thm:resolvelogsing} for the existence of $\pi$.
    Then $\mathbb{D}(T)$ is the nef Cartier b-divisor realized by $(\pi,\{R\})$. Note that \eqref{eq:volDT} is immediate in this case.

    \textbf{Step~3}. We prove the assertion for a general Kähler current $T$.
    Next, take a closed smooth real $(1,1)$-form $\theta$ in the cohomology class of $T$ and write $T=\theta_{\varphi}$ for some $\varphi\in \PSH(X,\theta)$. Since $T$ is a Kähler current, we choose the quasi-equisingular approximation $(\varphi_j)_j$ of $\varphi$ so that, after discarding finitely many terms, each $\theta+\ddc\varphi_j$ is a Kähler current. This is possible by Demailly regularization, in the form used in the proof of \cref{thm:charIgoodasclosure}. We claim that
    \begin{equation}\label{eq:conv_bdiv_temp1}
        \mathbb{D}(\theta+\ddc \varphi_j)\to \mathbb{D}(\theta+\ddc \varphi).
    \end{equation}
    By definition of this convergence, we need to establish the following: for every modification $\pi\colon Y\rightarrow X$, we have
    \[
        \left\{\Reg\left(\pi^*(\theta+\ddc\varphi_j)\right)\right\}\to \left\{\Reg\left(\pi^*(\theta+\ddc\varphi)\right)\right\}.
    \]
    This follows immediately from \cref{thm:Lelongcont} if $\Sing (\pi^*T)$ has only finitely many components.

    Fix a Kähler form $\omega_Y$ on $Y$. We want to show that for any $\epsilon>0$, we can find $j_0>0$ so that when $j\geq j_0$,
\begin{equation}\label{eq:Regpivarphij_conv_temp1}
\left\{\Reg\left(\pi^*(\theta+\ddc\varphi_j)\right)\right\}\leq \left\{\Reg\left(\pi^*(\theta+\ddc\varphi)\right)\right\}+\epsilon\{\omega_Y\}.
\end{equation}
Write the divisorial part of $\pi^*\theta+\ddc \pi^*\varphi_j$ and $\pi^*\theta+\ddc \pi^*\varphi$ as
\[
\sum_{i=1}^{\infty}a_i^j [E_i],\quad \sum_{i=1}^{\infty}a_i [E_i].
\]
Then $a_i^j\leq a_i$.

We can find an integer $N>0$ large enough so that, in the pseudo-effective order of cohomology classes,
\[
\sum_{i=N+1}^{\infty}a_i \{E_i\}\leq \frac{\epsilon}{2}\{\omega_Y\}.
\]
By \cref{thm:Lelongcont}, we can take $j_0$ large enough so that for $j>j_0$,
\[
\sum_{i=1}^{N}(a_i-a_i^j)\{E_i\}\leq \frac{\epsilon}{2}\{\omega_Y\}.
\]
Then \eqref{eq:Regpivarphij_conv_temp1} follows, and \eqref{eq:conv_bdiv_temp1} is established.

As a consequence, $\mathbb{D}(T)$ is nef and the volume can be computed using \cref{lma:volbdivnet}:
\[
\vol\mathbb{D}(T)=\lim_{j\to\infty}\vol \mathbb{D}\left( \theta+\ddc\varphi_j \right)=\lim_{j\to\infty}\vol \left(\theta+\ddc\varphi_j\right)=\vol T.
\]
Hence, \eqref{eq:volDT} follows.
\end{proof}

Conversely, we want to realize nef b-divisors as $\mathbb{D}(T)$. We first prove a continuity result.
\begin{proposition}\label{prop:dscontbdiv}
Let $\theta$ be a closed real smooth $(1,1)$-form on $X$, let $(\varphi_i)_{i\in I}$ be a net in $\PSH(X,\theta)$, and let $\varphi\in \PSH(X,\theta)$. If $\varphi_i\xrightarrow{d_S}\varphi$, then
\begin{equation}\label{eq:Dcont1}
    \mathbb{D}\left(\theta_{\varphi_i}\right)\to \mathbb{D}\left(\theta_{\varphi}\right).
\end{equation}
\end{proposition}
\begin{proof}
    Fix a modification $\pi\colon Y\rightarrow X$.
    Fix a Kähler form $\omega$ on $X$.
    It suffices to establish the following:
    \begin{equation}\label{eq:Dcont2}
    \mathbb{D}\left(\theta_{\varphi_i}\right)_Y\to \mathbb{D}\left(\theta_{\varphi}\right)_Y.
    \end{equation}
    As a map from the pseudometric space $\PSH(X,\theta)$ to the finite-dimensional Euclidean space $\mathrm{H}^{1,1}(Y,\mathbb{R})$, the continuity of $\mathbb{D}(\bullet)_Y$ can be tested on sequences. Thus, without loss of generality, we may assume that $(\varphi_i)_i$ is a sequence and $I=\mathbb{Z}_{>0}$. It is enough to prove that every subsequence has a further subsequence for which \eqref{eq:Dcont3} holds.

    By \cref{cor:dSconv_bimero}, $\pi^*\varphi_i\xrightarrow{d_S}\pi^*\varphi$. Hence it suffices to prove the claim on $Y$, so we may assume that $\pi$ is the identity. Therefore, we are reduced to the following assertion:
    \begin{equation}\label{eq:Dcont3}
        \left\{ \Reg \theta_{\varphi_i} \right\}\to \left\{ \Reg \theta_{\varphi} \right\}.
    \end{equation}
    After adding a Kähler form to $\theta$, we may also assume that $\theta_{\varphi}$ is a Kähler current. By \cref{thm:contPI}, replacing all potentials by their $\mathcal{I}$-envelopes does not change the divisorial parts and preserves $d_S$-convergence. The masses are then bounded from below along a tail, so after passing to a subsequence \cref{prop:incanddec} supplies increasing and decreasing envelopes $\eta_i\leq \varphi_i\leq \psi_i$ which both converge to $\varphi$ in $d_S$. Since
    \[
        \mathbb{D}(\theta_{\eta_i})_X\leq \mathbb{D}(\theta_{\varphi_i})_X\leq \mathbb{D}(\theta_{\psi_i})_X,
    \]
    it is enough to treat the case where $(\varphi_i)_i$ is monotone.

    Write the divisorial parts of $\theta_{\varphi_i}$ and $\theta_{\varphi}$ as
    \[
        \sum_{m=1}^{\infty}a_m^i[D_m],\quad \sum_{m=1}^{\infty}a_m[D_m].
    \]
    By \cref{thm:Lelongcont}, $a_m^i\to a_m$ for each fixed $m$. If $(\varphi_i)_i$ is decreasing, then $a_m^i\leq a_m$ for all $i,m$. The finite-tail argument in the proof of \eqref{eq:Regpivarphij_conv_temp1} then gives
    \[
        \left\{\Reg\theta_{\varphi}\right\}\leq \left\{\Reg\theta_{\varphi_i}\right\}\leq \left\{\Reg\theta_{\varphi}\right\}+\epsilon\{\omega\}
    \]
    for all sufficiently large $i$.

    If $(\varphi_i)_i$ is increasing, then $a_m\leq a_m^i\leq a_m^1$. We first choose $N$ so that the tails of both $\sum_m a_m[D_m]$ and $\sum_m a_m^1[D_m]$ are at most $(\epsilon/3)\{\omega\}$ in the pseudo-effective order. The convergence of the finitely many coefficients $a_m^i$, $m\leq N$, then gives
    \[
        \left\{\Reg\theta_{\varphi_i}\right\}\leq \left\{\Reg\theta_{\varphi}\right\}\leq \left\{\Reg\theta_{\varphi_i}\right\}+\epsilon\{\omega\}
    \]
    for all sufficiently large $i$. This proves \eqref{eq:Dcont3}.
\end{proof}

\begin{lemma}\label{lma:pushIpo}
    Let $\pi\colon X\rightarrow Z$ be a proper bimeromorphic morphism from $X$ to a Kähler manifold $Z$. Consider non-divisorial closed positive $(1,1)$ currents $T,S$ on $X$ in the same cohomology class. If $T\preceq_{\mathcal{I}} S$, then $\pi_*T\preceq_{\mathcal{I}} \pi_*S$.
\end{lemma}
\begin{proof}
    By Hironaka's Chow lemma \cref{thm:HironakaChow} and \cref{lma:Ipartorderinv}, we may assume that $\pi$ is a modification.

    By \cref{lma:Siudec_pullpush},
    \[
    \pi^*\pi_*T=T+F_T,\quad \pi^*\pi_*S=S+F_S,
    \]
    where $F_T$ and $F_S$ are effective $\pi$-exceptional divisors. Since $\{T\}=\{S\}$, the exceptional divisors $F_T$ and $F_S$ have the same cohomology class. Factoring $\pi$ into blow-ups and using \cref{prop:blowup_cone} at each step, the classes of $\pi$-exceptional prime divisors are linearly independent in the relative cohomology; hence $F_T=F_S$. It follows that
    \[
    \pi^*\pi_*T=T+F_T\preceq_{\mathcal{I}}S+F_S=\pi^*\pi_*S.
    \]
    Our assertion then follows from \cref{lma:Ipartorderinv}.
\end{proof}
\begin{theorem}\label{thm:bigandnefreal}
    Each big and nef b-divisor $\mathbb{D}$ over $X$ can be realized as $\mathbb{D}(T)$ for some $T\in \mathbb{D}_X$. Furthermore, we may always assume that $T$ is $\mathcal{I}$-good.
\end{theorem}
Note that $T$ is not unique. The current $T$ is necessarily non-divisorial.

\begin{proof}
    Fix a big and nef b-divisor $\mathbb{D}$ over $X$.

    For each $\pi\colon Y\rightarrow X$, we take a current with minimal singularities $T_Y$ in $\mathbb{D}_Y$.
    We claim that $\mathbb{D}(\pi_*T_Y)$ coincides with $\mathbb{D}$ up to the level of $Y$: For any modification $\pi'\colon Z\rightarrow X$ dominated by $\pi$ through a morphism $g\colon Y\rightarrow Z$, we have
    \[
    \mathbb{D}_Z=\mathbb{D}(\pi_*T_Y)_Z.
    \]
    The notations are summarized in the following commutative diagram:
    \begin{equation}\label{eq:domi2}
    \begin{tikzcd}
Y \arrow[rr,"g"] \arrow[rd, "\pi"'] &   & Z \arrow[ld, "\pi'"] \\
                                 & X. &
\end{tikzcd}
    \end{equation}
    After unfolding the definitions, this means
    \[
        \Reg(\pi'^*\pi_*T_Y)\in \mathbb{D}_Z.
    \]
    Note that
    \[
        \Reg(\pi'^*\pi_*T_Y)=\Reg(\pi'^* \pi'_* g_*T_Y).
    \]
    Due to \cref{prop:cones_modif} and \cref{prop:cones_modif2}, we know that $\mathbb{D}_Y$ is modified nef and big. In particular, $T_Y$ is non-divisorial, hence so is $g_*T_Y$ by \cref{lma:ndvpush}. It follows from \cref{lma:Siudec_pullpush} that
    \[
    \Reg(\pi'^* \pi'_* g_*T_Y)=\Reg(g_*T_Y)=g_*T_Y\in \mathbb{D}_Z.
    \]
    Note that
    \begin{equation}\label{eq:vol_unif_bound}
    \vol T_Y\geq \vol \mathbb{D}>0.
    \end{equation}

    Next we claim that the $P$-singularity types of the net $(\pi_*T_Y)_Y$ form a decreasing net.

    To see this, let us fix a diagram as \eqref{eq:domi2}. We need to show that
    \[
    \pi_*T_Y\preceq_P \pi'_*T_Z.
    \]
    Since $T_Z$ has minimal singularities, we have $g_*T_Y\preceq_{\mathcal{I}} T_Z$. In particular, \cref{lma:pushIpo} guarantees that $\pi_*T_Y\preceq_{\mathcal{I}} \pi'_*T_Z$.
    But by \cref{cor:Igoodpush} and \cref{ex:ImodelIgood}, both $\pi_*T_Y$ and $\pi'_*T_Z$ are $\mathcal{I}$-good, so there is no difference between the $P$-preorder and the $\mathcal{I}$-preorder in this case. Our assertion follows.

    Next observe that the net $(\pi_*T_Y)_Y$ is decreasing with respect to the $P$-preorder and has a $d_S$-limit as a consequence of \eqref{eq:vol_unif_bound} and \cref{cor:completenessdS}.
    Take a closed positive $(1,1)$-current $T\in \mathbb{D}_X$ such that
    \[
    \pi_*T_Y\xrightarrow{d_S}T.
    \]
    It follows from \cref{prop:dscontbdiv} that
    \[
    \mathbb{D}(\pi_*T_Y)\to \mathbb{D}(T).
    \]
    Therefore, we conclude that
    \[
    \mathbb{D}(T)=\mathbb{D}.
    \]
    Thanks to \cref{thm:DT}, $\vol T>0$. Write $T=\theta+\ddc \varphi$ for some $\varphi\in \PSH(X,\theta)$, then
    \[
    T'\coloneqq \theta+\ddc P_{\theta}[\varphi]_{\mathcal{I}}
    \]
    is $\mathcal{I}$-good, non-divisorial and $\mathbb{D}(T')=\mathbb{D}(T)$.
\end{proof}

\begin{corollary}\label{cor:mn_level_nefb}
    Let $\mathbb{D}$ be a b-divisor over $X$. Then the following are equivalent:
    \begin{enumerate}
        \item\label{itm:cor:mn_level_nefb:1} $\mathbb{D}$ is nef;
        \item\label{itm:cor:mn_level_nefb:2} for each modification $\pi\colon Y\rightarrow X$, the class $\mathbb{D}_Y$ is modified nef.
    \end{enumerate}
\end{corollary}
\begin{proof}
    \ref{itm:cor:mn_level_nefb:1} $\implies$ \ref{itm:cor:mn_level_nefb:2}. Suppose that \ref{itm:cor:mn_level_nefb:1} holds. Fix a modification $\pi\colon Y\rightarrow X$. Then we need to show that $\mathbb{D}_Y$ is modified nef. Since modified nefness is a closed condition, after approximating $\mathbb{D}$ by nef Cartier b-divisors, we may assume that $\mathbb{D}$ itself is a nef Cartier b-divisor. We can then find a modification $\pi'\colon Z\rightarrow X$ dominating $\pi$ so that $\mathbb{D}$ is realized by a nef class $\alpha$ on $Z$. Then $\mathbb{D}_Y$ is nothing but the pushforward of $\alpha$, and hence modified nef thanks to \cref{prop:cones_modif}.

    \ref{itm:cor:mn_level_nefb:2} $\implies$ \ref{itm:cor:mn_level_nefb:1}. Suppose that \ref{itm:cor:mn_level_nefb:2} holds. We need to show that $\mathbb{D}$ is nef. Fix a Kähler form $\omega$ on $X$. It suffices to show that $\mathbb{D}+\epsilon\mathbb{D}(\omega)$ is nef for each $\epsilon>0$. After replacing $\mathbb{D}$ by the latter, we may further assume that $\mathbb{D}_Y$ is big for each modification $\pi\colon Y\rightarrow X$, and $\vol \mathbb{D}_Y$ has a uniform positive lower bound.
    In this case, the argument of \cref{thm:bigandnefreal} shows that $\mathbb{D}=\mathbb{D}(T)$ for some closed positive $(1,1)$-current in $\mathbb{D}_X$. Therefore, $\mathbb{D}$ is nef, thanks to \cref{thm:DT}.
\end{proof}

Let $\alpha$ be a modified nef class on $X$. We write $\mathcal{G}(\alpha)$ for the set of closed positive $(1,1)$-currents $T$ on $X$ with $T=\Reg T\in \alpha$ and $\vol T>0$.
\begin{theorem}\label{thm:bdiv_current_corr}
    There is a natural bijection from $\mathcal{G}(\alpha)/\sim_{\mathcal{I}}$ to the set of big and nef b-divisors $\mathbb{D}$ over $X$ with $\mathbb{D}_X=\alpha$.
\end{theorem}
\begin{proof}

    Given $T\in \mathcal{G}(\alpha)$, we associate the b-divisor $\mathbb{D}(T)$. It is big and nef due to \cref{thm:DT}. This map clearly descends to $\mathcal{G}(\alpha)/\sim_{\mathcal{I}}$.

    This map is surjective by \cref{thm:bigandnefreal}. Now we show that it is injective. Let $T,T'\in \mathcal{G}(\alpha)$. Assume that $\mathbb{D}(T)=\mathbb{D}(T')$, we want to show that $T\sim_{\mathcal{I}} T'$.

    Let $E$ be a prime divisor over $X$, it suffices to show that
    \begin{equation}\label{eq:nuTE}
    \nu(T,E)=\nu(T',E).
    \end{equation}
    We may assume that $E$ is not a prime divisor on $X$, as otherwise both sides vanish.

    Choose a sequence of blow-ups with smooth \emph{connected} centers
    \[
    Y\coloneqq X_k\rightarrow X_{k-1}\rightarrow \cdots \rightarrow X_0\coloneqq X
    \]
    so that $E$ is a prime divisor on $Y$, exceptional with respect to $X_k\rightarrow X_{k-1}$. Denote the composition by $\pi\colon Y\rightarrow X$.
    Thanks to \cref{prop:blowup_cone},
    \[
    \mathrm{H}^{1,1}(X_k,\mathbb{R})=\mathrm{H}^{1,1}(X_{k-1},\mathbb{R})\oplus \mathbb{R}\{E_k\},
    \]
    where $E_k=E$ is the exceptional divisor of $X_k\rightarrow X_{k-1}$.

    By induction,
    \[
    \mathrm{H}^{1,1}(Y,\mathbb{R})=\mathrm{H}^{1,1}(X,\mathbb{R})\oplus \bigoplus_{i=1}^{k} \mathbb{R} \{E_i\},
    \]
    where $E_i$ is the exceptional divisor of $X_i\rightarrow X_{i-1}$. Now by \cref{lma:Siudec_pullpush},
    \begin{equation}\label{eq:RegpistarT}
    \Reg \pi^*T=\pi^*T-\sum_{i=1}^k \nu(T,E_i)[E_i].
    \end{equation}
    In particular, the cohomology class of $\Reg \pi^*T$ determines $\nu(T,E)$. Hence, \eqref{eq:nuTE} follows.
\end{proof}

\begin{corollary}\label{cor:nef_b_degenerate_1}
    The set of nef b-divisors over $X$ with root $\alpha$ can be naturally identified with
    \[
    \varprojlim_{\omega} \left( \mathcal{G}(\alpha+\omega)/\sim_{\mathcal{I}} \right),
    \]
    where $\omega$ runs over the directed set of Kähler forms on $X$ (with respect to the preorder of reverse domination), and given two Kähler forms $\omega\leq \omega'$ the transition map
    \[
    \mathcal{G}(\alpha+\omega)/\sim_{\mathcal{I}}\rightarrow \mathcal{G}(\alpha+\omega')/\sim_{\mathcal{I}}
    \]
    is induced by the map $\mathcal{G}(\alpha+\omega)\rightarrow \mathcal{G}(\alpha+\omega')$ sending $T$ to $T+\omega'-\omega$.
\end{corollary}
\begin{proof}
For each Kähler form $\omega$, \cref{thm:bdiv_current_corr} identifies $\mathcal G(\alpha+\omega)/\sim_{\mathcal I}$ with the set of big and nef b-divisors whose root on $X$ is $\alpha+\{\omega\}$. The transition map is exactly addition of $\mathbb D(\omega'-\omega)$. Passing to the inverse limit  gives precisely nef b-divisors with root $\alpha$.
\end{proof}

It is tempting to extend these results to more general currents, not necessarily non-divisorial ones. For this purpose, we introduce the following definition:
\begin{definition}\label{def:augmented_nefb}
    An \emph{augmented nef b-divisor}\index{augmented nef b-divisor} over $X$ is a pair $(\mathbb{D},D)$, where
    \begin{enumerate}
        \item\label{itm:def:augmented_nefb:1} $\mathbb{D}$ is a nef b-divisor over $X$;
        \item\label{itm:def:augmented_nefb:2} $D$ is a formal sum
        \[
        D=\sum_{E} c_E E,\quad c_E\in \mathbb{R}_{\geq 0},
        \]
        where $E$ runs over the set of prime divisors on $X$,
    \end{enumerate}
    such that the following condition is satisfied:
    \[
    \sum_{E}c_E \{E\}
    \]
    is convergent as a sum in $\mathrm{H}^{1,1}(X,\mathbb{R})$.

    The \emph{cohomology class} of $(\mathbb{D},D)$ is defined as
    \[
    \mathbb{D}_X+\sum_{E}c_E \{E\}\in \mathrm{H}^{1,1}(X,\mathbb{R}).
    \]
    The \emph{volume} of $(\mathbb{D},D)$ is defined as
    \[
    \vol (\mathbb{D},D)=\vol \mathbb{D}.
    \]
\end{definition}

Recall that we defined $\mathcal{Z}_+(X,\alpha)$ as the set of closed positive $(1,1)$-currents $T\in \alpha$ in \cref{def:spaceofcurrents}. We introduce a further notation here:
\[
\mathcal{Z}_+(X,\alpha)_{>0}\coloneqq \left\{T\in \mathcal{Z}_+(X,\alpha): \vol T>0\right\}.
\]
\index{$\mathcal{Z}_+(X,\alpha)_{>0}$}

\begin{corollary}\label{cor:Iequiv_curr_b_div}
    There is a canonical bijection between the following sets:
    \begin{enumerate}
        \item\label{itm:cor:Iequiv_curr_b_div:1} The set $\mathcal{Z}_+(X,\alpha)_{>0}/\sim_{\mathcal{I}}$;
        \item\label{itm:cor:Iequiv_curr_b_div:2} the set of augmented nef b-divisors over $X$ with positive volume with cohomology class $\alpha$.
    \end{enumerate}
\end{corollary}
More precisely, given a current $T\in \mathcal{Z}_+(X,\alpha)_{>0}$, we associate $(\mathbb{D},D)$ as follows:
\[
\mathbb{D}=\mathbb{D}(T),\quad D=\sum_E \nu(T,E)E,
\]
where $E$ runs over the set of prime divisors on $X$.
\begin{proof}
    This is a direct consequence of \cref{thm:bdiv_current_corr} and Siu's decomposition \cref{lma:Siudec}.
\end{proof}

In fact, \cref{cor:Iequiv_curr_b_div} can also be reformulated in an elementary manner, without referring to b-divisors at all. Suppose that we are given a big cohomology class $\alpha$ on $X$ and a non-negative real number $c_E$ for each prime divisor $E$ over $X$. A natural question is: When is there a closed positive $(1,1)$-current $T\in \alpha$ with positive volume such that $\nu(T,E)=c_E$ for every $E$? Then \cref{cor:Iequiv_curr_b_div} says that $T$ exists if and only if the following two conditions hold:
\begin{enumerate}
    \item\label{itm:chapter-bdiv:L472:1}
    \[
    \lim_{\pi\colon Y\rightarrow X} \vol\left( \pi^*\alpha-\sum_{E\subseteq Y} c_E \{E\} \right)>0\footnotemark,
    \]
    \footnotetext{In particular, implicitly, the sum $\sum_{E\subseteq Y} c_E \{E\}$ converges in $\mathrm{H}^{1,1}(Y,\mathbb{R})$.}
    where $\pi$ runs over all modifications of $X$;
    \item\label{itm:chapter-bdiv:L472:2} for each modification $\pi\colon Y\rightarrow X$, the class $ \pi^*\alpha-\sum_{E\subseteq Y} c_E \{E\}$ is modified nef.
\end{enumerate}

Similarly, we have the following generalization of \cref{cor:nef_b_degenerate_1}.
\begin{corollary}\label{cor:nef_b_degenerate_2}
    There is a canonical bijection between the following two sets:
    \begin{enumerate}
        \item\label{itm:cor:nef_b_degenerate_2:1} The set of
        \[
            \varprojlim_{\omega} \left( \mathcal{Z}_+(X,\alpha+\omega)_{>0}/\sim_{\mathcal{I}} \right),
        \]
        where $\omega$ runs over the directed set of Kähler forms on $X$;
    \item\label{itm:cor:nef_b_degenerate_2:2} the set of augmented nef b-divisors over $X$ with cohomology class $\alpha$.
    \end{enumerate}
\end{corollary}
\begin{proof}
Apply \cref{cor:nef_b_degenerate_1} to the non-divisorial parts and use Siu's decomposition on $X$ to record the divisorial part. The convergence condition in \cref{def:augmented_nefb} is exactly the condition that the roots of the non-divisorial parts converge to the prescribed cohomology class $\alpha$.
\end{proof}

\begin{corollary}\label{cor:nefbCequiv}
    Let $\mathbb{D}$ be a Cartier b-divisor over $X$. Then $\mathbb{D}$ is nef in the sense of \cref{def:nefCar} if and only if it is nef in the sense of \cref{def:nefb}.
\end{corollary}

\begin{proof}
    We only handle the nontrivial implication. Assume that $\mathbb{D}$ is nef in the sense of \cref{def:nefb}. We want to show that $\mathbb{D}$ is nef in the sense of \cref{def:nefCar}. After replacing $X$ by a realization of $\mathbb{D}$, we may assume that $\mathbb{D}$ is realized by $(X,\alpha)$ for some cohomology class $\alpha\in \mathrm{H}^{1,1}(X,\mathbb{R})$.

    We first assume that $\mathbb{D}$ is big. Choose a non-divisorial closed positive $(1,1)$-current $T$ on $X$ such that $\mathbb{D}=\mathbb{D}(T)$. Now $\mathbb{D}=\mathbb{D}(T)$ means that for each modification $\pi\colon Y\rightarrow X$, the current $\pi^*T$ is non-divisorial. In particular, $T$ has vanishing generic Lelong number along each prime divisor over $X$, see \eqref{eq:RegpistarT}. Equivalently, $T$ has vanishing Lelong number everywhere. By Demailly's regularization theorem, a pseudo-effective class carrying a positive current with vanishing Lelong number at every point is nef. Hence $\alpha=\{T\}$ is nef.

    In the general case, fix a Kähler form $\omega$ on the realization of $\mathbb{D}$. For every $\epsilon>0$, the Cartier b-divisor $\mathbb{D}+\epsilon\mathbb{D}(\omega)$ is big and nef in the sense of \cref{def:nefb}, hence is nef in the sense of \cref{def:nefCar} by the big case. Thus $\alpha+\epsilon\{\omega\}$ is nef for all $\epsilon>0$, and the closedness of the nef cone shows that $\alpha$ is nef.
\end{proof}

\begin{corollary}\label{cor:Dorder_Iorder}
    Let $T$ and $T'$ be non-divisorial closed positive $(1,1)$-currents on $X$. Suppose that $\{T\}=\{T'\}$. Then the following are equivalent:
    \begin{enumerate}
        \item\label{itm:cor:Dorder_Iorder:1} $\mathbb{D}(T)\leq \mathbb{D}(T')$;
        \item\label{itm:cor:Dorder_Iorder:2} $T\preceq_{\mathcal{I}}T'$.
    \end{enumerate}
\end{corollary}
\begin{proof}
    Let $E$ be a prime divisor over $X$. Choose a sequence of blow-ups with smooth connected centers as in the proof of \cref{thm:bdiv_current_corr}, so that $E$ appears as one of the exceptional prime divisors on the resulting model $\pi\colon Y\rightarrow X$. By \eqref{eq:RegpistarT}, the difference
    \[
        \mathbb{D}(T')_Y-\mathbb{D}(T)_Y
    \]
    is represented in the relative cohomology by
    \[
        \sum_i\left(\nu(T,E_i)-\nu(T',E_i)\right)\{E_i\}.
    \]
    Using \cref{prop:blowup_cone} inductively, or equivalently the negativity lemma for exceptional classes, pseudo-effectivity of this exceptional class is equivalent to the non-negativity of all coefficients. Hence $\mathbb{D}(T)\leq \mathbb{D}(T')$ is equivalent to
    \[
        \nu(T,E)\geq \nu(T',E)
    \]
    for every prime divisor $E$ over $X$, which is exactly $T\preceq_{\mathcal{I}}T'$ by \cref{prop:Iequivchar2}.
\end{proof}

\begin{corollary}\label{cor:nef_Cartierseq}
    Let $\mathbb{D}$ be a nef b-divisor over $X$. Then there is a decreasing sequence of nef and big Cartier b-divisors $\mathbb{D}_i$ over $X$ with limit $\mathbb{D}$.
\end{corollary}
\begin{proof}
    Take a Kähler form $\omega$ on $X$. By \cref{thm:bigandnefreal}, for each $i>0$, we can find a non-divisorial Kähler current $T_i\in \mathbb{D}_X+i^{-1}\{\omega\}$ such that
    \[
    \mathbb{D}(T_i)=\mathbb{D}+i^{-1}\mathbb{D}(\omega).
    \]
    We observe that $T_i$ and $T_{i+1}$ have the same $\mathcal{I}$-singularity type. Indeed, applying \cref{cor:Dorder_Iorder} to $T_i$ and $T_{i+1}+(i^{-1}-(i+1)^{-1})\omega$ shows that these two currents are $\mathcal{I}$-equivalent; adding the smooth form does not change the singularity type.
    Let $(T_i^j)_j$ be quasi-equisingular approximations of $T_i$ such that
    \begin{enumerate}
        \item\label{itm:chapter-bdiv:L544:1}
        $T_i^j$ is a Kähler current in $\mathbb{D}_X+i^{-1}\{\omega\}$ for $j\geq j_0(i)$, and
        \item\label{itm:chapter-bdiv:L544:2} for each fixed $j$, the singularity type of $T_i^j$ is independent of $i$.
    \end{enumerate}
    Note that \ref{itm:chapter-bdiv:L544:2} is possible by using the Bergman kernel construction of the quasi-equisingular approximations, since the $T_i$'s have the same multiplier ideals.

    For fixed $j$, the b-divisor
    \[
        \mathbb{B}_j\coloneqq \mathbb{D}(T_i^j)-i^{-1}\mathbb{D}(\omega)
    \]
    is independent of $i$. Since the quasi-equisingular approximations are decreasing in $j$, the sequence $(\mathbb{B}_j)_j$ is decreasing, and \cref{prop:dscontbdiv} gives $\mathbb{B}_j\to \mathbb{D}$. It suffices to take $\mathbb{D}_i=\mathbb{B}_{j_i}+i^{-1}\mathbb{D}(\omega)=\mathbb{D}(T_i^{j_i})$, where $j_i$ is chosen to be a sufficiently fast increasing sequence of positive integers with $j_i\geq j_0(i)$.
\end{proof}

\section{Properties of the intersection product}\label{sec:int_b_div}
Let $X$ be a connected compact Kähler manifold of dimension $n$. Recall that the mixed volume of currents is studied in \cref{sec:mix_vol}.

\begin{definition}\label{def:DF_ext2}
    Let $\mathbb{D}_1,\ldots,\mathbb{D}_n$ be big and nef b-divisors over $X$. Then we define their intersection as
    \index{b-divisor!intersection product}
    \[
    \left( \mathbb{D}_1,\ldots,\mathbb{D}_n \right)\coloneqq \vol (T_1,\ldots,T_n),
    \]
    \index{$\left( \mathbb{D}_1,\ldots,\mathbb{D}_n \right)$}%
    where $T_1,\ldots,T_n$ are closed positive $(1,1)$-currents in $\mathbb{D}_{1,X},\ldots,\mathbb{D}_{n,X}$ respectively such that $\mathbb{D}(T_i)=\mathbb{D}_{i}$.
\end{definition}
The definition makes sense thanks to \cref{thm:bigandnefreal}. It does not depend on the choices of $T_1,\ldots,T_n$ since they are uniquely defined up to $\mathcal{I}$-equivalence, as proved in \cref{thm:bdiv_current_corr}.

\begin{lemma}\label{lma:DFlinear}\leavevmode
    \begin{enumerate}
        \item\label{itm:lma:DFlinear:1} The product in \cref{def:DF_ext2} is symmetric in its $n$ variables.
        \item\label{itm:lma:DFlinear:2} Let $\mathbb{D}_1,\ldots,\mathbb{D}_n,\mathbb{D}'_1$ be big and nef b-divisors over $X$. Then
        \[
        \left(\mathbb{D}_1+\mathbb{D}_1',\ldots,\mathbb{D}_n \right)=\left(\mathbb{D}_1,\ldots,\mathbb{D}_n \right)+\left(\mathbb{D}_1',\ldots,\mathbb{D}_n \right).
        \]
        \item\label{itm:lma:DFlinear:3} Let $\mathbb{D}_1,\ldots,\mathbb{D}_n$ be big and nef b-divisors over $X$ and $\lambda\geq 0$. Then
        \[
        \left(\lambda \mathbb{D}_1,\ldots,\mathbb{D}_n \right)=\lambda \left(\mathbb{D}_1,\ldots,\mathbb{D}_n \right).
        \]
        \item\label{itm:lma:DFlinear:4} Let $\mathbb{D}_1,\ldots,\mathbb{D}_n$ and $\mathbb{D}'$ be big and nef b-divisors over $X$ such that $\mathbb{D}_1\leq \mathbb{D}'$. Then
        \[
            \left(\mathbb{D}_1,\ldots,\mathbb{D}_n  \right)\leq \left(\mathbb{D}',\mathbb{D}_2,\ldots,\mathbb{D}_n  \right).
        \]
    \end{enumerate}
\end{lemma}
\begin{proof}
    \ref{itm:lma:DFlinear:1}, \ref{itm:lma:DFlinear:2} and \ref{itm:lma:DFlinear:3} follow immediately from the corresponding results for mixed volumes, as proved in \cref{prop:vollinearlimit}.

    \ref{itm:lma:DFlinear:4} Take $\mathcal{I}$-good non-divisorial closed positive $(1,1)$-currents $T_1,\ldots,T_n$ and $T'$ so that $\mathbb{D}(T_i)=\mathbb{D}_i$ for all $i=1,\ldots,n$ and $\mathbb{D}(T')=\mathbb{D}'$.
    Furthermore, we may assume that the $T_i$'s and $T'$ are Kähler currents by the perturbation argument.

    Let $(T_{i}^j)_j$ be a quasi-equisingular approximation of $T_i$ for $i=2,\ldots,n$. It follows from \cref{thm:convdS} that
    \[
    \int_X T_1\wedge \cdots \wedge T_n=\lim_{j\to\infty} \int_X T_1\wedge T_2^j\wedge\cdots \wedge T_n^j.
    \]
    It suffices to show that for all $j\geq 1$,
    \[
    \int_X T_1\wedge T_2^j\wedge\cdots \wedge T_n^j\leq \int_X T'\wedge T_2^j\wedge\cdots \wedge T_n^j.
    \]
    Therefore, we have reduced to the case where $T_2,\ldots,T_n$ have analytic singularities. After a resolution, we may assume that they have log singularities along $\mathbb{Q}$-divisors. By \itemref{itm:prop:vollinearlimit:5}, we can further reduce to the case where $T_2,\ldots,T_n$ have bounded local potentials. Perturbing $T_2,\ldots,T_n$ by a Kähler form, we may further assume that $\{T_2\},\ldots, \{T_n\}$ are Kähler classes. By \itemref{itm:prop:vollinearlimit:4}, we finally reduce to the case where $T_2,\ldots,T_n$ are Kähler forms. In this case, the assertion follows from the pseudo-effectivity of $\{T'\}-\{T_1\}$.
\end{proof}

In general, we can extend the definition of the intersection product to nef b-divisors as follows:
\begin{definition}\label{def:DF_ext}
    Let $\mathbb{D}_1,\ldots,\mathbb{D}_n$ be nef b-divisors over $X$. Then we define their intersection as
    \index{b-divisor!intersection product}
    \begin{equation}\label{eq:mixed_product_limit_def}
        \left( \mathbb{D}_1,\ldots,\mathbb{D}_n \right)\coloneqq \lim_{\epsilon \to 0+}\Bigl( \mathbb{D}_1+\epsilon \mathbb{D}(\omega),\ldots,\mathbb{D}_n+\epsilon \mathbb{D}(\omega) \Bigr),
    \end{equation}
    \index{$\left( \mathbb{D}_1,\ldots,\mathbb{D}_n \right)$}%
    where $\omega$ is a Kähler form on $X$.

    Note that as $\epsilon\to 0+$, $\left( \mathbb{D}_1+\epsilon \mathbb{D}(\omega),\ldots,\mathbb{D}_n+\epsilon \mathbb{D}(\omega) \right)$ is decreasing, as a consequence of \itemref{itm:lma:DFlinear:4}. Similarly, by \itemref{itm:lma:DFlinear:4} again, the definition is independent of the choice of $\omega$.
\end{definition}

\begin{proposition}\label{prop:DFlinear}\leavevmode
    \begin{enumerate}
        \item\label{itm:prop:DFlinear:1} The product in \cref{def:DF_ext} is symmetric in its $n$ variables.
        \item\label{itm:prop:DFlinear:2} Let $\mathbb{D}_1,\ldots,\mathbb{D}_n,\mathbb{D}'_1$ be nef b-divisors over $X$. Then
        \[
        \left(\mathbb{D}_1+\mathbb{D}_1',\ldots,\mathbb{D}_n \right)=\left(\mathbb{D}_1,\ldots,\mathbb{D}_n \right)+\left(\mathbb{D}_1',\ldots,\mathbb{D}_n \right).
        \]
        \item\label{itm:prop:DFlinear:3} Let $\mathbb{D}_1,\ldots,\mathbb{D}_n$ be nef b-divisors over $X$ and $\lambda\geq 0$. Then
        \[
        \left(\lambda \mathbb{D}_1,\ldots,\mathbb{D}_n \right)=\lambda \left(\mathbb{D}_1,\ldots,\mathbb{D}_n \right).
        \]
        \item\label{itm:prop:DFlinear:4} Let $\mathbb{D}_1,\ldots,\mathbb{D}_n$ and $\mathbb{D}'$ be nef b-divisors over $X$ such that $\mathbb{D}_1\leq \mathbb{D}'$. Then
        \[
            \left(\mathbb{D}_1,\ldots,\mathbb{D}_n  \right)\leq \left(\mathbb{D}',\mathbb{D}_2,\ldots,\mathbb{D}_n  \right).
        \]
        \item\label{itm:prop:DFlinear:5} Let $T_1,\ldots,T_n$ be closed positive $(1,1)$-currents on $X$. Then
        \[
            \Bigl( \mathbb{D}(T_1),\ldots,\mathbb{D}(T_n) \Bigr)= \vol (T_1,\ldots,T_n).
        \]
    \end{enumerate}
\end{proposition}
\begin{proof}
    \ref{itm:prop:DFlinear:1}, \ref{itm:prop:DFlinear:2}, \ref{itm:prop:DFlinear:3} and \ref{itm:prop:DFlinear:4} These are immediate consequences of the corresponding results in the big and nef case established in \cref{lma:DFlinear}.

    \ref{itm:prop:DFlinear:5} By \eqref{eq:volmixedgeneral} and \eqref{eq:mixed_product_limit_def}, we may assume that $T_1,\ldots,T_n$ are Kähler currents. In this case, the assertion follows from the definitions.
\end{proof}

\begin{lemma}\label{lma:perturbDF}
    Let $\omega$ be a Kähler form on $X$. Fix a compact set $K\subseteq \mathrm{H}^{1,1}(X,\mathbb{R})$.
    Let $\mathbb{D}_1,\ldots,\mathbb{D}_n$ be nef b-divisors over $X$ such that $\mathbb{D}_{i,X}\in K$ for each $i=1,\ldots,n$. Then there is a constant $C$ depending only on $X,K,\{\omega\}$ such that for any $\epsilon\in [0,1]$, we have
    \[
    0\leq \Bigl( \mathbb{D}_1+\epsilon \mathbb{D}(\omega),\ldots,\mathbb{D}_n+\epsilon \mathbb{D}(\omega) \Bigr)-\left( \mathbb{D}_1,\ldots,\mathbb{D}_n \right)\leq C\epsilon.
    \]
\end{lemma}
\begin{proof}
    By \itemref{itm:lma:DFlinear:4}, the left-hand side is non-negative. For the upper bound, by the multi-linearity of the product proved in \cref{prop:DFlinear}, we can rewrite the difference as a sum of $2^n-1$ terms like
    \[
        \epsilon^{m} \Bigl( \mathbb{D}(\omega)^{m}, \mathbb{D}_{i_1},\ldots, \mathbb{D}_{i_{n-m}} \Bigr),
    \]
    where $J\subseteq\{1,\ldots,n\}$ is a non-empty subset of cardinality $m$, and $\{i_1,\ldots,i_{n-m}\}$ is the complement of $J$. We need to show that these coefficients are uniformly bounded by a constant depending only on the given data. For this purpose, we add $\mathbb{D}(\omega)$ to each $\mathbb{D}_i$ and then replace $K$ by $K+\{\omega\}$. Therefore, we reduce to the case where each $\mathbb{D}_i$ is nef and big, and hence $\mathbb{D}_i=\mathbb{D}(T_i)$ for some $\mathcal{I}$-good closed positive $(1,1)$-currents on $X$ with positive volumes. Then we need to bound
    \[
        \int_X \omega^m\wedge T_{i_1}\wedge \cdots \wedge T_{i_{n-m}}.
    \]
    For this purpose, let $D\geq 0$ be the minimal non-negative number so that $D\{\omega\}-\alpha$ is nef for any $\alpha\in K$. Then the integral is bounded by the cohomological intersection $\{\omega\}^m\cap (D\{\omega\})^{n-m}$. Our assertion follows.
\end{proof}

We first make a consistency check.
\begin{proposition}\label{prop:bdiv:L678}
    Suppose that $\mathbb{D}$ is a nef b-divisor over $X$. Then
    \[
    \left( \mathbb{D},\ldots,\mathbb{D} \right)=\vol \mathbb{D}.
    \]
\end{proposition}
\begin{proof}
    Using \cref{lma:perturbDF} and \cref{lma:volbdivnet}, we may easily reduce to the case where $\mathbb{D}$ is nef and big. In this case, take a non-divisorial closed positive $(1,1)$-current $T$ in $\mathbb{D}_X$ such that $\mathbb{D}(T)=\mathbb{D}$. Then we need to show that
    \[
    \vol \mathbb{D}=\vol T,
    \]
    which is proved in \cref{thm:DT}.
\end{proof}

\begin{proposition}\label{prop:bdiv:L692}
    Let $\mathbb{D}_1,\ldots,\mathbb{D}_n$ be nef b-divisors over $X$. Then
    \[
    \left(\mathbb{D}_1,\ldots,\mathbb{D}_n \right)\geq \prod_{i=1}^n \left( \vol \mathbb{D}_i\right)^{1/n}.
    \]
\end{proposition}
\begin{proof}
    We may assume that $\vol \mathbb{D}_i>0$ for each $i=1,\ldots,n$ since there is nothing to prove otherwise. In this case, our assertion follows from \cref{prop:logvol}.
\end{proof}

\begin{proposition}\label{prop:DFusc}
    The product in \cref{def:DF_ext} is upper semicontinuous in the following sense. Suppose that $(\mathbb{D}_i^j)_{j\in J}$ are nets of nef b-divisors over $X$ with limits $\mathbb{D}_i$ for each $i=1,\ldots,n$. Then
    \begin{equation}\label{eq:limsupint_prod}
    \varlimsup_{j\in J}\left(\mathbb{D}_1^j,\ldots,\mathbb{D}_n^j \right)\leq  \left( \mathbb{D}_1,\ldots,\mathbb{D}_n \right).
    \end{equation}
\end{proposition}
\begin{proof}
    \textbf{Step~1}. We first assume that the $\mathbb D_i^j$'s and the $\mathbb D_i$'s are all big.

    Choose $\mathcal I$-good non-divisorial currents $T_i^j$ and $T_i$ such that $\mathbb D(T_i^j)=\mathbb D_i^j$ and $\mathbb D(T_i)=\mathbb D_i$. Fix a Kähler form $\omega$ on $X$ and fix $\epsilon>0$. Since $\mathbb D_{i,X}^j\to\mathbb D_{i,X}$, for all sufficiently large $j$ the class
    \[
        \mathbb D_{i,X}+\epsilon\{\omega\}-\mathbb D_{i,X}^j
    \]
    is Kähler. Choose a Kähler form $\chi_i^j$ in this class and set $S_i^j=T_i^j+\chi_i^j$. Then $S_i^j$ lies in the fixed class $\mathbb D_{i,X}+\epsilon\{\omega\}$ and has the same divisorial Lelong numbers as $T_i^j$. Hence, by \cref{thm:usc_mix_vol} and \itemref{itm:prop:vollinearlimit:2},
    \[
        \varlimsup_{j\in J}\vol\Bigl(T_1^j,\ldots,T_n^j\Bigr)\leq \varlimsup_{j\in J}\vol\Bigl(S_1^j,\ldots,S_n^j\Bigr)\leq\vol\left(T_1+\epsilon\omega,\ldots,T_n+\epsilon\omega\right).
    \]
    Therefore,
    \[
        \varlimsup_{j\in J}\Bigl(\mathbb D_1^j,\ldots,\mathbb D_n^j\Bigr)\leq \Bigl(\mathbb D_1+\epsilon\mathbb D(\omega),\ldots,\mathbb D_n+\epsilon\mathbb D(\omega)\Bigr).
    \]
    Letting $\epsilon\to 0+$ gives \eqref{eq:limsupint_prod}.

    \textbf{Step~2}. Next we handle the general case.

    Take a Kähler form $\omega$ on $X$. Then by \cref{lma:perturbDF}, for any $\epsilon\in (0,1]$, we have
    \[
    \begin{aligned}
          \varlimsup_{j\in J}\left(\mathbb{D}_1^j,\ldots,\mathbb{D}_n^j \right)
        \leq & \varlimsup_{j\in J}\left(\mathbb{D}_1^j+\epsilon \mathbb{D}(\omega),\ldots,\mathbb{D}_n^j+\epsilon \mathbb{D}(\omega) \right)\\
        \leq & \left( \mathbb{D}_1+\epsilon \mathbb{D}(\omega),\ldots,\mathbb{D}_n+\epsilon \mathbb{D}(\omega) \right)\\
        \leq & \left( \mathbb{D}_1,\ldots,\mathbb{D}_n \right)+C\epsilon.
    \end{aligned}
    \]
    Since $\epsilon$ is arbitrary, our assertion follows.
\end{proof}

\begin{proposition}\label{prop:DFdecnet}
    The product in \cref{def:DF_ext} is continuous along decreasing nets in each variable. In other words, if $(\mathbb{D}_i^j)_{j\in J}$ ($i=1,\ldots,n$) are decreasing nets of nef b-divisors over $X$ with limits $\mathbb{D}_i$, then
    \[
    \lim_{j\in J} \left(\mathbb{D}_1^j,\ldots,\mathbb{D}_n^j\right)=\left(\mathbb{D}_1,\ldots,\mathbb{D}_n\right).
    \]
\end{proposition}
\begin{proof}
    This is a straightforward consequence of \itemref{itm:lma:DFlinear:4} and \cref{prop:DFusc}.
\end{proof}

\begin{remark}\label{rmk:bdiv:L749}
The algebraic nef b-divisors are introduced by Dang--Favre \cite{DF20}, as their intersection theory. It is a straightforward application of the results proved in this section that our transcendental intersection theory coincides with theirs in the algebraic setting.
\end{remark}

\begin{partbacktext}
\part{Applications}
\noindent
The three chapters of this part apply the singularity theory developed earlier to three rather different settings: big toric geometry, non-Archimedean pluripotential theory and partial Bergman kernels. They are largely independent of one another, and each uses a different part of the machinery built in the preceding chapters.

\cref{chap:toric_big} extends the toric dictionary of \cref{chap:toric_ample} from ample line bundles to big toric line bundles. In the big case, identifying the Newton body is no longer as simple as in the ample setting, and the behavior of potentials along toric strata becomes more delicate. The chapter establishes the toric-convex dictionary using the theory of partial Okounkov bodies developed in \cref{chap:Okounkov}.

\cref{chap:NApotential} develops a non-Archimedean pluripotential theory on compact K\"ahler manifolds. The basic objects are inverse systems of $\mathcal I$-model test curves, and the operations on test curves from \cref{sec:operationtc} provide the corresponding algebra of non-Archimedean potentials. The chapter constructs Duistermaat--Heckman measures and proves a comparison theorem identifying the resulting space with the Boucksom--Jonsson space of non-Archimedean psh metrics in the projective setting.

\cref{chap:Bergman} proves the convergence of partial Bergman kernels. Given a singular $\theta$-psh potential, a closed non-pluripolar weighted subset $(K,v)$ and a Bernstein--Markov measure, the normalized partial Bergman measure converges weakly to the corresponding partial equilibrium measure. The proof combines the envelope theory of \cref{chap:envelopes} and the asymptotic Riemann--Roch theorem \cref{thm:DXmain1}. It extends the theorem of Berman--Boucksom--Witt Nystr\"om \cite{BBWN11} to the singular relative setting considered here.

The prerequisites for the three chapters are different. The big toric chapter uses basic toric geometry, for which \cite{CLS11} and \cite{Fultor} are standard references. The non-Archimedean chapter is best read with the Boucksom--Jonsson theory in mind, especially \cite{BJ22}. The Bergman kernel chapter is closest in spirit to \cite{BBWN11}. These references are not logically required for following the arguments, but they provide useful context.
\end{partbacktext}

\motto{C'est l'harmonie des diverses parties, leur symétrie, leur heureux balancement; c'est en un mot tout ce qui y met de l'ordre, tout ce qui leur donne de l'unité, ce qui nous permet par conséquent d'y voir clair et d'en comprendre l'ensemble en même temps que les détails.

--- Henri Poincaré\footnote{Henri Poincaré (1854--1912) was a French mathematician, physicist, and philosopher of science. He is considered one of the greatest mathematicians of all time and a pioneer of several modern mathematical fields. He also played a key role in the development of special relativity, and was one of the first to understand the deep connection between mathematics and physics.}, \emph{L'avenir des mathématiques}}
\chapter{Toric pluripotential theory on big line bundles}
\label{chap:toric_big}

This chapter extends the toric pluripotential theory of \cref{chap:toric_ample} from ample toric line bundles to big toric line bundles. Let $L=\mathcal O_X(D)$ be a big toric line bundle on a smooth projective toric variety $X$, and let $P_D\subseteq M_{\mathbb R}$ be the polytope associated with the toric divisor $D$. In the ample case, $P_D$ is a smooth lattice polytope and its faces correspond directly to toric strata. This direct correspondence breaks down. This is the main new feature of the chapter.

The first goal is to recover a usable toric-convex dictionary in this more flexible setting. We begin by recalling the big toric setup, the Cartier data of $D$, and the stable base locus. We then compare the Newton body of a toric psh metric with its partial Okounkov body. The main computation, \cref{thm:toricpob}, identifies the partial Okounkov body of a toric singularity with its Newton body, up to the natural affine transformation determined by the chosen toric flag. This result is the bridge between the general Okounkov-body theory of \cref{chap:Okounkov} and the explicit convex geometry of toric metrics.

Once this comparison is established, several pluripotential-theoretic consequences follow. The Newton body computes non-pluripolar mass, describes Lelong numbers along toric divisors, and gives an explicit description of the potential with minimal singularities. In this way the abstract constructions of previous chapters become concrete convex-geometric statements in the toric setting.

The second half of the chapter develops the full toric pluripotential theory in a big class. Toric $\theta$-psh potentials are identified with proper convex functions on $N_{\mathbb R}$ satisfying the appropriate growth condition determined by $P_D$. Under this correspondence, sums, rooftop envelopes, suprema, decreasing limits and $P$-equivalence all have convex interpretations. We also describe toric subgeodesics, geodesics and the toric Ross--Witt Nystr\"om correspondence.

The final sections study restriction and intersection theory in the toric setting. The trace operator along toric strata is computed by slicing Newton bodies, and its volume formula becomes a statement about the Euclidean volume of faces and slices. We also prove a mixed-volume formula for Newton bodies and recover the augmented base locus from the convex data.

Throughout the chapter we use the notation and conventions of \cref{chap:toric_ample}. The partial Okounkov body machinery of \cref{chap:Okounkov} is used repeatedly.

\section{Toric setup}\label{sec:toricsetup2}
Let $T$ be a complex torus of dimension $n$ with character lattice $M$ and cocharacter lattice $N$. Some basic terminology is recalled in \cref{sec:toricsetup}.
Recall that $T_c$ is the compact torus contained in $T(\mathbb{C})$.

Consider a rational polyhedral fan $\Sigma$ in $N_{\mathbb{R}}$ corresponding to an $n$-dimensional smooth projective toric variety $X$.

Let
\[
D=\sum_{\rho\in \Sigma(1)} a_{\rho}D_{\rho}
\]
be a $T$-invariant integral big divisor on $X$. Let $P_D\subseteq M_{\mathbb{R}}$ be the following polytope\footnote{Note that $P_D$ is not necessarily a lattice polytope, see \cite[Example~10.5.4]{CLS11}. In fact, any rational polytope with positive volume can be realized in this manner.}
\begin{equation}\label{eq:bigPrep}
P_D=\left\{ m\in M_{\mathbb{R}}: \langle m,u_{\rho}\rangle \geq -a_{\rho}\quad \forall \rho\in \Sigma(1)\right\}.
\end{equation}
Since we have assumed that $D$ is big, $P_D$ is $n$-dimensional.

Let $L=\mathcal{O}_X(D)$.
Note that replacing $D$ by a linearly equivalent divisor amounts to replacing $P_D$ by an integral translation.

Recall that for each $\rho\in \Sigma(1)$, $u_{\rho}$ denotes the ray generator of $\rho$.
Let $\{m_{\sigma}\}_{\sigma\in \Sigma}$ denote the \emph{Cartier data} associated with $D$. In other words, for each $\sigma\in \Sigma$, $m_{\sigma}\in M$ satisfies that
\[
\langle m_{\sigma},u_{\rho} \rangle =-a_{\rho},\quad \forall \rho\in \sigma(1).
\]
The element $m_{\sigma}\in M$ is well-defined modulo
\begin{equation}\label{eq:Msigmadef}
    M(\sigma)\coloneqq \sigma^{\perp}\cap M,
\end{equation}
where
\[
\sigma^{\perp}\coloneqq \left\{ m\in M_{\mathbb{R}} : \langle m,u\rangle=0 \quad\forall u\in \sigma \right\}.
\]
Moreover, if $\tau$ is a face of $\sigma$, then
\begin{equation}\label{eq:Cartierdatalowdim}
m_{\sigma}\equiv m_{\tau} \mod M(\tau).
\end{equation}
See \cite[Theorem~4.2.8]{CLS11}. In particular, for an $n$-dimensional $\sigma\in \Sigma$, the element $m_{\sigma}$ is uniquely determined. In general, for $\sigma\in \Sigma(n)$, one need not have $m_{\sigma}\in P_D$. In fact, $m_{\sigma}\in P_D$ for all $\sigma\in \Sigma(n)$ if and only if $D$ is base-point free. See \cite[Theorem~6.1.7]{CLS11}.

Note that for any $n$-dimensional cone $\sigma$ in $\Sigma$ and any $\rho\in \sigma(1)$, we have
\begin{equation}\label{eq:mmsigpos}
    \langle m-m_{\sigma},u_{\rho} \rangle \geq 0,\quad \forall m\in P_D,
\end{equation}
as a consequence of \eqref{eq:bigPrep} and the choice of $m_\sigma$.

Recall that
\begin{equation}\label{eq:DUsigma}
D|_{U_{\sigma}}=\Div(\chi^{-m_{\sigma}})|_{U_{\sigma}}
\end{equation}
for all $\sigma\in \Sigma$, where $U_{\sigma}$ is the affine subvariety of $X$ corresponding to $\sigma$. See \cite[Proposition~4.1.2]{CLS11}.

Next consider a $T$-invariant irreducible subvariety $Y\subseteq X$. Since $X$ is smooth, so is $Y$. Let $\sigma$ be the cone in $\Sigma$ corresponding to $Y$.
We observe that $\sigma$ corresponds to a (possibly empty) face $Q_{\sigma}$ of $P_D$:
\begin{equation}\label{eq:Qsigma}
Q_{\sigma}=\left\{m\in P_D: \langle m,u_{\rho}\rangle=-a_{\rho}\quad \forall \rho\in \sigma(1) \right\}.
\end{equation}
The dimension of $\sigma$ is not necessarily equal to the codimension of $Q$ as we will see in \cref{ex:P2bl1pt}.

We have the following characterization of the base locus of $D$. We are not aware of a direct reference; the result is likely folklore.
\begin{proposition}\label{prop:baseloc}
    The base locus $\Bs(D)$ (with the reduced complex structure) of $D$ is a toric-invariant (possibly reducible) subvariety given by the union of $V(\tau)$, where $\tau$ runs over elements in $\Sigma$ satisfying the following condition:
    \begin{equation}\label{eq:arhononcrit}
        a_{\rho}+\langle m, u_{\rho}\rangle>0\quad \textup{for some }\rho\in \tau(1)
    \end{equation}
    for each $m\in M\cap P_D$.
\end{proposition}
Here $V(\tau)$ denotes the toric subvariety of $X$ corresponding to $\tau$. In other words,
\[
V(\tau)=\bigcap_{\rho\in \tau(1)} D_{\rho}.
\]
\begin{proof}
    Recall that
    \[
    \mathrm{H}^0(X,L)\cong \bigoplus_{m\in M\cap P_D} \mathbb{C}\chi^m.
    \]
    See \cite[Proposition~4.3.3]{CLS11} for example. So we only need to understand the common zeros of $D+\Div \chi^m$ for all $m\in M\cap P_D$. But we know that
    \[
    D+\Div \chi^m=\sum_{\rho\in \Sigma(1)} \left( a_{\rho}+\langle m, u_{\rho}\rangle \right) D_{\rho}.
    \]
    Our assertion follows.
\end{proof}

\begin{corollary}\label{cor:stablebl_toric}
    The stable base locus (with the reduced complex structure) of $D$ is a toric-invariant (possibly reducible) subvariety given by the union of the $V(\tau)$'s,  where $\tau$ runs over elements in $\Sigma$ satisfying \eqref{eq:arhononcrit} for each $m\in P_D$.
\end{corollary}
Geometrically, the condition means $Q_{\tau}=\varnothing$.
\begin{proof}
    It follows from \cref{prop:baseloc} that the stable base locus of $D$ is given by the union of the $V(\tau)$'s, where $\tau$ runs over elements of $\Sigma$ satisfying \eqref{eq:arhononcrit} for all $m\in M_{\mathbb{Q}}\cap P_D$. That is,
    \[
        M_{\mathbb{Q}}\cap \left\{ m\in P_D : a_{\rho}+\langle m, u_{\rho}\rangle=0\quad \forall \rho\in \tau(1) \right\}=\varnothing.
    \]
    Since the latter part is a rational polytope, this statement is equivalent to
    \begin{equation}\label{eq:PDcapvarnot}
    \left\{ m\in P_D : a_{\rho}+\langle m, u_{\rho}\rangle=0\quad \forall \rho\in \tau(1) \right\}=\varnothing.
    \end{equation}
    Our first assertion follows.
\end{proof}

We will keep two examples in mind.
\begin{example}\label{ex:P2}
    In this case, $\Sigma$ is the fan in \cref{fig:torP2} consisting of three $2$-dimensional cones $\sigma_1$, $\sigma_2$ and $\sigma_3$; three $1$-dimensional cones $\sigma_4$, $\sigma_5$ and $\sigma_6$; one $0$-dimensional cone $\sigma_0$.

    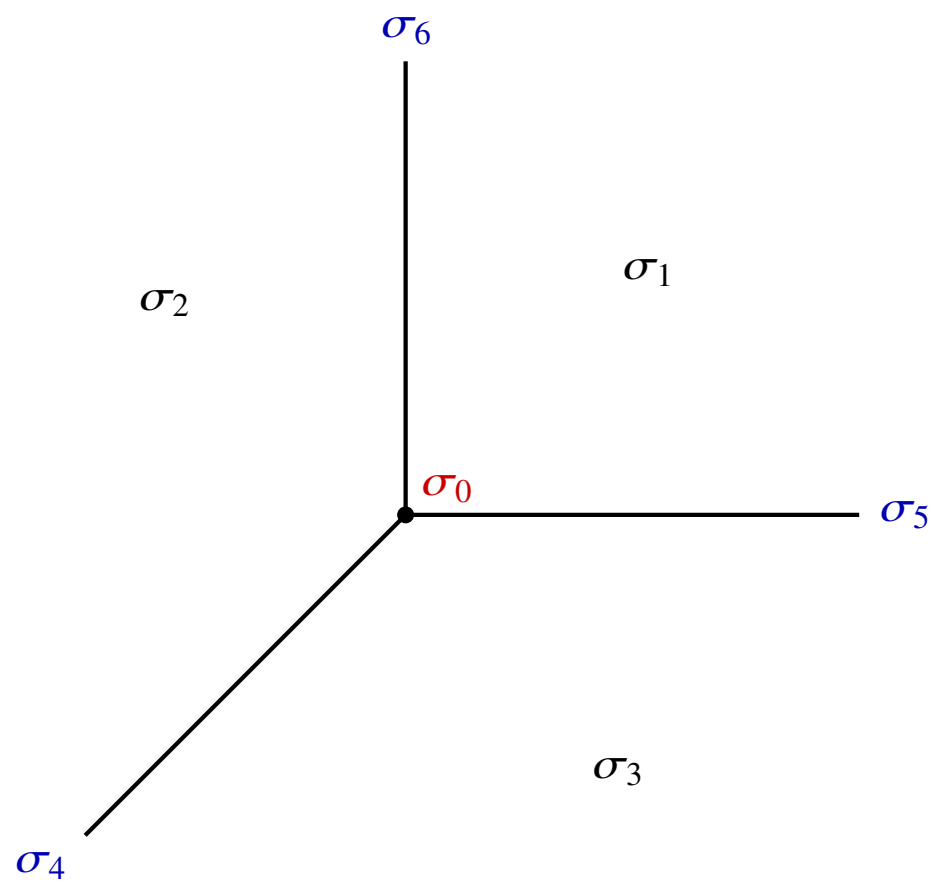
\begin{figure}[ht]
\centering
\begin{tikzpicture}[scale=1.0, line width=0.7pt, every node/.style={font=\small}]
  \draw[thick] (0,0) -- (3.0,0)   node[blue!70!black, right] {$\sigma_5$};
  \draw[thick] (0,0) -- (0,3.0)   node[blue!70!black, above] {$\sigma_6$};
  \draw[thick] (0,0) -- (-2.12,-2.12) node[blue!70!black, below left] {$\sigma_4$};
  \fill (0,0) circle (1.6pt);
  \node[red!80!black, anchor=south west, inner sep=2pt] at (0,0) {$\sigma_0$};
  \node at (1.6,1.6)   {$\sigma_1$};
  \node at (-1.6,1.4)  {$\sigma_2$};
  \node at (1.4,-1.7)  {$\sigma_3$};
\end{tikzpicture}
\caption{The fan $\Sigma$ of $\mathbb{P}^2$. The three rays $\sigma_4,\sigma_5,\sigma_6\in\Sigma(1)$ subdivide $N_{\mathbb{R}}$ into three full-dimensional cones $\sigma_1,\sigma_2,\sigma_3\in\Sigma(2)$; the origin is the unique cone of dimension zero.}
\label{fig:torP2}
\end{figure}

    The fan $\Sigma$ is just the fan of $X=\mathbb{P}^2$. Under the orbit-cone correspondence, we have
    \[
    \begin{aligned}
    D_{\sigma_1}=& \{[1:0:0]\},\quad D_{\sigma_2}=\{[0:1:0]\},\quad D_{\sigma_3}=\{[0:0:1]\},\\
    D_{\sigma_4}=& \{[0:X_1:X_2]:(X_1,X_2)\neq (0,0)\},\\
    D_{\sigma_5}=& \{[X_0:0:X_2]:(X_0,X_2)\neq (0,0)\},\\
    D_{\sigma_6}=& \{[X_0:X_1:0]:(X_0,X_1)\neq (0,0)\},\quad
    D_{\sigma_0}=\mathbb{P}^2.
    \end{aligned}
    \]
    In particular, $\Sigma(1)=\{\sigma_4,\sigma_5,\sigma_6\}$.
    We shall take
    \[
    D=D_{\sigma_4}.
    \]
    In other words,
    \[
    a_{\sigma_5}=a_{\sigma_6}=0,\quad a_{\sigma_4}=1.
    \]
    Note that the ray generators are given by
    \[
    u_{\sigma_4}=(-1,-1),\quad u_{\sigma_5}=(1,0),\quad u_{\sigma_6}=(0,1).
    \]
    It follows that
    \[
    P_D=\{m=(m_1,m_2)\in \mathbb{R}^2: m_1+m_2\leq 1, m_1\geq 0, m_2\geq 0\}.
    \]
    Therefore, $P_D$ is just the polytope in \cref{fig:polyP2}.
        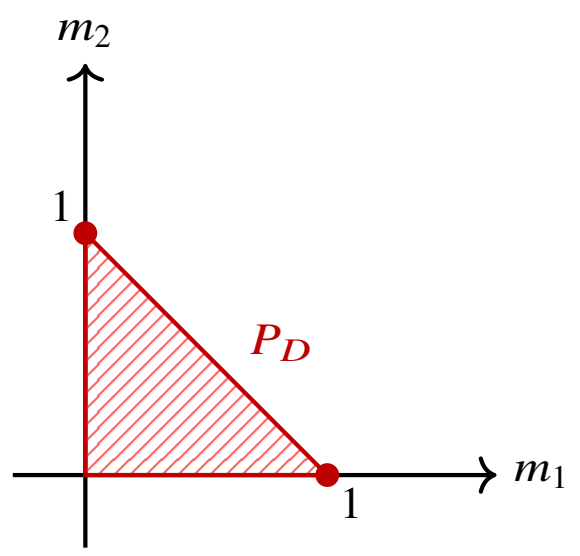
\begin{figure}[ht]
        \centering
        \begin{tikzpicture}[scale=1.6, line width=0.7pt, every node/.style={font=\small}]
          \draw[->, thick] (-0.3,0) -- (1.7,0) node[right] {$m_1$};
          \draw[->, thick] (0,-0.3) -- (0,1.7) node[above] {$m_2$};
          \fill[pattern=north east lines, pattern color=red!60] (0,0) -- (1,0) -- (0,1) -- cycle;
          \draw[red!75!black, thick] (0,0) -- (1,0) -- (0,1) -- cycle;
          \fill[red!75!black] (1,0) circle (1.4pt);
          \fill[red!75!black] (0,1) circle (1.4pt);
          \node[anchor=north west, inner sep=2pt] at (1,0) {$1$};
          \node[anchor=south east, inner sep=2pt] at (0,1) {$1$};
          \node[red!75!black, anchor=south west, inner sep=2pt] at (0.62,0.42) {$P_D$};
        \end{tikzpicture}
        \caption{The polytope $P_D$.}
        \label{fig:polyP2}
    \end{figure}
    In this case, the Cartier data for $2$-dimensional cones are given as follows:
    \[
    m_{\sigma_1}=(0,0),\quad m_{\sigma_2}=(1,0),\quad m_{\sigma_3}=(0,1);
    \]
    while the remaining Cartier data are determined by \eqref{eq:Cartierdatalowdim}.

    In this case, $L=\mathcal{O}_X(D)=\mathcal{O}_{\mathbb{P}^2}(1)$. Hence the line bundle $L$ is ample.

    We also observe that
    \[
    \begin{aligned}
    Q_{\sigma_1}=&\{(0,0)\},\quad Q_{\sigma_2}=\{(1,0)\},\quad Q_{\sigma_3}=\{(0,1)\},\\
    Q_{\sigma_4}=& \{(m_1,m_2):m_1\geq 0,m_2\geq 0,m_1+m_2=1\},\\
    Q_{\sigma_5}=& \{0\}\times [0,1],\quad Q_{\sigma_6}= [0,1]\times \{0\},\\
    Q_{\sigma_0}=& P_D.
    \end{aligned}
    \]
\end{example}

Next we give three non-ample examples.
\begin{example}\label{ex:P2bl1pt}
    Let $\Sigma$ be the fan shown in \cref{fig:polyP2b1p}.
    \begin{figure}[ht]
        \centering
        \begin{tikzpicture}[scale=1.0, line width=0.7pt, every node/.style={font=\small}]
          \draw[thick] (0,0) -- (3.0,0)   node[blue!70!black, right] {$\sigma_5$};
          \draw[thick] (0,0) -- (0,3.0)   node[blue!70!black, above] {$\sigma_6$};
          \draw[thick] (0,0) -- (-2.12,-2.12) node[blue!70!black, below left] {$\sigma_4$};
          \draw[thick] (0,0) -- (2.4,2.4) node[green!50!black, anchor=south west, inner sep=2pt, font=\small] {$\sigma_7$};
          \fill (0,0) circle (1.6pt);
          \node[red!80!black, anchor=north, inner sep=2pt] at (0,-0.30) {$\sigma_0$};
          \node[green!50!black] at (2.0,0.85) {$\sigma_1'$};
          \node[green!50!black] at (0.85,2.0) {$\sigma_1''$};
          \node at (-1.6,1.4) {$\sigma_2$};
          \node at (1.4,-1.7) {$\sigma_3$};
        \end{tikzpicture}
        \caption{The fan of $\Bl_0\mathbb{P}^2$. Adding the ray $\sigma_7$ generated by $(1,1)$ subdivides the cone $\sigma_1$ into two full-dimensional cones $\sigma_1',\sigma_1''$.}
        \label{fig:polyP2b1p}
    \end{figure}
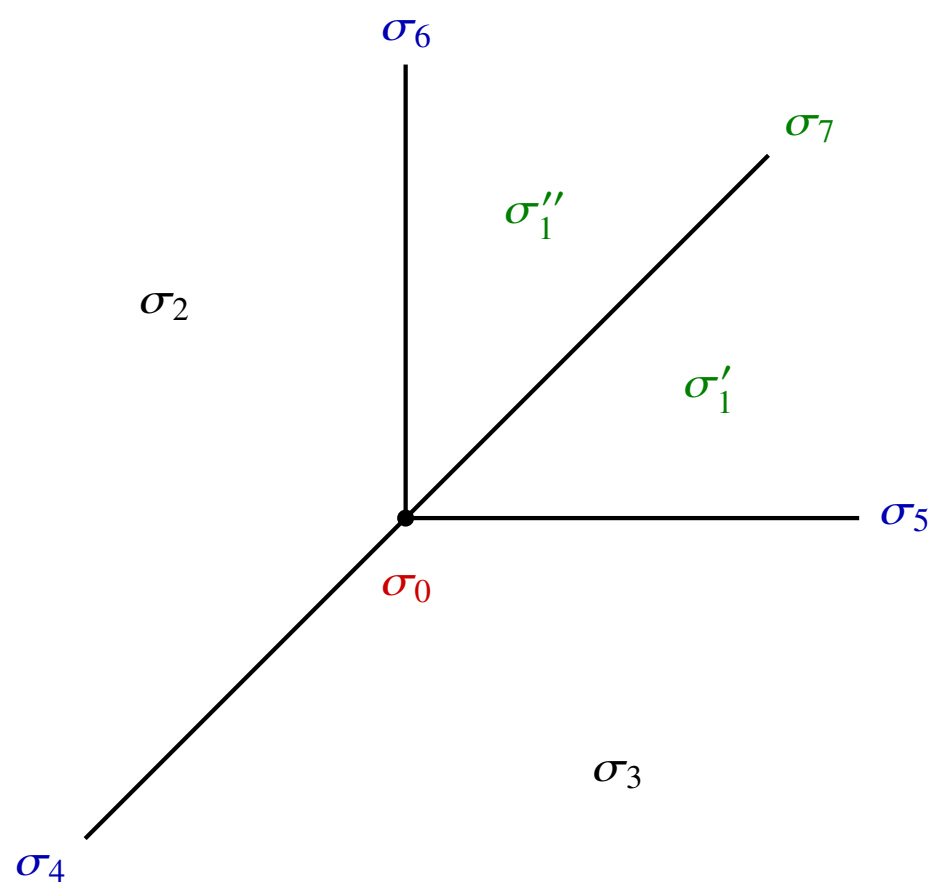
    Comparing with our previous example \cref{fig:torP2}, we have divided $\sigma_1$ from the middle, giving rise to two additional $2$-dimensional cones $\sigma_1'$ and $\sigma_1''$, and one additional $1$-dimensional cone $\sigma_7$.

    The corresponding $X=\Bl_0 \mathbb{P}^2$ is just the blow-up of $\mathbb{P}^2$ at the origin $0=[1:0:0]$ and hence $L=\pi^*\mathcal{O}_{\mathbb{P}^2}(1)$. Let $\pi\colon X\rightarrow \mathbb{P}^2$ denote the blow-up morphism.
    Let
    \[
    D=D_{\sigma_4}.
    \]
    Then $D$ is the pull-back of the divisor $D$ in \cref{ex:P2}. Note that $D$ is not ample, since it has degree $0$ on the exceptional divisor.

    In this case, we have
    \[
    \Sigma(1)=\{\sigma_4,\sigma_5,\sigma_6,\sigma_7\},
    \]
    and $D_{\sigma_7}$ is just the exceptional divisor.

    The corresponding ray generators are
    \[
    u_{\sigma_4}=(-1,-1),\quad u_{\sigma_5}=(1,0),\quad u_{\sigma_6}=(0,1),\quad u_{\sigma_7}=(1,1),
    \]
    while
    \[
    m_{\sigma_1'}=m_{\sigma_1''}=(0,0),\quad m_{\sigma_2}=(1,0),\quad m_{\sigma_3}=(0,1).
    \]
    Therefore, $P_D$ is the same as in \cref{fig:polyP2}.

    We also observe that
    \[
    \begin{aligned}
    Q_{\sigma_1'}=&\{(0,0)\},\quad Q_{\sigma_1''}=\{(0,0)\},\quad Q_{\sigma_7}=\{(0,0)\},\\
    Q_{\sigma_2}=&\{(1,0)\},\quad Q_{\sigma_3}=\{(0,1)\},\\
    Q_{\sigma_4}=& \{(m_1,m_2):m_1\geq 0,m_2\geq 0,m_1+m_2=1\},\\
    Q_{\sigma_5}=& \{0\}\times [0,1],\quad Q_{\sigma_6}= [0,1]\times \{0\},\\
    Q_{\sigma_0}=& P_D.
    \end{aligned}
    \]
\end{example}

\begin{example}\label{ex:F1big_non_nef}
    We continue with the surface in \cref{ex:P2bl1pt} with $X=\Bl_0\mathbb{P}^2$. Let $\pi\colon X\rightarrow \mathbb{P}^2$ be the blow-up morphism, and write
    \[
        H=D_{\sigma_4},\quad E=D_{\sigma_7}.
    \]
    Consider the toric Cartier divisor
    \[
        D=2D_{\sigma_4}+D_{\sigma_7}=2H+E.
    \]
    This gives a simple example of a big toric divisor which is not nef.

    Indeed, the intersection numbers on $X$ are
    \[
    H^2=1,\quad H\cdot E=0,\quad E^2=-1.
    \]
    Hence
    \[
    D\cdot E=(2H+E)\cdot E=-1<0,
    \]
    so $D$ is not nef.

    The inequalities defining $P_D$ are
    \[
    m_1\geq 0,\quad m_2\geq 0,\quad m_1+m_2\leq 2,\quad m_1+m_2\geq -1.
    \]
    The last inequality, corresponding to the exceptional ray $\sigma_7$, is redundant. This is a typical phenomenon in the non-ample situation.
    Thus
    \[
    P_D=\{m=(m_1,m_2)\in M_{\mathbb{R}}:m_1\geq 0,\ m_2\geq 0,\ m_1+m_2\leq 2\},
    \]
    as illustrated in \cref{fig:F1big_non_nef_polytope}.
    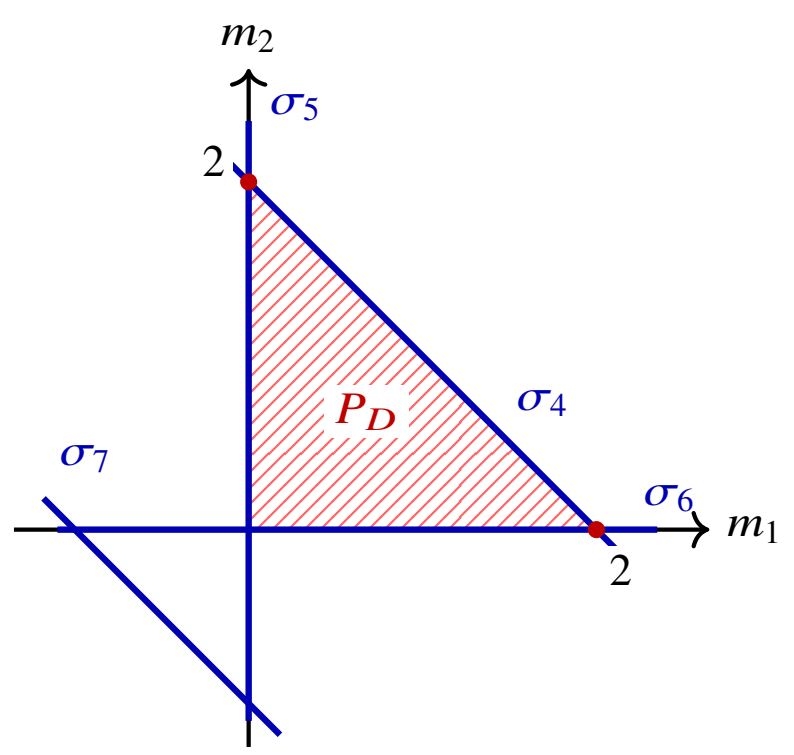
\begin{figure}[ht]
        \centering
        \begin{tikzpicture}[scale=1.15, line width=0.7pt, every node/.style={font=\small}]
          \draw[->, thick] (-1.35,0) -- (2.65,0) node[right] {$m_1$};
          \draw[->, thick] (0,-1.25) -- (0,2.65) node[above] {$m_2$};

          \fill[pattern=north east lines, pattern color=red!60] (0,0) -- (2,0) -- (0,2) -- cycle;
          \draw[red!75!black, thick] (0,0) -- (2,0) -- (0,2) -- cycle;
          \node[red!75!black, fill=white, inner sep=1.5pt] at (0.68,0.68) {$P_D$};

          \draw[blue!70!black, line width=1.2pt] (0,-1.10) -- (0,2.35);
          \draw[blue!70!black, line width=1.2pt] (-1.10,0) -- (2.35,0);
          \draw[blue!70!black, line width=1.2pt] (-0.18,2.18) -- (2.18,-0.18);
          \draw[blue!70!black, line width=1.2pt] (-1.18,0.18) -- (0.18,-1.18);

          \node[blue!70!black, anchor=south west, inner sep=2pt] at (0.03,2.28) {$\sigma_5$};
          \node[blue!70!black, anchor=south west, inner sep=2pt] at (2.18,0.03) {$\sigma_6$};
          \node[blue!70!black, anchor=south west, inner sep=2pt] at (1.45,0.58) {$\sigma_4$};
          \node[blue!70!black, anchor=south west, inner sep=2pt] at (-1.18,0.26) {$\sigma_7$};

          \fill[red!75!black] (2,0) circle (1.4pt);
          \fill[red!75!black] (0,2) circle (1.4pt);
          \node[fill=white, inner sep=1.5pt] at (2.14,-0.23) {$2$};
          \node[fill=white, inner sep=1.5pt] at (-0.20,2.12) {$2$};
        \end{tikzpicture}
        \caption{The polytope $P_D$ for $D=2D_{\sigma_4}+D_{\sigma_7}$. The blue lines are the four affine hyperplanes $a_{\rho}+\langle m,u_{\rho}\rangle=0$ for $\rho\in\Sigma(1)$.}
        \label{fig:F1big_non_nef_polytope}
    \end{figure}
\end{example}

\begin{example}\label{ex:F2_non_lattice_polytope}
    We give a surface example in which $P_D$ is not a lattice polytope. Let $\Sigma$ be the fan in $N_{\mathbb{R}}\simeq \mathbb{R}^2$ whose rays are generated, in cyclic order, by
    \[
    u_{\rho_1}=(1,0),\quad u_{\rho_2}=(0,1),\quad u_{\rho_3}=(-1,2),\quad u_{\rho_4}=(0,-1).
    \]
    The corresponding smooth projective toric surface is the Hirzebruch surface $\mathbb{F}_2$. Consider the toric Cartier divisor
    \[
    D=-D_{\rho_3}+2D_{\rho_4}.
    \]
    Thus
    \[
    a_{\rho_1}=a_{\rho_2}=0,\quad a_{\rho_3}=-1,\quad a_{\rho_4}=2.
    \]
    Hence
    \[
    P_D=\left\{m=(m_1,m_2)\in M_{\mathbb{R}}: m_1\geq 0,\ \frac{m_1+1}{2}\leq m_2\leq 2 \right\}.
    \]
    See \cref{fig:F2_non_lattice_polytope}.
    \begin{figure}[ht]
        \centering
        \begin{tikzpicture}[scale=1.05, line width=0.7pt, every node/.style={font=\small}]
          \draw[gray!25, thin, step=1] (0,0) grid (3,2);
          \draw[gray!50, dashed] (-0.18,0.5) -- (0,0.5);
          \draw[->, thick] (-0.35,0) -- (3.45,0) node[right] {$m_1$};
          \draw[->, thick] (0,-0.25) -- (0,2.45) node[above] {$m_2$};

          \fill[pattern=north east lines, pattern color=red!60] (0,0.5) -- (0,2) -- (3,2) -- cycle;
          \draw[red!75!black, thick] (0,0.5) -- (0,2) -- (3,2) -- cycle;
          \node[red!75!black, fill=white, inner sep=1.5pt] at (1.03,1.55) {$P_D$};

          \fill[red!75!black] (0,0.5) circle (1.4pt);
          \fill[red!75!black] (0,2) circle (1.4pt);
          \fill[red!75!black] (3,2) circle (1.4pt);
          \node[fill=white, inner sep=1.5pt, anchor=east] at (-0.02,0.5) {$\frac{1}{2}$};
          \node[fill=white, inner sep=1.5pt, anchor=east] at (-0.02,2) {$2$};
          \node[fill=white, inner sep=1.5pt, anchor=north] at (3,0) {$3$};
        \end{tikzpicture}
        \caption{A surface example where $P_D$ is not a lattice polytope.}
        \label{fig:F2_non_lattice_polytope}
    \end{figure}
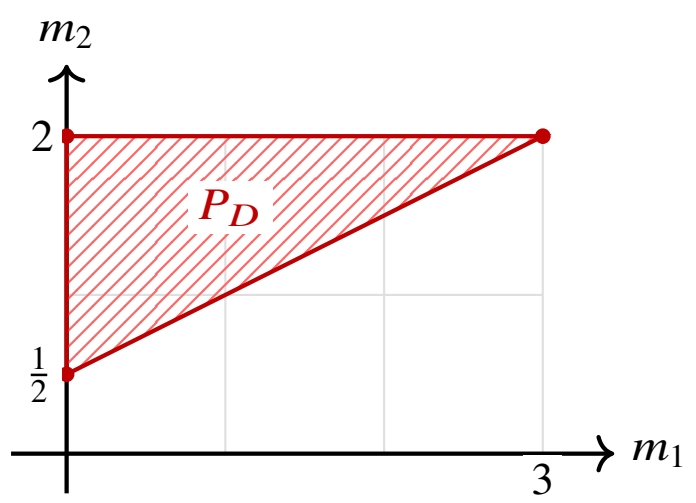
\end{example}

\section{Toric partial Okounkov bodies}\label{sec:toricpob}
We continue to use the notations in \cref{sec:toricsetup2}.

In order to study the toric-invariant singular plurisubharmonic metrics on $L$, we need to fix a reference toric-invariant smooth Hermitian metric, so that the psh metrics can be identified with quasi-psh functions.
Unlike the ample case studied in \cref{chap:toric_ample}, in the case of big line bundles, there does not seem to be a natural choice of a smooth Hermitian metric on $L$ similar to Guillemin's metric.

We shall fix a $T_c$-invariant Hermitian metric $h$ on $L$ so that $\theta=c_1(L,h)$.

The first observation is the following:
\begin{lemma}\label{lma:exist_Ftheta}
    There is a smooth function $F_{\theta}\colon N_{\mathbb{R}}\rightarrow \mathbb{R}$ such that
    \[
        \theta=\ddc \Trop^*F_{\theta}\quad\textup{on }T(\mathbb{C}).
    \]
\end{lemma}
\begin{proof}
    \textbf{Step~1}. We first prove the existence of $F_{\theta}$ for a specific choice of $h$.

    Write $D$ as the difference of two toric-invariant ample divisors, say $D_1-D_2$. Let $h_1$ and $h_2$ be the associated Guillemin's metrics on $D_1$ and $D_2$ respectively. We define $h=h_1\otimes h_2^{-1}$.

    In this case, our assertion follows since it holds in the case of Guillemin's metrics.

    \textbf{Step~2}. In general, fix $h_0$ as in Step~1. Then the general $h$ can be written as $h_0\exp(-g)$ for some smooth function $g$ on $X$. Since both metrics are $T_c$-invariant, so is $g$, and on the open torus we may write $g=\Trop^* r$ for some smooth function $r$ on $N_{\mathbb{R}}$.
    Hence it suffices to modify $F_{\theta}$ in Step~1 by $r$ to conclude.
\end{proof}
Note that $F_{\theta}$ is well-defined up to a linear term.

Next, we make an additional requirement on $F_{\theta}$ to fix the linear term. Let $s_D$ be a rational section of $L$ corresponding to $D$. Then $s_D$ is well-defined up to a non-zero multiple. By Lelong--Poincaré formula \cref{prop:LelongPoincare}, we have
\[
\ddc\left( \Trop^*F_{\theta}+\log |s_D|_h^2 \right)=0
\]
on $T(\mathbb{C})$. Therefore, $\Trop^*F_{\theta}+\log |s_D|_h^2$ is the tropicalization of a linear function. Hence, after adding a linear function to $F_{\theta}$, we can guarantee that
\begin{equation}\label{eq:normF0}
\Trop^*F_{\theta}+\log |s_D|_h^2=0
\end{equation}
from now on. Note that a different choice of $s_D$ means adding a constant to $F_{\theta}$.

Summarizing the situation, we have chosen the following data in addition to the fan $\Sigma$ so far: A divisor $D$, a Hermitian metric $h$ on $L$ and a rational section $s_D$ of $L$ subject to various conditions.

\subsection{Newton bodies}
\label{subsec:newton_big}
Let $\PSH_{\tor}(X,\theta)$ be the set of $T_c$-invariant functions in $\PSH(X,\theta)$.

\begin{definition}\label{def:toric-big:L416}
    A function $\varphi\in \PSH_{\tor}(X,\theta)$ can be written as
\[
\varphi|_{T(\mathbb{C})}=\Trop^*f
\]
for some unique function $f\colon N_{\mathbb{R}}\rightarrow [-\infty,\infty)$. Then we define $F_{\varphi}\colon N_{\mathbb{R}}\rightarrow \mathbb{R}$ as follows:
\begin{equation}\label{eq:Fvarphibig1}
F_{\varphi}=F_{\theta}+f.
\end{equation}
\index{$F_{\varphi}$}
\end{definition}
Observe that $F_{\varphi}$ is a convex function and takes finite values by \cref{lma:convextopsh}. In particular, $f$ is also real-valued.
Once $D$ and $h$ are fixed, $F_{\varphi}$ is well-defined up to a constant since $F_{\theta}$ is.

\begin{definition}\label{def:toric-big:L430}
    Let $\varphi\in \PSH_{\tor}(X,\theta)$, we define its \emph{Newton body}\index{Newton body} as
    \[
    \Delta(\theta,\varphi)\coloneqq \overline{\nabla F_{\varphi}(N_{\mathbb{R}})}\subseteq M_{\mathbb{R}}.
    \]
    \index{$\Delta(\theta,\varphi)$}
\end{definition}
Here the gradient is understood in the sense of Alexandrov, see \cref{def:subgrad} for the details.
Note that $\Delta(\theta,\varphi)$ is independent of the choice of $s_D$.

The Newton body $\Delta(\theta,\varphi)$ depends on the choice of $D$, not only on the associated line bundle $L$: A different choice of $D$ inducing the same line bundle corresponds to a translation of $\Delta(\theta,\varphi)$. We will see in a while (\cref{thm:toricpob}) that once $D$ is fixed $\Delta(\theta,\varphi)$ depends only on the current $\theta_{\varphi}$. Hence, the choice of $h$ is irrelevant.

For the moment, it is not clear if $\Delta(\theta,\varphi)\subseteq P_D$ or not. We shall prove this result using the theory of partial Okounkov bodies in the next section.

\begin{proposition}\label{prop:toricMAandrealMA2}
    Let $\varphi\in \PSH_{\tor}(X,\theta)$. Then
    \begin{equation}\label{eq:tropMAmea2}
    \Trop_* \left(\theta|_{T(\mathbb{C})}+\ddc\varphi|_{T(\mathbb{C})} \right)^n=\MA_{\mathbb{R}}\left(F_{\varphi}\right).
    \end{equation}
    In particular,
    \begin{equation}\label{eq:toricmass2}
        \int_X \theta_{\varphi}^n =n! \vol \Delta(\theta,\varphi)
    \end{equation}
\end{proposition}
\begin{proof}
    Let $F_0$ be a smooth convex function on $N_{\mathbb{R}}$ such that $\ddc \Trop^* F_0$ can be extended to a K\"ahler form on $X$. For example, Guillemin's construction \eqref{eq:G0def} with respect to a suitable Delzant polytope gives such an example.

    Write $\omega$ for a Kähler form on $X$ extending $\ddc\Trop^*F_0$. Then for any large enough $C>0$, $\theta+C\omega$ is a Kähler form. So we conclude from \cref{prop:toricMAandrealMA} that
    \[
    \Trop_* \left((\theta+C\omega)|_{T(\mathbb{C})}+\ddc\varphi|_{T(\mathbb{C})} \right)^n=\MA_{\mathbb{R}}(F_{\varphi}+CF_0).
    \]
    By the multi-linearity of the non-pluripolar product \itemref{itm:prop:npp1:6}, the left-hand side is a polynomial in $C$. By the expansion formula for the mixed real Monge--Amp\`ere measure \itemref{itm:thm:mixed_realMA:1}, the right-hand side is also a polynomial in $C$. We conclude that the same holds for $C=0$. Therefore, \eqref{eq:tropMAmea2} follows.

    \eqref{eq:toricmass2} is a direct consequence of \eqref{eq:tropMAmea2}.
\end{proof}

\subsection{Partial Okounkov bodies}\label{subsec:pobtorgeneral}

There are some canonical choices of smooth flags in the toric setting.

Since $X$ is smooth and projective, we can choose a full-dimensional cone $\sigma$ in $\Sigma$ with rays $\rho_1,\ldots,\rho_n\in \sigma(1)$ such that $u_{\rho_1},\ldots,u_{\rho_n}$ form a basis of $N$. Define
\[
Y_i=D_{\rho_1}\cap\cdots\cap D_{\rho_i},\quad i=1,\ldots,n.
\]
Then $Y_{\bullet}$ is a smooth flag on $X$. Let
\begin{equation}\label{eq:isoMZncanonical}
\Phi\colon M\rightarrow \mathbb{Z}^n,\quad m\mapsto \left(\langle m-m_{\sigma},u_{\rho_1} \rangle,\dots,\langle m-m_{\sigma},u_{\rho_n} \rangle \right).
\end{equation}
Then $\Phi$ is an isomorphism of lattices.
It induces an $\mathbb{Z}$-affine isomorphism
\[
\Phi_{\mathbb{R}}\colon M_{\mathbb{R}}\rightarrow \mathbb{R}^n.
\]

\begin{proposition}\label{prop:toricusualOko}
    We have
    \begin{equation}\label{eq:DeltakLtoric}
        k^{-1}\nu_{Y_{\bullet}}\left(\mathrm{H}^0(X,L^k)^{\times} \right)=\Phi_{\mathbb{R}} \left(P_{D}\cap k^{-1}M\right)
    \end{equation}
    for any $k\in \mathbb{Z}_{>0}$.
    In particular,
    \begin{equation}
        \Delta_{Y_{\bullet}}(L)=\Phi_{\mathbb{R}}(P_D).
    \end{equation}
\end{proposition}
Recall that $\Delta_{Y_{\bullet}}(L)\subseteq \mathbb{R}^n$ is the Okounkov body defined in \cref{def:Okobd}.
\begin{proof}
    We first reduce to the case where $D|_{U_{\sigma}}=0$. In fact, replacing $D$ by $D+\Div \chi^{m_{\sigma}}$ would result in changing $P_D$ to $P_D-m_{\sigma}$. So in view of \eqref{eq:DUsigma}, we may assume that $D|_{U_{\sigma}}=0$ and hence $m_{\sigma}=0$.

    Fix $k\in \mathbb{Z}_{>0}$. Let $s\in \mathrm{H}^0(X,L^k)$ be a non-zero toric eigen-section, say $\chi^m$ for some $m\in kP_D\cap M$. The zero-divisor of $s$ on $U_{\sigma}$ is given by
    \[
        \sum_{i=1}^n \langle m, u_{\rho_i} \rangle D_{\rho_i},
    \]
    see \cite[Proposition~4.1.2]{CLS11}.
    Therefore,
    \[
    \nu_{Y_{\bullet}}(s)=\Bigl(\langle m, u_{\rho_1} \rangle,\dots,\langle m, u_{\rho_n} \rangle \Bigr)=\Phi(m).
    \]
    Conversely, every section of $L^k$ is a finite linear combination of such characters. In the coordinates on $U_{\sigma}$ determined by $u_{\rho_1},\ldots,u_{\rho_n}$, the Okounkov valuation of a non-zero linear combination is the lexicographically first exponent among the characters with non-zero coefficient. Hence arbitrary sections give no valuations beyond the right-hand side of \eqref{eq:DeltakLtoric}, and \eqref{eq:DeltakLtoric} follows.
\end{proof}

\begin{example}\label{ex:OkounP2}
    Let us continue the example of $\mathbb{P}^2$ in \cref{ex:P2}. We use the same notations. Take $\sigma_1$ as our reference cone, and $\rho_1=\sigma_5$, $\rho_2=\sigma_6$. Then
    \[
    Y_1=\Bigl\{[X_0:0:X_2]:(X_0,X_2)\neq (0,0) \Bigr\},\quad Y_2=\Bigl\{[X_0:0:0]:X_0\neq 0\Bigr\}.
    \]
    The map $\Phi$ is given by
    \[
    \Phi(m_1,m_2)=(m_1,m_2).
    \]
    In this case, we see easily
    \[
    \Delta_{Y_{\bullet}}\left(\mathcal{O}_{\mathbb{P}^2}(1)\right)=P_D
    \]
    is the polytope in \cref{fig:polyP2}.
\end{example}

\begin{example}\label{ex:toric-big:L532}
    Let us continue the example of $\Bl_{0}\mathbb{P}^2$ in \cref{ex:P2bl1pt}. This time, let us take $\sigma_1'$ as our reference cone and $\rho_1=\sigma_5$, $\rho_2=\sigma_7$. Then $Y_1$ is just the strict transform of the line $\{[X_0:0:X_2]:X_0X_2\neq 0\}$ in $\mathbb{P}^2$, while $Y_2$ is the point $Y_1\cap E$, where $E$ is the exceptional divisor.

    In this case, the map $\Phi$ is given by
    \[
    \Phi(m_1,m_2)=(m_1,m_1+m_2).
    \]
    We find that
    \[
    \Delta_{Y_{\bullet}}\left(\Bl_{0}\mathbb{P}^2,\pi^*\mathcal{O}_{\mathbb{P}^2}(1)\right)
    \]
    is the polytope in \cref{fig:P2b1pOkoun}.

    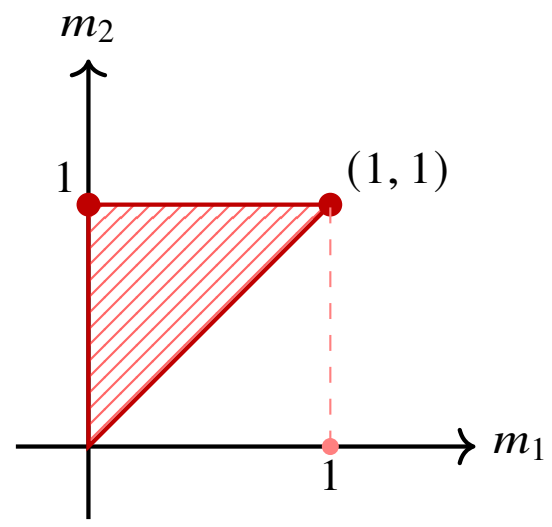
\begin{figure}[ht]
        \centering
        \begin{tikzpicture}[scale=1.6, line width=0.7pt, every node/.style={font=\small}]
          \draw[->, thick] (-0.3,0) -- (1.6,0) node[right] {$m_1$};
          \draw[->, thick] (0,-0.3) -- (0,1.6) node[above] {$m_2$};
          \fill[pattern=north east lines, pattern color=red!60] (0,0) -- (1,1) -- (0,1) -- cycle;
          \draw[red!75!black, thick] (0,0) -- (1,1) -- (0,1) -- cycle;
          \fill[red!75!black] (1,1) circle (1.4pt);
          \fill[red!75!black] (0,1) circle (1.4pt);
          \node[anchor=south west, inner sep=2pt] at (1,1) {$(1,1)$};
          \node[anchor=south east, inner sep=2pt] at (0,1) {$1$};
          \draw[red!50, dashed, thin] (1,0) -- (1,1);
          \fill[red!50] (1,0) circle (1.0pt);
          \node[anchor=north, inner sep=2pt] at (1,0) {$1$};
        \end{tikzpicture}
        \caption{The Okounkov body $\Delta_{Y_{\bullet}}(\Bl_{0}\mathbb{P}^2,\pi^*\mathcal{O}_{\mathbb{P}^2}(1))$.}
        \label{fig:P2b1pOkoun}
    \end{figure}
    Note that it differs from the polytope in \cref{ex:OkounP2}.\footnote{Although these examples are almost trivial, they did confuse me a lot at the beginning of 2023, when Kewei Zhang, Tamás Darvas and I were collaborating on \cite{DRWNXZ}. At that time, Kewei himself already proved the main theorem for a generic flag. I realized that some simple birational geometry would suffice to prove the same result for general flags.
    I persuaded myself and Kewei that the Okounkov bodies are always birationally invariant, and deduced some apparently wrong conclusions. I got no clue for a couple of weeks, then one day, on the noisy metro line~7 of Paris, I got nothing to do, so I said to myself: Why not compute the simplest toric examples? Then after a few minutes, all of a sudden, the whole picture became completely clear.}
\end{example}

\begin{lemma}\label{lma:toricBergmanapprox}
    Let $\varphi\in \PSH_{\tor}(X,\theta)$ and assume that $\theta_{\varphi}$ is a Kähler current. Then $\varphi$ admits a quasi-equisingular approximation $(\varphi_j)_j$ with the following additional properties:
    \begin{enumerate}
        \item\label{itm:lma:toricBergmanapprox:1} $\varphi_j\in \PSH_{\tor}(X,\theta)$ and $\theta_{\varphi_j}$ is a Kähler current for all $j\geq 1$;
        \item\label{itm:lma:toricBergmanapprox:2} for each $j\geq 1$, there are $q_j\in\mathbb Z_{>0}$ and toric eigen-sections
        \[
            s_{j,1},\ldots,s_{j,a_j}\in \mathrm{H}^0(X,L^{q_j})
        \]
        such that
        \[
            \varphi_j=\frac{1}{q_j}\log\sum_{\ell=1}^{a_j}|s_{j,\ell}|_{h^{q_j}}^2+\mathcal O(1).
        \]
    \end{enumerate}
\end{lemma}
\begin{proof}
    Apply Demailly's Bergman-kernel approximation theorem \cite[Theorem~13.21]{Dem12} to the singular Hermitian metric $h\exp(-\varphi)$ on $L$. Since $\theta_{\varphi}$ is a Kähler current, the theorem applies to the Hilbert spaces
    \[
        \mathrm{H}^0\bigl(X,L^q\otimes\mathcal I(q\varphi)\bigr)
    \]
    for $q\gg 1$, endowed with the $L^2$-norm induced by $h^q\exp(-q\varphi)$ and a fixed Kähler volume form. The usual decreasing extraction in Demailly's construction gives a quasi-equisingular approximation; taking $q$ sufficiently large at each step ensures that the approximants still have Kähler Chern currents.

    It remains to note that this construction is compatible with the torus action. Indeed, both $h$ and $\varphi$ are $T_c$-invariant, hence the above Hilbert spaces and their $L^2$-inner products are $T_c$-invariant. Since $T_c$ is compact and abelian, each of these finite-dimensional Hilbert spaces admits an orthonormal basis consisting of toric eigen-sections. With such a basis, the corresponding Bergman potential has the displayed form and is $T_c$-invariant. This proves the lemma.
\end{proof}

\begin{theorem}\label{thm:toricpob}
    Let $\varphi\in \PSH_{\tor}(X,\theta)_{>0}$. Then
    \begin{equation}\label{eq:toricOkounkovcomp}
    \Phi_{\mathbb{R}}\Bigl(  \Delta(\theta,\varphi)\Bigr)=\Delta_{Y_{\bullet}}(\theta,\varphi).
    \end{equation}
    In particular,
    \begin{equation}\label{eq:vol_Newton_general}
        n!\vol \Delta(\theta,\varphi)=\vol \theta_{\varphi}.
    \end{equation}
\end{theorem}
In particular, once $D$ is fixed, the Newton body $\Delta(\theta,\varphi)$ depends only on the current $\theta_{\varphi}$, not on the specific choices of $h$, $\varphi$ and $s_D$. It makes sense to write
\[
\Delta\left(\theta_{\varphi}\right)=\Delta(\theta,\varphi).
\]
A generalization of \eqref{eq:vol_Newton_general} to mixed volumes will be proved in \cref{thm:mix_vol_newbod}.
\begin{proof}
     We first reduce to the case where $D|_{U_{\sigma}}=0$. In fact, changing $D$ to $D+\Div \chi^{m_{\sigma}}$ would result in changing $F_{\theta}$ to $F_{\theta}-m_{\sigma}$. Hence, $F_{\varphi}$ changes to $F_{\varphi}-m_{\sigma}$. Therefore, $\Delta(\theta,\varphi)$ becomes $\Delta(\theta,\varphi)-m_{\sigma}$. Taking \eqref{eq:DUsigma} into consideration, we may assume that $m_{\sigma}=0$.

    \textbf{Step~1}. We first reduce to the case where $\theta_{\varphi}$ is a Kähler current.

    By \cref{lma:kahcurrentposmass}, we can find $\psi\in \PSH(X,\theta)$ such that $\psi\leq \varphi$ and $\theta_{\psi}$ is a Kähler current. Taking the average along $T_c$, we may assume that $\psi$ is $T_c$-invariant.

    For each $t\in (0,1)$, we let
    \[
    \varphi_t=(1-t)\psi+t\varphi.
    \]
    Suppose that Kähler current case is known. Then we get
    \[
    \Phi_{\mathbb{R}}\Bigl(  \Delta(\theta,\varphi_t)\Bigr)=\Delta_{Y_{\bullet}}(\theta,\varphi_t)
    \]
    for any $t\in (0,1)$.
    It follows from \cref{thm:convcomp} that
    \[
    \Phi_{\mathbb{R}}\Bigl(  \Delta(\theta,\varphi)\Bigr)\supseteq \Phi_{\mathbb{R}}\Bigl(  \Delta(\theta,\varphi_t)\Bigr)= \Delta_{Y_{\bullet}}(\theta,\varphi_t)
    \]
    for any $t\in (0,1)$. Thanks to \cref{thm:Okoucont}, we have
    \[
    \Phi_{\mathbb{R}}\Bigl(  \Delta(\theta,\varphi)\Bigr)\supseteq \Delta_{Y_{\bullet}}(\theta,\varphi).
    \]
    Comparing the volumes of both sides using \cref{prop:toricMAandrealMA2}, \eqref{eq:Okov}, and the $\mathcal{I}$-goodness of toric singularities \cref{ex:toricIgood}, we find that
    \[
    n!\vol \Phi_{\mathbb{R}}\Bigl(  \Delta(\theta,\varphi)\Bigr)=\int_X \theta_{\varphi}^n=\vol \theta_{\varphi}=n!\vol \Delta_{Y_{\bullet}}(\theta,\varphi).
    \]
    In particular, we conclude \eqref{eq:toricOkounkovcomp}.

    \textbf{Step~2}. We handle the case where $\theta_{\varphi}$ is a Kähler current.

    Let $(\varphi_j)_j$ be the quasi-equisingular approximation given by \cref{lma:toricBergmanapprox}.

    Now assume that the result is known for each $\varphi_j$. Then
    \[
    \Phi_{\mathbb{R}}\Bigl(  \Delta(\theta,\varphi_j)\Bigr)=\Delta_{Y_{\bullet}}(\theta,\varphi_j).
    \]
    In particular, thanks to \cref{prop:IcompimplyDeltacomp},
    \[
    \Phi_{\mathbb{R}}\Bigl(  \Delta(\theta,\varphi)\Bigr)\subseteq \Delta_{Y_{\bullet}}(\theta,\varphi_j)
    \]
    for each $j\geq 1$. It follows from \cref{thm:Okoucont} that
    \[
    \Phi_{\mathbb{R}}\Bigl(  \Delta(\theta,\varphi)\Bigr)\subseteq \Delta_{Y_{\bullet}}(\theta,\varphi).
    \]
    Comparing the volumes of both sides using \cref{prop:toricMAandrealMA2}, \eqref{eq:Okov}, and \cref{ex:toricIgood}, we conclude \eqref{eq:toricOkounkovcomp}.

    \textbf{Step~3}. It remains to prove the assertion for the analytic approximants appearing in Step~2.

    By \cref{lma:toricBergmanapprox}, these approximants may be chosen in the form
    \[
        \varphi=\frac{1}{q}\log\sum_{i=1}^a |s_i|_{h^q}^2+\mathcal O(1),
    \]
    where $q\in\mathbb Z_{>0}$ and $s_1,\ldots,s_a\in \mathrm{H}^0(X,L^q)$ are toric eigen-sections.

    Let $m_1,\ldots,m_a\in qP_D\cap M$ be the lattice points corresponding to $s_1,\ldots,s_a$.

    Observe that
    \[
    \begin{aligned}
    \Delta(\theta,\varphi)
    =&\overline{\nabla F_{\varphi}(N_{\mathbb{R}})}\\
    =&\left\{ m\in M_{\mathbb{R}} : F_{\varphi}(n)-\langle m,n \rangle \textup{ is bounded from below} \right\}\\
    =& \left\{ m\in M_{\mathbb{R}}:  q^{-1}\log \sum_{i=1}^a \mathrm{e}^{\langle m_i,n\rangle}-\langle m,n\rangle \textup{ is bounded below in } n \right\}\\
    =& q^{-1}\Conv \{m_1,\ldots,m_a\},
    \end{aligned}
    \]
    where we have applied \eqref{eq:normF0} on the second line and \cref{lma:polybdd2} on the third line.
    Since $m_1,\ldots,m_a\in qP_D$, the preceding formula also shows that
    \[
        \Delta(\theta,\varphi)\subseteq P_D.
    \]
    Thus, if $m\in k\Delta(\theta,\varphi)\cap M$, then $m\in kP_D\cap M$, and $\chi^m$ indeed determines a global section of $L^k$.

    In particular, by the same separation argument as in \cref{lma:polybdd}, for any $k\in \mathbb{Z}_{>0}$ and any $m\in k\Delta(\theta,\varphi)\cap M$, we have
    \[
    |\chi^m|^2 \mathrm{e}^{-k\varphi}
    \]
    is bounded from above on $T(\mathbb{C})$. In other words, the section $s$ of $L^k$ corresponding to $m$ satisfies
    \[
    s\in \mathrm{H}^0\left(X,L^k\otimes \mathcal{I}_{\infty}(k\varphi)\right).
    \]
    Therefore,
    \[
    \nu_{Y_{\bullet}}(s)=\Phi(m)\in k\Delta_k(\theta,\varphi),
    \]
    where $\Delta_k$ is defined in \cref{sec:PoB}.
    Hence,
    \[
    \Phi\Bigl(k\Delta(\theta,\varphi)\cap M\Bigr)\subseteq k\Delta_k(\theta,\varphi).
    \]
    Letting $k\to \infty$ and applying \cref{thm:HCP}, we find that
    \[
    \Phi_{\mathbb{R}}\Bigl(\Delta(\theta,\varphi)\Bigr) \subseteq \Delta_{Y_{\bullet}}(\theta,\varphi).
    \]
    Comparing the volumes of both sides using \cref{prop:toricMAandrealMA2} and \eqref{eq:Okov}, we conclude that the equality holds and \eqref{eq:toricOkounkovcomp} follows.
\end{proof}

The following two consequences are both due to Yi Yao.
\begin{corollary}\label{cor:toricLelong}
    Let $E$ be a $T$-invariant prime divisor on $X$ corresponding to a ray $\rho\in \Sigma(1)$. Then for any $\varphi\in \PSH_{\tor}(X,\theta)_{>0}$, we have
    \[
    \nu(\varphi,E)=\inf \left\{ \langle m-m_{\rho},u_{\rho}\rangle : m\in \Delta(\theta,\varphi) \right\}.
    \]
\end{corollary}
\begin{proof}
    This follows immediately from \cref{thm:toricpob} and \cref{thm:nuOk}. In fact, since $X$ is projective and smooth, there is always a $T$-invariant smooth flag $Y_{\bullet}$ with $Y_1=E$.
\end{proof}
This result seems new even in the ample setting. Intuitively, after taking \cref{thm:FVtheta} into consideration as well, in the ample case the generic Lelong number $\nu(\varphi,E)$ is the rescaled ``distance'' from $\Delta(\theta,\varphi)$ to the facet of $P_D$ corresponding to $E$.\footnote{Be cautious! In the big setting, in general, the condition $\langle m-m_{\rho},u_{\rho}\rangle=0$ does not necessarily define a facet of $P_D$. Hence the intuition fails.}

\begin{corollary}\label{cor:toricLelong2}
    Let $Y\subseteq X$ be a $T$-invariant subvariety corresponding to a cone $\sigma$ in $\Sigma$, and let $\varphi\in \PSH_{\tor}(X,\theta)_{>0}$. Then the following are equivalent:
    \begin{enumerate}
        \item\label{itm:cor:toricLelong2:1} $\nu(\varphi,Y)=0$;
        \item\label{itm:cor:toricLelong2:2} there is a point $m\in \Delta(\theta,\varphi)$ such that $\langle m-m_{\rho},u_{\rho}\rangle=0$ for any $\rho\in \sigma(1)$;
        \item\label{itm:cor:toricLelong2:3} we have
        \[
        \Delta(\theta,\varphi)\cap Q_{\sigma} \neq \varnothing.
        \]
    \end{enumerate}
\end{corollary}
Recall that $Q_{\sigma}$ is defined in \eqref{eq:Qsigma}.
\begin{proof}
    \ref{itm:cor:toricLelong2:2} $\iff$ \ref{itm:cor:toricLelong2:3}. This follows from the definition of $Q_{\sigma}$ in \eqref{eq:Qsigma}.

    \ref{itm:cor:toricLelong2:1} $\iff$ \ref{itm:cor:toricLelong2:2}.
    Let $\rho_1,\ldots,\rho_r$ be the rays of $\sigma$. Up to replacing $D$ by a translation, we may assume that $m_{\sigma}=0$. Hence, we may take $m_{\rho_i}=0$ for all $i$.

    Let $\pi\colon Z\rightarrow X$ be the blow-up of $X$ along $Y$.
    See \cite[Page~132]{CLS11} for the basic properties of the toric blow-up.
    Take the divisor $\pi^*D$ on $Z$.
    We choose the pull-back metric $\pi^*h$ on $\pi^*L$. Then $F_{\pi^*\theta}$ can be taken as $\pi^*F_{\theta}$ by \eqref{eq:normF0}. It follows that $\Delta(\theta,\varphi)=\Delta(\pi^*\theta,\pi^*\varphi)$. On the other hand, the ray corresponding to the exceptional divisor $E$ is generated by $u_{\rho_1}+\cdots+u_{\rho_r}$. Since $X$ is smooth, this vector is primitive.

    Recall that the support function of $\pi^*D$ is the same as the support function of $D$, see \cite[Proposition~6.2.7]{CLS11}.
    In particular, we can take the Cartier datum $m_{\rho}=m_{\sigma} \mod M(\rho)$, where $\rho$ is the ray corresponding to $E$.

    It follows from \cref{cor:toricLelong} and \cref{lma:blowupLelong} that
    \begin{equation}\label{eq:nuvarphiYtoric1}
    \nu(\varphi,Y)=\nu(\pi^*\varphi,E)=
    \inf \left\{ \langle m-m_{\sigma},u_{\rho_1}+\cdots+u_{\rho_r}\rangle:m\in \Delta(\theta,\varphi)\right\}.
    \end{equation}
    Our assertion follows in view of \eqref{eq:mmsigpos}.
\end{proof}
It follows from \eqref{eq:nuvarphiYtoric1} that
\[
\nu(\varphi,Y)\geq \sum_{i=1}^r \nu(\varphi,E_i),
\]
where the $E_i$'s are the prime divisors corresponding to the rays of $\sigma$. This inequality seems to be new as well.

The following consequence of \cref{thm:toricpob} is the key to the development of the toric pluripotential theory.
\begin{theorem}\label{thm:FVtheta}
    We have
    \begin{equation}\label{eq:FVthetaE}
    F_{V_{\theta}}\in \mathcal{E}^{\infty}(N_{\mathbb{R}},P_D).\footnotemark
    \end{equation}
    \footnotetext{Initially I was only able to show $F_{V_{\theta}}\in \mathcal{E}(N_{\mathbb{R}},P_D)$, the strengthened version was suggested by Robert Berman.}
    In particular,
    \begin{equation}\label{eq:toricminsingmass}
        \int_X  \theta_{V_{\theta}}^n=n!\vol P_D.
    \end{equation}
\end{theorem}
Recall that $\mathcal{E}^{\infty}$ is defined in \cref{def:convexPfunctions}. The equation \eqref{eq:FVthetaE} says
\begin{equation}\label{eq:FSupp_bdd_temp1}
F_{V_{\theta}}-\Supp_{P_D} \textup{ is bounded}.
\end{equation}
In particular,
\begin{equation}\label{eq:DeltaVthetaPD_temp1}
\Delta(\theta,V_{\theta})=P_D
\end{equation}
and hence the Newton bodies $\Delta(\theta,\varphi)$ are all contained in $P_D$.

\begin{proof}
    Taking $\varphi=V_{\theta}$ in \cref{thm:toricpob}, we find
    \[
    \Phi_{\mathbb{R}}\left(  \Delta(\theta,V_{\theta})\right)=\Delta_{Y_{\bullet}}(\theta,V_{\theta})=\Delta_{Y_{\bullet}}(L)=\Phi_{\mathbb{R}}(P_D),
    \]
    where we applied \cref{prop:toricusualOko} in the last equality.
    Therefore, \eqref{eq:DeltaVthetaPD_temp1} follows.
    In particular, \eqref{eq:toricminsingmass} follows from \cref{prop:toricMAandrealMA2}.

    Next we prove \eqref{eq:FSupp_bdd_temp1}. Take a smooth function $\chi\colon N_{\mathbb{R}}\rightarrow [0,1]$ so that
    \begin{enumerate}
        \item\label{itm:chapter-toricb:L672:1} $\Supp \chi$ is compact;
        \item\label{itm:chapter-toricb:L672:2} $\chi\equiv 1$ on a neighbourhood of $0$.
    \end{enumerate}
    Fix a $T_c$-invariant volume form $\mu$ on $T(\mathbb{C})$. Then $(\Trop^*\chi)\mu$ can be regarded as a $T_c$-invariant smooth non-negative measure on $X$. After choosing a suitable normalizing constant $C>0$, there is a unique solution
    \begin{equation}\label{eq:thetavarphinD_temp1}
    \theta_{\varphi}^n=C\left(\Trop^*\chi\right)\mu,\quad \sup_X \varphi=0,\quad \varphi\in \mathcal{E}^{\infty}(X,\theta).
    \end{equation}
    This is proved in \cite[Theorem~B]{BEGZ10}. The uniqueness of $\varphi$ further guarantees its $T_c$-invariance. Since $\varphi$ has minimal singularities, $F_{\varphi}\sim F_{V_{\theta}}$. Next, using \cref{prop:toricMAandrealMA2}, the Monge--Ampère equation \eqref{eq:thetavarphinD_temp1} implies
    \[
    \MA_{\mathbb{R}}\left(F_{\varphi}\right)=C \chi \Trop_*\mu.
    \]
    It follows from the regularity theorem \cite[Theorem~2.19]{BB13} and \eqref{eq:DeltaVthetaPD_temp1} that $F_{\varphi}\sim \Supp_{P_D}$. Hence, \eqref{eq:FSupp_bdd_temp1} follows.
\end{proof}

In particular, thanks to \cref{cor:toricLelong},
\begin{equation}\label{eq:nuVtheta_toric}
\nu(V_{\theta},E)=\min \left\{ \langle m-m_{\rho},u_{\rho}\rangle : m\in P_D \right\}.
\end{equation}
As an interesting consequence of \eqref{eq:nuVtheta_toric}, we have a geometric description of the divisorial Zariski decomposition of \cite{Bou02} in the toric setting: The negative part of $D$ is given by
\[
\sum_{\rho\in \Sigma(1)} \min \left\{ \langle m-m_{\rho},u_{\rho}\rangle : m\in P_D \right\} D_{\rho}.
\]
A prime divisor $D_{\rho}$ appears in the negative part if and only if the corresponding condition
\[
\langle m,u_{\rho} \rangle \geq -a_{\rho}
\]
is redundant in defining $P_D$ by \eqref{eq:bigPrep}. The modified nef part of $D$ is the $\mathbb{Q}$-divisor obtained after removing these redundancies.

\section{The pluripotential theory}\label{sec:toric_ppt_general}
We continue to use the notations in \cref{sec:toricsetup2}.

\begin{theorem}\label{thm:toricpshbig}
    There are canonical bijections\footnote{In the earlier version of this book, I required additional conditions in \ref{itm:thm:toricpshbig:2} and \ref{itm:thm:toricpshbig:3}, namely $F\preceq F_{V_{\theta}}$ and $G\succeq G_{V_{\theta}}$ respectively. These conditions can be removed, since $G_{V_{\theta}}$ is always bounded, as suggested by Robert Berman. See \cref{thm:FVtheta}.} between the following sets:
    \begin{enumerate}
        \item\label{itm:thm:toricpshbig:1} The set of $\varphi\in \PSH_{\tor}(X,\theta)$;
        \item\label{itm:thm:toricpshbig:2} the set of $F\in \mathcal{P}(N_{\mathbb{R}},P_D)$, and
        \item\label{itm:thm:toricpshbig:3} the set of closed proper convex functions $G\in \Conv(M_{\mathbb{R}})$ satisfying
        \[
        G|_{M_{\mathbb{R}}\setminus P_D}\equiv \infty.
        \]
    \end{enumerate}
\end{theorem}
The set $\mathcal{P}(N_{\mathbb{R}},P_D)$ is defined in \cref{def:convexPfunctions}.
As before, we write $F_{\varphi}$, $G_{\varphi}$ for the functions determined by this construction.

\begin{proof}
    The proof is similar to that of \cref{thm:toricpsh}, but due to its importance, we give the details. Again, the correspondence between \ref{itm:thm:toricpshbig:2} and \ref{itm:thm:toricpshbig:3} follows easily from \cref{prop:Fsuppchar}.

    Given $\varphi$, we can construct $F_{\varphi}$ in \ref{itm:thm:toricpshbig:2} as explained earlier in \eqref{eq:Fvarphibig1}. Conversely, suppose that $F\in \mathcal{P}(N_{\mathbb{R}},P_D)$, then $F\preceq F_{V_{\theta}}$ by \cref{thm:FVtheta}. Then
    \[
        \Trop^*\left(F-F_{\theta}\right)\in \PSH\left(T(\mathbb{C}),\theta|_{T(\mathbb{C})}\right)
    \]
    by \cref{lma:convextopsh}.
    Since $F\preceq F_{V_{\theta}}$, we see that $\Trop^*(F-F_{\theta})$ is bounded from above. It follows that Grauert--Remmert's extension theorem \cref{thm:GRexten} is applicable, and this function extends to a unique $\theta$-psh function $\varphi$. The uniqueness of the extension guarantees that $\varphi\in \PSH_{\tor}(X,\theta)$.

    The two maps are inverse to each other by construction.
\end{proof}
We fix a model potential $\phi\in \PSH_{\tor}(X,\theta)_{>0}$ with Newton body $\Delta(\theta,\phi)$.

A similar argument guarantees the following:
\begin{corollary}\label{cor:toric-big:L860}
    There is a canonical bijection between the following sets:
    \begin{enumerate}
        \item\label{itm:chapter-toricb:L739:1} The set of $\varphi\in \PSH_{\tor}(X,\theta;\phi)$,
        \item\label{itm:chapter-toricb:L739:2} the set of $F\in \mathcal{P}(N_{\mathbb{R}},\Delta(\theta,\phi))$, and
        \item\label{itm:chapter-toricb:L739:3} the set of closed proper convex functions $G\in \Conv(M_{\mathbb{R}})$ satisfying
        \[
        G|_{M_{\mathbb{R}}\setminus \Delta(\theta,\phi)}=\infty.
        \]
    \end{enumerate}
    Moreover, under these correspondences, we have the following bijections:
    \begin{enumerate}
        \item\label{itm:chapter-toricb:L748:1} The set $\mathcal{E}_{\tor}(X,\theta;\phi)$,
        \item\label{itm:chapter-toricb:L748:2} the set of $F\in \mathcal{E}(N_{\mathbb{R}},\Delta(\theta,\phi))$, and
        \item\label{itm:chapter-toricb:L748:3} the set of closed proper convex functions $G\in \Conv(M_{\mathbb{R}})$ satisfying
        \[
            \overline{\Int \{G<\infty\}}=\Delta(\theta,\phi).
        \]
    \end{enumerate}
\end{corollary}
Here the notations are defined as follows:
\[
\begin{aligned}
\PSH_{\tor}(X,\theta;\phi)\coloneqq & \{ \varphi\in \PSH_{\tor}(X,\theta):\varphi\preceq \phi \},\\
\mathcal{E}_{\tor}(X,\theta;\phi)\coloneqq & \mathcal{E}(X,\theta;\phi)\cap \PSH_{\tor}(X,\theta).
\end{aligned}
\]

The proofs of the following results are similar to the ample case studied in \cref{chap:toric_ample}. We omit the details.

\begin{proposition}\label{prop:toricpluscstbig}
    Given $\varphi\in \PSH_{\tor}(X,\theta)$ and $C\in \mathbb{R}$. We have
    \[
    F_{\varphi+C}=F_{\varphi}+C,\quad G_{\varphi+C}=G_{\varphi}-C.
    \]
\end{proposition}

\begin{proposition}\label{prop:toricrooftopbig}
    Given $\varphi,\psi\in \PSH_{\tor}(X,\theta)$, assume that $\varphi\land\psi\not\equiv -\infty$,
    then $\varphi\land \psi\in \PSH_{\tor}(X,\theta)$ and
    \[
    F_{\varphi\land \psi}=F_{\varphi}\land F_{\psi},\quad G_{\varphi\land \psi}=G_{\varphi}\lor G_{\psi}.
    \]
\end{proposition}
\begin{proposition}\label{prop:toricseqbig}
    Let $(\varphi_i)_{i\in I}$ be a family in $\PSH_{\tor}(X,\theta)$ uniformly bounded from above. Then $\uscsup[i\in I] \varphi_i\in \PSH_{\tor}(X,\theta)$ and
    \[
    F_{\uscsup[i\in I]\varphi_i}=\bigvee_{i\in I}F_{\varphi_i},\quad G_{\uscsup[i\in I]\varphi_i}=\cl \bigwedge_{i\in I} G_{\varphi_i}.
    \]
    Moreover, if $I$ is finite, then
    \[
        G_{\bigvee_{i\in I}\varphi_i}=\bigwedge_{i\in I} G_{\varphi_i}.
    \]

    Similarly, if $\{\varphi_i\}_{i\in I}$ is a decreasing net in $\PSH_{\tor}(X,\theta)$ such that $\inf_{i\in I}\varphi_i\not\equiv -\infty$, then $\inf_{i\in I}\varphi_i\in \PSH_{\tor}(X,\theta)$ and
    \[
    F_{\inf_{i\in I}\varphi_i}=\bigwedge_{i\in I}F_{\varphi_i},\quad G_{\inf_{i\in I}\varphi_i}=\bigvee_{i\in I}G_{\varphi_i}.
    \]
\end{proposition}

\begin{proposition}\label{prop:GPenvelopebig}
    Let $\varphi\in \PSH_{\tor}(X,\theta)$. Then $P_{\theta}[\varphi]\in \PSH_{\tor}(X,\theta)$ and
    \begin{equation}\label{eq:toricPenvbig}
    G_{P_{\theta}[\varphi]}(x)=
    \begin{cases}
        G_{V_{\theta}}(x),&\textup{if } x\in \Delta(\theta,\varphi);\\
        \infty,&\textup{otherwise}.
    \end{cases}
    \end{equation}
\end{proposition}
\begin{proof}
    By \eqref{eq:Penvsups},
    \[
    P_{\theta}[\varphi]=\uscsup[C\in \mathbb{R}]\bigl((\varphi+C)\land V_{\theta}\bigr).
    \]
    Hence $P_{\theta}[\varphi]\in \PSH_{\tor}(X,\theta)$ by \cref{prop:toricpluscstbig,prop:toricrooftopbig,prop:toricseqbig}. The same results give
    \[
    G_{P_{\theta}[\varphi]}=\cl\inf_{C\in \mathbb{R}}\left(G_{V_{\theta}}\lor (G_{\varphi}-C)\right).
    \]
    Put
    \[
    K=\inf_{C\in \mathbb{R}}\left(G_{V_{\theta}}\lor (G_{\varphi}-C)\right).
    \]
    If $x\in \Dom G_{\varphi}$, then $G_{\varphi}(x)-C\leq G_{V_{\theta}}(x)$ for all large enough $C$, so $K(x)=G_{V_{\theta}}(x)$. If $x\notin \Dom G_{\varphi}$, then $K(x)=\infty$. By \cref{prop:gradDomFstar}, we have $\Delta(\theta,\varphi)=\overline{\Dom G_{\varphi}}$. Moreover, \cref{thm:FVtheta} gives $F_{V_{\theta}}\sim \Supp_{P_D}$, hence $G_{V_{\theta}}\sim \chi_{P_D}$ by Legendre duality. In particular, $G_{V_{\theta}}$ is finite and continuous on $P_D$. Thus the lower semicontinuous regularization of $K$ is equal to $G_{V_{\theta}}$ on $\Delta(\theta,\varphi)$ and equal to $\infty$ outside $\Delta(\theta,\varphi)$.
\end{proof}
As a consequence, we have
\begin{corollary}\label{cor:PequivtorD}
    Let $\varphi,\psi\in \PSH_{\tor}(X,\theta)$. Then the following are equivalent:
    \begin{enumerate}
        \item\label{itm:cor:PequivtorD:1} $\varphi\preceq_P \psi$;
        \item\label{itm:cor:PequivtorD:2} $\varphi\preceq_{\mathcal{I}} \psi$;
        \item\label{itm:cor:PequivtorD:3} $\Delta(\theta,\varphi)\subseteq \Delta(\theta,\psi)$.
    \end{enumerate}
\end{corollary}
\begin{proof}
    \ref{itm:cor:PequivtorD:1} $\iff$ \ref{itm:cor:PequivtorD:2}. This follows from the $\mathcal{I}$-goodness of $\varphi$ and $\psi$, see \cref{ex:toricIgood}.

    \ref{itm:cor:PequivtorD:1} $\iff$ \ref{itm:cor:PequivtorD:3}. When $\varphi,\psi\in \PSH_{\tor}(X,\theta)_{>0}$, this follows from \cref{prop:GPenvelopebig}.

    In general, fix an ample toric-invariant divisor $H$ on $X$ and a toric-invariant Hermitian metric $h_H$ on $\mathcal{O}_X(H)$ such that $\omega_H\coloneqq \ddc h_H$ is a Kähler form. We perturb $L$ to $L+m^{-1}H$ for $m\in \mathbb{Z}_{>0}$, then thanks to \cref{lma:grad_img_sum}, we find that
    \[
    \Delta(\theta+m^{-1}\omega_H,\varphi)\coloneqq m^{-1}\Delta(m\theta+\omega_H,m\varphi)=\Delta(\theta,\varphi)+m^{-1}\Delta\left(\mathcal{O}_X(H)\right).
    \]
    Then since $\varphi,\psi\in \PSH_{\tor}(X,\theta+m^{-1}\omega_H)_{>0}$, we know that $\varphi\preceq_P \psi$ if and only if
    \[
    \Delta(\theta,\varphi)+m^{-1}\Delta\left(\mathcal{O}_X(H)\right)\subseteq \Delta(\theta,\psi)+m^{-1}\Delta\left(\mathcal{O}_X(H)\right)
    \]
    for any $m\in \mathbb{Z}_{>0}$. The latter condition is equivalent to $\Delta(\theta,\varphi)\subseteq \Delta(\theta,\psi)$.
\end{proof}
\begin{corollary}\label{cor:bij_conv_torpsh}
    There is a canonical bijection
    \[
    \PSH_{\tor}(X,\theta)/\sim_{P}=\PSH_{\tor}(X,\theta)/\sim_{\mathcal{I}}\rightarrow \left\{K\in \mathcal{K}_n: K \subseteq P_D \right\}.
    \]
    Here $\mathcal{K}_n$ denotes the set of convex bodies in $\mathbb{R}^n$, as in \cref{sec:cb}.
\end{corollary}
More precisely, given $\varphi\in \PSH_{\tor}(X,\theta)$, the corresponding convex body is $\Delta(\theta,\varphi)$.

\begin{proof}
    Thanks to \cref{cor:PequivtorD}, we have the equality
    \[
    \PSH_{\tor}(X,\theta)/\sim_{P}=\PSH_{\tor}(X,\theta)/\sim_{\mathcal{I}},
    \]
    and the map
    \[
    \PSH_{\tor}(X,\theta)/\sim_{\mathcal{I}}\rightarrow \left\{K\in \mathcal{K}_n: K \subseteq P_D \right\}
    \]
    is well-defined and injective.

    In order to see the surjectivity, take $K\in \mathcal{K}_n$, $K\subseteq P_D$. Let $\chi_K$ be the characteristic function of $K$ as defined in \cref{ex:charconvex}. Since $\chi_K$ satisfies Condition~\ref{itm:thm:toricpshbig:3} of \cref{thm:toricpshbig}, we can find $\varphi\in \PSH_{\tor}(X,\theta)$ so that $\chi_K=G_{\varphi}$. Then $\Delta(\theta,\varphi)=K$ due to \cref{prop:gradDomFstar}.
\end{proof}

Next we handle subgeodesics.
\begin{proposition}\label{prop:toricsubgeod}
    Let $\varphi_0,\varphi_1\in \PSH_{\tor}(X,\theta)$.
    There is a canonical bijection between the following sets:
    \begin{enumerate}
        \item\label{itm:prop:toricsubgeod:1} The set of $T_c$-invariant subgeodesics from $\varphi_0$ to $\varphi_1$;
        \item\label{itm:prop:toricsubgeod:2} the set of convex functions $F\colon N_{\mathbb{R}}\times (0,1)\rightarrow \mathbb{R}$ such that for each $r\in (0,1)$, the function
        \begin{equation}\label{eq:Ftdefsubgeo}
        F_r\colon N_{\mathbb{R}}\rightarrow \mathbb{R},\quad n\mapsto F(n,r)
        \end{equation}
        satisfies $F_r\to F_{\varphi_1}$ (resp. $F_r\to F_{\varphi_0}$) everywhere as $r\to 1-$ (resp. $r\to 0+$).
    \end{enumerate}
\end{proposition}

\begin{proof}
We begin with a subgeodesic $(\varphi_t)_{t\in (0,1)}$ from $\varphi_0$ to $\varphi_1$. Then we define $F\colon N_{\mathbb{R}}\times (0,1)\rightarrow \mathbb{R}$ as follows:
    \[
    F(n,t)=F_{\varphi_t}(n).
    \]
    Define $F_t$ as in \eqref{eq:Ftdefsubgeo}, we have
    \[
    \Trop^* F_t-\Trop^* F_{\theta}=\varphi_t,\quad t\in (0,1).
    \]
    Let $t_j\to 0+$. Since $\varphi_{t_j}\to\varphi_0$ in $L^1$, after passing to a subsequence we have almost-everywhere convergence on $X$. Restricted to the open torus and using the tropical coordinates, this gives $F_{t_j}\to F_{\varphi_0}$ almost everywhere on $N_{\mathbb R}$. Since the $F_{t_j}$'s are convex and locally uniformly bounded above, \cref{thm:seqconv} upgrades this to pointwise convergence everywhere. As the limit is independent of the subsequence, $F_t\to F_{\varphi_0}$ everywhere as $t\to0+$. Similarly, $F_t\to F_{\varphi_1}$ everywhere as $t\to 1-$.

    Next we show that $F$ is convex. Let $p_1\colon X\times S\rightarrow X$ be the projection, where
    \[
    S\coloneqq \{z\in \mathbb{C}:  \mathrm{e}^{-1}<|z|^2<1\}.
    \]
    Since $(\varphi_t)_{t\in(0,1)}$ is a subgeodesic, its complexification $\Phi$ is $p_1^*\theta$-psh. Recall that $\Phi$ is defined as
    \begin{equation}\label{eq:Phixz}
    \Phi(x,z)=\Trop^*\left( F_{-\log |z|^2}-F_{\theta} \right)(x).
    \end{equation}
    In particular, $\Psi\colon T(\mathbb{C})\times S\rightarrow \mathbb{R}$ defined by
    \[
    \Psi(x,z)\coloneqq \Phi(x,z)+\Trop^*F_{\theta}(x)=\Trop^*F_{-\log |z|^2}(x)
    \]
    is plurisubharmonic and $T_c\times S^1$-invariant.
    Fix a small enough $\epsilon>0$. We can find a decreasing sequence of $T_c\times S^1$-invariant plurisubharmonic functions $\Psi_i$ on $T(\mathbb{C})\times S_{\epsilon}$ converging to $\Psi$ everywhere, where
    \[
    S_{\epsilon}\coloneqq \{z\in \mathbb{C}: \mathrm{e}^{-1+\epsilon}<|z|^2<\mathrm{e}^{-\epsilon}\}.
    \]
    Let us write
    \[
    \Psi_i(x,z)=\Trop^*F_{i,-\log |z|^2}(x)
    \]
    for some $F_i\colon N_{\mathbb{R}}\times (\epsilon,1-\epsilon)\rightarrow \mathbb{R}$.

    The same computation as in \cref{lma:convextopsh} shows that $F_i$ is convex. It follows that $F$, as the decreasing limit of $F_i$, is also convex on $N_{\mathbb{R}}\times (\epsilon,1-\epsilon)$. Since $\epsilon>0$ is arbitrary, we conclude that $F$ is convex on $N_{\mathbb{R}}\times (0,1)$.

    Conversely, suppose that we are given $F$ in \ref{itm:prop:toricsubgeod:2}. We define $\Phi\colon T(\mathbb{C})\times S\rightarrow \mathbb{R}$ using \eqref{eq:Phixz}.
    The arguments in the previous part can be reversed to show that $\Phi$ is $p_1^*\theta|_{T(\mathbb{C})\times S}$-psh.

    By our assumption, for each $t\in (0,1)$, we have
    \begin{equation}\label{eq:Ftupptemp1}
    F_t\leq tF_{\varphi_1}+(1-t)F_{\varphi_0}\leq F_{V_{\theta}}+C
    \end{equation}
    for some constant $C\in \mathbb{R}$ independent of the choice of $t$.
    Therefore, $\Phi$ is bounded from above and hence, by \cref{thm:GRexten}, admits a unique $p_1^*\theta$-psh extension to $X\times S$, which we still denote by $\Phi$. We let
    \[
    \varphi_t(x)=\Phi(x,\mathrm{e}^{-t/2})
    \]
    for all $t\in (0,1)$ and $x\in X$. We claim that $(\varphi_t)$ is a subgeodesic from $\varphi_0$ to $\varphi_1$.

    For this purpose, we only need to show that $(\varphi_t)_{t\in (0,1)}$ has the correct boundary value. But from our assumption in \ref{itm:prop:toricsubgeod:2}, we know that as $t\to 0+$ (resp. $t\to 1-$), $\varphi_t\to \varphi_0$ (resp. $\varphi_t\to \varphi_1$) almost everywhere. In particular, $\sup_X \varphi_t\geq -C'$ for some large constant $C'>0$ independent of $t\in (0,1)$. Therefore, together with \eqref{eq:Ftupptemp1}, we deduce from \cref{prop:L1compa} that $\{\varphi_t\}_{t\in (0,1)}$ is a relatively compact family with respect to the $L^1$-topology. We need to show that each cluster point $\psi$ as $t\to 0+$ is equal to $\varphi_0$. But we already know that $\psi=\varphi_0$ almost everywhere. Hence we deduce $\psi=\varphi_0$ from \cref{prop:pshfuncdetdense}. As $t\to 0+$, we have $\varphi_t\xrightarrow{L^1}\varphi_0$.
    Similarly, as $t\to 1-$, we have $\varphi_t\xrightarrow{L^1}\varphi_1$.

    The two constructions are inverse to each other by construction.

\end{proof}

\begin{corollary}\label{cor:geodG}
    Let $\varphi_0,\varphi_1\in \PSH_{\tor}(X,\theta)$. Then there is a canonical bijection between the following sets:
    \begin{enumerate}
        \item\label{itm:cor:geodG:1} The set of $T_c$-invariant subgeodesics from $\psi_0$ to $\psi_1$, where $\psi_0,\psi_1\in \PSH_{\tor}(X,\theta)$ and $\psi_0\leq \varphi_0$, $\psi_1\leq \varphi_1$;
        \item\label{itm:cor:geodG:2} the set of closed proper convex functions $\Psi$ on $M_{\mathbb{R}}\times \mathbb{R}$ such that there is a closed proper convex function $G\in \Conv(M_{\mathbb{R}})$ satisfying
        \begin{equation}\label{eq:Psilower}
        G(m)+(s\lor 0)\geq \Psi(m,s)\geq G_{\varphi_0}(m)\lor \left( G_{\varphi_1}(m)+s \right)
        \end{equation}
        for all $m\in M_{\mathbb{R}}$ and $s\in \mathbb{R}$.
    \end{enumerate}
\end{corollary}
\begin{proof}
    Let us begin with a subgeodesic $(\psi_t)_{t\in (0,1)}$ as in \ref{itm:cor:geodG:1}. Let $F$ be the convex function as in \cref{prop:toricsubgeod}.
    We extend $F$ to a function $F\colon N_{\mathbb{R}}\times \mathbb{R}\rightarrow \mathbb{R}$ as follows: For any $n\in N_{\mathbb{R}}$, we define
\[
F(n,t)=\begin{cases}
    F_{\psi_t}(n),&\textup{if }0<t<1,\\
    F_{\psi_0}(n),& \textup{if }t=0,\\
    F_{\psi_1}(n),&\textup{if }t=1,\\
    \infty,&\textup{if } t>1 \textup{ or }t<0.
\end{cases}
\]
Then $F$ is a proper closed convex function on $N_{\mathbb{R}}\times \mathbb{R}$. Let $\Psi$ be the Legendre transform of $F$. Then $\Psi$ is a proper closed convex function on $M_{\mathbb{R}}\times \mathbb{R}$ by \cref{thm:Legendredual}.
By \eqref{eq:Legtrandef}, for any $m\in M_{\mathbb{R}}$ and $s\in \mathbb{R}$, we have
\[
\begin{aligned}
\Psi(m,s)=&\sup_{n\in N_{\mathbb{R}},t\in [0,1]} \left( \langle m,n\rangle + ts-F(n,t)\right)\\
=& \sup_{t\in [0,1]} \left( ts+F_t^*(m) \right).
\end{aligned}
\]
Therefore, the latter half of \eqref{eq:Psilower} follows. Next recall that
\[
\eta\coloneqq \inf_{t\in (0,1)}\psi_t\in \PSH_{\tor}(X,\theta),
\]
as follows from \cref{prop:convexsubgeod}. Therefore,
\[
\begin{aligned}
\Psi(m,s)=&\sup_{n\in N_{\mathbb{R}},t\in [0,1]} \left( \langle m,n\rangle + ts-F(n,t)\right)\\
\leq & \sup_{n\in N_{\mathbb{R}},t\in [0,1]} \left( \langle m,n\rangle + ts-F_{\eta}(n)\right)\\
=& \sup_{t\in [0,1]}ts+G_{\eta}(m)\\
=& (s\lor 0)+G_{\eta}(m).
\end{aligned}
\]

Conversely, let us begin with a function $\Psi$ as in \ref{itm:cor:geodG:2}. Let $F$ be the Legendre transform of $\Psi$.
We first observe that $F(n,t)=\infty$ for all $n\in N_{\mathbb{R}}$ and $t\not\in [0,1]$.
Indeed, the upper bound $\Psi\leq G+(s\lor 0)$ gives
\[
F(n,t)\geq G^*(n)+\sup_{s\in\mathbb{R}}\left(ts-(s\lor 0)\right),
\]
and the last supremum is infinite exactly when $t\notin[0,1]$.

For $t\in[0,1]$, the lower bound $\Psi\geq G_{\varphi_0}\lor(G_{\varphi_1}+s)$ implies
\[
\Psi(m,s)\geq (1-t)G_{\varphi_0}(m)+tG_{\varphi_1}(m)+ts.
\]
Consequently,
\[
F(n,t)\leq \left((1-t)G_{\varphi_0}+tG_{\varphi_1}\right)^*(n)
\leq (1-t)F_{\varphi_0}(n)+tF_{\varphi_1}(n)
\]
for all $t\in [0,1]$ and $n\in N_{\mathbb{R}}$.
On the other hand, the upper bound $\Psi\leq G+(s\lor 0)$ gives
\[
\begin{aligned}
F(n,t)=&\sup_{m\in M_{\mathbb{R}},s\in \mathbb{R}} \left(\langle m,n\rangle +ts-\Psi(m,s)\right)\\
\geq & \sup_{m\in M_{\mathbb{R}},s\in \mathbb{R}} \left(\langle m,n\rangle +ts-G(m)-(s\lor 0)\right)\\
=& G^*(n)+\sup_{s\in \mathbb{R}} \left( ts-(s\lor 0)\right)\\
=&
\begin{cases}
    G^*(n),&\text{if }t\in [0,1],\\
    \infty,&\text{otherwise}.
\end{cases}
\end{aligned}
\]
The convexity of $F$ gives one-sided boundary limits
\[
F_0(n)\coloneqq \lim_{t\to 0+}F(n,t),\quad
F_1(n)\coloneqq \lim_{t\to 1-}F(n,t).
\]
The estimates above give $F_0\leq F_{\varphi_0}$ and $F_1\leq F_{\varphi_1}$, and also show that $F_0,F_1\in \mathcal{P}(N_{\mathbb{R}},P_D)$. Let $\psi_0,\psi_1\in \PSH_{\tor}(X,\theta)$ be the corresponding potentials given by \cref{thm:toricpshbig}. Then $\psi_0\leq\varphi_0$ and $\psi_1\leq\varphi_1$. Applying \cref{prop:toricsubgeod} to the restriction of $F$ to $N_{\mathbb{R}}\times (0,1)$, we get a subgeodesic $(\psi_t)_{t\in (0,1)}$ from $\psi_0$ to $\psi_1$. Thus $(\psi_t)_{t\in (0,1)}$ satisfies \ref{itm:cor:geodG:1}.
The two operations are inverse to each other by construction.
\end{proof}

As an immediate corollary,
\begin{corollary}\label{cor:toricgeodgeneral}
    Let $\varphi_0,\varphi_1\in \PSH_{\tor}(X,\theta)\cap \PSH(X,\theta)_{>0}$.
    Then the following are equivalent:
    \begin{enumerate}
        \item\label{itm:cor:toricgeodgeneral:1} $\varphi_0\sim_P \varphi_1$;
        \item\label{itm:cor:toricgeodgeneral:2} there is a subgeodesic from $\varphi_0$ to $\varphi_1$.
   \end{enumerate}
    If these conditions are satisfied, let $(\varphi_t)_{t\in (0,1)}$ be the geodesic from $\varphi_0$ to $\varphi_1$. Then $\varphi_t\in \PSH_{\tor}(X,\theta)$ for all $t\in (0,1)$ and
    \begin{equation}\label{eq:Glinear}
    G_{\varphi_t}=(1-t)G_{\varphi_0}+tG_{\varphi_1}.
    \end{equation}
\end{corollary}
\begin{proof}
    The equivalence between \ref{itm:cor:toricgeodgeneral:1} and \ref{itm:cor:toricgeodgeneral:2} follows from a general fact \cref{thm:geodexi}. In the toric case, \ref{itm:cor:toricgeodgeneral:2} $\implies$ \ref{itm:cor:toricgeodgeneral:1} can also be argued more directly.

    Assume that these conditions are satisfied. Let $(\varphi_t)_{t\in (0,1)}$ be the geodesic from $\varphi_0$ to $\varphi_1$. Since the geodesic is the upper envelope of toric subgeodesics, $\varphi_t\in \PSH_{\tor}(X,\theta)$ for all $t\in (0,1)$. Let $\Psi'$ be the proper convex function on $M_{\mathbb{R}}\times \mathbb{R}$ defined by \cref{cor:geodG}. Then $\Psi'$ is the minimum of all $\Psi$ satisfying \eqref{eq:Psilower}. We claim that
    \begin{equation}\label{eq:Psiptemp1}
    \Psi'(m,s)=G_{\varphi_0}(m)\lor \left(G_{\varphi_1}(m)+s\right).
    \end{equation}
    It suffices to show that the right-hand side is proper, namely, $G_{\varphi_0}\lor G_{\varphi_1}$ is not identically $\infty$. But recall that by \cref{prop:convexsubgeod}, we have $\varphi_0\land \varphi_1\in \PSH(X,\theta)$. Therefore, by \cref{prop:toricrooftopbig},
    \[
    G_{\varphi_0}\lor G_{\varphi_1}=G_{\varphi_0\land \varphi_1}\not\equiv \infty.
    \]
    In particular, \eqref{eq:Psiptemp1} follows.

    Now by construction,
    \[
    G_{\varphi_t}(m)=\inf_{s\in \mathbb{R}} \left( \Psi'(m,s)-st\right)=(1-t)G_{\varphi_0}(m)+tG_{\varphi_1}(m)
    \]
    for all $t\in (0,1)$. So \eqref{eq:Glinear} follows.
\end{proof}

Let $\phi\in \PSH_{\tor}(X,\theta)\cap \PSH(X,\theta)_{>0}$ be a model potential. We write
\[
\mathcal{R}_{\tor}(X,\theta;\phi)\coloneqq \left\{ \ell\in \mathcal{R}(X,\theta;\phi): \ell_t\in \PSH_{\tor}(X,\theta)\quad\forall t\geq 0\right\}.
\]
\index{$\mathcal{R}_{\tor}(X,\theta;\phi)$}
Recall that $\mathcal{R}(X,\theta;\phi)$ is defined in \cref{def:geodraydef}.
\begin{corollary}\label{cor:geod_tor_char}
There is a canonical bijection between the following sets:
\begin{enumerate}
    \item\label{itm:cor:geod_tor_char:1} The set of $\ell\in \mathcal{R}_{\tor}(X,\theta;\phi)$;
    \item\label{itm:cor:geod_tor_char:2} The set of proper closed convex functions $g$ on $M_{\mathbb{R}}$ with $\overline{\Dom g}=\Delta(\theta,\phi)$.
\end{enumerate}
Moreover, given $\ell\in \mathcal{R}_{\tor}(X,\theta;\phi)$, then
\begin{equation}\label{eq:geod_ray_G}
G_{\ell_t}=G_{\phi}+tg,\quad \forall t>0.
\end{equation}
\end{corollary}
\begin{proof}
    First observe that given $g$ as in \ref{itm:cor:geod_tor_char:2}, \eqref{eq:geod_ray_G} indeed induces a geodesic ray in $\mathcal{R}_{\tor}(X,\theta;\phi)$ thanks to \cref{cor:toricgeodgeneral}. Therefore, we have a map from \ref{itm:cor:geod_tor_char:2} to \ref{itm:cor:geod_tor_char:1}.

    Conversely, given $\ell\in \mathcal{R}_{\tor}(X,\theta;\phi)$, it follows from \cref{cor:toricgeodgeneral} that $G_{\ell_t}$ is linear in $t\geq 0$ after restriction to $\Int \{G_{\phi}<\infty\}$. Let
    \[
    g'(m)=
    \begin{cases}
    G_{\ell_1}(m)-G_{\phi}(m),& m\in \Int \{G_{\phi}<\infty\};\\
    \infty,& m\in M_{\mathbb{R}}\setminus \Int \{G_{\phi}<\infty\}.
    \end{cases}
    \]
    Note that for $m\in \Int \{G_{\phi}<\infty\}$, we actually have
    \[
    g'(m)=\lim_{t\to\infty}t^{-1}G_{\ell_t}(m).
    \]
    Hence, $g'$ is a proper convex function on $M_{\mathbb{R}}$.
    Define $g=\cl g'$. Then $g$ is a proper closed convex function on $M_{\mathbb{R}}$, and $\overline{\Dom g}=\Delta(\theta,\phi)$.
    We have therefore a map from \ref{itm:cor:geod_tor_char:1} to \ref{itm:cor:geod_tor_char:2}.

    It remains to argue that this map is the converse of the preceding map from \ref{itm:cor:geod_tor_char:2} to \ref{itm:cor:geod_tor_char:1}. The non-trivial point is to verify \eqref{eq:geod_ray_G} holds for the $g$ we just constructed. Fix $t>0$. We need to show that
    \[
    G_{\ell_t}=G_{\phi}+tg.
    \]
    This holds on $\Int \{G_{\phi}<\infty\}$ by \cref{cor:toricgeodgeneral} and the definition of $g$, but then it holds everywhere thanks to \cref{prop:closedcvx_domain}. Our assertion follows.
\end{proof}

In particular, we can make the corresponding test curve more explicit.
\begin{corollary}\label{cor:RWN_toric}
    Let $\ell\in \mathcal{R}_{\tor}(X,\theta;\phi)$ and $g$ be the function as in \cref{cor:geod_tor_char}. Then
    \begin{equation}\label{eq:ellmax_temp1}
    \ell^*_{\max}=-\inf_{M_{\mathbb{R}}} g,
    \end{equation}
    and for each $\tau<-\inf_{M_{\mathbb{R}}} g$, the function $\ell^*_{\tau}$ is toric-invariant, and
    \[
    G_{\ell^*_{\tau}}(m)=\begin{cases}
        G_{\phi}(m),&\textup{if }g(m)\leq -\tau;\\
        \infty,&\textup{otherwise}.
    \end{cases}
    \]
    In particular, for such $\tau$,
    \[
    \Delta\left(\theta, \ell^*_{\tau}\right)=\left\{g\leq -\tau\right\}.
    \]
\end{corollary}
In other words,
\[
\left(\Delta\left(\theta, \ell^*_{\tau}\right)\right)_{\tau<-\inf_{M_{\mathbb{R}}} g}
\]
is the inverse Legendre transform of $-g$ using the terminology of \cref{def:LegOkoun}.

As a consequence, in the toric setting, the Ross--Witt Nyström correspondence \cref{thm:Legenbij} reduces essentially to \cref{thm:Okotestcurve}.
\begin{proof}
    Fix $\tau\in \mathbb{R}$. Then
    \[
    \ell^*_{\tau}=\inf_{t>0} \left(\ell_t-t\tau\right)
    \]
    is clearly toric-invariant.
    Therefore,
    \[
    F_{\ell^*_{\tau}}=\inf_{t>0} \left(F_{\ell_t}-t\tau\right).
    \]
    Fix $m\in M_{\mathbb{R}}$. We compute
    \[
    \begin{aligned}
    G_{\ell^*_{\tau}}(m)=& \sup_{n\in N_{\mathbb{R}}}\left( \langle m, n\rangle -F_{\ell^*_{\tau}}(n) \right)\\
    =& \sup_{t>0}\sup_{n\in N_{\mathbb{R}}}\left( \langle m, n\rangle - F_{\ell_t}(n) +t\tau\right)\\
    =&\sup_{t>0} \left(G_{\ell_t}(m)+t\tau \right)\\
    =& G_{\phi}(m)+\sup_{t>0} t\left(g(m)+\tau\right)\\
    =& \begin{cases}
        G_{\phi}(m),&\textup{if }g(m)+\tau\leq 0;\\
        \infty,&\textup{otherwise}.
    \end{cases}
    \end{aligned}
    \]
    Our assertions follow.
\end{proof}

Next we consider the trace operator studied in \cref{chap:trace}. We wish to understand the trace operator in the toric situation.
For this purpose, we will need to fix a $T$-invariant subvariety $Y\subseteq X$. Let $\sigma$ be the corresponding cone in $\Sigma$, and let $Q_{\sigma}\subseteq P_D$ be the subset defined in \eqref{eq:Qsigma}.
The cocharacter lattice of $Y$ is given by
\[
N(\sigma)\coloneqq N/N\cap \langle \sigma \rangle,
\]
where $\langle \sigma \rangle$ is the linear span of $\sigma$. See \cite[(3.2.6)]{CLS11}. In particular, we have a canonical identification of the character lattice $M(\sigma)$ of $Y$:
\[
M(\sigma)=\sigma^{\perp}\cap M,
\]
which is compatible with our previous notation \eqref{eq:Msigmadef}.
Let $i_{\sigma}\colon M(\sigma)\rightarrow M$ be the inclusion map.
Let $T_Y$ be the torus of $Y$. Then we have a natural surjection $q_T\colon T\rightarrow T_Y$. In particular, the tropicalization map
\[
\Trop\colon T(\mathbb{C})\rightarrow N_{\mathbb{R}}
\]
descends to the tropicalization map of $Y$:
\[
\Trop_Y\colon T_Y(\mathbb{C})\rightarrow N(\sigma)_{\mathbb{R}}.
\]
Choose a representative $m_{\sigma}\in M$ of the Cartier datum along $\sigma$. Then $m-m_{\sigma}\in M(\sigma)_{\mathbb{R}}$ for every $m\in Q_{\sigma}$.
We let
\[
D_Y=\sum_{\substack{\rho\in \Sigma(1)\\ \rho\not\preceq\sigma}} \left(a_{\rho}+\langle m_{\sigma},u_{\rho}\rangle\right)D_{\rho}|_Y,
\]
where $\rho\not\preceq\sigma$ means that $\rho$ is not a face of $\sigma$. Then $\mathcal{O}_Y(D_Y)=L|_Y$.
Equivalently, $D_Y$ is the divisor associated to the rational section
\[
s_{D_Y}\coloneqq (\chi^{m_{\sigma}}s_D)|_Y
\]
of $L|_Y$. Write $P_{D_Y}$ for the polytope of $D_Y$. With this choice, $Q_{\sigma}-m_{\sigma}$ is naturally contained in $P_{D_Y}$.

\begin{theorem}\label{thm:tracetor}
With the above choice of representative $m_{\sigma}$, let $\varphi\in \PSH_{\tor}(X,\theta)_{>0}$ satisfy $\nu(\varphi,Y)=0$. If $\Tr_Y^{\theta}(\varphi)$ is defined and $\vol (\theta|_Y,\Tr_Y^{\theta}(\varphi))>0$\footnote{Note that $\Tr_Y^{\theta}(\varphi)\in \PSH_{\tor}(Y,\theta|_Y)$.}, then
\[
\Delta\left(\theta|_Y,\Tr_Y^{\theta}(\varphi)\right)= \left(\Delta(\theta,\varphi)\cap Q_{\sigma}\right)-m_{\sigma}
\]
as subsets of $M(\sigma)_{\mathbb{R}}$.
\end{theorem}
Observe that the condition $\nu(\varphi,Y)=0$ means exactly that $\Delta(\theta,\varphi)\cap Q_{\sigma}\neq \varnothing$ by \cref{cor:toricLelong2}.

Since $Y$ itself is a smooth toric variety, the preceding constructions of $X$ all apply to $Y$. We briefly summarize the situation in \cref{Tab:corrXY}.
\begin{table}[ht]
\begin{center}
\begin{tabular}{|c|c|}
  \hline
Notions for $X$ & Notions for $Y$ \\
  \hline
 $N$ & $N(\sigma)$ \\
  \hline
 $M$ & $M(\sigma)$ \\
  \hline
 $\Sigma$ & $\Star(\sigma)$ \\
  \hline
$D$ & $D_Y$ \\
  \hline
  $L$ & $L|_Y$ \\
  \hline
  $h$ & $h|_Y$ \\
  \hline
  $\theta$ & $\theta|_Y$ \\
  \hline
  $\Trop$ & $\Trop_Y$ \\
  \hline
 $P_D$ & $P_{D_Y}$            \\
  \hline
  $s_D$ & $s_{D_Y}$            \\
  \hline
\end{tabular}
\caption{The correspondence between $X$ and $Y$}
\label{Tab:corrXY}
\end{center}
\end{table}

Recall that $\Star(\sigma)$ is the fan in $N(\sigma)_{\mathbb{R}}$ consisting of $\bar{\tau}$ for all faces $\tau\in \Sigma$ containing $\sigma$, where $\bar{\tau}$ is the image of $\tau$ in $N(\sigma)_{\mathbb{R}}$. See \cite[Proposition~3.2.7]{CLS11}.

\begin{proof}
    The idea of the proof is that since \cref{lma:slicepob,rmk:slicepob} describe how partial Okounkov bodies behave under restriction, and \cref{thm:toricpob} compares partial Okounkov bodies with Newton bodies, we can deduce the behavior of Newton bodies under restriction as well.

    First we note that by our assumption, $L|_Y$ is a big line bundle. In particular, if we set $r=\dim \sigma$, then $\dim Y=n-r$.

    We also note that $\Tr_Y^\theta(\varphi)$ is $T_{Y,c}$-invariant. Indeed, for any $g\in T_c$, let $g_Y$ be the induced automorphism of $Y$. Then $g_Y^*\Tr_Y(\varphi)$ is again a trace of $g^*\varphi=\varphi$ along $Y$. Hence $g_Y^*\Tr_Y(\varphi)\sim_P\Tr_Y(\varphi)$ by \cref{prop:traceunique}. Since $\theta|_Y$ is $T_{Y,c}$-invariant and the $P$-envelope depends only on the $P$-equivalence class, we get $g_Y^*\Tr_Y^\theta(\varphi)=\Tr_Y^\theta(\varphi)$. As the image of $T_c$ in $\Aut(Y)$ is $T_{Y,c}$, the assertion follows.

    For this purpose, let $\sigma^0$ be an $n$-dimensional cone of $\Sigma$ containing $\sigma$. The image $\overline{\sigma^0}$ in $N(\sigma)$ is then an $(n-r)$-dimensional cone of $\Star(\sigma)$. We shall use these cones as the reference cones while defining the partial Okounkov bodies.

    We list the rays in $\sigma^0(1)$ as follows:
    \begin{equation}\label{eq:rayssigma0}
        \rho_1,\ldots,\rho_n,
    \end{equation}
    where $\rho_1,\ldots,\rho_r\in \sigma(1)$
    and hence $\rho_{r+1},\ldots,\rho_n\not\in \sigma(1)$. In particular, the images
    \begin{equation}\label{eq:rayssigma0bar}
        \overline{\rho_{r+1}},\ldots,\overline{\rho_{n}}
    \end{equation}
    of the latter give a list of $\overline{\sigma^0}(1)$.

    We construct the flag $Y_{\bullet}$ on $X$ using the rays \eqref{eq:rayssigma0} and the flag $Z_{\bullet}$ on $Y$ using the rays \eqref{eq:rayssigma0bar}. Note that
    \[
        Z_i=Y_{r+i},
    \]
    where $i=1,\ldots,n-r$.

    Next we compute the Cartier data associated with $\overline{\sigma^0}$ for the divisor $D_Y$. By definition, $m_{\overline{\sigma^0}}\in M(\sigma)$ is the unique element satisfying
    \[
        \langle m_{\overline{\sigma^0}},u_{\overline{\rho_j}}\rangle=-a_{\rho_j}-\langle m_{\sigma},u_{\rho_j}\rangle
    \]
    for all $j=r+1,\ldots,n$.
    Since $m_{\sigma^0}-m_{\sigma}$ vanishes on $\sigma$, it lies in $M(\sigma)$, and the preceding formula gives
    \[
        m_{\overline{\sigma^0}}=m_{\sigma^0}-m_{\sigma}.
    \]

    Let $\Phi\colon M\rightarrow \mathbb{Z}^n$ and $\Psi\colon M(\sigma)\rightarrow \mathbb{Z}^{n-r}$ be defined as
    \[
    \begin{aligned}
    \Phi(m)=&\left( \langle m-m_{\sigma^0},u_{\rho_1} \rangle,\dots,\langle m-m_{\sigma^0},u_{\rho_n} \rangle \right)\\
    =& \left( \langle m,u_{\rho_1}\rangle,\dots,\langle m,u_{\rho_n}\rangle \right) + (a_{\rho_1},\ldots,a_{\rho_n}),\\
    \Psi(m)=&\left( \langle m-m_{\overline{\sigma^0}},u_{\overline{\rho_{r+1}}} \rangle,\dots,\langle m-m_{\overline{\sigma^0}},u_{\overline{\rho_n}} \rangle \right)
    \end{aligned}
    \]
    Observe that for $i=r+1,\ldots,n$, we have
    \[
        u_{\overline{\rho_i}}=u_{\rho_i}\mod N\cap \langle \sigma \rangle,
    \]
    so
    \[
    \begin{aligned}
    \Psi(m)=&\left( \langle m,u_{\rho_{r+1}}  \rangle,\ldots, \langle m,u_{\rho_{n}}  \rangle\right)\\
    &+\left(a_{\rho_{r+1}}+\langle m_{\sigma},u_{\rho_{r+1}}\rangle,\ldots,a_{\rho_n}+\langle m_{\sigma},u_{\rho_n}\rangle\right)
    \end{aligned}
    \]
    for $m\in M(\sigma)$. Therefore, we have a commutative diagram
    \[
    \begin{tikzcd}
	m_{\sigma}+M(\sigma)_{\mathbb{R}} \arrow[r, "\Phi_{\mathbb{R}}"] \arrow[d, "{m\mapsto m-m_{\sigma}}"']                         & \mathbb{R}^n \arrow[d, "L"] \\
	M(\sigma)_{\mathbb{R}} \arrow[r, "\Psi_{\mathbb{R}}"] & \mathbb{R}^{n-r},
	\end{tikzcd}
    \]
    where $L\colon \mathbb{R}^n\rightarrow \mathbb{R}^{n-r}$ is the map
    \[
    (b_1,\ldots,b_n)\mapsto (b_{r+1},\ldots,b_{n}).
    \]

    By \cref{thm:toricpob}, we have
    \[
    \Phi_{\mathbb{R}} \Bigl( \Delta(\theta,\varphi) \Bigr)= \Delta_{Y_{\bullet}}(\theta,\varphi),\quad \Psi_{\mathbb{R}} \Bigl( \Delta(\theta|_Y,\Tr_{Y}^{\theta}(\varphi)) \Bigr)= \Delta_{Z_{\bullet}}(\theta|_Y,\Tr_{Y}^{\theta}(\varphi)).
    \]
    By \cref{cor:toricLelong2}, the condition $\nu(\varphi,Y)=0$ is equivalent to
    \[
        \Delta(\theta,\varphi)\cap Q_\sigma\neq\varnothing.
    \]
    In particular, the first $r$ valuation coordinates of $\theta+\ddc\varphi$ along the flag $Y_\bullet$ vanish.
    By \cref{lma:slicepob} and \cref{rmk:slicepob},
    \[
    \begin{aligned}
    \Delta_{Z_{\bullet}}\Bigl(\theta|_Y,\Tr_{Y}^{\theta}(\varphi)\Bigr)=& L \Bigl(  \Delta_{Y_{\bullet}}(\theta,\varphi) \cap (\{0\}^r \times \mathbb{R}^{n-r})\Bigr)\\
    =& L\Bigl( \Phi_{\mathbb{R}}(\Delta(\theta,\varphi)) \cap (\{0\}^r \times \mathbb{R}^{n-r}) \Bigr)\\
    =& L\circ \Phi_{\mathbb{R}} \Bigl(\Delta(\theta,\varphi)\cap Q_{\sigma} \Bigr)\\
    =& \Psi_{\mathbb{R}}\Bigl((\Delta(\theta,\varphi)\cap Q_{\sigma})-m_{\sigma}\Bigr).
    \end{aligned}
    \]
    In the third line, we used $\Delta(\theta,\varphi)\subseteq P_D$ from \cref{thm:FVtheta}, the identity \eqref{eq:Qsigma}, and the fact that $m_{\sigma^0}-m_\sigma$ vanishes on $\sigma$.
    Since $\Psi_{\mathbb{R}}$ is an affine isomorphism, our assertion follows.
\end{proof}

\begin{corollary}\label{cor:tor_tr_vol}
    For any $\varphi\in \PSH_{\tor}(X,\theta)$, we have
    \begin{equation}\label{eq:volTr_tor}
    \vol \left(\theta|_Y+\ddc \Tr_Y^{\theta}(\varphi)\right)=(\dim Y)!\vol_{\dim Y} \Bigl(\Delta(\theta,\varphi)\cap Q_{\sigma}\Bigr).
    \end{equation}
\end{corollary}
The left-hand side of \eqref{eq:volTr_tor} is understood as $0$ if $\Tr_Y^{\theta}(\varphi)$ is not defined. On the right-hand side, we translate $\Delta(\theta,\varphi)\cap Q_{\sigma}$ by $-m_{\sigma}$ and compute the $\dim Y$-dimensional Lebesgue measure on $M(\sigma)_{\mathbb{R}}$, normalized so that a fundamental parallelepiped for $M(\sigma)$ has volume $1$.

\begin{proof}
    If $\nu(\varphi,Y)>0$, then by \cref{cor:toricLelong2}, $\Delta(\theta,\varphi)\cap Q_\sigma=\varnothing$. By convention the right-hand side of \eqref{eq:volTr_tor} is $0$, and the trace is not defined, so the left-hand side is $0$.

    Assume now that $\nu(\varphi,Y)=0$.
    Take a toric-invariant ample line bundle $H$ on $X$ and a toric-invariant Kähler metric $\omega\in c_1(H)$. For $\epsilon\in \mathbb{Q}_{>0}$, write $Q_{\sigma,\epsilon}$ for the face corresponding to $\sigma$ for the divisor $D+\epsilon H$, with the same convention on the Cartier datum along $\sigma$. Since $\theta_{\varphi}+\epsilon\omega$ is a toric Kähler current and $\nu(\varphi,Y)=0$, its trace along $Y$ is defined and has positive volume. Hence \cref{thm:tracetor} gives
    \[
    \vol \left(\theta|_Y+\epsilon\omega|_Y+\ddc \Tr_Y^{\theta+\epsilon\omega}(\varphi)\right)= (\dim Y)!\vol_{\dim Y} \Bigl(\Delta(\theta+\epsilon\omega,\varphi)\cap Q_{\sigma,\epsilon}\Bigr).
    \]
    By \cref{lma:grad_img_sum},
    \[
        \Delta(\theta+\epsilon\omega,\varphi)=\Delta(\theta,\varphi)+\epsilon P_H.
    \]
    We choose the Cartier datum of $D+\epsilon H$ along $\sigma$ to be the sum of the chosen datum of $D$ and $\epsilon$ times the chosen datum of $H$. With this convention, the face of $P_D+\epsilon P_H$ corresponding to $\sigma$ is $Q_\sigma+\epsilon Q_\sigma^H$. In particular,
    \[
        \Delta(\theta+\epsilon\omega,\varphi)\cap Q_{\sigma,\epsilon}=(\Delta(\theta,\varphi)\cap Q_{\sigma})+\epsilon Q^H_{\sigma},
    \]
    where $Q^H_{\sigma}$ denotes the corresponding face of the polytope of $H$. Indeed, the defining inequalities along the rays of $\sigma$ are non-negative on both summands, so equality for their sum forces equality on each summand. By \cref{ex:tracedefinedposmass}, the left-hand side converges to
    \[
        \vol\left(\theta|_Y+\ddc\Tr_Y^\theta(\varphi)\right),
    \]
    with the convention that this limit is $0$ when the trace in the limiting class is not defined. The right-hand side converges to
    \[
        (\dim Y)!\vol_{\dim Y}\Bigl(\Delta(\theta,\varphi)\cap Q_\sigma\Bigr).
    \]
    Hence \eqref{eq:volTr_tor} follows.
\end{proof}

As a corollary, we obtain an elegant characterization of the augmented base locus in the toric setting.
\begin{corollary}\label{cor:augbl}
    The augmented base locus (with the reduced complex structure) of $D$ is a toric-invariant (possibly reducible) subvariety given by the union of the $V(\tau)$'s, where $\tau$ runs over the elements in $\Sigma$ such that
    \[
    \dim Q_{\tau}<n-\dim \tau\footnotemark.
    \]
    \footnotetext{Here we understand that $\dim \varnothing=-\infty$.}
\end{corollary}
We call $n-\dim\tau$ the expected dimension of $Q_\tau$. An orbit closure $V(\tau)$ can also lie in the augmented base locus by being contained in a larger orbit closure whose face has deficient dimension.

\begin{proof}
Set $d=\dim V(\tau)$. The image of the restriction map of sections of $kD$ to $V(\tau)$ has dimension $\#(kQ_\tau\cap M)$. Indeed, as in \cref{prop:baseloc}, the monomial section indexed by $m\in kP_D\cap M$ restricts nontrivially precisely when its orders along the divisors indexed by $\tau(1)$ all vanish. On the dense orbit of $V(\tau)$ these sections have distinct characters, after translating the indices by $km_\tau$ into $M(\tau)$, so their restrictions are linearly independent.

Thus the restricted volume of $D$ along $V(\tau)$ vanishes exactly when $\dim Q_\tau<d$: if the face is nonempty, the lattice count in sufficiently divisible degrees grows with order $k^{\dim Q_\tau}$. By \cite[Corollary~5.9]{ELMNP09}, every orbit closure of zero restricted volume is contained in the augmented base locus, and every irreducible component of the augmented base locus has zero restricted volume. Since the augmented base locus is torus invariant, each of its irreducible components is an orbit closure. The claimed union follows.
\end{proof}

\begin{example}\label{ex:P2bl_nK}
    Let us consider the example \cref{ex:P2bl1pt} again, we consider the subvariety corresponding to $\sigma_7$. We have
    \[
    Q_{\sigma_7}=\Bigl\{m\in P_D: \langle m,(1,1) \rangle=0 \Bigr\}.
    \]
    The situation is explained in \cref{fig:Qsigma7}.
    \begin{figure}[ht]
        \centering
        \begin{tikzpicture}[scale=1.6, line width=0.7pt, every node/.style={font=\small}]
          \draw[->, thick] (-1.4,0) -- (1.6,0) node[right] {$m_1$};
          \draw[->, thick] (0,-1.4) -- (0,1.6) node[above] {$m_2$};
          \fill[red!10] (0,0) -- (1,0) -- (0,1) -- cycle;
          \draw[red!75!black, thick] (0,0) -- (1,0) -- (0,1) -- cycle;
          \node[red!75!black, anchor=south west, inner sep=2pt] at (0.62,0.42) {$P_D$};
          \draw[blue!70!black, thick] (-1.2,1.2) -- (1.2,-1.2);
          \node[blue!70!black, anchor=west, inner sep=2pt] at (1.30,-1.30) {$\langle m,(1,1)\rangle=0$};
          \fill[blue!70!black] (0,0) circle (2pt);
          \node[blue!70!black, anchor=south east, inner sep=2pt] at (0,0) {$Q_{\sigma_7}$};
        \end{tikzpicture}
        \caption{The face $Q_{\sigma_7}=\{m\in P_D:\langle m,(1,1)\rangle=0\}$ is the single point $\{(0,0)\}$, of dimension $0$ rather than the expected $1$.}
        \label{fig:Qsigma7}
    \end{figure}
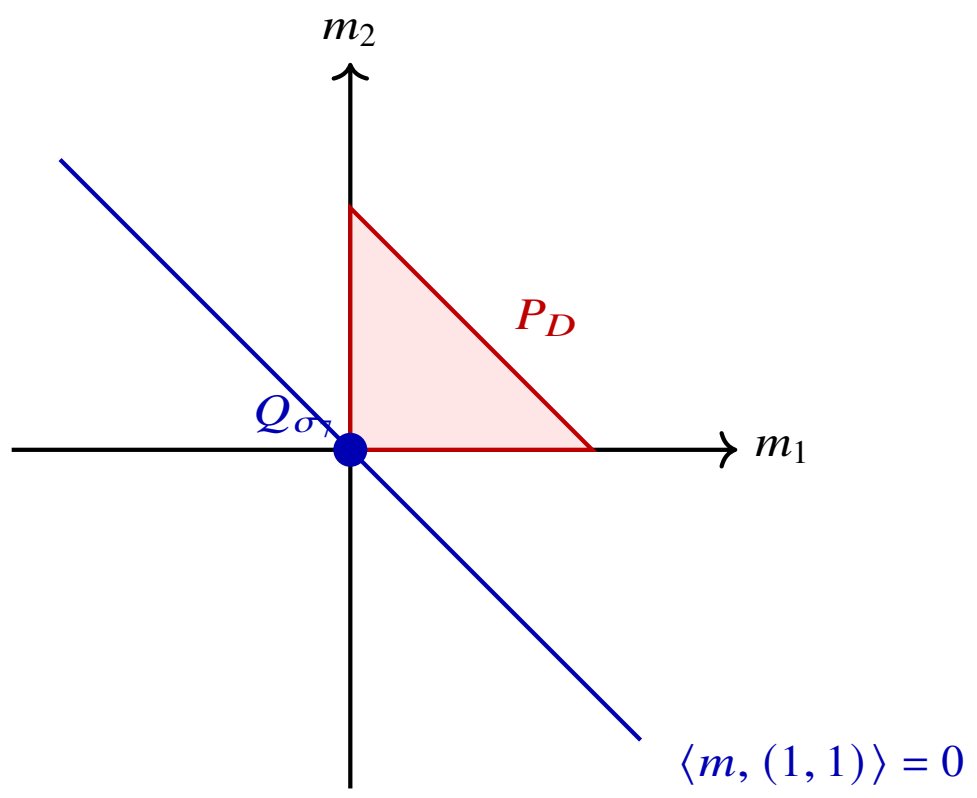
    We find that $\dim Q_{\sigma_7}=0$, but the expected dimension is $2-1=1$. So the corresponding divisor, namely the exceptional divisor on $\Bl_0\mathbb{P}^2$ is in the non-Kähler locus. Similarly, we can verify that the non-Kähler locus consists only of this divisor.
\end{example}

Finally, we compute the mixed volumes of Newton bodies.
\begin{theorem}\label{thm:mix_vol_newbod}
Take $T$-invariant big divisors $D_1,\ldots,D_n$ on $X$. Consider $T_c$-invariant closed positive $(1,1)$-currents $T_i$ in $c_1(\mathcal{O}_X(D_i))$ for $i=1,\ldots,n$. Then
\begin{equation}\label{eq:mixvol_pob}
n! \vol \Bigl( \Delta(T_1),\ldots,\Delta(T_n) \Bigr)=\vol \left(T_1,\ldots,T_n \right).
\end{equation}
\end{theorem}
The mixed volume of currents is studied in \cref{sec:mix_vol}.
\begin{proof}
    It suffices to prove the following assertion:
    \begin{equation}\label{eq:Delta_add}
    \Delta(T_1+T_2)=\Delta(T_1)+\Delta(T_2).\footnotemark
    \end{equation}
    \footnotetext{This is slightly problematic when $n=1$, since $T_2$ is not defined \emph{stricto sensu}. We take $D_2,T_2$ satisfying the same assumptions as $D_1,T_1$ in this case.}
    In fact, assuming this, by induction, for each tuple of non-negative integers $C_1,\ldots,C_n$, we have
    \[
    \Delta\left( \sum_{i=1}^n C_i T_i \right)=\sum_{i=1}^n C_i \Delta(T_i).
    \]
    Comparing the volumes, we have
    \[
    \vol \Delta\left( \sum_{i=1}^n C_i T_i \right)=\sum_{i_1,\ldots,i_n=1}^n C_{i_1}\cdots C_{i_n} \vol \left(\Delta\left(T_{i_1}\right),\ldots,\Delta\left(T_{i_n}\right)\right).
    \]
    On the other hand, by \cref{prop:toricMAandrealMA2}, the $\mathcal{I}$-goodness of toric singularities \cref{ex:toricIgood}, and \cref{prop:vollinearlimit}, we have
    \[
    n! \vol \Delta\left( \sum_{i=1}^n C_i T_i \right)=\vol \left( \sum_{i=1}^n C_i T_i \right)=\sum_{i_1,\ldots,i_n=1}^n C_{i_1}\cdots C_{i_n} \vol \left(T_{i_1},\ldots,T_{i_n}\right).
    \]
    Comparing the coefficients of the above two expressions, we conclude \eqref{eq:mixvol_pob}.

    It only remains to establish \eqref{eq:Delta_add}. For this purpose, we fix smooth $T_c$-invariant Hermitian metrics $h_1$ and $h_2$ on $\mathcal{O}_X(D_1)$ and $\mathcal{O}_X(D_2)$ respectively. Correspondingly, take $h_1\otimes h_2$ as the reference Hermitian metric on $\mathcal{O}_X(D_1+D_2)$. Write
    \[
    \theta_1=\ddc h_1,\quad \theta_2=\ddc h_2.
    \]
    We choose the rational section of $\mathcal O_X(D_1+D_2)$ to be $s_{D_1}\otimes s_{D_2}$. With this choice and with the reference metric $h_1\otimes h_2$, the normalization \eqref{eq:normF0} gives
    \[
        F_{\theta_1+\theta_2}=F_{\theta_1}+F_{\theta_2}.
    \]
    Thus
    \[
        F_{\varphi_1+\varphi_2}=F_{\varphi_1}+F_{\varphi_2}.
    \]
    Take $\varphi_i\in \PSH_{\tor}(X,\theta_i)$ with $T_i=\theta_i+\ddc \varphi_i$ for $i=1,2$. Then $F_{\varphi_1}+F_{\varphi_2}$ is a candidate for $F_{\varphi_1+\varphi_2}$; see \eqref{eq:Fvarphibig1}. Therefore, \eqref{eq:Delta_add} follows from the equality of closed gradient images proved in \cref{lma:grad_img_sum}.
\end{proof}

\motto{A good theorem lasts forever. Once proved, it will always stay proved, and other mathematicians are free to use it and build on it as they please, sometimes to great effect.

--- John Tate\footnote{John Torrence Tate Jr. (1925--2019), the grandfather of Dustin Clausen, was one of the greatest minds in the whole history of America. However, his aversion to publishing papers arguably impeded the progress of the development of mathematics to some extent. For example, his foundational work on rigid non-Archimedean geometry was written in 1962, but was not available to the public until 1971.}}
\chapter{Non-Archimedean pluripotential theory}
\label{chap:NApotential}

This chapter develops a non-Archimedean pluripotential theory for compact K\"ahler manifolds. In the projective setting, Boucksom and Jonsson constructed a space of non-Archimedean psh metrics on the Berkovich analytification of a polarized variety. Our approach starts instead from the analytic objects developed in the previous chapters: $\mathcal I$-model test curves and their operations. The resulting theory agrees with the Boucksom--Jonsson theory in the projective case.

The guiding idea is the Ross--Witt Nystr\"om correspondence. A geodesic ray can be encoded by a test curve, and the $\mathcal I$-model version of this construction retains precisely the multiplier-ideal information relevant to algebraic approximation. For a pseudo-effective class $\{\theta\}$, we define a non-Archimedean potential as an inverse system of $\mathcal I$-model test curves in the big classes $\{\theta+\omega\}$, where $\omega$ ranges over K\"ahler forms. The transition maps are given by the $\mathcal I$-envelope.

After the definition, we develop the basic operations on non-Archimedean potentials. Sums, translations, maxima, suprema, infima and rescalings are induced from the corresponding operations on test curves. We also define the non-Archimedean volume and prove its basic monotonicity and continuity properties.

The next part of the chapter constructs the Duistermaat--Heckman measure of a non-Archimedean potential based on the theory of Okounkov bodies developed in \cref{chap:Okounkov}. It is defined through Okounkov test curves and is shown to be independent of the chosen flag.

The final part compares this construction with the Boucksom--Jonsson theory. In the projective case, we define the analytification of our non-Archimedean potentials and prove that the resulting space is naturally identified with the Boucksom--Jonsson space of psh non-Archimedean metrics.

A related approach has been developed independently by Mesquita-Piccione \cite{MP24}, where a Berkovich-like analytification of a compact K\"ahler manifold is constructed. The present chapter is logically based on the test-curve formalism of \cref{chap:testcurve}, especially the $\mathcal I$-model theory of \cref{sec:Imodeltc} and the operations of \cref{sec:operationtc}. It is largely independent of \cref{chap:toric_big,chap:Bergman}, while the comparison with Boucksom--Jonsson uses the non-Archimedean theory recalled in \cref{sec:compBJ}.

\section{The definition of non-Archimedean metrics}\label{sec:nAm_def}
Let $X$ be a connected compact K\"ahler manifold of dimension $n$.
Let $\Kah(X)$ be the set of K\"ahler forms on $X$ with the partial order given as follows: We say $\omega\preceq \omega'$ if $\omega\geq \omega'$. Note that the partially ordered set $\Kah(X)$ is a directed set.

Let $\theta$ be a closed smooth real $(1,1)$-form representing a pseudo-effective cohomology class.
\begin{definition}\label{def:NApotential:L24}
    We define
    \[
        \PSH^{\NA}(X,\theta)=\varprojlim_{\omega\in \Kah(X)} \PSH^{\NA}(X,\theta+\omega)_{>0}\footnotemark
    \]
    \index{non-Archimedean metric}
    \index{$\PSH^{\NA}(X,\theta)$} \footnotetext{The cumbersome projective limit can be avoided if instead of relying on the language of quasi-plurisubharmonic functions, we use that of augmented nef b-divisors developed in \cref{def:augmented_nefb} instead. But given the fact that these two formulations are completely equivalent to each other by \cref{cor:nef_b_degenerate_2}, we just stick to the slightly more traditional language here.}
    in the category of sets, where the transition maps are given as follows: Suppose that $\omega,\omega'\in \Kah$ and $\omega\geq \omega'$. Then the transition map is defined in \cref{prop:transitionPI}:
    \begin{equation}\label{eq:PItransPSHNApositive}
        P_{\theta+\omega}[\bullet]_{\mathcal{I}}\colon \PSH^{\NA}(X,\theta+\omega')_{>0}\rightarrow \PSH^{\NA}(X,\theta+\omega)_{>0}.
    \end{equation}
    \index{$P_{\theta+\omega}[\bullet]_{\mathcal{I}}$}
    Note that the transition maps are compatible:
    \[
    P_{\theta+\omega''}[P_{\theta+\omega'}[\varphi]_\mathcal{I}]_\mathcal{I}=P_{\theta+\omega''}[\varphi]_\mathcal{I},\quad \forall \, \varphi\in \PSH^{\NA}(X,\theta+\omega)_{>0}
    \]
    for $\omega\leq\omega'\leq\omega''$, since both sides are $\mathcal{I}$-model and are $\mathcal{I}$-equivalent to $\varphi$.
\end{definition}
Recall that $\PSH^{\NA}(X,\theta)_{>0}$ is defined in \cref{def:Imodeltestcur}.

In general, when we denote an element in $\PSH^{\NA}(X,\theta)$ by $\Gamma$, its component in $\PSH^{\NA}(X,\theta+\omega)_{>0}$ ($\omega\in \Kah(X)$) will be written either as $\Gamma^{\theta+\omega}$ or as $P_{\theta+\omega}[\Gamma]_{\mathcal{I}}$.

Note that $\Gamma^{\theta+\omega}_{\max}$ is independent of the choice of $\omega\in \Kah(X)$. We denote this common value by $\Gamma_{\max}$\index{$\Gamma_{\max}$}.
\begin{remark}\label{rmk:NApotential:L47}
    Thanks to \cref{prop:Imodeltestcurveindeptheta}, for any other $\theta'$ representing $\{\theta\}$, we have a canonical bijection
    \[
    \PSH^{\NA}(X,\theta)\xrightarrow{\sim}\PSH^{\NA}(X,\theta').
    \]
    Moreover, these bijections satisfy the cocycle condition. If we view the set of closed real smooth $(1,1)$-forms representing $\{\theta\}$ as a category with a unique morphism between any two objects, then we can define
    \[
    \PSH^{\NA}(X,\{\theta\})=\varprojlim_{\theta} \PSH^{\NA}(X,\theta).
    \]
    This definition is independent of the choice of the explicit representative of the cohomology class $\{\theta\}$.

    However, given the fact that our notation is already quite heavy, we stick to the set $\PSH^{\NA}(X,\theta)$. The constructions below are independent of the choice of $\theta$ within its cohomology class.
\end{remark}

\begin{proposition}\label{prop:testcminftyPrela}
    Let $\Gamma\in \PSH^{\NA}(X,\theta)$. Then given $\omega,\omega'\in \Kah(X)$ with $\omega\geq \omega'$, we have
    \[
        P_{\theta+\omega}\left[ \Gamma^{\theta+\omega'}_{-\infty} \right]=P_{\theta+\omega}\left[ \Gamma^{\theta+\omega'}_{-\infty} \right]_{\mathcal{I}} = \Gamma^{\theta+\omega}_{-\infty}.
    \]
\end{proposition}
\begin{proof}
Since for any $\tau<\Gamma_{\max}$, the potential $\Gamma^{\theta+\omega'}_{\tau}$ is $\mathcal{I}$-good by \cref{ex:ImodelIgood}, it follows that
\[
P_{\theta+\omega}\left[ \Gamma^{\theta+\omega'}_{\tau} \right]=P_{\theta+\omega}\left[ \Gamma^{\theta+\omega'}_{\tau} \right]_{\mathcal{I}}=\Gamma^{\theta+\omega}_{\tau}
\]
for all $\tau<\Gamma_{\max}$.
Our assertion follows from \cref{prop:incnetmodel} and \cref{prop:incnetmodelPI}.
\end{proof}

\begin{proposition}\label{prop:NAposNAemb}
    There is a natural injective map
    \[
        \PSH^{\NA}(X,\theta)_{>0}\hookrightarrow \PSH^{\NA}(X,\theta),\quad \Gamma\mapsto \left(P_{\theta+\omega}[\Gamma]_{\mathcal{I}}\right)_{\omega\in \Kah(X)}.
    \]
\end{proposition}
In the sequel, we do not distinguish an element of $\PSH^{\NA}(X,\theta)_{>0}$ from its image in $\PSH^{\NA}(X,\theta)$. Note that given $\Gamma\in \PSH^{\NA}(X,\theta)_{>0}$, the value of $\Gamma_{\max}$ does not depend on whether we view it as an element of $\PSH^{\NA}(X,\theta)_{>0}$ or of $\PSH^{\NA}(X,\theta)$.
\begin{proof}
    The map is well-defined by the compatibility of the transition maps. It suffices to argue its injectivity. Suppose that $\Gamma,\Gamma'\in \PSH^{\NA}(X,\theta)_{>0}$ and
    \[
        P_{\theta+\omega}[\Gamma]_{\mathcal{I}}=P_{\theta+\omega}[\Gamma']_{\mathcal{I}}
    \]
    for some K\"ahler form $\omega$ on $X$. Then $\Gamma_{\max}=\Gamma'_{\max}$ and
    for any $\tau<\Gamma_{\max}$, we have
    \[
        \Gamma_{\tau}\sim_{\mathcal{I}}\Gamma'_{\tau}
    \]
    by \cref{prop:Icomparandenvelope}. It follows again from \cref{prop:Icomparandenvelope} that
    \[
    \Gamma_{\tau}=\Gamma'_{\tau}.
    \]
    This completes the proof.
\end{proof}

\begin{definition}\label{def:NApotential:L102}
    Let $\Gamma\in \PSH^{\NA}(X,\theta)$. We define its \emph{volume}\index{volume!volume of a non-Archimedean metric} as follows:
    \[
        \vol \Gamma\coloneqq \lim_{\omega\in \Kah(X)} \int_X \left(\theta+\omega+\ddc \Gamma^{\theta+\omega}_{-\infty}\right)^n\in [0,\infty).
    \]
    \index{$\vol \Gamma$}
\end{definition}
Observe that the net is decreasing, so the limit exists.

\begin{proposition}\label{prop:PSHNApvolGamma}
    Let $\Gamma\in \PSH^{\NA}(X,\theta)_{>0}$. Then
    \[
        \vol \Gamma=\int_X \left(\theta+\ddc \Gamma_{-\infty}\right)^n.
    \]
\end{proposition}
\begin{proof}
    This follows from \cref{prop:vol_limit_model}, \cref{cor:varphiperturbtheta} and \cref{prop:PpreGammainf}.
\end{proof}

\begin{lemma}\label{lma:NA_mass_drop_transition}
    Let $\Gamma\in \PSH^{\NA}(X,\theta)$, let $\omega$ be a Kähler form on $X$, and let $0<t\leq 1$. Then for every $\tau<\Gamma_{\max}$, we have
    \[
    \begin{aligned}
    &\int_X \left(\theta+t\omega+\ddc\Gamma^{\theta+t\omega}_{-\infty}\right)^n-\int_X \left(\theta+t\omega+\ddc\Gamma^{\theta+t\omega}_{\tau}\right)^n\\
    \leq{}&\int_X \left(\theta+\omega+\ddc\Gamma^{\theta+\omega}_{-\infty}\right)^n-\int_X \left(\theta+\omega+\ddc\Gamma^{\theta+\omega}_{\tau}\right)^n.
    \end{aligned}
    \]
\end{lemma}
\begin{proof}
    For $\sigma\in\{-\infty,\tau\}$, put
    \[
        S_t^\sigma=\theta+t\omega+\ddc\Gamma^{\theta+t\omega}_{\sigma}.
    \]
    By the transition compatibility and \cref{prop:testcminftyPrela,ex:ImodelIgood,prop:Ppresmass}, projecting $\Gamma^{\theta+t\omega}_{\sigma}$ to the class $\theta+\omega$ does not change the total mass computed in the larger class. Thus
    \[
    \int_X \left(\theta+\omega+\ddc\Gamma^{\theta+\omega}_{\sigma}\right)^n
    =
    \int_X \left(S_t^\sigma+(1-t)\omega\right)^n.
    \]
    Expanding the right-hand side and subtracting the two identities for $\sigma=-\infty$ and $\sigma=\tau$, we get
    \[
    \begin{aligned}
    &\int_X \left(\theta+\omega+\ddc\Gamma^{\theta+\omega}_{-\infty}\right)^n-\int_X \left(\theta+\omega+\ddc\Gamma^{\theta+\omega}_{\tau}\right)^n\\
    ={}&\sum_{a=0}^n \binom{n}{a}(1-t)^{n-a}
    \left(\int_X (S_t^{-\infty})^a\wedge\omega^{n-a}-\int_X (S_t^\tau)^a\wedge\omega^{n-a}\right).
    \end{aligned}
    \]
    Since $\Gamma^{\theta+t\omega}_{-\infty}\geq \Gamma^{\theta+t\omega}_{\tau}$, each summand is non-negative by the monotonicity theorem for mixed non-pluripolar products. Keeping only the term $a=n$ gives the desired inequality.
\end{proof}

\begin{lemma}\label{lma:positivemasspart}
The image of the canonical injection
\[
\PSH^{\NA}(X,\theta)_{>0}\hookrightarrow \PSH^{\NA}(X,\theta)
\]
is given by the set of $\Gamma\in \PSH^{\NA}(X,\theta)$ with positive volume.
\end{lemma}
\begin{proof}
By \cref{prop:PSHNApvolGamma}, every element in the image has positive volume.

Conversely, take $\Gamma\in \PSH^{\NA}(X,\theta)$ with positive volume. We want to construct $\Gamma'\in \PSH^{\NA}(X,\theta)_{>0}$ representing $\Gamma$.

Fix a K\"ahler form $\omega$ on $X$. Define
\begin{equation}\label{eq:Gamma'_tau}
    \Gamma'_{\tau}\coloneqq \lim_{k\to\infty}\Gamma^{\theta+k^{-1}\omega}_{\tau},\quad \tau<\Gamma_{\max}.
\end{equation}
We claim that it suffices to show
\begin{equation}\label{eq:theta_kinvomega_masslim_pos}
    \varliminf_{k\to\infty}\int_X \left(\theta+k^{-1}\omega+\ddc\Gamma^{\theta+k^{-1}\omega}_{\tau}\right)^n>0
\end{equation}
for some $\tau<\Gamma_{\max}$. If this holds, then the argument of \cref{lma:testcurvposmass} implies that the same holds for all $\tau<\Gamma_{\max}$. Then \cref{prop:vol_limit_model,prop:decnetmodelPI} guarantee that $\Gamma'\in \PSH^{\NA}(X,\theta)_{>0}$ and represents $\Gamma$.

It remains to argue \eqref{eq:theta_kinvomega_masslim_pos}.
Let $\epsilon=\vol \Gamma>0$.
Take $\tau<\Gamma_{\max}$ so that
\[
\int_X \left(\theta+\omega+\ddc\Gamma^{\theta+\omega}_{-\infty}\right)^n-\int_X \left(\theta+\omega+\ddc\Gamma^{\theta+\omega}_{\tau}\right)^n<\epsilon/2.
\]
Therefore, by \cref{lma:NA_mass_drop_transition}, for any $k\geq 1$,
\[
\int_X \left(\theta+k^{-1}\omega+\ddc\Gamma^{\theta+k^{-1}\omega}_{-\infty}\right)^n-\int_X \left(\theta+k^{-1}\omega+\ddc\Gamma^{\theta+k^{-1}\omega}_{\tau}\right)^n<\epsilon/2.
\]
By the definition of $\vol\Gamma$, the first integral in the preceding line tends to $\epsilon$ as $k\to\infty$. Hence
\[
    \varliminf_{k\to\infty}\int_X \left(\theta+k^{-1}\omega+\ddc\Gamma^{\theta+k^{-1}\omega}_{\tau}\right)^n\geq \epsilon/2>0,
\]
which proves \eqref{eq:theta_kinvomega_masslim_pos}.
\end{proof}

\begin{example}\label{ex:embedpshintona}
    Fix $\varphi\in \PSH(X,\theta)$. We can define an associated non-Archimedean metric
    \[
    \Gamma^{\varphi}\in \PSH^{\NA}(X,\theta)
    \]
    as follows:
    \begin{enumerate}
        \item\label{itm:ex:embedpshintona:1} $\Gamma^{\varphi}_{\max}=0$;
        \item\label{itm:ex:embedpshintona:2} for any $\omega\in \Kah(X)$ and any $\tau<0$, we set
        \[
        \Gamma^{\varphi,\theta+\omega}_{\tau}=P_{\theta+\omega}[\varphi]_{\mathcal{I}}.
        \]
    \end{enumerate}
    Such non-Archimedean metrics are called \emph{homogeneous}\index{homogeneous non-Archimedean metric} non-Archimedean metrics.

    Observe that
    \[
    \vol \Gamma^{\varphi}=\vol \theta_{\varphi}.
    \]
    See \cref{prop:vollinearlimit} and the footnote there.
\end{example}

\begin{definition}\label{def:PSHNAtrangeneral}
    Let $\omega$ be a closed real smooth positive $(1,1)$-form on $X$. We define the map
    \[
    P_{\theta+\omega}[\bullet]_{\mathcal{I}}\colon \PSH^{\NA}(X,\theta)\rightarrow \PSH^{\NA}(X,\theta+\omega)
    \]
    \index{$P_{\theta+\omega}[\bullet]_{\mathcal{I}}$}
    as follows: Given $\Gamma\in \PSH^{\NA}(X,\theta)$, we define $P_{\theta+\omega}[\Gamma]_{\mathcal{I}}$ as the element such that for any $\omega'\in \Kah(X)$, we have
    \[
    P_{\theta+\omega}[\Gamma]_{\mathcal{I}}^{\theta+\omega+\omega'}=\Gamma^{\theta+\omega+\omega'}.
    \]
\end{definition}
It is straightforward to check that under the identification of \cref{prop:NAposNAemb}, the map $P_{\theta+\omega}[\bullet]_{\mathcal{I}}$ extends the map \eqref{eq:PItransPSHNApositive}.

\begin{proposition}\label{prop:NApotential:L222}
    The maps $P_{\theta+\omega}[\bullet]_{\mathcal{I}}$ in \cref{def:PSHNAtrangeneral} together induce a bijection
    \begin{equation}\label{eq:PSHNAprojlimigeneral2}
    \PSH^{\NA}(X,\theta)\xrightarrow{\sim} \varprojlim_{\omega\in \Kah(X)}\PSH^{\NA}(X,\theta+\omega).
    \end{equation}
\end{proposition}
\begin{proof}
    It is a tautology that the maps $P_{\theta+\omega}[\bullet]_{\mathcal{I}}$ in \cref{def:PSHNAtrangeneral} are compatible with the transition maps. So the map \eqref{eq:PSHNAprojlimigeneral2} is well-defined. It is injective by the same argument as \cref{prop:NAposNAemb}. We argue the surjectivity.

    By unfolding the definitions, an object in the target of \eqref{eq:PSHNAprojlimigeneral2} is an assignment: With each $\omega\in \Kah(X)$, we associate a family $(\Gamma^{\omega,\omega'})_{\omega'\in \Kah(X)}$ satisfying:
    \begin{enumerate}
        \item\label{itm:chapter-NApotential:L242:1} $\Gamma^{\omega,\omega'}\in \PSH^{\NA}(X,\theta+\omega+\omega')_{>0}$ for each $\omega,\omega'\in \Kah(X)$;
        \item\label{itm:chapter-NApotential:L242:2} for each $\omega,\omega',\omega''\in \Kah(X)$ satisfying $\omega''\geq \omega'$, we have
        \[
        P_{\theta+\omega+\omega''}\left[\Gamma^{\omega,\omega'}\right]_{\mathcal{I}}=\Gamma^{\omega,\omega''};
        \]
        \item\label{itm:chapter-NApotential:L242:3} for each $\omega,\omega',\omega''\in \Kah(X)$ satisfying $\omega\leq \omega'$, we have
        \[
        P_{\theta+\omega'+\omega''}\left[\Gamma^{\omega,\omega''}\right]_{\mathcal{I}}=\Gamma^{\omega',\omega''}.
        \]
    \end{enumerate}
    The preimage of such an object is given by the family $(\Gamma^{\theta+\omega})_{\omega\in \Kah(X)}$ given by
    \[
        \Gamma^{\theta+\omega}=\Gamma^{\omega/2,\omega/2}.
    \]
    The fact that the image of $\Gamma$ is as expected follows directly from the definitions.
\end{proof}
With an almost identical argument involving \cref{prop:vol_limit_model}, we get
\begin{proposition}\label{prop:PSHNAreform1}
    The maps $P_{\theta+\omega}[\bullet]_{\mathcal{I}}$ in \cref{def:PSHNAtrangeneral} and the injective maps \cref{prop:NAposNAemb} together induce bijections
    \begin{equation}\label{eq:PSHNAprojlimigeneral}
    \PSH^{\NA}(X,\theta)\xrightarrow{\sim} \varprojlim_{\omega}\PSH^{\NA}(X,\theta+\omega)_{>0}\xrightarrow{\sim} \varprojlim_{\omega}\PSH^{\NA}(X,\theta+\omega),
    \end{equation}
    where $\omega$ runs over either the partially ordered set of all smooth closed real positive $(1,1)$-forms with positive volume\footnote{This partially ordered set is not directed.} on $X$ or $\Kah(X)$.
\end{proposition}

\begin{corollary}\label{cor:PSHNAbimero}
    Let $\pi\colon Y\rightarrow X$ be a proper bimeromorphic morphism from a compact K\"ahler manifold $Y$. Then $\pi^*$ induces a bijection
    \[
    \PSH^{\NA}(X,\theta)\xrightarrow{\sim} \PSH^{\NA}(Y,\pi^*\theta).
    \]
\end{corollary}
\begin{proof}
    This follows immediately from \cref{prop:PSHNAreform1}.
\end{proof}
It is immediate to verify that $\pi^*$ in \cref{cor:PSHNAbimero} extends the map \cref{prop:ETCIbimero}.

\section{Operations on non-Archimedean metrics}\label{sec:opnA}
Let $X$ be a connected compact K\"ahler manifold of dimension $n$ and $\theta,\theta',\theta''$ be closed real smooth $(1,1)$-forms on $X$ representing big cohomology classes.

This section relies heavily on \cref{sec:operationtc}. We shall use the notions introduced in that section without further explanations.

\begin{definition}\label{def:NApotential:L274}
    Let $\Gamma\in \PSH^{\NA}(X,\theta)$, $\Gamma'\in \PSH^{\NA}(X,\theta')$. We say $\Gamma\leq \Gamma'$ if for some $\omega\in \Kah(X)$, we have
    \index{non-Archimedean metric!partial order}
    \[
        \Gamma^{\theta+\omega} \leq \Gamma'{}^{\theta'+\omega}.
    \]
\end{definition}
This notion is independent of the choice of $\omega$ thanks to \cref{lma:testcurord1}.

Moreover, we have the following:
\begin{proposition}\label{prop:NApotential:L284}
    Let $\Gamma,\Gamma'\in \PSH^{\NA}(X,\theta)$ and $\omega$ be a closed smooth positive $(1,1)$-form on $X$. Then the following are equivalent:
    \begin{enumerate}
        \item\label{itm:chapter-NApotential:L296:1} $\Gamma\leq \Gamma'$;
        \item\label{itm:chapter-NApotential:L296:2} $P_{\theta+\omega}[\Gamma]_{\mathcal{I}}\leq P_{\theta+\omega}[\Gamma']_{\mathcal{I}}$.
    \end{enumerate}
\end{proposition}
\begin{proof}
    This follows immediately from \cref{lma:testcurord1}.
\end{proof}

Observe that this definition coincides with the corresponding definition in \cref{def:potestcurve} when $\Gamma,\Gamma'\in \PSH^{\NA}(X,\theta)_{>0}$.

\begin{proposition}\label{prop:vol_NAmetric_mono}
    Let $\Gamma,\Gamma'\in \PSH^{\NA}(X,\theta)$. Assume that $\Gamma\leq \Gamma'$. Then
    \[
    \vol \Gamma\leq \vol \Gamma'.
    \]
\end{proposition}
\begin{proof}
    It suffices to show that for any Kähler form $\omega$ on $X$, we have
    \[
    \vol \Gamma^{\theta+\omega}_{-\infty}\leq \vol \Gamma'{}^{\theta+\omega}_{-\infty},
    \]
    which is an immediate consequence of \cref{thm:mono}.
\end{proof}

\begin{definition}\label{def:sumNAmetrics}
    Let $\Gamma\in \PSH^{\NA}(X,\theta)$ and $\Gamma'\in \PSH^{\NA}(X,\theta')$. Then we define $\Gamma+\Gamma'\in \PSH^{\NA}(X,\theta+\theta')$ as the unique element such that for any $\omega\in \Kah(X)$, we have
    \index{non-Archimedean metric!sum}
    \[
    \left( \Gamma+\Gamma' \right)^{\theta+\theta'+2\omega}=\Gamma^{\theta+\omega}+\Gamma'^{\theta'+\omega}.
    \]
\end{definition}
This definition yields an element in $\PSH^{\NA}(X,\theta+\theta')$ by \cref{lma:testcursumcomp} and it extends the definition in \cref{def:sumtestcur} by \cref{lma:testcursumcomp} as well.

\begin{proposition}\label{prop:NApotential:L321}
    Let $\Gamma\in \PSH^{\NA}(X,\theta)$ and $\Gamma'\in \PSH^{\NA}(X,\theta')$. Suppose that $\omega,\omega'$ are two smooth closed positive $(1,1)$-forms on $X$. Then
    \[
    P_{\theta+\omega+\theta'+\omega'}[\Gamma+\Gamma']_{\mathcal{I}}=P_{\theta+\omega}[\Gamma]_{\mathcal{I}}+P_{\theta'+\omega'}[\Gamma']_{\mathcal{I}}.
    \]
\end{proposition}
\begin{proof}
    This is a direct consequence of \cref{lma:testcursumcomp}.
\end{proof}

\begin{proposition}\label{prop:NApotential:L331}
    The operation $+$ is commutative and associative: For any $\Gamma\in \PSH^{\NA}(X,\theta)$, $\Gamma'\in \PSH^{\NA}(X,\theta')$ and $\Gamma''\in \PSH^{\NA}(X,\theta'')$, we have
    \[
    \Gamma+\Gamma'=\Gamma'+\Gamma,\quad (\Gamma+\Gamma')+\Gamma''=\Gamma+(\Gamma'+\Gamma'').
    \]
\end{proposition}
\begin{proof}
    This is a direct consequence of \cref{prop:testcurvesumproperty}.
\end{proof}

\begin{definition}\label{def:NApotential:L341}
    Let $\Gamma\in \PSH^{\NA}(X,\theta)$ and $C\in \mathbb{R}$. We define $\Gamma+C\in \PSH^{\NA}(X,\theta)$ as the unique element such that for any $\omega\in \Kah(X)$, we have
    \index{non-Archimedean metric!translation}
    \[
    \left( \Gamma+C \right)^{\theta+\omega}=\Gamma^{\theta+\omega}+C.
    \]
\end{definition}
It is obvious from \cref{def:testcurvplusC} that $\Gamma+C\in \PSH^{\NA}(X,\theta)$. It is also obvious that this definition extends \cref{def:testcurvplusC}.

\begin{proposition}\label{prop:NApotential:L350}
    Let $\Gamma\in \PSH^{\NA}(X,\theta)$ and $C\in \mathbb{R}$. Suppose that $\omega$ is a smooth closed positive $(1,1)$-form on $X$. Then
    \[
    P_{\theta+\omega}[\Gamma]_{\mathcal{I}}+C=P_{\theta+\omega}[\Gamma+C]_{\mathcal{I}}.
    \]
\end{proposition}
\begin{proof}
    This is clear by definition.
\end{proof}

\begin{proposition}\label{prop:NAmetricplusC}
    Let $\Gamma\in \PSH^{\NA}(X,\theta)$, $\Gamma'\in \PSH^{\NA}(X,\theta')$ and $C,C'\in \mathbb{R}$. Then
    \begin{enumerate}
        \item\label{itm:prop:NAmetricplusC:1} $(\Gamma+\Gamma')+C=\Gamma+(\Gamma'+C)=(\Gamma+C)+\Gamma'$;
        \item\label{itm:prop:NAmetricplusC:2} $\Gamma+(C+C')=(\Gamma+C)+C'$.
    \end{enumerate}
\end{proposition}
\begin{proof}
    This is a direct consequence of \cref{prop:testcurveplusC}.
\end{proof}

\begin{proposition}\label{prop:vol_inv_plusC}
    Let $\Gamma\in \PSH^{\NA}(X,\theta)$ and $C\in \mathbb{R}$. Then
    \[
    \vol \Gamma=\vol (\Gamma+C).
    \]
\end{proposition}
\begin{proof}
    It suffices to show that for each Kähler form $\omega$ on $X$,
    \[
    \vol \Gamma^{\theta+\omega}_{-\infty}=\vol (\Gamma+C)^{\theta+\omega}_{-\infty},
    \]
    which is obvious.
\end{proof}

\begin{definition}\label{def:PSHNAlor}
    Let $\Gamma,\Gamma'\in \PSH^{\NA}(X,\theta)$, we define $\Gamma\lor \Gamma'\in \PSH^{\NA}(X,\theta)$ as the unique element such that for any $\omega\in \Kah(X)$, we have
    \index{non-Archimedean metric!maximum}
    \[
    \left( \Gamma\lor \Gamma' \right)^{\theta+\omega}=\Gamma^{\theta+\omega}\lor \Gamma'^{\theta+\omega}.
    \]
\end{definition}
It follows from \cref{lma:testcurlorPthetapomega} that $\Gamma\lor \Gamma'\in \PSH^{\NA}(X,\theta)$ and this definition extends the corresponding definition in \cref{def:testcurlor}.

\begin{proposition}\label{prop:NApotential:L394}
    Let $\Gamma,\Gamma'\in \PSH^{\NA}(X,\theta)$ and $\omega$ be a closed smooth positive $(1,1)$-form on $X$. Then
    \[
    P_{\theta+\omega}[\Gamma\lor \Gamma']_{\mathcal{I}}=P_{\theta+\omega}[\Gamma]_{\mathcal{I}}\lor P_{\theta+\omega}[\Gamma']_{\mathcal{I}}.
    \]
\end{proposition}
\begin{proof}
    This is a direct consequence of \cref{lma:testcurlorPthetapomega}.
\end{proof}

\begin{proposition}\label{prop:NApotential:L404}
    The operation $\lor$ is commutative and associative.
\end{proposition}
In particular, given a finite non-empty family $(\Gamma^i)_{i\in I}$ in $\PSH^{\NA}(X,\theta)$, we then define $\bigvee_{i\in I}\Gamma^i$ in the obvious way.
\begin{proof}
    This is a direct consequence of \cref{cor:testcurvlorprop}.
\end{proof}

\begin{definition}\label{def:NApotential:L413}
    Let $(\Gamma^i)_{i\in I}$ be a non-empty family in $\PSH^{\NA}(X,\theta)$. Assume that
    \index{non-Archimedean metric!usc supremum}
    \begin{equation}\label{eq:supPSHNAmaxfinite}
        \sup_{i\in I}\Gamma^i_{\max}<\infty.
    \end{equation}
    Then we define $\uscsup[i\in I]\Gamma^i\in \PSH^{\NA}(X,\theta)$ as the unique element such that for any $\omega\in \Kah(X)$, we have
    \[
        \left( \uscsup[i\in I]\Gamma^i \right)^{\theta+\omega}=\uscsup[i\in I]\Gamma^{i,\theta+\omega}.
    \]
\end{definition}
It follows immediately from \cref{lma:testcursupcompatible} that $\uscsup[i\in I]\Gamma^i\in \PSH^{\NA}(X,\theta)$ and this definition extends \cref{def:testcurvsupsgeneral}.
Moreover, this definition clearly extends \cref{def:PSHNAlor} as well.

\begin{proposition}\label{prop:NAsupproj}
    Let $(\Gamma^i)_{i\in I}$ be a non-empty family in $\PSH^{\NA}(X,\theta)$ satisfying \eqref{eq:supPSHNAmaxfinite}. Assume that $\omega$ is a closed smooth positive $(1,1)$-form on $X$. Then
    \[
    P_{\theta+\omega}\left[ \uscsup[i\in I] \Gamma^{i} \right]_{\mathcal{I}}= \uscsup[i\in I] P_{\theta+\omega}\left[ \Gamma^{i} \right]_{\mathcal{I}}.
    \]
\end{proposition}
\begin{proof}
    This is a direct consequence of \cref{lma:testcursupcompatible}.
\end{proof}

We also have a non-Archimedean version of Choquet's lemma.
\begin{proposition}\label{prop:NAChoquet}
    Let $(\Gamma^i)_{i\in I}$ be a non-empty family in $\PSH^{\NA}(X,\theta)$ satisfying \eqref{eq:supPSHNAmaxfinite}. Then there exists a countable subfamily $I'\subseteq I$ such that
    \[
        \uscsup[i\in I]\Gamma^i=\uscsup[i\in I']\Gamma^i.
    \]
\end{proposition}
\begin{proof}
    Fix a K\"ahler form $\omega_0$ on $X$. For each $m\geq 1$, \cref{prop:testcurvChoquet} gives a countable subfamily $I_m\subseteq I$ such that
    \[
    \uscsup[i\in I]P_{\theta+m^{-1}\omega_0}[\Gamma^i]_{\mathcal{I}}=\uscsup[i\in I_m]P_{\theta+m^{-1}\omega_0}[\Gamma^i]_{\mathcal{I}}.
    \]
    Set $I'=\bigcup_{m\geq 1}I_m$. This is countable. We claim that it works.

    Let $\omega'\in \Kah(X)$. Choose $m$ so large that $\omega'\geq m^{-1}\omega_0$. By \cref{prop:NAsupproj} and the compatibility of the transition maps,
    \[
    \begin{aligned}
    \uscsup[i\in I]P_{\theta+\omega'}[\Gamma^i]_{\mathcal{I}}
    =&P_{\theta+\omega'}\left[\uscsup[i\in I]P_{\theta+m^{-1}\omega_0}[\Gamma^i]_{\mathcal{I}}\right]_{\mathcal{I}}\\
    =&P_{\theta+\omega'}\left[\uscsup[i\in I']P_{\theta+m^{-1}\omega_0}[\Gamma^i]_{\mathcal{I}}\right]_{\mathcal{I}}\\
    =&\uscsup[i\in I']P_{\theta+\omega'}[\Gamma^i]_{\mathcal{I}}.
    \end{aligned}
    \]
    Thus the two non-Archimedean metrics have the same component in every class $\theta+\omega'$, and the assertion follows.
\end{proof}

\begin{proposition}\label{prop:supGammiotherprop2}
    Let $(\Gamma^i)_{i\in I}$ be a non-empty family in $\PSH^{\NA}(X,\theta)$ satisfying \eqref{eq:supPSHNAmaxfinite}. Let $C\in \mathbb{R}$. Then
    \[
    \uscsup[i\in I] (\Gamma^i+C)=\uscsup[i\in I] \Gamma^i+C.
    \]
    Suppose that $(\Gamma'^i)_{i\in I}$ is another family in $\PSH^{\NA}(X,\theta')$ satisfying \eqref{eq:supPSHNAmaxfinite}. Suppose that $\Gamma^i\leq \Gamma'^i$ for all $i\in I$. Then
    \[
    \uscsup[i\in I]\Gamma^i\leq \uscsup[i\in I]\Gamma'^i.
    \]
\end{proposition}
\begin{proof}
    This is an immediate consequence of \cref{prop:supGammiotherprop}.
\end{proof}

\begin{proposition}\label{prop:vol_NA_incnet_cont}
    Let $(\Gamma^i)_{i\in I}$ be an increasing net in $\PSH^{\NA}(X,\theta)$ satisfying \eqref{eq:supPSHNAmaxfinite}. Then
    \begin{equation}\label{eq:prop:vol_NA_incnet_cont}
    \vol \left( \uscsup[i\in I] \Gamma^i \right)=\lim_{i\in I} \vol \Gamma^i.
    \end{equation}
\end{proposition}
\begin{proof}
    The $\geq$ direction in \eqref{eq:prop:vol_NA_incnet_cont} is a direct consequence of \cref{prop:vol_NAmetric_mono}. It remains to prove the reverse inequality.

    Note that \eqref{eq:prop:vol_NA_incnet_cont} holds when $\vol \Gamma^i>0$ for each $i\in I$, as a consequence of \cref{prop:supsincnetteestcur}, \cref{cor:incnetdSconv} and \cref{thm:contvolu}.

    In particular, for each Kähler form $\omega$ on $X$, we have
    \[
    \vol \left( \uscsup[i\in I] \Gamma^{i,\theta+\omega} \right)=\lim_{i\in I} \vol \Gamma^{i,\theta+\omega}.
    \]
    For our purpose, we need to show that for any $\epsilon>0$, we can find $\omega$ so that
    \[
    \sup_{i\in I} \vol \Gamma^{i,\theta+\omega}<\sup_{i\in I} \vol \Gamma^{i}+\epsilon.
    \]
    We shall show that it is possible to choose $\omega$ so that the stronger statement holds:
    \[
    \vol \Gamma^{i,\theta+\omega}<\vol \Gamma^{i}+\epsilon,\quad \forall i\in I.
    \]
    Equivalently, we need to choose $\omega$ so that for any Kähler form $\omega'$ on $X$ dominated by $\omega$, we have
    \[
    \vol \Gamma^{i,\theta+\omega}_{-\infty}<\vol \Gamma^{i,\theta+\omega'}_{-\infty}+\epsilon/2,\quad \forall i\in I.
    \]
    Fix a K\"ahler form $\omega_0$ on $X$ and a K\"ahler form $\Omega$ such that $\Omega\geq \theta+\omega_0$. It is enough to take $\omega=\delta\omega_0$ for some $\delta\in (0,1)$. If $\omega'$ is dominated by $\omega$, then $\theta+\omega'\leq \Omega$. Set
    \[
        T_i\coloneqq \theta+\omega'+\ddc \Gamma^{i,\theta+\omega'}_{-\infty}.
    \]
    By the transition compatibility and \cref{prop:PpreGammainf}, the potential $\Gamma^{i,\theta+\omega}_{-\infty}$ has the same $P$-singularity type in the class $\theta+\omega$ as $\Gamma^{i,\theta+\omega'}_{-\infty}$. Hence its mass can be computed as
    \[
        \int_X\left(\theta+\omega+\ddc\Gamma^{i,\theta+\omega'}_{-\infty}\right)^n=\int_X \left(T_i+\omega-\omega'\right)^n.
    \]
    Therefore,
    \[
    \begin{aligned}
    \vol \Gamma^{i,\theta+\omega}_{-\infty}-\vol \Gamma^{i,\theta+\omega'}_{-\infty}
    =&\int_X \left(T_i+\omega-\omega'\right)^n-\int_X T_i^n\\
    =& \sum_{a=0}^{n-1} \binom{n}{a} \int_X T_i^a\wedge \left( \omega-\omega'\right)^{n-a}\\
    \leq &\sum_{a=0}^{n-1} \binom{n}{a} \int_X \Omega^a\wedge \omega^{n-a}.
    \end{aligned}
    \]
    In the last inequality we used the facts that $0\leq \omega-\omega'\leq \omega$ and $\theta+\omega'\leq \Omega$, and the monotonicity theorem \cref{thm:mono}.
    We can choose $\omega$ so that the right-hand side is less than $\epsilon/2$. Our assertion then follows.
\end{proof}

\begin{definition}\label{def:NApotential:L525}
    Let $(\Gamma^i)_{i\in I}$ be a decreasing net in $\PSH^{\NA}(X,\theta)$. Assume that
    \index{non-Archimedean metric!infimum of a decreasing net}
    \begin{equation}\label{eq:decnetcontition}
    \inf_{i\in I}\Gamma^i_{\max}>-\infty,
    \end{equation}
    then we define $\inf_{i\in I}\Gamma^i\in \PSH^{\NA}(X,\theta)$ as the unique element such that for each $\omega\in \Kah(X)$, we have
    \begin{equation}\label{eq:decnettestcurdef}
    \left( \inf_{i\in I}\Gamma^i \right)^{\theta+\omega}=\inf_{i\in I} \Gamma^{i,\theta+\omega}.
    \end{equation}
\end{definition}
We observe that
\[
\left( \inf_{i\in I}\Gamma^i \right)^{\theta+\omega}\in \PSH^{\NA}(X,\theta+\omega)_{>0}.
\]
This follows from \cref{prop:inftcistc}. Moreover, by \cref{lma:testcurinfcompatible}, we have $\inf_{i\in I}\Gamma^i\in \PSH^{\NA}(X,\theta)$, and this definition extends \cref{def:testcurinf}.

In general,
\[
\vol \left( \inf_{i\in I}\Gamma^i \right)\leq \lim_{i\in I} \vol \Gamma^i
\]
as a consequence of \cref{prop:vol_NAmetric_mono}. But the reverse inequality fails in general.

\begin{proposition}\label{prop:NApotential:L548}
    Let $(\Gamma^i)_{i\in I}$ be a decreasing net in $\PSH^{\NA}(X,\theta)$ satisfying \eqref{eq:decnetcontition}. Assume that $\omega$ is a closed smooth positive $(1,1)$-form on $X$. Then
    \[
    P_{\theta+\omega}\left[ \inf_{i\in I} \Gamma^{i} \right]_{\mathcal{I}}= \inf_{i\in I} P_{\theta+\omega}\left[ \Gamma^{i} \right]_{\mathcal{I}}.
    \]
\end{proposition}
\begin{proof}
    This follows from \cref{lma:testcurinfcompatible}.
\end{proof}

\begin{proposition}\label{prop:infGammiotherprop2}
    Let $(\Gamma^i)_{i\in I}$ be a decreasing net in $\PSH^{\NA}(X,\theta)$ satisfying \eqref{eq:decnetcontition}. Let $C\in \mathbb{R}$. Then
    \[
    \inf_{i\in I} (\Gamma^i+C)=\inf_{i\in I} \Gamma^i+C.
    \]
    Suppose that $(\Gamma'^i)_{i\in I}$ is another decreasing net in $\PSH^{\NA}(X,\theta')$ satisfying \eqref{eq:decnetcontition}. Suppose that $\Gamma^i\leq \Gamma'^i$ for all $i\in I$. Then
    \[
    \inf_{i\in I}\Gamma^i\leq \inf_{i\in I}\Gamma'^i.
    \]
\end{proposition}
\begin{proof}
    This is clear by definition.
\end{proof}

\begin{definition}\label{def:NApotential:L572}
    Let $\Gamma\in \PSH^{\NA}(X,\theta)$ and $\lambda\in \mathbb{R}_{>0}$. Then we define $\lambda \Gamma\in \PSH^{\NA}(X,\lambda\theta)$ as the unique element such that for any $\omega\in \Kah(X)$, we have
    \index{non-Archimedean metric!rescaling}
    \[
    \left( \lambda\Gamma \right)^{\lambda\theta+\omega}=\lambda \Gamma^{\theta+\lambda^{-1}\omega}.
    \]
\end{definition}
It follows immediately from \cref{lma:testcurvrescompatible} that $\lambda \Gamma\in \PSH^{\NA}(X,\lambda\theta)$ and this definition extends \cref{def:res}.

\begin{proposition}\label{prop:NApotential:L581}
    Let $\Gamma\in \PSH^{\NA}(X,\theta)$ and $\lambda\in \mathbb{R}_{>0}$. Then for any closed smooth positive $(1,1)$-form $\omega$ on $X$, we have
    \[
    P_{\lambda\theta+\omega}[\lambda\Gamma]_{\mathcal{I}}=\lambda P_{\theta+\lambda^{-1}\omega}[\Gamma]_{\mathcal{I}}.
    \]
\end{proposition}
\begin{proof}
    This follows immediately from \cref{lma:testcurvrescompatible}.
\end{proof}

\begin{proposition}\label{prop:resclacompat2}
    Let $\Gamma\in \PSH^{\NA}(X,\theta)$, $\Gamma'\in \PSH^{\NA}(X,\theta')$, $C\in \mathbb{R}$ and $\lambda,\lambda'>0$, we have
    \[
        \begin{aligned}
        \lambda(\Gamma+\Gamma')=&\lambda\Gamma+\lambda\Gamma',\\
        (\lambda\lambda')\Gamma=&\lambda(\lambda'\Gamma),\\
        \lambda(\Gamma+C)=&\lambda\Gamma+\lambda C.
        \end{aligned}
    \]
    Suppose that $(\Gamma^i)_{i\in I}$ is a non-empty family in $\PSH^{\NA}(X,\theta)$ satisfying \eqref{eq:supPSHNAmaxfinite}. Then
    \[
    \lambda \left( \uscsup[i\in I]\Gamma^i \right)=\uscsup[i\in I](\lambda \Gamma^i).
    \]
    If $(\Gamma^i)_{i\in I}$ is a decreasing net in $\PSH^{\NA}(X,\theta)$ satisfying \eqref{eq:decnetcontition}, then
    \[
    \lambda \left( \inf_{i\in I}\Gamma^i \right)=\inf_{i\in I}(\lambda \Gamma^i).
    \]
\end{proposition}
\begin{proof}
Everything except the last assertion follows from \cref{prop:resclacompat}. The last assertion follows directly from the definition.
\end{proof}

\begin{proposition}\label{prop:vol_NA_resc}
    Let $\Gamma\in \PSH^{\NA}(X,\theta)$ and $\lambda\in \mathbb{R}_{>0}$. Then
    \[
    \vol \left( \lambda \Gamma \right)=\lambda^n \vol \Gamma.
    \]
\end{proposition}
\begin{proof}
    This is clear by definition.
\end{proof}

\begin{definition}\label{def:NApotential:L623}
    Let $\Gamma\in \PSH^{\NA}(X,\theta)$. Let $Y\subseteq X$ be an irreducible analytic subset. We say that the trace operator of $\Gamma$ along $Y$ is \emph{well-defined} if
    \index{non-Archimedean metric!trace operator}
    \[
    \nu\left(\Gamma_{\tau}^{\theta+\omega},Y\right)=0
    \]
    for all sufficiently negative $\tau$ and every $\omega\in \Kah(X)$. We define
    \[
    \left(\Tr_Y(\Gamma)\right)_{\max}\coloneqq \sup\left\{ \tau<\Gamma_{\max}: \nu\left(\Gamma_{\tau}^{\theta+\omega},Y\right)=0\right\}.
    \]
    Here $\omega$ is any K\"ahler form on $X$; this number is independent of the choice of $\omega$ by the compatibility of the transition maps.

    In this case, we define $\Tr_Y(\Gamma)\in \PSH^{\NA}(\widetilde{Y},\theta|_{\widetilde{Y}})$\footnote{Here $\widetilde{Y}\rightarrow Y$ is the normalization of $Y$.} as the unique element such that for any $\omega\in \Kah(\widetilde{Y})$, the component
    \begin{equation}\label{eq:TrYinPSHNAg0}
    \Tr_Y(\Gamma)^{\theta|_{\widetilde{Y}}+\omega}\in \PSH^{\NA}(\widetilde{Y},\theta|_{\widetilde{Y}}+\omega)_{>0}
    \end{equation}
    is defined as follows:
    \begin{enumerate}
        \item\label{itm:eq:TrYinPSHNAg0:1} We let
        \begin{equation}\label{eq:tracemax}
        \left( \Tr_Y(\Gamma)^{\theta|_{\widetilde{Y}}+\omega} \right)_{\max}=\left(\Tr_Y(\Gamma)\right)_{\max};
        \end{equation}
        \item\label{itm:eq:TrYinPSHNAg0:2} for each $\tau\in \mathbb{R}$ less than $\left(\Tr_Y(\Gamma)\right)_{\max}$, we define
        \[
        \Tr_Y(\Gamma)^{\theta|_{\widetilde{Y}}+\omega}_{\tau}\coloneqq P_{\theta|_{\widetilde{Y}}+\omega}\left[ \Tr_Y^{\theta+\tilde{\omega}}\left(\Gamma^{\theta+\tilde{\omega}}_{\tau}\right) \right]_{\mathcal{I}},
        \]
        where $\tilde{\omega}$ is an arbitrary K\"ahler form on $X$ such that $\omega\geq \tilde{\omega}|_{\widetilde{Y}}$.
    \end{enumerate}
\end{definition}
It follows from \cite[Proposition~3.6]{GK20} that $\widetilde{Y}$ is a normal K\"ahler space and hence $\tilde{\omega}$ exists. We observe that the choice of the trace operator $\Tr_Y^{\theta+\tilde{\omega}}\left(\Gamma^{\theta+\tilde{\omega}}_{\tau}\right)$ is irrelevant since two different choices are $\mathcal{I}$-equivalent. Moreover, \eqref{eq:TrYinPSHNAg0} holds as a consequence of \cref{prop:traceunique} and \cref{prop:tracelinear}. It is therefore clear that $\Tr_Y(\Gamma)\in \PSH^{\NA}(\widetilde{Y},\theta|_{\widetilde{Y}})$. Here when $\widetilde{Y}$ is not smooth, we understand $\PSH^{\NA}(\widetilde{Y},\theta|_{\widetilde{Y}})$ as the set of non-Archimedean metrics on a resolution of $\widetilde{Y}$, which is independent of the choice of the resolution by \cref{prop:tracelinear}. See \cref{chap:unibranch} for more details.

\begin{proposition}\label{prop:NApotential:L655}
    Let $\pi\colon Y\rightarrow X$ be a proper bimeromorphic morphism from a compact K\"ahler manifold $Y$. Then all definitions in this section are invariant under pullback to $Y$.
\end{proposition}
The meaning is clear in most cases. In the case of the trace operator, this means the following: Suppose that $Z\subseteq X$ is an analytic subset and $\Gamma\in \PSH^{\NA}(X,\theta)$ has non-trivial restriction to $Z$. Suppose that $Z$ is not contained in the non-isomorphism locus of $\pi$ so that the strict transform $W$ of $Z$ is defined. If we write $\Pi\colon W\rightarrow Z$ for the restriction of $\pi$ and $\widetilde{\Pi}\colon \widetilde{W}\rightarrow \widetilde{Z}$ the strict transform of $\Pi$, then we have
\[
    \widetilde{\Pi}^* \Tr_Z(\Gamma)=\Tr_W(\pi^*\Gamma).
\]
The relevant notations are summarized in the following diagram:
\[
\begin{tikzcd}
\widetilde{W} \arrow[r] \arrow[d, "\widetilde{\Pi}"] & W \arrow[r] \arrow[d, "\Pi"] & Y \arrow[d, "\pi"] \\
\widetilde{Z} \arrow[r]                          & Z \arrow[r]                  & X.
\end{tikzcd}
\]
\begin{proof}
We only prove the assertion for the trace operator, as the other proofs are similar.

We shall use the notations above. Observe that for any closed positive smooth $(1,1)$-form $\omega$ on $X$ with positive mass, we have
\[
\left(\widetilde{\Pi}^* \Tr_Z(\Gamma)\right)_{\max}=\left( \Tr_Z(\Gamma)\right)_{\max}=\sup\left\{ \tau<\Gamma_{\max}: \nu\left(\Gamma^{\theta+\omega}_{\tau},Z\right)=0 \right\},
\]
and
\[
\begin{aligned}
\left(\Tr_W(\pi^*\Gamma)\right)_{\max}=&\sup\left\{ \tau<\Gamma_{\max}: \nu\left((\pi^*\Gamma)^{\pi^*\theta+\pi^*\omega}_{\tau},W\right)=0 \right\}\\
=&\sup\left\{ \tau<\Gamma_{\max}: \nu\left(\pi^*\Gamma^{\theta+\omega}_{\tau},W\right)=0 \right\}\\
=& \sup\left\{ \tau<\Gamma_{\max}: \nu\left(\Gamma^{\theta+\omega}_{\tau},Z\right)=0 \right\}.
\end{aligned}
\]
Here we implicitly used \cref{prop:PSHNAreform1}.
Therefore,
\[
\left(\widetilde{\Pi}^* \Tr_Z(\Gamma)\right)_{\max}=\left(\Tr_W(\pi^*\Gamma)\right)_{\max}.
\]
Let $\tau\in \mathbb{R}$ be less than this common value. Take a K\"ahler form $\omega_{\widetilde{Z}}$ (resp. $\omega_{\widetilde{W}}$) on $\widetilde{Z}$ (resp. $\widetilde{W}$).
Take a K\"ahler form $\omega_Y$ on $Y$ (resp. $\omega_X$ on $X$) such that
\[
\omega_{\widetilde{Z}}\geq \omega_X|_{\widetilde{Z}},\quad \omega_{\widetilde{W}} \geq \omega_Y|_{\widetilde{W}},\quad \omega_Y \geq \pi^*\omega_X.
\]
We want to show that, after projection to the common background $(\pi^*\theta)|_{\widetilde{W}}+\omega_{\widetilde{W}}$,
\[
P_{(\pi^*\theta)|_{\widetilde{W}}+\omega_{\widetilde{W}}}\left[\widetilde{\Pi}^*\left(\Tr_Z(\Gamma)^{\theta|_{\widetilde{Z}}+\omega_{\widetilde{Z}}}_\tau\right)\right]_{\mathcal I} \sim_P
\left(\Tr_W(\pi^*\Gamma)\right)^{(\pi^*\theta)|_{\widetilde{W}}+\omega_{\widetilde{W}}}_\tau.
\]

Unfolding the definitions, we reduce to
\[
\widetilde{\Pi}^* \Tr_Z^{\theta+\omega_X}\left( \Gamma^{\theta+\omega_X}_{\tau} \right)\sim_P \Tr_{W}^{\pi^*\theta+\omega_Y}\left( \left( \pi^*\Gamma\right)^{\pi^*\theta+\omega_Y}_{\tau}\right).
\]
Using \cref{prop:tracelinear}, this is equivalent to
\[
\widetilde{\Pi}^* \Tr_Z\left( \Gamma^{\theta+\omega_X}_{\tau} \right)\sim_P \Tr_{W}\left( \left( \pi^*\Gamma\right)^{\pi^*\theta+\pi^*\omega_X}_{\tau}\right).
\]
This is a consequence of \cref{lma:rescommpullback}.
\end{proof}

\section{Duistermaat--Heckman measures}\label{sec:DHmeasure}
Let $X$ be a connected compact K\"ahler manifold of dimension $n$ and $\theta$ be a closed real smooth $(1,1)$-form on $X$ representing a big cohomology class.

\begin{definition}\label{def:DHm}
Assume that $X$ admits a smooth flag $Y_{\bullet}$. Let $\Gamma\in \PSH^{\NA}(X,\theta)_{>0}$.
    The \emph{Duistermaat--Heckman measure}\index{Duistermaat--Heckman measure} $\DHm(\Gamma)$\index{$\DHm(\Gamma)$} of $\Gamma$ is defined as
    \[
    \DHm(\Gamma)\coloneqq n!\cdot \DHm\left(\Delta_{Y_{\bullet}}(\theta,\Gamma)\right).
    \]
\end{definition}
Recall that $\Delta_{Y_{\bullet}}(\theta,\Gamma)\in \TC(\Delta_{Y_{\bullet}}(\theta,\Gamma_{-\infty}))$ is the Okounkov test curve defined in \cref{thm:Okountescurvex}.
See \cref{def:DHmeasureOTC} for the definition of the Duistermaat--Heckman measure of an Okounkov test curve.

\begin{theorem}\label{thm:DHindep}
Assume that $X$ admits a smooth flag $Y_{\bullet}$.
    The Duistermaat--Heckman measure $\DHm(\Gamma)$ of $\Gamma\in \PSH^{\NA}(X,\theta)_{>0}$ in \cref{def:DHm} is independent of the choice of the smooth flag $Y_{\bullet}$. Furthermore, for any $m\in \mathbb{Z}_{>0}$, the following moment identity holds:
    \begin{equation}\label{eq:momentDHmtc1}
    \int_{\mathbb{R}}x^m\DHm(\Gamma)(x)=\Gamma_{\max}^m\vol \Gamma+m\int_{-\infty}^{\Gamma_{\max}}\tau^{m-1}\Bigl( \vol (\theta+\ddc \Gamma_{\tau})-\vol \Gamma \Bigr)\,\mathrm{d}\tau,
    \end{equation}
    and
    \begin{equation}\label{eq:momentDHmtc2}
        \int_{\mathbb{R}}\DHm(\Gamma)=\vol \Gamma.
    \end{equation}
\end{theorem}
For each $m\in \mathbb{N}$, $\int_{\mathbb{R}}x^m\DHm(\Gamma)(x)$ is indeed $n!$ times the $m$-th moment of the random variable $G[\Delta_{Y_{\bullet}}(\theta,\Gamma)]$, by \cref{def:DHm,def:DHmeasureOTC}.

\begin{proof}
We observe that \eqref{eq:momentDHmtc1} holds by the computations in \cref{prop:DHmoments} and \itemref{itm:thm:pobmain:1}. It is independent of the choice of the flag. Furthermore, \eqref{eq:momentDHmtc2} follows immediately from \itemref{itm:thm:pobmain:1}.

Assume first that $\Gamma$ is bounded.
Since the Duistermaat--Heckman measure has bounded support in this case (cf. \cref{thm:Okotestcurve}), we conclude that $\DHm(\Gamma)$ is uniquely determined.

By translating the parameter, it suffices to treat the case $\Gamma_{\max}=0$. Indeed, if $a=\Gamma_{\max}$ and
\[
    \widetilde\Gamma_\tau\coloneqq\Gamma_{\tau+a},
\]
then $\widetilde\Gamma_{\max}=0$, and the associated test function is translated by $a$. Consequently,
\[
    \DHm(\Gamma)=(x\mapsto x+a)_*\DHm(\widetilde\Gamma).
\]
Thus independence of the flag and the moment identity for $\widetilde\Gamma$ imply the corresponding assertions for $\Gamma$.

For each $\epsilon>0$, we define $\Gamma^{\epsilon}\in \PSH^{\NA}(X,\theta)_{>0}$ as follows:
\begin{enumerate}
        \item\label{itm:chapter-NApotential:L761:1} Let $\Gamma^{\epsilon}_{\max}=0$, and
        \item\label{itm:chapter-NApotential:L761:2} we set
        \[
        \Gamma^{\epsilon}_{\tau}=
        \begin{cases}
            \Gamma_{-\infty},&\text{ if }\tau\leq -\epsilon^{-1},\\
            P_{\theta}\left[ (1+\epsilon\tau)\Gamma_{\tau}-\epsilon\tau \Gamma_{-\infty} \right],&\text{ if }\tau\in \left(-\epsilon^{-1},0\right).
        \end{cases}
        \]
    \end{enumerate}
Let us check that $\Gamma^\epsilon$ is indeed an element of $\PSH^{\NA}(X,\theta)_{>0}$. The concavity and monotonicity in $\tau$ follow from \cref{prop:Pconc}. Since the $\Gamma_\tau$'s are $\mathcal I$-model, they are $\mathcal I$-good by \cref{ex:ImodelIgood}; hence
\[
(1+\epsilon\tau)\Gamma_\tau-\epsilon\tau\Gamma_{-\infty}
\]
is $\mathcal I$-good by \cref{prop:Igoodlinear}. Therefore its $P$-envelope is also its $\mathcal I$-envelope, and the corresponding level is $\mathcal I$-model. The value at $-\infty$ is
$\Gamma_{-\infty}$. Thus $\Gamma^\epsilon\in\PSH^{\NA}(X,\theta)_{>0}$.

For each fixed $\tau<0$, the argument of \cref{thm:Legenbij} Step~3.3 gives
\[
    \Gamma^\epsilon_\tau\searrow \Gamma_\tau \quad\textup{as }\epsilon\to 0+.
\]
By \cref{prop:vol_limit_model} and \cref{cor:decnetdS}, this convergence is in $d_S$. Hence \itemref{itm:thm:pobmain:4} gives
\[
    \Delta_{Y_\bullet}(\theta,\Gamma^\epsilon_\tau)\xrightarrow{d_{\Haus}}\Delta_{Y_\bullet}(\theta,\Gamma_\tau)
\]
for every $\tau<0$. Since
\[
    \Gamma^\epsilon_{-\infty}=\Gamma_{-\infty}
\]
for all $\epsilon$, the underlying convex body is fixed. Therefore \cref{lma:DHmconv} applies and gives
\[
    \DHm(\Gamma^\epsilon)\rightharpoonup\DHm(\Gamma).
\]
It follows that $\DHm(\Gamma)$ is independent of the choice of the flag.
\end{proof}

More generally, when $X$ does not admit a smooth flag, we may choose a modification $\pi\colon Y\rightarrow X$ so that $Y$ admits a flag. We define
\begin{equation}\label{eq:DHmgeneral}
\DHm(\Gamma)\coloneqq \DHm(\pi^*\Gamma).
\end{equation}
It follows from \itemref{itm:thm:pobmain:5} that this measure is independent of the choice of $\pi$.

\begin{proposition}\label{prop:contDH}
Let
\[
(\Gamma^i)_{i\geq 1}\subseteq \PSH^{\NA}(X,\theta)_{>0},\quad \Gamma\in \PSH^{\NA}(X,\theta)_{>0}.
\]
Assume one of the following conditions holds:
\begin{enumerate}
    \item\label{itm:prop:contDH:1} The sequence $(\Gamma^i)_{i\geq 1}$ is decreasing and $\Gamma=\inf_{i\geq 1}\Gamma^i$. Assume that
    \begin{equation}\label{eq:volGammalim}
    \vol \Gamma=\lim_{i\to\infty}\vol \Gamma^i.
    \end{equation}
    \item\label{itm:prop:contDH:2} The sequence $(\Gamma^i)_{i\geq 1}$ is increasing and $\Gamma=\uscsup[i\geq 1]\Gamma^i$.
\end{enumerate}
    Then
    \begin{equation}\label{eq:contDHm}
    \DHm(\Gamma^i)\rightharpoonup \DHm(\Gamma).
    \end{equation}
\end{proposition}
\begin{proof}
    We may assume that $X$ admits a smooth flag $Y_{\bullet}$.

    Assume \ref{itm:prop:contDH:1}. Let
\[
    \eta\coloneqq\inf_{i\geq 1}\Gamma^i_{-\infty}.
\]
Since $\Gamma_\tau\leq \Gamma^i_{-\infty}$ for every $i$ and every $\tau<\Gamma_{\max}$, we have
\[
    \Gamma_{-\infty}\leq \eta.
\]
By \cref{prop:vol_limit_model}, the decreasing limit $\eta$ is model and
\[
    \int_X\theta_\eta^n=\lim_{i\to\infty}\int_X\theta_{\Gamma^i_{-\infty}}^n=\lim_{i\to\infty}\vol\Gamma^i.
\]
Using \eqref{eq:volGammalim}, this equals
\[
    \vol\Gamma=\int_X\theta_{\Gamma_{-\infty}}^n.
\]
Since $\Gamma_{-\infty}\leq\eta$ and both potentials are model, \cref{prop:Pequivchar2} gives
\[
    \Gamma_{-\infty}=\eta.
\]
Thus
\[
    \Gamma_{-\infty}=\inf_{i\geq1}\Gamma^i_{-\infty}.
\]

We want to derive \eqref{eq:contDHm} from \cref{lma:DHmconv}.
    It boils down to prove the following: For any $\tau<\Gamma_{\max}$, we have
    \[
        \Delta_{Y_{\bullet}}(\theta,\Gamma^i_{\tau})\xrightarrow{d_{\Haus}} \Delta_{Y_{\bullet}}(\theta,\Gamma_{\tau}).
    \]
    For fixed $\tau<\Gamma_{\max}$, the potentials $\Gamma^i_\tau$ decrease to $\Gamma_\tau$. By \cref{prop:vol_limit_model}, their masses converge to the mass of $\Gamma_\tau$; hence $\Gamma^i_\tau\xrightarrow{d_S}\Gamma_\tau$ by \cref{cor:decnetdS}. Applying the continuity of the partial Okounkov body map, \itemref{itm:thm:pobmain:4}, gives
    \[
        \Delta_{Y_{\bullet}}(\theta,\Gamma^i_{\tau})\xrightarrow{d_{\Haus}}\Delta_{Y_{\bullet}}(\theta,\Gamma_{\tau}).
    \]

    The proof under the assumption \ref{itm:prop:contDH:2} is similar. We only need to apply \cref{lma:DHmconv2} instead of \cref{lma:DHmconv}.
\end{proof}

\begin{remark}\label{rmk:NApotential:L862}
    Let $\Gamma$ be induced by a test configuration $(\mathcal X,\mathcal L)$ of $(X,L)$ with $L$ ample, as in \cref{ex:RWN,rmk:PhongSturm}. The Duistermaat--Heckman measure of \cite[Section~3.2]{BHJ17} is normalized to have total mass one. By \eqref{eq:momentDHmtc2}, our measure has total mass $\vol\Gamma=(L^n)$, since $\Gamma_{-\infty}=V_\theta$. Thus $\DHm(\Gamma)$ is $(L^n)$ times their measure. See \cite[Remark~7.17]{Xia21} for the comparison via the induced filtration.
\end{remark}

\section{Comparison with Boucksom--Jonsson's theory}\label{sec:compBJ}

\subsection{A brief recap of Boucksom--Jonsson's theory}
\label{subsec:BJ_recap}

In this section, we briefly recall the non-Archimedean global pluripotential theory à la Boucksom--Jonsson \cite{BJ22}. Our presentation is selective, so the original paper remains the best companion for this section.

\subsubsection{Valuations}
\label{subsubsec:BJ_val}

Let $X$ be an irreducible reduced variety over $\mathbb{C}$ of dimension $n$. We recall the Berkovich analytification $X^{\An}$ of $X$ with respect to the trivial valuation on $\mathbb{C}$.

\begin{definition}\label{def:NApotential:L878}
    A (real-valued) \emph{valuation}\index{valuation} on $X$ (or a \emph{valuation} of $\mathbb{C}(X)$) is a map $v\colon \mathbb{C}(X)\rightarrow (-\infty,\infty]$ satisfying the following conditions:
\begin{enumerate}
    \item\label{itm:chapter-NApotential:L883:1} For $f\in \mathbb{C}(X)$, $v(f)=\infty$ if and only if $f=0$;
    \item\label{itm:chapter-NApotential:L883:2} For $f,g\in \mathbb{C}(X)$, $v(fg)=v(f)+v(g)$;
    \item\label{itm:chapter-NApotential:L883:3} For $f,g\in \mathbb{C}(X)$, $v(f+g)\geq v(f)\land v(g)$;
    \item\label{itm:chapter-NApotential:L883:4} For $c\in \mathbb{C}^{\times}$, $v(c)=0$.
\end{enumerate}
\end{definition}
The set of valuations on $X$ is denoted by $X^{\val}$. The center of a valuation $v$ is the scheme-theoretic point $c=c(v)$ of $X$ such that $v\geq 0$ on $\mathcal{O}_{X,c}$ and $v>0$ on the maximal ideal $\mathfrak{m}_{X,c}$ of $\mathcal{O}_{X,c}$. The center is unique if it exists, and it exists if $X$ is proper.

In the remainder of this section, we assume that $X$ is projective.

As a set, $X^{\An}$ is the set of \emph{semi-valuations} on $X$, in other words, real-valued valuations $v$ on irreducible reduced subvarieties $Y$ of $X$ that are trivial on $\mathbb{C}$. We call $Y$ the \emph{support} of the semi-valuation $v$.
In other words,
\[
X^{\An}=\coprod_{Y} Y^{\val}.
\]
We will write $v_{\triv}\in X^{\An}$ for the trivial valuation on $X$: $v_{\triv}(f)=0$ for any $f\in \mathbb{C}(X)^{\times}$.

We endow $X^{\An}$ with the coarsest topology such that
\begin{enumerate}
    \item\label{itm:chapter-NApotential:L902:1} for any Zariski open subset $U\subseteq X$, the subset $U^{\An}$ of $X^{\An}$ consisting of semi-valuations whose supports meet $U$ is open;
    \item\label{itm:chapter-NApotential:L902:2} for each Zariski open subset $U\subseteq X$ and each $f\in \mathrm{H}^0(U,\mathcal{O}_X)$ (here $\mathcal{O}_X$ is the sheaf of regular functions), the map $|f|\colon U^{\An}\rightarrow \mathbb{R}$ sending $v$ to $\exp(-v(f))$ is continuous.
\end{enumerate}
See \cite{Berk93} for more details.

We will be most interested in divisorial valuations. Recall that a \emph{divisorial valuation} on $X$ is a valuation of the form $t\ord_E$, where $t\in \mathbb{Q}_{>0}$ and $E$ is a prime divisor over $X$. The set consisting of all divisorial valuations on $X$ and the trivial valuation $v_{\triv}$ is denoted by $X^{\Div}$. When $\mathbb{Q}_{>0}$ is replaced by $\mathbb{R}_{>0}$, we can similarly define a space $X^{\Div}_{\mathbb{R}}$.

Given any coherent ideal $\mathfrak{a}$ on $X$ and any $v\in X^{\An}$, we define
\begin{equation}\label{eq:va}
v(\mathfrak{a})\coloneqq\min\{v(f):f\in \mathfrak{a}_{c(v)}\}\in [0,\infty],
\end{equation}
where $c(v)$ is the center of the valuation $v$ on $X$.

Given any valuation $v$ on $X$, the Gauss extension of $v$ is a valuation $\sigma(v)$ on $X\times \mathbb{A}^1$:
\begin{equation}\label{eq:Gext}
\sigma(v)\left(\sum_i f_i t^i\right)\coloneqq\min_i (v(f_i)+i).
\end{equation}
Here $t$ is the standard coordinate on $\mathbb{A}^1=\Spec \mathbb{C}[t]$. The key property is that when $v$ is a divisorial valuation, then so is $\sigma(v)$. See \cite[Lemma~4.2]{BHJ17}.

\subsubsection{Non-Archimedean plurisubharmonic functions}
\label{subsubsec:BJ_NApsh}

Let $X$ be an irreducible unibranch complex projective variety of dimension $n$ and $L$ be a holomorphic pseudo-effective $\mathbb{Q}$-line bundle on $X$. Through the GAGA morphism $X^{\An}\rightarrow X$ of ringed spaces, $L$ can be pulled back to an analytic line bundle $L^{\An}$ on $X$. The purpose of this section is to study the psh metrics on $L^{\An}$. We will follow the approach of \cite{BJ22}, which avoids the direct treatment of $L^{\An}$ itself.

\enlargethispage{16pt}
Following \cite[Definition~2.18]{BJ22}, we define $\mathcal{H}^{\gf}_{\mathbb{Q}}(L^{\An})$, the set of \emph{(rational) generically finite Fubini--Study functions} $\phi\colon X^{\An}\rightarrow [-\infty,\infty)$ that are of the following form:
    \begin{equation}
    \phi=\frac{1}{m}\max_j \{\log|s_j|+\lambda_j\}.
    \end{equation}
Here $m\in \mathbb{Z}_{>0}$ is an integer such that $L^m$ is a line bundle, the $s_j$'s are a finite collection of non-zero sections in $\mathrm{H}^0(X,L^m)$, and $\lambda_j\in \mathbb{Q}$. We follow the convention of Boucksom--Jonsson by writing $\log |s_j|(v)=-v(s_j)$.

\begin{definition}[{\cite[Definition~4.1]{BJ22}}]\label{def:BJpshmetric}
A \emph{plurisubharmonic metric}\index{plurisubharmonic metric} (or \emph{psh metric} for short) on $L^{\An}$ is a function $\phi\colon X^{\An}\rightarrow [-\infty,\infty)$ that is not identically $-\infty$, and is the pointwise limit of a decreasing net $(\phi_i)_{i\in I}$, where $\phi_i\in \mathcal{H}^{\gf}_{\mathbb{Q}}(L_i^{\An})$ for some $\mathbb{Q}$-line bundles $L_i$ on $X$ satisfying $c_1(L_i)\to c_1(L)$ in $\NS^1(X)_{\mathbb{R}}$.
\end{definition}

The set of psh metrics on $L^{\An}$ is denoted by $\PSH(L^{\An})$.
We endow $\PSH(L^{\An})$ with the topology of pointwise convergence on $X^{\Div}$. This topology is Hausdorff as functions in $\PSH(L^{\An})$ are completely determined by their restriction on $X^{\Div}$:
\begin{theorem}[{\cite[Theorem~4.22]{BJ22}}]\label{thm:pshdetbydiv}
Let $\phi \in\PSH(L^{\An})$ and $\psi:X^{\An}\rightarrow [-\infty,\infty)$ be an usc function. Assume that $\phi\leq \psi$ on $X^{\Div}$. Then the same holds on $X^{\An}$.
\end{theorem}

\begin{proposition}[{\cite[Theorem~4.7]{BJ22}}]\label{prop:BJclosedmax}
    Let $\phi,\phi'\in \PSH(L^{\An})$. Then so is their pointwise maximum $\phi\lor \phi'$.
\end{proposition}

\begin{proposition}\label{prop:BJdecseqstab}
    Let $H$ be an ample line bundle on $X$. Consider $\phi\in \PSH((L+H)^{\An})$.
    Assume that for each $m\in \mathbb{Z}_{>0}$, we have $\phi\in \PSH((L+m^{-1}H)^{\An})$. Then $\phi\in \PSH(L^{\An})$.
\end{proposition}
This is a special case of the decreasing stability property (PSH2), proved for unibranch $X$ in \cite[Corollary~8.4]{FHX26}.

Next we note that we may use sequences instead of nets in the definition of $\PSH(L^{\An})$:
\begin{theorem}[{\cite[Corollary~12.18]{BJ22}}]\label{thm:BJappseq}
Let $S$ be an ample line bundle on $X$. Let $\phi \in \PSH(L^{\An})$.
Then there is a sequence of rational numbers $\varepsilon_i\searrow 0$ and a decreasing sequence $\phi_i\in \mathcal{H}^{\gf}_{\mathbb{Q}}((L+\varepsilon_i S)^{\An})$ such that $\phi$ is the pointwise limit of $\phi_i$, as $i\to\infty$.
\end{theorem}

The space $\PSH(L^{\An})$ inherits most of the expected properties of (Archimedean) psh functions (\cite[Theorem~4.7]{BJ22}). Fang--Han--Xia proved the following compactness result:
\begin{theorem}[{\cite[Theorem~A and Proposition~2.3]{FHX26}}]\label{thm:NApshincreasing}
Assume that $X$ is unibranch. Then every bounded from above increasing net of elements in $\PSH(L^{\An})$ converges in $\PSH(L^{\An})$.
\end{theorem}
This is equivalent to the envelope property, formulated in \cite[Conjecture~5.14]{BJ22}: the regularized supremum of a nonempty bounded from above family of functions in $\PSH(L^{\An})$ lies in $\PSH(L^{\An})$. See \cite[Theorem~5.11]{BJ22} for the equivalence. The case when $X$ is smooth and $L$ is nef was proved in \cite{BJ22}. Boucksom--Jonsson subsequently established the case when $X$ is smooth and $L$ is pseudo-effective in \cite{BJ22b}.

\subsection{The analytifications}
\label{subsec:analytifications}
Let $X$ be a connected projective manifold of dimension $n$. Let $\theta$ be a closed smooth real $(1,1)$-form on $X$ representing a pseudo-effective cohomology class.

\subsubsection{The transcendental setting}
\label{subsubsec:analytic_trans}

\begin{definition}\label{def:NApotential:L975}
    For $\varphi\in \PSH(X,\theta)$, we define the \emph{analytification}\index{analytification} $\varphi^{\An}\colon X^{\An}\rightarrow [-\infty,0]$ as follows:
\begin{equation}\label{eq:varphianv}
\varphi^{\An}(v)\coloneqq -v(\varphi)=-\lim_{k\to\infty}\frac{1}{k}v\left(\mathcal{I}(k\varphi) \right).
\end{equation}
\end{definition}
By \cref{thm:multipsubadd} and Fekete's lemma, the limit in \eqref{eq:varphianv} exists.

Note that we can also write
\begin{equation}\label{eq:varphiAnalt}
\varphi^{\An}(v)=\inf_{k\in \mathbb{Z}_{>0}} -2^{-k}v\left( \mathcal{I}(2^k\varphi)\right).
\end{equation}

When $v=t\ord_E$ for some prime divisor $E$ over $X$, $\varphi^{\An}(v)=-t\nu(\varphi,E)$ by \cref{prop:Lelongnumfrommis}.

\begin{definition}\label{def:Gammaana}
    Let $\Gamma\in \PSH^{\NA}(X,\theta)$. We define the \emph{analytification}\index{analytification} $\Gamma^{\An}\colon X^{\Div}\rightarrow [-\infty,\infty)$ of $\Gamma$ as follows: For any $\omega\in \Kah(X)$, we define
    \begin{equation}\label{eq:psiandef}
        \Gamma^{\An}(v)\coloneqq \sup_{\tau\leq \Gamma_{\max}}\left(\Gamma_{\tau}^{\theta+\omega,\An}(v)+\tau\right).
    \end{equation}
\end{definition}
Clearly, \eqref{eq:psiandef} is independent of the choice of $\omega$.

Note that \eqref{eq:psiandef} can be equivalently written as
\[
\Gamma^{\An}(v)= \sup_{\tau\leq \Gamma_{\max}}\left(\Gamma_{\tau}^{\theta+\omega,\An}(v)+\tau\right)=\sup_{\tau\in \mathbb{R}}\left(\Gamma_{\tau}^{\theta+\omega,\An}(v)+\tau\right)
\]
with $(-\infty)^{\An}(v) = -\infty$ understood.

\begin{proposition}\label{prop:psiandPhi}
Let $\Gamma\in \PSH^{\NA}(X,\theta)_{>0}$ with $\Gamma_{\max}\leq 0$. Let $\Psi$ be the complexification of $\Gamma^*$. Then
\begin{equation}\label{eq:psiansigma}
\Gamma^{\An}(v)=-\sigma(v)(\Psi)\quad \forall v\in X^{\Div}.
\end{equation}
\end{proposition}
See \cref{def:subgeoray} for the definition of the complexification $\Psi\in \QPSH(X\times \Delta)$. Note that since $\Gamma_{\max}\leq 0$, by \cref{cor:reltestcursuplinear} and \cref{thm:GRexten}, $\Psi$ extends uniquely to a quasi-psh function on $X\times \Delta$.

\begin{proof}
Recall that
\[
\Psi(x,\delta)=\sup_{\tau\leq \Gamma_{\max}} (\Gamma_{\tau}(x)-\log|\delta|^2\tau)\quad \text{for }x\in X, \delta\in \Delta^*.
\]
By \eqref{eq:Gext}, we have $\sigma(v)(\log |\delta|^2)=1$ and $\sigma(v)(\Gamma_{\tau})=v(\Gamma_{\tau})$ for all $\tau\leq \Gamma_{\max}$. So we have
\[
\sigma(v)(\Gamma_\tau(x) - \log|\delta|^2 \tau) = v(\Gamma_\tau) - \tau.
\]
Lastly, since $\sigma(v)$ is a divisorial valuation on $X \times \Delta$, by \cref{cor:supsLelong}, we conclude \eqref{eq:psiansigma}.
\end{proof}

\begin{definition}\label{def:piecewiselincurv}
Let $N\in \mathbb{N}$, and $A_{0},\ldots,A_N$ be a finite collection of elements in $\PSH(X,\theta)$, and $\tau_0>\tau_1>\cdots>\tau_N$ be finitely many real numbers. Then the \emph{piecewise linear curve}\index{piecewise linear curve} $A=(A_{\tau})_{\tau\in \mathbb{R}}$ in $\PSH(X,\theta)\cup \{-\infty\}$ associated with these data is the affine interpolation of these data:
\begin{enumerate}
    \item\label{itm:def:piecewiselincurv:1} $A_{\tau_i}=A_i$ for $i=0,\ldots,N$;
    \item\label{itm:def:piecewiselincurv:2} $A_{\tau}=A_{\tau_N}$ for $\tau\leq \tau_N$;
    \item\label{itm:def:piecewiselincurv:3} for $t\in (0,1)$ and $i=0,\ldots,N-1$, we have
    \[
        A_{(1-t)\tau_i+t\tau_{i+1}}=(1-t)A_{\tau_i}+tA_{\tau_{i+1}};
    \]
    \item\label{itm:def:piecewiselincurv:4} $A_{\tau}\equiv -\infty$ for $\tau>\tau_0$.
\end{enumerate}
The \emph{analytification}\index{analytification} of $A$ is the function $A^{\An}\colon X^{\An}\rightarrow [-\infty,\infty)$ defined as follows:
\begin{equation}\label{eq:psianvlinear}
A^{\An}(v)\coloneqq \sup_{\tau\leq \tau_0} (A_\tau^{\An}(v)+\tau)=\max_{i=0,\ldots,N} \left(A_{\tau_i}^{\An}(v)+\tau_i\right)\quad \forall v\in X^{\An}.
\end{equation}
\end{definition}
We also say $A=(A_{\tau})_{\tau\leq \tau_0}$ is a piecewise linear curve in $\PSH(X,\theta)$.

\begin{remark}\label{rmk:pwlinhull}
Note that $\tau\mapsto A_{\tau}$ is upper semicontinuous, but not necessarily concave.
Let $A'_{\tau}$ be its pointwise upper concave envelope, followed by upper semicontinuous regularization in $X$ for each $\tau$. For $\tau\leq\tau_0$, set
\[
\Lambda_{\tau}\coloneqq\left\{(\lambda_0,\ldots,\lambda_N)\in[0,1]^{N+1}:\ \sum_{i=0}^N\lambda_i=1,\quad \sum_{i=0}^N\lambda_i\tau_i\geq\tau\right\}.
\]
Then
\begin{equation}\label{eq:pwlinhull-simplex}
A'_{\tau}=\max_{\lambda\in\Lambda_{\tau}}\sum_{i=0}^N\lambda_i A_i,\quad \tau\leq\tau_0,
\end{equation}
where summands with $\lambda_i=0$ are omitted, and $A'_{\tau}\equiv-\infty$ for $\tau>\tau_0$.
Let $e_0,\ldots,e_N$ denote the standard basis of $\mathbb R^{N+1}$. The maximum in \eqref{eq:pwlinhull-simplex} can be taken over the finitely many vertices of $\Lambda_{\tau}$. These are the vectors $e_i$ with $\tau_i\geq\tau$, together with
\[
\frac{\tau-\tau_j}{\tau_i-\tau_j}e_i+\frac{\tau_i-\tau}{\tau_i-\tau_j}e_j,
\quad i<j,\quad \tau_j<\tau<\tau_i.
\]
Thus $A'_{\tau}\in\PSH(X,\theta)$, and in particular
\[
A'_{\tau_0}=A_0,\qquad A'_{\tau}=\max_{0\leq i\leq N}A_i\quad\text{for }\tau\leq\tau_N.
\]

To justify the formula, denote its right-hand side by $B_{\tau}$. If $\sigma,\tau\leq\tau_0$ and $0<t<1$, combining a vector in $\Lambda_{\sigma}$ and a vector in $\Lambda_{\tau}$ with coefficients $1-t$ and $t$ gives a vector in $\Lambda_{(1-t)\sigma+t\tau}$, so $B$ is concave. The inclusion $\Lambda_{\tau}\subseteq\Lambda_{\sigma}$ for $\sigma\leq\tau$ shows that $B$ is decreasing, and its definition gives $B\geq A$. Fix $x$ outside the union of the polar sets $\{A_i=-\infty\}$, and let $f$ be any concave majorant of $\tau\mapsto A_{\tau}(x)$. Since $f$ is bounded below on the left half-line, it is decreasing. For every $\lambda\in\Lambda_{\tau}$, we have
\[
f(\tau)\geq f\left(\sum_i\lambda_i\tau_i\right)
\geq\sum_i\lambda_i f(\tau_i)
\geq\sum_i\lambda_i A_i(x).
\]
Thus $f(\tau)\geq B_{\tau}(x)$, which proves that $B$ agrees with the pointwise upper concave envelope outside these polar sets. Since $B_{\tau}$ is $\theta$-psh and majorizes the pointwise envelope everywhere, upper semicontinuous regularization gives $A'_{\tau}=B_{\tau}$ on $X$.

For each such $x$, the finite maximum defining $B_{\tau}(x)$ consists of a fixed collection of affine functions of $\tau$ on each interval $(\tau_{j+1},\tau_j)$. Hence it is convex as well as concave there, and therefore affine. Its endpoint values are given by continuity. This identity extends across the polar sets because both sides are $\theta$-psh. Hence $A'$ is itself the piecewise linear curve associated with the values $A'_{\tau_i}$ at the parameters $\tau_i$.
\end{remark}

\begin{lemma}\label{lma:psieqpsiprimeondiv}
Let $A$ be a piecewise linear curve in $\PSH(X,\theta)$, and let $A'$ be the regularized concave envelope in \cref{rmk:pwlinhull}. Define $\tilde A$ by
\[
    \tilde{A}^{\theta+\omega}_{\tau}\coloneqq P_{\theta+\omega}[A'_{\tau}]_{\mathcal I}, \quad \omega\in\Kah(X),\ \tau<\tau_0.
\]
Then $\tilde A\in\PSH^{\NA}(X,\theta)$, and
\begin{equation}\label{eq:psieqsiprime}
A^{\An}=\tilde{A}^{\An}\quad \text{on }X^{\Div}.
\end{equation}
\end{lemma}
Here $\tau_0$ is as in \cref{def:piecewiselincurv}.
\begin{proof}
We use the notation of \cref{rmk:pwlinhull}. By \cref{prop:PIconc,lma:PIenvmono1}, each component $\tilde A^{\theta+\omega}$ is concave and decreasing. Its slices are $\mathcal I$-model, and its constant left tail is $P_{\theta+\omega}[\max_i A_i]_{\mathcal I}$. This tail has positive mass: $\theta+\omega+\ddc\max_i A_i\geq\omega$, and taking the $\mathcal I$-envelope can only increase the mass. The cocycle property of the $\mathcal I$-projection gives compatibility of the components, so $\tilde A\in\PSH^{\NA}(X,\theta)$ with $\tilde A_{\max}=\tau_0$.

Fix $v\in X^{\Div}$ and put $a_i=A_i^{\An}(v)$. Analytification commutes with finite maxima and positive linear combinations on divisorial valuations. Since $\mathcal I$-projection preserves these valuations, \eqref{eq:pwlinhull-simplex} gives
\[
\left(\tilde A^{\theta+\omega}_{\tau}\right)^{\An}(v)+\tau
=\max_{\lambda\in\Lambda_{\tau}}\sum_i\lambda_i a_i+\tau
\leq\max_i(a_i+\tau_i),\quad \tau<\tau_0.
\]
Conversely, $A'_{\tau}\geq A_i$ whenever $\tau<\tau_i$. Letting $\tau\nearrow\tau_i$ gives the reverse inequality after taking the supremum over $\tau<\tau_0$. Now \eqref{eq:psianvlinear} proves \eqref{eq:psieqsiprime}.
\end{proof}

\subsubsection{The algebraic setting}
\label{subsubsec:analytic_alg}

Let $L$ be a $\mathbb{Q}$-line bundle on $X$ and $h$ be a Hermitian metric on $L$ with $\theta=c_1(L,h)$.

\begin{lemma}\label{lma:varphian}
For any $\varphi\in \PSH(X,\theta)$, we have $\varphi^{\An}\in \PSH(L^{\An})$.
\end{lemma}
\begin{proof}
After replacing $L$ with a sufficiently high power, we may assume that $L$ is a line bundle. Take a very ample line bundle $H$ on $X$.
By Siu's uniform global generation theorem \cite{Siu98}, \cite[Theorem~6.27]{Dem12} there exists $b>0$ large enough so that $H^{b}\otimes L^k\otimes \mathcal{I}(k\varphi)$ is globally generated for all $k>0$. Let $\{s_i\}_i$ be a finite set of global sections that generate the sheaf $H^{b}\otimes L^k\otimes \mathcal{I}(k\varphi)$. Then
\[
v(\mathcal{I}(k\varphi))=\min_i v(s_i).
\]
It follows that $v\mapsto -k^{-1}v\left(\mathcal{I}(k\varphi) \right)$ lies in $\mathcal{H}^{\gf}_{\mathbb{Q}}((L+\frac{b}{k}H)^{\An})$. Using \eqref{eq:varphiAnalt}, we conclude that $\varphi^{\An}\in \PSH(L^{\An})$.
\end{proof}

\begin{lemma}\label{lma:pltest}
Let $\Gamma$ be a piecewise linear curve in $\PSH(X,\theta)$.
Then $\Gamma^{\An}\in \PSH(L^{\An})$.
\end{lemma}
\begin{proof}  The result follows from \eqref{eq:psianvlinear}, \cref{prop:BJclosedmax} and \cref{lma:varphian}.
\end{proof}

\begin{lemma}\label{lma:istab}
Let $R$ be a commutative $\mathbb{C}$-algebra of finite type and $I$ be an ideal of $ R[t]$. If for any $a\in S^1$, $a^*I\subseteq I$, then $I$ is stable under the $\mathbb C^*$-action. Moreover, there are ideals $I_0\subseteq I_1\subseteq \cdots \subseteq I_m$ in $R$ so that
\begin{equation}\label{eq:Iexp}
    I=I_0 + I_1t+ \cdots + I_m (t^m).
    \end{equation}
\end{lemma}

\begin{proof}
It suffices to argue that $I$ can be expanded as in \eqref{eq:Iexp}.
To see this, assume that $a\in I$. We can write $a=a_0+a_1t+\cdots+a_mt^m$ with $a_i\in R$. Then our assumption implies that $\sum_i a_i \rho^i t^i\in I$ as well for all $\rho\in S^1$. So by the Lagrange interpolation formula, $a_it^i\in I$ for all $i$. Therefore, we can write $I$ as $I_0+I_1t+I_2t^2+\cdots$ for some ideals $I_0\subseteq I_1\subseteq \cdots$ in $R$. But as $R$ is noetherian, there is $m\geq 0$ so that $I_{m'}=I_m$ for $m'>m$. \eqref{eq:Iexp} follows.
\end{proof}

\begin{lemma}\label{lma:GAGA}
Let $X$ be a complex projective variety and $p:X\times \mathbb{C}\rightarrow X$ be the natural projection. Assume that $\mathcal{I}$ is an analytic coherent ideal sheaf on $X\times \mathbb{C}$. Assume that $\mathcal{I}|_{X \times \mathbb C^*}=p^*\mathcal{J}$ for some coherent ideal sheaf $\mathcal{J}$ on $X$. Then $\mathcal{I}$ is the analytification of an algebraic coherent ideal sheaf.
\end{lemma}

\begin{proof}
    Let $q:X\times (\mathbb{P}^1\setminus \{0\})\rightarrow X$ be the natural projection. Since $\mathbb C^* \subset \mathbb{P}^1\setminus \{0\}$, we can glue $q^*\mathcal{J}$ with $\mathcal{I}$ to get an analytic coherent ideal sheaf on $X\times \mathbb{P}^1$. By the GAGA principle, this ideal sheaf is necessarily algebraic, hence so is its restriction to $X\times \mathbb{C}$.
\end{proof}

Next we point out a version of Siu's uniform global generation lemma \cite{Siu88} that we will need in the proof of our next theorem:

\begin{lemma}\label{lma:Siu_global_uniform_generatedness} Let $L$ be a big line bundle on $X$ such that $c_1(L) = \{\theta\}$ and $\Phi \in \PSH(X\times \Delta, p_1^*\theta)$, where $\Delta$ is the unit disk. Let $G$ be an ample line bundle on $X$. Then there exists $k>0$, depending only on $X$ and $G$, such that $p_1^*(G^k \otimes L^m)\otimes \mathcal I(m \Phi)$ is globally generated for all $m \in \mathbb{Z}_{>0}$.
\end{lemma}

\begin{proof} The argument for this is exactly the same as the one in \cite[Lemma~5.6]{BBJ21} with Nadel's vanishing theorem replaced by the family version proved by Matsumura in \cite[Theorem~1.7]{Mat16}.
\end{proof}

\begin{proposition}\label{prop:PhiinducePSH}
Let $\phi\in \PSH(X,\theta)_{>0}$ be a model potential and $\ell\in \mathcal{R}(X,\theta;\phi)$ with $\sup_X \ell_1\leq 0$. Let $\Phi$ be the complexification of $\ell$. Then the function
\[
v\mapsto -\sigma(v)(\Phi) \quad \text{for } v\in X^{\Div}
\]
admits a unique extension to an element of $\PSH(L^{\An})$.
\end{proposition}
\begin{proof}
We may assume that $L$ is a line bundle. Observe that the extension is unique if it exists by \cref{thm:pshdetbydiv}.

Let $p_1\colon X\times \mathbb{C}\rightarrow X$ be the projection. For each $t>0$, we have $\ell_t\in \mathcal{E}(X,\theta;\phi)$, hence $\ell_t\sim_P\phi$ by \cref{thm:Pvarphidiffdef} and $\ell_t\sim_{\mathcal{I}}\phi$ by \cref{cor:PimpliesI}. Thanks to \cref{prop:pull-backmis} and \cref{lma:fibIeq}, for each $m\in \mathbb{Z}_{>0}$, we have
\[
\mathcal{I}(m\Phi)|_{X\times \Delta^*} = p_1^* \mathcal{I}(m\phi)|_{X\times \Delta^*}.
\]
In particular, $\mathcal{I}(m\Phi)$ admits an $S^1$-invariant coherent extension to $X\times \mathbb{C}$, obtained by gluing it with $p_1^*\mathcal{I}(m\phi)$ over $X\times \mathbb{C}^*$.

From \cref{lma:istab} and \cref{lma:GAGA}, we get that
\begin{equation}\label{eq:ImPhiexp}
\mathcal{I}(m\Phi)=\mathfrak{a}_0 +\mathfrak{a}_1 t+\cdots+\mathfrak{a}_{N-1} t^{N-1}+\mathfrak{a}_N (t^N)\,,
\end{equation}
where the $\mathfrak{a}_i$'s are coherent ideal sheaves on $X$.

Using \cref{lma:Siu_global_uniform_generatedness}, there exists an ample line bundle $T$ over $X$ such that $p_1^* (T\otimes L^m)\otimes \mathcal{I}(m\Phi)$ is globally generated, which is equivalent to $T\otimes L^m\otimes \mathfrak{a}_i$ being globally generated for all $i$.\footnote{In contrast with the case where $\phi$ is bounded, explored in \cite{BBJ21}, $\mathfrak{a}_N \neq \mathcal{O}_X$ in general.}

We define
\[
\varphi_m(v)\coloneqq -\frac{1}{m}\sigma(v)(\mathcal{I}(m\Phi))=-\frac{1}{m}\min_i (v(\mathfrak{a}_i)+i),\quad v\in X^{\Div}.
\]
From the right-hand side of the formula, $\varphi_m$ can be extended to an element in $\mathcal{H}^{\gf}_{\mathbb{Q}}((L+m^{-1}T)^{\An})$, which we denote by the same symbol.

For $v\in X^{\Div}$,
\[
-\sigma(v)(\Phi)=\lim_{m\to\infty}-\frac{1}{2^m}\sigma(v)(\mathcal{I}(2^m\Phi))=\lim_{m\to\infty}\varphi_{2^m}(v)
\]
and the right-hand side defines an element in $\PSH(L^{\An})$ by definition, since $\{\varphi_{2^m}\}_m$ is decreasing by the subadditivity of multiplier ideals \cref{thm:multipsubadd}.
\end{proof}

\begin{corollary}\label{cor:testcurvegivepshme}
Let $\Gamma\in \PSH^{\NA}(X,\theta)$. Then $\Gamma^{\An}$ defined in \cref{def:Gammaana} admits a unique extension to $\PSH(L^{\An})$.
\end{corollary}
The extension will be denoted by the same notation $\Gamma^{\An}$.

\begin{proof}
Observe that the extension is unique if it exists by \cref{thm:pshdetbydiv}. We may assume that $\Gamma_{\max}=0$ without loss of generality.

When $\Gamma\in \PSH^{\NA}(X,\theta)_{>0}$, our assertion follows from \cref{prop:PhiinducePSH} and \cref{prop:psiandPhi}.

In general, fix an ample line bundle $H$ on $X$ and a K\"ahler form $\omega\in c_1(H)$. Then we know that
\[
\Gamma^{\An}=\left(\Gamma^{m^{-1}\omega}\right)^{\An}\in \PSH((L+m^{-1}H)^{\An})
\]
for any $m\in \mathbb{Z}_{>0}$. Therefore, $\Gamma^{\An}\in \PSH(L^{\An})$ by \cref{prop:BJdecseqstab}.
\end{proof}

\subsection{The comparison theorem}
\label{subsec:NA_comparison}
Let $X$ be a connected projective manifold of dimension $n$. Let $L$ be a pseudo-effective $\mathbb{Q}$-line bundle on $X$ and $h$ be a Hermitian metric on $L$ with $\theta=c_1(L,h)$.

Thanks to \cref{cor:testcurvegivepshme}, we already have a map
\begin{equation}\label{eq:analymap}
    \PSH^{\NA}(X,\theta)\rightarrow \PSH(L^{\An}),\quad \Gamma\mapsto \Gamma^{\An}.
\end{equation}
We observe that for $\Gamma\in \PSH^{\NA}(X,\theta)$ and a K\"ahler form $\omega$ on $X$, we have
\begin{equation}\label{eq:proj_NA_an_invariance}
\left( P_{\theta+\omega}\left[  \Gamma \right]_{\mathcal{I}} \right)^{\An}=\Gamma^{\An}.
\end{equation}
Also observe that
\begin{equation}\label{eq:Gammamax}
    \Gamma_{\max}=\Gamma^{\An}(v_{\triv}),
\end{equation}

\begin{lemma}\label{lma:anaorder}
    The map \eqref{eq:analymap} is order preserving. Moreover, suppose that $\Gamma,\Gamma'\in \PSH^{\NA}(X,\theta)$ satisfies that $\Gamma^{\An}\leq \Gamma'^{\An}$, then $\Gamma\leq \Gamma'$.

    In particular, the map \eqref{eq:analymap} is injective.
\end{lemma}
\begin{proof}
    The map \eqref{eq:analymap} is order preserving by definition. Let us take $\Gamma,\Gamma'\in \PSH^{\NA}(X,\theta)$ with $\Gamma^{\An}\leq \Gamma'^{\An}$. Fix a K\"ahler form $\omega$ on $X$.

    Let $v\in X^{\Div}$ and $t \in \mathbb{Q}_{>0}$. Then, using \eqref{eq:psiandef} and the homogeneity of divisorial valuations, we notice that
\begin{equation}\label{eq:legendre_NA}
t \Gamma^\An \left(t^{-1} v \right) =  \sup_{\tau\in \mathbb{R}}\left( \left(\Gamma_{\tau}^{\theta+\omega}\right)^{\An}(v)+ t \tau\right).
\end{equation}
A similar equality holds for $\Gamma'$. Therefore, by \cref{cor:Leginvpartdef}, we have
\[
\left(\Gamma_{\tau}^{\theta+\omega}\right)^{\An}\leq \left(\Gamma'{}_{\tau}^{\theta+\omega}\right)^{\An}
\]
for all $\tau \in \mathbb{R}$. It follows that
\[
\Gamma_{\tau}^{\theta+\omega}\preceq_{\mathcal{I}} \Gamma'{}_{\tau}^{\theta+\omega}
\]
by \cref{prop:Iequivchar2}. Since both sides are $\mathcal{I}$-model, \cref{prop:Icomparandenvelope} gives
\[
\Gamma_{\tau}^{\theta+\omega}\leq \Gamma'{}_{\tau}^{\theta+\omega}
\]
for all $\tau\in \mathbb{R}$. Hence $\Gamma\leq \Gamma'$ by the definition of the order.
\end{proof}

\begin{lemma}\label{lma:gfinimage}
    Let $\phi \in \mathcal{H}^{\gf}_{\mathbb{Q}}(L^{\An})$. Then there is a piecewise linear curve $A$ in $\PSH(X,\theta)$ with $\phi=A^{\An}$. In particular, $\phi$ is in the image of \eqref{eq:analymap}.
\end{lemma}
In the construction below, the numbers $\tau_0,\ldots,\tau_N$ are rational, and the regularized concave envelope $A'_{\tau}$ is a generalized Fubini--Study metric for every rational $\tau\leq\tau_0$.
\begin{proof}
     Let us write
    \begin{equation}\label{eq:phiintominversemax}
        \phi=\frac{1}{m}\max_{i=1,\ldots,M} \left(\log|s_i|+\lambda_i\right),
    \end{equation}
where $m\in \mathbb{Z}_{>0}$, $s_{1},\ldots, s_M$ are a finite number of sections of $L^m$ and $\lambda_1,\ldots,\lambda_M\in \mathbb{Q}$.

Write $I_{\tau}$ for the set of $i$ such that $m^{-1}\lambda_i=\tau$. We denote the finitely many $\tau$ so that $I_{\tau}$ is non-empty as $\tau_0> \cdots > \tau_N$.
For each $i=0,\ldots,N$, we write
\[
A_{\tau_i}=\frac{1}{m}\max_{j\in I_{\tau_i}} \log |s_j|^2_{h^m}.
\]
We define $A$ as the piecewise linear curve associated with the $A_{\tau_i}$'s and the $\tau_i$'s. For every $v\in X^{\An}$,
\[
A^{\An}(v)=\max_{i=0,\ldots,N}\left(A_{\tau_i}^{\An}(v)+\tau_i\right)
=\frac{1}{m}\max_{j=1,\ldots,M}\left(\log|s_j|(v)+\lambda_j\right)=\phi(v),
\]
where we used the convention $\log|s_j|(v)=-v(s_j)$.

The final assertion follows from \cref{lma:psieqpsiprimeondiv} and \cref{thm:pshdetbydiv}.

For rational $\tau$, all vertices of $\Lambda_{\tau}$ in \cref{rmk:pwlinhull} are rational. Hence $A'_{\tau}$ is obtained from the Fubini--Study potentials $A_i$ by rational convex combinations and finite maxima, so it is a generalized Fubini--Study metric and has analytic singularities. Analytification commutes with these operations on all of $X^{\An}$, as follows by taking products and sums of the corresponding ideals. The estimate in the proof of \cref{lma:psieqpsiprimeondiv}, and rational parameters tending to each $\tau_i$ from below, therefore give
\begin{equation}\label{eq:gf-rational-hull}
\phi(v)=\sup_{\substack{\tau\in\mathbb{Q}\\\tau<\tau_0}}\left((A'_{\tau})^{\An}(v)+\tau\right),\quad v\in X^{\An}.
\end{equation}
\end{proof}

\begin{proposition}\label{prop:anpreinf}
    Let $(\Gamma_i)_{i\in I}$ be a decreasing net in $\PSH^{\NA}(X,\theta)$. Assume that \eqref{eq:decnetcontition} is satisfied. Then
    \[
        \left( \inf_{i\in I} \Gamma_{i} \right)^{\An}=\inf_{i\in I} \Gamma_{i}^{\An}.
    \]
\end{proposition}
\begin{proof}
    Take a K\"ahler form $\omega$ on $X$. We need to show that
    \[
    \left( \inf_{i\in I} P_{\theta+\omega}[\Gamma_i]_{\mathcal{I}} \right)^{\An}=\inf_{i\in I} \Gamma_{i}^{\An}.
    \]
    Therefore, after replacing $\theta$ by $\theta+\omega$, we may assume that $\Gamma_i\in \PSH^{\NA}(X,\theta)_{>0}$ for all $i\in I$ and $\inf_{i\in I}\Gamma_i\in \PSH^{\NA}(X,\theta)_{>0}$.
    Fix $v\in X^{\Div}$. By \cref{thm:pshdetbydiv}, it suffices to prove that
    \begin{equation}\label{eq:anpreinf1}
        \sup_{\tau\in \mathbb{R}} \left( \left(\inf_{i\in I}\Gamma_{i,\tau}\right)^{\An}(v)+\tau \right)= \adjustlimits\inf_{i\in I}\sup_{\tau\in \mathbb{R}} \left( \Gamma_{i,\tau}^{\An}(v)+\tau \right).
    \end{equation}
    Write $\Gamma=\inf_{i\in I}\Gamma_i$. We first compare the two sides inside the Legendre-type supremum. If $\tau<\Gamma_{\max}$, then by the definition of the infimum test curve,
    \[
        \Gamma_{\tau}=\inf_{i\in I}\Gamma_{i,\tau}.
    \]
    The net $(\Gamma_{i,\tau})_i$ decreases to $\Gamma_{\tau}$. By \cref{lma:testcurinfcompatible}, $\Gamma_{\tau}$ is stable under $P_{\theta}[\bullet]$.
    Hence \cref{prop:vol_limit_model} gives convergence of the masses, and \cref{cor:decnetdS} yields
    \[
        \Gamma_{i,\tau}\xrightarrow{d_S}\Gamma_{\tau}.
    \]
    Hence \cref{thm:Lelongcont} gives
    \[
        \Gamma_{\tau}^{\An}(v)=\inf_{i\in I}\Gamma_{i,\tau}^{\An}(v),\quad \tau<\Gamma_{\max}.
    \]
    If $\tau>\Gamma_{\max}$, then $\Gamma_{\tau}\equiv-\infty$, and also $\inf_i\Gamma_{i,\tau}\equiv-\infty$, since $\Gamma_{\max}=\inf_i(\Gamma_i)_{\max}$. Thus the same equality holds for $\tau>\Gamma_{\max}$ as well, with value $-\infty$. The possible endpoint $\tau=\Gamma_{\max}$ does not affect the supremum: for any test curve $\Lambda$ and any $v\in X^{\Div}$,
    \[
        \sup_{\tau\leq \Lambda_{\max}}\left(\Lambda_{\tau}^{\An}(v)+\tau\right)
        =\sup_{\tau< \Lambda_{\max}}\left(\Lambda_{\tau}^{\An}(v)+\tau\right),
    \]
    because $\Lambda_{\Lambda_{\max}}\leq \Lambda_{\tau}$ for $\tau<\Lambda_{\max}$, hence
    \[
        \Lambda_{\Lambda_{\max}}^{\An}(v)+\Lambda_{\max}
        \leq \Lambda_{\tau}^{\An}(v)+\Lambda_{\max}
    \]
    and we let $\tau\nearrow \Lambda_{\max}$.

    Consequently, in \eqref{eq:anpreinf1} we may replace $\left(\inf_i\Gamma_{i,\tau}\right)^{\An}(v)$ by $\inf_i\Gamma_{i,\tau}^{\An}(v)$ inside the supremum. The equality then follows from \cref{prop:Legendretranssup} applied to the closed convex functions $-\Gamma_{i,\bullet}^{\An}(v)$.
\end{proof}

\begin{theorem}\label{thm:DXZ}
    The map \eqref{eq:analymap} is an order preserving bijection.
\end{theorem}
\begin{proof}
    The map \eqref{eq:analymap} is an order preserving injection by \cref{lma:anaorder}. It remains to prove that it is surjective. Let $\phi\in \PSH(L^{\An})$. We want to construct $\Gamma\in \PSH^{\NA}(X,\theta)$ with $\Gamma^{\An}=\phi$.

    Let $H$ be an ample line bundle and $(\epsilon_i)_i$ be a decreasing sequence of rational numbers with limit $0$, $\phi_i\in \mathcal{H}^{\gf}_{\mathbb{Q}}((L+\epsilon_i H)^{\An})$ such that
    \[
    \phi=\inf_{i>0}\phi_i.
    \]
    The existence of these data is guaranteed by \cref{thm:BJappseq}. Fix a K\"ahler form $\omega\in c_1(H)$.

    Thanks to \cref{lma:gfinimage}, we can find $\Gamma^i\in \PSH^{\NA}(X,\theta+\epsilon_i \omega)$ with $(\Gamma^i)^{\An}=\phi_i$. It follows from \cref{lma:anaorder} that
    \[
    \Gamma^i\geq P_{\theta+\epsilon_i\omega}\left[ \Gamma^{i+1} \right]_{\mathcal{I}} \geq \Gamma^{i+1}.
    \]
    Therefore, for any $\omega'\in \Kah(X)$, after discarding finitely many $i$ so that $\omega'\geq \epsilon_i\omega$, the sequence $(P_{\theta+\omega'}\left[ \Gamma^{i} \right]_{\mathcal{I}})_i$ is decreasing. This does not change its infimum. We let
    \[
    \Gamma^{\theta+\omega'}=\inf_{i\gg 1}P_{\theta+\omega'}\left[ \Gamma^{i} \right]_{\mathcal{I}}\in \PSH^{\NA}(X,\theta+\omega').
    \]
    Note that the infimum is defined thanks to \eqref{eq:Gammamax}. It follows from \cref{prop:anpreinf} that
    \[
    \left(\Gamma^{\theta+\omega'} \right)^{\An}=\phi.
    \]
    From this, we conclude that for $\omega',\omega''\in \Kah(X)$ with $\omega'\leq \omega''$, we have
    \[
    P_{\theta+\omega''}\left[ \Gamma^{\theta+\omega'}\right]_{\mathcal{I}}=\Gamma^{\theta+\omega''}.
    \]
    Indeed, both sides have analytification $\phi$, by \eqref{eq:proj_NA_an_invariance} and \eqref{eq:analymap}. The injectivity part of \cref{lma:anaorder} gives our assertion.

    It follows that $(\Gamma^{\theta+\omega'})_{\omega'\in \Kah(X)}$ defines an element $\Gamma$ in $\PSH^{\NA}(X,\theta)$ and $\Gamma^{\An}=\phi$.
\end{proof}

\begin{theorem}\label{thm:NApotential:L1352}
    Under the bijection \cref{thm:DXZ}, the operations on the space $\PSH^{\NA}(X,\theta)$ defined in \cref{sec:opnA} all correspond to the corresponding operations on $\PSH(L^{\An})$ in Boucksom--Jonsson's theory.
\end{theorem}
The meaning is immediate for all operations except for the trace operator, and the proofs are elementary, as we have seen in \cref{prop:anpreinf} for the infimum operator. We shall only restate and prove the case of trace operators; the remaining arguments are recorded in \cite{Xia23Operations}.

\begin{theorem}\label{thm:comp_traceop}
    Let $Y\subseteq X$ be an irreducible analytic subset.
    Consider an element $\Gamma\in \PSH^{\NA}(X,\theta)$ with well-defined restriction to $Y$. Then
    \begin{equation}\label{eq:Tracena}
    \Tr_Y(\Gamma)^{\An}|_{Y^{\Div}}=\Gamma^{\An}|_{Y^{\Div}}.
    \end{equation}
\end{theorem}
Observe that there is a canonical identification $Y^{\Div}=\widetilde{Y}^{\Div}$. Recall that a generalized Fubini--Study metric is defined in \cref{def:gFSmetric}.

\begin{proof}
    We first reduce to the positive mass case. Choose an ample line bundle $H$ and a K\"ahler form $\omega_0\in c_1(H)$, and set $\widehat{\Gamma}=P_{\theta+\omega_0}[\Gamma]_{\mathcal{I}}$. Replacing $(L,\theta,\Gamma)$ by $(L+H,\theta+\omega_0,\widehat{\Gamma})$ does not change $\Gamma^{\An}$, by \eqref{eq:proj_NA_an_invariance}; it also only replaces $\Tr_Y(\Gamma)$ by the corresponding $\mathcal{I}$-projection on $\widetilde{Y}$, whose analytification is unchanged. Hence we may assume that $\Gamma\in \PSH^{\NA}(X,\theta)_{>0}$.
    Let $\phi=\Gamma^{\An}\in \PSH(L^{\An})$. By \cref{lma:non-trivrest}, $\phi(v_{Y,\triv})\neq -\infty$.

    Take an ample line bundle $S$ on $X$ and a K\"ahler form $\omega$ in $c_1(S)$.
    By \cref{thm:BJappseq}, there are rational numbers $0<\epsilon_i<1$ decreasing to $0$ and a decreasing sequence $\phi^i\in\mathcal{H}^{\gf}_{\mathbb{Q}}((L+\epsilon_i S)^{\An})$ with limit $\phi$.

    Take $\Gamma^{i}\in \PSH^{\NA}(X,\theta+\epsilon_i\omega)$ such that $\Gamma^{i,\An}=\phi^{i}$. By \cref{lma:anaorder}, $\Gamma^i\geq P_{\theta+\epsilon_i\omega}[\Gamma]_{\mathcal I}$. Thus $\Gamma^i$ has positive volume, and \cref{lma:positivemasspart} gives $\Gamma^i\in \PSH^{\NA}(X,\theta+\epsilon_i\omega)_{>0}$.

    Write $\widehat{\Gamma}^i\coloneqq P_{\theta+\omega}[\Gamma^i]_{\mathcal{I}}$ and $\widehat{\Gamma}\coloneqq P_{\theta+\omega}[\Gamma]_{\mathcal{I}}$. Applying \cref{prop:anpreinf} to the decreasing sequence $(\widehat{\Gamma}^i)_i$ and then using the injectivity in \cref{lma:anaorder}, we get $\inf_i\widehat{\Gamma}^i=\widehat{\Gamma}$. Hence, by the componentwise definition of this infimum, for any $\tau<\Gamma_{\max}$, we have
    \[
    \inf_{i\in \mathbb{N}}\widehat{\Gamma}^{i}_{\tau}=\widehat{\Gamma}_{\tau}.
    \]
    Now \cref{prop:vol_limit_model} and \cref{cor:decnetdS} give
    \[
        \widehat{\Gamma}^{i}_{\tau}\xrightarrow{d_S} \widehat{\Gamma}_{\tau}.
    \]
    for all $\tau<\Gamma_{\max}$.

    By \cref{lma:non-trivrest} again, each $\widehat{\Gamma}^{i}$ has non-trivial restriction to $Y$. By \cref{prop:tracedeclimit}, for any K\"ahler form $\omega'$ on $\widetilde{Y}$ satisfying $\omega'\geq \omega|_{\widetilde{Y}}$ we have
    \[
    \Tr_Y \left( \widehat{\Gamma}^{i}_{\tau} \right)\xrightarrow{d_S} \Tr_Y \left( \widehat{\Gamma}_{\tau}\right)
    \]
    for any $\tau<(\Tr_Y(\Gamma))_{\max}$. Thanks to \cref{thm:Lelongcont}, for every $c\ord_F\in Y^{\Div}$ and every such $\tau$,
    \[
    \Tr_Y \left( \widehat{\Gamma}_{\tau} \right)^{\An}(c\ord_F)=\inf_{i\geq 1}\Tr_Y \left( \widehat{\Gamma}^{i}_{\tau} \right)^{\An}(c\ord_F).
    \]
    Put $T_i=(\Tr_Y(\Gamma^i))_{\max}$ and $T=(\Tr_Y(\Gamma))_{\max}$. These endpoints satisfy $T_i\searrow T$. Indeed, $T_i\geq T$ by monotonicity. For $T<\tau<\Gamma_{\max}$, the $d_S$-convergence above and \cref{thm:Lelongcont} imply that $\nu(\widehat{\Gamma}^i_{\tau},Y)>0$ for all sufficiently large $i$, so $T_i\leq\tau$. If $T=\Gamma_{\max}$, use $T_i\leq\Gamma^i_{\max}\searrow\Gamma_{\max}$ instead. Thus the extended trace slices converge also for $\tau>T$, where they are eventually $-\infty$; the value at $T$ does not affect the Legendre supremum. Applying \cref{prop:Legendretranssup}, and using again that $\mathcal{I}$-projection does not change analytification, we get
    \[
    \Tr_Y(\Gamma)^{\An}(c\ord_F)=\inf_{i\geq 1}\Tr_Y(\Gamma^i)^{\An}(c\ord_F).
    \]
    Here the decreasing infimum of the conjugates is finite for every positive dual parameter, hence continuous there; the closure in \cref{prop:Legendretranssup} therefore does not change its value at $1$.
    Since also $\Gamma^{\An}=\inf_i \Gamma^{i,\An}$, it suffices to prove \eqref{eq:Tracena} with $\Gamma^i$ in place of $\Gamma$.

    In other words, we have reduced to the case where $\phi\in \mathcal{H}^{\gf}_{\mathbb{Q}}(L^{\An})$ and $L$ is big.

    Let $A$ and $A'$ be the curves constructed in \cref{lma:gfinimage}, write $\Gamma'_{\tau}=A'_{\tau}$, and set $T=(\Tr_Y(\Gamma))_{\max}$. By \cref{lma:psieqpsiprimeondiv,lma:anaorder}, the inverse system obtained from $A'$ is $\Gamma$. Hence $\Gamma_{\max}=\tau_0$ and
    \[
    \Gamma_{\tau}\sim_{\mathcal I}\Gamma'_{\tau},\qquad
    \Gamma_{\tau}=P_{\theta}[\Gamma'_{\tau}]_{\mathcal I},\quad \tau<\Gamma_{\max}.
    \]
    For rational $\tau<\Gamma_{\max}$, the potential $\Gamma'_{\tau}$ is a generalized Fubini--Study metric and has analytic singularities. Consequently, for rational $\tau<T$, its restriction to $\widetilde Y$ is not identically $-\infty$, while for rational $T<\tau<\Gamma_{\max}$ its analytification is $-\infty$ on $Y^{\An}$. The endpoint $T$ can be omitted since $\Gamma'$ is decreasing. Thus \eqref{eq:gf-rational-hull} gives, for any $c\ord_F\in\widetilde Y^{\Div}$,
    \[
    \begin{aligned}
    \phi|_{Y^{\An}}(c\ord_F)=&\sup_{\substack{\tau\in\mathbb Q\\\tau<T}} \left(\Gamma'_{\tau}{}^{\An}(c\ord_F)+\tau\right)\\
    =& \sup_{\substack{\tau\in\mathbb Q\\\tau<T}} \left(\left(\Gamma'_{\tau}|_{\widetilde{Y}}\right)^{\An}(c\ord_F)+\tau\right)\\
    =& \sup_{\substack{\tau\in\mathbb Q\\\tau<T}} \left(\Tr_Y(\Gamma'_{\tau})^{\An}(c\ord_F)+\tau\right)\\
     =& \sup_{\substack{\tau\in\mathbb Q\\\tau<T}} \left(\Tr_Y(\Gamma_{\tau})^{\An}(c\ord_F)+\tau\right)\\
    =& \Tr_Y(\Gamma)^{\An}(c\ord_F).
    \end{aligned}
    \]
    The third equality follows from \cref{ex:resanalyt}, and the fourth from \cref{prop:tracelinear} and $\Gamma_{\tau}\sim_{\mathcal I}\Gamma'_{\tau}$. For the last equality, the supremum defining the analytification of the trace curve can be taken over rational parameters: every real parameter can be approached from below, and the slices are decreasing. The $\mathcal I$-projections in the definition of the trace do not change divisorial valuations. It remains to justify the second line. Namely, for any generalized Fubini--Study metric $\varphi$ whose restriction to $\widetilde{Y}$ is not identically $-\infty$, we have
    \begin{equation}\label{eq:varphiAncordF}
    \varphi^{\An}(c\ord_F)=\left(\varphi|_{\widetilde{Y}}\right)^{\An}(c\ord_F).
    \end{equation}
    We can immediately reduce to the case where $\varphi$ is a Fubini--Study metric, and then to the case
    \[
    \varphi=\log |s|_{h_0}^2,
    \]
    where $s$ is a holomorphic section of $L$, not vanishing identically on $Y$, in which case \eqref{eq:varphiAncordF} is obvious.
\end{proof}

\begin{lemma}\label{lma:non-trivrest}
    Let $\Gamma\in \PSH^{\NA}(X,\theta)$ and $Y\subseteq X$ be an irreducible analytic subset. Then the following are equivalent:
    \begin{enumerate}
    \item\label{itm:lma:non-trivrest:1} $\Gamma^{\An}(v_{Y,\triv})\neq -\infty$;
    \item\label{itm:lma:non-trivrest:2} $\Gamma^{\An}|_{Y^{\An}}\not\equiv -\infty$;
    \item\label{itm:lma:non-trivrest:3} $\Gamma$ has a well-defined restriction to $Y$.
\end{enumerate}
\end{lemma}
Here $v_{Y,\triv}$ denotes the trivial valuation of $\mathbb{C}(Y)$.

\begin{proof}
The equivalence between \ref{itm:lma:non-trivrest:1} and \ref{itm:lma:non-trivrest:2} is a simple consequence of the maximum principle \cite[Lemma~1.4(i)]{BJ22}.

To see the equivalence between \ref{itm:lma:non-trivrest:1} and \ref{itm:lma:non-trivrest:3}, it suffices to observe that for any $\varphi\in \PSH(X,\theta)$,
\[
\varphi^{\An}(v_{Y,\triv})=
\begin{cases}
    -\infty,&\textup{if }\nu(\varphi,Y)>0;\\
    0,& \textup{if }\nu(\varphi,Y)=0.
\end{cases}
\]
This completes the proof.
\end{proof}

\motto{I speak twelve languages: English is the \emph{bestest}.

--- Stefan Bergman\footnote{Stefan Bergman (1895--1977), bearing a very Scandinavian name, was a Polish-American mathematician best known for his work in complex analysis, especially in several complex variables. He introduced the Bergman kernel, a fundamental concept in complex analysis that has influenced many areas of mathematics and theoretical physics.

Bergman was born in Poland (then part of the Russian Empire), and studied in Berlin. He fled Europe during World War II and eventually settled in the United States.}}
\chapter{Partial Bergman kernels}
\label{chap:Bergman}

The classical Bergman kernel of a Hermitian line bundle $L$ on a compact Kähler manifold $X$ encodes, through its asymptotic expansion as the tensor power $L^k$ goes to infinity, deep geometric information about the underlying metric. This chapter develops the corresponding \emph{partial} Bergman kernel. The word partial refers to the fact that we only use sections cut out by an ideal, typically the multiplier ideal of an auxiliary psh weight on a closed subset $K\subseteq X$.

The main result is the convergence theorem \cref{thm:pBMconvergence}. Given a $\theta$-psh potential $\varphi\in \PSH(X,\theta)_{>0}$, a closed non-pluripolar weighted subset $(K,v)$, and a Bernstein--Markov measure $\mu$, the normalized partial Bergman measures on $X$ converge weakly to the partial equilibrium measure of $(K,v,\varphi)$ as $k\to\infty$. This is the relative version of the theorem of Berman--Boucksom--Witt Nystr\"om \cite{BBWN11}, corresponding to the case $\varphi=V_\theta$. The assumptions are essentially the mildest ones under which the statement makes sense.

The first part of the chapter introduces the relative envelope operators. In \cref{sec:envrel} we define the relative $P$-envelope $P_{\theta,K}[\varphi](v)$ and the relative $\mathcal{I}$-envelope $P_{\theta,K}[\varphi]_{\mathcal{I}}(v)$ by replacing the constraint $\eta\leq 0$ in the usual $P$- and $\mathcal{I}$-envelopes with the constraint $\eta|_K\leq v$. We establish their basic properties, including monotonicity, mass-preservation, and the relation with the envelope operators of \cref{chap:envelopes}.

The second part develops the equilibrium theory associated with a weighted subset $(K,v)$. In \cref{sec:pem} we show that the Monge--Amp\`ere measure $\theta_{P_{\theta,K}[\varphi](v)}^n$ is carried by the contact locus $\{P_{\theta,K}[\varphi](v)|_K=v\}\subseteq K$, which justifies the name \emph{partial equilibrium measure}. The chapter concludes in \cref{sec:quant} with the proof of the weak convergence of partial Bergman measures towards this partial equilibrium.

\section{Partial envelopes}\label{sec:envrel}
Let $X$ be a connected compact Kähler manifold of dimension $n$ and $K \subseteq X$ be a closed non-pluripolar set.
Let $\theta$ be a smooth closed real $(1,1)$-form on $X$ representing a pseudo-effective cohomology class.

\begin{definition}\label{def:partial_env}
    Given a function $f\colon K\rightarrow [-\infty,\infty]$ and $\varphi\in \PSH(X,\theta)$, we introduce the \emph{relative $P$-envelope}\index{envelope!P-envelope@$P$-envelope!relative P-envelope@relative $P$-envelope} of $\varphi$ (with respect to $K,f,\theta$) as
    \begin{equation}\label{eq:PthetaKvarphiv}
    P_{\theta,K}[\varphi](f)\coloneqq  \uscsup \Bigl\{ \eta \in \PSH(X,\theta): \eta|_K \leq f  \textup{ quasi-everywhere and } \eta\preceq \varphi\Bigr\}.
    \end{equation}
    \index{$P_{\theta,K}[\varphi](f)$}

    Similarly, we define the \emph{relative $\mathcal{I}$-envelope}\index{envelope!I-envelope@$\mathcal{I}$-envelope!relative I-envelope@relative $\mathcal{I}$-envelope} of $\varphi$ (with respect to $K,f,\theta$) as
    \begin{equation}\label{eq:PthetaKvarphivI}
    P_{\theta,K}[\varphi]_{\mathcal{I}}(f)\coloneqq  \uscsup \Bigl\{ \eta \in \PSH(X,\theta): \eta|_K \leq f  \textup{ quasi-everywhere and } \eta\preceq_{\mathcal{I}} \varphi\Bigr\}.
    \end{equation}
\end{definition}
By $\eta|_K\leq f$ quasi-everywhere, we mean that there is a pluripolar set $S\subseteq X$ so that $\eta|_{K\setminus S}\leq f|_{K\setminus S}$. In general the upper envelope in \eqref{eq:PthetaKvarphiv} or \eqref{eq:PthetaKvarphivI} may be identically $\infty$.

Observe that when $f$ is bounded from below, neither envelope is identically $-\infty$.

The envelopes contain essentially all envelope operators that we have studied in this book:
\begin{example}\label{ex:rel_env_cases}
    When $K=X$ and $f=0$, we have
    \[
        P_{\theta,X}[\varphi](0)=P_{\theta}[\varphi],\quad P_{\theta,X}[\varphi]_{\mathcal{I}}(0)=P_{\theta}[\varphi]_{\mathcal{I}}.
    \]
    When $K=X$, we have
    \[
        P_{\theta,X}[\varphi](f)=P_{\theta}[\varphi](f).
    \]
    These definitions reduce to the envelopes studied in \cref{chap:envelopes}.
    On the other hand, when $\varphi=V_{\theta}$, we have
    \[
    \begin{aligned}
        P_{\theta,K}[V_{\theta}](f)
        &=P_{\theta,K}[V_{\theta}]_{\mathcal{I}}(f)\\
        &=\uscsup \left\{ \eta \in \PSH(X,\theta):
        \begin{array}{l}
        \eta|_K \leq f \textup{ quasi-everywhere}
        \end{array}
        \right\}.
    \end{aligned}
    \]
    In an even more special case where $\varphi=V_{\theta}$ and $K=X$, we find that
    \[
        P_{\theta,X}[V_{\theta}](f)=\uscsup \Bigl\{ \eta \in \PSH(X,\theta): \eta \leq f \textup{ quasi-everywhere}\Bigr\}.
    \]
    This is nothing but $P_{\theta}(f)$ as we studied in \cref{def:Pthetaf}.
\end{example}

\begin{proposition}\label{prop:partial_env}
    Let $f,g\colon K\rightarrow [-\infty,\infty]$ be functions and $\varphi,\psi\in \PSH(X,\theta)$. Then we have the following:
    \begin{enumerate}
        \item\label{itm:prop:partial_env:1} We have
        \[
        P_{\theta,K}[\varphi](f)\leq P_{\theta,K}[\varphi]_{\mathcal{I}}(f);
        \]
        \item\label{itm:prop:partial_env:2} assume that $f$ is bounded from above, then
        \[
        \begin{aligned}
            \left.P_{\theta,K}[\varphi](f)\right|_K\leq & f\textup{ quasi-everywhere},\\
             \left. P_{\theta,K}[\varphi]_{\mathcal{I}}(f)\right|_K\leq & f\textup{ quasi-everywhere};
        \end{aligned}
        \]
        \item\label{itm:prop:partial_env:3} suppose that $\varphi\preceq_{\mathcal{I}}\psi$, then
        \[
            P_{\theta,K}[\varphi]_{\mathcal{I}}(f)\leq P_{\theta,K}[\psi]_{\mathcal{I}}(f);
        \]
        \item\label{itm:prop:partial_env:4} if $f\leq g$, then
        \[
        P_{\theta,K}[\varphi](f)\leq P_{\theta,K}[\varphi](g),\quad P_{\theta,K}[\varphi]_{\mathcal{I}}(f)\leq P_{\theta,K}[\varphi]_{\mathcal{I}}(g);
        \]
        \item\label{itm:prop:partial_env:5} for any $t \in (0,1)$, we have
        \[
        \begin{aligned}
            P_{\theta,K}\left[t\varphi+(1-t)\psi\right](f)\geq & tP_{\theta,K}[\varphi](f)+(1-t)P_{\theta,K}[\psi](f),\\
            P_{\theta,K}\left[t\varphi+(1-t)\psi\right]_{\mathcal{I}}(f)\geq & tP_{\theta,K}[\varphi]_{\mathcal{I}}(f)+(1-t)P_{\theta,K}[\psi]_{\mathcal{I}}(f).
        \end{aligned}
        \]
    \end{enumerate}
\end{proposition}
\begin{proof}
    \ref{itm:prop:partial_env:1} This follows from $\eta\preceq \varphi\Rightarrow \eta\preceq_{\mathcal I}\varphi$.

    \ref{itm:prop:partial_env:2} This follows from \cref{prop:supsupstardiff} and \cref{rmk:L1cpt}.

    \ref{itm:prop:partial_env:3} This follows by inclusion of the defining candidate sets.

    \ref{itm:prop:partial_env:4} This follows by inclusion of the defining candidate sets.

    \ref{itm:prop:partial_env:5} If $\eta,\eta'\in \PSH(X,\theta)$ are such that
    \[
    \begin{split}
        &\eta|_K\leq f\textup{ quasi-everywhere},\quad \eta'|_K\leq f\textup{ quasi-everywhere},\\
        &\eta\preceq \varphi,\quad \eta'\preceq \psi,
    \end{split}
    \]
    then for $t\in (0,1)$,
    \[
        \left.\left(t\eta+(1-t)\eta'\right)\right|_K\leq f\textup{ quasi-everywhere},\quad t\eta+(1-t)\eta' \preceq t\varphi+(1-t)\psi.
    \]
    Therefore,
    \[
    t\eta+(1-t)\eta'\leq P_{\theta,K}\left[t\varphi+(1-t)\psi\right](f).
    \]
    Since $\eta$ and $\eta'$ are arbitrary, the first inequality follows. The second follows from the same argument.
\end{proof}

Similar to the usual envelopes studied in \cref{chap:envelopes}, the $P$-partial envelopes are well-behaved only when $\varphi$ has positive mass.

\begin{lemma}\label{lma:pos_mass_partial}
    Suppose that $f\colon K\rightarrow (-\infty,\infty]$ is bounded from below, and $\varphi\in \PSH(X,\theta)_{>0}$. Then
    \begin{equation}\label{eq:Ptheta_partial_pos_mass}
    P_{\theta,K}[\varphi](f)= \uscsup \Bigl\{ \eta \in \PSH(X,\theta)_{>0}: \eta|_K \leq f  \textup{ quasi-everywhere and } \eta\preceq \varphi\Bigr\}.
    \end{equation}
\end{lemma}
In other words, in this situation, in the definition \eqref{eq:PthetaKvarphiv}, we may restrict to only those $\eta$ with positive mass.
\begin{proof}
    By our assumption, we can find $C>0$ so that $\varphi|_K-C\leq f$. Then if $\eta$ is a candidate in \eqref{eq:PthetaKvarphiv}, $\eta\lor (\varphi-C)$ is also a candidate, and $\eta\lor (\varphi-C)\in \PSH(X,\theta)_{>0}$. Therefore, $\eta\lor (\varphi-C)$ is also a candidate for the right-hand side of \eqref{eq:Ptheta_partial_pos_mass}. The equality \eqref{eq:Ptheta_partial_pos_mass} follows.
\end{proof}

Note that in the most important cases, we can actually replace the condition $\eta\leq f$ quasi-everywhere in \eqref{eq:PthetaKvarphiv} by everywhere.
\begin{lemma}\label{lma:PKoutsidepps}
    Let $f\colon K\rightarrow (-\infty,\infty]$ be bounded from below, and $\varphi\in \PSH(X,\theta)_{>0}$.
    Then
    \begin{equation}\label{eq:PthetaK_new_def}
        \begin{split}
    P_{\theta,K}[\varphi](f)=&\uscsup \Bigl\{ \eta \in \PSH(X,\theta): \eta|_K \leq f  \textup{ and } \eta\preceq \varphi\Bigr\}\\
    =& \uscsup \Bigl\{ \eta \in \PSH(X,\theta)_{>0}: \eta|_K \leq f  \textup{ and } \eta\preceq \varphi\Bigr\}.
        \end{split}
    \end{equation}
\end{lemma}
\begin{proof}
Both $\geq$ directions in \eqref{eq:PthetaK_new_def} are trivial.

Conversely, let $\eta$ be a candidate in \eqref{eq:Ptheta_partial_pos_mass}.
Let $S\subseteq X$ be a pluripolar set so that $\eta|_{K\setminus S}\leq f|_{K\setminus S}$.

We now claim that we can choose $\chi\in \PSH(X,\theta)_{>0}$ with $\chi|_S\equiv -\infty$ and $\chi\leq \eta$.

Assume this for the moment. For any $\delta\in (0,1)$, we have
\[
(1-\delta)\eta|_K+\delta \chi|_K\leq f,\quad (1-\delta)\eta+\delta \chi\preceq \varphi.
\]
Hence, $(1-\delta)\eta+\delta \chi$ is dominated by the right-hand side of \eqref{eq:PthetaK_new_def}.
Letting $\delta \to 0+$, we conclude that $\eta$ is dominated by the right-hand side of \eqref{eq:PthetaK_new_def} and hence, \eqref{eq:PthetaK_new_def} follows.

    It remains to construct $\chi$.
    Since $\varphi$ has positive mass, the class $\{\theta\}$ is big. Fix a K\"ahler form $\omega_0$ on $X$ and choose $\rho\in \PSH(X,\theta)$ with $\theta_{\rho}\geq 2\epsilon_0\omega_0$ for some $\epsilon_0>0$. By \cref{thm:gloJosefson}, there is a quasi-psh function $q$ such that $q|_S\equiv -\infty$. Choose $A>0$ so that $A\omega_0+\ddc q\geq 0$. For $a>0$ small enough, $\rho+aq\in \PSH(X,\theta)$ and still has pole $-\infty$ on $S$. Replacing $\rho$ by $\rho+aq$, we may therefore assume that $\rho|_S\equiv -\infty$.
    Since $\int_X \theta_{\eta}^n>0$, fixing such a $\rho$ and taking $\epsilon\in (0,1)$ small enough, we have
    \[
    \int_X \theta^n_{\epsilon \rho+(1-\epsilon)V_{\theta}}+\int_X \theta_{\eta}^n>\int_X \theta_{V_{\theta}}^n.
    \]
    Thus, by \cref{prop:landfinitecond1} and \cref{thm:diamond}, we have
    \[
    \Bigl(\epsilon \rho+(1-\epsilon)V_{\theta}\Bigr) \land \eta \in \PSH(X,\theta)_{>0}.
    \]
    It suffices to define this function as $\chi$.
\end{proof}
\begin{corollary}\label{cor:Ppartbdd_bdabove}
    Suppose that $f\colon K\rightarrow \mathbb{R}$ is bounded, and $\varphi\in \PSH(X,\theta)_{>0}$. Then there is a constant $C=C(\sup_K f,K,X)>0$ such that
        \[
            P_{\theta,K}[\varphi](f)\leq C,\quad P_{\theta,K}[\varphi]_{\mathcal{I}}(f)\leq C.
        \]
\end{corollary}
\begin{proof}
    Since $K$ is non-pluripolar, \cref{rmk:L1cpt} then gives a constant $C>0$, depending only on $\sup_K f$, $K$ and $X$, such that every $\eta\in\PSH(X,\theta)$ satisfying $\eta|_K\leq f$ quasi-everywhere satisfies $\eta\leq C$ on $X$.

    This applies to the candidates defining both envelopes. Hence both candidate families are locally uniformly bounded above, so their upper semicontinuous suprema are genuine $\theta$-psh functions. By \itemref{itm:prop:partial_env:2}, both envelopes are dominated by $f$ quasi-everywhere on $K$. Applying \cref{rmk:L1cpt} once more gives
    \[
        P_{\theta,K}[\varphi](f)\leq C,\quad P_{\theta,K}[\varphi]_{\mathcal I}(f)\leq C.
    \]
\end{proof}

\begin{proposition}\label{prop:Penv_partial_monotone}
    Let $f\colon K\rightarrow (-\infty,\infty]$ be a function bounded from below and $\varphi,\psi\in \PSH(X,\theta)_{>0}$. Suppose that $\varphi\preceq_{P} \psi$. Then
    \[
    P_{\theta,K}[\varphi](f)\leq P_{\theta,K}[\psi](f).
    \]
\end{proposition}
\begin{proof}
    Let $\eta\in \PSH(X,\theta)_{>0}$ be such that $\eta|_K\leq f$ quasi-everywhere and $\eta\preceq \varphi$. Then $\eta\preceq_P \psi$, and thanks to \cref{cor:P_par_envelope_rel} for any $C>0$, we have $\eta\land (\psi+C)\in \PSH(X,\theta)_{>0}$. But then $\eta\land (\psi+C)$ is a candidate in the definition of $P_{\theta,K}[\psi](f)$, and we find
    \[
    \eta\land (\psi+C)\leq P_{\theta,K}[\psi](f).
    \]
    By \cref{cor:P_par_envelope_rel} again, we get
    \[
    \eta=P_{\theta}[\psi](\eta)\leq P_{\theta,K}[\psi](f).
    \]
    Our assertion follows now from \cref{lma:pos_mass_partial}.
\end{proof}

\begin{corollary}\label{cor:partial_P_new_def}
    Let $f\colon K\rightarrow (-\infty,\infty]$ be a function bounded from below and $\varphi \in \PSH(X,\theta)_{>0}$. Then
    \[
    \begin{aligned}
    P_{\theta,K}[\varphi](f)= & \uscsup \Bigl\{ \eta \in \PSH(X,\theta): \eta|_K \leq f  \textup{ quasi-everywhere and } \eta\preceq_P \varphi\Bigr\}\\
    =& \uscsup \Bigl\{ \eta \in \PSH(X,\theta)_{>0}: \eta|_K \leq f  \textup{ quasi-everywhere and } \eta\preceq_P \varphi\Bigr\}\\
    =& \uscsup \Bigl\{ \eta \in \PSH(X,\theta): \eta|_K \leq f  \textup{ and } \eta\preceq_P \varphi\Bigr\}\\
    =& \uscsup \Bigl\{ \eta \in \PSH(X,\theta)_{>0}: \eta|_K \leq f  \textup{ and } \eta\preceq_P \varphi\Bigr\}.
    \end{aligned}
    \]
\end{corollary}
\begin{proof}
    We first apply \cref{prop:Penv_partial_monotone},
    \[
    \begin{split}
    P_{\theta,K}[\varphi](f)= & P_{\theta,K}\left[P_{\theta}[\varphi]\right](f)\\
    =& \uscsup \Bigl\{ \eta \in \PSH(X,\theta): \eta|_K \leq f  \textup{ quasi-everywhere and } \eta\preceq_P \varphi\Bigr\}.
    \end{split}
    \]
    Next observe that by \cref{lma:pos_mass_partial}, the latter can also be written as
    \[
    \uscsup \Bigl\{ \eta \in \PSH(X,\theta)_{>0}: \eta|_K \leq f  \textup{ quasi-everywhere and } \eta\preceq_P \varphi\Bigr\}.
    \]
    Next we apply \cref{lma:PKoutsidepps}; the remaining statements follow.
\end{proof}

\begin{corollary}\label{cor:I-good_partial_env}
    Let $f\colon K\rightarrow (-\infty,\infty]$ be a function bounded from below and $\varphi \in \PSH(X,\theta)_{>0}$ be $\mathcal{I}$-good. Then
    \[
        P_{\theta,K}[\varphi]_{\mathcal{I}}(f)=P_{\theta,K}[\varphi](f).
    \]
\end{corollary}
\begin{proof}
    In fact, since $\varphi$ is $\mathcal{I}$-good, a potential $\eta\in \PSH(X,\theta)$ satisfies $\eta\preceq_P \varphi$ if and only if $\eta\preceq_{\mathcal{I}} \varphi$. Therefore, our assertion follows from \cref{cor:partial_P_new_def} and \eqref{eq:PthetaKvarphivI}.
\end{proof}

Next, we make the following observations about the singularity types of our envelopes:
\begin{lemma}\label{lma:same_sing_type} Let $f \colon K\rightarrow \mathbb{R}$ be a bounded function, and $\varphi\in \PSH(X,\theta)_{>0}$. We have
\begin{equation}\label{eq:PthetaKvarphiv_singtype}
P_{\theta,K}[\varphi](f)\sim P_{\theta}[\varphi],\quad P_{\theta,K}[\varphi]_{\mathcal{I}}(f)\sim P_{\theta}[\varphi]_{\mathcal{I}}.
\end{equation}
\end{lemma}
\begin{proof}
    We first prove the $P$-part of \eqref{eq:PthetaKvarphiv_singtype}.

    If $\eta\in \PSH(X,\theta)$ with $\eta\preceq \varphi$ and $\eta\leq 0$, then $\eta+\inf_K f$ is a candidate in the definition \eqref{eq:PthetaKvarphiv}, and hence,
    \[
    \eta +\inf_K f \leq P_{\theta,K}[\varphi](f).
    \]
    It follows that
    \[
        P_{\theta}[\varphi]+\inf_K f \leq P_{\theta,K}[\varphi](f).
    \]

    Conversely, let $\eta\in \PSH(X,\theta)$ be such that
    \[
        \eta|_K\leq f,\quad \eta\preceq \varphi.
    \]
    Since $K$ is non-pluripolar, by \cref{rmk:L1cpt} we can find $C=C(\sup_K f,K, X)>0$ so that $\eta\leq C$. Hence $\eta-C$ is a candidate in the definition of $P_{\theta}[\varphi]$, and hence
    \[
        \eta-C\leq P_{\theta}[\varphi].
    \]
    By \cref{lma:PKoutsidepps}, this implies that
    \[
        P_{\theta,K}[\varphi](f) \leq P_{\theta}[\varphi] + C.
    \]

    The proof of the $\mathcal{I}$-part is similar, but we spell it out since this comparison will be used below. If $\eta\in \PSH(X,\theta)$ satisfies $\eta\leq 0$ and $\eta\preceq_{\mathcal{I}}\varphi$, then $\eta+\inf_K f$ is a candidate in \eqref{eq:PthetaKvarphivI}. Hence
    \[
        P_{\theta}[\varphi]_{\mathcal{I}}+\inf_K f\leq P_{\theta,K}[\varphi]_{\mathcal{I}}(f).
    \]
    Conversely, let $\eta\in \PSH(X,\theta)$ be a candidate in \eqref{eq:PthetaKvarphivI}. Since $K$ is non-pluripolar and $f$ is bounded, the same compactness estimate gives a constant $C>0$, independent of $\eta$, such that $\eta\leq C$. Then $\eta-C\leq 0$ and $\eta-C\preceq_{\mathcal{I}}\varphi$. By \cref{cor:Ienvelopedef2},
    \[
        \eta-C\leq P_{\theta}[\varphi]_{\mathcal{I}}.
    \]
    Taking the upper envelope over all such $\eta$ gives
    \[
        P_{\theta,K}[\varphi]_{\mathcal{I}}(f)\leq P_{\theta}[\varphi]_{\mathcal{I}}+C.
    \]
    This completes the proof.
\end{proof}

\begin{corollary}\label{cor:PKtoPX}
    Let $f\colon K\rightarrow \mathbb{R}$ be a bounded function, and $\varphi\in \PSH(X,\theta)_{>0}$. Then
\begin{equation}\label{eq:PthetaKvarphiveq}
P_{\theta,K}[\varphi](f)=P_{\theta,X}[\varphi]\left( P_{\theta,K}[V_{\theta}](f) \right).
\end{equation}
\end{corollary}
\begin{proof}
Take any $\eta\in \PSH(X,\theta)$ with $\eta\preceq \varphi$ and $\eta|_K\leq f$ quasi-everywhere, then
\[
\eta\leq P_{\theta,K}[\varphi](f)\leq P_{\theta,K}[V_{\theta}](f).
\]
Therefore,
\[
\eta\leq P_{\theta,X}[\varphi]\left( P_{\theta,K}[V_{\theta}](f) \right).
\]
So we conclude the $\leq$ direction in \eqref{eq:PthetaKvarphiveq}.

For the reverse direction, take any candidate $\eta\in \PSH(X,\theta)$ for the right-hand side, so
\[
\eta\preceq \varphi,\quad \eta\leq P_{\theta,K}[V_{\theta}](f) \quad \textup{quasi-everywhere}.
\]
By \itemref{itm:prop:partial_env:2}, the intermediate envelope is at most $f$ quasi-everywhere on $K$. Combining the two inequalities outside the union of their pluripolar exceptional sets gives $\eta|_K\leq f$ quasi-everywhere. Hence
\[
\eta\leq P_{\theta,K}[\varphi](f).
\]
Taking the upper semicontinuous regularization of the supremum over these candidates gives the reverse inequality in \eqref{eq:PthetaKvarphiveq}.
\end{proof}

\section{Partial equilibrium measures}\label{sec:pem}
Let $X$ be a connected compact Kähler manifold of dimension $n$, and $\theta$ be a smooth closed real $(1,1)$-form on $X$ representing a pseudo-effective cohomology class.

\begin{definition}\label{def:weightedss}
A \emph{weighted subset}\index{weighted subset} of $X$ is a pair $(K,v)$ consisting of a closed non-pluripolar subset $K\subseteq X$ and a function $v\in C^0(K)$.
\end{definition}
The next result motivates our terminology to call the measures $\theta_{P_{\theta,K}[\varphi](f)}^n$ the \emph{partial equilibrium measures}\index{partial equilibrium measure} of our context:
\begin{lemma}\label{lma:balayage}
Let $(K,v)$ be a weighted subset of $X$, and $\varphi\in \PSH(X,\theta)_{>0}$.
Then
\begin{equation}\label{eq:intXminusKthetan_vanish}
\int_{X\setminus K}\theta_{P_{\theta,K}[\varphi](v)}^n=0,\quad \int_K \mathds{1}_{\left\{P_{\theta,K}[\varphi](v)|_K\neq v\right\}}\theta_{P_{\theta,K}[\varphi](v)}^n=0.
\end{equation}
More precisely, we have
\begin{equation}\label{eq:thetaPKuv}
\theta_{P_{\theta,K}[\varphi](v)}^n\leq \mathds{1}_{\left\{P_{\theta,K}[\varphi](v)|_K=P_{\theta,K}[V_{\theta}](v)|_K=v\right\}}\,\theta_{P_{\theta,K}[V_{\theta}](v)}^n.
\end{equation}
\end{lemma}
Here $\{P_{\theta,K}[\varphi](v)|_K\neq v\}$ and $\{P_{\theta,K}[\varphi](v)|_K=P_{\theta,K}[V_{\theta}](v)|_K=v\}$ are understood as consisting of points of $K$ satisfying these constraints.

\begin{proof}
If $n=0$, then $X=K$ is a point, both envelopes equal $v$, and their zeroth Monge--Ampère measures are the unit mass at that point. Both assertions follow. Assume henceforth that $n\geq 1$.

\textbf{Step~1}. We first prove \eqref{eq:intXminusKthetan_vanish} when $\varphi = V_{\theta}$.

Let $S\subseteq X$ be a proper analytic subset, such that $V_{\theta}$ is locally bounded on $X\setminus S$. Such $S$ exists by \cref{lma:Vthetalocbdd}.
Fix concentric coordinate balls $B\Subset B'\Subset X\setminus S$.

By \cref{prop:Choquet} and \cref{lma:PKoutsidepps}, we can take an increasing sequence $(\eta_j)_j$ in $\PSH(X,\theta)$ such that
\[
\eta_j\nearrow P_{\theta,K}[V_{\theta}](v)\textup{ a.e.},\quad \eta_j|_K\leq v\textup{ for all }j\geq 1.
\]
Note that $P_{\theta,K}[V_{\theta}](v)\sim V_{\theta}$ by \cref{lma:same_sing_type}. Take $C>0$ so that
\[
\left.\left( P_{\theta,K}[V_{\theta}](v) \right)\right|_K-C\leq v.
\]
Then we may replace $\eta_j$ by $\eta_j\lor (P_{\theta,K}[V_{\theta}](v)-C)$ and assume that $\eta_j\sim V_{\theta}$ for each $j$. In particular, $\eta_j$ is locally bounded on $X\setminus S$.

For each $j\geq 1$, let $\gamma_j\in \PSH(X,\theta)$ be the balayage of $\eta_j$ on $B$, defined in \cref{prop:bal_general}. Note that $(\gamma_j)_j$ is increasing due to \cref{prop:dom_prin_loc}.

For the first equality in \eqref{eq:intXminusKthetan_vanish},
we may assume that $B\Subset X\setminus (S \cup K)$. It suffices to show that $\theta_{P_{\theta,K}[V_{\theta}](v)}^n$ does not charge $B$.

Since $B$ lies outside $K$ by assumption, we also have $\gamma_j|_K=\eta_j|_K \leq v$ for all $j\geq 1$. In particular, $\gamma_j\leq P_{\theta,K}[V_{\theta}](v)$.

It follows that $\gamma_j \nearrow P_{\theta,K}[V_{\theta}](v)$ almost everywhere as well. Testing against non-negative continuous functions supported in $B$ and using \cref{thm:semicon}, we find that $\theta^n_{P_{\theta,K}[V_{\theta}](v)}$ does not charge $B$, as desired.

We next prove the second equality in \eqref{eq:intXminusKthetan_vanish}. The proof is essentially the same as that of \cref{lma:BT_freeenvelope}.

It suffices to show that $\theta^n_{P_{\theta,K}[V_{\theta}](v)}$ does not charge the set $K\cap \{P_{\theta,K}[V_{\theta}](v)|_K<v\}$.

Let $x\in (X\setminus S) \cap K$ satisfy $P_{\theta,K}[V_{\theta}](v)(x)<v(x)$. Choose $\epsilon>0$ with $P_{\theta,K}[V_{\theta}](v)(x)<v(x)-2\epsilon$. Choose $B$ and $B'$ centered at $x$, and a smooth local potential $g$ of $\theta$ on $B'$. By shrinking the balls we may assume that $\operatorname{osc}_{\overline B} g<\epsilon$ and that
\begin{enumerate}
    \item\label{itm:chapter-Bergman:L380:1} $P_{\theta,K}[V_{\theta}](v)|_{\overline B}<v(x)-2\epsilon$.
    \item\label{itm:chapter-Bergman:L380:2} $v|_{\overline B\cap K}>v(x)-\epsilon$.
\end{enumerate}
On $\partial B$, we have $\gamma_j+g=\eta_j+g\leq v(x)-2\epsilon+\sup_{\overline B} g$ by \ref{itm:chapter-Bergman:L380:1}. Applying the maximum principle to $\gamma_j+g$ and subtracting $g$ gives $\gamma_j|_B\leq v(x)-2\epsilon+\operatorname{osc}_{\overline B}g<v(x)-\epsilon$. By \ref{itm:chapter-Bergman:L380:2}, we have $\gamma_j|_{B\cap K}\leq v|_{B\cap K}$. On the other hand, $\gamma_j|_{K\setminus B}=\eta_{j}|_{K\setminus B}\leq v|_{K\setminus B}$. Therefore, $\gamma_j|_K\leq v$. It follows again that $\gamma_j\leq P_{\theta,K}[V_{\theta}](v)$. Therefore, we conclude as in the previous paragraph that $\gamma_j\nearrow P_{\theta,K}[V_{\theta}](v)$. Thus, testing against non-negative continuous functions supported in $B$, we conclude that $\theta^n_{P_{\theta,K}[V_{\theta}](v)}$ does not charge $B$ again by \cref{thm:semicon}.

\textbf{Step~2}. We handle the general case.

By \cref{cor:PKtoPX},
\[
P_{\theta,K}[\varphi](v) = P_{\theta,X}[\varphi]\left(P_{\theta,K}[V_{\theta}](v)\right).
\]
Now \cref{cor:Pvarphisupport} gives
\[
\begin{aligned}
\theta_{P_{\theta,K}[\varphi](v)}^n \leq & \mathds{1}_{\left\{P_{\theta,K}[\varphi](v) = P_{\theta,K}[V_{\theta}](v)\right\}} \theta_{P_{\theta,K}[V_{\theta}](v)}^n \\
\leq & \mathds{1}_{\left\{P_{\theta,K}[\varphi](v)|_K=P_{\theta,K}[V_{\theta}](v)|_K=v\right\}}\,\theta_{P_{\theta,K}[V_{\theta}](v)}^n,
\end{aligned}
\]
where in the second inequality we have used Step~1. Hence \eqref{eq:thetaPKuv} follows, and \emph{a fortiori} \eqref{eq:intXminusKthetan_vanish} follows as well.
\end{proof}

\begin{corollary}\label{cor:suppthetan}
Let $v\in C^0(K)$ and $\varphi\in \PSH(X,\theta)_{>0}$. Then
\begin{equation}\label{eq:thetandoesnotcharge1}
\begin{aligned}
\int_{(X\setminus K)\cup \left\{P_{\theta,K}[\varphi](v)<v\right\}}\theta^n_{P_{\theta,K}[\varphi](v)}=& 0,\\
\int_{(X\setminus K)\cup \left\{P_{\theta,K}[\varphi]_{\mathcal{I}}(v)<v\right\}}\theta^n_{P_{\theta,K}[\varphi]_{\mathcal{I}}(v)}=& 0.
\end{aligned}
\end{equation}
\end{corollary}

\begin{proof}
The first equation in \eqref{eq:thetandoesnotcharge1} follows from \cref{lma:balayage}. For the second, it suffices to observe that
\[
P_{\theta,K}[\varphi]_{\mathcal{I}}(v)=P_{\theta,K}\left[P_{\theta}[\varphi]_{\mathcal{I}}\right]_{\mathcal{I}}(v)=P_{\theta,K}\left[P_{\theta}[\varphi]_{\mathcal{I}}\right](v)
\]
by \cref{cor:I-good_partial_env}.
\end{proof}

\begin{definition}\label{def:partialequienergy}
Given $\varphi\in \PSH(X,\theta)_{>0}$, we define the \emph{partial equilibrium energy}\index{partial equilibrium energy} functional $\mathcal{E}_{\theta,K}^{\varphi}\colon C^0(K)\rightarrow \mathbb{R}$\index{$\mathcal{E}_{\theta,K}^{\varphi}$} as follows: given $v \in C^0(K)$, we set
\begin{equation}
\mathcal{E}_{\theta,K}^{\varphi}(v)\coloneqq E_{\theta}^{P_{\theta}[\varphi]_{\mathcal{I}}}\left(P_{\theta,K}[\varphi]_{\mathcal{I}}(v)\right).
\end{equation}
\end{definition}
Recall that the energy $E_{\theta}^{P_{\theta}[\varphi]_{\mathcal{I}}}$ functional is defined in \cref{def:MAenergy}.

Note that by \cref{lma:same_sing_type}, we have
\[
P_{\theta,K}[\varphi]_{\mathcal{I}}(v)\in \mathcal{E}^{\infty}\left(X,\theta;P_{\theta}[\varphi]_{\mathcal{I}}\right),
\]
so $\mathcal{E}_{\theta,K}^{\varphi}(v)\in \mathbb{R}$.

\begin{proposition}\label{prop:differential_P} Let $v,f \in C^0(K)$ and $\varphi \in \PSH(X,\theta)_{>0}$.
Then $\mathbb{R}\ni t \mapsto \mathcal{E}_{\theta,K}^{\varphi}(v+tf)$ is differentiable and
\begin{equation}\label{eq:ddtI}
\frac{\mathrm{d}}{\mathrm{d}t}\mathcal{E}_{\theta,K}^{\varphi}(v+tf) = \int_K f\, \theta^n_{P_{\theta,K}[\varphi]_{\mathcal{I}}(v+tf)}
\end{equation}
for all $t\in \mathbb{R}$.
\end{proposition}

\begin{proof}
We may assume that $\varphi$ is $\mathcal{I}$-model by replacing $\varphi$ by $P_{\theta}[\varphi]_{\mathcal{I}}$.

Note that it suffices to prove \eqref{eq:ddtI} at $t=0$, which is equivalent to
\begin{equation}\label{eq:to_prove_1}
\lim_{t\to 0}\frac{E_{\theta}^{\varphi}\left(P_{\theta,K}[\varphi]_{\mathcal{I}}(v+tf)\right)-E_{\theta}^{\varphi}\left(P_{\theta,K}[\varphi]_{\mathcal{I}}(v)\right)}{t}=\int_K f\,\theta^n_{P_{\theta,K}[\varphi]_{\mathcal{I}}(v)}.
\end{equation}
It suffices to prove the right derivative at $0$ for an arbitrary $f$; applying the same argument to $-f$ then gives the left derivative. In other words, we need to prove that
\begin{equation}\label{eq:to_prove_2}
\lim_{t\to 0+}\frac{E_{\theta}^{\varphi}\left(P_{\theta,K}[\varphi]_{\mathcal{I}}(v+tf)\right)-E_{\theta}^{\varphi}\left(P_{\theta,K}[\varphi]_{\mathcal{I}}(v)\right)}{t}=\int_K f\,\theta^n_{P_{\theta,K}[\varphi]_{\mathcal{I}}(v)}.
\end{equation}

We first reduce the proof of \eqref{eq:to_prove_2} to the case $f\leq 0$. Assume that this case is known. Choose $A\in\mathbb R$ such that $f-A\leq 0$ on $K$. Since for any $t>0$,
\[
P_{\theta,K}[\varphi]_{\mathcal I}\Bigl( v+t(f-A) \Bigr)= P_{\theta,K}[\varphi]_{\mathcal I}(v+tf)-tA,
\]
then we have
\[
    \int_K f\, \theta^n_{P_{\theta,K}[\varphi]_{\mathcal{I}}(v)}=\int_K (f-A)\, \theta^n_{P_{\theta,K}[\varphi]_{\mathcal{I}}(v)}+A\int_X \theta^n_{P_{\theta,K}[\varphi]_{\mathcal{I}}(v)}.
\]
By \eqref{eq:Ephiaddcst}, for any $t>0$,
\[
    E_{\theta}^{\varphi}\Bigl(P_{\theta,K}[\varphi]_{\mathcal{I}}\left(v+t(f-A)\right)\Bigr)=E_{\theta}^{\varphi}\left(P_{\theta,K}[\varphi]_{\mathcal{I}}(v+tf)\right)-tA\int_X\theta_{\varphi}^n.
\]
So the limit in \eqref{eq:to_prove_2} exists and is equal to the left-hand side of \eqref{eq:to_prove_2} with $f-A$ in place of $f$ plus $A\int_X\theta_{\varphi}^n$.
So \eqref{eq:to_prove_2} for $f-A$ implies the desired formula for $f$ in view of \cref{lma:same_sing_type} and the monotonicity theorem \cref{thm:mono}. Hence we may assume that $f\leq 0$.

Fix $t>0$.
Set
\[
\eta_t\coloneqq P_{\theta,K}[\varphi]_{\mathcal I}(v+tf),\quad \eta_0\coloneqq P_{\theta,K}[\varphi]_{\mathcal I}(v).
\]
By the comparison principle \eqref{eq:comp_princi1} and \cref{prop:cocycE1}, we find
\[
\begin{aligned}
E_{\theta}^{\varphi}(\eta_t)-E_{\theta}^{\varphi}(\eta_0)
={}&\frac{1}{n+1}\sum_{i=0}^n\int_X(\eta_t-\eta_0)\,
\theta_{\eta_t}^i\wedge\theta_{\eta_0}^{n-i}\\
\leq{}& \int_X(\eta_t-\eta_0)\,\theta_{\eta_0}^n.
\end{aligned}
\]
By \cref{lma:balayage} and \cref{cor:I-good_partial_env}, the measure $\theta_{\eta_0}^n$ is carried by $K\cap\{\eta_0=v\}$, and the measure $\theta_{\eta_t}^n$ is carried by $K\cap\{\eta_t=v+tf\}$.
On $K\cap\{\eta_0=v\}$, we have $\eta_t\leq v+tf=\eta_0+tf$, hence
\[
\eta_t-\eta_0\leq tf.
\]
On $K\cap\{\eta_t=v+tf\}$, we have $\eta_0\leq v$, hence
\[
\eta_t-\eta_0\geq tf.
\]
It follows that
\[
\int_X(\eta_t-\eta_0)\,\theta_{\eta_0}^n \leq t\int_K f\,\theta_{\eta_0}^n
\]
and
\[
\int_X(\eta_t-\eta_0)\,\theta_{\eta_t}^n \geq t\int_K f\,\theta_{\eta_t}^n.
\]
Thus, letting $t\to 0+$, we get the inequality
\[
\varlimsup_{t\to 0+}\frac{E_{\theta}^{\varphi}\left(P_{\theta,K}[\varphi]_{\mathcal{I}}(v+tf)\right)-E_{\theta}^{\varphi}\left(P_{\theta,K}[\varphi]_{\mathcal{I}}(v)\right)}{t}\leq\int_K f\,\theta^n_{P_{\theta,K}[\varphi]_{\mathcal{I}}(v)}.
\]
Similarly, for $t>0$, we have
\[
\begin{aligned}
&E_{\theta}^{\varphi}\left(P_{\theta,K}[\varphi]_{\mathcal{I}}(v+tf)\right)-E_{\theta}^{\varphi}\left(P_{\theta,K}[\varphi]_{\mathcal{I}}(v)\right)\\
\geq & \int_X \left(P_{\theta,K}[\varphi]_{\mathcal{I}}(v+tf)-P_{\theta,K}[\varphi]_{\mathcal{I}}(v)\right)\,\theta_{P_{\theta,K}[\varphi]_{\mathcal{I}}(v+tf)}^n\\
\geq & t\int_K f\,\theta_{P_{\theta,K}[\varphi]_{\mathcal{I}}(v+tf)}^n.
\end{aligned}
\]
Letting $t\to 0+$, we find
\[
\varliminf_{t\to 0+}\frac{E_{\theta}^{\varphi}\left(P_{\theta,K}[\varphi]_{\mathcal{I}}(v+tf)\right)-E_{\theta}^{\varphi}\left(P_{\theta,K}[\varphi]_{\mathcal{I}}(v)\right)}{t}\geq\int_K f\,\theta^n_{P_{\theta,K}[\varphi]_{\mathcal{I}}(v)}.
\]
Here we used the following elementary continuity consequence:
\[
P_{\theta,K}[\varphi]_{\mathcal{I}}(v)-t\|f\|_{C^0(K)}
\leq P_{\theta,K}[\varphi]_{\mathcal{I}}(v+tf)
\leq P_{\theta,K}[\varphi]_{\mathcal{I}}(v)
\]
for $t>0$, and hence $P_{\theta,K}[\varphi]_{\mathcal{I}}(v+tf) \nearrow P_{\theta,K}[\varphi]_{\mathcal{I}}(v)$ almost everywhere when $t\to 0+$.
Next we claim that
\[
    \int_K f\,\theta^n_{P_{\theta,K}[\varphi]_{\mathcal{I}}(v)}\leq \lim_{t\to 0+}  \int_K f\,\theta_{P_{\theta,K}[\varphi]_{\mathcal{I}}(v+tf)}^n.
\]
Choose a continuous extension $\tilde f\in C^0(X)$ of $f$. Since \cref{lma:balayage,cor:I-good_partial_env} show that all the measures appearing here are carried by $K$, the integrals over $K$ against $f$ are the same as the integrals over $X$ against $\tilde f$. Thus \cref{cor:incseqnppcont} applies.

We conclude \eqref{eq:to_prove_1}.
\end{proof}

\begin{lemma}\label{lma:global_env_approx} 
    Assume that $\{\theta\}$ represents a big cohomology class.
    Fix a Kähler form $\omega$ on $X$ and a weighted subset $(K,v)$ of $X$. There exists a uniformly bounded increasing sequence $(v^-_j)_j$ in $C^\infty(X)$ and a uniformly bounded decreasing sequence $(v^+_j)_j$ in $C^\infty(X)$, such that
    \begin{equation}\label{eq:unf_bd_vj}
    \sup_X |v_j^-|\leq C,\quad \sup_X |v_j^+|\leq C
    \end{equation}
    for some constant $C$ depending only on $\|v\|_{C^0(K)}$, $X$, $K$ and $\theta+\omega$. Furthermore, for any $\varphi \in \PSH(X,\theta)_{>0}$ and any $\delta\in [0,1]$ we have the following:
\begin{enumerate}
    \item\label{itm:lma:global_env_approx:1} We have
    \[
    P_{\theta+\delta \omega,X}[\varphi](v_j^+) \searrow P_{\theta+\delta \omega,K}[\varphi](v),
    \]
    and
    \[
    P_{\theta+\delta \omega,X}[\varphi](v^-_j) \nearrow P_{\theta+\delta \omega,K}[\varphi](v)
    \]
    almost everywhere.
    \item\label{itm:lma:global_env_approx:2} Suppose in addition that $\varphi$ is $\mathcal{I}$-good, then
    \[
        \lim_{j\to\infty}\mathcal{E}_{\theta,X}^{\varphi}(v_j^-)=\mathcal{E}_{\theta,K}^{\varphi}(v),\quad \lim_{j\to\infty}\mathcal{E}_{\theta,X}^{\varphi}(v_j^+)=\mathcal{E}_{\theta,K}^{\varphi}(v).
    \]
\end{enumerate}
\end{lemma}
\begin{proof}
    \textbf{Step~1}. First we construct $(v^-_j)_j$.

    Let
    \[
        C_{K,v} \coloneqq \sup\left\{\sup_X \eta : \eta \in \PSH(X,\theta+\omega),  \eta|_K \leq v\right\}.
    \]
    Since $K$ is non-pluripolar, we have $C_{K,v} \in \mathbb{R}$ by \cref{rmk:L1cpt}.
    Set
    \[
    M_{K,v}\coloneqq \max\{C_{K,v}+1,\sup_K v\}.
    \]
    Now define $\widetilde v \colon X \to \mathbb{R}$ as
    \[
    \widetilde{v}(x)=\left\{
    \begin{aligned}
        v(x),&\quad x\in K;\\
        M_{K,v},&\quad x\in X\setminus K.
    \end{aligned}
    \right.
    \]
    Since $M_{K,v}\geq \sup_K v$, the function $\widetilde v$ is lower semicontinuous. Therefore there exists an increasing and uniformly bounded sequence $(v^-_j)_j$ in $C^\infty(X)$ such that $v^-_j \nearrow \widetilde v$ and
    \[
        \|v_j^-\|_{C^0(X)}\leq \|\tilde{v}\|_{L^{\infty}(X)}+1.
    \]
    The right-hand side can clearly be bounded by a constant as in \eqref{eq:unf_bd_vj}.

    \textbf{Step~2}. We construct $(v_j^+)_j$.

    Since
\[
h \coloneqq P_{\theta + \omega,K}[V_{\theta+\omega}](v) \lor \left(\inf_K v -1\right)
\]
is usc, there exists a decreasing and uniformly bounded sequence $(v^+_j)_j$ in $C^\infty(X)$ such that $v^+_j \searrow h$. We may guarantee that
\[
\|v^+_j\|_{C^0(X)}\leq \|h\|_{L^{\infty}(X)}+1.
\]
Note that $\inf_K v-1\leq h\leq C_{K,v} \lor (\inf_K v-1)$. So \eqref{eq:unf_bd_vj} is satisfied.

\textbf{Step~3}. We finish the proof. Fix $\varphi\in \PSH(X,\theta)_{>0}$ and $\delta\in [0,1]$.

Note that \ref{itm:lma:global_env_approx:2} follows from \cref{cor:I-good_partial_env}, \cref{lma:same_sing_type}, \cref{prop:Econt_mono} and \ref{itm:lma:global_env_approx:1}. Let us prove \ref{itm:lma:global_env_approx:1}.

\textbf{Step~3.1}. We prove \ref{itm:lma:global_env_approx:1} for $(v^-_j)_j$.

Observe that $P_{\theta + \delta \omega,X}[\varphi](v^-_j)$ is increasing in $j\geq 1$, and
\[
P_{\theta + \delta \omega,X}[\varphi]\left(v^-_j\right)\leq P_{\theta + \delta \omega,K}[\varphi](v).
\]
We claim that
\begin{equation}\label{eq:Pthetapdeltaomega_temp1}
P_{\theta + \delta \omega,X}[\varphi]\left(v^-_j\right) \nearrow P_{\theta + \delta \omega,K}[\varphi](v)\quad \textup{a.e.}.
\end{equation}
Let $\eta\in \PSH(X,\theta+\delta\omega)$ be a potential such that $\eta\preceq \varphi$ and $\eta|_K\leq v$. Take some $\epsilon>0$. Then, since $\eta$ is upper semicontinuous, $\widetilde v$ is lower semicontinuous and $\eta < \epsilon+\widetilde v$, compactness gives $j_0>0$ such that $\eta < \epsilon+ v^-_j$ for $j \geq j_0$. Hence,
\[
\eta \leq \epsilon+P_{\theta + \delta \omega,X}[\varphi](v^-_j)
\]
for $j\geq j_0$.
Therefore,
\[
\eta\leq \uscsup[j\geq 1] P_{\theta + \delta \omega,X}[\varphi]\left(v^-_j\right)+\epsilon.
\]
Letting $\epsilon\to 0+$, we find
\[
\eta\leq \uscsup[j\geq 1] P_{\theta + \delta \omega,X}[\varphi]\left(v^-_j\right).
\]
Since $\eta$ is arbitrary, \eqref{eq:Pthetapdeltaomega_temp1} follows.

\textbf{Step~3.2}. We prove \ref{itm:lma:global_env_approx:1} for $(v^+_j)_j$.

Note that $P_{\theta + \delta \omega,X}[\varphi](v^+_j)$ is decreasing in $j$, and
\begin{equation}\label{eq:defchilimPthetadeltaomega}
\chi \coloneqq \lim_{j\to\infty} P_{\theta + \delta \omega,X}[\varphi]\left(v^+_j\right) \geq P_{\theta + \delta \omega,K}[\varphi](v).
\end{equation}
In particular, $\chi\in \PSH(X,\theta+\delta\omega)_{>0}$, since $P_{\theta + \delta \omega,K}[\varphi](v)$ has the same singularity types as $P_{\theta + \delta \omega}[\varphi]$ by \cref{lma:same_sing_type}.
Observe that
\[
P_{\theta + \delta \omega,X}[\varphi]\left(v_j^+\right)\leq v_j^+
\]
for all $j\geq 1$, since $v_j^+$ is continuous. Thus, $\chi\leq h$.
On the other hand, $h|_K\leq v$ quasi-everywhere. So in particular, $\chi|_{K}\leq v$ quasi-everywhere.
On the other hand,
\[
\chi\leq P_{\theta + \delta \omega,X}[\varphi]\left(v^+_1\right)\sim_P \varphi
\]
by \cref{lma:same_sing_type}. Hence, by \cref{cor:partial_P_new_def},
\[
\chi\leq P_{\theta + \delta \omega,K}[\varphi](v).
\]
Therefore, \eqref{eq:defchilimPthetadeltaomega} is an equality.
\end{proof}

\begin{proposition}\label{prop:conv_of_K_env}
    Let $(K,v)$ be a weighted subset of $X$. Let $(\varphi_j)_j$ be a sequence in $\PSH(X,\theta)_{>0}$ and $\varphi\in \PSH(X,\theta)_{>0}$. Then the following hold:
\begin{enumerate}
    \item\label{itm:prop:conv_of_K_env:1} If $(\varphi_j)_j$ is decreasing, and $\varphi_j\xrightarrow{d_S}\varphi$, then
    \begin{equation}\label{eq:PthetaKvarphij1}
    P_{\theta,K}[\varphi_j]_{\mathcal{I}}(v) \searrow P_{\theta,K}[\varphi]_{\mathcal{I}}(v),\quad P_{\theta,K}[\varphi_j](v) \searrow P_{\theta,K}[\varphi](v).
    \end{equation}
    \item\label{itm:prop:conv_of_K_env:2} If $\varphi_j \nearrow \varphi$ almost everywhere then
    \[
    P_{\theta,K}[\varphi_j]_{\mathcal{I}}(v) \nearrow P_{\theta,K}[\varphi]_{\mathcal{I}}(v)\textup{ a.e.},\quad P_{\theta,K}[\varphi_j](v) \nearrow P_{\theta,K}[\varphi](v) \textup{ a.e.}.
    \]
\end{enumerate}
\end{proposition}

\begin{proof}
Write $Q_{\theta,K}[\cdot]$ for either partial envelope and $Q_\theta[\cdot]$ for the corresponding global envelope. Set
\[
    \eta_j\coloneqq Q_{\theta,K}[\varphi_j](v),\quad
    \eta\coloneqq Q_{\theta,K}[\varphi](v).
\]
Suppose first that $(\varphi_j)_j$ is decreasing. By \cref{lma:same_sing_type}, $d_S(\eta_j,Q_\theta[\varphi_j])=0$ and $d_S(\eta,Q_\theta[\varphi])=0$. For the $\mathcal I$-envelope, \cref{thm:contPI} gives $Q_\theta[\varphi_j]\xrightarrow{d_S}Q_\theta[\varphi]$. For the ordinary envelope this follows from \cref{prop:ds0char} and $\varphi_j\xrightarrow{d_S}\varphi$. Hence $\eta_j\xrightarrow{d_S}\eta$.

By monotonicity, $\eta_j$ decreases to a $\theta$-psh function $\eta_\infty\geq\eta$. The masses of $\eta_j$ converge to the mass of $\eta$ by \cref{thm:convdS}. Since $\eta_j\geq\eta_\infty\geq\eta$, \cref{thm:mono} gives
\[
    \int_X\theta_{\eta_\infty}^n=\int_X\theta_\eta^n>0.
\]
By \cref{cor:decnetdS}, $\eta_j\xrightarrow{d_S}\eta_\infty$, so $\eta_\infty\sim_P\eta$ by \cref{prop:ds0char}. Also $\eta_\infty\leq v$ quasi-everywhere on $K$ since $\eta_\infty\leq\eta_1$. In the ordinary case $\eta_\infty\preceq_P\varphi$ by \cref{lma:same_sing_type}, so \cref{cor:partial_P_new_def} gives $\eta_\infty\leq\eta$. In the $\mathcal I$-case $\eta_\infty\preceq_{\mathcal I}\eta\preceq_{\mathcal I}\varphi$ by \cref{cor:PimpliesI} and \cref{lma:same_sing_type}; hence the same inequality follows from the defining candidates. This proves \ref{itm:prop:conv_of_K_env:1}.

Now suppose that $\varphi_j\nearrow\varphi$ almost everywhere. Then $\varphi_j\leq\varphi$ everywhere by \cref{prop:pshfuncdetdense}, and $\eta_j$ increases with $\eta_j\leq\eta$. By \cref{cor:incseqnppcont} the masses of $\varphi_j$ tend to the positive mass of $\varphi$. Thus \cref{lma:pathoenvelope} supplies $\epsilon_j\in(0,1)$ tending to zero and $\psi_j\in\PSH(X,\theta)_{>0}$ such that
    \[
        (1-\epsilon_j)\varphi+\epsilon_j\psi_j\leq\varphi_j.
    \]
Set $h_j=Q_{\theta,K}[\psi_j](v)$. Monotonicity and \itemref{itm:prop:partial_env:5} give
    \[
        (1-\epsilon_j)\eta+\epsilon_j h_j\leq\eta_j.
    \]
By \cref{cor:Ppartbdd_bdabove}, the $h_j$ are uniformly bounded from above. The potential $\psi_j-\sup_X\psi_j+\inf_K v$ is an ordinary candidate for $h_j$, so $\sup_X h_j\geq\inf_K v$. By \cref{prop:L1compa}, the $h_j$ have bounded $L^1$-norms. In particular, $\epsilon_j h_j\to 0$ in $L^1$. Along a subsequence this convergence holds almost everywhere, and the last displayed inequality gives $\eta\leq\sup_j\eta_j$ almost everywhere. This proves \ref{itm:prop:conv_of_K_env:2}.
\end{proof}

\section{Quantization of partial equilibrium measures}\label{sec:quant}

Let $X$ be a connected compact K\"ahler manifold of dimension $n$, let $L$ be a pseudo-effective line bundle on $X$, and let $h$ be a Hermitian metric on $L$. Set $\theta=c_1(L,h)$. Let $(T,h_T)$ be a Hermitian line bundle on $X$.

\subsection{Bernstein--Markov measures}
\label{subsec:BM}
Let $K\subseteq X$ be a closed non-pluripolar subset, and $v\colon K\rightarrow [-\infty,\infty)$ be a measurable function.

\begin{definition}\label{def:Bergman:L706}
    Let $\mu$ be a Borel probability measure on $K$.
We introduce the following functions on $\mathrm{H}^0(X,L^k \otimes T)$ ($k\geq 1$), with values possibly equal to $\infty$:
\index{$N^k_{v,\mu}$}
\index{$N^k_{v,K}$}
\begin{equation}\label{eq:Nkvmu_norm}
\begin{aligned}
N^k_{v,\mu}(s) \coloneqq & \left(\int_K h^k \otimes h_T(s,s) \mathrm{e}^{-kv} \,\mathrm{d}\mu\right)^{1/2},\\
N^k_{v,K}(s) \coloneqq & \sup_{K\setminus \{v=-\infty\}} \left( h^k \otimes h_T(s,s)\mathrm{e}^{-kv}\right)^{1/2}.
\end{aligned}
\end{equation}
\end{definition}
Observe that if $\mu$ puts no mass on $\{v=-\infty\}$ then we always have
\begin{equation}\label{eq:Nkinfcomp}
N^k_{v,\mu}(s)\leq N^k_{v,K}(s).
\end{equation}

We start with the following elementary observation:
\begin{lemma}\label{lma:mononorm}
Let $v_1,v_2\colon K\rightarrow [-\infty,\infty)$ be two measurable functions such that
\begin{enumerate}
    \item\label{itm:lma:mononorm:1} $v_1\leq v_2$;
    \item\label{itm:lma:mononorm:2} $\{v_1=-\infty\}=\{v_2=-\infty\}$.
\end{enumerate}
Then for any $s\in \mathrm{H}^0(X,L^k\otimes T)$ ($k\geq 1$), we have
\[
N_{v_1,K}^k(s)\geq N_{v_2,K}^k(s).
\]
\end{lemma}
\begin{proof}
    This follows immediately from the definition \eqref{eq:Nkvmu_norm}.
\end{proof}

\begin{definition}\label{def:BMmeasure}
Let $(K,v)$ be a weighted subset of $X$.
A Borel probability measure $\mu$ on $K$ is \emph{Bernstein--Markov}\index{Bernstein--Markov measure} with respect to $(K,v)$ if for each $\epsilon>0$, there is a constant $C_\epsilon>0$ such that
\begin{equation}\label{eq:BM}
     N^k_{v,K}(s) \leq C_\epsilon \mathrm{e}^{\epsilon k} N^k_{v,\mu}(s)
\end{equation}
for any $s \in \mathrm{H}^0(X,L^k \otimes T)$ and any $k\in \mathbb{Z}_{>0}$.
We write $\BM(K,v)$\index{$\BM(K,v)$} for the set of Bernstein--Markov measures with respect to $(K,v)$.
\end{definition}
As pointed out in \cite{BBWN11}, any smooth positive volume form of total mass one on $X$ is Bernstein--Markov with respect to $(X,v)$ for $v\in C^\infty(X)$.

\begin{proposition}\label{prop:BMNimplynorm}
Assume that $(K,v)$ is a weighted subset of $X$. Then
\begin{enumerate}
    \item\label{itm:prop:BMNimplynorm:1} $N^k_{v,K}$ is a norm on $\mathrm{H}^0(X,L^k\otimes T)$.
    \item\label{itm:prop:BMNimplynorm:2} For any $\mu\in \BM(K,v)$, $N^k_{v,\mu}$ is a norm on $\mathrm{H}^0(X,L^k\otimes T)$.
\end{enumerate}
\end{proposition}
\begin{proof}
\ref{itm:prop:BMNimplynorm:1}
As $v$ is bounded, $N^k_{v,K}$ is clearly finite on $\mathrm{H}^0(X,L^k\otimes T)$. In order to show that it is a norm, it suffices to show that for any $s\in \mathrm{H}^0(X,L^k\otimes T)$, $N^k_{v,K}(s)=0$ implies that $s=0$. In fact, we have $s|_K=0$, hence $s=0$ by the connectedness of $X$ and the assumption that $K$ is non-pluripolar.

\ref{itm:prop:BMNimplynorm:2}
As $v$ is bounded, clearly $N^k_{v,\mu}$ is finite and satisfies the triangle inequality. Non-degeneracy follows from the fact that $N^k_{v,K}$ is a norm and \eqref{eq:BM}.
\end{proof}

\subsection{Partial Bergman kernels}
\label{subsec:pBerg}
In this section, we fix a weighted subset $(K,v)$ of $X$ and $\mu\in \BM(K,v)$.

\begin{definition}\label{def:Bergman:L778}
    For any $\varphi\in \PSH(X,\theta)_{>0}$, we introduce the \emph{partial Bergman kernels}\index{partial Bergman kernels} of $\varphi$ (with respect to $(K,v)$) as follows: For any $k\geq 1$, we introduce
    \begin{equation}
    \begin{aligned}
     B^k_{v,\varphi, \mu}(x) \coloneqq \sup \left\{h^k \otimes h_T(s,s)\mathrm{e}^{-kv}(x):s \in \mathrm{H}^0\left(X,L^k \otimes T \otimes \mathcal{I}(k\varphi)\right),\right.\\
     \left. N^k_{v,\mu}(s) \leq 1 \right\},\quad x\in K.
     \end{aligned}
 \end{equation}
 \index{$B^k_{v,\varphi, \mu}$}
 We extend $B^k_{v,\varphi, \mu}$ to the whole $X$ by setting it to be $0$ outside $K$.

 The \emph{partial Bergman measures}\index{partial Bergman measures} of $\varphi$ (with respect to $(K,v)$) are defined as
 \begin{equation}
\beta^k_{v,\varphi,\mu} \coloneqq \frac{n!}{k^n} B^k_{v,\varphi,\mu} \,\mathrm{d}\mu
\end{equation}
\index{$\beta^k_{v,\varphi,\mu}$}
for each $k\geq 1$.
\end{definition}

\begin{remark}\label{rmk:Bergman:L797}
    Recall that in the definition of $\beta^k_{v,\varphi,\mu}$, the singular metric $\varphi$ only appears through the conditions $s \in \mathrm{H}^0(X,L^k \otimes T \otimes \mathcal{I}(k\varphi))$. The pointwise norm used to define $\beta^k_{v,\varphi,\mu}$ and $N^k_{v,\mu}$ is defined via $v$ not $\varphi$.

    This leads to the following seemingly paradoxical consequence, if $\psi\in \PSH(X,\theta)$ and $\varphi\preceq_{\mathcal{I}}\psi$, then
    \[
        \beta^k_{v,\varphi,\mu}\leq \beta^k_{v,\psi,\mu}\quad \forall\,k\geq 1.
    \]
    We shall apply this consequence without further mention in the sequel.
\end{remark}

Observe that
\[
    \int_X \beta^k_{v,\varphi,\mu}=\frac{n!}{k^n}h^0\left(X,T\otimes L^k \otimes \mathcal{I}(k\varphi)\right).
\]
In particular,
\begin{equation}\label{eq:betaint_conv}
    \lim_{k\to\infty}\int_X \beta^k_{v,\varphi,\mu}=\vol\left(\theta_{\varphi}\right)=\int_X \theta_{P_{\theta,K}[\varphi]_{\mathcal{I}}(v)}^n
\end{equation}
by \cref{thm:DXmain1} and \cref{lma:same_sing_type}. The goal of this section is to improve this result to the convergence of the Bergman kernels themselves:
\begin{theorem}\label{thm:pBMconvergence} Suppose that $\varphi \in \PSH(X,\theta)_{>0}$. Then
\begin{equation}\label{eq:pbkconvgeneral}
\beta^k_{v,\varphi,\mu} \rightharpoonup \theta_{P_{\theta,K}[\varphi]_{\mathcal{I}}(v)}^n,\quad \textup{as }k\to\infty.
\end{equation}
\end{theorem}

The proof will be given at the end of this subsection. We shall prove a few preliminaries at first.

Take a K\"ahler form $\omega$ on $X$ so that
\[
\int_X \omega^n = 1.
\]
\begin{proposition}\label{prop:smooth_weak_conv}
    Assume here that $K=X$.
    Let $\varphi \in \PSH(X,\theta)$ be a potential with analytic singularities such that $\theta_{\varphi}$ is a Kähler current. If $v \in C^\infty(X)$, then
\begin{equation}\label{eq:pbmconvanaly}
\beta^k_{v,\varphi,\omega^n} \rightharpoonup \theta_{P_{\theta,X}[\varphi]_{\mathcal{I}}(v)}^n=\theta_{P_{\theta,X}[\varphi](v)}^n\quad \textup{as }k\to\infty.
\end{equation}
\end{proposition}

\begin{proof} The equality part in \eqref{eq:pbmconvanaly} follows from \cref{cor:I-good_partial_env}. It remains to prove the weak convergence.

Let $\tau$ be the weak limit of a subsequence of $(\beta^k_{v,\varphi,\omega^n})_k$. In view of \eqref{eq:betaint_conv}, it suffices to show that
\begin{equation}\label{eq:muleq_temp1}
\tau\leq \theta_{P_{\theta,X}[\varphi](v)}^n.
\end{equation}

We can obtain an upper bound of $\tau$ via the observation:
\[
\beta^k_{v,\varphi,\omega^n} \leq \beta^k_{v,V_{\theta},\omega^n}.
\]
After replacing $\varphi$ by $V_{\theta}$, the corresponding problem is handled by \cite[Theorem~1.2]{Ber11}, giving us
\[
\beta^k_{v,V_{\theta},\omega^n} \rightharpoonup \theta_{P_{\theta,X}[V_{\theta}](v)}^n=\mathds{1}_{\left\{v = P_{\theta,X}[V_{\theta}](v)\right\}}\,\theta_v^n \quad \textup{as }k\to\infty,
\]
where the last identity is due to \cref{prop:DNTenvelope2}.

We obtain that
\begin{equation}\label{eq:Bergmanmeasure}
\tau \leq \mathds{1}_{\left\{v = P_{\theta,X}[V_{\theta}](v)\right\}}\,\theta_v^n.
\end{equation}
Now observe that it suffices to show that
\begin{equation}\label{eq:muP0_temp1}
\tau\left\{P_{\theta,X}[\varphi](v) < v\right\}=0.
\end{equation}
Assuming this, then \eqref{eq:Bergmanmeasure} further gives
\[
\tau \leq \mathds{1}_{\left\{P_{\theta,X}[\varphi](v)= v\right\}}\,\theta_v^n= \theta_{P_{\theta,X}[\varphi](v)}^n,
\]
where the equality follows from \cref{prop:DNTenvelope2}, and \eqref{eq:muleq_temp1} follows.

It now remains to argue \eqref{eq:muP0_temp1}. Fix a Kähler form $\omega$ on $X$. Observe that
\[
P_{\theta +\epsilon \omega,X}[\varphi](v) \sim \varphi
\]
for all $\epsilon\geq 0$ by \cref{lma:same_sing_type} and \cref{prop:analysingcompPandPI}. The family $P_{\theta+\epsilon\omega,X}[\varphi](v)$ is decreasing as $\epsilon\to0+$. Let
\[
    \chi\coloneqq \inf_{\epsilon>0}P_{\theta+\epsilon\omega,X}[\varphi](v).
\]
Then $\chi\in\PSH(X,\theta)$ since $\chi\geq P_{\theta,X}[\varphi](v)$. Moreover $\chi\leq v$ and $\chi\preceq\varphi$. Hence
\[
    \chi\leq P_{\theta,X}[\varphi](v).
\]
Thus, equality holds, and
\[
    P_{\theta+\epsilon\omega,X}[\varphi](v)\searrow P_{\theta,X}[\varphi](v).
\]
In other words,
\[
\left\{P_{\theta,X}[\varphi](v) < v\right\}=\bigcup_{\epsilon\in \mathbb{Q}_{>0}}\left\{P_{\theta+\epsilon\omega,X}[\varphi](v) < v\right\}.
\]
So in order to show \eqref{eq:muP0_temp1}, we may fix $\epsilon>0$ and argue that
\begin{equation}\label{eq:muP0_temp2}
\tau\left\{P_{\theta+\epsilon\omega,X}[\varphi](v) < v\right\}=0.
\end{equation}

Let $k\geq 1$, and $s \in \mathrm{H}^0(X,L^k \otimes T \otimes \mathcal I(k\varphi))$ be a section such that $N^k_{v,\omega^n}(s) \leq 1$. Then by \cite[Lemma~4.1]{Ber11}, there exists $C>0$ such that
\[
h^k\otimes h_T(s,s)\mathrm{e}^{-kv} \leq B^k_{v,\varphi,\omega^n} \leq B^k_{v,V_{\theta},\omega^n} \leq k^n C.
\]
This implies that
\[
\frac{1}{k}\log h^k \otimes h_T(s,s) \leq v + \frac{\log C }{k} + n \frac{\log k}{ k }.
\]
We define $\varphi_k$ as in \cref{prop:Bergman_approx}. Take $\alpha_k\nearrow 1$ as in \cref{prop:Bergman_approx}.
After adding a constant to $\varphi$, we may assume that $\varphi\leq 0$ on $X$. After replacing $\alpha_k$ by a subsequence, we also have $\alpha_k\varphi\in \PSH(X,\theta+\epsilon\omega)$ for all $k>0$, and $\alpha_k\varphi\searrow\varphi$.
Then
\[
\frac{1}{k}\log h^k \otimes h_T(s,s)\preceq \varphi_k\preceq \alpha_k\varphi.
\]
We notice that $\frac{1}{k}\log h^k \otimes h_T(s,s) \in \PSH(X, \theta + \epsilon \omega)$ for all $k \geq k_0(\epsilon)$.
Indeed,
\[
    \frac{1}{k}\log h^k\otimes h_T(s,s) \in \PSH\left(X,\theta+ k^{-1} c_1(T,h_T)\right)
\]
by \cref{prop:LelongPoincare}.

In particular,
\[
\frac{1}{k}\log h^k \otimes h_T(s,s)-\frac{\log C}{k}-n\frac{\log k}{k}\leq P_{\theta+\epsilon\omega,X}[\alpha_k \varphi](v).
\]
Now taking supremum over all candidates $s$, we obtain that
\begin{equation}\label{eq:smooth_Berg_est}
B^k_{v,\varphi,\omega^n} \leq C k^n \mathrm{e}^{k \left(P_{\theta +\epsilon \omega,X}[\alpha_k \varphi](v)-v\right)},\quad k \geq k_0.
\end{equation}

Now we prove \eqref{eq:muP0_temp2}. First observe that $\alpha_k\varphi\xrightarrow{d_S}\varphi$ by \cref{lma:dSsmallmult}. Thus, by \cref{prop:conv_of_K_env},
\[
P_{\theta +\epsilon \omega,X}[\alpha_k \varphi](v) \searrow P_{\theta +\epsilon \omega,X}[\varphi](v)
\]
and we get that
\[
\bigcup_{k=1}^{\infty}\left\{P_{\theta +\epsilon \omega,X}[\alpha_k \varphi](v) < v\right\}=\left\{P_{\theta +\epsilon \omega,X}[\varphi](v) < v\right\}.
\]
    It suffices to show that $\tau$ does not put mass on the set $\{P_{\theta +\epsilon \omega,X}[\alpha_m \varphi](v) < v\}$ for each fixed $m$. Fix a positive rational $\delta>0$, it suffices to show that $\tau$ puts no mass on $F\coloneqq \{P_{\theta +\epsilon \omega,X}[\alpha_m \varphi](v) < v-\delta\}$.
    Note that $F$ is open, and \eqref{eq:smooth_Berg_est} gives $\beta^\ell_{v,\varphi,\omega^n}(F)\to 0$. Therefore, $\tau$ puts no mass on $F$ as well, proving \eqref{eq:muP0_temp2}.
\end{proof}

\begin{definition}\label{def:Bergman:L942}
    Take $k\geq 1$ and $\varphi\in \PSH(X,\theta)$, let $\Nm \mathrm{H}^0(X,L^k \otimes T \otimes \mathcal{I}(k\varphi))$ be the space of norms on the vector space $\mathrm{H}^0(X,L^k \otimes T \otimes \mathcal{I}(k\varphi))$.

    Let $\mathcal{L}_{k,\varphi} \colon \Nm \mathrm{H}^0(X,L^k \otimes T \otimes \mathcal{I}(k\varphi)) \to \mathbb R$ be the \emph{partial Donaldson functional}\index{partial Donaldson functional}:
\begin{equation}
\mathcal{L}_{k,\varphi}(H) = \frac{n!}{k^{n+1}} \log \frac{\vol\{s \in \mathrm{H}^0\left(X,L^k \otimes T \otimes \mathcal{I}(k\varphi)\right) : H(s) \leq 1\}}{\vol \left\{s\in \mathrm{H}^0\left(X,L^k \otimes T \otimes \mathcal{I}(k\varphi)\right) : N^k_{0,\omega^n}(s) \leq 1\right\}},
\end{equation}
\index{$\mathcal{L}_{k,\varphi}$}
where $\vol$ is a Lebesgue measure on $\mathrm{H}^0(X,L^k \otimes T \otimes \mathcal{I}(k\varphi))$. The normalization does not matter.
\end{definition}

\begin{proposition}\label{prop:quant_I_smooth}
Let $w,w' \in C^0(X)$ and $\varphi \in \PSH(X,\theta)$ be a potential with analytic singularities such that $\theta_{\varphi}$ is a Kähler current. Then
\begin{equation}\label{eq:LdiffonXsmoothmeasure}
\lim_{k\to\infty} \left(\mathcal{L}_{k,\varphi}(N^k_{w,\omega^n}) - \mathcal{L}_{k,\varphi}(N^k_{w',\omega^n}) \right)= \mathcal{E}_{\theta,X}^{\varphi}(w)-\mathcal{E}_{\theta,X}^{\varphi}(w').
\end{equation}
In particular,
\begin{equation}\label{eq:LdiffonXsmoothmeasure2}
\lim_{k\to\infty}\mathcal{L}_{k,\varphi}(N^k_{w,\omega^n})=\mathcal{E}_{\theta,X}^{\varphi}(w)\,.
\end{equation}
\end{proposition}

\begin{proof}
First observe that by \cref{prop:BMNimplynorm}, for any $k\geq 1$, $N^k_{w,\omega^n}$ and $N^k_{w',\omega^n}$ are both norms, hence the expressions inside the limit in \eqref{eq:LdiffonXsmoothmeasure} make sense.

We first assume that $w,w'\in C^\infty(X)$. To start, we make the following observation:
\[
\begin{aligned}
\mathcal{L}_{k,\varphi}(N^k_{w,\omega^n}) - \mathcal{L}_{k,\varphi}(N^k_{w',\omega^n})&= \int_0^1 \frac{\mathrm{d}}{\mathrm{d}t} \mathcal{L}_{k,\varphi}(N^k_{w' + t(w-w'),\omega^n})\,\mathrm{d}t\\
& = \int_0^1 \int_X (w-w') \,\beta^k_{w' + t(w-w'),\varphi,\omega^n}\, \mathrm{d}t.
\end{aligned}
\]
The second equality follows from differentiating the Gram matrix of the Hilbert norm $N^k_{w'+t(w-w'),\omega^n}$:
\[
    \frac{\mathrm d}{\mathrm dt}\mathcal L_{k,\varphi}
    \left(N^k_{w'+t(w-w'),\omega^n}\right)=\int_X (w-w')\,\beta^k_{w'+t(w-w'),\varphi,\omega^n}.
\]

By \cref{prop:smooth_weak_conv}, we have
\[
\lim_{k\to\infty}\int_X (w-w') \,\beta^k_{w' + t(w-w'),\varphi,\omega^n}=\int_X (w-w')\,\theta_{P_{\theta,X}[\varphi](w' + t(w-w'))}^n.
\]
By \cref{thm:DXmain1}, we have
\[
\bigl|\int_X (w-w')\beta^k_{w' + t(w-w'),\varphi,\omega^n}\bigr| \leq C \sup_X |w - w'|.
\]
Hence, by the dominated convergence theorem we obtain that
\[
\begin{aligned}
\lim_{k\to\infty} \left(\mathcal{L}_{k,\varphi}(N^k_{w,\omega^n}) - \mathcal{L}_{k,\varphi}(N^k_{w',\omega^n})\right) &= \int_0^1 \int_X (w-w') \,\theta_{P_{\theta,X}[\varphi](w' + t(w-w'))}^n \,\mathrm{d}t\\
& = \mathcal{E}_{\theta,X}^{\varphi}(w)-\mathcal{E}_{\theta,X}^{\varphi}(w'),
\end{aligned}
\]
where in the last line we have used \cref{prop:differential_P}.

This proves \eqref{eq:LdiffonXsmoothmeasure} for smooth weights.

For general $w,w'\in C^0(X)$, choose smooth functions $w_j,w'_j$ converging uniformly to $w,w'$, respectively. If $\|w_j-w\|_{C^0}\leq \delta_j$, then
\[
\mathrm{e}^{-k\delta_j/2}N^k_{w,\omega^n}\leq N^k_{w_j,\omega^n}\leq \mathrm{e}^{k\delta_j/2}N^k_{w,\omega^n},
\]
and hence
\[
\left|\mathcal{L}_{k,\varphi}(N^k_{w_j,\omega^n})-\mathcal{L}_{k,\varphi}(N^k_{w,\omega^n})\right|
\leq \frac{n!}{k^n}h^0\left(X,L^k\otimes T\otimes \mathcal{I}(k\varphi)\right)\delta_j.
\]
The same estimate holds for $w'_j$ and $w'$. Moreover, \cref{prop:differential_P} shows that the energy side is Lipschitz with respect to the uniform norm of the weight. Letting first $k\to\infty$ in the smooth case, and then $j\to\infty$, proves \eqref{eq:LdiffonXsmoothmeasure} for continuous weights. Equation \eqref{eq:LdiffonXsmoothmeasure2} again follows by taking $w'=0$.
\end{proof}

\begin{lemma}\label{lma:BML}
Let $\varphi\in \PSH(X,\theta)$.
Let $(K,v)$ be a weighted subset of $X$. Let $\mu\in \BM(K,v)$. Then
\begin{equation}\label{eq:Bern_Mark_implies}
\lim_{k\to\infty} \big( \mathcal{L}_{k,\varphi}(N^k_{v,K}) - \mathcal{L}_{k,\varphi}(N^k_{v,\mu})\big) = 0\,.
\end{equation}
\end{lemma}
\begin{proof}
By \eqref{eq:Nkinfcomp} and \eqref{eq:BM}, for every $\epsilon>0$ and all $k\geq 1$,
\[
N^k_{v,\mu}\leq N^k_{v,K}\leq C_\epsilon \mathrm{e}^{\epsilon k}N^k_{v,\mu}.
\]
The corresponding unit balls therefore satisfy inclusions with radii differing by at most the factor $C_\epsilon \mathrm{e}^{\epsilon k}$. Taking logarithmic volumes and multiplying by $n!k^{-n-1}$ gives
\[
\left|\mathcal{L}_{k,\varphi}(N^k_{v,K})-\mathcal{L}_{k,\varphi}(N^k_{v,\mu})\right|
\leq \frac{2n!}{k^{n+1}} h^0\left(X,L^k\otimes T\otimes \mathcal{I}(k\varphi)\right)(\log C_\epsilon+\epsilon k).
\]
The right-hand side has limsup at most $2\epsilon\vol(\theta_\varphi)$ by \cref{thm:DXmain1}. Letting $\epsilon\to 0+$ proves \eqref{eq:Bern_Mark_implies}.
\end{proof}

\begin{corollary}\label{cor:Ninfdifflim}
Let $w \in C^0(X)$, $\varphi \in \PSH(X,\theta)$ be a potential with analytic singularities such that $\theta_\varphi$ is a Kähler current. Then
\[
\lim_{k\to\infty}\mathcal{L}_{k,\varphi}(N^k_{w,X})=\mathcal{E}_{\theta,X}^{\varphi}(w).
\]
\end{corollary}
\begin{proof}
For smooth $w$, the smooth volume form $\omega^n$ is Bernstein--Markov for $(X,w)$. The assertion then follows from \cref{lma:BML,prop:quant_I_smooth}. For continuous $w$, choose smooth $w_j\to w$ uniformly. The norm comparison
\[
    \mathrm{e}^{-k\|w_j-w\|_{C^0}/2} N^k_{w,X}
    \leq N^k_{w_j,X}
    \leq \mathrm{e}^{k\|w_j-w\|_{C^0}/2} N^k_{w,X}
\]
gives an error in $\mathcal{L}_{k,\varphi}$ bounded by $\|w_j-w\|_{C^0}$ times the normalized dimension $n!k^{-n}h^0(X,L^k\otimes T\otimes\mathcal{I}(k\varphi))$, which is bounded by \cref{thm:DXmain1}. The energy is continuous for uniform convergence of weights by \cref{prop:differential_P}. Letting first $k\to\infty$ and then $j\to\infty$ proves the assertion.
\end{proof}

It would be helpful to consider the following auxiliary functions:
\[
\begin{aligned}
    P'_{\theta,K}[\varphi](v)\coloneqq & \sup \Bigl\{ \eta \in \PSH(X,\theta): \eta|_K \leq v  \textup{ and } \eta\preceq \varphi\Bigr\},\\
    P'_{\theta,K}[\varphi]_{\mathcal{I}}(v)\coloneqq & \sup \Bigl\{ \eta \in \PSH(X,\theta): \eta|_K \leq v  \textup{ and } \eta\preceq_{\mathcal{I}} \varphi\Bigr\}.
\end{aligned}
\]
\index{$P'_{\theta,K}[\varphi](v)$}\index{$P'_{\theta,K}[\varphi]_{\mathcal{I}}(v)$}

We note the following maximum principles that follow from the above definitions:
\begin{lemma}\label{lma:max_princ} Let $v\in C^0(K)$, $\varphi\in \PSH(X,\theta)$, and $\eta \in \PSH(X,\theta)$.
Assume that $\eta \preceq \varphi$. Then
\begin{equation}\label{eq:max_princ}
\sup_K (\eta-v) = \sup_{\{\eta\neq -\infty\}} \left(\eta - P'_{\theta,K}[\varphi](v)\right)=\sup_{\left\{P'_{\theta,K}[\varphi](v)\neq -\infty\right\}} \left(\eta - P'_{\theta,K}[\varphi](v)\right).
\end{equation}
\end{lemma}
\begin{proof}
We prove the first equality first. We write $S=\{\eta=-\infty\}$.

By definition, $P'_{\theta,K}[\varphi](v)|_K\leq v$, so
\[
\left.\left(\eta-P'_{\theta,K}[\varphi](v)\right)\right|_{K\setminus S}\geq \eta|_{K\setminus S}-v|_{K\setminus S}\,.
\]
This implies that
\[
\sup_K (\eta-v) \leq \sup_{X\setminus S} \left(\eta - P'_{\theta,K}[\varphi](v)\right).
\]

Conversely, observe that $\sup_K(\eta-v)>-\infty$ as $K$ is non-pluripolar. It follows that $\eta-\sup_K (\eta-v)$ is a candidate in the definition of $P'_{\theta,K}[\varphi](v)$. Hence,
\begin{equation}\label{eq:etaminussupK_temp1}
\eta-\sup_K (\eta-v)\leq P'_{\theta,K}[\varphi](v).
\end{equation}
This implies that
\[
\sup_K (\eta-v) \geq \sup_{X\setminus S} \left(\eta - P'_{\theta,K}[\varphi](v)\right),
\]
finishing the proof of the first identity in \eqref{eq:max_princ}.

In order to prove the latter identity, observe that
\[
\left\{P'_{\theta,K}[\varphi](v) =-\infty\right\} \subseteq S
\]
as a consequence of \eqref{eq:etaminussupK_temp1}. On $\{P'_{\theta,K}[\varphi](v)\neq -\infty\}\cap S$, we have $\eta-P'_{\theta,K}[\varphi](v)=-\infty$. Hence the last equality of \eqref{eq:max_princ} also follows.
\end{proof}

\begin{proposition}\label{prop:quant_I_algebraic_BM} Let $\varphi\in \PSH(X,\theta)$ be a potential with analytic singularities such that $\theta_{\varphi}$ is a Kähler current. Let $(K,v)$, $(K',v')$ be two weighted subsets of $X$.
Then
\begin{equation}\label{eq:LkdiffconvtoI}
\lim_{k\to\infty} \left( \mathcal{L}_{k,\varphi}(N^k_{v,K}) - \mathcal{L}_{k,\varphi}(N^k_{v',K'}) \right) = \mathcal{E}_{\theta,K}^{\varphi}(v)-\mathcal{E}_{\theta,K'}^{\varphi}(v').
\end{equation}
In particular,
\begin{equation}\label{eq:Lkconv}
\lim_{k\to\infty}\mathcal{L}_{k,\varphi}(N^k_{v,K})=\mathcal{E}_{\theta,K}^{\varphi}(v).
\end{equation}
\end{proposition}

\begin{proof}
First observe that by \cref{prop:BMNimplynorm}, for any $k>0$, $N^k_{v,K}$ and $N^k_{v',K'}$ are both norms, hence the expressions inside the limit in \eqref{eq:LkdiffconvtoI} make sense. Moreover, \eqref{eq:Lkconv} is just a special case of \eqref{eq:LkdiffconvtoI} for $K'=X$ and $v'= 0$.

To prove \eqref{eq:LkdiffconvtoI} it is enough to show that for any fixed $w \in C^\infty(X)$ we have
\begin{equation}\label{eq:inproofLdiffsmw}
\lim_{k\to\infty} \big( \mathcal{L}_{k,\varphi}(N^k_{v,K}) - \mathcal{L}_{k,\varphi}(N^k_{w,X}) \big) = \mathcal{E}_{\theta,K}^{\varphi}(v) - \mathcal{E}_{\theta,X}^{\varphi}(w)\,.
\end{equation}

For $\epsilon \in (0,1)$ small enough, $\theta_{(1-\epsilon)\varphi}$ is still a Kähler current. Let us fix such $\epsilon$, along with an arbitrary $\epsilon' \in (0,1)$.

Let $(v^-_{j})_j$ and $(v^+_{j})_j$ be the sequences of smooth functions constructed in \cref{lma:global_env_approx} for the data $(K,v)$.

By \cref{prop:Bergman_approx} there exists $k_0(\epsilon,\epsilon') \in \mathbb N$ such that
\[
\frac{1}{k} \log h^k \otimes h_T (s,s) \preceq (1-\epsilon)\varphi,
\]
and $\frac{1}{k} \log h^k \otimes h_T (s,s) \in \PSH(X,\theta + \epsilon'\omega)$ for any $s \in \mathrm{H}^0(X,T\otimes L^k \otimes \mathcal{I}(k\varphi))$, as long as $k \geq k_0(\epsilon,\epsilon')$.

In particular, \cref{lma:max_princ} gives that
\[
\begin{aligned}
N^k_{P'_{\theta + \epsilon'\omega,K}[(1-\epsilon)\varphi](v),X}(s)=& N^k_{v,K}(s),\\
N^k_{P'_{\theta + \epsilon'\omega,X}[(1-\epsilon)\varphi](v_{j}^-),X}(s)=& N^k_{v_{j}^-,X}(s),\\
N^k_{P'_{\theta + \epsilon'\omega,X}[(1-\epsilon)\varphi](v_{j}^+),X}(s)=& N^k_{v^+_{j},X}(s).
\end{aligned}
\]
As
\[
P'_{\theta + \epsilon'\omega,X}[(1-\epsilon)\varphi](v_{j}^-) \leq P'_{\theta + \epsilon'\omega,K}[(1-\epsilon)\varphi](v) \leq P'_{\theta + \epsilon'\omega,X}[(1-\epsilon)\varphi](v_{j}^+),
\]
where the second inequality follows from \cref{lma:global_env_approx} and the fact that $P'=P$ for the global continuous obstacles $v_j^\pm$. By \cref{lma:mononorm} we have
\[
N^k_{v^+_{j},X}(s) \leq N^k_{v,K}(s) \leq N^k_{v_{j}^-,X}(s), \quad s \in \mathrm{H}^0(X,T \otimes L^k  \otimes \mathcal{I}(k\varphi)), k \geq k_0(\epsilon,\epsilon').
\]
Composing with $\mathcal{L}_{k,\varphi}$ we arrive at
\[
\mathcal{L}_{k,\varphi}(N^k_{v^-_{j},X})\leq  \mathcal{L}_{k,\varphi}(N^k_{v,K})\leq  \mathcal{L}_{k,\varphi}(N^k_{v^+_{j},X})\,, \quad k \geq k_0(\epsilon,\epsilon').
\]
For any $j>0$, by \cref{cor:Ninfdifflim} we get
\[
\begin{aligned}
\mathcal{E}_{\theta,X}^{\varphi}(v^-_{j}) - \mathcal{E}_{\theta,X}^{\varphi}(w)=& \lim_{k\to\infty} \left(\mathcal{L}_{k,\varphi}(N^k_{v^-_{j},X})-\mathcal{L}_{k,\varphi}(N^k_{w,X})\right)\\
\leq & \varliminf_{k\to\infty}\left(\mathcal{L}_{k,\varphi}(N^k_{v,K})- \mathcal{L}_{k,\varphi}(N^k_{w,X})\right)\\
\leq & \varlimsup_{k\to\infty}\left(\mathcal{L}_{k,\varphi}(N^k_{v,K})- \mathcal{L}_{k,\varphi}(N^k_{w,X})\right)\\
\leq & \lim_{k\to\infty} \left(\mathcal{L}_{k,\varphi}(N^k_{v^+_{j},X})-\mathcal{L}_{k,\varphi}(N^k_{w,X})\right)\\
=&  \mathcal{E}_{\theta,X}^{\varphi}(v^+_{j}) - \mathcal{E}_{\theta,X}^{\varphi}(w).
\end{aligned}
\]
Using \cref{lma:global_env_approx}, we can let $j \to \infty$ to arrive at
\[
\begin{split}
\mathcal{E}_{\theta,K}^{\varphi}(v)  -\mathcal{E}_{\theta,X}^{\varphi}(w) &\leq \varliminf_{k\to\infty}\left(\mathcal{L}_{k,\varphi}(N^k_{v,K})- \mathcal{L}_{k,\varphi}(N^k_{w,X})\right)\\
&\leq\varlimsup_{k\to\infty}\left(\mathcal{L}_{k,\varphi}(N^k_{v,K})- \mathcal{L}_{k,\varphi}(N^k_{w,X})\right)\\
& \leq \mathcal{E}_{\theta,K}^{\varphi}(v)  -\mathcal{E}_{\theta,X}^{\varphi}(w).
\end{split}
\]
Hence, \eqref{eq:inproofLdiffsmw} follows.
\end{proof}

\begin{corollary} \label{cor:LktoI}
Let $\varphi \in \PSH(X,\theta)$ be a potential with analytic singularities such that $\theta_{\varphi}$ is a Kähler current.
Then
\[
\lim_{k\to\infty}  \mathcal{L}_{k,\varphi}(N^k_{v,\mu})= \mathcal{E}_{\theta,K}^{\varphi}(v).
\]
\end{corollary}
\begin{proof}
Our claim follows from \cref{prop:quant_I_algebraic_BM} and \cref{lma:BML}.
\end{proof}

\begin{proposition}\label{prop:weakconvana}
Let $\varphi \in \PSH(X,\theta)$ be a potential with analytic singularities such that $\theta_{\varphi}$ is a Kähler current. Then
\[
\beta^k_{v,\varphi, \mu} \rightharpoonup \theta_{P_{\theta,K}[\varphi]_{\mathcal{I}}(v)}^n=\theta_{P_{\theta,K}[\varphi](v)}^n\quad \textup{as }k\to\infty.
\]
\end{proposition}

\begin{proof}
For $w\in C^0(X)$, let $f_k,g\colon \mathbb{R}\rightarrow \mathbb{R}$ ($k\geq 1$) be defined as
\[
f_k(t)\coloneqq \mathcal{L}_{k,\varphi}(N^k_{v+tw,\mu}),\quad g(t)\coloneqq \mathcal{E}_{\theta,K}^{\varphi}(v+tw).
\]
Since $\mu$ is a probability measure, $N^k_{v+tw,\mu}\leq N^k_{v+tw,K}$, so the corresponding unit balls give
\[
    f_k(t)\geq \mathcal{L}_{k,\varphi}(N^k_{v+tw,K}).
\]
By \cref{prop:quant_I_algebraic_BM}, the right-hand side converges to $g(t)$ for each $t$. Also, $f_k(0)\to g(0)$ by \cref{cor:LktoI}, using $\mu\in\BM(K,v)$. The function $f_k$ is concave by H\"older's inequality.
For each $k$, the function $f_k$ is differentiable and
\[
    f_k'(0)=\int_X w\,\beta^k_{v,\varphi,\mu}.
\]
This follows by differentiating the Gram matrix of the Hilbert norm $N^k_{v+tw,\mu}$. The function $g$ is differentiable at $0$ by \cref{prop:differential_P}, with
\[
    g'(0)=\int_K w\,\theta^n_{P_{\theta,K}[\varphi]_{\mathcal I}(v)}.
\]
For $t>0$, concavity bounds $f_k'(0)$ below by $(f_k(t)-f_k(0))/t$ and above by $(f_k(0)-f_k(-t))/t$. The preceding bounds therefore give
\[
    \frac{g(t)-g(0)}{t}\leq \varliminf_{k\to\infty} f_k'(0)
    \leq \varlimsup_{k\to\infty} f_k'(0)\leq \frac{g(0)-g(-t)}{t}.
\]
Letting $t\searrow 0$ gives $f_k'(0)\to g'(0)$,
which is equivalent to $\beta^k_{v, \varphi,\mu} \rightharpoonup \theta_{P_{\theta,K}[\varphi](v)}^n$, by \cref{prop:differential_P}.
\end{proof}

\begin{proposition} \label{prop:mainKahcurr}
Suppose that $\varphi \in \PSH(X,\theta)$ is such that $\theta_{\varphi}$ is a Kähler current.
Let $(K,v)$ be a weighted subset of $X$ and $\mu\in \BM(K,v)$. Then
\begin{equation}\label{eq:pbkconv1}
\beta^k_{v,\varphi,\mu} \rightharpoonup \theta_{P_{\theta,K}[\varphi]_{\mathcal{I}}(v)}^n\quad \textup{as } k\to\infty.
\end{equation}
\end{proposition}
\begin{proof} Let $\tau$ be the weak limit of a subsequence of $\beta^k_{v,\varphi,\mu}$. We claim that
\begin{equation}\label{eq:inproofmuleq}
\tau\leq \theta_{P_{\theta,K}[\varphi]_{\mathcal{I}}(v)}^n.
\end{equation}
Observe that this claim implies the conclusion. In fact, by \cref{thm:DXmain1}, we have equality of the total masses, so equality holds in \eqref{eq:inproofmuleq}. As $\tau$ is an arbitrary cluster point of the sequence $(\beta^k_{v,\varphi,\mu})_k$, we get \eqref{eq:pbkconv1}.

It remains to prove \eqref{eq:inproofmuleq}.
Let $(\varphi_j)$ be a quasi-equisingular approximation of $\varphi$ in $\PSH(X,\theta)$. We may assume that $\theta_{\varphi_j}$ is a K\"ahler current for all $j\geq 1$.
By \cref{cor:quasi-equichar}, we know that
\[
\varphi_j\xrightarrow{d_S}P_{\theta}[\varphi]_{\mathcal{I}}.
\]
The quasi-equisingular approximation is decreasing, so applying the decreasing case of \cref{prop:conv_of_K_env} to the limit $P_{\theta}[\varphi]_{\mathcal{I}}$, and using the definition of the relative $\mathcal{I}$-envelope, we get
\[
P_{\theta,K}\left[\varphi_j\right]_{\mathcal{I}}(v)\searrow P_{\theta,K}[\varphi]_{\mathcal{I}}(v).
\]
In particular,
\begin{equation}\label{eq:inproofeqmass}
\lim_{j\to\infty} \int_X \theta_{P_{\theta,K}\left[\varphi_j\right]_{\mathcal{I}}(v)}^n= \int_X  \theta_{P_{\theta,K}[\varphi]_{\mathcal{I}}(v)}^n.
\end{equation}
Observe that
\[
\beta^k_{v,\varphi,\mu} \leq \beta^k_{v,\varphi_j,\mu}
\]
for any $k \geq 1$.
As $\mu\in \BM(K,v)$,
by \cref{prop:weakconvana},
\[
\tau \leq \theta_{P_{\theta,K}[\varphi_j]_{\mathcal{I}}(v)}^n,
\]
for any $j \geq 1$ fixed.
By \eqref{eq:inproofeqmass} and the semicontinuity theorem \cref{thm:semicon}, we have
\[
\theta_{P_{\theta,K}\left[\varphi_j\right]_{\mathcal{I}}(v)}^n\rightharpoonup \theta_{P_{\theta,K}[\varphi]_{\mathcal{I}}(v)}^n
\]
as $j\to\infty$. Therefore, \eqref{eq:inproofmuleq} follows.
\end{proof}

\begin{proof}[Proof of \cref{thm:pBMconvergence}]
By \cref{lma:same_sing_type} and \cref{prop:Ienvprojection}, the three potentials below have the same multiplier ideals. Hence
\[
\begin{aligned}
\mathrm{H}^0\left(X,L^k \otimes T \otimes \mathcal{I}(k\varphi)\right) &= \mathrm{H}^0\left(X,L^k \otimes T \otimes \mathcal I(kP_{\theta}[\varphi]_{\mathcal{I}})\right)\\
 &= \mathrm{H}^0\left(X,L^k \otimes T \otimes \mathcal I\left(kP_{\theta,K}[\varphi]_{\mathcal{I}}(v)\right)\right).
 \end{aligned}
\]
This allows us to replace $\varphi$ with $P_{\theta,K}[\varphi]_{\mathcal{I}}(v)$.

By \cref{lma:kahcurrentposmass}, there exists a sequence $(\varphi_j)_j$ in $\PSH(X,\theta)$ such that $\varphi_j\nearrow \varphi$ a.e.\ and $\theta_{\varphi_j}$ is a Kähler current for each $j\geq 1$. This gives
\[
\beta^k_{v,\varphi_j,\mu} \leq \beta^k_{v,\varphi,\mu}.
\]
Let $\tau$ be the weak limit of a subsequence of $(\beta^k_{v,\varphi,\mu})_k$. Then by \cref{prop:mainKahcurr},
\[
\theta_{P_{\theta,K} [\varphi_j]_{\mathcal{I}}(v)}^n\leq \tau.
\]
By \cref{prop:conv_of_K_env} and \cref{cor:incseqnppcont}, we have
\[
\theta_{P_{\theta,K} [\varphi_j]_{\mathcal{I}}(v)}^n \rightharpoonup \theta_{P_{\theta,K} [\varphi]_{\mathcal{I}}(v)}^n.
\]
Hence,
\begin{equation}\label{eq:inproofmulower}
\theta_{P_{\theta,K} [\varphi]_{\mathcal{I}}(v)}^n\leq \tau.
\end{equation}
A comparison of total masses using \eqref{eq:betaint_conv} gives that equality holds in \eqref{eq:inproofmulower}. As $\tau$ is arbitrary, we obtain \eqref{eq:pbkconvgeneral}.
\end{proof}

\begin{remark}\label{rmk:Bergman:L1290}
    The results in this chapter can also be reformulated as the large deviation principle of a determinantal point process on $X$ using the Gärtner--Ellis theorem exactly as in \cite{Ber14}. We do not pursue this reformulation here.
\end{remark}

\Extrachap{Comments}
\label{chap:history}

\section*{A brief history}

Here we recall the origin of various results.

\textbf{\cref{chap:pshfunctions}}.

The notion of plurisubharmonic functions was introduced by Lelong \cite{Lel45} and Oka \cite{Oka42}, based on F.~Riesz's theory of subharmonic functions \cite{Rie26}. See \cite{Bre72} for an excellent introduction to the early history of the subject. We refer to \cite{Bre65} for the foundations of potential theory and \cite{GZ17} for the pluripotential theory.

The study of Josefson-type theorems was initiated by Josefson \cite{Jos78}, who proved \cref{thm:Josefson}.
The global Josefson theorem \cref{thm:gloJosefson} was due to Vu \cite{Vu19}. In the projective setting, it was due to Dinh--Sibony \cite{DS06} and in the K\"ahler setting, it was established by Guedj--Zeriahi \cite{GZ05}.

The extension theorem \cref{thm:GRexten} was proved in \cite{GR56}. In fact, they proved a more general version for complex spaces, see \cref{thm:GRexten2}. For some related important extension theorems, see \cite{Shi72, Wan22, CGZ13}.

\cref{prop:Kis} was due to Kiselman \cite{Kis78}.

The idea of fine topology was due to H.~Cartan \cite{Car41}.
The plurifine topology was introduced by Fuglede during the Séminaire d'analyse de Lelong--Dolbeault--Skoda of the year 1983/1984 \cite{Lel86}, but often wrongly attributed to Bedford--Taylor. The key result \cref{thm:pf_basis} was claimed in Bedford--Taylor's work \cite[Theorem~2.3]{BT87} without proof. The first rigorous proof was given by El~Marzguioui--Wiegerinck \cite{EMW06}. A weaker result was proved earlier in \cite[Theorem~4.8.7]{Kl91}.

Results in \cref{subsec:pftmfd} are certainly well-known and are already implicitly used in the literature. I could not find the proofs in the literature, so all details are presented.

The semicontinuity theorem \cref{thm:Siusemi} was due to Siu \cite{Siu74}.

Multiplier ideal sheaves were introduced by Nadel \cite{Nad90}.
The idea of \cref{thm:valuativemulti} first appeared in the ground-breaking work of Boucksom--Favre--Jonsson \cite{BFJ08}.

The openness conjecture was formulated and studied by Demailly--Koll\'ar \cite{DK01} and later proved by Berndtsson \cite{Bern13b}. The strong openness \cref{thm:strongopen} was first established by Guan--Zhou \cite{GZ15}.

\cref{lma:IandIinf} was due to \cite[Proposition~4.1.6]{Dem15}.

The notions of positivity of currents \cref{def:pos_curr} are due to Lelong \cite{Lel68} and Harvey--Knapp \cite{HK74}, although it is often mistakenly attributed to Demailly.

The $L^2$ extension theorems such as \cref{thm:OT} were first studied by Ohsawa--Takegoshi \cite{OT87}. Since then, it has become a very active area. It is almost impossible to give a complete list of the various generalizations and improvements.

\textbf{\cref{chap:npp}}

The Monge--Amp\`ere operators for bounded plurisubharmonic functions were introduced by Bedford--Taylor \cite{BT76, BT82}. The non-pluripolar product was due to Bedford--Taylor \cite{BT87}, Guedj--Zeriahi \cite{GZ07} and Boucksom--Eyssidieux--Guedj--Zeriahi \cite{BEGZ10}.

\cref{ex:1Dnpp} was claimed in \cite{BEGZ10} right after Definition~1.1 without proof. The first proof was essentially due to El Kadiri in \cite{EK24}.

The semicontinuity theorem \cref{thm:semicon} and the domination principle \cref{thm:dom_prin_1} were due to \cite{DDNL18fullmass}. The monotonicity theorem was due to \cite{WN19, DDNL18mono}.

The key lemma \cref{lma:pathoenvelope} was proved in \cite{DDNLmetric}. \cref{thm:concentration} was due to \cite{DDNLsurv}.

\textbf{\cref{chap:envelopes}}

\cref{lma:rooftopMA} was proved in \cite[Lemma~3.7]{DDNL18mono}.

The notion of the $P$-envelope is due to Ross--Witt Nystr\"om \cite{RWN14} based on the ideas of Rashkovskii--Sigurdsson \cite{RS05}.

\cref{thm:Pvarphisupport} is due to \cite[Theorem~3.8]{DDNL18mono}. The diamond inequality \cref{thm:diamond} and \cref{prop:landfinitecond1} are due to \cite[Theorem~5.4]{DDNLmetric}. Most results in \cref{subsec:fullmass} are simple generalizations of the corresponding results in \cite{DDNL18big, DDNL18fullmass}.

The study of finite energy potentials was initiated by Cegrell \cite{Ceg98} in the local setup; see also his general definition of the complex Monge--Amp\`ere operator \cite{Ceg04}. The global setup was first studied by Guedj--Zeriahi \cite{GZ07}. The variational approach to degenerate complex Monge--Amp\`ere equations in big cohomology classes was developed by Boucksom--Eyssidieux--Guedj--Zeriahi \cite{BEGZ10} and Berman--Boucksom--Guedj--Zeriahi \cite{BBGZ13}. For the complex Monge--Amp\`ere equation and pluripotential methods for weak solutions, see Ko{\l}odziej \cite{Kol98, Kol05} and the references therein.

The $\mathcal{I}$-envelope was introduced by Darvas--Xia \cite{DX22}, inspired by the works of Dano Kim \cite{Kim15} and Boucksom--Favre--Jonsson \cite{BFJ08}. The notion of $\mathcal{I}$-model singularities was first formulated in the explicit way in \cite{DX22} in 2020, although it was already essentially known in Boucksom--Jonsson's work. In fact, they correspond exactly to the homogeneous non-Archimedean potentials assuming that the relevant masses do not vanish. A less explicit equivalent formulation of $\mathcal{I}$-model potentials also appeared in \cite{Dem15}.
A few months later, the same notion was rediscovered by Trusiani \cite{Tru22}.

\cref{thm:DNTcontact} is due to \cite{DNT21}.

\textbf{\cref{chap:rays}}

The study of geodesics in the space of K\"ahler metrics was initiated by Mabuchi \cite{Mab87}, Semmes \cite{Sem92} and Donaldson \cite{Don99}. The existence of geodesics with smooth boundary data was established by Chen \cite{Chen00}.

The metric geometry and Mabuchi completion of finite-energy classes in the K\"ahler case were developed by Darvas \cite{Dar15, Dar17}. The notion of weak geodesics was studied in detail by Darvas \cite{Da17} in the K\"ahler case.

The case of general big classes was partly handled in \cite{DDNL18fullmass}, \cite{DDNL18big}.
However, the key fact that the geodesics between two full mass potentials have the correct limit at the end points does not seem to have been proved in any references. We give a proof in \cref{prop:perronenvissubgeod}. We also extend the relevant results to the relative setting.

In the special case $\phi=V_\theta$, \cref{prop:energylinear} was proved in \cite[Theorem~3.12]{DDNL18fullmass}, while the sequential case of \cref{prop:supsgeod} was proved in \cite[Proposition~3.3]{DDNL18fullmass}.

Previously, \cref{prop:geodsupsublinear} and \cref{prop:raysuplinear} were only known in the K\"ahler case.

Most results in \cref{sec:metE1} are simple extensions of \cite{DDNL18big}.

\textbf{\cref{chap:toric_ample}}

Classical toric K\"ahler geometry is rooted in the works of Atiyah \cite{Ati82}, Delzant \cite{Del88}, Guillemin \cite{Gui94}, Abreu \cite{Abr98} and Donaldson \cite{Don02}. The pluripotential-theoretic toric framework used here was first written down by Berman--Berndtsson \cite{BB13} and Coman--Guedj--Sahin--Zeriahi \cite{CGSZ19}.

The beautiful theorem \cref{thm:Yao} was first proved by Yi Yao, who did not publish the result. Later on, a new proof was found by Botero--Burgos Gil--Holmes--de Jong \cite{BBGHdJ21}. We chose to present the approach of Yao, which integrates naturally with our framework.

\Cref{prop:Etoric} is an extension of \cite[Proposition~2.9]{BB13}.

\textbf{\cref{chap:comp}}

The notion of $P$-preorder is new, as well as most results in \cref{sec:PIpartialorder}.

The $d_S$-pseudometric was introduced in \cite{DDNLmetric}. The basic properties are proved in \cite{DDNLmetric} and \cite{Xia21}.

\cref{ex:BBJ} was due to Berman--Boucksom--Jonsson \cite{BBJ21}.

\cref{thm:dSadditivity} can be viewed as a vast generalization of \cite{Cao15}, where a question raised by Demailly was confirmed.

\cref{thm:Lelongcont} is proved in \cite{Xia22}. \cref{thm:equising_cond_general} and \cref{thm:contvolu} appear to be new. These results appeared previously in the form of lecture notes.

\textbf{\cref{chap:Igood}}

The notion of $\mathcal{I}$-good singularities was due to \cite{DX21}. The name \emph{$\mathcal{I}$-good} was chosen in \cite{Xia22}.

The key results \cref{thm:charIgoodasclosure} and \cref{thm:DXmain1} are due to \cite{DX21, DX22}. Previously, Rashkovskii \cite{Ra13} also studied the local approximation problem with a different motivation.

There are some further examples of $\mathcal{I}$-good singularities provided by \cite{BBGHdJ21} with applications in the theory of modular forms in \cite{BBGHdJ22}.

\textbf{\cref{chap:trace}}

The trace operator was introduced in \cite{DX24}. Here we present a different point of view. \cref{thm:rest_volume} was proved in \cite{DX24}.

The analytic Bertini theorem \cref{thm:Bert} was proved in \cite{XiaBer}, based on the works of Matsumura--Fujino \cite{FM21} and Fujino \cite{Fuj23}. A weaker result was established by Meng--Zhou \cite{MZ23}.

The local analogue of the analytic Bertini theorem is proved in \cite{XiaBerII}.

\textbf{\cref{chap:testcurve}}

The technique of test curves originates from \cite{RWN14}. It was generalized by Darvas--Di Nezza--Lu \cite{DDNL18big}, \cite{DX21}, \cite{DZ22} and \cite{DXZ23}.
We give the full details of the proofs.

Test curves in \cref{def:testcur} are called \emph{maximal test curves} in the literature, a terminology which I do not like. I prefer to call the usual notion of test curves in the literature \emph{sub-test curves}.

\cref{prop:mis_RWN} was first proved by He--Testorf--Wang in \cite{HTW23}. \cref{prop:His_mass_tc} was due to Hisamoto \cite{His16}.

\cref{rmk:max_filt} was a folklore result. I am unaware of any written proof in the literature before our paper \cite{DX22}. Finski \cite{Fin22b} also gave a different proof with different techniques later on.

\cref{def:maximalgeoray} was not the original definition of maximal geodesic rays of Berman--Boucksom--Jonsson in \cite{BBJ21}. One of the first major applications of our theory was this pluripotential-theoretical characterization of maximal geodesic rays, as proved in our very first paper \cite{DX22}.

\cref{ex:RWN} is essentially due to Ross--Witt Nystr\"om \cite{RWN14}. The Phong--Sturm geodesics were first constructed in \cite{PS07}. For the related link between test configurations, filtrations and Okounkov bodies, see Boucksom--Chen \cite{BC11} and Witt Nystr\"om \cite{WN12}. For filtrations and generalized test configurations in K-stability, see Sz\'ekelyhidi \cite{Sze15}.

Results in \cref{sec:operationtc} are easy generalizations of the results proved in \cite{Xia23Operations}.

\textbf{\cref{chap:Okounkov}}

The algebraic theory of partial Okounkov bodies was developed in \cite{Xia21}.
The transcendental Okounkov body was first defined by Deng \cite{Deng17} as suggested by Demailly. The volume identity was proved in \cite{DRWNXZ}. The transcendental theory of partial Okounkov bodies is new.

\textbf{\cref{chap:bdiv}}

The notion of b-divisors was introduced by Shokurov \cite{Sho93} in the context of minimal model programs.

The applications of b-divisors in pluripotential theory began with \cite{BFJ09}. The intersection theory of nef b-divisors was introduced by Dang--Favre \cite{DF20}. The technique of singularity b-divisors was introduced in \cite{XiaPPT} in 2020. The general form first appeared in \cite{Xia22}. One year later, a special case was rediscovered in \cite{BBGHdJ21}.

The current chapter reproduces a large part of \cite{Xiabdiv}. For further developments, we refer to \cite{Xiabdiv2}.

\cref{ex:volbdivnetstrict} is reproduced from \cite[Example~3.3]{DF20}.

\textbf{\cref{chap:toric_big}}

The whole chapter appears to be new. The study of toric pluripotential theory on big line bundles was made possible by the development of partial Okounkov bodies. The key result is \cref{thm:FVtheta}.

Most results in this chapter resulted from discussions with Yi Yao.

\textbf{\cref{chap:NApotential}}

Most results in this chapter come from \cite{Xia23Operations}. Results from \cref{sec:DHmeasure} are new, although the main idea was already contained in \cite{Xia21}. The terminology of Duistermaat--Heckman measures goes back to the classical theorem of Duistermaat--Heckman \cite{DH82}; see also their addendum \cite{DH83}. In the setting of K\"ahler geometry, the Duistermaat--Heckman measures were studied in detail by Boucksom--Hisamoto--Jonsson \cite{BHJ17}.

Foundational non-Archimedean global pluripotential theory was developed by Boucksom--Favre--Jonsson \cite{BFJ16b} and, over a trivially valued field, by Boucksom--Jonsson \cite{BJ22}.

\cref{thm:DXZ} is due to \cite{DXZ23}.
An alternative approach to the transcendental theory is due to Mesquita-Piccione \cite{MP24}. For further progress in this direction, we refer to \cite{MPWN25, WN25}.

Special cases of the results in this section have been applied to study K-stability, see \cite{XiaPPT}, \cite{DZ22}, \cite{DXZ23} and \cite{DR22}. In \cite{DX22}, we established the bijective correspondence between a class of $\mathcal{I}$-model test curves and the maximal geodesic rays in the sense of \cite{BBJ21}.

\textbf{\cref{chap:Bergman}}

The study of asymptotic expansions of Bergman kernels in K\"ahler geometry was initiated by Tian \cite{Ti90} with later improvements due to Bouch\'e \cite{Bou90}, Ruan \cite{Rua98}, Catlin \cite{Cat99}, Zelditch \cite{Zel98}, Lu \cite{Lu00} among many others.

The special case of \cref{thm:pBMconvergence} without the prescribed singularity $\varphi$ was due to Berman--Boucksom--Witt Nystr\"om, see \cite{BB10}, \cite{BBWN11}. The general case is due to \cite{DX21}.

\appendix
\chapter{Convex functions and convex bodies}
\label[appendix]{chap:convex}

This appendix collects the basic facts about convex functions and convex bodies that are used throughout the book. Our standard reference is \cite{Roc70}. All results stated for convex functions apply equally to concave functions, after replacing the function by its negative.

\section{The notion of convex functions}
\label{sec:convex_notion}
Let $N$ be a real vector space of finite dimension.

\begin{definition}
    Let $F\colon N\rightarrow [-\infty,\infty]$ be a function. The \emph{epigraph}\index{epigraph} of $F$ is defined as the following set
    \[
    \epi F\coloneqq \left\{ (n,r)\in N\times \mathbb{R}: r\geq F(n) \right\}.
    \]
\end{definition}
\begin{definition}\label{def:convexfun}
    A \emph{convex function}\index{convex function} on $N$ is a function $F\colon N\rightarrow [-\infty,\infty]$ such that the epigraph $\epi F$ is a convex subset of $N\times \mathbb{R}$.

    The \emph{effective domain}\index{convex function!effective domain} of $F$ is the set
    \[
    \Dom F\coloneqq \left\{ n\in N : F(n)<\infty \right\}.
    \]
    A convex function $F$ on $N$ such that $\Dom F\neq \varnothing$ and $F(n)\neq -\infty$ for all $n\in N$ is said to be \emph{proper}\index{convex function!proper convex function}.

    The set of convex functions on $N$ is denoted by $\Conv(N)$. The subset of proper convex functions is denoted by $\Conv^{\prop}(N)$.
\end{definition}

The following characterization of convex functions is well-known.
\begin{lemma}\label{lma:charconvex}
    Let $F\colon N\rightarrow [-\infty,\infty]$. Then $F$ is convex if and only if the following condition holds: suppose that $n,r\in N$ and $a,b\in \mathbb{R}$ such that $a>F(n)$, $b>F(r)$, then for any $t\in (0,1)$, we have
    \[
    F(tn+(1-t)r)<ta+(1-t)b.
    \]
\end{lemma}
See \cite[Theorem~4.2]{Roc70} for the proof.

\begin{example}\label{ex:charconvex}
    Let $A\subseteq N$ be a convex subset. Then the \emph{characteristic function}\index{characteristic function}\index{$\chi_A$} $\chi_A\colon N\rightarrow \{0,\infty\}$ of $A$ is defined by
    \[
    \chi_A(n)\coloneqq
    \left\{
    \begin{aligned}
        0,&\quad n\in A;\\
        \infty,&\quad n\not\in A.
    \end{aligned}
    \right.
    \]
    The function $\chi_A$ lies in $\Conv(N)$.
\end{example}

\begin{example}\label{ex:suppfun}
    Let $M$ be the dual vector space of $N$ and $P\subseteq M$ be a convex subset. The \emph{support function}\index{support function} $\Supp_P\in \Conv(N)$\index{$\Supp_P$} of $P$ is defined as follows:
    \[
        \Supp_P(n)\coloneqq \sup\{ \langle m,n \rangle: m\in P \}.
    \]
\end{example}

It is well-known that convexity is preserved by a number of natural operations. We recall a few to fix the notation.

\begin{definition}\label{def:convex:L60}
    Let $F_1,\ldots,F_m\in \Conv^{\prop}(N)$ ($m\in \mathbb{Z}_{>0}$). We define their \emph{infimal convolution}\index{infimal convolution} $F_1 \square \cdots \square F_m\in \Conv(N)$\index{$F_1 \square \cdots \square F_m$} as follows:
    \[
    F_1 \square \cdots \square F_m(n)\coloneqq \inf\left\{ \sum_{i=1}^m F_i(n_i):n_i\in N, \sum_{i=1}^m n_i=n \right\}.
    \]
\end{definition}
The fact $F_1 \square \cdots \square F_m\in \Conv(N)$ is proved in \cite[Theorem~5.4]{Roc70}. One should note that $F_1 \square \cdots \square F_m$ is not always proper.

\begin{proposition}\label{prop:supconv}
    Let $\{F_i\}_{i\in I}$ be a non-empty family in $\Conv(N)$. Then $\sup_{i\in I} F_i\in \Conv(N)$.
\end{proposition}
This follows from \cite[Theorem~5.5]{Roc70}. In particular, this allows us to introduce
\begin{definition}\label{def:LCE}
    Let $f\colon N\rightarrow [-\infty,\infty]$. The \emph{lower convex envelope}\index{lower convex envelope}\index{$\LCE$} of $f$ is defined as
    \[
    \LCE f\coloneqq \sup\{F\in \Conv(N): F\leq f \}.
    \]
    It follows from \cref{prop:supconv} that $\LCE f\in \Conv(N)$.
\end{definition}

\begin{definition}\label{def:convwedge}
Given a non-empty family $\{F_i\}_{i\in I}$ in $\Conv(N)$, we define
\index{convex function!rooftop}
\index{$\bigwedge_i F_i$}
\[
\bigwedge_{i\in I}F_i\coloneqq \LCE \left(\inf_{i\in I} F_i\right).
\]
When the family $I$ is finite, say $I=\{1,\ldots,m\}$, we also write
\[
F_1\land \cdots\land F_m=\bigwedge_{i\in I}F_i.
\]
\end{definition}
\begin{definition}\label{def:convlor}
   Given a non-empty family $\{F_i\}_{i\in I}$ in $\Conv(N)$, we define
\index{convex function!maximum}
\index{$\bigvee_i F_i$}
\[
\bigvee_{i\in I}F_i\coloneqq \sup_{i\in I} F_i.
\]
When the family $I$ is finite, say $I=\{1,\ldots,m\}$, we also write
\[
F_1\lor \cdots\lor F_m=\bigvee_{i\in I}F_i.
\]
\end{definition}
Recall that $\bigvee_{i\in I}F_i\in \Conv(N)$ by \cref{prop:supconv}.

\begin{proposition}\label{prop:concavhull}
    Let $F_1,\ldots,F_m\in \Conv^{\prop}(N)$. Then
    \[
    \begin{split}
        F_1\land \cdots\land F_m(x)=\inf\left\{\sum_{i=1}^m \lambda_i F_i(x_i):x_i\in \Dom(F_i),\right.\\
        \left.\lambda_i\in [0,1],\sum_{i=1}^m\lambda_i=1, \sum_{i=1}^m\lambda_i x_i=x \right\}.
    \end{split}
    \]
\end{proposition}
See \cite[Theorem~5.6]{Roc70} for the more general result.

\begin{lemma}\label{lma:convdecnet}
    Let $\{F_i\}_{i\in I}$ be a decreasing net in $\Conv(N)$. Then $\inf_{i\in I} F_i\in \Conv(N)$.
\end{lemma}
\begin{proof}
Write $F=\inf_{i\in I} F_i$.
    We shall apply the characterization in \cref{lma:charconvex}. Take $n,r\in N$, $a,b\in \mathbb{R}$ such that $a>F(n)$, $b>F(r)$ and $t\in (0,1)$. We need to show that
    \begin{equation}\label{eq:convtemp1}
    F(tn+(1-t)r)<ta+(1-t)b.
    \end{equation}
    By definition, there are $j_n,j_r\in I$ such that $a>F_i(n)$ for $i\geq j_n$ and $b>F_i(r)$ for $i\geq j_r$. Since $I$ is directed, we can find $j\in I$ dominating both $j_n$ and $j_r$, and then for any $i\geq j$, we have
    \[
    a>F_i(n),\quad b>F_i(r).
    \]
    It follows from \cref{lma:charconvex} that
    \[
    F_i(tn+(1-t)r)<ta+(1-t)b
    \]
    for any $i\geq j$. Since $F_i$ is decreasing in $i$, we conclude \eqref{eq:convtemp1}.
\end{proof}

\begin{definition}\label{def:convexclosure}
    Let $F\in \Conv(N)$\index{$\Conv(N)$}. The \emph{closure}\index{convex function!closure of a convex function}\index{$\cl F$} $\cl F\in \Conv(N)$ of $F$ is defined as follows: If $F(n)=-\infty$ for some $n\in N$, then $\cl F\coloneqq -\infty$. Otherwise, we define $\cl F$ as the lower semicontinuity regularization of $F$.

    A convex function $F\in \Conv(N)$ is \emph{closed}\index{convex function!closed convex function} if $F=\cl F$. In other words, $F\in \Conv(N)$ is closed if one of the following conditions holds:
    \begin{enumerate}
        \item\label{itm:def:convexclosure:1} $F\equiv -\infty$;
        \item\label{itm:def:convexclosure:2} $F\equiv \infty$;
        \item\label{itm:def:convexclosure:3} $F$ is proper and lower semicontinuous.
    \end{enumerate}
\end{definition}
\begin{proposition}
    Let $F\in \Conv(N)$ be a closed convex function. Then $F$ is the supremum of all affine functions lying below $F$.
\end{proposition}
See \cite[Theorem~12.1]{Roc70}.

\begin{theorem}\label{thm:closurepropconv}
    Let $F\in \Conv^{\prop}(N)$\index{$\Conv^{\prop}(N)$}. Then $\cl F$ is a closed proper convex function. Moreover, $\cl F$ agrees with $F$ except possibly on the relative boundary of $\Dom F$.
\end{theorem}
See \cite[Theorem~7.4]{Roc70}.

\begin{proposition}\label{prop:closedcvx_domain}
    Let $F,F'\in \Conv(N)$ be closed convex functions. Assume that
    \begin{enumerate}
        \item\label{itm:prop:closedcvx_domain:1} $\RelInt \Dom F=\RelInt \Dom F'$, and
        \item\label{itm:prop:closedcvx_domain:2} $F=F'$ on $\RelInt \Dom F$.
    \end{enumerate}
    Then $F=F'$.
\end{proposition}
This is a special case of \cite[Corollary~7.3.4]{Roc70}.

\begin{definition}\label{def:partialorderconv}
    Given $F,F'\in \Conv(N)$, we write $F\preceq F'$ if there is $C\in \mathbb{R}$ such that
    \index{convex function!singularity type}
    \index{$F\preceq F'$}
    \[
        F\leq F'+C.
    \]
    We say $F\sim F'$ if $F\preceq F'$ and $F'\preceq F$ both hold.
\end{definition}

\begin{theorem}\label{thm:seqconv}
    Let $C\subseteq N$ be an open subset. Let $(f_i)_{i>0}$ be a sequence of real-valued convex functions on $C$. Suppose that the sequence converges on a dense subset of $C$ and the limit is finite. Then the limit
    \[
    f(x)\coloneqq \lim_{i\to\infty} f_i(x)
    \]
    exists for all $x\in C$ and is convex on $C$. Moreover, the sequence $(f_i)_i$ converges uniformly to $f$ on each compact subset of $C$.
\end{theorem}
This is a special case of \cite[Theorem~10.8]{Roc70}.

\section{Legendre transform}\label{sec:Leg}
Let $N$ be a real vector space of finite dimension and $M$ be the dual vector space.
The pairing $M\times N\rightarrow \mathbb{R}$ will be denoted by $\langle \bullet, \bullet \rangle$.

\begin{definition}\label{def:Legendregeneral}
    Let $F\in \Conv(N)$ be a convex function. We define the \emph{Legendre transform}\index{Legendre transform!Legendre transform of a convex function} of $F$ as the function $F^*\in \Conv(M)$:
    \begin{equation}\label{eq:Legtrandef}
        F^*(m)\coloneqq \sup_{n\in N} \left( \langle m,n \rangle-F(n) \right)=\sup_{n\in \RelInt \Dom F} \left( \langle m,n \rangle-F(n) \right).
    \end{equation}
\end{definition}
The latter equality follows from \cite[Corollary~12.2.2]{Roc70}.

Recall the well-known Legendre--Fenchel duality \cite[Theorem~12.2]{Roc70}.
\begin{theorem}\label{thm:Legendredual}
    Let $F\in \Conv(N)$. Then $F^*$ is a closed convex function. The function $F^*$ is proper if and only if $F$ is.

    Moreover, we have $(\cl F)^*=F^*$ and
    \[
        F^{**}=\cl F.
    \]
\end{theorem}

\begin{example}\label{ex:suppfundual}
    Let $P\subseteq M$ be a closed convex subset. Then
    \[
    \Supp_P^*=\chi_P,\quad \chi_P^*=\Supp_P.
    \]
    See \cite[Theorem~13.2]{Roc70}.
\end{example}

The following special case will be useful to us in the sequel.
\begin{corollary}\label{cor:Leginvpartdef}
    Let $F\colon (0,\infty)\rightarrow [-\infty,\infty)$ be a convex function. If we define $G\colon \mathbb{R}\rightarrow (-\infty,\infty]$ by
    \[
        G(\tau)=\sup_{t>0} \Bigl( t\tau-F(t) \Bigr),
    \]
    then $G$ is a convex function and
    \begin{equation}\label{eq:Leginvpartdef}
    F(t)=G^*(t),\quad \forall t>0.
    \end{equation}
    Moreover,
    \begin{equation}\label{eq:LegratQ}
        G(\tau)=\sup_{t\in \mathbb{R}_{>0}} \Bigl( t\tau-F(t) \Bigr)=\sup_{t\in \mathbb{Q}_{>0}} \Bigl( t\tau-F(t) \Bigr).
    \end{equation}
\end{corollary}
\begin{proof}
    We distinguish two cases.

    First suppose that $F(t)=-\infty$ for some $t>0$. Then $F(t)=-\infty$ for all $t>0$ by the convexity of $F$. Our assertions are clear in this case.

    Next assume that $F(t)\neq -\infty$ for all $t>0$. In this case, \cref{thm:closurepropconv} guarantees that $F$ admits a closed proper extension $\tilde{F}\in \Conv(\mathbb{R})$ with
    \[
    \tilde{F}(t)=\infty,\quad \forall t<0.
    \]
    It follows from \eqref{eq:Legtrandef} that
    \[
    G(\tau)=\tilde{F}^*(\tau),\quad \forall \tau\in \mathbb{R}.
    \]
    Now \cref{thm:Legendredual} implies \eqref{eq:Leginvpartdef}. Finally \eqref{eq:LegratQ} follows from the continuity of $F$.
\end{proof}

\begin{proposition}\label{prop:Legnonconvex}
    Let $F\colon N\rightarrow [-\infty,\infty]$. Define $F^*\colon M\rightarrow [-\infty,\infty]$ by
    \[
    F^*(m)\coloneqq \sup_{n\in N} \left( \langle m,n \rangle -F(n) \right).
    \]
    Then
    \[
    F^*=(\cl \LCE F)^*.
    \]
\end{proposition}
See \cite[Corollary~12.1.1]{Roc70}.

\begin{definition}\label{def:subgrad}
    Let $F\in \Conv(N)$ and $n\in N$. An element $m\in M$ is a \emph{subgradient}\index{subgradient} (in the sense of Alexandrov) of $F$ at $n$ if
    \begin{equation}\label{eq:subgrad}
    F(n')\geq F(n)+\langle n'-n,m \rangle,\quad \forall n'\in N.
    \end{equation}

    The set of subgradients of $F$ at $n$ is denoted by $\nabla F(n)$.

    More generally, for any subset $E\subseteq N$, we write
    \[
        \nabla F(E)=\bigcup_{n\in E}\nabla F(n).
    \]
\end{definition}

\begin{definition}\label{def:convexPorder}
    Given $F,F'\in \Conv(N)$, we write $F\preceq_P F'$ if
    \index{convex function!P-singularity type@$P$-singularity type}
    \index{$F\preceq_P F'$}
    \[
        \overline{\nabla F(N)}\subseteq \overline{\nabla F'(N)}.
    \]
    We write $F\sim_P F'$ if $F\preceq_P F'$ and $F'\preceq_P F$.
\end{definition}

\begin{theorem}\label{thm:nablaF}
    Suppose that $F\in \Conv^{\prop}(N)$. Then the following hold:
    \begin{enumerate}
        \item\label{itm:thm:nablaF:1} For any $n\not\in \Dom F$, $\nabla F(n)=\varnothing$;
        \item\label{itm:thm:nablaF:2} for any $n\in \RelInt\Dom F$, $\nabla F(n)\neq \varnothing$; moreover, for any $n'\in N$, we have
        \[
            \partial_{n'} F(n)=\sup\left\{\langle n', m \rangle: m\in \nabla F(n)\right\};
        \]
        \item\label{itm:thm:nablaF:3} for $n\in N$, the set $\nabla F(n)$ is non-empty and bounded if and only if $n\in \Int \Dom F$.
    \end{enumerate}
\end{theorem}
For the proof, we refer to \cite[Theorem~23.4]{Roc70}.

\begin{proposition}\label{prop:gradDomFstar}
    Let $F\in \Conv^{\prop}(N)$. Then
    \[
    \nabla F(N)\subseteq \Dom F^*.
    \]
    If moreover $F$ is closed, we have
    \begin{equation}\label{eq:relintdomFstar}
        \RelInt \Dom F^*\subseteq \nabla F(N).
    \end{equation}
    In particular, if $F$ is a proper closed convex function on $N$, then
    \[
    \overline{\nabla F(N)}=\overline{\Dom F^*}.
    \]
\end{proposition}
\begin{proof}
    Suppose that $m\in \nabla F(n)$ for some $n\in N$. Then \eqref{eq:subgrad} holds. In particular,
    \[
    \langle m,n'\rangle -F(n')\leq \langle m,n\rangle-F(n).
    \]
    It follows that
    \[
    F^*(m)\leq  \langle m,n\rangle-F(n)<\infty.
    \]
    \eqref{eq:relintdomFstar} is proved in \cite[Corollary~23.5.1]{Roc70}. For the last assertion, it suffices to observe that $\overline{\RelInt \Dom F^*}=\overline{\Dom F^*}$.
\end{proof}

\begin{proposition}\label{prop:Legendretranssup}
    Let $\{F_i\}_{i\in I}$ be a non-empty family in $\Conv^{\prop}(N)$. Then
    \[
    \left(\bigwedge_{i\in I} F_i\right)^*=\bigvee_{i\in I} F_i^*,\quad \left(\bigvee_{i\in I} \cl F_i\right)^*=\cl \bigwedge_{i\in I} F_i^*.
    \]
    If $I$ is finite and $\overline{\Dom F_i}$ is independent of the choice of $i\in I$, then
    \[
    \left(\bigvee_{i\in I} F_i\right)^*=\bigwedge_{i\in I} F_i^*.
    \]
\end{proposition}
Recall that $\wedge$ is defined in \cref{def:convwedge} and $\lor$ in \cref{def:convlor}. See \cite[Theorem~16.5]{Roc70} for the proof.

\begin{proposition}\label{prop:sumLegendre}
    Let $F_1,\ldots,F_r\in \Conv^{\prop}(N)$ ($r\in \mathbb{Z}_{>0}$). Assume that
    \[
        \bigcap_{i=1}^r\RelInt \Dom(F_i)\neq \varnothing,
    \]
    then for any $m\in M$,
    \[
    \left( \sum_{i=1}^r F_i \right)^*(m)=\inf\left\{\sum_{i=1}^r F_i^*(m_i):m_1,\ldots,m_r\in M, \sum_{i=1}^r m_i=m  \right\}.
    \]
\end{proposition}
See \cite[Theorem~16.4]{Roc70} for the proof.

\begin{proposition}\label{prop:Fsuppchar}
    Let $P\subseteq M$ be a convex body\footnote{Here a convex body refers to a non-empty compact convex subset, not necessarily having non-empty interior.} and $F\in \Conv^{\prop}(N)$. The following are equivalent:
    \begin{enumerate}
        \item\label{itm:prop:Fsuppchar:1} $F\preceq \Supp_P$;
        \item\label{itm:prop:Fsuppchar:2} $\Dom F=N$ and $F^*|_{M\setminus P}\equiv \infty$;
        \item\label{itm:prop:Fsuppchar:3} $\Dom F=N$ and $\nabla F(N)\subseteq P$.
    \end{enumerate}
    Moreover, under these conditions,
    \begin{equation}\label{eq:Fsupequal}
    F(n)-\Supp_P(n)\leq F(0),\quad \forall n\in N.
    \end{equation}
\end{proposition}
\begin{proof}
    \ref{itm:prop:Fsuppchar:1} $\implies$ \ref{itm:prop:Fsuppchar:2}. It is clear that $\Dom F=N$ since $\Dom \Supp_P=N$.
    From $F\preceq \Supp_P$ and \cref{ex:suppfundual}, we know that
    \[
    \chi_P=\Supp_P^*\preceq F^*.
    \]
    So \ref{itm:prop:Fsuppchar:2} follows.

    \ref{itm:prop:Fsuppchar:2} $\implies$ \ref{itm:prop:Fsuppchar:3}. This follows from \cref{prop:gradDomFstar}.

    \ref{itm:prop:Fsuppchar:3} $\implies$ \ref{itm:prop:Fsuppchar:1}. The function $t\mapsto F(tn)$ is convex on $[0,1]$, hence absolutely continuous and differentiable almost everywhere. For almost every $t\in(0,1)$,
    \[
        \frac{\mathrm d}{\mathrm dt}F(tn)=\partial_n F(tn) \leq \Supp_P(n)
    \]
    by \cref{thm:nablaF} and the assumption $\nabla F(N)\subseteq P$. Integrating in $t$ gives
    \[
        F(n)-F(0)\leq \Supp_P(n).
    \]
    In particular, \eqref{eq:Fsupequal} also follows.
\end{proof}

\section{Classes of convex functions}
\label{sec:convex_classes}
Let $N$ be a real vector space of finite dimension and $M$ be the dual vector space.

We shall fix a convex body $P\subseteq M$.

The following classes are introduced in \cite{BB13}.
\begin{definition}\label{def:convexPfunctions}
    We define the set $\mathcal{P}(N,P)$ as the set of proper convex functions $F\in \Conv(N)$ such that $F\preceq \Supp_P$.
    \index{$\mathcal{P}(N,P)$}

    We define the set $\mathcal{E}^{\infty}(N,P)$ as the set of closed convex functions $F\in \Conv(N)$ such that $F\sim \Supp_P$.
    \index{$\mathcal{E}^{\infty}(N,P)$}

    We define the set $\mathcal{E}(N,P)$ as follows: Suppose that $\Int P=\varnothing$. Then $\mathcal{E}(N,P)\coloneqq \mathcal{P}(N,P)$; otherwise, let
    \index{$\mathcal{E}(N,P)$}
    \[
    \mathcal{E}(N,P)=\left\{F\in \mathcal{P}(N,P) : P=\overline{\nabla F(N)} \right\}.
    \]
    We define the set $\mathcal{E}^1(N,P)$ as the subset of $\mathcal{E}(N,P)$ consisting of $F\in \mathcal{E}(N,P)$ with
    \index{$\mathcal{E}^{1}(N,P)$}
    \[
    \int_P F^*\,\mathrm{d}\vol<\infty,
    \]
    where $\mathrm{d}\vol$ is any Lebesgue measure on $M$.
\end{definition}
Observe that for any $F\in \mathcal{P}(N,P)$, we have $\Dom F=N$ and $F$ is necessarily closed.

\begin{proposition}\label{prop:convex:L415}
We have
\[
\mathcal{E}^{\infty}(N,P)\subseteq \mathcal{E}^{1}(N,P)\subseteq \mathcal{E}(N,P) \subseteq \mathcal{P}(N,P).
\]
\end{proposition}
\begin{proof}
    When $\Int P=\varnothing$, the assertion is clear. We assume that $\Int P\neq \varnothing$. The second inclusion follows from the definition. We only handle the first inequality. Take $F\in \mathcal{E}^{\infty}(N,P)$. By definition, $F\sim \Supp_P$ and hence $F^*\sim \chi_P$. It follows that $P=\Dom F^*$.

    By \cref{prop:Fsuppchar}, we already know that
    \[
    \nabla F(N)\subseteq P=\Dom F^*.
    \]
    On the other hand, by \cref{prop:gradDomFstar}, we have
    \[
    \Int P\subseteq \nabla F(N).
    \]
    So it follows that
    \[
    P=\overline{\nabla F(N)}.
    \]
    It is clear that $F^*\sim \chi_P$ is integrable.
\end{proof}

\begin{proposition}\label{prop:convex:L440}
    For any $F\in \mathcal{E}^{\infty}(N,P)$, we have $F^*|_{M\setminus P}\equiv \infty$ and $F^*$ is bounded on $P$.
\end{proposition}
\begin{proof}
    From $F\sim \Supp_P$, we take the Legendre transform to get $F^*\sim \Supp_P^*=\chi_P$, where we applied \cref{ex:suppfundual}.
\end{proof}

\begin{definition}\label{def:pointwiseconv}
    We endow the topology of pointwise convergence on $\mathcal{P}(N,P)$.
    \index{pointwise convergence}
\end{definition}
Note that this topology coincides with the compact-open topology. See \cite[Theorem~10.8]{Roc70}.

\begin{proposition}\label{prop:regularizationconvex}
    Let $F\in \mathcal{P}(N,P)$. Then there is a decreasing sequence $F_j\in \mathcal{E}^{\infty}(N,P)\cap C^{\infty}(N)$ converging to $F$.
\end{proposition}
See \cite[Lemma~2.2]{BB13}.

We observe that the point $0\in N$ plays a special role since it does in the definition of the support function.
\begin{proposition}\label{prop:convex:L462}
    For any $F\in \mathcal{P}(N,P)$, we have
    \[
    \max_{N}(F-\Supp_P)=F(0).
    \]
\end{proposition}
\begin{proof}
    It follows from \eqref{eq:Fsupequal} that
    \[
    \sup_{N}(F-\Supp_P)\leq F(0).
    \]
    The equality is clearly obtained at $0\in N$.
\end{proof}

\begin{lemma}\label{lma:grad_img_sum}
    Let $P'\subseteq M$ be another convex body. Then for any $F\in \mathcal{P}(N,P)$ and $F'\in \mathcal{P}(N,P')$, we have
    \[
    F+F'\in \mathcal{P}(N,P+P').
    \]

    Moreover, let $F$ and $F'$ be real-valued convex functions on $N$. Assume that $\Dom F^*$ is bounded. Then
    \[
    \overline{\nabla F(N)}+\overline{\nabla F'(N)}=\overline{\nabla (F+F')(N)}.
    \]
    In particular, if $F\in \mathcal{P}(N,P)$, $F'\in \mathcal{P}(N,P')$ and
    \[
    \overline{\nabla F(N)}=P,\quad \overline{\nabla F'(N)}=P',
    \]
    then $F+F'\in \mathcal{E}(N,P+P')$.
\end{lemma}
\begin{proof}
    The former assertion follows immediately from the observation
    \[
    \Supp_{P+P'}=\Supp_P+\Supp_{P'}.
    \]

    We prove the gradient-image equality. In view of \cref{prop:gradDomFstar}, this means
    \begin{equation}\label{eq:overlineDom_temp1}
    \overline{\Dom (F+F')^*}=\overline{\Dom F^*}+\overline{\Dom F'{}^*}.
    \end{equation}
    It follows from \cref{prop:sumLegendre} that
    \[
    \Dom (F+F')^*=\Dom F^*+\Dom F'{}^*.
    \]
    Since $\overline{\Dom F^*}$ is compact, \eqref{eq:overlineDom_temp1} follows\footnote{In general, the Minkowski sum does not commute with the closure.}.

    Finally, under the assumptions in the last assertion, the former assertion gives $F+F'\in\mathcal{P}(N,P+P')$. If $\Int(P+P')=\varnothing$, this is exactly the definition of $\mathcal{E}(N,P+P')$. Otherwise, the gradient-image equality gives
    \[
    \overline{\nabla (F+F')(N)}=P+P',
    \]
    so again $F+F'\in\mathcal{E}(N,P+P')$.
\end{proof}

\section{Monge--Amp\`ere measures}
\label{sec:convex_MA}
Let $N$ be a free Abelian group of finite rank $n$ (i.e., a lattice) and $M$ be its dual lattice. There is a canonical Lebesgue-type measure on $M_{\mathbb{R}}$, denoted by $\mathrm{d}\vol$, normalized so that the smallest cubes in $M$ have volume $1$. Similarly, the canonical measure on $N_{\mathbb{R}}$ is normalized in the same way and is denoted by $\mathrm{d}\vol$ as well.

We will write
\[
N_{\mathbb{R}}=N\otimes_{\mathbb{Z}}\mathbb{R},\quad M_{\mathbb{R}}=M\otimes_{\mathbb{Z}}\mathbb{R}.
\]

\begin{definition}\label{def:realMA}
    Let $F$ be a finite convex function on $N_{\mathbb{R}}$. We define the \emph{real Monge--Ampère measure}\index{real Monge--Ampère measure} $\MA_{\mathbb{R}}F$ as the Alexandrov measure on $N_{\mathbb{R}}$ given as follows: for each Borel measurable set $E\subseteq N_{\mathbb{R}}$, define
    \begin{equation}\label{eq:MARFE}
        \MA_{\mathbb{R}}F(E)\coloneqq n!\int_{\nabla F(E)} \mathrm{d}\vol.
    \end{equation}
    The measurability and countable additivity in this definition are standard properties of the subgradient map; see for example \cite[Section~2.1]{Fig17}.
\end{definition}

\begin{proposition}\label{prop:smoothMAreal}
    Suppose that $F\in C^{1,1}(N_{\mathbb{R}})\cap \Conv(N_{\mathbb{R}})$. Fix an identification $N=\mathbb{Z}^n$. Then
    \[
    \MA_{\mathbb{R}}F=n!\cdot \det \nabla^2 F \,\mathrm{d}\vol.
    \]
\end{proposition}
See \cite[Example~2.2]{Fig17}.

\begin{proposition}\label{prop:convex:L541}
    Let $P\subseteq M_{\mathbb{R}}$ be a convex body and $F\in \mathcal{P}(N_{\mathbb{R}},P)$. Then $F\in \mathcal{E}(N_{\mathbb{R}},P)$ if and only if
    \begin{equation}\label{eq:cvxfullmass}
    \int_{N_{\mathbb{R}}}\MA_{\mathbb{R}}F=n! \vol P.
    \end{equation}
\end{proposition}
\begin{proof}

   Since $F\in\mathcal{P}(N_{\mathbb{R}},P)$, the function $F$ is closed and has domain $N_{\mathbb{R}}$. By \cref{prop:gradDomFstar}, $\overline{\nabla F(N_{\mathbb{R}})}=\overline{\Dom F^*}$, so this closure is convex. As $\nabla F(N_{\mathbb{R}})$ contains the relative interior of $\Dom F^*$, the definition of $\MA_{\mathbb{R}}$ shows that \eqref{eq:cvxfullmass} is equivalent to
   \[
   \vol \overline{\nabla F(N_{\mathbb{R}})}=\vol P.
   \]

   We first handle the case where $\Int P\neq \varnothing$.
   By \cref{prop:Fsuppchar}, $\overline{\nabla F(N_{\mathbb{R}})}\subseteq P$. Hence the latter is equivalent to
   \[
   \overline{\nabla F(N_{\mathbb{R}})}=P.
   \]
   Indeed, a proper closed subset of the full-dimensional convex body $P$ has complement containing a non-empty relatively open subset of $P$, hence has strictly smaller volume.

   Now assume that $\Int P=\varnothing$, then $\vol \overline{\nabla F(N_{\mathbb{R}})}=\vol P=0$ by \cref{prop:Fsuppchar}. The assertion is clear.
\end{proof}

\begin{theorem}\label{thm:realMAcont}
    Let $F,F_j\in \mathcal{P}(N_{\mathbb{R}},P)$ ($j\in \mathbb{Z}_{>0}$). Assume that $F_j\to F$. Then
    \[
    \MA_{\mathbb{R}}(F_j)\rightharpoonup \MA_{\mathbb{R}}(F)
    \]
    weakly.
\end{theorem}
See \cite[Proposition~2.6]{Fig17}.

There is a well-known comparison principle.
\begin{theorem}\label{thm:convcomp}
    Let $F,F'\in \mathcal{P}(N_{\mathbb{R}},P)$. Assume that $F\preceq F'$. Then
    \[
    \overline{\nabla F(N_{\mathbb{R}})}\subseteq \overline{\nabla F'(N_{\mathbb{R}})},\quad
    \int_{N_{\mathbb{R}}} \MA_{\mathbb{R}}(F)\leq \int_{N_{\mathbb{R}}} \MA_{\mathbb{R}}(F').
    \]
\end{theorem}
\begin{proof}
    It suffices to observe that $F'{}^*\preceq F^*$, and hence the first assertion follows from \cref{prop:gradDomFstar}. By the definition of the real Monge--Ampère measure,
    \[
        \int_{N_{\mathbb R}}\MA_{\mathbb R}F=n!\vol\overline{\nabla F(N_{\mathbb R})}.
    \]
    So the second assertion follows from the first.
\end{proof}

\begin{lemma}\label{lma:realMA_Legendre_pushforward}
    Let $F\in \mathcal{E}(N_{\mathbb{R}},P)$ and set $G=F^*$. Then the map $\nabla G$ pushes forward $n!\,\mathrm{d}\vol|_{P}$ to $\MA_{\mathbb R}(F)$. Equivalently, for every bounded Borel function $h$ on $N_{\mathbb R}$,
    \[
        \int_{N_{\mathbb R}} h\,\MA_{\mathbb R}(F)= n!\int_P h(\nabla G(m))\,\mathrm{d}\vol(m).
    \]
\end{lemma}
\begin{proof}
    If $\vol P=0$, the assertion is clear from \cref{prop:Fsuppchar} and the definition of $\MA_{\mathbb R}$. We may therefore assume that $P$ has non-empty interior.

    Since $F\in \mathcal P(N_{\mathbb R},P)$, we have $\Dom G\subseteq P$ by \cref{prop:Fsuppchar}. On the other hand, by the definition of $\mathcal E(N_{\mathbb R},P)$ and \cref{prop:gradDomFstar},
    \[
        P=\overline{\nabla F(N_{\mathbb R})}\subseteq \overline{\Dom G}.
    \]
    Hence $\overline{\Dom G}=P$. Since $\Dom G$ is convex, $G$ is finite on $\Int P$ and is equal to $\infty$ outside $P$. Let $D$ be the set of points in $\Int P$ at which $G$ is differentiable. By \cite[Theorem~25.5]{Roc70}, $P\setminus D$ has measure zero. By \cite[Theorem~23.5]{Roc70}, $m\in\nabla F(n)$ if and only if $n\in\nabla G(m)$. Hence for any Borel set $E\subseteq N_{\mathbb R}$,
    \[
        \{m\in D:\nabla G(m)\in E\}=D\cap \nabla F(E).
    \]
    It follows from \cref{def:realMA} that
    \[
        n!\vol \{m\in P:\nabla G(m)\in E\}=n!\vol \nabla F(E)=\MA_{\mathbb R}(F)(E).
    \]
    Our assertions follow.
\end{proof}

Finally, let us recall a mixed version of the real Monge--Ampère measure:
\begin{theorem}\label{thm:mixed_realMA}
    There is a symmetric map assigning to every $n$-tuple of finite convex functions $F_1,\ldots,F_n$ on $N_{\mathbb R}$ a positive Radon measure
    \[
        \MA_{\mathbb R}(F_1,\ldots,F_n)
    \]
    on $N_{\mathbb R}$ with the following properties:
    \begin{enumerate}
        \item\label{itm:thm:mixed_realMA:1} For every $q\in \mathbb Z_{>0}$, finite convex functions $F_1,\ldots,F_q$ on $N_{\mathbb R}$, and $\lambda_1,\ldots,\lambda_q\geq 0$,
        \[
            \MA_{\mathbb R}(\lambda_1F_1+\cdots+\lambda_qF_q)
            =\sum_{i_1,\ldots,i_n=1}^q
            \lambda_{i_1}\cdots\lambda_{i_n}
            \MA_{\mathbb R}(F_{i_1},\ldots,F_{i_n}).
        \]
        \item\label{itm:thm:mixed_realMA:2} For every finite convex function $F$ on $N_{\mathbb R}$,
        \[
            \MA_{\mathbb R}(F,\ldots,F)=\MA_{\mathbb R}(F).
        \]
        \item\label{itm:thm:mixed_realMA:3} The map $\MA_{\mathbb R}$ is additive and positively homogeneous of degree $1$ in each entry. In other words, for finite convex functions $F,G,F_2,\ldots,F_n$ on $N_{\mathbb R}$ and $\lambda,\mu\geq 0$,
        \[
            \MA_{\mathbb R}(\lambda F+\mu G,F_2,\ldots,F_n)
            =\lambda\MA_{\mathbb R}(F,F_2,\ldots,F_n)
            +\mu\MA_{\mathbb R}(G,F_2,\ldots,F_n).
        \]
    \end{enumerate}
\end{theorem}
This is \cite[Theorem~4.3(a),(b),(f)]{CLM22}, multiplied by $n!$ to match the normalization in \cref{def:realMA}. The positivity is part of the construction: it is the non-negativity of the mixed discriminant in the smooth case and follows in general by weak continuity. The mixed real Monge--Ampère measure was also studied in \cite{TW02,PR04}.

\section{Several lemmata}
\label{sec:convex_sep}

\begin{lemma}\label{lma:polybdd}
	Let $\alpha,\beta_1,\ldots,\beta_m\in \mathbb{Z}^n$. Let $\Delta$ be the polytope generated by $\beta_1,\ldots,\beta_m$. Then the following are equivalent:
	\begin{enumerate}
		\item\label{itm:lma:polybdd:1}
		      \begin{equation}\label{eq:zalpha}
			      |z^{\alpha}|^2\left(\sum_{i=1}^m |z^{\beta_i}|^2 \right)^{-1}
		      \end{equation}
		      is a bounded function on $\mathbb{C}^{*n}$.
		\item\label{itm:lma:polybdd:2} $\alpha\in \Delta$.
	\end{enumerate}
\end{lemma}
\begin{proof}
	\ref{itm:lma:polybdd:2} $\implies$ \ref{itm:lma:polybdd:1}. Write $\alpha=\sum_i t_i\beta_i$, where $t_i\in [0,1]$, $\sum_i t_i=1$. Then
	\[
 \begin{aligned}
			|z^{\alpha}|^2\left(\sum_{i=1}^m |z^{\beta_i}|^2 \right)^{-1}
            &=\prod_i |z^{\beta_i}|^{2t_i}\left(\sum_{i=1}^m |z^{\beta_i}|^2 \right)^{-1}\\
            &\leq \prod_i\left(\sum_{j=1}^m |z^{\beta_j}|^2\right)^{t_i}
            \left(\sum_{i=1}^m |z^{\beta_i}|^2 \right)^{-1}=1.
	  \end{aligned}
		\]

		\ref{itm:lma:polybdd:1} $\implies$ \ref{itm:lma:polybdd:2}. Assume that $\alpha\not\in \Delta$. By the strict separation theorem, there is $a=(a_1,\ldots,a_n)\in\mathbb{R}^n$ such that
        \[
        a\cdot\alpha>\max_i a\cdot\beta_i.
        \]
        Set
		\[
			z(t)\coloneqq (t^{a_1},\ldots,t^{a_n}).
		\]
		Then \eqref{eq:zalpha} evaluated at $z(t)$ grows like
        \[
        \frac{t^{2a\cdot \alpha}}{\sum_i t^{2a\cdot\beta_i}}
        \]
        as $t\to\infty$, and hence is not bounded.
\end{proof}
\begin{lemma}\label{lma:polybdd2}
	Let $m\geq 1$, $\beta_1,\ldots,\beta_m\in \mathbb{N}^n$, and $\beta\in \mathbb{R}^n$. Then the following are equivalent:
	\begin{enumerate}
		\item\label{itm:lma:polybdd2:1} $\log \sum_{i=1}^m \mathrm{e}^{x\cdot\beta_i}-x\cdot\beta$ is bounded from below.
		\item\label{itm:lma:polybdd2:2} $\beta$ is in the convex hull of the $\beta_i$'s.
	\end{enumerate}
\end{lemma}
\begin{proof}
		If $\beta=\sum_i t_i\beta_i$ with $t_i\in[0,1]$ and $\sum_i t_i=1$, then by the weighted arithmetic--geometric mean inequality,
        \[
        x\cdot\beta=\sum_i t_i\,x\cdot\beta_i
        \leq \log\sum_i t_i\mathrm{e}^{x\cdot\beta_i}
        \leq \log\sum_i \mathrm{e}^{x\cdot\beta_i}.
        \]
        This proves that the function in \ref{itm:lma:polybdd2:1} is bounded from below.

        Conversely, if $\beta$ does not lie in the convex hull of the $\beta_i$'s, choose $a\in\mathbb{R}^n$ such that
        \[
        a\cdot\beta>\max_i a\cdot\beta_i.
        \]
        Evaluating at $x=ta$ and letting $t\to\infty$, we find
        \[
        \log \sum_i \mathrm{e}^{t a\cdot\beta_i}-t\,a\cdot\beta\to-\infty,
        \]
        contradicting \ref{itm:lma:polybdd2:1}.
\end{proof}

\begin{lemma}\label{lma:conc_lim_der}
    Let $I\subseteq \mathbb{R}$ be an open interval containing $0$. Let $(f_k\colon I\rightarrow \mathbb{R})_{k\geq 1}$ be a sequence of convex functions, and let $g\colon I\to \mathbb{R}$ be a function. Assume that
    \[
        \varlimsup_{k\to\infty} f_k(t)\leq g(t)\quad \text{for every } t\in I,
    \]
    and that
    \[
        \lim_{k\to\infty} f_k(0)=g(0).
    \]
    If each $f_k$ and $g$ is differentiable at $0$, then
    \[
        \lim_{k\to\infty} f_k'(0)=g'(0).
    \]
\end{lemma}
\begin{proof}
    Since $f_k$ is convex and differentiable at $0$, for every $t\in I$,
    \[
        f_k(t)\geq f_k(0)+t f_k'(0).
    \]
    Taking the upper limit as $k\to\infty$ and using the assumptions gives
    \[
        \varlimsup_{k\to\infty} t f_k'(0)
        \leq g(t)-g(0).
    \]

    For $t>0$, this yields
    \[
        \varlimsup_{k\to\infty} f_k'(0)
        \leq \frac{g(t)-g(0)}{t}.
    \]
    Letting $t\to 0+$ and using the differentiability of $g$ at $0$, we get
    \[
        \varlimsup_{k\to\infty} f_k'(0)\leq g'(0).
    \]

    For $t<0$, the same inequality gives
    \[
        \varliminf_{k\to\infty} f_k'(0)
        \geq \frac{g(t)-g(0)}{t}.
    \]
    Letting $t\to 0-$ and again using the differentiability of $g$ at $0$, we obtain
    \[
        \varliminf_{k\to\infty} f_k'(0)\geq g'(0).
    \]
    Therefore
    \[
        \lim_{k\to\infty} f_k'(0)=g'(0),
    \]
    as claimed.
\end{proof}

\chapter{Pluripotential theory on unibranch spaces}
\label[appendix]{chap:unibranch}

This appendix explains how much of the pluripotential theory developed in the body of the book extends from compact K\"ahler manifolds to compact unibranch K\"ahler spaces. The unibranch hypothesis is the natural generality in which plurisubharmonic functions descend through resolutions and the basic envelope formalism remains meaningful. Some constructions, especially the trace operator of \cref{chap:trace}, require extra care in the singular setting; we record here the precise statements needed in the main text.

\section{Complex spaces}
\label{sec:unib_complex}

A complex space is assumed to be reduced, Hausdorff and paracompact in the whole book.

\begin{definition}\label{def:primdiv}
A \emph{prime divisor}\index{prime divisor over a complex space} over an irreducible complex space $Z$ is a connected smooth hypersurface $E\subseteq X'$, where $X'\rightarrow Z$ is a proper bimeromorphic morphism with $X'$ smooth. Such a morphism $X'\rightarrow Z$ is also called a \emph{resolution}\index{resolution} of $Z$. The \emph{center}\index{center of a prime divisor} of the prime divisor is defined as the image of $E$ in $Z$.

Two prime divisors $E_1\subseteq X'_1$ and $E_2\subseteq X'_2$ over $Z$ are \emph{equivalent} if there is a common resolution $X''\rightarrow Z$ dominating both $X_1'$ and $X_2'$ such that the strict transforms of $E_1$ and $E_2$ coincide.

The set $Z^{\Div}$ is the set of pairs $(c,E)$, where $c\in \mathbb{Q}_{>0}$ and $E$ is an equivalence class of a prime divisor over $Z$. For simplicity, we will denote the pair $(c,E)$ by $c\ord_E$, although one should not really think of this object as a valuation unless $Z$ is projective and irreducible.
\end{definition}
Note that a prime divisor on $Z$ does not always define a prime divisor over $Z$ if $Z$ is singular.

\begin{definition}\label{def:unibranch:L21}
    A complex space $X$ is \emph{unibranch}\index{unibranch complex space} if for all $x\in X$, the local ring $\mathcal{O}_{X,x}$ is unibranch.
\end{definition}
See \cite[\href{https://stacks.math.columbia.edu/tag/0BPZ}{Tag 0BPZ}]{stacks-project} for the notion of unibranch local rings. In \cite{CAS}, the same notion is known as \emph{local irreducibility}. We prefer the name \emph{unibranchness}, because the corresponding terminology in algebraic geometry is standard and goes back to Grothendieck \cite[(0.23.2.1)]{EGAIV1}.

The notion of unibranchness coincides with the usual one when $\mathcal{X}$ is an algebraic variety:
\begin{proposition}\label{prop:unibranch_gaga}
    Let $\mathcal{X}$ be a reduced separated scheme of finite type over $\mathbb{C}$, and let $\mathcal{X}^{\mathrm{An}}$ be the associated complex space.
    Then $\mathcal{X}^{\An}$ is unibranch in the sense above if and only if $\mathcal{X}$ is unibranch in the algebraic sense (see \cite[\href{https://stacks.math.columbia.edu/tag/0BQ2}{Tag 0BQ2}]{stacks-project}).
\end{proposition}
This is \cite[Remark~2.7]{Xia23Mabuchi} in the arXiv version.
\begin{proof}
    Let $x\in \mathcal{X}(\mathbb{C})$. Then we claim that $\mathcal{X}$ is unibranch at $x$ if and only if $\mathcal{X}^{\An}$ is unibranch at $x$. By GAGA, there is a natural morphism of ringed spaces $(\mathcal{X}^{\An},\mathcal{O}_{\mathcal{X}^{\An}})\rightarrow (\mathcal{X},\mathcal{O}_{\mathcal{X}})$.
	 Then since $\mathcal{O}_{\mathcal{X},x}$ is excellent, $\mathcal{X}$ is unibranch at $x$ if and only if $\widehat{\mathcal{O}_{\mathcal{X},x}}=\widehat{\mathcal{O}_{\mathcal{X}^{\An},x}}$ is integral (\cite[Théorème~18.9.1]{EGAIV4}). The latter implies that $\mathcal{O}_{\mathcal{X}^{\An},x}$ is unibranch by \cite[Chapitre~5, Exercise~2.8]{AC}. Running the same argument the other way round, we conclude that $\mathcal{O}_{\mathcal{X},x}$ is unibranch if and only if $\mathcal{O}_{\mathcal{X}^{\An},x}$ is.

    It remains to check that unibranchness at all closed points suffices on the algebraic side. Let $\nu\colon \widetilde{\mathcal{X}}\to\mathcal{X}$ be the finite normalization. The number of branches at a point is the number of points in its normalization fiber (see \cite[\href{https://stacks.math.columbia.edu/tag/0C3Z}{Section 0C3Z}]{stacks-project}). The locus of fibers with at least two geometric points is constructible by Chevalley's theorem, since it is the image of the complement of the diagonal in $\widetilde{\mathcal{X}}\times_{\mathcal{X}}\widetilde{\mathcal{X}}$. At closed points geometric points and ordinary points give the same count. If nonempty it contains a closed point, which proves the assertion.
\end{proof}

\begin{theorem}[Zariski's main theorem]\label{thm:Zariskimain}\index{theorem!Zariski's main theorem}
    Let $\pi\colon Y\rightarrow X$ be a proper bimeromorphic morphism between complex spaces. Assume that $X$ is unibranch. Then $\pi$ has connected fibers.
\end{theorem}
We follow \cite[Proof of Théorème~1.7]{Dem85}.
\begin{proof}
    Let $Z\subseteq X$ be the smallest nowhere dense analytic subset such that $\pi$ is an isomorphism over $X\setminus Z$.
    Assume that the fiber $F=\pi^{-1}(x)$ over some $x\in X$ is disconnected.
    Since $\pi$ is proper, $F$ is compact. Choose a decomposition $F=F_1\sqcup F_2$ into two non-empty disjoint compact open subsets of $F$, and choose disjoint open neighborhoods $W_1,W_2$ of $F_1,F_2$ in $Y$. Since proper maps are closed, after shrinking an open neighborhood $U$ of $x$, we may assume that
    \[
        \pi^{-1}(U)\subseteq W_1\cup W_2.
    \]
    Thus $\pi^{-1}(U)$ is disconnected.

    Since $X$ is unibranch, we may shrink $U$ further so that $U$ is irreducible. Then $U\setminus Z$ is connected (cf. \cite[Chapter~9, \S1.2]{CAS}). Since $\pi$ is bimeromorphic, every irreducible component of $Y$ maps onto an irreducible component of $X$. Hence $\pi^{-1}(Z)$ contains no irreducible component of $Y$ and is nowhere dense. On the other hand, write
    \[
        \pi^{-1}(U)=V_1\sqcup V_2
    \]
    as a union of two non-empty disjoint open subsets. Each $V_i$ meets
    $\pi^{-1}(U\setminus Z)$: otherwise $V_i$ would be a non-empty open subset of $Y$ contained in the nowhere dense analytic subset $\pi^{-1}(Z)$, which is impossible. Hence
    \[
        \pi^{-1}(U\setminus Z)
        =
        \bigl(V_1\cap \pi^{-1}(U\setminus Z)\bigr)
        \sqcup
        \bigl(V_2\cap \pi^{-1}(U\setminus Z)\bigr)
    \]
    is disconnected. This contradicts the fact that
    $\pi\colon \pi^{-1}(U\setminus Z)\rightarrow U\setminus Z$ is an isomorphism and $U\setminus Z$ is connected.
\end{proof}

\begin{definition}\label{def:modif}
    A \emph{modification}\index{modification!modification of a complex space} of a compact complex space $X$ is a finite composition of blow-ups with smooth centers.

    We say a modification $\pi'\colon Z\rightarrow X$ \emph{dominates}\index{modification!domination of modifications} another $\pi\colon Y\rightarrow X$ if there is a morphism $g\colon Z\rightarrow Y$ of complex spaces making the following diagram commutative:
    \begin{equation}\label{eq:domi}
    \begin{tikzcd}
Z \arrow[rr,"g"] \arrow[rd, "\pi'"'] &   & Y \arrow[ld, "\pi"] \\
                                 & X. &
\end{tikzcd}
    \end{equation}
\end{definition}

The modifications of $X$ together with the domination relation form a directed set.

\begin{theorem}[Hironaka's Chow lemma]\label{thm:HironakaChow}\index{theorem!Hironaka's Chow lemma}
    Suppose that $X$ is a compact complex space. Then every proper bimeromorphic morphism to $X$ can be dominated by a modification.
\end{theorem}
This follows from the proof of \cite[Corollary~2]{Hir75}.

\begin{theorem}\label{thm:res}
    Let $X$ be a compact complex space. Then there is a modification $\pi\colon Y\rightarrow X$ such that $Y$ is smooth.
\end{theorem}
See \cite{BM97, Wlo09}.

\begin{corollary}\label{cor:primerealization}
    Let $X$ be an irreducible compact complex space and $E$ be a prime divisor over $X$. Then there is a modification $\pi\colon Y\rightarrow X$ such that $Y$ is smooth and $E$ can be realized as a prime divisor on $Y$.
\end{corollary}
\begin{proof}
    Choose a smooth model $X'\rightarrow X$ on which $E$ is realized. By \cref{thm:HironakaChow}, this proper bimeromorphic morphism is dominated by a modification $Y\rightarrow X$. After further resolving the pair if necessary, we may assume that $Y$ is smooth and that the strict transform of $E$ on $Y$ is smooth. This strict transform realizes the same prime divisor over $X$.
\end{proof}

\section{Plurisubharmonic functions}
\label{sec:unib_psh}
Let $X$ be a complex space.

\begin{definition}\label{def:unibranch:L105}
    A function $\varphi\colon X\rightarrow [-\infty,\infty)$ is \emph{plurisubharmonic}\index{plurisubharmonic function!plurisubharmonic function on a complex space} if
	\begin{enumerate}
		\item\label{itm:appendix-unibranch:L107:1} $\varphi$ is not identically $-\infty$ on any irreducible component of $X$, and
		\item\label{itm:appendix-unibranch:L107:2} for any $x\in X$, there is an open neighborhood $V$ of $x$ in $X$, a domain $\Omega\subseteq \mathbb{C}^N$, a closed immersion $V\hookrightarrow \Omega$ and a plurisubharmonic function $\tilde{\varphi}\in \PSH(\Omega)$ such that $\varphi|_V=\tilde{\varphi}|_V$ under this immersion.
	\end{enumerate}
	The set of plurisubharmonic functions on $X$ is denoted by $\PSH(X)$.

    Similarly, if $\theta$ is a smooth closed real $(1,1)$-form on $X$\footnote{Here a \emph{smooth closed real $(1,1)$-form} on $X$ is one which, near every point, is the restriction of a smooth closed real $(1,1)$-form on an open subset of $\mathbb{C}^N$ under a local embedding of $X$; equivalently, $\theta$ has local smooth $\ddc$-potentials.}, then a function $\varphi\colon X\rightarrow [-\infty,\infty)$ is \emph{$\theta$-plurisubharmonic}\index{plurisubharmonic function!theta-psh function on a complex space@$\theta$-psh function on a complex space} if for any $x\in X$, there is an open neighborhood $V$ of $x$ in $X$, a domain $\Omega\subseteq \mathbb{C}^N$, a closed immersion $V\hookrightarrow \Omega$ and a smooth function $g$ on $\Omega$ such that $\theta=(\ddc g)|_V$ and $g|_V+\varphi|_V\in \PSH(V)$.
\end{definition}

\begin{theorem}[Fornaess--Narasimhan]\label{thm:FN}\index{theorem!Fornaess--Narasimhan theorem}
	Let $\varphi\colon X\rightarrow [-\infty,\infty)$ be a function. Assume that $\varphi$ is not identically $-\infty$ on any irreducible component of $X$. Then the following are equivalent:
	\begin{enumerate}
		\item\label{itm:thm:FN:1} $\varphi$ is psh;
		\item\label{itm:thm:FN:2} $\varphi$ is usc and for any morphism $f\colon \Delta\rightarrow X$ from the open unit disk $\Delta$ in $\mathbb{C}$ to $X$ such that $f^*\varphi$ is not identically $-\infty$, the pull-back $f^*\varphi$ is subharmonic.
	\end{enumerate}
\end{theorem}
See \cite{FN80}.

\begin{theorem}[Grauert--Remmert]\label{thm:GRexten2}\index{theorem!Grauert--Remmert extension theorem}
Let $X$ be a unibranch\footnote{Unibranchness is very important here. Otherwise, consider the case where $X$ is the union of two copies of $\mathbb{C}$ intersecting only at their origins, $Z$ is the common origin. If we set $\varphi\equiv 0$ on one punctured plane and $\varphi\equiv 1$ on the other, then it is clear that $\varphi$ cannot be extended to $X$. This leads to a few misleading statements in the modern literature. The problem is that in the early German literature, \emph{komplexer Raum} is assumed to be either normal or unibranch!} complex space, let $Z$ be an analytic subset in $X$, and let $\varphi\in \PSH(X\setminus Z)$. Then the function $\varphi$ admits an extension to $\PSH(X)$ in the following two cases:
    \begin{enumerate}
        \item\label{itm:thm:GRexten2:1} The set $Z$ has codimension at least $2$ everywhere.
        \item\label{itm:thm:GRexten2:2} The set $Z$ has codimension at least $1$ everywhere and $\varphi$ is locally bounded from above on an open neighborhood of $Z$.
    \end{enumerate}
    In both cases, the extension is unique and is given by
    \begin{equation}\label{eq:GRextvarphi}
    \varphi(x)=\varlimsup_{X\setminus Z\ni y\to x}\varphi(y),\quad x\in X.
    \end{equation}
\end{theorem}
The proof given below combines \cite[Théorème~1.7]{Dem85} and \cite{GR56}.\footnote{This theorem has a quite entangled history. The corresponding results for subharmonic functions are generally attributed to Brelot. In \cite{GR56}, they cited a paper of Brelot written in 1934, which I cannot find. But in 1949, on the very first issue of \emph{Annales de l'Institut Fourier}, Brelot published a paper \cite{Bre49} with a very similar title, studying the behavior of a subharmonic function on the punctured neighborhood of a point. The general theorem was due to Grauert and Remmert, see \cite{GR56}. Their original proof was quite complicated, due to the fact that resolution of singularities was not available at that time. Later on, in 1985, Demailly published the influential paper \cite{Dem85} and gave a simpler proof. Oddly enough, Demailly did not cite either Grauert--Remmert or Brelot. He did not even mention that this result was already proved by Grauert--Remmert. The paper \cite{Dem85} is so influential that in France few people know the existence of \cite{GR56}.}

\begin{proof}
    We first prove the uniqueness, which is a local problem on $X$. Let $\psi$ denote the function defined by the right-hand side of \eqref{eq:GRextvarphi}. If $\eta\in\PSH(X)$ is an extension of $\varphi$, then upper semicontinuity gives $\eta\geq\psi$. Conversely, take $z\in Z$ and a holomorphic map $f\colon \Delta\rightarrow X$ such that $f(0)=z$ and $f(\Delta)\not\subset Z$. From the subharmonicity of $f^*\eta$ and \eqref{eq:GRextvarphi}, we find that
    \[
    \eta(z)=f^*\eta(0)=\varlimsup_{\Delta^* \ni w\to 0}f^*\eta(w)\leq \varlimsup_{X\setminus Z\ni y\to z}\varphi(y)=\psi(z),
    \]
    where the second equality follows from \cref{thm:GRexten}. Thus $\eta=\psi$, and uniqueness follows.

    Having established the uniqueness of the extension, the existence also becomes a local problem. So we can freely shrink $X$ around any point in $Z$ in the proof.

    \ref{itm:thm:GRexten2:2}
    Let $\pi\colon Y\rightarrow X$ be a resolution of singularities.
    By the smooth version of the extension theorem \cref{thm:GRexten}, we know that $\pi^*\varphi$ admits a unique extension to a psh function on $Y$, which we denote by $\eta$.
    Note that all fibers of $\pi$ are connected since $X$ is unibranch, thanks to \cref{thm:Zariskimain}. Furthermore, the fibers are compact, since $\pi$ is proper.
    The restriction of $\eta$ to each fiber is psh or identically $-\infty$, hence is constant by the maximum principle applied to a resolution. It therefore descends to an upper semicontinuous function $\psi$ on $X$, and $\psi|_{X\setminus Z}=\varphi$.

    We verify that $\psi$ is plurisubharmonic using \cref{thm:FN}. Let $f\colon \Delta\rightarrow X$ be a holomorphic map. We assume that $f^*\psi\not\equiv -\infty$.
    It suffices to show that $f^*\psi$ is subharmonic at $0\in \Delta$. After replacing $\Delta$ by a finite ramified cover if necessary, the germ of $f$ lifts to $Y$. Indeed, consider the proper projection
    \[
        \Delta\times_X Y\longrightarrow \Delta.
    \]
    Choose an irreducible analytic curve in $\Delta\times_X Y$ passing through a point over $0$ and dominating $\Delta$; after shrinking $\Delta$ and normalizing this curve, we obtain a finite ramified cover of the disk. After replacing $\Delta$ by a smaller disk, there are an integer $k>0$ and a holomorphic map $f'\colon \Delta\to Y$ such that
    \[
        f(t^k)=\pi(f'(t)).
    \]
    Therefore, $\psi(f(t^k))=\eta(f'(t))$ near $0$. Since subharmonicity descends under the finite map $t\mapsto t^k$, it follows that $f^*\psi$ is subharmonic near $0$.

     \ref{itm:thm:GRexten2:1} By the local description of complex spaces \cite[Section~3.4]{CAS}, we may assume that there is a domain $\Omega\subseteq \mathbb{C}^n$, a finite $s$-sheet branched covering $\Phi\colon X\rightarrow \Omega$ with branched locus contained in a proper analytic subset $V\subseteq \Omega$. We may assume that $X$ is connected, $n\geq 1$.

     It suffices to show that $\varphi$ is locally bounded from above near $Z$, since then \ref{itm:thm:GRexten2:2} applies. Suppose that this fails. Then we can find $z\in Z$ and $x_i\in X\setminus Z$ ($i\geq 1$) converging to $z$ such that
     \[
     \lim_{i\to\infty} \varphi(x_i)=\infty.
     \]
    Using \eqref{eq:GRextvarphi} around each $x_i$, we may further assume that $x_i\in X\setminus (\Phi^{-1}(\Phi(Z)\cup V))$.

     Let $L$ be a complex line passing through $\Phi(z)$ intersecting $(\Phi(Z)\cup V)\cap \overline{B}$ only at $\Phi(z)$, where $B\Subset B'\Subset \Omega$ are two small open balls centered at $\Phi(z)$. Since $\Phi$ is finite, $\Phi(Z)$ still has codimension at least $2$ in $\Omega$. Hence, for each $i$, the set of directions of lines through $\Phi(x_i)$ meeting $B'\cap \Phi(Z)$ is contained in a proper analytic subset of the projective space of directions. We can therefore choose lines $L_i$ through $\Phi(x_i)$ converging to $L$ such that $L_i\cap (B'\cap \Phi(Z))=\varnothing$. Moreover, since $\Phi(x_i)\notin V$, no such $L_i$ is contained in $V$, and hence $L_i\cap (B'\cap V)$ is discrete.
     Then $\Phi$ restricts to a branched covering over $B'\cap L_i$ for all $i\geq 1$.
     Adding a constant to $\varphi$, we may assume that $\varphi|_{\Phi^{-1}(L\cap \partial B)}<0$. We can then find an open neighborhood $U$ of $\Phi^{-1}(L\cap \partial B)$ so that $\varphi|_U<0$. For large $i$, we have $\Phi^{-1}(L_i\cap \partial B)\subseteq U$ and $\Phi(x_i)\in B$. Put
    \[
        C_i\coloneqq \Phi^{-1}(L_i\cap B').
    \]
    Thus $C_i$ is a one-dimensional complex space contained in $X\setminus Z$. Since $\Phi|_{C_i}$ is finite over $L_i\cap B'$, the inverse image $C_i\cap\Phi^{-1}(\overline B)$ is compact, and its boundary is contained in $\Phi^{-1}(L_i\cap\partial B)$.

    Passing to the normalization of each irreducible component of $C_i$, the pull-back of $\varphi$ is subharmonic. The one-dimensional maximum principle therefore gives
    \[
        \varphi(x_i)\leq \sup_{\Phi^{-1}(L_i\cap \partial B)}\varphi<0,
    \]
    which is a contradiction.
\end{proof}

\begin{corollary}\label{cor:PSH}
    Let $\pi\colon Y\rightarrow X$ be a proper bimeromorphic morphism between compact K\"ahler spaces. Let $\theta$ be a smooth closed real $(1,1)$-form on $X$. Assume that $X$ is unibranch. Then the pull-back induces a bijection
    \[
    \pi^*\colon \PSH(X,\theta)\xlongrightarrow{\sim}\PSH(Y,\pi^*\theta).
    \]
\end{corollary}
\begin{proof}
    Injectivity follows from the surjectivity of $\pi$. For surjectivity, let $u\in \PSH(Y,\pi^*\theta)$. The assertion is local on $X$, so after subtracting a local potential of $\theta$ we may assume that $\theta=0$.

    By \cref{thm:Zariskimain}, the fibers of $\pi$ are connected. The restriction of $u$ to any compact fiber is psh, hence is constant by the maximum principle. Thus $u$ descends to a function $v$ on $X$. Since $\pi$ is proper, $v$ is upper semicontinuous.

    It remains to check that $v$ is psh. By \cref{thm:FN}, it suffices to test on holomorphic disks. Let $f\colon \Delta\rightarrow X$ be a holomorphic disk such that $f^*v\not\equiv -\infty$. As in the proof of \cref{thm:GRexten2}, after a finite ramified cover $\rho\colon \Delta\rightarrow \Delta$, the germ of $f\circ\rho$ lifts to $Y$. The pull-back $(f^*v)\circ\rho$ is therefore subharmonic, so $f^*v$ is subharmonic. Hence $v\in \PSH(X)$, and the general $\theta$-psh case follows by adding back the local potentials.
\end{proof}

\section{Extensions of the results in the smooth setting}
\label{sec:unib_ext}
Let $X$ be an irreducible unibranch compact K\"ahler space of dimension $n$. Let $\theta$ be a closed real smooth $(1,1)$-form on $X$. We say \emph{the cohomology class $\{\theta\}$} is big if for any proper bimeromorphic morphism $\pi\colon Y\rightarrow X$ from a compact K\"ahler manifold $Y$, $\{\pi^*\theta\}$ is big. This definition is independent of the choice of $\pi$.

The non-pluripolar products can be defined exactly as in \cref{chap:npp} and the results in that chapter hold \emph{mutatis mutandis}.

The results in \cref{chap:envelopes} can also be easily extended. The definition of the $P$-envelope remains unchanged. As for the $\mathcal{I}$-envelope, we define
\begin{definition}\label{def:unibranch:L212}
    Given $\varphi\in \PSH(X,\theta)$, we define $P_{\theta}[\varphi]_{\mathcal{I}}\in \PSH(X,\theta)$ as the unique element with the following property: If $\pi\colon Y\rightarrow X$ is a proper bimeromorphic morphism from a compact K\"ahler manifold $Y$, then
    \index{envelope!I-envelope@$\mathcal{I}$-envelope!on a unibranch space}
    \[
        \pi^*P_{\theta}[\varphi]_{\mathcal{I}}=P_{\pi^*\theta}[\pi^*\varphi]_{\mathcal{I}}.
    \]
\end{definition}
It follows from \cref{cor:PSH} and \cref{prop:Ienvelopebimero} that $P_{\theta}[\varphi]_{\mathcal{I}}$ is independent of the choice of $\pi$ and is well-defined.
The other results can be easily extended.

\cref{chap:rays} and \cref{chap:comp} can be extended without big changes. The only exception is \cref{thm:equising_cond_general}, where we do not have the notion of multiplier ideal sheaves. So we do not know how to extend this theorem.

\cref{chap:Igood} can be extended except for \cref{sec:volHermitianbig} for the same reason as above.

If $Y\subseteq X$ is an irreducible analytic subspace not contained in $X^{\Sing}$, then the trace operator of \cref{chap:trace} can be extended to this setting by using an embedded resolution. In general, due to the lack of Demailly regularization on singular spaces, we do not know how to define the trace operator.

\cref{chap:testcurve} can be easily extended.

\cref{chap:Okounkov} is easy to extend since the partial Okounkov bodies are bimeromorphically invariant in the sense of \cref{thm:pobmain}.

\cref{chap:bdiv} is unchanged, since we always take projective limits with respect to all models in that section.

\cref{chap:toric_big} can be extended using the toric resolution \cite[Theorem~11.1.9]{CLS11}.

\cref{chap:NApotential} can be extended except for the parts involving the trace operator.

\cref{chap:Bergman} can be easily extended by considering a resolution.

\chapter{Almost semigroups}
\label[appendix]{chap:almostsg}

This appendix develops the theory of \emph{almost semigroups} in $\mathbb{Z}^{n+1}$ that underpins the partial Okounkov body construction of \cref{chap:Okounkov}. Almost semigroups generalize the integral semigroups arising in the classical Lazarsfeld--Musta\c{t}\u{a} and Kaveh--Khovanskii constructions to the setting where the underlying graded family is no longer multiplicative. The main result is the construction of a convex body $\Delta(S)\subseteq \mathbb{R}^n$ out of an almost semigroup $S$. The corresponding statements were stated without proof in \cref{chap:Okounkov}; the proofs are given here.

\section{Convex bodies}\label{sec:cb}

Fix $n\in \mathbb{N}$.

\begin{definition}\label{def:convbodies}
    A \emph{convex body}\index{convex body} in $\mathbb{R}^n$ is a non-empty compact convex set.
\end{definition}
We allow a convex body to have empty interior.

We write $\mathcal{K}_n$ for the set of convex bodies in $\mathbb{R}^n$.

\begin{definition}\label{def:Hausdorffmetric}
The \emph{Hausdorff metric}\index{Hausdorff metric}\index{$d_{\Haus}$} between $K_1,K_2\in \mathcal{K}_n$ is given by
\[
	d_{\Haus}(K_1,K_2)\coloneqq \max\left\{\adjustlimits\sup_{x_1\in K_1}\inf_{x_2\in K_2}|x_1-x_2|, \adjustlimits\sup_{x_2\in K_2}\inf_{x_1\in K_1}|x_1-x_2|\right\}.
\]

\end{definition}

It is well-known that the metric space $(\mathcal{K}_n,d_{\Haus})$ is complete.
We will need the following fundamental theorem:
\begin{theorem}[Blaschke selection theorem]\label{thm:Blaschke}\index{theorem!Blaschke selection theorem}
    The metric space $(\mathcal{K}_n,d_{\Haus})$ is locally compact.
\end{theorem}
We refer to \cite[Theorem~1.8.7]{Sch14} for details.

\begin{theorem}\label{thm:contvol}
	The Lebesgue volume $\vol\colon \mathcal{K}_n\rightarrow \mathbb{R}_{\geq 0}$ is continuous.
\end{theorem}
See \cite[Theorem~1.8.20]{Sch14}.

\begin{theorem}\label{thm:Hausconvcond}
	Let $K_i, K\in \mathcal{K}_n$ ($i\in \mathbb{N}$). Then $K_i\xrightarrow{d_{\Haus}}K$ if and only if the following conditions hold:
	\begin{enumerate}
		\item\label{itm:thm:Hausconvcond:1} each point $x\in K$ is the limit of a sequence $x_i\in K_i$, and
		\item\label{itm:thm:Hausconvcond:2} the limit of any convergent sequence $(x_{i_j})_{j\in \mathbb{N}}$ with $x_{i_j}\in K_{i_j}$ lies in $K$, where $i_j$ is a strictly increasing sequence in $\mathbb{Z}_{>0}$.
	\end{enumerate}
\end{theorem}
See \cite[Theorem~1.8.8]{Sch14}.

For a convex body $K\subseteq \mathbb{R}^n$ and $\epsilon>0$, we write $K_\epsilon\coloneqq \{x\in K: B(x,\epsilon)\subseteq K\}$ for the \emph{inner $\epsilon$-shrinking} of $K$. (Compare with the outer thickening $K^\epsilon$ used in \cref{cor:DXmaintoric}; the two are different operations.)

\begin{lemma}\label{lma:latcvb}
		Let $K\in \mathcal{K}_n$ be a convex body with positive volume and $K'\in \mathcal{K}_n$. Assume that, for some $k\in \mathbb Z_{>0}$, $K'$ contains $K\cap(k^{-1}\mathbb Z)^n$. Then
    \[
        K'\supseteq K_{n^{1/2}/k}.
    \]
\end{lemma}
\begin{proof}
    Let $x\in K_{n^{1/2}k^{-1}}$. By assumption, the closed ball $B$ with center $x$ and radius $n^{1/2}k^{-1}$ is contained in $K$. Observe that $x$ can be written as a convex combination of points in $B\cap (k^{-1}\mathbb{Z})^n$, which are contained in $K'$ by assumption. It follows that $x\in K'$.
\end{proof}

Given a sequence of convex bodies $K_i$ ($i\in \mathbb{N}$), we set
\[
	\varliminf_{i\to\infty}K_i= \overline{\bigcup_{i=0}^{\infty}\bigcap_{j\geq i}K_j}.
\]

Suppose $K$ is the limit of a subsequence of $K_i$, we have
\begin{equation}\label{eq:liminflimsup}
	\varliminf_{i\to\infty}K_i\subseteq K.
\end{equation}
This is a simple consequence of \cref{thm:Hausconvcond}.

\begin{lemma}\label{lma:Hausdorffconvslice}
    Assume that $n>0$.
    Let $K\subseteq \mathbb{R}^n$ be a convex body. Let
    \[
    t_{\min}\coloneqq \min \{t\in \mathbb{R}:\{x_1=t\}\cap K\neq \varnothing\},\quad t_{\max}\coloneqq \max \{t\in \mathbb{R}:\{x_1=t\}\cap K\neq \varnothing\}.
    \]
    Then for $t\in [t_{\min},t_{\max}]$, the map
    \[
    t\mapsto \{y\in \mathbb{R}^{n-1}:(t,y)\in K\}
    \]
    is continuous with respect to the Hausdorff metric.
\end{lemma}
Here $x_1$ denotes the first coordinate in $\mathbb{R}^n$.
\begin{proof}
We may assume that $t_{\min}<t_{\max}$ as otherwise there is nothing to prove.

For each $t\in [t_{\min},t_{\max}]$, we write $K_t=\{y\in \mathbb{R}^{n-1}:(t,y)\in K\}$.
Let $t_j\to t$ be a convergent sequence in $[t_{\min},t_{\max}]$, we want to show that $K_{t_j}$ converges to $K_t$ with respect to the Hausdorff metric. Recall that this amounts to the following two assertions:
\begin{enumerate}
    \item\label{itm:appendix-Okounkov:L87:1} For each convergent sequence $y_j\in K_{t_j}$ with limit $y$, we have $y\in K_t$;
    \item\label{itm:appendix-Okounkov:L87:2} Given any $y\in K_t$, up to replacing $t_j$ by a subsequence, we can find $y_j\in K_{t_j}$ converging to $y$.
\end{enumerate}
The first assertion is obvious. Let us prove the second. Take $y\in K_t$. Up to replacing $t_j$ by a subsequence and taking the symmetry into account, we may assume that $t_j> t$ for all $j$. In particular, $t<t_{\max}$.

We can find a point $z=(z^1,z')\in K$ such that $z^1>t$ (for example, there is always such a point with $z^1=t_{\max}$). Replacing $t_j$ by a subsequence, we may assume that $t_j\in (t,z^1)$ for all $j$. Then it suffices to take
\[
y_j=\frac{z^1-t_j}{z^1-t}y+\frac{t_j-t}{z^1-t}z'.
\]
This completes the proof.
\end{proof}

\begin{lemma}\label{lma:Hausdorffconvcoordbounds}
    Assume that $n>0$. The maps $\mathcal{K}_n\rightarrow \mathbb{R}$ given by
    \[
        K\mapsto \min_{x\in K}x_1,\qquad K\mapsto \max_{x\in K}x_1
    \]
    are continuous with respect to the Hausdorff metric.
\end{lemma}
\begin{proof}
    We prove the assertion for the maximum; the proof for the minimum is the same.
    Let $K,L\in \mathcal{K}_n$ and set $\epsilon=d_{\Haus}(K,L)$.
    For every $x\in K$, there exists $y\in L$ with $|x-y|\leq \epsilon$, hence
    \[
        x_1\leq y_1+\epsilon\leq \max_{z\in L}z_1+\epsilon.
    \]
    Taking the maximum over $x\in K$ gives
    \[
        \max_{x\in K}x_1\leq \max_{z\in L}z_1+\epsilon.
    \]
    Interchanging $K$ and $L$, we get
    \[
        \left|\max_{x\in K}x_1-\max_{z\in L}z_1\right|\leq d_{\Haus}(K,L).
    \]
    This proves the continuity.
\end{proof}

\begin{lemma}\label{lma:intconvexset}
    Let $D_j\subseteq \mathbb{R}^n$ $(j\geq 1)$ be a decreasing sequence of bounded convex sets. Assume that $\vol \bigcap_j D_j>0$. Then
    \[
    \overline{\bigcap_{j=1}^{\infty}D_j}=\bigcap_{j=1}^{\infty}\overline{D_j}.
    \]
\end{lemma}
\begin{proof}
    The $\subseteq$ direction is clear. By convexity, it suffices to show that both sides have the same positive volume: indeed, a proper inclusion of closed convex sets with positive volume would cut off a non-empty relatively open convex piece by a separating hyperplane, hence would strictly decrease the volume. As the boundary of convex sets has zero Lebesgue measure, it follows that the volumes of both sides are equal to $\lim_{j\to\infty}\vol D_j$.
\end{proof}

\begin{definition}\label{def:almostsg:L138}
    Let $K,K'\in \mathcal{K}_n$, their \emph{Minkowski sum}\index{Minkowski sum} is given by
    \[
    K+K'\coloneqq \{x+x': x\in K,x'\in K'\}.
    \]
\end{definition}
\begin{proposition}
    The Minkowski sum $\mathcal{K}_n\times \mathcal{K}_n\rightarrow \mathcal{K}_n$ is continuous.
\end{proposition}
See \cite[Page~139]{Sch14}.
\begin{theorem}[Brunn--Minkowski]\label{thm:BrunnMin}\index{theorem!Brunn--Minkowski theorem}
    Let $K,K'\in \mathcal{K}_n$. Then for any $t\in (0,1)$, we have
    \[
    \vol\bigl((1-t)K'+tK\bigr)\geq \bigl(\vol K'\bigr)^{1-t}\bigl(\vol K\bigr)^t.
    \]
\end{theorem}
In other words, the volume is log concave. See \cite[Page~372]{Sch14}.

\section{The Okounkov bodies of almost semigroups}\label{sec:clo}
Fix an integer $n\geq 0$. Fix a closed convex cone $C\subseteq \mathbb{R}^{n}\times \mathbb{R}_{\geq 0}$ such that $C\cap \{x_{n+1}=0\}=\{0\}$. Here $x_{n+1}$ is the last coordinate of $\mathbb{R}^{n+1}$.
\subsection{Generalities on semigroups}
\label{subsec:asg_general}
Write $\hat{\mathcal{S}}(C)$ for the set of subsets of $C\cap \mathbb{Z}^{n+1}$ and $\mathcal{S}(C)$ for the set of subsemigroups $S\subseteq C\cap \mathbb{Z}^{n+1}$. For each $k\in \mathbb{N}$ and $S\in \hat{\mathcal{S}}(C)$, we write
\[
S_k\coloneqq \left\{x\in \mathbb{Z}^n: (x,k)\in S \right\}.
\]
Note that $S_k$ is a finite set by our assumption on $C$.

We introduce a pseudometric on $\hat{\mathcal{S}}(C)$ as follows:
\begin{equation}\label{eq:dsgdefapp}
d_{\sg}(S,S')\coloneqq \varlimsup_{k\to\infty} k^{-n} \left(|S_k|+|S'_k|-2|(S\cap S')_k|\right).
\end{equation}
Here $|\bullet|$ denotes the cardinality of a finite set.
\begin{lemma}\label{lma:dpsapp}
    The above defined $d_{\sg}$ is a pseudometric on $\hat{\mathcal{S}}(C)$.
\end{lemma}
\begin{proof}
    Only the triangle inequality needs to be argued. Take $S,S',S''\in \hat{\mathcal{S}}(C)$. We claim that for any $k\in \mathbb{N}$,
    \[
    |S_k|+|S'_k|-2|S_k\cap S'_k|+|S''_k|+|S'_k|-2|S''_k\cap S'_k|\geq |S_k|+|S''_k|-2|S_k\cap S''_k|.
    \]
    From this the triangle inequality follows. To argue the claim, we rearrange it to the following form:
    \[
        |S'_k|-|S_k\cap S'_k| \geq |S'_k\cap S''_k|-|S_k\cap S''_k|,
    \]
    which is obvious.
\end{proof}
Given $S,S'\in \hat{\mathcal{S}}(C)$, we say $S$ is equivalent to $S'$ and write $S\sim S'$ if $d_{\sg}(S,S')=0$. This is an equivalence relation by \cref{lma:dpsapp}.

\begin{lemma}\label{lma:dBilapp}
    Given $S,S',S''\in \hat{\mathcal{S}}(C)$, we have
    \[
    d_{\sg}(S\cap S'',S'\cap S'')\leq d_{\sg}(S,S').
    \]
    In particular, if $S^i,S'^i\in \hat{\mathcal{S}}(C)$ ($i\in \mathbb{N}$) and $S^i\to S$, $S'^i\to S'$, then
    \[
    S^i\cap S'^i\to S\cap S'.
    \]
\end{lemma}
\begin{proof}
    Observe that for any $k\in \mathbb{N}$,
    \[
    |S_k\cap S_k''|-|S_k\cap S_k'\cap S_k''|\leq |S_k|-|S_k\cap S_k'|.
    \]
    The same holds if we interchange $S$ with $S'$. It follows that
    \[
    |S_k\cap S_k''|+|S_k'\cap S_k''|-2|S_k\cap S_k'\cap S_k''|\leq |S_k|+|S'_k|-2|S_k\cap S_k'|.
    \]
    The first assertion follows.

    Next we compute
    \[
    \begin{split}
    d_{\sg}(S^i\cap S'^i,S\cap S')\leq & d_{\sg}(S^i\cap S'^i,S^i\cap S')+d_{\sg}(S^i\cap S',S\cap S')\\
    \leq & d_{\sg}(S'^i,S')+d_{\sg}(S^i,S)
    \end{split}
    \]
    and the second assertion follows.
\end{proof}

The volume of $S\in \mathcal{S}(C)$ is defined as
\[
\vol S\coloneqq \lim_{k\to\infty} (ka)^{-n}|S_{ka}|=\varlimsup_{k\to\infty} k^{-n}|S_k|,
\]
where $a$ is a sufficiently divisible positive integer.
The existence of the limit and its independence from $a$ both follow from the more precise result \cite[Theorem~2]{KK12}, a special case of which can be stated as follows:
\begin{theorem}[Kaveh--Khovanskii]\label{thm:KK12growth}
    Let $S\in \mathcal{S}(C)$. Then, for every sufficiently divisible positive integer $a$, the limit
    \[
        \lim_{k\to\infty}(ka)^{-n}|S_{ka}|
    \]
    exists, is independent of $a$, and is equal to $\varlimsup_{k\to\infty} k^{-n}|S_k|$.
    If moreover $S$ generates $\mathbb{Z}^{n+1}$ as an Abelian group, and if $C(S)$ denotes the closed convex cone generated by $S\cup\{0\}$, then
    \[
        \lim_{k\to\infty} k^{-n}|S_k|
        =
        \vol \{x\in \mathbb{R}^n:(x,1)\in C(S)\}.
    \]
    In particular, for every convex body $K\subseteq \mathbb{R}^n$,
    \[
        \lim_{k\to\infty} k^{-n}\#\left(kK\cap \mathbb{Z}^n\right)=\vol K.
    \]
    The same assertion holds with $\mathbb{Z}^n$ replaced by any lattice, provided the Lebesgue measure is normalized by that lattice.
\end{theorem}
This is the part of \cite[Theorem~2]{KK12} used in this book.
\begin{lemma}\label{lma:vollipapp}
    Let $S,S'\in \mathcal{S}(C)$. Then
    \[
    |\vol S-\vol S'|\leq d_{\sg}(S,S').
    \]
\end{lemma}
\begin{proof}
    Choose a sufficiently divisible positive integer $a$. Taking the limsup in the definition of $d_{\sg}$ along $k=ma$ gives
    \[
        d_{\sg}(S,S')\geq \vol S+\vol S'-2\vol(S\cap S').
    \]
    It follows that $\vol S-\vol S'\leq d_{\sg}(S,S')$ and $\vol S'-\vol S\leq d_{\sg}(S,S')$.
\end{proof}
We define $\overline{\mathcal{S}}(C)$ as the closure of $\mathcal{S}(C)$ in $\hat{\mathcal{S}}(C)$ with respect to the topology defined by the pseudometric $d_{\sg}$.
By \cref{lma:vollipapp}, $\vol\colon \mathcal{S}(C)\rightarrow \mathbb{R}$ admits a unique $1$-Lipschitz extension to
\begin{equation}\label{eq:volexapp}
\vol\colon \overline{\mathcal{S}}(C)\rightarrow \mathbb{R}.
\end{equation}
\begin{lemma}\label{lma:volcompaapp}
    Suppose that $S,S'\in \overline{\mathcal{S}}(C)$ and $S\subseteq S'$. Then
    \[
    \vol S\leq \vol S'.
    \]
\end{lemma}
\begin{proof}
    Take sequences $S^j,S'^j$ in $\mathcal{S}(C)$ such that $S^j\to S$, $S'^j\to S'$. By \cref{lma:dBilapp}, after replacing $S^j$ by $S^j\cap S'^j$, we may assume that $S^j\subseteq S'^j$ for each $j$. Then our assertion follows easily.
\end{proof}

\subsection{Okounkov bodies of semigroups}
\label{subsec:asg_okoun}
Given $S\in \hat{\mathcal{S}}(C)$, we will write $C(S)\subseteq C$ for the closed convex cone generated by $S\cup \{0\}$. Moreover, for each $k\in \mathbb{Z}_{>0}$, we define
\[
\Delta_k(S)\coloneqq \Conv \left\{k^{-1}x\in \mathbb{R}^n:x\in S_k \right\}\subseteq \mathbb{R}^n.
\]
Here $\Conv$ denotes the convex hull.
\begin{definition}\label{def:almostsg:L284}
    Let $\mathcal{S}'(C)$ be the subset of $\mathcal{S}(C)$ consisting of semigroups $S$ such that $S$ generates $\mathbb{Z}^{n+1}$ (as an Abelian group).
    \index{semigroup!$\mathcal{S}'(C)$}
    \index{$\mathcal{S}'(C)$}
\end{definition}
Note that for any $S\in \mathcal{S}'(C)$, the cone $C(S)$ has full dimension (i.e., the topological interior is non-empty). Given a full-dimensional subcone $C'\subseteq C$, the semigroup $C'\cap \mathbb{Z}^{n+1}$ lies in $\mathcal{S}'(C)$.

This class behaves well under intersections:
\begin{lemma}\label{lma:intersecS'app}
    Let $S,S'\in \mathcal{S}'(C)$. Assume that $\vol (S\cap S')>0$. Then $S\cap S'\in \mathcal{S}'(C)$, and
    \[
    C(S\cap S')=C(S)\cap C(S').
    \]
\end{lemma}
The lemma obviously fails if $\vol (S\cap S')=0$.
\begin{proof}
    We first observe that the cone $C(S)\cap C(S')$ has full dimension since otherwise $\vol(S\cap S')=0$. Take a full-dimensional subcone $C'$ in $C(S)\cap C(S')$ such that $C'$ intersects the boundary of $C(S)\cap C(S')$ only at $0$. It follows from \cite[Theorem~1]{KK12} that there is an integer $N>0$ such that any $x\in \mathbb{Z}^{n+1}\cap C'$ with Euclidean norm at least $N$ lies in $S\cap S'$. Therefore, $S\cap S'\in \mathcal{S}'(C)$.
    Since $C'$ can be chosen around any ray in the interior of $C(S)\cap C(S')$, the same argument also shows that
    \[
    C(S\cap S')=C(S)\cap C(S').
    \]
    This completes the proof.
\end{proof}

We recall the following definition from \cite{KK12}.
\begin{definition}\label{def:Okokkapp}
 Given $S\in \mathcal{S}'(C)$, its \emph{Okounkov body}\index{Okounkov body!Okounkov body of an almost semigroup} is defined as follows
\[
\Delta(S)\coloneqq \left\{x\in \mathbb{R}^n:(x,1)\in C(S) \right\}.
\]
\end{definition}

\begin{theorem}\label{thm:HausOkounapp}
	For each $S\in \mathcal{S}'(C)$, we have
	\begin{equation}\label{eq:volWvolDeltaapp}
		\vol S=\lim_{k\to\infty} k^{-n}|S_k|=\vol \Delta(S)>0.
	\end{equation}
    Moreover, as $k\to\infty$,
	\begin{equation}\label{eq:HausconvDeltaGLSapp}
		\Delta_{k}(S)\xrightarrow{d_{\Haus}}\Delta(S).
	\end{equation}
\end{theorem}
This is essentially proved in \cite[Lemma~4.8]{WN14}, which itself follows from a theorem of Khovanskii \cite{Kho92}.
\begin{proof}
	The equalities \eqref{eq:volWvolDeltaapp} follow from \cref{thm:KK12growth}.



	It remains to prove \eqref{eq:HausconvDeltaGLSapp}. By the argument of \cite[Lemma~4.8]{WN14}, for any compact set $K\subseteq \Int \Delta(S)$, there is $k_0>0$ such that for any $k\geq k_0$, $\alpha\in K\cap (k^{-1}\mathbb{Z})^n$ implies that $\alpha\in \Delta_k(S)$.

	To conclude, fix $\epsilon>0$. Choose $\delta>0$ such that $K\coloneqq \Delta(S)_{\delta}$ has positive volume and $d_{\Haus}(K,\Delta(S))<\epsilon/2$. Applying \cref{lma:latcvb} to $K$, we find that $\Delta_k(S)$ contains $K_{n^{1/2}k^{-1}}$ for all large $k$. After increasing $k$, we may also assume that $d_{\Haus}(K,K_{n^{1/2}k^{-1}})<\epsilon/2$. Since $\Delta_k(S)\subseteq \Delta(S)$, this implies \eqref{eq:HausconvDeltaGLSapp}.
\end{proof}
\begin{corollary}\label{cor:distapp}
    Let $S,S'\in \mathcal{S}'(C)$. Assume that $\vol (S\cap S')>0$. Then we have
    \[
    d_{\sg}(S,S')=\vol(S)+\vol(S')-2\vol(S\cap S').
    \]
\end{corollary}
\begin{proof}
    This is a direct consequence of \cref{lma:intersecS'app} and \eqref{eq:volWvolDeltaapp}.
\end{proof}

\begin{lemma}\label{lma:regularizatapp}
    Given $S\in \mathcal{S}'(C)$, we have $S\sim \Reg(S)$.
\end{lemma}
Recall that the regularization $\Reg(S)$ of $S$ is defined as $C(S)\cap \mathbb{Z}^{n+1}$.
\begin{proof}
    Since $S$ and $\Reg(S)$ have the same Okounkov body,
    we have $\vol S=\vol \Reg(S)$ by \cref{thm:HausOkounapp}. By \cref{cor:distapp} again,
    \[
    d_{\sg}(\Reg(S),S)=\vol \Reg(S)-\vol S=0.
    \]
    This completes the proof.
\end{proof}

\begin{lemma}\label{lma:Deltaindclassapp}
	Let $S,S'\in \mathcal{S}'(C)$. Assume that $d_{\sg}(S,S')=0$. Then $\Delta(S)=\Delta(S')$.
\end{lemma}
\begin{proof}
Observe that $\vol (S\cap S')>0$, as otherwise
\[
d_{\sg}(S,S')\geq \vol S+\vol S'>0,
\]
which is a contradiction.

It follows from \cref{lma:intersecS'app} that $S\cap S'\in \mathcal{S}'(C)$. It suffices to show that $\Delta(S)=\Delta(S\cap S')$.
In fact, suppose that this holds, since $\vol \Delta(S')=\vol S'=\vol S=\vol \Delta(S)$, the inclusion $\Delta(S')\supseteq \Delta(S\cap S')=\Delta(S)$ is an equality.

After replacing $S'$ by $S\cap S'$, we still have $d_{\sg}(S,S')=0$ by \cref{lma:dBilapp}, and now $S'\subseteq S$.
Hence \cref{cor:distapp} gives
    \[
        0=d_{\sg}(S,S')=\vol S-\vol S',
    \]
so $\vol S=\vol S'$. Therefore
    \[
        \vol\Delta(S)=\vol\Delta(S')>0
    \]
by \eqref{eq:volWvolDeltaapp}. Therefore, they are equal.
\end{proof}

\begin{lemma}\label{lma:Sprimeintapp}
    Suppose that $S^i\in \mathcal{S}'(C)$ is a decreasing sequence such that
    \[
    \lim_{i\to\infty}\vol S^i>0.
    \]
    Then there is $S\in \mathcal{S}'(C)$ such that $S^i\to S$.
\end{lemma}
In general, one cannot simply take $S=\bigcap_i S^i$. For example, consider the sequence $S^i=S^1\cap \{x_{n+1}\geq i\}$.
\begin{proof}
    By \cref{lma:regularizatapp}, we may replace $S^i$ by its regularization and assume that $S^i=C(S^i)\cap \mathbb{Z}^{n+1}$. We define
    \[
    S=\left( \bigcap_{i=1}^{\infty}C(S^i) \right) \cap \mathbb{Z}^{n+1}.
    \]
    The convex bodies $\Delta(S^i)$ form a decreasing sequence and satisfy
    \[
    \lim_{i\to\infty}\vol\Delta(S^i)=\lim_{i\to\infty}\vol S^i>0.
    \]
    By continuity from above for Lebesgue measure,
    \[
    \vol\bigcap_i\Delta(S^i)=\lim_{i\to\infty}\vol\Delta(S^i)>0.
    \]
    Set
    \[
        K\coloneqq \bigcap_{i=1}^{\infty}C(S^i).
    \]
    The cone $K$ is full-dimensional. As noted above, $K\cap\mathbb Z^{n+1}$ generates $\mathbb Z^{n+1}$ as an Abelian group; hence $S=K\cap\mathbb Z^{n+1}$ lies in $\mathcal S'(C)$. Moreover, the closed convex cone generated by $S$ is exactly $K$. Therefore
    \[
        \Delta(S)=\{x:(x,1)\in K\}=\bigcap_{i=1}^{\infty}\{x:(x,1)\in C(S^i)\}=\bigcap_{i=1}^{\infty}\Delta(S^i).
    \]
    By \cref{cor:distapp} and \cref{thm:HausOkounapp}, we can compute the distance
    \[
    d_{\sg}(S,S^i)=\vol S^i-\vol S=\vol \Delta (S^i)-\vol \Delta (S),
    \]
    which tends to $0$ by construction.
\end{proof}

\subsection{Okounkov bodies of almost semigroups}\label{subsec:Okobalmosgapp}
\begin{definition}\label{def:SpCbarapp}
    We define $\overline{\mathcal{S}'(C)}_{>0}$ to be the set of elements in the closure of $\mathcal{S}'(C)$ in $\hat{\mathcal{S}}(C)$ with positive volume. An element in $\overline{\mathcal{S}'(C)}_{>0}$ is called an \emph{almost semigroup}\index{almost semigroup} in $C$.
    \index{$\overline{\mathcal{S}'(C)}_{>0}$}
\end{definition}
Recall that the volume here is defined in \eqref{eq:volexapp}.

Our goal is to prove the following theorem:
\begin{theorem}\label{thm:Okocontapp}
	The Okounkov body map $\Delta\colon \mathcal{S}'(C) \rightarrow \mathcal{K}_n$
    as defined in \cref{def:Okokkapp} admits a unique continuous extension
	\begin{equation}\label{eq:Deltagensgapp}
		\Delta \colon \overline{\mathcal{S}'(C)}_{>0}\rightarrow \mathcal{K}_n.
	\end{equation}
     Moreover, for any $S\in \overline{\mathcal{S}'(C)}_{>0}$, we have
	\begin{equation}\label{eq:volWfinalapp}
		\vol S=\vol \Delta(S).
	\end{equation}
\end{theorem}

\begin{proof}
The uniqueness of the extension is clear once existence is known. Moreover, \eqref{eq:volWfinalapp} follows easily from \cref{thm:HausOkounapp} and \cref{thm:contvol} by continuity.
It remains to argue the existence of the continuous extension. We first construct an extension and prove its continuity.

\textbf{Step~1}. We construct the desired map \eqref{eq:Deltagensgapp}. Let $S\in \overline{\mathcal{S}'(C)}_{>0}$. We wish to construct a convex body $\Delta(S)\in \mathcal{K}_n$.

    Let $S^i\in \mathcal{S}'(C)$ be a sequence that converges to $S$ such that
    \[
    d_{\sg}(S^i,S^{i+1})\leq 2^{-i}.
    \]
    For each $i,j\geq 0$, we introduce
    \[
    S^{i,j}=S^i\cap S^{i+1}\cap \cdots \cap S^{i+j}.
    \]
    Then by \cref{lma:dBilapp},
    \[
    d_{\sg}(S^{i,j},S^{i,j+1})\leq 2^{-i-j}.
    \]
    Take $i_0>0$ large enough so that, for $i\geq i_0$, we have $\vol S^i>2^{-1}\vol S$ and $2^{2-i}< \vol S$. Then
    \[
        \vol S^i-\vol S^{i,j}\leq d_{\sg}(S^{i,0},S^{i,1})+d_{\sg}(S^{i,1},S^{i,2})+\cdots+d_{\sg}(S^{i,j-1},S^{i,j})\leq 2^{1-i}.
    \]
    It follows that $\vol S^{i,j}> 2^{-1}\vol S-2^{1-i}>0$ whenever $i\geq i_0$. In particular, by \cref{lma:intersecS'app}, $S^{i,j}\in \mathcal{S}'(C)$ for $i\geq i_0$.

    By \cref{lma:Sprimeintapp}, for $i\geq i_0$, there exists $T^i\in \mathcal{S}'(C)$ such that $S^{i,j}\to T^i$ as $j\to\infty$. Moreover,
    \[
    d_{\sg}(T^i,S)=\lim_{j\to\infty} d_{\sg}(S^{i,j},S)\leq \lim_{j\to\infty} d_{\sg}(S^{i,j},S^i)+d_{\sg}(S^i,S)\leq 2^{1-i}+d_{\sg}(S^i,S).
    \]
    Therefore, $T^i\to S$. We then define
    \[
    \Delta(S)\coloneqq \overline{\bigcup_{i=i_0}^{\infty} \Delta(T^i)}.
    \]
    In other words, we have defined
    \[
	\Delta(S)\coloneqq \varliminf_{i\to\infty}\Delta(S^i).
    \]
    Indeed, by the construction in the proof of \cref{lma:Sprimeintapp} and the cone equality in \cref{lma:intersecS'app},
    \[
    \Delta(T^i)=\bigcap_{j\geq 0}\Delta(S^{i,j})
    =\bigcap_{j\geq i}\Delta(S^j).
    \]
    Thus the bodies $\Delta(T^i)$ are increasing in $i$ and are all contained in the compact slice $\{x\in \mathbb{R}^n:(x,1)\in C\}$. In particular, $\Delta(S)$ is a convex body and
    \[
    \vol\Delta(S)=\lim_{i\to\infty}\vol\Delta(T^i)=\lim_{i\to\infty}\vol T^i=\vol S.
    \]
    This is an honest limit: if $\Delta$ is the limit of a subsequence of $\Delta(S^i)$, then $\Delta(S)\subseteq \Delta$ by \eqref{eq:liminflimsup}. Moreover,
    \[
    \vol\Delta=\lim_{m\to\infty}\vol\Delta(S^{i_m})=\lim_{m\to\infty}\vol S^{i_m}=\vol S.
    \]
    Comparing the volumes, we find that equality holds. So by \cref{thm:Blaschke},
    \begin{equation}\label{eq:deltawtempapp}
	\Delta(S)=\lim_{i\to\infty}\Delta(S^i).
    \end{equation}

    Next we claim that $\Delta(S)$ as defined above does not depend on the choice of the sequence $S^i$. In fact, suppose that $S'^i\in \mathcal{S}'(C)$ is another sequence satisfying the same conditions as $S^i$.
    By \cref{lma:dBilapp}, we have
    \[
        R^i=S^{i+1}\cap S'^{i+1}\longrightarrow S\cap S=S.
    \]
    Since $\vol S>0$ and the volume is continuous on $\overline{\mathcal S}(C)$, we have $\vol R^i>0$ for all large $i$. Thus $R^i\in\mathcal S'(C)$ by \cref{lma:intersecS'app}. Moreover,
    \[
        d_{\sg}(R^i,R^{i+1}) \leq d_{\sg}(S^{i+1},S^{i+2}) + d_{\sg}(S'^{i+1},S'^{i+2}).
    \]
    Passing simultaneously to subsequences of $S^i$ and $S'^i$, we may assume
    \[
        d_{\sg}(S^{i+1},S^{i+2})+d_{\sg}(S'^{i+1},S'^{i+2})\leq 2^{-i}.
    \]
    Then the construction above applies to $R^i=S^{i+1}\cap S'^{i+1}$. Since $R^i\subseteq S^{i+1}$ and $R^i\subseteq S'^{i+1}$, every Hausdorff limit of $\Delta(R^i)$ is contained in the corresponding limits of $\Delta(S^i)$ and $\Delta(S'^i)$. The volumes of all three limits are $\vol S$, hence the inclusions are equalities. The same is true with $S'^i$ in place of $S^i$. So we conclude that $\Delta(S)$ as in \eqref{eq:deltawtempapp} does not depend on the choices we made.

    \textbf{Step~2}. It remains to prove the continuity of $\Delta$ defined in Step~1.
    Suppose that $S^i\in \overline{\mathcal{S}'(C)}_{>0}$ is a sequence with limit $S\in \overline{\mathcal{S}'(C)}_{>0}$.
    We want to show that
	\begin{equation}\label{eq:temp5app}
		\Delta(S^i)\xrightarrow{d_{\Haus}}\Delta(S).
	\end{equation}

    We first reduce to the case where $S^i\in \mathcal{S}'(C)$.
    By \eqref{eq:deltawtempapp}, for each $i$, we can choose $T^i\in \mathcal{S}'(C)$ such that $d_{\sg}(S^i,T^i)<2^{-i}$ and $d_{\Haus}(\Delta(S^i),\Delta(T^i))<2^{-i}$. If we have shown $\Delta(T^i)\xrightarrow{d_{\Haus}}\Delta(S)$, then \eqref{eq:temp5app} follows immediately.

    Next we reduce to the case where $d_{\sg}(S^i,S^{i+1})\leq 2^{-i}$. Thanks to \cref{thm:Blaschke}, in order to prove \eqref{eq:temp5app}, it suffices to show that each subsequence of $\Delta(S^i)$ admits a subsequence that converges to $\Delta(S)$. Hence, we may pass to a sufficiently fast subsequence and reduce to the required case.

    After these reductions, \eqref{eq:temp5app} is nothing but \eqref{eq:deltawtempapp}.
\end{proof}

\begin{corollary}\label{cor:Okocompapp}
    Suppose that $S,S'\in \overline{\mathcal{S}'(C)}_{>0}$ with $S\subseteq S'$. Then
    \begin{equation}\label{eq:Deltacontainapp}
    \Delta(S)\subseteq \Delta(S').
    \end{equation}
\end{corollary}
\begin{proof}
    Let $S^j,S'^j\in \mathcal{S}'(C)$ be elements such that $S^j\to S$, $S'^j\to S'$. Then it follows from \cref{lma:dBilapp} that $S^j\cap S'^j \to S$. Since $\vol$ is continuous, for large $j$, $S^j\cap S'^j$ has positive volume and hence lies in $\mathcal{S}'(C)$ by \cref{lma:intersecS'app}. We may therefore replace $S^j$ by $S^j\cap S'^j$ and assume that $S^j\subseteq S'^j$. Hence, \eqref{eq:Deltacontainapp} follows from the continuity of $\Delta$ proved in \cref{thm:Okocontapp}.
\end{proof}

\begin{remark}\label{rmk:almostsg:L542}
    The construction of $\Delta$ is independent of the choice of $C$ in the following sense. If $C'$ is another cone satisfying the same assumptions as $C$ and $C'\supseteq C$, then the Okounkov body map $\Delta \colon \overline{\mathcal{S}'(C')}_{>0}\rightarrow \mathcal{K}_n$ extends the corresponding map \eqref{eq:Deltagensgapp}. We will use this fact without further comment.
\end{remark}

\backmatter
\printindex

\bibliographystyle{sgptalpha}
\bibliography{ComplexGeometry}


\end{document}